%% file: NewFSZ.tex
\documentclass[11pt]{amsbook}

\usepackage{amssymb} 
\usepackage{amscd}
\usepackage{verbatim} 

\usepackage{amssymb}
\usepackage{graphicx}

\usepackage{pictex}
\usepackage{color}

\makeindex

\begin{document}


\newtheorem{theorem}{Theorem}[chapter] 
\newtheorem{corollary}[theorem]{Corollary}
\newtheorem{conjecture}[theorem]{Conjecture} 
\newtheorem{lemma}[theorem]{Lemma} 
\newtheorem{proposition}[theorem]{Proposition} 

\theoremstyle{definition}
\newtheorem{definition}[theorem]{Definition}
\newtheorem{example}[theorem]{Example}
\newtheorem*{theorem*}{Theorem}

\theoremstyle{remark}
\newtheorem{remark}[theorem]{Remark}
\newtheorem{exercise}[theorem]{Exercise}

\numberwithin{equation}{chapter} 


\def\vA{{\vec{A}}}
\def\vB{{\vec{B}}}
\def\vC{{\vec{C}}}
\def\vD{{\vec{D}}}
\def\vE{{\vec{E}}}
\def\vF{{\vec{F}}}
\def\vG{{\vec{G}}}
\def\vH{{\vec{H}}}
\def\vI{{\vec{I}}}
\def\vJ{{\vec{J}}}
\def\vK{{\vec{K}}}
\def\vL{{\vec{L}}}
\def\vM{{\vec{M}}}
\def\vP{{\vec{P}}}
\def\vQ{{\vec{Q}}}
\def\vR{{\vec{R}}}
\def\vS{{\vec{S}}}
\def\vT{{\vec{T}}}
\def\vU{{\vec{U}}}
\def\vW{{\vec{W}}}
\def\vX{{\vec{X}}}
\def\vY{{\vec{Y}}}
\def\vZ{{\vec{Z}}}

\def\va{{\vec{a}}}
\def\vb{{\vec{b}}}
\def\vc{{\vec{c}}}
\def\vd{{\vec{d}}}
\def\ve{{\vec{e}}}
\def\vf{{\vec{f}}}
\def\vg{{\vec{g}}}
\def\vh{{\vec{h}}}
\def\vi{{\vec{i}}}
\def\vj{{\vec{j}}}
\def\vk{{\vec{k}}}
\def\vl{{\vec{\ell}}}
\def\vm{{\vec{m}}}
\def\vn{{\vec{n}}}
\def\vo{{\vec{o}}}
\def\vp{{\vec{p}}}
\def\vq{{\vec{q}}}
\def\vr{{\vec{r}}}
\def\vs{{\vec{s}}}
\def\vt{{\vec{t}}}
\def\vu{{\vec{u}}}
\def\vw{{\vec{w}}}
\def\vx{{\vec{x}}}
\def\vy{{\vec{y}}}
\def\vz{{\vec{z}}}

\def\oORIG{{\overline{o}}}
\def\vORIG{{\vec{0}}}
\def\v1{{\vec{1}}}
\def\vbb1{\vec{{\mathbf 1}}}

\def\valpha{{\vec{\alpha}}}
\def\vbeta{{\vec{\beta}}}
\def\vdelta{{\vec{\delta}}}
\def\vgamma{{\vec{\gamma}}}
\def\veps{{\vec{\varepsilon}}}
\def\vtheta{{\vec{\theta}}}
\def\vrho{{\vec{\rho}}}
\def\vsigma{{\vec{\sigma}}}
\def\vchi{{\vec{\chi}}}

\def\vDelta{{\vec{\Delta}}}

\def\BFO{{\mathbf 0}}
\def\BF1{{\mathbf 1}}
\def\BFb{{\mathbf b}}
\def\BFA{{\mathbf A}}
\def\BFI{{\mathbf I}}
\def\BFJ{{\mathbf J}}
\def\BFD{{\mathbf \tTheta}}
\def\BFN{{\mathbf \nabla}}
\def\BFt{{\mathbf{\widehat{w}}}}
\def\BFj{{\mathbf j}}

\def\bE{{\mathbb E}}
\def\bF{{\mathbb F}}
\def\bU{{\mathbb U}}

\def\AA{{\mathbb A}}
\def\BB{{\mathbb B}}
\def\CC{{\mathbb C}}
\def\DD{{\mathbb D}}
\def\EE{{\mathbb E}}
\def\FF{{\mathbb F}}
\def\GG{{\mathbb G}}
\def\HH{{\mathbb H}}
\def\II{{\mathbb I}}
\def\JJ{{\mathbb J}}
\def\KK{{\mathbb K}}
\def\LL{{\mathbb L}}
\def\MM{{\mathbb M}}
\def\NN{{\mathbb N}}
\def\OO{{\mathbb O}}
\def\PP{{\mathbb P}}
\def\QQ{{\mathbb Q}}
\def\RR{{\mathbb R}}
\def\SS{{\mathbb S}}
\def\TT{{\mathbb T}}
\def\UU{{\mathbb U}}
\def\VV{{\mathbb V}}
\def\WW{{\mathbb W}}
\def\XX{{\mathbb X}}
\def\YY{{\mathbb Y}}
\def\ZZ{{\mathbb Z}}

\def\Bfone{{\mathbb 1}}
\def\Bfone{{\mathbb 1}}

\def\cA{{\mathcal A}}
\def\cB{{\mathcal B}}
\def\cC{{\mathcal C}}
\def\cD{{\mathcal D}}
\def\cE{{\mathcal E}}
\def\cF{{\mathcal F}}
\def\cG{{\mathcal G}}
\def\cH{{\mathcal H}}
\def\cI{{\mathcal I}}
\def\cJ{{\mathcal J}}
\def\cK{{\mathcal K}}
\def\cL{{\mathcal L}}
\def\cM{{\mathcal M}}
\def\cN{{\mathcal N}}
\def\cO{{\mathcal O}}
\def\cP{{\mathcal P}}
\def\cQ{{\mathcal Q}}
\def\cR{{\mathcal R}}
\def\cS{{\mathcal S}}
\def\cT{{\mathcal T}}
\def\cU{{\mathcal U}}
\def\cV{{\mathcal V}}
\def\cW{{\mathcal W}}
\def\cX{{\mathcal X}}
\def\cY{{\mathcal Y}}
\def\cZ{{\mathcal Z}}

\def\JacNer{{{\it J}_{\rm Ner}}}

\def\ha{{\widehat{a}}}
\def\hb{{\widehat{b}}}
\def\hc{{\widehat{c}}}
\def\hd{{\widehat{d}}}
\def\he{{\widehat{e}}}
\def\hf{{\widehat{f}}}
\def\hg{{\widehat{g}}}
\def\hh{{\widehat{h}}}
\def\hi{{\widehat{i}}}
\def\hj{{\widehat{j}}}
\def\hk{{\widehat{j}}}
\def\hl{{\widehat{l}}}
\def\hm{{\widehat{m}}}
\def\hn{{\widehat{n}}}
\def\ho{{\widehat{o}}}
\def\hp{{\widehat{p}}}
\def\hq{{\widehat{q}}}
\def\hr{{\widehat{r}}}
\def\hs{{\widehat{s}}}
\def\that{{\widehat{t}}}
\def\hu{{\widehat{u}}}
\def\hv{{\widehat{v}}}
\def\hw{{\widehat{w}}}
\def\hx{{\widehat{x}}}
\def\hy{{\widehat{y}}}
\def\hz{{\widehat{z}}}

\def\hA{{\widehat{A}}}
\def\hB{{\widehat{B}}}
\def\hC{{\widehat{C}}}
\def\hD{{\widehat{D}}}
\def\hE{{\widehat{E}}}
\def\hF{{\widehat{F}}}
\def\hG{{\widehat{G}}}
\def\hH{{\widehat{H}}}
\def\hI{{\widehat{I}}}
\def\hJ{{\widehat{J}}}
\def\hK{{\widehat{K}}}
\def\hL{{\widehat{L}}}
\def\hM{{\widehat{M}}}
\def\hN{{\widehat{N}}}
\def\hO{{\widehat{O}}}
\def\hP{{\widehat{P}}}
\def\hQ{{\widehat{Q}}}
\def\hR{{\widehat{R}}}
\def\hS{{\widehat{S}}}
\def\hT{{\widehat{T}}}
\def\hU{{\widehat{U}}}
\def\hV{{\widehat{V}}}
\def\hW{{\widehat{W}}}
\def\hX{{\widehat{X}}}
\def\hY{{\widehat{Y}}}
\def\hZ{{\widehat{Z}}}

\def\hcO{{\widehat{\mathcal O}}}
\def\hcC{{\widehat{\mathcal C}}}
\def\hcF{{\widehat{\mathcal F}}}
\def\hcP{{\widehat{\mathcal P}}}
\def\hcZ{{\widehat{\mathcal Z}}}
\def\hLambda{{\widehat{\Lambda}}}
\def\hdelta{{\widehat{\delta}}}
\def\htheta{{\widehat{\theta}}}
\def\halpha{{\widehat{\alpha}}}
\def\hsigma{{\widehat{\sigma}}}
\def\hGamma{{\widehat{\Gamma}}}
\def\hphi{{\widehat{\varphi}}}
\def\hPhi{{\widehat{\Phi}}}

\def\hhf{{\widehat{f}^{\prime}}}
\def\hhg{{\widehat{g}^{\prime}}}

\def\BCg{{B_{\cC_v^g}}}
\def\sGamma{{\scriptscriptstyle{\Gamma}}}
\def\Sg{S_{\gamma}}

\def\abar{{\overline{a}}}
\def\bbar{{\overline{b}}}
\def\cbar{{\overline{c}}}
\def\dbar{{\overline{d}}}
\def\ebar{{\overline{e}}}
\def\fbar{{\overline{f}}}
\def\gbar{{\overline{g}}}
\def\hbar{{\overline{h}}}
\def\ibar{{\overline{i}}}
\def\jbar{{\overline{j}}}
\def\kbar{{\overline{k}}}
\def\lbar{{\overline{l}}}
\def\mbar{{\overline{m}}}
\def\nbar{{\overline{n}}}
\def\obar{{\overline{o}}}
\def\pbar{{\overline{p}}}
\def\qbar{{\overline{q}}}
\def\rbar{{\overline{r}}}
\def\sbar{{\overline{s}}}
\def\tbar{{\overline{t}}}
\def\ubar{{\overline{u}}}
\def\vbar{{\overline{v}}}
\def\wbar{{\overline{w}}}
\def\xbar{{\overline{x}}}
\def\ybar{{\overline{y}}}
\def\zbar{{\overline{z}}}

\def\Abar{{\overline{A}}}
\def\Bbar{{\overline{B}}}
\def\Cbar{{\overline{C}}}
\def\Dbar{{\overline{D}}}
\def\Ebar{{\overline{E}}}
\def\Fbar{{\overline{F}}}
\def\Gbar{{\overline{G}}}
\def\Hbar{{\overline{H}}}
\def\Ibar{{\overline{I}}}
\def\Jbar{{\overline{J}}}
\def\Kbar{{\overline{K}}}
\def\Lbar{{\overline{L}}}
\def\Mbar{{\overline{M}}}
\def\Nbar{{\overline{N}}}
\def\Obar{{\overline{O}}}
\def\Pbar{{\overline{P}}}
\def\Qbar{{\overline{Q}}}
\def\Rbar{{\overline{R}}}
\def\Sbar{{\overline{S}}}
\def\Tbar{{\overline{T}}}
\def\Ubar{{\overline{U}}}
\def\Vbar{{\overline{V}}}
\def\Wbar{{\overline{W}}}
\def\Xbar{{\overline{X}}}
\def\Ybar{{\overline{Y}}}
\def\Zbar{{\overline{Z}}}

\def\cAbar{{\overline{\cA}}}
\def\bphi{{\overline{\varphi}}}
\def\bPhi{{\overline{\Phi}}}
\def\kvbar{{\overline{k}_v}}
\def\sigmabar{{\overline{\sigma}}}
\def\gammabar{{\overline{\gamma}}}
\def\Gammabar{{\overline{\Gamma}}}
\def\Deltabar{{\overline{\Delta}}}
\def\partialbar{{\overline{\partial}}}
\def\ellbar{{\overline{\ell}}}
\def\cCbar{{\overline{\mathcal C}}}
\def\cZbar{{\overline{\mathcal Z}}}
\def\cPbar{{\overline{\mathcal P}}}
\def\QQbar{{\widetilde{\mathbb Q}}}
\def\rhobar{{\overline{\rho}}}
\def\deltabar{{\overline{\delta}}}

\def\inftyGamma{\infty_{\scriptscriptstyle{\Gamma}}}
\def\inftybar{{\widetilde{\scriptstyle{\infty}}}}

\def\cka{{\check{a}}}
\def\ckb{{\check{b}}}
\def\ckc{{\check{c}}}
\def\ckd{{\check{d}}}
\def\cke{{\check{e}}}
\def\ckf{{\check{f}}}
\def\ckg{{\check{g}}}
\def\ckh{{\check{h}}}
\def\cki{{\check{i}}}
\def\ckj{{\check{j}}}
\def\ckk{{\check{j}}}
\def\ckl{{\check{l}}}
\def\ckm{{\check{m}}}
\def\ckn{{\check{n}}}
\def\cko{{\check{o}}}
\def\ckp{{\check{p}}}
\def\ckq{{\check{q}}}
\def\ckr{{\check{r}}}
\def\cks{{\check{s}}}
\def\ckt{{\check{t}}}
\def\cku{{\check{u}}}
\def\ckv{{\check{v}}}
\def\ckw{{\check{w}}}
\def\ckx{{\check{x}}}
\def\cky{{\check{y}}}
\def\ckz{{\check{z}}}

\def\ckA{{\check{A}}}
\def\ckB{{\check{B}}}
\def\ckC{{\check{C}}}
\def\ckD{{\check{D}}}
\def\ckE{{\check{E}}}
\def\ckF{{\check{F}}}
\def\ckG{{\check{G}}}
\def\ckH{{\check{H}}}
\def\ckI{{\check{I}}}
\def\ckJ{{\check{J}}}
\def\ckK{{\check{K}}}
\def\ckL{{\check{L}}}
\def\ckM{{\check{M}}}
\def\ckN{{\check{N}}}
\def\ckO{{\check{O}}}
\def\ckP{{\check{P}}}
\def\ckQ{{\check{Q}}}
\def\ckR{{\check{R}}}
\def\ckS{{\check{S}}}
\def\ckT{{\check{T}}}
\def\ckU{{\check{U}}}
\def\ckV{{\check{V}}}
\def\ckW{{\check{W}}}
\def\ckX{{\check{X}}}
\def\ckY{{\check{Y}}}
\def\ckZ{{\check{Z}}}

\def\ckcS{{\check{\mathcal{S}}}}
\def\ckPhi{{\check{\Phi}}}

\def\fA{{\mathfrak A}}
\def\fB{{\mathfrak B}}
\def\fC{{\mathfrak C}}
\def\fD{{\mathfrak D}}
\def\fE{{\mathfrak E}}
\def\fF{{\mathfrak F}}
\def\fG{{\mathfrak G}}
\def\fH{{\mathfrak H}}
\def\fI{{\mathfrak I}}
\def\fJ{{\mathfrak J}}
\def\fK{{\mathfrak K}}
\def\fL{{\mathfrak L}}
\def\fM{{\mathfrak M}}
\def\fN{{\mathfrak N}}
\def\fO{{\mathfrak O}}
\def\fP{{\mathfrak P}}
\def\fQ{{\mathfrak Q}}
\def\fR{{\mathfrak R}}
\def\fS{{\mathfrak S}}
\def\fT{{\mathfrak T}}
\def\fU{{\mathfrak U}}
\def\fV{{\mathfrak V}}
\def\fW{{\mathfrak W}}
\def\fX{{\mathfrak X}}
\def\fY{{\mathfrak Y}}
\def\fZ{{\mathfrak Z}}

\def\fa{{\mathfrak{a}}}
\def\fb{{\mathfrak{b}}}
\def\fm{{\mathfrak{m}}}
\def\fp{{\mathfrak{p}}}

\def\ta{{\widetilde{a}}}
\def\tb{{\widetilde{b}}}
\def\tc{{\widetilde{c}}}
\def\td{{\widetilde{d}}}
\def\te{{\widetilde{e}}}
\def\tf{{\widetilde{f}}}
\def\tg{{\widetilde{g}}}
\def\th{{\widetilde{h}}}
\def\ti{{\widetilde{i}}}
\def\tj{{\widetilde{j}}}
\def\tk{{\widetilde{j}}}
\def\tl{{\widetilde{l}}}
\def\tm{{\widetilde{m}}}
\def\tn{{\widetilde{n}}}
\def\to{{\widetilde{o}}}
\def\tp{{\widetilde{p}}}
\def\tq{{\widetilde{q}}}
\def\tr{{\widetilde{r}}}
\def\ts{{\widetilde{s}}}
\def\tt{{\widetilde{t}}}
\def\tu{{\widetilde{u}}}
\def\tv{{\widetilde{v}}}
\def\tw{{\widetilde{w}}}
\def\tx{{\widetilde{x}}}
\def\ty{{\widetilde{y}}}
\def\tz{{\widetilde{z}}}

\def\tA{{\widetilde{A}}}
\def\tB{{\widetilde{B}}}
\def\tC{{\widetilde{C}}}
\def\tD{{\widetilde{D}}}
\def\tE{{\widetilde{E}}}
\def\tF{{\widetilde{F}}}
\def\tG{{\widetilde{G}}}
\def\tH{{\widetilde{H}}}
\def\tI{{\widetilde{I}}}
\def\tJ{{\widetilde{J}}}
\def\tK{{\widetilde{K}}}
\def\tL{{\widetilde{L}}}
\def\tM{{\widetilde{M}}}
\def\tN{{\widetilde{N}}}
\def\tO{{\widetilde{O}}}
\def\tP{{\widetilde{P}}}
\def\tQ{{\widetilde{Q}}}
\def\tR{{\widetilde{R}}}
\def\tS{{\widetilde{S}}}
\def\tT{{\widetilde{T}}}
\def\tU{{\widetilde{U}}}
\def\tV{{\widetilde{V}}}
\def\tW{{\widetilde{W}}}
\def\tX{{\widetilde{X}}}
\def\tY{{\widetilde{Y}}}
\def\tZ{{\widetilde{Z}}}

\def\talpha{{\widetilde{\alpha}}}
\def\tbeta{{\widetilde{\beta}}}
\def\tdelta{{\widetilde{\delta}}}
\def\tDelta{{\widetilde{\Delta}}}
\def\tphi{{\widetilde{\varphi}}}
\def\trho{{\widetilde{\rho}}}
\def\tsigma{{\widetilde{\sigma}}}
\def\ttheta{{\widetilde{\theta}}}
\def\tvartheta{{\widetilde{\vartheta}}}
\def\tomega{{\widetilde{\omega}}}
\def\tLambda{{\widetilde{\Lambda}}}
\def\tlambda{{\widetilde{\lambda}}}
\def\tcR{{\widetilde{\mathcal R}}}
\def\tcB{{\widetilde{\mathcal B}}}
\def\tGamma{{\widetilde{\Gamma}}}
\def\tPhi{{\widetilde{\Phi}}}
\def\tOmega{{\widetilde{\Omega}}}
\def\tTheta{{\widetilde{\Theta}}}
\def\tUpsilon{{\widetilde{\Upsilon}}}
\def\tEE{{\widetilde{\mathbb E}}}

\def\ulgamma{{\underline{\gamma}}}
\def\ulG{{\underline{G}}}

\def\blacksquare{{\Box}}

\newcommand{\PL}{\operatorname{PL}}
\newcommand{\ord}{\operatorname{ord}}
\newcommand{\Aut}{\operatorname{Aut}}
\newcommand{\Gal}{\operatorname{Gal}}
\newcommand{\val}{\operatorname{val}}
\newcommand{\supp}{\operatorname{supp}}
\renewcommand{\div}{\operatorname{div}}
\newcommand{\Div}{{\operatorname{Div}}}
\newcommand{\BDiv}{\operatorname{\bf Div}}
\newcommand{\Det}{\operatorname{Det}}
\newcommand{\vol}{\operatorname{vol}}
\newcommand{\covol}{\operatorname{covol}}
\renewcommand{\dim}{{\operatorname{dim}}}
\renewcommand{\exp}{{\operatorname{exp}}}
\renewcommand{\vol}{{\operatorname{vol}}}
\renewcommand{\deg}{{\operatorname{deg}}}
\newcommand{\lcm}{\operatorname{lcm}}
\newcommand{\Sym}{\operatorname{Sym}}
\newcommand{\Spm}{\operatorname{Spm}}
\newcommand{\Proj}{\operatorname{Proj}}
\newcommand{\BProj}{\operatorname{\bf Proj}}
\newcommand{\Spec}{\operatorname{Spec}}
\newcommand{\CPA}{\operatorname{CPA}}
\newcommand{\GL}{\operatorname{GL}}
\newcommand{\PGL}{\operatorname{PGL}}
\newcommand{\LCM}{\operatorname{LCM}}
\newcommand{\Meas}{\operatorname{Meas}}
\renewcommand{\Re}{\operatorname{Re}}
\renewcommand{\Im}{\operatorname{Im}}
\newcommand{\Res}{\operatorname{Res}}
\newcommand{\Tr}{\operatorname{Tr}}
\newcommand{\id}{\operatorname{id}}
\newcommand{\Jac}{\operatorname{Jac}}
\newcommand{\Char}{\operatorname{char}}
\newcommand{\an}{\operatorname{an}}
\newcommand{\nat}{\operatorname{nat}}
\newcommand{\sep}{\operatorname{sep}}
\newcommand{\insep}{\operatorname{insep}}
\newcommand{\cl}{\operatorname{cl}}
\newcommand{\Adj}{\operatorname{Adj}}
\newcommand{\Pic}{\operatorname{Pic}}
\newcommand{\BPic}{\operatorname{\bf Pic}}
\newcommand{\BHilb}{\operatorname{\bf Hilb}}
\newcommand{\CH}{\operatorname{CH}}
\newcommand{\Hom}{\operatorname{Hom}}
\newcommand{\Max}{\operatorname{Max}}
\newcommand{\Ker}{\operatorname{Ker}}
\newcommand{\Ann}{\operatorname{Ann}}
\newcommand{\Mod}{\operatorname{Mod}}
\newcommand{\radius}{\operatorname{rad}}
\newcommand{\diam}{\operatorname{diam}}
\newcommand{\Band}{\operatorname{Band}}
\newcommand{\Block}{\operatorname{Block}}
\newcommand{\Sp}{\operatorname{Sp}}
\newcommand{\sn}{\operatorname{sn}}
\newcommand{\diag}{\operatorname{diag}}

\newcommand{\ii}[1]{{\bf #1}}

\newcommand{\Zh}{\operatorname{Zh}}
\newcommand{\BDV}{\operatorname{BDV}}
\newcommand{\SH}{\operatorname{SH}}

\newcommand{\BerkA}{{\bf A}}
\newcommand{\BerkB}{{\bf B}}
\newcommand{\BerkE}{{\bf E}}
\newcommand{\BerkF}{{\bf F}}
\newcommand{\BerkK}{{\bf K}}
\newcommand{\BerkS}{{\bf S}}
\newcommand{\BerkU}{{\bf U}}
\newcommand{\BerkX}{{\bf X}}
\newcommand{\BerkY}{{\bf Y}}
\newcommand{\ThCap}{\operatorname{Cap}}

\newcommand{\into}{\hookrightarrow}     
\newcommand{\onto}{\twoheadrightarrow}  
\def\isomap{{\buildrel \sim\over\longrightarrow}} 
\newcommand{\dirlim}{\varinjlim}
\newcommand{\invlim}{\varprojlim}

\def\Berk{{\rm Berk}}
\def\Gauss{{\rm Gauss}}
\def\dL{{d}}
\def\nz{{m}}

\def\IN{{\mathop{\rm In}}}
\def\OUT{{\mathop{\rm Out}}}

\newcommand{\Fpbar}{{\overline{\FF}_p}}
\newcommand{\RL}{{\rm RL}}


\def\rL{{\rm L}}
\def\rF{{\rm F}}

\def\dL{{d}}
\def\nz{{m}}

\def\Sg{S_{\gamma}}

\def\BA{{\bf{A}}}
\def\BP{{\bf{P}}}

\def\bj{{\bf j}}
\def\bfj{{\bf{j}}}
\def\crJ{{\bf J}}

\def\blacksquare{\Box}

\def\<{{\langle }}
\def\>{{\rangle }}

\def\<<{{\langle \! \langle}}
\def\>>{{\rangle \! \rangle}} 

\def\({(\!(}
\def\){)\!)}

\def\[{[\![}
\def\]{]\!]}

\def\lH{\{\!|}
\def\rH{|\!\}}

\def\ldB{[\![}
\def\rdB{]\!]}

\def\Jp{{\dot{+}}}
\def\Jm{{\dot{-}}}
\def\ap{{\ddot{+}}}
\def\am{{\ddot{-}}}

\newcommand{\fixme}[1]{{\sf $\spadesuit\spadesuit$ FIXME! [#1]}} 

\newenvironment{notation}[0]{%
  \begin{list}%
    {}%
    {\setlength{\itemindent}{0pt}
     \setlength{\labelwidth}{4\parindent}
     \setlength{\labelsep}{\parindent}
     \setlength{\leftmargin}{5\parindent}
     \setlength{\itemsep}{0pt}
     }%
   }%
  {\end{list}}

\frontmatter

\title[The Fekete-Szeg\"o Theorem]{The Fekete-Szeg\"o Theorem \\ with Local Rationality Conditions \\
on Algebraic Curves} 

\author[Robert Rumely]{Robert Rumely}
       



\address{ \ }

\address{\centerline{{\rm March 6, 2012 \ \ \ \ \ \ \ \ }}} 

\address{ \ }

\address{\centerline{Abstract: \ \ \ \ \ \  }}

\address{{ \rm Let $K$ be a number field or a function field in one variable over a finite field, 
and let $\tK^{\sep}$ be a separable closure of $K$.  Let $\cC/K$ be a smooth, complete, connected curve.  
We prove a strong theorem of Fekete-Szeg\"o type for adelic sets $\EE = \prod_v E_v$ on $\cC$ which 
satisfy local rationality conditions at finitely many places $v$ of $K$, 
showing that under appropriate conditions there are infinitely many points 
in $\cC(\tK^{\sep})$ whose conjugates all belong to $E_v$ at each place $v$.  
We give several variants of the theorem, including two for Berkovich curves,
and we provide examples illustrating the theorem on the projective line and on 
elliptic curves, Fermat curves, and modular curves.}}

\address{ \ }

\address{{\rm AMS Subject Classifications (2010):  
Primary 11G30, 14G40, 14G05; Secondary 31C15}}

\address{{\rm Keywords and Phrases:  Fekete-Szeg\"o theorem, Local rationality conditions, 
Algebraic curves, Adelic capacity theory, Cantor capacity, Global fields, Conjugate sets of points}}

\address{ \ } 

\address{\ \ \ \ \ \ \ \ \ \ \ \ \ \ \ \ \ \ \ \ \ \ \ \ \ \ \ \ \ Department of Mathematics} 

\address{\ \ \ \ \ \ \ \ \ \ \ \ \ \ \ \ \ \ \ \ \ \ \ \ \ \ \ \ \ University of Georgia} 

\address{\ \ \ \ \ \ \ \ \ \ \ \ \ \ \ \ \ \ \ \ \ \ \ \ \ \ \ \ \ Athens, Georgia  30602-7403, USA}

\address{\ \ \ \ \ \ \ \ \ \ \ \ \ \ \ \ \ \ \ \ \ \ \ \ \ \ \ \ \ {\rm rr@math.uga.edu}}

\address{ \ }

\address{\rm{This work was supported in part by NSF grants DMS 95-000892, DMS 00-70736, DMS 03-00784, and DMS 06-01037. 
Any opinions, findings and conclusions or recommendations expressed in this material are those of the author 
and do not necessarily reflect the views of the National Science Foundation.} } 





\thanks{To Cherilyn, who makes me happy.}




\maketitle

\tableofcontents

\input{NewFSZIntro}

\mainmatter

\input{NewFSZChap1}

\input{NewFSZChap2}

\input{NewFSZChap3}

\input{NewFSZChap4}

\input{NewFSZChap5}

\input{NewFSZChap6}

\input{NewFSZChap7}
\input{NewFSZChap8}

\input{NewFSZChap9}
\input{NewFSZChap10}

\input{NewFSZChap11}

\appendix

\input{NewFSZAppA} 
\input{NewFSZAppB} 
\input{NewFSZAppC}

\input{NewFSZAppD}

\backmatter

\input{NewFSZBiblio}
\printindex



\end{document}

%% file: NewFSZIntro.tex
\chapter*{Introduction}
\label{IntroChapter}


One of the gems of mid-twentieth century mathematics was Raphael Robinson's theorem 
on totally real algebraic integers in a closed interval $[a,b]$:  

\vskip .08 in
\noindent{\bf Theorem (Robinson \cite{Rob1}, 1964).} 
\index{Robinson, Raphael}  
{\it Let $a < b \in \RR$.  If $b-a > 4$, 
then there are infinitely many totally real algebraic integers whose
\index{algebraic integer!totally real}  
conjugates all belong to the interval $[a,b]$.
If $b-a < 4$, there are only finitely many.}

\vskip .08 in
\noindent{Four years later,} he gave a criterion for the existence of 
totally real units in $[a,b]$: 
\index{units!totally real}

\vskip .08 in
\noindent{\bf Theorem (Robinson \cite{Rob2}, 1968).}
{\it Suppose $0 < a < b \in \RR$ satisfy the conditions }
\begin{eqnarray}
\log(\frac{b-a}{4}) \ > \ 0 \ , \qquad \qquad \qquad \quad \label{FRobC1A}\\
\log(\frac{b-a}{4}) \cdot \log(\frac{b-a}{4ab})  -  
        \Big(\log(\frac{\sqrt{b}+\sqrt{a}}{\sqrt{b}-\sqrt{a}})\Big)^2 
       \ > \ 0 \ .  \label{FRobC2A} 
\end{eqnarray}
{\it Then there are infinitely many totally real units $\alpha$ 
whose conjugates all belong to $[a,b]$.
If either inequality is reversed, there are only finitely many.}  

\vskip .08 in
David Cantor's ``Fekete-Szeg\"o theorem with splitting conditions'' on the projective line (\cite{Can1}, 1980)
\index{Fekete-Szeg\"o theorem with splitting conditions}\index{Cantor, David}\index{Fekete, Michael}\index{Szeg\"o, G\'abor}
formulated Robinson's theorems adelically and set them in a potential-theoretic framework.   
In this work we generalize the Fekete-Szeg\"o theorem 
with splitting conditions\index{Fekete-Szeg\"o theorem with splitting conditions} 
to algebraic curves.  Below we state the theorem, recall some history, and outline its proof.  

\vskip .1 in
Let $K$ be a global field, that is, 
a number field or a finite extension of $\FF_p(T)$ for some prime $p$.  
Let $\tK$ be a fixed algebraic closure of $K$, 
and let $\tK^{\sep} \subseteq \tK$ be the separable closure of $K$.
\label{`SymbolIndextKsep'}
We will write $\Aut(\tK/K)$ for the group of automorphisms $\Aut(\tK/K) \cong \Gal(\tK^{\sep}/K)$.    
Let $\cM_K$ be the set of all places of $K$.
For each $v \in \cM_K$, 
let $K_v$ be the completion of $K$ at $v$, 
let $\tK_v$ be an algebraic closure of $K_v$, 
and let $\CC_v$ be the completion of $\tK_v$. 
We will write $\Aut_c(\CC_v/K_v)$ for the group of continuous automorphisms of $\CC_v/K_v$;  
thus $\Aut_c(\CC_v/K_v) \cong \Aut(\tK_v/K_v) \cong \Gal(\tK_v^{\sep}/K_v)$.

Let $\cC/K$ be a smooth, geometrically integral, projective curve. 
Given a field $F$ containing $K$, put $\cC_F = \cC\times_K\Spec(F)$
and write $\cC(F)$ for the set of $F$-rational points $\Hom_F(\Spec(F),\cC_F)$; 
write $F(\cC)$ for its function field.  
When $F = K_v$, \label{`SymbolIndexcCv'} write $\cC_v$ for $\cC_{K_v}$. 

Let $\fX = \{x_1, \ldots, x_m\}$ be a finite, galois-stable set 
points of $\cC(\tK)$, and let $\EE = \EE_K = \prod_{v \in \cM_K} E_v$
be a $K$-rational adelic set for $\cC$, that is, a product of sets
$E_v \subset \cC_v(\CC_v)$ such that each $E_v$ is stable under $\Aut_c(\CC_v/K_v)$.  
\label{`SymbolIndexEv'}
For each $v$, fix an embedding $\tK \hookrightarrow \CC_v$ over $K$, 
inducing an embedding $\cC(\tK) \hookrightarrow \cC_v(\CC_v)$.  
In this way $\fX$ can be regarded as a subset of $\cC_v(\CC_v)$:  
since $\fX$ is galois-stable, 
its image is well-defined, independent of the choice of embedding.  
The same is true for any other galois-stable set of points in $\cC(\tK)$,
for instance, the set of $\Aut(\tK/K)$-conjugates of a given point
$\alpha \in \cC(\tK)$. 

We will call a set $E_v \subset \cC_v(\CC_v)$ an $\RL$-domain (`Rational Lemniscate Domain')
\index{$\RL$-domain|ii}  
if there is a nonconstant rational function $f_v(z) \in \CC_v(\cC_v)$
such that $E_v = \{ z \in \cC_v(\CC_v) : |f_v(z)|_v \le 1 \}$.  This terminology is due to Cantor.
\index{Cantor, David}   
By combining (\cite{Fies}, Satz 2.2)  with (\cite{RR1}, Corollary 4.2.14),
one sees that a set is an $\RL$-domain if and only if it is a strict affinoid subdomain of $\cC_v(\CC_v)$, 
\index{$\RL$-domain} 
\index{affinoid!affinoid domain} 
in the sense of rigid analysis. 

Fix an embedding $\cC \hookrightarrow \PP^N_K = \PP^N/\Spec(K)$ for an appropriate $N$, 
and equip $\PP^N_K$ with a system of homogeneous coordinates. 
For each nonarchimedean $v$, this data determines a model $\fC_v/\Spec(\cO_v)$.  
There is a natural metric $\|x,y\|_v$ on $\PP^N_v(\CC_v)$:  
the chordal distance associated to the Fubini-Study metric, 
\index{Fubini-Study metric}\index{chordal distance} 
if $v$ is archimedean;
the $v$-adic spherical metric, if $v$ is nonarchimedean
\index{spherical metric} 
(see \S\ref{Chap3}.\ref{SphericalMetricSection} below).
The metric $\|x,y\|_v$ induces the $v$-topology on $\cC_v(\CC_v)$.  
Given $a \in \cC_v(\CC_v)$ and $r > 0$, we write   
$B(a,r)^- = \{z \in \cC_v(\CC_v) : \|z,a\|_v < r \}$
and $B(a,r) = \{z \in \cC_v(\CC_v) : \|z,a\|_v \le r \}$
for the corresponding `open' and `closed' balls.    

\begin{definition} \label{Xtrivial} 
If $v$ is a nonarchimedean place of $K$, 
a set $E_v \subset \cC_v(\CC_v)$ will be called {\em $\fX$-trivial}  
\index{$\fX$-trivial} 
if $\fC_v$ has good reduction at $v$, if
\index{good reduction} 
the points of $\fX$ specialize to distinct points$\pmod v$,
and if $E_v = \cC_v(\CC_v) \backslash  \bigcup_{i=1}^m B(x_i,1)^-$.
\end{definition}
\index{$\fX$-trivial}

\begin{definition} \label{Compatible} 
An adelic set $\EE = \prod_{v \in \cM_K} E_v \subset \prod_{v \in \cM_K} \cC_v(\CC_v)$ \label{`SymbolIndexEEv'}
will be called {\em compatible with $\fX$} if the following conditions hold:
\index{compatible with $\fX$}   

$(1)$  Each $E_v$ is bounded away from $\fX$ in the $v$-topology$;$

$(2)$  For all but finitely many $v$, $E_v$ is $\fX$-trivial.
\index{$\fX$-trivial}
\end{definition}  

If $E_v$ is $\fX$-trivial, 
\index{$\fX$-trivial}
it consists of all points of $\cC_v(\CC_v)$ 
which are $\fX$-integral at $v$ for the model $\fC_v$, 
i.e. which specialize to points complementary to $\fX \pmod v$.
If $E_v$ is $\fX$-trivial, it an $\RL$-domain and is stable under $\Aut_c(\CC_v/K_v)$.
\index{$\fX$-trivial}
\index{$\RL$-domain} 
The property of compatibility is 
independent of the embedding $\cC \hookrightarrow \PP^N_K$
and the choice of coordinates on $\PP^N_K$.

There is a potential-theoretic measure of size for the adelic set 
$\EE$ relative to the set of global points $\fX$: 
the Cantor capacity $\gamma(\EE,\fX)$, defined in (\ref{GlobalCapacityDef}) below. 
\index{capacity!Cantor capacity} 
Our main result is:     

\begin{theorem}[Fekete-Szeg\"o Theorem with Local Rationality Conditions, 
producing points in $\EE$]\label{aT1} 
\index{Fekete-Szeg\"o theorem with LRC!producing points in $\EE$|ii} 
\index{Fekete, Michael}\index{Szeg\"o, G\'abor}
 
Let $K$ be a global field, 
and let $\cC/K$ be a smooth, geometrically integral, projective curve.
Let $\fX = \{x_1, \ldots, x_m\} \subset \cC(\tK)$ 
be a finite set of points stable under $\Aut(\tK/K)$, and let
$\EE = \prod_v E_v \subset \prod_v \cC_v(\CC_v)$ be an adelic set compatible with $\fX$.
\index{compatible with $\fX$}  
Let $S \subset \cM_K$ be a finite set of places $v$, containing all archimedean $v$,
such that $E_v$ is $\fX$-trivial for each $v \notin S$.
\index{$\fX$-trivial}

Assume that $\gamma(\EE,\fX) > 1$.  Assume also that $E_v$ has the following form, for each $v \in S$:   

$(A)$ If $v$ is archimedean and $K_v \cong \CC$, 
then $E_v$ is compact, and is a finite union of sets $E_{v,\ell}$, 
each of which is the closure of its $\cC_v(\CC)$-interior 
and has a piecewise smooth boundary;\index{boundary!piecewise smooth}\index{closure of $\cC_v(\CC)$ interior} 

$(B)$ If $v$ is archimedean and $K_v \cong \RR$, then $E_v$ is compact, stable under complex conjugation, 
and is a finite union of sets $E_{v,\ell}$, where each $E_{v,\ell}$ is either 

\quad $(1)$ the closure of its $\cC_v(\CC)$-interior and has a piecewise smooth boundary, 
or\index{boundary!piecewise smooth}\index{closure of $\cC_v(\CC)$ interior}  

\quad $(2)$ is a compact, connected subset of $\cC_v(\RR)$; 

$(C)$ If $v$ is nonarchimedean, then $E_v$ is stable under $\Aut_c(\CC_v/K_v)$ 
and is a finite union of sets $E_{v,\ell}$, where each $E_{v,\ell}$ is either 

\quad $(1)$ an $\RL$-domain or a ball $B(a_\ell,r_\ell)$, or
\index{$\RL$-domain} 

\quad $(2)$ is compact and has the form $\cC_v(F_{w_\ell}) \cap B(a_\ell,r_\ell)$ 
for some finite separable extension $F_{w_\ell}/K_v$ in $\CC_v$, and some ball $B(a_\ell,r_\ell)$.  

Then there are infinitely many points $\alpha \in \cC(\tK^{\sep})$ such that for each $v \in \cM_K$, 
the $\Aut(\tK/K)$-conjugates of $\alpha$ all belong to $E_v$.  
\end{theorem} 

Note that for a given $v$, the extensions $F_{w_\ell}/K_v$ need not be galois,
the sets $E_{v,\ell}$ may overlap, and sets $E_{v,\ell}$ of more than one type 
(intervals, sets with nonempty interior, $\RL$-domains, balls, compact sets) may occur.      
The main content of the theorem is the satisfiability of the local rationality conditions
(the fact that the $E_{v,\ell}$ can be taken to be subsets of the $\cC_v(F_{w,\ell})$ 
and the conjugates belong to $E_v$, for each $v$);
the Fekete-Szeg\"o theorem without local rationality conditions,
\index{Fekete-Szeg\"o theorem}\index{Fekete, Michael}\index{Szeg\"o, G\'abor} 
which constructs points $\alpha$ whose conjugates belong to arbitrarily small $\cC_v(\CC_v)$ 
neighborhoods of $E_v$, was proved in (\cite{RR1}, Theorem 6.3.2). 
 In \S\ref{ExamApp}.\ref{FFSeparabilitySection} we provide examples due to 
Daeshik Park,\index{Park, Daeshik}\index{Fekete-Szeg\"o theorem with LRC!need for separability hypotheses in} 
showing the need for the hypothesis of separability for the extensions $F_{w_\ell}/K_v$ in (C2) of the theorem.

Suppose that in the theorem, for each $v \in S$ we have $E_v \subset \cC_v(K_v)$.  
Then for each $v \in S$, the conjugates of $\alpha$ belong to $\cC_v(K_v)$, 
which means that $v$ splits completely in $K(\alpha)$.  
In this case, we speak of ``the Fekete-Szeg\"o theorem with splitting conditions''.
\index{Fekete-Szeg\"o theorem with splitting conditions} 
\index{Fekete, Michael}\index{Szeg\"o, G\'abor}

Often it is the corollaries of a theorem, which are weaker but easier to apply, 
that are most useful.  The following consequence of the Fekete-Szeg\"o theorem 
with local rationality conditions 
strengthens Laurent Moret-Bailly's theorem on ``Incomplete Skolem Problems'' 
(\cite{MB2} Th\'eor\`eme 1.3, p.182) for curves, 
but does not require evaluating capacities.\index{Moret-Bailly, Laurent}\index{Incomplete Skolem Problems} 

\begin{corollary}[Fekete-Szeg\"o with LRC, for Incomplete Skolem Problems] \label{FSZiv} 
Let $K$ be a global field, 
and let $\cA/K$ be a geometrically integral $($possibly singular$)$ affine curve, 
embedded in $\AA^N$ for some $N$. 
Let $z_1, \ldots, z_N$ be the coordinates on $\AA^N;$  
given a place $v$ of $K$ and a point $P \in \AA^N(\CC_v)$, 
write $\|P\|_v = \max(|z_1(P)|_v, \ldots, |z_N(P)|_v)$.  

Fix a place $v_0$ of $K$, and let $S \subset \cM_K \backslash \{v_0\}$ be a finite set of places containing 
all archimedean $v \ne v_0$.     
For each $v \in S$, let a nonempty set $E_v \subset \cA_v(\CC_v)$
satisfying condition $(A)$, $(B)$ or $(C)$ of Theorem $\ref{aT1}$ be given,    
and put $\EE_S = \prod_{v \in S} E_v$.  
Assume that for each $v \in \cM_K \backslash (S \cup \{v_0\})$ there is a 
point $P \in \cA(\CC_v)$ with $\|P\|_v \le 1$. 
Then there is a constant $C = C(\cA,\EE_S,v_0)$
such that there are infinitely many points $\alpha \in \cA(\tK^{\sep})$ for which  

$(1)$ for each $v \in S$, all the conjugates of $\alpha$ in $\cA_v(\CC_v)$
belong to $E_v;$   

$(2)$ for each $v \in \cM_K \backslash (S \cup \{v_0\})$, 
all the conjugates of $\alpha$ in $\cA_v(\CC_v)$ satisfy $\|\sigma(\alpha)\|_{v_0} \le 1;$   

$(3)$ for $v = v_0$, all the conjugates of $\alpha$ in $\cA_{v_0}(\CC_{v_0})$ satisfy $\|\sigma(\alpha)\|_{v_0} \le C$.      
\end{corollary}  

\smallskip
In Chapter \ref{Variants} below, we will give several variants of Theorem \ref{aT1},  
including one involving ``quasi-neighborhoods'' analogous to the classical theorem of Fekete and Szeg\"o,
\index{quasi-neighborhood}\index{Fekete-Szeg\"o theorem}\index{Fekete, Michael}\index{Szeg\"o, G\'abor} 
one for more general sets $\EE$ using the inner Cantor capacity $\gammabar(\EE,\fX)$,
and two for sets on Berkovich curves.\index{Berkovich!curve}
Theorem \ref{aT1}, Corollary \ref{FSZiv}, and the variants in Chapter \ref{Variants} will be proved in
Chapter \ref{Chap4}.

\medskip
\section*{Some History} \label{Sec0A}

 
The original theorem of Fekete and Szeg\"o (\cite{F-SZ}, 1955) said  
\index{Fekete-Szeg\"o theorem} 
\index{Fekete, Michael}
\index{Szeg\"o, G\'abor}
that if $E \subset \CC$ is a compact set, stable under complex conjugation, 
with logarithmic capacity $\gamma_{\infty}(E) > 1$, 
then every neighborhood $U$ of $E$ contains infinitely many conjugates sets  
of algebraic integers.  (The neighborhood $U$ is needed to `fatten' 
sets like a circle $E = C(0,r)$ with transcendental radius $r$, 
which contain no algebraic numbers.)

A decade later Raphael Robinson\index{Robinson, Raphael} 
gave the generalizations of the 
Fekete-Szeg\"o theorem\index{Fekete-Szeg\"o theorem}\index{Fekete, Michael}\index{Szeg\"o, G\'abor} 
for totally real algebraic integers\index{algebraic integer!totally real} 
and totally real units stated above.\index{units!totally real}  
Independently, Bertrandias gave an adelic generalization of the Fekete-Szeg\"o 
theorem\index{Bertrandias, Fran\c{c}oise}\index{Fekete, Michael}\index{Szeg\"o, G\'abor}
concerning algebraic integers with conjugates near sets $E_p$
at a finite number of $p$-adic places as well as the archimedean place (see Amice \cite{Am}, 1975).  

In the 1970's David Cantor carried out an investigation of capacities\index{Cantor, David} 
on $\PP^1$ dealing with all three themes:  incorporating local rationality conditions, 
requiring integrality with respect to multiple poles, and formulating the theory adelically.   
In a series of papers  
culminating with (\cite{Can3}, 1980), he introduced  
the Cantor capacity $\gamma(\EE,\fX)$, 
which he called  the extended transfinite diameter.
\index{transfinite diameter!extended} 

Cantor's capacity $\gamma(\EE,\fX)$ is defined by means of a minimax property 
\index{minimax property} 
\index{capacity!Cantor capacity}  
which encodes a finite collection of linear inequalities;   
its definition is given in (\ref{GlobalCapacityDef}) below.  
The points in $\fX$ will be called the {\it poles} for the capacity.
In the special case where $\cC = \PP^1$ and $\fX = \{0,\infty\}$, Cantor's conditions
are equivalent those in Robinson's unit theorem.  
Among the applications Cantor gave in\index{Cantor, David} (\cite{Can3}) 
were generalizations of the P\'olya-Carlson theorem and Fekete's 
theorem,\index{P\'olya-Carlson theorem}\index{Fekete's theorem}\index{Fekete, Michael}\index{Szeg\"o, G\'abor}
and the Fekete-Szeg\"o theorem with splitting 
conditions.\index{Fekete-Szeg\"o theorem with splitting conditions} 
Unfortunately, as noted in (\cite{RR3}), 
the proof of the satisfiability of the splitting conditions had gaps.   

In the 1980's the author (\cite{RR1}) extended Cantor's theory to curves  
of arbitrary genus, and proved the Fekete-Szeg\"o theorem on curves, 
\index{Fekete-Szeg\"o theorem} {\it without} splitting conditions. 
As an application he obtained a local-global principle for the
existence of algebraic integer points on absolutely irreducible 
\index{algebraic integer} 
affine algebraic varieties (\cite{RR5}), which had been conjectured
by Cantor and Roquette (\cite{C-R}).  
\index{Cantor, David}
\index{Roquette, Peter} 

Moret-Bailly and Szpiro recognized that the
\index{Moret-Bailly, Laurent} 
\index{Szpiro, Lucien}  
theory of capacities (which imposes conditions at all places) 
was stronger than was needed for the existence of integral points.  
They reformulated the local-global principle 
in scheme-theoretic language as an ``Existence Theorem'' 
for algebraic integer points, and gave a much simpler proof. 
\index{algebraic integer}  
Moret-Bailly subsequently gave far-reaching generalizations
\index{Moret-Bailly, Laurent} 
of the Existence Theorem (\cite{MB1}, \cite{MB2}, \cite{MB3}), which allowed
imposition of $F_w$-rationality conditions at a finite number of places, 
for a finite galois extension $F_w/K_v$,  
and applied to algebraic stacks as well as schemes.  
However, the method required that there be at least one place $v_0$ where
no conditions are imposed. Roquette, Green, and Pop (\cite{R-G-P}) 
\index{Roquette, Peter}\index{Green, Barry}\index{Pop, Florian} 
independently proved the Existence Theorem
with $F_w$-rationality conditions, and Green, Matignon, and Pop (\cite{G-M-P}) 
have given very general conditions on the base field $K$ for such theorems to hold.  
Rumely (\cite{RR5}), van den Dries (\cite{vdD}), Prestel-Schmidt (\cite{P-S}),
and others\index{Rumely, Robert}\index{van den Dries, Lou}\index{Prestel, Alexander}\index{Schmidt, Joachim}
have given applications of these results to decision procedures in mathematical logic.  

Recently Tamagawa\index{Tamagawa, Akio} (\cite{Tam}) proved an extension of the Existence Theorem in characteristic $p$, 
which produces points that are unramified outside $v_0$ and the places where the $F_w$-rationality conditions are imposed.  
 
\smallskip
The Fekete-Szeg\"o theorem  with local rationality conditions constructs algebraic numbers 
\index{Fekete-Szeg\"o theorem with LRC} 
satisfying conditions at {\it all} places.  At its core it is analytic
in character, while the Existence Theorem is algebraic. 
The proof of the Fekete-Szeg\"o theorem involves a process 
called ``patching'', which takes an 
\index{patching} 
initial collection of local functions $f_v(z) \in K_v(\cC)$ 
with poles supported on $\fX$ and roots in $E_v$ for each $v$, 
and constructs  
a global function $G(z) \in K(\cC)$ (of much higher degree) with poles
supported on $\fX$, whose roots belong to $E_v$ for all $v$. 
In his doctoral thesis, Pascal Autissier (\cite{Aut}) gave a reformulation of the patching
\index{Autissier, Pascal}\index{patching} 
process in the context of Arakelov theory.
\index{Arakelov theory}   
   
In (\cite{RR2}, \cite{RR3}) the author proved the Fekete-Szeg\"o 
\index{Fekete-Szeg\"o theorem with splitting conditions} 
theorem with splitting conditions for sets 
$\EE$ in $\PP^1$, when $\fX = \{\infty\}$.  Those papers 
developed a method for carrying out the patching process in the $p$-adic  
compact case, and introduced a technique for patching together archimedean 
and nonarchimedean polynomials over number fields.
\index{patching argument!earlier results}   

When $\cC = \PP^1/K$, 
with $K$ a finite extension of $\FF_p(T)$, the Fekete-Szeg\"o theorem with splitting conditions
\index{Fekete-Szeg\"o theorem with splitting conditions}
was established in the doctoral thesis of Daeshik Park (\cite{DPark}).  
\index{Park, Daeshik} 

\vskip .1 in
\section*{A Sketch of the Proof of the Fekete-Szeg\"o Theorem} \label{Sec0B}

In outline, the proof of the classical Fekete-Szeg\"o theorem (\cite{F-SZ}, 1955)
\index{Fekete-Szeg\"o theorem} 
is as follows.  Let a compact set $E \subset \CC$
and a complex neighborhood $U$ of $E$ be given.  
Assume $E$ is stable under complex conjugation, 
and has logarithmic capacity $\gamma_{\infty}(E) > 1$.
For simplicity, assume also that the boundary\index{boundary!piecewise smooth} of $E$ is piecewise smooth 
and the complement of $E$ is connected.   

Under these assumptions, 
there is a real-valued function $G(z,\infty;E)$,
called the Green's function of $E$ respect to $\infty$,  
\index{Green's function|ii} 
which is continuous on $\CC$, $0$ on $E$, 
harmonic and positive in $\CC \backslash E$, 
and has the property that 
$G(z,\infty;E)-\log(|z|)$ is bounded as $z \rightarrow \infty$.  
(We write $\log(x)$ for $\ln(x)$.)\index{$\log(x)$!same as natural logarithm|ii}\index{logarithm!$\log(x) = \ln(x)$|ii}  
The theorem on removable singularities for harmonic 
functions shows that the {\em Robin constant}, defined by  
\index{Robin constant}
\begin{equation*}
V_{\infty}(E) \ = \ \lim_{z \rightarrow \infty} G(z,\infty;E) - \log(|z|) \ ,  
\end{equation*} 
\nobreak
exists.  By definition $\gamma_{\infty}(E) = e^{-V_{\infty}(E)}$; 
our assumption that $\gamma_{\infty}(E) > 1$ means $V_{\infty}(E) < 0$.  
It can be shown that $V_{\infty}(E)$ is the minimum possible value of the `energy integral'
\index{energy integral!classical|ii} 
\begin{equation*}
I_{\infty}(\nu) \ = \ \iint_{E \times E} -\log(|z-w|) \, d\nu(z) d\nu(w)
\end{equation*} 
as $\nu$ ranges over all probability measures supported on $E$.  
There is a unique probability measure $\mu_{\infty}$ on $E$,
called the equilibrium distribution of $E$ with respect to $\infty$,
\index{equilibrium distribution!classical|ii}\index{Robin constant!classical|ii} 
for which 
\begin{equation*}
V_{\infty}(E) \ = \ \iint_{E \times E} -\log(|z-w|) \, d\mu_{\infty}(z) d\mu_{\infty}(w) \ . 
\end{equation*}
The Green's function is related to the equilibrium distribution by 
\begin{equation*}
G(z,\infty;E) - V_{\infty}(E) \ = \ \int_E \log(|z-w|) \, d\mu_{\infty}(w) \ . 
\end{equation*}\index{Green's function!of a compact set|ii}

Because of its uniqueness, the measure $\mu_{\infty}$ is stable under complex conjugation.  
Taking a suitable discrete approximation 
$\mu_N = \frac{1}{N} \sum_{i=1}^N \delta_{x_i}(z)$ to $\mu_{\infty}$, 
stable under complex conjugation, 
one obtains a monic polynomial $f(z) = \prod_{i=1}^N (z-x_i) \in \RR[z]$
such that $\frac{1}{N}\log(|f(z)|)$ 
approximates $G(z,\infty,E) - V_{\infty}(E)$
very well outside $U$.  If the approximation is good enough, 
then since $V_{\infty}(E) < 0$, there will be 
an $\varepsilon > 0$ such that $\log(|f(z)|) > \varepsilon$ 
outside $U$.  

One then uses the polynomial $f(z) \in \RR[z]$ 
to construct a monic polynomial $G(z) \in \ZZ[z]$ of much higher degree, 
which has properties similar to those of $f(z)$.  
The construction is as follows.  By adjusting the
coefficients of $f(z)$ to be rational numbers and using continuity, 
\index{coefficients $A_{v,ij}$}
one first obtains  a polynomial $f(z) \in \QQ[z]$.  
and an $R > 1$ such that $|f(z)| \ge R$ outside $U$.  
For suitably chosen $n$, the multinomial theorem implies that
\index{multinomial theorem}
$f(z)^n$ will have a pre-designated number of high-order coefficients in $\ZZ$. 
\index{coefficients $A_{v,ij}$!high-order} 
By successively modifying the remaining coefficients 
\index{coefficients $A_{v,ij}$} 
of $G^{(0)}(z) := P(z)^n$ from highest to lowest order, 
writing $k = mN+r$ and adding 
$\delta_k \cdot z^r P(z)^m$ to change $a_k z^k$ with $a_k \in \RR$ 
to $(a_k + \delta_k)z^k$  with $a_k + \delta_k \in \ZZ$ (the ``patching'' process), 
\index{patching}
one obtains the desired polynomial $G(z) = G^{(n)}(z) \in \ZZ[z]$. 
One uses the polynomials $\delta_k z^r \phi(z)^m$ in patching, 
rather than simply the monomials $\delta_k z^k$, 
in order to control the sup-norms $\|z^r \phi(z)^m\|_E$.  
Each adjustment changes all the coefficients of order $k$ and lower,
\index{coefficients $A_{v,ij}$!high-order} 
but leaves the higher coefficients unchanged.  
Using a geometric series estimate to show that $|G(z)| > 1$ outside $U$, 
one concludes that $G(z)$ has all its roots in $U$.  
The algebraic integers produced by the classical Fekete-Szeg\"o theorem 
\index{Fekete-Szeg\"o theorem} 
\index{algebraic integer}
are the roots of $G(z)^{\ell} - 1$ for $\ell = 1, 2, 3, \ldots$.

\vskip .1 in
The proof of the Fekete-Szeg\"o theorem with local rationality conditions follows the same pattern, 
\index{Fekete-Szeg\"o theorem with LRC} 
but with many complications.  These arise from working on curves 
of arbitrary genus, from arranging that the zeros avoid the   
finite set $\fX = \{x_1, \ldots, x_m\}$ instead of a single point, 
from working adelically, and from imposing the local rationality conditions.  

We will now sketch the proof in the situation where $E_v \subset \cC_v(K_v)$ 
for each $v \in S$.  The proof begins reducing the theorem to a setting where one 
is given a $\cC_v(\CC_v)$-neighborhood $U_v$ of $E_v$ for each $v$, with $U_v = E_v$ if $v \notin S$,
and one must construct points $\alpha \in \cC(\tK^{\sep})$ whose conjugates belong to  
$U_v \cap \cC_v(K_v)$ for each $v \in S$, and to $U_v$ for each $v \notin S$.
The strategy is to construct rational functions $G(z) \in K(\cC)$ with poles supported on $\fX$, 
whose zeros have the property above.  

One first constructs an `initial approximating function'
\index{initial approximating functions $f_v(z)$}
$f_v(z) \in K_v(\cC)$ for each $v \in S$.  Each $f_v(z)$ has 
poles supported on $\fX$ and zeros in $U_v$, with the zeros in 
$\cC_v(K_v)$ if $v \in S$.  All the $f_v(z)$ have the same degree $N$, and they have the property
that outside $U_v$ the logarithms $\log_v(|f(z)|_v$ closely approximate a weighted sum of 
Green's functions $G(z,x_i;E_v)$.  The weights are determined by $\EE$ 
\index{Green's function}
and $\fX$, through the definition of the Cantor capacity.  
\index{capacity!Cantor capacity} 

The construction of the initial approximating functions
\index{initial approximating functions $f_v(z)$} 
is one of the hardest parts of the proof.  
When working on curves of positive genus, one cannot simply take a discrete 
approximation to the equilibrium distribution, 
\index{equilibrium distribution} 
but must arrange that the divisor whose zeros come from that
approximation and whose poles have the 
prespecified orders on the points in $\fX$, is principal. 
For places $v \in S$ there are additional constraints.   
When $K_v \cong \RR$ and $E_v \subset \cC_v(\RR)$,  
one must assure that $f_v(z)$ is real-valued and 
oscillates between large positive and negative values on $E_v$ 
(a property like that of Chebyshev polynomials,
\index{Chebyshev polynomial} 
first exploited by Robinson).  
\index{Robinson, Raphael} 
In this work, we give a general 
potential-theoretic construction of oscillating functions.  
When $K_v$ is nonarchimedean and $E_v \subset \cC_v(K_v)$,   
one must arrange that the zeros of $f_v(z)$ belong to $U_v \cap \cC_v(K_v)$ 
and are uniformly distributed with respect to a certain 
generalized equilibrium measure.     
\index{equilibrium distribution} 
Both cases are treated by constructing a nonprincipal 
divisor with the necessary properties, and then carefully moving some of 
its zeros to obtain a principal divisor.  In this construction,  
the `canonical distance function' $[x,y]_{\zeta}$,\index{canonical distance!$[z,w]_{\zeta}$}
introduced in (\cite{RR1}, \S2.1), 
plays an essential role:   given a divisor $D$ of degree $0$, 
the canonical distance tells what the $v$-adic absolute of a function 
with divisor $D$ `would be', if such a function were to exist.  

A further complication is that for archimedean $v$, 
one must arrange that the leading coefficients of the
\index{coefficients $A_{v,ij}$!leading} 
Laurent expansions of $f_v(z)$ at the points $x_i \in \fX$ have 
a property of `independent variability'.
\index{logarithmic leading coefficients!independent variability of archimedean}
\index{independent variability!of logarithmic leading coefficients}
When $K_v \cong \CC$, this 
was established in (\cite{RR1}) by using a convexity property of harmonic 
functions.  When $K_v \cong \RR$, we prove it by a continuity argument ultimately 
resting on the Brouwer Fixed Point theorem.
\index{Brouwer Fixed Point theorem}  

\vskip .1 in 
Once the initial approximating functions $f_v(z)$ have been constructed, 
\index{initial approximating functions $f_v(z)$}
we modify them to obtain `coherent approximating functions' $\phi_v(z)$ with specified leading coefficients,
\index{coherent approximating functions $\phi_v(z)$} 
using global considerations.
We then use the $\phi_v(z)$  to construct `initial patching functions' 
$G_v^{(0)}(z) \in K_v(\cC)$\index{patching functions, initial $G_v^{(0)}(z)$!construction of} 
of much higher degree which still have their zeros in $U_v$ 
(and in $\cC_v(K_v)$, for $v \in S$).  The $G_v^{(0)}(z)$ are obtained  
by raising the $\phi_v(z)$ to high powers, 
or by composing them with Chebyshev polynomials or generalized 
\index{Chebyshev polynomial}\index{Stirling polynomial} 
Stirling polynomials if $v \in S$. 
(This idea goes back to Cantor \cite{Can3}.) 
\index{Cantor, David}   

We next ``patch'' the functions $G_v^{(0)}(z)$,\index{patching functions, initial $G_v^{(0)}(z)$}   
inductively constructing $K_v$-rational functions $(G_v^{(k)}(z))_{v \in S}$,
$k = 1, 2, \ldots, n$,\index{patching functions, $G_v^{(k)}(z)$ for $1 \le k \le n$!constructed by patching}  
for which more and more of the high order Laurent coefficients 
\index{coefficients $A_{v,ij}$!Laurent} 
(relative to the points in $\fX$) are $K$-rational and independent of $v$.  
In the patching process,
\index{patching} 
we take care that the roots of 
$G_v^{(k)}(z)$\index{patching functions, $G_v^{(k)}(z)$ for $1 \le k \le n$!roots are confined to $E_v$} belong to $U_v$ for all $v$,
and belong to $\cC_v(K_v)$ for each $v \in S$.  
In then end we obtain a global $K$-rational function 
$G^{(n)}(z) = G_v^{(n)}(z)$\index{patching functions, $G_v^{(k)}(z)$ for $1 \le k \le n$!$G_v^{(n)}(z) = G^{(n)}(z)$ is independent of $v$} 
independent of $v$, which ``looks like'' $G_v^{(0)}(z)$ 
at each $v \in S$.\index{patching functions, initial $G_v^{(0)}(z)$} 

\vskip .05 in
The patching process has two aspects, global and local.
\index{patching}

\vskip .05 in
The global aspect concerns achieving $K$-rationality for $G(z)$,
while assuring that its roots remain outside the balls $B_v(x_i,1)^-$ 
for the infinitely many $v$ where $E_v$ is $\fX$-trivial.  
\index{$\fX$-trivial}
It is necessary to carry out the patching process in a galois-invariant way. 
\index{patching argument!global}
For this, we construct an $\Aut(\tK/K)$-equivariant
basis for the space of functions in $K(\CC)$ with poles supported on $\fX$, 
\index{basis!$L$-rational}
and arrange that when the functions $G_v^{(k)}(z)$ 
are expanded\index{patching functions, $G_v^{(k)}(z)$!expansion of}
relative to this basis, their coefficients are equivariant 
\index{coefficients $A_{v,ij}$} 
under $\Aut_c(\CC_v/K_v)$.  

The most delicate step involves patching the
leading coefficients:  one must arrange that they be $S$-units
\index{patching!leading coefficients}
\index{coefficients $A_{v,ij}$!leading}   
(the analogue of monicity in the classical case). 
The argument can succeed only if 
the orders of the poles of the $f_v(z)$ at the $x_i$
lie in a prescribed ratio to each other.  
The existence of such a ratio is intimately related 
to the fact that $\gamma(\EE,\fX) > 1$, 
and is at the heart of the definition of the Cantor capacity, 
\index{capacity!Cantor capacity}
as will be explained below.  

The remaining coefficients must be patched to be $S$-integers.  
As in the classical case, patching the high-order coefficients  
\index{coefficients $A_{v,ij}$!high-order}\index{patching!high-order coefficients}
presents special difficulties.  In general there are both
archimedean and nonarchimedean places in $S$.  
It is no longer possible to use continuity and the multinomial theorem 
\index{multinomial theorem}
as in the classical case;  instead, we use a phenomenon of `magnification'\index{magnification argument}
at the archimedean places, first applied in (\cite{RR3}), together with
a phenomenon of `contraction' at the nonarchimedean places.
In the function field case, additional complications arise from inseparability issues.
A different method is used to 
patch the high order coefficients than in the number field case:  
\index{coefficients $A_{v,ij}$!high-order} 
in the construction of initial patching functions, 
\index{patching functions, initial $G_v^{(0)}(z)$}
we arrange that the high order coefficients are all $0$, 
\index{coefficients $A_{v,ij}$!high-order} 
and that the patching process for the leading coefficients preserves this property. 
\index{patching!leading coefficients}  
\index{coefficients $A_{v,ij}$!leading}     

\vskip .05 in
The local aspect of the patching process consists of giving 
`confinement arguments'
\index{patching argument!local}
showing how to keep the roots of the $G_v^{(k)}(z)$ 
in the sets $E_v$,\index{patching functions, $G_v^{(k)}(z)$!roots are confined to $E_v$} 
while modifying the Laurent coefficients.  
\index{coefficients $A_{v,ij}$!Laurent} 
Four confinement arguments are required, 
\index{confinement argument}
corresponding to the cases $K_v \cong \CC$, 
$K_v \cong \RR$ with $E_v \subset \cC_v(\RR)$, 
$K_v$ nonarchimedean with $E_v$ being an  RL-domain, and $K_v$ nonarchimedean with 
$E_v \subset \cC_v(K_v)$.  
The confinement arguments in first and third case are adapted from (\cite{RR1}), 
and those in the second and fourth case are generalizations of those in (\cite{RR3}).  
The fourth case involves locally 
expanding the functions $G_v^{(k)}(z)$\index{patching functions, $G_v^{(k)}(z)$!expansion of}
as $v$-adic power series, and extending the Newton polygon\index{Newton Polygon}
construction in (\cite{RR3}) from polynomials to power series.  A crucial step
involves moving apart roots which have come close to each other.  This requires 
the theory of the Universal Function developed in Appendix \ref{AppC}, 
and the local action\index{local action of the Jacobian}
of the Jacobian developed Appendix \ref{AppD}.
\index{Jacobian variety}

\vskip .1 in
\section*{The Definition of the Cantor Capacity} \label{Sec0C}
\index{capacity!Cantor capacity|ii}

We next discuss the 
Cantor capacity $\gamma(\EE,\fX)$, 
which is treated more fully in (\cite{RR1}, \S5.1).      
Our purpose here is to explain its meaning 
and its role in the proof of the Fekete-Szeg\"o theorem.
\index{Fekete-Szeg\"o theorem with LRC} 
First, we will need some notation.

\smallskip
If $v$ is archimedean, write $\log_v(x) = \ln(x)$. If $v$ is nonarchimedean,
let $q_v$ be the order of the residue field of $K_v$,\index{$q_v$!definition of|ii} 
and write $\log_v(x)$ for the logarithm 
to the base $q_v$.\index{$\log_v(x)$!definition of|ii}\index{logarithm!definition of $\log_v(x)$|ii} 
  
Let $q_v = e$ if $K_v \cong \RR$ and $q_v = e^2$ if $K_v \cong \CC$.  

Define normalized absolute values on the $K_v$ by letting 
$|x|_v = |x|$ if $v$ is archimedean, and taking $|x|_v$ to be the 
the modulus of additive Haar measure if $v$ is nonarchimedean.  
\index{Haar measure}
For $0 \ne \kappa \in K$, the product formula reads
\begin{equation*} 
\sum_v \log_v(|\kappa|_v) \log(q_v) \ = \ 0 \ .
\end{equation*}
Each absolute value has a unique extension to $\CC_v$, 
which we continue to denote by $|x|_v$.
  
\smallskip
For each $\zeta \in \cC_v(\CC_v)$, 
the {\it canonical distance} $[z,w]_{\zeta}$ \index{canonical distance!$[z,w]_{\zeta}$}
on $\cC_v(\CC_v) \backslash \{\zeta\}$ 
(constructed in \S2.1 of \cite{RR1}) 
plays a role in the definition of $\gamma(\EE,\fX)$ similar to the role of 
the usual absolute value $|z-w|$ on $\PP^1(\CC) \backslash \{\infty\}$
for the classical logarithmic capacity $\gamma(E)$. 
\index{capacity!logarithmic} 
The canonical distance is a symmetric, real-valued, non-negative function
\index{canonical distance!$[z,w]_{\zeta}$}
of $z, w \in \cC_v(\CC_v)$,  
with $[z,w]_{\zeta} = 0$ if and only if $z = w$.  
For each $w$, it has a ``simple pole'' as $z \rightarrow \zeta$.    
It is uniquely determined up to scaling by a constant.  
The constant can be specified by choosing a 
uniformizing parameter\index{uniformizing parameter!normalizes canonical distance|ii}  
$g_{\zeta}(z) \in \CC_v(\cC)$ at $z = \zeta$,
\index{canonical distance!$[z,w]_{\zeta}$!normalization of|ii}   
and requiring that 
\begin{equation} \label{FDDA}
\lim_{z \rightarrow \zeta} \ [z,w]_{\zeta} \cdot |g_{\zeta}(z)|_v \ = \ 1 
\end{equation} 
for each $w$.  One definition of the canonical distance is that for each $w$,
\index{canonical distance!$[z,w]_{\zeta}$}
\index{canonical distance!constructed!directly}
\begin{equation*} 
 [z,w]_{\zeta} \ = \ \lim_{n \rightarrow \infty} |f_n(z)|_v^{1/\deg(f_n)}
\end{equation*}
where the limit is taken over any sequence of functions 
$f_n(z) \in \CC_v(\cC)$ having poles only at $\zeta$ 
whose zeros approach $w$, normalized so that 
\begin{equation*}
\lim_{z \rightarrow \zeta} |f_n(z) g_{\zeta}(z)^{\deg(f_n)}|_v \ = \ 1 \ .
\end{equation*}  
A key property of $[z,w]_{\zeta}$ is that it can be used to factor the 
absolute value of a rational function in terms of its divisor:  for each  
$f(z) \in \CC_v(\cC)$, there is a constant $C(f)$ such that 
\begin{equation*}
|f(z)|_v \ = \ C(f) \cdot \prod_{x \ne \zeta} [z,x]_{\zeta}^{\ord_x(f)}  
\end{equation*}  
for all $z \ne \zeta$.  For this reason, it is `right' kernel for use in arithmetic potential theory.
\index{potential theory!arithmetic}  

\vskip .1 in
The Cantor capacity is defined in terms of Green's functions $G(z,x_i;E_v)$.
\index{Green's function} 
\index{capacity!Cantor capacity} 
We first introduce the Green's function for compact sets 
$H_v \subset \cC_v(\CC_v)$, where there is a potential-theoretic construction 
like the one in the classical case. 
Suppose  $\zeta \notin H_v$.  
For each probability measure $\nu$ supported on $H_v$, 
consider the {\it energy integral}\index{energy integral|ii}  
\begin{equation*}
I_{\zeta}(\nu) \ = \ \iint_{H_v \times H_v} - \log_v([z,w]_{\zeta}) \, 
                d\nu(z) d\nu(w) \ .
\end{equation*}
Define the {\it Robin constant} 
\index{Robin constant|ii}
\begin{equation} \label{FDD1} 
V_{\zeta}(H_v) \ = \ \inf_{\nu} I_{\zeta}(\nu) \ .  
\end{equation}
It can be shown that  
either $V_{\zeta}(H_v) < \infty$ for all $\zeta \notin E_v$, 
or $V_{\zeta}(H_v) = \infty$ for all $\zeta \notin E_v$
(see Lemma \ref{ALemA1}).  In the first case we say that $H_v$
has {\em positive inner capacity}, and the second case that it has 
\index{capacity!inner|ii}
{\em inner capacity $0$}.
  
If $H_v$ has positive inner capacity,
\index{capacity!inner} 
there is a unique probability measure 
$\mu_{\zeta}$ on $H_v$ which achieves the infimum in (\ref{FDD1}).  
It is called the equilibrium distribution of $H_v$ with respect to $\zeta$.
\index{equilibrium distribution}  
\index{equilibrium distribution} 
We define the Green's function by 
\index{Green's function}
\begin{equation} \label{FDD3} 
G(z,\zeta;H_v) 
\ = \ V_{\zeta}(H_v) + \int_{H_v} \log_v([z,w]_{\zeta}) \, d\mu_{\zeta}(w) \ .
\end{equation} 
It is non-negative 
and has a logarithmic pole as $z \rightarrow \zeta$.  
If $H_v$ has inner capacity $0$, 
\index{capacity!inner}
we put $G(z,\zeta;H_v) = \infty$ for all $z, \zeta$.  

The Green's function is symmetric for  $z, \zeta \notin H_v$, 
\index{Green's function!properties of}
and is monotone decreasing in the set $H_v$: 
\index{Green's function!monotonic} 
for compact sets  $H_v \subset H_v^{\prime}$, 
and $z, \zeta \notin E_v^{\prime}$
\begin{equation} \label{GreenCompactIneq}
G(z,\zeta;H_v) \ \ge \ G(z,\zeta;H_v^{\prime}) \ .  
\end{equation}
If $H_v$ has positive inner capacity,
\index{capacity!inner}
then for each neighborhood $U \supset H_v$,  
and each $\varepsilon > 0$, 
by taking a suitable discrete approximation to $\mu_{\zeta}$, one sees that 
there are an $N > 0$ and a function $f_v(z) \in \CC_v(\cC)$ of degree $N$, 
with zeros in $U$ and a pole of order $N$ at $\zeta$, such that
\begin{equation*}
|G(z,\zeta;H_v) - \frac{1}{N} \log_v(|f_v(z)|_v)| \ < \ \varepsilon
\end{equation*}
for all $z \in \cC_v(\CC_v) \backslash (U \cup \{\zeta\})$.  

In \cite{RR1}, Green's functions $G(z,\zeta;E_v)$ 
\index{Green's function}
are defined for compact sets $E_v$ in the archimedean case, 
and by a process of taking limits, 
for `algebraically capacitable' sets in the nonarchimedean case.
Algebraically capacitable sets include all sets that are finite unions of 
\index{algebraically capacitable}
compact sets and affinoid sets;  see (\cite{RR1}, Theorem 4.3.11). 
\index{affinoid!affinoid domain} 
In particular, the sets $E_v$ in Theorem \ref{aT1} are algebraically capacitable.
\index{algebraically capacitable}

\smallskip
We next define local and global `Green's matrices'.
Let $L/K$ be a finite normal extension containing $K(\fX)$.  
For each place $v$ of $K$ and each $w$ of $L$ with $w|v$, 
after fixing an isomorphism $\CC_w \cong \CC_v$,
we can pull back $E_v$ to a set $E_w \subset \cC_w(\CC_w)$. 
The set $E_w$ is independent of the isomorphism chosen, since $E_v$ is 
stable under $\Aut_c(\CC_v/K_v)$.  If we identify $\cC_v(\CC_v)$ and
$\cC_w(\CC_w)$, then for $z, \zeta \notin E_v$
\begin{equation} \label{FDDB} 
G(z,\zeta;E_w) \log(q_w) \ = \ [L_w:K_v] \cdot G(z,\zeta;E_v) \log(q_v) \ .
\end{equation}

For each $x_i \in \fX$, 
fix a global uniformizing parameter\index{uniformizing parameter!normalizes Robin constant|ii} $g_{x_i}(x) \in L(\cC)$
and use it to define the upper Robin constants $V_{x_i}(E_w)$ for all $w$.
\index{Robin constant!upper}  
For each $w$, let the `local upper Green's matrix' be 
\index{Green's matrix!local!upper local}
\begin{equation} \label{LocalGreensMatrixDef} 
\Gamma(E_w,\fX) \ = \
      \left( \begin{array}{cccc} 
                V_{x_1}(E_w)   & G(x_1,x_2;E_w) & \cdots & G(x_1,x_m;E_w) \\
                G(x_2,x_1;E_w) &  V_{x_2}(E_w)  & \cdots & G(x_2,x_m;E_w) \\
                    \vdots     &   \vdots       & \ddots &     \vdots     \\
                G(x_m,x_1;E_w) & G(x_m,x_2;E_w) & \cdots &  V_{x_m}(E_w) 
             \end{array} \right)   \ .                                     
\end{equation}
Symmetrizing over the places of $L$, 
define the `global Green's matrix' by\index{$q_v$!weights $\log(q_v)$ in $\Gamma(\EE,\fX)$|ii}  
\index{Green's matrix!global}\label{`SymbolIndexGlobalGreen'}
\begin{equation} \label{GlobalGreensMatrixDef} 
\Gamma(\EE,\fX) 
\ = \ \frac{1}{[L:K]} \sum_{w \in \cM_L} \Gamma(E_w,\fX) \log(q_w) \ . 
\end{equation} 
If $\EE$ is compatible with $\fX$,\index{compatible with $\fX$} 
the sum defining $\Gamma(\EE,\fX)$ is finite.  
\index{Green's matrix!global}
By the product formula, $\Gamma(\EE,\fX)$ 
is independent of the choice of the $g_{x_i}(z)$.  By (\ref{FDDB}) it is
independent of the choice of $L$.   

The global Green's matrix is symmetric and non-negative off the diagonal.
\index{Green's matrix!global} 
Its entries are finite if and only if each $E_v$ has positive inner capacity. 
\index{capacity!inner}          

\smallskip
Finally, for each $K$-rational $\EE$ compatible with $\fX$,\index{compatible with $\fX$} we define 
the {\it Cantor capacity} to be 
\label{`SymbolIndexCantorCap'} 
\index{capacity!Cantor capacity|ii}
\begin{equation} \label{GlobalCapacityDef} 
\gamma(\EE,\fX) \ = \ e^{-V(\EE,\fX)} \ , 
\end{equation} 
where $V(\EE,\fX) = \val(\Gamma(\EE,\fX))$ is the value of 
\index{Green's matrix!global}
$\Gamma(\EE,\fX)$ as a matrix game.
\index{value of $\Gamma$ as a matrix game|ii}   
Here, for any $m \times m$ real-valued matrix $\Gamma$, 
\begin{equation} \label{aF3}
\val(\Gamma) \ = \ \ \max_{\vs \in \cP^m} \ \min_{\vr \in \cP^m}
                       \  \phantom{}^t\vs \, \Gamma \vr 
\end{equation} \label{`SymbolIndexValGamma'}
where $\cP^m = \{ \phantom{}^t(s_1,\ldots,s_m) \in \RR^m : 
s_ 1, \ldots, s_m \ge 0, \sum s_i = 1 \}$ 
is the set of $m$-dimensional `probability vectors'.\label{`SymbolIndexcPm'} 
Clearly $\gamma(\EE,\fX) > 0$ 
if and only if each $E_v$ has positive inner capacity.
\index{capacity!inner}

\vskip .1 in
The hidden fact behind the definition is that $\val(\Gamma)$ 
is a function of matrices which, for symmetric real matrices  $\Gamma$ 
which are non-negative off the diagonal, 
is negative if and only if $\Gamma$ is
negative definite:  this is a consequence of Frobenius' Theorem\index{Frobenius' Theorem|ii}  
(see (\cite{RR1}, p.328 and p.331) and (\cite{G}, p.53).  
Thus, $\gamma(\EE,\fX) > 1$ if and only if $\Gamma(\EE,\fX)$ 
\index{Green's matrix!global}
is negative definite.\index{Green's matrix!negative definite} 

If $\Gamma(\EE,\fX)$ is negative definite,  
\index{Green's matrix!global}
there is a unique probability vector $\hs = {}^t(\hs_1, \ldots, \hs_m)$
such that 
\begin{equation} \label{hsUnique}
\Gamma(\EE,\fX) \, \hs \ = \ 
\left( \begin{array}{c} \hV \\ \vdots \\ \hV  \end{array} \right) 
\end{equation}
has all its coordinates equal.  From the definition of
$\val(\Gamma)$, it follows that $\hV = V(\EE,\fX) < 0$.   
For simplicity, assume in what follows that $\hs$ has rational coordinates 
(in general, this fails;  overcoming the failure is a major difficulty in the proof). 

The probability vector $\hs$ determines the relative orders of the poles
of the function $G(z)$ constructed in the Fekete-Szeg\"o theorem. 
\index{Fekete-Szeg\"o theorem with LRC} 
The idea is that the initial 
local approximating functions $f_v(z)$ should have polar
divisor $\sum_{i=1}^m N\hs_i (x_i)$ for some $N$, 
and be such that for each $v$, 
outside the given neighborhood $U_v$ of $E_v$
\begin{equation*}
\frac{1}{N} \log_v(|f_v(z)|_v) \ = \ \sum_{j=1}^m G(z,x_j;E_v) \hs_j \ .
\end{equation*} 
(At archimedean places, this will only hold asymptotically as  
$z \rightarrow x_i$, for each $x_i$.) 
The fact that the coordinates of $\Gammabar(\EE,\fX) \hs$ are equal
means it is possible to scale the $f_v(z)$
so that in their Laurent expansions at $x_i$, the leading
coefficients $c_{v,i}$ satisfy  
\begin{equation*}
\sum_v \log_v(|c_{v,i}|_v) \log(q_v) \ = \ 0 
\end{equation*}
compatible with the product formula, allowing the patching process to begin.  
Reversing this chain of ideas 
lead Cantor to his definition of the capacity. 
\index{Cantor, David} 
\index{capacity!Cantor capacity}  

\medskip
For readers familiar with intersection theory, we remark that an Arakelov-like 
adelic intersection theory for curves was constructed in (\cite{RR6}). 
The arithmetic divisors in that theory include all pairs $\cD = (D,\{G(z,D;E_v)\}_{v \in \cM_K})$ 
where $D = \sum_{i=1}^m s_i (x_i)$ is a $K$-rational divisor on $\cC$ with real coefficients  
and $G(z,D;E_v) = \sum_{i=1} s_i G(z,x_i;E_v)$. If $\vs = \hs$ is the  
probability vector constructed in (\ref{hsUnique}), then relative to that intersection theory 
\begin{equation*}
V(\EE,\fX) \ = \ ^{t}\vs \Gamma(\EE,\fX) \vs \ = \ \cD \cdot \cD \ < \ 0 \ .
\end{equation*} 
As noted by Moret-Bailly, this says that the Fekete-Szeg\"o Theorem with local rationality conditions can be viewed as
a kind of arithmetic contractibility theorem. 
\index{Moret-Bailly, Laurent}

\vskip .1 in
\section*{Outline of the Manuscript} \label{Sec0E}
In this section we outline the contents and main ideas of the work.

\smallskip

This Introduction, and Chapters \ref{Variants} and \ref{ExamApp}, are expository, 
intended to give perspective on the 
Fekete-Szeg\"o theorem.\index{Fekete-Szeg\"o theorem with LRC}   
In the Introduction we have recalled history, sketched the proofs of the classical 
Fekete-Szeg\"o theorem and Theorem \ref{aT1}, 
and defined the Cantor capacity.\index{Fekete-Szeg\"o theorem}\index{capacity!Cantor capacity} 
In Chapter \ref{Variants} we state six variants of the theorem, which extend it in different directions.  
These include a version producing points in 
`quasi-neighborhoods' of $\EE$, generalizing the classical 
Fekete-Szeg\"o theorem;\index{Fekete-Szeg\"o theorem with LRC!for quasi-neighborhoods}\index{quasi-neighborhood} 
a version producing points in $\EE$ under weaker conditions than those of Theorem \ref{aT1};  
a version which imposes ramification 
conditions at finitely many primes outside $S$;  
a version for algebraically capacitable sets which expresses 
\index{algebraically capacitable}
the Fekete/Fekete-Szeg\"o dichotomy  in terms of the global Green's matrix $\Gamma(\EE,\fX)$; 
and two versions\index{Green's matrix!global}
for Berkovich curves. 
\index{Berkovich!curve}  

In Chapter \ref{ExamApp} we give numerical examples\index{numerical!examples} illustrating the theorem on 
$\PP^1$, elliptic curves, Fermat curves, and modular curves.\index{elliptic curve}\index{Fermat curve}\index{modular curve}
 We begin by proving several formulas for capacities and Green's functions of archimedean 
and nonarchimedean sets, aiming\index{Green's function!examples}
 to collect formulas useful for applications and going beyond those tabulated in 
(\cite{RR1}, Chapter 5). In the archimedean case, we give formulas for capacities and Green's functions
\index{Green's function} 
of  one, two, and arbitrarily many intervals in $\RR$.  The formulas for two intervals involve 
classical theta-functions,\index{theta-functions} and those for multiple intervals (due to Harold Widom)\index{Widom, Harold} 
involve hyperelliptic integrals.
In the nonarchimedean case we give a general algorithm for computing capacities of compact sets.
We determine the capacities and Green's functions of rings of integers, groups of units, and bounded tori 
\index{Green's function!examples}
in local fields.  We also give the first known computation of a capacity of a nonarchimedean set 
where the Robin constant is not a rational number.
\index{capacity}

In the global case, we give numerical criteria\index{numerical!criteria} for the existence/non-existence of infinitely many
algebraic integers and units satisfying various geometric conditions.  The existence of such criteria, 
for which the prototypes are Robinson's theorems\index{Robinson, Raphael} for totally real algebraic integers and units, 
is one of the attractive features of the subject.  In applying a general theorem like the Fekete-Szeg\"o 
theorem with local rationality conditions, it is often necessary to make clever reductions in order
\index{Fekete-Szeg\"o theorem with LRC} 
to obtain interesting results, and we have tried to give examples illustrating some of the reduction 
methods that can be used.  

Our results for elliptic curves include a complete determination of the capacities 
(relative to the origin) of the integral points on Weierstrass models and N\'eron models.
\index{N\'eron model!of elliptic curve}\index{Weierstrass equation}     
Our results for Fermat curves are based on McCallum's determination of the special fibre 
\index{McCallum, William}\index{Fermat curve}
for a regular model of the Fermat curve $\cF_p$ over $\QQ_p(\zeta_p)$.  They show how the geometry
of the model (in particular the number of `tame curves' in the special fibre) is reflected
\index{tame curve}
in the arithmetic of the curve.  Our results for the modular curves $X_0(p)$ use 
the Deligne-Rapoport model.\index{modular curve}\index{modular curve!Deligne-Rapoport model} 
In combination, they illustrate a general principle that it is
usually possible to compute nonarchimedean local capacities on a curve of higher genus, 
if a regular model of the curve is known.    

\smallskip 
Beginning with Chapter \ref{Chap3}, we develop the theory rigorously. 

Chapter \ref{Chap3} covers notation, conventions, and foundational material about capacities 
and Green's functions used throughout the work.  
\index{Green's function}
An important notion is the $(\fX,\vs)$-canonical distance $[z,w]_{\fX,\vs}$.
\index{canonical distance!$[z,w]_{\fX,\vs}$}  
Given a curve $\cC/K$ and a place $v$ of $K$, we will be 
interested in constructing rational functions $f \in \cC_v(\CC_v)$ whose poles are supported on a finite set 
$\fX = \{x_1, \ldots, x_m\}$ and whose polar divisor is proportional to $\sum_{i=1}^m s_i (x_i)$, 
where $\vs = (s_1, \ldots, s_m)$ is a fixed probability vector. 
\label{`SymbolIndexvs'} 
The $(\fX,\vs)$-canonical distance enables to treat $|f(z)|_v$
\index{canonical distance!factorization property}  
like the absolute value of a polynomial, factoring it in terms of the zero divisor of $f$ as
\begin{equation*}
|f(z)|_v \ = \ C(f) \cdot \prod_{\text{zeros $\alpha_i$ of $f$}} [z,\alpha_i]_{\fX,\vs} \ .  
\end{equation*} 
Furthermore, the product on the right -- which we call an $(\fX,\vs)$-pseudopolynomial -- 
is defined and continuous even for divisors which are not principal.
This lets us separate analytic and algebraic issues in the construction of $f$.  

Put $L = K(\fX) = K(x_1, \ldots, x_m)$, and let $L^{\sep}$ be the separable closure of $K$ in $L$.
\label{`SymbolIndexLsep'}
Another important technical tool from Chapter \ref{Chap3} are the $L$-rational and $L^{\sep}$-rational
bases.  These are multiplicatively finitely generated sets of functions 
which can be used to expand rational functions with poles supported on $\fX$, 
much like the monomials $1, z,z^2, \ldots$ can be used to expand polynomials.
As their names indicate, the functions in the $L$-rational basis are defined over $L$, 
and those in the $L^{\sep}$-rational basis are defined over $L^{\sep}$.  The construction arranges
that the transition matrix between the two bases is block diagonal, 
hence has bounded norm at each place $w$ of $L$.
\index{$L$-rational basis!transition matrix block diagonal}   

In Chapter \ref{Chap4} we state a version of
the Fekete-Szeg\"o theorem with local rationality conditions for ``$K_v$-simple sets''
\index{$K_v$-simple!set}  
(Theorem \ref{aT1-B}),
\index{Fekete-Szeg\"o theorem with LRC} 
\index{$K_v$-simple!set} 
and we reduce Theorem \ref{aT1}, Corollary \ref{FSZiv}, 
and the variants stated in Chapter \ref{Variants} to it.  
The rest of the manuscript (Chapters \ref{Chap5} -- \ref{Chap11} and Appendices \ref{AppA} -- \ref{AppD}) 
is devoted to the proof of Theorem \ref{aT1-B}.

Chapters \ref{Chap5} and \ref{Chap6} 
construct the ``initial approximating functions''  needed for Theorem \ref{aT1-B}. 
\index{initial approximating functions $f_v(z)$} 
Four constructions are needed:  
for archimedean sets $E_v \subset \cC_v(\CC)$ 
when the ground field is $\CC$ and $\RR$, 
and for nonarchimedean sets $E_v \subset \cC_v(\CC_v)$ which are $\RL$-domains   
\index{$\RL$-domain} or are compact.  
The first and third were done in (\cite{RR1}); the second and fourth 
are done here.

The probability vector $\vs$ ultimately used in the construction is determined by $\EE$ and $\fX$, 
through the global Green's matrix $\Gamma(\EE,\fX)$.  This means that for each $E_v$,
the local constructions 
must be carried out in a uniform way for all $\vs$.  In Appendix \ref{AppA} we develop 
potential theory with respect to the kernel $[z,w]_{\fX,\vs}$.  It turns out that there are 
\index{canonical distance!$[z,w]_{\fX,\vs}$!potential theory for}
$(\fX,\vs)$-capacities, $(\fX,\vs)$-Green's functions, and $(\fX,\vs)$-equilibrium distributions
\index{capacity!$(\fX,\vs)$}\index{potential theory!$(\fX,\vs)$} 
\index{equilibrium distribution!$(\fX,\vs)$} 
\index{Green's function!$(\fX,\vs)$} 
with properties analogous to the corresponding objects in classical potential theory.\index{potential theory!classical}  
The initial approximating functions are $(\fX,\vs)$-functions 
\index{initial approximating functions $f_v(z)$}\index{$(\fX,\vs)$-function}
whose normalized logarithms $\deg(f)^{-1} \log_v(|f(z)|_v)$ closely 
approximate the $(\fX,\vs)$-Green's function outside a neighborhood of $E_v$, 
and whose zeros are roughly equidistributed   
like the $(\fX,\vs)$-equilibrium distribution.

Chapter \ref{Chap5} deals with the construction of initial approximating functions $f(z) \in \RR(\cC_v)$
\index{initial approximating functions $f_v(z)$}
when the ground field $K_v$ is $\RR$, 
for galois-stable sets $E_v \subset \cC_v(\CC)$ which are finite unions
of intervals in $\cC_v(\RR)$ and closed sets in $\cC_v(\CC)$ with piecewise smooth boundaries.  
The desired functions  must oscillate with large magnitude on the real intervals. 
The construction has two parts:  a potential-theoretic part carried out in Appendix \ref{AppB}, 
which constructs `$(\fX,\vs)$ pseudo-polynomials' whose absolute value behaves like that of 
a Chebyshev polynomial, and an algebraic part which involves adjusting the divisor of the
\index{Chebyshev polynomial}  
pseudo-polynomial to make it principal.  The first part of the argument requires subdividing the 
real intervals into `short' segments, where the notion of shortness depends only on the deviation 
of the canonical distance $[z,w]_{\fX,\vs}$ from $|z-w|$ in local coordinates, 
\index{canonical distance!$[z,w]_{\fX,\vs}$}  
\index{canonical distance!archimedean!deviation from absolute value}  
\index{canonical distance!determines `shortness'}
and is uniform over compact sets.  
The second part of the argument uses a variant of the Brouwer Fixed Point theorem. 
\index{Brouwer Fixed Point theorem} 
An added difficulty involves assuring that the `logarithmic leading coefficients' of $f$ are 
\index{logarithmic leading coefficients} 
independently variable over a range independent of $\vs$,  which is needed as an input 
to the global patching process in Chapter \ref{Chap7}.  
\index{patching argument!global}

Chapter \ref{Chap6} deals with the construction of initial approximating functions $f \in K_v(\cC_v)$
\index{initial approximating functions $f_v(z)$}
when the ground field $K_v$ is a nonarchimedean local field, and the sets $E_v$ 
are galois-stable finite unions of balls in $\cC_v(F_{w,i})$, 
for fields $F_{w,i}$ are which are finite separable extensions of $K_v$. 
Again the construction has two parts:  an analytic part, 
which constructs an $(\fX,\vs)$ pseudo-polynomial by transporting Stirling polynomials
\index{Stirling polynomial}\index{pseudopolynomial!$(\fX,\vs)$}
for the rings of integers of the $F_{w,i}$ to the balls, and an algebraic part, 
which involves moving some of the roots of the pseudo-polynomial to make its divisor principal.
When $\cC_v$ has positive genus $g$, 
this uses an action of a neighborhood of the origin in $\Jac(\cC)(\CC_v)$ on $\cC_v(\CC_v)^g$ 
constructed in Appendix \ref{AppD}.

Chapter \ref{Chap7} contains the global patching argument for Theorem \ref{aT1-B}, 
\index{patching argument!global} 
which breaks into two cases:  when $\Char(K) = 0$,
and when $\Char(K) = p > 0$.  The two cases involve different difficulties.  When $\Char(K) = 0$,
the need to patch archimedean and nonarchimedean initial approximating functions together 
\index{initial approximating functions $f_v(z)$}
is the main constraint, and the most serious bottleneck involves patching the leading coefficients.
\index{coefficients $A_{v,ij}$!leading}\index{patching!leading coefficients}   
The ability to independently adjust the logarithmic leading coefficients for the archimedean
\index{logarithmic leading coefficients} 
initial approximating functions allows us to accomplish this.  When $\Char(K) = p > 0$, the leading
coefficients are not a problem, but separability/inseparability issues drive the argument. 
\index{coefficients $A_{v,ij}$!leading} 
These are dealt with by simultaneously monitoring the patching process with respect 
to the $L$-rational and $L^{\sep}$-rational bases from Chapter \ref{Chap3}.    

Chapters \ref{Chap8} -- \ref{Chap11} contain the local patching arguments needed for Theorem \ref{aT1-B}.
\index{patching argument!local} 
Chapter \ref{Chap8} concerns the case when $K_v \cong \CC$, 
Chapter \ref{Chap9} concerns the case when $K_v \cong \RR$, Chapter \ref{Chap10} concerns the
nonarchimedean case for $\RL$-domains, and Chapter \ref{Chap11} 
concerns the nonarchimedean case for compact sets.
\index{$\RL$-domain} 
Each  provides geometrically increasing bounds for the amount the coefficients  
can be varied, while simultaneously confining the movement of the roots,
as the patching proceeds from high order to low order coefficients.   

Chapter \ref{Chap8} gives the local patching argument when $K_v \cong \CC$.  
The aim of the construction is to confine the roots of the function to a prespecified neighborhood 
\index{patching argument!local}\index{confinement argument} 
$U_v$ of $E_v$, while providing the global patching construction with increasing freedom 
\index{patching argument!global} 
in modify the coefficents relative to the $L$-rational basis, as the degree
of the basis functions goes down.  For the purposes of the patching argument, the coefficients
\index{coefficients $A_{v,ij}$!high-order}\index{patching argument!local}  
\index{coefficients $A_{v,ij}$!middle} 
\index{coefficients $A_{v,ij}$!low-order} 
are grouped into `high-order', `middle' and `low-order'.  The construction begins by raising
the initial approximating function to a high power $n$.  A `magnification argument', 
\index{initial approximating functions $f_v(z)$}\index{magnification argument}
similar to the ones in (\cite{RR2}) and (\cite{RR3}), is used to gain the freedom needed to patch 
the high-order coefficients.

Chapter \ref{Chap9} gives the local patching argument when $K_v \cong \RR$. 
Here the construction must simultaneously confine the roots to a set $U_v$ which is the union 
\index{confinement argument}\index{patching argument!local}  
of $\RR$-neighborhoods of the components of $E_v$ in $\cC_v(\RR)$, 
and $\CC$-neighborhoods of the other components.  We call such a set a `quasi-neighborhood' 
\index{quasi-neighborhood|ii}
of $E_v$.   The construction is similar to the one over $\CC$, except that it begins by 
composing the initial approximating function with a Chebyshev polynomial of degree $n$.
\index{initial approximating functions $f_v(z)$}
\index{Chebyshev polynomial}   
Chebyshev polynomials have the property that they oscillate with large magnitude on a real interval,
\index{Chebyshev polynomial} 
and take a family of confocal ellipses in the complex plane to ellipses. 
Both properties are used in the confinement argument.
\index{confinement argument} 

Chapter \ref{Chap10} gives the local patching construction when $K_v$ is nonarchimedean
and $E_v$ is an $\RL$-domain.  The construction again begins by raising the initial approximating
\index{initial approximating functions $f_v(z)$}
function to a power $n$, and to facilitate patching the high-order coefficients, 
\index{coefficients $A_{v,ij}$!high-order}\index{patching argument!local}  
we require that $n$ be divisible by a high power of the residue characteristic $p$.  
If $K_v$ has characteristic $0$, this makes the high order coefficients be $p$-adically small;
\index{coefficients $A_{v,ij}$!high-order} 
if $K_v$ has characteristic $p$, it makes them vanish (apart from the leading coefficients), 
\index{coefficients $A_{v,ij}$!leading} 
so they do not need to be patched at all.

Chapter \ref{Chap11} gives the local patching construction when $K_v$ is nonarchimedean and $E_v$
is compact.  This case is by far the most intricate, and 
begins by composing the initial approximating\index{initial approximating functions $f_v(z)$}
function with a Stirling polynomial.  If $K_v$ has characteristic $0$, 
\index{Stirling polynomial}
this makes the high order coefficients be $p$-adically small; if $K_v$ has characteristic $p$,
\index{coefficients $A_{v,ij}$!high-order}\index{patching argument!local}   
it makes them vanish.  The confinement argument\index{confinement argument} 
generalizes those in (\cite{RR2}, \cite{RR3}),  
and the roots are controlled by tracking their positions within ``$\psi_v$-regular sequences''. 

A $\psi_v$-regular sequence is a finite sequence of roots which are $v$-adically spaced like 
\index{regular sequence!$\psi_v$-regular sequence}\index{$\psi_v$-regular sequence|see{regular sequence}}
an initial segment of the integers, viewed as embedded in $\ZZ_p$ (see Definition \ref{DDefCPC1}) 
The local rationality of each root is preserved by an argument 
involving Newton polygons\index{Newton Polygon} for power series.    
In the initial stages, confinement of the roots depends on the fact
\index{confinement argument}\index{patching argument!local}  
that the Stirling polynomial factors completely over $K_v$.  Some roots may move quite close close to others 
\index{Stirling polynomial}  
in early steps of the patching process, and the the middle part of argument involves an extra 
step of separating roots, first used in (\cite{RR2}). This is accomplished by multiplying 
the partially patched function with a carefully chosen rational function whose zeros and poles
are very close in pairs.  This function is obtained by specializing the 
`Universal Function' constructed in Appendix \ref{AppC}, which parametrizes all functions of given degree 
by means of their roots and poles and value at a normalizing point.    

Appendix A develops potential theory\index{potential theory!$(\fX,\vs)$} with respect the kernel $[z,w]_{\fX,\vs}$, 
paralleling the classical development of potential theory\index{potential theory!classical} over $\CC$ given in (\cite{Ts}).
There are $(\fX,\vs)$-equilibrium distributions, potential functions, transfinite diameters, Chebyshev constants,  
\index{equilibrium distribution!$(\fX,\vs)$} 
\index{transfinite diameter!$(\fX,\vs)$} 
\index{Chebyshev constant!$(\fX,\vs)$}
\index{equilibrium distribution!$(\fX,\vs)$} 
and capacities with the same properties as in the classical theory. 
\index{capacity!$(\fX,\vs)$}
A key result is Proposition \ref{BPropF1}, which asserts that `$(\fX,\vs)$-Green's functions', obtained 
\index{Green's function!$(\fX,\vs)$}
by subtracting an `$(\fX,\vs)$-potential function'from an `$(\fX,\vs)$-Robin constant',
\index{Robin constant!$(\fX,\vs)$}\index{potential function!$(\fX,\vs)$} 
are given by linear combinations of the Green's functions constructed in (\cite{RR1}). 
Other important results are Lemmas \ref{bL3} and \ref{ALJ29}, 
which provide uniform upper and lower bounds for the mass 
the $(\fX,\vs)$-equilibrium distribution can place on a subset, independent of $\vs$;  
\index{mass bounds} 
and Theorem \ref{BFThm2A}, which shows that nonarchimedean 
$(\fX,\vs)$-Green's functions and equilibrium distributions
can be computed using linear algebra.
  
Appendix B constructs archimedean local oscillating functions for short intervals, 
\index{short@`short' interval}  
and gives the potential-theoretic input 
for the construction of the initial approximating functions over $\RR$ in Chapter 5.
\index{initial approximating functions $f_v(z)$}  
In classical potential theory, the equality of the transfinite diameter, Chebyshev constant, 
\index{transfinite diameter}\index{potential theory!classical} 
\index{Chebyshev constant}  
and logarithmic capacity of a compact set $E \subset \CC$ is shown by means of a `rock-paper-scissors'
\index{capacity!logarithmic}\index{rock-paper-scissors argument} 
argument proving in a cyclic fashion that each of the three quantities is greater than or equal to the next.
Here, a rock-paper-scissors argument is used to prove Theorem \ref{ATF19}, 
which says that the probability measures associated to the roots of weighted Chebyshev polynomials 
\index{Chebyshev polynomial} 
for a set $E_v$ converge to the $(\fX,\vs)$-equilibrium measure of $E_v$.  
  
Appendix C studies the `universal function' of degree $d$ on a curve, 
\index{universal function} 
used in Chapter \ref{Chap11}.  
We give two constructions for it, 
one by the author using the theory of the Picard scheme,
\index{Picard scheme}
the other by Robert Varley using Grauert's theorem. 
\index{Varley, Robert}\index{Grauert's theorem}
We then use local power series parametrizations, together with a compactness argument,
to obtain uniform bounds for the
change in the norm of a function outside a union of balls containing its divisor,
if its zeros and poles are moved a distance at most $\delta$ (Theorem \ref{BMT2}).
We thank Varley for permission to include his construction here.  
\index{Varley, Robert}

Appendix D shows that in the nonarchimedean case, if the genus $g$ of $\cC$ is positive, 
then at generic points of $\cC_v(\CC_v)^g$ there is an action of a neighborhood 
of the origin of the Jacobian on $\cC_v(\CC_v)^g$,
\index{local action of the Jacobian}\index{Jacobian variety}
which makes $\cC_v(\CC_v)^g$ into a local principal homogeneous space.\index{principal homogeneous space}  
This is used in Chapters \ref{Chap6} and \ref{Chap11} in adjusting non-principal divisors to make them principal.  
The action is obtained by considering the canonical map $\cC_v^g(\CC_v) \rightarrow \Jac(\cC)(\CC_v)$,
which is locally an isomorphism outside a set of codimension $1$, pulling back the formal group of 
\index{formal group}
the Jacobian, and using properties of power series in several variables.  
Theorem \ref{BKeyThm1} gives the most general form of the action.

\section*{Thanks} \label{ThanksSection}
The author thanks  Pete Clark, Will Kazez, Dino Lorenzini, Ted Shifrin, and Robert Varley for useful 
conversations\index{Varley, Robert}\index{Shifrin, Ted}
during\index{Clark, Pete L.}\index{Lorenzini, Dino}\index{Kazez, William}  
the course of the investigation.  He also thanks the Institute Henri Poincar\'e in Paris, 
where the work was begun  
during the Special Trimestre on Diophantine Geometry in 1999,\index{Institute Henri Poincar\'e} 
and the University of Georgia,\index{University of Georgia} where most of the research was carried out.

The author gratefully acknowledges the National Science Foundation's support of this project 
through grants DMS 95-000892, DMS 00-70736, DMS 03-00784, and DMS 06-01037.\index{National Science Foundation|ii} 
Any opinions, findings and conclusions or recommendations expressed in this material are those of the author 
and do not necessarily reflect the views of the National Science Foundation.\index{National Science Foundation!disclaimer} 

\newpage
\section*{Symbol Table} \label{Sec0F}
 
\ Below are some symbols used throughout the work. 
See \S\ref{Chap3}.\ref{NotationSection}, \S\ref{Chap3}.\ref{AssumptionsSection}
for more conventions.

\vskip .2in
\hskip -.3in\begin{tabular}{lllll}
Symbol & & Meaning & & Defined\\ 
\hline
$K$ & & a global field & & p. \pageref{`SymbolIndexK'} \\
$\cC$ & & a smooth, projective, connected curve over $K$ 
                     & & p. \pageref{`SymbolIndexcC'} \\ 
$g = g(\cC)$ & & the genus of $\cC$ & & p. \pageref{`SymbolIndexgenus'} \\
$\tK$ & & a fixed algebraic closure of $K$ & & p. \pageref{`SymbolIndextK'} \\
$\tK^{\sep}$ & & the separable closure of $K$ in $\tK$ & & p. \pageref{`SymbolIndextKsep'} \\                   
$K_v$ & & the completion of $K$ at a place $v$ & & p. \pageref{`SymbolIndexKv'} \\
$\cO_v$ & & the ring of integers of $K_v$ & & p. \pageref{`SymbolIndexOv'} \\
$\tK_v$ & & a fixed algebraic closure of $K_v$ & & p. \pageref{`SymbolIndextKv'} \\
$\CC_v$ & & the completion of $\tK_v$ & & p. \pageref{`SymbolIndexCCv'} \\
$\Aut(\tK/K)$ & & the group of continuous automorphisms of $\tK/K$ & & p. \pageref{`SymbolIndexAuttK'} \\
$\Aut_c(\CC_v/K_v)$ & & the group of continuous automorphisms of $\CC_v/K_v$ 
                     & & p. \pageref{`SymbolIndexAutc'} \\
$\fX = \{x_1, \ldots, x_m\}$ & & a finite, $\Aut_c(\tK/K)$-stable set of points of $\cC(\tK)$ 
                     & & p. \pageref{`SymbolIndexfX'} \\
$\vs = (s_1, \ldots, s_m)$ & & a probability vector weighting the points in $\fX$ 
                     & & p. \pageref{`SymbolIndexvs'} \\                    
$L = K(\fX)$ & & the field $K(x_1, \ldots, x_m)$ & &  p. \pageref{`SymbolIndexL'} \\
$L^{\sep}$ & & the separable closure of $K$ in $L$ & & p. \pageref{`SymbolIndexLsep'} \\
$\cC_v$           & & the curve $\cC \times_K \Spec(K_v)$ 
                      & & p. \pageref{`SymbolIndexcCv'}  \\  
$\cCbar_v$         & & the curve $\cC_v \times_{K_v} \Spec(\CC_v)$ 
                      & & p. \pageref{`SymbolIndexcCbarv'}  \\    
$\|z,w\|_v$ & & the chordal distance or spherical metric on $\cC_v(\CC_v)$ & & p. \pageref{`SymbolIndexChordalDist'}ff \\
$\|f\|_{E_v}$ & & the sup norm $\sup_{z \in E_v} |f(z)|_v$ & & p. \pageref{`SymbolIndexSupNorm'} \\
$D(a,r)$ & & the `closed disc' $\{z \in \CC_v : |z-a|_v \le r \}$ & & p. \pageref{`SymbolIndexClosedDisc'} \\
$D(a,r)^-$ & & the `open disc' $\{z \in \CC_v : |z-a|_v < r \}$ & & p. \pageref{`SymbolIndexOpenDisc'} \\
$B(a,r)$ & & the `closed ball' $\{z \in \cC_v(\CC_v) : \|z,a\|_v \le r \}$ & & p. \pageref{`SymbolIndexClosedBall'} \\
$B(a,r)^-$ & & the `open ball' $\{z \in \cC_v(\CC_v) : \|z,a\|_v < r \}$ & & p. \pageref{`SymbolIndexOpenBall'} \\
$q_v$ & & the order of the residue field of $K_v$, if $v$ is nonarchimedean & & p. \pageref{`SymbolIndexqv'}ff \\
$w_v$       & & the distinguished place of $L$ over a place $v$ of $K$  & & p. \pageref{`SymbolIndexw_v'} \\
$\val_v(x)$ & & the exponent of the largest power of $q_v$ dividing $x \in \NN$ & & p. \pageref{`SymbolIndexValw'} \\
$\log_v(x)$ & & the logarithm to the base $q_v$, when $v$ is nonarchimedean  & & p. \pageref{`SymbolIndexLogv'} \\
$\ord_v(z)$ & & the exponential valuation $-\log_v(|z|_v)$, for $z \in \CC_v$  & & p. \pageref{`SymbolIndexOrdv'} \\
$\log(x)$ & & the natural logarithm $\ln(x)$ & & p. \pageref{`SymbolIndexLog'} \\                
$\zeta$             & & a point of $\cC_v(\CC_v)$ & & p. \pageref{`SymbolIndexZeta'} \\
$g_\zeta(z)$        & & a fixed uniformizing parameter at $\zeta$ & & p. \pageref{`SymbolIndexgZeta'} \\
$[z,w]_{\zeta}$ & & the canonical distance with respect to  $\zeta \in \cC_v(\CC_v)$ 
                     & & p. \pageref{`SymbolIndexCanD'} \\
$[z,w]_{(\fX,\vs)}$ & & the $(\fX,\vs)$-canonical distance on $\cC_v(\CC_v)$ 
                     & & p. \pageref{`SymbolIndexCanDfXvs'} \\
$E_v$               & & a subset of $\cC_v(\CC_v)$ & & p. \pageref{`SymbolIndexEv'} \\
$\cl(E_v)$ & & the topological closure of $E_v$ & & p. \pageref{`SymbolIndexcl'} \\
$\gamma_{\zeta}(E_v)$ & & the capacity of $E_v$ with respect to $\zeta$ and $g_\zeta(z)$ 
                    & & p. \pageref{`SymbolIndexCapZeta'} \\ 
$V_{\zeta}(E_v)$ & & the Robin constant of $E_v$ with respect to $\zeta$ and $g_\zeta(z)$ 
                    & & p. \pageref{`SymbolIndexVZeta'} \\     
$G(z,\zeta;E_v)$ & & the Green's function of $E_v$  & & p. \pageref{`SymbolIndexGreen'} \\ 
$\val(\Gamma)$ & & the value of $\Gamma \in M_n(\RR)$ as a matrix game & & p. \pageref{`SymbolIndexValGamma'} \\ 
$\EE = \prod_v E_v$ & & an adelic set in $\prod_v \cC_v(\CC_v)$ & & p. \pageref{`SymbolIndexEEv'} \\
$\Gamma(\EE,\fX)$ & & the global Green's matrix of $\EE$ relative to $\fX$ 
                     & & p. \pageref{`SymbolIndexGlobalGreen'}  \\                      
$\gamma(\EE,\fX)$ & & the global Cantor capacity of $\EE$ with respect to $\fX$  
                     & & p. \pageref{`SymbolIndexCantorCap'}   
\end{tabular}

\hskip -.3in\begin{tabular}{lllll}
Symbol & & Meaning & & Defined\\ 
\hline
$\cC_v^{\an}$         & & the Berkovich analytification of $\cCbar_v$ 
                      & & p. \pageref{`SymbolIndexcCan'} \\   
$\BerkE_v$         & & a subset of $\cC_v^{\an}$ & & p. \pageref{`SymbolIndexBerkEv'} \\  
$V_{\zeta}(\BerkE_v)^{\an}$ & & the Robin constant of $\BerkE_v$ with respect to $\zeta$ and $g_\zeta(z)$ 
                    & & p. \pageref{`SymbolIndexBerkVZeta'} \\     
$G(z,\zeta;\BerkE_v)^{\an} $ & & the Thuillier Green's function of $\BerkE_v$ 
                    & & p. \pageref{`SymbolIndexBerkGreen'} \\ 
$\gamma(\EE,\fX)^{\an}$ & & the global capacity of a Berkovich set $\EE = \prod_v \BerkE_v$ relative to $\fX$  
                     & & p. \pageref{`SymbolIndexBerkCantorCap'}  \\
$\cP^m = \cP^m(\RR)$ & & the set of probability vectors $\vs = (s_1, \ldots, s_m) \in \RR^m$ 
                   & & p. \pageref{`SymbolIndexcPm'} \\ 
$\cP^m(\QQ)$ & & the set of  probability vectors with rational coefficients  & & p. \pageref{`SymbolIndexcPmQ'} \\
$J$ & & $J = 2g+1$ if $\Char(K) = 0$, a power of $p$ if $\Char(K) = p > 0$ & & p. \pageref{`SymbolIndexJ'}ff \\
$\varphi_{ij}(z), \varphi_\lambda(z)$ & & functions in the $L$-rational basis & & p. \pageref{`SymbolIndexLRat'} \\
$\tphi_{ij}(z), \tphi_\lambda(z)$ & & functions in the $L^{\sep}$-rational basis 
                 & & p. \pageref{`SymbolIndexLsepRat'}   \\
$\Lambda_0$ & & number of low-order elements in the $L$- and $L^{\sep}$-rational bases 
                                   & & p. \pageref{`SymbolIndexLambda0'} \\
$\Lambda$ & & number of basis elements deemed low-order in patching & & p. \pageref{`SymbolIndexLambda'} \\                 
$f_v(z)$ & & an initial approximating function & & p. \pageref{`SymbolIndexInitialApproxf'} \\
$c_{v,i}$ & & the leading coefficient of $f_v(z)$ at $x_i$ & & p. \pageref{`SymbolIndexcvi'} \\
$\Lambda_{x_i}(f_v,\vs)$ & & the logarithmic leading coefficient of $f_v(z)$ at $x_i$ & & p. \pageref{`SymbolIndexLogLeadf'} \\
$\Lambda_{x_i}(E_v,\vs)$ & & the logarithmic leading coefficient of the Green's function of $E_v$ 
                                   & & p. \pageref{`SymbolIndexLogLeadGreen'} \\
$\phi_v(z)$ & & a coherent approximating function & & p. \pageref{`SymbolIndexCoherentApproxf'} \\
$\tc_{v,i}$ & & the leading coefficient of $\phi_v(z)$ at $x_i$ & & p. \pageref{`SymbolIndextcvi'} \\
$\cI$ & & the index set $\{(i,j) \in \ZZ^2 : 1 \le i \le m, j \ge 0\}$ & & p. \pageref{`SymbolIndexI'} \\
$\prec_N$ & & the order on $\cI$ determining how coefficients are patched & & p. \pageref{`SymbolIndexPrecN'} \\
$\Band_N(k)$ & & `Bands' of indices in $\cI$ for the order $\prec_N$ & & p. \pageref{`SymbolIndexBandN'} \\
$\Block_N(i,j)$ & & the `Galois orbit' of the index $(i,j) \in \cI$ & & p. \pageref{`SymbolIndexBlockN'} \\
$G_v^{(k)}(z)$ & & the patching function at $v$ in stage $k$ of the patching process & & p. \pageref{`SymbolIndexGvk'} \\
$A_{v,ij}, A_{v,\lambda}$ & & the coefficients of $G_v^{(k)}(z)$ relative to the $L$-rational basis 
                          & & p. \pageref{`SymbolIndexAvij'} \\
$\tA_{v,ij}, \tA_{v,\lambda}$ & & the coefficients of $G_v^{(k)}(z)$ relative to the $L^{\sep}$-rational basis 
                          & & p. \pageref{`SymbolIndextAvij'} \\
$\Delta^{(k)}_{v,ij}, \Delta^{(n)}_{v,\lambda}$ 
             & & the changes in the coefficients of $G_v^{(k)}(z)$ in stage $k$ of patching
                          & & p. \pageref{`SymbolIndexDeltavij'} \\
$\vartheta_{v,ij}^{(k)}(z), \vartheta_{v,\lambda}^{(n)}$ & & compensating functions for stage $k$ of patching  
                          & & p. \pageref{`SymbolIndexComp'} \\
$\psi_v(k)$ & & the basic well-distributed sequence for the ring $\cO_v$ 
                         & & p. \pageref{`SymbolIndexWellDistv'}  \\
$S_{n,v}(z)$ & & the Stirling polynomial $\prod_{k=0}^{n-1} (z-\psi_v(k))$ for $\cO_v$ 
                        & & p. \pageref{`SymbolIndexStirlingv'}  \\                          
$E(a,b)$ & & the filled ellipse $\{z = x+iy \in \CC: x^2/a^2 + y^2/b^2 \le 1 \}$ & & p. \pageref{`SymbolIndexFilledEllipse'} \\
$T_n(z)$ & & the Chebyshev polynomial of degree $n$ for $[-2,2]$ & & p. \pageref{`SymbolIndexChebPoly'}  \\
$T_{n,R}(z)$    & & the Chebyshev polynomial of degree $n$ for $[-2R,2R]$ & & p. \pageref{`SymbolIndexChebPolyR'}  \\
$\FF_p[[t]]$ & & the ring of formal power series over $\FF_p$ & & p. \pageref{`SymbolIndexFpsbt'} \\
$\FF_p((t))$ & & the field of formal Laurent series over $\FF_p$ & & p. \pageref{`SymbolIndexFprbt'} \\
$\Jac(\cC_v)$          & & the Jacobian of a curve $\cC_v$ with genus $g > 0$ & & p. \pageref{`SymbolIndexJacobian'}  \\ 
$\JacNer(\cC_v)$     & & the N\'eron model of $\cC_v$ & & p. \pageref{`SymbolIndexNeron'} 
\end{tabular}

%% file: NewFSZChap1.tex
\chapter{Variants} \label{Variants} 

In this chapter we give six variants of  Theorem \ref{aT1}, strengthening it in different directions.
\index{Fekete-Szeg\"o theorem with LRC}     
Theorem \ref{aT1}, Corollary \ref{FSZiv} and the variants stated here will 
be reduced to yet another variant (Theorem \ref{aT1-B}) in Chapter \ref{Chap4}, 
and we will spend most of the paper proving the theorem in that form.  





\medskip
Our first variant is similar to the original theorem of Fekete and Szeg\"o (\cite{F-SZ}).
\index{Fekete-Szeg\"o theorem} 
In that theorem the sets $E_v \subset \CC$ were compact, 
and the conjugates of the algebraic integers produced were required to lie 
in arbitrarily small open neighborhoods $U_v$ of the $E_v$. 
In Theorem \ref{aT1-A1} below, 
we lift the assumption of compactness 
and replace the Cantor capacity with {\it inner Cantor Capacity} $\gammabar(\EE,\fX)$,
\index{capacity!Cantor capacity}
\index{capacity!inner Cantor capacity} 
which is defined for arbitrary adelic sets.  We also replace the neighborhoods 
$U_v$ with ``quasi-neighborhoods'', which are finite unions of open sets in $\cC_v(\CC_v)$ and open sets in 
\index{quasi-neighborhood}
$\CC_v(F_w)$, for algebraic extensions $F_w/K_v$ in $\CC_v$.  

\smallskip
The inner Cantor capacity $\gammabar(\EE,\fX)$ is similar to 
\index{capacity!inner Cantor capacity|ii}
Cantor capacity except that it is defined in terms of upper Green's functions $\Gbar(z,x_i;E_v)$. 
\index{Green's function!upper}  
Here, we briefly recall the definitions of $\Gbar(z,x_i;E_v)$ and $\gammabar(\EE,\fX)$
and some of their properties; 
they are studied in detail in \S\ref{Chap3}.\ref{UpperGreenSection} 
and \S\ref{Chap3}.\ref{CantorCapacitySection} below.

Upper Green's functions are gotten by taking decreasing limits of Green's functions 
\index{Green's function!upper}
of compact sets.  For an arbitrary $E_v \subset \cC_v(\CC_v)$, 
if $\zeta \notin E_v$ the upper Green's function is
\begin{equation} \label{FGHK} 
\Gbar(z,\zeta;E_v) \ = \ 
\inf_{\substack{ H_v \subset E_v \\ \text{$H_v$ compact} }} G(z,\zeta;H_v) \ .
\end{equation}
If $\zeta$ is not in the closure of $E_v$, the upper Robin constant $\Vbar_{\zeta}(E_v)$ is finite 
\index{Robin constant!upper}
and is defined by 
\begin{equation} \label{UpGreenLimit}
\Vbar_{\zeta}(E_v) \ = \ \lim_{z \rightarrow \zeta} \Gbar(z,\zeta;E_v) 
+ \log_v(|g_{\zeta}(z)|_v) \ ,
\end{equation}  
where $g_{\zeta}(z)$ is the uniformizer from (\ref{FDDA}).  
By (\ref{GreenCompactIneq}), if $E_v$ is compact then by (\cite{RR1}, Theorem 4.4.4) 
$\Gbar(z,\zeta;E_v) = G(z,\zeta;E_v)$ and $\Vbar_{\zeta}(E_v) = V_{\zeta}(E_v)$. 
For nonarchimedean $v$, if $E_v$ is algebraically capacitable in the sense of (\cite{RR1}), 
\index{algebraically capacitable}
then  $\Gbar(z,\zeta;E_v) = G(z,\zeta;E_v)$ and $\Vbar_{\zeta}(E_v) = V_{\zeta}(E_v)$.
The upper Green's function is symmetric and nonnegative: 
\index{Green's function!upper}
for all $z, \zeta \notin E_v$,  $\Gbar(z,\zeta;E_v) = \Gbar(\zeta,z;E_v) \ge 0$.
It has functoriality properties under pullbacks and base extension similar to those of $G(z,\zeta;E_v)$.

\smallskip
Now assume that each $E_v$ is stable under $\Aut_c(\CC_v/K_v)$, and that
$\EE = \prod_v E_v$ is compatible with $\fX$.\index{compatible with $\fX$}
Let $L/K$ be a finite normal extension containing $K(\fX)$. 
For each place $v$ of $K$ and each place $w$ of $L$ with $w|v$, 
after fixing an isomorphism $\CC_w \cong \CC_v$,
we can pull back $E_v$ to a set $E_w \subset \cC_w(\CC_w)$,  
which is independent of the isomorphism chosen.
If we identify $\cC_v(\CC_v)$ with $\cC_w(\CC_w)$, then for $z, \zeta \notin E_v$
\begin{equation} \label{FDDB2} 
\Gbar(z,\zeta;E_w) \log(q_w) \ = \ [L_w:K_v] \cdot \Gbar(z,\zeta;E_v) \log(q_v) \ .
\end{equation}

For each $x_i \in \fX$, 
fix a global uniformizing parameter\index{uniformizing parameter!normalizes Robin constant} $g_{x_i}(x) \in L(\cC)$
and use it to define the upper Robin constants $\Vbar_{x_i}(E_w) $ for all places $w$ of $L$.
\index{Robin constant!upper}  
For each $w$, the `local upper Green's matrix' is 
\index{Green's matrix!local!upper local|ii} 
\begin{equation*}
\Gammabar(E_w,\fX) \ = \
      \left( \begin{array}{cccc} 
                \Vbar_{x_1}(E_w)   & \Gbar(x_1,x_2;E_w) & \cdots & \Gbar(x_1,x_m;E_w) \\
                \Gbar(x_2,x_1;E_w) &  \Vbar_{x_2}(E_w)  & \cdots & \Gbar(x_2,x_m;E_w) \\
                    \vdots     &   \vdots       & \ddots &     \vdots     \\
                \Gbar(x_m,x_1;E_w) & \Gbar(x_m,x_2;E_w) & \cdots &  \Vbar_{x_m}(E_w) 
             \end{array} \right)   \ ,                                     
\end{equation*}
and the `global upper Green's matrix' is 
\index{Green's matrix!global!upper global|ii}  
\begin{eqnarray*}
\Gammabar(\EE,\fX) 
& = & \frac{1}{[L:K]} \sum_{w \in \cM_L} \Gammabar(E_w,\fX) \log(q_w) \ . 
\end{eqnarray*} 
Since $\EE$ is compatible with $\fX$,\index{compatible with $\fX$} all but finitely many of the $\Gammabar(E_w,\fX)$ are $0$. 
By the product formula, $\Gammabar(\EE,\fX)$ 
is independent of the choice of the $g_{x_i}(z)$.  By (\ref{FDDB}) it is
independent of the choice of $L$.    
It is symmetric and non-negative off the diagonal; 
its entries are finite if and only if each $E_v$ has positive inner capacity.  
\index{capacity!inner}         

\smallskip
For each $K$-rational $\EE$ compatible with $\fX$,\index{compatible with $\fX$} 
the {\it inner Cantor capacity} is 
\index{capacity!inner Cantor capacity|ii}  
\begin{equation*}
\gammabar(\EE,\fX) \ = \ e^{-\Vbar(\EE,\fX)} \ , 
\end{equation*} 
where $\Vbar(\EE,\fX) = \val(\Gammabar(\EE,\fX))$ is the value of 
$\Gammabar(\EE,\fX)$ as a matrix game.  
When the sets $E_v$ are compact or algebraically capacitable, the inner Cantor capacity 
\index{algebraically capacitable}
\index{capacity!inner Cantor capacity}
coincides with the Cantor capacity $\gamma(\EE,\fX)$ defined in (\cite{RR1}).
\index{capacity!Cantor capacity}  
It reduces to the classical logarithmic capacity
\index{capacity!logarithmic} 
when $\cC = \PP^1/\QQ$, $\fX = \infty$, and all the nonarchimedean
$E_v$ are trivial.  

The reason the inner Cantor capacity is the appropriate capacity 
\index{capacity!inner Cantor capacity}
to use in the Fekete-Szeg\"o theorem is that 
\index{Fekete-Szeg\"o theorem with LRC!for quasi-neighborhoods} 
one of the initial reductions in the proof 
is to replace each $E_v$ which is not 
$\fX$-trivial by a compact set $H_v \subset E_v$.  
\index{$\fX$-trivial}
Since the Green's function is a limit of Green's functions of
\index{Green's function}
compact sets, this can be done in such a way that  
$\Gamma(\EE,\fX)$ remains negative definite.  
\index{Green's matrix!global}\index{Green's matrix!negative definite}

\begin{definition} \label{QuasiNeighborhood} 
Let $v$ be a place of $K$.   
A set $U_v \subset \cC_v(\CC_v)$ will be called a 
{\em quasi-neighborhood} 
\index{quasi-neighborhood|ii}
if there are open sets 
$U_{v,0}, U_{v,1}, \ldots, U_{v,D}$ in $\cC_v(\CC_v)$ and algebraic 
extensions $F_{w_1}/K_v, \ldots, F_{w_D}/K_v$  in $\CC_v$ 
(possibly of infinite degree) such that
\begin{equation*} 
U_v \ = \ U_{v,0} \cup \bigcup_{\ell=1}^D \big(U_{v,\ell} \cap \cC_v(F_{w_\ell})\big) \ .
\end{equation*} 
We allow the possibility that one or more of the $U_{v,\ell}$ are empty.
We will say that $U_v$ is {\em $K_v$-symmetric} if it is 
\index{$K_v$-symmetric!quasi-neighborhood|ii}
stable under $\Aut_c(\CC_v/K_v)$, 
and that it is {\em separable} if each $F_{w_\ell}/K_v$ is separable.  
If $U_v$ contains a set $E_v$,
we will say that $U_v$ is a quasi-neighborhood of $E_v$.  
\index{quasi-neighborhood}
\end{definition} 

Equivalently, a quasi-neighborhood $U_v \subset \cC_v(\CC_v)$ is the union 
\index{quasi-neighborhood}
finitely many sets, each of which is either open in $\cC_v(\CC_v)$ or 
is open in $\cC_v(F_{w_\ell})$ for some algebraic extension $F_{w_\ell}/K_v$ in $\CC_v$.   
Note that these sets need not be disjoint.  For example, take $\cC = \PP^1$ 
and identify $\PP^1(\CC_v)$ with $\CC_v \cup \infty$.  Suppose $v$ is nonarchimedean;  
let $F_{w_1},\ldots, F_{w_D}$ be algebraic extensions of $K_v$ contained in $\CC_v$, 
and let $\cO_{w_1}, \ldots, \cO_{w_D}$ be their rings of integers.  
Then the set $U_v = \cO_{v_1} \cup \cdots \cup \cO_{w_D}$ is a quasi-neighborhood of the origin
\index{quasi-neighborhood}
in $\PP^1(\CC_v)$.

\smallskip
If $\EE = \prod_v E_v \subseteq \prod_v \cC_v(\CC_v)$ is an adelic set, 
we will say that a set $\UU = \prod_v U_v \subseteq \prod_v \cC_v(\CC_v)$
is a {\em $K$-rational separable quasi-neighborhood of $\EE$} if each $U_v$ is a 
\index{quasi-neighborhood!separable|ii}
separable quasi-neighborhood of $E_v$, stable under $\Aut_c(\CC_v/K_v)$.

\begin{theorem} [FSZ with LRC for Quasi-neighborhoods] \label{aT1-A1}
\index{quasi-neighborhood}
\index{Fekete-Szeg\"o theorem with LRC!for quasi-neighborhoods|ii}
Let $K$ be a global field, 
and let $\cC/K$ be a smooth, connected, projective curve.
Let $\fX = \{x_1, \ldots, x_m\} \subset \cC(\tK)$ 
be a finite set of points stable under $\Aut(\tK/K)$, and let 
$\EE = \prod_v E_v \subset \prod_v \cC_v(\CC_v)$
be an adelic set compatible with $\fX$,\index{compatible with $\fX$} 
such that each $E_v$ is stable under $\Aut_c(\CC_v/K_v)$. 
 
Suppose $\gammabar(\EE,\fX) > 1$.
Then for any $K$-rational separable quasi-neighborhood $\UU$ of $\EE$, there 
\index{quasi-neighborhood!separable}
are infinitely many points $\alpha \in \cC(\tK^{\sep})$ 
such that for each $v \in \cM_K$,
the $\Aut(\tK/K)$-conjugates of $\alpha$ all belong to $U_v$.  
\end{theorem} 

\smallskip
Our next variant is a stronger, but more technical, version of Theorem \ref{aT1}, 
which requires that the points produced have all their conjugates in $\EE$.
It uses the inner capacity,
\index{capacity!inner} 
and weakens the conditions on the sets $E_v$ needed for local rationality conditions.  

\smallskip
Write $\cl(E_v)$ for the closure of $E_v$ in $\cC_v(\CC_v)$.\label{`SymbolIndexcl'}  
If $v$ is an archimedean place of $K$, and a set $E_v \subset \cC_v(\CC)$ 
and a subset $E_v^{\prime} \subset E_v$ are given,   
we will say that a point $z_0 \in E_v$ is {\em analytically accessible} from $E_v^{\prime}$ if 
for some $r > 0$, there is a non-constant analytic map $f : D(0,r)^- \rightarrow \cC_v(\CC)$ with
$f(0) = z_0$, such that $f((0,r)) \subset E_v^{\prime}$.  
(See Definition \ref{AnalyticallyAccessible}.) 

\begin{theorem}[Strong FSZ with LRC, producing points in $\EE$] \label{aT1-A}  
\index{Fekete-Szeg\"o theorem with LRC!Strong form|ii}
Let $K$ be a global field, 
and let $\cC/K$ be a smooth, geometrically integral projective curve.
Let $\fX = \{x_1, \ldots, x_m\} \subset \cC(\tK)$ 
be a finite set of points stable under $\Aut(\tK/K)$, and let
$\EE = \prod_v E_v \subset \prod_v \cC_v(\CC_v)$ be an adelic set compatible with $\fX$,\index{compatible with $\fX$} 
such that each $E_v$ is stable under $\Aut_c(\CC_v/K_v)$.  
Let $S \subset \cM_K$ be a finite set of places $v$, containing all archimedean $v$,
such that $E_v$ is $\fX$-trivial for each $v \notin S$.
\index{$\fX$-trivial}

Assume that $\gammabar(\EE,\fX) > 1$.  
Assume also that for each $v \in S$,  
there is a $($possibly empty$)$ $\Aut_c(\CC_v/K_v)$-stable Borel subset 
$e_v \subset \cC_v(\CC_v)$ of inner capacity $0$ such that 
\index{capacity!inner}

$(A)$ If $v$ is archimedean and $K_v \cong \CC$, 
then each point of $\cl(E_v) \backslash e_v$ is analytically accessible 
from the $\cC_v(\CC)$-interior of $E_v$.

$(B)$ If $v$ is archimedean and $K_v \cong \RR$, 
then each point of $\cl(E_v) \backslash e_v$ is  
 
\quad $(1)$ analytically accessible from the $\cC_v(\CC)$-interior of $E_v$, or 

\quad $(2)$ is an endpoint of an open segment contained in $E_v \cap \cC_v(\RR)$.   

$(C)$ If $v$ is nonarchimedean, then $E_v$ is the disjoint union of $e_v$ 
and finitely many sets $E_{v,1}, \ldots, E_{v,D_v}$, where each $E_{v,\ell}$ is 

\quad $(1)$ open in $\cC_v(\CC_v)$, or 

\quad $(2)$ of the form $U_{v,\ell} \cap \cC_v(F_{w_\ell})$, where $U_{v,\ell}$ is open 
in $\cC_v(\CC_v)$ and  $F_{w_\ell}$ is a separable algebraic extension of $K_v$ 
contained in $\CC_v$ $($possibly of infinite degree$)$.  

Then there are infinitely many points $\alpha \in \cC(\tK^{\sep})$ such that for each $v \in \cM_K$, 
the $\Aut(\tK/K)$-conjugates of $\alpha$ all belong to $E_v$.   
\end{theorem} 

Note that if $v$ is archimedean, then the set $e_v$ in Theorem \ref{aT1-A}
can be taken to belong to $\partial E_v$,
since trivially each point of the $\cC_v(\CC)$-interior of $E_v$ or the 
$\cC_v(\RR)$-interior of $E_v \cap \cC_v(\RR)$ is analytically accessible.  Any countable set has 
inner capacity $0$, so the conditions in Theorem \ref{aT1} imply those in Theorem \ref{aT1-A}.
\index{capacity!inner} 

If $v$ is nonarchimedean, note that $\RL$-domains and balls $B(a,r)^-$, $B(a,r)$,
\index{$\RL$-domain} 
are both open and closed in the $\cC_v(\CC_v)$-topology.  
Thus if $E_v$ is a finite union of sets which are $\RL$-domains, open or closed balls, 
\index{$\RL$-domain} 
or their intersections with $\cC_v(F_{w,i})$ for separable algebraic 
extensions $F_{w,i}/K_v$ in $\CC_v$, then the theorem applies with $e_v = \phi$.

For an example of an archimedean set satisfying the conditions of Theorem \ref{aT1-A} 
but not Theorem \ref{aT1}, take $K = \QQ$, $\cC = \PP^1$, and let $v$ be the archimedean place of $\QQ$.
Identify $\PP^1(\CC)$ with $\CC \cup \infty$, 
and take $E_v = \{0\} \cup \big(\bigcup_{n=2}^{\infty} D(2/n,1/n^2)\big)$.  
Then each point of $E_v \backslash \{0\}$
is analytically accessible from $E_v^0$.   
For an example where the conditions of Theorem \ref{aT1-A} fail, 
let $E_v$ be the union of a circle $C(0,r)$ 
and countably many pairwise disjoint discs $D(a_i,r_i)$ contained in $D(0,r)^-$  
chosen in such a way that each point of $C(0,r)$ is a limit point of those discs.

\smallskip
Our third variant is a version of Theorem \ref{aT1} which adds side conditions concerning ramification.    
It says that at a finite number of places outside $S$  
we can require that the algebraic numbers produced are ramified or unramified, ``for free''.

\begin{theorem}[FSZ with LRC and Ramification Side Conditions] \label{FSZi} 
\index{Fekete-Szeg\"o theorem with LRC!and Ramification Side Conditions|ii}
Let $K$ be a global field, 
and let $\cC/K$ be a smooth, connected, projective curve.
Let $\fX = \{x_1, \ldots, x_m\} \subset \cC(\tK)$ 
be a finite, Galois-stable set of points, and let
$\EE = \prod_v E_v \subset \prod_v \cC_v(\CC_v)$
be an adelic set compatible with $\fX$,\index{compatible with $\fX$}
such that each $E_v$ is stable under $\Aut_c(\CC_v/K_v)$.  

Let $S, S^{\prime}, S^{\prime \prime} \subset \cM_K$ be finite $($possibly empty$)$ 
sets of places of $K$ which are pairwise disjoint, such that the places in
$S^{\prime} \cup S^{\prime \prime}$ are nonarchimedean.  
Assume that $\gammabar(\EE,\fX) > 1$, and that 

$(A)$ for each $v \in S$, the set $E_v$ satisfies the conditions of Theorem $\ref{aT1}$ 
or Theorem $\ref{aT1-A}$.
 
$(B)$ for each $v \in S^{\prime}$, either $E_v$ is $\fX$-trivial, 
\index{$\fX$-trivial}
or $E_v$ is a finite union of 
closed isometrically parametrizable balls $B(a_i,r_i)$ 
\index{isometrically parametrizable ball}
whose radii belong to the value group of $K_v^{\times}$
and whose centers belong to an unramified extension of $K_v$;   
 
$(C)$ for each $v \in S^{\prime\prime}$, 
either $E_v$ is $\fX$-trivial and $E_v \cap \cC_v(K_v)$ is nonempty, 
\index{$\fX$-trivial}
or $E_v$ is a finite union of closed and$/$or open isometrically parametrizable 
\index{isometrically parametrizable ball}
balls $B(a_i,r_i)$, $B(a_j,r_j)^-$ with centers in $\cC_v(K_v)$.

Then there are infinitely many  points $\alpha \in \cC(\tK^{\sep})$ such that 

$(1)$ for each $v \in \cM_K$, the $\Aut(\tK/K)$-conjugates of $\alpha$ 
all belong to $E_v$;   


$(2)$ for each $v \in S^{\prime}$, each place of $K(\alpha)/K$ above $v$ is unramified over $v$; 

$(3)$ for each $v \in S^{\prime \prime}$, 
each place of $K(\alpha)/K$ above $v$ is totally ramified over $v$.
  
\end{theorem} 

\vskip .1 in
Our fourth variant involves a partial converse to the Fekete-Szeg\"o theorem, known as Fekete's theorem, 
\index{Fekete-Szeg\"o theorem}
\index{Fekete's theorem}
which asserts that if $\gamma(\EE,\fX) < 1$ then for a sufficiently small neighborhood 
$\UU$ of $\EE$,  there are only finitely many points $\alpha \in \cC(\tK)$ whose conjugates
all belong to $\UU$.  Fekete's theorem on curves is proved in (\cite{RR1}, Theorem 6.3.1).
However, Fekete's theorem requires a different notion of capacity 
than we have been using here: it concerns the ``outer capacity'' $\underline{\gamma}(\EE,\fX)$,
rather than the inner capacity $\gammabar(\EE,\fX)$. 
\index{capacity!inner Cantor capacity} 

Extending the definition of algebraic capacitability in (\cite{RR1}) to both archimedean and nonarchimedean sets, we will say that $E_v$ {\em algebraically capacitable} if it is closed
\index{algebraically capacitable|ii}
in $\cC_v(\CC_v)$ and $\gammabar_{\zeta}(E_v) = \underline{\gamma}_{\zeta}(E_v)$ 
for each $\zeta \notin E_v$.  If each $E_v$ is algebraically capacitable, 
\index{algebraically capacitable}
then $\gammabar(\EE,\fX)$ and  $\underline{\gamma}(\EE,\fX)$ are equal, 
and coincide with the capacity $\gamma(\EE,\fX)$ in (\cite{RR1}).  Here
\index{capacity!Cantor capacity}
\begin{equation*}
\gammabar_{\zeta}(E_v) \ = \ 
\sup_{\substack{ H_v \subset E_v \\ \text{$H_v$ compact} }} \gamma(H_v) \ , \qquad 
\underline{\gamma}_{\zeta}(E_v) 
        \ = \ \inf_{\substack{ U_v \supset E_v \\ \text{$U_v$ a $PL_{\zeta}$-domain}}} 
\gammabar(U_v) \ .
\end{equation*}
A set $U_v$ is a $PL_{\zeta}$-domain if there is a nonconstant rational function 
$f(z) \in \CC_v(\cC)$, whose only poles are at $\zeta$, 
for which $U_v = \{z \in \cC_v(\CC_v) : |f(z)|_v \le 1\}$.                               
In the nonarchimedean case, the compatibility of this definition with the one given in 
(\cite{RR1}, p.259) follows from (\cite{RR1}, Propositions 4.3.1 and 4.3.16).  
In (\cite{RR1}), algebraic capacitability was not defined in the archimedean case,  
but all archimedean sets were required to be compact. 

If $v$ is archimedean, it follows from (\cite{RR1}, Proposition 3.3.3) 
that every compact set is algebraically capacitable.  
\index{algebraically capacitable}
If $v$ is nonarchimedean, it is shown in (\cite{RR1}, Theorem 4.3.13) that 
 any set $E_v$ which can be expressed as a finite combination of unions and intersections 
of compact sets and $\RL$-domains, is algebraically capacitable.
\index{algebraically capacitable}
\index{$\RL$-domain}


Assuming algebraic capacitability for the sets $E_v$, the following result   
describes the dichotomy provided by Fekete's theorem and 
\index{Fekete's theorem}
the Fekete-Szeg\"o theorem in terms of the Green's matrix $\Gamma(\EE,\fX)$. 
\index{Green's matrix!global}\index{Fekete-Szeg\"o theorem} 
Recall (see \cite{RR1}, \S5.1) that $\gamma(\EE,\fX) > 1 $ 
if and only if $\Gamma(\EE,\fX)$ is negative definite,\index{Green's matrix!negative definite} 
and that $\gamma(\EE,\fX) < 1$ if and only if when the rows and columns of $\Gamma(\EE,\fX)$
are permuted to bring $\Gamma(\EE,\fX)$ into block diagonal form, 
then some eigenvalue of each block is positive.     

\begin{theorem}[Fekete/Fekete-Szeg\"o with LRC for Algebraically Capacitable Sets] \label{FSZii}
\index{Fekete-Szeg\"o theorem with LRC!for algebraically capacitable sets|ii} 
\index{algebraically capacitable}
Let $K$ be a global field 
and let $\cC/K$ be a smooth, connected, projective curve.
Let $\fX = \{x_1, \ldots, x_m\} \subset \cC(\tK)$ 
be a finite, galois-stable set of points, and let
$\EE = \prod_v E_v \subset \prod_v \cC_v(\CC_v)$
be an adelic set compatible with $\fX$.\index{compatible with $\fX$} 

Assume that each $E_v$ is algebraically capacitable and stable under $\Aut_c(\CC_v/K_v)$.  Then 
\index{algebraically capacitable} 

\smallskip
$(A)$  If all the eigenvalues of $\Gamma(\EE,\fX)$ are non-positive 
$($that is, $\Gamma(\EE,\fX)$ is either negative definite\index{Green's matrix!negative definite} 
or negative semi-definite$)$,  
let $\UU = \prod_v U_v$ be a separable $K$-rational quasi-neighborhood of $\EE$ 
\index{quasi-neighborhood!separable}
such that there is at least one place $v_0$ where $E_{v_0}$ is compact and the 
quasi-neighborhood $U_{v_0}$ properly contains $E_{v_0}$. If $v_0$ is archimedean, assume also 
\index{quasi-neighborhood}
that $U_{v_0}$ meets each component of $\cC_{v_0}(\CC) \backslash E_{v_0}$ containing a point of $\fX$.  
Then there are infinitely many points $\alpha \in \cC(\tK^{\sep})$ such that 
all the conjugates of $\alpha$ belong to $\UU$.  

\smallskip
$(B)$ If some eigenvalue of $\Gamma(\EE,\fX)$ is positive $($that is, $\Gamma(\EE,\fX)$ is 
either indefinite, nonzero and positive semi-definite, or positive definite$)$,
there is an adelic neighborhood $\UU$ of $\EE$ such that only 
finitely many points $\alpha \in \cC(\tK)$ have all their conjugates in $\UU$.
\end{theorem}

Finally, we formulate two Berkovich versions of the Fekete-Szeg\"o Theorem 
\index{Fekete-Szeg\"o theorem with LRC!Berkovich} 
with local rationality conditions. 

\smallskip   
For each nonarchimedean place $v$ of $K$, 
let $\cC_v^{\an}$ be the Berkovich analytic space associated to 
\label{`SymbolIndexcCan'}
\index{Berkovich!analytic space} 
$\cC_v \times_{K_v} \Spec(\CC_v)$ (see \cite{Berk}).  This is a locally ringed space 
whose underlying topological space is a compact, path connected Hausdorff space 
with $\cC_v(\CC_v)$ as a dense subset;  it has the same sheaf 
of functions as the rigid analytic space\index{rigid analytic space}
 associated to $\cC_v \times_{K_v} \Spec(\CC_v)$. 
In his doctoral thesis, Amaury Thuillier (\cite{Th}) constructed 
\index{Thuillier, Amaury}\index{potential theory!on Berkovich curves}
a potential theory on $\cC_v^{\an}$  which includes a $dd^c$ operator, 
harmonic functions, subharmonic functions, capacities, and Green's functions.
\index{Green's function!Berkovich}  
When $\cC \cong \PP^1$, Baker and Rumely (\cite{B-R}) 
\index{Baker, Matthew}\index{Rumely, Robert} 
constructed a similar theory in an elementary way. 

In what follows, 
we assume familiarity with Berkovich analytic spaces and Thuillier's theory.
\index{Berkovich!analytic space}  
\index{Thuillier, Amaury} 
\index{Berkovich, Vladimir} 
For each compact, non-polar subset $\BerkE_v \subset \cC_v^{\an}$ 
\label{`SymbolIndexBerkEv'}
and each $\zeta \in \cC_v^{\an} \backslash \BerkE_v $, Thuillier (\cite{Th}, Th\'eor\`eme 3.6.15) 
\index{Thuillier, Amaury}
has constructed a Green's function $g_{\zeta,\BerkE_v}(z)$ which is non-negative, 
\index{Green's function!Thuillier}
vanishes on $\BerkE_v$ except possibly on a set of capacity $0$, is subharmonic in $\cC_v^{\an}$,
\index{capacity!inner}
harmonic in $\cC_v^{\an} \backslash (\BerkE_v \cup \{\zeta\})$, 
and satisfies the distributional equation $dd^c g_{\zeta,\BerkE_v} = \mu - \delta_\zeta$ 
where $\mu$ is a probability measure supported on $K$. We will write 
$G(z,\zeta;\BerkE_v)^{\an}$ for $g_{\zeta,\BerkE_v}(z)$, and regard it as a function of 
two variables. 
\label{`SymbolIndexBerkGreen'}
By Proposition \ref{BerkGreenPropertiesProp} below, 
for all $z, \zeta \in \cC_v^{\an} \backslash \BerkE_v$ with $z \ne \zeta$, 
\begin{equation*} 
G(z,\zeta;\BerkE_v)^{\an} \ = \ G(\zeta,z;\BerkE_v)^{\an} \ , 
\end{equation*} 
and for each $\zeta \in \cC_v(\CC_v) \backslash \BerkE_v$, the Robin constant
\index{Robin constant!Berkovich|ii}
\label{`SymbolIndexBerkVZeta'}
\begin{equation*} 
V_{\zeta}(\BerkE_v)^{\an} \ = \ \lim_{\substack{z \rightarrow \zeta \\ z \in \cC_v^{\an}}} 
G(z,\zeta;\BerkE_v)^{\an} + \log(|g_{\zeta}(z)|_v) 
\end{equation*} 
exists. 
The group $\Aut_c(\CC_v/K_v)$ acts on $\cC_v^{\an}$ in a natural way, 
and for all $\sigma \in \Aut_c(\CC_v/K_v)$
\begin{equation*}
G(\sigma(z),\sigma(\zeta);\sigma(\BerkE_v))^{\an} \ = \ G(z,\zeta;\BerkE_v)^{\an} \ .
\end{equation*}
By Proposition \ref{BerkCompatibilityProp} below, the Green's functions $G(z,\zeta;\BerkE_v)^{\an}$ 
\index{Green's function!Berkovich}
and the functions $G(z,\zeta;E_v)$ from this work are compatible up to a normalizing factor, 
in the sense that if $E_v \subset \cC_v(\CC_v)$ is algebraically capacitable 
\index{algebraically capacitable}
(in particular, if $E_v$ is a finite union of $\RL$-domains and compact sets), 
\index{$\RL$-domain} 
and if $\BerkE_v $ is the closure of $E_v$ in $\cC_v^{\an}$ for the Berkovich topology, 
\index{Berkovich!topology}
then for all $z, \zeta \in \cC_v(\CC_v) \backslash E_v$, 
\begin{equation*}
G(z,\zeta;\BerkE_v)^{\an} \ = \ G(z,\zeta;E_v) \log(q_v) \ .
\end{equation*} 

If $v$ is an archimedean place of $K$, we take $\cC_v^{\an}$ to be the Riemann surface
\index{Riemann surface}  
$\cC_v(\CC)$, and for a set $\BerkE_v = E_v \subset \cC_v(\CC)$ we put 
$G(z,\zeta;\BerkE_v)^{\an} = G(z,\zeta;E_v)$ and $V_{\zeta}(\BerkE_v)^{\an} = V_{\zeta}(E_v)$.  

Let $\fX = \{x_1,\ldots, x_m\} \subset \cC(\tK)$ be a finite, galois-stable set of points.  We will 
now define the notion of a {\em compact Berkovich adelic set compatible with $\fX$}.
\index{Berkovich!adelic set}\index{compatible with $\fX$!Berkovich set compatible with $\fX$}   
For each place $v$ of $K$, let $\BerkE_v \subset\cC_v^{\an}$ 
be a compact, nonpolar set disjoint from $\fX$.
\index{nonpolar set|ii}
(A Berkovich set is nonpolar if and only if it has positive capacity:  
\index{capacity}
see (\cite{Th}, \S3.4.2 and Theorem 3.6.11).)  
We will say that $\BerkE_v$ is $\fX$-trivial if $v$ is nonarchimedean and  
\index{$\fX$-trivial}
the model $\fC_v/\Spec(\cO_v)$ from Definition \ref{Xtrivial} has good reduction, 
\index{good reduction} 
the points of $\fX$ specialize to distinct points in the special fibre $r_v(\fC_v)$, 
and $\BerkE_v$ consists of all points $z \in \cC_v^{an}$ 
whose specialization $r_v(z) \in r_v(\fC_v)$ is distinct from 
$\{r_v(x_1), \ldots, r_v(x_m)\}$.  Equivalently, $\BerkE_v$ is $\fX$-trivial if 
\index{$\fX$-trivial}
it is the closure of the $\fX$-trivial set $E_v = \cC_v(\CC_v) \backslash (\bigcup_{i=1}^m B(x_i,1)^-)$ 
in $\cC_v(\CC_v)$.  Then 
\begin{equation*} 
\EE \ := \ \prod_v \BerkE_v \ \subset \ \prod_v \cC_v^{\an}
\end{equation*}
is a compact Berkovich adelic set compatible with 
$\fX$\index{compatible with $\fX$!Berkovich set compatible with $\fX$}  
\index{Berkovich!adelic set} 
if each $\BerkE_v$ satisfies the conditions above, and $\BerkE_v$ is $\fX$-trivial for all but
\index{$\fX$-trivial}
finitely many $v$.

If $\EE$ is a compact Berkovich adelic set 
compatible with $\fX$,\index{compatible with $\fX$!Berkovich set compatible with $\fX$} 
\index{Berkovich!adelic set} 
we define the local and global Green's matrices $\Gamma(\BerkE_w,\fX)^{\an}$ and $\Gamma(\EE,\fX)^{\an}$
\index{Green's matrix!local!local Berkovich}
\index{Green's matrix!global!global Berkovich}
as in (\ref{LocalGreensMatrixDef}), (\ref{GlobalGreensMatrixDef}), replacing
$G(z,\zeta;E_v)$ by $G(z,\zeta;\BerkE_v)^{\an}$ and $V_{\zeta}(E_v)$  by $V_{\zeta}(\BerkE_v)^{\an}$, 
but omitting the weights $\log(q_v)$ at nonarchimedean places.    
We then define the global Robin constant $V(\EE,\fX)^{\an}$ using the minimax formula (\ref{aF3})
\index{Robin constant!global}\index{minimax property} 
taking $\Gamma = \Gamma(\EE,\fX)^{\an}$, and the global capacity by 
\label{`SymbolIndexBerkCantorCap'}
\index{capacity!global}
\begin{equation*}
\gamma(\EE,\fX)^{\an} \ = \ e^{-V(\EE,\fX)^{\an}} \ .
\end{equation*} 

We will call a set 
\begin{equation*} 
\UU \ = \ \prod_v \BerkU_v \ \subset \ \prod_v \cC_v^{\an} 
\end{equation*}
a {\em Berkovich adelic neigbhorhood} of $\EE$ if $\BerkU_v$ contains $\BerkE_v$ for each $v$, 
\index{Berkovich!adelic neighborhood} 
and either $\BerkU_v$ is 
an open set in $\cC_v^{\an}$, or $\BerkE_v$ is $\fX$-trivial and $\BerkU_v = \BerkE_v$.
\index{$\fX$-trivial}
We will call $\UU$ a {\em separable Berkovich quasi-neigbhorhood} of $\EE$ 
\index{quasi-neighborhood!Berkovich}
\index{Berkovich!quasi-neighborhood} 
if $\BerkU_v$ contains $\BerkE_v$ for each $v$,
and either $\BerkU_v$ is the union of a Berkovich open set and finitely many 
\index{Berkovich!open set} 
open sets in $\cC_v(F_w)$ for finite separable extensions $F_w/K_v$, 
or $\BerkE_v$ is $\fX$-trivial and $\BerkU_v = \BerkE_v$. We will say that $\UU$ is {\em $K$-rational}
\index{$\fX$-trivial}
if each $U_v$ is stable under $\Aut_c(\CC_v/K_v)$.

The following is the Berkovich analogue of Theorem \ref{aT1}:

\begin{theorem}[Berkovich FSZ with LRC, producing points in $\EE$] 
\index{Fekete-Szeg\"o theorem with LRC!Berkovich|ii}
\label{aT1-B1} \ 

Let $K$ be a global field, 
and let $\cC/K$ be a smooth, geometrically integral, projective curve.
Let $\fX = \{x_1, \ldots, x_m\} \subset \cC(\tK)$ 
be a finite set of points stable under $\Aut(\tK/K)$, and let
$\EE = \prod_v \BerkE_v \subset \prod_v \cC_v^{\an}$ be a $K$-rational Berkovic adelic 
set compatible with $\fX$.\index{compatible with $\fX$!Berkovich set compatible with $\fX$}   
Let $S \subset \cM_K$ be a finite set of places $v$, containing all archimedean $v$,
such that $\BerkE_v$ is $\fX$-trivial for each $v \notin S$.
\index{$\fX$-trivial}

Assume that $\gamma(\EE,\fX) > 1$.  
Assume also that $\BerkE_v$ has the following form, for each $v \in S$:   

$(A)$ If $v$ is archimedean and $K_v \cong \CC$, 
then $\BerkE_v$ is compact, and is a finite union of sets $E_{v,i}$, 
each of which is the closure of its $\cC_v(\CC)$-interior and has a piecewise smooth boundary;
\index{boundary!piecewise smooth}\index{closure of $\cC_v(\CC)$ interior} 

$(B)$ If $v$ is archimedean and $K_v \cong \RR$, then $\BerkE_v$ is compact, 
stable under complex conjugation, 
and is a finite union of sets $E_{v,\ell}$, where each $E_{v,\ell}$ is either 

\quad $(1)$ the closure of its $\cC_v(\CC_)$-interior and has a piecewise smooth boundary, or
\index{boundary!piecewise smooth}\index{closure of $\cC_v(\CC)$ interior} 

\quad $(2)$ is a compact, connected subset of $\cC_v(\RR)$; 

$(C)$ If $v$ is nonarchimedean, then $\BerkE_v$ is compact, stable under $\Aut_c(\CC_v/K_v)$, 
and is a finite union of sets $E_{v,\ell}$, where each $E_{v,\ell}$ is either 

\quad $(1)$ a strict closed Berkovich affinoid, or
\index{Berkovich!strict closed affinoid} 
\index{affinoid!Berkovich affinoid} 

\quad $(2)$ is a compact subset of $\cC_v(\CC_v)$ 
and has the form $\cC_v(F_{w_\ell}) \cap B(a_\ell,r_\ell)$ 
for some finite separable extension $F_{w_\ell}/K_v$ in $\CC_v$, and some ball $B(a_\ell,r_\ell)$.  

Then there are infinitely many points $\alpha \in \cC(\tK^{\sep})$ such that for each $v \in \cM_K$, 
the $\Aut(\tK/K)$-conjugates of $\alpha$ all belong to $\BerkE_v$.  
\end{theorem} 

Finally, we give a Berkovich version of the Fekete-Szeg\"o Theorem 
with local rationality conditions for quasi-neighborhoods, 
\index{quasi-neighborhood!Berkovich}
generalizing Theorem \ref{aT1-A1} and (\cite{B-R}, Theorem 7.48): 

\begin{theorem} [Berkovich Fekete/FSZ with LRC for Quasi-neighborhoods] \label{aT1-B2}
\index{Fekete-Szeg\"o theorem with LRC!for Berkovich quasi-neighborhoods|ii}
Let $K$ be a global field, 
and let $\cC/K$ be a smooth, connected, projective curve.
Let $\fX = \{x_1, \ldots, x_m\} \subset \cC(\tK)$ 
be a finite set of points stable under $\Aut(\tK/K)$, and let 
$\EE = \prod_v \BerkE_v \subset \prod_v \cC_v^{\an}$
be a compact Berkovich adelic set 
compatible with $\fX$,\index{compatible with $\fX$!Berkovich set compatible with $\fX$} 
\index{Berkovich!adelic set} 
such that each $\BerkE_v$ is stable under $\Aut_c(\CC_v/K_v)$. 

$(A)$ If $\gamma(\EE,\fX)^{\an} < 1$, 
there is a $K$-rational Berkovich neighborhood $\UU = \prod_v \BerkU_v$ of $\EE$
\index{Berkovich!neighborhood} 
such that there are only finitely many points of $\cC(\tK)$ 
whose $\Aut(\tK/K)$-conjugates are all contained in $\BerkU_v$, 
for each $v \in \cM_K$.

$(B)$ If $\gamma(\EE,\fX)^{\an} > 1$,
then for any $K$-rational 
separable Berkovich quasi-neighborhood $\UU$ of $\EE$, there 
\index{quasi-neighborhood!Berkovich}
\index{Berkovich!quasi-neighborhood} 
are infinitely many points $\alpha \in \cC(\tK^{\sep})$ 
such that for each $v \in \cM_K$,
the $\Aut(\tK/K)$-conjugates of $\alpha$ all belong to $\BerkU_v$.  
\end{theorem}

%% file: NewFSZChap2.tex
\chapter{Examples and Applications} \label{ExamApp}

In this chapter we illustrate the Fekete-Szeg\"o theorem with local rationality conditions. 
\index{Fekete-Szeg\"o theorem with LRC}
We first apply it on $\PP^1$, using it to construct algebraic integers and algebraic 
units satisfying various conditions.  We then apply it on elliptic curves\index{elliptic curve}, 
Fermat curves\index{Fermat curve}, and modular curves\index{modular curve}.  


\section{ Local capacities and Green's functions of Archimedean Sets} \label{ArchLocalExamples}  

Suppose $K_v = \RR$ or $K_v = \CC$. 
In this section we give formulas for local capacities and Green's functions of sets in $\PP^1(\CC)$ 
\index{Green's function!examples!archimedean}
which arise naturally in arithmetic applications.  
Some involve closed formulas, others require numerical computations.\index{numerical!computations}  
Most of the formulas appear in the literature; only a few are new.  
Further examples, mainly concerning sets in $\CC$ with geometric symmetry, 
are given in (\cite{RR1}, pp. 348-351). 

\smallskip 
For archimedean sets, 
the most effective way of determining capacities is by ``guessing'' the Green's function: 
\index{Green's function!archimedean!guessing} 
given $E$ and $\zeta \notin E$, if a function can be found which is continuous, $0$ on $E$,
and harmonic in the complement of $E$ except for a positive logarithmic pole at $\zeta$,
then by the maximum modulus principle, it must be the Green's function. 
\index{Green's function!archimedean!characterization of|ii} 
Then, given a uniformizing parameter\index{uniformizing parameter!normalizes Robin constant} $g_{\zeta}(z)$, 
the Robin constant and capacity of $E$ with respect to $\zeta$ can be read off by 
\index{Robin constant}
\index{capacity}
\begin{equation} \label{FVCap}
V_{\zeta}(E) \ = \ \lim_{z \rightarrow \zeta} G(z,\zeta;E) + \log(|g_{\zeta}(z)|) \ , \quad
\gamma_{\zeta}(E) = e^{-V_{\zeta}(E)} \ .
\end{equation}
For the sets we are dealing with here, which are compact unions of continua, 
the upper Green's function $G(z,\zeta;E)$ coincides 
with the usual Green's function $G(z,\zeta;E)$.
\index{Green's function}
\index{Green's function!upper}   

In the discussion below, we will identify $\PP^1(\CC)$ with $\CC \cup \{\infty\}$.  When $\zeta = \infty$, 
we take $g_{\zeta}(z) = 1/z$;  when $\zeta \in \CC$, we take $g_{\zeta}(z) = z - \zeta$. 

\medskip 

{\bf The Disc.} The most basic example is when $E$ is the disc $D(0,r) \subset \CC$. 
Here
\begin{eqnarray}
G(z,\infty;E) \ = \ \log^+(|z/r|) 
\ = \ \left\{ \begin{array}{ll} \log(|z/r|) & \text{if $|z| > r$} \\
                                0           & \text{if $|z| \le r$} \end{array} \right.  \ . 
                                \label{FDisc1} 
\end{eqnarray}
\index{Green's function!examples!archimedean}\index{examples!archimedean!disc}
Computing capacities relative to the parameter $g_{\infty}(z) = 1/z$, we find 
\begin{eqnarray}
\qquad  V_{\infty}(E) = \lim_{z \rightarrow \infty} G(z,\infty;E) - \log(|z|) = -\log(r) 
\ , \label{FDisc2} \\ 
\gamma_{\infty}(E) = e^{-V_{\infty}(E)} = r \ . \qquad \qquad \qquad \notag   
\end{eqnarray}
\index{Robin constant!examples!archimedean}
By applying a linear fractional transformation, 
one can find the Green's function of $D(0,r)$
\index{Green's function!examples!archimedean}
with respect to an arbitrary point $\zeta \in \CC$: 
\begin{equation}
G(z,\zeta;E) \ = \ \log^+\left(\left|
       \frac{r^2-\overline{\zeta}z}{r(z-\zeta)}\right|\right) \ . \label{FDisc3}
\end{equation}
\index{Green's function!examples!archimedean}
Computing capacities relative to $g_{\zeta}(z) = z-\zeta$, one has
\begin{eqnarray}
V_{\zeta}(E)  \ = \ \lim_{z \rightarrow \zeta} G(z,\infty;E) + \log(|z-\zeta|) 
    \ = \ \log(\frac{|\zeta|^2-r^2}{r}) \ , \label{FDisc4} \\
\gamma_{\zeta}(E)  = e^{-V_{\zeta}(E)} = \frac{r}{|\zeta|^2-r^2} \ . 
            \qquad \qquad \qquad \quad \label{FDisc5}
\end{eqnarray}      

\smallskip
{\bf The Segment.}\index{examples!archimedean!one segment}  
Another basic example is when $E$ is a segment $[a,b] \subset \RR$.  
Choosing the branch of $\sqrt{z}$ which is positive on the positive real axis 
and cut along the negative real axis, 
the map $z \mapsto w = \sqrt{(z-a)/(z-b)}$ takes $\PP^1(\CC) \backslash [a,b]$ to the right halfplane and takes $\infty$ to $1$;  then $w \mapsto (w+1)/(w-1)$ takes the right halfplane 
to the exterior of the unit disc, and takes $1$ to $\infty$.  It follows that 
\begin{equation}
G(z,\infty;E) \ = \ -\log^+ \left(\left| \frac{\sqrt{(z-a)/(z-b)}-1}
                       {\sqrt{(z-a)/(z-b)}+1} \right| \right) \ .  \label{F1Seg1}
\end{equation} 
\index{Green's function!examples!archimedean}
For an arbitrary $\zeta \in \CC$, a similar computation (see \cite{Can3}, p.165) gives 
\begin{equation} \label{F1Seg2}
G(z,\zeta;E) \ = \ -\log^+ \left(\left| \frac{\sqrt{(z-a)/(z-b)}-{\sqrt{(\zeta-a)/(\zeta-b)}}}
                      {\sqrt{(z-a)/(z-b)}+\overline{\sqrt{(\zeta-a)/(\zeta-bB)}}}
             \right| \right) \ .
\end{equation} 
\index{Green's function!examples!archimedean}
With $g_{\infty}(z) = 1/z$ and $g_{\zeta}(z) = z-\zeta$ for $\zeta \in \CC \backslash E$, one finds 
\begin{eqnarray} 
V_{\infty}(E) \ = \ -\log((b-a)/4) \ ,\  \quad \gamma_{\infty}(E) \ = \ (b-a)/4 \ , 
\qquad \label{F1Seg3} \\
\gamma_{\zeta}(E) = e^{-V_{\zeta}(E)} = 
    \frac{b-a}{4 \cdot \Re(\sqrt{(\zeta-a)|\zeta-a| \cdot (\overline{\zeta-b})|\zeta-b|})} 
                                                  \ \quad \label{F1Seg4}
\end{eqnarray} 

When $\zeta = \infty$ there is another expression for $G(z,\infty,E)$ 
which makes its geometric behavior clearer.  
For simplicity, assume $E = [-2r,2r]$ where $0 < r \in \RR$.  
It is well known, and easy to verify, that the Joukowski map\index{Joukowski map} 
\begin{equation} \label{rJoukowski}
z \ = \ J_r(w) \ = \ w + \frac{r^2}{w}
\end{equation}\index{Green's function!examples!archimedean}
maps $\CC \backslash D(0,r)$ conformally onto $\CC \backslash [-2r,2r]$. 
For each $R > r$, 
it takes the circle $C(0,R)$ parametrized by $w = R \cos(\theta) + i R \sin(\theta)$ 
to the ellipse $\cE(R+\frac{r^2}{R},R-\frac{r^2}{R})$ parametrized by 
\begin{equation} \label{REllipseParam}
z \ = \ x+iy \ = \ (R+\frac{r^2}{R}) \cos(\theta) + i (R-\frac{r^2}{R}) \sin(\theta) 
\ = \ J_r(R \cos(\theta) + i R \sin(\theta)) \ .
\end{equation} 
It maps the circle $C(0,R)$ in a $2-1$ manner to the interval $[-2r,2r]$, 
and takes $\infty $ to $\infty$.

The function $G_r(z) = \log(|J_r^{-1}(z)|/r)$ is harmonic on $\CC \backslash E$, 
with a logarithmic pole at $\infty$;  it has a continuous extension to $\CC$ 
which takes the value $0$ on $E$.  By the characterization of Green's functions,
\index{Green's function!archimedean!characterization of}
$G(z,\infty;[-2r,2r]) = G_r(z)$.  Thus, for each $R > r$,
\index{Green's function!examples!archimedean}
\begin{equation} \label{EllipseLevelCurve} 
\{z \in \CC : G(z,\infty;[-2r,2r]) = \log(R/r)\} \ = \ \cE(R+\frac{r^2}{R},R-\frac{r^2}{R}) \ .
\end{equation} 

\smallskip
{\bf Two segments.}\index{examples!archimedean!two segments}  
When $E = [a,b] \cup [c,d] \subset \RR$, there are closed formulas for the Green's function and capacity. 
\index{Green's function!examples!archimedean}
\index{capacity}  
When the segments have the same length, $G(z,\infty;E)$ and $\gamma_{\infty}(E)$ 
are given by elementary formulas.  
In general, they can be expressed in terms of theta-functions.\index{theta-functions!classical}

First suppose $E = [-b,-a] \cup [a,b] \subset \RR$. 
Put $f(z) = z^2$;  then $f^*((\infty)) = 2(\infty)$ and $E = f^{-1}([a^2,b^2])$.  
By the pullback formula for Green's functions (see (\ref{GreenPullbackF}) below),
\index{Green's function!pullback formula for} 
\begin{equation} 
G(z,\infty;E)    
\ = \ -\frac{1}{2} \log\left(\left| \frac{\sqrt{(z^2-a^2)/(z^2-b^2)}-1}
                                {\sqrt{(z^2-a^2)/(z^2-b^2)}+1} \right| \right) \ . 
\label{F2Seg1}        
\end{equation}\index{Green's function!examples!archimedean}
Using this, we find 
\begin{eqnarray}
\qquad  V_{\infty}(E) = \frac{1}{2}\log(4/(b^2-a^2)) \ , \quad 
\gamma_{\infty}(E) = \frac{\sqrt{b^2-a^2}}{2} \ , \label{F2Seg2}  \\
G(0,\infty;E) \ = \ G(\infty,0,E) \ = \ \frac{1}{2} \log(\frac{b+a}{b-a}) \ . 
\quad \label{F2Seg3} 
\end{eqnarray}
Similarly, when $\zeta = 0$, pulling back $[1/b^2,1/a^2]$ by $f(z) = 1/z^2$, we get  
\begin{eqnarray}
G(z,0;E) \ = \ \frac{1}{2} \log\left(\left| \frac{\sqrt{(z^2-b^2)/(z^2-a^2)}+b/a}
              {\sqrt{(z^2-b^2)/(z^2-a^2)}-b/a} \right| \right)\ , 
               \qquad  \label{F2Seg4} \\
V_0(E) = \frac{1}{2}\log(4a^2b^2/(b^2-a^2)) \ ,
\quad \gamma_{0}(E) = \frac{\sqrt{b^2-a^2}}{2ab} \ . \qquad \label{F2Seg5} 
\end{eqnarray}\index{Green's function!examples!archimedean}
 
\smallskip
Before dealing with a general set $E = [a,b] \cup [c,d] \subset \RR$, and arbitrary $\zeta$, 
it will be useful to recall some of the properties 
of classical theta-functions\index{theta-functions!classical}
(see Shimura, \cite{ShimuraTheta}, \cite{WW}).\index{Shimura, Goro} 
For $u \in \CC$, $\tau \in \fH = \{ \Im(z) > 0\}$, 
and $r, s \in \RR$, write $e(z) = e^{2 \pi i z}$ and put 
\begin{equation}
\theta(u,\tau;r,s) \ = \ \sum_{n \in \ZZ} e(\frac{1}{2}(n+r)^2\tau + (n+r)(u+s)) \ . 
            \label{FCTheta1} 
\end{equation}
Because of the quadratic dependence on $n$ in (\ref{FCTheta1}), 
the series defining $\theta(u,\tau;r,s)$ converges very rapidly.  
$\theta(u,\tau;r,s)$ is continuous in all four variables 
and is jointly holomophic in $u$ and $\tau$.  

We will be particularly interested in $\theta(u,\tau;\frac{1}{2},\frac{1}{2})$.  
When $r,s \in \{0,1/2\}$, the functions $\theta(u,\tau;r,s)$ appear in the literature wih several names.
Our notation follows Krazer-Prym and Shimura; 
\index{Shimura, Goro}     
in the notation of Riemann and Mumford (respectively Whittaker-Watson \cite{WW}),
\index{Riemann, Bernhard} 
\index{Mumford, David}  
\begin{eqnarray*} 
\theta(u,\tau;{\textstyle{\frac{1}{2}}},{\textstyle{\frac{1}{2}}}) 
= \theta_{11}(u,\tau) = \vartheta_1(\pi u|\tau)\ , & \quad &
\theta(u,\tau;{\textstyle{\frac{1}{2}}},0) = \theta_{10}(u,\tau) = \vartheta_2(\pi u|\tau)\ , \\
\theta(u,\tau;0,0) = \theta_{00}(u,\tau) = \vartheta_3(\pi u|\tau) \ , & \quad & 
\theta(u,\tau;0,{\textstyle{\frac{1}{2}}}) = \theta_{01}(u,\tau) = \vartheta_4(\pi u|\tau)\ .
\end{eqnarray*}
In the notation of Courant-Hilbert (\cite{C-H}), 
\index{Courant, Richard} 
\index{Hilbert, David}
$\theta(u,\tau;0,\frac{1}{2}) = \theta_0(u)$ and   
$\theta(u,\tau;\frac{1}{2},\frac{1}{2}) = \theta_1(u)$. 

Considering $\theta(u,\tau;r,s)$ as a function of $u$ and using the definition, 
one sees that for all $a, b \in \ZZ$ 
\begin{equation}
\theta(u+za+b,\tau;r,s) = e(rb-as) \cdot e(-\frac{1}{2} a^2\tau - au) \cdot \theta(u,z;r,s) \ .
          \label{FCTheta3}  
\end{equation} 
Applying the Argument Principle, it follows that $\theta(u,\tau,r,s)$ has a simple zero
\index{Argument Principle} 
in each period parallelogram for the lattice\index{period lattice} $\langle 1,\tau \rangle \subset \CC$; 
the zero occurs at $u \equiv (\frac{1}{2}-r)\tau+ (\frac{1}{2}-s) \pmod{\langle 1,\tau \rangle}$
(see \cite{WW}, p.465-466, and \cite{ShimuraTheta}, formula (11), p.675).   

Again using the definitions, one sees that  $\theta(u,\tau;\frac{1}{2},\frac{1}{2})$ is 
an odd function of $u$, 
that $\theta(u + \frac{1}{2},\tau;\frac{1}{2},\frac{1}{2}) = -\theta(u,\tau;\frac{1}{2},0)$, 
and if $\tau$ is pure imaginary, 
then $\theta(\overline{u},\tau;\frac{1}{2},\frac{1}{2}) = 
\overline{\theta(u,\tau;\frac{1}{2},\frac{1}{2})}$.  
Similarly $\theta(u,\tau;\frac{1}{2},0)$ is 
an even function of $u$, and if $\tau$ is pure imaginary, 
then $\theta(\overline{u},\tau;\frac{1}{2},0) = 
\overline{\theta(u,\tau;\frac{1}{2},0)}$. 

With these facts, one can check that if $\tau$ is pure imaginary, 
then for each $M \in \CC$ with $\Re(M) \notin \frac{1}{2}\ZZ$, the function  
\begin{equation} \label{FFFunc} 
\cG(u) \ = \ \frac{\theta(u-M,\tau;\frac{1}{2},\frac{1}{2})}
      {\theta(u+\overline{M},\tau;\frac{1}{2},\frac{1}{2})} 
\end{equation} 
satisfies $|\cG(u+\tau)| = |\cG(u+1)| = |\cG(u)|$, 
and if $\Re(u) = 0$ or if $\Re(u) = \frac{1}{2}$ then $|\cG(u)| = 1$. 
It has simple zeros at points $u \equiv M \pmod{\langle 1, \tau \rangle}$, 
simple poles at $u \equiv -\overline{M} \pmod{\langle 1, \tau \rangle}$, 
and no other zeros or poles.  

\smallskip
Now consider a set $E = [a,b] \cup [c,d] \subset \RR$, where $a < b < c < d$.  
We will give a (multivalued, periodic) conformal mapping of $\CC \backslash E$ onto a 
vertical strip, which will enable us to express $G(z,\zeta;E)$ in terms of the function
$\cG(u)$ in (\ref{FFFunc}).  
We follow Akhiezer (\cite{Akh}) and Falliero and Sebbar (\cite{Fall}, \cite{FS}), 
\index{Falliero, Th\'ere\`ese} 
\index{Sebbar, Ahmad} 
but obtain a different expression for the capacity. 
\index{capacity}  

First, put 
\begin{equation}
w \ = \ T(z) \ = \ \sqrt{\frac{(z-a)}{(z-b)} \frac{(d-b)}{(d-a)}} \ . \label{FTFunc} 
\end{equation} 
where $\sqrt{z}$ is positive for $z > 0$ 
and is slit along the negative real axis.  
$T(z)$ maps $\CC \backslash E$ conformally onto the right half-plane 
with the segment $[1,1/k]$ removed, where
\begin{equation} \label{FkDef} 
k \ = \ \frac{1}{T(c)} \ = \ \sqrt{\frac{(c-b)}{(c-a)} \frac{(d-a)}{(d-b)}} \ .
\end{equation}
$T(z)$ takes $a \mapsto 0$, $b \mapsto \infty$, $d \mapsto 1$, and $c \mapsto 1/k$.  Since
the linear fractional transformation $F(z) = (z-a)(d-b)/(z-b)(d-a)$ maps $\RR \cup {\infty}$ to itself and 
preserves the cyclic order of $a, b, c, d$, ones sees that $T(c) > 1$ and $0 < k < 1$. 
Note that $1/k^2$ is the crossratio $(a,b;c,d)$.  

Follow $T(z)$ with the elliptic integral 
\begin{equation} \label{FSDef} 
u \ = \ S(w) \ = \ \int_{0}^w \frac{dx}{\sqrt{(1-x^2)(1-k^2x^2)}} \ .
\end{equation}
Here $S(w)$ is the Schwarz-Christoffel map which sends the upper half-plane to the
\index{Schwarz-Christoffel map} 
rectangle with corners $\pm K, \pm K + iK^{\prime}$, where
\begin{eqnarray}
K & = & \int_{0}^1 \frac{dx}{\sqrt{(1-x^2)(1-k^2x^2)}}  \label{FKdef} \\
iK^{\prime} & = & \int_{1}^{1/k} \frac{dx}{\sqrt{(1-x^2)(1-k^2x^2)}} \label{FKprimeDef} 
\end{eqnarray} 
and $K, K^{\prime} > 0$.  It takes the imaginary axis to itself, and sends $0 \mapsto 0$,
$1 \mapsto K$, $1/k \mapsto K+iK^{\prime}$, and $\infty \mapsto iK^{\prime}$.  
By the Schwarz Reflection Principle, $S(w)$ extends to a multivalued holomorphic function
\index{Schwarz Reflection Principle}\index{multivalued holomorphic function} 
taking  $\{\Re(w) > 0\} \backslash [1,1/k]$ 
to the vertical strip $0 < \Re(u) < K$, with period $2iK^{\prime}$.  
The inverse function to $S(w)$ is the Jacobian elliptic function $w = \sn(u,k)$
\index{Jacobian elliptic function}    
(see \cite{WW}, \S 22, and \cite{Neh}, \S VI.3). 

Now let $\tau = iK^{\prime}/K$.  
Fix $\zeta \notin E$; put $u = S(T(z))$, $M = M(\zeta) = S(T(\zeta))$.  
Scaling $u \mapsto v = u/(2K)$ takes $0 < \Re(u) < K$ to the strip $0 < \Re(v) < 1/2$, 
with $2iK^{\prime} \mapsto \tau$.  We claim that 
\begin{equation} \label{FGreen} 
G(z,\zeta;E) \ = \ -\log\left(\left|\frac{\theta(\frac{u-M}{2K},\tau;\frac{1}{2},\frac{1}{2})}
                {\theta(\frac{u+\overline{M}}{2K},\tau;\frac{1}{2},\frac{1}{2})}\right|\right)
\end{equation} 
\index{Green's function!examples!archimedean}
Indeed, by our discussion of theta-functions,\index{theta-functions!classical}
the function on the right has the properties characterizing $G(z,\zeta;E)$:  
it is well-defined and continuous, 
vanishes on $E$, is harmonic on $\PP^1(\CC) \backslash (E \cup {\zeta})$, 
and has a positive logarithmic pole as $z \rightarrow \zeta$.   
This formula is one given by Falliero and Sebbar (\cite{Fall}; \cite{FS}, p.416). 
\index{Falliero, Th\'ere\`ese}   
\index{Sebbar, Ahmad}  

Numerically,\index{numerical!computations} $K$ and $K^{\prime}$ can be 
found using the hypergeometric function 
\begin{equation*}
F(a,b,c;z) \ = \ 1 + \frac{a \cdot b}{1 \cdot c} z 
     + \frac{a(a+1) \cdot b(b+1)}{1 \cdot 2 \cdot c(c+1)}z^2 + \cdots \ 
\end{equation*}
with $K = \frac{1}{2} \pi F(\frac{1}{2},\frac{1}{2},1;k^2)$, 
$K^{\prime} =  \frac{1}{2} \pi F(\frac{1}{2},\frac{1}{2},1;1-k^2)$ 
(see \cite{WW}, pp.499, 501); then $\tau = i K^{\prime}/K$. 
Another way to determine $\tau$, $K$ and $K^{\prime}$ 
is by first solving for $q = e^{i \pi \tau}$ using the relation
\begin{equation} \label{FkF2}
\frac{(c-b)}{(c-a)} \frac{(d-a)}{(d-b)} = 
k^2 = \frac{\theta(0,\tau,\frac{1}{2},0)^4}{\theta(0,\tau,0,0)^4} = 
\frac{16(q^{1/4} + q^{9/4} + q^{25/4} + \cdots)^4}{(1 + 2q^4 + 2q^9 + \cdots)^4} 
\end{equation}
and then using the formulas 
\begin{equation} \label{FKForm} 
K = \frac{1}{2} \pi \theta(0,\tau,0,0)^2\ , 
\qquad K^{\prime} = -i \tau K \ . 
\end{equation} 
Finally, $M$ can be determined by solving
\begin{equation} \label{FMFind} 
T(\zeta) \ = \ \sn(M,k) 
   \ = \ \frac{1}{k} \frac{\theta(M/2K,\tau;\frac{1}{2},\frac{1}{2})}
                          {\theta(M/2K,\tau;0,\frac{1}{2})} \ .
\end{equation}
(See \cite{WW}, pp.492, 501.)

We now determine the capacity of $E$.  If $\zeta = \infty$,
\index{capacity} 
put $\hz = 1/z$;  otherwise put $\hz = z - \zeta$.  Then as $\hz \rightarrow 0$,
we have $z \rightarrow \zeta$, $w \rightarrow T(\zeta)$, and $u \rightarrow M$.
Using (\ref{FGreen}), it follows that 
\index{Robin constant!examples!archimedean}
\begin{eqnarray}
V_{\zeta}(E) & = & \lim_{\hz \rightarrow 0} G(z,\zeta;E) + \log(|\hz|) \notag \\
& = & -\log\left(\left|\frac{\frac{d}{du}\theta(0,\tau;\frac{1}{2},\frac{1}{2})}
   {\theta(\frac{M+\overline{M}}{2K},\tau;\frac{1}{2},\frac{1}{2})} \cdot \frac{1}{2K} 
                \cdot \frac{dw}{du}(T(\zeta)) \cdot \frac{dT}{d\hz}(0)\right|\right)   \label{F2SegV} 
\end{eqnarray} 
The last two terms can be computed in terms of $\zeta$ and $a,b,c,d$; 
the expression can then be simplified using the Jacobi identity 
\begin{equation} \label{FJacId} 
\frac{d}{du}\theta(0,\tau;\frac{1}{2},\frac{1}{2}) 
= \pi \theta(0,\tau;0,0)\theta(0,\tau;\frac{1}{2},0) \theta(0,\tau;0,\frac{1}{2}) 
\end{equation}
(see \cite{WW}, p.470), together with (\ref{FKForm}) and (\ref{FkF2}).  
If $\zeta = \infty$ one obtains
\begin{equation} \label{F2SegV1}
\gamma_{\infty}(E) = e^{-V_{\infty}(E)} = 
  \frac{\sqrt[4]{(c-a)(c-b)(d-a)(d-b)}}
     {2 \left| \frac{\theta(\Re(M(\infty))/K,\tau;\frac{1}{2},\frac{1}{2})}
               {\theta(0,\tau;0,\frac{1}{2})} \right|} 
              \ ;
\end{equation}
if $\zeta \in \CC \backslash E$, then 
\begin{equation} \label{F2SegV2}
\gamma_{\zeta}(E) = 
 \frac{\sqrt[4]{(c-a)(c-b)(d-a)(d-b)}}
      {2 \left| \frac{\theta(\Re(M(\zeta))/K,\tau;\frac{1}{2},\frac{1}{2})}
               {\theta(0,\tau;0,\frac{1}{2})} \right| \cdot 
                      |(\zeta-a)(\zeta-b)(\zeta-c)(\zeta-d)|^{1/2} } \ .
\end{equation}  
Numerical examples\index{numerical!examples} confirm the compatibility of (\ref{F2SegV1}) with (\ref{F2Seg2}).  
However the formula of Akhiezer reported in (\cite{FS}, p.422) seems to be incorrect.
\index{Akhiezer, Naum}   

\smallskip
{\bf Three Segments.}\index{examples!archimedean!three segments} 
When $E = [a_1,b_1] \cup [a_2,b_2] \cup [a_3,b_3] \subset \RR$ and $\zeta = \infty$,  
Th\'ere\`se Falliero has given formulas for the Green's function 
\index{Falliero, Th\'ere\`ese} 
and capacity of $E$ using theta-functions of genus $2$;  for these, we refer the reader to 
\index{Green's function!examples!archimedean}\index{theta-functions!of genus two} 
\index{capacity} 
Falliero (\cite{Fall}) and Falliero-Sebbar  (\cite{FS}). 
\index{Falliero, Th\'ere\`ese}  
\index{Sebbar, Ahmad} 

\medskip  

{\bf Multiple segments.}\index{examples!archimedean!multiple segments}
 
When $E = [a_1,b_1] \cup [a_2,b_2] \cup \cdots \cup [a_n,b_n] \subset \RR$ 
with $a_1 < b_1 < a_2 < b_2 < \cdots < a_n < b_n$ and $n$ arbitrary, 
Harold Widom (\cite{Wid}, pp.224ff) has given formulas for  $G(z,\zeta;E)$ and $V_{\zeta}(E)$ 
\index{Robin constant!examples!archimedean}
\index{Widom, Harold} 
which we recall below.  

Let $q(z) =  \prod_{j=1}^n (z-a_j)(z-b_j)$, and let $q(z)^{1/2}$  
be the branch of $\sqrt{\prod_{j=1}^n (z-a_j)(z-b_j)}$ on $\CC \backslash E$  
which is positive as $z \rightarrow \infty$ along the real axis.  
This branch is well-defined throughout $\CC \backslash E$ and 
positive on $\RR \backslash E$.
For each $x \in E$, the limiting values of $q(z)^{1/2}$
 as $z \rightarrow x+i0^+$ and $z \rightarrow x + i0^-$ are 
pure imaginary, and are negatives of each other.  For convenience, we will extend $q(z)^{1/2}$ to $E$ 
by defining $q(x)^{1/2} = q(x+i0^-)^{1/2}$ when $x \in E$.  Thus, $q(z)^{1/2}$ is pure imaginary on $E$.

\smallskip
First take $\zeta = \infty$.  Fix a point $z_0 \in E$, 
and let $h(z) = h_0 + h_1 z + \cdots + h_{n-1} z^{n-1}\in \RR(z)$ be a polynomial 
of degree $\le n-1$ with real coefficients.  
Consider the multiple-valued function 
\begin{equation*}
G_h(z) \ = \ \int_{z_0}^z h(w)/q(w)^{1/2} \, dw 
\end{equation*} 
on $\CC$, where the integral is taken over any path from $z_0$ to $z$ which is disjoint from $E$ 
except for one or both of its endpoints.   
Since $G_h(z)$ has pure imaginary periods\index{pure imaginary periods} around $\infty$ 
and around each component $[a_j,b_j]$ of $E$, 
the function $\Re(G_h(z))$ is well-defined and continuous,  
and constant on each component of $E$.  
Since $G(z)$ has a holomorphic branch
in a neighborhood of each point $w \in \CC \backslash E$, 
$\Re(G_h(z))$ is harmonic in $\CC \backslash E$.

Clearly $\Re(G_h(z)) \equiv 0$ on the component of $E$ containing $z_0$.  If 
\begin{equation} \label{WidomGapEq}
\int_{b_j}^{a_{j+1}} h(x)/q(x)^{1/2} \, dx \ = \ 0
\end{equation}  
for each `gap' $(b_j,a_{j+1})$, then $\Re(G_h(z)) \equiv 0$ on $E$. 
If in addition $h(z)$ is monic, there is a number $V \in \RR$ 
such that $\Re(G_h(z))$ is asymptotic to $\log(|z|) + V$ 
as $z \rightarrow \infty$.   
In this setting, the characterization of Green's functions shows that
\index{Green's function!archimedean!characterization of}
\begin{equation} \label{WidomGreen} 
G(z,\infty;E) \ = \ \Re(G_h(z)) \ .
\end{equation} 
\index{Widom, Harold}
We will now show that such an $h(z)$ exists.

Following Widom, put $A_{jk} = \int_{b_j}^{a_{j+1}} x^k/q(x)\, dx$
for  $j = 1, \ldots, n-1$, $k = 0, \ldots, n-1$.  Note that $1/q(z)^{1/2}$ has singularities at $b_j$
and $a_{j+1}$ of order $z^{-1/2}$, so each $A_{jk}$ is finite and belongs to $\RR$.
We claim that there is a unique solution $h_0, h_1, \ldots, h_{n-1}$ to the 
system of linear equations  
\begin{equation} \label{WidomSystem} 
\left\{ \begin{array}{rl} \sum_{k=0}^{n-1} A_{jk} h_k  \ = \ 0 \ \ & 
                              \text{for $j = 1, \ldots, n-1$\ ,} \\
                          h_{n-1} \ = \ 1 \ . & \text{  }
        \end{array} \right. 
\end{equation}   
If $h_0, \ldots, h_{n-1}$ satisfy (\ref{WidomSystem}), and $h(z)$ is the corresponding polynomial,  
then $h(z)$ is monic and the conditions (\ref{WidomGapEq}) hold.
Thus $G(z,\infty;E) = \Re(G_h(z))$.  Solving the system (\ref{WidomSystem}) 
is called the `Jacobi Inversion Problem'.\index{Jacobi Inversion problem}  

To see that (\ref{WidomSystem}) has a unique solution, 
it suffices to show that the  $n \times n$ matrix associated to the system has rank $n$, 
or equivalently, that $h_0 = \cdots = h_{n-1} = 0$ is the only solution to the corresponding homogeneous
system.  Let $h_0, \ldots, h_{n-1} \in \RR$ be any solution to the homogeneous system, and let $h(z)$
be the corresponding polynomial.  Then $G_h(z)$ is harmonic on $\CC \backslash E$, vanishes on $E$, 
and remains bounded as  $z \rightarrow \infty$ since $h_{n-1}=0$, 
so it extends to a function on $\PP^1(\CC) \backslash E$ harmonic at $\infty$.
By the maximum principle for harmonic functions, $G_h(z) \equiv 0$.  
Restricting $G_h(z)$ to $\RR$, differentiating, 
and using the Fundamental Theorem of Calculus,\index{Fundamental Theorem of Calculus}
 we see that 
$h(z)/q(z)^{1/2} \equiv 0$ on $\RR \backslash E$.  Since $q(z)$ is nonzero except at the endpoints of $E$, 
it follows that $h(z) \equiv 0$, and hence that $h_0 = \cdots = h_{n-1} = 0$.  

We remark that $h(z)$ has one zero in each gap $(b_j,a_{j+1})$, and no other zeros. 
Indeed $G(x,\infty;E)$ vanishes at $b_j$ and $a_{j+1}$, 
and is real-valued and differentiable on $(b_j,a_{j+1})$, 
so by Rolle's Theorem there is a point 
\index{Rolle's theorem}
$x_j^* \in (b_j,a_{j+1})$ where $G^{\prime}(x_j^*,\infty;E) = 0$.  
The argument above shows that $h(x_j^*) = 0$.
Since $h(z)$ has degree $n-1$, it has a unique zero in each gap, and these are its only zeros. 
Thus, $h(z) = \prod_{j=1}^{n-1} (z-x_j^*)$, and $h(z)$ has constant sign on each component of $E$.

\smallskip
Now consider the case when $\zeta \in \CC \backslash E$.  Again, fix $z_0 \in E$.
We claim that for a suitable polynomial
$h(z) = h_0 + h_1 z + \cdots h_{n-1} z^{n-1} \in \CC[z]$, we have 
\begin{equation} \label{WidomGreenZeta}
G(z,\zeta;E) \ = \ \Re(G_h(z)) \ , 
\end{equation}
where now 
\begin{equation} \label{WidomGhZeta}
G_h(z) \ = \ \int_{z_0}^z \frac{h(w)}{q(w)^{1/2}(w-\zeta)} \, dw \ .
\end{equation} 
\index{Green's function!examples!archimedean}
Here $h(z)$ must be chosen so that the following properties are satisfied: 
\begin{enumerate}
\item The periods of $G_h(z)$ are pure imaginary,\index{pure imaginary periods|ii} so $\Re(G_h(z))$ is well-defined. 
\item The value of $\Re(G_h(z))$ is $0$ on each segment $[a_j,b_j]$.
\item $\Re(G_h(z))$ has a singularity of type $-\log(|z-\zeta|)$ at $\zeta$.
\end{enumerate} 
Let $\Gamma_j$ be a loop about $[a_j,b_j]$, traversed counterclockwise.  
Using Cauchy's theorem,\index{Cauchy's theorem} one sees that
\begin{equation*} 
\int_{\Gamma_j} \frac{h(w)}{q(w)^{1/2}(w-\zeta)} \, dw 
\ = \ 2 \int_{a_j}^{b_j} \frac{h(x)}{q(x)^{1/2}(x-\zeta)} \, dx \ .
\end{equation*} 
Thus for the periods of $G_h(z)$ about the intervals $[a_j,b_j]$ to be pure imaginary,\index{pure imaginary periods} we need
\begin{equation} \label{abPeriods0} 
\Re\big(\int_{a_j}^{b_j} \frac{h(x)}{q(x)^{1/2}(x-\zeta)} \, dx\big) \ = \ 0 \qquad 
\text{for $j = 1, \ldots, n$\ .}
\end{equation}
Let $\varepsilon > 0$ be small enough that the circle $C(\zeta,\varepsilon)$ is disjoint from $E$.
Since the differential $h(w)\, dw/\big(q(q)^{1/2}(w-\zeta)\big)$ is holomorphic at $\infty$, 
applying Cauchy's theorem\index{Cauchy's theorem} on the domain $\PP^1(\CC) \backslash (E \cup \{\zeta\})$  
we obtain    
\begin{equation*} 
\sum_{j=1}^n \int_{\Gamma_j} \frac{h(w)}{q(w)^{1/2}(w-\zeta)} \, dw 
+ \int_{C(\zeta,\varepsilon)} \frac{h(w)}{q(w)^{1/2}(w-\zeta)} \, dw \ = \ 0 \ .
\end{equation*}
Hence if the conditions (\ref{abPeriods0}) hold, the period of $G_h(z)$ 
about $\zeta$ (which is $2 \pi i h(\zeta)/q(\zeta)^{1/2})$ is pure imaginary as well.\index{pure imaginary periods}  

Under the conditions (\ref{abPeriods0}), $\Re(G_h(z))$ is well-defined, 
harmonic in $\PP^1(\CC) \backslash (E \cup \{\zeta\})$, and constant on each segment $[a_j,b_j]$.
Clearly its value on the segment containing $z_0$ is $0$.  For it to be identically $0$ on $E$, 
we need 
\begin{equation} \label{GapPeriod0}  
\Re\big( \int_{b_j}^{a_{j+1}} \frac{h(x)}{q(x)^{1/2}(x-\zeta)} \, dx \big) \ = \ 0 
\qquad \text{for $j = 1, \ldots, n-1$ \ .} 
\end{equation} 

Finally, for $\Re(G_h(z))$ to have a singularity of type $-\log(|z-\zeta|)$ at $\zeta$, we need 
\begin{equation*} 
-1 \ = \  \Res_{w = \zeta}\big(\frac{h(w)}{q(w)^{1/2}(w-\zeta)}\big)  
\ = \ \frac{h(\zeta)}{q(\zeta)^{1/2}}  \ . 
\end{equation*}
Since the period of $G_h(z)$ 
about $\zeta$ is imaginary\index{pure imaginary periods} we automatically have $\Im(h(\zeta)/q(\zeta)^{1/2}) = 0$, 
and it is enough to require 
\begin{equation} \label{ZetaPeriod0} 
-1 \ = \ \Re\big(\frac{h(\zeta)}{q(\zeta)^{1/2}}\big)  \ . 
\end{equation}

Writing $h_k = c_k + d_k i$ for $k = 0, \ldots, n-1$, with $c_k, d_k \in \RR$, 
the conditions (\ref{abPeriods0}), (\ref{GapPeriod0}) and (\ref{ZetaPeriod0}) represent 
a system of $2n$ linear equations with real coefficients in $2n$ real unknowns.  To show that it has a unique
solution, it is enough to show that the only solution to the corresponding homogeneous system is 
the trivial one.  

Suppose that $h(z)$ arises from a solution to the homogeneous system.  Then $\Re(G_h(z))$ is harmonic
in $\PP^1(\CC) \backslash (E \cup \{\zeta\})$ and extends to a function harmonic at $\zeta$, 
with boundary values $0$ on $E$.  By the Maximum Principle, $\Re(G_h(z)) \equiv 0$.  
\index{Maximum principle!for harmonic functions}
Differentiating, and using the Fundamental Theorem of Calculus\index{Fundamental Theorem of Calculus} 
on horizontal segments, we see that 
\begin{equation*} 
\Re\big(\frac{h(z)}{q(z)^{1/2}(z-\zeta)}\big) 
\ = \ \frac{\partial}{\partial x}\big(\Re(G_h(z))\big) \ \equiv \ 0 
\end{equation*} 
on $\CC \backslash (E \cup \{\zeta\})$. If real part of an analytic function is identically $0$,
that function is constant, so we must have 
\begin{equation*}
\frac{h(z)}{q(z)^{1/2}(z-\zeta)} \ = \ C 
\end{equation*}
for some purely imaginary constant $C$.  
However, $q(z)^{1/2}(z-\zeta)$ is not a polynomial, so this can hold only if $C = 0$. Thus $h(z) \equiv 0$,
which means that $c_1 = d_1 = \cdots = c_n = d_n = 0$. 

\smallskip
When $\zeta = \infty$, choosing $z_0 \ne 0$ and 
noting that $\log(|z|) = \Re\big(\int_{z_0}^{\infty} 1/w \, dw \big) - \log(|z_0|)$, 
Widom gives a formula for the Robin constant equivalent to 
\index{Robin constant!examples!archimedean}
\index{Widom, Harold} 
\begin{eqnarray*}
V_{\infty}(E) & = & \lim_{z \rightarrow \infty} \big( G(z,\infty;E) - \log(|z|)\big) \\
& = & \Re\big(\int_{z_0}^{\infty} \frac{h(w)}{q(w)^{1/2}} - \frac{1}{w} \, dw \big) \ + \ \log(|z_0|) \ .
\end{eqnarray*}  
Similarly, when $\zeta \in \CC \backslash E$, 
\begin{eqnarray*}
V_{\zeta}(E) & = & \lim_{z \rightarrow \zeta} \big( G(z,\zeta;E) + \log(|z-\zeta|)\big) \\
 & = & \Re\big(\int_{z_0}^{\infty} \frac{h(w)}{q(w)^{1/2}(w-\zeta)} + \frac{1}{w-\zeta} \, dw \big) 
\ - \ \log(|z_0-\zeta|) \ .
\end{eqnarray*} 

\smallskip
When $\zeta = \infty$, there is a more illuminating formula for $V_{\infty}(E)$.
\index{Robin constant!examples!archimedean}   
Put $c = (a_1 + b_n)/2$ and $r = (b_n-a_1)/4$, so 
$E \subset [a_1,b_n] = [c-2r,c+2r]$.  We claim that 
\begin{eqnarray}  
V_{\infty}(E) & = & -\log(\frac{b_n-a_1}{4}) \ + \ \sum_{j=1}^{n-1} \int_{b_j}^{a_{j+1}} G(x,\infty;E) \, 
\frac{1}{\pi} \frac{dx}{\sqrt{4r^2 - (x-c)^2}} \ , \label{VWidom} \\
\gamma_{\infty}(E) & = & e^{-V_{\infty}(E)} \ = \ \frac{b_n-a_1}{4} \cdot 
\prod_{j=1}^{n-1} e^{-\int_{b_j}^{a_{j+1}} G(x,\infty;E) \, \frac{1}{\pi} \frac{dx}{\sqrt{4r^2 - (x-c)^2}}} \ .
\notag 
\end{eqnarray}  
Readers familiar with capacity theory will recognize $dx/(\pi \sqrt{4r^2 - (x-c)^2})$ 
\index{capacity} 
as the equilibrium distribution of $[a_1,b_n]$ relative to $\infty$.

To derive (\ref{VWidom}), assume for simplicity that $c = 0$, so $[a_1,b_n] = [-2r,2r]$;  
this can always be arranged by a translation.  Note that since 
$G(z,\infty;[-2r,2r]) \sim \log(|z|) + V_{\infty}([-2r,2r])$ as $z \rightarrow \infty$, 
and $V_{\infty}([-2r,2r]) = V_{\infty}([a_1,b_n]) = -\log((b_n-a_1)/4)$, we have 
\index{Robin constant!examples!archimedean}
\begin{eqnarray} 
V_{\infty}(E) & := & \lim_{z \rightarrow \infty} G(z,\infty;E) - \log(|z|) \notag \\
       & = & \lim_{z \rightarrow \infty} \big( G(z,\infty;E) - G(z,\infty;[-2r,2r]) \big) 
                  - \log(\frac{b_n-a_1}{4}) \ . \label{WF1}
\end{eqnarray}
The function $g(z) := G(z,\infty;E) - G(z,\infty;[-2r,2r])$ is harmonic in 
$\CC \backslash [-2r,2r]$ and bounded as $z \rightarrow \infty$;
hence it extends to a function harmonic in $\PP^1(\CC) \backslash [-2r,2r]$, with
\begin{equation} \label{WF2}
g(\infty) \ = \ \lim_{z \rightarrow \infty} \big( G(z,\infty;E) - G(z,\infty;[-2r,2r]) \big) \ .
\end{equation}   
Let the Joukowski map $z = J_r(w) = w + r^2/w$ be as in (\ref{rJoukowski}). 
\index{Joukowski map}  
For each $R > r$, parametrize the ellipse $\cE(R+r^2/R,R-r^2/R)$ by 
$z = J_r(R \cos(\theta) + i R \sin(\theta))$ 
as in (\ref{REllipseParam}).
Let $\cD_R = \PP^1(\CC) \backslash D(0,R)$, 
and let $\cE_R$ be the connected component of $\PP^1(\CC) \backslash \cE(R+r^2/R,R-r^2/R)$
containing $\infty$. 
The map $J_r(w)$ gives a conformal equivalence from $\cD_R$ to $\cE_R$, and takes $\infty$ to $\infty$.  
Thus $H(w) := g(J_r(w))$ is harmonic in $\cD_R$, 
and $H(\infty) \ = \ g(\infty)$. 
By the mean value theorem for harmonic functions,
\begin{eqnarray*}
g(\infty) & = & H(\infty) \ = \ 
\frac{1}{2\pi} \int_0^{2 \pi} H(R \cos(\theta) + i R \sin(\theta)) \, d \theta \\
& = &  \frac{1}{2\pi} \int_0^{2 \pi}
g((R+\frac{r^2}{R}) \cos(\theta) + i (R -\frac{r^2}{R})\sin(\theta)) \, d \theta \ . 
\end{eqnarray*}
Since $\cE(R+r^2/R,R-r^2/R)$ is the level curve $\log(R/r)$ 
for $G(z,\infty;[-2r,2r])$ (see (\ref{EllipseLevelCurve})), it follows that  
\begin{equation*} 
g(\infty) \ = \  \frac{1}{2\pi} \int_0^{2 \pi} 
G((R+\frac{r^2}{R}) \cos(\theta) + i (R -\frac{r^2}{R})\sin(\theta)),\infty;E) \, d \theta 
\ - \ \log(R/r) \ .
\end{equation*}
Since $G(z,\infty;E) = G_h(z)$ is continuous on $\CC$,
letting $R \rightarrow r$ we see that  
\begin{equation*}
g(\infty) \ = \  \frac{1}{2\pi} \int_0^{2 \pi} G(2r \cos(\theta),\infty;E) \, d \theta \ .
\end{equation*} 
Finally, making the change of variables $x = 2r \cos(\theta)$ yields  
\begin{equation} \label{WF5}
g(\infty) \ = \ \int_{a_1}^{b_n} G(x,\infty;E) \, \frac{1}{\pi} \frac{dx}{\sqrt{4r^2 - x^2}}
\ = \ \sum_{j=1}^{n-1} \int_{b_j}^{a_{j+1}} G(x,\infty;E) \, \frac{1}{\pi} \frac{dx}{\sqrt{4r^2 - x^2}} \ .
\end{equation} 
Combining (\ref{WF1}), (\ref{WF2}) and (\ref{WF5}) gives (\ref{VWidom}).

\smallskip
{\bf The Real Projective Line.}\index{examples!archimedean!real projective line}  
If $E = \PP^1(\RR)$, 
the components of its complement in $\PP^1(\CC)$ are the upper and lower half-planes.  
Fix $\zeta \notin E$.
Using the characterization of the Green's function,
\index{Green's function!archimedean!characterization of}
it is easy to check that if $z$ and $\zeta$ belong to the same component of 
$\PP^1(\CC) \backslash \PP^1(\RR)$, then 
\begin{equation} \label{FGHP} 
G(z,\zeta;E) \ = \ -\log\left(\left|\frac{z-\zeta}{z-\overline{\zeta}}\right|\right) \ .
\end{equation}
\index{Green's function!examples!archimedean}
If $z$ and $\zeta$ are not in the same component, then $G(z,\zeta;E) = 0$.  
Taking $g_{\zeta}(z) = z-\zeta$ and using (\ref{FGHP}) we obtain
\index{Green's function!examples!nonarchimedean}   
\begin{eqnarray}
V_{\zeta}(E) & = & \lim_{z \rightarrow \zeta} 
-\log\left(\left|\frac{z-\zeta}{z-\overline{\zeta}}\right|\right) + \log(|z-\zeta|) 
\ = \ \log(2|\Im(\zeta)|) \ , \label{FVHP1}     \\
\gamma_{\zeta}(E) & = & \frac{1}{2|\Im(\zeta)|} \ . \label{FVHP2}
\end{eqnarray} 

\smallskip
{\bf The Disc with Opposite Radial Arms.}\index{examples!archimedean!disc with arms}
Take $L_1, L_2 \ge 0$, and let $E$ be the union of $D(0,R)$  
with the segment $[-L_1-R,R+L_2]$;  
thus $E$ is a disc with opposite radial arms of length $L_1$, $L_2$.   We claim that 
\begin{equation} \label{DiscOppRadial} 
\gamma_{\infty}(E) \ = \ \frac{1}{4} \left(2R + 
\frac{R^2 + R L_1 + L_1^2}{R+L_1} + \frac{R^2 + R L_2 + L_2^2}{R+L_2}\right) \ .
\end{equation}

To see this, first take $R = 1$. Put $a_1 = 1+L_1$, $a_2 = 1+L_2$;
then $E = D(0,1) \cup [-a_1,a_2]$.   Let $w = \varphi(z) = (z-1)^2/z$.  Then
$\varphi$ is the composite of the maps $z \rightarrow t = -1/(z+1)$, 
$t \rightarrow u = t+1/2$, $u \rightarrow v = u^2$,
and $v \rightarrow w = -1/(v-1/4)$.   
Using standard properties of conformal maps one sees that $\varphi(z)$ 
maps $\CC \backslash E$  conformally onto $\CC \backslash [A,B]$, where 
\begin{equation*} 
A = -\frac{(a_1+1)^2}{a_1} \ , \qquad B = \frac{(a_2-1)^2}{a_2} \ .
\end{equation*} 
Clearly $\varphi(\infty) = \infty$.  
Since $\lim_{z \rightarrow \infty} \log(|w|/|z|) = 0$, it follows that 
\begin{equation} \label{TwoRadial2} 
\gamma_{\infty}(E) \ = \ \gamma_{\infty}([A,B]) \ = \ \frac{B-A}{4}  
\ = \ \frac{(a_1 a_2 + 1)(a_1 + a_2)}{4 a_1 a_2} \ . 
\end{equation} 
In the general case, put $a_1 = 1 + L_1/R$, $a_2 = 1 + L_2/R$, 
and scale (\ref{TwoRadial2}) by $R$;  
after simplification, one obtains (\ref{DiscOppRadial}).   
The expression (\ref{TwoRadial2}) appears in (\cite{Kob}, p.82).   

For the set $E = D(0,1) \cup [-a_1,a_2]$ discussed above, 
and for $z, \zeta \notin E$, one has 
\begin{equation*}
G(z,\zeta;E) \ = \ G(\varphi(z),\varphi(\zeta);[A,B])
\end{equation*}
\index{Green's function!examples!archimedean}
where the Green's function of $[A,B]$ is given by (\ref{F1Seg2}).
This can be used to find $\gamma_{\zeta}(E)$ for any $\zeta \in \CC \backslash E$.

\smallskip
{\bf Two Concentric Circles.}\index{examples!archimedean!two concentric circles}
 Fix $r > 1$, and let $E$ be the union of
the circles $|z| = 1/r$ and $|z| = r$.  The complement of $E$ has three components.
If $z$ and $\zeta$ belong to different components, then $G(z,\zeta;E) = 0$.  
If they belong to the outer component, then $G(z,\zeta;E) = G(z,\zeta;D(0,r))$, 
while if they belong to the inner component then $G(z,\zeta;E) = G(1/z,1/\zeta;D(0,r))$. 

$G(z,\zeta;E)$ is also known when $z$ and $\zeta$ belong to the annular region between 
the circles.  Courant and Hilbert (\cite{C-H}, pp. 386--388) derive a formula for it
\index{Courant, Richard} 
\index{Hilbert, David}
using the Schwarz Reflection Principle:  define $q$ by $q^{1/2} = 1/r$, 
\index{Schwarz Reflection Principle}
and suppose $1/r < |z|, |\zeta| < r$.  Courant and Hilbert show that
\index{Courant, Richard} 
\index{Hilbert, David}
 $G(z,\zeta;E) = -\log(|f_{\zeta}(z)|)$, where
\index{Green's function!examples!archimedean}   
\begin{equation} \label{DFfunc}  
f_{\zeta}(z) = |z|^{-\log(|\zeta|)/\log(q)} \cdot 
       \frac{ q^{1/4} (\sqrt{\frac{z}{\zeta}}-\sqrt{\frac{\zeta}{z}}\ )  
                    \prod_{n=1}^{\infty} (1-q^{2n} \frac{z}{\zeta}) (1 - q^{2n} \frac{\zeta}{z}) }
     { \prod_{n=1}^{\infty} (1 - q^{2n-1} \overline{\zeta}z) (1 - q^{2n-1} 
                           \frac{1}{\overline{\zeta}z} ) }
\end{equation}  
Recalling the product expansions of the theta functions, they note that the second term 
is a quotient of two theta functions, leading to the following expression:  
writing $\tau = 2i \log(r)/\pi$,  $z = e^{2 \pi i u}$ and $\zeta = e^{2 \pi i \alpha}$,
then for $1/r < |z|, |\zeta| < r$, 
\begin{equation} \label{FGreenRing} 
G(z,\zeta;E) \ = \ -\frac{\log(|z|) \log(|\zeta|)}{2\log(r)} 
- \log\left( \left| \frac{\theta(u-\alpha,\tau;\frac{1}{2},\frac{1}{2})}
                    {\theta(u-\overline{\alpha},\tau;0,\frac{1}{2})} \right| \right) \ .
\end{equation}  
\index{Green's function!examples!archimedean}
Here we have corrected a minor error in Courant-Hilbert, who state (\ref{DFfunc}) 
\index{Courant, Richard} 
\index{Hilbert, David}
for positive real $\zeta$, and omit the conjugate on $\zeta$ in generalizing;  this changes their
$\theta(u+\alpha,\tau;0,\frac{1}{2})$ to $\theta(u-\overline{\alpha},\tau;0,\frac{1}{2})$.    
Using (\ref{FGreenRing}), we obtain
\index{Green's function!examples!archimedean}   
\begin{eqnarray} \label{FRingV} 
V_{\zeta}(E) & = & \lim_{z \rightarrow \zeta} G(z,\zeta;E) + \log(|z-\zeta|) \\
             & = & -\frac{(\log(|\zeta|))^2}{2\log(r)} 
                         + \log(|\theta(\alpha-\overline{\alpha},\tau;0,\frac{1}{2})|) \notag \\
             &  & \qquad \qquad  - \log(|\frac{d}{du}\theta(0,\tau;\frac{1}{2},\frac{1}{2})| 
                       + \log(|\frac{dz}{du}(\alpha)|) \notag \\
             & = & -\frac{(\log(|\zeta|))^2}{2\log(r)} + 
            \log\left(\left|\frac{2 \zeta \cdot \theta(\alpha-\overline{\alpha},\tau;0,\frac{1}{2})}
            {\theta(0,\tau;0,0)\theta(0,\tau;\frac{1}{2},0)\theta(0,\tau;0,\frac{1}{2})}
                 \right| \right) \notag
\end{eqnarray}     
where we have used Jacobi's identity\index{Jacobi identity} (\ref{FJacId}) to simplify
$\frac{d}{du}\theta(0,\tau;\frac{1}{2},\frac{1}{2})$.  When $\zeta = 1$, this becomes
\index{Robin constant!examples!archimedean}
\begin{equation} \label{FRingV1}
V_1(E) \ = \ \log\left(\frac{2}{|\theta(0,\tau;0,0)\theta(0,\tau;\frac{1}{2},0)|}\right) ,
\quad \gamma_1(E) = \frac{|\theta(0,\tau;0,0)\theta(0,\tau;\frac{1}{2},0)|}{2} \ .
\end{equation} 

\medskip
{\bf Sets arising in Polynomial Dynamics.} 

{\bf Julia Sets.}\index{examples!archimedean!Julia sets}
 Let $\varphi(x) = a_0 + a_1 x + \cdots + a_d x^d \in \CC[x]$  be 
a polynomial of degree $d \ge 2$.  
By definition, the {\em filled Julia set} $\cK_\varphi$ of $\varphi(x)$ 
\index{Julia set!filled}
is the set of all $z \in \CC$ 
whose forward orbit $z, \varphi(z), \varphi(\varphi(z)), \ldots$ under $\varphi$ remains bounded;
the {\em Julia set} is its boundary 
$\cJ_\varphi = \partial \cK_\varphi$.\index{boundary!of filled Julia set}\index{Julia set} 
Let $\varphi^{(n)} = \varphi \circ \varphi \circ \cdots \circ \varphi$ be the $n$-fold iterate. 
For each sufficiently large $R$, 
we have $D(0,R) \supset \varphi^{-1}(D(0,R)) \supset (\varphi^{(2)})^{-1}(D(0,R)) \supset \cdots \supset \cK_\varphi$,
and 
\begin{equation*}
\cK_\varphi \ = \ \bigcap_{n=1}^{\infty} (\varphi^{(n)})^{-1}(D(0,R)) \ .
\end{equation*} 
As in (\cite{SM}, p.147), for each $z \in \CC$ we have 
\begin{equation} \label{JuliaGreen}
G(z,\infty; \cJ_\varphi) \ = \ G(z,\infty; \cK_\varphi) 
\ = \ \lim_{n \rightarrow \infty} \frac{1}{d^n} \log^+(|\varphi^{(n)}(z)|) 
\end{equation}
(the `escape velocity' of $z$),\index{escape velocity} and 
\index{Robin constant!examples!archimedean}
\begin{equation} \label{JuliaCapacity} 
V_{\infty}(\cJ_\varphi) = V_{\infty}(\cK_\varphi)  = \frac{\log(|a_d|)}{d-1} \ , \qquad 
\gamma_{\infty}(\cJ_\varphi) = \gamma_{\infty}(\cK_\varphi) = |a_d|^{-1/(d-1)} \ .
\end{equation} 
\index{Green's function!examples!archimedean}

The proofs of (\ref{JuliaGreen}) and (\ref{JuliaCapacity}) are simple. 
It is easy to see that $\varphi^{(n)}(z)$ has degree $d^n$ and leading coefficient 
$a_0^{d^{n-1} + d^{n-2} + \cdots + d + 1}$.   
By the characterization of Green's functions it follows that  
\index{Green's function!archimedean!characterization of}
$G(z,\infty;(\varphi^{(n)})^{-1}(D(0,R))) = d^{-n} \log^+(|\varphi^{(n)}(z)|)$, and that 
\begin{equation*}
V_{\infty}((\varphi^{(n)})^{-1}(D(0,R))) \ = \ \frac{d^{n-1} + d^{n-2} + \cdots + d + 1}{d^n} \log(|a_0|) 
\ = \ \frac{1-1/d^{n}}{d-1} \log(|a_0|) \ .
\end{equation*}
The Green's functions $G(z,\infty;(\varphi^{(n)})^{-1}(D(0,R))$ decrease monotonically to $G(z,\infty;\cK_\varphi)$, 
\index{Green's function!examples!archimedean}
\index{Green's function!monotonic}
and the convergence is uniform outside any neighborhood of $\cK_\varphi$, 
so (\ref{JuliaGreen}) and (\ref{JuliaCapacity}) follow.   

\smallskip
{\bf The Mandelbrot Set.}\index{examples!archimedean!Mandelbrot set} 
Each quadratic polynomial $\varphi_c(x) = x^2 + c \in \CC[x]$ has $0$ as a critical point.  
The {\em Mandelbrot set $\cM$} is the set of all $c \in \CC$ for which the forward orbit 
\index{Mandelbrot set}  
$0, \varphi_c(0), \varphi_c^{(2)}(0), \ldots$ remains bounded; equivalently, 
$\cM$ is the set of all $c \in \cC$ for which $0$ belongs to the filled Julia set of $\varphi_c(x)$.
\index{Julia set!filled} 
It is easy to see that $\varphi_c(0) = c$, $\varphi_c^{(2)}(0) = c^2 - c$, and $\varphi_c^{(3)}(0) = c^4 - 2c^3 + c^2 -c$;  
in general $P_n(c) := \varphi_c^{(n+1)}(0)$ is a monic polynomial of degree $2^n$  
in $\ZZ[c]$. 
It can be shown that 
$D(0,2) \supset P_1^{-1}(D(0,2)) \supset P_2^{-1}(D(0,2)) \supset \cdots \supset \cM$,
and that 
\begin{equation} \label{MandelDef}
\cM \ = \ \bigcap_{n=1}^{\infty} P_n^{-1}(D(0,2)) \ ;
\end{equation}
see (\cite{SM}, p.158).  
By arguments like those for Julia sets,  
\index{Julia set} 
for each $c \in \CC$ we have 
\begin{equation} \label{MandelbrotGreen}
G(c,\infty;\cM) \ = \ \lim_{n \rightarrow \infty} \frac{1}{2^n} \log^+(|P_n(c)|) 
\ = \ \lim_{n \rightarrow \infty} \frac{1}{2^n} \log^+(\varphi_c^{(n+1)}(0)) \ ,
\end{equation} 
\index{Green's function!examples!archimedean}
and 
\index{Robin constant!examples!archimedean}
\begin{equation} \label{MandelbrotCapacity}
V_{\infty}(\cM) = 0 \ , \qquad \gamma_{\infty}(\cM) = 1 \ .
\end{equation}

\section{ Local capacities and Green's functions of Nonarchimedean Sets} \label{NonarchLocalExamples} 
 
In this section, $K_v$ will be a nonarchimedean local field.  
Identify $\PP^1(\CC_v)$ with $\CC_v \cup \{\infty\}$.
There are two methods of determining the Green's function 
\index{Green's function!computing nonarchimedean}
for sets $E_v \subset \PP^1(\CC_v)$:  
by using the pullback formula for Green's functions,  for noncompact sets; 
\index{Green's function!pullback formula for}
or by guessing the equilibrium distribution based on symmetry,
for compact sets. We are aided by the fact that
the capacity is monotonic under containment of sets.  
\index{capacity}  

The pullback formula for Green's functions is as follows. 
\index{Green's function!pullback formula for|ii} 
Let $\cC_1, \cC_2/\CC_v$ be smooth, complete curves,  
and let $f : \cC_1 \rightarrow \cC_2$ be a nonconstant rational map.  
Suppose $E_v \subset \cC_2(\CC_v)$ is an algebraically capacitable set of positive capacity.
\index{algebraically capacitable}
\index{capacity $> 0$}  
Fix $\zeta \in \cC_2(\CC_v) \backslash E_v$ and write $f^{*}((\zeta)) = \sum_{j=1}^m m_k (\xi_j)$
as a divisor. Then for each $z \in \cC_1(\CC_v)$,  
\begin{equation} \label{GreenPullbackF}
G(f(z),\zeta;E_v) \ = \ \sum_{j=1}^m m_k G(z,\xi_j;f^{-1}(E_v) \ .
\end{equation}
This holds for both nonarchimedean and archimedean sets (see \cite{RR1}, Theorems 3.2.9, 4.4.19).

\smallskip 
{\bf The Closed Disc.}\index{examples!nonarchimedean!closed disc}   
If $E_v = D(a,R) = \{z \in \CC_v : |z-a|_v \le R\}$ then
\begin{equation} \label{FGNAD} 
G(z,\zeta;E_v) \ = \ 
\left\{ \begin{array}{ll} 
\log_v^+(|z-a|_v/R) & \text{if $\zeta = \infty$ \ ,} \\
\log_v^+\left( \frac{|\zeta-a|_v}{R} \cdot \left|\frac{z-a}{z-\zeta}\right|_v \right)
               & \text{if $\zeta \in \CC_v \backslash D(a,R)$ \ .} 
    \end{array} \right. 
\end{equation}   
\index{Green's function!examples!nonarchimedean}
The first formula is essentially the definition of the Green's function  
as given by Cantor (\cite{Can3});  
\index{Cantor, David} 
the second follows from the first, by applying the pullback formula (\ref{GreenPullbackF}) 
to the map $f(z) = (z-a)/(z-\zeta)$ which takes $D(a,R)$ to $D(0,R/|\zeta-a|_v)$ 
and takes $\zeta$ to $\infty$. 

Taking $g_{\infty}(z) = 1/z$, 
and $g_{\zeta}(z) = z-\zeta$ if $\zeta \in \CC_v \backslash D(a,R)$,
we have 
\index{Robin constant!examples!nonarchimedean}
\begin{eqnarray}
V_{\infty}(E_v)  & = & -\log_v(R) \ , \qquad \qquad  \  
\gamma_{\infty}(E_v) \ = \ q_v^{-V_{\infty}(E_v)} \ = \ R \ ; \label{FVNAD1}\\
V_{\zeta}(E_v) & = &  -\log_v(R/|\zeta-a|_v^2) \ , \quad 
\gamma_{\zeta}(E_v) \ = \ R/|\zeta-a|_v^2 \ . \label{FVNAD2}
\end{eqnarray}  

\smallskip
{\bf The Open Disc.}\index{examples!nonarchimedean!open disc}  
If $E_v = D(a,R)^- = \{z \in \CC_v : |z|_v < R \}$,
formulas (\ref{FGNAD}), (\ref{FVNAD1}), and (\ref{FVNAD2}) 
for the Green's function, Robin constant and capacity remain valid. 
\index{Robin constant!examples!nonarchimedean}
\index{Green's function!examples!nonarchimedean}
\index{capacity}  

If $\zeta = \infty$, 
this is because for any $R_1 < R$ we have $D(a,R_1) \subset D(a,R)^- \subset D(a,R)$, 
and hence
\begin{eqnarray*}
G(z,\infty;D(a,R_1)) & \le & G(z,\infty;D(a,R)^-) \ \le \ G(z,\infty;D(a,R)) \ , \\
\gamma_{\infty}(D(a,R_1)) & \le & \gamma_{\infty}(D(a,R)^-) \ \le \ \gamma_{\infty}(D(a,R)) \ .
\end{eqnarray*}
\index{Green's function!examples!nonarchimedean}
Taking a limit as $R_1 \rightarrow R$, it follows from formulas 
(\ref{FGNAD}) and (\ref{FVNAD1}) that $G(z,\infty;D(a,R)^-) = G(z,\infty;D(a,R))$
and $\gamma_{\infty}(D(a,R)^-) = \gamma_{\infty}(D(a,R))$.  

If $\zeta \in \CC_v \backslash D(a,R)^-$, 
we can reduce to the case where $\zeta = \infty$
by applying the map $f(z) = (z-a)/(z-\zeta)$ and using the pullback formula (\ref{GreenPullbackF}). 
Thus (\ref{FGNAD}), (\ref{FVNAD1}), and (\ref{FVNAD2}) hold when $E_v = D(a,R)^-$,
for any $\zeta \notin E_v$. 

\smallskip
{\bf The Punctured Disc.}\index{examples!nonarchimedean!punctured disc}  
Suppose $E_v = D(a,R) \backslash (\bigcup_{i=1}^m D(a_i,R_i)^-)$, 
where $a_1, \ldots, a_m \in D(0,R)$ and $R_i \le R$ for each $i$. 
For each $\zeta \notin D(a,R)$, the Green's function and capacity are still given by 
\index{Green's function!examples!nonarchimedean}
\index{capacity} 
(\ref{FGNAD}), (\ref{FVNAD1}), and (\ref{FVNAD2}).  
Indeed, for any fixed $a_0 \in E_v$, 
we have $D(a_0,R)^- \subset E_v \subset D(a_0,R)$, 
so the result follows from the previous case.  

If $\zeta$ belongs to one of the ``holes'' $D(a_i,R_i)^-$, 
then $D(a_i,R_i)^- = D(\zeta,R_i)^-$ and by applying $f(z) = 1/(z-\zeta)$ 
and using the pullback formula (\ref{GreenPullbackF}), we find that 
\index{Green's function!pullback formula for} 
\index{Robin constant!examples!nonarchimedean}
\begin{eqnarray} 
G(z,\zeta;E_v) & = & G(\frac{1}{z-\zeta},\infty;D(0,\frac{1}{R_i})) 
\ = \ \log_v^+\left(\frac{R_i}{|z-\zeta|_v}\right) \ , \label{FGNAD2} \\
V_{\zeta}(E_v) & = & \log_v(R_i) \ , 
\quad \gamma_{\zeta}(E_v) \ = \ 1/R_i \ . \label{FVNAD3} 
\end{eqnarray} 

\smallskip 
{\bf The Ring of Integers $\cO_w$}.\index{examples!nonarchimedean!ring of integers $\cO_w$}  
We next determine the Green's function of the ring of integers of a finite extension $F_w/K_v$ in $\CC_v$.
\index{Green's function!examples!nonarchimedean}

\begin{proposition} \label{OwProp}
Let $F_w/K_v$ be a finite extension in $\CC_v$, with ramification index $e = e_{w/v}$ and residue degree
$f = f_{w/v}$.  Take  $E_v = \cO_w$, the ring of integers of $F_w$.  Given $z \in \CC_v$, put 
\begin{equation*}
r \ = \ \|z,\cO_w\|_v \ = \ \min_{x \in \cO_w} |z - x|_v \ . 
\end{equation*}  
Let $M = \lfloor -e \log_v(r) \rfloor$ and $\langle -e \log_v(r) \rangle$ 
be the integer and fractional parts of $-e \log_v(r)$, respectively.  Then 
\begin{equation} \label{FOwGreen}
G(z,\infty; \cO_w) \ = \ \left\{ \begin{array}{ll} 
             0  & \text{if $z \in \cO_w$,} \\
             \frac{1}{e} \frac{1}{q_v^f-1} \frac{1}{q_v^{fM}} 
                   - \langle -e \log_v(r) \rangle \frac{1}{e} \frac{1}{q_v^{f(M+1)}} 
              & \text{if $z \notin \cO_w$, $|z|_v \le 1$,} \\
             \frac{1}{e} \frac{1}{q_v^f-1} + \log_v(|z|_v) 
                 &   \text{if $|z|_v > 1$.}   
                                 \end{array} \right.
\end{equation}    
\index{Green's function!examples!nonarchimedean}                          
and if capacities are computed relative to the uniformizer $g_{\infty}(z) = 1/z$ then 
\index{Robin constant!examples!nonarchimedean}
\begin{equation} \label{FOwV} 
V_{\infty}(\cO_w) \ = \ \frac{1}{e} \frac{1}{q_v^f-1} \ , \qquad \gamma_{\infty}(\cO_w) \ = \ q_v^{-1/(e(q_v^f-1))} 
\end{equation}  

For any coset $a + b \cO_w$  where $a \in \CC_v$, $b \in \CC_v^{\times}$,  
\begin{equation*}
G(z,\infty;a+b\cO_w) \ = \ G((z-a)/b,\infty;\cO_w) \ ,
\end{equation*} 
\index{Robin constant!examples!nonarchimedean} 
so that $V_{\infty}(a+b \cO_w) = -\log_v(|b|_v) + V_{\infty}(\cO_w)$ 
and $\gamma_{\infty}(a + b \cO_w) = |b|_v \cdot \gamma_\infty(\cO_w)$.   
In particular, if $\pi_w$ is a generator for the maximal ideal of $\cO_w$, then 
\index{Robin constant!examples!nonarchimedean}
\begin{equation} \label{FOwV2}
V_{\infty}(a+\pi_w^m \cO_w) \ = \ \frac{1}{e} \frac{1}{q_v^f-1} + \frac{m}{e} \ ,
\qquad \gamma_{\infty}(a+\pi_w^m \cO_w) \ = \ q_v^{-m/e - 1/(e(q_v^f-1))} \ .
\end{equation} 
\end{proposition}  

\begin{proof}  See (\cite{RR1}, Example 5.2.17).  
The equilibrium distribution of $\cO_w$ is the additive Haar measure $\mu$ for $F_w$, normalized 
\index{Haar measure}
so that $\mu(\cO_w) = 1$ (see \cite{RR1}, p.212).  It follows that if we write $q_w = q_v^f$, 
and put  $M = \lfloor - e \log_v(r) \rfloor$ if $|z|_v \le 1$, $M = -1$ if $|z|_v > 1$,
then the potential function\index{potential function!of $\cO_w$} is given by  
\begin{eqnarray*}
u_{\cO_w}(z,\infty) & = & \int_{\cO_w} -\log_v(|z-x|_v) d\mu(x) \\
      & = & \sum_{k=0}^M \frac{k}{e} \cdot \frac{q_w-1}{q_w^{k+1}} 
                    + 
                    \sum_{k = M+1}^\infty (-\log_v(r)) \cdot \frac{q_w-1}{q_w^{k+1}} \\
      & = & \frac{1}{e} \frac{1}{q_w-1} \cdot [1- \frac{M+1}{q_w^M} + \frac{1}{q_w^{M+1}}] 
                     - \log_v(r) \cdot \frac{1}{q_w^{M+1}} 
\end{eqnarray*} 
The potential function\index{potential function}   is invariant under translation by $\cO_w$, 
\index{Robin constant!examples!nonarchimedean}
so $V_{\infty}(\cO_w) = u_{\cO_w}(0,\infty) = 1/(e(q_w-1))$. 
The expression (\ref{FOwGreen}) 
is obtained by simplifying $G(z,\infty;\cO_w) = 1/(e(q_w-1)) - u_{\cO_w}(z,\infty)$. 
(Compare \cite{RR1}, Example 4.1.24, p.212). 

The assertions about cosets follow easily.
\end{proof}         

We now recall a general procedure for computing capacities 
of finite disjoint unions of nonarchimedean sets (for more details, 
see Theorem \ref{BFThm2A} and Corollary \ref{BFCor2} of Appendix \ref{AppA}, or see \cite{RR1}, p.354).
\index{algorithm!to compute nonarchimedean capacities} 
\index{Green's function!computing nonarchimedean}
\index{Robin constant!computing nonarchimedean}

Let $\cC_v/K_v$ be a curve.  Suppose $E_v = \bigcup_{i=1}^N E_{v,i} \subset \cC_v(\CC_v)$ 
is a finite disjoint union of compact sets $E_{v,i}$ with positive inner capacity, 
\index{capacity!inner} 
and that $\zeta \in \cC_v(\CC_v)$ is such that the canonical distance $[z,w]_{\zeta}$ 
\index{canonical distance!$[z,w]_{\zeta}$}  
\index{canonical distance!nonarchimedean!constant on disjoint balls}  
(see \S\ref{Chap3}.\ref{CanonicalDistanceSection}) is constant on $E_{v,i} \times E_{v,j}$, 
for each $i \ne j$. 
For each $i$, let $\mu_{\zeta,i}$ be the equilibrium distribution of $E_{v,i}$ 
(see \S\ref{Chap3}.\ref{CompactGreenSection}).  Then each $i$, 
the potential function $u_{E_{v,i}}(z,\zeta) = \int_{E_{v,i}} -\log([z,w]_{\zeta} d\mu_{\zeta,i}(w)$ 
and Green's function 
\index{Green's function!computing nonarchimedean}\index{potential function}  
\index{Robin constant!computing nonarchimedean}   
$G(z,\zeta;E_{v,i}) = V_{\zeta}(E_v) - u_{E_{v,i}}(z,\zeta)$ are constant for $z \in E_{v,j}$, 
for each $j \ne i$.  

We now show that we can compute  $G(z,\zeta;E_v)$ and $V_{\zeta}(E_v)$ 
\index{Robin constant!computing nonarchimedean}
in terms of the potential functions $u_{E_{v,i}}(z,\zeta)$.
Let capacities be defined in terms of the uniformizer $g_{\zeta}(z)$.  
For each $E_{v,i}$, put 
\begin{equation*} 
W_{ii} \ = \ V_{\zeta}(E_{v,i}) 
\end{equation*}   
and for each $i \ne j$ let $W_{ij}$ be
the value that $u_{E_{v,i}}(z,\zeta)$ assumes on $E_{v,j}$.  Consider the system of
$N+1$ linear equations in the variables $V$, $s_1, \ldots, s_N$:   
\begin{eqnarray} 
1 & = & 0 V + s_1 + s_2 + \cdots + s_N \ , \label{FSolve}\\
0 & = & V - W_{i1}s_1 -W_{i2} s_2 - \ldots - W_{iN} s_N , \notag \\
& & \qquad \qquad \text{for $i = 1, \ldots, N$.} \notag 
\end{eqnarray} 

We claim that this system of equations has a unique solution, for which $s_1, \ldots, s_N > 0$; 
and for this solution, we have
\index{Robin constant!computing nonarchimedean}
\begin{eqnarray} \label{FUse}
V_{\zeta}(E_v) & = & V \ , \label{FUse1}  \\
G(z;\zeta;E_v) & = & \sum s_i G(z,\zeta;E_{v,i}) + \sum s_i W_{ii} - V \ .
                \label{FUse2}
\end{eqnarray}   

To see this, let $\mu$ be the equilibrium distribution of $E_v$ with respect to $\zeta$, 
and put $\hs_i = \mu(E_{v,i})$ for each $i$.  
Then $\hs_i > 0$: otherwise, $\mu$ would be supported on $E_v \backslash E_{v,i}$ 
and then $u_{\zeta}(z,E_v) = u_{\zeta}(z,E_v \backslash E_{v,i})$.  By (\cite{RR1}, Corollary 4.1.12)
we would have $u_{\zeta}(z,E_v \backslash E_{v,i}) <  V_{\zeta}(E_v\backslash E_{v,i}) = V_{\zeta}(E_v)$
for all $z \in E_{v,i}$, contradicting that $u_{\zeta}(z;E_v)$ takes the 
value $V_{\zeta}(E_v)$ for all $z \in E_v$ except possibly on a set of inner capacity $0$
\index{capacity!inner} 
(\cite{RR1}, Theorem 4.1.11).  Consider the probability measure $\mu_i = \hs_i^{-1} \mu|_{E_{v,i}}$,
and put $u_i(z,\zeta) = \int_{E_{v,i}} -\log_v([z,w]_{\zeta}) \, d\mu_i(w)$;  by our hypothesis on the
canonical distance, $u_i(z,\zeta)$ is constant on $E_j$, for each $j \ne i$.   Then
\index{canonical distance!$[z,w]_{\zeta}$}  
\index{canonical distance!nonarchimedean!constant on disjoint balls}  
\begin{eqnarray*}
u_{E_v}(z,\zeta) & = & \int_{E_v} -\log_v([z,w]_{\zeta}) \, d\mu(z) \\
                  & = & \sum_{i=1}^r \int_{E_{v,i}} -\log_v([z,w]_{\zeta}) \, d\mu(z) 
                  \ = \ \sum_{i=1}^r \hs_i u_i(z,\zeta) \ . 
\end{eqnarray*}  
For each $i$, since $u_{E_v}(z,\zeta)$ and the $u_j(z,\zeta)$ for $j \ne i$  
are constant on $E_{v,i}$ except possibly on a set of inner capacity $0$ , it follows that
\index{capacity!inner}
$u_i(z,\zeta)$ is constant on $E_{v,i}$ except possibly on a set of inner capacity $0$.  
Since this property characterizes the equilibrium potential, it follows that $\mu_i$
must be the equilibrium distribution of $E_{v,i}$ with respect to $\zeta$. 
Thus there are unique weights $\hs_1, \ldots, \hs_N > 0$ with $\sum_{i=1}^N \hs_i = 1$, for which 
\begin{equation} \label{FMun} 
u_{E_v}(z,\zeta) \ = \ \sum_{i=1}^N \hs_i u_{E_{v,i}}(z,\zeta) \ .
\end{equation} 

Evaluating (\ref{FMun}) at a generic point of each $E_{v,i}$, we see that 
$V = V_{\zeta}(E_v)$ and $\hs_1, \ldots, \hs_N$ are a solution to the system (\ref{FSolve})
\index{Robin constant!computing nonarchimedean}
with each $\hs_i > 0$.  
Conversely, any solution to (\ref{FSolve}) gives a system of weights 
for which $\mu = \sum s_i \mu_i$.  
The uniqueness of the equilibrium distribution (\cite{RR1}, Theorem 4.1.22) shows that
$s_1, \ldots, s_N$, and in turn $V$, are unique.  Thus $s_i = \hs_i$ for each $i$, 
and $V_{\zeta}(E_v) = V$. 
\index{Robin constant!computing nonarchimedean}
Since $G(z,\zeta;E_v) = V_{\zeta}(E_v) - u_{E_v}(z,\zeta)$, 
formula (\ref{FUse}) follows.  

\smallskip
{\bf The Group of Units  $\cO_w^{\times}$.}\index{examples!nonarchimedean!group of units $\cO_w^{\times}$}  
Using the machinery above, we will now determine the Green's function 
\index{Green's function!examples!nonarchimedean}
and the capacity of the set $\cO_w^{\times}$, relative to the point $\infty$.
\index{capacity}

\begin{proposition} \label{OwCross}  
Let $F_w/K_v$ be a finite extension, with ramification index $e = e_{w/v}$ 
and residue degree $f = f_{w/v}$.  Let $\cO_w^{\times}$ be the group of units of $\cO_w$.
For $z \in \CC_v$, put $r_0 = \min_{x \in \cO_w^{\times}} |z-x|_v$,
$M_0 = \lfloor - e \log_v(r_0)\rfloor$; 
note that $r_0 = |z|_v$ if $|z|_v > 1$.     
Then  
\begin{equation} \label{FOwCrossGreen}
G(z,\infty; \cO_w^{\times}) \ = \ \left\{ \begin{array}{ll} 
             0  & \text{if $z \in \cO_w^{\times}$,} \\  
           \frac{q_v^f}{e(q_v^f-1)^2} \cdot \frac{1}{q_v^{fM_0}} 
                    - \langle -e \log_v(r_0) \rangle \frac{1}{e} \frac{1}{q_v^{fM_0}}  
                &  \text{if $0 < r_0 \le 1$}, \\       
           \frac{q_v^f}{e(q_v^f-1)^2} + \log_v(|z|_v) 
                 &   \text{if $|z|_v > 1$.}   
                                 \end{array} \right.
\end{equation}    
\index{Green's function!examples!nonarchimedean}                             
If capacities are computed relative to the uniformizer $g_{\infty}(z) = 1/z$ then 
\begin{equation} \label{FOwCrossV} 
V_{\infty}(\cO_w^{\times}) \ = \ \frac{1}{e} \frac{1}{q^f-1} \big( 1+ \frac{1}{q^f-1} \big)
\ = \ \frac{q_v^f}{e(q_v^f-1)^2} \ .
\end{equation}         
\end{proposition} 

\begin{proof}  Put $N = q_v^f-1$ and let $a_1, \ldots, a_N$ be coset representatives 
for the nonzero classes in $\cO_w/\pi_w \cO_w$.  Then 
\begin{equation*}
\cO_w^{\times} \ = \ \bigcup_{i=1}^{N} (a_i + \pi_w \cO_w)
\end{equation*}
is a decomposition of the type needed to compute $G(z,\infty;\cO_w^{\times})$
in terms of the $G(z,\infty;a_i+\pi \cO_w)$.  Applying the last part of 
Proposition \ref{OwProp}, solving the system (\ref{FSolve})
and simplifying (\ref{FUse1}), (\ref{FUse2}) gives the result.  
Here $W_{ij} = 0$ if $i \ne j$ and each $W_{ii} = q_v^f/(e(q_v^f-1))$,   
giving $V = q_v^f/(e(q_v^f-1)^2)$ and $s_i = 1/(q_v^f-1)$ for each $i$.  
\end{proof}

\begin{corollary} \label{OwN} Let $K_v$ be nonarchimedean, 
and let $\pi_v$ be a uniformizer for the maximal ideal of $\cO_v$.  
Suppose $a_1, \ldots, a_N$ are representives for distinct cosets
of $\cO_v/\pi_v \cO_v$, and put $E_v = \cup_{i=1}^N (a_i + \pi_v \cO_v)$.  
Then
\index{Robin constant!examples!nonarchimedean}   
\begin{equation*}
V_{\infty}(E_v) \ = \frac{q_v}{N(q_v-1)}
\end{equation*}
\end{corollary} 

\begin{proof} The proof is similar to Proposition \ref{OwCross};  
with $s_i = 1/N$ for $i = 1, \ldots, N$.
\end{proof} 
  
\smallskip  
{\bf The punctured $\cO_v$-disc.}\index{examples!nonarchimedean!punctured $\cO_v$-disc} 
Next we determine the capacity of a union of cosets of $\cO_v^{\times}$, 
relative to the point $\zeta = \infty$.  This computation has important theoretical
consequences:  it is used in the proof Proposition \ref{IdentifyGreenProp}, 
which plays a key role in the reduction of Theorem \ref{aT1} to Theorem \ref{aT1-B}. 
\index{Robin constant!examples!nonarchimedean}   
  
\begin{proposition} \label{FmCosets}  
Put $E_{v,m} = \bigcup_{k=0}^m \pi_v^k \cO_v^{\times}$, and take $\zeta = \infty$.  Then 
\begin{eqnarray}
V_{\infty}(E_{v,m}) & = & 
\frac{1}{q_v-1} + \frac{1}{(q_v-1)^2 (1 + q_v^2 + q_v^4 + \cdots + q_v^{2m})} \label{FmV} \\
G(0,\infty;E_{v,m}) & = & \frac{q_v^{m+1}}{(q_v-1)^2 (1 + q_v^2 + q_v^4 + \cdots + q_v^{2m})}
\ , \label{FmG0} 
\end{eqnarray}\index{Robin constant!examples!nonarchimedean}   
and for each $k = 0, \ldots, m$ the mass of $\pi_v^k \cO_v^{\times}$ under the
equilibrium distribution $\mu_m$ of $E_{v,m}$ with respect to $\infty$ is 
\begin{equation} \label{FmWt}
\mu_m(\pi_v^k \cO_v^{\times}) \ = \ 
\frac{q_v^k + q_v^{2m+1-k}}{1 + q_v + q_v^2 + q_v^3 + \cdots + q_v^{2m+1}} \ .
\end{equation}  
\end{proposition} 
\index{Green's function!examples!nonarchimedean}

\begin{proof} Write $V_m = V_{\infty}(E_{v,m})$.  By Proposition \ref{OwCross},
we have 
\index{Robin constant!examples!nonarchimedean}
\begin{equation} \label{FV0}
V_0 \ = \  \frac{q_v}{(q_v-1)^2} \ = \ \frac{1}{q_v-1} + \frac{1}{(q_v-1)^2 } \ .
\end{equation} 
We will prove (\ref{FmV}) by induction on $m$.  
Note that $E_{v,m} = \pi_v E_{v,m-1} \cup \cO_v^{\times}$.  
For $z \in \pi_v E_{v,m-1}$ and $w \in \cO_v^{\times}$,
$-\log_v([z,w]_{\infty}) = -\log_v(|z-w|_v) = 0$, independent of $z, w$. 
Hence $u_{\pi_v E_{v,m-1}}(z,\infty) = 0$ if $z \in \cO_w^{\times}$, 
and $u_{\cO_w^{\times}}(z,\infty) = 0$ if $z \in \pi_v E_{v,m-1}$.     
Furthermore, by the scaling property of the capacity, 
\index{capacity!scaling property of}
\index{Robin constant!computing nonarchimedean}
$V_{\infty}(\pi_v E_{v,m}) = V_{\infty}(E_{v,m-1}) + 1 = V_{m-1}+1$.  
By (\ref{FSolve}), there are numbers $s_{1,m}, s_{2,m} > 0$ for which 
\begin{equation} \label{FSys1}
\left\{ \begin{array}{ccc} 
         1 & = & s_{1,m} + s_{2,m} \\
         V_m & = & (V_{m-1} + 1) \cdot s_{1,m} + 0 \cdot s_{2,m} \\
         V_m & = & 0 \cdot s_{1,m} + V_0 \cdot s_{2,m} 
           \end{array} \right. 
\end{equation} 
Solving (\ref{FSys1}) for $V_m$ and inserting (\ref{FV0}) leads to the recursion 
\begin{equation*}
V_m \ = \ \frac{q_v(1+V_{m-1})}{q_v + (q_v-1)^2 V_{m-1}} 
\end{equation*}
whose solution is easily seen to be (\ref{FmV}).

Once the $V_m$ are known, one sees that 
\begin{equation} \label{FSVals} 
s_{1,m} \ = \ \frac{q_v (1 + q_v + \cdots +q_v^{2m-1})}{1+q_v+ \cdots q_v^{2m+1}}  \ , 
\quad s_{2,m} \ = \ \frac{1 + q_v^{2m+1}}{1+q_v+ \cdots q_v^{2m+1}} \ .
\end{equation}

To obtain (\ref{FmG0}), note that since $u_{\cO_w^{\times}}(0,\infty) = 0$, 
we have  $u_{E_{v,m}}(0,\infty) \ = \ s_{1,m} \cdot (1+ u_{E_{v,m-1}}(0,\infty))$.  
Thus recursively 
\begin{equation} \label{FSRec} 
u_{E_{v,m}}(0,\infty) 
\ = \ s_{1,m} + s_{1,m} s_{1,m-1} + \cdots + s_{1,m}s_{1,m-1} \cdots s_{1,1} \ .
\end{equation} 
One gets (\ref{FmG0}) by inserting (\ref{FmV}), (\ref{FSVals}) and (\ref{FSRec}) in 
the formula 
\begin{equation*}
G(0,\infty;E_{v,m}) \ = \ V_m - u_{E_{v,m}}(0,\infty)
\end{equation*} 
\index{Robin constant!computing nonarchimedean}
and simplifying. 
Finally, the weights of the cosets $\pi_v^k \cO_w^{\times}$ 
under the equilibrium distribution $\mu_m$ can be found by using
\begin{eqnarray*}
\mu_m(\pi_v^k \cO_w^{\times}) & = & s_{1,m} \mu_{m-1}(\pi_v^{k-1} \cO_w^{\times}) \ = \ \cdots \\
           & = &  s_{1,m} s_{1,m-1} \cdots s_{1,m-k+1} \cdot \mu_{m-k}(\cO_w^{\times}) 
\end{eqnarray*} 
where $\mu_{m-k}(\cO_w^{\times}) = s_{2,m-k}$.  
Using (\ref{FSVals}), and simplifying, yields (\ref{FmWt}).
Once the weights $\mu_m(\pi_v^k \cO_v^{\times})$ are known, 
the value of $G(z,\infty;E_{v,m})$ can be found for any $z$.  
\end{proof} 

{\bf The union of two rings of integers.}\index{examples!nonarchimedean!non-rational Robin constant}  
Let $F_w$ be the unique unramified quadratic extension of $K_v$, 
and let $F_u$ be a totally ramified quadratic extension.    
We will compute the capacity of the set $E_v = \cO_w \cup \cO_u$ with respect to $\infty$.
\index{capacity}
This is the only nonarchimedean set known to the author whose Robin constant 
\index{Robin constant!examples!nonarchimedean|ii}
can be computed explicitly, and is not rational.    
The importance is not the result itself, but the method, 
which uses a partial self-similarity\index{partial self-similarity} of $E_v$ with itself,
and can be applied to non-disjoint unions of much more general sets.  

\begin{proposition} \label{TwoOvSets}  Fix a nonarchimedean local field $K_v$.  
Let $F_w/K_v$ be the unique unramified quadratic extension,
and let $F_u/K_v$ be a totally ramified quadratic extension.  
Put $E_v = \cO_w \cup \cO_u$ and let 
\begin{eqnarray*}
A & = & 2q_v^4+2q_v^3-4q_v^2+2q_v-2 \ , \\
B & = & q_v^4+2q_v^3-2q_v^2+2q_v-1 \ , \\ 
D & = & q_v^8+4 q_v^7 + 8q_v^6 +12q_v^5 + 18q_v^4 + 12q_v^3 + 8q_v^2 + 4 q_v +1 \ .
\end{eqnarray*} 
Then
\index{Robin constant!examples!nonarchimedean|ii}
\begin{equation} \label{TwoOvFormula} 
V_{\infty}(E_v) \ = \ \frac{-B + \sqrt{D}}{2A} \ .
\end{equation} 
\end{proposition} 

Below are some numerical examples\index{numerical!examples} when $K_v = \QQ_p$, for small primes $p$.  
We give the values of $V_{\infty}(\cO_w)$ and $V_{\infty}(\cO_u)$ for comparison.
\index{Robin constant!examples!nonarchimedean}

\vskip .1 in
\centerline{
\begin{tabular}{|c|c|c|c|c|c|} 
  \hline \hline        & $q_v = 2 $     & $q_v = 3$      & $q_v=5$       & $q_v = 7$    &  $q_v = 11$   \\                                                  \hline \hline 
    $V_{\infty}(E_v)$  & $.2750820518$  & $.1060035774$  & $.0366954968$ & $.0188065868$ & $.0077456591$ \\                                                  \hline
  $V_{\infty}(\cO_w)$  & $.3333333333$  & $.1250000000$  & $.0416666666$ & $.0208333333$ & $.0083333333$ \\                                                  \hline
  $V_{\infty}(\cO_u)$  & $.5000000000$  & $.2500000000$  & $.1250000000$ & $.0833333333$ & $.0500000000$ \\                                                  \hline \hline
\end{tabular}
}

\vskip .1 in
\noindent{It} can be shown that as $q_v \rightarrow \infty$, 
then $V_{\infty}(E_v) = 1/q_v^2-1/q_v^3+O(1/q_v^4)$.

\begin{proof}[Proof of Proposition \ref{TwoOvSets}]
Let $\pi = \pi_v$ be a generator for the maximal ideal of $\cO_v$,  
and write $q = q_v$. Then $\#(\cO_v/\pi \cO_v) = q$; let $\gamma_1, \ldots, \gamma_q$ 
be coset representatives for $\cO_v/\pi \cO_v$.  
Put $E_{0,i} = \gamma_i + \pi E_v = \gamma_i + \pi (\cO_w  \cup \cO_u)$, for $i = 1, \ldots, q$. 
There are $q^2-q$ cosets of $\cO_w/\pi \cO_w$ which do not contain elements of $\cO_v$;  
let these be $E_{1,j} = \alpha_j + \pi \cO_w$, for $j = 1, \ldots, q^2-q$.  Similarly,
there are $q^2-q$ cosets of $\cO_u/\pi \cO_u$ which do not contain elements of $\cO_v$;  let these
be $E_{2,k} = \beta_k + \pi \cO_u$, for $k = 1, \ldots, q^2-q$.  
Then the sets $E_{0,i}$, $E_{1,j}$ and $E_{2,k}$ are pairwise disjoint 
(in fact, they are contained in pairwise disjoint cosets $a + \pi \hcO_v$, 
where $\hcO_v = D(0,1)$ is the ring of integers of $\CC_v$), and we can write 
\begin{equation*} 
E_v \ = \ \Big( \bigcup_{i=1}^{q} E_{1,i} \Big) \cup 
          \Big( \bigcup_{j=1}^{q^2-q} E_{2,j} \Big) \cup
          \Big( \bigcup_{k=1}^{q^2-q} E_{3,k} \Big) \ .
\end{equation*}   

Let $\mu$ be the equilibrium distribution of $E_v$ with respect to $\infty$, 
and put $w_{0,i} = \mu(E_{0,i})$, $w_{1,j} = \mu(E_{1,j})$, $w_{2,k} = \mu(E_{2,k})$ for all $i$, $j$, $k$.
Then 
\begin{equation} \label{F3SumEq} 
u_{\infty}(z,E_v) \ = \ \sum_{i=1}^{q} w_{1,i} u_{\infty}(z,E_{0,i}) + 
\sum_{j=1}^{q^2-q} w_{0,j} u_{\infty}(z,E_{1,j}) + \sum_{k=1}^{q^2-q} w_{2,k} u_{\infty}(z,E_{2,k}) \ .
\end{equation}  
\index{Robin constant!computing nonarchimedean}
Let $V = V_{\infty}(E_v)$ be the (as yet unknown) Robin constant of $E_v = \cO_w \cup \cO_u$,
\index{Robin constant!examples!nonarchimedean} 
and let $V_1 = V_{\infty}(\cO_w)$, $V_2 = V_{\infty}(\cO_u)$.  
Since $E_v \subset D(0,1)$, we must have $V \ge 0$.  By Proposition \ref{OwProp} 
\begin{equation} \label{FV1V2Vals}  
V_1 = \frac{1}{q^2-1} \ , \qquad V_2 = \frac{1}{2(q-1)} \ .
\end{equation} 

In general, for any compact set $\tE \subset \CC_v$ of positive capacity, 
\index{capacity}
we have $V_{\infty}(a + \pi \tE) = V_{\infty}(\tE) + 1$ for each $a \in \CC_v$. 
If $\tE \subset D(a,r)$, then $u_{\infty}(z,\tE) = -\log_v(|z-a|_v)$ 
for all $z \notin D(a,r)$.  
It follows that for each $E_{0,i}$, one has  $u_{\infty}(z,E_{0,i}) = V + 1$  on $E_{0,i}$.
On the $q-1$ cosets $E_{2,k}$ contained in $\gamma_i + \sqrt{\pi} \hcO_v$,
one has $u_{\infty}(z,E_{0,i}) = 1/2$.  
On the other $q^2-2q+1$ cosets $E_{2,k}$ and the other $q-1$ cosets $E_{0,i^{\prime}}$,
as well as all the cosets $E_{1,j}$, one has $u_{\infty}(z,E_{0,i}) = 0$. 
For each $E_{1,j}$, one has $u_{\infty}(z,E_{1,j}) = V_1 + 1$ on $E_{1,j}$, 
\index{Robin constant!computing nonarchimedean}
and $u_{\infty}(z,E_{1,j}) = 0$ on all the $E_{0,i}$, all the $E_{2,j}$ and all the $E_{1,j^{\prime}}$
distinct from $j$.  
For each $E_{2,k}$, one has  $u_{\infty}(z,E_{2,k}) = V_2 + 1$  on $E_{2,k}$. 
There are $q-2$ other cosets $E_{2,k^{\prime}}$ and one coset $E_{1,j}$ 
contained in $\beta_k + \sqrt{\pi} \hcO_v$.  
On those cosets we have $u_{\infty}(z,E_{2,k}) = 1/2$. 
On the remaining $q^2-2q+1$ cosets $E_{2,k^{\prime}}$ and on all the cosets $E_{1,j}$, 
one has $u_{\infty}(z,E_{2,k}) = 0$.

\smallskip

Evaluating $u_{\infty}(z,E_v)$ on each of the sets $E_{r,s}$ in turn 
yields a system of $2q^2-q$ equations satisfied by $V$ and the $w_{r,s}$.   
Since $\mu$ and $V = V_{\infty}(E_v)$ are unique, 
\index{Robin constant!computing nonarchimedean}
these equations uniquely determine the $w_{r,s}$.  
Hence for any permutation $\sigma$ of the sets $E_{r,s}$   
which takes sets of type $r = 0, 1, 2$ to sets of the same type,   
and which preserves distances between corresponding pairs of sets, 
we must have $w_{r,\sigma(s)} = w_{r,s}$ for all $r, s$.  It is easy to see that there are 
enough permutations satisfying these conditions to assure that 
there are $w_0$, $w_1$, $w_2$ such that for all $i$, $j$, $k$
\begin{equation*}
w_{0,i} = w_0 \ , \qquad w_{1,j} = w_1 \ , \qquad w_{2,k} = w_2 \ .
\end{equation*}  

We can now determine $V$.  From $\mu(E_v) = 1$, we obtain the mass equation 
\begin{equation*} 
1 \ = \  (q) \cdot w_0 + (q^2-q) \cdot w_1 + (q^2-q) \cdot w_2  \ .
\end{equation*} 
Evaluating $u_{\infty}(z,E_v)$ on the sets $E_{0,i}$, $E_{1,j}$ and $E_{2,k}$  
gives the equations  
\begin{eqnarray*}
    V & = &  w_0 \cdot (V + 1) +  w_2 \cdot (q-1) \cdot (1/2) \ , \\
    V & = &  w_1 \cdot (V_1 + 1) \ , \\
    V & = &   w_0 \cdot (1/2) + w_2 \cdot ((V_2 + 1) + (q-2) \cdot (1/2)) \ .    
\end{eqnarray*} 
\index{Robin constant!computing nonarchimedean}
Treating this as a linear system in $w_0$, $w_1$, $w_2$, 
solving it in terms of $V$, $V_1$, $V_2$, 
and inserting the resulting values in the mass equation leads to    
\begin{equation*} 
1 \ = \ (q) \frac{V (\frac{1}{2} + V_2)  }{(1+V)(V_2 + \frac{q}{2})- \frac{q-1}{4}}
           + (q^2-q) \frac{V}{1+V_1} 
           + (q^2-q) \frac{ V^2 + \frac{1}{2} V} {(1+V)(V_2 + \frac{q}{2})- \frac{q-1}{4}} \ . 
\end{equation*} 
Clearing denominators and using the values for $V_1$, $V_2$ from (\ref{FV1V2Vals}) 
yields a quadratic equation in $V$.  
Its unique non-negative root (simplified using Maple)\index{Maple computations} 
is the one in (\ref{TwoOvFormula}). 
\end{proof}  
 
 \smallskip
\section{ Global Examples on $\PP^1$} \label{PP1Examples}

As will be seen, capacity theory provides a ``calculus'' 
\index{capacity theory}
for answering certain types of questions about algebraic integers and units.
Note that $\alpha \in \QQbar$ is an algebraic integer
\index{algebraic integer}  
if and only if its conjugates all 
satisfy $|\sigma(\alpha)|_v \le 1$ for all nonarchimedean $v$, and it is a unit
if and only if $|\sigma(\alpha)|_v = 1$ for all nonarchimedean $v$.  

\smallskip
{\bf Algebraic Integers}. 
The following example is a trivial application of capacity theory, 
\index{capacity theory}
but appears hard to prove without it.
 
\begin{example} \label{MandelbrotApplication}\index{examples!global!Mandelbrot set}
Let $\cM$ be the Mandelbrot set. Then 
\index{Mandelbrot set}

$(A)$ There are infinitely many algebraic integers whose conjugates all belong to $\cM$.

$(B)$ For each number $B > 0$ there are 
only finitely many algebraic integers $\alpha$ whose conjugates all belong to $\cM$,
and some prime $p \le B$ 
splits completely in $\QQ(\alpha)$. 
Indeed, there is a neighborhood $U = U(B)$ of $\cM$ with this property. 

$(C)$ On the other hand, for each neighborhood $U$ of $\cM$ in $\CC$, there is a number $C = C(U)$ 
such that for each prime $p > C$, there are infinitely many algebraic integers $\alpha$
such that all the conjugates of $\alpha$ belong to $U$, and $p$ splits completely in $\QQ(\alpha)$.  
\end{example}  

\begin{proof}
Take $K = \QQ$, $\cC = \PP^1$, and $\fX = \{\infty\}$. 

Part (A) is well known.  Indeed, put $\varphi_c(z) = z^2 + c$ and for each integer $n \ge 1$ 
put $P_n(c) = \varphi_c^{(n+1)}(0)$, as in the discussion preceding (\ref{MandelDef}).  Then $P_n(c)$
is a monic polynomial in $\ZZ[c]$ of degree $2^n$.  If $\alpha$ is a root of $P_n(c) = 0$, 
then $z=0$ is periodic for $\varphi_\alpha(z)$ (with period dividing $n+1$) since $\varphi_\alpha^{(n+1)}(0) = 0$.  
The same is true for all the $\Gal(\QQbar/\QQ)$-conjugates of $\alpha$, 
so $\alpha$ is an algebraic integer whose conjugates all belong to $\cM$. 
\index{algebraic integer}  
 
There are many ways to see that as a collection, the $P_n(c)$ have infinitely many distinct roots.
For example, note that $c = 0$ is the only number such that $0$ is periodic for 
$\varphi_c(z)$ with period $1$.  Taking $n = p-1$ where $p$ is prime, 
we obtain $2^n-1$ values of $c$ such that $0$ is periodic for $\varphi_c(z)$ 
with exact period $p$\,; 
thus there are infinitely many algebraic integers whose conjugates all belong to $\cM$. 

\smallskip
For part (B), fix a prime $p$, let  $E_\infty = \cM$,  $E_p = \ZZ_p$,
and let $E_q = D(0,1) \subset \CC_q$ for each prime $q \ne p$.  
Put $\EE = \prod_v E_v$.  Then $\EE$ is algebraically capacitable, and 
\index{algebraically capacitable}
\begin{equation*} 
\gamma(\EE,\fX) \ = \ \gamma_{\infty}(\cM) \cdot \gamma_{\infty}(E_p) 
\ = \ p^{-1/(p-1)} \ < \ 1 \ . 
\end{equation*} 
By Theorem \ref{FSZii} there is an 
adelic neigbhorhood $\UU = \UU_p = \prod_v U_{p,v}$ of $\EE$ such that there are 
only finitely many $\alpha \in \QQbar$ which have all their conjugates in $\UU_p$.
Each algebraic integer $\alpha$ such that $p$ splits completely in $\QQ(\alpha)$ is such a number.
\index{algebraic integer} 
  
Given $B > 0$, put $U = U(B) = \cap_{p \le B} U_{p,\infty} \subset \CC$.  
Then $U(B)$ has the desired properties. 

\smallskip
For part (C),
let $U \subset \CC$ be any neighborhood of $\cM$.  
By enlarging $\cM$ within $U$ (for example by choosing a point $a \in (U \cap \RR) \backslash \cM$
and adjoining a suitably small disc $D(a,r)$) 
we can obtain a set $\cM_U \subset U$ 
which has $\gamma_{\infty}(\cM_U) > 1$ and is stable under complex conjugation.  

Fix a prime $p$, and take $E_{\infty} = \cM_U$,  $E_p = \ZZ_p$,
and $E_q = D(0,1) \subset \CC_p$ for each prime $q \ne p$.  
Put $\EE = \prod_v E_v$.  By (\ref{FVNAD1}) and (\ref{FOwV}), 
the capacity of $\EE$ with respect to $\fX$ is 
\index{capacity}
\begin{equation*} 
\gamma(\EE,\fX) \ = \ \gamma_{\infty}(\cM_U) \cdot p^{-1/(p-1)} \ .  
\end{equation*}   
It follows that if $p$ is sufficiently large, then $\gamma(\EE,\fX) > 1$. 
Put $U_{\infty} = U$, $U_p = \ZZ_p$, and $U_q = D(0,1)$ for each prime $q \ne p$.  
By Theorem \ref{aT1-A1}, there are  
infinitely many $\alpha \in \QQbar$ whose conjugates belong to $U_v$ for each $v$.  
Each such $\alpha$ is an algebraic integer whose archimedean conjugates belong to $U$ 
\index{algebraic integer} 
and whose conjugates in $\CC_p$ belong to $\ZZ_p$, so $p$ splits completely in $\QQ(\alpha)$.
\end{proof} 

We remark that the same assertions hold for the Julia set of a monic polynomial $g(x) \in \ZZ[x]$ 
\index{Julia set} 
with degree $d > 1$ (see (\ref{JuliaCapacity})).  
In this case, each repelling periodic point for $g(x)$ is an algebraic integer 
\index{algebraic integer} 
whose conjugates belong to the Julia set.

\smallskip
Our next result, originally formulated by Cantor (\cite{Can3}) and proved in (\cite{RR2}),
\index{Cantor, David} 
generalizes the classical theorem of Robinson (\cite{Rob1}) 
which was the prototype for the Fekete-Szeg\"o theorem with local rationality conditions.
\index{Fekete-Szeg\"o theorem with LRC}  

\begin{example}  \label{RealAlgInt}\index{examples!global!segment $[a,b]$} 
Let $\cQ$ be a finite set of primes of $\QQ$, and let $[a,b] \subset \RR$.  If 
\begin{equation*}
b-a \ > \ 4  \cdot \prod_{q \in \cQ} q^{1/(q-1)} \ ,
\end{equation*} 
then there are infinitely many algebraic integers $\alpha$ 
whose conjugates all belong to $[a,b]$ and for which the primes in $\cQ$ split completely in 
$\QQ(\alpha)$;  if $b-a < 4 \cdot \prod_{q \in \cQ} q^{1/(q-1)}$ there are only finitely many.
\end{example}

\begin{proof} Take $K = \QQ$, $\cC = \PP^1$, and $\fX = \{\infty\}$.  Put $E_{\infty} = [a,b]$,
$E_q = \ZZ_q$ for $q \in \cQ$, and $E_p = \hcO_p$ for finite $p \notin \cQ$.  
Then for $\EE = E_{\infty} \times \prod_p E_p$, 
using $g_{\infty} = 1/z$ to compute the local capacities, formulas (\ref{F1Seg3}) and (\ref{FOwV}) give   
\begin{equation*}
\gamma(\EE,\fX) \ = \ \prod_{p,\infty} \gamma_{\infty}(E_p) \ = 
\ \frac{b-a}{4} \cdot \prod_{q \in \cQ} q^{-1/(q-1)} \ .
\end{equation*}
Thus the result follows from Theorem \ref{FSZii}.
\end{proof}   
 
\vskip .05 in
Over an arbitrary number field, we have the following generalization of Example \ref{RealAlgInt},   
motivated by a result of Moret-Bailly (\cite{MB2}, Th\'eor\`eme 1.3, p.182).
\index{Moret-Bailly, Laurent} 

\begin{example}\index{examples!global!of Moret-Bailly type}  \label{AlgInt} 
Let $K$ be a number field, with $r_1$ real places and $r_2$ complex places.
Write $n = [K:\QQ]$, and let $\cQ$ be a finite set of nonarchimedean places of $K$.  
For each $v \in \cQ$, let $F_w/K_v$ be a finite galois extension, with ramification index
$e_{w/v}$ and residue degree $f_{w/v}$.  
If
\begin{equation} \label{FRsize} 
R^n \ > \ 2^{r_1} \prod_{v \in \cQ} q_v^{1/(e_{w/v}(q_v^{f_{w/v}}-1))}
\end{equation}
then there are infinitely many algebraic integers $\alpha$  
whose archimedean conjugates belong to $D(0,R)$ at each $v$ where $K_v \cong \CC$,
to $[-R,R]$ at each $v$ where $K_v \cong \RR$, 
and are such that for each $v \in \cQ$ all the conjugates in $\CC_v$ belong to $\cO_{F_w}$.
  
If $R^n$ is less than the bound in $(\ref{FRsize})$, there are only finitely 
many such algebraic integers. 
\end{example}

\begin{proof}   
For each complex archimedean $v$, put $E_v = D(0,R)$;  then $\gamma_{\infty}(E_v) = R$. 
For each real archimedean $v$, put $E_v = [-R,R] \subset \RR$; then $\gamma_{\infty}(E_v) = R/2$. 
For each nonarchimedean $v \in \cQ$, put $E_v = \cO_w$, and write $e = e_{w/v}$, $f = f_{w/v}$; 
then  $\gamma_{\infty}(E_v) = q_v^{-1/e(q_v^{f}-1)}$ by (\ref{FOwV}).
For all other nonarchimedean $v$, put $E_v = D(0,1)$, and put $\EE = \prod_v E_v$. 
 By our convention about weights and absolute values in 
the complex archimedean case, 
\begin{eqnarray*}
\gamma(\EE,\{\infty\}) & = & \prod_{\text{real $v$}} \gamma_{\infty}(E_v) 
              \cdot \prod_{\text{complex $v$}} \gamma_{\infty}(E_v)^2  
              \cdot \prod_{\text{finite $v$}} \gamma_{\infty}(E_v) \\
                       & = & R^n \cdot 2^{-r_1} \cdot \prod_{v \in \cQ} q_v^{-1/(q_v-1)} \ .
\end{eqnarray*}
Again the result follows from Theorem \ref{FSZii}. 
\end{proof} 

\vskip .05 in
In general, the behavior is not known in the extremal case when 
$R^n = 2^{r_1} \prod_{v \in S} q_v^{1/(q_v-1)}$.  
When $K$ is totally real, and $S$ is empty, we have $R =2$, 
and the roots of Chebyshev polynomials belong to $[-2,2]$.
\index{Chebyshev polynomial} 
When $K$ is totally complex, $R = 1$ and the roots of unity belong to $D(0,1)$.  
Thus in these two cases there are infinitely many algebraic integers whose conjugates 
satisfy the required conditions.  There are no known examples where there are only finitely
many.  

\medskip
Algebraic numbers satisfing various arithmetic conditions, with controlled
archimedean conjugates, can be constructed by imposing appropriate geometric conditions
on the sets $E_v$.  The following (admittedly contrived) example illustrates 
some of the possibilities. 

\begin{example}\index{examples!global!contrived} \label{Contrived}  
For any $\varepsilon > 0$, there are infinitely many algebraic integers $\alpha$ such that

$(1)$ each archimedean conjugate  $\sigma(\alpha)$ is real 
and satisfies  $0 < \sigma(\alpha) < 12 \sqrt{5} + \varepsilon$;

$(2)$ the primes $\fp_v$ above $2$ in $\QQ(\alpha)$ have residue degree $1$, 
and $|\alpha|_v = 1$ at each $v|2$.  

$(3)$ the prime $3$ is unramified in $\QQ(\alpha)$, and $\ord_v(\alpha) = 1$ at all $v$ above $3$;

$(4)$ the prime $5$ splits completely in $\QQ(\alpha)$, 
and $\alpha$ is a quadratic nonresidue at each $v|5$;  

$(5)$ for all primes $\fp_v$ of $\QQ(\alpha)$ above $7$, we have $\alpha \equiv -1 \pmod{\fp_v}$.  
\end{example} 

\begin{proof}  Take $K = \QQ$ and let $L > 0$ be a parameter.  
Put $E_{\infty} = [0,L] \subset \RR$,  
$E_2 = D(1,1)^-$, 
$E_3 = D(2/3,1/3)^-$, 
$E_5 = \ZZ_5 \cap (D(2,1)^- \cup D(3,1)^-)$,  
$E_7 = D(-1,1)^-$.  Put $E_p = D(0,1)$ for all other primes $p$.  
Then $\gamma_{\infty}(E_{\infty}) = L/4$.   
As seen in the discussion of capacities of nonarchimedean open discs, 
$\gamma_{\infty}(E_2) = \gamma_{\infty}(D(1,1)) = \gamma_{\infty}(D(0,1)) = 1$ and 
$\gamma_{\infty}(E_3) = \gamma_{\infty}(D(2/3,1/3)) = \gamma_{\infty}(D(0,1/3)) = 1/3$.  
Corollary \ref{OwN} shows that $\gamma_{\infty}(E_5) = 5^{-2/4}$.  
At $p = 7$, we have $E_7 \subset D(0,1)$ so $\gamma_{\infty}(E_7) \le 1$;  
on the other hand for each totally ramified finite extension $F_w/\QQ_7$ we have
$-1+\pi_w \cO_w \subset E_7$, and Proposition \ref{OwProp} shows that 
$\gamma_{\infty}(-1+\pi_w \cO_w) = 7^{-1/e_w} \cdot 7^{-1/6e_w}$. 
Letting $e_w \rightarrow \infty$, we see that $\gamma_{\infty}(E_7) = 1$.
As noted in Theorem \ref{FSZi}, the condition that the primes above $2$ have residue degree $1$, 
and the primes above $3$ be unramified, 
can be imposed `for free'.  
For the Fekete-Szeg\"o theorem \ref{FSZi} to be applicable,
\index{Fekete-Szeg\"o theorem with LRC} 
we need $L/4 \cdot 1/3 \cdot 5^{-1/2} > 1$.
\end{proof}  

\medskip
The following example illustrates a case in which some of the $E_v$ are unions 
of ``different types'' of sets, with overlaps.  

\begin{example}\index{examples!global!overlapping sets} \label{Intersection sets}  Take $K = \QQ$;  
let $\fX = \{\infty\}$, put $E_{\infty} = D(0,1) \cup [1, 1 + L]$ where $L \ge 0$,
and put $E_3 = \cO_{v_1} \cup \cO_{v_2}$ where  $\cO_{v_1}$ is the 
ring of integers of the unramified extension $L_{v_1} = \QQ_3(\sqrt{-1})$, 
and $\cO_{v_2}$ is the  ring of integers of the totally ramified extension 
$L_{v_2} = \QQ_3(\sqrt{-3})$.  
For each finite prime $p \ne 2$, let $E_p = \hcO_v$ be the $\fX$-trivial set.
\index{$\fX$-trivial}
  
By formula (\ref{DiscOppRadial}) with $L_1 = 0$ and $L_2 = L$, 
and formula (\ref{TwoOvFormula}) with $q_v = 3$, we have 
\begin{eqnarray*} 
\gamma(\EE,\{\infty\}) & = & 3^{-(-61+\sqrt{6481})/184} 
\cdot \Big(1 + \frac{L^2}{4(1+L)}\Big) \\
& \cong & 0.89000685 + 0.22251713 L^2/(1+L) \ .
\end{eqnarray*} 
By Theorems \ref{aT1} and \ref{FSZii}, 
if $L > 0.99240793$ then there are infinitely many 
algebraic integers whose archimedean conjugates all lie in  
$D(0,1) \cup [1,1+L]$ and whose $\CC_3$-conjugates all lie in 
$\cO_{v_1} \cup \cO_{v_2}$, while if $L < 0.99240792$
there are only finitely many.  
\end{example} 

\smallskip
Our last result in this section 
is a continuation of an example of Cantor (\cite{Can3}, p.167).
\index{Cantor, David} 
Suppose that instead of constructing algebraic integers, one has a rational 
function $f(x)$ and is interested in constructing numbers $\alpha$
for which $f(\alpha)$ is an algebraic integer.
\index{algebraic integer} 

For instance, let $f(x) = 1/(1+x^2)$.  Using Fekete's Theorem, Cantor showed
\index{Fekete's theorem} 
\index{Cantor, David} 
that there are only finitely many totally real $\alpha$ 
for which $f(\alpha)$ is an algebraic integer;  indeed, $\alpha = 0$ 
\index{algebraic integer} 
and $\alpha = \infty$ are the only such points.  

Suppose, however, that we were willing to accept numbers $\alpha$, all of whose 
conjugates had a small imaginary part.  How large would the imaginary parts have 
to be to guarantee the existence of infinitely many solutions?

\begin{example}\index{examples!global!Cantor example continued} \label{FunctionIntVals}  
Take $f(x) = 1/(1+x^2)$.  
Suppose $T > 3/4$.  
Then there are infinitely many  $\alpha \in \widetilde{\QQ}$, 
all of whose conjugates satisfy $|\Im(\sigma(\alpha))| < T$,
for which $f(\alpha)$ is an algebraic integer.
\index{algebraic integer} 
However, if $T < 3/4$, there are only finitely many.  
\end{example} 

\begin{proof} Take $K = \QQ$, and let $\fX = \{i,-i\}$, the set of poles of $f(x)$.
Put $g_i(z) = z-i$, $g_{-i}(z) = z+i$, and let $L = \QQ(i)$.  Fix $T > 0$.  

At the archimedean place of $\QQ$, take 
\begin{equation*}
E_{\infty} \ = \ \{z \in \CC : -T \le \Im(z) \le T \} \cup \{\infty\} \ .
\end{equation*}
At each finite prime $p$, put $E_p = f^{-1}(\hcO_p)$.  It is easy to see that
$E_2 = \PP^1(\CC_2) \backslash B(1,1)^-$, while for each $p \ge 3$, 
$E_p = \PP^1(\CC_p) \backslash (B(i,1)^- \cup B(-i,1)^-)$ where the two balls
are distinct.  Put $\EE = E_{\infty} \times \prod_p E_p$.   

To compute the Green's matrices we must make a base change to $L$.   Recall that 
$\Gamma(\EE,\fX) = [L:K]^{-1} \Gamma(\EE_L,\fX)$.
\index{Green's matrix!global}

There is one archimedean place of $L$, which we will denote $w_\infty$.
By (\ref{FVHP1}), we have $V_i(E_{w_\infty}) = V_{-i}(E_{w_\infty}) = \log(2(1-T))$,
\index{Robin constant!examples!archimedean}
while $G(-i,i;E_{w_\infty}) = G(i,-i;E_{w_\infty}) = 0$ since $i$ and $-i$ belong
to distinct components of $\PP^1(\CC) \backslash E_{\infty}$.  
Since $L_{w_\infty} \cong \CC$, we have $\log(q_{w_\infty}) = \log(e^2) = 2$. 
There is one place $w_2$ of $L$ above $2$;  $L_{w_2}/\QQ_2$ is totally ramified. 
Fixing an isomorphism $\CC_{w_2} \cong \CC_2$, identify $E_{w_2}$ with $E_2$.  
Then $V_i(E_{w_2}) = V_{-i}(E_{w_2}) = 0$, 
\index{Robin constant!examples!nonarchimedean}
while $G(-i,i;E_{w_2}) = G(i,-i;E_{w_2}) = -\log_2(|i-(-i)|_{w_2}) = 2$. 
We have $\log(q_{w_2}) = \log(2)$.
For all other places $v$ of $L$ the Green's matrices are trivial.  Thus
\begin{equation} \label{FFGreen} 
\Gamma(\EE,\fX) 
\ = \       \left( \begin{array}{cc} 
                 \log(2(1-T)) &  \log(2) \\
                  \log(2) & \log(2(1-T)) 
            \end{array} \right) \ .          
\end{equation}
\index{Green's matrix!global}
By definition $\gamma(\EE,\fX) = \exp(-\val(\Gamma(\EE,\fX)))$, where
\begin{equation*}
\val(\Gamma(\EE,\fX)) 
\ = \ \min_{\vr \in \fp} \max_{\vs \in \fp} {}^t\vr \Gamma(\EE,\fX) \vs 
\ = \ \max_{\vr \in \fp} \min_{\vs \in \fp} {}^t\vr \Gamma(\EE,\fX) \vs
\end{equation*}
\index{Green's matrix!global}\index{value of $\Gamma$ as a matrix game}
is the value of $\Gamma(\EE,\fX)$ as a matrix game;
here $\fp$ is the set of probability vectors in $\RR^2$.  
If we take $\vs = {}^t(\frac{1}{2},\frac{1}{2})$ then both entries 
of $\Gamma(\EE,\fX)\vs$ are equal to $\frac{1}{2}\log(4(1-T))$;  
combining the ``mini-max'' and ``maxi-min'' expressions for 
\index{minimax property} 
$\val(\Gamma(\EE,\fX))$ shows that 
$\val(\Gamma(\EE,\fX)) = \frac{1}{2} \log(4(1-T))$.  

Thus $\gamma(\EE,\fX) > 1$ iff $T > 3/4$,
while $\gamma(\EE,\fX) < 1$ iff $T < 3/4$, 
and the result follows from the Fekete-Szeg\"o theorem \ref{FSZii}.
\index{Fekete-Szeg\"o theorem with LRC} 
\end{proof} 

\smallskip
{\bf Algebraic Units.} As in Example \ref{Contrived}, 
the condition that an algebraic integer $\alpha$ be a
\index{algebraic integer} 
unit at finitely many specified primes can be imposed `for free'.  
However, if we want global algebraic units, the construction must assure
that they avoid $0$ (and $\infty$) at all nonarchimedean $v$.  
This can be accomplished by using the capacity relative to two points 
\index{capacity!Cantor capacity}
$\fX = \{0,\infty\}$.  

Below is the theorem of Robinson (\cite{Rob2}) cited in the Introduction,\index{Robinson, Raphael} 
which was originally proved without using capacity theory.\index{capacity theory}
The fact that Robinson's conditions arise naturally in the context of 
capacities was first recognized by Cantor (\cite{Can3}):  
\index{Cantor, David}  

\begin{example}[{\bf Robinson}]\index{examples!global!Robinson's unit theorem} \label{RobUnitThm} 
Suppose $0 < a < b \in \RR$ satisfy the conditions
\begin{eqnarray}
\log(\frac{b-a}{4}) \ > \ 0 \ , \qquad \qquad \qquad \quad \label{FRobC1}\\
\log(\frac{b-a}{4}) \cdot \log(\frac{b-a}{4ab})  -  
        \Big(\log(\frac{\sqrt{b}+\sqrt{a}}{\sqrt{b}-\sqrt{a}})\Big)^2 
       \ > \ 0 \ .  \label{FRobC2} 
\end{eqnarray}
Then there are infinitely many totally real units $\alpha$ 
whose conjugates all belong to $[a,b]$.  
\end{example} 

\begin{proof}  We follow Cantor (\cite{Can3}, p.166).
\index{Cantor, David} 
Take $K = \QQ$, $\cC = \PP^1$, and $\fX = \{0,\infty\}$.  Put
$E_{\infty} = [a,b]$, and put $E_p = D(0,1) \backslash D(0,1)^-$ for each 
finite prime $p$.  Each nonarchimedean $E_p$ is the `$\fX$-trivial'  
\index{$\fX$-trivial}
set in $\PP^1(\CC_p)$, so we can take $\EE = E_{\infty} \times \prod_p E_p$.  

Let the uniformizing parameters\index{uniformizing parameter!normalizes capacity} 
used to compute capacities be $g_0(z) = z$, 
$g_{\infty}(z) = 1/z$. 
By formulas (\ref{F1Seg3}), (\ref{F1Seg1}),  at the archimedean place w
e have $V_{\infty}([a,b]) = \log(4/(b-a))$ and
\index{Robin constant!examples!archimedean}
$G(0,\infty;[a,b]) = \log((\sqrt{b}+\sqrt{a})/(\sqrt{b}-\sqrt{a}))$. Pulling back by $1/z$,
we have $G(z,0;[a,b]) = G(1/z,\infty;[1/b,1/a])$.  In view of our choices of the
uniformizing parameters,\index{uniformizing parameter!normalizes Robin constant} 
this yields $V_0([a,b]) = V_{\infty}([1/b,1/a]) = \log(4ab/(b-a))$.  
At each finite prime $p$, one sees easily that 
$V_0(E_p) = V_{\infty}(E_p) = G(0,\infty;E_p) = 0$.
\index{Robin constant!examples!nonarchimedean}
Thus
\begin{equation} \label{FGamma}
\Gamma(\EE,\fX) \ = \ \Gamma(E_{\infty},\fX) \ = \ 
\left( \begin{array}{cc} 
         \log(\frac{4ab}{b-a}) & \log(\frac{\sqrt{b}+\sqrt{a}}{\sqrt{b}-\sqrt{a}}) \\  
         \log(\frac{\sqrt{b}+\sqrt{a}}{\sqrt{b}-\sqrt{a}}) & \log(\frac{4}{b-a}) 
               \end{array} \right) \ .
\end{equation}
\index{Green's matrix!global}
The conditions (\ref{FRobC1}), (\ref{FRobC2}) 
in the Theorem are simply the determinant inequalities on the  
minors of $\Gamma(\EE,\fX)$,\index{Determinant Criterion!for negative definiteness} 
necessary and sufficient for it to be negative definite.\index{Green's matrix!negative definite}
Hence the result follows from the Fekete-Szeg\"o theorem \ref{FSZii}.
\index{Fekete-Szeg\"o theorem with LRC} 
\end{proof}   

\smallskip
In the next result, we bound the size of the units and their reciprocals,
as well as imposing conditions at nonarchimedean places.

\begin{example}\index{examples!global!units in two segments} \label{UnitThm2} 
There are infinitely many totally real algebraic units $\alpha$ whose archimedean conjugates 
belong to $[-r,-1/r]\cup[1/r,r]$, if  
\begin{equation*}
r \ > \ 1 + \sqrt{2} \ .
\end{equation*} 
More generally, let $\cQ$ be a finite set of primes, 
and put $A = \prod_{q \in \cQ} q^{q/(q-1)^2}$.  
Then there are infinitely many totally real algebraic units $\alpha$ 
for which the primes $q \in \cQ$ split completely in $\QQ(\alpha)$, 
and whose archimedean conjugates belong to $[-r,-1/r]\cup[1/r,r]$ if
\begin{equation} \label{FUnit2}
r \ > \ A^2 +  \sqrt{A^4+1} \ . 
\end{equation} 
If the opposite inequality holds, there are only finitely many.  
\end{example}
            
\begin{proof}  Take $K = \QQ$, $\cC = \PP^1$, and $\fX = \{0,\infty\}$.  
Let the uniformizers be $g_0(z) = z$, $g_{\infty}(z) = 1/z$ as before.  
Take $r \ge 1$ and put $E_{\infty} = [-r,-1/r] \cup [1/r,r] \subset \RR$. 
For each $q \in \cQ$, put $E_q = \ZZ_q^{\times}$.  For all other primes $p$, 
put $E_p = D(0,1) \backslash D(0,1)^- \subset \CC_p$, 
then let $\EE = E_{\infty} \times \prod_p E_p$.              

By formulas (\ref{F2Seg2}), (\ref{F2Seg3}) and (\ref{F2Seg5}), we have 
\begin{equation*}
\Gamma(E_{\infty},\fX) 
          \ = \ \left( \begin{array}{cc} 
               \frac{1}{2} \log(\frac{4r^2}{r^4-1}) & \frac{1}{2} \log(\frac{r^2+1}{r^2-1}) \\ 
               \frac{1}{2} \log(\frac{r^2+1}{r^2-1}) & \frac{1}{2} \log(\frac{4r^2}{r^4-1})
                       \end{array} \right) \ .
\end{equation*}\index{Green's function!examples!archimedean}                    
For primes $q \in \cQ$, formulas (\ref{FOwCrossGreen}) and  (\ref{FOwCrossV}) 
give $V_{\infty}(E_q) = G(0,\infty;E_q) = q/(q-1)^2$.  Pulling back by $1/z$ and using 
\index{Robin constant!examples!nonarchimedean}
that $\ZZ_p^{\times}$ is stable under taking reciprocals, we have 
$G(z,0;E_q) = G(1/z,\infty;E_q)$ and hence $V_{0}(E_q) = G(\infty,0;E_q) = q/(q-1)^2$
\index{Green's function!examples!nonarchimedean}
as well.  Thus 
\begin{equation} \label{FqinQ} 
\Gamma(E_q,\fX) 
          \ = \ \left( \begin{array}{cc} 
               q/(q-1)^2 & q/(q-1)^2 \\ 
               q/(q-1)^2 & q/(q-1)^2 \end{array} \right) \ .
\end{equation}  
For all other $p$, $\Gamma(E_p,\fX)$ is the $0$ matrix.  
\index{Green's matrix!local}
Hence 
\begin{eqnarray*}                    
\Gamma(\EE,\fX) & = & \Gamma(E_{\infty},\fX) + \sum_p \Gamma(E_p,\fX) \log(p) \\
      & = & \frac{1}{2}\left( \begin{array}{cc} 
             \log(\frac{4A^2r^2}{r^4-1}) &  \log( A^2 \frac{r^2+1}{r^2-1}) \\ 
             \log( A^2 \frac{r^2+1}{r^2-1}) & \log(\frac{4A^2 r^2}{r^4-1})
                       \end{array} \right) \ .
\end{eqnarray*}
\index{Green's matrix!global}
Take $\vs = {}^t(\frac{1}{2},\frac{1}{2})$.  Then $\Gamma(\EE,\fX)\vs$ has equal entries  
\begin{equation*}
V \ = \ \frac{1}{2} \log \left( \frac{4 A^4 r^2}{(r^2-1)^2} \right) \ .  
\end{equation*} 
By the definition of the value of a matrix\index{value of $\Gamma$ as a matrix game}
game, it follows that $V(\EE,\fX) := \val(\Gamma(\EE,\fX)) = V$.  
Since $r \ge 1$, $\gamma(\EE,\fX) = e^{-V(\EE,\fX)} = (r^2-1)/(2A^2r)$. 
It is easy to see that $\gamma(\EE,\fX) > 1$ if and only if condition (\ref{FUnit2}) holds, 
and that $\gamma(\EE,\fX) < 1$ if and only if the opposite inequality holds.
Hence the result follows from the Theorem \ref{FSZii}.
\end{proof}  

\smallskip
If $r = 1+\sqrt{2}$ in the first part of Example \ref{UnitThm2}, 
then there are infinitely many units whose conjugates
lie in $[-r,-1/r]\cup[1/r,r]$.  Note that this set is the pullback of $[-2,2]$
by  $f(z) = z-1/z$.  For each $n \ge 1$, let $T_n(x)$ denote the
Chebyshev polynomial of degree $n$.  It is well known that $T_n(x)$ is a monic polynomial 
\index{Chebyshev polynomial} 
with integer coefficients, whose roots are simple and belong to the interval $[-2,2]$.
\index{algebraic integer} 
Put $P_n(z) = z^n T_n(z-1/z)$.  Then $P_n(z)$ is monic with integer coefficients, and has
constant coefficient $(-1)^n$.  Thus the roots of the $P_n(z)$ are the units we need. 

\smallskip
Next we give an $S$-unit analogue of Example \ref{UnitThm2}.
By a trick, 
we are able to require that the $S$-units constructed be totally $p$-adic, 
while their archimedean conjugates all have absolute value $1$:

\begin{example}\index{examples!global!$S$-unit analogue} \label{UnitThm3} 
Let $k = \QQ$ and fix a $($nonarchimedean$)$ prime $p$.   
Let $\cQ$ be a finite set of nonarchimedean primes of $\QQ$, disjoint from $\{p\}$,
and put $A = \prod_{q \in \cQ} q^{q/(q-1)^2}$ as in Example $\ref{UnitThm2}$.
Suppose $0 < m \in \ZZ$ is such that
\begin{equation} \label{FmCond}
m \log(p) \ > \ 2 \log(A) + \frac{p \, (p^{2m}+1)}{(p-1)(p^{2m+1}-1)} \log(p) \ .
\end{equation}
Then there are infinitely many numbers\, $\alpha \in \widetilde{\QQ}$ 
for which the primes\, $q \in \cQ$ split completely in $\QQ(\alpha)$, 
which are units at all nonarchimedean places $v$ of\, $\QQ(\alpha)$ not above $p$,  
whose archimedean conjugates all satisfy $|\sigma(\alpha)| = 1$, 
and whose conjugates in $\CC_p$ all belong to $\QQ_p$ 
and satisfy $|\ord_p(\sigma(\alpha))| \le m$.  

If the opposite inequality to $(\ref{FmCond})$ holds, there are only finitely many.
\end{example} 

\begin{proof}  Take $K = \QQ$, and let $\fX = \{0,\infty\}$.  Let the uniformizing parameters
be $g_0(z) = z$, $g_{\infty}(z) = 1/z$ as usual.\index{uniformizing parameter!normalizes Robin constant} 

The proof makes use of two $\QQ$-rational adelic sets, 
which we will denote $\EE$ and $\EE^{\prime}$.  
To construct $\EE$, let $E_{\infty} = C(0,1)$, the unit circle.  
For each $q \in \cQ$, put $E_q = \ZZ_q^{\times}$, and put 
\begin{equation*}
E_p \ = \ \{ x \in \QQ_p : -m \le \ord_p(x) \le m \} \ = \ p^{-m} E_{p,2m}
\end{equation*}
where $E_{p,2m} = \bigcup_{k=0}^{2m} p^k \ZZ_p^{\times}$ is as in Proposition \ref{FmCosets}.
For all other finite primes $q$ take $E_q = \hcO_q^{\times} = D(0,1) \backslash D(0,1)^-$, 
the $\fX$-trivial set in $\PP^1(\CC_q)$.   
\index{$\fX$-trivial}
Set $\EE = E_{\infty} \times \prod_{q \ne \infty} E_q$.   

To construct $\EE^{\prime}$, first choose a square-free integer $d < 0$ which satisfies
$d \equiv 1 \pmod{8}$ and $d \equiv 1 \pmod{q}$ for each $q \in \cQ \cup \{p\}$.  
Thus the primes in $\cQ \cup \{p\}$ split completely in the quadratic imaginary field 
$F = \QQ(\sqrt{d})$.  Let
\begin{equation*}
f(x) \ = \ \frac{x - \sqrt{d}}{x+\sqrt{d}} \ ,
\end{equation*}
and for each prime $q$ (archimedean or nonarchimedean) put $E_q^{\prime} = f^{-1}(E_q)$.  
Then $E_{\infty}^{\prime} = \PP^1(\RR)$, 
while for each $q \in \cQ \cup \{p\}$ we have $E_q^{\prime} \subset \QQ_q$.  
For all other primes $q$, 
$E_q^{\prime}$ is the $\RL$-domain in $\PP^1(\CC_q)$ gotten by omitting two open discs 
\index{$\RL$-domain} 
centered on $\pm \sqrt{d}$;  for all but finitely many $q$ these discs are disjoint 
and have radius $1$.  Note that for each $q$, the set $E_q^{\prime}$  
is stable under $\Aut_c(\CC_q/\QQ_q)$.  If $q \in \cQ \cup \{p,\infty\}$ this is trivial;  
for all other $q$, note that for each $\sigma \in \Aut_c(\CC_q/\QQ_q)$, either 
$\sigma(f)(x) = f(x)$ or $\sigma(f)(x) = 1/f(x)$.  
Since $E_q = \hcO_q^{\times}$  is stable under inversion and $\Aut_c(\CC_q/\QQ_q)$, 
it follows that $x \in E_q^{\prime}$ if and only if $\sigma(x) \in E_q^{\prime}$.
     
Set $\EE^{\prime} = E_{\infty}^{\prime} \times \prod_{q \ne \infty} E_q^{\prime}$,
and take $\fX^{\prime} = \{\sqrt{d},-\sqrt{d}\}$.   We claim that
\begin{equation*}
\Gamma(\EE^{\prime},\fX^{\prime}) \ = \ \Gamma(\EE,\fX) \ .
\end{equation*}\index{Green's matrix!global}
This follows by pulling back using $f(x)$:  let $\EE_F$, $\EE_F^{\prime}$
be the $F$-rational adelic sets obtained from $\EE$, $\EE^{\prime}$ by base change 
(see \cite{RR1}, \S5.1).  Then $\Gamma(\EE_F,\fX) = [F:\QQ]\cdot \Gamma(\EE,\fX)$
and $\Gamma(\EE_F^{\prime},\fX^{\prime}) = [F:\QQ] \cdot \Gamma(\EE^{\prime},\fX^{\prime})$ 
(\cite{RR1}, p.326, formula (9)).  On the other hand $f(x)$ is rational over $F$, 
so by (\cite{RR1}, p.335, formula (16)) and the fact that $\deg(f) = 1$,  
\begin{equation*} 
\Gamma(\EE_F^{\prime},\fX^{\prime}) \ = \  \Gamma(\EE_F,\fX) \ .
\end{equation*}  
This establishes the claim. 

By the Fekete-Szeg\"o theorem \ref{aT1}, if $\Gamma(\EE^{\prime},\fX^{\prime})$ 
\index{Fekete-Szeg\"o theorem with LRC}\index{Green's matrix!negative definite} 
is negative definite, there are infinitely many algebraic numbers whose archimdean 
conjugates belong to $E_{\infty}^{\prime} = \PP^1(\RR)$ and 
whose $q$-adic conjugates belong to $E_q^{\prime}$ for all nonarchimedean $q$.  
The images of these numbers under $f(x)$ will be the numbers $\alpha$ in the Example.

Hence it suffices to show that $\Gamma(\EE^{\prime},\fX^{\prime}) = \Gamma(\EE,\fX)$ 
\index{Green's matrix!global}\index{Green's matrix!negative definite}
is negative definite under condition (\ref{FmCond}).  We have 
\begin{equation*} 
\Gamma(\EE,\fX) \ = \ \Gamma(E_{\infty},\fX) + \sum_{q \ne \infty} \Gamma(E_q,\fX) \log(q) \ .
\end{equation*}
\index{Green's matrix!global}\index{Green's matrix!negative definite} 
Here $\Gamma(E_{\infty},\fX)$ is the $0$ matrix.  For each $q \in \cQ$, $\Gamma(E_q,\fX)$ is 
\index{Green's matrix!local}
the same as in (\ref{FqinQ}) in the proof of Example \ref{UnitThm2}.  To express 
$\Gamma(E_p,\fX)$, write 
\index{Green's matrix!local}
\begin{eqnarray*} 
B & = & \frac{1}{p-1} + \frac{1}{(p-1)^2(1 + p^2 + p^4 + \cdots + p^{4m})} \ , \\
C & = & \frac{p^{2m+1}}{(p-1)^2(1 + p^2 + p^4 + \cdots + p^{4m})} \ .
\end{eqnarray*}
By Proposition \ref{FmCosets} and the scaling property of the capacity,
\index{capacity!scaling property of}
$V_{\infty}(E_p) = V_{\infty}(p^{-m}E_{p,2m}) = -m + B$.  The map $f(z) = 1/z$
\index{Robin constant!examples!nonarchimedean}
stabilizes $E_p$ but takes $0$ to $\infty$, 
so our choice of uniformizing parameters\index{uniformizing parameter!normalizes Robin constant} gives 
$V_0(E_p) = V_{\infty}(E_p)$.  Finally  $f(z) = p^m z$ takes 
$0 \mapsto 0$, $\infty \mapsto \infty$, and $E_p$ to $E_{p,2m}$, 
so by the pullback formula (\ref{GreenPullbackF}) and Proposition (\ref{FmCosets}),
\index{Green's function!pullback formula for}
$G(0,\infty;E_p) = G(0,\infty;E_{p,2m}) = C$.   Thus 
\begin{equation*}
\Gamma(E_p,\fX)  \ = \ \left( \begin{array}{cc} 
                                     -m+B & C \\
                                     C & -m+B 
                                 \end{array} \right) \ .
\end{equation*}
For all other primes $q$, $\Gamma(E_q,\fX)$ is the $0$ matrix, so
\index{Green's matrix!local}
\begin{equation*} 
\Gamma(\EE,\fX) \ = \ \left( \begin{array}{cc} 
                    \log(A) + (-m+B) \log(p) & \log(A) + C \log(p) \\
                    \log(A) + C \log(p) & \log(A) + (-m+B) \log(p)
                                 \end{array} \right) \ .
\end{equation*}   
\index{Green's matrix!global}
              
Since $\Gamma(\EE,\fX)$ has equal row sums, as in the proof of Example \ref{UnitThm2} it follows
that $\val(\Gamma(\EE,\fX)) = \frac{1}{2}\big(-m \log(p) + 2\log(A) + (B+C) \log(p)\big)$. 
Simplifying, we have $\val(\Gamma(\EE,\fX)) < 0$, and hence 
\index{Green's matrix!global}
$\gamma(\EE,\fX) > 1$, if and only if (\ref{FmCond}) holds; similarly, $\gamma(\EE,\fX) < 1$ if and only
if the opposite inequality holds.  Thus the result follows from the Fekete-Szego theorem \ref{FSZii}.
\index{Fekete-Szeg\"o theorem with LRC} 
\end{proof}

For instance, if $\cQ = \{2\}$ in Example \ref{UnitThm3}, then for $p =3$ we need $m \ge 4$;
for $5 \le p \le 17$ we can take $m = 2$, and for $p \ge 19$ we can take $m = 1$.  
If $\cQ = \{2,3\}$ in Example \ref{UnitThm3}, then for $p =5$ we need $m \ge 7$;  
for $7 \le p \le 11$ we  can take $m = 5$;  for $13 \le p \le 23$ we can take $m = 4$;
for $29 \le p \le 109$ we can take $m = 3$; for $113 \le p \le 11673$ we can take $m = 2$;
and for $p \ge 11677$ we can take $m = 1$.  

Note that the $S$-units constructed in Example \ref{UnitThm3} are not roots of unity, 
because there are only finitely many roots of unity $\zeta_n$ for which $p$ splits completely
in $\QQ(\zeta_n)$.

\smallskip
In the next example, a limit argument allows us deal with a situation where one 
of the points in $\fX$ belongs to $E_{\infty}$:  

\begin{example}\index{examples!global!example with $E_{\infty} \cap \fX \ne \phi$} \label{AUnitThm}  
Let $A > 0$.  If $A \ge 4$, then there are infinitely many units whose conjugates
all lie in $[0,1] \cup [A,A+1] \subset \RR$.  
If $A < 4$, there are only finitely many.   
\end{example}

\begin{proof}  Write $E = [0,1] \cup [A,A+1]$.  
If $A \le 1$, then $E \subset [0,2]$, so $\gamma_{\infty}(E) < 1/2$.  
Otherwise, $E$ is a translate of $[-(A+1)/2,-(A-1)/2] \cup [(A-1)/2,(A+1)/2]$ 
and so $\gamma_{\infty}(E) = \sqrt{A}/2$ by formula (\ref{F2Seg2}).  
Thus if $A < 4$, we have $\gamma_{\infty}(E) < 1$, and Fekete's theorem \ref{FSZii}(B) shows 
\index{Fekete-Szeg\"o theorem with LRC} 
there are only finitely many  {\em algebraic integers}, and in particular
finitely many units, whose conjugates lie in $E$.  

\smallskip
\vskip .05 in
Next suppose $A = 4$. We will explicitly construct infinitely many units whose conjugates
lie in $[0,1] \cup [4,5]$.  
To do so, note that $[0,1] \cup [4,5]$ is the pullback of $[-2,2]$
by  $f(z) = z^2 - 5 z + 2$. Let $T_n(x)$ denote the
Chebyshev polynomial of degree $n$.  As before, $T_n(x)$ is a monic polynomial 
\index{Chebyshev polynomial} 
with integer coefficients whose roots are simple and belong to $[-2,2]$.
It oscillates $n$ times between $\pm 2$ on $[-2,2]$;   in particular,  $T_n(2) = 2$.  
Furthermore, $T_n(x)$ is an even function if $n$ is even, 
and is an odd function if $n$ is odd.  Consider the polynomials $Q_n(z) = T_n(f(z))$.  
They are monic with integer coefficients, and have all their roots in $[0,1] \cup [4,5]$.
Unfortunately, $Q_n(z)$ has constant coefficient $Q_n(0) = T_n(2) = 2$.  However, if $n$ is odd,
then $Q_n(z)$ has $f(z)$ as a factor, so $P_n(z) := Q_n(z)/f(z)$ has constant coefficient $1$.  
Thus the roots of the $P_n(z)$ for odd $n$ are the required units.   

\smallskip
Finally, suppose $A > 4$, and let $0 < \varepsilon < 1$.  
Consider the set $E_{\varepsilon,A} = [\varepsilon,1] \cup [A,A+1]$.  
Take $K = \QQ$, $\cC = \PP^1$, $\fX = \{0,\infty\}$.  Let $E_{\infty} = E_{\varepsilon,A}$,
and for each finite prime $p$ let $E_p = D(0,1) \backslash D(0,1)^- \subset \CC_p$
be the $\fX$-trivial set.  Put $\EE_{\varepsilon,A} = E_{\varepsilon,A} \times \prod_p E_p$,
\index{$\fX$-trivial} 
and take $g_0(z) = z$, $g_{\infty}(z) = 1/z$ as before.  Then
\index{Robin constant!examples!nonarchimedean}
\begin{equation*}
\Gamma(\EE_{\varepsilon,A},\fX) \ = \   
        \left( \begin{array}{cc} 
                 V_0(E_{\varepsilon,A})  &  G(\infty,0;E_{\varepsilon,A}) \\
                 G(0,\infty;E_{\varepsilon,A}) & V_{\infty}(E_{\varepsilon,A}) 
               \end{array} \right) \ .
\end{equation*}  
Formula (\ref{FGreen}) expresses the Green's function of two intervals in terms of a quotient of
\index{Green's function!examples!archimedean} 
two theta functions.  These theta-functions\index{theta-functions!classical} 
and their parameters vary continuously with 
$\varepsilon$ and $A$, hence the Green's function varies continuously as well.  
Letting $\varepsilon \rightarrow 0$, formulas (\ref{F2Seg2}), (\ref{FGreen}), 
and (\ref{F2SegV2}) show that 
\index{Robin constant!examples!archimedean}
\begin{eqnarray*}
\lim_{\varepsilon \rightarrow 0} V_{\infty}(E_{\varepsilon,A}) & = & V_{\infty}(E) 
\ = \ -\frac{1}{2}\log(A/4) \ < \ 0 \ , \\ 
\lim_{\varepsilon \rightarrow 0} G(\infty,0;E_{\varepsilon,A}) & = & 
    \lim_{\varepsilon \rightarrow 0} G(0,\infty;E_{\varepsilon,A})  
      \ = \ G(0,\infty;E) \ = \ 0 \ , \\
\lim_{\varepsilon \rightarrow 0} V_{0}(E_{\varepsilon,A}) & = & -\infty \ . 
\end{eqnarray*} 
Thus for all sufficiently small $\varepsilon > 0$, $\Gamma(\EE_{\varepsilon,A},\fX)$ is negative
definite, and the Fekete-Szeg\"o Theorems \ref{aT1} and \ref{FSZii} yield the result.
\index{Fekete-Szeg\"o theorem with LRC} 
\end{proof}

As a whimsical side note, we remark that an argument similar to the one in 
Example \ref{AUnitThm} shows that $E = [0,1] \cup [A,A+.001]$ 
contains infinitely many conjugate sets of units if $A \ge 30.19249489$, 
but only finitely many if  $0 < A < 30.19249488$.  
This is obtained by using Maple\index{Maple computations} 
to evaluate $V_{\infty}(E)$ in formula (\ref{F2SegV1}).
\index{Robin constant!examples!archimedean}

\medskip
It is also possible to use the Fekete-Szeg\"o theorem to construct units 
\index{Fekete-Szeg\"o theorem with LRC} 
whose conjugates {\em globally omit} residue classes.  
View $\PP^1/\QQ$ as the generic fibre of $\PP^1_{\ZZ}/\Spec(\ZZ)$.  
Given $\alpha, \beta \in \PP^1(\widetilde{\QQ})$, 
we will say that $\alpha$ is integral with respect to $\beta$ if the Zariski closures of 
$\alpha$ and $\beta$ in $\PP^1_{\ZZ}$ do not meet.  If $\beta_1, \ldots, \beta_N$ 
are the conjugates of $\beta$, this is equivalent to requiring that for
every prime $p$, all the conjugates $\sigma(\alpha)$ in $\PP^1(\CC_p)$ belong to   
$\PP^1(\CC_p) \backslash (\bigcup_{i=1}^N B(\beta_i,1)^-)$. 

Recall that the (absolute, logarithmic) Weil height of a number $\alpha \in \widetilde{\QQ}$ is 
\index{Weil!height} 
\begin{equation*}
h(\alpha) \ = \ \frac{1}{[F:\QQ]} \sum_{\text{$w$ of $F$}} \log_w^+(|\alpha|_w) \log(q_w) \ , 
\end{equation*}
where $F$ is any finite extension $\QQ(\alpha)$.  The height is independent of the field
used to compute it.  By definition, $h(\infty) = 0$.  The points of $\PP^1(\widetilde{\QQ})$ 
with $h(\alpha) = 0$ are precisely $0$, $\infty$, and the roots of unity.  

\smallskip
To put the following result in context, we remark that in (\cite{BIR}) the authors show that
if $h(\beta) \ne 0$, there are only finitely many roots of unity which are integral with 
respect to $\beta$.   

\begin{example}\index{examples!global!regarding Ih's conjecture} \label{OmitThm}
Take $K = \QQ$, and let $\beta \in \widetilde{\QQ}$.  
Then there are infinitely many algebraic units $\eta \in \QQbar$
which are integral with respect to $\beta$.   
For any $\varepsilon > 0$, these units can be required to have Weil height  
\index{Weil!height}  
$h(\eta) < \varepsilon$.
\end{example}  

\begin{proof}
Put $L = \QQ(\beta_1, \ldots, \beta_N)$,
where $\beta_1, \ldots, \beta_N$ are the conjugates of $\beta$ over $\QQ$,
and take $\fX = \{\infty, 0, \beta_1,\ldots, \beta_N\}$.  
Let $g_{\infty}(z) = 1/z$, $g_0(z) = z$, and $g_{\beta_i}(z) = z-\beta_i$ for each $i$. 
In the discussion below, we will assume that $\beta \ne 0$.  
If $\beta = 0$, then we are merely asking for units with height $h(\eta) < \varepsilon$, 
and the argument carries through in a simplified form. 

Viewing the $\beta_i$ as embedded in $\CC$, 
let $r > 1$ be any number small enough that $r < \min_{|\beta_i| > 1}(|\beta_i|)$, 
and then let $\rho > 0$ be any number small enough that the discs 
$D(\beta_i, \rho)$ for $i =1, \ldots, N$ and $D(0,\delta)$
are pairwise disjoint and do not meet the circles $\{|z| = r\}$.  
(Note that $r$ and $\rho$ are independent of choice of the embedding of the $\beta_i$.).  
Eventually we will let $\rho \rightarrow 0$,
and then let $r \rightarrow 1$.  Put 
\begin{equation*}
E_{\infty} \ = \ \Big( D(0,r) \cup \big(\bigcup_{|\beta_i| > r} C(\beta_i,\rho)\big)\Big) \backslash 
\Big( D(0,\rho)^- \cup \bigcup_{|\beta_i| < r} D(\beta_i,\rho)^-\Big) \ .
\end{equation*}
Thus, $E_{\infty}$ consists of a disc $D(0,r)$ with tiny holes deleted around $0$ and the 
$\beta_i \in D(0,r)$, together with tiny circles adjoined around the 
$\beta_i \notin D(0,r)$.  By construction $E_{\infty}$ 
is stable under complex conjugation. 

For each finite prime $p$, regarding the $\beta_i$ as embedded in $\CC_p$, put 
\begin{equation*} 
E_p \ = \ \PP^1(\CC_p) \backslash 
      \big(B(\infty,1)^- \cup B(0,1)^- \cup \bigcup_{i=1}^N B(\beta_i,1)^-\big) \ .
\end{equation*}   
Then $E_p$ is an $\RL$-domain, stable under $\Aut_c(\CC_p/\QQ_p)$, and for all but finitely many
\index{$\RL$-domain} 
$p$ \, it is $\fX$-trivial.
\index{$\fX$-trivial}

Put $\EE = E_{\infty} \times \prod_p E_p$.  To compute $\Gamma(\EE,\fX)$,
\index{Green's matrix!global} 
we must first make a base change to the field $L$, 
over which the $\beta_i$ are rational. By definition
\begin{equation} \label{FBaseChange}
\Gamma(\EE,\fX) \ = \ \frac{1}{[L:K]} \, \Gamma(\EE_L,\fX) 
\ = \ \frac{1}{[L:K]} \, \sum_{\text{places $w$ of $L$}} \Gamma(E_w,\fX) \log(q_w) 
\end{equation}
\index{Green's matrix!global}
where $\EE_L = \prod_{\text{$w$ of $L$}} E_w$.  
Here, for each place $w$ of $L$, if $w$ lies over $p$, 
then $E_w$ is gotten by choosing an embedding $\sigma : L \hookrightarrow \CC_p$ which induces $v$, 
extending $\sigma$ to an isomorphism $\overline{\sigma} : \CC_w \rightarrow \CC_p$, 
and setting $E_w = \overline{\sigma}^{-1}(E_p)$.  
Basically, $E_w$ is the same as $E_p$,   
but the way the $\beta_i$ are embedded depends $w$. 
 
We now compute the matrices $\Gamma(\EE_w,\fX)$.  
First suppose $w|\infty$.  By construction, each point of $\fX$ belongs to 
a different connected component of $\PP^1(\CC_w) \backslash E_w$, 
so $\Gamma(E_w,\fX)$ is a diagonal matrix.  
\index{Green's matrix!local}
\index{Robin constant!examples!archimedean}
We have $V_{\infty}(E_w) = -\log(r)-\delta(\rho)$ where $\delta(\rho) > 0$ 
and $\delta(\rho) \rightarrow 0$ as $\rho \rightarrow 0$, 
while $V_{\beta_i}(E_w) = -\log(\rho)$ for each $i$.

Next suppose $w$ is nonarchimedean.  Since $E_w$ is obtained by deleting a finite
number of open discs of radius $1$ from $\PP^1(\CC_w)$, 
one of which is $B(\infty,1)^-$, we have  $V_{\infty}(E_w) = V_{\infty}(D(0,1)) = 0$. 
\index{Robin constant!examples!nonarchimedean}
The other entries of $\Gamma(E_w,\fX)$ will not matter to us:
\index{Green's matrix!local}  
$\Gamma(E_w,\fX)$ is an $(N+2)\times(N+2)$ matrix whose $V_{\infty}(E_w)$ entry is $0$.
For all but finitely many $w$, $E_w$ is $\fX$-trivial and $\Gamma(E_w,\fX)$ is the $0$ matrix.
\index{Green's matrix!local}
\index{$\fX$-trivial}

By definition $\Gamma(\EE_F,\fX) = \sum_w \Gamma(E_w,\fX) \log(q_w)$;  
for archimedean $w$, $q_w = e$ if $L_w \cong \RR$, 
while $q_v = e^2$ if $L_w \cong \CC$. By (\ref{FBaseChange})  
\begin{equation*}
\Gamma(\EE,\fX) \ = \ \left( \begin{array}{cccc} 
                        -\log(r)-\delta(\rho) & A_{12} & \cdots & A_{1,N+2} \\
              
                        A_{21} & -\log(\rho) + A_{22} & \cdots & A_{2,N+2}   \\
                        \vdots & \vdots & \ddots & \vdots \\
                        A_{N+2,1} & A_{N+1,2} & \cdots & -\log(\rho) + A_{N+2,N+2} 
                              \end{array} \right) 
\end{equation*}\index{Green's matrix!global}                            
where the $A_{ij}$ do not depend on $r$ or $\rho$, and $A_{ij} = A_{ji}$ for all $i, j$.
By the determinant criterion for negative definiteness from linear algebra 
(see for example \cite{RR1}, Proposition 5.1.8, p.331),\index{Determinant Criterion!for negative definiteness} 
for each fixed $r$  
if $\rho$ is sufficiently small then $\Gamma(\EE,\fX)$ is negative definite.  Thus
\index{Green's matrix!global}\index{Green's matrix!negative definite}
for any neighborhood $U$ of $E_{\infty}$ the Fekete-Szeg\"o theorem \ref{FSZii} 
\index{Fekete-Szeg\"o theorem with LRC} 
produces infinitely many units $\eta$ whose archimedean conjugates all lie in $U$, 
and whose nonarchimedean conjugates avoid the balls $B(\beta_i,1)^-$ 
at all places $w$ of $L$.

To see why the numbers $\eta$ can be assumed to have arbitrarily small height
requires some understanding of the proof of the Fekete-Szeg\"o theorem 
\index{Fekete-Szeg\"o theorem} 
(either Theorem 6.3.2 of \cite{RR1}, or Theorem \ref{aT1-B} in this work).
We will now sketch the argument, assuming the reader is loosely familiar with the proof.  

Fix $r > 1$, and let $\rho$ be small enough that 
$\Gamma(\EE,\fX)$ is negative definite.\index{Green's matrix!negative definite}  Then there is a probability vector
\index{Green's matrix!global}
$\vs = {}^t(s_1, \ldots, s_{N+2})$ for which $\Gamma(\EE,\fX) \vs$ has
\index{Green's matrix!global} 
all its entries equal to $V(\EE,\fX)$.  
These $s_i$ are essentially the relative orders of the poles
of the initial patching functions $G_v^{(0)}(z)$ at the points $x_i \in \fX$.  
\index{patching functions, initial $G_v^{(0)}(z)$}
As $\rho \rightarrow 0$, we will have $s_1 \rightarrow 1$ and 
$s_2, \ldots, s_{N+2} \rightarrow 0$  since the first row of $\Gamma(\EE,\fX)$ 
\index{Green's matrix!global}
(and hence $V(\EE,\fX)$) remains bounded but the diagonal entries 
in the other rows approach $\infty$.  

The archimedean initial local patching function $G_{\infty}^{(0)}(z)$ 
\index{patching functions, initial $G_v^{(0)}(z)$}
is chosen so that the discrete probability density of its zeros approximates 
$\sum_{i=1}^{N+2} s_i \mu_i$, 
where $\mu_i$ is the equilibrium distribution of $E_v$
with respect to $x_i$.  
Here each $\mu_i$ is a probability measure supported on the boundary of the
\index{boundary!equilibrium distribution supported on}
component of $\PP^1(\CC_v) \backslash E_v$ containing $x_i$.  
As $\rho \rightarrow 0$, 
the amount of mass which $\mu_1$ (corresponding to $x_1 = \infty$) places on the 
circles $C(\beta_i,\rho)$ goes to $0$. Thus the proportion of the zeros 
of $G_{\infty}^{(0)}(z)$ which lie near near $C(0,r)$ goes to $1$.  The remaining 
zeros all lie near the circles $C(\beta_i,\rho)$.  
If $U$ is chosen small enough that each $C(\beta_i,\rho)$ outside $D(0,r)$ 
lies in a separate component of $U$, then the patching process preserves
\index{patching argument!local} 
the number of zeros which lie in each component.  
Thus the final patched function $G^{(n)}(z)$, 
whose zeros are numbers constructed by the Fekete-Szeg\"o theorem,
\index{Fekete-Szeg\"o theorem with LRC}  
has the same number of zeros in each component of $U$ as the initial function $G_{\infty}^{(0)}(z)$.  

Since the zeros of $G^{(n)}(z)$ (in our instance) are algebraic units, 
the only contribution to their height is from archimedean places. 
By the discussion above, that contribution approaches $\log(r)$ 
as $\rho \rightarrow 0$.  So, if we first let $\rho \rightarrow 0$,
and then let $r \rightarrow 1$, the Fekete-Szeg\"o theorem produces 
\index{Fekete-Szeg\"o theorem with LRC} 
numbers whose heights approach $0$.  
\end{proof} 

\smallskip
Our final example constructs units which avoid the residue class of $1$ at every prime, 
and whose archimedean conjugates all lie very close to the circle $C(0,r)$ or the circle $C(0,1/r)$
(so $|\log(|\sigma(\alpha)|)| \approx \log(r)$), for suitable $r$.     

\begin{example}\index{examples!global!units near circles} \label{ThmOmit1}
Let $r$ satisfy $1 < r < 2.96605206$.  Then for any $\varepsilon > 0$, 
there are infinitely many units $\alpha$ whose conjugates all satisfy 
\begin{equation*}
||\sigma(\alpha)| - r| \ < \ \varepsilon \quad \text{or} \quad  
||\sigma(\alpha)| - 1/r| \ < \ \varepsilon \ ,
\end{equation*}
and are such that $\sigma(\alpha) \not\equiv 1 \pmod{\fp}$ 
for each prime $\fp$ of $\cO_{\QQ(\sigma(\alpha))}$.   If $r > 2.96605207$,
there are only finitely many.
\end{example}             
          
\begin{proof} 
Take $K = \QQ$, $\cC = \PP^1$, and $\fX = \{0,1,\infty\}$.  
Let $E_{\infty} = C(0,r) \cup C(0,1/r)$,
and for each finite prime let $E_p$ be the $\fX$-trivial set 
\index{$\fX$-trivial}
\begin{equation*}
E_p \ = \ \PP^1(\CC_p) \backslash (B(0,1)^- \cup B(1,1)^- \cup B(\infty,1)^-) \ .
\end{equation*} 
Put $\EE = E_{\infty} \times \prod_p E_p$, and take $g_0(z) = z$, $g_1(z) = z-1$,
$g_{\infty}(z) = 1/z$.

Note that $0$, $1$ and $\infty$ belong to different components of
$\PP^1(\CC) \backslash E_{\infty}$.  
Then $V_{\infty}(E_{\infty}) = V_0(E_{\infty}) = -\log(r)$ by formula (\ref{FDisc2}), 
\index{Robin constant!examples!archimedean}
while $V_1(E_{\infty})$ is given by (\ref{FRingV1}) with $\tau = 2 i \log(r)/\pi$.  
At each nonarchimedean place, $\Gamma(E_p,\fX)$ is the $0$ matrix. 
\index{Green's matrix!local}  
Hence $\Gamma(\EE,\fX) = \Gamma(E_{\infty},\fX)$ is the diagonal matrix
\index{Green's matrix!global}
\begin{equation} \label{FG3}
\Gamma(\EE,\fX) \ = \ \left( \begin{array}{ccc} 
                               -\log(r) & 0 & 0 \\
       0 & -\log(\frac{|\theta(0,\tau;0,0)\theta(0,\tau;\frac{1}{2},0)|}{2}) & 0 \\
                               0 & 0 & -\log(r) 
                             \end{array} \right)
\end{equation}
\index{Green's matrix!global}
Clearly $\Gamma(\EE,\fX)$ is negative definite\index{Green's matrix!negative definite} if and only if 
the middle term is negative.  A computation with Maple\index{Maple computations} yields the result.  
\end{proof}          
  
\section{ Function Field Examples concerning Separability} \label{FFSeparabilitySection}

In this section, we will take $K = \FF_p(t)$ where $p$ is a prime, 
$\FF_p$ is the finite field with $p$ elements, 
and $t$ is transcendental over $\FF_p$.  
We give three examples showing the need for the separability
hypothesis in Theorem \ref{aT1}.C.2, and in Theorems \ref{aT1-A1}, \ref{aT1-A}, and \ref{FSZii}. 
This was discovered by Daeshik Park in his doctoral thesis (\cite{DPark}).
\index{Park, Daeshik}
\index{Fekete-Szeg\"o theorem with LRC!need for separability hypotheses in}
\index{examples!for function fields!concerning separability}

\smallskip
Let $v_0$, $v_1$ and $v_{\infty}$ be the valuations of $\FF_p(t)$ for which $v_0(t) = 1$, $v_1(t-1) = 1$,
and $v_{\infty}(1/t) = 1$, respectively.  
For each of the corresponding places, the residue field is $\FF_p$, 
and we have $\cO_{v_0} \cong \FF_p[[t]]$,\label{`SymbolIndexFpsbt'} $\cO_{v_1} \cong \FF_p[[t-1]]$, 
and $\cO_{v_\infty} \cong \FF_p[[\frac{1}{t}]]$.  

\medskip
Our first example, which is due to Park, concerns a set where all the 
\index{Park, Daeshik}
hypotheses of the Fekete-Szeg\"o theorem \ref{aT1} are satisfied except for 
\index{Fekete-Szeg\"o theorem with LRC} 
separability of the extension $F_{w_0}/K_{v_0}$, yet the conclusion of the theorem fails for $r$ 
in a certain range.  

\begin{example} \label{FFSepExample1}  Let $K = \FF_p(t)$, and let $\cC = \PP^1/K$.  Identify 
$\PP^1(\CC_v)$ with $\CC_v \cup \{\infty\}$, and take $\fX = \{\infty\}$.  
Put $F_{w_0} = K_{v_0}(t^{1/p}) = \FF_p((t^{1/p}))$,\label{`SymbolIndexFprbt'}
 so that $\cO_{w_0} = \FF_p[[t^{1/p}]]$
and $F_{w_0}/K_{v_0}$ is purely inseparable.  Fix a place
$v_2 \in \cM_K$ distinct from  $v_0, v_1, v_{\infty}$ and 
define an adelic set $\EE = \EE(r) = \prod_{v \in \cM_K} E_v$ by putting 
$E_{v_0} = \cO_{w_0}$, taking  
$E_{v_1} = \cO_{v_1}$, $E_{v_\infty} = \cO_{v_\infty}$, $E_{v_2} = D(0,r)$,  
and letting $E_v = D(0,1)$ be $\fX$-trivial for all $v \ne v_0, v_1, v_2, v_{\infty}$.  
\index{$\fX$-trivial}
Then $\EE$ is compatible with $\fX$\index{compatible with $\fX$}  and satisfies all the hypotheses of Theorem $\ref{aT1}$ 
apart from the inseparability of $F_{w_0}/K_{v_0}$, and 
\begin{equation*}
\gamma(\EE,\fX) \ = \ r \cdot p^{-\frac{2+1/p}{p-1}}  \ . 
\end{equation*}
However, if 
\begin{equation*} 
p^{\frac{2+1/p}{p-1}} \ < \ r \ < \ p^{\frac{3}{p-1}}  
\end{equation*} 
then $\gamma(\EE,\fX) > 1$, 
yet there are only finitely many numbers in $\tK$ whose conjugates belong to $E_v$ for each $v \in \cM_K$.
\end{example}  

\begin{proof} By Proposition \ref{OwProp} we have 
$\gamma_{\infty}(E_{v_1}) = \gamma_{\infty}(E_{v_\infty}) = p^{-1/(p-1)}$.  
The extension $F_{w_0}/K_{v_0}$ is totally ramified, so by the same Proposition,  
$\gamma_\infty(E_{v_0}) = p^{-1/(p(p-1))}$.  By (\ref{FVNAD1}) it follows that 
\begin{equation*}
\gamma(\EE,\fX) 
\ = \ p^{-\frac{1}{p(p-1)}} \cdot p^{-\frac{1}{p-1}} \cdot p^{-\frac{1}{p-1}} \cdot r 
\ = \ r \cdot p^{-\frac{2+1/p}{p-1}} \ .
\end{equation*}

Suppose $\alpha \in \tK$ has all its conjugates in $E_v$, for each $v \in \cM_K$.
Recall that $K_{v_0} = \FF_p((t))$ and $K_{v_1} = \FF_p((t-1))$ are separable over $K = F_p(t)$ 
(see Grothendieck, \cite{Gro}, EGA IV, 7.8.3ii, or Matsumura \cite{Mat}, Proposition 28.M, p.207).
\index{Grothendieck, Alexander}\index{$K_v$ is separable over $K$} 
Since the conjugates of $\alpha$ in $\CC_{v_1}$ all belong to 
$E_{v_1} = \cO_{v_1}$, $\alpha$ must be separably algebraic over $K$.  
On the other hand each element of $F_{w_0} \backslash K_{v_0}$ 
is purely inseparable over $K_{v_0}$, so the only elements of $\cO_{w_0}$ 
which can be separably algebraic over $K$ are those in $\cO_{v_0}$.  
Thus the conjugates of $\alpha$ in $\CC_{v_0}$ (which \`a priori belong to $E_{v_0} = \cO_{w_0}$)  
must actually belong to $\cO_{v_0}$.  
It follows that the conjugates of $\alpha$ belong to 
\begin{equation*}
\EE^{\prime} \ = \ \cO_{v_0} \times \cO_{v_1} \times \cO_{v_\infty} \times D(0,r) \times 
\prod_{v \ne v_0, v_1, v_2, v_{\infty}} E_v 
\end{equation*} 
whose capacity is $\gamma(\EE^{\prime},\fX)  = r \cdot p^{-\frac{3}{p-1}}$.
\index{capacity!Cantor capacity}

Each local set occurring in $\EE^{\prime}$ is algebraically capacitable, so by Fekete's theorem
\index{Fekete's theorem} 
\index{algebraically capacitable} 
(\cite{RR1}, Theorem 6.3.1, p.414), if $\gamma(\EE^{\prime},\fX) < 1$ there is an adelic
neighborhood $\UU$ of $\EE^\prime$ which contains only finitely many conjugate sets of numbers in 
$\tK$.  In particular, there are only finitely many $\alpha \in \tK$
which have all their conjugates in $\EE^{\prime}$.  

When 
\begin{equation*} 
p^{\frac{2+1/p}{p-1}} \ < \ r \ < \ p^{\frac{3}{p-1}}  
\end{equation*}
we have $\gamma(\EE,\fX) > 1$ but $\gamma(\EE^{\prime},\fX) < 1$, so all the hypotheses of
Theorem \ref{aT1} hold for $\EE$ except for the inseparability of $F_{w_0}/K_{v_0}$, 
yet the conclusion of Theorem \ref{aT1} fails.  
\end{proof}

Our next example provides sets $\EE$ of arbitrarily large capacity, where all the 
\index{capacity!Cantor capacity}
hypotheses of Theorem \ref{aT1} are satisfied except for the separability condition, 
yet there are no $\alpha \in \tK$ with all their conjugates in $\EE$.

\begin{example} \label{FFSepExample2}  Let $K = \FF_p(t)$, and let $\cC = \PP^1/K$.  Identify 
$\PP^1(\CC_v)$ with $\CC_v \cup \{\infty\}$, and take $\fX = \{\infty\}$.  
Again put $F_{w_0} = K_{v_0}(t^{1/p}) = \FF_p((t^{1/p}))$, so that $F_{w_0}/K_{v_0}$
is purely inseparable.  
Fix a place $v_2 \in \cM_K$ distinct from  $v_0$ and $v_1$, and 
define $\EE = \EE(r) = \prod_{v \in \cM_K} E_v$ by putting 
\begin{equation*}
E_{v_0} \ = \ t^{-1/p} + \cO_{w_0} \ = \ B(t^{-1/p},1) \cap F_{w_0} \ , 
\end{equation*} 
taking  $E_{v_1} = \cO_{v_1}$, $E_{v_2} = D(0,r)$,
and letting $E_v = D(0,1)$ be $\fX$-trivial for all $v \ne v_0, v_1, v_2$.  
\index{$\fX$-trivial}
Then $\EE$ is compatible with $\fX${compatible with $\fX$} and satisfies all the hypotheses of Theorem $\ref{aT1}$ 
apart from the inseparability of $F_{w_0}/K_{v_0}$, and  
\begin{equation*}
\gamma(\EE,\fX) \ = \ r \cdot p^{-\frac{1+1/p}{p-1}}  \ . 
\end{equation*}
If $r > p^{(1+1/p)/(p-1)}$
we have $\gamma(\EE,\fX) > 1$, yet there are no $\alpha \in \tK$ 
whose conjugates all belong to $\EE$.
\end{example}  

\begin{proof} The argument is the same as that in Example \ref{FFSepExample1}, except that 
in this case $E_{v_0} \subset F_{w_0} \backslash K_{v_0}$, so there are no $\alpha \in \tK$
whose conjugates in $\CC_{v_1}$ belong to $E_{v_1}$ and whose conjugates in $\CC_{v_0}$ 
belong to $E_{v_0}$.
\end{proof}  

In the previous examples, the conclusion of Theorem \ref{aT1} failed because of interactions
between places of $K$ where the extensions $F_w/K_v$ were separable and inseparable.  
In our last example, the conclusion of Theorem \ref{aT1} fails because of interaction between two places 
where $F_v/K_v$ is inseparable, with different degrees of inseparability. 

\begin{example} \label{FFSepExample3}  Let $K = \FF_p(t)$, and let $\cC = \PP^1/K$.  Identify 
$\PP^1(\CC_v)$ with $\CC_v \cup \{\infty\}$, and take $\fX = \{\infty\}$.  
Put $F_{w_0} = K_{v_0}(t^{1/p}) = \FF_p((t^{1/p}))$, so that $F_{w_0}/K_{v_0}$
is purely inseparable of degree $p$, and put $F_{w_1} = K_{v_1}((t-1)^{1/p^2}) = \FF_p(((t-1)^{1/p^2}))$,
so that $F_{w_1}/K_{v_1}$ is purely inseparable of degree $p^2$.  
Fix a place $v_2 \in \cM_K$ distinct from  $v_0$ and $v_1$, and 
define $\EE = \EE(r) = \prod_{v \in \cM_K} E_v$ by putting 
\begin{eqnarray*}
E_{v_0} \ = \ t^{-1/p} + \cO_{w_0} & = & B(t^{-1/p},1) \cap F_{w_0} \ , \\ 
E_{v_1} \ = \ (t-1)^{-1/p^2} + \cO_{w_1} & = & B((t-1)^{-1/p^2},1) \cap F_{w_1} \ ,
\end{eqnarray*} 
taking  $E_{v_2} = D(0,r)$, 
and letting $E_v = D(0,1)$ be $\fX$-trivial for all $v \ne v_0, v_1, v_2$.  
\index{$\fX$-trivial}
Then $\EE$ is compatible with $\fX${compatible with $\fX$} and satisfies all the hypotheses of Theorem $\ref{aT1}$ 
apart from the inseparability of $F_{w_0}/K_{v_0}$ and $F_{w_1}/K_{v_1}$, and 
\begin{equation*}
\gamma(\EE,\fX) \ = \ r \cdot p^{-\frac{p+1}{p^2(p-1)}}  \ . 
\end{equation*}
If $r > p^{(p+1)/(p^2(p-1)})$
we have $\gamma(\EE,\fX) > 1$, but there are no $\alpha \in \tK$ 
whose conjugates all belong to $\EE$.
\end{example}   

\begin{proof}  The argument is similar to that in Examples \ref{FFSepExample1} and \ref{FFSepExample2},
except that here each element of $\tK \cap E_{v_0}$ satisfies $[K(\alpha) : K]^{\insep} = p$, 
while each element of $\tK \cap E_{v_1}$ satisfies $[K(\alpha) : K]^{\insep} = p^2$. 
\end{proof}
  
\section{ Examples on Elliptic Curves} \label{EllipticExamplesSection}
\index{elliptic curve}

\smallskip
{\bf Capacities of Archimedean Sets on Elliptic Curves.}

It is difficult to find capacities of archimedean sets on curves of positive genus, 
but explicit formulas for some sets can be obtained 
using pullbacks from $\PP^1$.\index{examples!elliptic curves!archimedean pullback sets}
  
Let $K_v$ be $\RR$ or $\CC$, and suppose $\cE_v/K_v$ 
is defined by a Weierstrass equation
\index{Weierstrass equation!for specific elliptic curves}
\begin{equation} \label{FWe1}
y^2 + a_0 xy + a_3 y \ = \ x^3 + a_2 x^2 + a_4 x + a_6 \ .
\end{equation}
We will compute capacities of sets relative to the origin $\oORIG = \infty$, 
using $z = x/y$ as the uniformizing parameter.\index{uniformizing parameter!normalizes capacity}

Let $f \in K_v(\cE_v)$ be a rational function of degree $d > 0$ whose only poles
are at $\oORIG$,  and let $H \subset \CC$ be a compact set of positive capacity.
\index{capacity}  
Take $E_v = f^{-1}(H) \subset \cE_v(\CC)$.   By the pullback formula (\ref{GreenPullbackF}), 
\index{Green's function!examples!elliptic curve} 
\index{Green's function!pullback formula for}
\begin{equation*}
G(p,\oORIG;E_v) \ = \ \frac{1}{d} G(f(p),\infty;H) \ .
\end{equation*} 
Assume that capacities of sets in $\CC$ are computed using the 
standard uniformizing parameter\index{uniformizing parameter!normalizes capacity} at $\infty$, 
and that $\lim_{p \rightarrow \oORIG} |f(p) \cdot z(p)^d| = A$.  Then 
\index{Robin constant!examples!elliptic curve}
\begin{equation} \label{FEVF} 
V_{\oORIG}(E_v) \ = \ \lim_{p \rightarrow \oORIG} G(p,\oORIG;E_v) + \log(|z(p)|) 
                \ =  \frac{1}{d} (V_{\infty}(H) + \log(A)) \ .
\end{equation} 
In particular, taking $f(p) = x(p)$, then 
\begin{equation} \label{FEVXF} 
V_{\oORIG}(E_v) \ = \ \frac{1}{2} V_{\infty}(H) \ .
\end{equation}
\index{Robin constant!examples!elliptic curve}
If $f(p) = y(p)$, then
\begin{equation} \label{FEVYF} 
V_{\oORIG}(E_v) \ = \ \frac{1}{3} V_{\infty}(H) \ .
\end{equation}
For example, if $E_v = \{p \in \cE_v(\CC) : |y(p)| \le R\} = y^{-1}(D(0,R))$, then
$V_{\oORIG}(E_v)  =  -\frac{1}{3} \log(R)$. 

Now assume $K_v = \RR$;  
we will compute some capacities of sets $E_v \subset \cE_v(\RR)$.
Completing the square on the left side of (\ref{FWe1}), we get 
\begin{equation} \label{FWe2}
(y + \frac{1}{2}a_0 x + \frac{1}{2}a_3)^2 \ = 
\ x^3 + (a_2 + \frac{1}{4}a_0^2) x^2 + (a_4 + \frac{1}{2} a_0a_3) x 
      + (a_6 + \frac{1}{4}a_3^2) \ .
\end{equation}
Let $g(x)$ be the polynomial on the right side of (\ref{FWe2}).  
Then $g(x)$ has either one or three real roots.  

If $g(x) = (x-a)(x-b)(x-c)$ with $a < b < c$, then $\cE_v(\RR)$ has two components, 
the bounded loop $x^{-1}([a,b])$ and the unbounded loop $x^{-1}([c,\infty])$.  
If $E_v \subset \cE_v(\RR)$ is the bounded loop, then by formula (\ref{F1Seg3})
\begin{equation} \label{FEVBL1}
V_{\oORIG}(E_v) \ = \ -\frac{1}{2}\log(\frac{b-a}{4}) \ .
\end{equation}
If $T > c$ and $E_v = \{p \in \cE(\RR) : x(p) \le T\}$, 
then $E_v = x^{-1}([a,b] \cup [c,T])$ and by formula (\ref{F2SegV1})  
\index{Robin constant!examples!elliptic curve}
\begin{equation} \label{FEVBL2}
V_{\oORIG}(E_v) \ = \ -\frac{1}{2}\log\left( 
  \frac{\sqrt[4]{(c-a)(c-b)(T-a)(T-b)}}
     {2 \left| \frac{\theta(\Re(M(\infty))/K,\tau;\frac{1}{2},\frac{1}{2})}
               {\theta(0,\tau;0,\frac{1}{2})} \right|}\right) 
              \ .
\end{equation} 

If $g(x)$ has only one real root, $x = c$, then $\cE_v(\RR) = x^{-1}([c,\infty])$
has one component.  If $T > c$ and we take 
$E_v = \{p \in \cE(\RR) : x(p) \le T\}$, 
then $E_v = x^{-1}([c,T])$ and by formula (\ref{F1Seg3})
\index{Robin constant!examples!elliptic curve}
\begin{equation} \label{FEVUL1}
V_{\oORIG}(E_v) \ = \ -\frac{1}{2}\log(\frac{T-c}{4}) \ .
\end{equation}

\smallskip
{\bf Capacities of Nonarchimedean Sets on Elliptic Curves.}   

In this subsection we will compute the capacities of certain sets of integral
points on N\'eron models and Weierstrass models. 
\index{N\'eron model!of elliptic curve}

\smallskip
\begin{theorem} \label{EllipticNonarch}\index{Kodaira classification of elliptic curves}
Suppose $K_v$ is nonarchimedean.  
Let $\cE/K_v$ be an elliptic curve, and let $\cE_{\cN}/\Spec(\cO_v)$ be its Neron model.
Let $\oORIG$ be the origin of $\cE$, and let 
$E_v \subset \cE(K_v)$ be the set of $K_v$-rational points 
which do not specialize $\pmod{v}$ to the origin of the special fibre $\cE_{\cN,v}$. 
Equivalently, if $\cE_{\cW}$ is the affine model of $\cE$ 
defined by a minimal Weierstrass equation for $\cE$, then $E_v = \cE_{\cW}(\cO_v)$.
\index{Weierstrass equation!minimal}\index{examples!elliptic curves!nonarchimedean Kodiara fibres}  

Write $k_v$ for the residue field of $\cO_v$, and let $q_v$ be its order.   
Let  $g_{\oORIG}(z) \in K_v(\cE)$ be a uniformizing parameter\index{uniformizing parameter!normalizes Robin constant}
which specializes $\pmod{v}$ to a uniformizer at the origin of $\cE_{\cN,v}$,
so $g_{\oORIG}(z)$ is a local coordinate function which defines the 
formal group at the origin of $\cE$;  
\index{formal group}
for example, take $g_{\oORIG}(z) = x/y$, in terms of the standard coordinates
on a minimal Weierstrass model $\cE_{\cW}$. 
\index{Weierstrass equation!minimal}   
 
Then the local Robin constant 
\index{Robin constant!examples!elliptic curve|ii}
\begin{equation*} 
V_{\oORIG}(E_v) \ := \ \lim_{z \rightarrow \infty} G(z,\oORIG;E_v) + \log_v(|g_{\oORIG}(z)|_v)
\end{equation*} 
is given by the following formulas, according to the reduction type of $\cE$:  
\index{Green's function!examples!elliptic curve}

$(A)$  Type $I_0:$  Good reduction.  If $\#(\cE_{\cN,v}(k_v)) = N$, then     
\index{Robin constant!examples!elliptic curve} 
\begin{equation} 
V_{\oORIG}(E_v) \ = \ \frac{q_v}{(N-1)(q_v-1)} \ . \label{FEGood}
\end{equation}   

$(B)$  Type $I_1:$  Nodal reduction, one component,  
\index{Robin constant!examples!elliptic curve}
\begin{equation}
V_{\oORIG}(E_v) \ = \ \frac{q_v}{(q_v-2)(q_v-1)}\ ;\label{FENodal}
\end{equation} 
here we assume $q_v > 2:$ \ if $q_v = 2$, then $E_v$ is empty.  
        
$(C)$  Type $I_n$, $n \ge 2:$  
Multiplicative reduction, a loop of $n$ lines.

\qquad $(C1)$  Split multiplicative reduction.  
Let $\{P_k(x)\}_{k \ge 0}$ be the polynomials defined recursively by  
$P_0(x) = 1$, $P_1(x) = x$, $P_k(x) = x P_{k-1}(x) - P_{k-2}(x)$ for $k \ge 2$, 
so  
\begin{equation}
P_k(x) \ = \ \frac{1}{2^{k+1}} \cdot
\frac{(x+\sqrt{x^2+4})^{k+1}-(x-\sqrt{x^2+4})^{k+1}}{\sqrt{x^2+4}} \ \ . \label{FECP}
\end{equation} 
Then 
\begin{equation}
V_{\oORIG}(E_v)  \ = \ \frac{q_v P_{n-1}(q_v + \frac{1}{q_v})}
{(q_v^2-q_v+2) P_{n-1}(q_v+\frac{1}{q_v}) - 2q_v P_{n-2}(q_v+\frac{1}{q_v}) - 2q_v} \ .
\label{FESplitMult}
\end{equation} 

\qquad $(C2)$  Non-split multiplicative reduction, $n$ odd: 
one component with rational points, 
\index{Robin constant!examples!elliptic curve}
\begin{equation}
V_{\oORIG}(E_v) \ = \ \frac{1}{q_v-1}\ . \label{FENonSplitMult1}
\end{equation}
                       
\qquad $(C3)$  Non-split multiplicative reduction, $n$ even: 
two components with rational points, 
\index{Robin constant!examples!elliptic curve}
\begin{equation}
     V_{\oORIG}(E_v) \ = \ 
\frac{n q_v^2 + 4q_v-n}{(q_v-1)(n q_v^2 + 8q_v - n +4)} \ . \label{FENonSplitMult2}
\end{equation} 

$(D)$ Type $II:$  Cuspidal reduction, one component,
\index{Robin constant!examples!elliptic curve}
\begin{equation}
V_{\oORIG}(E_v) \ = \ \frac{q_v}{(q_v-1)^2} \ . \label{FEII}
\end{equation}

$(E)$ Type $III:$  Two lines tangent at a point,  
\index{Robin constant!examples!elliptic curve}
\begin{equation}
     V_{\oORIG}(E_v) \ = \ \frac{q_v(q_v+1)}{(q_v-1)(q_v^2+2q_v-1)} \ .  
\label{FEIII}
\end{equation} 

$(F)$ Type $IV:$  Three lines meeting transversely at a point.

\qquad $(F1)$ One $k_v$-rational component,
\begin{equation} 
V_{\oORIG}(E_v) \ = \ \frac{q_v}{(q_v-1)^2}\ .  \label{FEIV1}
\end{equation}

\qquad $(F2)$ All three components $k_v$-rational, 
\index{Robin constant!examples!elliptic curve} 
\begin{equation}
     V_{\oORIG}(E_v) \ = \ \frac{q_v}{(q_v+1)(q_v-1)} \ . \label{FEIV3}
\end{equation} 

$(G)$  Type $I_0^*:$  Four lines of multiplicity $1$ 
meeting a line of multiplicity $2$ at distinct points.   
  
\qquad $(G1)$ One $k_v$-rational component of multiplicity $1$, 
\index{Robin constant!examples!elliptic curve} 
\begin{equation}
     V_{\oORIG}(E_v) \ = \ \frac{q_v}{(q_v-1)^2} \ . \label{FEI0*1}
\end{equation} 

\qquad $(G2)$ Two $k_v$-rational components of multiplicity $1$,  
\index{Robin constant!examples!elliptic curve}
\begin{equation}
     V_{\oORIG}(E_v) \ = \ \frac{1}{q_v-1} \ . \label{FEI0*2}
\end{equation} 
 
\qquad $(G3)$ Four $k_v$-rational components of multiplicity $1$,  
\begin{equation}
     V_{\oORIG}(E_v) \ = \ \frac{q_v(2 q_v-1)}{(q_v-1)(2 q_v^2 + 1)} \ . 
\label{FEI0*4}
\end{equation} 
 
$(H)$  Type $I_n^*$, $n \ge 1:$  Two lines of multiplicity $1$ 
at each end of a chain of $n+1$ lines of multiplicity $2$. 

\qquad $(H1)$ Two $k_v$-rational components of multiplicity $1$ $($adjacent$)$,
\index{Robin constant!examples!elliptic curve}
\begin{equation}
     V_{\oORIG}(E_v) \ = \ \frac{1}{q_v-1} \ . \label{FEIn*2}
\end{equation} 
 
\qquad $(H2)$ Four $k_v$-rational components of multiplicity $1$,  
\index{Robin constant!examples!elliptic curve}
\begin{equation}
 V_{\oORIG}(E_v) \ = \ 
\frac{q_v((n+2)q_v^2 - (n-1)q_v - 1)}{(q_v-1)((n+2)q_v^3-(n-2) q_v^2 + q_v + 1)} \ .
\label{FEIn*4}
\end{equation}

$(I)$ Type $IV^*:$  Three lines of multiplicity $1$, 
each meeting a line of multiplicity $2$,
which in turn meets a line of multiplicity $3$.   

\qquad $(I1)$  One $k_v$-rational component of multiplicity $1$, 
\index{Robin constant!examples!elliptic curve}
\begin{equation}  
V_{\oORIG}(E_v)  =  \frac{q_v}{(q_v-1)^2} \ . \label{FEIV*1}
\end{equation}

\qquad $(I2)$ Three $k_v$-rational components of multiplicity $1$, 
\begin{equation}
     V_{\oORIG}(E_v) \ = \ \frac{q_v(2q_v -1)}{(q_v-1)(2q_v^2 - q_v + 1)} \ .
     \label{FEIV*3}
\end{equation} 

$(J)$ Type $III^*:$  
A chain of lines with multiplicities $1-2-3-4-3-2-1$ 
with another line of multiplicity $2$ meeting the component of multiplicity $4$,   
\index{Robin constant!examples!elliptic curve}
\begin{equation}
     V_{\oORIG}(E_v) \ = \ \frac{q_v(3q_v-1)}{(q_v-1)(3q_v^2 - 2 q_v + 1)} \ .
     \label{FEIII*}
\end{equation}
 
$(K)$ Type $II^*:$  A chain of lines with multiplicities $1-2-3-4-5-6-4-2$  
with another line of multiplicity $3$ meeting the component of multiplicity $6$,  
\index{Robin constant!examples!elliptic curve}
\begin{equation} 
V_{\oORIG}(E_v) \ = \ \frac{q_v}{(q_v-1)^2} \ . \label{FEII*}
\end{equation}
\end{theorem} 
 
The proof requires a formula for the canonical distance
\index{canonical distance!nonarchimedean!intersection theory formula for}  
in terms of intersection theory, derived for minimal models in (\cite{RR1}, \S2.4)
and for `well-adjusted' models in (\cite{CR}).  
In proving Theorem \ref{EllipticNonarch}, 
we will need the intersection theory formula for an arbitrary model.
\index{well-adjusted model}   

Let $\cC_v/K_v$ be a smooth, connected, projective curve, 
and let $\fC_v/\Spec(\cO_v)$ be any regular model of $\cC_v$.  
Given a point $p \in \cC_v(\CC_v)$ (respectively, a divisor $D$ on $\cC_v$), 
write $(p)$ (respectively $\cl(D)$) for its closure in $\fC_v$.  Let $F_1, \ldots, F_m$
be the irreducible components of the special fibre of $\fC_v$.  
Recall that the $m \times m$ intersection matrix $(F_i \cdot F_j)$ 
is symmetric and negative semidefinite, with rank $m-1$.  I
ts kernel consists of vectors which are multiples of the special fibre,
and its image consists of all vectors orthogonal to the special fibre.  In particular, any vector
$\sum a_i F_i$ supported on components of multiplicity $1$ in the special fibre, 
for which $\sum a_i = 0$, belongs to the image.   

If $f \in K_v(\cC_v)$ is a nonzero rational function, 
write $\div_{\cC_v}(f)$ for its divisor on $\cC_v$,
and $\div_{\fC_v}(f)$ for its divisor on $\fC_v$;  then there are integers
$c_1, \ldots, c_m$ for which 
\begin{equation} \label{FDivf}
\div_{\fC_v}(f) \ = \ \cl(\div_{\cC_v}(f)) + \sum_{j=1}^m c_j F_j \ .
\end{equation} 
If $a \in \cC_v(K_v)$, then
\begin{equation} \label{FCD1} 
-\log_v(|f_v(a)|_v) \ = \ \ord_v(f(a)) \ = \ (a) \cdot \div_{\fC_v}(f) \ .
\end{equation}   

Suppose $\zeta \in \cC_v(K_v)$.  Then $\zeta$ specializes to a nonsingular 
closed point on the special fibre of $\fC_v$, and we can choose a uniformizer 
$g_{\zeta}(z) \in K_v(\cC)$ in such a way such that for all $t \in \cC_v(K_v)$ 
which specialize to that same closed point, 
\begin{equation*}
-\log_v(|g_{\zeta}(t)|_v) \ = \ \ord_v(g_{\zeta}(t)) \ = \ (t) \cdot (\zeta) 
\end{equation*} 
Normalize the canonical distance $[x,y]_{\zeta}$ so that 
\index{canonical distance!$[z,w]_{\zeta}$}  
\index{canonical distance!$[z,w]_{\zeta}$!normalization of}  
\begin{equation} \label{FCDNorm}  
\lim_{y \rightarrow \zeta} \, [x,y]_{\zeta} |g_{\zeta}(y)|_v \ = \ 1 \ .  
\end{equation} 

Let $a \ne b$ be points of $\cC_v(\CC_v) \backslash \{\zeta\}$.  
By (\cite{RR1}, Proof of Uniqueness for Theorem 2.1.1, p.57),  
$[a,b]_{\zeta}$  is given by 
\begin{equation} \label{FCD2} 
-\log_v([a,b]_{\zeta}) \ = \ \lim_{f} \frac{1}{\deg(f)}(- \log_v(|f(a)|_v)) 
\end{equation}
where the limit is taken over any sequence of functions whose only poles are at $\zeta$, 
whose zeros approach $a$, and which are normalized so that 
$\lim_{z \rightarrow \zeta} |f(z)|_v \cdot |g_{\zeta}(z)^{\deg(f}|_v = 1$.   
Consider $\div_{\fC_v}(f)$ for such an $f$.  
After relabeling the $F_i$, we can assume that $(\zeta)$ specializes to $F_1$.  
The normalization (\ref{FCDNorm}) determines the constant $c_1$ in (\ref{FDivf}):
if the zeros of $f$ are $b_1, \ldots, b_n$ then $c_1 = -\sum_{i=1}^n (\zeta) \cdot (b_i)$. 
The remaining $c_j$ are determined by the equations $F_i \cdot \div_{\fC_v}(f) = 0$, 
for $i = 2, \ldots, m$.  Put $\hat{c}_j = (c_j-c_1)/\deg(f)$.    

Combining (\ref{FCD1}) and (\ref{FCD2}), passing to the limit in $f$,
and using the asymptotic stability of the various terms in the intersection products, 
we obtain the intersection theory formula for the canonical distance: 
\index{canonical distance!nonarchimedean!intersection theory formula for|ii}
\index{Intersection Theory formula!for the canonical distance|ii}  

\begin{proposition} \label{CanonDistForm} 
Let $\cC_v/K_v$ be a smooth, connected, projective curve. 
Fix a regular model $\fC_v/\Spec(\cO_v)$ of $\cC_v$, and let   
$F_1, \ldots, F_m$ be the irreducible components of the special fibre of $\fC_v$.  
If $\zeta \in \cC_v(K_v)$ and the canonical distance $[x,y]_{\zeta}$ 
\index{canonical distance!$[z,w]_{\zeta}$}\index{canonical distance!$[z,w]_{\zeta}$!normalization of}  
is normalized as in $(\ref{FCDNorm})$, 
then for distinct $a, b \in \cC_v(K_v) \backslash \{\zeta\}$   
\begin{equation} \label{FCD3} 
-\log_v([a,b]_{\zeta}) \ = \ (a) \cdot (b) - (a) \cdot (\zeta) - (b) \cdot (\zeta)
+ \sum_{j=1}^m \hat{c}_j \, F_j \cdot (a)
\end{equation} 
where $\hat{c}_1, \ldots, \hat{c}_m \in \QQ$ are uniquely determined by the equations
\begin{equation} \label{FCD4} 
\left\{ \begin{array}{ll}
\sum_{j=1}^m \hat{c}_j \, F_i \cdot F_j \ = \ F_i \cdot (\zeta) - F_i \cdot (b)  
                   & \text{for $i = 1, \ldots, m$;} \\
\hat{c}_1 = 0 \ , & \ 
\end{array} \right.
\end{equation} 
if \, $\zeta$ specializes to $F_1$. 
\end{proposition} 
 
In (\ref{FCD4}), the numbers $\hat{c}_j$ depend only on the components to which $b$
and $\zeta$ specialize, and if $a$ specializes to $F_k$, then in (\ref{FCD3})
\begin{equation} \label{FCD5} 
\sum_{j=1}^m \hat{c}_j \, F_j \cdot (a) \ = \ \hat{c}_k \ .
\end{equation} 
If $b$ specializes to $F_{\ell}$, we will write $j_{\zeta}(F_k,F_{\ell})$ for $\hat{c}_k$.
It is easily seen that $j_{\zeta}(F_k,F_{\ell}) \ge 0$, 
and that $j_{\zeta}(F_k,F_{\ell}) = j_{\zeta}(F_{\ell},F_k)$.  
If the model $\fC_v$ is projective, and if $\|x,y\|_v$ is the spherical metric\index{spherical metric} 
on $\cC_v$ determined by the projective embedding of $\fC_v$, then 
\begin{equation} \label{FCD6}
(a) \cdot (b) \ = \ -\log_v(\|a,b\|_v) \ .
\end{equation}
Thus (\ref{FCD3}) can be rewritten 
\begin{equation} \label{FCD7} 
-\log_v([a,b]_{\zeta}) \ = \ -\log_v\left(\frac{\|a,b\|_v}{\|a,\zeta\|_v\|b,\zeta\|_v}\right) 
+ j_{\zeta}(F_k,F_{\ell}) \ . 
\end{equation} 

We now apply this to potential functions.  
Suppose $E_v \subset \cC_v(K_v)$ is compact with positive capacity, 
\index{capacity}
and that all the points of $E_v$ specialize to the same component $F_{\ell}$.  
Suppose in addition that no point of $E_v$ specializes to the
same closed point of the special fibre as $\zeta$.  If $\mu$ is the equilibrium distribution
of $E_v$ with respect to $\zeta$, and if $a \in \cC_v(K_v)$ specializes to $F_k$, 
then by (\ref{FCD7}), 
\begin{equation} \label{FCD8}
u_{E_v}(a,\zeta) \ = \ 
\int_{E_v} -\log_v(\|a,b\|_v) d\mu(b) + j_{\zeta}(F_k,F_{\ell}) -\log_v(\|a,\zeta\|_v) \ .
\end{equation} 
Note that if $a, b \in \cC_v(K_v)$ specialize to different closed points of the special fibre,
then $-\log_v(\|a,b\|_v) = 0$.  If they specialize to the same closed point (which is necessarily 
nonsingular on the special fibre), and if we fix a $K_v$-rational isometric parametrization 
\index{isometric parametrization}
of the ball $B(b,1)^-$, then $-\log_v(\|a,b\|_v) = -\log_v(|a^{\prime}-b^{\prime}|_v)$ 
where $a^{\prime}, b^{\prime} \in K_v$ correspond to $a, b$ under the isometric parametrization.
\index{isometric parametrization}
In particular, if $E_v = \cC_v(K_v) \cap B(b_0,1)^-$, then the integral appearing in (\ref{FCD8}) 
is the same as the one studied in Proposition \ref{OwProp}, and  
\index{Robin constant!computing nonarchimedean}
\begin{equation} \label{FCD9}
V_{\zeta}(E_v) \ = \ 1 + \frac{1}{q_v-1} + j_{\zeta}(F_{\ell},F_{\ell}) \ .
\end{equation} 
For $a \in \cC_v(K_v) \backslash (B(b_0,1)^- \cup B(\zeta,1)^-)$ specializing to $F_k$, 
\begin{equation} \label{FCD10}
u_{E_v}(a,\zeta) \ = \ j_{\zeta}(F_k,F_{\ell}) \ .\
\end{equation}   

If $E_v$ consists of points belonging to several balls in a single component $F_{\ell}$, 
the averaging procedure used in Corollary \ref{OwN} applies.  
Thus, if $E_v = \cC_v(K_v) \cap (\bigcup_{i=1}^M B(b_i,1)^-)$ 
where $b_1, \ldots, b_M \in \cC_v(K_v)$ specialize to distinct closed points of $F_{\ell}$ 
(and $\zeta$ does not specialize to any of those points), then 
\index{Robin constant!computing nonarchimedean}
\begin{equation} \label{FCD11}
V_{\zeta}(E_v) \ = \ \frac{q_v}{M(q_v-1)} + j_{\zeta}(F_{\ell},F_{\ell}) \ , 
\end{equation}  
while for $a \in \cC_v(K_v) \backslash (\bigcup_{i=1}^M B(b_i,1)^- \cup B(\zeta,1)^-)$ 
specializing to $F_k$, 
\begin{equation} \label{FCD12}
u_{E_v}(a,\zeta) \ = \ j_{\zeta}(F_k,F_{\ell}) \ .\
\end{equation}   

Finally, if $E_v$ has points belonging to several components, 
we can find $V_{\zeta}(E_v)$ by solving the system of equations (\ref{FSolve}) 
\index{Robin constant!computing nonarchimedean}
for the potential functions of the sets $E_{v,\ell}$, 
where $E_{v,\ell} \subset E_v$ is the set of points specializing to $F_{\ell}$.  

\smallskip
\begin{proof}[Proof of Theorem \ref{EllipticNonarch}]  
As might be expected, the proof involves considering the various reduction types individually.   
One must solve the system of equations discussed above, in each case.  

\vskip .05 in
In cases (A), (B), (C2), (D), (F1), (G1), (I1) and (K), only the identity component of the
special fibre has rational points.  Since $E_v = \cC_v(K_v) \backslash B(\oORIG,1)^-$, 
we can apply (\ref{FCD11}) with $M = \#\cE_0(k_v)-1$ and $j_{\oORIG}(F_1,F_1) = 0$,
where $F_1 = \cE_0$ is the identity component. 
In case B), $\#\cE_0(k_v) = q_v-1$;  in case C2), $\#\cE_0(k_v) = q_v+1$; 
and in cases D)--K), $\#\cE_0(k_v) = q_v$.  

\vskip .05 in
In cases (E), (F2), (G2), (G3), (I2) and  (J), 
where the special fibre has a fixed number of components,
the computations are similar except for details.   
For each, one first solves the system of equations (\ref{FCD4}) 
to determine the numbers $j_{\oORIG}(F_k,F_{\ell})$; then finds the 
potential functions $u_{E_{v,\ell}}(z,\zeta)$ corresponding to the various components $F_{\ell}$,
using (\ref{FCD11}) and (\ref{FCD12}) and taking into account the number of $k_v$-rational 
closed points on each component; and finally solves the system (\ref{FSolve}) 
to find $V_{\oORIG}(E_v)$, using the potential functions $u_{E_{v,\ell}}(z,\zeta)$.
\index{Robin constant!computing nonarchimedean}
The computations were carried out using Maple.\index{Maple computations}

We will illustrate the method in case (F2), where $\cE$ has Type IV additive reduction
and the special fibre consists of three components meeting transversely at a point, 
each component being $k_v$-rational.
Let these components be $F_1 = \cE_0$, $F_2$, and $F_3$. Each has $q_v$ $k_v$-rational 
closed points, so $E_{v,1}$ is formed from $q_v-1$ balls, 
while $E_{v,2}$ and $E_{v,3}$ each have $q_v$ balls.   
Trivially 
\begin{equation*} 
j_{\oORIG}(F_1,F_1) = j_{\oORIG}(F_2,F_1) = j_{\oORIG}(F_3,F_1) = 0 \ .
\end{equation*} 
To find the $j_{\oORIG}(F_i,F_2)$, note that each $F_i^2= -2$, while  
$F_i \cdot F_j = 1$ if $i \ne j$,
and solve the system (\ref{FCD4}) which reads 
\begin{equation*} 
\left\{ \begin{array}{l} 
\hat{c}_1 \cdot (-2) + \hat{c}_2 \cdot 1 + \hat{c}_3 \cdot 1 = 1 \\
\hat{c}_1 \cdot 1 + \hat{c}_2 \cdot (-2) + \hat{c}_3 \cdot 1 = -1 \\
\hat{c}_1 \cdot 1 + \hat{c}_2 \cdot 1 + \hat{c}_3 \cdot (-2) = 0  \\
\hat{c}_1 = 0 
      \end{array} \right.
\end{equation*} 
giving $j_{\oORIG}(F_1,F_2) = \hat{c}_1 = 0$, $j_{\oORIG}(F_2,F_2) = \hat{c}_2 = 2/3$,
$j_{\oORIG}(F_3,F_2) = \hat{c}_3 = 1/3$.  
Similarly $j_{\oORIG}(F_1,F_3) = 0$, $j_{\oORIG}(F_2,F_3) = 1/3$,
$j_{\oORIG}(F_3,F_3) = 2/3$.  The potential functions $u_{E_{v,i}}(z,\oORIG)$ 
are then given by (\ref{FCD11}) and (\ref{FCD12}), 
with $V_{\oORIG}(E_{v,1}) = q_v/(q_v-1)^2$ 
\index{Robin constant!computing nonarchimedean}
and $V_{\oORIG}(E_{v,2}) = V_{\oORIG}(E_{v,3 }) = 2/3 + 1/(q_v-1)$. 

To find $V_{\oORIG}(E_v)$, solve the system (\ref{FSolve}) which reads 
\begin{equation*}
\left\{ \begin{array}{l}
1 \ = \ 0 V + s_1 + s_2 + s_3 \\
0 \ = \ V - \frac{q_v}{(q_v-1)^2} s_1 - 0 s_2 - 0 s_3 \\
0 \ = \ V - 0 s_1 - (\frac{2}{3} + \frac{1}{q_v-1}) s_2 - \frac{1}{3} s_3 \\
0 \ = \ V - 0 s_1 - \frac{1}{3} s_2 - (\frac{2}{3} + \frac{1}{q_v-1}) s_3 
      \end{array} \right. 
\end{equation*}   
giving $V = V_{\oORIG}(E_v) = q_v/(q_v^2-1)$, and $s_1 = (q_v-1)/(q_v+1)$,
$s_2 = s_3 = 1/(q_v+1)$.  If necessary, the weights $s_1, s_2, s_3$ could be used to find 
$u_{E_v}(z,\zeta)$ for any $z \in \cE(\CC_v)$.
\index{Robin constant!computing nonarchimedean}

\vskip .05 in
The remaining cases (C1), (C3) and (H), where the number of components 
depends on $n$, must be treated separately.  

\vskip .05 in
First consider case (C3), non-split multiplicative reduction with $n = 2N$.  
Among the $n$ components $\cE_0, \ldots, \cE_{2N-1}$ (listed cyclically around the loop), 
only $\cE_0$ and $\cE_N$ have $k_v$-rational points (each with $q_v+1$), 
while the other components have none.
By (\cite{RR1}, p.96), $j_{\oORIG}(\cE_{\ell},\cE_{\ell}) = \ell - \ell^2/n$, 
so $j_{\oORIG}(\cE_0,\cE_0) = 0$ and $j_{\oORIG}(\cE_N,\cE_N) = n/4$.  Thus
\index{Robin constant!computing nonarchimedean}
\begin{equation*}
V_{\oORIG}(E_{v,0}) = \frac{1}{q_v-1} \ , 
\quad  V_{\oORIG}(E_{v,N}) = \frac{n}{4} + \frac{q_v}{q_v^2-1} \ .
\end{equation*}
The equations (\ref{FSolve}) read 
\begin{equation*}
\left\{ \begin{array}{l}
1 \ = \ 0 V + s_1 + s_2 \\
0 \ = \ V - \frac{1}{q_v-1} s_1 - 0 s_2  \\
0 \ = \ V - 0 s_1 - (\frac{n}{4} + \frac{q_v}{q_v^2-1}) s_2 
      \end{array} \right. 
\end{equation*}   
giving 
\begin{eqnarray*}
V = V_{\oORIG}(E_v) = \frac{nq_v^2+4q_v-n}{(q_v-1)(nq_v^2 + 8q_v-n +4)} \ , \quad \ \\ 
s_1 = \frac{nq_v^2 + 4 q_v - n}{nq_v^2 + 8q_v - n +4} \ , \quad 
s_2 = \frac{4q_v + 4 }{nq_v^2 + 8q_v - n +4} \ . 
\end{eqnarray*}  

\smallskip
Next, consider case (H):  Type $I_n^*$ additive reduction, $n \ge 1$.  
Let $F_1$, $F_2$, $F_3$, and $F_4$ be
the four components of multiplicity $1$, and let $G_1, \ldots, G_{n+1}$ be the components of
multiplicity $2$, listed sequentially along the chain; assume $F_1$ and $F_2$ meet $G_1$, 
and $F_3$ and $F_4$ meet $G_{n+1}$, with $F_1 = \cE_0$ being the identity component.  

We first determine the numbers $j_{\oORIG}(F_k,F_{\ell})$ and $j_{\oORIG}(G_i,F_{\ell})$.

Trivially $j_{\oORIG}(F_k,F_1) = j_{\oORIG}(G_i,F_1) = 0$ for all $k$ and $i$.  
For $F_2$, $F_3$, and $F_4$ the equations (\ref{FCD4}) can be solved recursively. 
For $F_2$, one finds in turn $j_{\oORIG}(F_1,F_2) = 0$,
$j_{\oORIG}(G_1,F_2) = 1$, $j_{\oORIG}(F_2,F_2) = 1$, 
then $j_{\oORIG}(G_i,F_2) = 1$ for $i = 2, \ldots, n+1$, 
and finally $j_{\oORIG}(F_3,F_2) = j_{\oORIG}(F_4, F_2) =1/2$. For $F_3$, one finds
$j_{\oORIG}(F_1,F_3) = 0$, $j_{\oORIG}(G_1,F_3) = 1$, $j_{\oORIG}(F_2,F_3) = 1/2$,
then $j_{\oORIG}(G_i,F_2) = (i+1)/2$ for $i = 2, \ldots, n+1$, and finally
$j_{\oORIG}(F_3,F_3) = 1 + n/4$, $j_{\oORIG}(F_4,F_3) = 1/2 + n/4$.  For $F_4$,
the values are the same as for $F_3$, except that $j_{\oORIG}(F_3,F_4) = 1/2 + n/4$
and $j_{\oORIG}(F_4,F_4) = 1 + n/4$.  

In subcase (H1), $F_1$ and $F_2$ are $k_v$-rational but $F_3$ and $F_4$ are not.  
The computation is identical to the one in case G2), and one gets $V_{\oORIG}(E_v) = 1/(q_v-1)$.  
\index{Robin constant!computing nonarchimedean}

In subcase (H2), all of $F_1, F_2, F_3, F_4$ are $k_v$-rational.  
Each has $q_v$ $k_v$-rational closed points.  The equations (\ref{FSolve}) read
\begin{equation*}
\left\{ \begin{array}{l}
1 \ = \ 0 V + s_1 + s_2 + s_3 + s_4\\
0 \ = \ V - \frac{q_v}{(q_v-1)^2} s_1 - 0 s_2 - 0 s_3 - 0 s_4\\
0 \ = \ V - 0 s_1 - (1 + \frac{1}{q_v-1}) s_2 - \frac{1}{2} s_3 - \frac{1}{2} s_4\\
0 \ = \ V - 0 s_1 - \frac{1}{2} s_2 - (1 + \frac{n}{4} + \frac{1}{q_v-1}) s_3 
                                    - (\frac{1}{2} + \frac{n}{4}) s_4 \\
0 \ = \ V - 0 s_1 - \frac{1}{2} s_2 - (\frac{1}{2} + \frac{n}{4}) s_3 
                                    - (1 + \frac{n}{4} + \frac{1}{q_v-1}) s_4
      \end{array} \right. 
\end{equation*}   
and Maple\index{Maple computations} gives 
\index{Robin constant!computing nonarchimedean}
\begin{equation*} 
V_{\oORIG}(E_v) \ = \ V 
\ = \ \frac{q_v[(n+2)q_v^2 - (n-1)q_v - 1]}{(q_v-1)[(n+2)q_v^3-(n-2) q_v^2 + q_v + 1]} \ .
\end{equation*}                                  

\vskip .05 in
Case (C1), split multiplicative reduction with $n \ge 2$ components, is the most difficult.  
Let the components (listed cyclically around the loop) be $\cE_0, \ldots, \cE_{n-1}$, 
where $\cE_0$ is the identity component, and let $\bigcup_{i=0}^{n-1} E_{v,i}$ 
be the corresponding decomposition of $E_v$.  There are $q_v-1$ $k_v$-rational points on each
$\cE_i$, so $E_{v,0}$ consists of $q_v-2$ balls and all the other $E_{v,i}$ consist of $q_v-1$
balls.  Put 
\begin{equation*} 
\hE_v \ = \ \bigcup_{i=1}^{n-1} E_{v,i} \ ,
\end{equation*} 
so $E_v = E_{v,0} \cup \hE_v$.  
We will first find  $V_{\oORIG}(\hE_v)$, and then use it to find $V_{\oORIG}(E_v)$. 
\index{Robin constant!computing nonarchimedean}

For this, we will need a lemma.

\begin{lemma}[{\bf Cantor's Lemma}] \label{CantorLemma}\index{Cantor's Lemma|ii} 
\index{Cantor, David} 
Suppose $A \in M_k(\RR)$ is symmetric and negative definite.  
Let $\vbb1$ be the row vector $(1,1,\ldots,1) \in \RR^k$, 
and consider the matrix $B \in M_{k+1}(\RR)$ given in block form by 
\begin{equation*} 
B \ = \ \left( \begin{array}{cc} 
                            0     & \vbb1 \\
                        {}^t\vbb1 &   A 
               \end{array} \right) \ .
\end{equation*}
Then $B$ is invertible, and if $\vb = - \vbb1 A^{-1}$ 
and $\alpha^{-1} = - \vbb1 (A^{-1})\, {}^t\vbb1$, 
then 
\begin{equation} \label{FC1} 
B^{-1} \ = \ \left( \begin{array}{cc} 
                        \alpha & \alpha \vb \\
                        \alpha \, {}^t\vb & A^{-1} + \alpha \, {}^t\vb \cdot \vb  
               \end{array} \right) \ .
\end{equation} 
\end{lemma}

\begin{proof} See (\cite{Can3}, Lemma 3.2.3) or (\cite{RR3}, p.406). 
The proof is a block by block verification that the matrix $C$ in (\ref{FC1}) satisfies $CB = I$. 
\end{proof}   

\vskip .05 in
\index{Robin constant!computing nonarchimedean}
To find $V_{\oORIG}(\hE_v)$, let $A \in M_{n-1}(\RR)$ be the matrix 
$(-j_{\oORIG}(\cE_k,\cE_{\ell}))_{1 \le k,\ell \le n-1}$.  
By the equations (\ref{FCD4}) defining the $j_{\oORIG}(\cE_k,\cE_{\ell})$  
and the fact that $j_{\oORIG}(\cE_0,\cE_{\ell}) = 0$ for each $\ell$, it follows that 
$A$ is inverse to the tridiagonal matrix
\begin{equation*}
\Delta \ = \ \left(\begin{array}{cccccc}
                     -2  &   1    &    0   &   0    & \cdots &  0     \\
                      1  &  -2    &    1   &   0    & \cdots &  0     \\
                      0  &   1    &   -2   &   1    & \cdots &  0     \\
                  \vdots & \vdots & \ddots & \ddots & \ddots & \vdots \\
                      0  &   0    & \cdots &   1    &  -2    &  1     \\
                      0  &   0    & \cdots &   0    &  1     & -2     \\
                      \end{array} \right) \ .
\end{equation*} 
It is well known (and easy to check) that $\Delta$ is negative definite, 
so $A$ is also negative definite.  

Let $B$ be as in Lemma \ref{CantorLemma}.  Then $\alpha = 1/2$ 
and $\vb = (1, 0, \ldots, 0, 1)$ in the formula for $B^{-1}$ in that Lemma.   
Put $Q = q_v/(q_v-1)^2$ and let
\begin{equation*}
B_Q \ = \ \left( \begin{array}{cc} 
                           0        &    \vbb1                \\
                        {}^t\vbb1 & A - Q I_{n-1} 
               \end{array}  \right) 
\end{equation*} 
where $I_{n-1}$ is the $(n-1) \times (n-1)$ identity matrix.  Then the system of 
equations (\ref{FSolve}) determining $V_{\oORIG}(\hE_v)$ reads 
\begin{equation*} 
B_Q  \left( \begin{array}{c} V \\ s_1 \\ \vdots \\s_{n-1} \end{array} \right) 
\ = \ \left( \begin{array}{c} 1 \\ 0 \\ \vdots \\ 0 \end{array} \right) \ . 
\end{equation*}
Left-multiplying by $B^{-1}$ yields the simpler system 
\index{Robin constant!computing nonarchimedean}
\begin{equation*}  
\left( \begin{array}{ccccccc}
            1    &   -\frac{1}{2}Q  &   0     &   0    & \cdots &   0      & -\frac{1}{2}Q    \\
            0    & 1 + \frac{3}{2}Q &  -Q     &   0    & \cdots &   0      & -\frac{1}{2}Q    \\
            0    &      -Q          & 1+2Q    &  -Q    & \cdots &   0      &    0             \\
          \vdots &     \vdots       &  \ddots & \ddots & \ddots & \vdots   &  \vdots          \\
            0    &      0           &  \cdots &   0     &  -Q   &  1+2Q    &   -Q             \\ 
            0    &  -\frac{1}{2}Q   &   0     & \cdots  &   0   &   -Q     & 1 + \frac{3}{2}Q
        \end{array} \right) 
 \left( \begin{array}{c} V \\ s_1 \\ s_2 \\ \vdots \\ s_{n-2} \\s_{n-1} \end{array} \right) 
\ = \ \left( \begin{array}{c} 1/2 \\ 1/2 \\ 0 \\ \vdots \\ 0 \\ 1/2 \end{array} \right) \ . 
\end{equation*} 
We now solve for $V$ using Cramer's rule.\index{Cramer's Rule}  It will be useful to write 
\begin{equation*}
P_k(x) \ = \ \det \left( \begin{array}{ccccc} 
                               x &     -1   &    0   & \cdots  &   0    \\
                              -1 &      x   &   -1   & \cdots  &   0    \\  
                          \vdots & \ddots   & \ddots & \ddots  & \vdots \\  
                               0 &  \cdots  &   -1   &   x     &  -1    \\
                               0 &  \cdots  &    0   &   -1    &   x      
                          \end{array} \right) \ , 
\end{equation*} 
so $P_0(x) = 1$, $P_1(x) = x$, and in general $P_k(x) = x P_{k-1}(x) - P_{k-2}(x)$.
Solving the linear recurrence yields the formula for $P_k(x)$ in Theorem 
\ref{EllipticNonarch}.  Cramer's rule\index{Cramer's Rule} then gives 
\index{Robin constant!computing nonarchimedean}
\begin{equation} \label{FVCramer} 
V_{\oORIG}(\hE_v) \ = \ V  \ = \ 
\frac{\frac{1}{2} P_{n-1}(2 + \frac{1}{Q})}
     { P_{n-1}(2+\frac{1}{Q}) - P_{n-2}(2+\frac{1}{Q}) - 1} \ .
\end{equation}      

With $V_{\oORIG}(\hE_v)$ in hand, we can use the decomposition $E_v = E_{v,0} \cup \hE_v$
\index{Robin constant!computing nonarchimedean}
to find $V_{\oORIG}(E_v)$. The equations (\ref{FSolve}) read  
\begin{equation*}
\left\{ \begin{array}{l}
1 \ = \ 0 V + s_0 + \hat{s}_1 \\
0 \ = \ V - \frac{q_v}{(q_v-2)(q_v-1)} s_0 - 0 \hat{s}_1  \\
0 \ = \ V - 0 s_0 - V_{\oORIG}(\hE_v) \hat{s}_1 
      \end{array} \right. 
\end{equation*}   
giving 
\begin{equation} \label{FC1B} 
V_{\oORIG}(E_v) \ = \ \frac{q_v P_{n-1}(2+\frac{1}{Q})}
      {(q_v^2-q_v+2) P_{n-1}(2+\frac{1}{Q}) - 2q_v P_{n-2}(2+\frac{1}{Q}) - 2q_v} \ .
\end{equation}
If $q_v = 2$, then $E_{v,0}$ is empty and $V_{\oORIG}(E_v) = V_{\oORIG}(\hE_v)$. 
\index{Robin constant!computing nonarchimedean}
A quick check shows that (\ref{FC1B}) remains valid even in this case.    
Since $2+1/Q = q_v + 1/q_v$, this completes the proof of Theorem \ref{EllipticNonarch}.      
\end{proof} 

\smallskip
The next proposition gives the capacity of the set 
\index{capacity}
$E_v = \{ P \in \cE_v(K_v) : \	|x(P)|_v \le q_v^k \}$
for an elliptic curve $\cE_v/K_v$ in Weierstrass normal form,
\index{Weierstrass equation} 
\begin{equation}
y^2 + a_1 xy + a_3 y = x^3 + a_2 x^2 + a_4 x + a_6 \ ,  \label{FWeier} 
\end{equation} 
whose coefficients belong to $\cO_v$.
Let $\Delta$ be its discriminant, and let $\Delta_0$ be the discriminant
of a minimal Weierstrass equation.\index{Weierstrass equation!minimal}  Let $\pi_v$ be a
generator for the maximal ideal of $\cO_v$.    
Then there is an integer $m \ge 0$ for which 
\begin{equation*}
\Delta \ = \ \pi_v^{12m} \Delta_0 \ . \label{FvDelta}
\end{equation*} 
This is the number of times that Step 11 in Tate's algorithm 
\index{Tate, John} 
\index{Tate's algorithm}
(replacing $x$ by $\pi_v^2 x^{\prime}$ and  $y$ by $\pi_v^3 y^{\prime}$;   
see \cite{Sil2}, pp.364-368) is executed in computing 
the N\'eron model and a minimal Weierstrass equation for $\cE_v$.
\index{N\'eron model!of elliptic curve}
\index{Weierstrass equation!minimal}

\begin{proposition} \label{Weierstrass}\index{examples!elliptic curves!nonarchimedean Weierstrass equations}
Let $v$ be a nonarchimedean place of $K$, and  
let $\cE_v/K_v$ be the elliptic curve defined by the Weierstrass equation 
with integral coefficients $(\ref{FWeier})$.   

Let $q_v$ be the order of the residue field of $K_v$, 
and let $Q_v = V_{\oORIG}(E_v)$ be the
number associated to the N\`eron model of $\cE_v$ in Theorem $\ref{EllipticNonarch}$.
For each integer $\ell \ge 0$ put
\index{Robin constant!examples!elliptic curve}
\begin{equation} \label{FVm}  
V_{\ell}(q_v,Q_v) =  
   \frac{q_v[Q_v(q_v-1)(q_v^{2\ell-1}+1)+ (q_v^{2\ell}-1)]}
            {(q_v-1)[Q_v(q_v-1)(q_v^{2\ell}-1)+(q_v^{2\ell+1} + 1)]} \ . 
\end{equation}  

Suppose $m \ge 0$ is the number of times Step $11$ of 
Tate's algorithm is executed
\index{Tate's algorithm}
in computing the N\'eron model of $\cE_v$.
\index{N\'eron model!of elliptic curve}  
Let $k \ge -m$ be an integer, and put 
$E_{v,k} = \{p \in \cE_v(K_v) : |x(p)|_v \le q_v^{2k}\}$.   
Let $z = x/y$ be the standard uniforming parameter 
at the origin $\oORIG = \infty$ of $\cE_v$.  
Then, computing capacities relative to $z$, 
\begin{equation}
V_{\oORIG}(E_{v,k}) \ = \ -k+V_{m+k}(q_v,Q_v) \ . \label{FVWeier} 
\end{equation} 
\end{proposition} 
\index{Robin constant!examples!elliptic curve}

\begin{proof}  Let 
\begin{equation}
y_0^2 + a_{1,0} x_0y_0 + a_{3,0} y_0 
= x_0^3 + a_{2,0} x_0^2 + a_{4,0} x_0 + a_{6,0}   \label{FWeierMin} 
\end{equation}
be a minimal Weierstrass equation\index{Weierstrass equation!minimal}  for $\cE_v$. Then $z_0 = x_0/y_0$ 
can be used as the uniformizing parameter\index{uniformizing parameter!normalizes Robin constant} 
in Theorem \ref{EllipticNonarch}. 
For each $\ell \ge 0$, put 
\begin{equation} \label{FEell} 
E_v^{(\ell)} \ = \ \{P \in \cE_v(K_v) : |x_0(P)|_v \le q_v^{2\ell} \} \ , 
\end{equation} 
and put $z_{\ell} = \pi_v^{-\ell}z_0$.  Then $z_m$ is the uniformizing parameter 
$z = x/y$ in the Proposition, and $E_{v,k} = E_v^{(m+k)}$.  
    
Let   $V_0 = Q_v$, and recursively define $V_1, V_2, \ldots$ by 
requiring that $V_{\ell}$ be determined by the system of equations 
\begin{equation} \label{FVSystem}
\left\{ \begin{array}{l} 
            1 = 0V_{\ell} + s_1 + s_2 \\
            0 = V_{\ell} - \frac{q_v}{(q_v-1)^2} s_1 - 0 s_2 \\
            0 = V_{\ell} - 0 s_1 - (1+V_{\ell-1}) s_2
         \end{array} \right. 
\end{equation} 
\index{Robin constant!examples!elliptic curve}
so that  
\begin{equation} \label{FVQRec1} 
V_{\ell} \ = \ \frac{q_v(1+V_{\ell-1})}{q_v+(q_v-1)^2(1+V_{\ell-1})} 
\quad \text{for $\ell \ge 1$ \ .} 
\end{equation}
An easy induction shows that $V_{\ell} = V_{\ell}(q_v,Q_v)$ 
is given by (\ref{FVm}).

We claim that $V_{\ell}$ is the Robin constant of the set $E_v^{(\ell)}$
\index{Robin constant!examples!elliptic curve} 
relative to the uniformizing parameter $z_{\ell}$.\index{uniformizing parameter!normalizes Robin constant}  
To see this, note that $E_v^{(0)}$ coincides
with the set $E_v$ attached to the N\'eron model of $\cE_v$ in Theorem \ref{EllipticNonarch},
\index{N\'eron model!of elliptic curve}
and for each $\ell \ge 1$  
\begin{equation*}
E_v^{(\ell)} \ = \ E_v^{(\ell-1)} \bigcup \{ P \in \cE_v(K_v) : |x_0(P)| = q_v^{2\ell}\} \ .
\end{equation*} 
By the intersection theory formula for canonical distance,\index{Intersection Theory formula!for the canonical distance}  
\index{canonical distance!nonarchimedean!intersection theory formula for}  
this decomposition satisfies the conditions needed to find 
$V_{\oORIG}(E_v^{(\ell)})$ using a system (\ref{FSolve}). 
 By Theorem \ref{EllipticNonarch}, $V_{\oORIG}(E_v^{(0)}) = Q_v$ when 
 \index{Robin constant!examples!elliptic curve}
capacities are computed relative to $z_0$.  Assume that  
$V_{\oORIG}(E_v^{(\ell-1)}) = V_{\ell-1}$ when the capacity is computed relative
\index{capacity}
to $z_{\ell-1}$.  Relative to $z_{\ell}$ it is $1+V_{\ell-1}$.  
Hence the system of equations (\ref{FSolve}) for finding the capacity of 
\index{capacity}
$V_{\oORIG}(E_v^{\ell})$ relative to $z_{\ell}$ 
is exactly (\ref{FVSystem}), and our claim holds by induction. 

If the Robin constant of $E_v^{(m+k)}$ relative to $z_{m+k}$ is $V_{m+k}$,
\index{Robin constant!examples!elliptic curve}
then relative to $z_m$ it is $-k+V_{m+k}$.  This yields the result.
\end{proof}                                    

\smallskip
{\bf Global Examples on Elliptic Curves.}

\vskip .05 in
In the following examples, $N$ is the conductor of $\cE$. We take $k = \QQ$ and 
consider elliptic curves $\cE/\QQ$ defined by Weierstrass equations.\index{Weierstrass equation}     
If $p$ is a prime, by $\cE(\ZZ_p)$ or $\cE(\hcO_p)$ we mean the corresponding integral points 
on the affine curve defined by the given equation.   

\begin{example}[N = 50] \label{EllEx1}\index{examples!elliptic curves!global, Cremona $50(A1)$}
Let $\cE/\QQ$ be the elliptic curve defined by the Weierstrass equation
\index{Weierstrass equation!for specific elliptic curves}   
$y^2 + xy + y= x^3-x-2$, curve $50(A1)$ in Cremona's tables.  Then for any 
$T \ge 41.898861528$ there are infinitely many points $\alpha \in \cE(\QQbar)$ 
whose archimedean conjugates belong to $\cE(\RR)$ and satisfy $x(\alpha) < T$, 
whose conjugates in $\cE(\CC_3)$ all belong to $\cE(\ZZ_3)$, 
whose conjugates in $\cE(\CC_5)$ all belong to $\cE(\ZZ_5)$, 
and whose conjugates in $\cE(\CC_p)$ belong to $\cE(\hcO_p)$, for all primes $p \ne 3,5$.
  
If $T \le 41.898861527$, there are only finitely many $\alpha \in \cE(\QQbar)$
satisfying these conditions.   
\end{example} 

\begin{proof}  The given Weierstrass equation is minimal; after completing the
\index{Weierstrass equation!minimal}
square on the left side it becomes 
\begin{equation} \label{FComplete} 
(y+\frac{1}{2}x+\frac{1}{2})^2 \ = \ x^3 + \frac{1}{2} x^2 - \frac{1}{2}x - \frac{7}{4} \ .
\end{equation}  
At $p =2$ it has reduction type $I_1$;  at $p =3$ it has good reduction, and at $p = 5$ 
\index{good reduction}
it has reduction type $IV$ with $3$ rational components 
(see Cremona \cite{Cremona}, p.93).  
\index{Cremona, John} 
We will compute capacities\index{uniformizing parameter!normalizes capacity} 
with respect to the uniformizing parameter $z = x/y$.   

The real locus $\cE(\RR)$ consists of one unbounded loop 
$x^{-1}([\alpha,\infty])$,  where $\alpha \cong 1.256458778$ is the unique real root of the
polynomial on the right side of (\ref{FComplete}).  
Take $E_{\infty} = x^{-1}([\alpha,T])$ where $T > \alpha$.  
By formula (\ref{FEVUL1})
\index{Robin constant!examples!elliptic curve}
\begin{equation} \label{FReal1}
V_{\oORIG}(E_{\infty}) \ = \ -\frac{1}{2} \ln(\frac{T-\alpha}{4})  \ .
\end{equation}
At the prime $p =2$, the set $\cE(\ZZ_2)$ is empty (see Case B of Theorem \ref{EllipticNonarch})
so we cannot impose splitting and integrality conditions simultaneously;  
we require integrality by taking 
\begin{equation*}
E_2 \ = \ \cE(\hcO_2) \ = \ \{P \in \cE(\CC_2) : |x(P)|_2 \le 1\} \ .  
\end{equation*} 
\index{Robin constant!examples!elliptic curve}
so $V_{\oORIG}(E_p) = 0$. At $p=3$, where $\cE$ has good reduction, we take $E_3 = \cE(\ZZ_3)$. 
\index{good reduction} 
A simple check shows that $\cE \pmod{3}$ has $N = 3$ points rational over $\FF_3$.   
By formula (\ref{FEGood}) $V_{\oORIG}(E_3) = 3/4$.  At $p = 5$, we take $E_5 = \cE(\ZZ_5)$. 
\index{Robin constant!examples!elliptic curve} 
By formula (\ref{FEIV3}), $V_{\oORIG}(E_5) = 5/24$.  For $p > 5$, take $E_p = \cE(\hcO_p)$.  
Since the given model of $\cE$ and the parameter $z$ have good reduction at $p$,
\index{good reduction} 
$V_{\oORIG}(E_p) = 0$.  

Let $\EE = \prod_{p,\infty} E_p$, and take $\fX = \{\oORIG\}$.  Then  
\begin{equation*}
V(\EE,\fX) \ = \ 
-\frac{1}{2} \ln(\frac{T-\alpha}{4}) + \frac{3}{4} \ln(3) + \frac{5}{24} \ln(5) \ .
\end{equation*}   
 Maple\index{Maple computations} shows that the value 
of $T$ for which $V(\EE,\fX) = 0$ satisfies
\begin{equation*}
41.898861527 \ < \ T \ < \ 41.898861528 \ , 
\end{equation*} 
and the Fekete-Szeg\"o theorems \ref{aT1} and \ref{FSZii} yield the result. 
\index{Fekete-Szeg\"o theorem with LRC} 
\end{proof}
\index{elliptic curve}

\begin{example}[N = 32]\index{examples!elliptic curves!global, non-minimal Weierstrass equation} \label{EllEx2}
Let $\cE/\QQ$ be the elliptic curve defined by the non-minimal Weierstrass equation 
\index{Weierstrass equation!for specific elliptic curves}
\index{Weierstrass equation!non-minimal}
$y^2 = x^3-256x$.  There are infinitely many points $\alpha \in \cE(\QQbar)$ 
whose archimedean conjugates all belong to the bounded real loop in $\cE(\RR)$,
whose conjugates in $\cE(\CC_2)$ all belong to $\cE(\ZZ_2)$, and whose conjugates
in $\cE(\CC_p)$ belong to $\cE(\hcO_p)$, for each $p \ge 3$.  
\end{example} 
\index{elliptic curve}

\begin{proof}  The given Weierstrass equation is not minimal;  
the minimal equation is $y^2 = x^3-x$ (curve $32(A2)$ in Cremona's tables
\index{Cremona, John}  
\cite{Cremona}).  We will use the parameter $z = x/y$ in computing capacities
with respect to $\oORIG = \infty$.  

The bounded real loop is $E_{\infty} = x^{-1}([-16,0])$, 
for which formula (\ref{FEVBL1}) gives
\index{Robin constant!examples!elliptic curve}
\begin{equation*}
V_{\oORIG}(E_{\infty}) \ = \ -\frac{1}{2}\ln(\frac{16}{4}) \ = \ -\ln(2) \ .
\end{equation*}
Take $E_2 = \cE(\ZZ_2)$.  The curve  
$y^2 = x^3-x$ has Kodaira reduction type $III$ at $p =2$ 
(Cremona \cite{Cremona}, p.91).  
\index{Cremona, John} 
In passing from $y^2 = x^3-256x$ to $y^2 = x^3-x$ we have $m = 2$.  
By formulas (\ref{FEIII}) and (\ref{FVWeier}), $V_{\oORIG}(E_2) = 106/107$.
For all other primes $p$, take $E_p = \cE(\hcO_p)$, the trivial set 
with respect to $\oORIG$.  The model of $\cE$ given by $y^2 = x^3-256x$
and the parameter $z$ have good reduction outside $2$, so $V_{\oORIG}(E_p) = 0$.
\index{Robin constant!examples!elliptic curve} 
\index{good reduction} 

Let $\EE = \prod_{p,\infty} E_p$, and take $\fX = \{\oORIG\}$.  Then  
\begin{equation*}
V(\EE,\fX) \ = \ -\ln(2) + \frac{106}{107} \ln(2)  \ < \ 0 \ , 
\end{equation*}  
so the Fekete-Szeg\"o theorems \ref{aT1} and \ref{FSZii} yield the result.
\index{Fekete-Szeg\"o theorem with LRC} 
\end{proof}
\index{elliptic curve}

\begin{example}[N = 48]\index{examples!elliptic curves!global, Cremona $48(A3)$} \label{EllEx3}
Let $\cE/\QQ$ be the elliptic curve defined by the Weierstrass equation 
\index{Weierstrass equation!for specific elliptic curves}
$y^2 = x^3+x^2-24x+36$, curve $48(A3)$ in Cremona's tables.  Then for any 
$T \ge 28.890384202$ there are infinitely many points $\alpha \in \cE(\QQbar)$ 
whose archimedean conjugates belong to $\cE(\RR)$ and satisfy $x(\alpha) < T$, 
whose conjugates in $\cE(\CC_2)$ all belong to $\cE(\ZZ_2)$, 
whose conjugates in $\cE(\CC_3)$ all belong to $\cE(\ZZ_3)$, 
and whose conjugates in $\cE(\CC_p)$ belong to $\cE(\hcO_p)$, for $p \ge 5$.
  
If $T \le 28.890384201$, there are only finitely many $\alpha \in \cE(\QQbar)$
satisfying the conditions above.   
\end{example} 

\begin{proof}  The given Weierstrass equation is minimal; 
\index{Weierstrass equation!minimal}
it factors as $y^2 = (x+6)(x-2)(x-3)$.  At $p =2$ it has reduction type $I_2^*$,
with $4$ components rational over $k_2$;  at $p =3$ it has split multiplicative
reduction of type $I_4$ (see Cremona \cite{Cremona}, p.93).
\index{Cremona, John}   
We will compute capacities\index{uniformizing parameter!normalizes capacity}
 with respect to the uniformizing parameter $z = x/y$.   

The real locus $\cE(\RR)$ consists of the bounded loop 
$x^{-1}([-6,2])$, whose Robin constant is 
\index{Robin constant!examples!elliptic curve}
$-\frac{1}{2} \ln(2)$, together with the unbound loop $x^{-1}([3,\infty])$.  
Take $E_{\infty} = x^{-1}([-6,2] \cup [3,T])$ where $T \ge 3$.  
By formula (\ref{FEVBL2}) 
\index{Robin constant!examples!elliptic curve}
\begin{equation} \label{FReal2}
V_{\oORIG}(E_{\infty}) \ = \ f_{\infty}(T) \ := \ 
  -\frac{1}{2} \ln\left( \frac{\sqrt[4]{9(T+6)(T-2)}}
     {2 \left| \frac{\theta(\Re(M(\infty))/K,\tau;\frac{1}{2},\frac{1}{2})}
               {\theta(0,\tau;0,\frac{1}{2})} \right|}\right)  \ .
\end{equation}
Take $E_2 = \cE(\ZZ_2)$.  By formula (\ref{FEIn*4}), 
$V_{\oORIG}(E_2) = 26/35$.  Similarly, take $E_3 = \cE(\ZZ_3)$.
\index{Robin constant!examples!elliptic curve}
By formula (\ref{FESplitMult}) with $n = 4$, we have 
 $V_{\oORIG}(E_3) = 123/238$.  For $p > 3$, take $E_p = \cE(\hcO_p)$.  
Since the given model of $\cE$ and the parameter $z$ have good reduction at $p$,
\index{good reduction} 
$V_{\oORIG}(E_p) = 0$.  

Let $\EE = \prod_{p,\infty} E_p$, and take $\fX = \{\oORIG\}$.  Then  
\begin{equation*}
V(\EE,\fX) \ = \ 
f_{\infty}(T) + \frac{26}{35} \ln(2) + \frac{123}{238} \ln(3) \ .
\end{equation*}   
  
If we had taken $E_{\infty}$ to be the real loop $x^{-1}([-6,2])$, 
then by the Fekete-Szeg\"o theorem there would 
\index{Fekete-Szeg\"o theorem with LRC} 
be only finitely many $\alpha \in \cE(\QQ)$
whose conjugates meet the given conditions.  Maple\index{Maple computations} shows that the value 
of $T$ for which $V(\EE,\fX) = 0$ satisfies
\begin{equation*}
28.890384201 \ < \ T \ < \ 28.890384202 \ , 
\end{equation*} 
and the Fekete-Szeg\"o theorems \ref{aT1} and \ref{FSZii} yield the result.
\index{Fekete-Szeg\"o theorem with LRC}  
\end{proof}

\smallskip
\begin{example}[N = 360]\index{examples!elliptic curves!global, Cremona $360(E4)$} \label{EllEx4}
Let $\cE/\QQ$ be the elliptic curve defined by the Weierstrass equation
\index{Weierstrass equation!for specific elliptic curves} 
$y^2 = x^3 + 117x + 918$, curve $360(E4)$ in Cremona's tables.  Then for any 
$R \ge 142.388571238$ there are infinitely many $\alpha \in \cE(\QQbar)$ 
whose archimedean conjugates satisfy $|y(\alpha)| \le R$, 
whose conjugates in $\cE(\CC_2)$ all belong to $\cE(\ZZ_2)$, 
whose conjugates in $\cE(\CC_3)$ all belong to $\cE(\ZZ_3)$, 
whose conjugates in $\cE(\CC_5)$ all belong to $\cE(\ZZ_5)$, 
and whose conjugates in $\cE(\CC_p)$ belong to $\cE(\hcO_p)$, for all primes $p  > 5$.
  
If $R \le 142.388571237$, there are only finitely many $\alpha \in \cE(\QQbar)$
meeting these conditions.   
\end{example} 
\index{elliptic curve}

\begin{proof}  The given Weierstrass equation is minimal.
\index{Weierstrass equation!minimal}
At $p =2$ it has reduction type $III^*$, with $2$ rational components;  
at $p =3$ it has reduction type $I_0^*$, with $2$ rational components;  
and at $p = 5$ it has non-split multiplicative reduction ($n = 4$) with $2$ rational components 
(see Cremona \cite{Cremona}, p.133).  
\index{Cremona, John}   
We compute capacities\index{uniformizing parameter!normalizes capacity}
 with respect to the uniformizing parameter $z = x/y$.   
 
Take $E_{\infty} = y^{-1}(D(0,R))$, where $R > 0$.  
By formula (\ref{FEVYF}), 
$V_{\oORIG}(E_{\infty}) \ = \ -\frac{1}{3} \ln(R)$.  
\index{Robin constant!examples!elliptic curve}
At $p =2$, take $E_2 = \cE(\ZZ_2)$.  By formula (\ref{FEIII*}), $V_{\oORIG}(E_2) = 10/9$.  
At $p=3$,  take $E_3 = \cE(\ZZ_3)$.  By formula (\ref{FEGood}) $V_{\oORIG}(E_3) = 1/2$.  
At $p = 5$, take $E_5 = \cE(\ZZ_5)$.  By formula (\ref{FENonSplitMult2}), 
$V_{\oORIG}(E_5) = 29/140$.  For $p > 5$, take $E_p = \cE(\hcO_p)$.  
\index{Robin constant!examples!elliptic curve}
Since the given model of $\cE$ and the parameter $z$ have good reduction at $p$,
\index{good reduction} 
$V_{\oORIG}(E_p) = 0$.  

Let $\EE = \prod_{p,\infty} E_p$, and take $\fX = \{\oORIG\}$.  Then   
\begin{equation*}
V(\EE,\fX) \ = \ 
-\frac{1}{3} \ln(R) + \frac{10}{9} \ln(2) + \frac{1}{2} \ln(3) + \frac{29}{140} \ln(5) \ .
\end{equation*}   
The value of $R$ for which $V(\EE,\fX) = 0$ satisfies
\begin{equation*}
142.388571237 \ < \ R \ < \ 142.388571238 \ , 
\end{equation*} 
and the Fekete-Szeg\"o theorems \ref{aT1} and \ref{FSZii} yield the result.
\index{Fekete-Szeg\"o theorem with LRC}  
\end{proof}


\section{ The Fermat Curve} \label{FermatCurveSection}
\index{Fermat curve}\index{examples!Fermat curve}

Let $p \ge 3$ be an odd prime, and let $\zeta = e^{2 \pi i/p}$.   
In this section we will apply the Fekete-Szeg\"o theorem 
with local rationality conditions to the Fermat curve 
\index{Fekete-Szeg\"o theorem with LRC} 
\index{Fermat curve}
\begin{equation} \label{FFer1}  
\cF^p : X^p + Y^p  \ = \ Z^p \ , 
\end{equation}  
taking the ground field to be $K = \QQ$.  
Let $\fF^p/\Spec(\ZZ)$ be the corresponding scheme.  
To obtain a nontrivial set $\EE$, we make use of William McCallum's description 
(\cite{McC}) of a regular model for $\fF^p_{v_p} := \fF^p \times \Spec(\cO_{L,v_p})$,  
where $L = \QQ(\zeta)$ and $v_p$ is the place of $L$ over $p$.
The author thanks Dino Lorenzini for suggesting this example.
\index{Lorenzini, Dino} 

Writing $x = X/Z$, $y = Y/Z$, let the part of $\cF^p$ in
the coordinate patch $Z \ne 0$ be the affine curve 
\begin{equation} \label{FFer2} 
\cF^{p,0} : \ x^p + y^p  \ = \ 1 \ .
\end{equation}  
Let $\fX = \{\xi_1, \ldots, \xi_p\}$ 
be the set of points at infinity, 
where $\xi_k = (1:-\zeta^k:0)$, 
and take  $\EE  = \prod_v E_v$ where the sets $E_v$ are as follows:  
for the archimedean place, let 
\begin{equation*} 
E_{\infty} \ = \ x^{-1}(D(0,R)) 
         \ = \ \{z \in \cF^p(\CC) : |x(z)| \le R \} \ .
\end{equation*}
At the place $p$, take $E_p = \cF^{p,0}(\cO_{L,v_p})$, 
and for all the other nonarchimedean places $q$, take $E_q$ to be the $\fX$-trivial 
\index{$\fX$-trivial}
set $E_q = \cF^{p,0}(\hcO_q)$.  

Let $z$ vary over $\cF^p(\CC)$.  Writing (\ref{FFer2}) in the form  
\begin{equation*}
\prod_{k=1}^p (\frac{y}{x} + \zeta^k) \ = \ \big(\frac{1}{x}\big)^p
\end{equation*} 
we see that as $z \rightarrow \xi_k$, 
then $(y/x) + \zeta^k$ vanishes to order $p$;  
at each $\xi_k \in \fX$ we will take
the local uniformizing parameter\index{uniformizing parameter!normalizes Robin constant} to be 
\begin{equation*} 
g_{\xi_k}(z) = \frac{1}{x(z)} \ .
\end{equation*} 

\smallskip
{\bf The Green's matrix at the archimedean place.} 
\index{Green's matrix!local}
Since $E_{\infty}$ and $|1/x(z)|$ are invariant under 
the automorphisms of $\cF^p$  given by 
\begin{equation*} 
(X:Y:Z) \ \mapsto \ (\zeta^k X: \zeta^{\ell} Y:Z) \ ,
\end{equation*} 
while the $\xi_k$ are permuted by those automorphisms, 
there are numbers $A$, $B$ such that $G(\xi_k,\xi_{\ell};E_\infty) = A$ and 
$V_{\xi_k}(E_{\infty}) = B$ for all $k \ne \ell$.  Thus the
archimedean local Green's matrix is 
\index{Robin constant!examples!elliptic curve}
\index{Green's matrix!local}
\begin{equation} \label{FFer6} 
\Gamma(E_{\infty},\fX) 
        \ = \ \left( \begin{array}{cccc} 
                       B & A & \cdots & A \\
                       A & B & \cdots & A \\
                       \vdots & \vdots & \ddots & \vdots \\
                       A & A & \cdots & B \end{array} \right) 
\end{equation}  
Although we are unable to determine the numbers $A$, $B$ explicitly, 
we will see below that 
\begin{equation} \label{FFer3} 
(p-1)A + B \ = \ -\log(R) \ .
\end{equation}
This relation will enable us to determine the capacities we need.  

For the divisor $(\infty)$ on $\PP^1$, we have $x^{-1}((\infty)) = (\xi_1)+ \cdots + (\xi_p)$, 
so the pullback formula (\ref{GreenPullbackF}) shows that
\index{Green's function!pullback formula for}
for each $z \in \cF^p(\CC)$, 
\begin{equation*}
G(x(z),\infty;D(0,R)) \ = \ \sum_{\ell=1}^p G(z,\xi_\ell,E_{\infty}) \ .
\end{equation*} 
Since $G(w,\infty;D(0,R)) = \log^+(|w/R|)$ in $\PP^1$, for each $\xi_k$ we have 
\begin{eqnarray}
-\log(R) & = & \lim_{z \rightarrow \xi_k} G(x(z),\infty;D(0,R)) + \log(|1/x(z)|) \notag \\ 
         & = & V_{\xi_k}(E_{\infty}) + \sum_{\ell \ne k} G(\xi_k,\xi_{\ell};E_{\infty}) \ , 
              \label{FFer4}
\end{eqnarray} 
and (\ref{FFer3}) follows.
\index{Robin constant!examples!Fermat curve}

\smallskip
{\bf The Green's matrix at the place $p$.}
\index{Green's matrix!local}
Put $L = \QQ(\zeta)$, and let $v_p$ be the unique place of $L$ above $p$;  
thus  $\cO_{L,v_p} \cong \ZZ_p[\zeta]$.  Put $\pi_{v_p} = 1-\zeta$.
The residue field $k_v = \cO_{L,v_p}/\pi_{v_p}\cO_{L,v_p}$ 
is isomorphic to  $\FF_p$.  
Write $\cF^p_{v_p} = \cF^p \times_{\QQ} \Spec(L_{v_p})$
and $\fF^p_{v_p} = \fF^p \times_{\ZZ} \Spec(\cO_{L,v_p})$.  

McCallum (\cite{McC}, see Theorem 3, p.59; Diagram 3, p.69)
\index{McCallum, William} 
has determined a regular model for $\cF^p_{v_p}$.  Put 
\begin{equation*} 
\phi(x,y) \ = \ \frac{(x+y)^p-x^p-y^p}{p} \ .
\end{equation*}
Then $\phi(x,y)$ is a polynomial with integer coefficients, 
divisible by $xy(x+y)$.  Let $\widetilde{\FF}_p$ be the algebraic closure
of $\FF_p$;  McCallum notes that $\phi(x,-y) \pmod{p}$ 
has a factorization over $\widetilde{\FF}_p$ of the form 
\begin{equation*}
xy(x-y) \cdot \prod_{i} (x-\alpha_i y)^2 \cdot \prod_{j} (x-\beta_j y)
\end{equation*} 
in which the $\alpha_i, \beta_j \in \widetilde{\FF}_p$ are distinct,
the $\alpha_i$ belong to $\FF_p \backslash \{0,1\}$, 
and the $\beta_j$ belong to $\widetilde{\FF}_p \backslash \FF_p$.
\index{Fermat curve!McCallums's regular model for}

McCallum shows that there is a regular model $\fG^p_{v_p}/\Spec(\cO_{L,v_p})$, 
gotten by blowing up $\fF^p_{v_p}$,
whose geometric special fibre has the configuration shown in Figure \ref{figure:FermatFibre}. 
The components $\rL_0$, $\rL_1$, $\rL_{\infty}$, $\rL_{\alpha_i}$ and $\rL_{\beta_j}$ 
meeting $\rL$ are indexed by the irreducible factors of $\phi(x,-y) \pmod{p}$, 
and for each $\alpha_i$
there are $p$ components $\rF_{\alpha_i,k}$ meeting $\rL_{\alpha_i}$. 
All components are nonsingular and isomorphic to $\PP^1/\widetilde{\FF}_p$, 
and all intersections are transverse: 

\begin{figure}
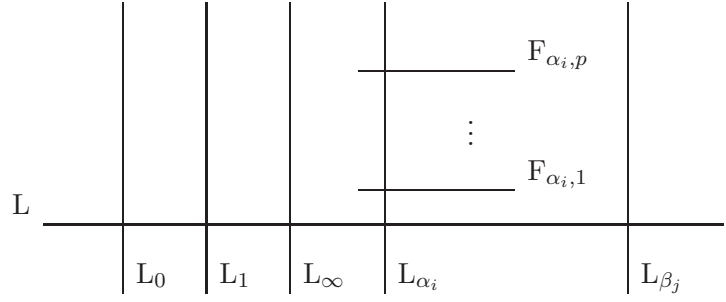
 
\centerline{
\begin{tabular}{lccc|cc|cc|cc|ccccc|cc}
  \ & \ & \ & \ &     \ & \ &     \ & \ &            \ & \ &    \   &  \ &  \  &         \         &   \  &  \  &  \    \\ 
  \ & \ & \ & \ &     \ & \ &     \ & \ &            \ & \ &    \   &  \ &  \  & $\rF_{\alpha_i,p}$  &   \  &  \  &  \    \\   
                                               \cline{10-13}                                                     
  \ & \ & \ & \ &     \ & \ &     \ & \ &            \ & \ &    \   &  \ &   \ &         \         &   \  &  \  &  \    \\
  \ & \ & \ & \ &     \ & \ &     \ & \ &            \ & \ &    \ & \vdots & \ &         \         &   \  &  \  &  \    \\ 
  \ & \ & \ & \ &     \ & \ &     \ & \ &            \ & \ &    \   &  \ &   \ & $\rF_{\alpha_i,1}$  &   \  &  \  &  \    \\   
                                               \cline{10-13}
$\rL$ & \ & \ & \ &     \ & \ &     \ & \ &            \ & \ &    \   &  \ &   \ &         \         &   \  &  \  &  \    \\
            \cline{2-17}\cline{2-17}\cline{2-17}
  \ & \ & \ & \ &     \ & \ &     \ & \ &            \ & \ &    \   &  \ &   \ &         \         &   \  &  \  &   \   \\
  \ & \ & \ & \ & $\rL_0$ & \ & $\rL_1$ & \ & $\rL_{\infty}$ & \ & $\rL_{\alpha_i}$
                                                                 &  \ &   \ &         \         &   \  & $\rL_{\beta_j}$  & \
\end{tabular} 
}
\caption{Fermat Curve Special Fibre}\index{Fermat curve special fibre} \label{figure:FermatFibre} 
\end{figure}

\noindent{The} components $\rL$, $\rL_0$, $\rL_1$,  
$\rL_{\infty}$, $\rL_{\alpha_i}$, and $\rF_{\alpha_i,j}$ 
are rational over $k_v = \FF_p$; 
each $\rL_{\beta_j}$ is rational over $\FF_p(\beta_j)$.  
Furthermore $\rL$ has multiplicity $p$ and self-intersection $-1$; 
$\rL_0$, $\rL_1$, and $\rL_{\infty}$ have multiplicity $1$ and self-intersection $-p$;  
the  $\rL_{\alpha_i}$ have multiplicity $2$ and self-intersection $-p$;  
the $\rL_{\beta_j}$ have multiplicity $1$ and self-intersection $-p$;  and the 
$\rF_{\alpha_i,k}$ have multiplicity $1$ and self-intersection $-2$.

The points of $\cF^p_{v_p}(L_{v_p})$ 
specialize to the $k_v$-rational closed points of 
the $k_v$-rational multiplicity $1$ components, 
which are not intersection points of components.  
There are $p$ such points on each of 
$\rL_0$, $\rL_1$, and the components $\rF_{\alpha_i,k}$.
Each such point lifts to a subset 
of $E_p$ isomorphic to $\pi_{v_p} \cO_{L,v_p}$, 
and $E_p = \cF^{p,0}(\cO_{L,v_p})$ is the union of those subsets.  
On the other hand, $\xi_1, \ldots, \xi_p$ 
specialize to distinct $k_v$-rational closed points of  $\rL_{\infty}$.   
 
Using Proposition \ref{CanonDistForm} and the above description of $E_p$, 
we can determine $G(z,\xi_k;E_p)$ for each $k$.  Since the computations 
\index{Green's function!examples!Fermat curve}
are somewhat tedious, and the methods are the same as those in the proof 
of Theorem \ref{EllipticNonarch}, we only give the final result:  
if $n_p$ is the number of components $\rL_{\alpha_i}$, and if 
\begin{equation}
V \ = \ \frac{1}{p} + \frac{2p-1}{(2n_p+2)p-n_p} \ , \label{FFer7}
\end{equation} 
then for points $z \in \cF^p(L_{v_p})$ specializing to $\rL_{\infty}$, 
\begin{equation*}
G(z;\xi_k;E_p) \ = \ 
 \frac{1}{p-1} \cdot \Big(V + \log_{v_p}\big(((z) \cdot (\xi_k))_{\fG^p_{v_p}}
        \big)\Big)
\end{equation*} 
where $((z) \cdot(\xi_k))_{\fG^p_{v_p}}$ is the intersection number 
of the closures of $z$ and $\xi_k$ in the model $\fG^p_{v_p}$.  
The factor $1/(p-1)$ appears because the ramification index of $L_{v_p}/\QQ_p$ is $p-1$. 
By analyzing the blowups in the construction of $\fG^p_{v_p}$, and writing
$z \equiv_{\rL_{\infty}} \xi_k$ if $z$ and $\xi_k$ specialize to the same
closed point of $\rL_{\infty}$, one further sees that 
\begin{equation*} 
\log_{v_p}\big(((z) \cdot(\xi_k))_{\fG^p_{v_p}}\big) = \left\{ \begin{array}{cl}
0 & \text{ if $z \not \equiv_{\rL_{\infty}} \xi_k$} \\
\log_{v_p}( |x(z)|_{v_p}))-1 
 & \text{ if $z \equiv_{\rL_{\infty}} \xi_k$} 
 \end{array} \right. \ .
\end{equation*} 
Since $g_{\xi_k}(z) = 1/x(z)$ for each $k$, 
the local Green's matrix at $p$ is 
\index{Green's matrix!local}
\begin{equation} \label{FFer8} 
\Gamma(E_p,\fX) 
        \ = \ \frac{1}{p-1} \left( \begin{array}{cccc} 
                       V-1 & V & \cdots & V \\
                       V & V-1 & \cdots & V \\
                       \vdots & \vdots & \ddots & \vdots \\
                       V & V & \cdots & V-1 \end{array} \right)  \ .
\end{equation}  

\smallskip
{\bf The Global Green's Matrix.}  
\index{Green's matrix!global}
For each prime $q$ of $\QQ$ with $q \ne p$,
the model $\fF^p$ has good reduction at $q$, 
\index{good reduction}
the points $\xi_k$ specialize to distinct points of the special fibre, 
and the function $1/x(z)$ specializes to a nonconstant function $\pmod{q}$.
Since $E_q$ is $\fX$-trivial, $\Gamma(E_q,\fX)$ is the zero matrix.
\index{Green's matrix!local}

Thus the Global Green's matrix is 
\index{Green's matrix!global}
\begin{equation*}
\Gamma(\EE,\fX) = \Gamma(E_{\infty},\fX) + \Gamma(E_p,\fX) \log(p) \ .
\end{equation*}  
When $\vs = {}^t(\frac{1}{p}, \cdots, \frac{1}{p}) \in \cP^p(\RR)$, 
entries of $\Gamma(\EE,\fX) \vs$ are all equal, so using (\ref{FFer3}) 
we conclude that
\begin{eqnarray*}
V(\EE,\fX) & = & \frac{1}{p} \Big( \big(B+(p-1)A\big) 
       + \frac{1}{p-1}\big(pV-1\big)\Big) \\
          & = & \frac{1}{p} \big(-\log(R) + \frac{pV-1}{p-1} \big)
\end{eqnarray*} 
Thus, by (\ref{FFer7}) and the Fekete-Szeg\"o theorems \ref{aT1} and \ref{FSZii} 
\index{Fekete-Szeg\"o theorem with LRC} 
we obtain:

\begin{theorem}  \label{FermatTheorem} Let $p$ be an odd prime.  Then on the affine
Fermat curve $x^p+y^p = 1$, \index{Fermat curve} if 
\begin{equation}
R \ > \ p^{ \frac{p(2p-1)}{(p-1)^2((2n_p+2)p-n_p)}} \ , \label{FFer9}
\end{equation} 
there are infinitely many integral points $\alpha$ 
whose $p$-adic conjugates are all rational over $L_{v_p}$ 
and whose archimedean conjugates satisfy $|x(\sigma(\alpha))| < R$.

If the inequality $(\ref{FFer9})$ is reversed, there are only finitely many.  
\end{theorem} 

For small primes, $n_p$ can be computed using Maple.\index{Maple computations} 
For $p = 2$ and $p = 5$, we have $n_p = 0$;  for all primes with $5 < p < 75$ 
except $p = 59$, we have $n_p = 2$;  for $p = 59$ we have $n_p = 13$.  
Below are some examples for the critical value of $R$: 

\vskip .1 in
\centerline{
\begin{tabular}{||c|c|ll||} 
  \hline \hline  $p$ & $n_p$ & \text{critical $R$} &  \\ \hline \hline 
                 $3$ &  $0$  & $3^{5/8}$ & $\cong 1.987013346$ \\ \hline
                 $5$ &  $0$  & $5^{9/32}$ & $\cong 1.572480664$ \\ \hline 
                 $7$ &  $2$  & $7^{91/1440}$ & $\cong 1.130851299$ \\ \hline
                 $53$ & $2$  & $53^{5565/854464}$ & $\cong 1.026195152$ \\ \hline
                 $59$ & $13$ & $59^{6903/5513596}$ & $\cong 1.005118113$ \\ \hline
                 $61$ & $2$  & $61^{7381/1310400}$ & $\cong 1.023425196$ \\ \hline
                 $73$ & $2$  & $73^{10585/2260224}$ & $\cong 1.020296147$ \\ 
  \hline \hline
\end{tabular}}
\vskip .1 in
\noindent{As} McCallum remarks, 
$n_p$ is the number of ``tame curves''\index{tame curve|ii} $C_s$ 
for which $\Jac(C_s)$ is isogenous to a factor of $\Jac(\cF_p)$.  For $p = 59$ 
the abnormally large number of tame curves means the critical value of $R$
is unusually small.  It would be interesting to know if there are 
other phenomena related to this. 
\index{McCallum, William}


\section{ The Modular Curve $X_0(p)$}  \label{ModularCurveSection}

In this section we will give an example applying the Fekete-Szeg\"o theorem 
with local rationality conditions to the modular curve $X_0(p)/\QQ$, 
\index{Fekete-Szeg\"o theorem with LRC}\index{modular curve!$X_0(p)$}\index{examples!modular curve $X_0(p)$}
where $p \ge 5$ is prime.   
The author thanks Pete Clark for help with this.  
\index{Clark, Pete L.} 

As is well known, $X_0(p)$ is the compactification of the moduli space
for pairs $(E,C)$ consisting of an elliptic curve 
and a cyclic subgroup of order $p$.  
As a Riemann surface, $X_0(p)(\CC)$ is 
\index{Riemann surface} 
gotten from $\Gamma_0(p) \backslash \fH$ 
by adjoining the `cusps' $c_0$ and $c_{\infty}$;\index{cusps}  
here  $\fH$ is the complex upper half plane and $\Gamma_0(p)$ 
is the congruence subgroup 
\begin{equation*}
\Gamma_0(p) \ = \ 
\Big\{ \left( \begin{array}{cc} a & b \\ c & d \end{array} \right) \in SL_2(\ZZ) : 
                                              c \equiv 0 \pmod{p} \Big\} \ .
\end{equation*}  
The function field of $X_0(p)/\QQ$ is $\QQ(j(z),j(pz))$ 
where $j(z)$ is the modular function\index{modular function $j(z)$} 
\begin{equation*}
j(z) \ = \ \frac{1728 g_2^3}{g_2^3 - 27 g_3^2} \ = \ \frac{1}{q} + 744 + 196884 q + \cdots \ .
\end{equation*} 
Here $X = j(z)$ and $Y = j(pz)$ satisfy the ``Modular Equation'' 
$\Phi(X,Y) = 0$,\index{Modular Equation|ii} where 
\begin{equation*}
\Phi(X,Y) \ = \ -(X^p-Y)(Y^p-X) + \sum_{\max(i,j) \le p} a_{ij} X^i Y^j \ \in \ \ZZ[X,Y] 
\end{equation*} 
and each $a_{ij}$ is divisible by $p$.  The genus of $X_0(p)$ is  
\begin{equation*}
g_p \ = \ \left\{ \begin{array}{ll} 
              (p-13)/12 & \text{if $p \equiv 1 \pmod{12}$,} \\
              (p-5)/12 & \text{if $p \equiv 5 \pmod{12}$,} \\
              (p-7)/12 & \text{if $p \equiv 7  \pmod{12}$,} \\
              (p+1)/12 & \text{if $p \equiv 11 \!\! \pmod{12}$,} 
                  \end{array} \right.  
\end{equation*}

Deligne-Rapoport\index{Deligne-Rapoport model|ii}  determined a regular model $\fM_0(p)/\Spec(\ZZ)$ for $X_0(p)$.  
It can be described as follows (see \cite{MR}, Theorem 1.1, p.175).  
First, consider the projective normalization $M_0(p)$ of $\Spec(\ZZ[X,Y]/(\Phi(X,Y))$.  
It is smooth outside the points corresponding to supersingular elliptic curves\index{supersingular points|ii} 
in characteristic $p$ with  $j \ne 0, 1728$; its special fibre at $p$ has two components, 
each isomorphic to $\PP^1$, which meet transversely at the supersingular points.  
These components will be denoted $Z_0$ and $Z_{\infty}$; 
the reduction of $j$ (that is, $X$) is a coordinate function on $Z_{\infty}$.  
If $p \equiv 2 \pmod{3}$ then $j = 0$ is supersingular in the fibre at $p$, and 
$M_0(p)$ has a singularity of type $A_3$ at the corresponding point;  
if $p \equiv 3 \pmod{4}$ then $j = 1728$ is supersingular in the fibre at $p$ and $M_0(p)$ 
has a singularity of type $A_2$ at the corresponding point.\index{supersingular points}  

The model $\fM_0(p)$ is gotten by resolving these singularities, 
introducing a chain of two components $F_1$, $F_2$ in the first case, 
and a single component $G$ in the second.  
The special fibre of $\fM_0(p)$ is reduced, 
and all its components are rational over $\FF_p$.
There are $m = g_p+1$ supersingular points,\index{supersingular points} 
each of which is rational over $\FF_{p^2}$.
The components $Z_0$ and $Z_{\infty}$ have self-intersection $-m$;  
the components $F_1$, $F_2$, and $G$ (if present) have self-intersection $-2$.
The cusps\index{cusps} are rational over $\QQ$; 
$c_0$ specializes to $Z_0$, and $c_{\infty}$ specializes to  $Z_{\infty}$.
Their images are not supersingular,\index{supersingular points}
and are the points ``at infinity'' on those components.  
 
\smallskip
We will take $\fX = \{c_{\infty},c_0\}$ to be the set of cusps,\index{cusps}   
and we will take $\EE = \prod_v E_v$, where 
\begin{equation*}
E_{\infty} \ = \ j^{-1}(D(0,R)) \ = \ \{ z \in X_0(p)(\CC) : |j(z)| \le R \} \ .
\end{equation*}
and where $E_p$ is the set of points of $X_0(p)(\QQ_p)$ specializing to
the `ordinary' (i.e non-supersingular and non-cuspidal) points\index{ordinary point|ii} of $Z_{\infty}$. 
For all the other nonarchimedean places $q$, we will take $E_q$ 
to be the $\fX$-trivial set 
\begin{equation*}
E_q \ = \fM_0(p)(\CC_p) \backslash (B(c_0,1)^- \bigcup B(c_{\infty},1)^-) \ . 
\end{equation*}  
We will take the local uniformizing parameters\index{uniformizing parameter!normalizes Robin constant} 
to be $g_{c_{\infty}}(z) = 1/j(z)$,
$g_{c_0}(z) = 1/j(pz)$.   

This set $\EE$ is chosen mainly because we can do explicit computations with it,
rather than for its intrinsic interest.  
However, it illustrates nicely how arithmetic 
and geometric information about a curve enter into capacities.

\smallskip
{\bf The Green's matrix at the archimedean place.}
\index{Green's matrix!local} 

Let $\cD = \{z \in \fH : -1/2 \le \Re(z) \le 1/2, |z| \ge 1 \}$ be the standard
closed fundamental domain for $SL_2(\ZZ)$.  As a function from $\fH$ to $\CC$,
$j(z)$ maps this region conformally onto $\CC$, 
taking the ray from $i$ to $\infty$ along the imaginary axis
to the real interval $[1728,\infty)$, with $j(i) = 1728$;  
the circular arc at the bottom of $\cD$ 
\index{arc!circular}
to the real interval $[0,1728]$ (covering it twice), 
with $j(e^{\pi i /3}) = j(e^{2 \pi i/3}) = 0$; 
and the vertical sides of $\cD$  
to the real interval $[-\infty,0]$.  
It also takes the part of the imaginary axis from $0$ to $i$ to $[1728,\infty)$.
A fundamental domain for $\Gamma_0(p)$ is given by  
\begin{equation*}
\cD(p) \ = \ \cD \cup \Big( \bigcup_{k=-(p-1)/2}^{(p-1)/2} f_k(\cD) \Big) 
\end{equation*}
where $f_k(z) = -1/(z+k)$.  
Under the quotient $\Gamma_0(p) \backslash \fH$, 
the image of the circular arc at the bottom of $\cD$ 
\index{arc!circular}
separates $X_0(p)(\CC)$ into two components,
one containing $c_{\infty}$ and the other containing $c_0$.  On the other hand, 
the image of the imaginary axis joins the cusps\index{cusps} $c_0$, $c_{\infty}$.
  
By our choice of $E_{\infty}$ and discussion above, 
it follows that when $j(z)$ is viewed as a map from $X_0(p)(\CC)$ to $\PP^1(\CC)$, 
if $R \ge 1728$ then $X_0(p)(\CC) \backslash E_{\infty}$ 
has two connected components, while if  $R < 1728$ it has one component.  

As a divisor $j^{-1}((\infty)) = p(c_0) + (c_{\infty})$, 
so the pullback formula (\ref{GreenPullbackF}) gives
\index{Green's function!pullback formula for}
\begin{equation*} 
G(j(z),\infty;D(0,R)) \ = \  p G(z,c_0;E_{\infty}) + G(z,c_{\infty};E_{\infty}) \ .
\end{equation*} 
Since $1/j(z)$ is the uniformizing parameter\index{uniformizing parameter!normalizes Robin constant} 
at $c_{\infty}$, it follows that
\index{Green's function!examples!modular curve}
\index{Robin constant!examples!modular curve}
\begin{equation*}
-\log(R) \ = \ p G(c_{\infty},c_0;E_{\infty}) + V_{c_{\infty}}(E_{\infty}) \ .
\end{equation*} 
Similarly, since $\lim_{z \rightarrow c_0} j(z)^p/j(pz) = 1$, and 
since $1/j(pz)$ is the uniformizing parameter at $c_0$, 
\begin{equation*}
-\log(R) \ = \ G(c_0,c_{\infty};E_{\infty}) + p V_{c_0}(E_{\infty}) \ .
\end{equation*} 
Hence, writing  $B(R) = G(c_{\infty},c_0;E_{\infty}) = G(c_0,c_{\infty};E_{\infty})$, 
the archimedean local Green's matrix is 
\index{Green's matrix!local}
\begin{equation} \label{FMod1} 
\Gamma(E_{\infty},\fX) \ = \ 
\left( \begin{array}{cc} -\log(R) & 0 \\
                         0 & -\frac{1}{p}\log(R) \end{array} \right) 
 + B(R) \left( \begin{array}{cc} -p & 1 \\
                         1 & -1  \end{array} \right) \ .
\end{equation} 
Here $B(R) = 0$ if $R \ge 1728$, while $B(R) > 0$ if $R < 1728$. 
It will turn out that $R > 1728$ in the situation of interest to us;  
however, note that in any case the second matrix in (\ref{FMod1}) 
is negative semi-definite.    

\smallskip
{\bf The Green's matrix at the place $p$.} 
\index{Green's matrix!local}
Using Proposition \ref{CanonDistForm} and the definition of $E_p$ as the
set of points of $X_0(\QQ_p)$ specializing to ordinary points\index{ordinary point} on 
the component $Z_{\infty}$, we can determine $G(z,c_{\infty};E_p)$ 
\index{Green's function!examples!modular curve}
and $G(z,c_0;E_p)$.  

Let $\cN_p$ be the number of $\FF_p$-rational ordinary points\index{ordinary point} on $Z_{\infty}$.
For $z \in X_0(p)(\CC_p)$, 
write $z \equiv_{Z_{0}} c_{0}$ if $z$ specializes to same point of $Z_0$ as $c_0$, 
and write  $z \equiv_{Z_{\infty}} c_{\infty}$ specializes to the same point of
$Z_{\infty}$ as $c_{\infty}$.  
Using Proposition \ref{CanonDistForm} 
and the methods in the proof of Theorem \ref{EllipticNonarch}, 
we find that 
\begin{eqnarray*} 
G(z,c_{\infty};E_p) & = & \left\{ \begin{array}{ll}
    \frac{1}{\cN_p} \frac{p}{p-1} & 
        \qquad \quad \, \text{if $z \equiv_{Z_{0}} c_{0}$,} \\
    \frac{1}{\cN_p} \frac{p}{p-1} + \log_p(|j(z)|_p) & 
         \qquad \quad \, \text{if $z \equiv_{Z_{\infty}} c_{\infty}$,}
    \end{array} \right. \\
G(z,c_0;E_p) & = & \left\{ \begin{array}{ll}
    \frac{1}{\cN_p} \frac{p}{p-1} & 
           \text{if $z \equiv_{Z_{\infty}} c_{\infty}$,} \\
    \frac{1}{\cN_p} \frac{p}{p-1} + \frac{12}{p-1} + \log_p(|j(pz)|_p) & 
          \text{if $z \equiv_{Z_{0}} c_{0}$. }
    \end{array} \right.
\end{eqnarray*}
\index{Green's function!examples!modular curve}
Here the number $12/(p-1)$ is actually the quantity 
$j_{c_0}(Z_{\infty},Z_{\infty})$ in the notation of (\ref{FCD7}), 
obtained by solving the equations (\ref{FCD4}) relating components. 
Since the special fibre of $\fM_0(p)$ at $p$ has different configurations  
according as $p \equiv 1, 5, 7, 11 \pmod{12}$, 
it is somewhat surprising that the same value arises in all cases.

It follows that the local Green's matrix at $p$ is 
\index{Green's matrix!local}
\begin{equation} \label{FMod2} 
\Gamma(E_p,\fX) \ = \ 
\left( \begin{array}{cc} 
    \frac{1}{\cN_p} \frac{p}{p-1} & \frac{1}{\cN_p} \frac{p}{p-1} \\
    \frac{1}{\cN_p} \frac{p}{p-1} & \frac{1}{\cN_p} \frac{p}{p-1} + \frac{12}{p-1} 
       \end{array} \right) 
 \ .
\end{equation} 
The number $\cN_p$ can be expressed in terms of the class number $h(-p)$
of the ring of integers of $\QQ(\sqrt{-p})$.  Indeed, $\cN_p = p - n_{ss}(\FF_p)$, 
where $n_{ss}(\FF_p)$ is the number of $\FF_p$-rational supersingular\index{supersingular points} 
points on $Z_{\infty}$.  It is known (see for example \cite{PClark}, pp.75-76) 
that
\begin{equation*}
n_{ss}(\FF_p) \ = \ \frac{h^{\prime}(-p) + h^{\prime}(-4p)}{2}
\end{equation*}
where $h^{\prime}(D)$ is the class number of the quadratic order of discriminant
$D$ if there is such an order, and is $0$ otherwise.  
Using the formula relating class numbers of 
orders in quadratic fields to those of the maximal orders 
(see \cite{LE}, Theorem 7, p.95), 
this simplifies to $n_{ss}(\FF_p) = c_p h(-p)$, where 
\begin{equation} \label{FMod3} 
c_p \ = \ \left\{ \begin{array}{ll}
        1/2 & \text{if $p \equiv 1 \pmod{8}$,} \\
        2   & \text{if $p \equiv 3 \pmod{8}$,} \\
        1/2 & \text{if $p \equiv 5 \pmod{8}$,} \\
        1   & \text{if $p \equiv 7 \pmod{8}$.} 
                  \end{array} \right.
\end{equation}
Thus $\cN_p = p - c_p h(-p)$.   
It is known that $\cN_p$ is always positive, so $E_p$ is nonempty.

\smallskip
{\bf The Global Green's Matrix.}  
\index{Green's matrix!global}
For each prime $q$ of $\QQ$ with $q \ne p$,
the model $\fM_0(p)$ has good reduction at $q$, 
\index{good reduction}
the cusps\index{cusps} $c_{\infty}, c_0$ specialize to distinct points of the special fibre, 
and the uniformizing parameters\index{uniformizing parameter!normalizes Robin constant} $g_{c_{\infty}}(z)$ and $g_{c_0}(z)$
specialize to nonconstant functions $\pmod{q}$.
Since $E_q$ is $\fX$-trivial, $\Gamma(E_q,\fX)$ is the zero matrix.
\index{Green's matrix!local}

Suppose for the moment that $R \ge 1728$;  this assumption will be justified below.  
Then $B(R) = 0$, and the global Green's matrix $\Gamma(\EE,\fX)$ is 
\index{Green's matrix!global}
\begin{eqnarray*}
\left( \begin{array}{cc} 
    -\log(R)+ \frac{1}{\cN_p} \frac{p}{p-1} \log(p) &  
              \frac{1}{\cN_p} \frac{p}{p-1} \log(p) \\
     \frac{1}{\cN_p} \frac{p}{p-1} \log(p) & 
        -\frac{1}{p}\log(R) + (\frac{1}{\cN_p} \frac{p}{p-1} + \frac{12}{p-1}) \log(p) 
            \end{array} \right) \ .
\end{eqnarray*}  
By the minimax definition of $V(\EE,\fX)$ (see formula (\ref{MinMax1}) below)\index{minimax property}
\begin{equation*} 
V(\EE,\fX) \ = \ \min_{\vs \in \cP^2(\RR)} \max_i \big(\Gamma(\EE,\fX) \vs \big)_i \ .
\end{equation*}
\index{Green's matrix!global}
Thus $V(\EE,\fX) < 0$ if and only if for some  
$\vs = {}^t(s_1,s_2) \in \cP^2(\RR)$,  
\begin{equation*} 
\left\{ \begin{array}{c}
\qquad s_1 \big(-\log(R)+ \frac{1}{\cN_p} \frac{p}{p-1} \log(p)\big)  
     + s_2 \big(\frac{1}{\cN_p} \frac{p}{p-1} \log(p)\big) \ < \ 0 \ , \\
s_1 \big(\frac{1}{\cN_p} \frac{p}{p-1} \log(p) \big) 
   + s_2 \big(-\frac{1}{p}\log(R) + \frac{1}{\cN_p} \frac{p}{p-1} 
                  + \frac{12}{p-1} \log(p)\big) \ < \ 0 \ .
        \end{array} \right. 
\end{equation*}
Equivalently, $V(\EE,\fX) < 0$ if and only if for some $s \in \RR$ 
with $0 < s < 1$, 
\begin{equation}  \label{FMod4}
\left\{ \begin{array}{c}
\qquad \frac{1}{s} \cdot \big(\frac{1}{\cN_p} \frac{p}{p-1}\big)  \ < \log_p(R) \ , \\
\frac{1}{1-s} \cdot \big(\frac{1}{\cN_p} \frac{p^2}{p-1}\big)
    + \frac{12p}{p-1}  \ < \ \log_p(R) \ . 
        \end{array} \right. 
\end{equation}  
The left side of the first inequality in (\ref{FMod4}) is decreasing with $s$,
while that in the second inequality is increasing, so the 
extremal value of $R$ is obtained when they are equal.  Solving, 
and using the Fekete-Szeg\"o Theorems \ref{aT1} and \ref{FSZii}, one obtains 
\index{Fekete-Szeg\"o theorem with LRC} 

\begin{theorem} \label{ModularEx} Let $p \ge 5$ be a prime, 
and consider the Deligne-Rapoport model\index{Deligne-Rapoport model}  $\fM_0(p)$ for modular curve $X_0(p)/\QQ$.
\index{modular curve!$X_0(p)$|ii}\index{modular curve!Deligne-Rapoport model|ii}\index{Deligne-Rapoport model}  
Put $\cN_p = p-c_p h(-p)$, where $c_p$ is as in $(\ref{FMod3})$ 
and $h(-p)$ is the class number of $\QQ(\sqrt{-p})$.  Then if 
\begin{equation} \label{FMod5}
R \ > \ p^{ \frac{1+p+12 \cN_p + \sqrt{(1 + p + 12 \cN_p)^2 - 48 \cN_p}}{2 \cN_p} 
              \cdot \frac{p}{p-1}} \ ,
\end{equation} 
there are infinitely many $\alpha \in X_0(p)(\QQbar)$ 
whose archimedean conjugates satisfy $|j(\sigma(\alpha))| < R$, 
whose $p$-adic conjugates all belong to $X_0(p)(\QQ_p)$ 
and specialize$\pmod p$ to ordinary points\index{ordinary point} in $Z_{\infty}$, 
and whose conjugates in $X_0(\CC_q)$ specialize$\mod q$ to non-cuspidal points
of $\fM_0(p)$,
for all $q \ne p$. 

If the inequality $(\ref{FMod5})$ is reversed, there are only finitely many.  
\end{theorem} 

Note that the right side of (\ref{FMod5}) is greater than $p^6$, 
and for $p \ge 5$ this is at least $15625$.  By the second part of
the theorem and the monotonicity of the sets $j^{-1}(D(0,R))$, 
the first part cannot hold for any $R < 15625$.  
This validates our assumption that $R > 1728$.

%% file: NewFSZChap3.tex
\chapter{Preliminaries} \label{Chap3}  

In this chapter we systematically lay out notation, conventions, 
and foundational material used in the rest of the paper.  

This work can be regarded as a sequel to the author's monograph 
``Capacity theory on algebraic curves''
\index{capacity theory} 
(\cite{RR1}), and we recall several results from that work.  
In particular, we consider spherical metrics\index{spherical metric} on $\PP^N$, 
the `canonical distance' $[z,w]_{\zeta}$ on an algebraic curve, 
\index{canonical distance!$[z,w]_{\zeta}$}  
sets of capacity $0$, upper Green's functions, the inner Cantor capacity,
\index{capacity!inner Cantor capacity} 
\index{Green's function!upper}
and the $L$-rational basis for algebraic functions on a curve with poles
\index{basis!$L$-rational}\index{$L$-rational basis} 
supported on a finite set $\fX$.

\section{ Notation and Conventions} \label{NotationSection} 

Throughout the paper, we write $\log(x)$ 
for $\ln(x)$.\index{$\log(x)$!same as natural logarithm|ii}\index{logarithm!$\log(x) = \ln(x)$|ii}
\label{`SymbolIndexLog'}  
 
If $K$ is a number field, $\tK$ will be a fixed algebraic closure of $K$,
\label{`SymbolIndexK'} \label{`SymbolIndextK'}
and $\tK^{\sep}$ will be the separable closure of $K$ in $\tK$.  We write $\Aut(\tK/K)$
for the group of automorphisms $\Aut(\tK/K) \cong \Gal(\tK^{\sep}/K)$. 
\label{`SymbolIndexAuttK'}
Given a place $v$ of $K$, 
let $K_v$ be the completion of $K$ at $v$, 
\label{`SymbolIndexKv'}
let $\tK_v$ be an algebraic closure of $K_v$, let
\label{`SymbolIndextKv'}
$\tK_v^{\sep}$ be the separable closure of $K_v$ in $\tK_v$,
and let $\CC_v$ be the completion of $\tK_v$. 
\label{`SymbolIndexCCv'}
If $v$ is nonarchimedean, write $\hcO_v$ for the ring of integers of $\CC_v$.  
Let $\Aut_c(\CC_v/K_v) \cong \Aut(\tK_v/K_v) \cong \Gal(\tK_v^{\sep}/K_v)$ be the
\label{`SymbolIndexAutc'}
group of continuous automorphisms of $\CC_v$ fixing $K_v$.

If $v$ is archimedean, let $|z|_v = |z|$ be the usual absolute value
on $\RR$ or $\CC$ for which the triangle inequality holds.
For $0 < x \in \RR$, write $\log_v(x) = \ln(x)$.  
If $K_v \cong \RR$, put $q_v = e$;  if $K_v \cong \CC$, put $q_v = e^2$.
\index{$\log_v(x)$!definition of|ii}\index{$q_v$!definition of|ii}
\index{logarithm!definition of $\log_v(x)$|ii}\label{`SymbolIndexqv'} 
 
If $v$ is nonarchimedean, let $|x|_v$ be the absolute value 
on $K_v$ given by the modulus of additive Haar measure.  
\index{Haar measure}
Let $\cO_v$ be the ring of integers of $K_v$,   
\label{`SymbolIndexOv'}
let $\pi_v$ be a uniformizer for  $\cO_v$,  
and let $k_v = \cO_v/\pi_v \cO_v$ be the residue field.  
Put  $q_v = \#(k_v)$\index{$q_v$!definition of|ii}  
Then $|\pi_v|_v = 1/q_v$.     
If $\Char(K) = p > 0$, then $q_v = p^{f_v}$ is a power of $p$.
If $\Char(K) = 0$ and $p$ is the rational prime under $v$, let $e_v$ and $f_v$ 
be the absolute ramification index and residue degree;  then $q_v = p^{f_v}$
and $|p|_v = (1/p)^{-[K_v:\QQ_p]} = 1/q_v^{e_v}$.  
This absolute value has a unique extension to $\CC_v$, which we will
continue to write as $|x|_v$. 
Let $\log_v(x)$ be the logarithm to the base $q_v$,\label{`SymbolIndexLogv'} 
\index{$\log_v(x)$!definition of|ii}\index{logarithm!definition of $\log_v(x)$|ii}
and let $\ord_v(x)$ be the additive valuation on $\CC_v$ associated to $|x|_v$.\label{`SymbolIndexOrdv'}  
Then $\ord_v(x) = - \log_v(|x|_v) \in \QQ$ for all $x \in \CC_v^{\times}$.

Let $\cM_K$ be the set of all places of $K$.  For $0 \ne \kappa \in K$, 
the product formula reads
\begin{equation*}
\sum_{v \in \cM_K} \log_v(|\kappa|_v) \log(q_v) \ = \  0 \ .
\end{equation*}
We will sometimes need to write the product formula multiplicatively.
To this end, define weights $D_v = \log_v(q_v)$,\index{weights!in the product formula|ii} 
so that $D_v = 1$ unless $K_v \cong \CC$, in which case $D_v = 2$.  Then 
\begin{equation*}
\prod_{v \in \cM_K} |\kappa|_v^{D_v} \ = \ 1 \ . 
\end{equation*}
This combination of normalized absolute values 
and weights is made to preserve compatibility 
with the literature in analysis concerning capacities.  

If $L/K$ is a finite extension, we use similar conventions in
defining normalized absolute values $|x|_w$ for places $w$ of $L$, 
as well as $q_w$, $D_w$ and $\log_w(x)$.  Thus for $0 \ne \lambda \in L$,
\begin{center}
$\displaystyle{\sum_{w \in \cM_L} \log_w(|\lambda|_w) \log(q_w) \ =  \ 0}$ \ , 
\qquad $\displaystyle{\prod_{w \in \cM_L} |\lambda|_w^{D_w}  \ = \ 1}$ \ . \\
\end{center}
If $v$ is a place of $K$ and $w$ is a place of $L$ over $K$,
then on  $\CC_w \cong \CC_v$
we have the (extended) absolute values $|x|_v$ and $|x|_w$.
For archimedean places, $|x|_w = |x|_v = |x|$
for all $x \in \CC_w \cong \CC_v \cong \CC$.  For nonarchimedean places,
$|x|_w = |x|_v^{[L_w:K_v]}$ for  $x \in \CC_w \cong \CC_v$.  
For all places,  $|x|_w^{D_w} \ = \ |x|_v^{[L_w:K_v]D_v}$.

\vskip .1 in

Let $F$ be a field.  
By a `variety' $\cV/F$ we mean a separated scheme $\cV$ 
of finite type over $\Spec(F)$.  
By a `curve' $\cC/F$, we mean a smooth, projective, 
connected scheme of dimension $1$ over $\Spec(F)$.
If $\cV/F$ is a variety and $L$ is a field containing $F$, we write
$\cV_L = \cV \times_F \Spec(L)$ and let 
$\cV(L) = \Hom_F(\Spec(L),\cV) \cong \Hom_L(\Spec(L),\cV_L)$ 
be the set of $L$-rational points of $\cV$.  If $\cV_L$ is irreducible,
we write  $L(\cV)$ for its function field.

If $K$ is a global field, $v$ is a place of $K$, and $\cV/K$ is a variety, 
we abbreviate  $\cV_{K_v}$ by $\cV_v$.  
Note that $\cV_v(\CC_v) \cong \cV(\CC_v)$ since $K$
is embedded in $K_v$ and $\CC_v$. 
If $E_v \subset \cV_v(\CC_v)$ is a nonempty set, then for each $f \in \CC_v(\cV)$, 
we write $\|f\|_{E_v} = \sup_{z \in E_v} |f(z)|_v$ for its sup norm.\label{`SymbolIndexSupNorm'}  

\section{ Basic Assumptions} \label{AssumptionsSection}

         Throughout the paper we will assume that:  

\vskip .1 in
         $K$ is a global field and $\cC/K$ is a curve (smooth, projective, and geometrically connected).
\label{`SymbolIndexcC'} 
Fix an embedding $\cC \hookrightarrow \PP^N_K = \PP^N/\Spec(K)$ for
an appropriate $N$, and equip $\PP^N_K$ with a system of 
homogeneous coordinates.\index{homogeneous coordinates!choice of|ii} 
For each $v$, let $\|x,y\|_v$ be the $v$-adic metric\index{spherical metric} associated to 
this embedding (see \S\ref{Chap3}.\ref{SphericalMetricSection}  below).    
For each nonarchimedean $v$, the choice of homogeneous coordinates yields gives an integral structure 
$\PP^N/\Spec(\cO_v)$.  By taking the closure of $\cC_v$ in $\PP^N/\Spec(\cO_v)$, we obtain
a model $\fC_v/\Spec(\cO_v)$.  

\vskip .05 in
         $\fX = \{x_1, \ldots, x_n\} \subset \cC(\tK)$ is a 
finite set of global algebraic points, stable under $\Aut(\tK/K)$.  
\label{`SymbolIndexfX'}
We will call the points in $\fX$ {\em poles}.

\vskip .05 in
$L = K(\fX)$, so $L/K$ is a finite normal extension.  For each place $v$ of $K$,\label{`SymbolIndexL'} 
fix\index{embedding $\iota_v : \tK \hookrightarrow \CC_v$|ii} 
an embedding $\iota_v : \tK \hookrightarrow \CC_v$ over $K$. This induces 
a distinguished place $w_v$\index{distinguished place $w_v$|ii}\index{$w_v$|see{distinguished place $w_v$}} 
of $L$ lying over $v$,\label{`SymbolIndexw_v'} 
and an embedding $\fX \hookrightarrow \cC_v(\CC_v)$.  In this way, 
we regard $\fX$ as a subset of $\cC_v(\CC_v)$, for each $v$.\index{$\fX$, viewed as embedded in $\cC_v(\CC_v)$|ii}  

\vskip .05 in
          
$E_v \subset \cC_v(\CC_v)$ is a nonempty set, for each place $v$ of $K$. 
A set $E_v$ will be called {\em $\fX$-trivial} $($for the model $\fC_v)$ 
\index{$\fX$-trivial}
if $v$ is nonarchimedean, $\fC_v$ has good reduction at $v$, the balls
\index{good reduction} 
$B(x_i,1)^- = \{z \in \cC_v(\CC_v) : \|z,x_i\|_v < 1\}$ are pairwise disjoint, 
and 
\begin{equation*}
E_v \ = \ \cC_v(\CC_v) \backslash \bigcup_{i=1}^m B(x_i,1)^- \ .  
\end{equation*}  
Equivalently, $E_v$ is $\fX$-trivial if $\fC_v$ has good reduction at $v$,
\index{$\fX$-trivial}
\index{good reduction}  
the points in $\fX$ $($identified with points of $\cC_v(\CC_v)$ 
and extended to sections of $\fC_v \otimes_{\cO_v} \Spec(\hcO_v))$ 
specialize to distinct points $\pmod v$, and $E_v$ is precisely  
the set of points of $\cC_v(\CC_v)$ which specialize to points complementary 
to $\fX \pmod v$.

We will assume the sets $E_v$ satisfy the following conditions:   

\vskip .03 in
(1)  Each $E_v$ is stable under $\Aut_c(\CC_v/K_v)$.  

(2)  Each $E_v$ is bounded away from $\fX$ in the $v$-topology.

(3)  For all but finitely many $v$,  $E_v$ is $\fX$-trivial. 
\index{$\fX$-trivial}
 
\vskip .03 in
\noindent{Here,} property (3) is independent of the choice of 
the embedding $\cC \hookrightarrow \PP^N/K$ and the choice of homogeneous coordinates on $\PP^N/K$.  
Note that $E_v$ may be closed, open, or neither.  
However, an $\fX$-trivial set is both open and closed. 
\index{$\fX$-trivial} 
 
\vskip .05 in
Put  $\EE = \prod_v E_v$.    
We will call $\EE$ a {\em $K$-rational adelic set}.   
By our assumptions,  $\EE$ is compatible with $\fX$, \index{compatible with $\fX$} in the
terminology of the Introduction.

\vskip .05 in

$U_v \subset \cC_v(\CC_v)$ is an open set containing $E_v$, 
for each place $v$ of $K$.
We do not assume that $U_v$ is bounded away from $\fX$   
or stable under $\Aut_c(\CC_v/K_v)$,
though we will reduce to that situation in the proof of the Fekete-Szeg\"o theorem. 
\index{Fekete-Szeg\"o theorem with LRC}  
The set $\UU = \prod_v U_v$ will be called an {\em adelic neighborhood} of $\EE$.

\vskip .05 in
For each $x_i \in \fX$, fix a rational function $g_{x_i}(z) \in K(x_i)(\cC)$ 
with a simple zero at $x_i$.  We require that the choices be made so that 
$g_{\sigma(x_i)}(z) = \sigma(g_{x_i})(z)$
for all $\sigma \in \Aut(\tK/K)$.  
These uniformizing parameters\index{uniformizing parameter|ii}\index{uniformizing parameter!galois equivariant system}
will be used for normalizations throughout the paper. 

\medskip
We will often deal with objects on which $\Aut(\tK/K)$ acts:  
for instance points in $\cC(\tK)$ and functions in $\tK(\cC)$.
If $L/K$ is galois, the points of $\cC(L)$
fixed by $\Aut(\tK/K)$ are the $K$-rational points, and the 
functions in $L(\cC)$ fixed by $\Aut(\tK/K)$ are the $K$-rational
functions.  Likewise, if $v$ is a place of $K$, we will often deal with objects
on which $\Aut_c(\CC_v/K_v)$ acts.   
We will use the following terminology 
(due to Cantor \cite{Can3}): 

\begin{definition}[$K$-symmetric, $K_v$-symmetric] \label{DefD1}
\index{$K$-symmetric|ii} 
\index{$K_v$-symmetric|ii}
Let $K$ be a global field, 
and let $Y$ be a collection of objects on which $\Aut(\tK/K)$ acts.  
If an element $y_0 \in Y$ is fixed by that action, 
we will say that $y_0$ is $K$-symmetric.  
If a subset $Y_0 \subset Y$ is stable under $\Aut(\tK/K)$, we will say that $Y_0$ is $K$-symmetric. 
 
Let $v$ be a place of $K$, and let $Y$ be a set on which $\Aut_c(\CC_v/K_v)$ acts. If an element $y_0 \in Y$ is fixed by that action, 
we will say that $y_0$ is $K_v$-symmetric.
Likewise, if a subset $Y_0 \subset Y$ is stable under $\Aut_c(\CC_v/K_v)$, we will say that $Y_0$ is $K_v$-symmetric.
\end{definition}


Here is an important example.  

Let $\fX = \{x_1, \ldots, x_m\} \subset \cC(\tK)$ be as above. 
Define a permutation representation of $G = \Aut(\tK/K)$ in $S_m$ by 
$x_{\sigma(i)} = \sigma(x_i)$ for each $\sigma \in G$.   
There is an induced action of $G$ on $\RR^m$:   
given $\vs = \phantom{}^t (s_1, \ldots, s_m) \in \RR^m$, 
put $\sigma(\vs) = \phantom{}^t (s_{\sigma(1)}, \ldots, s_{\sigma(m)})$. 
If $s_{\sigma(i)} = s_i$ for all $\sigma \in G$ and all $i$, 
we call $\vs$ a {\em $K$-symmetric vector}.  
Similarly,  $G$ acts on $M_m(\RR)$ by simultaneously permuting
the rows and columns.  
These two actions are compatible:  for $\vs \in \RR^m$ 
and $\Gamma \in M_m(\RR)$, we have 
$\sigma(\Gamma \vs) = \sigma(\Gamma) \sigma(\vs)$.   
If $\Gamma_{ij} = \Gamma_{\sigma(i),\sigma(j)}$ 
for all $\sigma \in G$ and all $i$, $j$, 
we call $\Gamma$ a {\em $K$-symmetric matrix}.  

\smallskip
For another example, let $Y$ be the collection of all finite sets of functions in $\tK(\cC)$.  
Then the set of uniformizing parameters\index{uniformizing parameter} $\{g_{x_i}(z)\}_{x_i \in \fX}$ 
chosen above is $K$-symmetric, since $\sigma(g_{x_i})(z) = g_{\sigma(x_i)}(z)$ for all $i$. 

\medskip
The following fact is well known, but we do not have a convenient reference for it. 

\begin{proposition} \label{SepClosureDense}  Let $K$ be a global field,
and let $\cC/K$ be a smooth, projective, connected curve.  Then 
$\cC_v(\tK_v^{\sep})$ is dense in $\cC_v(\CC_v)$ under the $v$ topology, 
for each place $v$ of $K$.  
\end{proposition}\index{$\cC(\tK_v^{\sep})$ is dense in $\cC_v(CC_v)$}   

\begin{proof}
We first show that $\tK_v^{\sep}$ is dense in $\CC_v$.  Since $\CC_v$ is the completion
of $\tK_v$, it suffices to show that $\tK_v^{\sep}$ is dense in $\tK_v$.  
If $\Char(K) = 0$ there is nothing to prove.  Assume $\Char(K) = p > 0$,
and take $\alpha \in \tK_v \backslash \tK_v^{\sep}$.  
Let $P(x) \in K_v[x]$ be the minimal polynomial of $\alpha$.  
Since $\alpha$ is inseparable over $K_v$, 
there is a polynomial $Q(x) \in K_v[x]$ such that $P(x) = Q(x^p)$.  
For each $0 \ne b \in K_v$, put $P_b(x) = P(x)+b\,x$.  
Then $P_b^{\prime}(x) = b$, so $(P_b(x),P_b^{\prime}(x)) = 1$, 
and $P_b(x)$ has distinct roots.  
In particular, its roots all belong to $\tK_v^{\sep}$.  
Fix $\varepsilon > 0$, and let $b \rightarrow 0$.  
By the continuity of roots of polynomials in $\CC_v[x]$ 
under variation of the coefficients (\cite{LANT}, p.44), 
if $|b|_v$ is sufficiently small, then $P_b(x)$ has a root $\alpha_b$
with $|\alpha_b - \alpha|_v < \varepsilon$.  
 
We next show that $\cC_v(\tK_v^{\sep})$ is dense in $\cC_v(\CC_v)$.  
Since $\cC/K$ is smooth, so is $\cC_v/K_v$, hence the 
function field $K_v(\cC_v)$ is separably generated over $K_v$ 
(see \cite{RR1}, p.21).  Since $K_v(\cC_v)/K_v$ is finitely generated
and has transcendence degree $1$, there is an $f \in K_v(\cC)$ 
such that $K_v(\cC_v)$ is finite and separable over $K_v(f)$.  
By the Primitive Element theorem,\index{Primitive Element theorem} 
there is a $g \in K_v(\cC)$ such that $K_v(\cC) = K_v(f,g)$.  
Let $F(x,y) \in K_v[x,y]$ be a nonzero polynomial of minimal degree 
for which $F(f,g) = 0$.  Regarding $F(x,y)$ as a polynomial in $y$ with
coefficients in $K_v[x]$, write
\begin{equation*}
F(x,y) \ = \ a_0(x) y^n + a_1(x) y^{n-1} + \cdots + a_n(x) \ .
\end{equation*} 
Let $R(x)$ be the resultant of $F$ and $\frac{\partial F}{\partial y}$.  
It is not the zero polynomial, since $g$ is separable over $K_v(f)$.  
There are finitely many values of $x \in \CC_v$ for which $R(x) = 0$;  
for all other $x$, the polynomial $F_x(y) = F(x,y)$ has $n$ distinct roots.  

Let $\cC_{v,1}/K_v$ be the projective closure of the 
plane curve defined by $F(x,y) = 0$.
There is a birational morphism $Q : \cC_v \rightarrow \cC_{v,1}$ defined over $K_v$. 
Let $S_1 \subset \cC_{v,1}(\CC_v)$ be the set consisting of all singular points, 
 branch points of $Q$, points ``at infinity'', 
and points where $R(x) = 0$.  Put $S = Q^{-1}(S_1) \subset \cC_v(\CC_v)$;  
both $S_1$ and $S$ are finite.  The map $Q$ 
induces a topological isomorphism from 
$\cC_v(\CC_v) \backslash S$ onto $\cC_{v,1}(\CC_v) \backslash S_1$  
which takes $\cC_v(\tK_v^{\sep}) \backslash S$ onto 
$\cC_{v,1}(\tK_v^{\sep}) \backslash S_1$.  To show that $\cC_v(\tK_v^{\sep})$ is
dense in $\cC_v(\CC_v)$, it suffices to show that 
$\cC_{v,1}(\tK_v^{\sep}) \backslash S_1$ 
is dense in $\cC_{v,1}(\CC_v) \backslash S_1$.  

Identify points of $\cC_{v,1}(\CC_v) \backslash S_1$ 
with solutions to $F(x,y) = 0$ in $\CC_v^2$, 
and fix $P = (b,c) \in \cC_{v,1}(\CC_v) \backslash S_1$.  
Then $c$ is a root of $F_b(y) = a_0(b) y^n + \cdots + a_0(b)$. 
Since $\tK_v^{\sep}$ is dense in $\CC_v$, there is a sequence 
$b_1, b_2, \cdots \in \tK_v^{\sep}$ converging to $b$.  
We can assume that none of the $b_i$ is the $x$-coordinate of a point in $S_1$.   
By the continuity of the roots of polynomials, 
there are $c_1, c_2, \cdots \in \tK_v^{\sep}$ 
such that each $c_i$ is a root of $F_{b_i}(y)$,  
and the points $P_i = (a_i,b_i)$  converge to $P$.  
\end{proof} 


\section{ The $L$-rational and $L^{\sep}$-rational bases} \label{LRationalBasisSection}

\index{$L$-rational basis|(}\index{$L^{\sep}$-rational basis|(}
Let $L = K(\fX)$, as above.  Given finite extension $F/K$, 
let $F^{\sep}$ be the separable closure of $K$ in $F$.  
In this section we will construct a $K$-symmetric (that is, $\Aut(\tK/K)$-equivariant) 
set of $L$-rational functions which we will use  
to expand functions in $\tK(\cC)$ with poles supported on $\fX$.
We will call this the $L$-rational basis. 
\index{basis!$L$-rational}\index{$L$-rational basis} 
To deal with separability issues in the global patching construction, 
\index{patching argument!global} 
we also construct a related set of $L^{\sep}$-rational functions, which we call the 
 $L^{\sep}$-rational basis. 
\index{basis!$L^{\sep}$-rational}\index{$L^{\sep}$-rational basis} 
When $\Char(K) = 0$, the $L$-rational and $L^{\sep}$-rational bases are the same,  
but when $\Char(K) = p > 0$ they are different. 

Given a $\tK$-rational divisor $D$ on $\cC$, put 
\begin{equation*}
\tGamma(D) \ = \ H^{0}(\cC_{\tK},\cO_{\cC_{\tK}}(D))
        \ = \ \{ f \in \tK(\cC) : \div(f) + D \ge 0 \} \ .
\end{equation*}
A theorem of Weil (see for example Lang \cite{L}, Theorem 5, p.174) 
\index{Weil, Andr\'e}\index{Lang, Serge} 
asserts that if $D$ is rational over a finite extension $F/K$, 
then $\tGamma(D)$ has a basis consisting of $F$-rational functions.
In scheme-theoretic terms, this comes from the faithful flatness of 
$\Spec(\tK)/\Spec(F)$.

\smallskip
We first give the construction when $\Char(K) = 0$.
Let $g \ge 0$ be the genus of $\cC$, and put $J = 2g+1$.\label{`SymbolIndexJ'} 
\label{`SymbolIndexgenus'}

Taking $D_0 = \sum_{i=1}^m J \cdot (x_i)$, 
put $\Lambda_0 = \dim_{\tK}(\tGamma((D_0))) \ge 1$.\label{`SymbolIndexLambda0'}
Noting that $D_0$ is rational over $K$, choose an arbitrary $K$-rational basis
$\{\varphi_1, \ldots, \varphi_{\Lambda_0}\}$ for $\tGamma(D_0)$.
Next, choose a representative $x_i$ from each $\Aut(L/K)$-orbit in $\fX$. 
For each $j = J, \ldots, 2J$, choose a $K(x_i)$-rational function
$\varphi_{i,j}(z) \in \tGamma(j \cdot (x_i))$ with a pole of exact order $j$ at $x_i$, 
normalized in such a way that\index{uniformizing parameter!normalizes $L$-rational basis|ii}
\begin{equation*}
\lim_{z \rightarrow x_i} \varphi_{i,j}(z) \cdot g_{x_i}(z)^j \ = \ 1 \ ,
\end{equation*}
and taking $\varphi_{i,2J}(z) = \varphi_{i,J}(z)^2$.  
Such functions exist by the Riemann-Roch Theorem.
\index{Riemann-Roch theorem} 
For each $\sigma(x_i)$ in the orbit of $x_i$ 
and each $j = J, \ldots, 2J$, put $\varphi_{\sigma(i),j} = \sigma(\varphi_{i,j})$. 
 
For each $j \ge 2J+1$, 
we can uniquely write  $j = \ell \cdot J + r$, where
$\ell$ and $r$ are integers with $\ell > 0$ and $J+1 \le r \le 2J$.  Define
\begin{equation} \label{FDQ1} 
   \varphi_{i,j}(z) = (\varphi_{i,J}(z))^{\ell} \cdot \varphi_{i,r}(z) .
\end{equation}
Then  $\varphi_{i,j}$ is rational over $K(x_i)$, 
has a simple pole of order $j$ at $x_i$ and no other poles, 
and is normalized so that
\begin{equation} \label{aF1}
\lim_{z \rightarrow x_i} \varphi_{i,j}(z) \cdot g_{x_i}(z)^j \ = \ 1 \ . 
\end{equation}
For all $\sigma \in \Aut(\tK/K)$ and all $(i,j)$, we have 
$\varphi_{\sigma(i),j} = \sigma(\varphi_{i,j})$.   

In this way we obtain a multiplicatively finitely generated, 
$K$-symmetric set of functions
\begin{equation}
\{\varphi_\lambda  : \lambda = 1 \ldots, \Lambda_0 \} \cup 
\{ \varphi_{i,j} : i = 1, \ldots, m,  j \ge J+1 \}  \label{FLRat}  
\end{equation}  
defined over $L$, which we will call the {\em $L$-rational basis} 
\label{`SymbolIndexLRat'}
\index{basis!$L$-rational|ii}\index{$L$-rational basis!definition of|ii} 
(note that the functions $\varphi_{i,J}$ do not belong to the $L$-rational basis,
\index{basis!$L$-rational}  
though they were used in constructing it).  
Each $f \in \tK(\cC)$ with poles supported on $\fX$ can be uniquely expanded
as a $\tK$-linear combination of the $\varphi_{k}$ and the $\varphi_{i,j}$, 
and if $f$ has a pole of order $n_i$ at each $x_i$,
only the $\varphi_{\lambda}$ and the $\varphi_{i,j}$ with
$j \le n_i$ are required in the expansion.  
Similarly, for each $v$, each $f \in \CC_v(\cC)$  with poles supported on $\fX$ can be uniquely expanded
as a $\CC_v$-linear combination of the $\varphi_{k}$ and the $\varphi_{i,j}$.   

When $\Char(K) = 0$, the $L^{\sep}$-rational basis 
\index{basis!$L^{\sep}$-rational}\index{$L$-rational basis!definition of|ii} 
will be the same as the $L$-rational basis.  However, we will write it as  
\index{$L^{\sep}$-rational basis!definition of|ii} 
\begin{equation}
\{\tphi_\lambda  : \lambda = 1 \ldots, \Lambda_0 \} \cup 
\{ \tphi_{i,j} : i = 1, \ldots, m,  j \ge J+1 \} \ , \label{FLsepRat}  
\end{equation}  
with $\tphi_{\lambda} = \varphi_{\lambda}$ and $\tphi_{ij} = \varphi_{ij}$
for all $\lambda$, $i$, $j$.  

\smallskip
Next suppose $\Char(K) = p > 0$. 
For each $x_i$, let $[K(x_i):K]^{\insep} = [K(x_i):K(x_i)^{\sep}]$ be the inseparable degree
of the extension $K(x_i)/K$.   
Let $J = p^A$ be the least power of $p$ 
such that $p^A \ge \max\big(2g+1,\max_i([K(x_i):K]^{\insep})\big)$.  
Put $D_0 = \sum_{i=1}^m J \cdot (x_i)$ and let $\Lambda_0 = \dim_{\tK}(D_0) \ge 1$. 
By construction, $D_0$ is $K$-rational.

We first construct the $L^{\sep}$-rational basis, 
\index{basis!$L^{\sep}$-rational} 
then we use it to construct the $L$-rational basis.
\index{basis!$L$-rational} 
To assure $L^{\sep}$-rationality, we must relax the condition that 
each $\tphi_{i,j}$ has a pole of exact order $j$ at $x_i$.  
 
Fix a $K$-rational basis $\{\tphi_1, \ldots, \tphi_{\Lambda_0}\}$ for $\tGamma(J \cdot D)$.
Then, for each $i = 1, \ldots, m$, consider the divisors $D_{i,J} = J \cdot (x_i)$ 
and $D_{i,2J} = 2J  \cdot (x_i)$.  
Both are rational over $K(x_i)^{\sep}$, 
so $\tGamma(D_{i,J})$ and $\tGamma(D_{i,2J})$ have $K(x_i)^{\sep}$-rational bases.  
By the Riemann-Roch theorem, there is a $K(x_i)^{\sep}$-rational function $\tphi_{i,J}$ 
\index{Riemann-Roch theorem} 
with a pole of order precisely $J$ at $(x_i)$.  
Since $J$ is divisible by $[K(x_i):K]^{\insep}$, the function 
$g_{x_i}(z)^J$ is rational over $K(x_i)^{\sep}$, and we can normalize $\tphi_{i,J}$
so that\index{uniformizing parameter!normalizes $L^{\sep}$-rational basis|ii} 
\begin{equation*}
\lim_{z \rightarrow x_i} \tphi_{i,J}(z) \cdot g_{x_i}(z)^J \ = \ 1 \ .
\end{equation*} 
Again by the Riemann-Roch theorem, $\dim_{\tK}(\tGamma(D_{i,2J})/\tGamma(D_{i,J})) = J$.
\index{Riemann-Roch theorem} 
Choose $K(x_i)^{\sep}$-rational functions $\tphi_{i,J+1}, \ldots, \tphi_{i,2J} \in \tGamma(D_{i,2J})$ 
in such a way that $\tphi_{i,2J} = (\tphi_{i,J})^2$ and the images of 
$\tphi_{i,J+1}, \ldots, \tphi_{i,2J}$ in $\tGamma(D_{i,2J})/\tGamma(D_{i,J})$ form a basis
for that space. For each $j > 2J$, we can uniquely write $j = \ell \cdot J + r$ with $\ell, r \in \ZZ$,
$\ell \ge 1$, and $J + 1\le r \le 2J$;  
put 
\begin{equation} \label{FMultF1}
\tphi_{i,j} \ = \ \tphi_{i,J}^{\ell} \cdot \tphi_{i,r} \ .
\end{equation} 
  
Thus, for each $j \ge J$, $\tphi_{i,j}$ is rational over $K(x_i)^{\sep}$.
For an index $j > J$ not divisible by $J$, 
the function $\tphi_{i,j}$ has a pole of order at most $J \cdot \lceil j/J \rceil$ at $x_i$,
but its pole will not in general have exact order $j$.  

We will require that the $\tphi_{i,j}$ for different $x_i$ be chosen 
in a $\Gal(L^{\sep}/K)$-equivariant way, so  
$\tphi_{\sigma(i),j} = \sigma(\tphi_{i,j})$ for each $\sigma \in \Aut(\tK/K)$.
The collection of $L^{\sep}$-rational functions 
\begin{equation}
\{\tphi_\lambda(z) : \lambda = 1 \ldots, \Lambda_0 \} \cup 
\{ \tphi_{i,j}(z) : i = 1, \ldots, m,  j \ge J+1 \}  \label{FLRatSep} 
\end{equation}  
will be called the {\em $L^{\sep}$-rational basis}.
\label{`SymbolIndexLsepRat'}  
\index{basis!$L^{\sep}$-rational|ii}\index{$L^{\sep}$-rational basis!definition of|ii} 
By construction, it is $K$-symmetric and 
multiplicatively finitely generated.  Note that although $\tphi_{i,J}$ was used in constructing
the $L^{\sep}$-rational basis, it is not an element of the basis. 
\index{basis!$L^{\sep}$-rational}  

Each $f \in \tK(\cC)$ with poles supported on $\fX$ can be uniquely expanded
in terms of the $L^{\sep}$-rational basis,
\index{basis!$L^{\sep}$-rational}  
and if $f(z)$ has a pole of order $n_i$ at each $x_i$,
only the $\tphi_{\lambda}$ and the $\tphi_{i,j}$ with
$j \le J \cdot \lceil n_i/J \rceil$ are required in the expansion.  
Similarly, for each $v$, each  $f \in \CC_v(\cC)$ can uniquely be expanded 
as a $\CC_v$-linear combination of the $\tphi_{k}$ and the $\tphi_{i,j}$.

\smallskip
We now use the $L^{\sep}$ rational basis to construct the $L$-rational basis.  
\index{basis!$L^{\sep}$-rational} 
\index{basis!$L$-rational}\index{$L$-rational basis!definition of|ii} 

Put $\varphi_\lambda = \tphi_\lambda$ for $\lambda = 1, \ldots, \Lambda_0$.  
For each $x_i$, put $\varphi_{i,J} = \tphi_{i,J}$ and $\varphi_{i,2J} = \tphi_{i,J}^2$.  
By the Riemann-Roch theorem, for each $j = J+1, \ldots, 2J-1$ 
\index{Riemann-Roch theorem} 
there is a $K(x_i)$-rational function $\varphi_{i,j}(z) \in \tGamma(j \cdot (x_i))$ 
with a pole of exact order $j$ at $x_i$. We will choose $\varphi_{i,j}(z)$
to be a $K(x_i)$-rational linear combination of $\tphi_{i,J+1}, \ldots, \tphi_{i,2J}$, 
normalized so that\index{uniformizing parameter!normalizes $L$-rational basis} 
\begin{equation*}
\lim_{z \rightarrow x_i}  \varphi_{i,j}(z) \cdot g_{x_i}(z)^j \ = \ 1 \ ;
\end{equation*} 
this is possible since $\tphi_{i,J+1}, \ldots, \tphi_{i,2J}$ span 
$\tGamma(2J \cdot (x_i))/\tGamma(J \cdot (x_i))$.  
We will require that for distinct $i$, the $\varphi_{i,j}$ be chosen so that 
$\varphi_{\sigma(i),j}= \sigma(\varphi_{i.j})$ for all $\sigma \in \Aut(L/K)$. 
For each $j > 2J$, we can uniquely write
$j = \ell \cdot J + r$, where $0 \le \ell \in \ZZ$ and $J+1 \le r \le 2J$, and we define 
\begin{equation} \label{FMultF2}
\varphi_{i,j}(z) \ = \ (\varphi_{i,J}(z))^{\ell} \cdot \varphi_{i,r}(z) .
\end{equation}
Then  $\varphi_{i,j}$ is rational over $K(x_i)$, 
has a simple pole of order $j$ at $x_i$ and no other poles, 
and is normalized so that
\begin{equation} \label{aF1b}
\lim_{z \rightarrow x_i}  \varphi_{i,j}(z) \cdot g_{x_i}(z)^j \ = \ 1 \ .
\end{equation}
For all $\sigma \in \Aut(\tK/K)$ and all $(i,j)$, we have 
$\varphi_{\sigma(i),j} = \sigma(\varphi_{i,j})$. 

\smallskip
The construction has a number of consequences, which we record for future use.
  
\begin{proposition}[Uniform Transition Coefficients] \label{TransitionProp} { \ } 
Let $\Char(K)$ be arbitrary.  Then\index{$L$-rational basis!uniform transition coefficients|ii} 

$(A)$ $\tphi_\lambda = \varphi_\lambda$ for $\lambda = 1, \ldots, \Lambda_0$,

$(B)$ For each $i = 1, \ldots, m$ and each $\ell \ge 2$, 
$\varphi_{i,\ell J} = \tphi_{i,\ell J} = \tphi_{i,J}^\ell$ is $K(x_i)^{\sep}$-rational
and belongs to both the $L$-rational and $L^{\sep}$-rational bases; 
it has a pole of exact order $\ell J$ at $x_i$, and is normalized 
so that\index{uniformizing parameter!normalizes $L$-rational basis} 
\begin{equation*}
\lim_{z \rightarrow x_i}  \varphi_{i,\ell J}(z) \cdot g_{x_i}(z)^{\ell J} \ = \ 1 \ .
\end{equation*}
For each $j \ge J+1$ we have 
$\varphi_{i,\ell J}(z) \cdot \varphi_{ij}(z) = \varphi_{i,\ell J+j}(z)$
and $\tphi_{i,\ell J}(z) \cdot \tphi_{i,j}(z) = \tphi_{i,\ell J +j}(z)$.

$(C)$ For each $i = 1, \ldots, m$, there is an invertible matrix 
$\tcB_i = (\tcB_{i,jk})_{1 \le j, k \le J}$ with coefficients in $K(x_i)$ such that 
for each $\ell \ge 1$, and each $j = 1, \ldots, J$, 
\begin{equation} \label{FForward}
\varphi_{i,\ell J + j} \ = \ \sum_{k=1}^J \tcB_{i,jk} \cdot \tphi_{i, \ell J + k} \ .
\end{equation} 
Likewise, put $\cB_i = \tcB_i^{-1}$ and write 
$\cB_i = (\cB_{i,jk})_{1 \le j, k \le J}$. Then the $\cB_{i,jk}$ belong to $K(x_i)$, and  
for each $\ell \ge 1$, and each $j = 1, \ldots, J$, 
\begin{equation} \label{FBackward}
\tphi_{i,\ell J + j} \ = \ \sum_{k=1}^J \cB_{i,jk} \cdot \varphi_{i, \ell J + k} \ .
\end{equation} 
\end{proposition}

\begin{proof}
When $\Char(K) = 0$, the proposition is trivial, since the $L$-rational
and $L^{\sep}$-rational bases coincide.
\index{basis!$L$-rational}\index{basis!$L^{\sep}$-rational} 
 
When $\Char(K) = p$, 
the low-order basis functions $\varphi_1, \ldots, \varphi_{\Lambda_0}$ coincide with 
$\tphi_1, \ldots, \tphi_{\Lambda_0}$.  
The high-order basis functions $\varphi_{i,j}$ and $\tphi_{i,j}$ 
are closely related as well. For each $i$ and each $\ell \ge 1$, 
the function $\varphi_{i,\ell J} = \tphi_{i,\ell J} = \tphi_{i,J}^\ell$
is defined over $K(x_i)^{\sep}$;  for $\ell \ge 2$, it belongs to both the 
$L$-rational and $L^{\sep}$-rational bases. 
By construction, the functions $\varphi_{i,\ell J + 1}, \ldots, \varphi_{i,(\ell+1)J}$
and $\tphi_{i,\ell J + 1}, \ldots, \tphi_{i,(\ell+1)J}$, 
and $\tphi_{i,\ell J + 1}, \ldots, \tphi_{i,(\ell+1)J}$
are $K(x_i)$-rational linear combinations of each other.
Since $\varphi_{i,J} = \tphi_{i,J}$, (\ref{FMultF1}) and (\ref{FMultF2}) 
show that for each block of $J$ functions we have the same transition coefficients.
\index{$L$-rational basis!transition matrix block diagonal|ii}  
\end{proof}

For each place $v$ of $K$, fix an embedding of $\tK$ into $\CC_v$, 
and use it to identify functions in $\tK(\cC)$ with functions in $\CC_v(\cC)$.  

\begin{corollary}[Uniform Comparison of Expansion Coefficients] \label{UnifBoundednessCor}  
Let $\Char(K)$ be arbitrary.  For each place $v$ of $K$, 
there are constants $B_v, \tB_v > 0$ 
$($with $B_v = \tB_v = 1$ for all but finitely many $v)$, 
such that for any $f \in \CC_v(\cC)$ with poles supported on $\fX$, 
if we expand $f$ using the $L$-rational and $L^{\sep}$-rational bases as 
\begin{equation*} 
f  =  \sum_{i,j} A_{i,j} \varphi_{i,j} + \sum_{\lambda} A_{\lambda} \varphi_{\lambda} \ , \quad 
f  =  \sum_{i,j} \tA_{i,j} \tphi_{i,j} 
     + \sum_{\lambda} \tA_{\lambda} \tphi_{\lambda} \ , 
\end{equation*} 
then 
\begin{equation*} 
\max_{i,j,\lambda}(|A_{i,j}|_v, |A_{\lambda}|_v) \ \le \
\tB_v \cdot \max_{i,j,\lambda}(|\tA_{i,j}|_v, |\tA_{\lambda}|_v)  
\end{equation*}
and
\begin{equation*}
\max_{i,j,\lambda}(|\tA_{i,j}|_v, |\tA_{\lambda}|_v) \ \le \
B_v \cdot \max_{i,j,\lambda}(|A_{i,j}|_v, |A_{\lambda}|_v) \ .
\end{equation*} 
\end{corollary}\index{$L$-rational basis!comparability of expansion coefficients|ii}  

\begin{corollary}[Rationality Properties of Expansion Coefficients] \label{CoeffRationalityCor}
Let $\Char(K)$ be arbitrary.  Suppose $f \in K(\cC)$ is a $K$-rational function with poles supported 
on $\fX$.  When $f$ is expanded in terms of the $L$-rational and $L$-rational bases as
\begin{equation} \label{fExpansions}
f  =  \sum_{i=1}^m \sum_{j=J+1}^{N_i} A_{i,j} \varphi_{i,j} + \sum_{\lambda = 1}^{\Lambda_0} A_{\lambda} \varphi_{\lambda} \ , \quad 
f  =   \sum_{i=1}^m \sum_{j=J+1}^{\tN_i} \tA_{i,j} \tphi_{i,j} 
     + \sum_{\lambda = 1}^{\Lambda_0} \tA_{\lambda} \tphi_{\lambda} \ , 
\end{equation}
then each $A_{ij}$ belongs to $K(x_i)$, each $A_{\lambda}$ belongs to $K$, each $\tA_{ij}$ belongs to $K(x_i)^{\sep}$, each $\tA_{\lambda}$ belongs to $K$, and the $A_{ij}$, $A_\lambda$, $\tA_{ij}$, 
and $\tA_{\lambda}$ are $K$-symmetric.  

Similarly, for each place $v$ of $K$, if $f \in K_v(\cC)$ is a $K_v$-rational function with poles
supported on $\fX$, when $f$ is expanded in terms of the $L$-rational and $L$-rational bases as
in $(\ref{fExpansions})$, then each $A_{ij}$ belongs to $K_v(x_i)$, each $A_{\lambda}$ belongs to $K_v$, each $\tA_{ij}$ belongs to $K_v(x_i)^{\sep}$, each $\tA_{\lambda}$ belongs to $K_v$,
and the $A_{ij}$, $A_\lambda$, $\tA_{ij}$, 
and $\tA_{\lambda}$ are $K_v$-symmetric.  
\index{$K_v$-symmetric!set of functions}
\end{corollary}\index{$L$-rational basis!rationality of expansion coefficients|ii} 

\begin{proof} We only give the proof in the global case, since the local case is similar.
Let $f \in K(\cC)$ be a $K$-rational function with poles supported on $\fX$.

First consider the expansion of $f$ in terms of the $L^{\sep}$-rational basis.  
\index{basis!$L^{\sep}$-rational} 
Since $f$, the $\tphi_{ij}$, and the $\tphi_{\lambda}$ are all defined over $L^{\sep}$, 
the $\tA_{ij}$ and $\tA_{\lambda}$ belong to $L^{\sep}$.  Since 
the $\tphi_{ij}$ are defined over $K(x_i)^{\sep}$ and are galois-equivariant, 
and the $\tphi_{\lambda}$ are defined over $K$, it follows from invariance of $f$
under $\Gal(L^{\sep}/K)$ that each $\tA_{ij}$ belongs to $K(x_i)^{\sep}$ 
and each $\tA_{\lambda}$ belongs to $K$, 
and as a collection, the $\tA_{ij}$ and $\tA_{\lambda}$ are $K$-symmetric. 

Next consider the expansion of $f$ in terms of the $L$-rational basis. 
\index{basis!$L$-rational} 
When the $L^{\sep}$-rational basis is expressed in terms of the $L$-rational basis, 
\index{basis!$L^{\sep}$-rational} 
for each $(i,j)$, $\tphi_{ij}$ is a $K(x_i)$-linear combination of the $\varphi_{ik}$, 
and for each $\lambda$, $\phi_{\lambda} = \tphi_{\lambda}$.
It follows that for each $(i,j)$,
$A_{ij}$ is a $K(x_i)$-linear combination of the $\tA_{ik}$, 
and for each $\lambda$, $A_{\lambda} = \tA_{\lambda}$. 
Since each $\tA_{ik}$ is $K(x_i)^{\sep}$-rational, it follows that each $A_{ij}$ is $K(x_i)$-rational.
Similarly, each $A_{\lambda}$ is $K_v$-rational, and as a collection the $A_{ij}$ and $A_{\lambda}$ 
are $K$-symmetric.   
\end{proof}

\begin{corollary}[Good Reduction Almost Everywhere] \label{GoodRedAECor}
There is a finite set $S$ of places of $K$, 
such that for each $v \notin S$, $\fC_v$ has good 
reduction at $v$ and each of $\varphi_\lambda$, $\varphi_{ij}$, $\tphi_\lambda$, 
and $\tphi_{i,j}$  specializes to a well defined\index{$L$-rational basis!has good reduction almost everywhere|ii}  
non-constant function on $\fC_v \pmod v$.  
\end{corollary}\index{$L$-rational basis!has good reduction almost everywhere} 

\begin{proof} This follows from the fact that the $L$-rational and $L^{\sep}$-rational bases are   
multiplicatively finitely generated.\index{$L$-rational basis!multiplicatively finitely generated} 
\end{proof}   

For each place $v$ of $K$, let $U_v$ be the neighborhood of $E_v$ chosen 
in \S\ref{Chap3}.\ref{AssumptionsSection}.  
For any $\varphi \in \CC_v(\cC)$ with poles supported on $\fX$,
let $\|\varphi\|_{U_v} = \sup_{x \in U_v} (|\varphi(x)|_v)$ be the
sup norm.  

\begin{proposition}[Uniform Growth Bounds] \label{GBL}  
Suppose each $U_v$ is bounded away from $\fX$ in the $v$-topology, 
and that $U_v$ is $\fX$-trivial, for all but finitely many $v$. 
\index{$\fX$-trivial}\index{$L$-rational basis!uniform growth bounds|ii} 
Then for each $v$, there is a constant $C_v > 0$ such that 
\begin{equation*}
\left\{ \begin{array}{ll}
  \|
  \varphi_{k} \|_{U_v}, \ \|\tphi_{\lambda} \|_{U_v} \le C_v&
                \text{for all $\lambda = 1, \ldots, \Lambda_0$ \ ,} \\
  \| \varphi_{i,j} \|_{U_v}, \ \|\tphi_{i,j} \|_{U_v}\le C_v^j &
                \text{for all $i$ and all $j > J$ \ .}
         \end{array} \right.
\end{equation*}
Moreover, for all but finitely many $v$,
we can take $C_v = 1$.
\end{proposition}

\begin{proof}  Since the $L$-rational and $L^{\sep}$-rational bases are multiplicatively
\index{$L$-rational basis!multiplicatively finitely generated} 
finitely generated, the proposition is immediate from the construction and our assumption that $U_v$ is bounded away from $\fX$ and is $\fX$-trivial
\index{$\fX$-trivial}
for all but finitely many $v$.  
\end{proof} 

\smallskip
For an example comparing the $L$-rational and $L^{\sep}$-rational bases 
when $\Char(K) = p > 0$, 
let $K = \FF_p(t)$ where $t$ is transcendental over $\FF_p$.
Take $\cC = \PP^1/K$ and identify $\PP^1$ with $\AA^1 \cup \{\infty\}$, 
using $z$ as the standard coordinate function on $\AA^1$.  
Take $\fX = \{x_1\}$, where $x_1 = t^{1/p}$ in affine coordinates.  
Then $L = K(\fX) = K(t^{1/p})$ and $L^{\sep} = K$.  
Choose $g_{x_1} = z - t^{1/p}$, noting that $(g_{x_i})^p = z^p - t$
is $K$-rational.  
 
Then $J = p$, $D_0 = p \cdot (x_1)$,
 and $\Lambda_0 = \dim_{\tK}(\tGamma(D_0)) = p+1$.   
We can take the low-order part of the $L$-rational and $L^{\sep}$
rational bases to be 
\begin{equation*}
\big\{\varphi_1, \ldots, \varphi_{p+1} \big\} 
\ = \ \big\{\tphi_1, \ldots, \tphi_{p+1} \big\} 
\ = \ \big\{1, \frac{1}{z^p-t}, \frac{z}{z^p-t}, \ldots, \frac{z^{p-1}}{z^p-t} \big\}  \ . 
\end{equation*} 
For the high-order part of the $L$-rational basis  
\index{basis!$L$-rational} 
we can take $\varphi_{1,j} = 1/(z-t^{1/p})^j$ for $j = p, \ldots, 2p-1$,   
and for the the high-order part of the $L^{\sep}$-rational basis 
\index{basis!$L^{\sep}$-rational}  
we can take $\tphi_{1,p} = 1/(z^p-t)$ and $\tphi_{i,j} = z^{2p-j}/(z^p-t)^2$ 
for $j = p+1, \ldots, 2p$. 
Thus in general for $j > p$, 
if we write $j = \ell \cdot p - s$ with $0 \le s < p$, 
then 
\begin{equation*} 
\varphi_{1,j} = \frac{1}{(z - t^{1/p})^j} \quad \text{and} 
\quad \tphi_{1,j} \ = \ \frac{z^{s}}{(z^p-t)^{\ell}} \ .
\end{equation*} 
Observe that $\tphi_{1,j}$ has a pole of order $p \lceil j/p \rceil$ 
at $x_1$.  
For each $\ell > 1$, 
$\varphi_{1,\ell p} = \tphi_{1, \ell p} = 1/(z^p-t)^\ell$ is $K$-rational,
with a pole of exact order $\ell p$ at $x_1$. 
\index{$L$-rational basis|)}\index{$L^{\sep}$-rational basis|)} 


\section{ The Spherical Metric and Isometric Parametrizability} \label{SphericalMetricSection}

Consider $\PP^N/K$, equipped with a system of homogeneous coordinates
$x_0, \ldots,x_N$. Write $\AA^N_k$ for affine patch on which $x_k \ne 0$.   
There is a natural metric
$\|z,w\|_v$ on $\PP^N_v(\CC_v)$ called the $v$-adic spherical metric\index{spherical metric|ii}
(see \cite{RR1}, \S 1.1):

Write $z = (z_0:\cdots:z_N)$, $w = (w_0:\cdots:w_N)$.
If $v$ is archimedean, and we fix an isomorphism $\CC_v \cong \CC$,
then $\|z,w\|_v$ is the chordal distance associated
to the Fubini-Study metric:  explicitly, 
\index{Fubini-Study metric}\index{chordal distance|ii}\label{`SymbolIndexChordalDist'}
\begin{equation*}
\|z,w\|_v \ = \ \frac{ \sqrt{ \sum_{0 \le i < j \le N} |z_i w_j - w_i z_j|^2 } }
                         { \sqrt{ \sum_i |z_i|^2} \sqrt{\sum_j |w_j|^2} } \ .
\end{equation*}
It has the following geometric interpretation.  Let $\varphi$ be the length
of the geodesic from $z$ to $w$ under the 
Fubini-Study metric on $\PP^N_v(\CC)$,
\index{Fubini-Study metric|ii}
so $0 \le \varphi \le \pi$.  Then $\|z,w\|_v = \sin(\varphi/2)$, the length
of a chord subtending an central arc of measure $\varphi$ 
in a circle of diameter $1$ (see \cite{RR1}, p.26).  
When $N=1$, and if we write $\PP^1(\CC) = \CC \cup \{\infty\}$,
identifying the affine patch $\AA_0^1(\CC)$ with $\CC$, then for $z = (1:z_1)$ and $w = (1:w_1)$,  
\begin{equation*}
\|z,w\|_v \ = \ \frac{|z_1 - w_1|}
            {\sqrt{1+|z_1|^2} \sqrt{1+|w_1|^2}} 
\end{equation*}
is the usual chordal distance\index{chordal distance|ii}  on $\PP^1$.

If $v$ is nonarchimedean, then $\|z,w\|_v$ is defined by
\begin{equation*}
\|z,w\|_v \ = \ \frac{\max_{0 i \le j \le N} |z_i w_j - w_i z_j|_v}
                         {\max_i(|z_i|_v) \max_j(|w_j|_v)} \ .
\end{equation*}                         
If $z$, $w$ belong to the affine patch $\AA_0^N$ and are
scaled so that $z_0 = w_0 = 1$, then 
\begin{equation*}
\|z,w\|_v \ = \ \frac{  \max_{1 \le i \le N} |z_i - w_i|_v   }
       { (\max(1,\max_{1 \le i \le N} |z_i|_v)
         (\max(1,\max_{1 \le j \le N} |w_j|_v) } \ .
\end{equation*}
In particular, if $z,w \in \AA_0^N(\hcO_v)$,
then $\|z,w\|_v = \max_{1 \le i \le N} (|z_i-w_i|_v)$.

If $v$ is archimedean, then $\|z,w\|_v$ is invariant under the
action of the unitary group $U(N+1,\CC)$ on $\PP^N_v(\CC)$.  
If $v$ is nonarchimedean, then $\|z,w\|_v$ in invariant under the
action of $\GL(N+1,\hcO_v)$ on $\PP^N_v(\CC_v)$.  

\smallskip
Clearly  $0 \le \|z,w\|_v \le 1$ for all $z, w \in \PP^N_v(\CC_v)$, 
with $\|z,w\|_v = 0$ if and only if $z = w$.  
Furthermore, $\|z,w\|_v = \|w,z\|$.  If $v$ is archimedean,
$\|z,w\|_v$ satisfies the triangle inequality;  if $v$ is nonarchimedean,
it satisfies the ultrametric inequality (\cite{RR1}, p. 26).  
In particular, $\|z,w\|_v$ is a metric on $\PP^N_v(\CC_v)$.   

A deeper fact is that for each $\zeta \in \PP^N_v(\CC_v)$, 
the function 
\begin{equation} \label{FNCC1} 
\[z,w\]_{\zeta} \ := \ \frac{\|z,w\|_v}{\|z,\zeta\|_v \|w,\zeta\|_v}
\end{equation}
is a metric on $\PP^N_v(\CC_v) \backslash \{\zeta\}$.  
If $v$ is archimedean,
$\[z,w\]_{\zeta}$ satisfies the triangle inequality;  
if $v$ is nonarchimedean,
it satisfies the ultrametric inequality.  
For this, see (\cite{RR1}, Theorem 2.5.1, p.122).

\vskip .1 in 
If $\cC/K$ is a smooth curve and $\iota : \cC \hookrightarrow \PP^N$
is a projective embedding, then we get an induced metric $\|z,w\|_v$
on $\cC_v(\CC_v)$, for each $v$.

It can be shown that $\|z,w\|_v$ is a Weil distribution for the diagonal
\index{Weil!distribution} 
divisor on $\cC_v \times \cC_v$ (\cite{RR1}, Theorem 1.1.1, p.27).
From this, it follows that if $\iota_1 : \cC \hookrightarrow \PP^{N_1}$ and
$\iota_2 : \cC \hookrightarrow \PP^{N_2}$ are two embeddings,
and $\|z,w\|_{v,1}$, $\|z,w\|_{v,2}$ are the corresponding metrics 
on $\cC_v(\CC_v)$,
they are equivalent: there are constants $C_1, C_2 > 0$ such that
\index{spherical metric!from different embeddings comparable|ii}
for all $z, w \in \cC_v(\CC_v)$,  
\begin{equation}  \label{FBB1}
C_1 \|z,w\|_{v,1} \ \le \ \|z,w\|_{v,2} \ \le \  C_2 \|z,w\|_{v,1} \ .
\end{equation}

We will call any such metric $\|z,w\|_v$ on $\cC_v(\CC_v)$ a spherical metric\index{spherical metric|ii}.  
We use the following notation for `discs' in $\CC_v$, and `balls' in $\cC_v(\CC_v)$:  
\begin{eqnarray*}
\begin{array}{ll} 
D(a,r)^- = \{z \in \CC_v : |z-a|_v < r\}, &  
D(a,r) = \{z \in \CC_v : |z-a|_v \le r\} \ ; \\
B(a,r)^- = \{z \in \cC_v(\CC_v) : \|z,a\|_v < r\}, &   
B(a,r) = \{z \in \cC_v(\CC_v) : \|z,a\|_v \le r\}.
\end{array} 
\end{eqnarray*}
\label{`SymbolIndexClosedDisc'}\label{`SymbolIndexOpenDisc'}
\label{`SymbolIndexClosedBall'}\label{`SymbolIndexOpenBall'}

If $\zeta \in \cC_v(\CC_v)$, and $g_{\zeta}(z) \in \CC_v(\cC)$
\label{`SymbolIndexZeta'}
\label{`SymbolIndexgZeta'}
is a uniformizer at $\zeta$, then there is a constant
$C_{\zeta} > 0$ such that
\begin{equation} \label{FMBA1}
\lim_{z \rightarrow \zeta}
  \frac{ |g_{\zeta}(z)|_v} {\|z,\zeta\|_v} \ = \ C_{\zeta} \ .
\end{equation}
This follows from the nonsingularity of the curve $\cC$, and definition of
$\|z,w\|_v$ in terms of local coordinate functions.

\begin{definition} \label{IsoParamDef} Let $v \in \cM_K$ be nonarchimedean.  
An open ball $B(a,r)^- \subset \cC_v(\CC_v)$ is 
{\em isometrically parametrizable} if it is contained 
\index{isometrically parametrizable ball|ii}
in some affine patch $\AA^N_k$, and there are power series 
$\lambda_1(z), \ldots, \lambda_N(z) \in \CC_v[[z]]$ 
converging on the disc $D(0,r)^-$ 
such that the map $\Lambda: D(0,r)^- \rightarrow B(a,r)^-$ 
given in affine coordinates on $\AA^N_k$ by 
$\Lambda(z) = (\lambda_1(z), \ldots, \lambda_N(z))$ 
is a surjective isometry:  
$\Lambda(D(0,r)^-) = B(a,r)^-$ and 
for all $x,y \in D(0,r)^-$, 
\begin{equation*}
\|\Lambda(x),\Lambda(y)\|_v \ = \ |x-y|_v \ . 
\end{equation*} 
We call the map $\Lambda$ an {\em isometric parametrization}.
\index{isometric parametrization} 
If $F_u \subseteq \CC_v$ is a field such that each $\lambda_i(z) \in F_u[[z]]$,
we say that $\Lambda$ is $F_u$-rational.

A closed ball $B(a,r)$ will be called {\em isometrically parametrizable}
\index{isometrically parametrizable ball|ii} 
if it is contained
in an isometrically parametrizable open ball $B(a,r_1)^-$ for some $r_1 > r$.   
\end{definition} 

\vskip .1 in

If $v$ is nonarchimedean, and $\cC$ is embedded in $\PP^N$, 
then all sufficiently small balls with respect 
to the corresponding spherical metric $\|z,w\|_v$\index{spherical metric}
are isometrically parametrizable:  
\index{isometrically parametrizable ball|ii}

\begin{theorem} \label{IsoParamThm} Let $v \in \cM_K$ be nonarchimedean.  
Then there is a number $0 < R_v \le 1$, depending only on $v$ and the embedding 
$\cC \rightarrow \PP^N$, such that each ball $B(a,r)^- \subset \cC_v(\CC_v)$ 
with $0 < r \le R_v$ is isometrically parametrizable. 
\index{isometrically parametrizable ball}    
If $\cC$ has good reduction at $v$ for the given embedding, we can take $R_v = 1$.
\index{good reduction}

If $a \in \cC_v(\CC_v)$ and $0 < r \le R_v$, 
then for any point $a_0 \in B(a,r)^-$ and any complete field $F_u \subseteq \CC_v$ 
such that $a_0 \in \cC_v(F_u)$ and $K_v \subseteq F_u$, 
there is an $F_u$-rational isometric parametrization $\Lambda : D(0,r)^- \rightarrow B(a,r)^-$  
\index{isometric parametrization}   
with $\Lambda(0) = a_0$.

For any isometric parametrization $\tLambda : D(0,r)^- \rightarrow B(a,r)^-$,
\index{isometric parametrization}
there is an index $i_0$ such that for all $x, y \in D(0,r)^-$ 
\begin{equation} \label{FOneComp} 
\|\tLambda(x),\tLambda(y)\|_v \ = \ |\tlambda_{i_0}(x)-\tlambda_{i_0}(y)|_v \ = \ |x-y|_v \ .
\end{equation} 
Furthermore, if $F_u$ is any field such that $\tLambda$ is $F_u$-rational, 
and if $L_w \subset \CC_v$ is a complete field containing $F_u$,
then  $\tLambda(D(0,r)^- \cap L_w) = B(a,r)^- \cap \cC_v(L_w)$. 
\end{theorem} 

\begin{proof} The existence of the number $R_v$ and the existence of isometric parametrizations
\index{isometric parametrization} 
with the specified properties relative to $a_0$ and $F_u$ are proved in (\cite{RR1}, Theorem 1.2.3, p.31).  
The fact that we can take $R_v = 1$ when $\cC$ has good reduction at $v$ is proved in 
\index{good reduction}
(\cite{RR1}, Corollary 1.2.4, p.39).

\smallskip
Now let $\tLambda = (\tlambda_1, \ldots, \tlambda_N) : D(0,r)^- \rightarrow B(a,r)^-$ be 
an arbitrary isometric parametrization.  
\index{isometric parametrization}
We first show that there is an index $i_0$ for which (\ref{FOneComp}) holds.

Put $a_0 = \tLambda(0)$.  After replacing $a$ by $a_0$ and changing coordinates by a translation,
we can assume that $a = a_0 = \vORIG$.  
For each $i$, write 
\begin{equation*}
\tlambda_i(z) \ = \ \sum_{n=1}^{\infty} a_{i,n} z^n \ .
\end{equation*}  
Since $\lambda_i$ converges on $D(0,r)^-$, 
for each $R$ with $0 < R < r$ we have $\lim_{n \rightarrow \infty} |a_{i,n}|R^n = 0$. 
Furthermore, by the Maximum Modulus Principle for power series, if $R \in |\CC_v^{\times}|_v$ then  
\index{Maximum principle!nonarchimedean!for power series}
\begin{equation} \label{FMMPF}
\|\tlambda_i\|_{D(0,R)} \ = \ \max_n |a_{i,n}|_v R^n \ .
\end{equation} 

Since $\tLambda$ is an isometric parametrization with $\tLambda(0) = \vORIG$, 
\index{isometric parametrization}
for each $i$ and each $R$ with $0 < R < r$, if $x \in D(0,R)$ we have 
\begin{equation*} 
|\tlambda_i(x)|_v = |\tlambda_i(x) - \tlambda_i(0)|_v \ \le \ \|\tLambda(x),\tLambda(0)|_v 
\ = \ |x-0|_v \ \le \ R \ . 
\end{equation*}
This means that $\|\tlambda_i\|_{D(0,R)} \le R$, for each $i$ and $R$.  
Similarly, for all $x, y \in D(0,r)^-$
\begin{equation} \label{Flamineq}  
|\tlambda_i(x) - \tlambda_i(y)|_v \ \le \ \|\tLambda(x),\tLambda(y)\|_v 
\ = \ |x-y|_v \ . 
\end{equation}
On the other hand for each $R$ with $0 < R < r$, and each $x \in D(0,R)$, we have 
\begin{equation*}
\max_{1 \le i \le N} |\tlambda_i(x)|_v \ = \ \|\tLambda(x),\tLambda(0)\|_v \ = \ |x-0|_v \ = \ |x|_v \ .
\end{equation*} 
Letting $|x|_v$ approach $R$ and using the Pigeon-hole Principle,\index{Pigeon-hole Principle}  
we see for each $R$ there is some $i$ for which $\|\tlambda_i\|_{D(0,R)} = R$.  

Take a sequence $0 < R_1 < R_2 < \cdots < r$ with $\lim_{\ell \rightarrow \infty} R_{\ell} = r$, 
such that each $R_{\ell} \in |\CC_v^{\times}|_v$.  By the Pigeon-hole Principle,\index{Pigeon-hole Principle}
there is an $i_0$ such that $\|\tlambda_{i_0}\|_{D(0,R_{\ell})} = R_{\ell}$ for infinitely many $\ell$.  
After replacing $\{R_{\ell}\}_{\ell \ge 1}$ by a subsequence, we can assume this holds for all $\ell$.
For notational convenience, relabel the coordinates so that $i_0 =1$.

For each $\ell$, (\ref{FMMPF}) shows that 
$|a_{1,n}|_v R_{\ell}^n \le \|\tlambda_1\|_{D(0,R_{\ell})} = R_{\ell}$ for each $n$,
with equality for some $n$.  
Let $n(\ell)$ be the maximal index for which $|a_{1,n}|_v R_{\ell}^n = R_{\ell}$.  We claim that 
$n_{\ell} = 1$ for each $\ell$.  Suppose to the contrary that $n_{\ell} \ge 2$ for some $\ell$.  
Since the function $f_{\ell}(R) = |a_{1,n_{\ell}}|_v R^{n_{\ell}} - R$ is convex upward for $R > 0$, 
is negative for small positive $R$, and satisfies $f_{\ell}(R_{\ell}) = 0$, 
it must be positive for $R > R_{\ell}$.  
Hence for each $R \in |C_v^{\times}|_v$ with $R_{\ell} < R < r$ we would have 
\begin{equation*} 
\|\tlambda_1\|_{D(0,R)} \ \ge \ |a_{i,n_{\ell}}|_v R^{n_{\ell}} \ > \ R \ , 
\end{equation*}  
contradicting $\|\tlambda_1\|_{D(0,R)} \le R$.  Thus 
$|a_{1,1}|_v R_{\ell} = R_{\ell}$ and $|a_{1,n}|_v R_{\ell}^n < R_{\ell}$ for each $n \ge 2$.   
Letting $\ell \rightarrow \infty$ and using the convexity of $|a_{1,n}|_v R^n$ 
for $n \ge 2$, we see that for $0 < R < r$ 
\begin{equation} \label{FCoeffIneqs}
\left\{ \begin{array}{l}
\text{$|a_{1,1}|_v = 1$\ ,} \\  
\text{$|a_{1,n}|_v R^{n-1} < 1$ \  for $n \ge 2$\ .}
        \end{array} \right. 
\end{equation} 
Take $x, y \in D(0,r)^-$, and choose $R$ with $\max(|x|_v,|y|_v) < R < r$.  Then  
\begin{eqnarray}
|\tlambda_1(x)-\tlambda_1(y)|_v & = & |\sum_{n=1}^{\infty} a_{1,n}(x^n-y^n)|_v \notag \\
& = & |x-y|_v \cdot |a_{1,1} + \sum_{n=2}^{\infty} a_{1,n}(\sum_{k=0}^{n-1}x^k y^{n-1-k})|_v 
\ = \ |x-y|_v \ , \label{FFullEq}
\end{eqnarray}
where the last step uses (\ref{FCoeffIneqs}) and the ultrametric inequality.
Combining (\ref{Flamineq}) and (\ref{FFullEq}) yields (\ref{FOneComp}).

Now let $F_u$ be any field over which $\tLambda(z)$ is rational.  
Since $a_{1,1} \ne 0$, under composition of power series $\tlambda_1(z)$ 
has a formal inverse $\tlambda_1^{-1}(z) \in F_u[[z]]$.  By (\ref{FCoeffIneqs}) and a simple recursion, 
$\tlambda_1^{-1}(z)$ converges on $D(0,r)^-$: for each $x \in D(0,r)^-$,   
\begin{equation*}
\tlambda_1^{-1}(\tlambda_1(x)) \ = \ \tlambda_1(\tlambda_1^{-1}(x)) \ = \ x \ .
\end{equation*} 
Thus $\tlambda_1$ and $\tlambda_1^{-1}$ induce inverse isometries from $D(0,r)^-$ onto itself. 
  
If $F_u \subseteq L_w \subseteq \CC_v$, then $\tLambda$ is $L_w$-rational.
Suppose $L_w$ is complete.  In this case $a_0 := \tLambda(0) \in \cC_v(L_w)$, 
and the initial reductions allowing us to assume $a_0 = \vORIG$ 
do not affect the $L_w$-rationality of $\tLambda$.   
Clearly $\tLambda(D(0,r)^- \cap L_w) \subseteq B(a,r)^- \cap \cC_v(L_w)$.  For the opposite
containment, take $b \in B(a,r)^- \cap \cC_v(L_w)$.  Write $b = (b_1, \ldots, b_N)$.  
Then $b_1 \in D(0,r)^- \cap L_w$;  put $x_1 = \tlambda_1^{-1}(b_1)$.  
Since $\tlambda_1^{-1}$ is $L_w$-rational, it follows that $x_1 \in D(0,r)^- \cap L_w$.  
We claim that $\tLambda(x_1) = b$.  In fact this is immediate, since there is some $x \in D(0,r)^-$
for which $\tLambda(x) = b$;  hence by (\ref{FOneComp}), 
\begin{equation*}
\|b,\tLambda(x_1)\|_v \ = \ \|\tLambda(x),\tLambda(x_1)\| \ = \ |b_1-\tlambda_1(x_1)|_v \ = \ 0 \ .
\end{equation*}
\vskip -.25 in 
\end{proof} 


\section{ The Canonical Distance and the $(\fX,\vs)$-Canonical Distance}
        \label{CanonicalDistanceSection}
\index{canonical distance!$[z,w]_{\zeta}$|(}    
\index{canonical distance!$[z,w]_{\fX,\vs}$|(}    

     Let $v$ be a place of $K$. Consider the usual distance function 
$|z-w|_v$  on $\CC_v = \PP^1(\CC_v) \backslash \{\infty\}$, which has the property 
that for each nonzero rational function $f \in \CC_v(\PP^1)$,  
if  $\div(f) = \sum m_i(a_i)$,  
then there is a constant $C(f)$ such that for all $z \in \CC_v$  
\begin{equation} \label{FGAB1}
|f(z)|_v \ = \ C(f) \cdot  \prod_{a_i \ne \infty} |z-a_i|_v^{m_i} \ .
\end{equation} 

Fix $\zeta \in \cC_v(\CC_v)$.
In (\cite{RR1}, \S 2) a `canonical distance' $[z,w]_{\zeta}$ on 
\label{`SymbolIndexCanD'}
\index{canonical distance!$[z,w]_{\zeta}$|ii}
$\cC_v(\CC_v) \backslash \{\zeta\}$ was introduced, generalizing $|z-w|_v$.  
The canonical distance $[z,w]_{\zeta}$
is a symmetric, nonnegative real-valued function of   
$z, w \in \cC_v(\CC_v) \backslash \{\zeta\}$, and is unique
up to scaling by a constant.  Its existence is shown in (\cite{RR1}, Theorem 2.1.1).    
It can be normalized by specifying a 
uniformizing parameter\index{uniformizing parameter!normalizes canonical distance} $g_{\zeta}(z)$, 
in which case it is characterized by the following three properties
(see \cite{RR1}, Theorem 2.1.1, p.57, and Corollary 2.1.2, p.69):
\index{canonical distance!factorization property|ii} 
\index{canonical distance!continuity of|ii}
\index{canonical distance!$[z,w]_{\zeta}$!normalization of|ii}   

\begin{enumerate}
\item {\rm (Continuity):} $[z,w]_{\zeta}$ is jointly continuous in $z$ and $w$.

\item {\rm (Factorization):}  
Let $0 \ne f(z) \in \CC_v(\cC)$
have divisor $\div(f) = \sum m_i(a_i)$.  Then there is a constant $C(f)$
such that for all $z \in \cC_v(\CC_v) \backslash \{\zeta\}$, 
\begin{equation} \label{FCan1}
|f(z)|_v \ = \ C(f) \cdot \prod_{a_i \ne \zeta} [z,a_i]_{\zeta}^{m_i} \ .
\end{equation}

\item {\rm (Normalization):} \ 
    For each $w \in \cC_v(\CC_v) \backslash \{\zeta\}$,
\begin{equation*}
\lim_{z \rightarrow \zeta} \ [z,w]_{\zeta} \cdot |g_{\zeta}(z)|_v \ = \ 1 \ .
\end{equation*}
\end{enumerate}
Two other important properties of the canonical distance are as follows: 
\index{canonical distance!symmetry of|ii}   
\index{canonical distance!galois equivariance|ii}  

\begin{proposition} \label{CanonProp} 
For all $z, w \in \cC_v(\CC_v) \backslash \{\zeta\}$  

\ \, $(4)$  {\rm (Symmetry):} $[z,w]_{\zeta} = [w,z]_{\zeta}$   

\ \, $(5)$ {\rm (Galois equivariance):} 
For each \ $\sigma \in \Aut_c(\CC_v/K_v)$, 
\begin{equation*}  
[\sigma(z),\sigma(w)]_{\sigma(\zeta)} = [w,z]_{\zeta} 
\end{equation*}
if $[x,y]_{\zeta}$ and $[x,y]_{\sigma(\zeta)}$ are normalized compatibly 
$($e.g. if $g_{\sigma(\zeta)}(z) = \sigma(g_{\zeta})(z))$.
\end{proposition}   
  
\begin{proof} 
The canonical distance can be defined directly using 
rational functions.  Fix $\zeta \in \cC_v(\CC_v)$ and fix a uniformizing
parameter $g_{\zeta}(z)$.  For each $w \ne \zeta$,  
there is a sequence of functions $f_n(z) \in \CC_v(\cC)$ 
having poles only at $\zeta$ and whose zeros approach $w$ in the $v$-topology. 
This follows from the `Jacobian Construction
Principle' of (\cite{RR1}, Theorem 1.3.1, p.48),
\index{Jacobian variety}
\index{Jacobian Construction Principle} 
and depends on the fact that the residue field of $\CC_v$ is the algebraic
closure of the prime field $\FF_p$ and its 
valuation group is the same as that of $\tK_v$.  
Fixing a uniformizing parameter\index{uniformizing parameter!normalizes canonical distance} $g_{\zeta}(z)$, 
normalize the $f_n(z)$ so that
\begin{equation*}
\lim_{z \rightarrow \zeta} \ |f_n(z) \cdot g_{\zeta}(z)^{\deg(f_n)}|_v 
           \ = \ 1 \ .
\end{equation*}
Then one can define the canonical distance by 
\index{canonical distance!constructed!directly}
\begin{equation*}
[z,w]_{\zeta} \ = \ \lim_{n \rightarrow \infty} \ |f_n(z)|_v^{1/\deg(f_n)} \ ; 
\end{equation*}
see (\cite{RR1}, Theorem 2.1.1, pp.57-58).  
The limit is independent of the sequence $\{f_n\}$
and the convergence is uniform outside each ball $B(\zeta,r)^-$, 
with $r > 0$. 

The fact that $[z,w]_{\zeta}$ is symmetric in $z$ and $w$ 
is proved in (\cite{RR1}, Theorem 2.1.1, p.57). Its galois equivariance
follows immediately from the galois equivariance of functions.
\end{proof}

\smallskip
The fact that $[z,w]_{\zeta}$ can be approximated by absolute values 
of rational functions is the reason it is the kernel which appears 
in arithmetic potential theory\index{potential theory!arithmetic}, and is the key to the proof of the 
Fekete-Szeg\"o theorem.   
\index{Fekete-Szeg\"o theorem with LRC} 

\vskip .1 in
Several alternate constructions of the canonical distance are given in \cite{RR1}.
To clarify its relation with other objects in arithmetic
geometry, we recall two of them:    

\vskip .1 in
First, the canonical distance is intimately related to 
N\'eron's local height pairing.
\index{canonical distance!constructed!using N\'eron's pairing}
\index{N\'eron's local height pairing|ii}
Recall that N\'eron's pairing $\langle \cdot, \cdot \rangle_v$ is a 
continuous, real-valued, $\Aut_c(\CC_v/K_v)$-equivariant bilinear
function defined for pairs of divisors on $\cC_v(\CC_v)$ of degree $0$
with coprime support, 
having the property that for each $0 \ne f \in \CC_v(\cC)$
and each $a, b \in \cC_v(\CC_v)$ disjoint from the support of $\div(f)$,
\begin{equation*}
\langle\div(f),(a)-(b)\rangle_v \ = \ -\log_v(|f(a)/f(b)|_v) \ .
\end{equation*}
(Here we adopt the normalization of N\'eron's pairing used in \cite{RR1},
which differs from N\'eron's  normalization by a factor $-1/\log(q_v)$.)  
N\'eron originally defined his pairing only for $K_v$-rational divisors: 
the fact that it can be extended
to a galois-equivariant pairing on $\CC_v$-rational divisors
follows from its continuity and invariance under base change;
see (\cite{RR1}, pp.74-76).  

The canonical distance, normalized as in (\ref{FGAM1}),
can be defined using N\'eron's pairing
\index{N\'eron's local height pairing}
by coalescing the poles of $\langle(z)-(t),(w)-(\zeta)\rangle_v$;  
see (\cite{RR1}, \S2.2):  
\begin{equation*}
-\log_v([z,w]_{\zeta}) \ = \ \lim_{t \rightarrow \zeta} \ 
         \langle (z)-(t),(w)-(\zeta)\rangle_v + \log_v(|g_{\zeta}(t)|_v) \ .
\end{equation*}
Conversely, N\'eron's pairing can be recovered from the canonical distance.
\index{N\'eron's local height pairing}
Suppose $D_1 = \sum m_i(a_i)$ and $D_2 = \sum n_j(b_j)$ are
divisors of degree $0$ with coprime support.
Take $\zeta$ distinct from the $a_i$, $b_j$.  Then
\begin{equation*}
\langle D_1,D_2\rangle_v \ = \ - \sum_{i,j} m_i n_j \log_v([a_i,b_j]_{\zeta}) \ .
\end{equation*}

\vskip .1 in
Second, the canonical distance can be expressed in terms of 
Arakelov functions.  If $v$ is archimedean, identify $\CC_v$ with 
\index{canonical distance!constructed!using Arakelov functions}
\index{Arakelov functions|ii}
$\CC$ and let $\(z,w\)_v$ be an Arakelov function on 
$\cC_v(\CC) \times \cC_v(\CC)$.  Then
\begin{equation} \label{JJB1}
[z,w]_{\zeta} \ = \  \frac{\(z,w\)_v}{\(z,\zeta\)_v\(w,\zeta\)_v}
\end{equation}
is a canonical distance function:  see (\cite{RR1}, \S2.3).  
If $v$ is nonarchimedean, functions $\(z,w\)_v$ on
$\cC_v(\CC_v) \times \cC_v(\CC_v)$ for which (\ref{JJB1}) holds
can be constructed using intersection theory and the semistable 
model theorem;  see (\cite{RR1}, \S2.4) and (\cite{CR}, \S2).  
They will also be called Arakelov functions.   

For each $v$, the function $\(z,w\)_v$ is bounded, continuous, symmetric, 
and vanishes only on the diagonal. In (\cite{RR1}, \S2.3, \S2.4) it is shown 
there is a constant $C_v \ge 1$ (depending on the choice of the 
spherical metric $\|z,w\|_v$\index{spherical metric}) 
such that for all $z, w \in \cC_v(\CC_v)$ we have
\begin{equation} \label{FJKLM}
1/C_v \cdot \|z,w\|_v \ \le \ \(z,w\)_v \ \le \ C_v \|z,w\|_v
\end{equation} 

From (\ref{JJB1}), one obtains the following `change of pole' formula
for the canonical distance: 
\index{canonical distance!change of pole formula} 
for any $\xi, \zeta \in \cC_v(\CC_v)$, there is a constant $C_{\xi,\zeta}$ 
such that for all $z, w \ne \xi, \zeta$, 
\begin{equation} \label{FChangeOfSingularity} 
[z,w]_{\xi} = C_{\xi,\zeta} \cdot \frac{[z,w]_{\zeta}}{[z,\xi]_{\zeta} [w,\xi]_{\zeta}} \ .
\end{equation} 
Using (\ref{JJB1}) or (\ref{FChangeOfSingularity}), 
one easily derives the following alternate form of (\ref{FCan1}): 
if $\deg(f) = N$ and we write $\div(f) = \sum_{i=1}^N (\alpha_i) - \sum_{i=1}^N (\xi_i)$,
listing the zeros and poles of $f$ with multiplicities, then there is a constant $\tC(f)$ such that 
for all $z \in \cC_v(\CC_v) \backslash \{\xi_1, \ldots, \xi_N\}$, 
\begin{equation} \label{FCan1a}
|f(z)|_v \ = \ \tC(f) \cdot  \prod_{i=1}^N [z,\alpha_i]_{\xi_i} \ .
\end{equation}

From (\ref{JJB1})  and the fact that 
\begin{equation*}
\[z,w\]_{\zeta} \ = \  \frac{\|z,w\|_v}{\|z,\zeta\|_v \|w,\zeta\|_v}
\end{equation*}
is a metric on $\cC_v(\CC_v)$ (see (\ref{FNCC1})), 
it follows that $[z,w]_{\zeta}$ satisfies a weak triangle inequality:  
there is a constant $B_v$ 
such that for each $\zeta$ and all 
$z, w, p \in \cC_v(\CC_v) \backslash \{\zeta\}$, 
\begin{equation*}
[z,w]_{\zeta} \ \le \ B_v \cdot ([z,p]_{\zeta} + [p,w]_{\zeta}) \ .
\end{equation*}
This property justifies calling $[z,w]_{\zeta}$ a `distance'.    
However, it seems not to be very important in practice, and examples show that
one cannot always take $B_v = 1$ (see \cite{RR1}, p.128). 

If $v$ is nonarchimedean and $\cC$ has good reduction at $v$ for the
\index{good reduction}
projective embedding which induces $\|z,w\|_v$, then $\|z,w\|_v$ is 
an Arakelov function.  Further, if $g_{\zeta}(z) \in K(\cC)$ is a uniformizing
parameter at $\zeta$, then for all but finitely many $v$  
\begin{equation} \label{FGAM2} 
|g_{\zeta}(z)|_v \ = \ \|z,\zeta\|_v 
\end{equation} 
on the ball $B(\zeta,1)^-$.   
Hence, if the canonical distances are normalized as in (\ref{FGAM1}), 
then for all but finitely many $v$,  
\index{canonical distance!constructed!in good reduction case}
\begin{equation*}
[z,w]_{\zeta} \ = \ \frac{\|z,w\|_v}{\|z,\zeta\|_v \|w,\zeta\|_v} \ .
\end{equation*} 
See \cite{RR1}, pp. 90-92. 

\medskip
For most of this work, we will be interested in the case
where $\zeta$ belongs to the $K$-symmetric set  $\fX = \{x_1, \ldots, x_m\}$ 
in the Fekete-Szeg\"o theorem.   Let $g_{x_i}(z) \in K(\cC)$ 
\index{Fekete-Szeg\"o theorem with LRC} 
be the fixed global uniformizing 
parameter chosen in \S\ref{Chap3}.\ref{AssumptionsSection}.  
For each $v$ we will normalize $[z,w]_{x_i}$ so that
\begin{equation} \label{FGAM1}
\lim_{z \rightarrow x_i} \, [z,w]_{x_i} \cdot |g_{x_i}(z)|_v \ = \ 1 \ .
\end{equation}

\medskip
A mild generalization of the canonical
distance, which we call the $(\fX,\vs)$-canonical distance, will play an important role in this work.
Given a probability vector $\vs \in \cP^m$, define 
\index{canonical distance!$[z,w]_{\fX,\vs}$|ii}
\label{`SymbolIndexCanDfXvs'} 
\begin{equation} \label{FXsDef}  
[z,w]_{\fX,\vs} \ = \ \prod_{i=1}^m ([z,w]_{x_i})^{s_i} \ , 
\end{equation} 
where the $[z,w]_{x_i}$ are normalized as in (\ref{FGAM1}).  
The case of interest
is where $\fX \subset \cC(\tK) \subset \cC_v(\CC_v)$ 
is the set of global algebraic points in the Fekete-Szeg\"o theorem,
\index{Fekete-Szeg\"o theorem with LRC} 
and the uniformizing parameters\index{uniformizing parameter!normalizes canonical distance|ii} 
are the ones chosen in \S\ref{Chap3}.\ref{AssumptionsSection}

\vskip .1 in
If $v$ is archimedean, and we identify $\CC_v$ with $\CC$, 
then $\cC_v(\CC)$ is a Riemann surface.  By a coordinate patch 
\index{Riemann surface}
on $\cC_v(\CC)$, we mean a simply connected open\index{simply connected} set $U \subset \cC_v(\CC)$ 
for which there is a chart $\varphi : U \rightarrow \CC$
giving an isomorphism of $U$ 
with an open set $U^{\prime} \subset \CC$.
Given $z, w \in U$, by abuse of notation we will write
$|z-w|$ for $|\varphi(z)-\varphi(w)|$.  

If $v$ is nonarchimedean, 
recall from Theorem \ref{IsoParamThm} that there is an $R_v > 0$
such that each ball $B(a,r)^- \subset \cC_v(\CC_v)$ with $r \le R_v$
is isometrically parametrizable by power series; 
\index{isometrically parametrizable ball} 
if $\varphi : B(a,r)^- \rightarrow D(0,r)^-$
is the inverse map to an isometric parametrization, 
\index{isometric parametrization}
then for all  $z, w \in B(a,r)^-$ 
we have $\|z,w\|_v = |\varphi(z)-\varphi(w)|_v$.

The following result, which is an immediate consequence of 
(\cite{RR1}, Proposition 2.1.3, p.69), asserts that  
$-\log_v([z,w]_{\fX,\vs})$ is `harmonic in $z$ except for 
logarithmic singularities at $w$ and the $x_i \in \fX$', 
and varies continuously with $w$.  It will be used in 
developing potential theory\index{potential theory!$(\fX,\vs)$} for the kernel $-\log_v([z,w]_{\fX,\vs})$.

\begin{proposition} \label{APropA2}
Let $\cC/K$ be a curve, and $v$ a place of $K$. 
Fix $\fX$.  

$(A)$ If $v$ is archimedean, 

\begin{enumerate} 
\item If $U$ and $V$ are disjoint open sets not meeting $\fX$, 
then for each $\vs$, $-\log([z,w]_{\fX,\vs})$ is continuous 
on $U \times V$ and is harmonic in each variable separately.  

\item On any coordinate patch $U \subset \cC_v(\CC)$
not containing not meeting $\fX$, 
there are continuous, real-valued functions
$\eta_{U,x_j}(z,w)$ on $U \times U$, 
harmonic in each variable separately, such that for all $z,w \in U$ and $\vs$, 
\begin{equation*}
-\log([z,w]_{\fX,\vs}) \ = \ -\log(|z-w|) + \sum_{j=1}^m s_j \eta_{U,x_j}(z,w) \ .
\end{equation*}

\item If $U$ is a coordinate patch containing exactly one point $x_i \in \fX$, 
and $V$ is a coordinate patch disjoint from $U$ and $\fX$, 
then there are continuous real-valued functions
$\eta_{U,V,x_j}(z,w)$ on $U \times V$,   
harmonic in each variable separately, such that for all
$z \in U$ and $w \in V$ and all $\vs$,
\begin{equation*}
-\log([z,w]_{\fX,\vs}) \ = \ s_i \log(|z-x_i|) 
        + \sum_{j=1}^m s_j \eta_{U,V,x_j}(z,w) \ .
\end{equation*}

\end{enumerate} 

$(B)$  If  $v$ is nonarchimedean, 

\begin{enumerate}

\item If \, $U = B(a,r)^-$ and $V = B(b,s)^-$ 
are isometrically parametrizable balls disjoint from each other 
\index{isometrically parametrizable ball}
and from $\fX$, then $-\log_v([z,w]_{\fX,\vs})$ is constant on $U \times V$.
More precisely, there are constants 
$\eta_{U,V,x_j} \in \QQ$ such that for all $\vs$ and all $z \in U$, $w \in V$, 
\begin{equation*}
-\log_v([z,w]_{\fX,\vs}) \ =  \ \sum_{j=1}^m s_j \eta_{U,V,x_j} \ .
\end{equation*}
 
\item 
If \, $U = B(a,r)^-$ is an isometrically parametrizable ball not containing 
\index{isometrically parametrizable ball}
any points of $\fX$,
then there are constants $\eta_{U,x_j} \in \QQ$ 
such that for all $z, w \in U$  and all $\vs$
\begin{equation*}
-\log_v([z,w]_{\fX,\vs}) \ = \ -\log_v(\|z,w\|_v) + \sum_{j=1}^m s_j \eta_{U,x_j} \ .
\end{equation*}

\item 
If \, $U = B(a,r)^-$ is an isometrically parametrizable ball 
\index{isometrically parametrizable ball}
containing exactly one point $x_i \in \fX$, 
and $V = B(b,s)^-$ is an isometrically parametrizable ball disjoint 
\index{isometrically parametrizable ball}
from $U$ and $\fX$, then there are constants 
$\eta_{U,V,x_j} \in \QQ$ 
such that for all $z \in U$ and $w \in V$, and all $\vs$, 
\begin{equation*}
-\log_v([z,w]_{\fX,\vs}) \ = \ s_i \log_v(\|z,x_i\|_v) 
                 + \sum_{j=1}^m s_j \eta_{U,V,x_j} \ .
\end{equation*}

\end{enumerate} 
\end{proposition}
\index{canonical distance!$[z,w]_{\fX,\vs}$|)}  
\index{canonical distance!$[z,w]_{\zeta}$|)}  

\section{ $(\fX,\vs)$-Functions and $(\fX,\vs)$-Pseudopolynomials} 
\index{pseudopolynomial!$(\fX,\vs)$}  
    \label{XSfunctionSection}

Fix a place $v$ of $K$.\index{$(\fX,\vs)$-function|(}  
Let $\fX = \{x_1, \ldots, x_m\} \subset \cC(\tK)$ 
be the $K$-symmetric set from \S\ref{Chap3}.\ref{AssumptionsSection}, and let the canonical distances
$[z,w]_{x_i}$ be normalized as in (\ref{FGAM1}), where the uniformizing
parameters $g_{x_i}(z)$ are the ones from \S\ref{Chap3}.\ref{AssumptionsSection}.

\begin{definition} \label{FXvsFunctionDef} Suppose $\vs \in \cP^m \cap \QQ^m$.   
By an {\em $(\fX,\vs)$-function} we mean a rational function $f(z) \in \CC_v(\cC)$
whose poles are supported on $\fX$, such that if $N = \deg(f)$, then $f(z)$ 
has a pole of exact order $N s_i$ at each $x_i \in \fX$.\index{$(\fX,\vs)$-function|ii}  
\end{definition}  

\begin{definition} \label{PseudoPolyDef} Let $\vs \in \cP^m$ be arbitrary.  
By an {\em $(\fX,\vs)$-pseudopolynomial} $($or simply a {\em pseudopolynomial}$)$ 
\index{pseudopolynomial|ii}  
we mean a function $P : \cC_v(\CC_v) \rightarrow [0,\infty]$ of the form 
\begin{equation} \label{FPPD} 
P(z) \ = \ C \cdot \prod_{k=1}^N \, [z,\alpha_k]_{\fX,\vs} \ .
\end{equation}
where $C > 0$ is a constant 
and $\alpha_1, \ldots, \alpha_N \in \cC_v(\CC_v) \backslash \fX$.
We will call $\alpha_1, \ldots, \alpha_N$ the {\em roots} of $P(z)$.
If $C = 1$ and we wish to emphasize that fact, we will say that $P(z)$
is {\em monic}.  We call 
\begin{equation*}
\nu(z) \ = \  \frac{1}{N} \sum_{k=1}^N \delta_{\alpha_k}(z)
\end{equation*}
the  probability measure associated to $P(z)$. 
\end{definition} 

\vskip .1 in 
In the proof of the Fekete-Szeg\"o theorem, $(\fX,\vs)$-functions occur naturally.
\index{Fekete-Szeg\"o theorem with LRC} 
\index{$(\fX,\vs)$-function}
The reason for introducing the $(\fX,\vs)$-canonical distance 
\index{canonical distance!$[z,w]_{\fX,\vs}$}  
is that it allows us to view the absolute value of an $(\fX,\vs)$-function
as an $(\fX,\vs)$-pseudopolynomial, factoring it in the form 
\index{pseudopolynomial!$(\fX,\vs)$|ii}  
\begin{equation} \label{FGHB} 
|f(z)|_v \ = \ C(f) \cdot \prod_{k=1}^N \, [z,\alpha_k]_{\fX,\vs}
\end{equation} 
where $\alpha_1, \ldots, \alpha_N$ are the zeros of $f(z)$, listed with
multiplicities.  This follows by an easy symmetrization argument:  
suppose $f(z)$ is an $(\fX,\vs)$-function,\index{$(\fX,\vs)$-function} and  
let $\xi_{1}, \ldots, \xi_N$ be the points $x_1, \ldots, x_m$ 
listed according to their multiplicities in $\div(f)$.  Thus, each $x_i$ occurs $N s_i$ times.    
For each permutation $\pi$ of $\{1, \ldots, N\}$, by (\ref{FCan1a}) there is 
a constant $C(f,\pi)$ such that 
$|f(z)|_v  =  C(f,\pi) \cdot \prod_{k = 1}^N [z,\alpha_k]_{\xi_{\pi(k)}}$ 
for all $z \in \cC_v(\CC_v) \backslash \fX$. 
Taking the product over all $\pi$, and then extracting the $(N!)^{th}$ root, 
gives (\ref{FGHB}).  

Note that a pseudopolynomial $P(z)$ makes sense even when its roots 
\index{pseudopolynomial}  
are not the zeros of an $(\fX,\vs)$-function\index{$(\fX,\vs)$-function} $f(z)$, 
but it agrees with $|f(z)|_v$ (up to a multiplicative constant) 
when such a function exists.  
Furthermore, $P(z)$ varies continuously with its roots.  
This allows us to investigate absolute values of 
$(\fX,\vs)$-functions\index{$(\fX,\vs)$-function} 
without worrying about questions of principality, which will play
a key role in the construction of the initial local approximating functions
in \S\ref{Chap5} and \S\ref{Chap6} below.\index{$(\fX,\vs)$-function|)}  

\section{ Capacities} \label{CapacitySection}

In this section we define sets of capacity $0$ and sets of positive capacity,
\index{capacity $= 0$}
and we introduce several numerical measures of capacity.
\index{capacity $> 0$}


\medskip
Fix a place $v$ of $K$.

\begin{definition} \label{FCapDef}
If $H$ is a compact subset of $\cC_v(\CC_v)$, 
we will say $H$ has {\em positive capacity} 
\index{capacity $> 0$|ii}
if there is a positive measure $\nu$ supported on $H$ for which
\begin{equation*}
I(\nu) \ := \ \iint_{H \times H} -\log_v(\|z,w\|_v) \, d\nu(z) d\nu(w) 
\ < \ \infty \ .
\end{equation*}
If $I(\nu) = \infty$ for all positive measures $\nu$ on $H$, 
we say that $H$ has {\em capacity $0$}.
\index{capacity $= 0$|ii}  
\end{definition} 

By (\ref{FBB1}) the property of having positive capacity or capacity $0$
\index{capacity $> 0$} 
\index{capacity $= 0$}
is independent of the choice of the spherical metric\index{spherical metric}.
Clearly it suffices to test $I(\nu)$ only for probability measures.  

\vskip .1 in
We next define the capacity of a compact set relative to a point. 
\index{capacity|ii}

Fix  $\zeta \in \cC_v(\CC_v)$, 
and fix a uniformizing parameter\index{uniformizing parameter!normalizes canonical distance} $g_{\zeta}(z)$, 
giving a normalization of the canonical distance $[z,w]_{\zeta}$.  
\index{canonical distance!$[z,w]_{\zeta}$}  
Let $H \subset \cC_v(\CC_v) \backslash \{\zeta\}$ be compact.   
Given a probability measure $\nu$ supported on $H$, 
we define the {\em energy integral of $\nu$ with respect to $\zeta$}\index{energy integral|ii} by  
\begin{equation*} 
I_{\zeta}(\nu) \ = \
      \iint_{H \times H} -\log_v([z,w]_{\zeta}) \, d\nu(z) d\nu(w) \ .
\end{equation*}
If $H$ is nonempty, the {\em Robin constant} of $H$ with respect to $\zeta$ is 
\index{Robin constant!of compact set|ii}
\label{`SymbolIndexVZeta'}
\begin{equation*}
V_{\zeta}(H) 
\ = \ \inf_{\substack{\text{prob meas $\nu$} \\ \text{on $H$}}} I_{\zeta}(\nu) \ ,
\end{equation*}  
where the infimum is taken over all probability measures supported on $H$.  
If $H$ is empty, we put $V_{\zeta}(H) = \infty$.  
The {\em capacity} of $H$ with respect to $\zeta$ to be
\index{capacity} 
\label{`SymbolIndexCapZeta'}
\begin{equation*}
\gamma_{\zeta}(H) \ = \ q_v^{-V_{\zeta}(H)}
\end{equation*} 
Thus, $\gamma_{\zeta}(H) > 0$ 
if and only if $I_{\zeta}(\nu) < \infty$ 
for some probability measure $\nu$ supported on $H$.

\vskip .1 in
Likewise, given a probability vector $\vs \in \cP^m$,
if $H \subset \cC_v(\CC_v) \backslash \fX$ is compact, 
then for any probability measure $\nu$ on $H$, we define the $(\fX,\vs)$-energy 
\begin{equation*} 
I_{\fX,\vs}(\nu) \ = \
      \iint_{H \times H} -\log_v([z,w]_{\fX,\vs}) \, d\nu(z) d\nu(w) \ .
\end{equation*}
We put $V_{\fX,\vs}(H) = \inf_{\nu} I_{\fX,\vs}(\nu)$, 
and define the $(\fX,\vs)$-capacity 
\index{capacity!$(\fX,\vs)$|ii}
\begin{equation*} 
\gamma_{\fX,\vs}(H) \ = \ q_v^{-V_{\fX,\vs}(H)} \ .
\end{equation*} 
  
\vskip .1 in
The following lemma shows that for a given compact set $H$, 
either $\gamma_{\zeta}(H) > 0$ for all $\zeta \notin H$, 
or $\gamma_{\zeta}(H) = 0$ for all $\zeta \notin H$.  

\begin{lemma} \label{ALemA1}
Let $H \subset \cC_v(\CC_v)$ be compact.    
Then the following are equivalent:  

\qquad $(1)$  $H$ has capacity $0$;
\index{capacity $= 0$}

\qquad $(2)$ For some $\zeta \in \cC_v(\CC_v) \backslash H$, $\gamma_{\zeta}(H) = 0$;  

\qquad $(3)$ For each $\zeta \in \cC_v(\CC_v) \backslash H$, $\gamma_{\zeta}(H) = 0$. 

If $H \subset \cC_v(\CC_v) \backslash \fX$, 
these are equivalent to 

\qquad $(4)$  For some $\vs \in \cP^m$, $\gamma_{\fX,\vs}(H) = 0$;  

\qquad $(5)$  For each $\vs \in \cP^m$, $\gamma_{\fX,\vs}(H) = 0$.  
\end{lemma}

\begin{proof}  

If $v$ is archimedean, the set $H = \cC_v(\CC)$ is compact.  
In this case $H$ has positive capacity, and the lemma is vacuously true.
\index{capacity $> 0$}
If $v$ is nonarchimedean, then $\cC_v(\CC)$ is not compact.  
Now suppose $H \ne \cC_v(\CC)$.  If $\zeta \notin H$, 
then $\|x,\zeta\|_v$ is uniformly bounded away from $0$ for $x \in H$, 
because $H$ is compact. Thus the lemma follows from (\ref{JJB1}) and (\ref{FJKLM}).
\end{proof}

We next define the inner capacity and the outer capacity of a set, 
\index{capacity!inner|ii}\index{inner capacity|ii}
\index{capacity!outer|ii}\index{outer capacity|ii}
relative to a point.

\begin{definition}  \label{InnerCapDef}
For an arbitrary set $E_v \subset \cC_v(\CC_v)$,
   we say $E_v$ has {\em positive inner capacity} 
   if there is some compact set $H \subset E_v$ with positive capacity.  
   If every compact set $H \subset E_v$ has capacity $0$, 
   we say that $E_v$ has {\em inner capacity $0$}.  
\end{definition}

For each $\zeta \in \cC_v(\CC_v)$, 
we define the inner capacity $\gammabar_{\zeta}(E_v)$ by
\index{capacity!inner|ii}
\begin{equation} \label{InnerCapacityDef}
\gammabar_{\zeta}(E_v) 
  \ = \ \sup_{\substack{ H \subset E_v \backslash \{\zeta\} \\ \text{$H$ compact}}}
 \gamma_{\zeta}(H) \ . 
\end{equation} 
Thus $0 \le \gammabar_{\zeta}(E_v) \le \infty$.  
When $E_v$ is compact and $\zeta \notin E$, 
clearly $\gammabar_{\zeta}(E_v) = \gamma_{\zeta}(E_v)$.  

Sets of inner capacity $0$ are ``negligible'' for many purposes
\index{capacity!inner}
in potential theory.\index{potential theory}  Each countable set $E_v$ has inner capacity $0$, because
\index{capacity!inner}
a probability measure supported on a compact subset of $E_v$ 
necessarily consists of point masses.   On the other hand, 
any set $E_v$ which contains a nonempty open subset of $\cC_v(\CC_v)$,
or a nonempty open subset of $\cC_v(L_w)$ for some finite extension $L_w/K_v$, 
or a continuum (if $v$ is archimedean), 
has positive inner capacity.  
\index{capacity $> 0$}
This follows from (\cite{RR1}, Proposition 3.1.3, p.137) 
and (\cite{RR1}, Example 4.1.24, p.212).

\smallskip
To define the outer capacity, we will need the notion of a $\PL_{\zeta}$-domain.
\index{capacity!outer|ii}
\index{$\PL_{\zeta}$-domain|ii}  

\begin{definition} If $\zeta \in \cC_v(\CC_v)$, a {\em $\PL_{\zeta}$-domain} is a set
of the form 
\begin{equation*} 
U = \{z \in \cC_v(\CC_v) : |f(z)|_v \le 1 \} \ , 
\end{equation*} 
where $f(z) \in \CC_v(\cC)$ is a nonconstant function 
whose only poles are at $\zeta$.
\end{definition} 

Fix $\zeta \in \cC_v(\CC_v)$, and fix a uniformizer $g_{\zeta}(z)$.  
If $U$ is a $\PL_{\zeta}$-domain, let $f \in \CC_v(\cC)$ 
be a function which defines it.  Write $N = \deg(f)$ and define
\begin{equation*} 
V_{\zeta}(U) \ = \ 
\lim_{z \rightarrow \infty} \frac{1}{N} \log_v(f(z) \cdot g_{\zeta}(z)^N)
\end{equation*} 
We then put 
\begin{equation} \label{OuterCapDef} 
\gamma_{\zeta}(U) \ = \ q_v^{-V_{\zeta}(U)} \ .
\end{equation} 
Using (\cite{RR1}, Theorem 3.2.2 and Proposition 4.3.1) 
one sees that this definition is independent of the choice 
of $f$ defining $U$.  If $v$ is archimdean, 
a $\PL_{\zeta}$-domain is compact, and (\cite{RR1}, Theorem 3.2.2) shows 
that the two definitions we have given for $\gamma_{\zeta}(U)$ coincide.  If
$v$ is nonarchimedean, a $\PL_{\zeta}$-domain $U$ is never compact;  however,
by (\cite{RR1}, Proposition 4.3.1) $\gammabar_{\zeta}(U) = \gamma_{\zeta}(U)$.  

For an arbitrary set $E_v$, if $\zeta \notin E_v$, we define the outer capacity to be
\index{capacity!outer|ii}
\begin{equation*} 
\ulgamma_{\zeta}(E_v) 
  \ = \ \inf_{\substack{ U \supset E_v \\ \text{$U$ a $\PL_{\zeta}$-domain}}}
               \gamma_{\zeta}(U) \ . 
\end{equation*} 
Trivially $\gammabar_{\zeta}(E_v) \le \ulgamma_{\zeta}(E_v)$.  

\begin{definition} \label{AlgCapacitabilityDef}  
Let $E_v \subset \cC_v(\CC_v)$ be arbitrary.  
If $\zeta \notin E_v$, and $\gammabar_{\zeta}(E_v) = \ulgamma_{\zeta}(E_v)$,
we say that $E_v$ is {\em algebraically capacitable with respect to $\zeta$}.
\index{algebraically capacitable!with respect to $\zeta$|ii}  
If $E_v$ is algebraically capacitable with respect to every $\zeta \notin E_v$, 
we simply say say that {\em algebraically capacitable}.  
\index{algebraically capacitable|ii}
\end{definition} 

If $E_v$ is algebraically capacitable with respect to $\zeta$,
\index{algebraically capacitable!with respect to $\zeta$} 
we define its capacity $\gamma_{\zeta}(E_v)$ to be
\index{capacity}
\begin{equation} \label{AlgCapacityDef} 
\gamma_{\zeta}(E_v) \ = \ \gammabar_{\zeta}(E_v) \ = \ \ulgamma_{\zeta}(E_v) \ .
\end{equation} 

In (\cite{RR1}) algebraic capacitability was only defined for 
\index{algebraically capacitable}
sets $E_v$ at nonarchimedean places.
If $v$ is nonarchimedean, it is shown in (\cite{RR1}, Theorem 4.3.13) 
that compact sets, $\RL$-domains, and finite unions of them 
\index{$\RL$-domain} 
are algebraically capacitable. 
\index{algebraically capacitable} 

If $v$ is archimedean, then each $\RL$-domain is compact, 
\index{$\RL$-domain} 
and it follows from (\cite{RR1}, Propositions 3.1.17 and 3.3.3) 
that each compact archimedean set is algebraically capacitable.
\index{algebraically capacitable}

\smallskip
\noindent{\bf Remark.}  
The reason for introducing the notion of algebraic capacitability 
is that the inner capacity turns out to be the `right' notion of capacity
\index{capacity!inner}
for the Fekete-Szeg\"o theorem, whereas the outer capacity is the 
\index{Fekete-Szeg\"o theorem} 
\index{capacity!outer}
right notion for Fekete's theorem (see Theorem \ref{FSZii}). 
\index{Fekete's theorem}  
This is because the initial reductions in the proof of the 
Fekete-Szeg\"o theorem involve replacing an arbitrary set $E_v$ 
with a compact subset whose capacity is arbitrarily near 
\index{capacity}
$\gammabar_{x_i}(E_v)$,  for each $x_i \in \fX$.  Likewise, the initial 
reductions in the proof of Fekete's theorem involve replacing $E_v$
with an algebraically defined neighborhood of itself.  

Thus, algebraic capacitability is the hypothesis which makes a set permissible 
\index{algebraically capacitable}
in both theorems.  In this work, we are primarily interested in the 
Fekete-Szeg\"o theorem, and in stating the most general versions of the theorem
\index{Fekete-Szeg\"o theorem with LRC} 
we use the inner capacity. 
\index{capacity!inner}


\section{ Green's functions of Compact Sets} \label{CompactGreenSection} 

In this section we define and study the Green's functions $G(z,\zeta;H_v)$ 
of compact sets.
\index{Green's function|ii}

\vskip .05 in
Let $H_v \subset \cC_v(\CC_v)$ be compact.  
We first define $G(z,\zeta;H_v)$ when $\zeta \notin H_v$. 
Fix a uniformizing parameter\index{uniformizing parameter!normalizes canonical distance} $g_{\zeta}(z)$, 
which determines the normalization of $[z,w]_{\zeta}$. 
If $H_v$ has positive inner capacity, so $V_{\zeta}(H) < \infty$,
\index{capacity $> 0$} 
there is a unique probability measure
$\mu_{\zeta}$ supported on $H_v$ 
for which $I_{\zeta}(\mu_{\zeta}) = V_{\zeta}(H_v)$.    
It is called the {\it equilibrium distribution} of $H_v$ relative to $\zeta$.

In the archimedean case, the existence of $\mu_{\zeta}$ 
is shown in (\cite{RR1}, p.137), 
and its uniqueness in (\cite{RR1}, Theorem 3.1.12, p.145);  
in the nonarchimedean case, its existence is shown in (\cite{RR1}, p.190), 
and its uniqueness in (\cite{RR1}, Theorem 4.1.22, p.211). 
 
The {\em potential function}\index{potential function|ii} $u_{H_v}(z,\zeta)$ is defined by  
\begin{eqnarray}
u_{H_v}(z,\zeta) & = & \int_{H_v} -\log_v([z,w]_{\zeta}) \, d\mu_{\zeta}(w) \ .
                                        \label{FBUM1} 
\end{eqnarray}
Since $V_{\zeta}(H_v) 
= \int_{H_v \times H_v} -\log_v([z,w]_{\zeta}) \, d\mu_{\zeta}(z) d\mu_{\zeta}(w)$, 
clearly $V_{\zeta}(H_v)-u_{H_v}(z,\zeta)$ is independent of the normalization
of $[z,w]_{\zeta}$.  
\index{potential function|ii} 

\begin{definition} \label{GreenDef1} 
Let $H_v \subset \cC_v(\CC_v)$ be compact, and fix $\zeta \notin H_v$.   
If $H_v$ has positive inner capacity, we define its Green's function
\label{`SymbolIndexGreen'}
\index{capacity $> 0$}
with respect to $\zeta$ to be  
\begin{equation} \label{FMUD1} 
G(z,\zeta;H_v) \ = \ V_{\zeta}(H_v) - u_{H_v}(z,\zeta) 
\end{equation}
for all $z \in \cC_v(\CC_v)$.  
If $H_v$ is compact and has inner capacity $0$,
\index{capacity $= 0$}
we put $G(z,\zeta;H_v) \equiv \infty$.
\end{definition} 
\index{Green's function!of a compact set|ii} 

\vskip .1 in
{\bf Remark.}  
This definition of the Green's function for a compact set
differs from the one in (\cite{RR1}).
In (\cite{RR1}, p.277) both `upper' and `lower' Green's functions 
\index{Green's function!upper}
\index{Green's function!lower} 
$\Gbar(z,\zeta;H_v)$ and $\ulG(z,\zeta;H_v)$ were defined.   
By (\cite{RR1}, Theorems 3.1.9 and 4.1.11), if $H_v$ is algebraically
capacitable (and in particular if $H_v$ is compact) 
then $\Gbar(z,\zeta;H_v)$ and $\ulG(z,\zeta;H_v)$ 
agree for all $z \notin H_v$;  however $\ulG(z,\zeta;H_v) = 0$ 
for all $z \in H_v$ while $\Gbar(z,\zeta;H_v)$ may be positive 
on a subset $e \subset H_v$ of inner capacity $0$.  
\index{capacity $= 0$}

Our $G(z,\zeta;H_v)$ is the same as the upper Green's function 
$\Gbar(z,\zeta;H_v)$, whereas in (\cite{RR1}) $G(z,\zeta;H_v)$ 
was defined to be the lower Green's function $\ulG(z,\zeta;H_v)$.
\index{Green's function|ii}
\index{Green's function!lower|ii}
\index{Green's function!upper|ii}

We have made the change in order to simplify notation, 
and because of the author's conviction
that the choice of $G(z,\zeta;H_v)$ made in (\cite{RR1}) should have been reversed:
$\Gbar(z,\zeta;H)$ carries more information than $\ulG(z,\zeta;H)$,
and is easier to work with.  

\vskip .1 in

The following proposition describes the main properties of the 
Green's function.    
\index{Green's function!properties of|ii}

\begin{proposition} \label{UpGreenProp1} 
Let $H_v \subset \cC_v(\CC_v)$ be a compact set of positive inner capacity, 
\index{capacity $> 0$}
and fix $\zeta \notin H_v$. 
Then $G(z,\zeta;H_v)$ has the following properties:  
for each $\zeta \in \cC_v(\CC_v) \backslash H_v$, 

$(1)$ $G(z,\zeta;H_v) \ge 0$ \ for all $z \in \cC_v(\CC_v)$.

$(2)$  If $v$ is nonarchimedean, then $G(z,\zeta;H_v) > 0$ for all 
         $z \notin H_v$.  
          If $v$ is archimedean, then $G(z,\zeta;H_v) > 0$
         on the connected component of $\cC_v(\CC) \backslash H_v$ 
         containing $\zeta$, and is $0$ on all other components. 
         Furthermore $G(z,\zeta;H_v) = G(z,\zeta;\widehat{H}_v)$, 
         where $\widehat{H}_v = \cC_v(\CC) \backslash D_{\zeta}$ 
         and $D_{\zeta}$ is the connected component of 
         $\cC_v(\CC) \backslash H_v$ containing $\zeta$.  
  
$(3)$ $G(z,\zeta;H_v) = 0$ \ for $z \in H_v$, except 
         on a $($possibly empty$)$ exceptional set $e_v \subset H_v$ 
   \index{exceptional set}
         of inner capacity $0$, which is an $F_{\sigma}$ set  
 \index{capacity!inner}
         and in particular is Borel measurable.  
                 If $v$ is archimedean, 
         $e_v$ is contained in the boundary $\partial H_v$\label{`SymbolIndexBoundary'} 
         and in fact is contained in the `outer boundary' $\partial D_{\zeta}$.
         \index{boundary!exceptional set contained in} 
         Each point of $H_v$ which belongs to a continuum in $H_v$
         is non-exceptional.  In particular, if $H_v$ is a union of continua, 
         the exceptional set is empty.   
\index{exceptional set}      
       
$(4)$ $G(z,\zeta;H_v)$ is continuous for each $z \notin H_v$,
         and at each $z_0 \in H_v$ where $G(z_0,\zeta;H_v) = 0$.
         If $v$ is archimedean, then $G(z,\zeta;H_v)$ is harmonic
         on $\cC_v(\CC) \backslash (H_v \cup \{\zeta\})$ and 
         subharmonic on $\cC_v(\CC) \backslash \{\zeta\}$.    
         
$(5)$    $G(z,\zeta;H_v)$ is upper semi-continuous everywhere:  
\index{semi-continuous!Green's function is upper semi-continuous|ii} 
         in fact, for each $z_0 \in \cC_v(\CC_v)$, 
         \begin{equation*}
         \limsup_{z \rightarrow z_0} \ G(z,\zeta;H_v) 
                     \ = \ G(z_0,\zeta;H_v) \ ,  
         \end{equation*} 
         and if $v$ is nonarchimedean and $z_0 \in H_v$ or if 
         $v$ is archimedean and $z_0 \in \partial D_{\zeta}$, then  
         \begin{equation*}  
         \limsup_{\substack{ z \rightarrow z_0 \\ z \notin H_v }} \ 
                 G(z,\zeta;H_v) \ = \ G(z_0,\zeta;H_v) \  .
         \end{equation*}
      
$(6)$  If $v$ is archimedean, then on any coordinate patch $U$ 
       with $\zeta \in U \subset D_{\zeta}$, 
       there is a harmonic function $\eta_{H_v,\zeta}(z)$ such that 
       $G(z,\zeta;H_v) = -\log_v(|z-\zeta|) + \eta_{H_v,\zeta}(z)$
       on $U$.  
       If $v$ is nonarchimedean, 
       then for any isometrically parametrizable ball 
 \index{isometrically parametrizable ball}
       $B(\zeta,r)^- \subset \cC_v(\CC_v) \backslash H_v$, 
       there is a constant $\eta_{H_v,\zeta}$ such that       
       $G(z,\zeta;H_v) = -\log_v(\|z,\zeta\|_v) + \eta_{H_v,\zeta}$ 
       on $B(\zeta,r)^-$.          
\end{proposition}

\begin{proof} Parts (1)--(4) follow from 
(\cite{RR1}, Lemma 3.1.2, Theorem 3.1.7, Lemma 3.1.8, and Theorem 3.1.9) 
in the archimedean case, and (\cite{RR1}, Lemma 4.1.9, Theorem 4.1.11,
and Corollary 4.1.12) in the nonarchimedean case.  
  Part (5) is contained in (\cite{RR1}, Lemma 3.1.2)
in the archimedean case, and follows from Proposition \ref{APropA2}.B(2)
and the definition of $u_{H_v}(z,\zeta)$ as an integral,
in the nonarchimedean case.   
Part (6) is immediate from the definition.  
\end{proof} 

An important fact is that the Robin constant can be read off from the upper Green's
function:  if $g_{\zeta}(z)$ is the 
uniformizing parameter\index{uniformizing parameter!normalizes canonical distance} 
determining the normalization of $[z,w]_{\zeta}$, then 
\index{Robin constant!upper|ii}
\begin{equation} \label{VReadoff}
    \lim_{z \rightarrow \zeta} 
         \ G(z,\zeta;H_v) + \log_v(|g_{\zeta}(z)|_v) 
             \ = \ V_{\zeta}(H_v) \ .
\end{equation}      
This follows trivially from the definition of $G(z,\zeta;E_v)$ 
in terms of the potential function.\index{potential function}  
However, it emphasizes the fact that the 
Robin constant depends on the choice of the uniformizing parameter,\index{uniformizing parameter!normalizes Robin constant|ii}
while $G(z,\zeta;E_v)$ is absolute.  

\vskip .05 in
The Green's function is decreasing as a function of  $H_v$:
\index{Green's function!monotonic|ii}

\begin{lemma} \label{CompactP1}
Let $H_v \subset H_v^{\prime}$ be compact sets in $\cC_v(\CC_v)$, 
and suppose $\zeta \notin H_v^{\prime}$.   

Then $G(z,\zeta;H_v) \ge G(z,\zeta;H_v^{\prime})$ 
for all $z$.
\end{lemma}

\begin{proof}
In the nonarchimedean case this is (\cite{RR1}, Proposition 4.1.21, p.209).

In the archimedean case (using our notation) 
it is shown in (\cite{RR1}, Lemma 3.2.5, p.157) 
that $\ulG(z,\zeta;H_v) \ge \ulG(z,\zeta;H_v^{\prime})$ for all $z$.  
By the discussion above, it follows that 
$G(z,\zeta;H_v) \ge G(z,\zeta;H_v^{\prime})$
except possibly on a set of capacity $0$ contained in $\partial H_v^{\prime}$. 
\index{capacity $= 0$}  
However, if $z_0 \in \partial H_v^{\prime}$, 
then by Proposition \ref{UpGreenProp1}(4), 
\begin{equation*} 
G(z_0,\zeta;H_v^{\prime}) 
\ = \ \limsup_{\substack{z \rightarrow z_0 \\ z \notin H_v^{\prime} }} 
        G(z,\zeta;H_v^{\prime}) \\
\ \le \ \limsup_{\substack{ z \rightarrow z_0 \\ z \notin H_v^{\prime} }} 
        G(z,\zeta;H_v) 
                \ \le \ G(z_0,\zeta;H_v) \ .
\end{equation*}                 
\end{proof}

The following result seems intrinsically obvious, 
but requires a surprising amount of work to prove.
It will be used in the proof of Theorem \ref{aT1-A1}.   

\begin{proposition} \label{InnerSeqProp} 
Fix $v$, and let $E_v \subset \cC_v(\CC_v)$ be a compact 
set of positive capacity.  
\index{capacity $> 0$}
Let $H_{v,1} \subseteq H_{v,2} \subseteq \cdots \subseteq E_v$
be an exhaustion of $E_v$ by an increasing sequence of compact sets.  
Then for each $\zeta \notin E_v$,
we have  $\lim_{n \rightarrow \infty} V_{\zeta}(H_{v,n}) \ = \ V_{\zeta}(E_v)$,
and for each $z \ne \zeta$
\begin{equation*}
\lim_{n \rightarrow \infty} G(z,\zeta;H_{v,n}) \ = \ G(z,\zeta;E_v) \ .
\end{equation*} 
\end{proposition} 
\index{Green's function|ii} 

\begin{proof} Fix $\zeta \in \cC_v(\CC_v) \backslash E_v$.  

Without loss, we can assume that each $H_{v,n}$ has positive capacity. 
\index{capacity $> 0$}
By the definition of the Robin constant, we have
\index{Robin constant!monotonicity of}
$V_{\zeta}(H_{v,1}) \ge \ V_{\zeta}(H_{v,2}) \ge \cdots \ge V_{\zeta}(E_v)$. 
Put 
\begin{equation*} 
\hV \ = \ \lim_{n \rightarrow \infty} V_{\zeta}(H_{v,n}) \ .
\end{equation*} 
For each $n$, let $\mu_n$ be 
the equilibrium distribution of $H_{v,n}$ with respect to $\zeta$, 
and let $u_n(z,\zeta) = u_{H_{v,n}}(z,\zeta)$ be the potential function.\index{potential function} 
By definition 
\begin{equation*} 
G(z,\zeta;H_{v,n}) \ = \ V_{\zeta}(H_{v,n}) - u_n(z,\zeta) \ .
\end{equation*} 
By Proposition \ref{UpGreenProp1}(3) there is an $F_{\sigma}$-set 
$e_n \subset H_{v,n}$ with inner capacity $0$ 
\index{capacity $= 0$}
such that $G(z,\zeta;H_{v,n}) = 0$, or equivalently 
$u_n(z,\zeta) \ = \ V_{\zeta}(H_{v,n})$, 
for all $z \in H_{v,n} \backslash e_n$.  
By Lemma \ref{CompactP1}, the functions
$G(z,\zeta;H_{v,n})$ are nonnegative and decreasing with $n$. Put 
\begin{equation*} 
G_{\zeta}(z) \ = \ \lim_{n \rightarrow \infty} G(z,\zeta;H_{v,n}) \ . 
\end{equation*} 

Similarly, $G(z,\zeta;E_v) = V_{\zeta}(E_v) - u_{E_v}(z,\zeta)$,
and there is an $F_{\sigma}$-set $e_0 \subset E_v$ of inner capacity $0$
\index{capacity $= 0$}
such that $G(z,\zeta;E_v) = 0$, or equivalently $u_{E_v}(z,\zeta) = V_{\zeta}(E_v)$,
for all $z \in E_v \backslash e_0$.  Let $e = \bigcup_{n=0}^{\infty} e_n$.  
By (\cite{RR1}, Propositions 3.1.15 and 4.1.14) the union of countably many 
Borel sets of inner capacity $0$ itself has inner capacity $0$,
\index{capacity $= 0$} 
so $e$ has inner capacity $0$.  For each 
$z \in E_v \backslash e$, we have 
\begin{equation*}
G_{\zeta}(z) \ = \ G(z,\zeta;E_v) \ = \ 0 \ .
\end{equation*}
We will show, successively, that $\hV = V_{\zeta}(E_v)$, that the $\mu_n$ 
converge weakly to $\mu$, and that $G_{\zeta}(z) = G(z,\zeta;E_v)$ for each 
$z \ne \zeta$.  

\smallskip
By the discussion above, the potential function\index{potential function}  
$u_{E_v}(z,\zeta)$ is identically
equal to $V_{\zeta}(E_v)$ on $E_v \backslash e$.  
Since a Borel set of inner capacity $0$ must have mass $0$ 
\index{capacity $= 0$}
for any positive Borel measure whose potential function\index{potential function}  
is bounded above (\cite{RR1}, Lemmas 3.1.4 and 4.1.7), 
for each $n$ the Fubini-Tonelli theorem gives\index{Fubini-Tonelli theorem}
\begin{eqnarray} 
\int_{E_v} u_n(z,\zeta) \, d\mu(z)  & = & 
\iint_{E_v \times H_{v,n}} -\log_v([z,w]_{\zeta}) \, d\mu_n(w) d\mu(z) \notag \\ 
& = & \int_{H_{v,n}} u_{E_v}(w,\zeta) \, d\mu_n(w) \ = \ V_{\zeta}(E_v) \ .
\label{FFT1A} 
\end{eqnarray} 
On the other hand, pointwise for each $z_0 \in E_v \backslash e$,  we have 
\begin{eqnarray*} 
\lim_{n \rightarrow \infty} u_n(z_0,\zeta) 
& = & \lim_{n \rightarrow \infty}  
             \big(V_{\zeta}(H_{v,n}) - G(z,\zeta;H_{v,n}) \big) 
 \ = \ \hV - G_{\zeta}(z_0) \ = \ \hV \ .
\end{eqnarray*} 

Since $E_v$ is bounded away from $\zeta$, there is a constant $B_1 > -\infty$
such that $u_n(z,\zeta) \ge B_1$  on $E_v$, for all $n$.
On the other hand, since $u_{H_{v,n}}(z,\zeta) \le V_{\zeta}(H_{v,n})$ for all $z$
and the $V_{\zeta}(H_{v,n})$ are finite and decreasing with $n$, 
there is a $B_2 < \infty$ 
such that $u_{H_{v,n}}(z,\zeta) \le B_2$ on $E_v$, for all $n$. 
From (\ref{FFT1A}), the fact that $\mu(e) = 0$, 
and the Dominated Convergence Theorem,\index{Dominated Convergence theorem} it follows that 
\begin{equation} \label{FPAD1A} 
V_{\zeta}(E_v) \ = \ 
\lim_{n \rightarrow \infty} \int_{E_v \backslash e} u_n(z,\zeta) \, d\mu(z) 
\ = \ \int_{E_v \backslash e} \hV \, d\mu(z) \ = \ \hV \ .
\end{equation}  

\smallskip
We next show that the $\mu_n$ converge weakly to $\mu$.  
Let $\widehat{\mu}$ be any weak limit of a subsequence of $\{\mu_n\}_{n \ge 1}$.  
After passing to that subsequence, we can assume that the $\mu_n$ converge
weakly to $\widehat{\mu}$.  We will show that $\widehat{\mu} = \mu$. 
 
For each $M \in \RR$, write
\begin{equation*}
-\log_v^{(M)}([z,w]_{\zeta}) \ = \ \min(M,-\log_v([z,w]_{\zeta}) \ .
\end{equation*}  
Since $-\log_v^{(M)}([z,w]_{\zeta}) \le -\log_v([z,w]_{\zeta})$, 
for each $n$ and each $M$ we have  
\begin{eqnarray} 
\lefteqn{\iint -\log_v^{(M)}([z,w])_{\zeta} \, d\mu_n(z) d\mu_n(w)} & & \notag \\
& & \qquad \le \ \iint -\log_v([z,w])_{\zeta} \, d\mu_n(z) d\mu_n(w) 
\ = \ I_{\zeta}(\mu_n) \ = \ V_{\zeta}(H_{v,n}) \ . \label{FNUG1A} 
\end{eqnarray} 
A standard argument shows that the measures $\mu_n \times \mu_n$ 
converge weakly to $\mu \times \mu$ on $E_v \times E_v$.  
Since the functions $-\log_v^{(M)}([z,w]_{\zeta})$ are continuous on 
$E_v \times E_v$, for each $M$ we have
\begin{equation*} 
\lim_{n \rightarrow \infty} 
   \iint -\log_v^{(M)}([z,w])_{\zeta} \, d\mu_n(z) d\mu_n(w) 
\ = \  \iint -\log_v^{(M)}([z,w])_{\zeta} \, d\widehat{\mu}(z) d\widehat{\mu}(w) \ .
\end{equation*} 
Combining this with (\ref{FPAD1A}) and (\ref{FNUG1A}) shows that for each $M$ 
\begin{equation*} 
\iint -\log_v^{(M)}([z,w])_{\zeta} \, d\widehat{\mu}(z) d\widehat{\mu}(w) 
\ \le \ \lim_{n \rightarrow \infty} V_{\zeta}(H_{v,n})
\ = \ \hV \ = \ \ V_{\zeta}(E_v) \ . 
\end{equation*}
On the other hand, by the Monotone Convergence Theorem\index{Monotone Convergence theorem} 
\begin{equation*} 
I_{\zeta}(\widehat{\mu}) \ = \ \lim_{M \rightarrow \infty} 
\iint -\log_v^{(M)}([z,w])_{\zeta} \, d\widehat{\mu}(z) d\widehat{\mu}(w) \ .
\end{equation*} 
Thus $I_{\zeta}(\widehat{\mu}) \le V_{\zeta}(E_v)$.
However, by the definition of $V_{\zeta}(E_v)$ we must have 
$I_{\zeta}(\widehat{\mu}) \ge V_{\zeta}(E_v)$. 
Hence $I_{\zeta}(\widehat{\mu}) = V_{\zeta}(E_v)$, and the uniqueness of the
equilibrium measure gives $\widehat{\mu} = \mu$.  

\smallskip
Lastly, we show that $G_{\zeta}(z) = G(z,\zeta;E_v)$ for each $z \ne \zeta$.
For each $M$, put
\begin{equation*}
u_n^{(M)}(z,\zeta)  =  \int -\log_v^{(M)}([z,w]_{\zeta}) \, d\mu_n(w) \ ,
\quad 
u_{E_v}^{(M)}(z,\zeta)  =  \int -\log_v^{(M)}([z,w]_{\zeta}) \, d\mu_n(w) \ .
\end{equation*} 
Since the kernels $-\log_v^{(M)}([z,w]_\zeta)$ are increasing with $M$,
the Monotone Convergence Theorem\index{Monotone Convergence theorem} shows that for each $z \ne \zeta$ 
\begin{equation*}
\lim_{M \rightarrow \infty} u_n^{(M)}(z,\zeta) \ = \ u_n(z,\zeta) \ , 
\qquad 
\lim_{M \rightarrow \infty} u_{E_v}^{(M)}(z,\zeta) \ = \ u_{E_v}(z,\zeta) \ ,
\end{equation*}
where the limits are increasing.

Now fix $z_0 \ne \zeta$.  
For each $M$, since $-\log_v^{(M)}([z_0,w]_\zeta)$ is continuous
on $E_v$ as a function of $w$ and the $\mu_n$ converge weakly to $\mu$, we have 
\begin{equation*}
\lim_{n \rightarrow \infty} u_n^{(M)}(z_0,\zeta)  \ = \ u_{E_v}^{(M)}(z_0,\zeta) \ .
\end{equation*} 
Hence, for each $\varepsilon > 0$, there is an $N = N(M,\varepsilon)$ 
such that for all $n \ge N(M,\varepsilon)$, we have 
$u_n^{(M)}(z_0,\zeta) > u_{E_v}^{(M)}(z_0,\zeta) - \varepsilon$.
By the monotonicity of $u_n^{(M)}(z_0,\zeta)$ in $M$, 
for each $n \ge N(M,\varepsilon)$ and each $M_1 > M$ we have 
\begin{equation*} 
u_n^{(M_1)}(z_0,\zeta) \ \ge \ u_n^{(M)}(z_0,\zeta) 
\ > \ u_{E_v}^{(M)}(z_0,\zeta) - \varepsilon \ .
\end{equation*}
Letting $M_1 \rightarrow \infty$ and then letting $M \rightarrow \infty$,
we see that for all sufficiently large $n$ 
\begin{equation*} 
u_n(z_0,\zeta) \ \ge \ u_{E_v}(z_0,\zeta) - \varepsilon \ .
\end{equation*}
Consequently 
\begin{eqnarray*}
G_{\zeta}(z_0) & = & \lim_{n \rightarrow \infty} G(z_0,\zeta;H_{v,n}) 
     \ \,  = \ \, \lim_{n \rightarrow \infty} V_{\zeta}(H_{v,n}) - u_n(z_0,\zeta) \\
     & \le & V_{\zeta}(E_v) - u_{E_v}(z_0,\zeta) + \varepsilon 
     \ \, = \ \, G(z_0,\zeta;E_v) + \varepsilon \ .
\end{eqnarray*} 
Letting $\varepsilon \rightarrow 0$, 
we see that $G_{\zeta}(z_0) \le G(z_0,\zeta;E_v)$.  
On the other hand, for each $n$ the monotonicity of Green's functions
\index{Green's function!monotonic}
shows that $G(z_0,\zeta;H_{v,n}) \ge G(z_0,\zeta;E_v)$, and hence trivially
$G_{\zeta}(z_0) \ge G(z_0,\zeta;E_v)$.  
It follows that $G_{\zeta}(z_0) = G(z_0,\zeta;E_v)$.  

This complete the proof. 
\end{proof} 

\section{ Upper Green's functions} \label{UpperGreenSection} 

In this section we introduce the upper 
Green's function's $\Gbar(z,\zeta;E_v)$ of an arbitrary set. 

\begin{definition} \label{GreenDef2} 
Given an arbitrary set $E_v \subset \cC_v(\CC_v)$,
for each  $z, \zeta \in \cC_v(\CC_v)$ we define the upper Green's function by
\index{Green's function!upper|ii}
\begin{equation} \label{FMUD2}
\Gbar(z,\zeta;E_v) \ = \
      \inf_{\substack{ H_v \subset E_v \backslash \{\zeta\} 
                \\ \text{$H_v$ compact} }}
             G(z,\zeta;H_v) \ ,
\end{equation}
and $($fixing a uniformizing parameter\index{uniformizing parameter!normalizes Robin constant} $g_{\zeta}(z)$)  
we define the upper Robin constant by  
\index{Robin constant!upper|ii}            
\begin{equation} \label{FMUD3} 
\Vbar_{\zeta}(E_v) \ = \
   \inf_{\substack{ H_v \subset E_v \backslash \{\zeta\} \\ \text{$H_v$ compact} }}
                V_{\zeta}(H_v) \ . 
\end{equation}            
\end{definition} 

If $E_v$ has inner capacity $0$ then  
\index{capacity $= 0$}
$\Gbar(z,\zeta;E_v) \equiv \infty$.
If $E_v$ has positive inner capacity, then for each $z \ne \zeta$,
\index{capacity $> 0$} 
$\Gbar(z,\zeta;E_v)$ is finite and non-negative, 
while $G(\zeta,\zeta;E_v) = \infty$.
By Lemma \ref{CompactP2}, for a compact set $H_v$ 
we have $\Gbar(z,\zeta;H_v) = G(z,\zeta;H_v)$.

\smallskip
The following results will be needed later.

\begin{lemma} \label{DecreasingSeqLemma} 
Let $E_v \subset \cC_v(\CC_v)$ have positive inner capacity.  
\index{capacity $> 0$}
Then for each $\zeta \in \cC_v(\CC_v)$, there is an increasing 
sequence of compact sets 
$H_{v,1} \subseteq H_{v,2} \subseteq \cdots \subset E_v \backslash \{\zeta\}$ 
such that $\Vbar_{\zeta}(E_v) = \lim_{n \rightarrow \infty} V_{\zeta}(H_{v,n})$ 
and for each $z \in \cC_v(\CC_v) \backslash \{\zeta\}$,
\begin{equation*}
\Gbar(z,\zeta;E_v) \ = \ \lim_{n \rightarrow \infty} G(z,\zeta;H_{v,n}) 
\end{equation*} 
\end{lemma}

\begin{proof} Fix $\zeta$, and let $A$ be the collection of all compact sets 
$K \subset E_v \backslash \{\zeta\}$ with positive inner capacity.  Consider the 
\index{capacity $> 0$}
family $\{G(z,\zeta;K)\}_{K \in A}$.  The set $\cC_v(\CC) \backslash \{\zeta\}$
is a separable metric space, and by Proposition \ref{UpGreenProp1}.5 
each $G(z,\zeta;K)$ for $K \in A$ is upper semi-continuous 
\index{semi-continuous!Green's function is upper semi-continuous} 
in $\cC_v(\CC) \backslash \{\zeta\}$.  By 
a well-known property of upper semi-continuous functions (see \cite{Kl}, Lemma 2.3.2)
there is a countable sequence $\{K_n\}_{n \ge 1}$ of sets in $A$ such that for each 
$z \in \cC_v(\CC) \backslash \{\zeta\}$.  
\begin{equation*}
\Gbar(z,\zeta;E_v) \ := \ \inf_{K \in A} G(z,\zeta;K) 
\ = \ \inf_{n \ge 1} G(z,\zeta;K_n) \ .
\end{equation*}
Likewise, since $\Vbar_{\zeta}(E_v) = \inf_{K \in A} V_{\zeta}(K)$, there is a 
sequence  $\{K_n^{\prime}\}_{n \ge 1}$ of sets in $A$ such that 
$\Vbar_{\zeta}(E_v) = \lim_{n \rightarrow \infty} V_{\zeta}(K_n^{\prime})$.

For each $n \ge 1$ put 
$H_{v,n} = (\bigcup_{i=1}^n K_i) \cup (\bigcup_{i=1}^n K_i^{\prime})$.  
By the monotonicity of Green's functions of compact sets, for each $n$ we have 
\index{Green's function!monotonic}
$G(z,\zeta;K_n) \ge G(z,\zeta;H_{v,n}) \ge G(z,\zeta;H_{v,n+1})$;
similarly $V_{\zeta}(K_n^{\prime}) \ge V_{\zeta}(H_{v,n}) \ge V_{\zeta}(H_{v,n+1})$.
It follows that 
\begin{equation*} 
\Gbar(z,\zeta;E_v) \ \le \ \lim_{n \rightarrow \infty} G(z,\zeta;H_{v,n}) 
               \ \le \ \inf_{n \ge 1} G(z,\zeta;K_n) \ = \ \Gbar(z,\zeta;E_v) 
\end{equation*} 
and that $\Vbar_{\zeta}(E_v) = \lim_{n \rightarrow \infty} V_{\zeta}(H_{v,n})$.   
\end{proof} 

\begin{lemma}  \label{CompactP2}

If $E_v$ is compact, $e_v$ has inner capacity $0$, and $\zeta \notin E_v$, then 
\index{capacity $= 0$}
\begin{equation} \label{FCP2B} 
\Gbar(z,\zeta;E_v \backslash e_v) \ = \ G(z,\zeta;E_v) \ .
\end{equation} 
If $E_v \subset \cC_v(\CC_v)$ is arbitrary  
and $e_v$ has inner capacity $0$, then for each $\zeta \in \cC_v(\CC_v)$,
\index{capacity $= 0$}  
\begin{equation} \label{FCP2A}
\Gbar(z,\zeta;E_v \backslash e_v) \ = \ \Gbar(z,\zeta;E_v) \ .                
\end{equation}
\end{lemma}  

\begin{proof}  We first prove (\ref{FCP2B}).  
Suppose $E_v$ is compact and $\zeta \notin E_v$.
If $E_v$ has inner capacity $0$ the result holds trivially, 
\index{capacity $= 0$}
so we can assume that $E_v$ has positive inner capacity.  
\index{capacity $> 0$}

By (\cite{RR1}, Corollaries 3.1.16 and 4.1.15), we have 
$\Vbar_{\zeta}(E_v \backslash e_v) = V_{\zeta}(E_v)$. 
(Note that in (\cite{RR1}), our inner capacity $\gammabar(E_v)$ is denoted $\ulgamma(E_v)$.)  
By Lemma \ref{DecreasingSeqLemma} there is an increasing sequence of compact sets
$H_{v,1} \subseteq H_{v,2} \subseteq \cdots \subseteq E_v \backslash e_v$
such that $\Vbar_{\zeta}(E_v \backslash e_v) 
= \lim_{n \rightarrow \infty} V_{\zeta}(H_{v,n})$ 
and for each $z \in \cC_v(\CC_v) \backslash \{\zeta\}$,
\begin{equation*}
\Gbar(z,\zeta;E_v \backslash e_v) \ = \ 
\lim_{n \rightarrow \infty} G(z,\zeta;H_{v,n}) \ .
\end{equation*} 
Let $\mu_n$ be the equilibrium distribution of $H_{v,n}$ with respect to $\zeta$,
and let $\mu$ be the equilibrium distribution of $E_v$ with respect to $\zeta$.  
By the discussion above we have 
\begin{equation*} 
\lim_{n \rightarrow \infty} V_{\zeta}(H_{v,n}) \ = \ V_{\zeta}(E_v) \ .
\end{equation*} 
Now the same argument as in the part of the proof of 
Proposition \ref{InnerSeqProp} after formula (\ref{FPAD1A}) shows that the
$\mu_n$ converge weakly to $\mu$, 
and that $\Gbar(z,\zeta;E_v \backslash e_v) = G(z,\zeta;E_v)$. 

\smallskip
We can now deduce (\ref{FCP2A}) formally.  Let $E_v \subset \cC_v(\CC_v)$
and $\zeta \in \cC_v(\CC_v)$ be arbitrary.  By definition, 
\begin{equation*} 
\Gbar(z,\zeta;E_v)  = 
      \inf_{\substack{ H_v \subset E_v \backslash \{\zeta\} 
                \\ \text{$H_v$ compact} }}
             G(z,\zeta;H_v) \ , \quad \ 
\Gbar(z,\zeta;E_v \backslash e_v)  = 
  \inf_{\substack{ H_v^{\prime} \subset (E_v \backslash e_v) \backslash \{\zeta\} 
                \\ \text{$H_v^{\prime} $ compact} }}
             G(z,\zeta;H_v^{\prime}) \ .
\end{equation*}         
On the other hand, by what has been shown above,
for each compact subset $H_v \subset E_v$ we have 
\begin{equation*}
G(z,\zeta;H_v) \ = \ \Gbar(z,\zeta;H_v \backslash e_v) \ = \ 
 \inf_{\substack{ H_v^{\prime} \subset (H_v \backslash e_v) \backslash \{\zeta\} 
                \\ \text{$H_v^{\prime} $ compact} }}
             G(z,\zeta;H_v^{\prime}) \ .
\end{equation*} 
It follows that $\Gbar(z,\zeta;E_v) = \Gbar(z,\zeta;E_v \backslash e_v)$.  
\end{proof}

The upper Green's function has the following properties:    

\begin{proposition} \label{GreenProp} 
\index{Green's function!upper|ii} 

Let $E_v \subset \cC_v(\CC_v)$.  Then 

$(1)$ {\rm (Finiteness):}  
$\Gbar(z,\zeta;E_v)$ is valued in $[0,\infty]$.  If $E_v$ has inner capacity $0$,
\index{capacity!inner} 
then $\Gbar(z,\zeta;E_v) \equiv \infty$.  If $E_v$ has positive inner capacity,
then $\Gbar(z,\zeta;E_v)$ is finite and upper semi-continuous 
\index{semi-continuous!Green's function is upper semi-continuous|ii} 
on $\cC_v(\CC_v) \backslash \{\zeta\}$.  If $v$ is archimedean
and $E_v$ has positive inner capacity, then $\Gbar(z,\zeta;E_v)$ is subharmonic on 
$\cC_v(\CC) \backslash \{\zeta\}$.    

$(2)$ {\rm (Symmetry):} \ $\Gbar(z,\zeta;E_v)  =  \Gbar(\zeta,z;E_v)$ 
for all $z, \zeta \in \cC_v(\CC_v)$ . 

$(3)$ {\rm (Galois equivariance):}  
$\Gbar(\sigma(z),\sigma(\zeta);\sigma(E_v))  =  \Gbar(z,\zeta;E_v)$
 for all $z, \zeta \in \cC_v(\CC_v)$ and all $\sigma \in \Aut_c(\CC_v/K_v)$.  
In particular, if $E_v$ is stable under $\Aut_c(\CC_v/K_v)$, then 
\begin{equation*}
\Gbar(\sigma(z),\sigma(\zeta);E_v) \ = \ \Gbar(z,\zeta;E_v) \ .
\end{equation*}

$(4)$ {\rm (Approximation):}  Let $\fX = \{x_1, \ldots, x_m\} \subset \cC_v(\CC_v)$ 
be a finite set of points, none of which belongs to the closure of $E_v$, 
and fix $\varepsilon > 0$.  Then there is a compact set $H_v \subset E_v$
such that for all $x_i, x_j \in \fX$,
\begin{equation*}
\left\{ \begin{array}{cl}
    \Gbar(x_i,x_j;E_v) \ \le \ \Gbar(x_i,x_j;H_v) 
                \ \le \ \Gbar(x_i,x_j;E_v) + \varepsilon
            & \text{if $i \ne j$\ ;} \\ 
    \Vbar_{x_i}(E_v) \ \le \ \Vbar_{x_i}(H_v) \ \le \ \Vbar_{x_i}(E_v) + \varepsilon
            & \text{if $i = j$\ .} 
         \end{array}  \right.
\end{equation*} 

$(5)$ {\rm (Base Change):}  Let $L_w/K_v$ be a finite extension with 
ramification index $e_{w/v}$ $($take $e_{w/v} = 1$ if  $v$ is archimedean$)$.
Fix an isomorphism $\iota_{w/v} : \CC_w \rightarrow \CC_v$ and put 
$E_w = \iota_{w/v}^{-1}(E_v) \subset \cC_w(\CC_w)$.  
Then 
\begin{equation*}
\Gbar(z,\zeta;E_w) \ = \ e_{w/v} \cdot \Gbar(z,\zeta;E_v) \ .
\end{equation*}

$(6)$ {\rm (Pullback):}\index{Green's function!pullback formula for|ii}  Let $\cC_1, \cC_2/K_v$ be curves and let 
$f : \cC_1 \rightarrow \cC_2$ be a nonconstant rational map.  
Given $\zeta \in \cC_2(\CC_v)$,  let its pullback divisor be 
$f^*((\zeta)) = \sum_i m_i (\xi_i)$.  
Then for each $E_v \subset \cC_2(\CC_v)$ and each $z \in \cC_1(\CC_v)$
\begin{equation*}
\Gbar(f(z),\zeta;E_v) \ = \ \sum_i m_i \Gbar(z,\xi_i;f^{-1}(E_v)) \ .
\end{equation*}
\end{proposition} 

\vskip .05 in
\begin{proof} 
In \cite{RR1} these properties are established for Green's functions of
\index{Green's function!of a compact set}
compact sets $H_v$, with $z, \zeta \notin H_v$.  Using Definition \ref{GreenDef2},
they carry over to arbitrary sets by taking limits:  

(1) The finiteness properties in assertion (1) are immediate from the definition
and Proposition \ref{UpGreenProp1}.  If $E_v$ has positive inner capacity, then
\index{capacity!inner} 
Lemma \ref{DecreasingSeqLemma} and Proposition \ref{UpGreenProp1}.
$G(z,\zeta;H_v)$ is a decreasing limit of upper semi-continuous functions
\index{semi-continuous!Green's function is upper semi-continuous} 
on $\cC_v(\CC_v) \backslash \{\zeta\}$, hence is itself upper semi-continuous.  
If in addition $v$ is archimedean, then  $\Gbar(z,\zeta;E_v)$ is the limit of a
decreasing sequence of subharmonic functions which is bounded below, 
so it is subharmonic (see \cite{Kl}, Theorem 2.6.ii).   

(2) By Lemma \ref{CompactP2}, since any finite set has inner capacity $0$,
\index{capacity!inner}
for each fixed $z_0, \zeta_0$ we have 
\begin{equation*} 
\Gbar(z_0,\zeta_0;E_v) \ = \ \Gbar(z_0,\zeta;E_v \backslash \{z_0,\zeta\}) \ = \ 
      \inf_{\substack{ H_v \subset E_v \backslash \{z_0,\zeta_0\} \ ,
                \\ \text{$H_v$ compact} }}
             G(z_0,\zeta_0;H_v) \ ,
\end{equation*} 
and a similar formula holds for $\Gbar(\zeta_0,z_0;E_v)$.  
By (\cite{RR1}, Theorem 3.2.7 and Theorem 4.4.14), 
for each compact $H_v \subset E_v \backslash \{z_0,\zeta_0\}$ we have  
$G(z_0,\zeta_0;H_v) = G(\zeta_0,z_0;H_v)$.

(3) By the Galois-equivariance of $[z,w]_{\zeta}$ (Proposition \ref{CanonProp} (B))
and formula (\ref{FMUD1}), it follows that for compact sets $H_v$ 
and $z, \zeta \notin H_v$, 
we have $G(\sigma(z),\sigma(\zeta),\sigma(H_v)) = G(z,w;H_v)$.   
By taking limits, it follows that in general 
\begin{equation} \label{FGAL1}  
\Gbar(\sigma(z),\sigma(\zeta);\sigma(E_v)) \ = \ \Gbar(z,\zeta;E_v) \ .
\end{equation}

(4) The approximation property is immediate from Definition \ref{GreenDef2} 
and Proposition \ref{APropM2}. 

(5) The base change property follows immediately from our normalizations of  
$|x|_v$,$|x|_w$, $\log_v(z)$ and $\log_w(z)$ (see \S\ref{Chap3}.\ref{NotationSection})) 

(6) The pullback formula for is established for compact sets in 
(\cite{RR1}, Theorem 3.2.9 and Theorem 4.4.19).  Fix $z_0$ and $\zeta$, 
and put $S = f^{-1}(\{z_0,\zeta\})$.  As $H_v$ runs over compact subsets of 
$E_v \backslash \{z_0,\zeta\}$, then $f^{-1}(H_v)$ runs over compact subsets
of $f^{-1}(E_v) \backslash S$.  These sets are cofinal in the compact subsets 
of $f^{-1}(E_v) \backslash S$, 
since for any compact $H_v^{\prime} \subset f^{-1}(E_v) \backslash S$, its image 
$f(H_v^{\prime}) \subset E_v \backslash \{z_0,\zeta\}$ is compact,
and $H_v^{\prime} \subset f^{-1}(f(H_v^{\prime})$.  Hence, the pullback
formula follows in general by applying Definition \ref{GreenDef2} to $E_v$
and using Lemma \ref{CompactP2} to $f^{-1}(E)$ and $S$.    
\end{proof}  

We will now show that just as for compact sets, 
the Robin constant can be read off from the upper Green's function, 
\index{Green's function!upper}
\index{Robin constant!upper|ii}
using the chosen uniformizing parameter:\index{uniformizing parameter!normalizes Robin constant}  

\begin{proposition} \label{APropM2}
If $E_v$ has positive inner capacity and is bounded away from $\zeta$,
\index{capacity!inner}
then the upper Robin constant is finite, and 
\begin{equation} \label{FMUD4}
\Vbar_{\zeta}(E_v) \ = \ \lim_{z \rightarrow \zeta} \Gbar(z,\zeta;E_v) 
                              + \log_v(|g_{\zeta}(z)|_v) \ .
\end{equation}    
\end{proposition}

\begin{proof}
In the archimedean case, fix a neighborhood $U$ of $\zeta$
in $\cC_v(\CC_v) \backslash E_v$
such that the uniformizer $g_{\zeta}(z)$ has no zeros
or poles in $U$ except at $\zeta$, and put $F = \cC_v(\CC_v) \backslash U$.
Then $F$ is a compact set of positive capacity and $E_v \subset F$,
\index{capacity $> 0$}
so $\Gbar(z,\zeta,H_v) \ge \Gbar(z,\zeta;F)$ for each $H_v \subset E_v$.
For each compact $H_v \subset E_v$ of positive capacity,
\index{capacity $> 0$}
the function $\eta_{H_v}(z) = \Gbar(z,\zeta;H_v)+\log_v(|g_{\zeta}(z)|_v)$
is harmonic on $U \backslash \{\zeta\}$ and is bounded in a neighborhood
of $\zeta$, so it extends to a harmonic function on $U$.
Similarly $\eta_F(z) = \Gbar(z,\zeta;F)+\log_v(|g_{\zeta}(z)|_v)$ extends
to a harmonic function on $U$.  By Lemma \ref{DecreasingSeqLemma},
there is an increasing sequence of compact sets $\{H_{v,n}\}$ contained in $E_v$
such that $\Gbar(z,\zeta;E_v) = \lim_{n \rightarrow \infty} G(z,\zeta;H_{v,n})$
and $\Vbar_{\zeta}(E_v) = \lim_{n \rightarrow \infty} V_{\zeta}(H_{v,n})$.
It follows that $\{\eta_{H_{v,n}}(z) \}$ is a decreasing sequence of 
functions harmonic in $U$, which is bounded below by $\eta_F(z)$.  
By Harnack's Principle, 
\index{Harnack's Principle}
$\eta_{E_v}(z) := \lim_{n \rightarrow \infty} \eta_{H_{v,n}}(z)$
is harmonic in $U$ and 
$\Gbar(z,\zeta;E_v) = \eta_{E_v}(z) - \log_v(|g_{\zeta}(z)|_v)$  in
$U \backslash \{\zeta\}$.  
Since $\eta_{H_{v,n}}(\zeta) = V_{\zeta}(H_{v,n})$ for each $v$,  
our assertion follows.

In the nonarchimedean case, let $r_0$ be as in Proposition \ref{APropA2}(B),
and take $0 < r < r_0$ small enough
that $r < \|\zeta,E_v\|_v := \inf_{z \in E_v} (\|z,\zeta\|_v)$.
By Proposition \ref{APropA2}(B) and the definition (\ref{FMUD1}),
there is a constant $C_{\zeta}$ such that for each
compact $H_v \subset E_v$ we have
\begin{equation*}
\Gbar(z,\zeta;H_v) \ = \  C_{\zeta} + V_{\zeta}(H_v) - \log_v(\|z,\zeta\|_v)
\end{equation*}
for all $z \in (B(\zeta,r) \backslash \{\zeta\})$.  For each such $z$ the
monotonicity of the Green's functions of compact sets shows that the
\index{Green's function!of a compact set}
values $\Gbar(z,\zeta;H_v)$ form a directed set, bounded below by $0$,
hence convergent.  Thus for each $z \in (B(\zeta,r) \backslash \{\zeta\})$
\begin{eqnarray*}
\Gbar(z,\zeta;E_v) & := & \inf_{H_v \subset E_v} \Gbar(z,\zeta;H_v) \\
               & = & C_{\zeta} + \big(\inf_{H_v \subset E_v} V_{\zeta}(H_v))
                      - \log_v(\|z,\zeta\|_v) \big)
\end{eqnarray*}
and it follows that
\begin{equation*}
\Vbar_{\zeta}(E_v)
  \ = \ \lim_{z \rightarrow_{\zeta}} \Gbar(z,\zeta;E_v) + \log_v(|g_{\zeta}(z)|_v)
  \ = \ \inf_{H_v \subset E_v} V_{\zeta}(H_v) \ .
\end{equation*}
\vskip -.25 in
\end{proof}

For nonarchimedean $v$, the Green's functions and Robin constants
\index{Green's function!properties of|ii} 
\index{Robin constant!properties of|ii}
 of `nice' sets take on values in $\QQ$.

\begin{proposition} \label{GreenRationalProp} 
Let $v$ be nonarchimedean.  Let $E_v \subset \cC_v(\CC_v)$.  
Suppose that either

$(A)$  $E_v = \bigcup_{\ell=1}^D \big(B(a_\ell,r_\ell) \cap \cC_v(F_{w_\ell})\big)$ where 
$B(a_1,r_1), \ldots, B(a_D,r_D) \subset \cC_v(\CC_v)$  
are pairwise disjoint isometrically parametrizable balls and 
\index{isometrically parametrizable ball}
$F_{w_1}, \ldots, F_{w_D}$ are finite extensions of $K_v$ in $\CC_v$,  
with $a_{\ell} \in \cC_v(F_{w_{\ell}})$ and $r_{\ell} \in |F_{w_{\ell}}^{\times}|_v$  
for each $\ell;$ or 

$(B)$  $E_v$ is an $\RL$-domain, that is $E_v = \{z \in \cC_v(\CC_v) : |f(z)|_v \le 1\}$\index{$\RL$-domain} 
for some nonconstant $f(z) \in \CC_v(\cC)$. 
\index{Green's function!nonarchimedean!takes on rational values|ii} 
\index{Robin constant!nonarchimedean!takes on rational values|ii}   

Then for each $\zeta \in \cC_v(\CC_v) \backslash E_v$, 
we have $\Vbar_{\zeta}(E_v) \in \QQ$, and $\Gbar(z,\zeta;E_v) \in \QQ$ for all $z \ne \zeta$.
\end{proposition}

\begin{proof}  
In case (A), the assertion is proved in Corollary \ref{BFCor3} of Appendix \ref{AppA}.

In case (B), the assertion follows from results proved in (\cite{RR1}).  
Suppose $E_v$ is an $\RL$-domain.
\index{$\RL$-domain} 
First note that by (\cite{RR1}, Theorem 4.3.3), each $\RL$-domain is algebraically
\index{$\RL$-domain} 
capacitable in the sense of (\cite{RR1}, Definition 4.3.2). 
By (\cite{RR1}, Theorem 4.4.4), this means that for $z \notin E_v$, 
the upper Green's function $\Gbar(z,\zeta;E_v)$ coincides with the lower
\index{Green's function!upper}
Green's function $\underline{G}(z,\zeta;E_v)$ defined in (\cite{RR1}, p.282).  
\index{Green's function!lower} 
On the other hand, since each point of $E_v$ has a neighborhood contained in $E_v$, 
trivially $\Gbar(z,\zeta;E_v) = 0$ for $z \in E_v$. 
Since $\Gbar(z,\zeta;E_v) \ge \underline{G}(z,\zeta;E_v) \ge 0$ for all $z$,
we conclude that $\Gbar(z,\zeta;E_v) = \underline{G}(z,\zeta;E_v)$ for all $z$.
 
Thus it suffices to work with the lower Green's function $\underline{G}(z,\zeta;E_v)$. 
\index{Green's function!lower} 
By (\cite{RR1}, Theorem 4.2.15) $E_v$ can be
uniquely written in the form 
\begin{equation*}
E_v \ = \ \bigcap_{j=1}^M D_j
\end{equation*}
where $D_1, \ldots, D_M$ are `$PL$-domains' with pairwise disjoint complements.
\index{$\PL$-domain} 
For each $\zeta \notin E_v$, there is a unique $D_j$ with $\zeta \notin D_j$,
and by (\cite{RR1}, Theorem 4.2.12) there is a function $h_\zeta(z) \in \CC_v(\cC)$,
whose only pole is at $\zeta$, such that 
\begin{equation*}
D_j \ = \ \{z \in \cC_v(\CC_v) : |h_\zeta(z)|_v \le 1 \} \ .
\end{equation*}  
By (\cite{RR1}, Corollary 4.2.13), $D_j$ is minimal among $PL$-domains containing $E_v$
\index{$\PL$-domain}
whose complement contains $\zeta$.  

By (\cite{RR1}, Proposition 4.4.1), if $\deg(h_{\zeta}) = N$, then 
\begin{equation*} 
\underline{G}(z,\zeta;E_v) \ = \ \underline{G}(z,\zeta;D_j) \ = \  
    \left\{ \begin{array}{ll} 
            \frac{1}{N} \log_v(|h_\zeta(z)|_v) & \text{if $z \notin D_j$ \ ,} \\
            0                                    & \text{if $z \in D_j$ \ .}  
            \end{array} \right. 
\end{equation*} 
It follows that $\Vbar_{\zeta}(E_v) \in \QQ$ 
and that $\Gbar(z,\zeta;E_v) \in \QQ$ for all $z \ne \zeta$.
\end{proof}

\smallskip
The next proposition plays a key role in the reduction of Theorem \ref{aT1} 
to Theorem \ref{aT1-B}.  We begin with a definition.  

\begin{definition} \label{AnalyticallyAccessible}
Let a set $E_v \subset \cC_v(\CC_v)$ and a subset $E_v^0 \subset E_v$ be given.
Let $z_0 $ be a point in $E_v$.   
If $v \in \cM_K$ is archimedean, 
we will say that $z_0$ is {\em analytically accessible}
 from $E_v^0$ if for some $r > 0$, there is a nonconstant analytic map 
$f : D(0,r) \rightarrow \cC_v(\CC)$ with
$f(0) = z_0$, such that $f((0,r]) \subset E_v^0$. 
\index{analytically accessible|ii} 

If $v$ is nonarchimedean, we will say that $z_0$ is analytically accessible
if there are an isometrically parametrizable ball $B(z_0,r)$
\index{isometrically parametrizable ball}
and a ($\CC_v$-rational) isometric parametrization $f : D(0,r) \rightarrow B(z_0,r)$
\index{isometric parametrization}
with $f(0) = z_0$, 
such that $f((\cO_v \cap D(0,r)) \backslash \{0\}) \subset E_v^0$. 
\end{definition} 

\begin{proposition} \label{IdentifyGreenProp}  
Let $E_v \subset \cC_v(\CC_v)$ be a compact set of positive capacity, 
\index{Green's function!of a compact set}
\index{capacity $> 0$}
and let a subset $E_v^0 \subset E_v$ be given.  Suppose there is a Borel subset
$e \subset E_v$ of inner capacity $0$ such that each point of 
\index{capacity $= 0$}
$E_v \backslash e$ is analytically accessible from $E_v^0$.  
Then for each $\zeta \in \cC_v(\CC_v) \backslash E_v$,
we have $\Vbar_{\zeta}(E_v^0) = V_{\zeta}(E_v)$, 
and $\Gbar(z,\zeta;E_v^0) = G(z,\zeta;E_v)$ for all $z \ne \zeta$.  
\end{proposition} 
\index{Green's function!identifying|ii}

\begin{proof} Fix $\zeta \notin E_v$. 
We begin by showing that $\Gbar(z,\zeta;E_v^0) = G(z,\zeta;E_v) = 0$ 
for each $z \in E_v \backslash e$.  

\smallskip 
First assume $v$ is archimedean. 
Recall (\cite{Kl}, p.53) that a set $Z \subset \cC_v(\CC)$ is said to be {\em thin}\index{thin set|ii}
at a point $z_0$ if either $z_0$ is not a limit point of $Z$, 
or there are a neighborhood $V$ of $z_0$ and a 
subharmonic function $u(z)$ on $V$ such that
\begin{equation*}
\limsup_{\substack{x \rightarrow z_0 \\ x \in Z \backslash \{z_0\} }} u(x) \ < \ u(z_0) \ .
\end{equation*} 
Fix $z_0 \in E_v \backslash e$.  
Since $z_0$ is analytically accessible from $E_v^0$, 
there is a nonconstant analytic map $f: D(0,r) \rightarrow \cC_v(\CC)$
with $f(0) = z_0$ such that $f((0,r]) \subseteq E_v^0$. 
Without loss we can assume $r$ is small enough that $\zeta \notin f(D(0,r))$.
By (\cite{Kl},Corollary 4.8.5), $f([0,r])$ is not thin at $z_0$.  On the other hand,\index{thin set}
for each $0 < \varepsilon < r$, the set $H_{\varepsilon} := f([\varepsilon,r])$ 
is a compact continuum contained in $E_v^0$.  By Proposition \ref{UpGreenProp1}.3,
$G(z,\zeta;H_{\varepsilon})$ is identically $0$ on $H_{\varepsilon}$.  Thus
$\Gbar(z,\zeta;E_v^0)$ is identically $0$ on $f((0,r])$.  
By Proposition \ref{UpGreenProp1}.5,  $\Gbar(z,\zeta;E_v^0)$ is subharmonic 
and non-negative on $f(D(0,r))$, so 
\begin{equation*} 
0 \ \le \ G(z_0,\zeta;E_v) \ \le \ \Gbar(z_0,\zeta;E_v^0) \ \le \ 
\limsup_{\substack{x \rightarrow z_0 \\ x \in f([0,r]) \backslash \{z_0\} }}
\Gbar(x,\zeta;E_v^0) \ = \ 0 \ .
\end{equation*} 
Thus $G(z_0,\zeta;E_v) = \Gbar(z_0,\zeta;E_v^0) = 0$.

Next suppose $v$ is nonarchimedean.  Given $z_0 \in E_v \backslash e$, 
let $f : D(0,r) \rightarrow B(z_0,r) \subset \cC_v(\CC_v)$ be an isometric
parametrization with $f(0) = z_0$, such that 
$f((D(0,r)\backslash \{0\}) \cap \cO_v) \subseteq E_v^0$. Without loss we 
can assume $r$ is small enough that $\zeta \notin B(z_0,r)$, 
and that $r \in |K_v^{\times}|$.   
For each $0 < \varepsilon < r$, 
put $K_{\varepsilon} = (D(0,r) \backslash D(0,\varepsilon)^-) \cap \cO_v$ and 
put $H_{\varepsilon} = f(K_{\varepsilon}) \subset E_v^0$.  
Since $B(z_0,r)$ is an isometrically parametrizable ball disjoint from $\zeta$,
\index{isometrically parametrizable ball}
by Proposition \ref{APropA2}.B.2, there is a constant $C$ such that 
$[z,w]_{\zeta} = C\|z,w\|_v$ for all $z, w \in B(z_0,r)$.  Pulling this back to
$D(0,r)$, we see that $[f(x),f(y)]_{\zeta} = C|x-y|$ for all $x, y \in D(0,r)$,

Regarding $K_{\varepsilon}$ as a subset of $\CC_v = \PP^1(\CC_v) \backslash \{\infty\}$ and considering the definitions of $G(z,\infty;K_{\varepsilon})$ 
and $G(z,\zeta;H_{\varepsilon})$, it follows that 
$G(f(x),\zeta;H_{\varepsilon}) = G(x,\infty;K_{\varepsilon})$
for each $x \in D(0,r)$.
By the explicit computation in Proposition \ref{FmCosets}, 
together with a simple scaling argument, we have 
\begin{equation*}
\lim_{\varepsilon \rightarrow 0^+} G(z_0,\infty;H_{\varepsilon}) \ = \ \lim_{\varepsilon \rightarrow 0^+} G(0,\infty;K_{\varepsilon}) \ = \ 0 \ .
\end{equation*} 
Since
\begin{equation*} 
0 \ \le \ G(z_0,\zeta;E_v) \ \le \ \Gbar(z_0,\zeta;E_v^0) 
\ \le \ \lim_{\varepsilon \rightarrow 0} G(z_0,\zeta;H_{\varepsilon}) \ = \ 0 \ , 
\end{equation*} 
once again we see that $G(z_0,\zeta;E_v) = \Gbar(z_0,\zeta;E_v^0) = 0$.

\smallskip
We next show that $\Vbar_{\zeta}(E_v^0) = V_{\zeta}(E_v)$. 
The argument is very similar to the one in Proposition \ref{InnerSeqProp},
but since the context is somewhat different we give the details.  

Write $\hV = \Vbar_{\zeta}(E_v^0)$.  
By Proposition \ref{DecreasingSeqLemma}, 
there is an increasing sequence of compact sets 
$H_{v,1} \subseteq H_{v,2} \subseteq  \cdots \subseteq E_v^0$ such that 
$\lim_{n \rightarrow \infty} G(z,\zeta;H_{v,n}) = \Gbar(z,\zeta;E_v^0)$
for each $z \in \cC_v(\CC_v) \backslash \{\zeta\}$,
and $\lim_{n \rightarrow \infty} V_{\zeta}(H_{v,n}) = \hV$;   
we can assume without loss that each $H_{v,n}$ has positive capacity.
\index{capacity $> 0$}  
Write $\mu_n$ for the equilibrium distribution of $H_{v,n}$ with
respect to $\zeta$, and let $\mu$ be the equilibrium distribution
of $E_v$ with respect to $\zeta$.  After replacing $\{H_{v,n}\}_{n \ge 1}$
with a subsequence, if necessary, we can assume that the measures $\mu_n$
converge weakly to a probability measure $\widehat{\mu}$ on $E_v$. 

By what has been shown above,  
the potential function\index{potential function}  $u_{E_v}(z,\zeta) = V_{\zeta}(E_v) - G(z,\zeta;E_v)$ 
is identically equal to $V_{\zeta}(E_v)$ on $E_v \backslash e$.  
Since a set of inner capacity $0$
\index{capacity $= 0$}
has mass $0$ under any positive measure with a finite energy integral,
\index{energy integral} 
for each $n$ the Fubini-Tonelli theorem gives\index{Fubini-Tonelli theorem}
\begin{eqnarray} 
\int_{E_v} u_{H_{v,n}}(z,\zeta) \, d\mu(z)  & = & 
\iint_{E_v \times H_{v,n}} -\log_v([z,w]_{\zeta}) \, 
                      d\mu_n(w) d\mu(z) \notag \\ 
& = & \int_{H_{v,n}} u_{E_v}(w,\zeta) \, d\mu_n(w) \ = \ V_{\zeta}(E_v) \ .
\label{FFT1} 
\end{eqnarray} 
On the other hand, pointwise for each $z_0 \in E_v \backslash e$,  we have 
\begin{eqnarray*} 
\lim_{n \rightarrow \infty} u_{H_{v,n}}(z_0,\zeta) 
& = & \lim_{n \rightarrow \infty}  
             \big(V_{\zeta}(H_{v,n}) - G(z_0,\zeta;H_{v,n})\big) 
 \ = \ \hV - \Gbar(z_0,\zeta;E_v^0) \ = \ \hV \ .
\end{eqnarray*} 
Since $E_v$ is bounded away from $\zeta$, there is a constant $B_1 > -\infty$
such that $u_{H_{v,n}}(z,\zeta) \ge B_1$  on $E_v$, for all $n$.
On the other hand, since $u_{H_{v,n}}(z,\zeta) \le V_{\zeta}(H_{v,n})$ for all $z$
and the $V_{\zeta}(H_{v,n})$ are decreasing with $n$, there is a $B_2 < \infty$ 
such that $u_{H_{v,n}}(z,\zeta) \le B_2$ on $E_v$, for all $n$. 
By (\ref{FFT1}) and the Dominated Convergence Theorem,\index{Dominated Convergence theorem} it follows that 
\begin{equation*} 
V_{\zeta}(E_v) \ = \ 
\lim_{n \rightarrow \infty} \int_{E_v \backslash e} 
          u_{H_{v,n}}(z,\zeta) \, d\mu(z) \\
\ = \ \int_{E_v \backslash e} \hV \, d\mu(z) \ = \ \hV \ .
\end{equation*}  
Hence $\Vbar_{\zeta}(E_v^0) = \lim_{n \rightarrow \infty} V_{\zeta}(H_{v,n}) 
= \hV = V_{\zeta}(E_v)$.

\smallskip
The remainder of the proof is identical 
to the part of the proof of Proposition \ref{InnerSeqProp} after
formula (\ref{FPAD1A}).  
Using the Monotone Convergence Theorem,\index{Monotone Convergence theorem} 
one shows that $I_{\zeta}(\widehat{\mu}) = V_{\zeta}(E_v)$.
The uniqueness of the equilibrium distribution implies $\widehat{\mu} = \mu$, and
this in turn yields $\Gbar(z,\zeta;E_v^0) = G(z,\zeta;E_v)$ for all $z \ne \zeta$.  
\end{proof}


\section{Green's Matrices and the Inner Cantor Capacity} 
\index{capacity!inner Cantor capacity|ii}
\label{CantorCapacitySection} 

Let notations and assumptions be as in \S\ref{Chap3}.\ref{AssumptionsSection}.  
Thus, $K$ is a global field and $\cC/K$ is a curve.  
We are given a finite, galois-stable set of points 
$\fX = \{x_1, \ldots, x_m\} \subset \cC(\tK)$, 
and a $K$-rational adelic set $\EE = \prod_v E_v$ compatible with $\fX$:\index{compatible with $\fX$}
each $E_v \subset \cC_v(\CC_v)$ is nonempty, stable under 
$\Aut_c(\CC_v/K_v)$, and bounded away from $\fX$, with
$E_v$ being $\fX$-trivial for all but finitely many $v$.
\index{$\fX$-trivial} 
For each $x_i \in \fX$ we are given a 
uniformizing parameter\index{uniformizing parameter!normalizes Robin constant} $g_{x_i}(z) \in K(\cC)$, 
with $g_{\sigma(x_i)}(z) = \sigma(g_{x_i})(z)$ for each $\sigma \in \Aut(\tK/K)$.  

\smallskip
In this section we introduce the inner Cantor capacity, 
\index{capacity!inner Cantor capacity|ii} 
extending the definitions and results from (\cite{RR1}, \S5.3) 
to arbitrary $K$-rational sets $\EE$ compatible with $\fX$.\index{compatible with $\fX$}  
When each $E_v$ is algebraically capacitable (in particular, if each $E_v$ is compact or $\fX$-trivial)
\index{$\fX$-trivial}
\index{algebraically capacitable} 
then the upper Green's matrix $\Gammabar(\EE,\fX)$ and the inner Cantor capacity $\gammabar(\EE,\fX)$
\index{capacity!inner Cantor capacity} 
defined here coincide with the Green's matrix $\Gamma(\EE,\fX)$ and the Cantor capacity $\gamma(\EE,\fX)$
\index{Green's matrix!global}
\index{capacity!Cantor capacity} 
from (\cite{RR1}, \S5.3).   

\smallskip
To define the inner Cantor capacity $\gammabar(\EE,\fX)$, 
\index{capacity!inner Cantor capacity}
it is necessary to first make a base change to $L = K(\fX)$.  
For each place $v$ of $K$, and 
each place $w$ of $L$ over $v$, fix a continuous isomorphism 
$\iota_{w/v} : \CC_w \rightarrow \CC_v$, and put
$E_w = \iota_{w/v}^{-1}(E_v) \subset \cC_w(\CC_w)$. Since $E_v$
is stable under $\Aut_c(\CC_v/K_v)$, the set $E_w$ is independent of the choice 
of $\iota_{w/v}$.  The Green's functions $\Gbar(z,x_i;E_w)$ and Robin constants 
\index{Green's function!upper}
\index{Robin constant!upper}
$\Vbar_{x_i}(E_w)$ are defined the same way as the corresponding objects over 
$K$, but using the normalized absolute values $|x|_w$ and the normalized logarithms 
$\log_w(z)$, which means that 
\begin{equation*}
\Gbar(x_i,x_j;E_w) = e_{w/v} \Gbar(x_i,x_j;E_v)\ , 
        \quad \Vbar_{x_i}(E_w) = e_{w/v} \Vbar_{x_i}(E_v)
\end{equation*}
where $e_{w/v}$ is the ramification index (for archimedean places, our 
convention is that $e_{w/v} = 1$).  Put $\EE_L = \prod_w E_w$.  
        
For each place $w$ of $L$, we first define the {\em local upper Green's matrix} by
\begin{equation*}
\Gammabar(E_w,\fX) \ = \ \left( \begin{array}{cccc}
            \Vbar_{x_1}(E_w)   & \Gbar(x_1,x_2;E_w) & \hdots & \Gbar(x_1,x_m;E_w) \\ 
            \Gbar(x_2,x_1;E_w) & \Vbar_{x_2}(E_w)   & \hdots & \Gbar(x_2,x_m;E_w) \\
               \vdots      &   \vdots       & \ddots &   \vdots       \\
            \Gbar(x_m,x_1;E_w) & \Gbar(x_m,x_2;E_w) & \hdots &  \Vbar_{x_m}(E_w)  
                       \end{array} \right) \ .
\end{equation*}
For all but finitely many  $w$, each $g_{x_i}(z)$ has good reduction at $w$
\index{good reduction}
and $E_w$ is $\fX$-trivial; for such $w$, $\Gamma(E_w,\fX)$ is the zero matrix
\index{Green's matrix!local}
\index{$\fX$-trivial}
(see \cite{RR1}, Proposition 5.1.2).    
 
The global upper Green's matrix over $L$ is defined by 
\begin{equation} \label{aF2}
\Gammabar(\EE_L,\fX) \ = \ \sum_w \Gammabar(E_w,\fX) \log(q_w) \ .
\end{equation}
By our remarks above, this is actually a finite sum.  

The local and global upper Green's matrices over $K$ are then defined by
\begin{eqnarray*}
\Gammabar(E_v,\fX) \log(q_v) 
            & = & \frac{1}{[L:K]} \sum_{w|v} \Gammabar(E_w,\fX) \log(q_w) \ , \\
\Gammabar(\EE_K,\fX) & = & \frac{1}{[L:K]} \Gammabar(\EE_L,\fX) \ ,
\end{eqnarray*}
so that\index{$q_v$!weights $\log(q_v)$ in $\Gamma(\EE,\fX)$|ii}
\begin{equation*}
\Gammabar(\EE_K,\fX) \ = \ \sum_v \Gammabar(E_v,\fX) \log(q_v) \ .
\end{equation*} 
The entries of $\Gammabar(\EE_K,\fX)$ are finite if and only if 
each $E_v$ has positive innner capacity.  
\index{capacity $> 0$}
Clearly $\Gammabar(\EE_K,\fX)$ and the $\Gammabar(E_v,\fX)$  
are symmetric and non-negative off the diagonal;  
they are also $K$-symmetric in the sense
of Definition \ref{DefD1}, as shown by the following lemma:    

\begin{lemma}  \label{aL1} 
Let $\EE$, $\fX$, $L/K$, and the $g_{x_i}(z)$ be as above.  
For each $\sigma \in \Aut(L/K)$, each place $w$ of $L$, and each $i \ne j$,  
\begin{eqnarray*}
\Gbar(x_i,x_j;E_w) & = & \Gbar(x_{\sigma(i)},x_{\sigma(j)},E_{\sigma(w)}) \ , \\
\Vbar_{x_i}(E_w) & = & \Vbar_{x_{\sigma(i)}}(E_{\sigma(w)}) \ .
\end{eqnarray*}
\end{lemma}  

\begin{proof}  This is essentially a tautology.  

We can view $\cC(L)$ as embedded on the diagonal in $\oplus_{w|v} \cC_w(L_w)$.  
Each $\sigma \in \Aut(L/K)$ acts on $L \otimes_K K_v$ through its action on $L$.  
Using the canonical isomorphism $L \otimes_K K_v \cong \oplus_{w|v} L_w$,
(which holds when $\Char(K) = p > 0$, as well as when $\Char(K) = 0$; see (\cite{RR1}, p.321))
this action can also be described by a collection of isomorphisms  
$\tau_{\sigma,w} : L_w \rightarrow L_{\sigma(w)}$.  

The action of $\sigma$ on $\fX \subset \cC(L)$ is described globally 
by the permutation representation $\sigma(x_i) = x_{\sigma(i)}$. 
When this is combined with the semilocal description of its
action  on $\oplus_{w|v} \cC_w(L_w)$, after identifying the points
in $\fX$ with their images in $\cC_w(L_w)$ and $\cC_{\sigma(w)}(L_{\sigma(w)})$,  
we find that $\tau_{\sigma,w}$ takes $x_i \in \cC_w(L_w)$ to 
$x_{\sigma(i)} = \sigma(x_i) \in \cC_{\sigma(w)}(L_{\sigma(w)})$.  

Extend each $\tau_{\sigma,w}$ to a continuous isomorphism 
$\overline{\tau}_{\sigma,w} : \CC_w \rightarrow \CC_{\sigma(w)}$.  
Under this isomorphism $E_w \subset \cC_w(\CC_w)$ is taken to $E_{\sigma(w)}$ 
since both are pullbacks of $E_v$, which is stable under $\Aut_c(\CC_v/K_v)$.  
Viewing $\overline{\tau}_{\sigma,w}$ as an identification, we then have 
$\Gbar(z,x_j,E_w) 
= \Gbar(\overline{\tau}_{\sigma,w}(z),x_{\sigma(j)};E_{\sigma(w)})$
for all $z \in \cC_w(\CC_w)$.  Both assertions in the lemma now follow;  
the second uses Proposition \ref{APropM2} and the fact that the local uniformizers
$g_{x_i}(z)$ are $K$-symmetric. 
\end{proof}
 
\vskip .1 in
We can now define the inner Cantor capacity.  
\index{capacity!inner Cantor capacity|ii} 
Let the set of  $m$-dimensional real {\em probability vectors} be 
\begin{equation*}
\cP^m \ = \ \cP^m(\RR) \  = \ 
\{ \vs \in \RR^m : \text{$\sum s_i = 1$, $s_i \ge 0$ for each $i$} \} \ ,
\end{equation*}
The global upper Robin constant $\Vbar(\EE_K,\fX)$ is the value of
\index{Robin constant!upper global} 
$\Gammabar(\EE_K,\fX)$ as a matrix game:
\begin{eqnarray}  
\Vbar(\EE_K,\fX) & = & \val(\Gammabar(\EE_K,\fX)) 
     \ := \  \max_{\vs \in \cP^m} \ \min_{\vr \in \cP^m} 
          \  \phantom{}^t \vs \Gammabar(\EE_K,\fX) \vr \label{MinMax} \\
    & = & \min_{\vs \in \cP^m} \max_i (\Gammabar(\EE_K,\fX) \vs)_i \ .
                    \label{MinMax1} 
\end{eqnarray}
The equality of (\ref{MinMax}) and (\ref{MinMax1}), along with many other similar
expressions, follows from the Fundamental theorem of Game Theory 
(see \cite{RR1}, p.327).  The inner Cantor capacity is then defined by 
\index{capacity!inner Cantor capacity|ii} 
\begin{equation*}
\gammabar(\EE_K,\fX) \ = \ e^{-\Vbar(\EE_K,\fX)} \ .
\end{equation*}
Clearly $\gammabar(\EE,\fX) > 0$ if and only if each $E_v$ has positive inner
capacity.

The inner Cantor capacity has properties like those of the classical
Cantor capacity  in (\cite{RR1}).  All of these are formal consequences of 
\index{capacity!Cantor capacity} 
properties of matrices and Green's functions, so the proofs given in (\cite{RR1})
\index{Green's function}
carry over without change.  The most important properties are as follows:    

\begin{proposition} \label{NegDef} 
$\gammabar(\EE,\fX) > 1$ iff \ $\Gammabar(\EE,\fX)$ is negative definite.\index{Green's matrix!negative definite} 
\end{proposition} 

\begin{proof}  See (\cite{RR1}, Proposition 5.1.8, p.331).  
The quantity $\val(\Gamma)$ is a statistic of symmetric, real-valued matrices 
such that $\val(\Gamma) < 0$ iff $\Gamma$ is negative definite.\index{Green's matrix!negative definite}
\end{proof} 

\begin{proposition} \label{UniqueSVec}
If \, $\gammabar(\EE,\fX) > 1$, then there is a unique probability vector 
$\hs = \phantom{}^t(\hs_1,\ldots,\hs_m) \in \cP^m(\RR)$   
for which 
\begin{equation*}
\Gammabar(\EE,\fX) \, \hs 
\ = \ \left( \begin{array}{c} \hV \\ \vdots \\ \hV \end{array} \right) 
\end{equation*} 
has all of its entries equal.  In this situation, $\hV = \Vbar(\EE,\fX) < 0$ and $\hs$
is $K$-symmetric, and $\hs_i > 0$ for each $i = 1, \ldots, m$.  
\end{proposition} 

\begin{proof}  The proof is the same as (\cite{RR1}, Theorem 5.1.6, p.328), 
and goes back to (\cite{Can3}).  However, 
because this proposition plays a 
key role in the proof of the Fekete-Szeg\"o theorem, 
\index{Fekete-Szeg\"o theorem with LRC} 
we give the argument here.  

For brevity, write $\Gamma = \Gammabar(\EE,\fX)$.  By Proposition \ref{NegDef}, 
$\Gamma$ is negative definite, so all of its eigenvalues are negative.\index{Green's matrix!negative definite}  
Choose $\alpha > 0$ large enough that  $\Gamma + \alpha I$ is positive definite, where
$I$ is the $m \times m$ identity matrix.  
Then each entry of $\Gamma + \alpha I$ is non-negative, 
and $0 < \lambda_i < \alpha$ for each eigenvalue $\lambda_i$ 
of $\Gamma + \alpha I$.    
It follows that the series 
\begin{eqnarray*}
-\Gamma^{-1} & = & \left(\alpha I - (\Gamma + \alpha I)\right)^{-1} 
     \ = \ \alpha^{-1}
       \left( I - \left(\frac{\Gamma + \alpha I}{\alpha}\right)\right)^{-1} \\
             & = & \alpha^{-1} \sum_{k = 0}^{\infty} 
                            \left(\frac{\Gamma + \alpha I}{\alpha}\right)^k
\end{eqnarray*}
converges, so $-\Gamma^{-1}$ has only non-negative entries.  

Set 
\begin{equation*}
\vs^{\prime} \ = 
\ -\Gamma^{-1} \left( \begin{array}{c} 1 \\ \vdots \\ 1 \end{array} \right) \ .
\end{equation*}    
By construction, each entry of $\vs^{\prime}$ is non-negative;  if some entry were $0$, 
it would mean that every entry in the corresponding row of $-\Gamma^{-1}$  were $0$, 
and hence $\Gamma^{-1} \cdot \Gamma = I$ would not be possible.  Scale $\vs^{\prime}$
so as obtain a probability vector $\hs$.  Then all the entries of $\hs$ are positive,
and we have  
\begin{equation} \label{FDM} 
\Gamma \, \hs 
\ = \ \left( \begin{array}{c} \hV \\ \vdots \\ \hV \end{array} \right) 
\end{equation}  
for some $\hV$.  The minimax inequality (\ref{MinMax}) defining $\val(\Gamma)$ then
\index{minimax property} 
shows $\hV = \Vbar(\EE,\fX) < 0$.  

The uniqueness of $\hs$ and $\hV$ satisfying (\ref{FDM}) are clear;  
since $\Gammabar(\EE,\fX)$ is $K$-symmetric,
it follows from the uniqueness that $\hs$ must be $K$-symmetric as well.
\index{$K$-symmetric!probability vector}
\end{proof} 

\vskip .1 in
Define the set of rational probability vectors to be 
\begin{equation*}\label{`SymbolIndexcPmQ'}
 \cP^m(\QQ) \  = \ \cP^m(\RR) \cap \QQ^m \ .
\end{equation*}
Unfortunately,  $\hs$ need not be rational.  This causes major technical
difficulties in the proof of the Fekete-Szeg\"o theorem.  To overcome it, 
\index{Fekete-Szeg\"o theorem with LRC} 
we show that in the number field case the ``logarithmic leading coefficients''
\index{logarithmic leading coefficients!independent variability of archimedean}  
of the archimedean initial\index{independent variability!of logarithmic leading coefficients}
approximating functions can be ``macroscopically independently varied'', 
which allows us to replace $\hs$ 
with a suitable rational $\vs$ closely approximating $\hs$.  
In the function field case, we show that if the sets $E_v$ satisfy appropriate
hypotheses, then all entries in $\Gammabar(\EE,\fX)$ 
are rational multiples of $\log(p)$, and in that case $\hs$ belongs to $\cP^m(\QQ)$.

\vskip .05 in
We conclude this section by noting some properties of 
the inner Cantor capacity.    
\index{capacity!inner Cantor capacity} 
\index{capacity!inner Cantor capacity!functoriality properties}

\begin{proposition} \label{Functoriality} { \ } 

$(1)$  {\rm (Base Change):} 
If $M/K$ is any finite extension, then $\Gammabar(\EE_M,\fX) = [M:K] \Gammabar(\EE_K,\fX)$,
$\Vbar(\EE_M,\fX) = [M:K] \Vbar(\EE_K,\fX)$, and
\begin{equation*}
\gammabar(\EE_M,\fX) \ = \ \gammabar(\EE_K,\fX)^{[M:K]}  \ . 
\end{equation*} 

$(2)$  {\rm (Pullback):}  Let $\cC_1, \cC_2/K$ be curves, and let
$f : \cC_1 \rightarrow \cC_2$ be a nonconstant rational map defined over $K$.  
If $\EE$ is a $K$-rational adelic set on $\cC_2$, compatible with $\fX \subset \cC_2(\tK)$,
then\index{compatible with $\fX$} 
\begin{equation*}
\gammabar(f^{-1}(\EE),f^{-1}(\fX)) \ = \ \gammabar(\EE,\fX)^{1/\deg(f)} \ .  
\end{equation*}  
\end{proposition} 

\begin{proof} These follow from the corresponding properties of
upper Green's functions.  
\index{Green's function!upper}
\index{Green's function!properties of}
See the proofs of (\cite{RR1}, Theorems 5.1.13 and 5.1.14, p.333). 
\end{proof} 


\section{Newton Polygons of Nonarchimedean Power Series}  \label{NewtonPolygonSection} 
\index{Newton Polygon|(}
\index{Newton Polygon|ii}

In this section we recall some some facts about Newton polygons
of nonarchimedean power series needed for the proof of Theorem \ref{CompactThm}. 

\vskip .05 in 
Given $r > 0$, write $D(0,r) = \{z \in \CC_v : |z|_v \le r \}$ for the
`closed' disc of radius $r$ in $\CC_v$, and $D(0,r)^- = \{z \in \CC_v : |z|_v < r \}$
for the `open' disc.

The Newton polygon of a power series
$f(Z) = \sum_{k=0}^{\infty} c_k Z^k \in \CC_v[[Z]]$ is the lower convex hull
of the set of points $\{ (k,\ord_v(c_k)) \}$, with a vertical side above
the point corresponding to the first nonzero coefficient.  If  $f(Z)$ is
a polynomial of degree $m$, its Newton polygon is also considered 
to have a vertical side above $(m,\ord_v(c_m))$.\index{Newton Polygon|ii} 

A power series of the form $h(Z) = 1 + \sum_{k=1}^{\infty} b_k Z^k \in \CC_v[[Z]]$ 
will be called a {\em unit power series} for $D(0,r)$ 
if $h(Z)$ is a unit in the ring of power series converging on $D(0,r)$.
This holds if and only if $|b_k|_v < 1/r^k$ for each $k \ge 1$,  
in which case $|h(z)|_v = 1$ for all $z \in D(0,r)$, and the series 
$h(Z)^{-1} = 1 + \sum_{k=1}^{\infty} b_k^{\prime} Z^k$ with $h(Z) \cdot h(Z)^{-1} = 1$
also satisfies $|b_k^{\prime}|_v < 1/r^k$ for each $k \ge 1$.  
If $h(Z) \in K_v[[Z]]$, then $h(Z)^{-1} \in K_v[[Z]]$.  

\begin{lemma} \label{BLem1} 
   Suppose $f(Z) = \sum_{k=0}^{\infty} c_k Z^k \in \CC_v[[Z]]$ converges on a disc $D(0,r)$,
where $r > 0$ belongs to the value group of $\CC_v^{\times}$.  Then 

$(A)$  $f(Z)$ has finitely many zeros in  $D(0,r)$, 
      and  $f(Z)$ can be factored as $f(Z) = g_r(Z) \cdot h_r(Z)$ 
      where $g_r(Z)$ is a polynomial having the same roots $($with multiplicities$)$ 
      as $f(Z)$ in $D(0,r)$, and where $h_r(Z) = 1 + \sum_{k=1}^{\infty} b_k Z^k$ 
      is a unit power series for $D(0,r)$. 
      In particular, $|b_k|_v < 1/r^k$ for all $k \ge 1$.       
      If $f(Z) \in F_w[[Z]]$ for some finite extension $F_w/K_v$, 
      then $g_r(Z) \in F_w[Z]$ and $h_r(Z) \in F_w[[Z]]$. 
      
$(B)$    If $f(Z)$ has $m \ge 0$ roots in $D(0,r)$ $($counted with multiplicities$)$, 
      then $(m,\ord_v(c_m))$ is a vertex of the Newton polygon of $f(Z)$, and the part 
of the Newton polygon of $f(Z)$ on and to the left of \ $(m,\ord_v(c_m))$ coincides with 
the Newton polygon of $g_r(Z)$.      
\end{lemma} 

\begin{proof}  This is well known;  we sketch the proof.  

The fact that $f(Z)$ has finitely many zeros in $D(0,r)$, 
and the existence of the factorization $f(Z) = g_r(Z) \cdot h_r(Z)$, 
follow from the Weierstrass Preparation Theorem and a change of variables
(see \cite{BGR}, p.201), or from Hensel's Lemma (as in \cite{Ar}, \S 2.5).  
The assertions about $F_w$-rationality hold because the Weierstrass Preparation 
Theorem and Hensel's Lemma are valid over any complete nonarchimedean field.  

To establish the relation between the Newton polygons of
$f(Z)$ and $g_r(Z)$, take $\beta \in \CC_v$ with $|\beta|_v = r$.  After replacing $f(Z)$ 
by $f(Z/\beta)$, which translates  
the Newton polygons of both $f(Z)$ and $g_r(Z)$ upwards
by the line $y = x \log_v(r)$, we can assume that $r = 1$.  
Suppose $f(Z)$ has $m$ roots in $D(0,1)$
and write $g(Z) = g_1(Z)$, $h(Z) = h_1(Z)$.   
 
Consider the factorization $f(Z) = g(Z) h(Z)$.  Write 
\begin{equation*}
f(Z)  =  \sum_{\ell = 0}^{\infty} c_k Z^k, \quad
g(Z)  =  \sum_{\ell=0}^{m} a_\ell Z^\ell \ ,
\end{equation*}
and expand 
\begin{equation*}
h(Z) \ = \ 1 + \sum_{j=1}^{\infty} b_j Z^j,
\qquad h(Z)^{-1} \ = \ 1 + \sum_{j=1}^{\infty} b_j^{\prime} Z^j \ .
\end{equation*}
Here $\ord_v(b_j), \ord_v(b_j^{\prime}) > 0$ for all $j \ge 1$.
By the Weierstrass Factorization theorem,\index{Weierstrass Factorization theorem}  $|a_m|_v = |c_m|_v$.  
After dividing through by $a_m$, we can assume that $g(Z)$ is monic 
and that $\ord_v(c_k) \ge 0$ for all $k$. 

Using that $g(Z) = f(Z) h(Z)^{-1}$, we see that if  $J$ is the
smallest index for which $c_J \ne 0$, then $a_k = c_k = 0$ for $k < J$,
while $a_J = c_J$.  By hypothesis, $\ord_v(a_m) = \ord_v(c_m) = 0$.  
For each $\ell$ with $J < \ell < m$, since $\ord_v(b_j^{\prime}) > 0$ for all $j \ge 1$
and $\ord_v(c_k) \ge 0$ for all $k$,  
\begin{eqnarray}
\ord_v(a_\ell) &=& \ord_v(c_\ell + c_{\ell-1} b_1^{\prime} + \ldots + c_0 b_\ell^{\prime})
            \label{BFG1} \\
   & \ge & \min(\ord_v(c_\ell), \ord_v(c_{\ell-1}), \ldots, \ord_v(c_0)) \ .
           \notag
\end{eqnarray}
If $(\ell,\ord_v(c_\ell))$ is a corner of the Newton polygon of $f(Z)$, then necessarily
$\ord_v(c_\ell) < \ord_v(c_0), \ldots, \ord_v(c_{\ell-1})$
and it follows from (\ref{BFG1}) that $\ord_v(a_\ell) = \ord_v(c_\ell)$.
Since $\ord_v(c_k) \ge 0$ for all $k$, the Newton polygon of $g(Z)$
lies on or below the initial part of the Newton polygon of $f(Z)$.
Applying the same arguments to $f(Z) = g(Z) h(Z)$, 
we see that initial part of the Newton polygon of $f(Z)$
lies on or below the Newton polygon of $g(Z)$.
\end{proof}  

For a polynomial, the absolute values of its roots and the slopes of the sides of its Newton polygon 
determine each other (see \cite{Ar}, \S 2.5):  if the Newton polygon
has a side with slope $M$, whose projection on the horizontal axis has length
$S$, then $f(Z)$ has exactly $S$ roots $\alpha$ for which $\log_v(|\alpha|_v) = M$.
In this correspondence, the vertical ray above the vertex
$(k,\ord_v(c_k))$ corresponding to the first nonzero coefficient 
is deemed to have projection length $k$.  

The correspondence holds for power series as well.  
If the initial segment of the Newton
polygon is a vertical ray above $(k,\ord_v(c_k))$, that ray is deemed to
have slope $-\infty$ and projection length $k$.  If $f(Z)$ is a polynomial,
the vertical side above its rightmost vertex is deemed to have infinite projection length.
There is also a special case when the radius of convergence \ $r$ belongs
to the value group of $\CC_v^{\times}$ and $f(Z)$ converges in $D(0,r)$.
In that situation the Newton polygon has a terminal ray\index{terminal ray!of Newton polygon} of slope $\log_v(r)$
which can have at most a finite number of vertices on it,
and the last such vertex is deemed to be the right endpoint
of the rightmost side of finite length:

\begin{proposition} \label{BProp2}
   Suppose  $f(Z) = \sum_{k=0}^{\infty} c_k Z^k \in \CC_v[[Z]]$ has radius of convergence
$r > 0$.  Then

$(A)$ The roots of $f(Z)$ correspond to the sides of the Newton polygon of
      $f(Z)$ of finite length in the same way as for a polynomial:
      for each finite length side with slope $M$ and projection length $S$,
      $f(Z)$ has exactly $S$ roots $\alpha_{M,1}, \ldots, \alpha_{M,S}$
      $($listed with multiplicity$)$ for which $\log_v(|\alpha_{M,i}|_v) = M$.

$(B)$ If $f(Z) \in F_w[[Z]]$, where $F_w/K_v$ is a finite extension,
       then given a side of the Newton polygon with slope $M$ 
      and projection length $S$, the polynomial $\prod_{i=1}^S (Z-\alpha_{M,i})$ 
      belongs to $F_w[Z]$.  In particular, if $S = 1$,
      the unique associated root is rational over $F_w$.  

$(C)$  The Newton polygon of $f(Z)$ is completely determined by the absolute value 
       of the first nonzero coefficient of $f(Z)$, the absolute values of the roots of $f(Z)$, 
       and the radius of convergence $r$.     
\end{proposition}
     
\begin{proof}  
This too is well known.  The domain of convergence of $f(Z)$ 
is either $D(0,r)$ or $D(0,r)^-$.  Exhausting it by discs $D(0,r_1)$ where $0 < r_1 \le r$
belongs to  $|\CC_v^{\times}|_v$, the roots of $f(Z)$ are accounted for by the roots of the polynomials $g_{r_1}(Z)$ in Lemma \ref{BLem1}. 
By that Lemma and properties of Newton polygons of polynomials proved in (\cite{Ar}, \S 2.5), 
assertions (A) and (B) hold.  

For (C), note that the absolute value of the first nonzero coefficient determines the location
of the first corner of the Newton polygon.  The absolute values of the nonzero
roots determine the lengths and slopes of the sides of the Newton polygon of
finite length.  If there are infinitely many roots, the Newton polygon is 
completely determined;  if there are only finitely many roots, 
the Newton polygon has a terminal ray\index{terminal ray!of Newton polygon} with slope $\log_v(r)$.
\end{proof}

The following concept will play an important role in nonarchimedean constructions  
throughout the paper:  

\begin{definition}  \label{ScaledIsometryDef} 
Let $K_v$ be nonarchimedean.  Let $B(a,\rho) \subset \cC_v(\CC_v)$ be an 
isometrically parametrizable ball and let $D(b,r) \subset \CC_v$ be a disc.
\index{isometrically parametrizable ball}  
We will call a map $f : B(a,\rho) \rightarrow D(b,r)$ 
a {\em scaled isometry} if $f$ is a $1-1$ correspondence  
satisfying $|f(z_1)-f(z_2)|_v = (r/\rho) \|z_1,z_2\|_v$ for all $z_1, z_2 \in B(a,\rho)$. 
Given two discs $D(a,\rho)$, $D(b,r)$, or two balls $B(a,\rho)$, $B(b,r)$, 
we define scaled a  isometry $f : D(a,\rho) \rightarrow D(b,r)$ 
or $f : B(a,\rho) \rightarrow B(b,r)$ in a similar way. 
\end{definition} 
\index{scaled isometry|ii}

We next give a criterion for a map defined by a power series to induce a scaled isometry.
\index{scaled isometry!power series map induces|ii} 

\begin{proposition} \label{ScaledIsometryProp} 
Suppose $f(Z) = \sum_{n=0}^{\infty} c_n Z^N \in \CC_v[[Z]]$ converges on a disc $D(0,r)$, with $r > 0$. 
If $f$ has a single zero in $D(0,r)$ $($counted with multiplicity$)$, 
then $f$ induces a scaled isometry from $D(0,r)$ onto $D(0,R)$ where $R = |c_1|_v r$.  
\index{scaled isometry}
If there is a finite extension $F_w/K_v$ in $\CC_v$ such that $f(Z) \in F_w[[Z]]$, 
then $f$ maps $F_w \cap D(0,r)$ onto $F_w \cap D(0,R)$.

If $H(Z) = \sum_{n=0}^{\infty} b_n Z^n$ is another power series which converges on $D(0,r)$, 
and $H$ has no zeros in $D(0,r)$, then $|H(z)|_v$ takes the constant value $B = |b_0|$ for 
all $z \in D(0,r)$, and $f \cdot H$ induces a scaled isometry from $D(0,r)$ 
\index{scaled isometry}
onto $D(0,BR)$.  If there is a finite extension $F_w/K_v$ in $\CC_v$ such that 
$f(Z), H(Z)\in F_w[[Z]]$, then $f \cdot H$ maps $F_w \cap D(0,r)$ onto $F_w \cap D(0,BR)$.
\end{proposition}

\begin{proof} Since $f$ converges on $D(0,r)$ and has a single zero there, 
the Newton polygon of $f(Z)$ lies on or above the line $y = (x-1) \log_v(r) + \ord_v(c_1)$, 
\index{Newton Polygon}
and the points $(n,\ord_v(c_n))$ for $n \ge 2$ lie strictly above that line.  
Thus $|c_0|_v \le |c_1|_v r$ and $|c_n|_v < |c_1|_v/r^{n-1}$ for $n \ge 2$.

By the ultrametric inequality, for each $a \in D(0,r)$ we have $|f(a)|_v \le |c_1|_v r$.  
On the other hand, for each $b \in D(0,|c_1|_v r)$ the Newton polygon of $f(Z)-b$ has the same 
properties as that of $f(Z)$, so there is a unique point $a \in D(0,r)$ such that $f(a) = b$.  
Thus $f$ induces a $1-1$ correspondence from $D(0,r)$ onto  $D(0,|c_1|_v r)$. 
In particular, $|c_1|_v = R/r$. 
To see that $f$ is a scaled isometry, note that for all $z, w \in D(0,r)$ 
\index{scaled isometry}
\begin{equation*}
|f(z)-f(w)|_v 
       \ = \ |z-w|_v \cdot |c_1 + \sum_{n=2}^{\infty} c_n \Big(\sum_{k=0}^{n-1} z^k w^{n-1-k}\Big)|_v 
       \ = \ \ |c_1|_v \cdot |z-w|_v \ ,
\end{equation*}
where the last step follows from the ultrametric inequality using $|c_n|_v <|c_1|_v/r^{n-1}$.

If $f(Z) \in F_w[[Z]]$ for some finite extension $F_w/K_v$, clearly $f$ maps $F_w \cap D(0,r)$
into $F_w \cap D(0,|c_1|_v r)$.  On the other hand, if $b \in F_w \cap D(0,|c_1|_v r)$ 
then $f(Z)-b \in F_w[[Z]]$ and so the unique solution to $f(a) - b = 0$ 
belongs to $F_w \cap D(0,r)$ by Proposition \ref{BProp2}(B).  
Thus $f$ maps $F_w \cap D(0,r)$ onto $F_w \cap D(0,|c_1|_v r)$.

If $H(Z) = \sum_{n=0}^\infty b_n Z^n \in \CC_v[[Z]]$ is another power series 
which converges on $D(0,r)$, but has no zeros there, then $H(Z) = b_0 \cdot h(Z)$
where $h(Z)$ is a unit power series.  It follows that $|H(z)|_v = |b_0|_v$ for all $z \in D(0,R)$,
and $|b_n|_v < B/r^n$ for all $n \ge 1$.  
If we write $f(Z) H(Z) = \sum_{n=0}^{\infty} a_n Z^n$, 
then  $|a_1|_v = |b_0 c_1 + b_1 c_0|_v = BR/r$ since $|b_0 c_1|_v = BR/r$ 
while $|c_0 b_1|_v < BR/r$.  We can now apply the previous discussion to $f \cdot H$. 
\end{proof}

\index{Newton Polygon|ii}
\index{Newton Polygon|)} 

\section{Stirling Polynomials} 
\label{StirlingPolynomialSection}

In this section we will consider $v$-adic Stirling polynomials, 
which play an important role in the construction of the nonarchimedean initial patching functions.
\index{patching functions, initial $G_v^{(0)}(z)$}   
They are also used in `degree-raising' arguments in the global patching construction.
\index{patching argument!global}\index{degree-raising}        

\vskip .05 in
Let $F_w/K_v$ be a finite extension in $\CC_v$.  
Let $e_w = e_{w/v}$ be its ramification index 
and $f_w = f_{w/v}$ its residue degree, so $e_w f_w = [F_w:K_v]$.  

Fix a prime element $\pi_w$ for $\cO_w$,
and let $\ord_w(x)$ be the valuation on $\CC_v$ for which $\ord_w(\pi_w) = 1$;  
thus $\ord_w(x) = e_w \ord_v(x)$. 
Put  $q = q_w = q_v^{f_w}$, 
and write $\log_w(x)$ for the logarithm to the base $q_w$.
 
We first construct an explicit uniformly 
distributed sequence of points in $\cO_w$. Let 
$\psi_w(0), \ldots, \psi_w(q-1)$ be the Teichm\"uller representatives\index{Teichm\"uller representatives|ii} 
for the cosets of  $\cO_w/\pi_w \cO_w \cong \FF_q$, with $\psi_w(0) = 0$.  
Thus when $\Char(K) = 0$, the representatives are $0$ and the $q-1^{st}$ roots of unity in $F_w$.
When $\Char(K) = p > 0$, so $F_w \cong \FF_q((\pi_w))$, they are the 
the elements of $\FF_q$.  This is needed for the global patching argument  
(see Lemma \ref{QTLem} and Proposition \ref{SnCharpForm}).
\index{patching argument!global}

Extend $\{ \psi_w(k) \}_{0 \le k < q}$  
to a sequence $\{ \psi_w(k) \}_{0 \le k < \infty}$ as follows:  
for each $k \ge q$, write $k$ in base $q$ as
\begin{equation*}
k \ = \ \sum_{i=0}^N d_i(k) q^i
\end{equation*}
where $N = \lfloor \log_v(k) \rfloor$ and $0 \le d_i(k) \le q-1$
for each $i$, then put
\begin{equation*}
\psi_w(k) \ = \ \sum_{i=0}^N \psi_w(d_i(k)) \pi_w^i \ \in \ \cO_w \ .
\end{equation*}
For each $m$, the sequence $\{ \psi_w(k) \}_{0 \le k < q^m}$
is a system of coset representatives for $\cO_w/\pi_w^m\cO_w$.  

Define a function $\val_w(k)$\label{`SymbolIndexValw'} for integers $k \ge 0$ by letting 
$\val_w(k)$ be the smallest $i$ for which $d_i(k) \ne 0$, if $k > 0$, 
and putting $\val_w(0) = \infty$.  For each $k$ it follows that  
\begin{equation} \label{BGF0}
\ord_w(\psi_w(k)) \ = \ \val_w(k) \ . 
\end{equation}
Similarly, for all $k  \ne \ell$,
\begin{equation} \label{BGF1}
\ord_w(\psi_w(k) - \psi_w(\ell)) \ = \ \val_w(|k-\ell|) \ .
\end{equation}
This is because if $\val_w(|k-\ell|) = j$, then the digits $d_{i}(k)$ and $d_i(\ell)$ 
coincide for $i < j$, while $d_j(k) \ne d_j(\ell)$.
It follows that for each $n > 0$, if $0 \le k, \ell < n$ and $k \ne \ell$, 
then by (\ref{BGF1})
\begin{equation} \label{BGGG1}
\ord_w(\psi_w(k)-\psi_w(\ell)) \ < \  \log_w(n) \ ,
\end{equation}
since $\val_w(|k-\ell|) \le \log_w(|k-\ell|) < \log_w(n)$.

\begin{definition} \label{BGD1}
The {\em basic well-distributed sequence for $\cO_w$} is 
$\{ \psi_w(k) \}_{0 \le k < \infty}$. 
\index{basic well-distributed sequence in $\cO_w$|ii}
\label{`SymbolIndexWellDistv'}
The {\em Stirling polynomial of degree $n$ for $\cO_w$} is
\index{Stirling polynomial!for $\cO_w$|ii} 
\label{`SymbolIndexStirlingv'}
\begin{equation} \label{FST1} 
S_{n,w}(z) \ = \ \prod_{k=0}^{n-1} (z-\psi_w(k)) \ .
\end{equation}
\end{definition}

The polynomials $S_{n,w}(z)$ were first studied by Polya (\cite{P}), 
and were used by Cantor in (\cite{Can3}).
The following proposition, which is similar to (\cite{RR3}, Lemma 8.7), summarizes their main properties:

\begin{proposition} \label{BGProp1}  Let $F_w/K_v$ be a finite, separable extension in $\CC_v$.  
Let $S_{n,w}(z)$ be the Stirling polynomial of degree $n$ for $\cO_w$,
\index{Stirling polynomial!for $\cO_w$}
and let $S_{n,w}^{\prime}(z)$ be its derivative. Then

$(A)$ If $0 \le i, j < n$ and $i \ne j$ then $|\psi_v(i)-\psi_v(j)|_v > 1/n$.  

$(B)$ For each $k$, $0 \le k < n$, we have 
\begin{equation} \label{BGF4}
\frac{n}{e_w(q_v^{f_w}-1)} - \frac{1}{e_w}\big(2\log_w(n) + 3\big) \ < \
    \ord_v \big( S_{n,w}^{\prime}(\psi_w(k))\big) \ < \ \frac{n}{e_w(q_v^{f_w}-1)} \ . 
\end{equation}

$(C)$ Fix $x \in \CC_v$.  If $0 \le J < n$ is such that
$|x - \psi_w(J)|_v = \min_k(|x-\psi_w(k)|_v)$,
then
\begin{equation} \label{BGF5}
 \ord_v(S_{n,w}(x)) \ < \ \frac{n}{e_w(q_v^{f_w}-1)} + \ord_v(x - \psi_w(J)) \ .
\end{equation}
If $x \notin D(0,1)$, then $\ord_v(S_{n,w}(x)) = n \ord_v(x)$.

\end{proposition}

\begin{proof}
Assertion (A) is a reformulation of (\ref{BGGG1}).  

To prove assertion (B), it will be convenient to work over $F_w$ rather than $K_v$.
Fix $ 0 \le k < n$, and note that by (\ref{BGF1}), if  $j \ne k$, 
then $\ord_w(\psi_w(k)-\psi_w(j)) = \val_w(|k-j|)$.  
This leads to a generalization of the well-known formula 
\begin{equation*} 
\sum_{\ell=1}^k \ord_p(\ell) \ =  \ord_p(k!) \ = \ \sum_{m \ge 1} \left\lfloor \frac{k}{p^m} \right\rfloor 
           \ = \ \frac{k}{p-1} - \frac{1}{p-1} \sum_{i \ge 0} a_i \ ,
\end{equation*} 
where the $a_i$ are the base $p$ digits of $k$. 
Writing $d_i = d_i(k)$, and $q = q_w$, for each $m \ge 1$ 
there are exactly $\lfloor k/q^m \rfloor$ 
integers $1 \le \ell \le k$ with $\val_w(\ell) \ge m$.  Hence 
\begin{eqnarray}
\sum_{\ell=1}^k \val_w(\ell) & = & \sum_{m \ge 1} 
                   \left\lfloor \frac{k}{q^m} \right\rfloor \notag \\
      & = &  (d_1 + d_2 q + d_3 q^2 + \cdots) + (d_2 + d_3 q + \cdots) + \cdots 
                                  \notag \\
      & = & d_1 \cdot \frac{q-1}{q-1} + d_2 \cdot \frac{q^2-1}{q-1} + 
                          d_3 \cdot \frac{q^3-1}{q-1} + \cdots \notag \\
      & = & \frac{d_0 + d_1 q + d_2 q^2 + \cdots}{q-1} - 
                 \frac{d_0 + d_1 + d_2 + \cdots}{q-1} \notag \\
      & = &  \frac{k}{q-1} - \frac{1}{q-1} \sum_{i \ge 0} d_i(k) \ .
                             \label{FactorialFormula}
\end{eqnarray}      
Consequently
\begin{eqnarray} \ord_w\big(\prod^{n-1}_{\substack{ \ell = 0 \\ \ell \ne k }} 
     (\psi_w(k)-\psi_w(\ell))\big) 
  & = & \sum_{\ell = 1}^k \val_w(\ell) + \sum_{\ell = 1}^{n-k-1} \val_w(\ell) 
         \label{FNSI1} \\
  & = & \frac{n}{q_w-1} - \frac{\sum d_i(k) + \sum d_i(n-k-1) + 1}{q_w-1} \ . \notag  
\end{eqnarray}
In particular 
\begin{equation*} 
0 < \ \sum d_i(k) + \sum d_i(n-k-1) + 1 \ < \  2(q_w-1) ( \log_w(n) + 1) + 1
\end{equation*} 
and hence 
\begin{equation} \label{BGF4A}
\frac{n}{q_w-1} -  2 \log_w(n) - 3 \ < \
    \ord_w \big( S_{n,w}^{\prime}(\psi_w(k))\big) \ < \ \frac{n}{q_w-1} \ . 
\end{equation}
Since $\ord_w(x) = e_w \ord_v(x)$, this is equivalent to (\ref{BGF4}).
 
\smallskip     
For (C), fix $x$ and let $0 \le J < n$ be an index for which $|x-\psi_w(J)|_v$ is minimal.
For each $k \ne J$, we claim that $|x-\psi_w(k)|_v \ge |\psi_w(J)-\psi_w(k)|_v$.
Suppose to the contrary that $|x-\psi_w(k)|_v < |\psi_w(J)-\psi_w(k)|_v$.
By the ultrametric inequality, 
\begin{eqnarray*}
|x-\psi_w(J)|_v & = & \max(|x-\psi_w(k)|_v,|\psi_w(J)-\psi_w(k)|_v ) \\
                 & = & |\psi_w(J) - \psi_w(k)|_v \ > \ |x-\psi_w(k)|_v \ ,
\end{eqnarray*} 
contradicting our choice of $J$.  
Hence $|x-\psi_w(k)|_v \ge |\psi_w(J) - \psi_w(k)|_v$.  

It follows that 
$\ord_v(S_{n,w}(x)) \le \ord_v(x-\psi_w(J)) 
+ \ord_v(S_{n,w}^{\prime}(\psi_w(J)))$.  Thus the first assertion in (C) is a consequence of (B).  
The second is trivial, since if $x \notin D(0,1)$ then $|x-\psi_w(k)|_v = |x|_v$ for all $k$. 
\end{proof}

\begin{corollary} \label{SnOwMappingCor}
Let $F_w$ be a finite, separable extension of $K_v$ in $\CC_v$, and let $S_{n,w}(z)$
be the Stirling polynomial of degree $n$ for $\cO_w$.
\index{Stirling polynomial!for $\cO_w$}
 Given a radius $R$ satisfying
\begin{equation*} 
0 \ < \ R \ \le \ q_v^{-n/(e_w(q_v^{f_w}-1))} \cdot n^{-1/[F_w:K_v]} \ , 
\end{equation*} 
put $\rho_k = R/|S_{n,w}^{\prime}(\psi_w(k))|_v$ for $k = 0, \ldots, n-1$.

Then $S_{n,w}^{-1}\big(D(0,R)\big) = \bigcup_{k=0}^{n-1} D(\psi_w(k),\rho_k) \subset D(0,1)$,  
where the discs $D(\psi_w(k),\rho_k)$ are pairwise disjoint. 
For each $k$, $S_{n,w}(z)$ induces an $F_w$-rational scaled isometry from
\index{scaled isometry}
$D(\psi_w(k),\rho_k)$ onto $D(0,R)$, so that 
$S_{n,w}(F_w \cap D(\psi_w(k),\rho_k))= F_w \cap D(0,R)$ and 
\begin{equation} \label{BGFF1}
S_{n,w}^{-1}\big(F_w \cap D(0,R)\big) 
\ = \ \bigcup_{k=0}^{n-1} \big(F_w \cap D(\psi_w(k),\rho_k)\big) \ .
\end{equation} 
Moreover, if $R \in |F_w^{\times}|_v$ then  $\rho_k \in |F_w^{\times}|_v$ for each $k$.
\end{corollary}

\begin{proof}
\smallskip
Note that $\log_v(n^{-1/[F_w:K_v]}) = -\log_v(n)/(e_w f_w) = -\log_w(n)/e_w$.  
By the definition of $\rho_k$, our assumption on $R$, and Proposition \ref{BGProp1}(A), for each $k$ 
\begin{eqnarray}
-\log_v(\rho_k) & = & -\log_v(R) - \ord_v(S_{n,w}^{\prime}(\psi_w(k))) \notag \\
                & > & \frac{1}{e_w}\big(\frac{n}{q_w-1}  \  
         + \log_w(n)\big) - \frac{1}{e_w}\frac{n}{q_w-1} \ = \ \frac{1}{e_w} \log_w(n) \ . \label{FNZC1}
\end{eqnarray} 
On the other hand by (\ref{BGGG1}), for all $\ell \ne k$
\begin{equation*} 
\ord_v(\psi_w(k) - \psi_w(\ell)) \ = \ 
\frac{1}{e_w} \ord_w(\psi_w(k) - \psi_w(\ell)) \ < \ \frac{1}{e_w} \log_w(n) \ .
\end{equation*}    
Thus the discs $D(\psi_w(k), \rho_k)$ are pairwise disjoint and contained in $D(0,1)$.                  

Fix $k$ and expand $S_{n,w}(x)$ about $\psi_w(k)$ as
\begin{equation*}
S_{n,w}(z) \ = \ \sum_{\ell =1}^N b_{\ell}(z-\psi_w(k))^{\ell} 
\end{equation*}
where $b_1 = S_{n,w}^{\prime}(\psi_w(k))$.  
The definition of $\rho_k$ shows that $|b_1|_v \cdot \rho_k = R$.
By Proposition \ref{ScaledIsometryProp}, $S_{n,w}(z)$ induces an $F_w$-rational scaled isometry from 
\index{scaled isometry}
$D(\psi_w(k),\rho_k)$ onto $D(0,R)$, 
which takes $F_w \cap D(\psi_w(k), \rho_k)$ onto $F_w \cap D(0,R)$.  


Now let $k$ vary. As noted above, the discs
$D(\psi_w(k),\rho_k)$ are pairwise disjoint.  
Since $S_{n,w}(z)$ has degree $n$, 
for each $x \in D(0,R)$ the solutions to $S_{n,w}(z) = x$ 
in $\bigcup_{k=0}^{n-1} D(\psi_w(k),\rho_k)$
account for all the solutions in $\CC_v$. Hence 
\begin{equation*}
S_{n,w}^{-1}\big(D(0,R)\big) \ = \ \bigcup_{k=0}^{n-1} D(\psi_w(k),\rho_k) \ .
\end{equation*}
Similar considerations show that
\begin{equation*} 
S_{n,w}^{-1}(F_w \cap D(0,R)) \ = \ \bigcup_{k=0}^{n-1} \big(F_w \cap D(\psi_w(k),\rho_k)\big) \ .
\end{equation*}
Finally, note that if $R \in |F_w^{\times}|_v$, 
then $\rho_k = R/|S_{n,w}^{\prime}(\psi_w(k))|_v \in |F_w^{\times}|_v$ 
for each $k$.  
\end{proof}    

\smallskip

%% file: NewFSZChap4.tex
\chapter{Reductions} \label{Chap4}  


In this chapter we will formulate a simplified version of Theorem \ref{aT1}, 
which is the form of the theorem we will actually prove.  
After stating this theorem (Theorem \ref{aT1-B}), we will use it to deduce Theorem \ref{aT1}, 
Corollary \ref{FSZiv}, and the variants given in Chapter \ref{Variants}.

We begin with a definition.    

\begin{definition} \label{KvSimple} 
Let $v$ be a place of $K$.  A set $E_v \subset \cC_v(\CC_v)$
will be called {\em $K_v$-simple} if it is stable under $\Aut_c(\CC_v/K_v)$ and 
\index{$K_v$-simple!set|ii}  
is a union of finitely many pairwise disjoint, nonempty compact sets $E_{v,1}, \ldots, E_{v,D}$
such that:  

$(A)$ if $K_v \cong \CC$, then each $E_{v,\ell}$ is simply connected,\index{simply connected} 
has a piecewise smooth boundary, and is the closure of its $\cC_v(\CC)$-interior;
\index{boundary!piecewise smooth}\index{closure of $\cC_v(\CC)$ interior}\index{$K_v$-simple!$\CC$-simple|ii}\index{$\CC$-simple set|ii} 

$(B)$ if $K_v \cong \RR$, then each $E_{v,\ell}$ is either 

\quad $(1)$ a closed segment of positive length contained in $\cC_v(\RR)$, or 

\quad $(2)$ is disjoint from $\cC_v(\RR)$, and is simply connected,\index{simply connected} 
has a piecewise smooth boundary, and is the closure of its $\cC_v(\CC)$-interior;
\index{boundary!piecewise smooth}\index{closure of $\cC_v(\CC)$ interior}\index{$K_v$-simple!$\RR$-simple|ii}\index{$\RR$-simple set|ii}

$(C)$ if $K_v$ is nonarchimedean, then 

\quad $(1)$ there are finite separable extensions
$F_{w_1}, \ldots, F_{w_n}$ of $K_v$ contained in $\CC_v$, 
and pairwise disjoint isometrically parametrizable balls
\index{isometrically parametrizable ball}
$B(a_1,r_1), \ldots, B(a_D,r_D)$, 
such that $E_{v,\ell} = \cC_v(F_{w_\ell}) \cap B(a_\ell,r_\ell)$ for $\ell = 1, \ldots, D$. 

\quad $(2)$  The collection of balls $\{B(a_1,r_1), \ldots, B(a_D,r_D)\}$
is stable under $\Aut_c(\CC_v/K_v)$, and as $\sigma$ ranges over $\Aut_c(\CC_v/K_v)$,
each ball $B(a_\ell,r_\ell)$ has $[F_{w_\ell}:K_v]$ distinct conjugates.  
For each $\sigma$, if $\sigma(B(a_\ell,r_\ell)) = B(a_j,r_j)$, then 
$\sigma(F_{w_\ell}) = F_{w_j}$ and $\sigma(E_{v,\ell}) = E_{v,j}$. 
\end{definition} 
 
If $E_v$ is $K_v$-simple, we will call a decomposition 
\index{$K_v$-simple!set}  
$E_v = \bigcup_{\ell=1}^n E_{v,\ell}$ of the type in Definition \ref{KvSimple} 
a {$K_v$-simple decomposition}. 
\index{$K_v$-simple!decomposition|ii}  
If $v$ is archimedean, a $K_v$-simple set has a unique $K_v$-simple decomposition.
If $v$ is nonarchimedean, a $K_v$-simple decomposition can always be refined 
to another $K_v$-simple decomposition with smaller balls and more sets.

\begin{theorem}[FSZ with LRC for $K_v$-simple sets]
\label{aT1-B} \ 
\index{Fekete-Szeg\"o theorem with LRC!for $K_v$-simple sets|ii} 
\index{$K_v$-simple!set}  

Let $K$ be a global field, 
and let $\cC/K$ be a smooth, geometrically integral, projective curve.
Let $\fX = \{x_1, \ldots, x_m\} \subset \cC(\tK)$ 
be a finite set of points stable under $\Aut(\tK/K)$, and let
$\EE = \prod_v E_v \subset \prod_v \cC_v(\CC_v)$ be a 
$K_v$-rational adelic set compatible with $\fX$.\index{compatible with $\fX$}    
Let $S \subset \cM_K$ be a finite set of places $v \in \cM_K$ 
containing all archimedean $v$.      

Assume that $\gamma(\EE,\fX) > 1$, and that 

$(A)$ $E_v$ is $K_v$-simple for each $v \in S$,

$(B)$ $E_v$ is $\fX$-trivial for each $v \notin S$.  
\index{$\fX$-trivial}

\noindent{Then} 
there are infinitely many points $\alpha \in \cC(K^{\sep})$ 
such that for each $v \in \cM_K$, 
 the $\Aut(\tK/K)$-conjugates of $\alpha$ all belong to $E_v$. 
\end{theorem}

\medskip
We will now prove Theorems \ref{aT1}, \ref{aT1-A1}, \ref{aT1-A}, \ref{FSZi}, \ref{FSZii}, 
\ref{aT1-B1}, \ref{aT1-B2} and Corollary \ref{FSZiv}, assuming Theorem \ref{aT1-B}.  
For the convenience of the reader, we restate each theorem before proving it.

\newpage
\noindent{{\bf Theorem \text{\ref{aT1}}.}
(FSZ with LRC, producing conjugate points in $\EE$).}
\index{Fekete-Szeg\"o theorem with LRC!producing points in $\EE$|ii}

{\it
Let $K$ be a global field, 
and let $\cC/K$ be a smooth, geometrically integral, projective curve.
Let $\fX = \{x_1, \ldots, x_m\} \subset \cC(\tK)$ 
be a finite set of points stable under $\Aut(\tK/K)$, and let
$\EE = \prod_v E_v \subset \prod_v \cC_v(\CC_v)$ be an adelic set compatible with $\fX$.\index{compatible with $\fX$}  
Let $S \subset \cM_K$ be a finite set of places $v$, containing all archimedean $v$,
such that $E_v$ is $\fX$-trivial for each $v \notin S$.\index{$\fX$-trivial}

Assume that $\gamma(\EE,\fX) > 1$.  Assume also that $E_v$ has the following form, for each $v \in S$:   

$(A)$ If $v$ is archimedean and $K_v \cong \CC$, 
then $E_v$ is compact, and is a finite union of sets $E_{v,i}$, 
each of which is the closure of its $\cC_v(\CC)$-interior and has a 
piecewise smooth boundary;\index{boundary!piecewise smooth}\index{closure of $\cC_v(\CC)$ interior} 

$(B)$ If $v$ is archimedean and $K_v \cong \RR$, then $E_v$ is compact, stable under complex conjugation, 
and is a finite union of sets $E_{v,\ell}$, where each $E_{v,\ell}$ is either 

\quad $(1)$ the closure of its $\cC_v(\CC)$-interior and has a piecewise smooth boundary, or
\index{boundary!piecewise smooth}\index{closure of $\cC_v(\CC)$ interior}  

\quad $(2)$ is a compact, connected subset of $\cC_v(\RR)$; 

$(C)$ If $v$ is nonarchimedean, then $E_v$ is stable under $\Aut_c(\CC_v/K_v)$ 
and is a finite union of sets $E_{v,\ell}$, where each $E_{v,\ell}$ is either 

\quad $(1)$ an $\RL$-domain or a ball $B(a_\ell,r_\ell)$, or
\index{$\RL$-domain} 

\quad $(2)$ is compact and has the form $\cC_v(F_{w_\ell}) \cap B(a_\ell,r_\ell)$ 
for some finite separable extension $F_{w_\ell}/K_v$ in $\CC_v$, and some ball $B(a_\ell,r_\ell)$.  

Then there are infinitely many points $\alpha \in \cC(\tK^{\sep})$ such that for each $v \in \cM_K$, 
the $\Aut(\tK/K)$-conjugates of $\alpha$ all belong to $E_v$.  
}

\begin{proof}[Proof of Theorem \ref{aT1}, assuming Theorem \ref{aT1-B}]

{ \ }

The idea is to reduce Theorem \ref{aT1} to the case 
where the $E_v = \bigcup_{\ell=1}^{D_v} E_{v,\ell}$ are $K_v$-simple, 
\index{$K_v$-simple!set}  
and in particular where the $E_{v,\ell}$ are pairwise disjoint.  

Assume Theorem \ref{aT1-B}, 
and let $\EE = \prod_v E_v \subset \prod_v \cC_v(\CC_v)$
be an adelic set compatible with $\fX$\index{compatible with $\fX$}  
for which the hypotheses of Theorem \ref{aT1} hold.  
We will construct a new adelic set $\EE^{\prime} = \prod_v E_v^{\prime} \subset \EE$ 
such that the hypotheses of Theorem \ref{aT1-B} hold.
Let $S \subseteq \cM_K$ be a finite set of places 
containing all archimedean places and all nonarchimedean places where $E_v$
is not $\fX$-trivial. 
\index{$\fX$-trivial}

By hypothesis, $\gamma(\EE,\fX) > 1$.  
Let $\Gamma$ range over all symmetric matrices in $M_m(\RR)$.
By (\ref{MinMax1}) the value of $\Gamma$ 
as a matrix game\index{value of $\Gamma$ as a matrix game} 
is a continuous function of its entries, so there is an $\varepsilon > 0$ 
such that for any $\Gamma$ whose entries 
satisfy $|\Gamma(\EE,\fX)_{ij} - \Gamma_{ij}| < \varepsilon$ for all $i, j$,
\index{Green's matrix!global} 
we have $\val(\Gamma) < 0$.  Choose numbers $\varepsilon_v > 0$ for $v \in S$
such that $\sum_{v \in S} \varepsilon_v \log(q_v) < \varepsilon$.  
In constructing the sets $E_v^{\prime}$ for $v \in S$, in order 
to assure that $\gamma(\EE^{\prime},\fX) > 1$ it suffices to have 
\begin{equation} \label{XIneqs} 
\left\{ \begin{array}{ll} 
|G(x_i,x_j;E_v^{\prime}) - G(x_i,x_j;E_v)| \ < \ \varepsilon_v & 
                \text{for all $i \ne j$ \ ,} \\
|V_{x_i}(E_v^{\prime}) - V_{x_i}(E_v)| \ < \ \varepsilon_v &
                \text{for each $i$\ .} 
        \end{array} \right. 
\end{equation} 
\index{Green's function!properties of} 
\index{Robin constant!properties of} 

For each $v \notin S$, put $E_v^{\prime} = E_v$.  
Now suppose $v \in S$:    
\smallskip

{\bf Case 1.}  If $K_v \cong \CC$, then $\cC_v(\CC)$ is 
a Riemann surface.  Fix a triangulation $\cT$ of $\cC_v(\CC)$. 
\index{Riemann surface}\index{triangulation}  
Without loss we can assume that each edge 
of the triangulation is a smooth arc.\index{arc!smooth}
For each $\delta > 0$, let $\cT_{\delta}$ be a refinement of $\cT$  
such that that each edge of $\cT_{\delta}$ is a smooth arc and each triangle in $\cT_{\delta}$ 
\index{arc!smooth}
has diameter less than $\delta$ under the spherical distance $\|x,y\|_v$\index{spherical metric}. 

By assumption $E_v$ is the closure of its interior $E_v^0$,\index{closure of $\cC_v(\CC)$ interior} 
and its boundary\index{boundary!piecewise smooth} is a finite union of smooth arcs.  
In particular, each point of $\partial E_v$ is analytically accessible\index{analytically accessible} 
from $E_v^0$.  By Proposition \ref{IdentifyGreenProp}, 
this means that $\Gbar(z,x_i;E_v^0) = G(z,x_i;E_v)$ for each $\zeta \notin E_v$.

For each $\delta > 0$, 
let $X_\delta$ be the set of closed triangles in $\cT_{\delta}$ 
which are contained in $E_v^0$, and let $E_{v,\delta}$ 
be the union of the triangles in $X_\delta$.  As $\delta \rightarrow 0$, 
the sets $E_{v,\delta}$ exhaust $E_v^0$.  
Hence there is some $\delta = \delta_0$ such that   
\begin{equation*} 
\left\{ \begin{array}{ll} 
|G(x_i,x_j;E_{v,\delta_0}) - G(x_i,x_j;E_v)| \ < \ \varepsilon_v/2 & 
                \text{for all $i \ne j$ \ ,} \\
|V_{x_i}(E_{v,\delta_0}) - V_{x_i}(E_v)| \ < \ \varepsilon_v/2 &
                \text{for each $i$\ .} 
        \end{array} \right. 
\end{equation*} 
\index{Green's function!properties of} 
\index{Robin constant!properties of}

Let $U_v^{\prime} \subset E_{v,\delta_0}$ be the union of the interiors 
of the triangles in $X_{\delta_0}$.  

Each point of $E_{v,\delta_0}$ is analytically accessible\index{analytically accessible}  from $U_v^{\prime}$
in the sense of Definition \ref{AnalyticallyAccessible}.
By exhausting the interior of each triangle in $X_{\delta_0}$ by an 
increasing sequence of closed subtriangles  
and applying Proposition \ref{IdentifyGreenProp} again, 
we can find a compact set 
$E_v^{\prime} \subset U_v^{\prime}$ which 
is a finite union of closed triangles,
one contained in each connected component of $U_v^{\prime}$, such that
\begin{equation*} 
\left\{ \begin{array}{ll} 
|G(x_i,x_j;E_{v,\delta}) - G(x_i,x_j;E_v^{\prime})| \ < \ \varepsilon_v/2 & 
                \text{for all $i \ne j$ \ ,} \\
|V_{x_i}(E_{v,\delta}) - V_{x_i}(E_v^{\prime})| \ < \ \varepsilon_v/2 &
                \text{for each $i$\ .} 
        \end{array} \right. 
\end{equation*} 
\index{Green's function!properties of} 
\index{Robin constant!properties of} 
Thus, $E_v^{\prime} \subset E_v$, 
and $E_v^{\prime}$ satisfies (\ref{XIneqs}).  Moreover 
$E_v^{\prime}$ is compact and has finitely many connected components, each of which is simply connected,
\index{simply connected} 
has a piecewise smooth boundary\index{boundary!piecewise smooth} and is the closure of its interior.
\index{closure of $\cC_v(\CC)$ interior}  
Thus it is $\CC$-simple.\index{$\CC$-simple set}

\smallskip 
{\bf Case 2.}  If $K_v \cong \RR$, 
again choose a triangulation $\cT$ of $\cC_v(\CC)$.  
\index{triangulation} 
After making adjustments to $\cT$, if necessary, 
we can assume that $\cT$ is stable under 
complex conjugation, that each edge of $\cT$ is a smooth arc,
\index{arc!smooth} 
and that $\cC_v(\RR)$ is contained in the union of the edges of $\cT$.    
For each $\delta > 0$, let $\cT_{\delta}$ be a refinement of $\cT$  
with the properties in Case 1, such that each triangle in $\cT_{\delta}$ 
has diameter less than $\delta$ under the spherical distance $\|x,y\|_v$\index{spherical metric}. 

By assumption, $E_v$ is stable under complex conjugation, and 
a finite union of connected sets $E_{v,\ell}$ such that
each $E_{v,\ell}$ is either the closure of its $\cC_v(\CC)$-interior,\index{closure of $\cC_v(\CC)$ interior} 
or is contained in $\cC_v(\RR)$.  Without loss, we can assume that no $E_{v,\ell}$
is reduced to a point, since removing a finite set of points from $E_v$
does not change its capacity or Green's functions (Lemma \ref{CompactP2}).
\index{Green's function!properties of}  
\index{capacity} 

Let $E_{v,\CC}$ be the union of the sets $E_{v,\ell}$
which are closures of their $\cC_v(\CC)$-interiors, and let 
$E_{v,\RR} = E_v \backslash E_{v,\CC}$.  Then $E_{v,\CC}$ stable under 
complex conjugation and is the closure of its $\cC_v(\CC)$-interior,\index{closure of $\cC_v(\CC)$ interior} 
and $E_{v,\RR}$ is contained in $\cC_v(\RR)$.  Since $E_{v,\CC}$ is closed, 
no component of $E_{v,\RR}$ is reduced to a point. 

We now apply a modification of the argument from Case 1. 
Let $E_{v,\CC}^0$ be the $\cC_v(\CC)$-interior of $E_{v,\CC}$,
and let $E_{v,\CC}^{00} = E_{v,\CC}^0 \backslash \cC_v(\RR)$.
Let $E_{v,\RR}^0$ be the $\cC_v(\RR)$-interior of $E_{v,\RR}$, and   
put $E_v^1 = E_{v,\CC}^{00} \cup E_{v,\RR}^0$.  
Each point of $E_{v,\CC}$ is analytically accessible\index{analytically accessible}  from 
$E_{v,\CC}^{00}$ and each point of $E_{v,\RR}$
is analytically accessible from $E_{v,\RR}^0$,   
so each point of $E_v$ is analytically accessible from 
$E_v^1$.  By Proposition \ref{IdentifyGreenProp},
$\Gbar(z,x_i;E_v^1) = G(z,x_i;E_v)$ for each $\zeta \notin E_v$.  

For each $\delta > 0$, 
let $E_{v,\delta}$ be the union of the triangles in $\cT_{\delta}$  
contained in $E_{v,\CC}^{00}$, together with the edges of $\cT_{\delta}$
which are contained in $E_{v,\RR}^0$.  
Each $E_{v,\delta}$ is compact and stable under complex conjugation.  
Furthermore, each compact subset of $E_v^1$ is contained in some $E_{v,\delta}$, 
so there is a $\delta_0 > 0$ such that 
\begin{equation*} 
\left\{ \begin{array}{ll} 
|G(x_i,x_j;E_{v,\delta_0}) - G(x_i,x_j;E_v)| \ < \ \varepsilon_v/2 & 
                \text{for all $i \ne j$ \ ,} \\
|V_{x_i}(E_{v,\delta_0}) - V_{x_i}(E_v)| \ < \ \varepsilon_v/2 &
                \text{for each $i$\ .} 
        \end{array} \right. 
\end{equation*} 
\index{Green's function!properties of} 
\index{Robin constant!properties of}
  
Let $U_v^{\prime}$ be the union of the interiors of the triangles in 
$\cT_{\delta_0}$ which are contained in $E_{v,\delta_0}$, 
together the (real) interiors of the edges of $\cT_{\delta_0}$ 
which are contained in $E_{v,\RR}^0$.  Then $U_v^{\prime}$ 
is stable under complex conjugation, and is the disjoint union of
finitely many open triangles in $\cC_v(\CC)$ whose closures are 
disjoint from $\cC_v(\RR)$, 
together with finitely many open segments in $\cC_v(\RR)$.

Each point of $E_{v,\delta_0}$
is analytically accessible\index{analytically accessible}  from $U_v^{\prime}$.  
By exhausting the interior of each open triangle in $U_v^{\prime}$ by an 
increasing sequence of closed subtriangles and each open segment by an increasing
sequence of closed subintervals, 
and applying Proposition \ref{IdentifyGreenProp} again, 
we can find a compact set 
$E_v^{\prime} \subset U_v^{\prime}$ which 
is a finite disjoint union of closed triangles 
and closed subintervals of $\cC_v(\RR)$,
one contained in each component of $U_v^{\prime}$, such that
\begin{equation*} 
\left\{ \begin{array}{ll} 
|G(x_i,x_j;E_{v,\delta}) - G(x_i,x_j;E_v^{\prime})| \ < \ \varepsilon_v/2 & 
                \text{for all $i \ne j$ \ ,} \\
|V_{x_i}(E_{v,\delta}) - V_{x_i}(E_v^{\prime})| \ < \ \varepsilon_v/2 &
                \text{for each $i$\ .} 
        \end{array} \right. 
\end{equation*} 
\index{Green's function!properties of} 
\index{Robin constant!properties of}
By choosing the closed subtriangles appropriately, 
we can also arrange that $E_v^{\prime}$ is stable under complex conjugation.

Thus, $E_v^{\prime}$ is compact, $\RR$-simple,\index{$\RR$-simple set} and contained in $E_v$.  
By construction it satisfies (\ref{XIneqs}).  

\smallskip
{\bf Case 3.}  Suppose $K_v$ is nonarchimedean.  By assumption, $E_v$
is a finite union of $\RL$-domains and compact sets of the form
\index{$\RL$-domain} 
$B(a_\ell,r_\ell) \cap \cC_v(F_{w_\ell})$, where each $B(a_\ell,r_\ell)$ 
is isometrically parametrizable and each $F_{w_\ell}$ 
\index{isometrically parametrizable ball}
is a finite separable extension of $K_v$ contained in $\CC_v$.   
Since a finite union of $\RL$-domains is an $\RL$-domain 
\index{$\RL$-domain} 
(\cite{RR1}, Theorem 4.2.15) we can assume that there is at most one
$\RL$ domain in the decomposition. 

Our first goal is to reduce to the case where there are no $\RL$-domains
\index{$\RL$-domain} 
in the decomposition of $E_v$.  
Suppose to the contrary that there is an $\RL$-domain $U_v$. 
\index{$\RL$-domain}  
Let $E_v^{(1)}$ be the union of the compact sets $B(a_\ell,r_\ell) \cap \cC_v(F_{w_\ell})$
in the decomposition of $E_v$.
Both $U_v$ and $E_v^{(1)}$ are stable under $\Aut_c(\CC_v/K_v)$.  
In (\cite{RR1}, Theorem 4.3.11) it is shown that the union of an $\RL$-domain and a
\index{$\RL$-domain}  
compact set is algebraically capacitable.  In fact,
\index{algebraically capacitable}
the proof of that theorem shows there is a compact set $E_v^{(2)} \subset U_v$, 
which itself is a finite union of compact sets of the
form $B(a_\ell^{\prime},r_\ell^{\prime}) \cap \cC_v(F_{w_\ell}^{\prime})$, 
with $B(a_\ell^{\prime},r_\ell^{\prime})$ isometrically parametrizable 
\index{isometrically parametrizable ball}
and $F_{w_\ell}^{\prime}/K_v$ finite, such that  
\begin{equation*} 
\left\{ \begin{array}{ll} 
|G(x_i,x_j;E_v^{(1)} \cup E_v^{(2)} - G(x_i,x_j;E_v)| \ < \ \varepsilon_v/2 & 
                \text{for all $i \ne j$ \ ,} \\
|V_{x_i}(E_v^{(1)} \cup E_v^{(2)} - V_{x_i}(E_v)| \ < \ \varepsilon_v/2 &
                \text{for each $i$\ .} 
        \end{array} \right. 
\end{equation*} 
\index{Green's function!properties of} 
\index{Robin constant!properties of}
The centers of the balls $B(a_\ell^{\prime},r_\ell^{\prime})$ can be required
to belong to $\cC_v(\tK_v^{\sep})$, since $\tK_v^{\sep}$ is dense in $\CC_v$,
and the extensions $F_{w_\ell}$ can be required to be separable extensions of $K_v$,  
since all that is needed for the proof of (\cite{RR1}, Theorem 4.3.11) is that
the residue degree or the ramification index of $F_{w_\ell}/K_v$ can be taken 
arbitrarily large.  
The set $E_v^{(2)}$ need not be stable under $\Aut_c(\CC_v/K_v)$,
but by its form it has only finitely many conjugates, and each of these is 
also contained in $U_v$.  
By replacing $E_v^{(2)}$ with the union of its conjugates, 
and using the monotonicity of the Green's functions (Lemma \ref{CompactP1}),
\index{Green's function!monotonic} 
we can arrange that $E_v^{(1)} \cup E_v^{(2)}$ 
is stable under $\Aut_c(\CC_v/K_v)$.  

Thus we can assume that $E_v$ has no $\RL$-domains in its decomposition.
\index{$\RL$-domain} 
Write $E_v = \bigcup_{\ell=1}^{D_v} B(a_\ell,r_\ell) \cap \cC_v(F_{w_\ell})$, and 
fix a number  $r \in |K_v^{\times}|_v$ with $0 < r \le \min_\ell(r_\ell)$.  After   
replacing $B(a_\ell,r_\ell) \cap \cC_v(F_{w_\ell})$ by finitely many sets 
$B(a_{\ell j},r) \cap \cC_v(F_{w_\ell})$ for each $\ell$, 
we can also assume that all the $r_\ell$ are equal to $r$, and that $r$ belongs 
to the value group of $K_v^{\times}$.  

\smallskip  
Our next goal is to arrange that balls $B(a_\ell,r)$ are disjoint.  
Consider the sets $E_{v,\ell} = B(a_\ell,r) \cap \cC_v(F_{w_\ell})$.  
Without loss, we can assume that none of the $E_{v,\ell}$  
is properly contained in another. 
However, it is possible that the same ball $B(a_\ell,r)$ occurs 
with several different fields $F_{w_{\ell_j}}$;  given such a ball, let 
$\ell_1 = \ell, \ell_2, \ldots, \ell_t$ be the indices 
for which $B(a_{\ell_j},r) = B(a_\ell,r)$.  For each $j \ge 2$,  
the field $F_{w_\ell} \cap F_{w_{\ell_j}}$ is a proper subfield of $F_{w_\ell}$, 
and by a simple argument involving $K_v$-vector spaces,
there is an element $u_\ell \in F_{w_\ell}$ 
which does not belong to any $F_{w_\ell} \cap F_{w_{\ell_j}}$.  
Let $0 \ne \pi_v \in K_v$ be such that $|\pi_v|_v < 1$.  
After replacing $u_i$ by $1 + \pi_v^N u_\ell$ for sufficiently large $N$, 
we can assume that $u_\ell \in \cO_{w_\ell}^{\times}$.   
  
For each $\ell$, put 
\begin{equation*}
X_\ell  \ = \ E_{v,\ell} \backslash \Big( \bigcup_{j \ne \ell} E_{v,j} \Big) \ ,
\end{equation*} 
and let $X = \bigcup_{\ell=1}^{D_v} X_\ell$.  
Note that $X$ is stable under $\Aut_c(\CC_v/K_v)$. 

We claim that each point of $E_v$ is analytically accessible\index{analytically accessible}  from $X$ 
in the sense of Definition \ref{AnalyticallyAccessible}.
Clearly each $b \in X$ is analytically accessible\index{analytically accessible}  from $X$, since if $b \in X_\ell$
then there is a ball $B(b,s)$ for which $B(b,s) \cap \cC_v(F_{w_\ell}) \subseteq X_\ell$.  
Suppose that $b \in E_{v,j} \cap E_{v,\ell}$ for some $j \ne \ell$.  Then
 $b  \in B(a_j,r) \cap \cC_v(F_{w_j}) \cap \cC_v(F_{w_\ell})$, 
so $b \in \cC_v(F_{w_j} \cap F_{w_\ell})$.  Let $F_b$ be the field $K_v(b)$.  
By Theorem \ref{IsoParamThm}
there is an $F_b$-rational isometric parametrization
\index{isometric parametrization}
$f_b : D(0,r) \rightarrow B(a_\ell,r)$ with $f_b(0) = b$.  Consider the 
image of $D(0,r) \cap (\cO_v \backslash \{0\})$ 
under the isometric parametrization $\tf_b(z) := f_b(u_\ell z)$: we have 
\index{isometric parametrization}
\begin{equation*} 
\tf_b(D(0,r) \cap (\cO_v \backslash \{0\})) 
        \ = \ f_b(u_\ell \cdot (D(0,r) \cap (\cO_v \backslash \{0\}))) \ .
\end{equation*}   
Since $f_b$ is  $F_b$-rational, 
for each complete field $H$ with $F_b \subseteq H \subseteq \CC_v$, 
$f_b$ induces a $1-1$ correspondence between points of $D(0,r) \cap H$ and $B(a_i,r) \cap \cC_v(H)$.  
By our choice of $u_\ell$, we have 
$u_\ell \cdot (\cO_v \backslash \{0\}) \subset 
       \cO_{w_\ell} \backslash (\bigcup_{j=1}^t  F_{w_{\ell_j}})$. 
It follows that $\tf_b(D(0,r) \cap (\cO_v \backslash \{0\})) \subset X_\ell$.   
This establishes our claim.

By Proposition \ref{IdentifyGreenProp}, $G(z,\zeta;E_v) = \Gbar(z,\zeta;X)$
for each $\zeta \notin E_v$.  Hence there is a compact set $Y \subset X$ such that
\begin{equation*} 
\left\{ \begin{array}{ll} 
|G(x_i,x_j;Y) - G(x_i,x_j;E_v)| \ < \ \varepsilon_v/2 & 
                \text{for all $i \ne j$ \ ,} \\
|V_{x_i}(Y) - V_{x_i}(E_v)| \ < \ \varepsilon_v/2 &
                \text{for each $i$\ .} 
        \end{array} \right. 
\end{equation*} 
\index{Green's function!properties of} 
\index{Robin constant!properties of}
The set $Y$ thus constructed 
may not be stable under $\Aut_c(\CC_v/K_v)$, but it only has finitely
many conjugates and each of them is contained in $X$. 
By replacing $Y$ with the union of its conjugates, and using the monotonicity
of Green's functions, we can assume $Y$ is stable under $\Aut_c(\CC_v/K_v)$. 
\index{Green's function!monotonic} 

For each $\ell = 1, \ldots, D_v$, put $Y_\ell = Y \cap B(a_\ell,r) \cap \cC_v(F_{w_\ell})$.  
Then $Y_\ell$ is compact and $Y_\ell \subset X_i$.  Since the $Y_\ell$ are compact
and pairwise disjoint, there is a number $0 < R \in |K_v^{\times}|_v$ such that 
\begin{equation*}
\min_{\ell \ne j} \min_{z \in Y_\ell, w \in Y_j} \|z,w\|_v \ > \ R \ .
\end{equation*} 
Cover  $Y$ with finitely many balls $B(b_j,R)$, $j = 1, \ldots, N$,
for points $b_j \in Y$.  

We can now construct $E_v^{\prime}$. 
For each $j = 1, \ldots, N$, if $b_j \in Y_\ell$, set $F_{w_\ell}^{\prime} = F_{w_\ell}$, 
and put
\begin{equation*}   
E_v^{\prime} \ = \ \bigcup_{j=1}^N B(b_j,R) \cap \cC_v(F_{w_j}^{\prime}) \ . 
\end{equation*}
Since $\Aut_c(\CC_v/K_v)$ preserves the spherical distance\index{spherical metric!Galois equivariance of}, 
and stabilizes both $X$ and $Y$, it follows that $E_v^{\prime}$ is stable under $\Aut_c(\CC_v/K_v)$, 
By construction it is $K_v$-simple, contained in $E_v$, 
\index{$K_v$-simple!set}  
and satisfies (\ref{XIneqs}).
\end{proof}

\noindent{{\bf Corollary \text{\ref{FSZiv}}.} (Fekete-Szeg\"o with LRC for Incomplete Skolem Problems)} 
\index{Incomplete Skolem Problems}
\index{Fekete-Szeg\"o theorem with LRC!for Incomplete Skolem Problems|ii} 
{\it Let $K$ be a global field, 
and let $\cA/K$ be a geometrically integral $($possibly singular$)$ affine curve, 
embedded in $\AA^N$ for some $N$. 
Let $z_1, \ldots, z_N$ be the coordinates on $\AA^N;$  
given a place $v$ of $K$ and a point $P \in \AA^N(\CC_v)$, 
write $\|P\|_v = \max(|z_1(P)|_v, \ldots, |z_N(P)|_v)$.  

Fix a place $v_0$ of $K$, and let $S \subset \cM_K \backslash \{v_0\}$ be a finite set of places containing 
all archimedean $v \ne v_0$.     
For each $v \in S$, let a nonempty set $E_v \subset \cA_v(\CC_v)$
satisfying condition $(A)$, $(B)$ or $(C)$ of Theorem $\ref{aT1}$ be given,    
and put $\EE_S = \prod_{v \in S} E_v$.  
Assume that for each $v \in \cM_K \backslash (S \cup \{v_0\})$ there is a 
point $P \in \cA(\CC_v)$ with $\|P\|_v \le 1$. 
Then there is a constant $C = C(\cA,\EE_S,v_0)$
such that there are infinitely many points $\alpha \in \cA(\tK^{\sep})$ for which  

$(1)$ for each $v \in S$, all the conjugates of $\alpha$ in $\cA_v(\CC_v)$
belong to $E_v;$   

$(2)$ for each $v \in \cM_K \backslash (S \cup \{v_0\})$, 
all the conjugates of $\alpha$ in $\cA_v(\CC_v)$ satisfy $\|\sigma(\alpha)\|_{v_0} \le 1;$   

$(3)$ for $v = v_0$, all the conjugates of $\alpha$ in $\cA_{v_0}(\CC_{v_0})$ satisfy $\|\sigma(\alpha)\|_{v_0} \le C$.      
}

\smallskip
\begin{proof}[Proof of Corollary \ref{FSZiv}, using Theorem \ref{aT1}]

{ \ }

Let $\cAbar$ be the projective closure of $\cA$, and let $\cC/K$ be a desingularization of $\cAbar$.
Then $\cC$ is a smooth, geometrically integral, projective curve birational to $\cAbar$.    
Let $\pi : \cC \rightarrow \cAbar$ be the natural morphism;  
it is an isomorphism away from the finitely many preimages of the singular points.
For each $v \in \cM_K$, $\pi$ induces a map $\pi_v : \cC_v(\CC_v) \rightarrow \cAbar(\CC_v)$.  

For each $v \in \cM_K$, and each $0 < R \in \RR$, put $B_v(R) = \{x \in \cA(\CC_v) : \|x\|_v \le R \}$.
Let $\fY = \cAbar(\tK) \backslash \cA(\tK)$ be the set of points at infinity for $\cA$, 
and put $\fX = \pi^{-1}(\fY)$.  Then $\fX$ is finite and stable under $\Aut(\tK/K)$.  
For each $v \in \cM_K$, define a set $\tE_v \subset \cC_v(\CC_v)$ as follows.
If $v \in \cS$, put $\tE_v = \pi_v^{-1}(E_v)$.  If $v \in \cM_K \backslash (S \cup \{v_0\})$, 
(so in particular $v$ is nonarchimedean), put $\tE_v = \pi_v^{-1}(B_v(1))$.  By assumption $B_v(1)$ is nonempty,
so $\tE_v$ is nonempty and open.  Indeed, it is an $\RL$-domain,\index{$\RL$-domain} 
since when we regard each coordinate function $z_i(x)$ as an element of $K_v(\cC)$, 
it is the intersection of the finitely many
$\RL$-domains $\{x \in \cC_v(\CC_v) : |z_i(x)|_v \le 1\}$ (see \cite{RR1}, Theorem 4.2.15).\index{good reduction} 
For all but finitely many $v$, the curve $\cC_v$ and the functions $z_i(x)$ have good reduction $\pmod{v}$, 
and the points in $\fX$ specialize to distinct points $\pmod{v}$.  For such $v$, $\tE_v$ is the $\fX$-trivial set.      
Finally, at the place $v_0$, take a number $R > 0$ and put $\tE_{v_0} = \tE_{v_0}(R) = \pi_{v_0}^{-1}(B_{v_0}(R))$.
Clearly each $\tE_v$ is stable under $\Aut_c(\CC_v/K_v)$ and satisfies the conditions in Theorem \ref{aT1}.     
  
Let  
\begin{equation*}
\tEE(R) \ = \ \tE_{v_0}(R) \times \prod_{v \ne v_0} \tE_v \ .
\end{equation*} 
Then $\tEE(R)$ is an adelic set compatible with $\fX$.  
We claim that for all sufficiently large $R$, we have $\gamma(\tEE(R),\fX) > 1$.  
Indeed, fix a spherical metric\index{spherical metric} on $\cC_{v_0}(\CC_{v_0})$, 
and let $\varepsilon > 0$ be small enough that the balls $B_{v_0}(x_i,\varepsilon)$ for $x_i \in \fX$
are pairwise disjoint.  If $R$ is big enough then $\cC_{v_0}(\CC_{v_0}) \backslash E_{v_0}(R)$ will 
be contained in $\bigcup_{x_i \in \fX} B_{v_0}(x_i,\varepsilon)$,
and in that case the local Green's matrix $\Gamma_{v_0}(E_{v_0}(R),\fX)$ will be diagonal.
By letting $R \rightarrow \infty$ we can make the diagonal entries arbitrarily large and negative.   
If $R$ is sufficiently large the matrix $\Gamma(\tEE(R),\fX)$ will be negative definite.
\index{matrix!negative definite} 
Taking $C = R$ for such an $R$ and pushing forward the points produced by Theorem \ref{aT1} yields the Corollary.    
\end{proof} 

\noindent{{\bf Theorem \text{\ref{aT1-A1}}.} (FSZ with LRC for Quasi-neighborhoods)} 
\index{quasi-neighborhood}
\index{Fekete-Szeg\"o theorem with LRC!for quasi-neighborhoods|ii} 
{\it Let $K$ be a global field, 
and let $\cC/K$ be a smooth, connected, projective curve.
Let $\fX = \{x_1, \ldots, x_m\} \subset \cC(\tK)$ 
be a finite set of points stable under $\Aut(\tK/K)$, and let 
$\EE = \prod_v E_v \subset \prod_v \cC_v(\CC_v)$
be an adelic set compatible with $\fX$,\index{compatible with $\fX$} 
such that each $E_v$ is stable under $\Aut_c(\CC_v/K_v)$. 
 
Suppose $\gammabar(\EE,\fX) > 1$.
Then for any $K$-rational separable quasi-neighborhood $\UU$ of $\EE$, there 
\index{quasi-neighborhood!separable}
are infinitely many points $\alpha \in \cC(\tK^{\sep})$ 
such that for each $v \in \cM_K$,
the $\Aut(\tK/K)$-conjugates of $\alpha$ all belong to $U_v$.  
}
 
\begin{proof}[Proof of Theorem \ref{aT1-A1}, using Theorem \ref{aT1}] 

{ \ } 

The idea is to adjust the sets $E_v$ within their 
separable quasi-neighborhoods $U_v$ (see Definition \ref{QuasiNeighborhood}), 
\index{quasi-neighborhood!separable}
and reduce Theorem \ref{aT1-A1} to the case where the 
$E_v$ satisfy the conditions of Theorem \ref{aT1}.  
In particular, in the nonarchimedean case,
we may need to pass from infinite algebraic extensions $F_{w_\ell}/K_v$ to ones 
of finite degree.  

Assume Theorem \ref{aT1} is true.  Let $\EE = \prod_v E_v$
be a $K$-rational adelic set compatible with $\fX$,\index{compatible with $\fX$} 
and let $\UU = \prod_v U_v$ be a $K$-rational
separable quasi-neighborhood of $\EE$, for which the hypotheses of Theorem \ref{aT1-A1} hold.  
\index{quasi-neighborhood!separable}
We will construct a new adelic set $\EE^{\prime} = \prod_v E_v^{\prime} \subset \UU$ 
such that the hypotheses of Theorem \ref{aT1} hold for $\EE^{\prime}$.
Let $S \subseteq \cM_K$ be a finite set of places 
containing all archimedean places and all nonarchimedean places where $E_v$
is not $\fX$-trivial.  
\index{$\fX$-trivial}
For each $v \notin S$, put $E_v^{\prime} = E_v$.  

By hypothesis, $\gammabar(\EE,\fX) > 1$.  
Let $\Gamma$ range over all symmetric matrices in $M_m(\RR)$.
By (\ref{MinMax1}) there is an $\varepsilon > 0$ 
such that for any $\Gamma$ whose entries 
satisfy $\Gamma_{ij} < \Gamma(\EE,\fX)_{ij} + \varepsilon$ for all $i, j$,
\index{Green's matrix!global}
we have $\val(\Gamma) < 0$.  Choose numbers $\varepsilon_v > 0$ for $v \in S$
such that $\sum_{v \in S} \varepsilon_v \log(q_v) < \varepsilon$.  
In constructing the sets $E_v^{\prime}$ for $v \in S$,  
to assure that $\gamma(\EE^{\prime},\fX) > 1$ it suffices to arrange that 
\begin{equation} \label{XIneqs1} 
\left\{ \begin{array}{ll} 
G(x_i,x_j;E_v^{\prime})  \ < \  \Gbar(x_i,x_j;E_v) + \varepsilon_v & 
                \text{for all $i \ne j$ \ ,} \\
V_{x_i}(E_v^{\prime}) \ < \ \Vbar_{x_i}(E_v) + \varepsilon_v &
                \text{for each $i$\ .} 
        \end{array} \right. 
\end{equation} 
\index{Green's function!properties of} 
\index{Robin constant!properties of}

We will now construct the sets $E_v^{\prime}$ for $v \in S$. 

\smallskip
{\bf Case 1.}  If $K_v \cong \CC$, then $U_v \subset \cC_v(\CC)$ is 
open.  By Proposition \ref{GreenProp}(4), there is a compact set $H_v \subset E_v$  such that 
\begin{equation} \label{XIneqs1C} 
\left\{ \begin{array}{ll} 
|G(x_i,x_j;H_v) - \Gbar(x_i,x_j;E_v)| \ < \ \varepsilon_v & 
                \text{for all $i \ne j$ \ ,} \\
|V_{x_i}(H) - \Vbar_{x_i}(E_v)| \ < \ \varepsilon_v &
                \text{for each $i$\ .} 
        \end{array} \right. 
\end{equation} 
\index{Green's function!properties of} 
\index{Robin constant!properties of}
For each $a \in H_v$, there is an $r_a > 0$ such that the closed ball $B(a,r_a)$ is contained 
in $U_v$.  The corresponding open balls $B(a,r_a)^-$ cover $H_v$.  By compactness, 
finitely many of these balls, say $B(a_1,r_{a_1})^-, \ldots, B(a_n,r_{a_n})^-$
also cover $H_v$.  Let 
\begin{equation*} 
E_v^{\prime} \ = \ B(a_1,r_{a_1}) \cup \cdots \cup B(a_n,r_{a_n}) \ \subset \ U_v \ .
\end{equation*} 
Each $B(a_i,r_{a_i})$ is compact, has a smooth boundary,\index{boundary}\index{closure of $\cC_v(\CC)$ interior} 
and is the closure of its interior,
so $E_v^{\prime}$ satisfies the conditions of Theorem \ref{aT1}. 
Since $H_v \subset E_v^{\prime}$, the monotonicity of Green's functions (Lemma \ref{CompactP1})
\index{Green's function!monotonic}
shows that \ref{XIneqs1} holds.  

\smallskip 
{\bf Case 2.}  If $K_v \cong \RR$, then the quasi-neighborhood $U_v$ of $E_v$ 
\index{quasi-neighborhood}
is the union of a set $U_{v,0}$ open in $\cC_v(\CC)$ and a set $U_{v,1}$  open in $\cC_v(\RR)$.
Since $U_v$ is stable under complex conjugation, so is $U_{v,0}$.  
By Proposition \ref{GreenProp}(4), there is a compact set $H_v \subset E_v$  for which 
\begin{equation} \label{XIneqs1R} 
\left\{ \begin{array}{ll} 
|G(x_i,x_j;H_v) - \Gbar(x_i,x_j;E_v)| \ < \ \varepsilon_v/2 & 
                \text{for all $i \ne j$ \ ,} \\
|V_{x_i}(H_v) - \Vbar_{x_i}(E_v)| \ < \ \varepsilon_v/2 &
                \text{for each $i$\ .} 
        \end{array} \right. 
\end{equation}
\index{Green's function!properties of} 
\index{Robin constant!properties of}  

We will construct the set $E_v^{\prime}$ in two steps.
First, put $H_{v,1} = H_v \backslash U_{v,0}$;  
it is compact and contained in $U_{v,1}$. 
Next, for each $r > 0$, put 
$H_{v,0}(r) = \{x \in H_v : \|x,z\| \ge r \text{ for all $z \in H_{v,1}$} \}$,
and let 
\begin{equation*}
H_v(r) \ = \ H_{v,0}(r) \cup H_{v,1} \ .
\end{equation*}   
Then  $H_{v,0}(r)$ is contained in $U_{v,0}$, 
and $H_{v,0}(r)$ and $H_v(r)$ are compact. 

The sets $H_v(r)$ increase as $r \rightarrow 0$, 
and they form an exhaustion of $H_v$.   
Choose a sequence $r_1 > r_2 > \cdots > 0$ with 
$\lim_{m \rightarrow \infty} r_m = 0$.  By Proposition \ref{InnerSeqProp}, 
if we take $M$ sufficiently large and put $X_v = H_v(r_M)$, then 
\begin{equation} \label{XIneqs1RB} 
\left\{ \begin{array}{ll} 
|G(x_i,x_j;X_v) - G(x_i,x_j;H_v)| \ < \ \varepsilon_v/2 & 
                \text{for all $i \ne j$ \ ,} \\
|V_{x_i}(X_v) - V_{x_i}(H_v)| \ < \ \varepsilon_v/2 &
                \text{for each $i$\ .} 
        \end{array} \right. 
\end{equation}
\index{Green's function!properties of} 
\index{Robin constant!properties of}

For each $a \in H_{v,0}(r_M)$, there is an open ball $B(a,r)^-$ 
whose closure $B(a,r)$ is contained in $U_{v,0}$.  
Finitely many of these balls, 
say $B(a_1,r_1)^-, \ldots, B(a_{m_0},r_{m_0})^-$ cover $H_{v,0}(r_M)$.  
For each $a \in T_{v,1}$, 
there is an open interval $I_a \subset \cC_v(\RR)$ containing $x$, 
whose closure $\Ibar_a$ is contained in $U_{v,1}$.  
Finitely many of these intervals, 
say $I_{b_1}, \ldots, I_{b_{m_1}}$, cover $H_{v,1}$. 
Put  
\begin{equation*} 
E_v^{\prime} \ = \ \Big( \bigcup_{i=1}^{m_0} B(a_i,r_i) \Big) 
\cup \Big( \bigcup_{i=1}^{m_0} B(\overline{a_i},r_i) \Big) 
\cup \big( \Ibar_{b_1} \cup \ldots \cup \Ibar_{b_{m_1}} \big)
\ .
\end{equation*} 
Then $E_v^{\prime}$ is compact, stable under complex conjugation, 
and contained in $U_v$.
By construction, it satisfies the conditions of Theorem \ref{aT1}. 
Hence (\ref{XIneqs1R}) and (\ref{XIneqs1RB}), 
together with the monotonicity of Green's functions,
\index{Green's function!monotonic}
show that (\ref{XIneqs1}) holds.

\smallskip
{\bf Case 3.}  If $K_v$ is nonarchimedean, then the quasi-neighborhood $U_v$
\index{quasi-neighborhood}
has the form  
\begin{equation} \label{QuasiRep} 
U_v \ = \ U_{v,0} \cup \big(U_{v,1} \cap \cC_v(F_{w_1})\big) \cup \cdots 
\cup \big(U_{v,D} \cap \cC_v(F_{w_D})\big) 
\end{equation}
where $U_{v,0}, \ldots, U_{v,D}$ are open in $\cC_v(\CC_v)$ and each $F_{w_\ell}$
is a separable algebraic extension of $K_v$ (possibly of infinite degree) 
contained in $\CC_v$.  
By hypothesis, $U_v$ is $\Aut_c(\CC_v/K_v)$-stable and contains $E_v$.

In what follows, it will convenient to take $F_{w_0} = \CC_v$, 
though we will be careful to make a distinction 
between $F_{w_0}$ and the $F_{w_\ell}$ for $\ell \ge 1$,
which are separably algebraic over $K_v$. 
Define a {\em representation of $U_v$}\index{representation of $U_v$|ii} to be a collection of pairs $(U_{v,\ell},F_{w_\ell})$ 
such that (\ref{QuasiRep}) holds.  
We allow the possibility that some of the $U_{v,i} \cap \cC_v(F_{w_i})$ may be empty.
We will say that a set $U_{v,i}$ or a field $F_{w_i}$ {\em occurs}, 
if it is a component of one of the pairs.  

There are many representation\index{representation of $U_v$} for $U_v$.   
We begin by adjusting the given representation to make it easier to work with. 
First, for each finite separable extension $F_{w_\ell}/K_v$ which occurs, 
there are finitely many intermediate fields $K_v \subseteq G_u \subseteq F_{w_\ell}$.  
Adjoin each of the pairs $(U_{v,\ell},G_u)$ to the representation.\index{representation of $U_v$}  
Second, a given field may occur in several pairs.  
Replace those pairs with a single pair 
whose first component is the union of the sets in the original pairs.
In this way, we can assume that the $F_{w_\ell}$ are distinct, 
and that whenever a finite separable extension $F_{w_\ell}/K_v$ occurs, 
so do all of its subextensions.  In addition, 
if two finite separable extensions $F_{w_j} \subseteq F_{w_\ell}$ occur, then 
\begin{equation} \label{FDifSet1} 
U_{v,j} \cap \cC_v(F_{w_j}) \ \subseteq \ U_{v,\ell} \cap \cC_v(F_{w_\ell}) \ .
\end{equation}  

We will now construct a new representation\index{representation of $U_v$} 
$\{(W_{v,j},F_{w_j})\}_{0 \le j \le n}$ with the same fields $F_{w_j}$, 
but giving a different decomposition of $U_v$.
After reordering the pairs $(U_{v,\ell},F_{w_\ell})$, 
we can assume that $F_{w_\ell}/K_v$ is of infinite degree for $\ell = 0, \ldots, D_0$ 
and that $F_{w_\ell}/K_v$ is finite for $\ell = D_0+1, \ldots, D$. 
After reordering them further, we can also assume that for each $\ell$, 
$F_{w_\ell}$ is maximal among the $F_{w_j}$ with $j \ge \ell$ 
(under the partial order given by containment).   

For each $j \le D_0$, 
let $W_{v,j}$ be the union of all the isometrically parametrizable balls $B(a,r)$
\index{isometrically parametrizable ball}
such that  $a \in U_{v,j} \cap F_{w_j}$ and $B(a,r) \subset U_{v,j}$.  
Since isometrically parametrizable balls are 
cofinal in the neighborhoods of a given point, 
we have $W_{v,j} \cap \cC_v(F_{w_j}) = U_{v,j} \cap \cC_v(F_{w_j})$. 
For each $j \ge D_0+1$, note that if $k > j$, then $\cC_v(F_{w_k})$ is compact, 
hence closed in $\cC_v(\CC_v)$;  put 
\begin{equation} \label{FW2NN1}
W_{v,j} \ = \ U_{v,j} \backslash \ \Big(\bigcup_{k > j} \cC_v(F_{w_k}) \Big) \ .
\end{equation} 
This means that the sets  $W_{v,j} \cap \cC_v(F_{v,j})$ 
corresponding to finite extensions $F_{w_j}/K_v$  are pairwise disjoint. 
However, by the construction of the ordering, for each $F_{w_k}$ with 
$K_v \subseteq F_{w_k} \subsetneq F_{w_j}$ we must have $k > j$.   
Hence by (\ref{FDifSet1}), we still have 
\begin{equation} \label{NewQuasiRep} 
U_v \ = \ \bigcup_{j=0}^D \big(W_{v,j} \cap \cC_v(F_{w_j})\big) \ .
\end{equation} 

\smallskip
By Proposition \ref{GreenProp}(4), there is a compact set $H_v \subseteq E_v$  for which 
\begin{equation} \label{XIneqs1N} 
\left\{ \begin{array}{ll} 
|G(x_i,x_j;H_v) - \Gbar(x_i,x_j;E_v)| \ < \ \varepsilon_v/3 & 
                \text{for all $i \ne j$ \ ,} \\
|V_{x_i}(H_v) - \Vbar_{x_i}(E_v)| \ < \ \varepsilon_v/3 &
                \text{for each $i$\ .} 
        \end{array} \right. 
\end{equation} 
\index{Green's function!properties of} 
\index{Robin constant!properties of}
We aim to pass from $H_v$ (whose structure is completely unknown) 
to an $\Aut_c(\CC_v/K_v)$-stable set $E_v^{\prime}$ contained in 
$U_v \cap \cC_v(F_w)$ for some finite galois extension $F_w/K_v$, 
which satisfies (\ref{XIneqs1}) and has the form required by Theorem \ref{aT1}. 
This will be done in four steps, first shrinking $H_v$ to a disjoint union
of compact sets respecting the decomposition (\ref{NewQuasiRep}),
then enlarging those sets by means of a finite covering with balls, 
then shrinking them again to get a set contained in $\cC_v(F_w)$ for a finite galois   
extension $F_w/K_v$, and finally taking the union of its conjugates
to get an $\Aut_c(\CC_v/K_v)$-stable set.  

\smallskip
For the first step, define $T_{v,-1} = H_v$.
Put $H_{v,0} = H_v \cap W_{v,0}$ and put $T_{v,0} = H_v \backslash W_{v,0}$.   
Inductively, for $k = 1, \ldots, D$, 
put $H_{v,k} = T_{v,k-1} \cap W_{v,k}$ and 
$T_{v,k} = T_{v,k-1} \backslash W_{v,k}$.  
Then each $T_{v,k}$ is compact, and $T_{v,k-1} = H_{v,k} \cup T_{v,k}$.  
Since $H_v \subseteq W_{v,0} \cup \cdots \cup W_{v,n}$ it follows that 
$T_{v,D} = \phi$ and 
\begin{equation*} 
H_v \ = \ H_{v,0} \cup H_{v,1} \cup \cdots \cup H_{v,D} \ .
\end{equation*} 
By construction the $H_{v,k}$ are pairwise disjoint,
but they are not in general compact.  

If $D_0 + 1 \le k \le D$, 
we claim that $H_{v,k} \subseteq W_{v,k} \cap \cC_v(F_{w_k})$.  
To see this, note that by definition 
$H_{v,k} \subseteq H_v \backslash (W_{v,0} \cup \cdots \cup W_{v,k-1})$.
From (\ref{NewQuasiRep}) it follows that 
\begin{equation*} 
H_{v,k} \ \subseteq \ \bigcup_{j=k}^D \big(W_{v,j} \cap \cC_v(F_{w_j}) \big) \ .
\end{equation*}
However, also $H_{v,k} \subseteq W_{v,k}$, and by (\ref{FW2NN1})  
this means that $H_{v,k} \cap \cC_v(F_{w_j}) = \phi$ 
for $j = k+1, \ldots, D$.  
Hence $H_{v,k} \subseteq W_{v,k} \cap \cC_v(F_{w_k})$.      

For each $r > 0$ and each $k = 0, \ldots, n$, put 
\begin{equation*} 
H_{v,k}(r) \ = \ 
\{z \in T_{v,k-1} : \|z,a\|_v \ge r \text{ for each $a \in T_{v,k}$} \}
           \ \subseteq \ H_{v,k} \ ,
\end{equation*} 
and then put 
\begin{equation*}
H_v(r) \ = \ \bigcup_{k=0}^D H_{v,k}(r) \ .
\end{equation*} 
Since $T_{v,k-1}$ and $T_{v,k}$ are compact, each $H_{v,k}(r)$ is compact,
and $H_v(r)$ is compact.  

The sets $H_v(r)$ increase monotonically as $r$ decreases.  
For each $z \in H_v$ there are an index $k$ such that $z \in H_{v,k}$
and an $r > 0$ such that $\|z,a\|_v > r$ for all $a \in T_{v,k}$.  
Thus the $H_v(r)$ form an exhaustion of $H_v$.  
Choose a sequence $r_1 > r_2 > \cdots > 0$ with 
$\lim_{m \rightarrow \infty} r_m = 0$.  
By Proposition \ref{InnerSeqProp}, if we take $M$ sufficiently large
and put $X_v = H_v(r_M)$, then 
\begin{equation} \label{XIneqs1NAA} 
\left\{ \begin{array}{ll} 
|G(x_i,x_j;X_v) - G(x_i,x_j;H_v)| \ < \ \varepsilon_v/3 & 
                \text{for all $i \ne j$ \ ,} \\
|V_{x_i}(X_v) - V_{x_i}(H_v)| \ < \ \varepsilon_v/3 &
                \text{for each $i$\ .} 
        \end{array} \right. 
\end{equation}
\index{Green's function!properties of} 
\index{Robin constant!properties of}
Combining (\ref{XIneqs1N}) and (\ref{XIneqs1NAA}) gives 
\begin{equation} \label{XIneqs1NBB} 
\left\{ \begin{array}{ll} 
G(x_i,x_j;X_v) \ <  \ \Gbar(x_i,x_j;E_v)| + 2\varepsilon_v/3 & 
                \text{for all $i \ne j$ \ ,} \\
V_{x_i}(X_v) \ < \ \Vbar_{x_i}(E_v)+  2\varepsilon_v/3 &
                \text{for each $i$\ .} 
        \end{array} \right. 
\end{equation}
\index{Green's function!properties of} 
\index{Robin constant!properties of}
For each $k = 0, \ldots, D$ put $X_{v,k} = H_{v,k}(r_M) \subset W_{v,k}$.  
Then the $X_{v,k}$ are compact and pairwise disjoint, and
\begin{equation*} 
X_v \ = \ \bigcup_{k=0}^D X_{v,k} \ .
\end{equation*} 
By the continuity of the spherical distance\index{spherical metric!continuity of}, 
there is an $R > 0$ such that for all $k \ne \ell$,
and all $x \in X_{v,k}$, $y \in X_{v,\ell}$, we have $\|x,y\|_v > R$.

\smallskip 
For the second step, we will enlarge $X_v$ to a set $Y_v$ 
which is the union of finitely many isometrically parametrizable balls
\index{isometrically parametrizable ball}
and compact sets, as follows.  

If $x \in X_{v,0}$, then $x \in W_{v,0}$, and by the definition of $W_{v,0}$ 
there is an isometrically parametrizable ball $B(a,r) \subseteq U_{v,0}$ 
\index{isometrically parametrizable ball}
which contains $x$.     
By Proposition \ref{SepClosureDense},  
$\cC_v(\tK_v^{\sep})$ is dense in $\cC_v(\CC_v)$,
so we can assume that $a \in \cC_v(\tK_v^{\sep})$.  Since $X_{v,0}$ is compact,
finitely many such balls $B(a_{0,1},r_{0,1}), \ldots, B(a_{0,m_0},r_{0,m_0})$ 
cover $X_{v,0}$.  Without loss, we can assume that $r_{0,\ell} < R$ 
and $r_{0,\ell} \in |K_v^{\times}|_v$ for each $\ell$;  by construction each 
$a_{0,\ell} \in \cC_v(\tK_v^{\sep})$.    

If $1 \le k \le D_0$ and $x \in X_{v,k}$, 
then by the definition of $W_{v,k}$ there is an isometrically parametrizable 
\index{isometrically parametrizable ball}
ball $B(a,r)$, with $a \in W_{v,k} \cap \cC_v(F_{w_k})$ 
and $B(a,r) \subset U_{v,k}$, such that $x \in B(a,r)$.
By compactness, finitely many such balls 
$B(a_{k,1},r_{k,1}), \ldots, B(a_{k,m_k},r_{k,m_k})$ 
cover $X_{v,k}$.  Without loss, we can assume that 
$r_{k,\ell} < R$ and $r_{k,\ell} \in |K_v^{\times}|_v$ for each $k,\ell$.
By construction each $a_{k,\ell} \in \cC_v(F_{w_k})$.   

If $D_0 + 1 \le k \le D$, and $x \in X_{v,k}$, 
then $x \in H_{v,k}$ and by the discussion above we have 
$x \in W_{v,k} \cap \cC_v(F_{w_k})$.     
As $W_{v,k}$ is open, there is an isometrically parametrizable
\index{isometrically parametrizable ball}
ball centered at $x$ for which $B(x,r) \subset W_{v,k}$.  
By the properties of an isometric parametrization, 
\index{isometric parametrization}
$B(x,r) \cap \cC_v(F_{w_k})$ is open in $\cC_v(F_{w_k})$.  
Since $X_{v,k}$ is compact and contained in $\cC_v(F_{w_k})$,
there are finitely many of these balls, say  
$B(a_{k,1},r_{k,1}), \ldots, B(a_{k,m_k},r_{k,m_k})$ for which 
\begin{equation} \label{FMADBAD}
X_{v,k} \ \subseteq \ 
 \bigcup_{\ell=1}^{m_k} \big( B(a_{k,\ell},r_{k,\ell}) \cap \cC_v(F_{w_k}) \big) \ .
\end{equation}  
Again we can assume that $r_{k,\ell} < R$ 
and $r_{k,\ell} \in |K_v^{\times}|_v$ for each $\ell$;  
by construction each $a_{k,\ell} \in \cC_v(F_{w_k})$. 
The right side of (\ref{FMADBAD}) is contained in 
$W_{v,k} \cap \cC_v(F_{w_k})$, hence in $U_v$.    

Put 
\begin{equation*} 
Y_v \ = \ \Big( \bigcup_{k=0}^{D_0} 
\Big( \bigcup_{\ell = 1}^{m_k} B(a_{k,\ell},r_{k,\ell}) \Big) \Big) \cup 
\Big( \bigcup_{k = D_0+1}^D 
\Big( \bigcup_{\ell = 1}^{m_k} B(a_{k,\ell},r_{k,\ell}) \cap \cC_v(F_{w_k}) \Big) \Big) \ .
\end{equation*} 
Here we have purposely omitted the intersection with $\cC_v(F_{w_k})$ 
for the balls $B(a_{k,\ell},r_{k,\ell})$ with $k \le D_0$.  
This means that $Y_v$ need not be contained in $U_v$.
However, $Y_v$ contains $X_v$. By the monotonicity of Green's functions, 
\index{Green's function!monotonic}
together with (\ref{XIneqs1NBB}), this gives  
\begin{equation} \label{XIneqs1NCC} 
\left\{ \begin{array}{ll} 
\Gbar(x_i,x_j;Y_v) \ < \ \Gbar(x_i,x_j;E_v) + 2\varepsilon_v/3 & 
                \text{for all $i \ne j$ \ ,} \\
\Vbar_{x_i}(Y_v) \ < \ \Vbar_{x_i}(E_v)+  2\varepsilon_v/3 &
                \text{for each $i$\ .} 
        \end{array} \right. 
\end{equation}  
 
\smallskip
For the third step, we will cut $Y_v$ down to a set $Z_v$ 
contained in $U_v \cap \cC_v(F_w)$, where $F_w$ is a suitable 
finite galois extension of $K_v$.   
By construction, $Y_v$ is the union of finitely many 
isometrically parametrizable balls whose radii belong to $|K_v^{\times}|_v$,\
\index{isometrically parametrizable ball}
and finitely many compact sets.    
By (\cite{RR1}, Theorems 4.2.16 and 4.3.11), it is algebraically capacitable. 
\index{algebraically capacitable} 
However, the proof of (\cite{RR1}, Theorem 4.3.11) gives more:  
the Green's function of $Y_v$ is the limit of Green's functions 
\index{Green's function!of a compact set}
of compact sets contained in finite extensions of $K_v$.  
An examination of the proof shows these sets can be chosen to lie in $U_v$. 
 
Explicitly, this comes out as follows.   
Put $\tF_{w_0} = \tK_v^{\sep}$, and put $\tF_{w_k} = F_{w_k}$ 
for $1 \le k \le D_0$.  For each $k = 0, \ldots, D_0$, 
choose an exhaustion of $\tF_{w_k}$ 
by an increasing sequence of finite separable extensions $\tF_{w_k,j}/K_v$: 
\begin{equation*} 
\tF_{w_k,1} \ \subseteq \ \tF_{w_k,2} \ \subseteq \ \tF_{w_k,3} \ \subseteq \ 
\cdots \ \subseteq \ \tF_{w_k} \ .
\end{equation*} 
Without loss, we can assume $\tF_{w_k,1}$ is large enough 
that each of $a_{k,1}, \ldots, a_{k,m_k}$ belongs to $\cC_v(\tF_{w_k,1})$.  
For each $j = 1,2,3 \cdots$, put
\begin{equation*} 
Z_{v,j} \ = \ \ \Big( \bigcup_{k=0}^{D_0} \big( \bigcup_{\ell = 1}^{m_k} 
B(a_{k,\ell},r_{k,\ell}) \cap \cC_v(\tF_{w_k,j}) \big) \Big) \cup 
\Big( \bigcup_{k = D_0+1}^D \big( \bigcup_{\ell = 1}^{m_k} 
B(a_{k,\ell},r_{k,\ell}) \cap \cC_v(F_{w_k}) \big) \Big) \ .
\end{equation*} 
The sets $Z_{v,j}$ play the same role as the sets $H_j$ in (\cite{RR1}, p.270),
though they are constructed somewhat differently. 
In (\cite{RR1}, p.269) the finite extensions $\tF_{w_k,j}$ (denoted $L_j$ there)
were made to exhaust the algebraic closure $\tK_v$, 
but all that is needed is the fact that for each $k = 0, \ldots, D_0$, 
as $j \rightarrow \infty$ either the ramification index $e_{w/v,k,j}$ 
or the residue degree $f_{w/v,k,j}$ of $\tF_{w_k,j}/K_v$ grows arbitrarily large.
This means that as $j \rightarrow \infty$, 
\begin{equation*}
\frac{1}{e_{w/v,k,j}} \cdot \frac{1}{q_v^{f_{w/v,k,j}}-1} \ \longrightarrow \ 0 \ .
\end{equation*}   
Hence the proof of (\cite{RR1}, Theorem 4.3.11) shows that 
if we take $J$ large enough and put $Z_v = Z_{v,J}$  then 
\begin{equation} \label{XIneqs1NDD} 
\left\{ \begin{array}{ll} 
|G(x_i,x_j;Z_v) - \Gbar(x_i,x_j;Y_v)| \ < \ \varepsilon_v/3 & 
                \text{for all $i \ne j$ \ ,} \\
|V_{x_i}(Z_v) - \Vbar_{x_i}(Y_v)| \ < \ \varepsilon_v/3 &
                \text{for each $i$\ .} 
        \end{array} \right. 
\end{equation} 
\index{Green's function!properties of} 
\index{Robin constant!properties of}
Furthermore, if $F_w$ is the galois closure of the composite of the fields 
$\tF_{w_0,J}, \ldots, \tF_{w_{D_0},J}$ and  $F_{w_{D_0+1}}, \ldots, F_{w_n}$, 
then $F_w$ is a finite galois extension of $K_v$ and 
$Z_v  \subset U_v \cap \cC_v(F_w)$.

\smallskip
For the last step, 
let $E_v^{\prime}$ be the union of the $\Gal(F_w/K_v)$-conjugates of $Z_v$.  
Then $E_v^{\prime}$ is compact.  
It is a finite union of sets of the kind in Theorem \ref{aT1}.C(2).  
Since $U_v$ is stable under $\Aut_c(\CC_v/K_v)$, 
$E_v^{\prime}$ is contained in $U_v$. 
From (\ref{XIneqs1NCC}), (\ref{XIneqs1NDD}), 
and the monotonicity of Green's functions, it follows that (\ref{XIneqs1}) holds:  that is,
\index{Green's function!monotonic} 
\begin{equation*}  
\left\{ \begin{array}{ll} 
G(x_i,x_j;E_v^{\prime}) \ < \ \Gbar(x_i,x_j;E_v) + \varepsilon_v & 
                \text{for all $i \ne j$ \ ,} \\
V_{x_i}(E_v^{\prime}) \ < \ \Vbar_{x_i}(E_v) +  \varepsilon_v &
                \text{for each $i$\ .} 
        \end{array} \right. 
\end{equation*}
\index{Green's function!properties of} 
\index{Robin constant!properties of}
This completes the proof.
\end{proof} 

\noindent{{\bf Theorem \text{\ref{aT1-A}}.} (Strong FSZ with LRC, producing points in $\EE$).}  
\index{Fekete-Szeg\"o theorem with LRC!Strong form|ii}
{\it Let $K$ be a global field, 
and let $\cC/K$ be a smooth, geometrically integral projective curve.
Let $\fX = \{x_1, \ldots, x_m\} \subset \cC(\tK)$ 
be a finite set of points stable under $\Aut(\tK/K)$, and let
$\EE = \prod_v E_v \subset \prod_v \cC_v(\CC_v)$ be an adelic set compatible with $\fX$,\index{compatible with $\fX$} 
such that each $E_v$ is stable under $\Aut_c(\CC_v/K_v)$.  
Let $S \subset \cM_K$ be a finite set of places $v$, containing all archimedean $v$,
such that $E_v$ is $\fX$-trivial for each $v \notin S$.
\index{$\fX$-trivial}

Assume that $\gammabar(\EE,\fX) > 1$.  
Assume also that for each $v \in S$,  
there is a $($possibly empty$)$ $\Aut_c(\CC_v/K_v)$-stable Borel subset 
$e_v \subset \cC_v(\CC_v)$ of inner capacity $0$ such that 
\index{capacity $= 0$} 

$(A)$ If $v$ is archimedean and $K_v \cong \CC$, 
then each point of $\cl(E_v) \backslash e_v$ is analytically accessible\index{analytically accessible}  
from the $\cC_v(\CC)$-interior of $E_v$.

$(B)$ If $v$ is archimedean and $K_v \cong \RR$, 
then each point of $\cl(E_v) \backslash e_v$ is  
 
\quad $(1)$ analytically accessible\index{analytically accessible}  from the $\cC_v(\CC)$-interior of $E_v$, or 

\quad $(2)$ is an endpoint of an open segment contained in $E_v \cap \cC_v(\RR)$.   

$(C)$ If $v$ is nonarchimedean, then $E_v$ is the disjoint union of $e_v$ 
and finitely many sets $E_{v,1}, \ldots, E_{v,D_v}$, where each $E_{v,\ell}$ is 

\quad $(1)$ open in $\cC_v(\CC_v)$, or 

\quad $(2)$ of the form $U_{v,\ell} \cap \cC_v(F_{w_\ell})$, where $U_{v,\ell}$ is open 
in $\cC_v(\CC_v)$ and  $F_{w_\ell}$ is a separable algebraic extension of $K_v$ 
contained in $\CC_v$ $($possibly of infinite degree$)$.  

Then there are infinitely many points $\alpha \in \cC(\tK^{\sep})$ such that for each $v \in \cM_K$, 
the $\Aut(\tK/K)$-conjugates of $\alpha$ all belong to $E_v$.   
}

\begin{proof}[Proof of Theorem \ref{aT1-A}, using Theorem \ref{aT1-A1}]  

{ \  }  

Assume Theorem \ref{aT1-A1}, 
and let $\EE = \prod_v E_v \subset \prod_v \cC_v(\CC_v)$
be an adelic set compatible with $\fX$\index{compatible with $\fX$}  
for which the hypotheses of Theorem \ref{aT1-A} hold.  
In this case, apart from a set of inner capacity $0$, 
\index{capacity $= 0$} 
each $E_v$ is itself a separable quasi-neighborhood. 
\index{quasi-neighborhood!separable}

We will construct new adelic sets $\EE^{\prime} = \prod_v E_v^{\prime}$ 
and $\UU^{\prime} = \prod_v U_v^{\prime}$, 
with $\EE^{\prime} \subseteq \UU^{\prime} \subseteq \EE$, 
such that the hypotheses of Theorem \ref{aT1-A1} hold for $\EE^{\prime}$
and $\UU^{\prime}$.  In fact, we will have $\EE^{\prime} = \UU^{\prime}$  
and $\gammabar(\EE^{\prime}, \fX) = \gammabar(\EE,\fX)$.  
Let $S \subseteq \cM_K$ be a finite set of places 
containing all archimedean places and all nonarchimedean places where $E_v$
is not $\fX$-trivial.  
\index{$\fX$-trivial}
For each $v \notin S$, put $E_v^{\prime} = U_v^{\prime} = E_v$.
  
For each archimedean $v \in S$ such that $K_v \cong \CC$, 
the set $\cl(E_v)$ is compact and there is a Borel 
subset $e_v \subset \cC_v(\CC)$ of inner capacity $0$ 
\index{capacity $= 0$} 
such that each point of $\cl(E_v) \backslash e_v$
is analytically accessible\index{analytically accessible}  from the $\cC_v(\CC)$-interior $E_v^{0}$ of $E_v$. 
Put $E_v^{\prime} = U_v^{\prime} = E_v^0$.  
Since $U_v^{\prime}$ is open, it is a quasi-neighborhood of $E_v^{\prime}$.
\index{quasi-neighborhood}    
By Proposition \ref{IdentifyGreenProp},
for each $\zeta \notin E_v$ we have $\Gbar(z,\zeta;E_v^{\prime}) = G(z,\zeta;E_v)$. 

For each archimedean $v \in S$ such that $K_v \cong \RR$, 
the set $\cl(E_v)$ is compact and there is a Borel 
subset $e_v \subset \cC_v(\CC)$ of inner capacity $0$ 
\index{capacity $= 0$} 
such that each point of $\cl(E_v) \backslash e_v$
is analytically accessible\index{analytically accessible}  from the $\cC_v(\CC)$-interior $E_v^{0}$ of $E_v$  
or from the $\cC_v(\RR)$-interior $E_v^1$ of $E_v \cap \cC_v(\RR)$. 
Put $E_v^{\prime} = U_v^{\prime} := E_v^0 \cup E_v^1$.  
Since $E_v$ is stable under complex conjugation, so are $E_v^0$ and $E_v^{\prime}$.  
Again $U_v^{\prime}$ is a quasi-neighborhood of $E_v^{\prime}$, 
\index{quasi-neighborhood}    
and by Proposition \ref{IdentifyGreenProp},
for each $\zeta \notin E_v$ we have $\Gbar(z,\zeta;E_v^{\prime}) = G(z,\zeta;E_v)$.

For each nonarchimedean $v \in S$, 
$E_v$ is the disjoint union of 
a Borel subset $e_v$ of inner capacity $0$ 
\index{capacity $= 0$} 
and sets $E_{v,1}, \ldots, E_{v,n}$, each of which is either open in $\cC_v(\CC_v)$
or of the form $U_{v,i} \cap \cC_v(F_{w_i})$ for some separable algebraic
extension $F_{w_i}/K_v$ contained in $\cC_v(\CC_v)$.  
Put $E_v^{\prime} = U_v^{\prime} = E_{v,1} \cup \cdots \cup E_{v,n}$.  
By assumption, $e_v$ is stable under $\Aut_c(\CC_v/K_v)$, 
hence so is $E_v^{\prime}$.  
By its construction, $U_v^{\prime}$ 
is a quasi-neighborhood of $E_v^{\prime}$.  
\index{quasi-neighborhood}
Finally, by Lemma \ref{CompactP2}, removing 
a set of inner capacity $0$ from a set does not change its Green's functions,
\index{Green's function!properties of}
\index{capacity $= 0$}  
so for each $\zeta \notin E_v$
we have $\Gbar(z,\zeta;E_v^{\prime}) = \Gbar(z,\zeta;E_v)$ 
\end{proof}

\noindent{{\bf Theorem \text{\ref{FSZi}}.} (FSZ with LRC and Ramification Side Conditions).}
\index{Fekete-Szeg\"o theorem with LRC!and Ramification Side Conditions|ii}  
{\it
Let $K$ be a global field, 
and let $\cC/K$ be a smooth, connected, projective curve.
Let $\fX = \{x_1, \ldots, x_m\} \subset \cC(\tK)$ 
be a finite, galois-stable set of points, and let
$\EE = \prod_v E_v \subset \prod_v \cC_v(\CC_v)$
be an adelic set compatible with $\fX$,\index{compatible with $\fX$}
such that each $E_v$ is stable under $\Aut_c(\CC_v/K_v)$.  

Let $S, S^{\prime}, S^{\prime \prime} \subset \cM_K$ be finite $($possibly empty$)$ 
sets of places of $K$ which are pairwise disjoint, such that the places in
$S^{\prime} \cup S^{\prime \prime}$ are nonarchimedean.  
Assume that $\gammabar(\EE,\fX) > 1$, and that 

$(A)$ for each $v \in S$, the set $E_v$ satisfies the conditions of Theorem $\ref{aT1}$ 
or Theorem $\ref{aT1-A}$.
 
$(B)$ for each $v \in S^{\prime}$, either $E_v$ is $\fX$-trivial, 
\index{$\fX$-trivial}
or $E_v$ is a finite union of 
closed isometrically parametrizable balls $B(a_i,r_i)$ 
\index{isometrically parametrizable ball}
whose radii belong to the value group of $K_v^{\times}$
and whose centers belong to an unramified extension of $K_v$;   
 
$(C)$ for each $v \in S^{\prime\prime}$, 
either $E_v$ is $\fX$-trivial and $E_v \cap \cC_v(K_v)$ is nonempty, 
\index{$\fX$-trivial}
or $E_v$ is a finite union of closed and$/$or open isometrically parametrizable 
\index{isometrically parametrizable ball}
balls $B(a_i,r_i)$, $B(a_j,r_j)^-$ with centers in $\cC_v(K_v)$.

Then there are infinitely many  points $\alpha \in \cC(\tK^{\sep})$ such that 

$(1)$ for each $v \in \cM_K$, the $\Aut(\tK/K)$-conjugates of $\alpha$ 
all belong to $E_v$;   


$(2)$ for each $v \in S^{\prime}$, each place of $K(\alpha)/K$ above $v$ is unramified over $v$; 

$(3)$ for each $v \in S^{\prime \prime}$, 
each place of $K(\alpha)/K$ above $v$ is totally ramified over $v$.
}  

\begin{proof}[Proof of Theorem \ref{FSZi}, using Theorem \ref{aT1-A}]


For each $v \in S^{\prime} \cup S^{\prime\prime}$, 
the hypotheses in (B) and (C) will enable us to replace 
the given set $E_v$ by sets of the form $\cC_v(F_w) \cap E_v$ 
for suitably chosen finite galois extensions $F_w/K_v$, 
which are unramified if $v \in S^{\prime}$ and are totally ramified 
if $v \in S^{\prime\prime}$, in such a way that 
we still have  $\gammabar(\EE,\fX) > 1$.
 
More precisely, we claim that we can choose the extensions $F_w/K_v$
so that the new Green's matrix is arbitrarily near the old one.
Since $\gammabar(\EE,\fX) > 1$ if and only if $\overline{\Gamma}(\EE,\fX)$ is 
negative definite,\index{Green's matrix!negative definite} the new Green's matrix will be negative definite if it is
sufficiently close to the old one.  Thus, the theorem reduces to Theorem \ref{aT1-A}. 

The claim is a consequence of explicit formulas for the Robin constants 
\index{Robin constant!examples!nonarchimedean}
and Green's functions of the sets in question, derived in (\cite{RR1}, pp.353-359) 
\index{Green's function!examples!nonarchimedean}
and stated 
in (\ref{FSolve}) of this work.  
Fix a nonarchimedean place $v$, and fix a spherical metric\index{spherical metric} on $\cC_v(\CC_v)$.   
Let $H_v = B(a,r) \subset \cC_v(\CC_v)$ be a closed isometrically 
parametrizable ball, and take $\zeta \in \cC_v(\CC_v) \backslash B(a,r)$.  
By (\cite{RR1}, Theorem 4.3.15, p.274) isometrically parametrizable balls are
\index{isometrically parametrizable ball}
algebraically capacitable, and by the proof of (\cite{RR1}, Theorem 4.4.4) 
\index{algebraically capacitable}
their Green's functions and upper Green's functions coincide.
\index{Green's function}
\index{Green's function!upper}  
Fix a normalization for the canonical distance $[z,w]_{\zeta}$.  Then there is
\index{canonical distance!$[z,w]_{\zeta}$}  
an $R > 0$ such that $B(a,r) = \{z \in \cC_v(\CC_v) : [z,a]_{\zeta} \le R\}$, 
and in terms of this $R$ (see \cite{RR1}, p.357) 
we have  
\index{Robin constant!examples!nonarchimedean} 
\index{Green's function!examples!nonarchimedean}
\begin{eqnarray}
V_{\zeta}(H_v) & = & -\log_v(R) \ , \label{FGGG1} \\
G(z,\zeta;H_v) & = & \left\{ \begin{array}{ll} 
                0 & \text{if $z \in H_v$,} \\       
                \log_v([z,a]_{\zeta}/R) & \text{if $z \notin B(a,r)$.}
                             \end{array} \right. \notag 
\end{eqnarray}
It will be useful to define $u_{H_v}(z,\zeta) = V_{\zeta}(H_v) - G(z,\zeta;H_v)$, so  
\begin{equation}
u_{H_v}(z,\zeta) \ = \ \left\{ \begin{array}{ll} 
                -\log_v(R) & \text{if $z \in H_v$,} \\       
                -\log_v([z,a]_{\zeta}) & \text{if $z \notin B(a,r)$.}
                             \end{array} \right. \label{FGGG2} 
\end{equation} 
Furthermore, if $F_w/K_v$ is a finite
extension with ramification index $e$ and residue degree $f$, 
if $q_v$ is the order of the residue field of $K_v$, 
and if $a \in \cC_v(F_w)$ and $r$ belongs to the value group of $F_w^{\times}$, 
then for the set $\cC_v(F_w) \cap H_v$ we have (see \cite{RR1}, p.358) 
\begin{eqnarray}
V_{\zeta}(\cC_v(F_w) \cap H_v) & = & -\log_v(R) + \frac{1}{e(q_v^f-1)} \ , \label{FGGG3} \\
u_{\cC_v(F_w) \cap H_v}(z,\zeta)   
          & = & \left\{ \begin{array}{ll} 
    -\log_v(R) + \frac{1}{e(q_v^f-1)} & \text{if $z \in \cC_v(F_w) \cap H_v$,} \\       
    -\log_v([z,a]_{\zeta}) & \text{if $z \notin B(a,r)$.}
                             \end{array} \right. \label{FGGG4} 
\end{eqnarray} 

More generally, if $H_v = \cup_{i=1}^N B(a_i,r_i)$ is a finite union
of closed, pairwise disjoint isometrically parametrizable balls, write $H_{v,i} = B(a_i,r_i)$;
\index{isometrically parametrizable ball}
then $V_{\zeta}(H_v)$ and $G(z,\zeta;H_v)$ can be determined 
\index{Robin constant!computing nonarchimedean} 
by solving the following system of equations for $V, s_1, \ldots, s_N$
(see \cite{RR1}, p.359):  
\begin{equation}
\left\{ \begin{array}{l} 
            1 \ = \ 0V + s_1 + s_2 + \cdots + s_N \ ,\\
            0 \ = \ V - \sum_{i=1}^N s_i u_{H_{v,i}}(a_j,\zeta) 
                      \quad \text{for $j = 1, \ldots, N$}
        \end{array} \right. \label{FGGS1} 
\end{equation} 
By (\cite{RR1}, Proposition 4.2.7) the solution is unique, 
and has $s_1, \ldots, s_N > 0$;  in terms of it, 
\index{Robin constant!computing nonarchimedean}
\begin{eqnarray*} 
V_{\zeta}(H_v) & = & V \ ,\\
G(z,\zeta;H_v) & = & V - \sum_{i=1}^N s_i u_{H_{v,i}}(z,\zeta) \ .  
\end{eqnarray*}
Similarly, if $F_w/K_v$ is a finite
extension with ramification index $e$ and residue degree $f$, 
and if the $a_i$ belong to $\cC_v(F_w)$ 
and the $r_i$ belong to the value group of $F_w^{\times}$, 
then we can determine $V_{\zeta}(\cC_v(F_w) \cap H_v)$ 
and $G(z,\zeta;\cC_v(F_w) \cap H_v)$ by solving for 
$V_w, s_{1,w}, \ldots, s_{N,w}$ in the following system of equations:
\index{Robin constant!computing nonarchimedean} 
\begin{equation}
\left\{ \begin{array}{l} 
        1 \ = \ 0 V_w + s_{1,w} + s_{2,w} + \cdots + s_{N,w} \ ,\\
   0 \ = \ V_w - \sum_{i=1}^N s_{i,w} 
                     u_{\cC_v(F_w) \cap H_{v,i}}(a_j,\zeta) 
                      \quad \text{for $j = 1, \ldots, N$}
        \end{array} \right.  \label{FGGS2}
\end{equation} 
By the existence and uniqueness of the equilibrium measure for $\cC_v(F_w) \cap H_v$
(\cite{RR1}, Theorem 3.1.12), 
again the solution is unique, with $s_{1,w}, \ldots, s_{N,w} > 0$;  and 
\index{Robin constant!computing nonarchimedean}
\begin{eqnarray*} 
V_{\zeta}(\cC_v(F_w) \cap H_v) & = & V_w  \ ,\\
G(z,\zeta;\cC_v(F_w) \cap H_v) & = & 
     V_w - \sum_{i=1}^N s_{i,w} u_{\cC_v(F_w) \cap H_{v,i}}(z,\zeta) \ .  
\end{eqnarray*}  
Comparing the systems (\ref{FGGS1}) and (\ref{FGGS2}), 
as $F_w$ passes through a sequence of extensions
for which $1/(e(q_v^f-1)) \rightarrow 0$,
then the $V_{\zeta}(\cC_v(F_w) \cap H_v)$ converge to $V_{\zeta}(H_v)$ and the 
$G(z,\zeta; \cC_v(F_w) \cap H_v)$ converge (uniformly) to $G(z,\zeta;H_v)$.

\smallskip
We now apply this to the sets $E_v$ for $v \in S^{\prime} \cup S^{\prime \prime}$ in the theorem.
  
First suppose $E_v = \bigcup_{i=1}^N B(a_i,r_i)$ is a finite union of closed 
isometrically parametrizable balls.  Without loss, we can assume the $B(a_i,r_i)$
\index{isometrically parametrizable ball}
are pairwise disjoint.  

If $v \in S^{\prime}$, there is a finite unramified extension
$F_w^{\prime}/K_v$ with $a_1, \ldots, a_N \in \cC_v(F_w^{\prime})$, and the $r_i$
belong to the value group of $K_v^{\times}$.  Letting $F_w$ pass through all finite unramified
extensions of $K_v$ containing $F_w^{\prime}$, for all $x_i \ne x_j \in \fX$ 
we can make $V_{x_i}(\cC_v(F_w) \cap E_v)$ arbitrarily near $V_{x_i}(E_v)$, and we can make
\index{Robin constant!properties of}
the $G(x_i,x_j;\cC_v(F_w) \cap E_v)$ arbitrarily near the $G(x_i,x_j;E_v)$.  

Similarly if $v \in S^{\prime\prime}$, then $a_1, \ldots, a_N \in \cC_v(K_v)$.  Letting $F_w/K_v$
pass through a sequence of finite, galois, totally ramified extensions for which 
$e_{w/v} \rightarrow \infty$ (for example, cyclotomic $p$-extensions, 
where $p$ is the residue characteristic of $K_v$), 
we obtain the same conclusion as before.  If $v \in S^{\prime\prime}$ and 
$E_v = (\bigcup_{i=1}^N B(a_i,r_i)) \cup (\bigcup_{j=N+1}^{N+M} B(a_j,r_j)^-)$, 
then by exhausting the open balls with closed balls $B(a_j,r_j^{\prime})$ and taking a limit 
as the $r_j^{\prime} \rightarrow r_j$, we are reduced to the previous case.  
Note that any compact  $H_v \subset E_v$ is contained in a set of the form 
$(\cup_{i=1}^N B(a_i,r_i)) \cup (\cup_{j=N+1}^{N+M} B(a_j,r_j^{\prime}))$. 
Since the upper Green's function $\Gbar(z,w;E_v)$ is by definition the (pointwise) limit of the
\index{Green's function!upper}
upper Green's functions $\Gbar(z,w;H_v)$ for compact $H_v \subset E_v$, and 
upper Green's functions are monotonic under containment, $\Gbar(z,w;E_v)$ is the 
\index{Green's function!monotonic}
limit of the Green's functions for the unions of closed balls discussed above, and hence also  
of the Green's functions $G(z,w;\cC_v(F_w) \cap E_v)$ as $F_w$ passes through finite,
galois, totally ramified extensions of $K_v$.    

Next, consider the case where $E_v$ is $\fX$-trivial.\index{$\fX$-trivial} 
The $\fX$-triviality implies that $\cC_v$ has good
reduction at $v$ and the points of $\fX$ specialize to distinct points $\pmod{v}$. 
Furthermore, if $\|x,y\|_v$ is the spherical metric\index{spherical metric} associated to the given embedding of $\cC_v$,
then the canonical distance (up to scaling by a constant) is given by 
$[x,y]_{\zeta} = \|x,y\|_v/(\|x,\zeta\|_v \|y,\zeta\|_v)$
\index{canonical distance!$[z,w]_{\zeta}$}  
(see \cite{RR1}, p.91).  Relative to this normalization of the canonical distance,
$\gamma_{\zeta}(E_v) = 1$ for each $\zeta \notin E_v$. 

First suppose $v \in S^{\prime}$.  
Let $k_v$ be the residue field of the ring of integers of $K_v$, 
and let $\overline{k}_v$ be its algebraic closure.  Then $\overline{k}_v$ is the residue
field of $K_v^{nr}$, the maximal unramified algebraic extension of $K_v$.  
Write $\overline{\cC}_v$ for the reduction of $\cC_v \pmod{v}$.  
By Hensel's lemma, each point of $\overline{\cC}_v(\overline{k}_v)$ lifts
to a point in $\cC_v(K_v^{nr})$.  Put $r = 1/q_v \in |K_v^{\times}|_v$.  
Then for arbitrarily large $N$ 
we can find galois-stable sets of the form $H_v(N) = \bigcup_{i=1}^N B(a_i,r) \subset E_v$, 
where each $a_i \in\cC_v(K_v^{nr})$ and distinct $a_i$ specialize to distinct points $\pmod{v}$.  
For each $\zeta \notin E_v$, and each ball $B(a_i,r)$,  we have 
$B(a_i,r) = \{z \in \cC_v(\CC_v) : [z,a]_{\zeta} \le 1/q_v\}$,  so by (\ref{FGGG2}) 
\begin{eqnarray*}
u_{B(a_i,r)}(z,\zeta) \ = \ \left\{ \begin{array}{ll} 
                1 & \text{if $z \in B(a_i,r)$,} \\       
                -\log_v([z,a_i]_{\zeta}) & \text{if $z \notin B(a_i,r)$.}
                             \end{array} \right. 
\end{eqnarray*} 
In particular, $u_{B(a_i,r)}(a_j,\zeta) = 0$ for each $j \ne i$.  
Inserting this in (\ref{FGGS1}), 
we find that $V_{\zeta}(H_v(N)) = 1/N$, and that if $z \notin E_v$ then 
v\index{Robin constant!nonarchimedean}
$G(z,\zeta;H_v(N)) = G(z,\zeta;E_v) + 1/N$.  Thus 
by replacing $E_v$ with $H_v(N)$  for a sufficiently large $N$, 
we are reduced to a previous case.  

If $v \in S^{\prime\prime}$, 
then by hypothesis $\cC_v(K_v) \cap E_v$ is nonempty;  fix $a \in \cC_v(K_v) \cap E_v$. 
Consider the open ball $B(a,1)^- \subset E_v$.  Exhausting it by closed balls $B(a,r)$,
and noting that $B(a,r) = \{z \in \cC_v(\CC_v) : [z,a]_{\zeta} \le r\}$ 
for each $\zeta \notin E_v$, it follows by (\ref{FGGG1}) that
$\gammabar_{\zeta}(B(a,1)^-) = 1$.  By (\cite{RR1}, Lemma 4.4.7), 
$\Gbar(z,\zeta;B(a,1)^-) = G(z,\zeta;E_v)$ for all $z, \zeta \notin E_v$.  
Thus, by replacing $E_v$ with $B(a,1)^-$, 
again we are reduced to a case considered before.                         
\end{proof}

\noindent{{\bf Theorem \text{ \ref{FSZii}}.} 
(Fekete/Fekete-Szeg\"o with LRC for Algebraically Capacitable Sets).}
\index{Fekete-Szeg\"o theorem with LRC!for algebraically capacitable sets|ii} 
\index{algebraically capacitable} 
{\it  
Let $K$ be a global field 
and let $\cC/K$ be a smooth, connected, projective curve.
Let $\fX = \{x_1, \ldots, x_m\} \subset \cC(\tK)$ 
be a finite, galois-stable set of points, and let
$\EE = \prod_v E_v \subset \prod_v \cC_v(\CC_v)$
be an adelic set compatible with $\fX$.\index{compatible with $\fX$} 

Assume that each $E_v$ is algebraically capacitable and stable under $\Aut_c(\CC_v/K_v)$.  Then  
\index{algebraically capacitable}

\smallskip
$(A)$  If all the eigenvalues of $\Gamma(\EE,\fX)$ are non-positive 
\index{Green's matrix!global}
$($that is, $\Gamma(\EE,\fX)$ is either negative definite\index{Green's matrix!negative definite} or negative semi-definite$)$,  
let $\UU = \prod_v U_v$ be a separable $K$-rational quasi-neighborhood of $\EE$ 
\index{quasi-neighborhood!separable}
such that there is at least one place $v_0$ where $E_{v_0}$ is compact and the 
quasi-neighborhood $U_{v_0}$ properly contains $E_{v_0}$. If $v_0$ is archimedean, assume also 
\index{quasi-neighborhood}
that $U_{v_0}$ meets each component of $\cC_{v_0}(\CC) \backslash E_{v_0}$ containing a point of $\fX$.  
Then there are infinitely many points $\alpha \in \cC(\tK^{\sep})$ such that 
all the conjugates of $\alpha$ belong to $\UU$.  

\smallskip
$(B)$ If some eigenvalue of $\Gamma(\EE,\fX)$ is positive $($that is, $\Gamma(\EE,\fX)$ is 
\index{Green's matrix!global}
either indefinite, nonzero and positive semi-definite, or positive definite$)$,
there is an adelic neighborhood $\UU$ of $\EE$ such that only 
finitely many points $\alpha \in \cC(\tK)$ have all their conjugates in $\UU$.
}

\begin{proof}[Proof of Theorem \ref{FSZii}, using Theorem \ref{aT1-A1}]

Since each $E_v$ is algebraically capacitable, 
\index{algebraically capacitable}
we have $\gammabar(\EE,\fX) = \gamma(\EE,\fX)$.  
Recall that a symmetric matrix $\Gamma \in M_k(\RR)$ with non-negative off-diagonal entries 
is called {\em irreducible}\index{matrix!irreducible} (\cite{RR1}, p.328) if the graph on the set $\{1, \ldots, k\}$, 
for which there is an edge between $i$ and $j$ iff $\Gamma_{i,j} > 0$, is connected. 
By (\cite{RR1}, Lemma 5.1.7, p.328) if $\Gamma$ is irreducible,\index{irreducible matrix|ii}  
then $\val(\Gamma)$ is positive, $0$, or negative, according as the largest eigenvalue 
of $\Gamma$ is positive, $0$, or negative.  
 
By re-ordering the elements of $\fX$ if necessary, 
we can bring $\Gamma(\EE,\fX)$ to block-diagonal form
\index{Green's matrix!global}
$\diag(\Gamma_1, \ldots, \Gamma_r)$, where each $\Gamma_i$ is irreducible. 
Note that $\Gamma(\EE,\fX)$ is negative definite\index{Green's matrix!negative definite} if and only if each $\Gamma_i$ is negative
definite, and is negative semi-definite if and only if each $\Gamma_i$ is negative definite
or negative semi-definite. 

If $\Gamma(\EE,\fX)$ is negative definite, 
\index{Green's matrix!global}\index{Green's matrix!negative definite}
then $\gamma(\EE,\fX) > 1$, and the result follows from Theorem \ref{aT1-A1}.
If each $\Gamma_i$ is negative semi-definite, then  
by enlarging the set $E_{_0}$ 
within its quasi-neighborhood $U_{v_0}$ (keeping $E_{v_0}$ stable under $\Aut_c(\CC_{v_0}/K_{v_0}$), 
\index{quasi-neighborhood}
we can decrease all the 
diagonal entries of $\Gamma(\EE,\fX)$, 
while either decreasing or leaving unchanged each off-diagonal entry.    
This makes $\Gamma(\EE,\fX)$ negative definite,\index{Green's matrix!negative definite} and we can again apply Theorem \ref{aT1-A1}.
\index{Green's matrix!global}

If some $\Gamma_i$ has a positive eigenvalue, let $\fX^{\prime}$ be the subset of $\fX$ 
consisting of all $x_{\ell} \in \fX$ corresponding to blocks which have positive eigenvalues.
Enlarge each $E_v$ which is $\fX$-trivial to a set $E_v^{\prime}$ which is $\fX^{\prime}$-trivial,
\index{$\fX$-trivial}
and let $\EE^{\prime}$ be corresponding adelic set.  Since the action of $\Aut(\tK/K)$ on $\fX$ 
permutes the $x_{\ell}$ and hence the blocks $\Gamma_i$, 
the sets  $\fX^{\prime}$ and $\EE^{\prime}$ are galois stable.  
By (\cite{RR1}, Lemma 5.1.7 and Theorem 5.1.6, p.328),  
$\val(\Gamma(\EE^{\prime},\fX^{\prime})) > 0$ 
and hence $\gamma(\EE^{\prime},\fX^{\prime}) < 1$.  The result now follows from Fekete's theorem 
\index{Fekete's theorem} 
applied to $\EE^{\prime}$ and $\fX^{\prime}$  (\cite{RR1}, Theorem 6.3.1, p.414).      
\end{proof} 

We now prepare for the proofs of Theorems \ref{aT1-B1} and \ref{aT1-B2}.  In the following, 
we assume familiarity with Berkovich analytic spaces (see \cite{Berk}) 
\index{Berkovich!analytic space} 
and Thuillier's potential theory\index{potential theory!on Berkovich curves} on Berkovich curves (\cite{Th}).
\index{Thuillier, Amaury}
\index{Berkovich!curve}
\index{Berkovich, Vladimir}    
For nonarchimedean places $v$, 
Thuillier established the compatibility of capacities 
\index{Thuillier, Amaury}
for sets in $\cC_v^{\an}$,  defined by him in (\cite{Th}), with capacities for sets in $\cC_v(\CC_v)$, 
as defined by Rumely in (\cite{RR1}) and used in this work (see  (\cite{Th}, Appendix 5.1)).
\index{Rumely, Robert} 
However, he did not explicitly state the compatibility of Green's functions.
\index{Green's function!Berkovich}  
Before proving Theorems \ref{aT1-B1} and \ref{aT1-B2}, we establish this.

Recall that for each compact, nonpolar\index{nonpolar set} subset $\BerkE_v \subset \cC_v^{\an}$ 
and each $\zeta \in \cC_v^{\an} \backslash \BerkE_v $, Thuillier (\cite{Th}, Th\'eor\`eme 3.6.15) has 
\index{Thuillier, Amaury}
constructed a Green's function $g_{\zeta,\BerkE_v}(z)$ which is non-negative,
\index{Green's function!Thuillier} 
vanishes on $\BerkE_v$ except possibly on a set of capacity $0$, is subharmonic in $\cC_v^{\an}$,
\index{capacity $= 0$} 
harmonic in $\cC_v^{\an} \backslash (\BerkE_v \cup \{\zeta\})$, 
and satisfies the distributional equation $dd^c g_{\zeta,\BerkE_v} = \mu - \delta_\zeta$ 
where $\mu$ is a probability measure supported on $K$.  We write 
$G(z,\zeta;\BerkE_v)^{\an}$ for $g_{\zeta,\BerkE_v}(z)$, regarding it as a function of 
two variables.  

By abuse of language, 
we write $\cC_v^{\an}$ for the topological space underlying the ringed space $\cC_v^{\an}$.
Following Thuillier, let $I(\cC_v^{\an}) := \cC_v^{\an} \backslash \cC_v(\CC_v)$ be the 
\index{Thuillier, Amaury}
set of non-classical points of $\cC_v^{\an}$.  Our first proposition shows that 
Berkovich Green's functions have properties analogous to those of classical Green's functions.   
\index{Berkovich!Green's function}  
\index{Green's function!Berkovich}

\begin{proposition} \label{BerkGreenPropertiesProp}
\index{Green's function!Berkovich!properties|ii}
\index{Robin constant!Berkovich!properties of|ii}
Let $K$ be a global field, and let $\cC/K$ be a smooth, connected, projective curve.  
Let $v$ be a nonarchimedean place of $K$, and let $\cC_v^{\an}$ be the Berkovich analytification 
\index{Berkovich!analytification} 
of $\cC_v \times_{K_v} \Spec(\CC_v)$.  
Let $\BerkE_v \subsetneq \cC_v^{an}$ be a proper compact, nonpolar subset of $\cC_v^{an}$.  Then
\index{nonpolar set} 

$(A)$ For each $\zeta \in \cC_v(\CC_v) \backslash \BerkE_v$, 
if we fix a uniformizer $g_{\zeta}(z)$ at $\zeta$,
then the Robin constant  
\index{Robin constant!Berkovich|ii}
\begin{equation*} 
V_{\zeta}(\BerkE_v)^{\an} \ = \ \lim_{\substack{z \rightarrow \zeta \\ z \in \cC_v^{\an}}} G(z,\zeta;\BerkE_v)^{\an} + \log(|g_{\zeta}(z)|_v) 
\end{equation*} 
is well defined and finite.

$(B)$  For all $x, y \in \cC_v^{an} \backslash \BerkE_v$ with $x \ne y$,  
\begin{equation*}
G(x,y;\BerkE_v)^{\an} \ = \ G(y,x;\BerkE_v)^{\an} \ .
\end{equation*} 

$(C)$ Let $\BerkE_{v,1} \subseteq \BerkE_{v,2} $ be nonpolar, proper compact subsets of $\cC_v^{an}$. 
\index{nonpolar set} 
Then for each $y \in \cC_v^{\an} \backslash \BerkE_{v,2}$, for all $x \in \cC_v^{\an}$ with $x \ne y$ we have
\begin{equation*}
G(x,y;\BerkE_{v,1})^{\an} \ \ge \ G(x,y;\BerkE_{v,2})^{\an} \ .
\end{equation*} 
This also holds when $y \in I(\cC_v^{\an}) \backslash \BerkE_{v,2}$ and $x = y$.  
For each $\zeta \in \cC_v(\CC_v) \backslash \BerkE_{v,2}$,
\begin{equation*} 
V_{\zeta}(\BerkE_{v,1})^{\an} \ \ge \ V_{\zeta}(\BerkE_{v,2})^{\an} \ .
\end{equation*} 
\index{Robin constant!Berkovich!properties of|ii}

$(D)$  Let $\BerkK_1 \supseteq \BerkK_2 \supseteq \cdots \supseteq \BerkK_n \cdots \supseteq \BerkE_v$
be a descending sequence of compact sets with $\bigcap_{n=1}^\infty \BerkK_n = \BerkE_v$. 
Then for all $x, y \in \cC_v^{\an} \backslash \BerkE_v$ such that $x \ne y$,
or such that $x = y \in I(\cC_v^{\an})$,  \begin{equation} \label{GDxyLim} 
\lim_{n \rightarrow \infty} G(x,y;\BerkK_n)^{\an} \ = \ G(x,y;\BerkE_v)^{\an} \ ,
\end{equation}
and for each $\zeta \in \cC_v(\CC_v) \backslash \BerkE_v$,
\begin{equation} \label{VDzzLim} 
\lim_{n \rightarrow \infty} V_{\zeta}(\BerkK_n)^{\an} \ = \ V_{\zeta}(\BerkE_v)^{\an} \ .
\end{equation}
\index{Robin constant!Berkovich!properties of|ii}

$(E)$  Let  $\BerkK_1 \subseteq \BerkK_2 \subseteq \cdots \subseteq \BerkK_n \cdots \subseteq \BerkE_v$
be an ascending sequence of compact sets with $\bigcup_{n=1}^{\infty} \BerkK_n = \BerkE_v$.  
Then for all $x, y \in \cC_v^{\an} \backslash \BerkE_v$ such that $x \ne y$, 
or such that $x = y \in I(\cC_v^{\an})$,   
\begin{equation} \label{GAxyLim} 
\lim_{n \rightarrow \infty} G(x,y;\BerkK_n)^{\an} \ = \ G(x,y;\BerkE_v)^{\an} \ ,
\end{equation}
and for each $\zeta \in \cC_v(\CC_v) \backslash \BerkE_v$,
\begin{equation} \label{VAzzLim} 
\lim_{n \rightarrow \infty} V_{\zeta}(\BerkK_n)^{\an} \ = \ V_{\zeta}(\BerkE_v)^{\an} \ .
\end{equation}
\index{Robin constant!Berkovich!properties of|ii}

$(F)$  For each $\sigma \in \Aut_c(\CC_v/K_v)$, 
and all $x, y \in \cC_v^{an} \backslash \BerkE_v$ with $x \ne y$,
\begin{equation*} 
G(\sigma(x),\sigma(y);\sigma(\BerkE_v))^{\an} \ = \ G(x,y;\BerkE_v)^{\an} \ ,
\end{equation*} 
and for each $\zeta \in \cC_v(\CC_v) \backslash \BerkE_v$, if the uniformizer $g_{\sigma(\zeta)}(z)$
is taken to be $\sigma(g_{\zeta})(z)$, then 
\index{Robin constant!Berkovich!properties of|ii}
\begin{equation*} 
V_{\sigma(\zeta)}(\sigma(\BerkE_v))^{\an} \ = \ V_{\zeta}(\BerkE_v)^{\an} \ .
\end{equation*} 
\end{proposition}  

\begin{proof}  We first prove (A).  When $\zeta \in \cC_v(\CC_v) \backslash \BerkE_v$, 
the existence and finiteness of the limit defining $V_{\zeta}(\BerkE_v)^{an}$
\index{Robin constant!Berkovich!properties of}
follows from the construction of $g_{\zeta,\BerkE_v}$:  
\index{Green's function!Thuillier} 
see the proof of (\cite{Th}, Th\'eor\`eme 3.6.15), 
noting that if $V$ is the Berkovich closure of a suitably 
\index{Berkovich!closure} 
small isometrically parametrizable ball $B(\zeta,r)$ and $y$ is its unique boundary point, 
\index{isometrically parametrizable ball}
the restriction of $\log(|g_{\zeta}|_v)$ to $V$ 
satisfies $dd^c\log(|g_{\zeta}|_v) = \delta_{\zeta} - \delta_y$ and thus coincides, 
up to an additive constant, with the function $g_{\zeta,y}(z)$ from (\cite{Th}, Lemma 3.4.14).  

\smallskip
We next establish the diagonal case in (D) and (E). 
Let  $\{\BerkK_n\}_{n \ge 1}$ be a sequence of compact sets such that 
 $\BerkK_1 \subseteq \BerkK_2 \subseteq \cdots \subseteq \BerkK_n 
\cdots \subseteq \BerkE_v$
and $\bigcup_{n =1}^{\infty} \BerkK_n = \BerkE_v$, 
or $\BerkK_1 \supseteq \BerkK_2 \supseteq \cdots \supseteq \BerkK_n \cdots \supseteq \BerkE_v$
and $\bigcap_{n=1}^\infty \BerkK_n = \BerkE_v$.    

First suppose $x = y \in I(\cC_v^{\an}) \backslash \BerkE_v$, 
and put $\Omega = \cC_v^{\an} \backslash \{y\}$.  
For each compact nonpolar  $\BerkK \subset \Omega$,
\index{nonpolar set} 
let $C(\BerkK,\Omega)$ be the capacity defined in (\cite{Th}, \S3.6).
\index{capacity!Thuillier} 
The construction in (\cite{Th}, Th\'eor\`eme 3.6.15) 
shows that $g_{y,\BerkE_v}(y) = C(\BerkE_v,\Omega)^{-1}$ 
and $g_{y,\BerkK_n}(y) = C(\BerkK_n,\Omega)^{-1}$ for all $n$.
By(\cite{Th}, Proposition 3.6.8, parts (ii) and (iv)), we have 
\begin{equation} \label{ThCapLim111} 
C(\BerkE_v,\Omega) \ = \ \lim_{n \rightarrow \infty} C(\BerkK_n,\Omega) \ , 
\end{equation} 
so (\ref{GDxyLim}) and (\ref{GAxyLim}) hold when $x = y \in I(\cC_v^{\an})$.  

A word is in order concerning the proof of (\cite{Th}, Proposition 3.6.8).  For an arbitrary subset 
$\BerkA \subset \Omega$, Thuillier defines
\index{Thuillier, Amaury}
\begin{equation*} 
C^*(\BerkA,\Omega) \ = \ \sup_{\text{compact $\BerkK \subseteq \BerkA$}} C(\BerkK,\Omega) \ ,
\end{equation*} 
then shows that $C^*(\cdot,\Omega)$ is a Choquet capacity on subsets of $\Omega$. 
\index{capacity!Choquet capacity|ii}\index{Choquet capacity|ii} 
Recall that this means $C^*(\cdot,\Omega)$ is an increasing set function, valued in $\RR_{\ge 0}$
with $C(\Omega,\phi) = 0$, such that 
\begin{enumerate} 
\item For each increasing sequence $\BerkA_1 \subseteq \BerkA_2 \subseteq \cdots $ 
of arbitrary sets in $\Omega$, 
if $A = \bigcup_{n=1}^{\infty} \BerkA_n$ 
then $C^*(\Omega,A) = \lim_{n \rightarrow \infty} C^*(\Omega,A_n)$.
\item For each decreasing sequence of {\em open} sets
 $\BerkU_1 \supseteq \BerkU_2 \supseteq \cdots$ in $\Omega$, 
if $A = \bigcap_{n=1}^{\infty} \BerkU_n$ 
then $C^*(\Omega,A) = \lim_{n \rightarrow \infty} C^*(\Omega,U_n)$.
\end{enumerate}  
In (\cite{Th}, Proposition 3.6.8(ii)) Thuillier shows that for a compact subset $\BerkK \subset \Omega$,
\index{Thuillier, Amaury} 
one has $C(\BerkK,\Omega) = C^*(\BerkK,\Omega)$, and he deduces (\ref{ThCapLim111})
for descending sequences of compact sets from this.  
His proof uses that $\BerkK = \bigcap_{n=1}^{\infty} \BerkU_n$ 
for a decreasing sequence of strict open affinoids 
$\BerkU_1 \supseteq \BerkU_2 \supseteq \cdots \supseteq \BerkK$
(see the first line on (\cite{Th}, p.112)).   For a compact set on a Berkovich curve 
\index{Berkovich!curve} 
over an arbitrary complete valued field $k$ this is not does not always hold 
(for example, if $k$ has an uncountable residue field, 
take $\BerkK$ to consist of a single type II point) but it does hold when $k = \CC_v$.  
This is because $\CC_v$ has a countable dense set, hence so do $\cC_v(\CC_v)$ and $\cC_v^{\an}$.  
Consequently if $\BerkK \subseteq \cC_v^{\an}$ is compact, then $\cC_v^{\an} \backslash \BerkK$
has countably many components.  By (\cite{Th}, Proposition 2.2.3), each component can be exhausted
by an increasing sequence of strict closed affinoids.  Using a diagonalization argument, one sees   
that $\BerkK$ is the intersection of a decreasing sequence of strict open affinoids.     

Next suppose $\zeta \in \cC(\CC_v) \backslash \BerkE_v$.  
Let $t$ be a tangent vector at $\zeta$,    
and let $\ThCap_{y,t}(\BerkK) = \|t\|_{\BerkK}^c$ be the function  
defined in (\cite{Th}, Corollary 3.6.19). By (\cite{Th}, Th\'eor\`eme 3.6.20), $\ThCap_{y,t}(\BerkK)$ 
induces a Choquet capacity on subsets of $\cC_v^{\an} \backslash \{y\}$.   
\index{capacity!Choquet capacity} 
Thus if $\BerkK_1 \subseteq \BerkK_2 \subseteq \cdots \subseteq \BerkK_n 
\cdots \subseteq \BerkE_v$
and $\bigcup_{n =1}^{\infty} \BerkK_n = \BerkE_v$, 
or $\BerkK_1 \supseteq \BerkK_2 \supseteq \cdots \supseteq \BerkK_n \cdots \supseteq \BerkE_v$
and $\bigcap_{n=1}^\infty \BerkK_n = \BerkE_v$, then    
\begin{equation} \label{ThCapLim222}
\lim_{n \rightarrow \infty} \ThCap_{y,t}(\BerkK_n) \ = \ \ThCap_{y,t}(\BerkE_v) \ .
\end{equation}  
Now fix a uniformizing parameter\index{uniformizing parameter!normalizes Robin constant} $g_{\zeta}(z)$, and choose $t$ 
so that $\langle t,g_{\zeta} \rangle = 1$.  By the discussion on (\cite{Th}, p.175), 
for each nonpolar compact $\BerkK \subset \cC_v^{\an} \backslash \BerkE_v$,   
\index{nonpolar set}  
\index{Robin constant!Berkovich!properties of}
\begin{equation*} 
V_{\zeta}(\BerkK)^{\an} \ = \ -\log(\ThCap_{y,t}(\BerkK)) \ .
\end{equation*} 
This yields (\ref{VDzzLim}) and (\ref{VAzzLim}).  

\smallskip
We next prove a special case of (D).  Suppose 
$\BerkK_1 \supseteq \BerkK_2 \supseteq \cdots \supseteq \BerkK_n \cdots \supseteq \BerkE_v$
is a descending sequence of compact sets with $\bigcap_{n=1}^{\infty} \BerkK_n = \BerkE_v$, 
and that in addition $\partial \BerkK_n \subset I(\cC_v^{\an})$ for each $n$.
Fix $x, y \in \cC_v^{\an} \backslash \BerkE_v$ with $x \ne y$.  After omitting finitely many $\BerkK_n$ 
from the sequence, we can assume that $x, y \notin \BerkK_1$.  

If $x$ and $y$ belong to distinct 
components of $\cC_v^{\an} \backslash \BerkE_v$, they belong to distinct components of 
$\cC_v^{\an} \backslash \BerkK_n$ for all $n$, so $G(x,y;\BerkK_n)^{\an} = G(x,y;\BerkE_v)= 0$
for all $n$, and the result is trivial.  Assume they belong to the same component $U$ of  
$\cC_v^{\an} \backslash \BerkE_v$.  For all sufficiently large $n$, they belong to the same 
component $U_n$ of $\cC_v^{\an} \backslash \BerkK_n$, and without loss we can assume they
belong to $U_n$ for all $n$.  For each $n$, put 
\begin{equation} \label{Harmonic1} 
h_n(z) \ = \ G(z,y;\BerkE_v)^{\an} - G(z,y;\BerkK_n)^{\an} \ ,
\end{equation} 
\index{Robin constant!Berkovich!properties of}
taking $h_n(y) = \lim_{z \rightarrow y} G(z,y;\BerkE_v)^{\an} - G(z,y;\BerkK_n)^{\an}  
= V_y(\BerkE_v)^{\an} - V_y(\BerkK_n)^{\an}$ if $y \in \cC_v(\CC_v)$.  Then $h_n(z)$ 
is harmonic on $U_n$ in the sense of (\cite{Th}, \S2.3).  
Note that $\partial U_n \subset \partial \BerkK_n \subset I(\cC_v^{\an})$. 
By (\cite{Th}, Propositions 3.1.19 and 3.1.20),
$G(z,y;\BerkK_n)$ is continuous at each point of $\partial U_n$ and vanishes on $\partial U_n$, 
so for each $p \in \partial U_n$ 
\begin{equation*}
\liminf_{\substack{z \rightarrow p \\  z \in U_n}} h_n(z) \ \ge \ 0 \ .
\end{equation*}
Hence the maximum principle for harmonic functions (\cite{Th}, Proposition 3.1.1)  
shows that $h_n(z) \ge 0$ for all $z \in U_n$, and in particular for all $z \in U_1$.  
Similarly, we see that $h_1(z) \ge h_2(z) \ge \cdots \ge 0$ for all $z \in U_1$.  
 
To conclude the argument, we apply Harnack's Principle (\cite{Th}, Proposition 3.1.2) 
\index{Harnack's Principle!Berkovich Harnack's Principle}
to the functions $h_n(z)$ on $U_1$.  
By the diagonal case of $(D)$ shown above, we have $\lim_{n \rightarrow \infty} h_n(y) = 0$.  
It follows from Harnack's Principle that the $h_n(z)$ converge uniformly to $0$ 
\index{Harnack's Principle!Berkovich Harnack's Principle}
on compact subsets of $U_1$, and consequently 
\begin{equation} \label{ThGreenLimit} 
\lim_{n \rightarrow \infty} G(x,y;\BerkK_n)^{\an} \ = \ G(x,y;\BerkE_v)^{\an} \ .
\end{equation} 

\smallskip
We can now prove (B), the symmetry of $G(x,y;\BerkE_v)$.  
Fix $x, y \in \cC_v^{\an} \backslash \BerkE_v$ with $x \ne y$.  If $x$ and $y$ belong 
to distinct components of $\cC_v^{\an} \backslash \BerkE_v$ then trivially 
$G(x,y;\BerkE_v)^{\an} = 0 = G(y,x;\BerkE_v)^{\an}$, so we can assume they belong to the same
component $U$.  Put $\BerkE_v^{U} = \cC_v^{\an} \backslash U$.  The characterization of 
Green's functions in (\cite{Th}, Th\'eor\`eme 3.6.15)
\index{Green's function!Berkovich!characterization of}
shows that $G(z,y;\BerkE_v^{U})^{\an} = G(z,y;\BerkE_v)^{\an}$ for all $z \in U$;  
in particular $G(x,y;\BerkE_v^{U})^{\an} = G(x,y;\BerkE_v)^{\an}$.

By (\cite{Th}, Proposition 2.2.23), there is an exhaustion of $U$ by an increasing sequence of 
strict open affinoids domains $V_1 \subset V_2 \subset \cdots \subset U$ 
\index{affinoid!Berkovich affinoid} 
with $\partial V_n \subset V_{n+1}$ for each $n$;  without loss, we can assume that $x, y \in V_1$.
By the definition of an affinoid, $\partial V_n \subset I(\cC_v^{\an})$ for each $n$.
\index{affinoid!Berkovich affinoid}
Let $g_y^{V_n}(z)$ be the Green's function of the domain $V_n$ 
\index{Green's function!Thuillier}
defined in (\cite{Th}, Proposition 3.3.7(ii)).  Then $g_y^{V_n}(z)$
is smooth, vanishes on $\partial V_n$, and satisfies the distributional equation 
$dd^c g_y^{V_n} = \mu_y^{V_n} - \delta_y$ where $\mu_y^{V_n}$ is a probability measure supported on 
$\partial V_n$.  Put $\BerkK_n = \cC_v^{an} \backslash V_n$. 
Again by the characterization of Green's functions in (\cite{Th}, Th\'eor\`eme 3.6.15),
\index{Green's function!Berkovich!characterization of}
for all $z \in V_n$ we have $G(z,y;\BerkK_n)^{an} = g_y^{V_n}(z)$. 

Here $\BerkK_1 \supset \BerkK_2 \supset \cdots \supset \BerkE_v^U$
is a descending sequence of compact sets with $\bigcap_{n=1}^{\infty} \BerkK_n = \BerkE_v^U$,
and $\partial \BerkK_n = \partial V_n \subset I(\cC_v^{\an})$ for each $n$.
By the special case of (E) shown above, we have 
\begin{equation} \label{ThGreenLimit3}  
\lim_{n \rightarrow \infty} G(x,y;\BerkK_n)^{\an} \ = \ G(x,y;\BerkE_v^{U})^{\an} 
\ = \ G(x,y;\BerkE_v)^{\an} \ .
\end{equation}
A similar formula holds with $x$ and $y$ interchanged.
By (\cite{Th}, Corollary 3.3.9(i)), for each $n$ we have $g_x^{V_n}(y) = g_y^{V_n}(x)$. 
\index{Green's function!Thuillier}
Combining these facts shows that $G(x,y;\BerkE_v)^{\an} = G(y,x;\BerkE_v)^{\an}$. 

\smallskip
Part (C), the monotonicity of $G(x,y,\BerkE_v)^{\an}$,  follows by a related argument.
Let $\BerkE_{v,1} \subseteq \BerkE_{v,2}$ be nonpolar, proper compact sets of $\cC_v^{\an}$.
\index{nonpolar set}   
Fix $x, y \in \cC_v^{\an}$ with $x \ne y$ and $y \in \cC_v^{\an} \backslash \BerkE_{v,2}$, 
and let $U$ be the component of $\cC_v^{\an} \backslash \BerkE_{v,2}$ containing $y$.

First suppose $x \in U$.  Putting $\BerkE_{v,2}^U = \cC_v^{\an} \backslash \BerkE_{v,2}$, 
let $\BerkK_1 \supseteq \BerkK_2 \supseteq \cdots \supseteq \BerkE_{v,2}^U$
be a descending sequence of compact sets with $\bigcap_{n=1}^{\infty} \BerkK_n = \BerkE_{v,2}^U$ 
and $\partial \BerkK_n = \partial V_n \subset I(\cC_v^{\an})$ for each $n$.  
By (\ref{ThGreenLimit3}) we have  
\begin{equation*}
\lim_{n \rightarrow \infty} G(x,y;\BerkK_n)^{\an} \ = \ G(x,y;\BerkE_{v,2})^{\an} \ .
\end{equation*} 
On the other hand, the same argument that gave (\ref{Harmonic1}) shows that
\begin{equation*}
G(x,y;\BerkE_{v,1}) \ \ge \ G(x,y;\BerkK_n)^{\an} 
\end{equation*}
for each $n$.  Thus $G(x,y;\BerkE_{v,1})^{\an} \ \ge \ G(x,y;\BerkE_{v,2})^{\an}$.

Now suppose $x \notin U$. 
By the characterization of Green's functions in (\cite{Th}, Th\'eor\`eme 3.6.15)
\index{Green's function!Berkovich!characterization of}
we have $G(z,y;\BerkE_v)^{\an} = G(z,y;\partial U)^{\an}$ for all $z$.    
First assume $x \notin \partial U$. 
Then $x$ and $y$ belong to distinct components of $\cC_v^{\an} \backslash \partial U$, 
and trivially $G(x,y;\BerkE_{v,1})^{\an} \ge 0 = G(x,y;\partial U) = G(x,y;\BerkE_{v,2})^{\an}$. 
Last, assume $x \in \partial U$.  Since $G(z,y;E_{v,1})$ and $G(z,y;E_{v,2})$ are subharmonic,
necessarily they are upper semi-continuous (see \cite{Th}, D\'efinition 3.1.5). 
\index{semi-continuous!Green's function is upper semi-continuous}  
By what has been shown above, 
\begin{equation*}
G(x,y;\BerkE_{v,2}) \ \ge \ \lim_{\substack{z \rightarrow x \\ z \in U}} G(z,y;\BerkE_{v,2}) 
\ \ge \ \lim_{\substack{z \rightarrow x \\ z \in U}} G(z,y;\BerkE_{v,1}) \ = \ G(x,y;\BerkE_{v,1}) \ .
\end{equation*}  
This establishes the desired inequality in all cases. 

Finally, consider the Robin constants.  If $\zeta \in \cC_v^{\an} \backslash \BerkE_{v,2}$, it follows that 
\index{Robin constant!Berkovich}
\index{Robin constant!Berkovich!properties of}
\begin{eqnarray*}
V_{\zeta}(\BerkE_{v,1})^{\an} 
& = & \lim_{z \rightarrow \zeta} G(x,\zeta;\BerkE_{v,1})^{\an} + \log(|g_{\zeta}(z)|_v) \\
& \ge & \lim_{z \rightarrow \zeta} G(x,\zeta;\BerkE_{v,2})^{\an} + \log(|g_{\zeta}(z)|_v) 
\ = \ V_{\zeta}(\BerkE_{v,2})^{\an} \ .
\end{eqnarray*}

\smallskip
We can now prove (D) in full generality.  
Since we have already established the diagonal case,
we only consider the non-diagonal case.  
Let $\{\BerkK_n\}_{n \ge 1}$ be a sequence of compact sets with 
$\BerkK_1 \supseteq \BerkK_2 \supseteq \cdots \supseteq \BerkK_n \cdots \supseteq \BerkE_v$ such 
that $\bigcap_{n=1}^{\infty} \BerkK_n = \BerkE_v$, 
and fix $x, y \in \cC_v^{\an} \backslash \BerkE_v$ with $x \ne y$. 
For all sufficiently large $n$ we have $x, y \notin \BerkK_n$, so after omitting finitely many 
$\BerkK_n$ we can assume without loss that $x, y \notin \BerkK_1$.

Let $U$ be the component of $\cC_v^{\an} \backslash \BerkE_v$ containing $y$. 
If $x \notin U$, then $x$ and $y$ belong to distinct components of $K_n$ for all $n$,   
so $G(x,y;\BerkK_n) = 0 = G(x,y;\BerkE_v)$ for all $n$, and (\ref{GDxyLim}) is trivial. 
Suppose $x \in U$.  After omitting finitely many $K_n$, we can assume that $x$ and $y$ belong
to the same component $U_1$ of $\cC_v^{\an} \backslash \BerkK_1$.    
For each $n$, put $h_n(z) = G(z,y;E_v) - G(z,y;\BerkK_n)$, 
taking $h_n(y) = V_y(\BerkK_n) - V_y(\BerkE_v)$ if $y \in \cC_v(\CC_v)$.
\index{Robin constant!Berkovich!properties of}
Then $h_n(z)$ is harmonic in $U_1$, and by part (C), it is non-negative.  
By the diagonal case of (D) shown above, we have $\lim_{n \rightarrow \infty} h_n(y) = 0$.  
Hence Harnack's Principle (\cite{Th}, Proposition 3.1.2)
\index{Harnack's Principle!Berkovich Harnack's Principle}
gives that as $n \rightarrow \infty$, then $h_n(z) \rightarrow 0$ uniformly on compact subsets of $U_1$. 
In particular  
\begin{equation*} 
\lim_{n \rightarrow \infty} G(x,y;\BerkK_n)^{\an} \ = \ G(x,y;\BerkE_v)^{an} \ .
\end{equation*}  

\smallskip
The nondiagonal case of (E) follows by a similar argument, 
but uses the topology of $\cC_v^{\an}$ in a stronger way.
Let $\{\BerkK_n\}_{n \ge 1}$ be a sequence of compact sets with 
$\BerkK_1 \subseteq \BerkK_2 \subseteq \cdots \subseteq \BerkK_n \cdots \subseteq \BerkE_v$ such 
that $\bigcup_{n=1}^{\infty} \BerkK_n = \BerkE_v$, 
and fix $x, y \in \cC_v^{\an} \backslash \BerkE_v$ with $x \ne y$. 

Let $\Gamma_{x,y}$ be the union of all paths connecting $x$ and $y$ in $\cC_v^{\an}$.    
Recall that there is a finite subgraph $\BerkS$ of $\cC_v^{\an}$, 
called its {\em skeleton},\index{skeleton of a Berkovich curve}  
such that there is retraction $\tau : \cC_v^{\an} \rightarrow \BerkS$ 
(see (\cite{Th}, Th\'eor\`eme 2.210)).  Each component of $\cC_v^{\an} \backslash \BerkS$ is a tree.  
It follows that $\Gamma_{x,y}$ is a graph with finitely many edges.  

If $x$ and $y$ belong to distinct components of $\cC_v^{\an} \backslash \BerkE_v$, 
they belong to distinct components of $\Gamma_{x,y} \backslash \BerkE_v$.  
Since $\Gamma_{x,y}$ has finite connectivity, 
there is a finite subset $P \subset \BerkE_v \cap \Gamma_{x,y}$ which disconnects $x$ from $y$. 
For all sufficiently large $n$ we have $P \subset \BerkK_n$, and for such $n$ it follows that
$G(x,y;\BerkK_n) = 0 = G(x,y;\BerkE_v)$, so (\ref{GAxyLim}) is trivial.  

Suppose $x$ and $y$ belong to the same component $U$ of $\cC_v^{\an} \backslash \BerkE_v$.  
For each $n$, put $h_n(z) = G(z,y;K_n) - G(z,y;\BerkE_v)$, 
taking $h_n(y) = V_y(\BerkE_v) - V_y(\BerkK_n)$ if $y \in \cC_v(\CC_v)$.  
\index{Robin constant!Berkovich!properties of}
Then $h_n(z)$ is harmonic in $U$, and by part (C), it is non-negative.  
By the diagonal case of (E) shown above, we have $\lim_{n \rightarrow \infty} h_n(y) = 0$.  
Hence Harnack's Principle gives that as $n \rightarrow \infty$,
\index{Harnack's Principle!Berkovich Harnack's Principle} 
then $h_n(z) \rightarrow 0$ uniformly on compact subsets of $U$, and again we conclude that  
\begin{equation*} 
\lim_{n \rightarrow \infty} G(x,y;\BerkK_n)^{\an} \ = \ G(x,y;\BerkE_v)^{an} \ .
\end{equation*}  
  
\smallskip
Part (F), the functoriality of $G(x,y;\BerkE_v)^{\an}$ and $V_{\zeta}(\BerkE_v)^{\an}$ 
\index{Robin constant!Berkovich!properties of}
under $\Aut_c(\CC_v/K_v)$, is immediate from the definition of the action of 
$\sigma \in \Aut_c(\CC_v/K_v)$ on $\cC_v^{\an}$ through its action on $\CC_v$, 
and the characterization of $g_{y,\BerkE_v}(z)$ in (\cite{Th}, Th\'eor\`eme 3.6.15). 
\end{proof}

The following proposition establishes the compatibility of our Green's functions 
and Berkovich Green's functions.  
\index{Green's function} 
\index{Berkovich!Green's function} 
\index{Green's function!Berkovich}
Note that in (\cite{Th}, Th\'eor\`eme 5.1.2), Thuillier has already established the compatibility
\index{Thuillier, Amaury}
of our capacities with his. Throughout (\cite{Th}) Thuillier uses the natural logarithm $\log(x)$, 
\index{Thuillier, Amaury}
while in $v$-adic constructions we use the logarithm $\log_v(x)$ to the base $q_v$ 
in order to have $\log_v(|z|_v) \in \QQ$ for $z \in \CC_v^{\times}$.  This gives rise to 
a factor of $\log(q_v)$ in comparing our Green's functions and his.

\begin{proposition}[Compatibility of Green's Functions] \label{BerkCompatibilityProp}
\index{Green's function!Berkovich!compatible with classical|ii}
Let $K$ be a global field, and let $\cC/K$ be a smooth, connected, projective curve.  
Let $v$ be a nonarchimedean place of $K$, and let $\cC_v^{\an}$ be the Berkovich analytification
\index{Berkovich!analytification}  
of $\cC_v \times_{K_v} \Spec(\CC_v)$.  

Suppose $E_v \subsetneq \cC_v(\CC_v)$ is an algebraically capacitable set with positive capacity,
\index{algebraically capacitable} 
\index{capacity $> 0$} 
and $\BerkE_v$ is its closure in $\cC_v^{\an}$.  Then $\BerkE_v$ is a proper compact, 
nonpolar subset of $\cC_v^{\an}$, 
\index{nonpolar set} 
and for all $z, \zeta \in \cC_v(\CC_v) \backslash E_v$, 
\begin{equation*} 
G(z,\zeta;\BerkE_v)^{\an} \ = \ G(z,\zeta;E_v) \log(q_v)\ , \qquad 
V_{\zeta}(\BerkE_v)^{\an} \ = \ V_{\zeta}(E_v) \log(q_v) \ .
\end{equation*}
\index{Robin constant!Berkovich!compatible with classical|ii} 
\end{proposition} 
 
\begin{proof}  We begin by considering two special cases:  compact sets and $\PL_{\zeta}$-domains.     

First, let $E_v \subset \cC_v(\CC_v)$ be a compact set with positive capacity.  
\index{capacity $> 0$} 
Since $E_v$ is compact and the restriction of the topology on $\cC_v^{\an}$ to $\cC_v(\CC_v)$
is the usual $v$-adic topology, $E_v$ coincides with its Berkovich closure $\BerkE_v$.
\index{Berkovich!closure} 
Fix a point $\zeta \in \cC_v(\CC_v) \backslash E_v$ and a 
uniformizing parameter\index{uniformizing parameter!normalizes canonical distance} $g_\zeta(z)$, 
and let $[z,w]_{\zeta}$ be the canonical distance  
\index{canonical distance!$[z,w]_{\zeta}$}  
normalized so that $\lim_{z \rightarrow \zeta} [z,w]_{\zeta} \cdot |g_\zeta(z)|_v = 1$ 
for each $w \ne \zeta$ (see \S\ref{Chap3}.\ref{CanonicalDistanceSection}).
By definition, our Robin constant is 
\index{Robin constant!of compact set}
\begin{equation*}
V_{\zeta}(E_v) \ = \ \inf_{\nu} \iint -\log_v([x,y]_{\zeta}) \, d\nu(x) d\nu(y)
\end{equation*} 
and our capacity is 
\begin{equation*}
\gamma_{\zeta}(E_v) \ = \ q_v^{-V_{\zeta}(E_v)} \ .
\end{equation*}
Let $\mu_{\zeta}$\index{equilibrium distribution} 
be the equilibrium distribution of $E_v$ with respect to $\zeta$:  the unique probability measure
on $E_v$ which minimizes the energy integral\index{energy integral} 
$I_{\zeta}(\nu) = \inf_{\nu} \iint -\log_v([x,y]_{\zeta}) \, d\nu(x) d\nu(y)$  
(see \S\ref{Chap3}.\ref{CompactGreenSection}).
Then the potential function\index{potential function!takes constant value a.e. on $E_v$}  
\begin{equation*}
u_{E_v}(z,\zeta) \ = \ \int -\log_v([z,w]_{\zeta}) \, d\mu_{\zeta}(w) 
\end{equation*} 
satisfies $u_{E_v}(z,\zeta) \le V_{\zeta}(E_v)$ for all $z \in \cC_v(\CC_v) \backslash \{\zeta\}$, 
takes the value $V_{\zeta}(E_v)$ on $E_v \backslash e_v$ 
where $e_v$ is an $F$-sigma set of inner capacity $0$, 
\index{capacity $= 0$} 
and has $u_{E_v}(z,\zeta) < V_{\zeta}(E_v)$ for all $z \in \cC_v(\CC_v) \backslash E_v$.  
In (\ref{FMUD1}) we have defined 
\begin{equation*}
G(z,\zeta;E_v) \ = \ V_{\zeta}(z) - u_{E_v}(z,\zeta) 
\ = \ V_{\zeta}(E_v) + \int \log_v([z,w]_{\zeta}) \, d\mu_{\zeta}(w)  \ .
\end{equation*} 

In (\cite{Th}, Th\'eor\`eme 5.1.5), 
Thuillier shows there is a unique extension of the canonical distance to a function 
\index{canonical distance!$[z,w]_{\zeta}$} 
\index{canonical distance!$[z,w]_{\zeta}^{\an}$|ii} 
\index{Thuillier, Amaury}
on $(\cC_v^{\an} \backslash \{\zeta\}) \times (\cC_v^{\an} \backslash \{\zeta\})$, 
which we will denote $[z,w]_{\zeta}^{\an}$.  
The function $[z,w]_{\zeta}^{\an}$ is continuous, symmetric, 
and satisfies  $dd^c \log([z,w]_{\zeta}^{\an}) = \delta_w - \delta_\zeta$   
for each $w \ne \zeta$.  Noting that $\log(x) = \log_v(x) \log(q_v)$, define 
\begin{equation*}
g_{\mu_\zeta}^{\an}(z) \ = \ V_{\zeta}(E_v) \log(q_v) + \int \log([z,w]_{\zeta}^{\an}) \, d\mu_{\zeta}(w)  
\end{equation*}
\index{Green's function!Thuillier}
for $z \in \cC_v^{\an} \backslash \{\zeta\}$.  By arguments like those in (\cite{Th},Proposition 3.4.16), 
$g_{\mu_\zeta}^{\an}(z)$ is subharmonic on $\cC_v^{\an} \backslash \{\zeta\}$,
harmonic on $\cC_v^{\an} \backslash (E_v \cup \{\zeta\})$, 
and satisfies $dd^c g_{\mu_\zeta}^{\an} = \mu_{\zeta} - \delta_{\zeta}$.  
Clearly $g_{\mu_{\zeta}}^{\an}(z) = G(z,\zeta;E_v) \log(q_v)$ 
for $z \in \cC_v(\CC_v) \backslash \{\zeta\}$.  In particular, it vanishes on $\BerkE_v = E_v$
except possibly on the set $e_v$. By (\cite{Th}, Th\'eor\`eme 3.6.11 and Th\'eor\`eme 5.1.2),
the set $e_v$ is polar.\index{polar set}  The characterization of Green's functions in (\cite{Th}, Th\'eor\`eme 3.6.15),
\index{Green's function!Berkovich!characterization of}
shows that $G(z,\zeta;\BerkE_v)^{an} = g_{\zeta,\BerkE_v}(z) = g_{\mu_{\zeta}}^{\an}(z)$.
Thus for all $z \in \cC_v(\CC_v) \backslash \{\zeta\}$,
\begin{equation*} 
G(z,\zeta;\BerkE_v)^{\an} \ = \ G(z,\zeta;E_v) \log(q_v) \ .
\end{equation*}
It follows that $V_{\zeta}(\BerkE_v)^{\an} \ = \ V_{\zeta}(E_v) \log(q_v)$.  
\index{Robin constant!Berkovich!compatible with classical} 

\smallskip
Next, let $E_v \subsetneq \cC_v(\CC_v)$ be a $\PL_{\zeta}$-domain 
in the sense of (\cite{RR1}, Definition 4.2.6):
there is a nonconstant $f(z) \in \CC_v(\cC)$ 
having poles only at $\zeta$, such that $E_v = \{z \in \cC_v(\CC_v) : |f(z)|_v \le 1\}$.  
Given a function $f(z)$ defining $E_v$ as a $\PL_{\zeta}$-domain,  
\begin{equation*}
G(z,\zeta;E_v) \ = \ \left\{ \begin{array}{ll} 
                 \frac{1}{\deg}(f) \log_v(|f(z)|_v) & \text{if $z \in \cC_v(\CC) \backslash E_v$ \ ,} \\
                 0 & \text{if $z \in E_v$ \ ,} 
                              \end{array} \right.
\end{equation*} 
where $\log_v(x)$ is the logarithm to the base $q_v$;  
by (\cite{RR1}, Proposition 4.4.1), $G(z,\zeta;E_v)$ is independent of the choice of $f$.   
The closure of $E_v$ in $\cC_v^{\an}$ is  $\BerkE_v = \{z \in \cC_v^{\an} : |f(z)|_v \le 1\}$ 
and by the discussion on (\cite{Th}, p.175) 
\begin{equation*}
G(z,\zeta;\BerkE_v)^{\an} \ = \ \left\{ \begin{array}{ll} 
                 \frac{1}{\deg}(f) \log(|f(z)|_v) & \text{if $z \in \cC_v^{\an} \backslash \BerkE_v$ \ ,} \\
                 0 & \text{if $z \in \BerkE_v$ \ ,} 
                            \end{array} \right.
\end{equation*} 
where $\log(x) = \ln(x)$.  It follows that for all $z \in \cC_v(\CC_v) \backslash \{\zeta\}$ 
\begin{equation} \label{GCompatibility} 
G(z,\zeta;\BerkE_v)^{\an} \ = \ G(z,\zeta;E_v) \log(q_v) \ .
\end{equation} 
\index{Green's function!Berkovich!compatible with classical} 
and that 
\begin{equation} \label{VCompatibility}
V_{\zeta}(\BerkE_v)^{\an} \ = \ V_{\zeta}(E_v) \log(q_v) \ ,
\end{equation}
\index{Robin constant!Berkovich!compatible with classical}  

\smallskip
We can now deal with the general case.  Let $E_v \subsetneq \cC_v(\CC_v)$ be an algebraically 
capacitable set with positive capacity, and let $\BerkE_v$ be its closure in $\cC_v^{\an}$.  
\index{capacity $> 0$} 
Note that $E_v$ is closed in $\cC_v(\CC_v)$ by (\cite{RR1}, Proposition 4.3.15).
Since the topology on $\cC_v^{\an}$ restricts to $v$-adic topology on $\cC_v(\CC_v)$,
this implies that $\BerkE_v \cap \cC_v(\CC_v) = E_v$.  In particular, $\BerkE_v$ 
is a proper subset of $\cC_v^{\an}$.  It is clearly compact, and it is nonpolar since $E_v$ 
\index{nonpolar set} 
contains compact subsets of $\cC_v(\CC_v)$ with positive capacity.  
\index{capacity $> 0$}

Fix $\zeta \in \cC_v(\CC_v) \backslash E_v$.  By (\cite{RR1}, Definition 4.3.2) 
we have 
\begin{equation*}
\inf_{\text{compact $K \subseteq E_v$}} V_{\zeta}(K) = V_{\zeta}(E_v) \ , \qquad 
\sup_{\text{$\PL_{\zeta}$-domains $U \supseteq E_v$}} V_{\zeta}(U) = V_{\zeta}(E_v) \ .
\end{equation*} 
\index{Robin constant!Berkovich!compatible with classical} 
Since the union of finitely compact sets is compact, 
and the intersection of finitely many $\PL_{\zeta}$-domains is a $\PL_{\zeta}$-domain 
(see (\cite{RR1}, Corollary 4.2.13)), there are an ascending sequence of compact sets 
$K_1 \subseteq K_2 \subseteq \cdots \subseteq E_v$ with 
\begin{equation} \label{VzLimit1}
\lim_{n \rightarrow \infty} V_{\zeta}(K_n) \ = \ V_{\zeta}(E_v) \ ,
\end{equation}  
and a descending sequence of $\PL_{\zeta}$-domains $U_1 \supseteq U_2 \supseteq \cdots \supseteq E_v$ with 
\begin{equation} \label{VzLimit22}
\lim_{n \rightarrow \infty} V_{\zeta}(U_n) \ = \ V_{\zeta}(E_v) \ .
\end{equation} 
\index{Robin constant!Berkovich!compatible with classical} 
By (\cite{RR1}, Lemma 4.4.7 and Definition 4.4.12), 
for each $z \in \cC_v(\CC_v) \backslash (E_v \cup \{\zeta\})$ we have 
\begin{equation} \label{GnLimit1}  
\lim_{n \rightarrow \infty} G(z,\zeta;K_n) = G(z,\zeta;E_v) \ , \qquad 
\lim_{n \rightarrow \infty} G(z,\zeta;U_n) = G(z,\zeta;E_v) \ .
\end{equation} 

By the compatibility of Green's functions for compact sets and $\PL_{\zeta}$-domains shown above,
\index{Green's function!Berkovich!compatible with classical}
if $\BerkK_n$ and $\BerkU_n$ are the closures of $K_n$ and $U_n$ in $\cC_v^{\an}$ respectively, 
then for all $z \in \cC_v(\CC_v) \backslash \{\zeta\}$ and all $n$,  
\begin{equation} \label{GnEquality2}
G(z,\zeta;\BerkK_n)^{\an} = G(z,\zeta;K_n) \log(q_v) \ , 
\qquad \ G(z,\zeta;\BerkU_n)^{\an} = G(z,\zeta;\BerkU_n) \log(q_v) \ .
\end{equation} 
Clearly 
\begin{equation*}
\BerkK_1 \ \subseteq \ \BerkK_2 \ \subseteq \ \cdots \ \subseteq \ \BerkE_v \ 
\subseteq \cdots \ \subseteq \BerkU_2 \ \subseteq \ \BerkU_1 \ ,
\end{equation*}
so by the monotonicity of Green's functions proved in Proposition \ref{BerkGreenPropertiesProp}(C), 
\index{Green's function!Berkovich!monotonic}
for all $z \in \cC_v^{\an} \backslash \{\zeta\}$
\begin{eqnarray} 
G(z,\zeta;\BerkK_1)^{\an} & \ge & G(z,\zeta;\BerkK_2)^{\an} \ \ge \ \cdots \notag \\
& \ge & G(z,\zeta;\BerkE_v)^{\an} 
\ \ge \cdots \ \ge \ G(z,\zeta;\BerkU_2)^{\an} \ \ge \ G(z,\zeta;\BerkU_1)^{\an} \ . \label{GnInequality3}
\end{eqnarray} 
Combining (\ref{GnLimit1}), (\ref{GnEquality2}) and (\ref{GnInequality3}) 
shows that for each $z \in \cC_v(\CC_v) \backslash (E_v \cup \{\zeta\})$ we have 
\begin{equation*} 
G(z,\zeta;\BerkE_v)^{\an} \ = \ G(z,\zeta;E_v) \log(q_v) \ .
\end{equation*} 
In a similar way, from (\ref{VzLimit1}), (\ref{VzLimit22}), 
the compatibility of Robin constants for compact sets and $\PL_{\zeta}$ domains, 
\index{Robin constant!Berkovich!monotonicity of}
and the monotonicity of Robin constants proved in Proposition \ref{BerkGreenPropertiesProp}(C),
\index{Green's function!Berkovich!monotonic}
we see that 
\begin{equation*}
V_{\zeta}(\BerkE_v)^{\an} \ = \ V_{\zeta}(E_v) \log(q_v) \ . 
\end{equation*}
\end{proof}

We can now prove Theorems \ref{aT1-B1} and \ref{aT1-B2}. 

\smallskip
\noindent{{\bf Theorem \text{\ref{aT1-B1}}.}
(Berkovich FSZ with LRC, producing points in $\EE$) 
\index{Fekete-Szeg\"o theorem with LRC!Berkovich|ii}

{\it 
Let $K$ be a global field, 
and let $\cC/K$ be a smooth, geometrically integral, projective curve.
Let $\fX = \{x_1, \ldots, x_m\} \subset \cC(\tK)$ 
be a finite set of points stable under $\Aut(\tK/K)$, and let
$\EE = \prod_v \BerkE_v \subset \prod_v \cC_v^{\an}$ be a $K$-rational Berkovich adelic 
set compatible with $\fX$.\index{compatible with $\fX$!Berkovich set compatible with $\fX$}  
Let $S \subset \cM_K$ be a finite set of places $v$, containing all archimedean $v$,
such that $\BerkE_v$ is $\fX$-trivial for each $v \notin S$.
\index{$\fX$-trivial}

Assume that $\gamma(\EE,\fX) > 1$.  
Assume also that $\BerkE_v$ has the following form, for each $v \in S$:   

$(A)$ If $v$ is archimedean and $K_v \cong \CC$, 
then $\BerkE_v$ is compact, and is a finite union of sets $E_{v,\ell}$, 
each of which is the closure of its $\cC_v(\CC)$-interior and has a piecewise smooth 
boundary;\index{boundary!piecewise smooth}\index{closure of $\cC_v(\CC)$ interior} 

$(B)$ If $v$ is archimedean and $K_v \cong \RR$, then $\BerkE_v$ is compact, 
stable under complex conjugation, 
and is a finite union of sets $E_{v,\ell}$, where each $E_{v,\ell}$ is either 

\quad $(1)$ the closure of its $\cC_v(\CC)$-interior and has a piecewise smooth boundary, or
\index{boundary!piecewise smooth}\index{closure of $\cC_v(\CC)$ interior}  

\quad $(2)$ is a compact, connected subset of $\cC_v(\RR)$; 

$(C)$ If $v$ is nonarchimedean, then $\BerkE_v$ is compact, stable under $\Aut_c(\CC_v/K_v)$, 
and is a finite union of sets $E_{v,\ell}$, where each $E_{v,\ell}$ is either 

\quad $(1)$ a strict closed Berkovich affinoid, or
\index{affinoid!Berkovich affinoid}
\index{Berkovich!affinoid} 

\quad $(2)$ is a compact subset of $\cC_v(\CC)$ 
and has the form $\cC_v(F_{w_\ell}) \cap B(a_\ell,r_\ell)$ 
for some finite separable extension $F_{w_\ell}/K_v$ in $\CC_v$, and some ball $B(a_\ell,r_\ell)$.  

Then there are infinitely many points $\alpha \in \cC(\tK^{\sep})$ such that for each $v \in \cM_K$, 
the $\Aut(\tK/K)$-conjugates of $\alpha$ all belong to $\BerkE_v$.  
}

\begin{proof}[Proof of Theorem \ref{aT1-B1}, using Theorem \ref{aT1}.] 

 For each $v \in \cM_K$, 
put $E_v^0 = \BerkE_v \cap \cC_v(\CC_v)$.  By the hypotheses of Theorem \ref{aT1-B1}, 
$E_v^0$ is algebraically capacitible and satisfies the hypotheses of Theorem \ref{aT1}.  
Those hypotheses in turn show that the Berkovich closure of $E_v^0$ is $\BerkE_v$, 
\index{Berkovich!closure} 
so by Proposition \ref{BerkCompatibilityProp} for all $x_i \ne x_j \in \fX$ we have 
\begin{equation*}
G(x_i,x_j;E_v^0)  \ = \  G(x_i,x_j;\BerkE_v)^{\an} \ , \qquad 
V_{x_i}(E_v^0) \ = \ V_{x_i}(\BerkE_v)^{\an} \ .
\end{equation*}   

Put $\EE^0 = \prod_v E_v^0$. 
Then $\EE^0$ is a $K$-rational adelic set compatible with $\fX$\index{compatible with $\fX$} 
and has $\gamma(\EE^0,\fX) > 1$. By Theorem \ref{aT1-B1} 
there are infinitely many points $\alpha \in \cC(\tK^{\sep})$ such that for each $v \in \cM_K$, 
the $\Aut(\tK/K)$-conjugates of $\alpha$ all belong to $E_v$, hence $\BerkE_v$.   
\end{proof} 

\noindent{{\bf Theorem \text{\ref{aT1-B2}}.} (Berkovich Fekete/FSZ with LRC for Quasi-neighborhoods).}
\index{quasi-neighborhood!Berkovich} 
\index{Fekete-Szeg\"o theorem with LRC!for Berkovich quasi-neighborhoods|ii}  
{\it  
Let $K$ be a global field, 
and let $\cC/K$ be a smooth, connected, projective curve.
Let $\fX = \{x_1, \ldots, x_m\} \subset \cC(\tK)$ 
be a finite set of points stable under $\Aut(\tK/K)$, and let 
$\EE = \prod_v \BerkE_v \subset \prod_v \cC_v^{\an}$
be a compact Berkovich adelic set 
compatible with $\fX$,\index{compatible with $\fX$!Berkovich set compatible with $\fX$} 
\index{Berkovich!adelic set} 
such that each $\BerkE_v$ is stable under $\Aut_c(\CC_v/K_v)$. 

$(A)$ If $\gamma(\EE,\fX)^{\an} < 1$, 
there is a $K$-rational Berkovich neighborhood $\UU = \prod_v \BerkU_v$ of $\EE$
\index{Berkovich!neighborhood} 
such that there are only finitely many points of $\cC(\tK)$ 
whose $\Aut(\tK/K)$-conjugates are all contained in $\BerkU_v$, 
for each $v \in \cM_K$.

$(B)$ If $\gamma(\EE,\fX)^{\an} > 1$,
then for any $K$-rational 
separable Berkovich quasi-neighborhood $\UU$ of $\EE$, there 
\index{quasi-neighborhood!Berkovich}
\index{Berkovich!quasi-neighborhood} 
are infinitely many points $\alpha \in \cC(\tK^{\sep})$ 
such that for each $v \in \cM_K$,
the $\Aut(\tK/K)$-conjugates of $\alpha$ all belong to $\BerkU_v$.  
}

\begin{proof}[Proof of Theorem \ref{aT1-B2}, using Theorem \ref{aT1-A1}.]
\ 

We first prove (A).  Suppose $\gamma(\EE,\fX)^{\an} < 1$.  
We begin by enlarging $\EE = \prod_v \BerkE_v$ 
to a set $\FF = \prod_v \BerkF_v$ with $\gamma(\FF,\fX) < 1$,  
such that $\BerkF_v$ is a strict closed affinoid for each nonarchimedean $v$.   Let $\varepsilon > 0$
\index{affinoid!strict closed affinoid} 
be small enough that if $\Gamma \in M_n(\RR)$ is a symmetric $n \times n$ matrix whose entries differ
from those of $\Gamma(\EE,\fX)^{\an}$ by at most $\varepsilon$, then $\val(\Gamma) > 1$.
\index{Green's matrix!global}\index{value of $\Gamma$ as a matrix game}
Fix a nonempty finite set of places $S$ of $K$  
containing all archimedean places and all nonarchimedean places where $\BerkE_v$ is not $\fX$-trivial, 
\index{$\fX$-trivial}
and choose a set of numbers $\{\varepsilon_v\}_{v \in S}$ with $\varepsilon_v > 0$ for each $v$ 
and $\sum_{v \in S} \varepsilon_v = \varepsilon$.   

If $v$ is archimedean, put $\BerkF_v = \BerkE_v$;  
likewise if $v \notin S$, so $\BerkE_v$ is $\fX$-trivial, put $\BerkF_v = \BerkE_v$.
\index{$\fX$-trivial}  
Suppose $v \in S$ is nonarchimedean.  
By hypothesis $\BerkE_v$ is compact, nonpolar, and stable under $\Aut_c(\CC_v/K_v)$.
\index{nonpolar set}   
As noted in the proof of Proposition \ref{BerkGreenPropertiesProp}, 
\index{Green's function!Berkovich!properties}
the fact that $\CC_v$ has a countable dense set means there is a descending sequence of 
strict closed affinoids $\BerkK_1 \supseteq \BerkK_2 \supseteq \cdots \supseteq \BerkE_v$  
with $\bigcap_{n=1}^{\infty} \BerkK_n = \BerkE_v$.  By Proposition \ref{BerkGreenPropertiesProp}(D),
if $n$ is large enough, then for all $x_i, x_j \in \fX$ with $i \ne j$ we have 
\begin{equation*} 
|G(x_i,x_j;\BerkE_v)^{\an} - G(x_i,x_j;\BerkK_n)^{\an}| \ < \ \varepsilon_v \ ,
\end{equation*} 
and for each $x_i \in \fX$ 
\begin{equation*} 
|V_{x_i}(\BerkE_v)^{\an} - V_{x_i}(\BerkK_n)^{\an}| \ < \ \varepsilon_v \ .
\end{equation*} 
\index{Robin constant!Berkovich} 
Fix such an $n$.  Since $\tK_v$ is dense in $\CC_v$, the strict closed affinoid $\BerkK_n$ can be 
\index{affinoid!strict closed affinoid}
defined by equations in $\tK_v$ and has only finitely many distinct conjugates under $\Aut_c(\CC_v/K_v)$.  
Put
\begin{equation*}
\BerkF_v \ = \ \bigcap_{\sigma \in \Aut_c(\CC_v/K_v)} \sigma(\BerkK_n) \ .
\end{equation*} 
Since the intersection of finitely many strict closed affinoids is again a strict closed affinoid, 
\index{affinoid!strict closed affinoid}
$\BerkF_v$ is a strict closed affinoid with $\BerkE_v \subseteq \BerkF_v \subseteq \BerkK_n$.
\index{affinoid!strict closed affinoid} 
By construction it is stable under $\Aut_c(\CC_v/K_v)$.   
The monotonicity of Green's functions in Proposition \ref{BerkGreenPropertiesProp}(C) shows that 
\index{Green's function!Berkovich!monotonic}
\begin{equation*} 
|G(x_i,x_j;\BerkE_v)^{\an} - G(x_i,x_j;\BerkF_v)^{\an}| \ < \ \varepsilon_v \ ,
\end{equation*} 
and for each $x_i \in \fX$ 
\begin{equation*} 
|V_{x_i}(\BerkE_v)^{\an} - V_{x_i}(\BerkF_v)^{\an}| \ < \ \varepsilon_v \ .
\end{equation*} 
\index{Robin constant!Berkovich} 

We now reduce to the classical case.  For each $v$, put $F_v^0 = \BerkF_v \cap \cC_v(\CC_v)$.
Thus, if $v \in S$ is archimedean, then $F_v^0 = \BerkF_v$;  if $v \in S$ is nonarchimedean, 
then $F_v^0$ is an $\RL$-domain whose closure in $\cC_v^{\an}$ is $\BerkF_v$, and if $v \notin S$
\index{$\RL$-domain} 
then $F_v^0$ is $\fX$-trivial and again its closure in $\cC_v^{\an}$ is $\BerkF_v$.  In particular,
\index{$\fX$-trivial}
each $F_v^0$ is algebraically capacitable and stable under $\Aut^c(\CC_v/K_v)$.
\index{algebraically capacitable}  
Set $\FF^0 = \prod_v F_v^0$.  By Proposition \ref{BerkCompatibilityProp}, the Green's matrices 
$\Gamma(\FF^0,\fX)$ and $\Gamma(\FF,\fX)^{\an}$ coincide.  
Our choice of $\varepsilon$ and the $\varepsilon_v$ shows that 
the entries of $\Gamma(\FF^0,\fX)$ differ from those of $\Gamma(\FF,\fX)^{\an}$
by at most $\varepsilon$.  Hence $\val(\FF^0,\fX) > 1$, and $\gamma(\FF^0,\fX) < 1$. 

By (\cite{RR1}, Theorem 6.2.1), there is a function $f(z) \in K(\cC)$ with poles supported on $\fX$, 
such that for each $v \in S$
\begin{equation*} 
F_v^0 \ \subset \ \{z \in \cC_v(\CC_v) :  |f(z)|_v < 1 \} \ ,
\end{equation*} 
and for each $v \notin S$ 
\begin{equation*} 
F_v^0 \ \subseteq \ \{z \in \cC_v(\CC_v) :  |f(z)|_v \le 1 \} \ .
\end{equation*}
For each $v \in S$, put $\BerkU_v = \{z \in \cC_v^{\an} :  |f(z)|_v < 1 \}$, 
and for each $v \notin S$, let $\BerkU_v = \BerkF_v \subseteq \{z \in \cC_v^{\an} :  |f(z)|_v \le 1 \}$
be the $\fX$-trivial set.  Then $\UU = \prod_v \BerkU_v$ is a $K$-rational Berkovich adelic neighborhood 
\index{$\fX$-trivial}
\index{Berkovich!adelic neighborhood} 
of $\EE$.  

We claim that there are only finitely many points of $\cC(\tK)$ whose 
$\Aut(\tK/K)$-conjugates belong to $\BerkU_v$ for each $v$.  Indeed, if $\alpha$ is such a point,
then $|N_{K(\alpha)/K}(f(\alpha))|_v < 1$ for each $v \in S$, and $|N_{K(\alpha)/K}(f(\alpha))|_v \le 1$
for all $v$, so 
\begin{equation*} 
\prod_v |N_{K(\alpha)/K}(f(\alpha))|_v \ < \ 1 \ .
\end{equation*}
By the Product Formula, we must have $f(\alpha) = 0$.  Since $f$ has only finitely many zeros, 
the conclusion follows. 

\smallskip
We now turn to the proof of (B).  
We are given a compact Berkovich adelic set $\EE = \prod_v \BerkE_v$ with $\gamma(\EE,\fX)^{\an} > 1$
\index{Berkovich!adelic set}   
and a $K$-rational separable Berkovich quasi-neighborhood $\UU = \prod_v \BerkU_v$ of $\EE$. 
\index{quasi-neighborhood!Berkovich}
\index{Berkovich!quasi-neighborhood} 
In this case, we will reduce the result to Theorem \ref{aT1-A1}  
by first shinking $\EE$, then enlarging it within $\UU$, 
and finally cutting back to classical points,  
obtaining a classical set $\FF^0 = \prod_v F_v^0 \subset \prod_v \cC_v(\CC_v)$ 
with a $K$-rational separable quasi-neighborhood $U_v^0 = \prod_v U_v^0 \subset \prod_v \cC_v(\CC_v)$ 
\index{quasi-neighborhood!Berkovich}
satisfying the conditions of Theorem \ref{aT1-A1}.

Since $\gamma(\EE,\fX)^{\an} > 1$, 
we have $\val(\Gamma(\EE,\fX)^{\an}) < 1$.   Let $\varepsilon > 0$
\index{Green's matrix!global}\index{value of $\Gamma$ as a matrix game}
be small enough that if $\Gamma \in M_n(\RR)$ is a symmetric $n \times n$ matrix whose entries differ
from those of $\Gamma(\EE,\fX)^{\an}$ by at most $\varepsilon$, then $\val(\Gamma) < 1$. 
\index{Green's matrix!global!global Berkovich}\index{value of $\Gamma$ as a matrix game}
Again fix a nonempty finite set of places $S$ of $K$  
containing all archimedean places and all nonarchimedean places where $\BerkE_v$ is not $\fX$-trivial, 
\index{$\fX$-trivial}
and choose a set of numbers $\{\varepsilon_v\}_{v \in S}$ with $\varepsilon_v > 0$ for each $v$ 
and $\sum_{v \in S} \varepsilon_v = \varepsilon$.   

If $v$ is archimedean, put $F_v^0 = \BerkE_v$, and let $U_v^0 = \BerkU_v$;  
if $v \notin S$, so $\BerkE_v$ is $\fX$-trivial, put $F_v^0 = U_v^0 = \BerkE_v \cap \cC_v(\CC_v)$, 
\index{$\fX$-trivial}
so $F_v^0$ and $U_v^0$ are the classical $\fX$-trivial sets.   
\index{$\fX$-trivial}
For each nonarchimedean $v \in S$, the separable Berkovich quasi-neighborhood $\BerkU_v$ of $\BerkE_v$
\index{quasi-neighborhood!Berkovich}
\index{Berkovich!quasi-neighborhood} 
can be written as 
\begin{equation*}
\BerkU_v \ = \ \BerkU_{v,0} \cup \big(U_{v,1} \cap \cC_v(F_{w,1})\big) \cup \cdots \cup 
\big(U_{v,N} \cap \cap(F_{w,N})\big) \ , 
\end{equation*} 
where $\BerkU_{v,0} \subset \cC_v^{\an}$ is a Berkovich open set,
\index{Berkovich!open set}  
$U_{v,1}, \ldots, U_{v,N} \subset \cC_v(\CC_v)$ are classical open sets, 
and $F_{w,1}, \ldots, F_{w,N}$ are separable algebraic extensions of $K_v$. 
By hypothesis, $\BerkU_v$ is stable under $\Aut_c(\CC_v/K_v)$, 
which means that $\BerkU_{v,0}$ is stable under $\Aut_c(\CC_v/K_v)$ as well. 
Put 
\begin{equation*}
Y_v \ = \ \BerkE_v \backslash \BerkU_{v,0} \ \subset \ \cC_v(\CC_v) \ . 
\end{equation*} 
Since $Y_v$ is compact as a subset of $\cC_v^{\an}$, it is compact as a subset of $\cC_v(\CC_v)$. 
 
For $n = 1, 2, 3$, choose a finite open cover of $Y_v$ by balls $B(x_1,1/n)^-, \ldots, B(x_{k_n},1/n)^-$ 
and let $\BerkB(x_1,1/n)^-, \ldots, \BerkB(x_{k_n},1/n)^-$ be the corresponding Berkovich open sets. 
\index{Berkovich!open set}  
Put 
\begin{equation*} 
\BerkK_n \ = \ 
\Big(\BerkE_v \backslash (\BerkB(x_1,1/n)^- \cup \ldots \cup \BerkB(x_{k_n},1/n)^- \Big) \cup Y_v \ . 
\end{equation*} 
(If $Y_v$ is empty, take $\BerkK_n = \BerkE_v$ for each $n$.)  
Then $\BerkK_1 \subseteq \BerkK_2 \subseteq \cdots \subseteq \BerkK_n \cdots \subseteq \BerkE_v$ 
is an ascending sequence of compact sets with $\bigcup_{n=1}^{\infty} \BerkK_n = \BerkE_v$, 
so by Proposition \ref{BerkGreenPropertiesProp}(E) there is an $n$
\index{Green's function!Berkovich!properties}
such that for all $x_i \ne x_j \in \fX$ with $i \ne j$ we have  
\begin{equation} \label{BerkGI1}
|G(x_i,x_j;\BerkE_v)^{\an} - G(x_i,x_j;\BerkK_n)^{\an}| \ < \ \varepsilon_v \ ,
\qquad 
|V_{x_i}(\BerkE_v)^{\an} - V_{x_i}(\BerkK_n)^{\an}| \ < \ \varepsilon_v \ .
\end{equation} 
\index{Robin constant!Berkovich!properties} 
Fix such an $n$ and put  
\begin{equation*} 
\BerkX_v \ = \ \BerkK_n \backslash 
\Big(\BerkE_v \backslash \big(\BerkB(x_1,1/n)^- \cup \ldots \cup \BerkB(x_{k_n},1/n)^- \big) \Big) \ .
\end{equation*} 
Then $\BerkX_v$ is compact, $\BerkX_v \subset \BerkU_{v,0}$, and $\BerkK_n = \BerkX_v \cup Y_v$.

Since strict closed affinoids are cofinal in the closed neighborhoods of a compact Berkovich set,
\index{Berkovich!compact set}  
which in turn are cofinal in the open neighborhoods of the set,  
there is a strict closed Berkovich affinoid $\BerkA_v$ 
\index{Berkovich!strict closed affinoid} 
\index{affinoid!strict closed affinoid}
with $\BerkX_v \subseteq \BerkA_v \subset \BerkU_{v,0}$.
As noted above, each strict closed Berkovich affinoid 
\index{Berkovich!strict closed affinoid} 
\index{affinoid!strict closed affinoid}
has finitely many conjugates under $\Aut_c(\CC_v/K_v)$.
Since the union of finitely many strict closed Berkovich affinoids 
\index{Berkovich!strict closed affinoid} 
is either a strict closed Berkovich affinoid or is all of $\cC_v^{\an}$ 
\index{Berkovich!strict closed affinoid} 
\index{affinoid!strict closed affinoid}
(see \cite{Th}, Corollaire 2.1.17), 
after replacing $\BerkA_v$ with the union of its conjugates (which are contained in $\BerkU_{v,0}$), 
we can assume that $\BerkA_v$ is stable under $\Aut_c(\CC_v/K_v)$. 
The intersection $A_v = \BerkA_v \cap \cC_v(\CC_v)$ is a $K_v$-rational 
closed affinoid in the sense of rigid analysis.
\index{affinoid!strict closed affinoid}
Since each rigid analytic\index{rigid analytic space} strict closed affinoid is an $\RL$-domain
\index{$\RL$-domain} 
\index{affinoid!strict closed affinoid}
(see (\cite{Fies}, Satz 2.2) and (\cite{RR1}, Corollary 4.2.14), 
or Corollary \ref{RLeqAffinoid} of Appendix \ref{AppC} below),   
there is a function $f(z) \in \CC_v(\cC_v)$ such that 
\begin{equation*} 
A_v \ = \ \{z \in \cC_v(\CC_v) : |f(z)|_v \le 1 \} \ , \quad 
\BerkA_v \ = \ \{z \in \cC_v^{\an} : |f(z)|_v \le 1 \}  \ .
\end{equation*}
Since $\tK_v$ is dense in $\CC_v$, we can assume that $f(z) \in \tK_v(\cC_v)$, 
and after replacing it with its norm to $K_v$, that $f(z) \in K_v(\cC)$.  
Put $\BerkF_v = \BerkA_v \cup X_v$.  Then $\BerkF_v \subset \BerkU_v$.  
Since $\BerkK_n \subseteq \BerkF_v$, 
the monotonicity of Green's functions in Proposition \ref{BerkGreenPropertiesProp}(C) 
\index{Green's function!Berkovich!monotonic}
shows that for all $x_i \ne x_j \in \fX$ we have 
\begin{equation} \label{BerkGI2}
G(x_i,x_j;\BerkK_n)^{\an}  \ \ge \  G(x_i,x_j;\BerkF_v)^{\an} \ , \qquad 
V_{x_i}(\BerkK_n)^{\an} \ \ge \ V_{x_i}(\BerkF_v)^{\an} \ .
\end{equation}
\index{Robin constant!Berkovich!properties} 

Finally, put 
\begin{equation*} 
F_v^0 \ = \ A_v \cup X_v \ = \ \BerkF_v \cap \cC_v(\CC_v) \ , \qquad 
U_v^0 \ = \ \BerkU_v \cap \cC_v(\CC_v) \ . 
\end{equation*}
Then $F_v^0$ and $U_v^0$ are stable under $\Aut_c(\CC_v/K_v)$, and $F_v^0 \subset U_v^0$.  Furthermore   
$U_v^0 = \big(\BerkU_{v,0} \cap \cC_v(\CC_v)\big) \cup \big(U_{v,1} \cap \cC_v(F_{w,1})\big) \cup 
\cdots \cup \big(U_{v,N} \cap \cC_v(F_{w,N})\big)$ 
so $U_v^0$ is a $K$-rational separable quasi-neighborhood of $F_v^0$. 
\index{quasi-neighborhood!separable}
Since $F_v^0$ is the union of an $\RL$-domain and a compact set, 
\index{$\RL$-domain} 
its closure in $\cC_v^{\an}$ is $\BerkF_v$.  
By (\cite{RR1}, Theorem 4.3.11) it is algebraically capacitable.
\index{algebraically capacitable} 
By Proposition \ref{BerkCompatibilityProp}, for all $x_i \ne x_j \in \fX$ we have 
\begin{equation} \label{BerkGI3}   
G(x_i,x_j;F_v^0) \ = \ G(x_i,x_j;\BerkF_v)^{\an} \ , 
\qquad V_{x_i}(F_v^0) \ = \ V_{x_i}(\BerkF_v)^{\an} \ .
\end{equation} 

Globalizing, take $\FF^0 = \prod_v F_v^0$ and $\UU^0 = \prod_v U_v^0$.  
Then $\FF^0$ is a $K$-rational adelic set in $\prod_v \cC_v(\CC_v)$, 
and $\UU^0$ is a $K$-rational separable quasi-neighborhood of $\FF^0$.
\index{quasi-neighborhood!separable} 
Since $\gamma(\FF,\fX)^{\an} > 1$, 
by (\ref{BerkGI1}), (\ref{BerkGI2}) and (\ref{BerkGI3}) we have $\gamma(\FF^0,\fX) > 1$ as well.  
By Theorem \ref{aT1-A1}, there are infinitely many points $\alpha \in \cC(\tK^{\sep})$ 
whose $\Aut(\tK/K)$-conjugates belong to $U_v^0$ (hence $\BerkU_v$) for each $v$.    
\end{proof}

%% file: NewFSZChap5.tex
\chapter{Initial Approximating Functions:  Archimedean Case} \label{Chap5}

Throughout this section $v$ will be an archimedean place of $K$, 
so $K$ is a number field and $K_v \cong \RR$ or $K_v \cong \CC$. 
Thus $\cC_v(\CC)$ is a connected, compact Riemann surface.
\index{Riemann surface} 
In this section we will construct the archimedean initial approximating functions
\index{initial approximating functions $f_v(z)$!archimedean|ii} 
needed for the proof of Theorem \ref{aT1-B}.  When $K_v \cong \RR$, the construction 
uses results about oscillating pseudopolynomials proved in Appendix \ref{AppB}.  

\smallskip
In Theorem \ref{aT1-B}, 
we are given a compact, $K_v$-simple set $E_v \subset \cC_v(\CC)$, 
\index{$K_v$-simple!set}\index{$\CC$-simple set}  
of positive inner capacity, which is disjoint from $\fX$.  
\index{capacity $> 0$}    
If $K_v \cong \CC$, so $E_v$ is $\CC$-simple, this means that $E_v$ 
is a finite union of pairwise disjoint compact sets $E_{v,i}$, 
each of which 
\begin{enumerate} 
\item is simply connected, has a piecewise smooth boundary,\index{$K_v$-simple!$\CC$-simple}
\index{boundary!piecewise smooth}\index{closure of $\cC_v(\CC)$ interior}\index{simply connected} 
and is the closure of its interior. 
\end{enumerate} 
If $K_v \cong \RR$,\index{$\RR$-simple set} so $E_v$ is $\RR$-simple, 
then $E_v$ is stable under complex conjugation
and is a finite union of pairwise disjoint compact sets $E_{v,i}$, 
where each $E_{v,i}$ either

\begin{enumerate} 
\item is a closed subinterval of $\cC_v(\RR)$ with positive length; or 

\item is disjoint from $\cC_v(\RR)$ and is simply connected,\index{simply connected} 
has a piecewise smooth boundary, and is the closure of its interior.
\index{boundary!piecewise smooth}\index{closure of $\cC_v(\CC)$ interior}
\index{$K_v$-simple!$\RR$-simple}\index{$\RR$-simple set} 
\end{enumerate} 

\noindent{Since} 
$E_v$ is compact we can use the usual Green's functions and Robin constant 
\index{Green's function!of a compact set}
$G(z,x_i;E_v)$ and $V_{x_i}(E_v)$,  instead of the upper ones
$\Gbar(z,x_i;E_v)$ and $\Vbar_{x_i}(E_v)$.  
\index{Green's function!upper}
\index{Robin constant!upper} 
   
If $K_v \cong \CC$, let $E_v^0$ be the interior of $E_v$;  if $K_v \cong \RR$,
let $E_v^0$ be the the union of the $\cC_v(\CC)$-interiors of the components $E_{v,i}$ disjoint from $\cC_v(\RR)$, 
together with the $\cC_v(\RR)$-interiors of the components $E_{v,i}$ contained in $\cC_v(\RR)$.  
We call $E_v^0$ the ``quasi-interior'' of $E_v$.\index{quasi-interior|ii} 

\smallskip
Fix $\varepsilon_v > 0$.  In constructing the approximating functions, we first 
c
replace $E_v$ with a $K_v$-simple set $\tE_v \subset E_v^0$ 
\index{$K_v$-simple!set}  
such that for each $x_i \ne x_j \in \fX$ 
\begin{equation*}
|V_{x_i}(\tE_v) - V_{x_i}(E_v)| \ < \ \varepsilon_v \ , 
\quad 
|G(x_i,x_j;\tE_v) - G(x_i,x_j;E_v)| \ < \ \varepsilon_v \ .
\end{equation*} 
\index{Green's function!properties of} 
\index{Robin constant!properties of}
Substituting $\tE_v$ for $E_v$ gives us ``freedom of movement'' 
in the constructions below.
Next, let $U_v \subset \cC_v(\CC)$ be an open set  
such that $U_v \cap E_v = E_v^0$.  
After shrinking $U_v$ if necessary, we can assume it is bounded away from $\fX$,
and that its connected components are in one to one correspondence with those of $E_v$.
If $K_v \cong \RR$, we can assume it is stable under complex conjugation as well. 

Let $\vs = (s_1, \ldots, s_m)$ be a $K_v$-symmetric probability vector
\index{$K_v$-symmetric!probability vector}
with rational coefficients.  
By a $K_v$-rational $(\fX,\vs)$-function,\index{$(\fX,\vs)$-function!$K_v$-rational|ii} 
we mean a function $f(z) \in K_v(\cC_v)$, whose poles are supported on $\fX$, 
such that if $N_i$ is the order of the pole
of $f$ at $x_i$ and $N = \deg(f)$, then $\frac{1}{N}(N_1, \ldots, N_m) = \vs$. 

\vskip .1 in
The initial approximating functions $f_v(z)$ will be 
\index{initial approximating functions $f_v(z)$!archimedean}
$K_v$-rational $(\fX,\vs)$-functions\index{$(\fX,\vs)$-function!$K_v$-rational} 
having several properties:

First,  $\frac{1}{N} \log_v(|f_v(z)|_v)$ will 
closely approximate $\sum_{i=1}^m s_i G(z,x_i;\tE_v)$ outside $U_v$. 
 
Second, $f_v(z)$ will have all its zeros in $E_v^0$, 
and satisfy $\{z \in \cC_v(\CC) : |f_v(z)|_v \le 1\} \subset U_v$.  
If $K_v \cong \RR$, we also require that $f_v(z)$ have
a property like that of Chebyshev polynomials, oscillating between large 
\index{Chebyshev polynomial} 
positive and negative values on $E_v^0 \cap \cC_v(\RR)$.    
This means that when $f(z)$ is perturbed slightly, 
its zeros continue to belong to $E_v^0$. 

Third, for global aspects of the proof of the Fekete-Szeg\"o theorem, 
\index{Fekete-Szeg\"o theorem with LRC}
we need to be able to independently vary the 
{\em logarithmic leading coefficients|independent variability of archimedean} of $f_v(z)$ 
\index{logarithmic leading coefficients|ii}
at the points in $\fX$. 
We define the logarithmic leading coefficient 
of $f_v(z)$ at $x_i$ to be\label{`SymbolIndexLogLeadf'} 
\begin{equation*}
\Lambda_{x_i}(f_v,\vs) \ = \ \lim_{z \rightarrow x_i} 
          \left(\frac{1}{N} \log_v(|f_v(z)|_v) + s_i \log_v(|g_{x_i}(z)|_v)\right) \ .
\end{equation*}
Similarly, we define\label{`SymbolIndexLogLeadGreen'} 
\begin{eqnarray*}
\Lambda_{x_i}(\tE_v,\vs) & = & \lim_{z \rightarrow x_i} 
     \left(\big(\sum_{j = 1}^m s_j G(z,x_j;\tE_v)\big) + s_i\log_v(|g_{x_i}(z)|_v)\right) \\
      & = & s_i V_{x_i}(\tE_v) + \sum_{j \ne i} s_j G(x_i,x_j;\tE_v) \ .
\end{eqnarray*}
\index{Robin constant!archimedean}
We will require that for pre-specified numbers $\beta_1, \ldots, \beta_m$ 
belonging to an interval $[-\delta_v,\delta_v]$ 
depending only on $\tE_v$ and $U_v$,  
\begin{equation*}
           \Lambda_{x_i}(f_v,\vs) \ = \ \Lambda_{x_i}(\tE_v,\vs) + \beta_i \ .       
\end{equation*}
Here the $\beta_i$ must be $K_v$-symmetric, but otherwise can be chosen arbitrarily. 
\index{$K_v$-symmetric!set of numbers} 
This ``independent variability of the logarithmic leading
coefficients'' is needed to deal with the problem that the probability vector 
\index{logarithmic leading coefficients|ii}
\index{logarithmic leading coefficients!independent variability of archimedean|ii}
\index{independent variability!of logarithmic leading coefficients|ii}
$\hat{s}$ for which $\Gamma(\EE,\fX) \hat{s}$ has equal entries 
\index{Green's matrix!global}
(constructed in \S\ref{Chap3}.\ref{CantorCapacitySection}), 
may not have rational entries.  


\section{ The Approximation Theorems}  \label{ArchApproxThmSection}

There are two cases to consider in constructing the initial 
approximating functions: when $K_v \cong \CC$, and when $K_v \cong \RR$.  
\index{initial approximating functions $f_v(z)$!archimedean} 
The case when $K_v \cong \CC$ follows from results \cite{RR1}:
\index{initial approximating functions $f_v(z)$!when $K_v \cong \CC$|ii}
\index{initial approximation theorem!when $K_v \cong \CC$|ii} 

\begin{theorem}  \label{CThm1}  Suppose $K_v \cong \CC$. 
Let $E_v \subset \cC_v(\CC)$ be a $\CC$-simple\index{$\CC$-simple set} 
set which is disjoint from $\fX$ 
and has positive capacity, 
\index{capacity $> 0$} 
and let $U_v = E_v^0$ be the $\cC_v(\CC)$-interior of $E_v$. 
Fix $\varepsilon_v > 0$.  

Then there is a compact set $\tE_v \subset U_v$ 
composed of a finite union of analytic arcs, 
with  $\cC_v(\CC) \backslash \tE_v$ connnected,
which has the following properties:

$(A)$ For each $x_i \in \fX$ 
\begin{equation} \label{FCApprox1}
|V_{x_i}(\tE_v) - V_{x_i}(E_v)| \ < \ \varepsilon_v \ ,
\end{equation} 
\index{Robin constant!archimedean}
and for all $x_i, x_j \in \fX$ with $x_i \ne x_j$, 
\begin{equation} \label{FCApprox2} 
|G(x_i,x_j;\tE_v) - G(x_i,x_j;E_v)| \ < \ \varepsilon_v \ .
\end{equation}
   
$(B)$ There is a $\delta_v > 0$ such that for any probability vector
$\vs = {}^t(s_1, \ldots, s_m)$ with rational entries,  
and any $\vbeta = {}^t(\beta_1, \ldots, \beta_m) \in [-\delta_v,\delta_v]^m$,
there is an integer $N_v \ge 1$ such that for  
each positive integer  $N$ divisible by $N_v$, 
there exists an $(\fX,\vs)$-function\index{$(\fX,\vs)$-function!$K_v$-rational} 
$f_v(z) \in K_v(\cC_v)$ of degree $N$, satisfying

\quad $(1)$  for each $x_i \in \fX$, 
\begin{equation} \label{FCApprox3} 
\Lambda_{x_i}(f_v,\vs) \ = \ \Lambda_{x_i}(\tE_v,\vs) + \beta_i \ .
\end{equation}

\quad $(2)$ $\{ z \in \cC_v(\CC) : |f_v(z)|_v \le 1 \}$ is contained in $E_v^0;$
in particular, all the zeros of $f_v(z)$ belong to $E_v^0$.
\end{theorem} 

\begin{proof}  We first construct the set $\tE_v$.  Since $E_v$ is the closure
of its interior, and its boundary is a finite union of smooth arcs, each point of
\index{boundary!piecewise smooth} 
$E_v$ is analytically accessible from $E_v^0$.  
Hence Proposition \ref{IdentifyGreenProp} shows that
$\Gbar(z,x_i;E_v^0) = G(z,x_i,E_v)$ for each $x_i$.  
By the monotonicity of Green's functions, 
\index{Green's function!monotonic}
there is a compact set $E_v^{*} \subset E_v^0$ such that for all $x_i \ne x_j$ 
\begin{equation} \label{C1eps2}
|V_{x_i}(E_v^{*}) - V_{x_i}(E_v)| \ < \ \varepsilon_v/2 \ , 
\quad 
|G(x_i,x_j;E_v^{*}) - G(x_i,x_j;E_v)| \ < \ \varepsilon_v/2 \ .
\end{equation}
\index{Robin constant!archimedean}
The remainder of the construction is a combination of results from \cite{RR1}.
By (\cite{RR1}, Proposition 3.3.2) there is a compact set $\tE_v \subset U_v$ 
which is a finite union of analytic arcs, with $\cC_v(\CC) \backslash \tE_v$
connected, such that for each $x_i$.  
\begin{equation} \label{C2eps2}
|V_{x_i}(\tE_v) - V_{x_i}(E_v^{*})| \ < \ \varepsilon_v/2 \ , 
\quad 
|G(x_i,x_j;\tE_v) - G(x_i,x_j;E_v^{*})| \ < \ \varepsilon_v/2 \ .
\end{equation}
\index{Green's function!archimedean}
\index{Robin constant!properties of}
\index{Robin constant!archimedean}
The set $\tE_v$ is obtained by first covering $E_v^{*}$ 
with a finite collection of closed discs contained in $U_v$, 
then taking the union of the boundaries of those discs, 
and finally cutting short intervals out of each boundary arc to obtain a 
\index{arc!smooth}\index{boundary!arc}
set such that $\cC_v(\CC) \backslash \tE_v$  connected.  From (\ref{C1eps2}) 
and (\ref{C2eps2}) we obtain (\ref{FCApprox1}) and (\ref{FCApprox2}).

The existence of a number $\delta_v > 0$, 
and for each rational probability vector $\vs$ and each $\vbeta \in (-\delta_v,\delta_v)^m$,   
the existence of an integer $N_v \ge 1$ 
and an $(\fX,\vs)$-function\index{$(\fX,\vs)$-function!$K_v$-rational} 
$f_{v,0}(z) \in K_v(\cC_v)$ 
of degree $N_v$, with the properties in the theorem, 
is proved in (\cite{RR1}, Theorem 3.3.7).
After shrinking $\delta_v$, one can replace the conditions $|\beta_i| < \delta_v$ 
in (\cite{RR1}, Theorem 3.3.7) with $|\beta_i| \le \delta_v$.
We remark that the independent variability of 
the logarithmic leading coefficients is based on a convexity 
\index{logarithmic leading coefficients!independent variability of archimedean|ii}
\index{independent variability!of logarithmic leading coefficients|ii}
argument using that $-\log([z,w]_{\zeta}$ is 
is everywhere harmonic in $z$, 
apart from logarithmic singularities when $z = w$ or $z = \zeta$      
(see \cite{RR1}, Lemma 3.3.9).

Note that properties (B1) and (B2) are preserved 
when $f_{v,0}(z)$ is raised to a power.
Given an arbitrary multiple $N = kN_v$ 
we can obtain the approximating function of degree $N$
\index{initial approximating functions $f_v(z)$!when $K_v \cong \CC$} 
by putting $f_v(z) = f_{v,0}(z)^k$. 
\end{proof}


\vskip .1 in
When $K_v \cong \RR$, the approximation theorem we need is as follows.
Note that if $f(z) \in \RR(\cC_v)$, then $f$ is real valued on $\cC_v(\RR)$.  
Given a number $M > 0$, we say that $f$ {\em oscillates $k$ times between $\pm M$}  
on an interval $I \subset \cC_v(\RR)$ 
if it varies $k$ times from $-M$ to $M$, or from $M$ to $-M$, on $I$.  
In particular, it has at least $k$ zeros in $I$.  Conversely,
if it has exactly $k$ zeros in $I$ and oscillates $k$ times between $\pm M$,
then each of those zeros is simple.
\index{initial approximating functions $f_v(z)$!when $K_v \cong \RR$|ii}
\index{initial approximation theorem!when $K_v \cong \RR$|ii}   

\begin{theorem}  \label{RThm2}  Suppose $K_v \cong \RR$.  
Let $E_v$ be a compact $\RR$-simple\index{$\RR$-simple set} set which is disjoint from $\fX$ 
and has positive capacity,
\index{capacity $> 0$} 
and let $E_v^0$ be the quasi-interior of $E_v$.  
\index{quasi-interior} 
Fix a $\cC_v(\CC)$-open set $U_v$ such that $U_v \cap E_v = E_v^0$,  
and which is stable under complex conjugation and bounded away from $\fX$.
Take $\varepsilon_v > 0$.    

Then there is a $\RR$-simple\index{$\RR$-simple set} compact set $\tE_v \subset E_v^0$   
such that $\cC_v(\CC) \backslash \tE_v$ is connnected,
which has the following properties:

$(A)$ For each $x_i \in \fX$ 
\begin{equation} \label{FRRApprox4}
|V_{x_i}(\tE_v) - V_{x_i}(E_v)| \ < \ \varepsilon_v \ ,
\end{equation} 
\index{Robin constant!archimedean}
and for all $x_i, x_j \in \fX$ with $x_i \ne x_j$, 
\begin{equation} \label{FRRApprox5} 
|G(x_i,x_j;\tE_v) - G(x_i,x_j;E_v)| \ < \ \varepsilon_v \ .
\end{equation}
\index{Green's function!archimedean} 
   
$(B)$ Given $0 < \cR_v < 1$, there is a 
$\delta_v > 0$ $($depending on $\tE_v$, $U_v$, $\varepsilon_v$, and $\cR_v)$ 
such that for each $K_v$-symmetric probability vector
\index{$K_v$-symmetric!probability vector}
$\vs = {}^t(s_1, \ldots, s_m)$ with rational entries, 
and for each $K_v$-symmetric 
$\vbeta = {}^t(\beta_1, \ldots, \beta_m) \in [-\delta_v,\delta_v]^m$,  
\index{$K_v$-symmetric!vector}
there is an integer $N_v \ge 1$  such that 
for each positive integer  $N$ divisible by $N_v$, 
there is an $(\fX,\vs)$-function\index{$(\fX,\vs)$-function!$K_v$-rational} 
$f_v(z) \in K_v(\cC_v)$ of degree $N$  which satisfies 

\quad $(1)$  For each $x_i \in \fX$, 
\begin{equation} \label{FRRApprox7} 
\Lambda_{x_i}(f_v,\vs) \ = \ \Lambda_{x_i}(\tE_v,\vs) + \beta_i \ .
\end{equation}

\quad $(2)$ $\{z \in \cC_v(\CC) : |f_v(z)| \le 1 \} \ \subset \ U_v$\ .

\quad $(3)$  All the zeros of $f_v(z)$ belong to $E_v^0$,   
and if $E_{v,i}$ is a component of $E_v$ contained in $\cC_v(\RR)$  
and $f_v(z)$ has $N_i$ zeros in $E_{v,i}$, 
then $f_v(z)$ oscillates $N_i$ times between $\pm \cR_v^N$ on $E_{v,i}$.  

\end{theorem} 

The proof of Theorem \ref{RThm2} will occupy the rest of this chapter.   
For notational convenience we identify $K_v$ with $\RR$, and write $\log_v(x) = \log(x)$.   


\section{ Outline of the Proof of Theorem $\ref{RThm2}$}  \label{OutlineSection} 

In this section we sketch the ideas behind the proof of Theorem \ref{RThm2}.
In \S\ref{Chap5}.\ref{IndependenceSection} we establish an independence lemma, 
and in \S\ref{Chap5}.\ref{RealCaseProofSection} we give the details of the proof.
\index{initial approximating functions $f_v(z)$!construction when $K_v \cong \RR$!outline|(}    

\vskip .1 in
     
If the $\RR$-simple\index{$\RR$-simple set} set $E_v = \bigcup_{i=1}^n E_{v,i}$ has no components 
in $\cC_v(\RR)$, Theorem \ref{RThm2} 
follows from results in (\cite{RR1}).     
For the remainder of the discussion below, assume that some $E_{v,i}$ is contained in $\cC_v(\RR)$.    
By standard potential-theoretic arguments,  
we can construct a $K_v$-simple compact set $E_v^{*} \subset E_v^0$ such that for each $x_i \ne x_j$ 
\index{$K_v$-simple!set} 
\begin{equation*} 
|V_{x_i}(E_v^{*})-V_{x_i}(E_v)| < \varepsilon_v \ , \quad  
|G(x_i,x_j;E_v^{*})-G(x_i,x_j;E_v)| < \varepsilon_v \ . 
\end{equation*}\index{Robin constant!archimedean}\index{Green's function!archimedean}
In doing so, we can arrange that each of the intervals making up $E_v^{*} \cap \cC_v(\RR)$
is ``short'', in a sense to be made precise later.  
\index{short@`short' interval}  
   
To assist in constructing $(\fX,\vs)$-functions\index{$(\fX,\vs)$-function} 
with prescribed logarithmic leading coefficients,
\index{logarithmic leading coefficients} 
we next adjoin a finite number of short intervals to $E_v^{*}$, which can be `wiggled' inside $E_v^0$. 
\index{short@`short' interval}    
We will show that there are finitely many points 
$t_1, \ldots, t_d \subset (E_v^0 \cap \cC_v(\RR)) \backslash E_v^{*}$ 
and a number $h > 0$, such that if $\tE_{v,\ell} = [t_{\ell}-h,t_{\ell}+h]$ \
for $\ell = 1, \ldots, d$ (the intervals are defined 
using local coordinates at the points $t_{\ell}$) then the set 
\begin{equation*}
\tE_v \ := \ E_v^{*} \cup (\bigcup_{\ell=1}^d \tE_{v,\ell})
\end{equation*}
meets the needs of the theorem, in particular satisfying 
$|V_{x_i}(\tE_v)-V_{x_i}(E_v)| < \varepsilon_v$ for each $x_i$,
and  $|G(x_i,x_j;\tE_v)-G(x_i,x_j;E_v)| < \varepsilon_v$ 
for each $x_i \ne x_j$.  
\index{Robin constant!archimedean}
\index{Green's function!archimedean}
The points $t_1, \ldots, t_d$ must be in ``general position'', 
in a sense to be described in \S\ref{Chap5}.\ref{IndependenceSection}, 
and  $h$ must be small enough that the intervals $\tE_{v,\ell}$
are contained in $(E_v^0 \cap \cC_v(\RR)) \backslash \fX$ and 
are disjoint from $E_v^{*}$ and each other.  

\vskip .1 in
We next construct the functions $f_v(z)$.  The first part of the construction
is purely potential-theoretic, and is carried out in Appendices \ref{AppA} and \ref{AppB}.    

Let each $[z,w]_{x_i}$ be normalized so 
$\lim_{z \rightarrow x_i} [z,w]_{x_i} \cdot |g_{x_i}(z)|_v = 1$.  
As in (\ref{FXsDef}), given a probability vector $\vs \in \cP^m$, we define the 
{\em $(\fX,\vs)$-canonical distance} to be
\index{canonical distance!$[z,w]_{\fX,\vs}$|ii} 
\begin{equation*}
[z,w]_{\fX,\vs} = \prod_{i=1}^{m} ([z,w]_{x_i})^{s_i} \ .
\end{equation*}
There is a potential theory\index{potential theory!$(\fX,\vs)$} for the $(\fX,\vs)$-canonical distance
similar to the one for the usual canonical distance (see Section\ref{XSPotTheorySection} of Appendix \ref{AppA}):
\index{canonical distance!$[z,w]_{\zeta}$}
\index{canonical distance!$[z,w]_{\fX,\vs}$!potential theory for}

Let $H_v \subset \cC_v(\CC) \backslash \fX$ be any compact set with positive capacity.
\index{capacity $> 0$} 
For each probability measure $\nu$ supported on $H_v$, 
define the $(\fX,\vs)$-energy by
\begin{equation*}
I_{\fX,\vs}(\nu) 
\ = \ \iint_{H_v \times H_v} -\log([z,w]_{\fX,\vs}) \ d\nu(z) d\nu(w)  
\end{equation*}
and the $(\fX,\vs)$-Robin constant by  
\begin{equation} \label{FXV1} 
V_{\fX,\vs}(H_v) \ = \ \inf_{\nu} I_{\fX,\vs}(\nu) \ . 
\end{equation} 
\index{Robin constant!archimedean!archimedean $(\fX,\vs)$}
By Theorem \ref{ATE10B}, there is a unique probability measure $\mu_{\fX,\vs}$ which achieves 
the infimum in (\ref{FXV1});  it will be called the $(\fX,\vs)$-equilibrium
measure of $H_v$.  By the same theorem, the $(\fX,\vs)$-potential function\index{potential function!$(\fX,\vs)$}  
\begin{equation*}
u_{\fX,\vs}(z) \ = \ \int_{H_v} -\log([z,w]_{\fX,\vs}) \, d\mu_{\fX,\vs}(w)
\end{equation*}
satisfies $u_{\fX,\vs}(z) \le V_{\fX,\vs}(H_v)$  for all $z$, 
with $u_{\fX,\vs}(z) = V_{\fX,\vs}(H_v)$ on $H_v$.
\index{Robin constant!archimedean!archimedean $(\fX,\vs)$}
 
By Proposition \ref{BPropF1}, the $(\fX,\vs)$-Green's function 
$G_{\fX,\vs}(z;H_v) := V_{\fX,\vs}(H_v) - u_{\fX,\vs}(z)$ 
\index{Green's function!$(\fX,\vs)$}
\index{Robin constant!archimedean!archimedean $(\fX,\vs)$}
can be decomposed as 
\begin{equation} \label{FGreen1} 
G_{\fX,\vs}(z;H_v) \ = \ \sum_{i=1}^m s_i G(z,x_i;H_v) 
\end{equation}
and the $(\fX,\vs)$-equilibrium measure of $H_v$ is given by
\begin{equation*}
\mu_{\fX,\vs} \ = \ \sum_{i=1}^m s_i \mu_{x_i}
\end{equation*}
where $\mu_{x_i}$ is the equilibrium measure of $H_v$ with respect to $x_i$.  

\vskip .1 in
Recall that an $(\fX,\vs)$-function\index{$(\fX,\vs)$-function} 
$f(z) \in K_v(\cC_v)$ of degree $N$ 
is a function with  polar divisor $\sum_{i=1}^m Ns_i(x_i)$. 
If the zeros of $f(z)$ are $\alpha_1, \ldots, \alpha_N$ 
(listed with multiplicities), then  an easy symmetrization argument 
shows there is a constant $C$ such that 
\begin{equation} \label{FGHB1} 
|f(z)|_v \ = \ C \cdot \prod_{k=1}^N [z,\alpha_k]_{\fX,\vs} 
\end{equation}  
for all $z \in \cC_v(\CC) \backslash \fX$:  
let $\xi_{1}, \ldots, \xi_N$ be the points $x_1, \ldots, x_m$ 
listed according to their multiplicities in $\div(f)$.  
For each permutation $\pi$ of $\{1, \ldots, N \}$, by (\ref{FCan1a}) there is 
a constant such that 
$|f(z)|_v  =  C(\pi) \cdot \prod_{k = 1}^N [z,\alpha_k]_{\xi_{\pi(k)}}$. 
Taking the product over all $\pi$, and then extracting $(N!)^{th}$ roots, 
gives (\ref{FGHB1}).  

This motivates the definition of an {\em $(\fX,\vs)$-pseudopolynomial} 
\index{pseudopolynomial|ii} 
\index{pseudopolynomial!$(\fX,\vs)$}
(usually we will just say pseudopolynomial).  
Given a constant $C$ and 
points $\alpha_1, \ldots, \alpha_N \in \cC_v(\CC) \backslash \fX$, 
the associated pseudopolynomial is the non-negative real valued function 
\begin{equation*} 
P(z) \ = \ P_{\valpha}(z) \ = \ C \cdot \prod_{k=1}^N [z,\alpha_k]_{\fX,\vs} \ .
\end{equation*} 
We write $N = \deg(P)$.  
We call the $\alpha_k$ the {\em roots} of $P$, 
we call $\div(P) := \sum_{i=1}^N (\alpha_i) - \sum_{i=1}^m N s_i (x_i)$ 
the {\em divisor} of $P$, and we call 
\begin{equation*}
\nu(z) \ = \ \nu_P(z) \ := \ \frac{1}{N} \sum_{k=1}^N \delta_{\alpha_k}(z)
\end{equation*}
the {\em probability measure associated to $P$}.   
Note that $P_{\valpha}(z)$ makes sense even when the $\alpha_k$
are not the zeros of an $(\fX,\vs)$-function\index{$(\fX,\vs)$-function} $f(z)$, 
but it agrees with $|f(z)|_v$ (up to a multiplicative constant) when such a function exists.  
Furthermore, $P_{\valpha}(z)$ varies continuously with its roots.  
This allows us to investigate absolute values of 
$(\fX,\vs)$-functions\index{$(\fX,\vs)$-function} 
with prescribed zeros, without worrying about principality of the divisors. 

For each $x_i \in \fX$, we define the logarithmic leading coefficient of 
$P(z)$ at $x_i$ to be
\begin{equation} \label{FQLogLead}
\Lambda_{x_i}(P,\vs) \ = \ \lim_{z \rightarrow x_i} \frac{1}{N} \log(P(z)) 
         + s_i \log(|g_{x_i}(z)|_v) \ . 
\end{equation}

\smallskip
We now apply this to the $\RR$-simple\index{$\RR$-simple set} set $H_v = \tE_v$.  
A detailed study of pseudopolynomials is carried out in Appendix \ref{AppB}. 
\index{pseudopolynomial} 
There it is shown that if the components of $\tE_v$ contained in $\cC_v(\RR)$ 
are sufficiently short (the precise meaning of ``short'' is given in 
Definition \ref{ShortnessDef}, 
\index{short@`short' interval}  
in terms of the canonical distance functions $[z,w]_{x_i}$ 
\index{canonical distance!determines `shortness'}
\index{canonical distance!$[z,w]_{\zeta}$}
relative to the $x_i \in \fX$), 
then by potential-theoretic methods one can show the existence of 
pseudopolynomials which behave like absolute values of classical Chebyshev polynomials, 
\index{pseudopolynomial} 
\index{Chebyshev polynomial} 
and have large oscillations on $\tE_v$. 

Let $D > d$ be the number of components of $\tE_v$, and label those 
in $E_v^{*}$ as $\tE_{v,d+1}, \ldots, \tE_{v,D}$, 
so that $\tE_v = \bigcup_{\ell=1}^D \tE_{v,\ell}$.
For each $\ell = 1, \ldots, D$, put $\sigma_{\ell} = \mu_{\fX,\vs}(\tE_{v,\ell})$, 
and put $\vec{\sigma} = (\sigma_1, \ldots, \sigma_D)$.  
The following is a specialization of Theorem \ref{bT3} of Appendix \ref{AppB}, 
formulated using the notation of this section. We write $\zbar$ for the complex conjugate of $z$.    

\begin{theorem} \label{PotPart}
Suppose $K_v \cong \RR$. 
Assume that $\fX$ is stable under complex conjugation, 
and that $\tE_v \subset \cC_v(\CC) \backslash \fX$ is $K_v$-simple, 
\index{$K_v$-simple!set}   
with components $\tE_{v,1}, \ldots, \tE_{v,D}$. 
Assume also that each component $\tE_{v,\ell}$ 
contained in $\cC_v(\RR)$ is a ``short interval'' relative to $\fX$  
\index{short@`short' interval}  
in the sense of Definition $\ref{ShortnessDef}$.  
  
Fix a $K_v$-symmetric probability vector $\vs \in \cP^m$.  
\index{$K_v$-symmetric!probability vector}
For each $\ell = 1, \ldots, D$ put $\sigma_{\ell} = \mu_{\fX,\vs}(\tE_{v,\ell})$ 
and let $\vsigma = (\sigma_1, \ldots, \sigma_D)$. 
Given a $K_v$-symmetric vector $\vn \in \NN^D$ 
\index{$K_v$-symmetric!vector}
write $N = N_{\vn} = \sum_{\ell} n_{\ell}$.
Then there are a collection of $(\fX,\vs)$-pseudopolynomials 
\index{pseudopolynomial!$(\fX,\vs)$} 
\begin{equation*}
\big\{Q_{\vn}(z)\big\}_{\vn \in \NN^D} \ , 
\end{equation*} 
and numbers $0 < R_{\vn} \le 1$, with the following properties$:$ 

\vskip .03 in
$(A)$ For each $\vn$, $Q_{\vn}$ satisfies $\|Q_{\vn}\|_{\tE_v} = 1$,  
with $Q_{\vn}(z) = Q_{\vn}(\zbar)$ for all $z \in \cC_v(\CC)$.
The roots of $Q_{\vn}$ all belong to $\tE_v$, 
with $n_{\ell}$ roots in each $\tE_{v,\ell}$.
For each $\tE_{v,\ell}$ which is a short interval
\index{short@`short' interval}   
the roots of $Q_{\vn}$ in $\tE_{v,\ell}$ are distinct,  
and $Q_{\vn}$ varies $n_{\ell}$ times from $R_{\vn}^N$ to $0$ to $R_{\vn}^N$ 
on $\tE_{v,\ell}$. 

\vskip .03 in
$(B)$ Let  $\{\vn_k\}_{k \in \NN}$ be a sequence  
with $N_{\vn_k} \rightarrow \infty$ 
and $\vn_k/N_{\vn_k} \rightarrow \vec{\sigma}$.  
Then $\lim_{k \rightarrow \infty} R_{\vn_k}  = 1$, and 
the discrete measures $\omega_{\vn_k}$ associated to the $Q_{\vn_k}$
converge weakly to the equilibrium distribution $\mu_{\fX,\vs}$ of $\tE_v$.  
For each neighborhood $\tU_v$ of $\tE_v$,     
the functions $\frac{1}{N_{\vn_k}}\log(Q_{\vn_k}(z))$ 
converge uniformly to $G_{\fX,\vs}(z,\tE_v)$ 
on $\cC_v(\CC_v) \backslash (\tU_v \cup \fX)$,
and for each $x_i \in \fX$,   
\begin{equation*} 
\lim_{k \rightarrow \infty} \Lambda_{x_i}(Q_{\vn_k},\vs) 
\ = \ \Lambda_{x_i}(\tE_v,\vs) \ .    
\end{equation*}     
\end{theorem}    
          
Theorem \ref{PotPart} is the potential-theoretic input to the construction.
We will call the $(\fX,\vs)$-pseudopolynomials $Q_{\vn}(z)$ given by  
\index{pseudopolynomial!special|ii} 
Theorem \ref{PotPart} {\em special pseudopolynomials} for $\tE_v$.  
\vskip .1 in

However, we want to construct $(\fX,\vs)$-functions,\index{$(\fX,\vs)$-function} 
not just pseudopolynomials, with large oscillations on $\tE_v$.  
\index{pseudopolynomial} 
The second part of the construction addresses this.  

\vskip .1 in
Let $\cR_v$ be as in Theorem \ref{RThm2}. 
Applying Theorem \ref{PotPart}, for an appropriate $\vn$  
we obtain a pseudopolynomial $Q(z) = Q_{\vn}(z)$ 
\index{pseudopolynomial} 
which varies $n_{\ell}$ times from $\cR_v^N$ to $0$ to $\cR_v^N$ on 
each real component $\tE_{v,\ell}$ of $\tE_v$.     
If $\div(Q) := \sum_{k=1}^N (\alpha_k) - \sum_{i=1}^m N s_i (x_i)$ 
were principal, then since it is $K_v$-symmetric 
\index{$K_v$-symmetric!divisor}
there would be an $(\fX,\vs)$-function\index{$(\fX,\vs)$-function} $f(z) \in \RR(\cC_v)$ 
with $|f(z)| = Q(z)$ for all $z$.  
Moreover, by Theorem \ref{PotPart} the roots of $Q(z)$ in $\cC_v(\RR)$
are simple, so each time $Q(z)$ varies from $\cR_v^N$ to $0$ to $\cR_v^N$
on a real component of $\tE_v$, the 
function $f(z)$ oscillates   
from $\cR_v^N$ to $-\cR_v^N$, or from $-\cR_v^N$ to $\cR_v^N$. 

Of course, it is unreasonable to expect $Q(z)$ to 
have a principal divisor.  We must assume from the start that 
$\vs = (s_1, \ldots, s_m) \in \QQ^m$ and that $N \vs \in \ZZ^m$, 
but still $\div(Q)$ will generally not be principal.

Our plan is to modify  $Q(z)$  by scaling it and 
``sliding some of its roots along $\cC_v(\RR)$'' 
to make $\div(Q)$ principal. 
In this process, some of the roots may move outside $\tE_v$, 
but they will remain inside $E_v^0$.

In the decomposition $\tE_v = E_v^{*} \cup (\bigcup_{\ell = 1}^d \tE_{v,\ell})$,  
write  $\tE_{v,\ell} = [t_{\ell}-h,t_{\ell}+h]$ 
in suitable local coordinates.  
Suppose we are given a number $r > 0$
small enough that $[t_{\ell}-h-r,t_{\ell}+h+r]$
is contained in $U_v$ and in the coordinate patch of $\tE_{v,\ell}$, 
for each $\ell = 1, \ldots, d$.  Suppose 
we are also given real numbers $\varepsilon_1, \ldots, \varepsilon_d$ 
with each $|\varepsilon_{\ell}| \le r$.  
Then for each $z \in \tE_{v,\ell}$, it makes sense to speak of the point 
$z+\varepsilon_{\ell}$ in terms of the given local coordinates.
Let $\varepsilon_0$ be another small real number, and 
put $\veps = (\varepsilon_0, \ldots, \varepsilon_d)$. 
If $Q(z) = C \cdot \prod_{k =1}^N [z,\alpha_k]_{\fX,\vs}$, 
define 
\begin{equation}\label{bF6}
Q^{\veps}(z) \ = \ \exp(\varepsilon_0)^N \cdot C 
   \cdot \prod_{\alpha_k \in E_v^{*}} [z,\alpha_k]_{\fX,\vs}
   \cdot \prod_{\ell=1}^d \, \prod_{\substack{ \alpha_k \in \tE_{v,\ell} }}
                               [z,\alpha_k + \varepsilon_{\ell}]_{\fX,\vs} \ .
\end{equation}
Using the continuity of $[z,w]_{\fX,\vs}$, one can choose $r$ 
so that if each $\varepsilon_{\ell}$ is permitted 
to vary over the interval $[-r,r]$, then $Q^{\veps}(z)$
oscillates $N$ times from $\cR_v^N$ to $0$ to $\cR_v^N$ on $U_v$.
(In the construction, $r$ will be chosen before $h$, 
and $h$ will be much smaller than $r$.)

\vskip .1 in
We claim that for sufficiently large $N$, one can choose  $\veps$
in such a way that $\div(Q^{\veps})$ becomes principal.   
If $\cC_v$ has genus $g = 0$, there is nothing to prove.                              
If $g > 0$, consider how $\div(Q)$ 
changes when $Q$ is replaced by $Q^{\veps}$.  
Let $\Jac(\cC_v) \cong \text{Div}^0(\cC_v)/P(\cC_v)$ 
be the Jacobian of $\cC_v$,
\index{Jacobian variety}
where $\text{Div}^0(\cC_v)$ is the set of ($\CC$-rational) 
divisors of degree $0$, and $P(\cC_v)$ 
is the subgroup of principal divisors.  Fixing a base point $p_0 \in \cC_v(\RR)$, 
there is a natural embedding $\varphi : \cC_v(\CC) \rightarrow \Jac(\cC_v)(\CC)$,  
$\varphi(p) = \cl((p)-(p_0))$.  Since $\cC_v$ is defined over $\RR$, 
the space of holomorphic differentials\index{differential!holomorphic} $H^1(\cC_v,\CC)$ 
has a basis consisting of real differentials\index{differential!real} $\omega_1, \ldots, \omega_g \in H^1(\cC_v,\RR)$.  
If $\cL \subset \CC^g$ is the corresponding period lattice,\index{period lattice} 
and we identify $\Jac(\cC_v)(\CC)$ with $\CC^g/\cL$, then 
\begin{equation*}
\varphi(p) \ = \ (\int_{p_0}^p \omega_1  , \ldots, 
                             \int_{p_0}^p \omega_g ) \ (\text{mod}\ \cL) \ .
\end{equation*}     
The embedding $\varphi$ induces a canonical surjective map
\begin{equation*}
\Phi : \cC_v(\CC)^g \rightarrow \Jac(\cC_v)(\CC) \ , \qquad  
\Phi((p_1, \ldots, p_g)) = \sum_{i=1}^g \varphi(p_i) \ .
\end{equation*}
This factors through the $g$-fold symmetric product $\text{Sym}^{(g)}(\cC_v)$ as 
\begin{equation*}
\Phi \ : \ \cC_v(\CC)^g \ \rightarrow \ \text{Sym}^{(g)}(\cC_v)
 \ \rightarrow \  \Jac(\cC_v)(\CC)
\end{equation*}
where the first map is finite of degree $g!$, and the second is a birational 
morphism.  Hence the image $\Phi(\cC_v(\RR)^g)$ contains
a $g$-dimensional open subset of $\Jac(\cC_v)(\RR)$.  On the other hand, $\Jac(\cC_v)(\RR)$
is a compact real Lie subgroup of $\Jac(\cC_v)(\CC)$ of dimension at most $g$.
Hence, $\Jac(\cC_v)(\RR)$ must have dimension exactly $g$, and the  
identity component of $\Jac(\cC_v)(\RR)$ must be isomorphic
to a real torus $\RR^g/\cL_0$ for some lattice $\cL_0 \subset \cL \bigcap \RR^g$, 
under our identification $\Jac(\cC_v)(\CC) \cong \CC^g/\cL$.  
The component group of $\Jac(\cC_v)(\RR)$ is an elementary abelian $2$-group,
and if $N \vs \in 2 \cdot \NN^m$,
one can arrange that $\div(Q)$ belongs to the identity component.
\index{Jacobian variety!structure of $\Jac(\cC_v)(\RR)$} 
 
Note that since $p_0$ and the $\omega_i$ are real, the integrals
$\int_{p_0}^p \omega_i(z) \, dz$ are real-valued on the component of $\cC_v(\RR)$ 
containing $p_0$. 

Now consider what happens when 
the roots $\alpha_k \in \tE_{v,\ell} = [t_{\ell} - h,t_{\ell} +h]$ 
are translated by $\varepsilon_{\ell}$.
In terms of the local coordinate at $t_{\ell}$, we can write
$\omega_1 = h_1(z) \, dz, \ldots, \omega_g = h_g(z) \, dz$.
The resulting change in $\varphi(\div(Q))$ is
\begin{eqnarray}  
 \lefteqn{ \sum_{\ell = 1}^d \sum_{\alpha_k \in \tE_{v,\ell}} 
       ( \int_{\alpha_k}^{\alpha_k + \varepsilon_{\ell}} \omega_1, 
     \ldots, \int_{\alpha_k}^{\alpha_k + \varepsilon_{\ell}} \omega_g)
                           \ (\text{mod}\ \cL) } \label{FF14}  \\
    & & \qquad  \cong \ \sum_{\ell =1}^d N \sigma_{\ell} \cdot 
        ( \int_{t_{\ell}}^{t_{\ell} + \varepsilon_{\ell}} h_1(z) \, dz,
     \ldots, \int_{t_{\ell}}^{t_{\ell} + \varepsilon_{\ell}} h_g(z) \, dz ) 
                            \ (\text{mod}\ \cL)  \notag \\
    & & \qquad \qquad 
    \cong \ \sum_{\ell = 1}^d N \sigma_{\ell} \varepsilon_{\ell} \cdot 
             (h_1(t_{\ell}), \ldots, h_g(t_{\ell})) \ (\text{mod}\ \cL) \ . \notag
\end{eqnarray}
However, it is more useful to consider the normalized lift  
$\widehat{\varphi} : \RR^{d+1} \rightarrow \RR^g$, 
\begin{eqnarray}
\widehat{\varphi}(\veps) & := & 
\frac{1}{N} \sum_{\ell = 1}^d \sum_{\alpha_k \in \tE_{v,\ell}} 
       ( \int_{\alpha_k}^{\alpha_k + \varepsilon_{\ell}} h_1(z) \, dz, 
     \ldots, \int_{\alpha_k}^{\alpha_k + \varepsilon_{\ell}} h_g(z) \, dz)  
                      \label{bF14} \\
   & \cong & \sum_{\ell=1}^d  
   \sigma_{\ell} \varepsilon_{\ell} \cdot (h_1(t_{\ell}), \ldots, h_g(t_{\ell})) \ ,
                     \notag                          
\end{eqnarray}
for which 
\begin{equation} \label{FNNorm}
\varphi(\div(Q^{\veps})) - \varphi(\div(Q)) \ = \ 
 N \widehat{\varphi}(\veps) \ (\text{mod}\ \cL_0) \ .
\end{equation}
We will show that for sufficiently large $N$, 
as $\veps$  varies over the ball $B(0,r) = \{\vx \in \RR^{d+1} : |\vx| \le r \}$, 
the image $\widehat{\varphi}(B(0,r))$ contains a fixed neighborhood 
of the origin in $\RR^g$.  
Hence the non-normalized change $N \widehat{\varphi}(\veps)$ 
varies over a region containing a fundamental domain for $\cL_0$,  
and $\veps$ can be chosen 
so $\div(Q^{\veps})$ is principal.  

\vskip .1 in
We also need to consider the change in the logarithmic leading coefficients
\index{logarithmic leading coefficients} 
produced by passing from $Q(z)$ to $Q^{\veps}(z)$.  
Because $\lim_{z \rightarrow x_i} [z,w]_{x_i} |g_{x_i}(z)|_v = 1$ for each $x_i$,   
\begin{eqnarray}
\Lambda_{x_i}(Q,\vs) 
    &=& \lim_{z \rightarrow x_i} \left( \frac{1}{N} \log(Q(z)) 
                           + s_i \log(|g_i(z)|) \right)  \notag \\
    &=& V_{\fX,\vs}(\tE_v) + \sum^m_{\substack{ j = 1 \\ j \ne i }}                         
         s_j \left( \sum_{k = 1}^N 
         \frac{1}{N} \log([x_i,\alpha_k]_{x_j}) \right) \ . 
                                    \label{bF12}
\end{eqnarray}
\index{Robin constant!archimedean!archimedean $(\fX,\vs)$}
For each $1 \le \ell \le d$, consider the contribution of the roots 
$\alpha_k \in \tE_{v,\ell} = [t_{\ell}-h,t_{\ell}+h]$ to (\ref{bF12}).  
Recall that $\sigma_{\ell} = \mu_{\fX,\vs}(\tE_{v,\ell})$.   
If $N$  is sufficiently large, $Q(z)$ has approximately $N \sigma_{\ell}$ \, 
roots in $\tE_{v,\ell}$.  
If $h$ is sufficiently small, the $\alpha_k$ belonging 
to $\tE_{v,\ell}$ can be viewed as being essentially equal to $t_{\ell}$, 
and will contribute approximately
\begin{equation*}
\sum^m_{\substack{ j = 1 \\ j \ne i }}                         
                s_j \log([x_i,t_{\ell}]_{x_j}) \, \sigma_{\ell} 
\end{equation*}
to (\ref{bF12}).  If these roots are translated by $\varepsilon_{\ell}$,
and the contributions from all $\tE_{v,\ell}$ are summed, along with the 
change due to scaling by $\exp(\varepsilon_0)^N$, we see that                   
\begin{equation}
  \Lambda_{x_i}(Q^{\veps},\vs) - 
                   \Lambda_{x_i}(Q,\vs)                   
   \ \cong \  \varepsilon_0 + \sum_{\ell=1}^d 
           \sum^m_{\substack{ j = 1 \\ j \ne i }}
    s_j \left. \left( \frac{d}{dt} \log([x_i,t_{\ell}]_{x_j}) \right)
                       \right|_{t_{\ell}}
                 \cdot \sigma_{\ell} \varepsilon_{\ell} \ . \label{bF13}
\end{equation}

\vskip .1 in
In the next section we will show that $t_1, \ldots, t_d$ can be chosen
so that the quantities (\ref{bF14}) and (\ref{bF13}) are independent.
Given this, a topological argument shows that there exist 
$\veps$ for which $\div(Q^{\veps})$ is principal
and such that the logarithmic leading coefficients
\index{logarithmic leading coefficients} 
$\Lambda_{x_i}(Q^{\veps},\vs)$ corresponding to 
distinct orbits $\fX_{\ell}$ can be specified arbitrarily, 
provided they are sufficently close to the  
$\Lambda_{x_{\ell}}(\tE_v,\vs)$.  
It will also be seen that the construction can be carried out
uniformly for all $\vs$.  
\index{initial approximating functions $f_v(z)$!construction when $K_v \cong \RR$!outline|)}                          
         

\section{ Independence} \label{IndependenceSection}

Let $\cC_v/\RR$ and $\fX$
be as in \S\ref{Chap5}.\ref{OutlineSection}.  If $m_1$ points in $\fX$ are real, and $2m_2$ 
points are in complex conjugate pairs, let $\fX_1, \ldots, \fX_{m_1+m_2}$
denote the corresponding orbits.
\index{initial approximating functions $f_v(z)$!construction when $K_v \cong \RR$!independence of differentials|(}        

Fix a $K_v$-symmetric, rational probability vector $\vs$, 
\index{$K_v$-symmetric!probability vector}
and consider formula (\ref{bF13}) giving
$\Lambda_{x_i}(Q^{\veps},\vs) - \Lambda_{x_i}(Q,\vs)$.  
Our first goal is to express the differential 
\begin{equation}  \label{bF15}
 \frac{d}{dt} \log([x_i,t]_{x_j}) \, dt  
\end{equation}
in terms of meromorphic differentials\index{differential!meromorphic} on the Riemann surface $\cC_v(\CC)$.
\index{Riemann surface} 

For each $x_i \ne x_j$, there is a multi-valued holomorphic function
$\Omega_{i,j}(z)$ on $\cC_v(\CC) \backslash \{x_i,x_j\}$
whose real part coincides with $\log([x_i,z]_{x_j})$,
and which has pure imaginary periods\index{pure imaginary periods} 
over all cycles of $\cC_v(\CC)$ and loops around $x_i$, $x_j$. 
(See \cite{SS}, or \cite{RR1}, pp. 64-65).  
In a given coordinate patch, 
$\Omega_{i,j}(z)$ is only defined up to a pure imaginary constant;  
however the differential $d\Omega_{i,j}(z) = G_{i,j}(z) \, dz$ 
is a globally well-defined\index{differential!of the third kind} 
differential of the third kind on $\cC_v(\CC)$, which is
holomorphic except for simple poles with residue $+1$ at $x_i$ and
residue $-1$ at $x_j$.  
  
Letting $\zbar$ be the complex conjugate of $z$, 
and writing $\tau$ for the permutation of $\{1, \ldots, m\}$  
such that $\overline{x_i} = x_{\tau(i)}$, one sees that for a suitable
normalization of $\Omega_{\tau(i),\tau(j)}(z)$, 
\begin{equation*}
\Omega_{\tau(i),\tau(j)}(\overline{z}) = \overline{\Omega_{i,j}(z)} \ .
\end{equation*} 
With this normalization, for each $t \in \cC_v(\RR)$,
\begin{equation*}
\log([x_i,t]_{x_j})
   = \frac{1}{2} \left( \Omega_{i,j}(t)
                 + \Omega_{\tau(i),\tau(j)}(t) \right) \ .
\end{equation*}
Assuming $t$ is the real part of the local coordinate function $z$, 
this means that at a point $t_{\ell} \in \cC_v(\RR)$, 
\begin{equation*}
\left.\left( \frac{d}{dt} \log([x_i,t]_{x_j}) \right) \right|_{t_{\ell}}
     \ = \ \frac{1}{2} \left( G_{i,j}(t_{\ell}) 
               + G_{\tau(i),\tau(j)}(t_{\ell}) \right) \ .
\end{equation*}

Note that if $x_i$ is complex, and $x_j = x_{\tau(i)}$, 
then $[x_i,t]_{x_j}$ is constant for $t \in \cC_v(\RR)$.
This is because $[x_i,z]_{x_j} \cdot [x_j,z]_{x_i}$
is constant on $\cC_v(\CC)$:  its logarithm is harmonic
everywhere except possibly at $x_i$, $x_j$;  but the
singularities at those points cancel so it is harmonic everywhere.  
On the other hand, applying $\tau$ to $[x_i,t]_{x_j}$, 
we have $[x_i,t]_{x_j} = [x_j,t]_{x_i}$.  
Combining these shows $[x_i,t]_{x_j}$ is constant.
Thus the terms with $j = \tau(i)$ (those for which $x_i$ and $x_j$ belong 
to the same orbit $\fX_a$) can be omitted from (\ref{bF13}).  

Given an orbit $\fX_a$, fix $x_i \in \fX_a$ and
define the differential\index{differential!meromorphic}
\begin{equation}  \label{bF16}
H_{\vs,a}(z) \, dz = \sum_{x_k \notin \fX_a} s_k \cdot
 \frac{1}{2} \left( G_{i,k}(z) \, dz + G_{\tau(i),\tau(k)}(z) \, dz \right) \ .
\end{equation}
The same differential is obtained if $x_i$
is replaced by $\tau(x_i)$.  Clearly $H_{\vs,a}(z) \, dz$
is holomorphic except at the points in $\fX$.  If $\fX_a = \{x_i\}$ is real
then $H_{\vs,a}(z) \, dz$ has a simple pole with residue $1-s_i$ at 
$x_i$;  if $\fX_a = \{x_i,x_{\tau(i)}\}$ is complex, 
it has simple poles with residue 
$\frac{1}{2} - s_i$ at $x_i$ and $x_{\tau(i)}$.  
In both cases it has a simple pole with residue $-s_k$ at each  
$x_k \notin \fX_a$.  Writing the $H_{\vs,a}(z)$ are in 
appropriate local coordinates, then for $x_i \in \fX_a$, 
formula (\ref{bF13}) becomes 
\begin{equation} \label{bF17}
  \Lambda_{x_i}(Q^{\veps},\vs) - 
                   \Lambda_{x_i}(Q,\vs)                   
   \ \cong \  \varepsilon_0 + \sum_{\ell=1}^{d} 
           H_{\vs,a}(t_{\ell}) \cdot \sigma_{\ell} \varepsilon_{\ell} \ .  
\end{equation}

We now ask about linear relations between the meromorphic differentials\index{differential!meromorphic} 
$H_{\vs,a}(z) \, dz$ and the holomorphic differentials\index{differential!holomorphic} 
$\omega_j = h_j(z) \, dz$.  Put $J = g + m_1 + m_2 - 1$, 
where $g$ is the genus of $\cC$.

\begin{proposition} \label{bP4} { \ }  
For each probability vector $\vs \in \cP^m$, the meromorphic differentials\index{differential!meromorphic}
\begin{equation*} 
H_{\vs,1}(z) \, dz\ , \ \ldots\ , \ H_{\vs,m_1 + m_2}(z) \, dz
\end{equation*} 
and the holomorphic differentials\index{differential!holomorphic} 
$\omega_1 = h_1(z) \, dz, \ldots, \omega_g = h_g(z) \, dz$ span a
vector space of dimension $J = g + m_1 + m_2 - 1$.  
If the orbits $\fX_1, \ldots, \fX_{m_1}$ are 
real and $\fX_{m_1 + 1}, \ldots, \fX_{m_1 + m_2}$ are complex, and if 
$x_i$ is a representative for $\fX_i$, $i = 1, \ldots, m_1+m_2$, 
then every linear relation among the $H_{\vs,i}(z) \, dz$ and the 
$\omega_j = h_j(z) \, dz$ is a consequence of the relation 
\begin{equation}  \label{bF18}
\sum_{i=1}^{m_1} s_i H_{\vs,i}(z) \, dz 
     + \sum_{i=m_1 + 1}^{m_1 + m_2} 2s_i H_{\vs,i}(z) \, dz \ = \ 0 \ .
\end{equation}
\end{proposition}

\begin{proof}  
Note that $J = 0$ iff $g = 0$ and $\fX$ consists of a single
point.  In that case $m_1 =1$, $m_2 =0$, and $H_{\vs,1}(z) dz \equiv 0$,  
so (\ref{bF18}) holds trivially.  Suppose $J > 0$, and let
\begin{equation} \label{bF19}
\sum_{i=1}^{m_1+m_2} c_i H_{\vs,i}(z) \, dz
       + \sum_{j=1}^g d_j h_j(z) \, dz \ = \ 0
\end{equation}
be an arbitrary relation. Considering the residues at the poles 
$x_1, \ldots, x_{m_1+m_2}$ and writing 
$^t\vec{c} = (c_1, \ldots, c_{m_1+m_2})$ we see that
\begin{equation*}
\left( \begin{array}{cccccc}
     1-s_1    & \hdots &    -s_1      &   -s_1       & \hdots &    -s_1     \\
    \vdots    & \ddots &  \vdots      &  \vdots      &        &   \vdots    \\
    -s_{m_1}  & \hdots & 1-s_{m_1}    & -s_{m_1}     & \hdots &  -s_{m_1}   \\
   -s_{m_1+1} & \hdots & -s_{m_1+1}   & \frac{1}{2} -s_{m_1 + 1} 
                                                     & \hdots &  -s_{m_1+1} \\
       \vdots &        &  \vdots      & \vdots       & \ddots & \vdots      \\
 -s_{m_1+m_2} & \hdots & -s_{m_1+m_2} & -s_{m_1+m_2} & \hdots & 
                                                 \frac{1}{2}-s_{m_1+m_2}
       \end{array} \right)       \cdot \vec{c} \ = \ \vec{0} \ .
\end{equation*}
Put $C = c_1 + \ldots + c_{m_1 + m_2}$.  The equations above imply  
$c_i = s_i C$ if $\fX_i$ is real, and $c_i = 2s_i C$ if $\fX_i$ is complex.  
Conversely, for any $C$, we obtain a solution by taking $c_i = s_i C$ when 
$\fX_i$ is real and $c_i = 2s_i C$ when $\fX_i$ is complex,   
since $\sum_{i=1}^{m_1} s_i + 
\sum_{i=m_1+1}^{m_1+m_2} 2s_i = 1$.  

Now consider the differential $\sum_{i=1}^{m_1+m_2} c_i H_{\vs,i}(z) \, dz$.
It is meromorphic\index{differential!meromorphic} with pure imaginary periods\index{pure imaginary periods} 
and at worst simple poles at the points in $\fX$.
We have just seen that it has residue $0$ at each $x_i$,
so it is everywhere holomorphic.  However, the only holomorphic differential\index{differential!holomorphic}
with pure imaginary periods\index{pure imaginary periods} is the $0$ differential, so   
\begin{equation*} 
\sum_{i=1}^{m_1+m_2} c_i H_{\vs,i}(z) \, dz \ = \ 0 \ .
\end{equation*} 
Inserting this in (\ref{bF19}), 
we see that $\sum_{j=1}^g d_j \omega_j = 0$.  
Since the $\omega_j$ are linearly independent, the $d_j$ must be $0$.   
(If $g = 0$, there are no holomorphic differentials, and the result is vacuously true.)
This yields the Proposition. 
\end{proof} 

\begin{corollary}  \label{bC1} 
For any open subinterval $I \subset \cC_v(\RR) \backslash \fX$, 
there exist points $t_1, \ldots, t_J \in I$ such that the matrix 
\begin{eqnarray} 
\cG_{\vs}(t_1, \ldots, t_J) \ := \
 \left( \begin{array}{cccccc}
   1   & H_{\vs,1}(t_1)       & & \hdots & & H_{\vs,1}(t_J) \\     
\vdots &   \vdots                 & &        & &            \vdots      \\
   1   & H_{\vs,m_1+m_2}(t_1) & & \hdots & & H_{\vs,m_1+m_2}(t_J) \\
   0   &   h_1(t_1)               & & \hdots & &     h_1(t_J)       \\
\vdots &   \vdots                 & &        & &      \vdots            \\
   0   &   h_g(t_1)               & & \hdots & &     h_g(t_J)  
     \end{array} \right)   \notag \\
                           \label{bF20}
\end{eqnarray} 
is nonsingular.        
\end{corollary}

\begin{proof}
After relabeling the $\fX_i$ if necessary, we can assume that $s_1 > 0$.
No nontrivial relation of the form
\begin{equation*}
\sum_{i=2}^{m_1+m_2} c_i H_{\vs,i}(z) \, dz
       + \sum_{j=1}^g d_j h_j(z) \, dz \ = \ 0
\end{equation*} 
can hold identically on $I$;  otherwise, it would hold identically on 
$\cC_v(\CC)$, contrary to Proposition \ref{bP4}.  Hence we can find 
$t_1, \ldots, t_J \in I$ such that submatrix of (\ref{bF20}) 
consisting of the last $J$ rows and columns is nonsingular. 
But then, row reducing (\ref{bF20}) by using the relation (\ref{bF18}) 
and replacing the first row with the corresponding linear combination of 
the first $m_1+m_2$ rows, that row becomes $(1, 0, \ldots, 0)$.  Hence
$\cG_{\vs}(t_1, \ldots, t_J)$ is nonsingular.  
\end{proof} 
\index{initial approximating functions $f_v(z)$!construction when $K_v \cong \RR$!independence of differentials|)} 


\section{ Proof of Theorem $\ref{RThm2}$} \label{RealCaseProofSection} 

In this section we will prove Theorem \ref{RThm2}. 
By assumption $K_v \cong \RR$. 
We are given a $K_v$-simple set $E_v$ which is bounded away from $\fX$ 
\index{$K_v$-simple!set} 
and has positive inner capacity, 
\index{capacity $> 0$} 
together with numbers $\varepsilon_v > 0$ and $0 < \cR_v < 1$, 
and an open set $U_v \subset \cC_v(\CC)$ such that  $E_v \cap U_v = E_v^0$.  
After shrinking $U_v$ we can assume it is bounded away from $\fX$,
and is stable under complex conjugation. 
Note that for each component $E_{v,\ell}$ of $E_v$ which not 
contained in $\cC_v(\RR)$, the interior $E_{v,\ell}^0$ is one of the 
connected components of $U_v$.   

The proof has several steps. 
We use the notation from \S\ref{Chap5}.\ref{OutlineSection}, \S\ref{Chap5}.\ref{IndependenceSection}.

\medskip{\bf Step 0.} 
\index{initial approximating functions $f_v(z)$!construction when $K_v \cong \RR$!Step 0: the case $E_v \cap \cC_v(\RR) = \phi$}  
If $E_v$ has no components contained in $\cC_v(\RR)$, 
then Theorem \ref{RThm2} follows by the same argument as Theorem \ref{CThm1}, 
using the assertions in (\cite{RR1}, Proposition 3.3.2)  
and (\cite{RR1}, Theorem 3.3.7) that deal with the case where $K_v \cong \RR$,
with $E_v$ and $\fX$ stable under complex conjugation,
and $\vbeta$ and $\vs$ being $K_v$-symmetric. 
\index{$K_v$-symmetric!vector} 

For the remainder of the proof, we will assume that at least one of 
the components of $E_v$ is a closed interval in $\cC_v(\RR)$.    

\medskip
\noindent{\bf Step 1.}
\index{initial approximating functions $f_v(z)$!construction when $K_v \cong \RR$!Step 1: reduction to short intervals|(}    
We first construct a $K_v$-simple set\index{$K_v$-simple!set}  
$E_v^{*} \subset E_v^0$, whose capacity is close to that of $E_v$,
\index{capacity}  
such that each real interval in $E_v^*$ is ``short'' 
in the sense of Definition \ref{ShortnessDef}. 
\index{short@`short' interval}  
 
\smallskip
 Since each point of $E_v$ is analytically accessible from $E_v^0$,
by Proposition \ref{IdentifyGreenProp} we have $G(z,x_i;E_v) = \Gbar(z,x_i,E_v^0)$
for each $x_i \in \fX$.  
Hence there is a compact subset $E_v^{**} \subset E_v^0$, such that for each 
$x_i, x_j \in \fX$ with $x_i \ne x_j$, 
\begin{eqnarray} \label{FK132} 
|V_{x_i}(E_v^{**}) - V_{x_i}(E_v)| &<& \varepsilon_v/2 
\quad \text{for each $i$} \ ,                         \\
|G(x_i,x_j;E_v^{**}) - G(x_i,x_j;E_v)| &<& \varepsilon_v/2 
\quad \text{for all $i \ne j$.}      \notag               
\end{eqnarray}
\index{Robin constant!properties of}
\index{Green's function!properties of}
Without loss we can assume that $E_v^{**}$ is $K_v$-simple (and in particular,
\index{$K_v$-simple!set} 
stable under complex conjugation). 
Indeed, since $E_v$ is $K_v$-simple, 
\index{$K_v$-simple!set}   
its quasi-interior $E_v^0$  has an exhaustion by $K_v$-simple sets.
\index{quasi-interior} 
\index{$K_v$-simple!set}  

 
To construct $E_v^{*}$, we will need a lemma.  
Fix a spherical metric $\|z,w\|_v$\index{spherical metric} on $\cC_v(\CC)$. 
For $p \in \cC_v(\CC)$ and $\delta > 0$, 
write $B(p,\delta)^- = \{ z \in \cC_v(\CC) : \|z,p\|_v < \delta \}$ .  
For $p \in \cC_v(\RR)$, 
let $I_p(\delta) = \{z \in \cC_v(\RR) : \|z,p\|_v < \delta \}$  
and let $\overline{I}_p(\delta) = \{z \in \cC_v(\RR) : \|z,p\|_v \le \delta \}$.  
     
\begin{lemma} \label{ALK30} 
Given a compact set $H \subset \cC_v(\CC) \backslash \fX$ 
and a number $\delta > 0$, let $H(\delta)$ be 
obtained from $H$ in any of the following ways:

$(A)$  $H(\delta) = \{ x \in \cC_v(\RR) : \|x,z\|_v \le h \ \text{for some $z \in H$} \}$;     

$(B)$  For some $p_1, \ldots, p_M \in H$, \quad 
    $H(\delta) = H \backslash \left( \bigcup_{k=1}^M B(p_k,\delta)^- \right)$;

$(C)$  For some $p_1, \ldots, p_M \in \cC_v(\RR) \backslash \fX$, \quad 
  $H(\delta) = H \bigcup \left( \bigcup_{k=1}^M \overline{I}_{p_k}(\delta) \right)$. 

\noindent{Then} for each $x_i \in \fX$, 
\begin{equation} \label{AFK129}
\lim_{\delta \rightarrow 0} \, V_{x_i}(H(\delta)) = V_{x_i}(H)  
\end{equation}
\index{Robin constant!archimedean}
and for each $x_i \ne x_j \in \fX$,  
\begin{equation} \label{AFK130}
\lim_{\delta \rightarrow 0} G(x_i,x_j;H(\delta)) 
          \ = \ G(x_i,x_j;H(\delta)) \ . 
\end{equation}
\end{lemma}

\begin{proof}  This follows from (\cite{RR1} , Corollary 3.1.16, p.149, 
Proposition 3.1.17, p.149, and Lemma 3.2.6, p.158).  
\end{proof} 

We next show that we can assume that each component of $E_v^{**}$ 
contained in $\cC_v(\RR)$ is a ``short'' interval.
\index{short@`short' interval|ii}  
Let $C(E_v,\fX)$ be the number gotten by taking $H = E_v$ 
in formula (\ref{FShort}), and put 
\begin{equation*}
B \ = \ B(E_v,\fX) \ = \ \min(1/C(E_v,\fX), 1/\sqrt{2C(E_v,\fX)}) \ .
\end{equation*}
Then any closed subinterval of $E_v^{**} \cap \cC_v(\RR)$,
with length at most $B$ under $\|z,w\|_v$,   
is ``short'' in the sense of Definition \ref{ShortnessDef}.
\index{short@`short' interval}    
Choose $p_1, \ldots, p_M \in E_v^{**} \cap \cC_v(\RR)$ 
such that $(E_v^{**} \cap \cC_v(\RR)) \backslash \{p_1, \ldots, p_M\}$ 
is composed of segments of length at most $B$. Lemma \ref{ALK30}(B) 
then shows that by deleting small open balls about the $p_k$, 
we can find a $K_v$-simple compact set $E_v^{*} \subset E_v^{**}$ 
\index{$K_v$-simple!set} 
such that each real interval in $E_v^{*}$ is ``short'', and   
\index{short@`short' interval}  
\begin{eqnarray} \label{FK133} 
|V_{x_i}(E_v^{*}) - V_{x_j}(E_v)| &<& \varepsilon_v 
\quad \text{for each $i$} \ ,                         \\
|G(x_i,x_j;E_v^{*}) - G(x_i,x_j;E_v)| &<& \varepsilon_v
\quad \text{for all $i \ne j$\ .}  \notag                   
\end{eqnarray}  
\index{Green's function!archimedean!properties of}
\index{Robin constant!archimedean!properties of} 
This set $E_v^{*}$ meets our needs.
\index{initial approximating functions $f_v(z)$!construction when $K_v \cong \RR$!Step 1: reduction to short intervals|)} 

\medskip
\noindent{\bf Step 2.}  The choice of $t_1, \ldots, t_d$.
\index{initial approximating functions $f_v(z)$!construction when $K_v \cong \RR$!Step 2: the choice of $t_1, \ldots, t_d$|(}  

\medskip
As in \S\ref{Chap5}.\ref{IndependenceSection}, put $J = g + m_1 + m_2 - 1$.  Then $J \ge 0$, 
with $J = 0$ if and only if $g = 0$ and $\fX = \{x_1\}$.  In that situation 
every divisor of degree $0$ is principal, and the variation in the logarithmic leading
coefficient at $x_1$ is accomplished by scaling alone 
(e.g. via $\varepsilon_0$ in (\ref{AFK147}) below).  
If $J = 0$, take $d = 0$ and ignore all constructions 
related to points $t_{\ell}$ in the rest of the proof.     

Assume now that $J \ge 1$.
Fix a closed interval $I \subset (E_v^0 \cap \cC_v(\RR)) \backslash E_v^{*}$
with positive length.  This interval will play an important role in the construction
below;  the points $t_i$, and the intervals we construct below, will belong to it.

Fix a local coordinate function $z$ on a neighborhood of $I$,
in such a way that $z$ is real-valued on $I$.  
Write the differentials\index{differential!meromorphic}\index{differential!holomorphic} 
$H_{\vs,i}(z) \, dz$ and $\omega_j(z) = h_j(z) \, dz$ 
from \S \ref{Chap5}.\ref{IndependenceSection} in terms of $z$.  
Translations of points, $z \mapsto z + \varepsilon$, will also be
understood relative to this coordinate.  

Let $\cP^m_v \subset \cP^m$ denote the set of  {\em $K_v$-symmetric} 
\index{$K_v$-symmetric!probability vector}
real probability vectors.  If $\vs_0 \in \cP^m_v$ and $\rho > 0$,
let $B(\vs_0,\rho) \subset \RR^m$ denote the open Euclidean ball about $\vs_0$ 
with radius $\rho$. 
 
For a given $\vs_0 \in \cP^m_v$, Corollary \ref{bC1} shows we can 
find points $t_1, \ldots, t_J$ in the interior of $I$ such that the 
matrix $\cG_{\vs_0}(\vec{t})$ defined there is nonsingular. 
Fixing  $\vec{t} = (t_1, \ldots, t_J)$, 
the function $w_{\vec{t}}(\vs) := \text{det}(\cG_{\vs})(\vec{t})$ 
is continuous for $\vs \in \cP^m_v$, so there is a $\rho = \rho(\vs_0) > 0$ 
such that $\text{det}(\cG_{\vs})(\vec{t})$ is nonsingular for all 
$\vs \in \cP^m_v \bigcap B(\vec{s_0},\rho)$.   Moreover, 
$|\text{det}(\cG_{\vs}(\vec{t})|$ is uniformly bounded away from $0$
if we restrict to $\vs \in \cP^m_v \bigcap B(\vec{s_0},\rho/2)$. 
Since $\cP^m_v$ is compact, we can cover it with a finite number of balls 
$B(\vec{s_0},\rho(\vec{s_0})/2)$.  Let $\cT = \{t_1, \ldots, t_d\}$ be the 
union of the sets $\{t_1, \ldots, t_J\}$ associated to these $\vs_0$.

The set $\tE_v$ will have the form 
\begin{equation}
E_v^{*}(h) \ := \ E_v^{*} \cup (\bigcup_{\ell=1}^d [t_{\ell}-h,t_{\ell}+h])  \label{AFK134}
\end{equation} 
for a suitably small $h$ which we will construct below, 
where the intervals $[t_{\ell}-h,t_{\ell}+h]$ (defined in terms of the local 
coordinate function $z$) are contained in the interior of $I$
and are pairwise disjoint.  

\vskip .1 in
However, before proceeding further, we note two facts concerning $\cT$: 
 
First, given a square matrix $\cG$ with real entries, 
denote its $L^2$ operator norm by 
\begin{equation*}
\|\cG\| \ = \ \max_{|\vec{x}| = 1} \, |\cG  \vec{x}| \ .
\end{equation*}
By the construction of $\cT$, there is a constant $B_1$ such that for each 
$\vs \in \cP^m_v$, we can find $J$ points in $\cT$ such that the 
corresponding matrix $\cG_{\vs}(\vec{t})$ from Corollary \ref{bC1}
is nonsingular, and 
\begin{equation} \label{AFK133}
\|\cG_{\vs}(\vec{t})^{-1}\| \ \le \ B_1 \ .
\end{equation} 

Second, for all sufficiently small $h$, there are \`a priori bounds\index{mass bounds|ii} on 
the relative mass which the $(\fX,\vs)$-equilibrium distribution of $E_v^{*}(h)$
gives to each segment $e_{\ell}(h) := [t_{\ell}-h,t_{\ell}+h]$. 
This is a consequence of potential-theoretic results in 
Appendix A.\ref{ArchMassBoundsSection}, as follows.       
   
If $I = [a,b]$, fix a number $r_0 > 0$ small enough that $5r_0$ is less 
than the minimum of the distances $|t_k-t_{\ell}|$ 
for $t_k \ne t_{\ell} \in \cT$, and 
such that $2r_0$ is less than the minimum of the distances to the endpoints,  
$|t_{\ell}-a|$, $|t_{\ell}-b|$, for each $t_{\ell} \in \cT$.  Thus the intervals 
$e_{\ell}(2r_0) \subset U_v$ are bounded away from each other, 
and are contained in the interior of $I$ so they are bounded away from $E_v^{*}$.  
We will also require that $r_0$ be small enough that each  
$e_{\ell}(2r_0)$ is ``short'' in the sense of Definition \ref{ShortnessDef},
\index{short@`short' interval}  
permitting the construction of oscillating pseudopolynomials.  

Consider $E_v^{*}(2r_0)$:  it is $K_v$-simple 
\index{$K_v$-simple!set} 
(so in particular $\cC_v(\CC) \backslash E_v^{*}(2r_0)$ is connected), 
and each of the segments $e_{\ell}(2r_0)$ is a component of $E_v^{*}(2r_0)$.  
Put 
\begin{equation}               \label{AFK135}    
B_2 \ = \ \min_{1 \le i \le m} \, \min_{1 \le \ell \le d} \, 
     \min_{z \in e_{\ell}(2r_0)} \, 
              G(z,x_i;E_v^{*}(2r_0) \backslash e_{\ell}(2r_0)) \ . 
\end{equation}
Then $B_2 > 0$, 
and  by the monotonicity of Green's functions, 
\index{Green's function!monotonic}
for each $0 < h \le 2r_0$ each $x_i$ and each $\ell$, 
we have $G(z,x_i;E_v^{*}(h) \backslash e_{\ell}(h)) \ge B_2$ on $e_{\ell}(h)$.      

By Proposition \ref{BPropF1}, $V_{\fX,\vs}(E_v^{*})$ is a continuous function of
\index{Robin constant!archimedean!archimedean $(\fX,\vs)$} 
$\vs$, so there is a finite upper bound $B_3$ for 
the values $V_{\fX,\vs}(E_v^{*})$ as $\vs$ ranges over the compact set $\cP^m_v$.  
Trivially 
\begin{equation*}
V_{\fX,\vs}(E_v^{*}(h)) \ \le \ 
    V_{\fX,\vs}(E_v^{*}(h) \backslash e_{\ell}(h)) 
\ \le \ V_{\fX,\vs}(E_v^{*}) \ \le \ B_3
\end{equation*}
for all $h$, $\ell$, and $\vs$.
\index{Robin constant!archimedean!archimedean $(\fX,\vs)$} 

By Lemma \ref{bL2}, there is a constant $A > 0$ such that for all $x_i \in \fX$, 
all $\ell = 1, \ldots, d$, and all sufficiently small $h > 0$, 
\begin{equation*}
-\log(h) - A \ < \ V_{x_i}(e_{\ell}(h)) \ < \ -\log(h) + A \ .
\end{equation*} 
\index{Robin constant!archimedean} 
In particular for all sufficiently small $h$, and all $x_i \in \fX$,
\begin{equation*}
V_{x_i}(e_{\ell}(h))  \ > \ V_{x_i}(E_v^{*}) \ > \ 
V_{x_i}(E_v^{*}(h) \backslash e_{\ell}(h)) \ ,
\end{equation*} 
validating the hypothesis of Lemma \ref{ALJ29}.
 
For a given $\vs \in \cP^m_v$, let $\mu_{\fX,\vs,h}$
be the $(\fX,\vs)$-equilibrium distribution of $E_v^{*}(h)$.  
Then by Lemma  \ref{bL3},\index{mass bounds|ii} there is a constant $C$ such that 
for each sufficiently small $h$,  
\begin{equation*}      
\mu_{\fX,\vs,h}(e_{\ell}(h)) 
    \ \le \ \frac {V_{\fX,\vs}(E_v^{*}(h)) + C}
                                 {V_{\fX,\vs}(e_{\ell}(h)) + C} 
    \ \le \ \frac {B_3 + C} {-\log(h) - A + C}  \ . 
\end{equation*}
\index{Robin constant!archimedean!bounds for} 
Likewise, by Lemma \ref{ALJ29}, for sufficiently small $h$, 
\begin{eqnarray*}      
\mu_{\fX,\vs,h}(e_{\ell}(h)) 
    &\ge& \frac {(B_2)^2}{ 2 (V_{\fX,\vs}(E_v^{*}(h) \backslash e_{\ell}(h)) + C)
                             (V_{\fX,\vs}(e_{\ell}(h)) + C + 2B_2) }   \\
    &\ge& \frac {(B_2)^2} {2 (B_3 + C)(-\log(h) + A + C + 2B_2)}  \ .                            
\end{eqnarray*} 
Put 
\begin{equation} \label{FDDef}
B_4  \ = \ \frac{3d \cdot (B_3+C)^2} {(B_2)^2} \ .
\end{equation}  
Then there is an $h_0 > 0$ such that for each $\vs \in \cP^m_v$, 
each $\ell = 1, \ldots, d$, and each $0 < h \le h_0$, 
\begin{equation}                               \label{AFK136} 
\frac {\mu_{\fX,\vs,h}(\bigcup_{k=1}^d e_k(h))} 
          {\mu_{\fX,\vs,h}(e_{\ell}(h))}
    \ \le \ B_4 \ .   
\end{equation}
\index{Robin constant!archimedean!bounds for} 

The parameter $h$ will be chosen in Step 4 below.
Given $0 < h \le h_0$, 
we will put  $\tE_v = E_v^{*}(h) = \bigcup_{i=1}^D \tE_{v,\ell}$, 
where $\tE_{v,\ell} = [t_\ell-h,t_\ell+h]$ for $\ell =1, \ldots, d$ and 
$\tE_{v,d+1}, \ldots, \tE_{v,D}$ are the components of $E_v^{*}$.  We can apply 
Theorem \ref{PotPart} to $E_v^{*}(h)$, constructing $(\vX,\vs)$-pseudopolynomials 
\index{pseudopolynomial!$(\fX,\vs)$} 
$Q(z)$ with large oscillations on the real components of $E_v^{*}(h)$.
\index{initial approximating functions $f_v(z)$!construction when $K_v \cong \RR$!Step 2: the choice of $t_1, \ldots, t_d$|)} 

\medskip
\noindent{\bf Step 3.}  The choice of $r$.
\index{initial approximating functions $f_v(z)$!construction when $K_v \cong \RR$!Step 3: the choice of $r$|(} 
\medskip

In later stages of the construction, 
we will need to move some of the roots of the 
pseudopolynomials $Q(z)$, in order to make their divisors principal and 
\index{pseudopolynomial} 
vary their logarithmic leading coefficients.\index{logarithmic leading coefficients} 
We now define a number $r$ which governs how far we 
can move the roots. 

Let $0 < \cR_v < 1$ be the oscillation bound required in the Theorem.  
Fix a number $\tcR_v$ with $\cR_v < \tcR_v < 1$, 
and put $\Delta_1 = \log(\tcR_v) - \log(\cR_v)$.
Recall that $U_v$ is the open set for which $E_v^0 = E_v \cap U_v$.
Consider the set $E_v^{*} \cup I$, which is contained in $U_v$.   
The Green's functions $G(z,x_i;E_v^{*} \cup I)$ are continuous, 
\index{Green's function!continuous on boundary}
and are positive in the complement of $E_v^{*} \cup I$.  
Since $\partial U_v$ is compact and disjoint from $E_v^{*} \cup I$, 
there is a $\Delta_2 > 0$ such that for each $x_i \in \fX$, 
and $z \notin E_v^{*} \cup I$, 
we have $G(z,x_i;E_v^{*} \cup I)  \ge \Delta_2$.  
For each $x_i \in \fX$, 
recall that there is a $\cC^{\infty}$ function $\eta_i(z,w)$ such that 
$\log([z,w]_{x_{\ell}}) = \log(|z-w|) + \eta_i(z,w)$ on $I \times I$.

\vskip .1 in
Let $r$ be a number in the range  
\begin{equation} \label{FRSiz1}
0 \ < \ r \ \le \ r_0 \ ,
\end{equation} 
(so in particular, the intervals $e_{\ell}(r)$ are pairwise disjoint, 
contained in $I$, and are ``short'' in the sense of Definition \ref{ShortnessDef}),
\index{short@`short' interval}  
which satisfies 
\begin{equation} \label{FRSiz2}
 r \ < \ \min(\Delta_1, \Delta_2)/3  \ ,
\end{equation} 
and which is small enough that for each $x_i \in \fX$,    
the following six conditions hold:     

\medskip
(1)  for each $t_{\ell} \in \cT$, and each $w \in E_v^{*}$, 
\begin{equation}                                        \label{AFK137}
\left| \max_{z \in e_{\ell}(r)} \log([z,w]_{x_i}) -  
  \min_{z \in e_{\ell}(r)} \log([z,w]_{x_i}) \right| 
                \ < \ \Delta_1/2 \ ;
\end{equation} 

(2)  for each $t_{\ell} \in \cT$, each $t_k \ne t_\ell \in \cT$, 
and each $w \in e_k(2r_0)$, 
\begin{equation}                                        \label{AFK138}
\left| \max_{z \in e_{\ell}(r)} \log([z,w]_{x_i}) -  
  \min_{z \in e_{\ell}(r)} \log([z,w]_{x_i}) \right| 
                 \ < \ \Delta_1/2 \ ;
\end{equation}  

(3)  for each $t_{\ell} \in \cT$,  
\begin{equation}                                        \label{AFK139}
\left| \max_{z,w \in e_{\ell}(2r)} \eta_i(z,w)-  
  \min_{z,w \in e_{\ell}(2r)} \eta_i(z,w) \right| \ < \ \Delta_1/2 \ ;
\end{equation} 

(4) for each $x_i \ne x_j \in \fX$,  
and for each $t_\ell \in \cT$ and each $z \in e_{\ell}(2r)$,
\begin{equation}                                       
\ \Big| \big(\log([x_i,z]_{x_j})
    - \log([x_i,t_{\ell}]_{x_j})\big) - (z-t_{\ell}) \cdot
  \left. \frac{d}{dt} \log([x_i,t]_{x_j}) 
                 \right|_{t=t_{\ell}} \Big|              
\ < \  \frac{r}{24 B_1 B_4 \sqrt{J+1}} \ ; \label{AFK140}
\end{equation} 

(5) for each $\omega_j = h_j(z) \, dz$, 
and for each $t_\ell \in \cT$ 
and each $z \in e_{\ell}(2r)$,  
\begin{equation}  \label{AFK141}  
\left| h_j(z) - h_j(t_{\ell}) \right| \ < \ \frac{1}{24 B_1 B_4 \sqrt{J+1}} \ ;
\end{equation}        

\medskip
(6) for all $z \in \partial U_v$ and all $w_1, w_2 \in E_v^{*} \cup I$ 
with $\|w_1, w_2\|_v \le r$, 
\begin{equation} \label{AFK141D} 
|\log([z,w_1]_{x_i}) - \log([z,w_2]_{x_i})| \ < \ \Delta_2/3 \ .
\end{equation} 

\medskip
Conditions (1), (2), and (3) hold for all sufficiently small $r$ 
by the continuity of the functions $\log([z,w]_{x_i})$ and  $\eta_i(z,w)$.  
Condition (4) holds for all sufficiently small $r$ since the functions 
$\log([x_i,t]_{x_j})$ are $\cC^{\infty}$ on the intervals
$e_{\ell}(2r_0)$. Condition (5) holds for all sufficiently small $r$ 
by the continuity of the $h_j(z)$.  Condition (6) holds for all sufficiently 
small $r$ since the functions $\log([z,w]_{x_i})$ are uniformly continuous
for $(z,w) \in \partial U \times (E_v^{*} \cup I)$.  

Fix such an $r$.  
By our assumptions on $E_v^{*}$ and $r$, for each $0 < h \le r$ 
we can apply Theorem \ref{PotPart} to the set $E_v^{*}(h)$,
constructing $(\fX,\vs)$-pseudopolynomials $Q(z)$.  
\index{pseudopolynomial!$(\fX,\vs)$}  
Recall if $\vs \in \cP^m_v$, then 
$G_{\fX,\vs}(z,E_v^{*}(h)) = \sum_{i=1}^m s_i G(z,x_i;E_v^{*}(h))$.

\begin{proposition} \label{APK31} 
Let $r$ be as above.  Given $0 < h \le r$, 
put $\tE_v = E_v^{*}(h) = E_v^{*} \cup \bigcup_{i=1}^d [t_{\ell}-h,t_{\ell}+h]$.
Let $\veps = (\varepsilon_0,\varepsilon_1, \ldots, \varepsilon_d) \in \RR^{d+1}$ 
be such that  $|\varepsilon_{\ell}| \le r$ for each $\ell$, 
and put $\tE_v(\veps) = E_v^{*} \cup \bigcup_{\ell=1}^d 
[t_{\ell}+\varepsilon_{\ell}-h,t_{\ell}+\varepsilon_{\ell}+h]$.  
 Write $\tE_{v,1}, \ldots, \tE_{v,D}$ for the components of $\tE_v$, where
$\tE_{v,\ell} = [t_{\ell}-h,m_{\ell}+h]$ for $\ell = 1, \ldots, d$, 
and $\tE_{v,d+1}, \ldots, \tE_{v,D}$
are the components of $E_v^{*}$.
Write $\tE_{v,1}(\veps), \ldots, \tE_{v,D}(\veps)$ 
for the corresponding components of $\tE_v(\veps)$. 

Fix $\vs \in \cP^m_v$.  Given $\vn = (n_1, \ldots, n_D) \in \NN^D$, 
let $Q(z) = Q_{\vn}(z) = C \cdot \prod_{k=1}^N [z,\alpha_k]_{\fX,\vs}$
be an $(\fX,\vs)$-pseudopolynomial with $n_\ell$ roots in $\tE_{v,\ell}$,
\index{pseudopolynomial!$(\fX,\vs)$}  
for $\ell = 1, \ldots, D$.    
Suppose that  $Q(z)$ oscillates $n_\ell$ times from $\tcR^N$ to $0$ to $\tcR_v^N$   
on each component $\tE_{v,\ell}$ of $\tE_v$ contained in $\cC_v(\RR)$
and that $|\frac{1}{N} \log(Q(z)) - G_{\fX,\vs}(z,\tE_v)| < r$
for each $z \in \partial U_v$. Put  
\begin{equation} \label{FNCh}
Q^{\veps}(z) \ := \ e^{N\varepsilon_0} \cdot C 
    \cdot \prod_{x_i \in E_v^{*}} [z,\alpha_k]_{\fX,\vs}
    \cdot \prod_{\ell=1}^d \, \prod_{\alpha_k \in e_{\ell}(r)} 
             [z,\alpha_k + \varepsilon_{\ell}]_{\fX,\vs} \ .
\end{equation}
Then $Q^{\veps}(z)$ has $n_{\ell}$ roots in each $\tE_{v,\ell}(\veps)$, 
$Q^{\veps}(z)$ oscillates $n_{\ell}$ times from $\cR^N$ to $0$ to $\cR_v^N$ 
on each component $\tE_{v,\ell}(\veps)$ contained in $\cC_v(\RR)$,
and $\{z \in \cC_v(\CC) : Q^{\veps}(z) \le 1\} \subset U_v$.
\end{proposition}

\begin{proof}  
By its construction, $Q^{\veps}(z)$ has $n_\ell$ roots in $\tE_{v,\ell}(\veps)$, 
for each $\ell = 1, \ldots, D$.   
Note that since $|\varepsilon_\ell| \le r < r_0$, 
if $\alpha_k \in [t_{\ell}-h,t_{\ell}+h]$, 
then $\alpha_k + \varepsilon_\ell \in [t_{\ell}-2r_0,t_{\ell}+2r_0]$.  
We first show that $Q^{\veps}(z)$ oscillates $n_{\ell}$ 
times from $\cR^N$ to $0$ to $\cR_v^N$ 
on each real component $\tE_{v,\ell}(\veps)$.  There are two cases to consider.
 
\smallskip
First suppose that $1 \le \ell \le d$, 
so that $\tE_{v,\ell} = [t_{\ell}-h,t_{\ell}+h]$ and 
$\tE_{v,\ell}(\veps) = [t_\ell+\varepsilon_\ell-h,t_\ell+\varepsilon_\ell+h]$.  
Let $z_0 \in [t_{\ell}-h,t_{\ell}+h]$ be a point where $Q(z_0) = \tcR_v^N$. 
By (\ref{AFK138}) and (\ref{AFK139}),   
\begin{eqnarray*}
\lefteqn{ \Big|\frac{1}{N}\log(Q^{\veps}(z_0+\varepsilon_{\ell})) 
        - \frac{1}{N}\log(Q(z_0))\Big| }   \\
 &\le& |\varepsilon_0| + \frac{1}{N} \sum_{\alpha_k \in E_v^{*}} \, \sum_{i=1}^m 
            s_i \big|\log([z_0+\varepsilon_{\ell},\alpha_k]_{x_i}) 
                              - \log([z_0,\alpha_k]_{x_i})\big|    \\  
 & & + \ \frac{1}{N} \sum^d_{\substack{ j=1 \\ j \ne \ell }} \, 
        \sum_{\alpha_k \in [t_j-h,t_j+h]} \, \sum_{i=1}^m 
       s_i \big|\log([z_0+\varepsilon_{\ell},\alpha_k+\varepsilon_j]_{x_i}) 
                 - \log([z_0,\alpha_k]_{x_i})\big|              \\
 & & + \ \frac{1}{N}  \sum_{\alpha_k \in [t_{\ell}-h,t_{\ell}+h]} \, \sum_{i=1}^m 
       s_i \big|\eta_i(z_0+\varepsilon_{\ell},\alpha_k+\varepsilon_{\ell})  
                 - \eta_i(z_0,\alpha_k)\big|                        \\  
& < & r + \frac{1}{N} \cdot N \cdot \sum_{i=1}^m s_i 
           \cdot r \ = \ 2r \ < \ \Delta_1 \ .
\end{eqnarray*}        
Since $\frac{1}{N} \log(Q(z_0)) = \log(\tcR_v)$ and 
$\log(\tcR_v) - \Delta_1 = \log(\cR_v)$, it follows that            
$\frac{1}{N} \log(Q^{\veps}(z_0+\varepsilon_\ell)) \ge \log(\cR_v)$.         
If $Q(z)$ takes the values $\tcR_v^N$, $0$, $\tcR_v^N$ at successive points 
$z_0$, $\alpha_k$, and $z_0^{\prime}$ of $[t_{\ell}-h,t_{\ell}+h]$, then 
$Q^{\veps}(z)$ takes values $\ge \cR_v^N$, $0$, $\ge \cR_v^N$ at 
the points $z_0+\varepsilon_\ell$, $\alpha_k+\varepsilon_\ell$, 
and $z_0^{\prime}+\varepsilon_\ell$
in $[t_{\ell}+\varepsilon_\ell-h,t_{\ell}+\varepsilon_\ell+h]$. 

\smallskip
Next suppose $\ell \ge d+1$, 
so $\tE_{v,\ell}(\veps) = \tE_{v,\ell}$ is a component of $E_v^{*}$  
contained in $\cC_v(\RR)$.  
Let $z_0 \in \tE_{v,\ell}$ be a point where $Q(z_0) = \tcR_v^N$.  
By (\ref{AFK137}), 
\begin{eqnarray*}
\lefteqn{ \Big|\frac{1}{N}\log(Q^{\veps}(z_0)) 
     - \frac{1}{N}\log(Q(z_0))\Big| } \\
  & \le & |\varepsilon_0| + \frac{1}{N} \sum_{\ell=1}^d \, 
       \sum_{\alpha_k \in e_{\ell}(h)} \, \sum_{i=1}^m 
            s_i \big|\log([z_0,\alpha_k]_{x_i}) 
                 - \log([z_0,\alpha_k+\varepsilon_\ell]_{x_i})\big| \\
& < & r + \frac{1}{N} \cdot N \cdot \sum_{i=1}^m s_i 
           \cdot r \ = \ 2 r \ < \ \Delta_1 \ .
\end{eqnarray*}
By the same argument as before, we conclude that 
$\frac{1}{N} \log(Q^{\veps}(z_0)) \ge \log(\tR_v)$.  Hence if $Q(z)$
takes the values $\tcR_v^N$, $0$, $\tcR_v^N$ at successive points 
$z_0$, $\alpha_k$, $z_0^{\prime}$ in $\tE_{v,\ell}$, then $Q^{\veps}(z)$ 
takes values $\ge \cR_v^N$, $0$, $\ge \cR_v^N$ at those points.   

\smallskip
Finally, we show that $\{z \in \cC_v(\CC) : Q^{\veps}(z) \le 1\} \subset U_v$.
The function 
$\frac{1}{N} \log(Q^{\veps}(z))$ is harmonic on 
$\cC_v(\CC) \backslash (\tE_v(\veps) \cup \fX)$.
For each $x_i \in \fX$, if $s_i > 0$ then 
$\lim_{z \rightarrow x_i} \frac{1}{N} \log(Q^{\veps}(z)) = \infty$,  
while if $s_i = 0$ then $\frac{1}{N}\log(Q^{\veps}(z))$
has a removable singularity at $x_i$.  Hence if we show that 
$\frac{1}{N}\log(Q^{\veps}(z)) > 0$ on $\partial U_v$, the maximum 
principle for harmonic functions will give 
$\frac{1}{N}\log(Q^{\veps}(z)) > 0$ outside $U_v$.  

Let $z_0$ be a point of $\partial U_v$.  
Since $\tE_v = E_v^{*}(h) \subset E_v^{*} \cup I$, 
the monotonicity of Green's functions 
\index{Green's function!monotonic}
shows that $G(z_0,x_i;\tE_v) \ge G(z_0,x_i;E_v^{*} \cup I) \ge \Delta_2$, 
and consequently 
\begin{equation*} 
G_{\fX,\vs}(z_0,\tE_v) 
\ = \ \sum_{i=1}^m s_i G(z_0,x_i;\tE_v) \ \ge \ \Delta_2 \ .
\end{equation*} 
Since 
$|G_{\fX,\vs}(z_0,\tE_v) - \frac{1}{N} \log(Q(z_0))| < r$ and $r < \Delta_2/3$,  
it follows that $\frac{1}{N} \log(Q(z_0)) > 2 \Delta_2/3$.  
On the other hand, by (\ref{AFK141D}) 
\begin{eqnarray*}
\lefteqn{ \Big|\frac{1}{N}\log(Q^{\veps}(z_0)) 
     - \frac{1}{N}\log(Q(z_0))\Big| } \\
  & \le & |\varepsilon_0| + \frac{1}{N} \sum_{\ell=1}^d \, 
       \sum_{\alpha_k \in e_{\ell}(h)} \, \sum_{i=1}^m 
            s_i \big|\log([z_0,\alpha_k]_{x_i}) 
                 - \log([z_0,\alpha_k+\varepsilon_\ell]_{x_i})\big| \\
& < & r + \frac{1}{N} \cdot N \cdot \sum_{i=1}^m s_i 
           \cdot r \ = \ 2 r \ \le \ 2 \Delta_2/3 \ .
\end{eqnarray*}
Hence $\frac{1}{N}\log(Q^{\veps}(z_0)) > 0$.  

Since $\frac{1}{N}\log(Q^{\veps}(z)) > 0$ on $\partial U_v$, 
the maximum principle for harmonic functions
(applied to $-\frac{1}{N} \log(Q^{\veps}(z))$ on 
$\cC_v(\CC) \backslash (U_v \cup \fX)$)
shows that $\frac{1}{N}\log(Q^{\veps}(z)) > 0$ 
for all $z \notin U_v$.  Thus  
$\{z \in \cC_v(\CC) : Q^{\veps}(z) \le 1\} \subset U_v$.
\end{proof} 
\index{initial approximating functions $f_v(z)$!construction when $K_v \cong \RR$!Step 3: the choice of $r$|)}    

\medskip
\noindent{\bf Step 4.}  The choice of $h$ and the construction of $\tE_v$.
\index{initial approximating functions $f_v(z)$!construction when $K_v \cong \RR$!Step 4: the construction of $\tE_v$|(} 
\medskip

We will now choose the number $h$, and hence determine the set $\tE_v$ in 
the Theorem.  Let  $0 < h < r$ be small enough that the following three conditions are satisfied:

\medskip
(1)  for the set $E_v^{*}(h) = E_v^{*} \cup \big(\bigcup_{j=1}^d e_{\ell}(h)\big)$,
\begin{equation}\label{AFK142}  
\left\{ \begin{array}{ll} 
|V_{x_i}(E_v^{*}(h)) - V_{x_i}(E_v^{*})| < \varepsilon_v/3 & 
   \text{for each $i$\ ,}                         \\ 
G(x_i,x_j;E_v^{*}(h)) - G(x_i,x_j;E_v^{*})| < \varepsilon_v/3 &
  \text{for all $i \ne j$\ .}  
   \end{array} \right.
\end{equation}
\index{Green's function!archimedean!properties of} 
\index{Robin constant!archimedean!properties of} 

(2)  for all $i \ne j$, all $t_{\ell} \in \cT$, all $\varepsilon_{\ell}$ with 
$|\varepsilon_\ell| \le r$, and all $z \in [t_{\ell}-h,t_{\ell} + h]$,  
\begin{equation}                                      \label{AFK143}
\left| \log([x_i,z+\varepsilon_{\ell}]_{x_j}) - 
            \log([x_i,t_{\ell}+\varepsilon_{\ell}]_{x_j}) \right| 
                 \ < \ \frac{r}{24 B_1 B_4 \sqrt{J+1}} \ .
\end{equation}

(3) $h \le h_0$, so the mass bounds (\ref{AFK136}) are valid.  

\medskip
Condition (1) holds for all small $h$ by Lemma \ref{ALK30}.C.  
Condition (2) holds for all small $h$ by the continuity of the 
$\log([x_i,t]_{x_j})$.  And, condition (3) clearly holds for all 
small $h$.                    

\medskip
Fix $h$ satisfying (1), (2) and (3), 
and write $\tE_{v,\ell} = e_{\ell}(h) =  [t_{\ell}-h,t_{\ell}+h]$, 
for $\ell = 1, \ldots, d$.  Put 
\begin{equation}     \label{AFK144} 
\tE_v \ = \ E_v^{*} \cup \big(\bigcup_{\ell=1}^d [t_{\ell}-h,t_{\ell}+h]\big) \ .
\end{equation}
Then $\tE_v \subset E_v^0$, and $\tE_v$ is $K_v$-simple;
\index{$K_v$-simple!set} 
in particular $\cC_v(\CC) \backslash \tE_v$ is connected.   
By (\ref{FK132}), (\ref{FK133}) and (\ref{AFK142}),
\begin{equation} \label{AFK145}  
\left\{ \begin{array}{ll} 
|V_{x_i}(\tE_v) - V_{x_i}(E_v)| < \varepsilon_v & 
   \text{for each $i$,}                         \\ 
G(x_i,x_j;\tE_v) - G(x_i,x_j;E_v)| < \varepsilon_v &
  \text{for all $i \ne j$.}  
   \end{array} \right.
\end{equation}
\index{Green's function!archimedean!properties of} 
\index{Robin constant!archimedean!properties of} 
as required by Theorem \ref{RThm2}.  
Moreover all the intervals making up $\tE_v$ are ``short'' in the sense of
Definition \ref{ShortnessDef}, and Proposition \ref{APK31} applies to $\tE_v$.   
\index{short@`short' interval}
\index{initial approximating functions $f_v(z)$!construction when $K_v \cong \RR$!Step 4: the construction of $\tE_v$|)}

\medskip
\noindent{\bf Step 5.}  The total change map.
\index{initial approximating functions $f_v(z)$!construction when $K_v \cong \RR$!Step 5: study of the total change map|(} 
\medskip

We now begin to address problem of constructing rational functions with 
large oscillations and specified logarithmic leading coefficients.
\index{logarithmic leading coefficients}   

Fix $\vs \in \cP^m_v \bigcap \QQ^m$.  
Consider a special $(\fX,\vs)$-pseudopolynomial for $\tE_v$,
\index{pseudopolynomial!special} 
\begin{equation}     \label{AFK146} 
Q(z) \ = \ Q_{\vec{n}}(z) \ = \ C 
           \cdot \prod_{k=1}^N [z,\alpha_k]_{\fX,\vs}  \ ,
\end{equation}   
as constructed in Theorem \ref{PotPart}.  
In particular $\|Q(z)\|_{\tE_v} = 1$, and for each $\ell = 1, \ldots, D$, 
if $\vn = (n_1, \ldots, n_D)$ then $Q(z)$ has $n_\ell$ roots in $\tE_{v,\ell}$. 
Recalling that $J = g+ m_1 + m_2 -1$, fix a set of $J$ points in $\cT$ 
such that the matrix $\cG_{\vs}(\vec{t})$ 
in Corollary \ref{bC1} is nonsingular and satisfies 
\begin{equation*} 
\|\cG_{\vs}(\vec{t})^{-1}\| \ \le \ B_1 \ , 
\end{equation*} 
as in (\ref{AFK133});  the construction of $\cT$ shows that this can be done.   
Without loss, suppose these points are $t_1, \ldots, t_J$, and put
$E_v^{00} = \tE_v \backslash (\bigcup_{\ell=1}^J \tE_{v,\ell})$.   
We will move the roots of $Q(z)$,
constructing a pseudopolynomial $Q^{\veps}(z)$ as in (\ref{FNCh}); 
\index{pseudopolynomial}  
however, instead of 
moving the roots in all the intervals $\tE_{v,\ell}, \ \ell = 1, \ldots, d$ 
we will only move the ones in $\tE_{v,1}, \ldots, \tE_{v,J}$.   
Equivalently,  take $\varepsilon_{\ell} = 0$ for $\ell = J+1, \ldots, d$ 
in Proposition \ref{APK31}. Let $\veps = 
(\varepsilon_0, \varepsilon_1, \ldots, \varepsilon_J) \in \RR^{J+1}$ 
be such that  $|\varepsilon_{\ell}| \le r$ for each $\ell = 0, \ldots, J$, 
and consider  
\begin{equation}   \label{AFK147}
Q^{\veps}(z) 
     \ = \  \exp(\varepsilon_0)^N \cdot C \cdot
     \prod_{\alpha_k \in E_v^{00}} [z,\alpha_k]_{\fX,\vs} 
         \cdot \prod_{\ell=1}^J \, \prod_{\alpha_k \in [t_{\ell}-h,t_{\ell}+h]} 
                  [z,\alpha_k + \varepsilon_\ell]_{\fX,\vs} \ .    
\end{equation}

\vskip .1 in
Let $\fX_1, \ldots, \fX_{m_1+m_2}$ be the orbits in $\fX$ under
complex conjugation, as in \S \ref{Chap5}.\ref{IndependenceSection}.  
If $x_i \in \fX_a$, then passing
from $Q(z)$ to $Q^{\veps}(z)$ produces the following change 
in the logarithmic leading coefficient at $x_i$:  
\index{logarithmic leading coefficients}   
\begin{eqnarray*}  
 \lambda_i(\veps) & := &\Lambda_{x_i}(Q^{\veps},\vs)
       - \Lambda_{x_i}(Q,\vs)       \notag \\
    &=& \varepsilon_0 + \sum_{\ell=1}^J \,  
          \sum_{\alpha_k \in [t_{\ell}-h,t_{\ell}+h]}  \frac{1}{N}
         \sum^m_{\substack{ j=1 \\ x_j \notin \fX_a }} s_k 
     \left( \log([x_i,\alpha_k+\varepsilon_\ell]_{x_j})
                - \log([x_i,\alpha_k]_{x_j}) \right) \ . 
\end{eqnarray*}
This is immediate if $\fX_a = \{x_i\}$;                                   
if $\fX_a = \{x_i,x_j\}$ consists of two points, 
it follows from the fact that $[x_i,t]_{x_j}$ is constant on $\cC_v(\RR)$, 
as shown in \S\ref{Chap5}.\ref{IndependenceSection}.    

Recalling that $Q(z)$ has $n_{\ell}$ roots in $\tE_{v,\ell}$, 
for each $a = 1, \ldots, m_1 + m_2$ fix a point $x_i \in \fX_a$ and put 
\begin{eqnarray*}
L_a(\veps) &=& \varepsilon_0 + \sum_{\ell=1}^J 
           H_{\vs,a}(t_{\ell}) \cdot \frac{n_{\ell}}{N} \varepsilon_{\ell}  \\ 
       &=&  \varepsilon_0 + \sum_{\ell=1}^J
           \sum^m_{\substack{ j = 1 \\ x_j \notin \fX_a }}
    s_j \left. \left( \frac{d}{dt} \log([x_i,t]_{x_j}) \right)
                       \right|_{t=t_{\ell}}
                 \cdot \frac{n_{\ell}}{N} \varepsilon_\ell      \ .               
\end{eqnarray*}
Then $L_a(\veps)$ is a linear map which approximates $\lambda_i(\veps)$,  and 
\begin{eqnarray*} 
\lefteqn{ | \lambda_i(\veps) - L_a(\veps) |   }\\       
   &\le& \frac{1}{N} \sum_{\ell=1}^J \, \sum_{\alpha_k \in [t_{\ell}-h,t_{\ell}+h]}  
         \sum^m_{\substack{ j=1 \\ x_j \notin \fX_a }} 
   s_j \left\{ 
          \left| \log([x_i,\alpha_k + \varepsilon_\ell]_{x_j}) -    
            \log([x_i,t_{\ell} + \varepsilon_\ell]_{x_j}) \right| 
                    \phantom{\frac{d}{dt}}  \right.  \\       
   & & \qquad \qquad \qquad \qquad \qquad \quad 
              +  \left| \log([x_i,\alpha_k]_{x_j}) -    
                    \log([x_i,t_{\ell}]_{x_j}) \right|    \\          
   & & \qquad \left. 
     + | (\log([x_i,t_{\ell}+\varepsilon_\ell]_{x_j})
           - \log([x_i,t_{\ell}]_{x_j}))
        - \left. \left( \frac{d}{dt} \log([x_i,t]_{x_j}) \right)
                       \right|_{t=t_{\ell}} \varepsilon_\ell | \right\} \ .
\end{eqnarray*}    
By (\ref{AFK140}) and (\ref{AFK143}), the magnitude of each term in absolute values 
on the right side is at most $r/(24B_1B_4\sqrt{J+1})$.  Hence
\begin{equation}    \label{AFK148}
 | \lambda_i(\veps) - L_a(\veps) |       
    \ \le \  \frac{n_1 + \ldots + n_J}{N} \cdot \frac{3r}{24B_1B_4\sqrt{J+1}} \ .
\end{equation} 

\vskip .1 in

Next, consider the ``normalized divisor change'' map 
$\widehat{\varphi}(\veps)$ (see (\ref{bF14})) 
induced by passing from $Q(z)$ to $Q^{\veps}(z)$.  
If $g = 0$, every divisor of degree $0$
is principal, and there are no nonzero holomorphic differentials;\index{differential!holomorphic}  in  that case 
ignore all constructions related to the $\omega_j = h_j(z) dz$ in the rest of the proof.
 
If $g \ge 1$, for each $j =1, \ldots, g$, the $j^{th}$ coordinate of $\widehat{\varphi}(\veps)$ is 
\begin{equation*}
\widehat{\varphi}_j(\veps) 
\ = \ \sum_{\ell=1}^J \, \sum_{\alpha_k \in [t_{\ell}-h,t_{\ell}+h]} 
          \frac{1}{N} \int_{\alpha_k}^{\alpha_k+\varepsilon_\ell} h_j(z) \, dz \ .
\end{equation*}
By the Mean Value Theorem for integrals,\index{Mean Value theorem for integrals} for each $\alpha_k$ there is an 
$\alpha_k^* \in [\alpha_k,\alpha_k+\varepsilon_{\ell}]$ such that
 $\int_{\alpha_k}^{\alpha_k+\varepsilon_{\ell}} h_j(z) dz \ = \  
            h_j(\alpha_k^*) \cdot \varepsilon_{\ell}$.  
Consequently, if 
\begin{equation*}
f_j(\veps) \ := \  
\sum_{\ell=1}^J h_j(t_{\ell}) \cdot \frac{n_{\ell}}{N} \varepsilon_{\ell} 
\end{equation*}
is the linear approximation to $\widehat{\varphi}_j(\veps)$ at the origin, then 
\begin{equation*}  
|\widehat{\varphi}_j(\veps) - f_j(\veps)| \ \le \ 
       \sum_{\ell=1}^J \, \sum_{\alpha_k \in [t_{\ell}-h,t_{\ell}+h]}          
       \frac{1}{N} \, |  h_j(\alpha_k^*) - h_j(t_{\ell})| \, |\varepsilon_{\ell}| \ .
\end{equation*}
From (\ref{AFK141}) and the fact that $|\varepsilon_{\ell}| \le r$ for each $\ell$, 
it follows that 
\begin{equation}   \label{AFK149} 
|\widehat{\varphi}_j(\veps) - f_j(\veps)| \ \le \ 
     \frac{n_1 + \ldots + n_J}{N} \cdot \frac{r}{24B_1B_4\sqrt{J+1}} \ .
\end{equation}

Let $B^{J+1}(0,r) = \{ (\varepsilon_0,\ldots, \varepsilon_J) \in \RR^{J+1} 
: |\veps| \le r \}$.  
Assuming $x_1, \ldots, x_{m_1+m_2}$ form a set of representatives 
for the orbits $\fX_1, \ldots, \fX_{m_1+m_2}$,   
we define the ``total change map'' $F^Q : B^{J+1}(0,r) \rightarrow \RR^{J+1}$ 
for $Q(z)$ by 
\begin{equation}       \label{AFK150}   
F^Q(\veps) \ = \ \left( \begin{array}{c} 
       \lambda_1(\veps)                    \\ 
                       \vdots                          \\ 
       \lambda_{m_1+m_2}(\veps)            \\ 
              \widehat{\varphi}_1(\veps)       \\                            \\
                       \vdots                          \\
              \widehat{\varphi}_g(\veps)  \end{array} \right) \ ,
\end{equation}
and note that its linear approximation
\begin{equation*}       
F^Q_0(\veps)
   \ := \ \left( \begin{array}{c} L_1(\veps) \\ \vdots \\ L_{m_1+m_2}(\veps) \\
                     f_1(\veps) \\ \vdots \\ f_g(\veps) \end{array} \right) 
  \ = \ \left( \begin{array}{c} 
     \varepsilon_0 + \sum_{{\ell} =1}^J 
         H_{\vs,1}(t_{\ell}) \cdot \frac{n_{\ell} }{N} \varepsilon_{\ell}         \\ 
                    \vdots                                           \\ 
     \varepsilon_0 + \sum_{{\ell} =1}^J 
         H_{\vs,m_1+m_2}(t_{\ell}) \cdot \frac{n_{\ell} }{N} \varepsilon_{\ell}   \\  
      \sum_{{\ell} =1}^J h_1(t_{\ell}) \cdot \frac{n_{\ell} }{N} \varepsilon_{\ell}     \\
                     \vdots                                          \\
      \sum_{{\ell} =1}^J h_g(t_{\ell}) \cdot \frac{n_{\ell} }{N} \varepsilon_{\ell}   
           \end{array} \right)
\end{equation*}     
can be decomposed as 
\begin{equation} \label{AFK151}
F^Q_0(\veps) \ = \ \cG_{\vs}(\vec{t}) \cdot \cN \cdot \veps 
\end{equation} 
where $\cN$ is the $(J+1) \times (J+1)$ diagonal matrix 
$\cN = \ \text{diag}(1,n_1/N,\ldots,n_J/N)$ and 
$\cG_{\vs}(\vec{t})$ is as in Corollary \ref{bC1}.

By (\ref{AFK148}), (\ref{AFK149}), (\ref{AFK150}) and (\ref{AFK151}), if 
\begin{equation}         \label{AFK152} 
 d(\veps) \ = \ F^Q(\veps) - F^Q_0(\veps) 
\end{equation} 
for $\veps \in B^{J+1}(0,r)$, 
then each coordinate function $d_i(\veps)$ of $d(\veps)$ satisfies 
\begin{equation}         \label{AFK153}
|d_i(\veps)| 
\ \le \ \frac{r}{8B_1B_4\sqrt{J+1}} \cdot \sum_{{\ell} =1}^J \frac{n_{\ell} }{N} \ .
\end{equation}  

Next consider the renormalized map $\cF^Q : B^{J+1}(0,r) \rightarrow \RR^{J+1}$
given by 
\begin{equation}       \label{AFK154}
\cF^Q(\veps) \ = \ \cN^{-1} \cG_{\vs}(\vec{t})^{-1} \cdot F^Q(\veps) \ , 
\end{equation}    
and the corresponding difference function  
\begin{equation}     \label{AFK155}
D(\veps) \ = \ \cF^Q(\veps) - \veps 
               \ = \ \cN^{-1} \cG_{\vs}(\vec{t})^{-1} \cdot d(\veps) \ .
\end{equation}
By construction, the operator norm 
$\|\cG_{\vs}(\vec{t})^{-1}\|$ is bounded by $B_1$, and clearly
\begin{equation*}
\|\cN^{-1}\| \ = \ \frac{1}{\min_{1 \le \ell \le J} n_{\ell}/N} \ .
\end{equation*}
It follows from (\ref{AFK152}) and (\ref{AFK153}) that 
for each $\veps \in B^{J+1}(0,r)$
\begin{eqnarray} 
|D(\veps)| 
      &\le& \frac{1}{\min_{1 \le \ell \le J} n_{\ell}/N} \cdot B_1 \cdot \sqrt{J+1} \cdot  
        \frac{r}{8B_1B_4\sqrt{J+1}} \cdot \sum_{\ell=1}^J n_{\ell}/N
                                                         \notag \\
          &=& \frac{r}{8} \cdot \frac{1}{B_4} \cdot
        \frac{ \sum_{\ell=1}^J n_{\ell}/N } {\min_{1 \le \ell \le J} n_{\ell}/N } \ .
                      \label{AFK156} 
\end{eqnarray}

\medskip
We now cite a topological fact  
related to the Brouwer Fixed Point Theorem.\footnote{The author 
thanks Ted Shifrin for helpful discussions concerning this argument.} 
\index{Shifrin, Ted} 
\index{Brouwer Fixed Point theorem|ii}  

\begin{proposition} \label{APK32}
Fix $0 < \eta < 1$ and suppose
$g : B^{J+1}(0,r) \rightarrow \RR^{J+1}$ is a continuous map such that for each 
$\vec{x}$ in $\partial B^{J+1}(0,r)$, 
\begin{equation*}
|g(\vec{x}) - \vec{x}| \ < \ \eta \cdot r \ .
\end{equation*}
Then $B^{J+1}(0,(1-\eta)r) \subset g(B^{J+1}(0,r))$.                    
\end{proposition}
       
\begin{proof}  
If $J = 0$ the result follows easily from the Intermediate Value theorem,\index{Intermediate Value theorem} 
so we can assume $J \ge 1$.  Suppose there were some 
$\vec{x}_0 \in B(0,(1-\eta) r)$ which did not belong to $g(B^{J+1}(0,r))$.
Then we could define a continuous map  
$g_0 : B^{J+1}(0,r) \rightarrow \partial B^{J+1}(0,r)$, such that $g_0(\vec{x})$
is the point where the ray from $\vec{x}_0$ through $g(\vec{x})$
meets $\partial B^{J+1}(0,r)$.  For each $\vec{x} \in \partial B^{J+1}(0,r)$,
$|g(\vec{x}) - \vec{x}| < \eta r$ while by hypothesis 
$|\vec{x}_0 - \vec{x}| > \eta r$.  Hence there would be a homotopy from 
$g_0|_{\partial B^{J+1}(0,r)}$  to the identity, such that for each $t \in [0,1]$ 
and each $\vec{x} \in \partial B^{J+1}(0,r)$, $g_t(\vec{x})$ is the point where
the ray from $\vec{x}_0$ through
$tg(\vec{x}) + (1-t)\vec{x}$ meets $\partial B^{J+1}(0,r)$.
Thus $g_0 : B^{J+1}(0,r) \rightarrow \partial B^{J+1}(0,r)$ would be a continuous map 
whose restriction to $\partial B^{J+1}(0,r)$ was homotopic to the identity.  
Now let $i : \partial B^{J+1}(0,r) \rightarrow B^{J+1}(0,r)$ be the inclusion, and
consider the induced map $(g_0 \circ i)_* = (g_0)_* \circ i_*$ on homology:
\begin{equation*}
H_J(\partial B^{J+1}(0,r)) \ \stackrel{i_*}{\longrightarrow} \ H_J(B^{J+1}(0,r))
          \ \stackrel{(g_0)_*}{\longrightarrow} \ H_J(\partial B^{J+1}(0,r)) \ .
\end{equation*}
Here $H_J(B^{J+1}(0,r)) = 0$ so $(g_0)_* \circ i_*$ is the $0$ map.  
On the other hand $H_J(\partial B^{J+1}(0,r)) \cong \ZZ$ and 
$g_0 \circ i$ is homotopic to the identity, so $(g_0 \circ i)_*$ 
is the identity map.  This is a contradiction, so $\vec{x}_0 \in g(B^{J+1}(0,r))$.
\end{proof}  

\medskip
Let $\vec{\sigma}$ be the vector of weights of the components of $\tE_v$ under
$\mu_{\fX,\vs}$.  For any special  pseudopolynomial
\index{pseudopolynomial!special}  
$Q_{\vec{n}}(z)$ for $\tE_v$ with $N = \sum n_i$ sufficiently large, 
and with $\vec{n}/N$ close enough to
$\vec{\sigma}$ that $|n_{\ell} /N - \mu_{\fX,\vs}(\tE_{v,\ell})|
  < \frac{1}{2} \mu_{\fX,\vs}(\tE_{v,\ell})$ for ${\ell}  = 1, \ldots, J$
we will have
\begin{equation} \label{AFK157}
\frac{ \sum_{{\ell} =1}^J n_{\ell} /N}{ \min_{1 \le \ell \le J}  n_{\ell} /N } 
\ < \ \frac{3/2 \cdot \mu_{\fX,\vs}(\bigcup_{k = 1}^J \tE_{v,k})}
           {1/2 \cdot \min_{1 \le \ell \le J} \mu_{\fX,\vs}(\tE_{v,\ell})}
\ = \ 3B_4.
\end{equation}
If (\ref{AFK157}) holds, then in (\ref{AFK156}) we will have
$|D(\veps)| < r/2$.  Applying Proposition \ref{APK32} to $\cF^Q(\veps)$ 
with $\eta = 1/2$, it follows that
\begin{equation} \label{AFK158}
B^{J+1}(0,r/2) \ \subset \ \cF^Q(B^{J+1}(0,r)) \ .
\end{equation}
\index{initial approximating functions $f_v(z)$!construction when $K_v \cong \RR$!Step 5: study of the total change map|)} 

\medskip
\noindent{\bf Step 6.}  The choice of $\delta_v$.
\index{initial approximating functions $f_v(z)$!construction when $K_v \cong \RR$!Step 6: the choice $\delta_v$|(} 
\medskip

Recalling that $\mu_{\fX,\vs} = \sum_{i=1}^m s_i \mu_i$ where
$\mu_i$ is the equilibrium distribution of $\tE_v$ with respect to $x_i$,
put
\begin{equation} \label{AFK159}
B_5\ = \ \min_{1 \le i \le m} \min_{1 \le \ell  \le d} 
            \mu_i(\tE_{v,\ell}) \ > \ 0 \ .
\end{equation}
Then for any $\vs \in \cP^m_v$, and any $\ell = 1, \ldots, d$,
\begin{equation} \label{AFK160}
\mu_{\fX,\vs}(\tE_{v,\ell}) \ \ge \ B_5\ .
\end{equation}

Henceforth we will assume $\vn$ and $N$ are such that 
\begin{equation} \label{AFK161}
\min_{1 \le \ell \le d} n_{\ell}/N  \ \ge \ B_5/2 \ .
\end{equation}
In the light of (\ref{AFK158}), (\ref{AFK154}), (\ref{AFK161}) and the construction
of $\cN$ and $\cG_{\vs}(\vec{t})$, for any such $\vn$, if $Q(z) = Q_{\vn}$ 
as in Theorem \ref{PotPart} then  
\begin{equation} \label{AFK162}
B(0,\frac{r B_5}{4B_1}) \ = \ \frac{B_5}{2B_1} \cdot B(0,r/2)
     \ \subset \ F^Q(B^{J+1}(0,r)) \ .
\end{equation}
This radius $\frac{r B_5}{4B_1}$ is uniform for all $\vs \in \cP^m_v$.

\vskip .1 in
Let the number $\delta_v > 0$ in Theorem \ref{RThm2} be such that
\begin{equation}  \label{AFK163}
\left( \prod_{i=1}^{m_1+m_2} [-2\delta_v,2\delta_v] \right) \times B^g(0,\delta_v)
 \ \subset \ B(0,\frac{r B_5}{4B_1}) \ ,
\end{equation}
where $B^g(0,\delta_v) = \{\vec{y} \in \RR^g : |\vec{y}| \le \delta_v \}$.  
Then if (\ref{AFK157}) and (\ref{AFK161}) hold, 
\begin{equation} \label{FBoth}
[-2\delta_v,2\delta_v]^{m_1+m_2} \times B^g(0,\delta_v) 
           \ \subset \ F^{Q}(B^{J+1}(0,r)) \ .
\end{equation}
\index{initial approximating functions $f_v(z)$!construction when $K_v \cong \RR$!Step 6: the choice $\delta_v$|)} 

\medskip
\noindent{\bf Step 7.}  Achieving principality and 
varying  the logarithmic leading coefficients.
\index{initial approximating functions $f_v(z)$!construction when $K_v \cong \RR$!Step 7: achieving principality|(} 
\index{logarithmic leading coefficients!independent variability of archimedean|ii} 

\vskip .1 in
As in \S\ref{Chap5}.\ref{OutlineSection}, 
let $\Jac(\cC_v)$ be the Jacobian of $\cC_v$,
\index{Jacobian variety}
and let $\omega_1, \ldots, \omega_g$ be a basis for the space of 
real holomorphic differentials\index{differential!real} of $\cC_v$.  
Then $\Jac(\cC_v)(\CC) \cong \CC^g/\cL$ where $\cL$ 
is the period lattice\index{period lattice} corresponding to $\omega_1, \ldots, \omega_g$.
Let $\varphi : \cC_v(\CC) \rightarrow \Jac(\cC_v)(\CC)$ be the canonical map
associated to a base point $p_0 \in \cC_v(\RR)$,
\begin{equation}    \label{AFK164} 
\varphi(p) = (\int_{p_0}^p \omega_1, \ldots, \int_{p_0}^p \omega_g)
\ (\text{mod $\cL$}) \ .
\end{equation}
and let
$\Phi : \cC_v(\CC)^g \rightarrow \Jac(\cC_v)(\CC)$ be the corresponding summatory map.

Since $\Jac(\cC_v)$ is nonsingular and defined over $\RR$, and its origin is
rational over $\RR$, the implicit function theorem shows that $\Jac(\cC_v)(\RR)$
is a $g$-dimensional real manifold.  Since $\Jac(\cC_v)$ is defined over $\RR$,
$\Jac(\cC_v)(\RR)$ is a compact Lie subgroup of $\Jac(\cC_v)(\CC)$.  Thus its identity
component must be isomorphic to a real torus $\RR^g/\cL_0$,
where $\cL_0 = \cL \bigcap \RR^g$ is a $g$-dimensional sublattice
of $\cL$.

Let $\Jac(\cC_v)(\RR)_0$ be the identity component of $\Jac(\cC_v)(\RR)$.

We claim that for each $p \in \cC_v(\RR)$, 
the point $2 \varphi(p)$ belongs to $\Jac(\cC_v)(\RR)_0$.  
To see this, suppose that in (\ref{AFK164}), 
$z = \varphi(p) \in \CC^g/\cL$ is obtained    
by integrating over some path $\gamma$ from $p_0$ to $p$.  
Applying complex conjugation to (\ref{AFK164}), we see that since
$p_0$, $p$ and the $\omega_i$ are fixed but 
$\gamma$ is taken to its conjugate path, $z \ (\text{mod $\cL$})$ is 
changed to $\overline{z} \ (\text{mod $\cL$)}$.  Hence, lifting $z$ to $\CC^g$, 
it must be that
\begin{equation} \label{AFK165}
\overline{z} \ = \ z - \ell
\end{equation}
for some $\ell \in \cL$.  We can find a basis
$\ell_1, \ldots, \ell_{2g}$ for $\cL$ such that $\ell_1, \ldots, \ell_g$ is
a basis for $\cL_0$.  Writing
\begin{eqnarray}
z &=& \sum_{i=1}^g a_i \ell_i + \sum_{i = 1}^g b_i \ell_{g+i}
   \quad \text{where $a_i, b_i \in \RR$} \ ,     \label{AFK166} \\
\ell &=& \sum_{i=1}^g c_i \ell_i + \sum_{i = 1}^g d_i \ell_{g+i}
   \quad \text{where $c_i, d_i \in \ZZ$} \ ,     \label{AFKW3}
\end{eqnarray}
it follows that
\begin{equation}  \label{AFK168}
\overline{z} = \sum_{i=1}^g a_i \ell_i + \sum_{i = 1}^g b_i
\overline{\ell}_{g+i} \ .
\end{equation}
But $\overline{\ell}_{g+i} \in \cL$ for each $i$, and since
$\overline{\ell}_{g+i} + \ell_{g+i} \in \RR^g$, for each $i$
\begin{equation}  \label{AFK169}
\overline{\ell}_{g+i} = -\ell_{g+i} + \sum_{k=1}^g M_{ik} \ell_k
\end{equation}
for certain integers $M_{ik}$.  Inserting (\ref{AFK166})--(\ref{AFK169}) 
in (\ref{AFK165}), we see that
$2b_i = d_i \in \ZZ$ for each $i = 1, \ldots, g$.  Thus, $2z \in \cL+\RR^g$,
which shows that $2 \varphi(p) \in \Jac(\cC_v)(\RR)_0$.

\vskip .1 in
Fix  $\vs \in \cP^m_v \cap \QQ^m$, 
and consider a special $(\fX,\vs)$-pseudopolynomial
\index{pseudopolynomial!special} 
\index{pseudopolynomial!$(\fX,\vs)$} 
\begin{equation*}
Q(z) \ = \ Q_{\vec{n}}(z) \ = \ 
C_{\vn} \cdot \prod_{k=1}^N [z,\alpha_k]_{\fX,\vs}
\end{equation*} 
with roots in $\tE_v$. Clearly roots belonging to the same segment in $\tE_v$
map to the same connected component of $\Jac(\cC_v)(\RR)$ under $\varphi$.
Hence if we arrange that each coordinate $n_{\ell}$ of $\vn$ is even, 
then $\sum_{i=1}^N \varphi(\alpha_k)$ will belong to $\Jac(\cC_v)(\RR)_0$.  
Likewise if $N = \sum n_{\ell}$ is such that $N \vs \in 2\ZZ^m$, then
$\varphi(\sum_{i=1}^m Ns_i(x_i)) \in \Jac(\cC_v)(\RR)_0$.  
If $S$ is the least common
denominator of the $s_i$, both of these conditions can
be assured by requiring  that
\begin{equation}                   \label{AFK170}
n_{\ell} \ \equiv \ 0 \ (\text{mod $2S$})  \quad \text{for all ${\ell}$} \ .
\end{equation}
Suppose $\vn/N$ is close enough to $\vec{\sigma}$ that (\ref{AFK157}) 
and (\ref{AFK161}) hold. Then 
\begin{equation*}
B^g(0,\delta_v) \ \subset \ \widehat{\varphi}(B^{J+1}(0,r)) \ .
\end{equation*} 
If in addition $N$ is sufficiently large, 
then $N \cdot B^g(0,\delta_v)$ will
contain a fundamental domain for $\RR^g/\cL_0$.

Since $\varphi(\div(Q^{\veps}) - \varphi(\div(Q)) = N \cdot \widehat{\varphi}(\veps)$,
it follows that if (\ref{AFK157}), (\ref{AFK161}), and (\ref{AFK170}) 
are satisfied and $N$ is sufficiently large, 
then $\veps \in B^{J+1}(0,r)$ can be chosen so that $Q^{\veps}(z)$ is principal.  
Furthermore, by (\ref{FBoth}) 
this can be done in such a way that for any $K_v$-symmetric 
\index{$K_v$-symmetric!probability vector}
$\vbeta^{\prime} = (\beta_1^{\prime}, \ldots, \beta_m^{\prime}) 
\in [-2\delta_v,2\delta_v]^m$, 
\begin{equation}   \label{AFK171}
\Lambda_{x_i}(Q^{\veps},\vs) \ = \ \Lambda_{x_i}(Q,\vs) + \beta_i^{\prime} 
            \quad \text{for $i = 1, \ldots, m$}\ .
\end{equation}
\index{initial approximating functions $f_v(z)$!construction when $K_v \cong \RR$!Step 7: achieving principality|)} 

\medskip
\noindent{\bf Step 8.}  The choice of the number $N_v$, and conclusion of the proof.
\index{initial approximating functions $f_v(z)$!construction when $K_v \cong \RR$!Step 8: the choice of $N_v$|(} 

\vskip .1 in
Let $\tcR_v$ be as in Proposition \ref{APK31}, 
let $r$ be the number constructed in Step 3,
and let $\delta_v$ be the number constructed in Step 7.   
Fix $\vs = (s_1, \ldots, s_m) \in \cP_v^m \cap \QQ^m$, 
and let $S$ be the least common denominator of the $s_i$.  

Take a sequence of $K_v$-symmetric vectors $\vn_k = (n_{k,1}, \ldots, n_{k,D}) \in \NN^D$ 
\index{$K_v$-symmetric!vector}
satisfying the following three properties:  

\vskip .03 in

$(1)$ for each $k$ and $\ell$, $n_{k,\ell} \equiv \ 0 \ (\text{mod $2S$})$ ; 

$(2)$ for each $k$, $\sum_{\ell} n_{k,\ell} = k \cdot 2S$ ;  

$(3)$ as $k \rightarrow \infty$, writing $N_k = \sum_{\ell} n_{k,\ell}$, 
we have $\vn_k/N_k \rightarrow \vec{\sigma}$.  

\vskip .03 in 
\noindent{Applying} Theorem \ref{PotPart} to $\{\vn_k\}_{k \in \NN}$, 
we obtain a sequence of $(\fX,\vs)$-pseudopolynomials $Q_{\vn_k}(z)$, 
whose roots belong to $\tE_v$, with $Q_{\vn_k}(z) = Q_{\vn_k}(\zbar)$ for all $z$ and $k$,
such that for all sufficiently large $k$,  
\vskip .03 in

$(4)$  $|\frac{1}{N_k} \log(Q_{\vn_k}(z)) - G_{\fX,\vs}(z,\tE_v)| < r$
for all $z \in \cC_v(\CC) \backslash (U_v \cup \fX)$\ ; 

$(5)$ $|\Lambda_{x_i}(Q_{\vn_k},\vs) - \Lambda_{x_i}(\tE_v,\vs)| \le \delta_v$ 
for each $x_i \in \fX$ ;  

$(6)$ $Q_{\vn_k}$ has $n_{k,\ell}$ roots in each component $\tE_{v,\ell}$ of $\tE_v$,
and on each component 

\quad \ contained in $\cC_v(\RR)$, $Q_{\vn_k}(z)$ varies $n_{k,\ell}$ times from $\tcR_v^{N_k}$ to 
$0$ to $\tcR_v^{N_k}$.
\vskip .03 in

\noindent{Likewise,} for all sufficiently large $k$, the set    
$N_k \cdot B^g(0,\delta_v)$ contains a fundamental domain for $\RR^g/\cL_0$, 
and (\ref{AFK157}) and (\ref{AFK161}) hold.  
Let $k_0$ be the least number such that these
conditions hold for all $k \ge k_0$, and put $N_v = N_{k_0} = k_0 \cdot 2S$.

Let $N > 0$ be a multiple of $N_v$, say $N = t \cdot N_v$ for some integer $t \ge 1$.  
Put $k = t \cdot k_0$, let $\vn = \vn_k$, and write $Q(z) = Q_{\vn_k}(z)$.  
Since $N$ is a multiple of $S$, 
the divisor $\div(Q) = \sum_{k=1}^N {\alpha_k} -  \sum_{i=1}^m N s_i (x_i)$ is integral.  
Since $\div(Q)$ is $K_v$-symmetric and $N s_i$ is even for each $i$, 
\index{$K_v$-symmetric!divisor}
the class of $\div(Q)$ in $\Jac(\cC_v)(\CC)$ belongs to $\Jac(\cC_v)(\RR)_0$.  
Let a $K_v$-symmetric vector $\vbeta = (\beta_1, \ldots, \beta_m) \in \RR^m$ with
\index{$K_v$-symmetric!vector}
each $|\beta_i| \le \delta_v$ be given.  Putting $\beta_{\ell}^{\prime} 
= \beta_{i} + \Lambda_{x_i}(\tE_v,\vs) - \Lambda_{x_i}(Q,\vs)$, 
we see that $|\beta_i^{\prime}| \le 2\delta_v$ for each $i$,
and that $\vbeta^{\prime} = (\beta_1^{\prime}, \ldots, \beta_m^{\prime})$
is $K_v$-symmetric. 
\index{$K_v$-symmetric!vector}
By (\ref{AFK171}) we can choose $\veps \in B^{J+1}(0,r)$ so that
$\div(Q^{\veps})$ is principal and
\begin{equation*}
\Lambda_{x_i}(Q^{\veps},\vs)
\ = \ \Lambda_{x_i}(Q,\vs) + \beta_{i}^{\prime}
\ = \ \Lambda_{x_i}(\tE_v,\vs) + \beta_i 
\end{equation*}
for each $i$. 

By Proposition \ref{APK31}, if $\tE_{v,\ell}$ is a component of $\tE_v$, 
then $Q^{\veps}$ has $n_{\ell}$ roots in the corresponding component 
$\tE_{v,\ell}(\veps)$ of $\tE_v(\veps)$. 
Furthermore, if $\tE_{v,\ell}$ is contained in $\cC_v(\RR)$, then 
$\tE_{v,\ell}(\veps)$ is contained in $\cC_v(\RR)$, and $Q^{\veps}(z)$ varies 
$n_{\ell}$ times from $\cR_v^N$ to $0$ to $\cR_v^N$ on $\tE_{v,\ell}(\veps)$. 

Since $\div(Q^{\veps})$ is principal and stable under complex conjugation, 
there is a rational function $f_v(z) \in \RR(\cC_v)$ with 
$|f_v(z)| = Q^{\veps}(z)$ for all $z$.  Thus 
$\Lambda_{x_i}(f_v,\vs) = \Lambda_{x_i}(\tE_v,\vs) + \beta_i$
for each $x_i \in \fX$.
Clearly $f_v(z)$ is real for all $z \in \cC_v(\RR)$, 
and $\div(f_v) = \div(Q^{\veps})$.  
Since $\tE_v(\veps)$ is contained in $E_v^0$, all the roots of $f_v$
belong to $E_v^0$. 
By Proposition \ref{APK31}, 
\begin{equation*}
\{ z \in \cC_v(\CC) : |f_v(z)| \le 1 \} \ \subset \ U_v \ .
\end{equation*} 

Let $E_{v,1}, \ldots, E_{v,k}$ be the 
components of the $K_v$-simple set $E_v$. 
\index{$K_v$-simple!set}  
For each $E_{v,i}$, put 
\begin{equation*}
N_i \ = \ \sum_{\tE_{v,\ell} \subset E_{v,i}} n_{\ell}
\ \ = \  \sum_{\tE_{v,\ell}(\veps) \subset E_{v,i}} n_{\ell} \ .
\end{equation*}
Then $f_v$ has $N_i$ roots in $E_{v,i}$, counted with multiplicities.
If $E_{v,i}$ is a component of $E_v$ contained in 
$\cC_v(\RR)$, then by construction $|f_v(z)| = Q^{\veps}(z)$ 
varies $N_i$ times from $\cR_v^N$ to $0$ to $\cR_v^N$ on $E_{v,i}$.  
Each of those oscillations accounts for at least one root of $f_v(z)$.  
Since $f_v$ has exactly $N_i$ roots in $E_{v,i}$, those roots
must be simple.  Since the only places $f_v(z)$ can change sign on $E_{v,i}$
are at the roots, each time $|f_v(z)|$ varies 
from $\cR_v^N$ to $0$ to $\cR_v^N$ on $E_{v,i}$, the function $f_v(z)$
varies from $\cR_v^N$ to $-\cR_v^N$, or from $-\cR_v^N$ to $\cR_v^N$.
Thus, the roots of $f_v$ in $E_{v,i}$ are simple and 
$f_v(z)$ oscillates $N_i$ times between $\pm \cR_v^N$ on $E_{v,i}$.  
\hfill $\blacksquare$
\index{initial approximating functions $f_v(z)$!construction when $K_v \cong \RR$!Step 8: the choice of $N_v$|)}

%% file: NewFSZChap6.tex
\chapter{Initial Approximating Functions:  Nonarchimedean Case}
\label{Chap6}   

In this section we will construct the nonarchimedean initial approximating functions
\index{initial approximating functions $f_v(z)$!nonarchimedean|ii}
needed for Theorem \ref{aT1-B}.  Because of the ultrametric inequality, 
the nonarchimedean theory has a different flavor from the archimedean theory:\, 
the constructions are more rigid, but at the same time more explicit.

\smallskip
As in \S\ref{Chap3}.\ref{AssumptionsSection}, let $K$ be a global field,
let $\cC/K$ be a smooth, connected, projective curve, 
and let $\fX = \{x_1, \ldots, x_m\} \subset  \cC(\tK)$ 
be a finite, $K$-symmetric set of points.  Put $L = K(\fX)$, 
\index{$K$-symmetric!set of points}
and let $\{g_{x_1}(z), \ldots, g_{x_m}(z)\}$
a system of uniformizing parameters,\index{uniformizing parameter!galois equivariant system} chosen in such a way that
$g_{x_i}(z) \in K(x_i)(\cC)$ for each $x_i$, with 
$g_{x_{\sigma(i)}}(z) = \sigma(g_{x_i})(z)$ for each 
$\sigma \in \Aut(\tK/K)$.

\smallskip
Let $v$ be a nonarchimedean place of $K$. 
In the statement of Theorem \ref{aT1-B}, 
we are given a $K_v$-symmetric 
set $E_v \subset \cC_v(\CC_v)$ which is bounded away from $\fX$ and has positive capacity, 
\index{capacity $> 0$} 
and a finite set of places $S$ of $K$.  
If $v \notin S$, then $E_v$ is $\fX$-trivial;  if $v \in S$, then 
\index{$\fX$-trivial}
$E_v$ is $K_v$-simple (Definition \ref{KvSimple}).  
\index{$K_v$-simple!set}  
This means that $E_v$ is compact, 
and there is a decomposition $E_v = \bigcup_{\ell=1}^D E_{v,\ell}$ with pairwise disjoint,  
nonempty compact sets $E_{v,\ell}$ such that 
\begin{enumerate}

\item There are 
finite separable extensions
$F_{w_1}, \ldots, F_{w_D}$ of $K_v$ contained in $\CC_v$,
and pairwise disjoint isometrically parametrizable balls
\index{isometrically parametrizable ball}
$B(a_1,r_1), \ldots, B(a_n,r_D)$, 
for which $E_{v,\ell} = \cC_v(F_{w_\ell}) \cap B(a_\ell,r_\ell)$ for $\ell = 1, \ldots, D$; 

\item  The set of balls $\{B(a_1,r_1), \ldots, B(a_D,r_D)\}$
is stable under $\Aut_c(\CC_v/K_v)$, and as $\sigma$ ranges over $\Aut_c(\CC_v/K_v)$,
each ball $B(a_\ell,r_\ell)$ has $[F_{w_\ell}:K_v]$ distinct conjugates.  
For each $\sigma$, if $\sigma(B(a_j,r_j)) = B(a_k,r_k)$, then 
$\sigma(F_{w_j}) = F_{w_k}$ and $\sigma(E_{v,j}) = E_{v,k}$. 
\end{enumerate} 
Both $\fX$-trivial sets and compact sets are algebraically capacitable
\index{$\fX$-trivial}
\index{algebraically capacitable}
(see \cite{RR1}, Theorems 4.3.13, 4.3.15), so $\Gbar(z,x_i;E_v) = G(z,x_i;E_v)$
and $\Vbar_{x_i}(E_v) = V_{x_i}(E_v)$ for each $x_i$;  throughout this section
we will write $G(z,x_i;E_v)$ and $V_{x_i}(E_v)$ for the Green's functions
\index{Green's function!nonarchimedean}
and Robin constants.  
\index{Robin constant!nonarchimedean}  

\vskip .05 in  

Let $\vs = (s_1, \ldots, s_m) \in \cP^m$ be a 
$K_v$-symmetric probability vector.  As in the archimedean case, 
\index{$K_v$-symmetric!probability vector}
the logarithmic leading coefficient of the Green's function of $E_v$ 
\index{Green's function!nonarchimedean}
\index{logarithmic leading coefficients}   
at $x_i$ is defined to be 
\begin{eqnarray} 
\Lambda_{x_i}(E_v,\vs) & = & \lim_{z \rightarrow x_i} 
   \Big(\sum_{j = 1}^m s_j G(z,x_j;E_v) + s_i\log_v(|g_{x_i}(z)|_v)\Big) \notag \\
   & = & s_i V_{x_i}(E_v) + \sum_{j \ne i} s_j G(x_i,x_j;E_v) \ . \label{FLL2}
\end{eqnarray}
\index{Robin constant!nonarchimedean}  
Likewise, if $\vs \in \cP^m(\QQ)$ and $f_v(z) \in K_v(\cC)$ is an
$(\fX,\vs)$-function\index{$(\fX,\vs)$-function!$K_v$-rational} of degree $N$,  
the logarithmic leading coefficient of $f_v(z)$ at $x_i$ is 
\index{logarithmic leading coefficients}   
\begin{equation} \label{FLL1}
\Lambda_{x_i}(f_v,\vs)  =  \lim_{z \rightarrow x_i} 
          \left(\frac{1}{N} \log_v(|f_v(z)|_v) + s_i \log_v(|g_{x_i}(z)|_v)\right) \ .
\end{equation}
Note that $|\CC_v^{\times}|_v = q_v^{\QQ}$ 
is dense in $\RR_{> 0}$ but is not equal to it. 
This means we cannot continuously vary the logarithmic 
leading coefficients $\Lambda_{x_i}(f_v,\vs)$ 
\index{logarithmic leading coefficients}   
as in the archimedean case;  this is one source of  
rigidity in the construction. 
  
\vskip .1 in
If $E_v$ is $\fX$-trivial, then it is an $\RL$-domain, 
\index{$\fX$-trivial}
\index{$\RL$-domain} 
and the initial local approximating function will be 
\index{initial approximating functions $f_v(z)$!for nonarchimedean $\RL$-domains}
an $(\fX,\vs)$-function $f_v(z) \in K_v(z)$ which defines $E_v$ as an $\RL$-domain: 
\index{$\RL$-domain}\index{$(\fX,\vs)$-function!$K_v$-rational} 
\begin{equation*}
E_v \ = \ \{ z \in \cC_v(\CC_v) : |f_v(z)|_v \le 1 \} \ .
\end{equation*}  
In this situation, $\sum_{i=1}^m s_i G(z,x_i;E_v) = \frac{1}{N} \log_v(|f(z)|_v)$
\index{Green's function!nonarchimedean}
for all $z \notin E_v$, and 
$\Lambda_{x_i}(f_v,\vs) = \Lambda_{x_i}(E_v,\vs)$
for each $x_i \in \fX$.

\vskip .05 in
If $E_v$ is $K_v$-simple, the initial local approximating function will be 
\index{initial approximating functions $f_v(z)$!for nonarchimedean $K_v$-simple sets}
\index{$K_v$-simple!set}  
an $(\fX,\vs)$-function\index{$(\fX,\vs)$-function!$K_v$-rational} $f_v(z) \in K_v(z)$ whose zeros
are distinct and belong to $E_v$, and are well-distributed relative 
to the $(\fX,\vs)$-equilibrium measure $\mu_{\fX,\vs}$ 
of a $K_v$-simple set $\tE_v \subseteq E_v$ constructed below
\index{$K_v$-simple!set}  
(see Appendix \ref{AppA} regarding $\mu_{\fX,\vs}$).    
We will show that for any $\beta_v > 0$ in $\QQ$, there exist functions having these 
properties and satisfying 
\begin{equation*}
\Lambda_{x_i}(f_v,\vs) \ = \  \Lambda_{x_i}(\tE_v,\vs) + \beta_v
\end{equation*}
for all $x_i \in \fX$.

In the $K_v$-simple case, the functions $f_v(z)$ will be contructed explicitly, 
\index{$K_v$-simple!set}  
by an argument generalizing the methods of (\cite{RR2}, \cite{RR3}).  
The idea is to first define an infinite subsequence of $E_v$ which is very
uniformly distributed relative to $\mu_{\fX,\vs}$, 
using a lemma of Balinski and Young \index{Balinski, Michael} \index{Young, H. Peyton}  
(\cite{BY}) originally proved for the purpose of apportioning seats 
in the US House of Representatives,\index{apportionment in House of Representatives}  
and to then take an initial segment of that sequence 
and modify it to become the zeros of $f_v(z)$.  The modification, which involves
moving a finite number of points so as to obtain a principal divisor,  
is carried out by using an action of a small neighborhood of the origin in 
$\Jac(\cC_v)(\CC_v)$ on a suitably generic polydisc in $\cC_v(\CC_v)^g$ 
which makes the polydisc into a principal homogeneous\index{principal homogeneous space}  
space for the neighborhood.  This action is studied in Appendix \ref{AppD}.      


\section{ The Approximation Theorems} \label{NonArchApproxSection}  

There are two cases to consider in constructing the initial local 
approximating functions: the $\RL$-domain case, and the compact case.  
\index{initial approximating functions $f_v(z)$!nonarchimedean}
\index{$\RL$-domain}  

\vskip .1 in
Recall that $\cC_v(\CC_v) \backslash E_v$ can be partitioned 
into equivalence classes called $\RL$-components, 
which play the same role as the connected components of the complement of $E_v$ 
in the archimedean case.  
The equivalence relation is defined by $z \equiv w$ iff $G(z,w;E_v) > 0$;  
\index{Green's function!nonarchimedean}
\index{$\RL$-component|ii} 
see (\cite{RR1}, Theorems 4.2.11 and 4.4.17, and Definition 4.4.18).  

For the $\RL$-domains constructed 
\index{$\RL$-domain} 
in the initial reductions of the Fekete-Szeg\"o theorem, 
\index{Fekete-Szeg\"o theorem with LRC}
each $\RL$-component of $\cC_v(\CC_v) \backslash E_v$ automatically 
\index{$\RL$-component} 
contains at least one point of $\fX$:  
If $E_v$ is a finite union of isometrically parametrizable balls, 
\index{isometrically parametrizable ball}
then by (\cite{RR1},Theorem 4.2.16 and Proposition 4.4.1(B)) 
the complement of $E_v$ consists of a single $\RL$-component.  
\index{$\RL$-component} 
If $E_v$ is $\fX$-trivial, then by definition $\cC_v$ has good reduction at $v$,
\index{$\fX$-trivial}
\index{good reduction} 
the points $x_i$ specialize to distinct points $\pmod v$, 
and $E_v = \cC_v(\CC_v) \backslash \bigcup_{i=1}^m B(x_i,1)^-$ 
relative to the spherical metric\index{spherical metric} on $\cC_v(\CC_v)$.   
In this case the $\RL$-components of $\cC_v(\CC_v) \backslash E_v$ are precisely
the balls $B(x_i,1)^-$ for the $x_i \in \fX$.  

The construction of the initial approximating functions
\index{initial approximating functions $f_v(z)$} 
when $E_v$ is an $\RL$-domain was already treated in (\cite{RR1}, \S4):
\index{$\RL$-domain}
\index{initial approximating functions $f_v(z)$!for nonarchimedean $\RL$-domains|ii}
\index{initial approximation theorem!for nonarchimedean $\RL$-domains|ii}  

\begin{theorem}  \label{RLThm}  Suppose $K_v$ is nonarchimedean, 
and let $E_v \subset \cC_v(\CC_v) \backslash \fX$ be an RL-domain 
such that each $\RL$-component of \ $\cC_v(\CC_v) \backslash E_v$ contains at
\index{$\RL$-component} 
least one point of $\fX$.   Let $\vs \in \cP^m(\QQ)$ be a 
$K_v$-symmetric rational probability vector, whose entries are all positive.
\index{$K_v$-symmetric!probability vector}

Then there is a integer $N_v$ with the following property.
For each positive integer $N$ divisible by $N_v$, 
there is an $(\fX,\vs)$-function\index{$(\fX,\vs)$-function!$K_v$-rational} 
$f_v(z) \in K_v(\cC)$ of degree $N$ which defines $E_v$ as an $\RL$-domain: 
\begin{equation} \label{FRLF1} 
E_v \ = \ \{z \in \cC_v(\CC_v) : |f_v(z)|_v \le 1 \} \ .
\end{equation} 
For any such function $f_v(z)$ we have 
\begin{equation*}
 \sum_{i=1}^m s_i G(z,x_i;E_v)
\ = \ \left\{ \begin{array}{ll} 
         \frac{1}{N} \log_v(|f_v(z)|_v) & \text{ if $z \notin E_v$ \ , } \\
         0 & \text{ if $z \in E_v$ \ . }
     \end{array} \right.
\end{equation*} 
\index{Green's function!nonarchimedean}
In particular, $\Lambda_{x_i}(f_v,\vs) = \Lambda_{x_i}(E_v,\vs)$ 
for each $x_i \in \fX$.

Furthermore, if $\Char(K_v) = p > 0$, we can require that for each $x_i$, the leading coefficient 
\index{coefficients $A_{v,ij}$!leading}   
$c_{v,i} = \lim_{z \rightarrow x_i} f_v(z) \cdot g_{x_i}(z)^{Ns_i}$ 
belongs to $K_v(x_i)^{\sep}$.   
\end{theorem}      

\begin{proof}  This is essentially (\cite{RR1}, Theorem 4.5.4, p.316).  
That theorem provides an integer $N_0$ 
and an $(\fX,\vs)$-function\index{$(\fX,\vs)$-function!$K_v$-rational}
$f_0(z) \in K_v(\cC)$ of degree $N_0$ for which (\ref{FRLF1}) holds.  Put $N_v = N_0$.
Given a multiple $N = k N_v$, we can take $f_v(z) = f_0(z)^k$.  

If $\Char(K_v) = p > 0$, put $p^B = \max_i([K_v(x_i):K_v]^{\insep})$, and take $N_v = p^B N_0$ 
instead.  Then if $N = k N_v$, we have $f_v(z) = (f_0(z)^{p^C})^k$.  Since $f_0$ is $K_v$-rational,
for each $i$ its leading coefficient at $x_i$ belongs to $K_v(x_i)$ (see Corollary \ref{CoeffRationalityCor}), 
and the leading coefficient\index{coefficients $A_{v,ij}$!leading}
of $f_0(z)^{p^C}$ at $x_i$ belongs to $K_v(x_i)^{\sep}$.   
\end{proof} 

Before stating the approximation theorem for the compact case,
we need a definition.  

\begin{definition} \label{CompatibleDef} 
Suppose $K_v$ is nonarchimedean.  
We will say that $K_v$-simple sets $E_v$ and $\tE_v$ are {\em compatible},
or that $\tE_v$ is {\em compatible with} $E_v$, if $E_v$ and $\tE_v$ 
\index{$K_v$-simple!set!compatible with another set|ii}  
are nonempty and have $K_v$-simple decompositions 
\begin{equation} \label{E_vH_vFormula}
E_v = \bigcup_{\ell=1}^n B(a_\ell,r_\ell) \cap \cC_v(F_{w_\ell}) \ , \quad 
\tE_v = \bigcup_{k=1}^{\tn} B(\ta_k,\tr_k) \cap \cC_v(\tF_{w_k}) \ , 
\end{equation} 
such that 
$\bigcup_{k=1}^{\tn} B(\ta_k,\tr_k) \subseteq \bigcup_{\ell=1}^n B(a_\ell,r_\ell)$, 
and whenever $B(\ta_k,\tr_k) \subseteq B(a_{\ell},r_{\ell})$ 
we have  $\tF_{w_k} = F_{w_\ell}$. 
We will call a  pair of $K_v$-simple decompositions (\ref{E_vH_vFormula})
satisfying the conditions above 
{\em compatible decompositions}.\index{$K_v$-simple!decomposition!compatible with another decomposition|ii} 
\end{definition}

\smallskip
In the compact case, the theorem we need is the following.  
It may appear somewhat strange at first reading.  The set $\tE_v \subset E_v$
plays an auxiliary role in the construction: 
we will initially replace $E_v$ by $\tE_v$.  
This reserves an `unused' part $E_v \backslash \tE_v$ of $E_v$,
and allows us to move some zeros from $\tE_v$
into $E_v$ in constructing the approximating function $f_v(z)$.  
\index{initial approximating functions $f_v(z)$!nonarchimedean}
The final approximation set $H_v$ is then defined in terms of $f_v(z)$;  
it is contained in $E_v$ but not in general in $\tE_v$. 
\index{initial approximating functions $f_v(z)$!for nonarchimedean $K_v$-simple sets|ii}
\index{initial approximation theorem!for nonarchimedean $K_v$-simple sets|ii}      
   
\begin{theorem}   \label{CompactThm} 
Suppose $K_v$ is nonarchimedean.  
Let $E_v$ be a compact $K_v$-simple set which is disjoint from $\fX$ 
\index{$K_v$-simple!set}  
and has positive capacity. Fix a $K_v$-simple decomposition 
\index{$K_v$-simple!decomposition}  
\index{capacity $> 0$} 
\begin{equation} \label{E_vSimple}
E_v \ = \ \bigcup_{\ell=1}^D B(a_\ell,r_\ell) \cap \cC_v(F_{w_\ell}) \ ,
\end{equation}     
and fix $\varepsilon_v > 0$. Then there is a compact, $K_v$-simple set $\tE_v \subseteq E_v$ 
\index{$K_v$-simple!set}  
compatible with $E_v$ such that 

$(A)$ For each $x_i \in \fX$ 
\begin{equation} \label{FRRApprox4B}
|V_{x_i}(\tE_v) - V_{x_i}(E_v)| \ < \ \varepsilon_v \ ,
\end{equation} 
\index{Robin constant!nonarchimedean}  
and for all $x_i, x_j \in \fX$ with $x_i \ne x_j$, 
\begin{equation} \label{FRRApprox5B} 
|G(x_i,x_j;\tE_v) - G(x_i,x_j;E_v)| \ < \ \varepsilon_v \ .
\end{equation}
   
$(B)$ For each $0 < \beta_v \in \QQ$  
and each $K_v$-symmetric probability vector
\index{$K_v$-symmetric!probability vector}
$\vs = {}^t(s_1, \ldots, s_m)$ with rational entries, 
there is an integer $N_v \ge 1$  such that 
for each positive integer  $N$ divisible by $N_v$, 
there is an $(\fX,\vs)$-function\index{$(\fX,\vs)$-function!$K_v$-rational}
 $f_v \in K_v(\cC_v)$ of degree $N$ such that
   
\quad $(1)$  For each $x_i \in \fX$, 
\begin{equation} \label{FRRApprox7B} 
\Lambda_{x_i}(f_v,\vs) \ = \ \Lambda_{x_i}(\tE_v,\vs) + \beta_v \ .
\end{equation}   
   
\quad $(2)$ The zeros $\theta_1, \ldots, \theta_N$ of $f_v$ are distinct and belong to $E_v$;

\quad $(3)$ $f_v^{-1}(D(0,1)) \subseteq \bigcup_{\ell = 1}^D B(a_{\ell},r_{\ell})$;

\quad $(4)$  There is a decomposition $f_v^{-1}(D(0,1)) = \bigcup_{h=1}^N B(\theta_h,\rho_h)$, 
where the balls $B(\theta_h,\rho_h)$ are pairwise disjoint and isometrically parametrizable. 
\index{isometrically parametrizable ball} 
For each $h = 1, \ldots, N$, if $\ell = \ell(h)$ is such that 
$B(\theta_h,\rho_h) \subseteq B(a_{\ell},r_{\ell}),$ put $F_{u_h} = F_{w_\ell}$.  
Then $\rho_h \in |F_{u_h}^{\times}|_v$
and $f_v$ induces an $F_{u_h}$-rational scaled isometry from $B(\theta_h,\rho_h)$ to $D(0,1)$,
\index{scaled isometry} 
with
\begin{equation*}
f_v\big( B(\theta_h,\rho_h) \cap \cC_v(F_{u_h})\big) \ = \ \cO_{F_{u_h}} \ , 
\end{equation*}      
such that  $|f_v(z_1)-f_v(z_2)|_v = (1/\rho_h) \|z_1,z_2\|_v$ for all $z_1, z_2 \in B(\theta_h, \rho_h)$. 

\quad $(5)$ The set $H_v := E_v \cap f_v^{-1}(D(0,1))$ is $K_v$-simple and compatible with $E_v$.
\index{$K_v$-simple!set} 
Indeed, 
\begin{equation} \label{E_vPBF}
H_v  \ = \ 
  \bigcup_{h=1}^N \big(B(\theta_h,\rho_h) \cap \cC_v(F_{u_h})\big) 
\end{equation}
is a $K_v$-simple decomposition of $H_v$ compatible with the $K_v$-simple
\index{$K_v$-simple!decomposition}
\index{$K_v$-simple!decomposition!compatible with another decomposition} 
decomposition $(\ref{E_vSimple})$. 

$(C)$ If $\Char(K_v) = p > 0$, then for each $x_i$ the leading coefficient 
\index{coefficients $A_{v,ij}$!leading}
$c_{v,i} = \lim_{z \rightarrow x_i} f_v(z) \cdot g_{x_i}(z)^{N s_i}$ belongs to $K_v(x_i)^{\sep}$. 
\end{theorem} 

The proof of Theorem \ref{CompactThm} will occupy the remainder of this chapter. 
\index{initial approximating functions $f_v(z)$!for nonarchimedean $K_v$-simple sets|(}   

\section{ Reduction to a Set $E_v$ in a Single Ball} 
\label{NonArchProofReductionsSection} 

In this section we will reduce proving Theorem \ref{CompactThm} 
to proving it over a finite separable extension $F_{w}/K_v$, 
in the case where $E_v = \cC_v(F_w) \cap B(a,r)$ 
for a single isometrically parametrizable ball.  
\index{isometrically parametrizable ball}
\index{initial approximating functions $f_v(z)$!construction when $K_v$ nonarchimedean!reduction to the case of a single ball|(} 

To do this, we first recall some facts about nonarchimedean $(\fX,\vs)$-capacities  
established in Appendix \ref{AppA}. 
Let $\vs = (s_1, \ldots, s_m) \in \cP^m$ be a probability vector.  
As in \S\ref{Chap3}.\ref{LRationalBasisSection}, define the $(\fX,\vs)$-canonical distance by
\index{canonical distance!$[z,w]_{\fX,\vs}$|ii}
\begin{equation*}
[z,w]_{\fX,\vs} \ = \ \prod_{i=1}^m ([z,w]_{x_i})^{s_i} \ ,
\end{equation*}
\index{canonical distance!$[z,w]_{\fX,\vs}$!normalization of}
where  the $[z,w]_{x_i}$ are normalized 
so that  $\lim_{z \rightarrow x_i} [z,w]_{x_i} \cdot |g_{x_i}(z)|_v = 1$.  

Given a compact set $H_v$ disjoint from $\fX$, define its $(\fX,\vs)$-Robin constant by 
\begin{equation} \label{FCapDef1} 
V_{\fX,\vs}(H_v) \ = \ \inf_{\nu} \iint_{H_v \times H_v} -\log_v([z,w]_{\fX,\vs}) \, 
                                      d\nu(z) d\nu(w) 
\end{equation} 
\index{Robin constant!nonarchimedean!nonarchimedean $(\fX,\vs)$}
where $\nu$ runs over all Borel probability measures supported on $H_v$,
and its $(\fX,\vs)$-capacity by 
\index{capacity!$(\fX,\vs)$}
\begin{equation*} 
\gamma_{\fX,\vs}(H_v) \ = \ q_v^{-V_{\fX,\vs}(H_v)} \ ,
\end{equation*}
where $q_v$ is the order of the residue field of $\cO_v$.  

By Theorem \ref{ATE10B}, if $H_v$ has positive capacity,\index{capacity $> 0$}  
there is a unique probability measure 
$\mu_{\fX,\vs}$ supported on $H_v$ for which the \ $\inf$ \  in (\ref{FCapDef1}) 
is achieved;  this measure is called the $(\fX,\vs)$-equilibrium distribution. 
The $(\fX,\vs)$-potential function is defined by\index{potential function!$(\fX,\vs)$}   
\begin{equation*}
u_{\fX,\vs}(z;H_v) \ = \ \int_{H_v} -\log_v([z,w]_{\fX,\vs}) \, d\mu_{\fX,\vs}(w) \ .
\end{equation*}
Here  $u_{\fX,\vs}(z;H_v) = V_{\fX,\vs}(H_v)$ for all $z \in H_v$ 
\index{Robin constant!nonarchimedean!nonarchimedean $(\fX,\vs)$}  
except possibly a set of inner capacity $0$;  
\index{capacity $= 0$} 
and $u_{\fX,\vs}(z;H_v) < V_{\fX,\vs}(H_v)$ for all $z \notin H_v$.   
By Proposition \ref{BPropF1}, $V_{\fX,\vs}(H_v)$ is a continuous function 
of $\vs \in \cP^m$, and the $(\fX,\vs)$-equilibrium distribution and the 
$(\fX,\vs)$-Green's function  $G_{\fX,\vs}(z,H_v) = V_{\fX,\vs}(H_v) - u_{\fX,\vs}(z;H_v)$
\index{Green's function!$(\fX,\vs)$} 
\index{Robin constant!nonarchimedean!nonarchimedean $(\fX,\vs)$}
can be decomposed in terms of the corresponding 
objects for the individual points $x_i$:  
\begin{eqnarray}
\mu_{\fX,\vs} & = & \sum_{i=1}^m s_i \mu_i \ , \label{FSum1} \\
G_{\fX,\vs}(z;H_v) & = & \sum_{i=1}^m s_i G(z,x_i;H_v) \ , \label{FSum2}
\end{eqnarray}\index{Green's function!nonarchimedean!$(\fX,\vs)$}
where $\mu_i$ is the equilibrium distribution of $H_v$ with respect to the point
$x_i$, and $G(z,x_i;H_v) = V_{x_i}(H_v) - u_{x_i}(z;H_v)$. 
\index{Robin constant!nonarchimedean!nonarchimedean $(\fX,\vs)$}  

If $H_v$ is $K_v$-simple, with the $K_v$-simple decomposition 
\index{$K_v$-simple!decomposition}  
\begin{equation} \label{H_vK_vSimple}
H_v \ = \ \bigcup_{\ell=1}^D \big(\cC_v(F_{w_\ell}) \cap B(a_\ell,r_\ell)\big) \ ,
\end{equation} 
then by Lemma \ref{BFLem2} and Proposition \ref{BFProp2} of Appendix \ref{AppA}, 
the exceptional set of inner capacity $0$, discussed above, is empty:\index{capacity $= 0$}   
$u_{\fX,\vs}(z;H_v) = V_{\fX,\vs}(H_v)$ for all $z \in H_v$. 
\index{Robin constant!nonarchimedean!nonarchimedean $(\fX,\vs)$}  
Moreover, the equilibrium
distribution $\mu_{\fX,\vs}$ can be described as follows. 
For each $\ell$, write $H_{v,\ell} = \cC_v(F_{w_\ell}) \cap B(a_\ell,r_\ell)$ and 
let $\sigma_\ell : D(0,r_{\ell}) \rightarrow B(a_{\ell},r_{\ell})$ 
be an $F_{w_\ell}$-rational isometric parametrization with $\sigma_{\ell}(0) = a_{\ell}$.
\index{isometric parametrization}   
Let $\mu_{\ell}^*$ be the pushforward 
of additive Haar measure on $F_{w_\ell} \cap D(0,r_{\ell})$ to $H_{v,\ell}$ 
\index{Haar measure} 
by $\sigma_{\ell}$, normalized to have mass $1$.  
By Corollary \ref{BFCor2} of Appendix \ref{AppA}, 
there are weights $w_{\ell}(\vs) > 0$, 
satisfying $\sum_{\ell = 1}^D w_{\ell}(\vs) = 1$, for which 
\begin{equation}\label{FMuDecomp}
\mu_{\fX,\vs} \ = \ \sum_{\ell=1}^D w_{\ell}(\vs) \mu_{\ell}^* \ .
\end{equation}  
The weights $w_{\ell}(\vs)$ are uniquely determined by the requirement that
\begin{equation} \label{FSumF}
u_{\fX,\vs}(z,H_v) \ = \ \sum_{\ell=1}^D w_{\ell}(\vs) u_{\fX,\vs}(z,H_{v,\ell})
\end{equation}
takes the same value $V_{\fX,\vs}(H_v)$ on each $H_{v,\ell}$.
\index{Robin constant!nonarchimedean!nonarchimedean $(\fX,\vs)$} 

\smallskip
This description of $u_{\fX,\vs}(z,H_v)$ leads to the following system of linear equations.  Writing 
$V = V_{\fX,\vs}(H_v)$ and $w_\ell = w_{\ell}(\vs)$, 
\index{Robin constant!computing nonarchimedean}   
and evaluating $u_{\fX,\vs}(z,H_v)$ at a generic point 
of each $H_{v,\ell}$, we have (see Theorem \ref{BFThm2A}):   
\begin{equation} \label{FSystem}
\left\{ \begin{array}{ccl}
        1  & = & 0 \cdot V + \sum_{\ell=1}^D w_{\ell} \ ,     \\   
0 & = & V + w_\ell \cdot \big(-V_{\fX,\vs}(H_{v,\ell})\big) 
           + \displaystyle{\sum^D_{\substack{ j=1 \\ j \ne \ell }}} 
            w_j  \cdot \log_v \big([a_\ell,a_j]_{\fX,\vs}\big)  \\
 &  & \qquad \qquad \text{for $\ell = 1, \ldots, D$.}  
          \end{array} \right.
\end{equation}\index{Robin constant!nonarchimedean!nonarchimedean $(\fX,\vs)$}
By Theorem \ref{BFThm2A} the system (\ref{FSystem}) has a unique solution. 
Since the $V_{\fX,\vs}(H_{v,\ell})$ and $[a_k,a_{\ell}]_{\fX,\vs}$ are continuous in $\vs$
(see Proposition \ref{BPropF1} and Proposition \ref{APropA2}(B1)), 
$V_{\fX,\vs}(H_v)$ and the $w_{\ell}(\vs)$ are continuous functions of $\vs$.  Since 
$w_{\ell}(\vs) = \mu_{\fX,\vs}(H_{v,\ell}) = \sum_{i=1}^m s_i \mu_i(H_{v,\ell})$ and $\mu_i(H_{v,\ell}) > 0$
for each $i$ and $\ell$,  
there is a constant $W_0 = W_0(H_v,\fX) > 0$ such that $w_{\ell}(\vs) \ge W_0$, uniformly for all 
$\ell$ and $\vs$. 

Since $H_v$ is $K_v$-simple, if $\vs \in \cP^m(\QQ)$ the coefficients of the system (\ref{FSystem})
\index{$K_v$-simple!set}  
are rational, so the solutions $V_{\fX,\vs}(H_v)$ and $w_\ell(\vs)$ are rational as well. 
\index{Robin constant!nonarchimedean!nonarchimedean $(\fX,\vs)$}
\index{Robin constant!nonarchimedean!takes on rational values|ii}
\index{Green's function!nonarchimedean!takes on rational values|ii}
\index{weights!for nonarchimedean equilibrium distribution are rational|ii}      
This fact is shown in Corollary \ref{BFCor2}, 
but because it is crucial to our construction, and because the ideas motivate 
later parts of the construction, we repeat the reasoning here. 
It involves computations with explicit examples.   

\smallskip
First suppose $\cC = \PP^1/K$. 
Identify $\PP^1_v(\CC_v)$ with $\CC_v \bigcup \{\infty\}$, and take $\fX = \{\infty\}$.  
Let $F_w/K_v$ be a finite, separable extension embedded in $\CC_v$.  
We do not assume $F_w/K_v$ is galois.  
Consider the set $H_v^0 = \cO_w$, where $\cO_w$ is the ring of integers of $F_w$.  
The canonical distance $[x,y]_\infty$ 
\index{canonical distance!$[z,w]_{\infty}$}
(with respect to the 
uniformizing parameter\index{uniformizing parameter!normalizes canonical distance} $g_{\infty}(z) = 1/z$)  
is just $|x-y|_v$, which is unchanged when $x$ and $y$ are replaced by $x-a$ and $y-a$.   
It follows that the equilibrium measure $\mu = \mu_{\infty}$ is translation invariant 
under the additive group of $\cO_w$, and hence must be additive Haar measure on $\cO_w$. 
\index{Haar measure}  
The Robin constant and potential function\index{potential function!of $\cO_w$}  
of $\cO_w$ can be computed explicitly (see Lemma \ref{BFLem2}, or \cite{RR1}, Example 4.1.24, p.212)  
and in particular, writing $q_w = q_v^{f_w}$ for the order of the residue field of $\cO_w$, 
we have  
\begin{eqnarray}  
u_{\infty}(z,\cO_w) & = & \left\{ \begin{array}{ll} 
                   \ \frac{1}{e_w(q_w-1)} & \text{if $z \in \cO_w$\ ,} \\
                   -\log_v(|z|_v) & \text{if $|z|_v > 1$\ .} 
                                  \end{array} \right. \label{FFirst}
\end{eqnarray}  

Next let $\cC/K$ and $\fX$ be arbitrary, 
and consider the set  $H_{v,\ell} = \cC_v(F_{w_\ell}) \cap B(a_\ell,r_\ell)$ 
where $B(a_\ell,r_\ell)$ is isometrically parametrizable, 
\index{isometrically parametrizable ball}
$F_{w_\ell}/K_v$ is a finite separable extension in $\CC_v$, 
$a_\ell \in \cC_v(F_{w_\ell})$, and $r_\ell \in |F_{w_\ell}^{\times}|_v$.  
Put $q_{w_\ell} = q_v^{f_{w_\ell}}$.  
Fix an $F_{w_\ell}$-rational isometric parametrization 
\index{isometric parametrization}  
$\sigma_{\ell}: D(0,r_\ell) \rightarrow B(a_\ell,r_\ell)$
with $\sigma_\ell(0) = a_\ell$, and use it to identify $B(a_\ell,r_\ell)$ with $D(0,r_\ell)$.  
For each $x_i$ there is a constant $A_{\ell,i} \in |\CC_v^{\times}|_v$ 
such that $[y,z]_{x_i} = A_{\ell,i} |y-z|_v$ for all $z,w \in B(a_\ell,r_\ell)$.  
It follows that for $y,z \in B(a,r)$, we have $[y,z]_{\fX,\vs} = C_{\fX,\vs}\|z,w\|_v$ where 
$C_{\fX,\vs,\ell} = \prod_i A_{\ell,i}^{s_i}$.  
By Proposition \ref{APropA2} the canonical distance is constant on pairwise disjoint
\index{canonical distance!nonarchimedean!constant on disjoint balls}
isometrically parametrizable balls disjoint from $\fX$, 
\index{isometrically parametrizable ball}
so for  $y \in B(a_\ell,r_\ell)$ and $z \notin B(a_\ell,r_\ell)$ we have 
$[z,y]_{\fX,\vs} = [z,a_\ell]_{\fX,\vs}$.  The equilibrium distribution of $H_{v,\ell}$ 
is the pushforward of additive Haar measure on $F_{w_\ell} \cap D(0,r)$.
\index{Haar measure}  
Hence, by (\ref{FFirst})
\begin{equation}  
u_{\fX,\vs}(z,H_{v,\ell}) \ = \ \left\{ \begin{array}{ll} 
         \ \frac{1}{e_{w_\ell}(q_{w_\ell}-1)} 
         - \log_v(C_{\fX,\vs,\ell} \cdot r_{\ell}) & \text{if $z \in H_{v,\ell}$\ ,} \\
   -\log_v([z,a_\ell]_{\fX,\vs}) & \text{if $z \notin B(a_\ell,r_\ell)$\ .} 
                                  \end{array} \right. \label{FSecond} 
\end{equation} 
In particular, 
$V_{\fX,\vs}(H_{v,\ell}) = 1/(e_{w_\ell}(q_{w_\ell}-1)) -\log_v(C_{\fX,\vs,\ell} \cdot r_{\ell})$.
\index{Robin constant!nonarchimedean!nonarchimedean $(\fX,\vs)$}
Note that $C_{\fX,\vs,\ell}$ belongs to $|\CC_v^{\times}|_v$ 
if all the $s_i$ are rational.  
Thus, if $\vs \in \cP^m(\QQ)$, then $V_{\fX,\vs}(H_{v,\ell}) \in \QQ$. 
\index{Robin constant!nonarchimedean!nonarchimedean $(\fX,\vs)$}
Likewise, by Proposition \ref{APropA2}(B1),  $-\log_v([a_\ell,a_j]_{\fX,\vs}) \in \QQ$
for all $\ell \ne j$. 

Hence if $\vs \in \cP^m(\QQ)$, the coefficients of the system (\ref{FSystem}) are rational. 
It follows that the solution is rational as well.        

\smallskip
We next apply the theory above to sets 
$\tE_v \subset E_v \subset \cC_v(\CC_v) \backslash \fX$, 
where $E_v$ and $\tE_v$ are $K_v$-simple and $\tE_v$ is compatible with $E_v$.
\index{$K_v$-simple!set} 
\index{$K_v$-simple!decomposition}     
Given a $K_v$-simple decomposition $E_v = \bigcup_{\ell = 1}^D \cC_v(F_{w_\ell}) \cap B(a_\ell,r_\ell)$,
write $E_{v,\ell} = \cC_v(F_{w_\ell}) \cap B(a_\ell,r_\ell)$, 
and put $\tE_{v,\ell} = \tE_v \cap B(a_\ell,r_\ell) \subset E_{v,\ell}$.  
Thus $E_v = \bigcup_{\ell=1}^D E_{v,\ell}$ and $\tE_v = \bigcup_{\ell=1}^D \tE_{v,\ell}$.

\begin{lemma} \label{ComparisonLemma} 
Let $E_v \subset \cC_v(\CC_v) \backslash \fX$ 
be a $K_v$-simple set,\index{$K_v$-simple!decomposition}\index{$K_v$-simple!set}  
with a $K_v$-simple decomposition $E_v = \bigcup_{\ell = 1}^D \cC_v(F_{w_\ell}) \cap B(a_\ell,r_\ell)$.

Let $\varepsilon_v > 0$ be given.  Then there is a $\delta_v > 0$ such that 
for any $K_v$-simple set $\tE_v \subset E_v$, if  
$|V_{x_i}(\tE_{v,\ell}) - V_{x_i}(E_{v,\ell})| < \delta_v$ for all $\ell = 1, \ldots, D$
\index{Robin constant!nonarchimedean}   
and $i = 1, \ldots, m$, then for each $x_i \in \fX$ we have 
\begin{equation} \label{FRRApprox4C}
|V_{x_i}(\tE_v) - V_{x_i}(E_v)| \ < \ \varepsilon_v \ ,
\end{equation} 
and for all $x_i, x_j \in \fX$ with $x_i \ne x_j$, 
\begin{equation} \label{FRRApprox5C} 
|G(x_i,x_j;\tE_v) - G(x_i,x_j;E_v)| \ < \ \varepsilon_v \ .
\end{equation}
\end{lemma} 

\begin{proof} 
By Theorem \ref{BFThm2A}, there are systems of equations analogous to (\ref{FSystem})
for the sets $E_v$ and $\tE_v$: 
write $V = V_{\fX,\vs}(E_v)$ and $w_{\ell} = w_{\ell}(\vs)$ for the Robin constant and weights  
\index{Robin constant!nonarchimedean!nonarchimedean $(\fX,\vs)$}
associated to $E_v$; let $\tV = V_{\fX,\vs}(\tE_v)$, $\tw_{\ell} = \tw_{\ell}(\vs)$ 
be the corresponding objects for $\tE_v$.  Then  
\begin{equation} \label{BFFM6A0}
\left\{ \begin{array}{ccl}
        1  & = & 0 \cdot V + \sum_{\ell=1}^D w_{\ell} \ ,     \\   
        0  & = & V + w_\ell \cdot \big(-V_{\fX,\vs}(E_{v,\ell})\big) 
           + \displaystyle{\sum^D_{\substack{j=1 \\ j \ne \ell }}} 
            w_j  \cdot \log_v \big([a_\ell,a_j]_{\fX,\vs}\big),   
             \quad \text{for $\ell = 1, \ldots, D$}  
          \end{array} \right.
\end{equation} 
\index{Robin constant!nonarchimedean!nonarchimedean $(\fX,\vs)$}
and 
\begin{equation} \label{BFFM6B}
\left\{ \begin{array}{ccl}
        1  & = & 0 \cdot \tV + \sum_{\ell=1}^D \tw_{\ell} \ ,     \\   
        0  & = & \tV + \tw_\ell \cdot \big(-V_{\fX,\vs}(\tE_{v,\ell})\big) 
           + \displaystyle{\sum^D_{\substack{ j=1 \\ j \ne \ell }}} 
                \tw_j  \cdot \log_v \big([a_\ell,a_j]_{\fX,\vs}\big),  
              \quad \text{for $\ell = 1, \ldots, D$\ .} 
          \end{array} \right.
\end{equation} 

Fix $x_i \in \fX$ and specialize to the case where $[x,y]_{\fX,\vs} = [x,y]_{x_i}$; 
that is, take $\vs = \vec{e}_i = (0,\ldots,1,\ldots,0)$.  
Note that the coefficients of the systems (\ref{BFFM6A0}), (\ref{BFFM6B}) 
are the same except for their diagonal terms.  
\index{Robin constant!nonarchimedean!nonarchimedean $(\fX,\vs)$}
Write $V_\ell = V_{\fX,\vs}(E_{v,\ell})$ and  
$\tV_{\ell} = V_{\fX,\vs}(\tE_{v,\ell})$, for $\ell = 1, \ldots, D$.

Henceforth we will regard $\tV_1, \ldots, \tV_D$ as variables, keeping $V_1, \ldots, V_D$ fixed.  
By Cramer's rule,\index{Cramer's Rule} the solution vectors $(V,w_1,\ldots,w_D)$ and $(\tV,\tw_1,\ldots,\tw_D)$ 
are continuous functions of the coefficients in (\ref{BFFM6A0}), (\ref{BFFM6B}). 
It follows that when $(\tV_1, \ldots, \tV_D) \longrightarrow (V_1, \ldots, V_D)$, then  
$(\tV,\tw_1, \ldots, \tw_D) \longrightarrow (V,w_1, \ldots, w_D)$.  In particular, 
$V_{x_i}(\tE_v) \longrightarrow V_{x_i}(E_v)$.  

Furthermore, for each $z \ne x_i$ we have 
\begin{equation} \label{GEq1} 
\left\{ \begin{array}{l}
     G(z,x_i;E_v) \ = \  V - u_{x_i}(x,E_v) 
             \ = \ V - \sum_{\ell=1}^D w_{\ell} \, u_{x_i}(z,E_{v,\ell}) \ , \\
     G(z,x_i;\tE_v) \ = \  \tV - u_{x_i}(x,\tE_v) 
               \ = \ \tV - \sum_{\ell=1}^D \tw_{\ell} \, u_{x_i}(z,\tE_{v,\ell}) \ . 
\end{array} \right. 
\end{equation} 
\index{Green's function!computing nonarchimedean}
Since $\fX$ is disjoint from $\bigcup_{\ell=1}^D B(a_{\ell},r_{\ell})$,
it follows from (\ref{FSecond}) that for each $x_j \ne x_i$ 
we have $u_{x_i}(x_j,E_{v,\ell}) = u_{x_i}(x_j,\tE_{v,\ell}) = -\log_v([x_j,a_{\ell}]_{x_i})$.
Thus, when $(\tV_1, \ldots, \tV_D) \longrightarrow (V_1, \ldots, V_D)$, 
we have $G(x_j,x_i;\tE_v) \longrightarrow G(x_j,x_i;E_v)$ as well.

Since $\fX$ is finite, the Lemma follows.  
\end{proof}

We now claim that in order to prove Theorem $\ref{CompactThm}$, 
it suffices to establish the following 
weaker version of the theorem for $E_{v,\ell}$ over $F_{w_\ell}$, 
for each $\ell = 1, \ldots, D$:

\begin{proposition} \label{K_vReductionProp}
Let $E_{v,\ell} = \cC_v(F_{w_\ell}) \cap B(a_{\ell},r_{\ell})$, 
where $F_{w_\ell}$ is a finite, separable extension of $K_v$ in $\CC_v$, 
$B(a_{\ell},r_{\ell}) \subset \cC_v(\CC_v) \backslash \fX$ 
is an isometrically parametrizable ball, and $a_{\ell} \in \cC_v(F_{w_\ell})$. 
\index{isometrically parametrizable ball}
Let $\varepsilon_{v,\ell} > 0$ be given.
Then there is a compact set $\tE_{v,\ell} \subseteq E_{v,\ell}$ for which

$(A)$ There are points $\alpha_{\ell,j} \in \cC_v(F_{w_\ell}) \cap B(a_\ell,r_\ell)$ 
and pairwise disjoint isometrically parametrizable balls 
\index{isometrically parametrizable ball}
$B(\alpha_{\ell,j},r_{\ell,j}) \subseteq B(a_\ell,r_\ell)$, for $j = 1, \ldots, d_\ell$, 
such that $\tE_{v,\ell}$ has the form  
\begin{equation} \label{FtEvLF} 
\tE_{v,\ell} \ = \ \bigcup_{j=1}^{d_\ell} \big(\cC_v(F_{w_\ell}) \cap B(\alpha_{\ell,j},r_{\ell,j})\big) 
\end{equation} 
and for each $x_i \in \fX$ 
\begin{equation} \label{FRRApprox4B1}
|V_{x_i}(\tE_{v,\ell}) - V_{x_i}(E_{v,\ell})| \ < \ \varepsilon_{v,\ell} \ .
\end{equation}
\index{Robin constant!nonarchimedean}    
   
$(B)$  For each $0 < \beta_v \in \QQ$  
and each $F_{w_\ell}$-symmetric probability vector
$\vs = {}^t(s_1, \ldots, s_m)$ with rational entries, 
there is an integer $N_{v,\ell} \ge 1$  such that 
for each positive integer  $N$ divisible by $N_{v,\ell}$, 
there is an $(\fX,\vs)$-function\index{$(\fX,\vs)$-function!$K_v$-rational}
 $f_{v,\ell} \in F_{w_\ell}(\cC_v)$ of degree $N$ such that
   
\quad $(1)$  For all $z \in \cC_v(\CC_v) \backslash (B(a_\ell,r_\ell) \cup \fX)$, 
\begin{equation} \label{FRRApprox7B1} 
\frac{1}{N} \log_v(|f_{v,\ell}(z)|_v) \ = \ G_{\fX,\vs}(z,\tE_{v,\ell}) + \beta_v \ .
\end{equation}   
   
\quad $(2)$ The zeros $\theta_1, \ldots, \theta_N$ of $f_{v,\ell}$ 
are distinct and belong to $E_{v,\ell}$ $($hence $\cC_v(F_{w_\ell}))$.

\quad $(3)$  
$f_{v,\ell}^{-1}(D(0,1)) = \bigcup_{h=1}^N B(\theta_h,\rho_h)$, 
where the balls $B(\theta_h,\rho_h)$ are pairwise disjoint and contained in $B(a_{\ell},r_{\ell})$. 
%
\end{proposition}   

Proposition \ref{K_vReductionProp} 
will be proved in Section \ref{Chap6}.\ref{NonArchProofSection}. 

\smallskip
\begin{proof}[Proof of Theorem \ref{CompactThm}, 
assuming Proposition \ref{K_vReductionProp}] { \ } 

Let $E_v \subset \cC_v(\CC_v) \backslash \fX$ be a $K_v$-simple set 
\index{$K_v$-simple!set}  
with the $K_v$-simple decomposition 
\index{$K_v$-simple!decomposition}  
\begin{equation} \label{E_vSimple1} 
E_v \ = \ \bigcup_{\ell=1}^D \big(\cC_v(F_{w_\ell}) \cap B(a_{\ell},r_{\ell})\big) \ . 
\end{equation} 
For each $\ell$, write $E_{v,\ell} = \cC_v(F_{w_\ell}) \cap B(a_{\ell},r_{\ell})$.   

By the definition of a $K_v$-simple decomposition, 
\index{$K_v$-simple!decomposition}  
the collection of balls $\{B(a_\ell,r_{\ell})\}_{1 \le \ell \le D}$ 
is stable under $\Aut_c(\CC_v/K_v)$ and for each $\ell$ the ball $B(a_{\ell},r_{\ell})$ 
has $[F_{w_\ell}:K_v]$ distinct conjugates.  
For each $\sigma \in \Aut_c(\CC_v/K_v)$ and each $\ell = 1, \ldots, D$, 
let $\sigma(\ell)$ be the index
such that $B(a_{\sigma(\ell)},r_{\sigma(\ell)}) = B(a_{\ell},r_{\ell})$.  
Then $F_{\sigma(\ell)} = \sigma(F_{\ell})$. 

We first construct the set $\tE_v$ in Theorem \ref{CompactThm}.  
Given $\varepsilon_v > 0$, let $\delta_v > 0$ be the
number given by Lemma \ref{ComparisonLemma} for $E_v$ 
and the $K_v$-simple decomposition (\ref{E_vSimple1}).  
\index{$K_v$-simple!decomposition}  
Suppose that under the action of $\Aut_c(\CC_v/K_v)$ on 
$B(a_1,r_1), \ldots, B(a_D,r_D)$, there are $T$ distinct orbits.    
We can assume without loss that  $B(a_1,r_1), \ldots, B(a_T,r_T)$ 
are representatives for the orbits.  For each $\ell = 1, \dots, T$, 
take $\varepsilon_{v,\ell} = \delta_v$ 
and let $\tE_{v,\ell} \subset E_{v,\ell}$ and $N_{v,\ell}$
be the $F_{w_\ell}$-simple set and number given for $E_{v,\ell}$ and $\varepsilon_{v,\ell}$, 
by Proposition \ref{K_vReductionProp}. 
Let $\tV_{\ell} = V_{\fX,\vs}(\tE_{v,\ell})$ be the 
$(\fX,\vs)$-Robin constant of $\tE_{v,\ell}$.  We define the sets 
\index{Robin constant!nonarchimedean!nonarchimedean $(\fX,\vs)$} 
$\tE_{v,\ell}$ for $T+1 \le \ell \le D$ by galois conjugacy:  given such an $\ell$, 
there are a $k$ with $1 \le k \le T$ and a $\sigma \in \Aut_c(\CC_v/K_v)$ such that 
$\ell = \sigma(k)$.  Put $\tE_{v,\ell} = \sigma(\tE_{v,k})$.  It is easy to see that 
$\tE_{v,\ell}$ is independent of the choice of $\sigma$, 
and that each $\tE_{v,\ell}$ is $F_{w_\ell}$-simple.  With this definition,
for each $1 \le \ell \le D$ and each $\sigma \in \Aut_c(\CC_v/K_v)$ we have 
$\tE_{v,\sigma(\ell)} = \sigma(\tE_{v,\ell})$ 
and $V_{\fX,\vs}(\tE_{v,\sigma(\ell)}) = V_{\fX,\vs}(\tE_{v,\ell})$.   Put 
\index{Robin constant!nonarchimedean!nonarchimedean $(\fX,\vs)$}  
\begin{equation*}
\tE_v \ = \ \bigcup_{\ell = 1}^D \tE_{v,\ell} \ .
\end{equation*} 
By construction $\tE_v$ is $K_v$-simple and compatible with $E_v$.  
\index{$K_v$-simple!set}\index{$K_v$-simple!set!compatible with another set}
By Proposition \ref{K_vReductionProp}, part (A) of Theorem \ref{CompactThm} 
holds for $\tE_v$ and $E_v$.  

\smallskip 
We next establish part (B).  Fix $0 < \beta_v \in \QQ$ 
and a $K_v$-symmetric probability vector $\vs \in \cP^m(\QQ)$.
\index{$K_v$-symmetric!probability vector}  
Let $\tV = V_{\fX,\vs}(\tE_v)$ and $\tw_\ell = \tw_\ell(\vs)$, for $\ell = 1, \ldots, D$, 
\index{Robin constant!nonarchimedean!nonarchimedean $(\fX,\vs)$} 
be the solutions to the system of equations (\ref{FSystem}) associated to $\tE_v$.  
Since  $\tE_v$, $\fX$ and $\vs$ are $K_v$-symmetric, 
\index{$K_v$-symmetric!probability vector}
we have $V_{\fX,\vs}(\tE_{v,\sigma(\ell)}) = V_{\fX,\vs}(\tE_{v,\ell})$ for all $\ell$ 
and $\sigma$, and each $\sigma$ permutes the equations (\ref{FSystem}). 
Hence the weights $\tw_\ell$ satisfy $\tw_{\sigma(\ell)} = \tw_\ell$
for all $\ell$ and all $\sigma$.    
 
Since each $\tE_{v,\ell}$ is $F_{w_\ell}$-simple, Corollary \ref{BFCor2} 
(applied to an $F_{w_\ell}$-simple decomposition of $\tE_{v,\ell}$)
shows that each $V_{\fX,\vs}(\tE_{v,\ell})$ is rational.  
\index{Robin constant!nonarchimedean!nonarchimedean $(\fX,\vs)$}
\index{Robin constant!nonarchimedean!takes on rational values} 
By Theorem \ref{BFThm2A}, $\tV = V_{\fX,\vs}(\tE_v)$ and the 
weights $\tw_1, \ldots, \tw_D$ are rational.
Write $\tw_{\ell} = P_{\ell}/Q_{\ell}$ with positive integers $P_{\ell}, Q_{\ell}$,
and let $Q$ be the least common multiple of $Q_1, \ldots, Q_T$.  
Write $V_{\fX,\vs}(\tE_v) = X_0/Y_0$ and 
$V_{\fX,\vs}(\tE_{v,\ell}) = X_{\ell}/Y_{\ell}$ with integers 
$X_{\ell}, Y_{\ell}$, and put $Y = \LCM(Y_0,Y_1,\ldots, Y_T)$.  

\smallskip
Put $\tN_v = \LCM(N_{v,1},\ldots,N_{v,T})$, and set $N_v = Y Q \tN_v$.  

Suppose $N$ is a multiple of $N_v$, say $N = k N_v$.  
For each $\ell = 1, \ldots, T$ put $n_{\ell} = \tw_{\ell} N$, noting that $n_{\ell} \in \NN$ 
and that $N_{v,\ell} | n_{\ell}$.  Let $f_{v,\ell}(z) \in F_{w_\ell}(\cC)$ 
be the $(\fX,\vs)$-function\index{$(\fX,\vs)$-function}
 of degree $n_{\ell}$ given by Proposition \ref{K_vReductionProp} 
for $\tE_{v,\ell}$, $E_{v,\ell}$, $\vs$, and $\beta_v$.  
For the remaining sets $\tE_{v,\ell}$ with $\ell = T+1,\ldots, D$,
define $f_{v,\ell}$ by conjugacy, so that if $\ell = \sigma(k)$ with $1 \le k \le T$
and $\sigma \in \Aut_c(\CC_v/K_v)$ then 
$f_{v,\ell} = (f_{v,k})^{\sigma} = \sigma \circ f_{v,k} \circ \sigma^{-1}$.  
It follows that for each $\ell = 1, \ldots, D$,
the function $f_{v,\ell}$ belongs to $F_{w_\ell}(\cC_v)$, has degree $\tw_{\ell} N$, 
and satisfies the conditions of Proposition \ref{K_vReductionProp} 
relative to $\tE_{v,\ell}$, $\vs$, and $\beta_v$.  In particular, 
for all $z \in \cC_v(\CC_v) \backslash \big(B(a_\ell,r_\ell) \cup \fX\big)$, we have   
\begin{equation} \label{FSingleApprox}
\frac{1}{n_{\ell}} \log_v(|f_{v,\ell}(z)|_v) \ = \ G_{\fX,\vs}(z,\tE_{v,\ell}) + \beta_v \ .
\end{equation}   
Clearly $f_{v,\sigma(\ell)} = (f_{v,\ell})^{\sigma}$ 
for all $\sigma \in \Aut_c(\CC_v/K_v)$ and all $\ell$. 

Note that for each $\ell$ we have 
$u_{\fX,\vs}(z,\tE_{\ell}) = V_{\fX,\vs}(\tE_{v,\ell}) -  u_{\fX,\vs}(z,\tE_{v,\ell})$, 
\index{Robin constant!nonarchimedean!nonarchimedean $(\fX,\vs)$}  
and that 
\begin{equation*}
G_{\fX,\vs}(z,\tE_v) 
\ = \ V_{\fX,\vs}(\tE_v) - \sum_{\ell = 1}^D \tw_{\ell} \, u_{\fX,\vs}(z,\tE_{v,\ell}) \ . 
\end{equation*}
Thus if we put $C = V_{\fX,\vs}(\tE_v) - \sum_{\ell=1}^D \tw_{\ell} V_{\fX,\vs}(\tE_{v,\ell})$, 
\index{Robin constant!nonarchimedean!nonarchimedean $(\fX,\vs)$} 
then 
\begin{equation} \label{FGXs1}
G_{\fX,\vs}(z,\tE_v) \ = \ C + \sum_{\ell=1}^D \tw_{\ell} \, G_{\fX,\vs}(z,\tE_{v,\ell}) \ .
\end{equation} 
\index{Green's function!nonarchimedean!$(\fX,\vs)$}   
By our choice of $N$, we have $N \cdot C \in \ZZ$.  
Since $G_{\fX,\vs}(z,\tE_v) = 0$ for all $z \in \tE_v$, 
and $G_{\fX,\vx}(z,\tE_{v,\ell}) \ge 0$ for all $z$ and all $\ell$, 
by evaluating both sides of (\ref{FGXs1})) at a point $z \in \tE_v$ we see that $C \le 0$.  

Let $\pi_v$ be a uniformizing element for the maximal ideal of $\cO_v$, and define 
\begin{equation} \label{NonarchFcnDef} 
f_v(z) \ = \ \pi_v^{-NC} \cdot \prod_{\ell=1}^D  f_{v,\ell}(z) \ . 
\end{equation} 
Since each $f_{v,\ell}$ is an $(\fX,\vs)$-function,\index{$(\fX,\vs)$-function} so is $f_v$.  
By construction $f_v$ is stable under $\Aut_c(\CC_v/K_v)$.  
Since each $F_{w_\ell}/K_v$ is separable, $f_v$ belongs to $K_v(\cC)$.  
It clearly has degree $N$.   

\smallskip
We will now show that $f_v$ satisfies the conditions of Theorem \ref{CompactThm} relative 
to $\tE_v$, $\vs$, and $\beta_v$ and the $K_v$-simple decomposition (\ref{E_vSimple1}).
\index{$K_v$-simple!decomposition}    
First, by our hypotheses on the $f_{v,\ell}(z)$, 
for each $z \notin \bigcup_{\ell = 1}^D B(a_{\ell},r_{\ell}) \cup \fX$ we have
\begin{eqnarray*}
\frac{1}{N} \log_v(|f_v(z)|_v) 
& = & \frac{1}{N} \Big(NC + \sum_{\ell=1}^D \log_v(|f_{v,\ell}(z)|_v)\Big) \\ 
& = & \frac{1}{N} \Big(NC + \sum_{\ell=1}^D \tw_{\ell}\, N \cdot
           \big( G_{\fX,\vs}(z,\tE_{v,\ell}) + \beta_v\big) \Big) \\
& = & G_{\fX,\vs}(z,\tE_v) + \beta_v \ ,
\end{eqnarray*}
using (\ref{FGXs1}) and the fact that $\sum_{\ell=1}^D \tw_{\ell} = 1$.  In particular, 
for each $x_i \in \fX$, 
\begin{equation*} 
\Lambda_{x_i}(f_v,\vs) \ = \ 
\lim_{z \rightarrow x_i} \big(G_{\fX,\vs}(z,\tE_v) + s_i\log_v(|g_i(z)|_v)\big) + \beta_v
\ = \ \Lambda_{x_i}(\tE_v,\vs) + \beta_v \ .
\end{equation*}  

Next, we claim that the zeros of $f_v$ are distinct and belong to $E_v$.  
Indeed, for each $\ell$ the zeros of $f_{v,\ell}$ are distinct and belong to $E_{v,\ell}$. 
This holds for $\ell = 1, \ldots, T$ by Proposition \ref{K_vReductionProp}, and 
for the remaining $\ell$ by conjugacy.  Since the sets $E_{v,\ell}$ are pairwise disjoint,
our claim follows.   

Fix $\ell$.  Recalling that $n_{\ell} = \deg(f_{v,\ell}) = \tw_{\ell} N$, 
let $\theta_{\ell,1}, \ldots, \theta_{\ell,n_{\ell}}$ be the zeros of $f_{v,\ell}$.  
By Proposition \ref{K_vReductionProp} there pairwise disjoint balls 
$B(\theta_{\ell,1},\rho_{\ell,1}), \ldots, B(a_{\ell,n_{\ell}},\rho_{\ell,n_{\ell}})$ 
contained in $B(a_{\ell},r_{\ell})$ such that 
\begin{equation*} 
f_{v,\ell}^{-1}(D(0,1)) \ = \ \bigcup_{j=1}^{n_{\ell}} B(\theta_{\ell,j},\rho_{\ell,j}) \ . 
\end{equation*} 
Here, the balls $B(\theta_{\ell,j},\rho_{\ell,j})$ 
are isometrically parametrizable since $B(a_\ell,r_\ell)$ is isometrically parametrizable.
\index{isometrically parametrizable ball}
The $\theta_{\ell,j}$ belong to  $\cC_v(F_{w_\ell})$ since they belong to $E_{v,\ell}$. 
By choosing an $F_{w_\ell}$-rational isometric parametrization 
\index{isometric parametrization}  
$\varphi_{\ell,j} : D(0,r_\ell) \rightarrow B(a_{\ell},r_{\ell})$ 
with $\varphi_{\ell,j}(0) = \theta_{\ell,j}$, 
expanding $f_{v,\ell}$ as a power series $c_0 + c_1 Z + \cdots$, 
and applying Proposition \ref{ScaledIsometryProp}, 
we see that $f_{v,\ell}$ induces a scaled isometry from $B(\theta_{\ell,j},\rho_{\ell,j})$ to $D(0,1)$.
\index{scaled isometry}
Here $c_0, c_1, \cdots \in F_{w_\ell}$ since $f_{v,\ell} \in F_{w_\ell}(\cC_v)$ 
and $\theta_{\ell,j} \in \cC_v(F_{w_\ell})$.  
Proposition \ref{ScaledIsometryProp} gives $|c_1|_v \cdot \rho_{\ell,j} = 1$, 
so $\rho_{\ell,j} = 1/|c_1|_v \in |F_{w_\ell}^{\times}|_v$.


On the other hand, the function  $H_{\ell}(z) = f_v(z)/f_{v,\ell}(z)$ 
is also $F_{w_\ell}$-rational, 
and its zeros and poles are disjoint from $B(a_{\ell},r_{\ell})$.
Hence there is a constant $B_\ell$ such that $|H_{\ell}(z)|_v = B_\ell$ 
for all $z \in B(a_{\ell},r_{\ell})$.  
Evaluating $H_{\ell}(z)$ at a point $z_\ell \in \tE_{v,\ell}$, 
and successively using $G_{\fX,\vs}(z_{\ell},\tE_{v,\ell}) = 0$,
(\ref{FGXs1}), and $G_{\fX,\vs}(z_{\ell},\tE_v) = 0$, we see that  
\begin{eqnarray*} 
\log_v(B_\ell) & = & |\pi_v^{-N C} \cdot \prod_{k \ne \ell} f_{v,k}(z_{\ell})|_v 
\ = \ N C + \sum_{k \ne \ell} \tw_k N \cdot \big(G_{\fX,\vs}(z_{\ell},\tE_{v,k}) + \beta_v \big) \\
& = &  NC + \Big(\sum_{k=1}^D \tw_k N \cdot G_{\fX,\vs}(z_{\ell},\tE_{v,k})\Big)
 + N \cdot(1-\tw_{\ell}) \beta_v \\
& = & N \cdot G_{\fX,\vs}(z_0,\tE_v) + N \cdot (1-\tw_{\ell}) \beta_v 
\ = \ N \cdot (1-\tw_{\ell}) \beta_v \ \ge \ 0 \ , 
\end{eqnarray*}  
so $B_\ell \ge 1$.  Since $B_\ell = |\pi_v^{-NC} \cdot H_{\ell}(z_\ell)|_v$, where 
$H_{\ell} \in F_{w_\ell}(\cC_v)$ and $z_{\ell} \in \cC_v(F_{w_\ell})$, 
we have $B_\ell \in |F_{w_\ell}^{\times}|_v$.
 
Choose an $F_{w_\ell}$-rational isometric parametrization of $B(\theta_{\ell,j},\rho_{\ell,j})$, 
\index{isometric parametrization}  
and expand $f_{v,\ell}$ and $H_{\ell}$ as power series.
By Proposition \ref{ScaledIsometryProp}, 
$f_v = f_{v,\ell} \cdot H_{\ell}$ induces an $F_{w_\ell}$-rational scaled isometry 
\index{scaled isometry}
from $B(\theta_{\ell,j},B_\ell^{-1}\rho_{\ell,j})$ onto $D(0,1)$
which maps $\cC_v(F_{w_\ell}) \cap B(\theta_{\ell,j},B_\ell^{-1}\rho_{\ell,j})$
onto $\cO_{w_\ell}$, for each $j$.  
Clearly $B_\ell^{-1}\rho_{\ell,j} \in |F_{w_\ell}^{\times}|_v$. 

Now let $\ell$ vary.  For each $\ell$ and each $j$ we have   
\begin{equation*}
B(\theta_{\ell,j},B_\ell^{-1}\rho_{\ell,j}) \ \subseteq \ 
B(\theta_{\ell,j},\rho_{\ell,j}) \ \subseteq \ B(a_{\ell},r_{\ell}) \ .
\end{equation*}
For a given $\ell$ the balls $B(\theta_{\ell,j},\rho_{\ell,j})$ are pairwise disjoint,    
so the balls $B(\theta_{\ell,j},B_\ell^{-1}\rho_{\ell,j})$ are pairwise disjoint.  
For different $\ell$, the balls $B(a_\ell,r_\ell)$ are pairwise disjoint, so in fact the balls
$B(\theta_{\ell,j},B_\ell^{-1}\rho_{\ell,j})$ are pairwise disjoint for all $\ell$ and $j$.
There are exactly $N = \sum_{\ell=1} n_{\ell} = \deg(f_v)$ such balls, so 
\begin{equation*}
f_v^{-1}(D(0,1)) \ = \ \bigcup_{\ell=1}^D \bigcup_{j=1}^{n_\ell} 
B(\theta_{\ell,j},B_\ell^{-1}\rho_{\ell,j}) \ \subseteq \ \bigcup_{\ell=1}^D B(a_{\ell},r_{\ell}) \ .
\end{equation*}     
It follows that $H_v := E_v \cap f_v^{-1}\big(D(0,1)\big)$ is $K_v$-simple, 
\index{$K_v$-simple!set}  
and has the $K_v$-simple decomposition
\index{$K_v$-simple!decomposition}  
\begin{equation*}
H_v \ = \ \bigcup_{\ell=1}^D \bigcup_{j=1}^{n_\ell} 
\cC_v(F_{w_\ell}) \cap B(\theta_{\ell,j},B_\ell^{-1}\rho_{\ell,j}) 
\end{equation*}  
which is compatible with the $K_v$-simple decomposition (\ref{E_vSimple1}) of $E_v$.
\index{$K_v$-simple!decomposition!compatible with another decomposition}  
This completes the proof of part (B).

\smallskip
Finally, suppose $\Char(K_v) = p > 0$.  
We will show that by modifying the construction above, 
we can arrange that the leading coefficients 
\index{coefficients $A_{v,ij}$!leading}
\begin{equation*}
c_{v,i} \ = \ \lim_{z \rightarrow x_i} f_v(z) \cdot g_{x_i}(z)^{N_i}
\end{equation*}
belong to $K_v(x_i)^{\sep}$, so that part (C) holds. 
 
Fix a positive, $K_v$-symmetric probability vector $\vs \in \QQ^m$ 
\index{$K_v$-symmetric!probability vector}
and a number $0 < \beta_v \in \QQ$ as before, 
and carry out the construction in part (B) for $E_v$, $\tE_v$,   
and $\vs$, but with $\beta_v$ replaced by $\beta_v/2$.  
Let $N_{v,0} > 0$ be the integer given by part (B) for $\beta_v/2$;
note that $N_{v,0} \cdot \beta_v/2 \in \ZZ$. 
Put $p^B = \max_{1 \le i \le m}([K_v(x_i):K_v]^{\insep})$.
After replacing $N_{v,0}$ by a multiple of itself if necessary,
we can assume that  
\begin{equation*}
p^B N_{v,0} \cdot  \frac{\beta_v}{2} \ \ge \ \frac{p^B}{q_v-1} + \log_v(p^B) + 2 \ .
\end{equation*} 
We will take $N_v = p^B N_{v,0}$.

Given a positive integer $N$ divisible by $N_v$, put $N_0 = N/p^B$.
Then $N_0$ is divisible by $N_{v,0}$;   
let $f_{v,0} \in K_v(\cC)$ be the $(\fX,\vs)$-function\index{$(\fX,\vs)$-function} 
of degree $N_0$ constructed in part (B)
for $E_v$ and $\tE_v$, relative to $\vs$ and $\beta_v/2$.  
Noting that $N \cdot \beta_v/2 \in \ZZ$, compose $f_{v,0}$ 
with the Stirling polynomial $S_{p^B,\cO_v}(z)$ and put 
\index{Stirling polynomial!for $\cO_v$}
\begin{equation} \label{FStirlingCompose}
f_v(z) \ = \ \pi_v^{-N \beta_v/2} \cdot S_{p^B,\cO_v}(f_{v,0}(z)) \ .
\end{equation}
For each $i$, the leading coefficient $c_{v,i,0}$ of $f_{v,0}$ at $x_i$ 
\index{coefficients $A_{v,ij}$!leading}
belongs to $K_v(x_i)$, since $f_{v,0} \in K_v(\cC)$ and $g_{x_i}$ is rational
over $K_v(x_i)$.  This means that the leading coefficient of $f_v(z)$ at $x_i$ is
\index{coefficients $A_{v,ij}$!leading}
\begin{equation*}
c_{v,i} \ = \ \pi_v^{-N \beta_v/2} \cdot c_{v,i,0}^{p^B} \ , 
\end{equation*} 
which belongs to $K_v(x_i)^{\sep}$.  Thus part (C) of Theorem \ref{CompactThm} holds for $f_v$.  

We now show that $f_v(z)$ continues to satisfy properties (B.1)--(B.4) of of Theorem \ref{CompactThm}.   
First, note that since $\Lambda_{x_i}(f_{v,0},\vs) = \Lambda_{x_i}(\tE_v,\vs) + \beta_v/2$, 
we have 
\begin{eqnarray*}
\Lambda_{x_i}(f_v,\vs) 
    & = & \frac{1}{N} \log_v(|c_{v,i}|_v) 
    \ = \ \frac{1}{N} \cdot N \beta_v/2 +  \frac{1}{N} \cdot p^B \log_v(|c_{v,i,0}|_v) \\
    & = & \beta_v/2 +  \Lambda_{x_i}(f_{v,0},\vs) \ = \ \Lambda_{x_i}(\tE_v,\vs) + \beta_v  \ .
\end{eqnarray*}
This proves property (B.1).

For property (B.2), recall that the zeros of $S_{p^B,\cO_v}(z)$ are distinct and belong to $\cO_v$.
For each $k = 1, \ldots, N_0$, the function $f_{v,0}$ induces a 
scaled isometry from $B(\theta_k,\rho_k)$ onto $D(0,1)$,
\index{scaled isometry}
which takes $B(\theta_k,\rho_k) \cap \cC_v(F_{w_\ell(k)})$ onto $\cO_{w_\ell(k)}$.  
The zeros of $f_v$ in $B(\theta_k,\rho_k)$, 
which we will denote $\theta_{k,j}$ for $j = 0, \ldots, p^B-1$,  
therefore belong to $B(\theta_k,\rho_k) \cap \cC_v(F_{w_\ell(k)}) \subseteq E_v$.  
Letting $k$ vary, we see that $f_v$ has $N_0 \cdot p^B = N$ zeros in $E_v$. 
Since $\deg(f_v) = N$, these are all the zeros of $f_v$. 
Thus the zeros of $f_v$ are distinct and belong to $E_v$. 

For property (B.3), note that by Corollary \ref{SnOwMappingCor}, if 
$0 < R < q_v^{-p^B/(q_v-1)} \cdot (p^B)^{-1}$, then the inverse image of $D(0,R)$ under $S_{p^B,\cO_v}(z)$
consists of $p^B$ disjoint discs contained in $D(0,1)$, centered on the roots of $S_{p^B,\cO_v}(z)$.
In addition, if $R$ belongs to the value group of $K_v^{\times}$, 
the radii of those discs belong to the value group of $K_v^{\times}$. 

Take $R = q_v^{-\eta}$, where 
\begin{equation*} 
\eta \ = \ \lceil \frac{p^B}{q_v-1} + \log_v(p^B) \rceil + 1 \ < \ \frac{p^B}{q_v-1} + \log_v(p^B) + 2 \ .
\end{equation*} 
By our choice of $N_v$, we have $|\pi_v^{-N \beta_v/2}|_v \cdot R > 1$.  
Thus $D(0,1) \subset D(0,|\pi_v^{-N \beta_v/2}|_v R)$, and the 
inverse image of $D(0,1)$ under $\pi_v^{-N \beta_v/2} S_{p^B,\cO_v}(z)$ consists of $p^B$ disjoint
discs in $D(0,1)$, centered on the roots of $S_{p^B,\cO_v}(z)$ and having radii in $|K_v^{\times}|_v$.  
Since $f_{v,0}(z)$ induces an $F_{w_\ell(k)}$-rational scaled isometry from $B(\theta_k,\rho_k)$
\index{scaled isometry}
onto $D(0,1)$ for each $k = 1, \ldots, N_0$, it follows that
\begin{equation*}
f_v^{-1}(D(0,1)) \ = \ \bigcup_{k=1}^N \bigcup_{j=0}^{p^B-1} B(\theta_{k,j},\rho_{k,j}) 
\end{equation*}
where the balls on the right are pairwise disjoint and isometrically parametrizable. 
\index{isometrically parametrizable ball}
Furthermore, $\rho_{k,j}$ belongs to the value group of $F_{w_\ell(k)}$, and   
$f_v(z)$ induces an $F_{w_\ell(k)}$-rational scaled isometry from $B(\theta_{k,j},\rho_{k,j})$
\index{scaled isometry}
onto $D(0,1)$, for all $k, j$.  This establishes property (B.3).  

It is clear from the discussion above that the set 
\begin{equation} \label{FpBSimpleDec}
H_v \ = \ E_v \cap f_v^{-1}(D(0,1)) \ = \ 
\bigcup_{k=1}^N \bigcup_{j=0}^{p^B-1} \big( B(\theta_{k,j},\rho_{k,j}) \cap \cC_v(F_{w_\ell(k)}) \big)
\end{equation}
is $K_v$-simple, and that (\ref{FpBSimpleDec}) is a $K_v$-simple decomposition 
\index{$K_v$-simple!set}  
\index{$K_v$-simple!decomposition!compatible with another decomposition} 
compatible with the $K_v$-simple decomposition (\ref{E_vSimple}) of $E_v$.  
This yields property (B.4), and completes the proof.
\end{proof}
\index{initial approximating functions $f_v(z)$!construction when $K_v$ nonarchimedean!reduction to the case of a single ball|)}   

\section{ Generalized Stirling Polynomials} 

In this section we construct Stirling polynomials for sets of the form 
\index{Stirling polynomial!generalized}
\index{initial approximating functions $f_v(z)$!construction when $K_v$ nonarchimedean!generalized Stirling polynomials|(} 
\begin{equation*} 
H_w \ = \ \bigcup_{\ell=1}^d F_w \cap D(a_{\ell},r_{\ell}) \ \subset \ \CC_v \ .
\end{equation*} 
These will play a key role in the proof of Theorem \ref{CompactThm}.
The idea is that within any isometrically parametrizable ball
\index{isometrically parametrizable ball}
$B(a,r)$ disjoint from $\fX$, the canonical distance $[x,y]_{\fX,\vs}$
\index{canonical distance!$[z,w]_{\fX,\vs}$}
is a multiple of the spherical distance $\|x,y\|_v$\index{spherical metric}.  Under an isometric
parametrization, $\|x,y\|_v$ pulls back to a multiple of the usual distance 
$|X-Y|_v$ on the disc $D(0,r) \subset \CC_v$.  
This means that potential-theoretic constructions on $B(a,r)$ 
relative to $[x,y]_{\fX,\vs}$ are essentially the same as potential-theoretic
constructions on $\CC_v$ relative to $[X,Y]_{\infty} = |X-Y|_v$.  

\smallskip

First suppose $H_w = F_w \cap D(a,r)$ for a single disc, 
where $a \in F_w$ and $r \in |F_w^{\times}|_v$.  
Fix $b \in F_w$ with $|b|_v = r$.  We can can obtain a well-distributed sequence in $H_w$ 
by composing the basic well-distributed sequence $\{\psi_w(k)\}_{k \ge 0}$ for $\cO_w$
\index{basic well-distributed sequence in $\cO_w$} 
with the affine map $a + bz$:  
for each integer $n \ge 1$, we define the Stirling polynomial $S_{n,H_w}(z)$ by 
\index{Stirling polynomial!generalized|ii}
\begin{equation*} 
S_{n,H_w}(z) \ = \ \prod_{k=0}^{n-1} \big(z - (a + b \psi_w(k)) \big) 
             \ = \ b^n \cdot S_{n,\cO_w}(\frac{z-a}{b}) \ .
\end{equation*} 

Now consider the general case:  suppose
\begin{equation*} 
H_w \ = \ \bigcup_{\ell = 1}^{d} \big(F_w \cap D(a_{\ell},r_{\ell})\big)
\end{equation*} 
where  $a_{\ell} \in F_w$, $r_{\ell} \in |F_w^{\times}|_v$ for each $\ell$,
and the discs $D(a_{\ell},r_{\ell})$ are pairwise disjoint.   
Put $H_{w_\ell} = F_w \cap D(a_\ell,r_\ell)$, so  $H_w = \bigcup_{\ell=1}^{d} H_{w_\ell}$.
By Corollary \ref{BFCor1} the potential function\index{potential function} of $H_{w_\ell}$ is given by
\begin{equation} \label{FuPotF} 
u_\infty(z,H_{w_\ell}) \ = \ \left\{
   \begin{array}{ll}   \frac{1}{e_w(q_w - 1)} - \log_v(r_{\ell}) 
                          & \text{if $z \in H_{w_\ell}$\ ,} \\
                      -\log_v(|z-a_{\ell}|_v) & \text{if $z \notin D(a_{\ell},r_{\ell})$ \ , }
   \end{array} \right. 
\end{equation}
and in particular $V_{\infty}(H_{w_\ell}) = \frac{1}{e_w(q_w - 1)} - \log_v(r_{\ell})$. 
\c  

Let $\mu_{\infty}$ be the equilibrium distribution of $H_w$ relative
to the point $\infty$, and let $V = V_{\infty}(H_w)$ be the Robin constant.
\index{Robin constant!nonarchimedean!nonarchimedean $(\fX,\vs)$}   
For each $\ell$, put $w_{\ell} = \mu_{\infty}(H_{w_\ell}) > 0$.  
As in \S\ref{Chap6}.\ref{NonArchProofReductionsSection}, by Corollary \ref{BFCor2}
the following system of linear equations uniquely determine $V$ and the $w_{\ell}$\, :
\begin{equation} \label{FSystem10}
\left\{ \begin{array}{ccl}
        1  & = & 0 \cdot V + \sum_{\ell=1}^d w_{\ell} \ ,     \\   
0 & = & V + w_\ell \cdot \big(\log_v(r_{\ell}) - \frac{1}{e_w(q_w - 1)}\big) 
           + \displaystyle{\sum^{d}_{\substack{ j=1 \\ j \ne \ell }}} 
            w_j  \cdot \log_v (|a_\ell-a_j|_v)  \\
 &  & \qquad \qquad \text{for $\ell = 1, \ldots, d$.}  
          \end{array} \right.
\end{equation}
Since the coefficients of this system are rational, $V$ and the $w_{\ell}$ belong to $\QQ$.
Write $w_{\ell} = P_{\ell}/Q_{\ell}$ with positive integers $P_{\ell}, Q_{\ell}$,
and put $Q = \LCM(Q_1, \ldots, Q_{d})$.   

For each $\ell$, fix an element $b_{\ell} \in F_w$ with $|b_{\ell}|_v = r_{\ell}$.  
Then the affine map $\varphi_{\ell}(z) = a_{\ell}+b_{\ell}z$ takes $\cO_w$ to $H_{w_\ell}$.  

Let $n$ be a positive integer divisible by $Q$,
and put $n_{\ell} = w_{\ell} \cdot n$ for $\ell = 1, \ldots, d$. 
Define the Stirling polynomial $S_{n,H_w}(z)$ by 
\index{Stirling polynomial!generalized|ii}
\begin{equation} \label{GenStirling} 
S_{n,H_w}(z) \ = \ \prod_{\ell=1}^{d} \prod_{k=0}^{n_{\ell}-1} 
\big(z-(a_{\ell} + b_{\ell} \psi_w(k))\big) 
\ = \ \prod_{\ell=1}^{d} S_{n_{\ell},H_{w_\ell}}(z) \ .
\end{equation}
Note that $Q$ and $S_{n,H_w}(z)$ 
depend on the decomposition $H_w = \bigcup_{\ell = 1}^{d} H_{w_\ell}$  
and the maps $\varphi_{\ell}(z) = a_\ell + b_\ell z$,  not just  $H_w$.
For the rest of this section we will assume these are fixed.

\medskip

The following proposition generalizes Proposition \ref{BGProp1}:

\begin{proposition} \label{BGProp1AA}  Let $F_w/K_v$ be a finite, separable extension in $\CC_v$.
Suppose $H_w = \bigcup_{\ell=1}^{d} \big(F_w \cap D(a_{\ell},r_{\ell})\big)$, 
where the discs $D(a_{\ell},r_{\ell})$ are pairwise disjoint,
and $a_{\ell} \in F_w$ and $r_{\ell} \in |F_w^{\times}|_v$ for each $\ell$.  
Let $0 < w_1, \ldots, w_{d} \in \QQ$ be the weights corresponding to 
the sets $H_{w_\ell} = F_w \cap D(a_{\ell},r_{\ell})$ by the system $(\ref{FSystem10})$,
and let $Q$ be the least common multiple of their denominators.  

For each positive integer $n$ divisible by $Q$, 
let $S_{n,H_w}(z) = \prod_{\ell=1}^d S_{n_\ell,H_{w_\ell}}(z)$ 
be the Stirling polynomial of degree $n$ for $H_w$,
\index{Stirling polynomial!generalized}
and let $S_{n,H_w}^{\prime}(z)$ be its derivative. 
Write $\vartheta_1, \ldots, \vartheta_n$ for its zeros, 
and for each $\ell$ let $\varphi_{\ell}(z) = a_\ell + b_\ell z : \cO_w \rightarrow H_{w_\ell}$ 
be the affine map used in defining $S_{n,H_w}(z)$.
Then

$(A)$ $\vartheta_1, \ldots, \vartheta_n$ are distinct and belong to $H_w$.
There is a constant $A > 0$, independent of $n$, such that for all $i \ne j$
\begin{equation} \label{FThetaDist} 
|\vartheta_i - \vartheta_j|_v \ > \ A/n \ .
\end{equation} 

$(B)$ For each $k$, $1 \le k \le n$, if $\vartheta_k \in D(a_\ell,r_\ell)$ then we have 
\begin{equation} 
    \ord_v \big( S_{n,H_w}^{\prime}(\vartheta_k)\big) \ < \ n \cdot V_{\infty}(H_w) - \ord_v(b_\ell) \ .
  \label{BGF4AA}  
  \index{Robin constant!nonarchimedean}   
\end{equation} 

$(C)$  Given $x \in \CC_v$, fix $1 \le J \le n$ with 
$|x - \vartheta_J|_v = \min_k (|x-\vartheta_k|_v)$. 
If $\vartheta_J \in H_{w_\ell}$ then   
\begin{equation} \label{BGF5AA}
\ord_v(S_{n,H_w}(x)) \ < \ n \cdot V_{\infty}(H_w) + \ord_v\Big(\frac{x - \vartheta_J}{b_{\ell}} \Big) \ .
\end{equation}
\index{Robin constant!nonarchimedean}   
If $x \in \CC_v \backslash \bigcup_{\ell = 1}^{d} D(a_\ell,r_\ell)$, then 
\begin{equation} \label{BGF5BB}
\ord_v(S_{n,H_w}(x)) \ = \ n \cdot \big( V_{\infty}(H_w) - G(x,\infty;H_w) \big) \ . 
\end{equation}

\end{proposition}

\begin{proof}
We first prove the result when $H_w = F_w \cap D(a,r) = a + b \cO_w$,  
with $a, b \in F_w$ and $|b|_v = r > 0$.  For notational convenience, we relabel 
the zeros as $\vartheta_0, \ldots, \vartheta_{n-1}$, with $\vartheta_k = a + b \psi_w(k)$ for $k = 0, \ldots, n-1$.  

In this case, part (A) holds with $A = r$, since 
\begin{equation*} 
|\vartheta_i - \vartheta_j|_v \ = \ |b|_v |\psi_w(i)-\psi_w(j)|_v \ > \  r/n
\end{equation*} 
by Proposition \ref{BGProp1}(A) 

For part (B), note that since $S_{n,H_w}(z) = b^n S_{n,\cO_w}((z-a)/b)$, 
we have $S_{n,H_w}^{\prime}(\vartheta_k) = b^{n-1} S_{n,\cO_w}^{\prime}(\psi_w(k))$ 
for $k = 0, \ldots, n-1$.  Since 
\begin{equation*}
V_{\infty}(H_w) \ = \ \frac{1}{e_w(q_w-1)} - \log_v(r)  \ = \ V_{\infty}(\cO_w) + \ord_v(b) \ ,
\end{equation*} 
\index{Robin constant!nonarchimedean}   
part (B) follows from Proposition \ref{BGProp1}(B). 

For part (C), observe that   
if $x \in D(a,r)$ and $|x-\vartheta_J|_v = \min_{0 \le k < n}(|x-\vartheta_k|_v)$, 
then for $X := (x-a)/b \in D(0,1)$ we have $|(x-\vartheta_k)/b|_v = |X - \psi_w(J)|_v = \min_{0 \le k < n}(|X-\psi_w(k)|_v)$.
Hence (\ref{BGF5AA}) follows from Proposition \ref{BGProp1}(C).  

If $x \notin D(a,r)$, 
then by the ultrametric inequality we have  $|x - \vartheta_k|_v = |x - a|_v$ for each $k$.  
Furthermore, $u_{\infty}(x,H_w) = -\log_v(|x-a|_v)$ by (\ref{FuPotF}). Thus
\index{Robin constant!nonarchimedean}   
\begin{equation*}
\ord_v(S_{n,H_w}(x)) \ = \ n \cdot u_\infty(x,H_w) 
\ = \ n \cdot \big(V_\infty(H_w) - G(x,\infty;H_w)\big) \ .
\end{equation*} 
This yields (\ref{BGF5BB}) since $G(x,\infty;H_w) > 0$.

\smallskip 
Now consider the general case, 
where $H_w = \bigcup_{\ell=1}^d H_{w_\ell} = \bigcup_{\ell=1}^d \big( F_w \cap D(a_{\ell},r_{\ell}) \big)$,
with finitely many pairwise disjoint discs $D(a_{\ell},r_{\ell})$ such that $a_{\ell} \in F_w$ and $r_{\ell} \in |F_w^{\times}|_v$.  

Fix $0 < n \in \NN$ divisible by $Q$, put $n_\ell = n \cdot w_\ell$ for each $\ell$, 
and consider $S_{n,H_w}(z) = \prod_{k=1}^n (z-\vartheta_k) = \prod_{\ell=1}^d S_{n_\ell,H_{w_\ell}}(z)$. 

For part (A), put $A = \min_{1 \le \ell \le d} (r_\ell)$.  Fix $n$ divisible by $Q$, 
and let $\vartheta_i \ne \vartheta_j$ be distinct roots of $S_{n,H_w}(z)$.  Necessarily $n \ge 2$.  
Since the balls $B(a_{\ell},r_{\ell})$ are pairwise disjoint, if $\vartheta_i$ and $\vartheta_j$ belong 
to distinct balls then $|\vartheta_i - \vartheta_j|_v > A > A/n$.  On the other hand, 
if $\vartheta_i$ and $\vartheta_j$ belong to the same ball $B(a_\ell,r_\ell)$,
then there are indices $k \ne h$ with $1 \le k, h \le n_\ell$ such that 
$\vartheta_i = a_\ell + b_\ell \psi_w(k)$ and $\vartheta_j = a_\ell + b_\ell \psi_w(h)$.  
In this case
\begin{equation*} 
|\vartheta_i-\vartheta_j|_v \ \ge \ |b_\ell|_v \cdot |\psi_w(k) - \psi_w(h)|_v 
\ > \ A / n_\ell \ \ge \ A/n \ .
\end{equation*}   

For part (B), first note that if $x \in D(a_{\ell},r_{\ell})$, 
then for each $j \ne \ell$ and each $\vartheta_k \in D(a_j,r_j)$, we have
$|x-\vartheta_k|_v = |a_\ell - a_j|_v$.  Hence 
\begin{equation} \label{FNqp1}
|S_{n,H_w}(x)|_v \ = \ |S_{n_\ell,H_{w_\ell}}(x)|_v \cdot 
            \prod_{\substack{j = 1 \\ j \ne \ell}}^d |a_{\ell} - a_j|_v^{n_j} \ .
\end{equation} 
Similarly, if $\vartheta_h \in D(a_{\ell},r_{\ell})$ then 
\begin{equation} \label{FNqp2}
|S_{n,H_w}^{\prime}(\vartheta_h)|_v \ = \ |S_{n_\ell,H_{w_\ell}}^{\prime}(\vartheta_h)|_v \cdot 
            \prod_{\substack{j = 1 \\ j \ne \ell}}^d |a_{\ell} - a_j|_v^{n_j} \ .
\end{equation} 
On the other hand, since $n_j = n \cdot w_j$, 
for each $\ell$ the equations (\ref{FSystem10}) give 
\index{Robin constant!nonarchimedean}   
\begin{equation} \label{FNqp3}
n \cdot V_\infty(H_w) \ = \ n_{\ell} V_{\infty}(H_{w_\ell}) 
        + \sum_{\substack{j=1 \\ j \ne \ell}}^d n_j \cdot \big(-\log_v(|a_\ell-a_j|_v)\big) \ .  
\end{equation} 
One obtains part (B) by combining (\ref{BGF4AA}) for $H_{w_\ell}$ with (\ref{FNqp2}) and (\ref{FNqp3}), 
and using $\log_w(n_{\ell}) \le \log_w(n)$.    

For part (C), if $\vartheta_J \in D(a_{\ell},r_{\ell})$ one obtains (\ref{BGF5AA})  
by combining  (\ref{BGF5AA}) for $H_{w_\ell}$ with (\ref{FNqp1}) and (\ref{FNqp3}).
If $x \in \CC_v \backslash \bigcup_{\ell=1}^d D(a_{\ell},r_{\ell})$, 
one obtains (\ref{BGF5BB}) by using (\ref{FuPotF}) and noting that  
\index{Robin constant!nonarchimedean}   
\begin{eqnarray*} 
-\log_v(|S_{n,H_w}(x)|_v) & = & - \sum_{\ell=1}^d n_{\ell} \log_v(|x-a_{\ell}|_v) 
   \ = \ n \cdot \Big( \sum_{\ell=1}^d w_{\ell} \, u_\infty(x,H_{w_\ell})\Big) \\ 
 & = & n \cdot u_\infty(x,H_w) \ = \ n \cdot \big(V_{\infty}(H_w) - G(x,\infty;H_w) \big) \ .
\end{eqnarray*} 
\end{proof}

\noindent{\bf Remark.} In Proposition \ref{BGProp1AA}(B), 
one can show that if $x \in \CC_v \backslash \bigcup_{\ell=1}^d D(a_\ell,r_\ell)$, then 
\index{Robin constant!nonarchimedean}   
\begin{equation*}
\ord_v(S_{n,H_w}(x)) \ < \  n \cdot \Big( V_{\infty}(H_w) - 
\min_{1 \le \ell \le d}(w_\ell) \cdot \frac{1}{e_w(q_w-1)} \Big) \ . 
\end{equation*}
One can also prove the following generalization of Corollary \ref{SnOwMappingCor}:
given an $R$ satisfying
\begin{equation} \label{FRBound}
0 \ < \ R \ \le \ q_v^{-n V_{\infty}(H_w)} \cdot n^{-1/[F_w:K_v]} \ , 
\end{equation} 
put $\rho_k = R/|S_{n,H_w}^{\prime}(\vartheta_k)|_v$ for $k = 1, \ldots, n$.
Then the discs $D(\vartheta_k,\rho_k)$ are pairwise disjoint and 
\begin{equation} \label{FInvFAA}
S_{n,H_w}^{-1}\big(D(0,R)\big) \ = \ \bigcup_{k=1}^{n} D(\vartheta_k,\rho_k)  
\ \subseteq \ \bigcup_{\ell = 1}^d D(a_{\ell},r_{\ell}) \ .
\end{equation}   

\index{initial approximating functions $f_v(z)$!construction when $K_v$ nonarchimedean!generalized Stirling polynomials|)}

\section{ Proof of Proposition $\ref{K_vReductionProp}$} 
\label{NonArchProofSection} 

In this section we will prove Proposition \ref{K_vReductionProp}, 
completing the proof of Theorem \ref{CompactThm}.  
The proof breaks into two cases, according as the genus $g(\cC_v) = 0$
or $g(\cC_v) > 0$.  
 
\medskip

Proposition \ref{K_vReductionProp}
concerns each subset $E_{v,\ell} = \cC_v(F_{w_\ell}) \cap B(a_\ell,r_\ell)$ individually. 
To simplify notation, we restate the Proposition, dropping the index $\ell$ and 
relabeling $E_{v,\ell}$ as $E_w$\ :

\vskip .15 in

Proposition \ref{K_vReductionProp}{\rm A}. \label{K_vReductionPropA}
{\em
Let $F_w/K_v$ be a finite, separable extension in $\CC_v$
and take $E_w = \cC_v(F_w) \cap B(a,r)$, where $a \in \cC_v(F_w)$ 
and $B(a,r) \subset \cC_v(\CC_v) \backslash \fX$ is an isometrically parametrizable ball.  
\index{isometrically parametrizable ball} 
Let $\varepsilon_w > 0$ be given.  Then there is a compact subset $\tE_w \subseteq E_w$ for which  

$(A)$ There are points $\alpha_j \in \cC_v(F_w) \cap B(a,r)$ 
and pairwise disjoint isometrically parametrizable balls
\index{isometrically parametrizable ball} 
$B(\alpha_1,r_1), \ldots, B(\alpha_d,r_d) \subseteq B(a,r)$,
such that $\tE_w$ has the form  
\begin{equation} \label{FtESetForm}
\tE_w \ = \ \bigcup_{j=1}^{d} \big(\cC_v(F_w) \cap B(\alpha_j,r_j)\big) 
\end{equation} 
and for each $x_i \in \fX$,  
\begin{equation} \label{FRRApprox4B1A}
|V_{x_i}(\tE_w) - V_{x_i}(E_w)| \ < \ \varepsilon_w\ .
\end{equation} 
   
$(B)$  For each $0 < \beta_w \in \QQ$  
and each $F_w$-symmetric probability vector
$\vs = {}^t(s_1, \ldots, s_m)$ with rational entries, 
there is an integer $N_w \ge 1$  such that 
for each positive integer  $N$ divisible by $N_w$, 
there is an $(\fX,\vs)$-function\index{$(\fX,\vs)$-function} 
$f_w \in F_w(\cC_v)$ of degree $N$ such that
   
\quad $(1)$  For all $z \in \cC_v(\CC_w) \backslash \big( B(a,r) \cup \fX \big)$, 
\begin{equation} \label{FRRApprox7B1A} 
\frac{1}{N} \log_v(|f_w(z)|_v) \ = \ G_{\fX,\vs}(z,\tE_w) + \beta_w \ .
\end{equation}   
   
\quad $(2)$ The zeros $\theta_1, \ldots, \theta_N$ of $f_w$ 
are distinct and belong to $E_w$ $($hence $\cC_v(F_w))$.

\quad $(3)$  
$f_w^{-1}(D(0,1)) = \bigcup_{h=1}^N B(\theta_h,\rho_h)$, 
where the balls $B(\theta_h,\rho_h)$ are pairwise disjoint and contained in $B(a,r)$.  
} 
\vskip .15 in

The following lemma will be helpful in proving Proposition \ref{K_vReductionProp}A:
\index{initial approximating functions $f_v(z)$!construction when $K_v$ nonarchimedean!reduction to finding a principle divisor|(} 

\begin{lemma} \label{ConstructionLemma}
Let $E_w = \cC_v(F_w) \cap B(a,r)$ be as Proposition $\ref{K_vReductionProp}{\rm A}$.  
Given $\varepsilon_w > 0$, suppose $\tE_w \subseteq E_w$ is a compact subset satisfying part $(A)$ 
of Proposition $\ref{K_vReductionProp}{\rm A}$. 
 
Then part $(B)$ of Proposition $\ref{K_vReductionProp}{\rm A}$ holds for $E_w$ and $\tE_w$
if for each $F_w$-symmetric 
$\vs \in \cP^m(\QQ)$, there are an integer $N_0 = N_0(\vs,\tE_w) > 0$, 
and constants $A = A(\vs,\tE_w) > 0$ and $B = B(\vs,\tE_w)$, such that 
$N_0 \vs \in \NN^m$ and for each sufficiently large 
integer $N$ divisible by $N_0$ there is a divisor $D_N$ of the form 
\begin{equation} 
D_N \ = \ \sum_{h=1}^N (\theta_h) - \sum_{i=1}^m Ns_i(x_i) 
\end{equation}
satisfying conditions $(1)$, $(2)$ and $(3)$ below:   

\vskip .05 in
$(1)$ $D_N$ is principal.

$(2)$ $\theta_1, \ldots, \theta_N$ are distinct and belong to $E_w$, 
and for all  $i \ne j$  
\begin{equation*}
\|\theta_i,\theta_j\|_v \ > \ A/N\ . 
\end{equation*}

$(3)$ The pseudopolynomial $Q_N(z) = \prod_{h=1}^N [z,\theta_h]_{\fX,\vs}$ 
\index{pseudopolynomial!$(\fX,\vs)$} 
with divisor $D_N$ has the following property$:$  for each $z \in B(a,r)$, 
if $\theta_J$ is such that $\|z,\theta_J\|_v = \min_{1 \le h \le N} \|z,\theta_h\|_v$, then 
\index{Robin constant!nonarchimedean!nonarchimedean $(\fX,\vs)$}   
\begin{equation*}
-\log_v(Q_N(z)) \ \le \ N \cdot V_{\fX,\vs}(\tE_w)  
                                -\log_v(\|z,\theta_J\|_v) + B \ . 
\end{equation*}
%
%
\end{lemma} 

\begin{proof}  Fix an $F_w$-symmetric probability vector $\vs \in \cP^m(\QQ)$
and a number $0 < \beta_w \in \QQ$.  Let $A > 0$, $B \ge 0$, and $0 < N_0 \in \NN$
be the numbers given by Lemma \ref{ConstructionLemma} for $E_w$, $\tE_w$, and $\vs$. 
As in \S\ref{Chap6}.\ref{NonArchProofReductionsSection},
there is a constant $C_{\fX,\vs} \in |\CC_v^{\times}|_v$ 
such that $[x,y]_{\fX,\vs} = C_{\fX,\vs} \|x,y\|_v$ for all $x, y \in B(a,r)$.  
Let $N_1$ be the least positive integer such that $C_{\fX,\vs}^{N_1} \in |F_w^{\times}|_v$.
By (\ref{FtESetForm}) and Corollary \ref{BFCor2}, $V_{\fX,\vs}(\tE_w) \in \QQ$, 
and by hypothesis, $\beta_w \in \QQ$;\index{Robin constant!nonarchimedean!nonarchimedean $(\fX,\vs)$}
let $N_2$ be the least common denominator for $V_{\fX,\vs}(\tE_w)$ and $\beta_w$.

Let $N_3$ be the smallest natural number such that for each $N \ge N_3$ divisible by $N_0$, 
there is a divisor $D_N = \sum (\theta_h) - \sum Ns_i(x_i)$ 
satisfying conditions (1), (2), and (3) of Lemma \ref{ConstructionLemma}. 
Let $N_4$ be the smallest natural number such that for each $N \ge N_4$ we have 
\begin{equation} \label{FNBigEnough} 
N \beta_w - B + \log_v(A) - \log_v(N) \ > \ 0 \ .
\end{equation}  
The number $N_w$ in part (B) of Proposition \ref{K_vReductionProp}A 
will be the least multiple of $N_0$, $N_1$, and $N_2$ which is greater than $N_3$ and $N_4$.  

Given a positive integer $N$ divisible by $N_w$, let $D_N$ be the corresponding divisor. 
Since $D_N$ is $F_w$-rational and principal,
there is a function $f_N \in F_w(\cC_v)$ with $\div(f_N) = D_N$.  By the factorization property
of the canonical distance, 
\index{canonical distance!$[z,w]_{\fX,\vs}$}\index{canonical distance!factorization property}
there is a constant $C$ such that 
$|f_N(z)|_v = C \prod_{h=1}^N [z,\theta_h]_{\fX,\vs}$ for all $z \in \cC_v(\CC_v)$.  

Fix a point $z_0 \in E_w \backslash \{\theta_1, \ldots, \theta_N\}$.  
Since $f_N$ is $F_w$-rational and $z_0 \in \cC_v(F_w)$, we have $f(z_0) \in F_w^{\times}$.  
Likewise, for each $h$ we have $\|z_0,\theta_h\|_v \in |F_w^{\times}|_v$, so 
\begin{equation*}
\prod_{h=1}^N [z_0,\theta_h]_{\fX,\vs} 
\ = \ C_{\fX,\vs}^N \cdot \prod_{h=1}^N \|z_0,\theta_h\|_v \ \in \ |F_w^{\times}|_v \ .
\end{equation*}
Since $|f(z_0)|_v = C \prod_{h=1}^N [z_0,\theta_h]_{\fX,\vs}$ 
it follows that $C \in |F_w^{\times}|_v$. Thus, after scaling $f_N$ by a suitable constant, 
we can assume that $C = 1$, and that $|f_N(z)|_v = \prod_{h=1}^N [z_0,\theta_h]_{\fX,\vs}$. 

By construction, $N \cdot (V_{\fX,\vs}(\tE_w) + \beta_w) \in \ZZ$. 
Let $C_N \in F_w^{\times}$ be such that\index{Robin constant!nonarchimedean!nonarchimedean $(\fX,\vs)$}
\begin{equation} \label{FC_NSize} 
|C_N|_v \ = \ q_v^{N V_{\fX,\vs}(\tE_w) + N \beta_w} \ ,
\end{equation}
and define the function $f_w$ in part (B) of Proposition \ref{K_vReductionProp}A 
to be $f_w = C_N \cdot f_N$. 
By assumption the zeros of $f_w$ are distinct and belong to $E_w$.
If $A$ is the constant given Proposition \ref{BGProp1AA}(A) for $H_w$, 
then for all pairs of distinct roots $\theta_i \ne \theta_j$ we have 
\begin{equation*}
\|\theta_i, \theta_j\|_v \ = \ |\vartheta_i - \vartheta_j|_v \ > \ A/N \ .
\end{equation*}
Thus property (B2) in Proposition \ref{K_vReductionProp}A holds.  

We next show that property (B3) holds.
Without loss, we can assume $A$ is small enough that $B(\theta_h,A/N) \subseteq B(a,r)$
for each $h$.  Thus $B(\theta_h,A/N)$ is isometrically parametrizable.  Fix an $F_w$-rational 
\index{isometrically parametrizable ball}
isometric parametrization $\varphi_h : D(0,A/N) \rightarrow B(\theta_h,A/N)$ 
\index{isometric parametrization}  
with $\varphi_h(0) = \theta_h$, and expand 
\begin{equation*}
f_w(\varphi_h(z)) \ = \ \sum_{k=1}^{\infty} c_{h,k} z^k  \ \in \ F_w[[z]] \ .
\end{equation*} 
The zeros of $f_w$ are distinct, so $c_{h,1} \ne 0$.  
Hence if $|z|_v$ is sufficiently small, then 
\begin{equation*} 
-\log_v(|f_w(\varphi_n(z))|_v) \ = \ -\log_v(|c_{h,1}|_v) - \log_v(|z|_v)\ .
\end{equation*}   
On the other hand, by condition (3)
\begin{eqnarray*}
-\log_v(|f_w(\varphi_n(z))|_v) & \le & (B + N V_{\fX,\vs}(\tE_w) - \log_v(|z|_v)) - (N V_{\fX,\vs}(\tE_w) + N \beta_w) \\
& = & B - N \beta_w - \log_v(|z|_v) \ ,  
\end{eqnarray*} 
which means that $|c_{h,1}|_v \ \ge \ q_v^{-B + N \beta_w}$ \ .
\index{Robin constant!nonarchimedean!nonarchimedean $(\fX,\vs)$}

By condition (2) of Lemma \ref{ConstructionLemma}, 
$\theta_h$ is the only zero of $f_w$ in $B(\theta_h,A/N)$.  Thus Proposition \ref{ScaledIsometryProp}
shows that $f_w$ induces an $F_w$-rational scaled isometry from 
\index{scaled isometry}
$B(\theta_h,A/N)$ onto $D(0,|c_{h,1}|_v A/N)$.  Since  (\ref{FNBigEnough}) holds, we have 
\begin{equation*} 
|c_{h,1}|_v A/N \ \ge \ q_v^{-B + N \beta_w + \log_v(A) - \log_v(N)} \ > \ 1 \ . 
\end{equation*}
Put $\rho_h = 1/|c_{h,1}|_v < A/N$;  then $f_w$ induces a scaled isometry from $B(\theta_h,\rho_h)$ onto $D(0,1)$, 
\index{scaled isometry}
and $\rho_h \in |F_w^{\times}|_v$. 
 
Again by condition (2), the balls $B(\theta_h,A/N)$ for $h = 1, \ldots, N$ are pairwise disjoint. 
Since $B(\theta_h,\rho_h) \subset B(\theta_h,A/N)$, the balls $B(\theta_h,\rho_h)$
are pairwise disjoint.  Since $f_w$ is a rational function of degree $N$, the $N$ balls $B(\theta_h,\rho_h)$ 
account for all the solutions to $f(z) = x$ with $x \in D(0,1)$, and it follows that 
\begin{equation*} 
f_w^{-1}(D(0,1)) \ = \ \bigcup_{h=1}^N B(\theta_h,\rho_h) \ \subset \ B(a,r) \ .
\end{equation*}

It remains to establish property (B1).  Fix $z \in \cC_v(\CC_v) \backslash B(a,r)$.
Since the canonical distance is constant on pairwise disjoint
\index{canonical distance!$[z,w]_{\fX,\vs}$}
\index{canonical distance!nonarchimedean!constant on disjoint balls}
isometrically parametrizable balls, we have $[z,\theta_h]_{\fX,\vs} = [z,a]_{\fX,\vs}$ for each $h$, 
\index{isometrically parametrizable ball}
and so 
\begin{equation*} 
\frac{1}{N}\log_v(|f_N(z)|_v) \ = \ \frac{1}{N} \sum_{h=1}^N \log_v([z,\theta_h]_{\fX,\vs}) 
\ = \  \log_v([z,a]_{\fX,\vs}) \ .
\end{equation*}
However, by (\ref{FSecond}), $u_{\fX,\vs}(z,\tE_w) = -\log_v([z,a]_{\fX,\vs})$.  
Since $f_w = C_N \cdot f_N$ and $G_{\fX,\vs}(z,\tE_w) = V_{\fX,\vs}(\tE_w) - u_{\fX,\vs}(z,\tE_w)$,
\index{Robin constant!nonarchimedean!nonarchimedean $(\fX,\vs)$}
it follows from (\ref{FC_NSize}) that
\begin{equation*} 
\frac{1}{N} \log_v(|f_w(z)|_v) \ = \ (V_{\fX,\vs}(\tE_w) + \beta_w) - u_{\fX,\vs}(z,\tE_w) 
\ = \ G_{\fX,\vs}(z,\tE_w) + \beta_w \ . 
\end{equation*} 
\vskip -.15 in
\end{proof}\index{initial approximating functions $f_v(z)$!construction when $K_v$ nonarchimedean!reduction to finding a principle divisor|)}

\begin{proof}[Proof of Proposition \ref{K_vReductionProp}{\rm A} when $g(\cC_v) = 0$.]
In this case, the proof is relatively easy.
We can take $\tE_w = E_w = \cC_v(F_w) \cap B(a,r)$, so (A) holds trivially. 
\index{initial approximating functions $f_v(z)$!construction when $K_v$ nonarchimedean!the proof when $g(\cC_v) = 0$|(}

For (B), let an $F_w$-symmetric probability vector $\vs \in \cP^m(\QQ)$ and a number $0 < \beta_w \in \QQ$ be given.
There is a constant $C_{\fX,\vs} \in |\CC_v^{\times}|_v$ such that $[z,w]_{\fX,\vs} = C_{\fX,\vs}\|z,w\|_v$ 
for all $z,w \in B(a,r)$.  
Fix an $F_w$-rational isometric parametrization $\varphi : D(0,r) \rightarrow B(a,r)$
\index{isometric parametrization}  
and let $H_w = F_w \cap D(0,r)$.  Then $\varphi(H_w) = E_w = \tE_w$.  

Let $A > 0$ be the constant given in part (A) 
of Proposition \ref{BGProp1AA} for $H_w$.
The idea for constructing the functions $f_w$ 
is to push forward the zeros of Stirling polynomials for $H_w$,
\index{Stirling polynomial!generalized}
and let them be the zeros of $f_w$.  
Let $0 < Q \in \ZZ$ be the number given for $H_w$ for Proposition \ref{BGProp1AA}.
Write $\beta_w = S/R$ with coprime integers $S, R$, where $R > 0$,  
and put $V = V_\infty(H_w)$. By Corollary \ref{BFCor2}, $V \in \QQ$;  
\index{Robin constant!nonarchimedean}
write $V = X/Y$ with coprime integers $X, Y$, where $Y > 0$.
Finally, for each $i = 1, \ldots, m$ write $s_i = A_i/B_i$ with coprime integers $A_i,B_i$, where $B_i > 0$,       
and set $N_0 = Q R Y \cdot \LCM(B_1, \ldots, B_m)$.

Suppose $0 < N \in \NN$ is a multiple of $N_0$.
Let $S_{N,H_w}(z) = \prod_{k=1}^N (z-\vartheta_k)$ 
be the Stirling polynomial of degree $N$ for $H_w$ constructed 
\index{Stirling polynomial!generalized}
in Proposition \ref{BGProp1AA}, and put $\theta_k = \varphi(\vartheta_k) \in E_w$ for $k = 1, \ldots, N$.  
Then  
\begin{equation*} 
D_N \ := \ \sum_{k=0}^{N-1} (\theta_k) - \sum_{i=1}^m N s_i(x_i) 
\end{equation*} 
is an $F_w$-rational divisor of degree $0$.  Since $\cC_v(F_w)$ is nonempty, 
$\cC_w = \cC_v \times_{K_v} \Spec(F_w)$ is $F_w$-isomorphic to $\PP^1/F_w$,  
and each divisor of degree $0$ is principal.  
Thus condition (1) in Lemma \ref{ConstructionLemma} holds.

By construction $\theta_1, \ldots, \theta_N$ are distinct and belong to $E_w$.
Since $\varphi : D(0,r) \rightarrow B(a,r)$ is an $F_w$-rational isometric parametrization,
\index{isometric parametrization}  
for all $i \ne j$ we have $\|\theta_i , \theta_j\|_v = |\vartheta_i-\vartheta_j|_v > A/N$.
Thus condition (2) in Lemma \ref{ConstructionLemma} holds.

To show condition (3), 
let $Q_N(z) = \prod [z,\theta_k]_{\fX,\vs}$ be the $(\fX,\vs)$-pseudo-polynomial associated with $D_N$.
Note that for $x, y \in D(0,r)$ we have $\|\varphi(x),\varphi(y)\|_v = |x-y|_v$, 
so $[\varphi(x),\varphi(y)]_{\fX,\vs} = C_{\fX,\vs}|x-y|_v$. Thus for $x \in D(0,r)$ 
\begin{equation*} 
Q_N(\varphi(x)) \ = \ \prod_{k=1}^N [\varphi(x),\theta_k]_{\fX,\vs} 
\ = \ C_{\fX,\vs}^N \prod_{k=1}^N |x-\vartheta_k|_v
\ = \ C_{\fX,\vs}^N \, |S_{N,H_w}(x)|_v \ ,
\end{equation*} 
and so
\begin{equation} \label{FQNNN1} 
- \log_v(Q_N(\varphi(x))) \ = \ -N \log_v(C_{\fX,\vs}) - \log_v(|S_{N,H_w}(x)|_v) \ .
\end{equation} 
Similarly
\begin{eqnarray}
V_{\fX,\vs}(\tE_w) & = & \inf_{\text{prob meas $\nu$ on $\tE_w$}} 
                   \iint_{\tE_w \times \tE_w} -\log_v([z,w]_{\fX,\vs}) \, d\nu(z) d\nu(w) \notag \\
                   & = & \inf_{\text{prob meas $\nu$ on $H_w$}} 
                   \iint_{H_w \times H_w} -\log_v(C_{\fX,\vs}|x-y|_v) \, d\nu(x) d\nu(y) \notag \\
                   & = & V_{\infty}(H_w) - \log_v(C_{\fX,\vs}) \ . \label{FQNNN2}
\end{eqnarray} 
\index{Robin constant!nonarchimedean!nonarchimedean $(\fX,\vs)$}
Let $B = \log_v(r)$.  Given $x \in D(0,r)$, let $\vartheta_J$ be the root of $S_N(z,H_w)$
for which $|x - \vartheta_J|_v$ is minimal.  
Using (\ref{FQNNN1}), (\ref{FQNNN2}) and Proposition \ref{BGProp1AA}(C), we obtain 
\begin{eqnarray*}
-\log_v(Q_N(\varphi(x))) 
& = & -N \log_v(C_{\fX,\vs}) - \log_v(|S_{N,H_w}(x)|_v) \\
& \le & - N \log_v(C_{\fX,\vs}) + N \cdot V_{\infty}(H_w) - \log_v(|x-\vartheta_J|_v) + B \\
& = & \ N \cdot V_{\fX,\vs}(\tE_w) - \log_v(\|\varphi(x),\theta_J\|_v) + B \ .
\end{eqnarray*}
\index{Robin constant!nonarchimedean}
This yields condition (3) in Lemma \ref{ConstructionLemma}.

Applying Lemma \ref{ConstructionLemma}, we obtain part (B) of Proposition \ref{K_vReductionProp}A.
\end{proof}\index{initial approximating functions $f_v(z)$!construction when $K_v$ nonarchimedean!the proof when $g(\cC_v) = 0$|)}

For the remainder of this section, we will assume that $g = g(\cC_v) > 0$.  
To prove Theorem \ref{CompactThm} when $g > 0$, 
we must first do some preparations. Given $\varepsilon_w > 0$, 
we will construct a subset $\tE_w \subset E_w$ of the form 
\begin{equation*} 
\tE_w \ = \  E_w \backslash \big( \bigcup_{j=1}^g B(\alpha_j,\rho_j) \big) 
\cup \big( \bigcup_{j=1}^g \cC_v(F_w) \cap B(\alpha_j,\trho_j) \big) \ , 
\end{equation*} 
with $0 < \trho_j < \rho_j$ for each $j$, 
such that $|V_{x_i}(E_w) - V_{x_i}(\tE_w)| < \varepsilon_w$ for each $x_i \in \fX$. 
\index{Robin constant!nonarchimedean}
That is, we first remove finitely many discs from $E_w$, 
and replace them with the $F_w$-rational points in smaller discs having the same centers. 
The idea is that given a suitable divisor $D_N^* = \sum_{k=1}^N (\theta_k^*) - \sum_{i=1}^N Ns_i(x_i)$ 
with $\theta_1^*, \ldots, \theta_N^* \in \tE_w$, we will be able to create a principal divisor 
$D_N = \sum_{k=1}^N (\theta_k) - \sum_{i=1}^N Ns_i(x_i)$ with $\theta_1, \ldots, \theta_N \in E_w$ 
by moving some of the $\theta_k^*$ into $E_w \backslash \tE_w$.    

The construction of $\tE_w$ is based on the following two facts, 
proved in the Appendices.
First, removing small balls from $E_w$ does not significantly change its capacity:  
\index{capacity} 

\begin{proposition} \label{BPropF3A1}  
Let $E_w \subset \cC_v(\CC_v) \backslash \fX$ be as in Proposition $\ref{K_vReductionProp}{\rm A}$.  
Fix  $\alpha_1, \ldots, \alpha_g \in E_w$.  
Then for each $\varepsilon_w > 0$, there is an $R_1 > 0$ such that for any compact set $\tE_w$ such that
\index{Robin constant!nonarchimedean}
\begin{equation*}
E_w \backslash \big(\bigcup_{j=1}^g B(\alpha_j,R_1)\big) \ \subseteq \ \tE_w \ \subseteq \ E_w \ ,
\end{equation*}
we have $|V_{x_i}(E_w) - V_{x_i}(\tE_w)| < \varepsilon_w$ for each $x_i \in \fX$.  
\end{proposition} 

\begin{proof}  
This is a special case of Proposition \ref{BPropF3} of Appendix \ref{AppA}.  
\end{proof}

Second, for each generic, sufficiently small 
polyball $\prod_{j=1}^g B(\alpha_j,\rho) \subset \cC_v(\CC_v)^g$, 
there is an action of a neighborhood of the origin in $\Jac(\cC_v)(\CC_v)$ 
which makes the polyball into a principal homogeneous space.\index{principal homogeneous space}
This action will enable us to move the points $\theta_k^*$ and obtain the principal divisor $D_N$.
\index{initial approximating functions $f_v(z)$!construction when $K_v$ nonarchimedean!the local action theorem|(}

Let $\JacNer(\cC_v)/\Spec(\cO_v)$ be the N\'eron model of $\Jac(\cC_v)$.
\index{N\'eron model!of Jacobian}  
By (\cite{BLR}, Theorem 1, p.153), $\JacNer(\cC_v)$ is quasi-projective.   
We regard it as embedded in $\PP^N/\Spec(\cO_v)$, for an appropriate $N$,  
and identify $\Jac(\cC_v)$ with generic fibre of $\JacNer(\cC_v)$ in $\PP^N/\Spec(K_v)$.

Let $\|x,y\|_{J,v}$ be the corresponding spherical metric\index{spherical metric!on Jacobian} on $\Jac(\cC_v)(\CC_v)$,
and let $O$ be the origin in $\Jac(\cC_v)(\CC_v)$.  
Then the unit ball $B_J(O,1)^- := \{z \in \Jac(\cC_v)(\CC_v) : \|z,O\|_{J,v} < 1\}$ is a subgroup.
Since $O$ is nonsingular on the special fibre of the N\'eron model,\index{N\'eron model!of Jacobian} 
there is a $K_v$-rational isometric parametrization
\index{isometric parametrization}   
\begin{equation} \label{JPsiDef1}
\Psi : D(\vORIG,1)^- \ \rightarrow \ B_J(O,1)^- 
\end{equation}
by power series converging on $D(\vORIG,1)^- \subset \CC_v^g$, 
taking $\vORIG$ to $O$ (Theorem \ref{IsoParamThm}).
Pulling the group action back to $D(\vORIG,1)^-$ using $\Psi$ 
yields the formal group of $\Jac(\cC_v)$. 
\index{formal group}

Put $\cCbar_v = \cC_v \times_{K_v} \Spec(\CC_v)$,  
\label{`SymbolIndexcCbarv'}
and let $\BPic_{\cCbar_v/\CC_v}$ be its Picard scheme.
\index{Picard scheme}
Given a divisor $D$ on $\cC_v(\CC_v)$, let $[D]$ be its 
class in $\BPic_{\cCbar_v/\CC_v}(\CC_v)$.
The identity component $\BPic_{\cCbar_v/\CC_v}^0$ 
is canonically isomorphic to $\Jac(\cC_v) \times_{K_v} \Spec(\CC_v)$, 
and we identify $\BPic_{\cCbar_v/\CC_v}^0(\CC_v)$ with $\Jac(\cC_v)(\CC_v)$. 
Given $\alpha \in \cC_v(\CC_v)$, 
the Abel map\index{Abel map|ii} $\BFj_{\alpha} : \cC_v(\CC_v) \rightarrow \Jac(\cC_v)(\CC_v)$ is defined by 
\begin{equation*}
\BFj_{\alpha}(x) \ = \ [(x)-(\alpha)] \ .
\end{equation*}
The Abel map is continuous for the $v$-topology, 
and if $\alpha \in \cC_v(F_w)$ for some finite extension $F_w/K_v$,
it is $F_w$-rational.
Given $\valpha = (\alpha_1, \ldots, \alpha_g) \in \cC_v(\CC_v)^g$, 
let $\crJ_{\valpha} : \cC_v(\CC_v)^g \rightarrow \Jac(\cC_v)(\CC_v)$
be the map 
\begin{equation*} 
\crJ_{\valpha}(\vx) \ = \ [\sum_{j=1}^g (x_j) - \sum_{j=1}^g (\alpha_j)] 
   \ = \ \sum_{j=1}^g \BFj_{\alpha_j}(x_j) \ .
\end{equation*}  
If $\alpha_1, \ldots, \alpha_g \in \cC_v(F_w)$, then $\crJ_{\valpha}$ is $F_w$-rational.
It is shown in Appendix \ref{AppD} that $\crJ_{\valpha}$ is nonsingular at $\valpha$ 
for a dense set of $\valpha$, and if $\crJ_{\valpha}$ is nonsingular at $\valpha$,  
then for each sufficiently small $\rho > 0$, the image of 
 $\prod_{j=1}^g B(\alpha_j,\rho)$ under $\crJ_{\valpha}$ is an open subgroup 
$W_{\valpha}(\rho)$ of $\Jac(\cC_v)(\CC_v)$. Furthermore, 
if $\Jp$ is addition in $\BPic_{\cCbar_v/\CC_v}(\CC_v)$, 
then $\prod_{j=1}^g B(\alpha_j,\rho)$
is a principal homogeneous space\index{principal homogeneous space} for $W_{\valpha}(\rho)$ under the action 
\begin{equation*}
w \ap \vx  \ = \ \crJ_{\valpha}^{-1}(w \Jp \crJ_{\valpha}(\vx)) \ ,  
\end{equation*}   
and if we write $[\vx]$ for $[(x_1)+\cdots+(x_g)]$, then 
\begin{equation*}
[w \ap \vx] \ = \ w \Jp [\vx] \ .
\end{equation*} 

Below are the properties of the action we will need;  for a more general statement, see Theorem \ref{BKeyThm1} 
of Appendix \ref{AppD}.  

\vskip .05 in

\begin{theorem}  \label{BKeyThm1A}\index{local action of the Jacobian|ii} 
Let $K_v$ be a nonarchimedean local field, 
and let $\cC_v/K_v$ be a smooth, projective, 
geometrically integral curve of genus $g > 0$.  
Then the points $\valpha = (\valpha_1, \ldots, \valpha_g) \in \cC_v(\CC_v)^g$ such that 
$\crJ_{\valpha} : \cC_v(\CC_v)^g \rightarrow \Jac(\cC_v)(\CC_v)$
is nonsingular at $\valpha$ are dense in $\cC_v(\CC_v)^g$ for the $v$-topology.
If $F_w/K_v$ is a finite extension in $\CC_v$ and $\cC_v(F_w)$ is nonempty, 
they are dense in $\cC_v(F_w)^g$.    
  
Fix such an $\valpha;$  then $\alpha_1, \ldots, \alpha_g$ are distinct, 
and for each $0 < \eta < 1$, there is a number $0 < R_2 < 1$ 
$($depending on $\valpha$ and $\eta)$ such that $B(\alpha_1,R_2), \ldots, B(\alpha_g,R_2)$ 
are pairwise disjoint and isometrically parametrizable,
\index{isometrically parametrizable ball}
the map $\crJ_{\valpha} : \cC_v(\CC_v)^g \rightarrow \Jac(\cC_v)(\CC_v)$ 
is injective on $\prod_{j=1}^g B(\alpha_j,R_2)$,  
and for each $0 < \rho \le R_2$ the following properties hold$:$ 

$(A)$ {\rm $($Subgroup$)$} 
The set $W_{\valpha}(\rho) := \crJ_{\valpha}\big(\prod_{j=1}^g B(\alpha_j,\rho)\big)$ 
is an open subgroup of $\Jac(\cC_v)(\CC_v)$.  

$(B)$ {\rm $($Limited Distortion$)$}  
For each $j = 1, \ldots g$, let $\varphi_j : D(0,\rho) \rightarrow B(\alpha_j,\rho)$
be an isometric parametrization with $\varphi_i(0) = \alpha_i$, and 
\index{isometric parametrization}  
let $\Phi_{\valpha} = (\varphi_1,\ldots,\varphi_g) : D(0,\rho)^g \rightarrow \prod_{j=1}^g B(\alpha_j,\rho)$ 
be the associated map.  Let $\Psi : D(\vORIG,1)^- \ \rightarrow \ B_J(O,1)^-$ 
be the isometric parametrization inducing the formal group,
\index{formal group}
\index{isometric parametrization}    
and let $L_{\valpha} : \CC_v^g \rightarrow \CC_v^g$ be the linear map 
$(\Psi^{-1} \circ \crJ_{\valpha} \circ \Phi_{\valpha})^{\prime}(\vORIG)$.  
 
Then $W_{\valpha}(\rho) = \Psi(L_{\valpha}(D(0,\rho)^g))$. 
If we give $D(0,\rho)^g$ the structure of an additive subgroup of $\CC_v^g$,
the map $\Psi \circ L_{\valpha}$ induces an isomorphism of groups 
\begin{equation} \label{FGrpCong}
D(0,\rho)^g/D(0,\eta \rho)^g \ \cong \ W_{\valpha}(\rho)/W_{\valpha}(\eta \rho) 
\end{equation} 
with the property that for each $\vx \in D(0,\rho)^g$, 
\begin{equation} \label{FCommute}
\crJ_{\valpha}(\Phi_{\valpha}(\vx)) \ \equiv \ \Psi(L_{\valpha}(\vx))  \pmod{W_{\valpha}(\eta \rho)} \ . 
\end{equation} 

$(C)$ {\rm $($Action$)$} There is an action $(\omega,\vx) \mapsto \omega \ap \vx$
of $W_{\valpha}(\rho)$ on $\prod_{j=1}^g B(\alpha_j,\rho)$ 
which makes $\prod_{j=1}^g B(\alpha_j,\rho)$ into a 
principal homogeneous space\index{principal homogeneous space} for $W_{\valpha}(\rho)$. 
If we restrict the domain of $\crJ_{\valpha}$ to $\prod_{j=1}^g B(\alpha_j,\rho)$, 
then $\omega \ap \vx = \crJ_{\valpha}^{-1}(\omega \Jp \crJ_{\valpha}(\vx))$. 
For each \, $\omega \in W_{\valpha}(\rho)$   
and $\vx \in \prod_{j=1}^g B(\alpha_j,\rho)$,
\begin{equation} \label{FNat1AA}
[\omega \ap \vx] \ = \ \omega \Jp [\vx] \ .
\end{equation}

$(D)$ {\rm $($Uniformity$)$} For each $\vbeta \in \prod_{i=1}^g B(\alpha_i,\rho)$, 
\begin{equation} \label{FCenter1AA}  
W_{\valpha}(\eta \rho) \, \ap \, \vbeta  \ = \ \prod_{j=1}^g B(\beta_j,\eta \rho) 
\quad \text{and} \quad 
\crJ_{\vbeta}\Big(\prod_{j=1}^g B(\beta_j,\eta \rho)\Big) \ = \ W_{\valpha}(\eta \rho)\ . 
\end{equation} 

$(E)$ {\rm $($Rationality$)$} If $F_w/K_v$ is a finite extension in $\CC_v$, 
and $\valpha \in \cC_v(F_w)^g$,  then   
\begin{eqnarray}
\crJ_{\valpha}(\prod_{j=1}^g (B(\alpha_j,\rho) \cap \cC_v(F_w))) \ = \ 
W_{\valpha}(\rho) \cap \Jac(\cC_v)(F_w) \ , \label{FRat1AA} \\
(W_{\valpha}(\rho) \cap \Jac(\cC_v)(F_w)) \, \ap \, \valpha 
\ = \ \prod_{j=1}^g \big(B(\alpha_j,\rho) \cap \cC_v(F_w)\big) \ . \label{FRat2AA}
\end{eqnarray}

$(F)$ {\rm $($Trace$)$} If $F_w/K_v$ is finite and separable, 
there is a constant $C = C(F_w,\valpha) > 0$, depending on $F_w$ and $\valpha$ but not on $\rho$, such that 
\begin{equation} \label{FTraceContainsAA} 
B_J(O,C \rho) \cap \Jac(\cC_v)(K_v) \ \subseteq \ 
\Tr_{F_w/K_v}\big(W_{\valpha}(\rho) \cap \Jac(\cC_v)(F_w)\big)  \ .
\end{equation} 
\end{theorem}

\begin{proof} 
This is a specialization of Theorem \ref{BKeyThm1} of Appendix \ref{AppD}. 
In particular, $R_2$ is the number $R$ from Theorem \ref{BKeyThm1}.
We restrict to the case where the vectors $\vr = (r_1, \ldots, r_g)$\ in that theorem have equal coordinates, 
so the condition $\eta \cdot \max(r_i) \le \min(r_i)$ is automatic. 
For $\vr = (\rho, \ldots,\rho) \in \RR^g$, the subgroups denoted $W_{\valpha}(\vr)$ 
in Appendix \ref{AppD} are written $W_{\valpha}(\rho)$ here.  
The balls denoted $B_{\cC_v}(a_j,r)$ in Appendix \ref{AppD}  
are the balls $B(\alpha_j,\rho) \subset \cC_v(\CC_v)$ here. 
\end{proof}
\index{initial approximating functions $f_v(z)$!construction when $K_v$ nonarchimedean!the local action theorem|)}

\begin{proof}[Proof of Proposition \ref{K_vReductionProp}A when $g = g(\cC_v) > 0$.] { \ }
\index{initial approximating functions $f_v(z)$!construction when $K_v$ nonarchimedean!the proof when $g(\cC_v) > 0$|(}

Let $F_w$ be a finite, separable extension of $K_v$ in $\CC_v$, 
and let $E_w = \cC_v(F_w) \cap B(a,r)$,  
where $a \in \cC_v(F_w)$ and $B(a,r) \subset \cC_v(\CC_v) \backslash \fX$
is isometrically parametrizable. 
\index{isometrically parametrizable ball}

We first construct the set $\tE_w$.  

Let $\varepsilon_w > 0$ be given.   
Since $E_w$ is open in $\cC_v(F_w)$,\index{local action of the Jacobian}  
by Theorem \ref{BKeyThm1A} we can choose $\alpha_1, \ldots, \alpha_g \in E_w$
such that the map $\crJ_{\valpha} : \cC_v(\CC_v)^g \rightarrow \Jac(\cC_v)(\CC)$ 
is injective at $\valpha = (\alpha_1, \cdots, \alpha_g)$.   
Let $R_1$ be the number given  by Proposition \ref{BPropF3A1} for $\valpha$ and $\varepsilon_w$.
Fix a uniformizer $\pi_w$ for the maximal ideal of $\cO_w$.  
If the residue characteristic of $F_w$ is odd, take $\eta = |\pi_w|_v$.  
If the residue characteristic of $F_w$ is $2$, take $\eta = |\pi_w^2|_v$.  
Let $R_2$ be the number given by Theorem \ref{BKeyThm1A} for $\valpha$ and $\eta$.\index{local action of the Jacobian}   
Without loss we can assume that $R_1$ and $R_2$ belong to $|F_w^{\times}|_v$.  
Set $\rho = \min(R_1,R_2,r)$, and put 
\begin{equation} \label{E_wContains}
\tE_w \ = \ \big(E_w \backslash \bigcup_{j=1}^g B(\alpha_j,\rho)\big) 
                         \cup \big(\bigcup_{j=1}^g (\cC_v(F_w) \cap B(\alpha_j,\eta \rho)) \big) \ .
\end{equation} 

By Theorem \ref{BKeyThm1A}(A),\index{local action of the Jacobian}  
$W_{\valpha}(\eta\rho)$ is an open subgroup of $\Jac(\cC_v)(\CC_v)$.  
For each $\tau \in \tE_w$ the continuity of the Abel map\index{Abel map}
 $j_\tau : \cC_v(\CC_v) \rightarrow \Jac(\cC_v)(\CC_v)$
shows there is a radius $r_{\tau} > 0$ 
such that $\BFj_\tau(B(\tau,r_{\tau})) \subseteq W_{\valpha}(\eta \rho)$. 
Without loss we can assume that $r_{\tau} \le \eta \rho$ and that $r_{\tau} \in |F_w^{\times}|_v$. 
If $\tau \in B(\alpha_j,\eta \rho)$ for some $j$, 
then Theorem \ref{BKeyThm1A}(D) shows that\index{local action of the Jacobian}  
$\BFj_\tau(B(\tau,\eta \rho)) \subseteq W_{\valpha}(\eta \rho)$, 
because if $\vbeta = (\alpha_1, \ldots, \alpha_{j-1},\tau, \alpha_{j+1}, \ldots, \alpha_g)$ 
and $\vx = (\alpha_1, \ldots, \alpha_{j-1},x,\alpha_{j+1}, \ldots, \alpha_g)$ then 
$\BFj_\tau(x) = \BFJ_\vbeta(\vx)$. 
 
Since $E_w \backslash \bigcup_{j=1}^g B(\alpha_j,\rho)$ is compact, 
we can choose $\alpha_{g+1}, \ldots, \alpha_d \in E_w \backslash \bigcup_{j=1}^g B(\alpha_j,\rho)$
such that the balls $B(\alpha_k,r_{\alpha_k})$ for $k = g+1, \ldots, d$
cover $E_w \backslash \bigcup_{j=1}^g B(\alpha_j,\rho)$. For notational simplicity, write $r_k$ for $r_{\alpha_k}$.
Since any two balls are either pairwise disjoint or one is contained in the other, 
we can assume that $B(\alpha_{g+1},r_{g+1}), \ldots, B(\alpha_d,r_d)$ are pairwise disjoint. 
They are also disjoint from $B(\alpha_1,\rho), \ldots, B(\alpha_g,\rho)$.
For $k = 1, \ldots, g$, put $r_k = \eta \rho$.
Then $B(\alpha_1,r_1), \ldots, B(\alpha_d,r_d)$ are pairwise disjoint and cover $\tE_w$, and 
\begin{equation} \label{tE_wSimple}
\tE_w \ = \ \bigcup_{j=1}^d \big(\cC_v(F_w) \cap B(\alpha_j,r_j)\big) 
\end{equation} 
is an $F_w$-simple decomposition of $\tE_w$.

Since $\rho \le R_1$ and (\ref{E_wContains}) holds, Proposition \ref{BPropF3A1} shows that  
$|V_{x_i}(E_w) - V_{x_i}(\tE_w)| < \varepsilon_w$ for each $x_i \in \fX$.  
\index{Robin constant!nonarchimedean} 
Thus $\tE_w$ satisfies part (A) of Proposition \ref{K_vReductionProp}A. 

\smallskip 
We next apply Lemma \ref{ConstructionLemma} to show that 
part (B) of Proposition \ref{K_vReductionProp}A holds.  

Fix an $F_w$-rational isometric parametrization $\varphi : D(0,r) \rightarrow B(a,r)$ with $\varphi(0) = a$.
\index{isometric parametrization}    
Let $\talpha_1, \ldots, \talpha_d \in F_w \cap D(0,r)$ be the points with $\varphi(\talpha_j) = \alpha_j$, 
and put   
\begin{equation*} 
\tH_w \ = \ \bigcup_{j=1}^d \big(F_w \cap D(\talpha_j,r_j)\big) \ .
\end{equation*}  
Then $\varphi(\tH_w) = \tE_w$.  

Let an $F_w$-symmetric probability vector $\vs \in \cP^m(\QQ)$
and a number $0 < \beta_w \in \QQ$ be given.  
To construct the $(\fX,\vs)$-functions\index{$(\fX,\vs)$-function} 
$f_w$ in Proposition \ref{K_vReductionProp}A, 
we begin with Stirling polynomials for $\tH_w$,
\index{Stirling polynomial!generalized}
as in the proof when $g(\cC_v) = 0$, but we then modify them.     

Let $J_w$ be the group $\Jac(\cC_v)(F_w)$, 
and let 
\begin{equation*}
J_w(\rho) = W_{\valpha}(\rho) \cap \Jac(\cC_v)(F_w) \ , \qquad 
J_w(\eta \rho) = W_{\valpha}(\eta \rho) \cap \Jac(\cC_v)(F_w) \ .
\end{equation*}   
By Theorem \ref{BKeyThm1A}(A),\index{local action of the Jacobian}  $J_w(\eta \rho)$ is open in $J_w$.
Since $J_w$ is compact, the quotient group $J_w/J_w(\eta \rho)$ is finite.
Let $I$ be its order.   
Let $0 < Q \in \ZZ$ be the number given for $\tH_w$ for Proposition \ref{BGProp1AA}, 
related to the construction of Stirling polynomials; 
\index{Stirling polynomial!generalized} 
in particular, if $w_j$ is the weight of $\cC_v(F_w) \cap B(\alpha_j,r_j)$
for the equilibrium distribution $\mu_{\fX,\vs}$ of $\tH_w$, then $0 < Q w_j \in \ZZ$ for each $j$.   
Write $\beta_w = S/R$ with coprime integers $S, R$, where $R > 0$,  
and put $V = V_\infty(H_w)$. By Corollary \ref{BFCor2}, $V \in \QQ$; 
\index{Robin constant!nonarchimedean} 
write $V = X/Y$ with coprime integers $X, Y$, where $Y > 0$.
Finally, for each $i = 1, \ldots, m$ write $s_i = A_i/B_i$ with coprime integers $A_i,B_i$, where $B_i > 0$,       
and set $N_0 = 4 I Q R Y \cdot \LCM(B_1, \ldots, B_m)$.   
      
Suppose $0 < N \in \NN$ is a multiple of $N_0$.
Let $S_{N,\tH_w}(z) = \prod_{k=1}^N (z-\vartheta_k)$ 
be the Stirling polynomial of degree $N$ for $\tH_w$ given by Proposition \ref{BGProp1AA}, 
\index{Stirling polynomial!generalized}
and put $\theta_k^* = \varphi(\vartheta_k) \in \tE_w$ for $k = 1, \ldots, N$.  
Then  
\begin{equation*} 
D_N^* \ := \ \sum_{k=0}^{N-1} (\theta_k^*) - \sum_{i=1}^m N s_i(x_i) 
\end{equation*} 
is an $F_w$-rational $(\fX,\vs)$-divisor of degree $0$ on $\cC_v$. 
Its positive part is supported on $\tE_w$ and its polar part is supported on $\fX$. 

\smallskip
Condition (1) of Lemma \ref{ConstructionLemma} may fail for $D_N^*$, since it need not be principal,  
but we claim that conditions (2) and (3) hold for it.  

For condition (2), note that by Proposition \ref{BGProp1AA}(A) applied to $\tH_w$, there is a number $\tA > 0$ 
such that $|\vartheta_i - \vartheta_j|_v > \tA/N$ for all $N$ and all $i \ne j$.
Since $\varphi : D(0,r) \rightarrow B(a,r)$ is an $F_w$-rational isometric parametrization,
\index{isometric parametrization}  
we have $\|\theta_i^* , \theta_j^*\|_v = |\vartheta_i-\vartheta_j|_v > \tA/N$ for all $i \ne j$.

For condition (3), 
let $Q_N^*(z) = \prod [z,\theta_k^*]_{\fX,\vs}$ be the $(\fX,\vs)$-pseudo-polynomial associated with $D_N^*$.
By Proposition \ref{APropA2}(B2), there is a number $C_{\fX,\vs} \in |\CC_v^{\times}|_v$ 
such that $[z,w]_{\fX,\vs} = C_{\fX,\vx} \|z,w\|_v$ for all $z, w \in B(a,r)$.
For $x, y \in D(0,r)$ we have $\|\varphi(x),\varphi(y)\|_v = |x-y|_v$, 
so $[\varphi(x),\varphi(y)]_{\fX,\vs} = C_{\fX,\vs}|x-y|_v$. Thus for $x \in D(0,r)$ 
\begin{equation*} 
Q_N^*(\varphi(x)) \ = \ \prod_{k=1}^N [\varphi(x),\theta_k^*]_{\fX,\vs} 
\ = \ C_{\fX,\vs}^N \prod_{k=1}^N |x-\vartheta_k|_v
\ = \ C_{\fX,\vs}^N \, |S_{N,\tH_w}(x)|_v \ ,
\end{equation*} 
and so
\begin{equation} \label{FQNNN1B} 
- \log_v(Q_N^*(\varphi(x))) \ = \ -N\log_v(C_{\fX,\vs}) -  \log_v(|S_{N,\tH_w}(x)|_v) \ .
\end{equation} 
As in (\ref{FQNNN2}), 
\begin{equation}
V_{\fX,\vs}(\tE_w)  \ = \ V_{\infty}(\tH_w) - \log_v(C_{\fX,\vs}) \ . \label{FQNNN2B}
\end{equation} 
\index{Robin constant!nonarchimedean}\index{Robin constant!nonarchimedean!nonarchimedean $(\fX,\vs)$}
Put $\tB = \max_{1 \le j \le d}(\log_v(r_j))$.  Given $z \in B(a,r)$, let $x \in D(0,r)$ be such that $\varphi(x) = z$, 
and let $\vartheta_J$ be the root of $S_N(z,\tH_w)$
for which $|x - \vartheta_J|_v$ is minimal.  
Using (\ref{FQNNN1B}), (\ref{FQNNN2B}) and Proposition \ref{BGProp1AA}(C), we obtain
\begin{eqnarray} 
-\log_v(Q_N^*(z)) 
& = & -N \log_v(C_{\fX,\vs}) - \log_v(|S_{N,\tH_w}(x))|_v \notag\\
& \le & - N \log_v(C_{\fX,\vs}) + N \cdot V_{\infty}(\tH_w) - \log_v(|x-\vartheta_J|_v) + \tB  \notag \\
& = & \ N \cdot V_{\fX,\vs}(\tE_w) - \log_v(\|z,\theta_J^*\|_v) + \tB \ . \label{FCond3}
\end{eqnarray}
Thus condition (3) holds.\index{Robin constant!nonarchimedean!nonarchimedean $(\fX,\vs)$}

\smallskip
We will now modify $D_N^*$ to obtain a principal divisor $D_N$ which satisfies 
all the conditions of Lemma \ref{ConstructionLemma}.

Put $\delta = [D_N^*]$.  We claim that $\delta \in J_w(\eta \rho)$.  
To see this, for each $k = 1, \ldots, N$ let $j(k)$ denote the index $1 \le j \le d$ 
for which $\theta_k^* \in B(\alpha_j,r_j)$.  For each $j = 1, \ldots, d$ put $P_j = N \cdot w_j$ 
and for each $i = 1, \ldots, m$ put $Q_i = N \cdot s_i$.  The numbers $P_j$ and $Q_i$ belong to $\NN$ since $N_0|N$, 
and we can write
\begin{equation*}
\delta \ = \ [D_N^*] \ = \ \sum_{k=1}^N [(\theta_k^*) - (\alpha_{j(k)})] + \sum_{j=1}^d P_j [(\alpha_j) - (a)] 
- \sum_{i=1}^m Q_i [(x_i) - (a)] \ .
\end{equation*}
By the construction of the balls $B(\alpha_j,r_j)$, 
for each $k$ we have $[(\theta_k^*) - (\alpha_{j(k)})] \in J_w(\eta \rho)$.
Since $\alpha_j$ and $a$ belong to $\cC_v(F_w)$, for each $j$ we have $[(\alpha_j) - (a)] \in J_w$.  
By our choice of $N_0$, each $P_j$ is divisible by $I$, so for each $j$ we have $P_j [(\alpha_j) - (a)] \in J_w(\eta \rho)$.
Similarly, each $Q_i$ is divisible by $I$, and $\vs$ and $\fX$ are $F_w$-symmetric,
so $\sum_{i=1}^m Q_i [(x_i) - (a)] \in J_w(\eta \rho)$. Thus $\delta \in J_w(\eta \rho)$.

Let $\ell$ be such that $|\pi_w^\ell|_v = \rho$.  If the residue characteristic of $F_w$ is odd, 
then by parts (B) and (E) of Theorem \ref{BKeyThm1A},\index{local action of the Jacobian}  
together with our choice of $\eta$, 
the group $J_w(\rho)/J_w(\eta \rho)$ is isomorphic to $(\pi_w^\ell \cO_w/\pi_w^{\ell+1} \cO_w)^g$.  
If the residue characteristic of $F_w$ is $2$, 
then $J_w(\rho)/J_w(\eta \rho)$ is isomorphic to $(\pi_w^\ell \cO_w/\pi_w^{\ell+2} \cO_w)^g$.  

Let $\Psi : D(\vORIG,0)^- \rightarrow B_J(O,1)^-$ be the $F_w$-rational isometric parametrization
\index{isometric parametrization}  
in Theorem \ref{BKeyThm1A}(B).\index{local action of the Jacobian}   
Using the $F_w$-rational isometric parametrization
\index{isometric parametrization}  
$\varphi : D(0,r) \rightarrow B(a,r)$, we get $F_w$-rational isometric parametrizations  
$\varphi_j : D(0,\rho) \rightarrow B(\alpha_j,\rho)$, for 
$j = 1, \ldots, g$, by setting $\varphi_j(z) = \varphi(\talpha_j + z)$ where $\varphi(\talpha_j) = \alpha_j$.  
Define $\Phi : D(0,\rho)^g \rightarrow \prod_{j=1}^g B(\alpha_j,\rho)$ by 
$\Phi(\vz) = (\varphi_1(z_1), \ldots, \varphi_g(z_g))$.
Then $\Psi((\pi_w^{\ell} \cO_w)^g) = J_w(\rho)$, and $\varphi_j(\pi_w^{\ell} \cO_w) = \cC_v(F_w) \cap B(\alpha_j,\rho)$
for each $j$.
Let $L_{\valpha} : \CC_v^g \rightarrow \CC_v^g$ be the $F_w$-rational linear map defined by  
$L_{\valpha} = (\Psi^{-1} \circ \BFJ_{\valpha} \circ \Phi)^{\prime}(\vORIG)$.
Let $\ap$ be the action of $W_{\valpha}(\rho)$ on $\prod_{j=1}^g B(\alpha_j,\rho)$ from Theorem \ref{BKeyThm1A}.  
By parts (C) and (E) of Theorem \ref{BKeyThm1A},
for each $\vz \in (\pi_w^\ell \cO_w)^g$ we have $\Psi(L_{\valpha}(\vz)),\, \BFJ_{\valpha}(\Phi(\vz)) \in J_w(\rho)$,
with $\Psi(L_{\valpha}(\vz)) \equiv\BFJ_{\valpha}(\Phi(\vz)) \pmod{J_w(\eta \rho)}$.

First suppose the residue characteristic of $F_w$ is odd.  
By our choice of $N_0$, at least two points $\theta_k^*$ 
belong to $B(\alpha_j,r_j)$ for each $j$.  
Recall that $r_j = \eta \rho$ for $j = 1, \ldots, g$.  
Without loss, we can assume that $\theta_1^*, \ldots, \theta_N^*$ 
are labeled in such a way that $\theta_j^*$ and $\theta_{g+j}^*$ 
belong to $B(\alpha_j,\eta \rho)$, for $j = 1, \ldots, g$. 

Fix an element $t \in \pi_w^{\ell} \cO_w$ whose image in $\pi_w^{\ell} \cO_w / \pi_w^{\ell+1} \cO_w$ 
is nonzero. Put  $\vt = (t, \cdots, t) \in (\pi_w^\ell \cO_w)^g$ and 
set $\Delta = \Psi(L_{\valpha}(\vt))$, 
then define $\theta_1, \ldots, \theta_{2g} \in \cC_v(F_w)$ by 
\begin{equation*}
(\theta_1, \ldots, \theta_g) \ = \ \Delta \ap (\theta_1^*, \ldots,\theta_g^*) \ , \qquad 
(\theta_{g+1}, \ldots, \theta_{2g}) \ = \ (-(\delta \Jp \Delta)) \ap (\theta_1^*, \ldots,\theta_g^*) \ .
\end{equation*}
Put $\theta_k = \theta_k^*$ for $k = 2g+1, \ldots, N$, and set 
\begin{equation*}
D_N \ = \ \sum_{k=0}^{N-1} (\theta_k) - \sum_{i=1}^m N s_i(x_i) \ .
\end{equation*}
The divisor $D_N$ is principal, since by (\ref{FNat1AA}) we have 
$[(\theta_1)+ \cdots + (\theta_g)]  =  \Delta \Jp [(\theta_1^*)+ \cdots + (\theta_g^*)]$ and  
$[(\theta_{g+1})+ \cdots + (\theta_{2g})] = (-(\delta \Jp \Delta)) \Jp [(\theta_{g+1}^*)+ \cdots + (\theta_{2g}^*)]$, which gives
\begin{eqnarray} 
[D_N] & = & [D_N^*] \, \Jp \, \Big([\sum_{k=1}^{g} (\theta_k)] \Jm [\sum_{k=1}^{g} (\theta_k^*)]\Big) \, \Jp \, 
\Big([\sum_{k=g+1}^{2g} (\theta_{k})] \Jm [\sum_{k=g+1}^{2g} (\theta_{k}^*)]\Big) \notag \\
& = & \delta \, \Jp \, \Delta \, \Jp \, (-(\delta \Jp \Delta)) \ = \ 0 \ . \label{DNPrincipal}
\end{eqnarray} 

We will now show that $\theta_1, \ldots, \theta_{2g}$ belong to $E_w \backslash \tE_w$, 
and that they are well-separated from 
each other and $\theta_{2g+1}, \ldots, \theta_N$.

Put $\vbeta = (\beta_1, \ldots, \beta_g) = (\varphi_1(t), \ldots, \varphi_g(t))$. 
Then $\beta_j \in \cC_v(F_w) \cap B(\alpha_j,\rho)$ for each $j$, 
and by Theorem \ref{BKeyThm1A}(B)\index{local action of the Jacobian} 
\begin{equation*} 
\crJ_{\valpha}(\vbeta) \ = \ \crJ_{\valpha}(\Phi(\vt)) 
\ \equiv \ \Psi(L_{\valpha}(\vt)) \ = \ \Delta \pmod{W_{\valpha}(\eta \rho)} \ .
\end{equation*} 
This means there is a $\delta^{\prime} \in W_{\valpha}(\eta \rho)$ 
such that $\crJ_{\valpha}(\vbeta) = \delta^{\prime} \Jp \Delta$.
It follows that $\delta^{\prime} \ap (\Delta \ap \valpha) = (\delta^{\prime} \Jp \Delta) \ap \valpha 
= \crJ_{\valpha}(\vbeta) \ap \valpha = \vbeta$, 
so  $\Delta \ap \valpha = (-\delta^{\prime}) \ap \vbeta$. 
Likewise, put $\delta^* = \crJ_{\valpha}((\theta_1^*, \ldots, \theta_g^*)) \in W_{\valpha}(\eta \rho)$,
so $(\theta_1^*, \ldots, \theta_g^*) = \delta^* \ap \valpha$.   
By properties of the action $\ap$ and Theorem \ref{BKeyThm1A}(D)\index{local action of the Jacobian}  
\begin{eqnarray*} 
(\theta_1, \ldots, \theta_g) & = & \Delta \ap (\theta_1^*, \ldots, \theta_g^*)  
\ = \ \Delta \ap (\delta^* \ap \valpha) \ = \ \delta^* \ap (\Delta \ap \valpha)
\ = \ \delta^* \ap ((-\delta^{\prime}) \ap \vbeta) \\
& = & (\delta^* \Jm \delta^{\prime}) \ap \vbeta 
\ \in \ \prod_{j=1}^g \big(\cC_v(F_w) \cap B(\beta_j,\eta \rho)\big) \ .
\end{eqnarray*}
Similarly, put $\vgamma = (\gamma_1, \ldots, \gamma_g) = (\varphi_1(-t), \ldots, \varphi_g(-t))$;  
then $\gamma_j \in \cC_v(F_w) \cap B(\alpha_j,\rho)$ for each $j$.    
Let $\delta^{**} = \crJ_{\valpha}(\theta_{g+1}^*, \ldots, \theta_{2g}^*) \in W_{\valpha}(\eta \rho)$.
By computations like those above one sees that there is a 
$\delta^{\prime \prime} \in W_{\valpha}(\eta \rho)$ 
such that $\crJ_{\valpha}(\vgamma) = \delta^{\prime \prime} \Jp (-\Delta)$, 
and that 
\begin{equation*}
(\theta_{g+1}, \ldots, \theta_{2g}) 
\ = \ (\delta^{**} \Jm \delta^{\prime \prime} \Jm \delta) \ap \vgamma 
\ \in \ \prod_{j=1}^g \big(\cC_v(F_w) \cap B(\gamma_j,\eta \rho)\big) \ .
\end{equation*}   

For each $j$ the map $\varphi_j : D(0,\rho) \rightarrow B(\alpha_j,\rho)$ is an isometric parametrization,
\index{isometric parametrization}   
so our choice of $t$ means that $\|\beta_j,\alpha_j\|_v = |t|_v = |\pi_w^{\ell}|_v = \rho$
and $\|\gamma_j,\alpha_j\|_v = |-t|_v = \rho$.  
Moreover $\|\beta_j,\gamma_j\|_v = |t-(-t)|_v = |2 t|_v = \rho$ since the residue characteristic is odd.
Thus the balls $B(\alpha_j, \eta \rho)$, $B(\beta_j, \eta \rho)$ and $B(\gamma_j,\eta \rho)$
are pairwise disjoint and contained in $B(\alpha_j, \rho)$.  Since $B(\alpha_1, \rho), \ldots, B(\alpha_g,\rho)$  
are pairwise disjoint, all the balls 
\begin{equation*}
B(\alpha_1, \eta \rho), \ldots, B(\alpha_g, \eta \rho), \ B(\beta_1, \eta \rho), \ldots, B(\beta_g, \eta \rho), \ 
B(\gamma_1, \eta \rho), \ldots, B(\gamma_g, \eta \rho)
\end{equation*} 
are pairwise disjoint.  Since $\tE_w = \big(E_w \backslash \bigcup_{j=1}^g B(\alpha_j, \rho)\big) 
\cup \big( \bigcup_{j=1}^g (\cC_v(F_w) B(\alpha_j,\eta \rho)\big)$,    
it follows that $\theta_1, \ldots, \theta_{2g} \in E_w \backslash \tE_w$. 
In addition, for each $j = 1, \ldots, 2g$,  if $1 \le k \le N$ and $k \ne j$ then 
\begin{equation} \label{FSep1}
\|\theta_j, \theta_k\|_v \ \ge \ \rho \ .  
\end{equation}

\smallskip
Next suppose the residue characteristic of $F_w$ is $2$.  
Since $\pi_w^\ell \cO_w/\pi_w^{\ell+2} \cO_w$ is an abelian $2$-group 
with at least $4$ elements, it either has a subgroup isomorphic to $\ZZ/4\ZZ$ or a subgroup isomorphic to 
$(\ZZ/2 \ZZ)^2$.   If it has a subgroup isomorphic to $\ZZ/4\ZZ$, choose an element $t \in \pi_w^{\ell} \cO_w$ whose image in 
$\pi_w^\ell \cO_w/\pi_w^{\ell+2} \cO_w$ has order $4$.  
Then the same construction as in the case of odd residue characteristic applies, 
but at the very end, in place of (\ref{FSep1}) one gets that for 
$j = 1, \ldots, 2g$ and $1 \le k \le N$ with $k \ne j$ 
\begin{equation} \label{FSep2}
\|\theta_j, \theta_k\|_v \ \ge \ |\pi_w|_v \cdot \rho \ .
\end{equation}   

If the residue characteristic is $2$ and $\pi_w^\ell \cO_w/\pi_w^{\ell+2} \cO_w$ has no elements of order $4$,
we modify the construction as follows.  
Let $t_1, t_2 \in \pi_w^\ell \cO_w$ be elements whose images in $\pi_w^\ell \cO_w/\pi_w^{\ell+2} \cO_w$ 
generate a subgroup isomorphic to $(\ZZ/2\ZZ)^2$.   
By our choice of $N_0$, at least three of the points $\theta_1^*, \ldots, \theta_N^*$ 
belong to $\cC_v(F_w) \cap B(\alpha_j,r_j)$ for each $j$.  
Without loss, we can assume that $\theta_1^*, \ldots, \theta_N^*$ 
are indexed in such a way that $\theta_j^*$, $\theta_{g+j}^*$ and $\theta_{2g+j}^*$ 
belong to $B(\alpha_j,\eta \rho)$ for $j = 1, \ldots, g$.
Put  $\vt_1 = (t_1, \cdots, t_1)$ and $\vt_2 = (t_2, \cdots, t_2)$, then 
set $\Delta_1 = \Psi(L_{\valpha}(\vt_1))$, $\Delta_2 = \Psi(L_{\valpha}(\vt_2))$.
Recall that $\delta = [D_N^*] \in J_w(\eta \rho)$.   
Define $\theta_1, \ldots, \theta_{3g} \in \cC_v(F_w)$ by 
\begin{eqnarray*}
& & (\theta_1, \ldots, \theta_g) \ = \ \Delta_1 \ap (\theta_1^*, \ldots,\theta_g^*) \ , \quad 
(\theta_{g+1}, \ldots, \theta_{2g}) \ = \ \Delta_2 \ap (\theta_{g+1}^*, \ldots,\theta_{2g}^*) \ , \\
& & \qquad \quad
(\theta_{2g+1}, \ldots, \theta_{3g}) \ = \ (-(\delta \Jp \Delta_1 \Jp \Delta_2)) \ap (\theta_{2g+1}^*, \ldots,\theta_{3g}^*) \ , 
\end{eqnarray*}
and put $\theta_k = \theta_k^*$ for $k = 3g+1, \ldots, N$. If we take 
\begin{equation*}
D_N \ = \ \sum_{k=0}^{N-1} (\theta_k) - \sum_{i=1}^m N s_i(x_i) \ ,
\end{equation*} 
an argument similar to the one before shows that $D_N$ is principal, 
that $\theta_1, \ldots, \theta_{3g}$ belong to $E_w \backslash \tE_w$, and that
for each $j = 1, \ldots, 3g$,  if $1 \le k \le N$ and $k \ne j$ then 
\begin{equation} \label{FSep3}
\|\theta_j, \theta_k\|_v \ \ge \ |\pi_w|_v \cdot \rho \ .  
\end{equation}  

\smallskip
Finally, we show that for sufficiently large $N$ divisible by $N_0$, 
the divisor $D_N$ satisfies conditions (1), (2), and (3) of Lemma \ref{ConstructionLemma}.  

We have already seen that $D_N$ is principal, so condition (1) holds.

For condition (2), let $A = \tA$ be the constant from Proposition \ref{BGProp1AA}(A) 
for the set $\tH_w$.  As shown above, for all $N$ divisible by $N_0$ 
and all $j \ne k$ with $1 \le j, \, k \le N$, we have $\|\theta_j^*,\theta_k^*\|_v \ge \tA/N$.
Suppose in addition that $N$ is large enough that $A/N < |\pi_w|_v \cdot \rho$.    
For $1 \le j \le 2g$ (resp. $3g$ in the third case),  
and all $k \ne j$ we have $\|\theta_j,\theta_k\|_v \ge |\pi_w|_v \cdot \rho >  A/N$ 
by (\ref{FSep1}), (\ref{FSep2}) and (\ref{FSep3}).  
By symmetry, this also holds for  $1 \le k \le 2g$ (resp. $3g$) and all $j \ne k$.
For $j, k > 2g$ (resp. $3g$ in the third case) with $j \ne k$
we have $\|\theta_j,\theta_k\|_v > A/N$ since $\theta_j = \theta_j^*$, $\theta_k = \theta_k^*$.  
Thus condition (2) holds. 

For condition (3), consider the $(\fX,\vs)$ pseudo-polynomials $Q_N(z) = \prod_{j=1}^N [z,\theta_j]_{\fX,\vs}$
and $Q_N^*(z) = \prod_{j=1}^N [z,\theta_j^*]_{\fX,\vs}$.  
We must show there is a constant $B$  such that for each $z \in B(a,r)$, if 
$J = J(z)$ is an index for which $\|z,\theta_J\|_v = \min_{1 \le j \le N} (\|z,\theta_j\|_v)$, then  
\begin{equation} \label{FNextQ}
-\log_v(|Q_N(z)|_v) \ \le \ N \cdot V_{\fX,\vs}(\tE_w) + B + \big(-\log_v(\|z,\theta_J\|_v)\big) \ .
\end{equation}\index{Robin constant!nonarchimedean!nonarchimedean $(\fX,\vs)$}
By (\ref{FCond3}), there is a constant $\tB$ such that for each $z \in B(a,r)$, 
if $K = K(z)$ is an index for which $\|z,\theta_K^*\|_v = \min_{1 \le j \le N} (\|z,\theta_j^*\|_v)$, then
\begin{equation*} 
-\log_v(Q_N^*(z)) \ \le \ \ N \cdot V_{\fX,\vs}(\tE_w) + \tB + \big(- \log_v(\|z,\theta_K^*\|_v)\big) \ .
\end{equation*} 
Throughout the discussion below, $J$ and $K$ will have this meaning.

Let $G$ be the number of roots of $Q_N^*(z)$ that were moved, i.e. $G = 2g$ or $G = 3g$.  Then 
\begin{equation*}
Q_N(z) \ = \ \prod_{j=1}^{G} \frac{[z,\theta_j]_{\fX,\vs}}{[z,\theta_j^*]_{\fX,\vs}} \cdot Q_N^*(z) \ .
\end{equation*} 
(If $z = \theta_h^*$ for some $h = 1, \ldots, G$, we regard the right side as defined by its limit
as $z \rightarrow \theta_h^*$.)

For all $z, w \in B(a,r)$ we have $[z,w]_{\fX,\vs} = C_{\fX,\vs} \|z,w\|_v$.  Hence if we set 
\begin{equation} \label{FDSum}
D(z) \ = \ \sum_{j=1}^G (-\log_v(\|z,\theta_j\|_v) + \sum_{j=1}^G \log_v(\|z,\theta_j^*\|_v) 
+ ( - \log_v(\|z,\theta_K^*\|_v)) \ ,
\end{equation} 
then to prove (\ref{FNextQ}) it will suffice to show there is a constant $\hB$ 
such that for all $z \in B(a,r)$ 
\begin{equation} \label{FDBound}
D(z) \ \le \  \hB + \big(-\log_v(\|z,\theta_J\|_v)\big)  \ .
\end{equation} 
(If $z = \theta_h^*$ for some $h = 1, \ldots, G$, then $\theta_K^* = \theta_h^*$, 
and we define $D(z)$ by the sum gotten by omitting the corresponding terms 
from the right side of (\ref{FDSum}); 
if $z = \theta_j$ for some $j$, we regard both sides of (\ref{FDSum}) and (\ref{FDBound}) as being $\infty$.)

We will prove (\ref{FDBound}) by considering cases.  
By (\ref{FSep1}), (\ref{FSep2}) and (\ref{FSep3}), the balls 
\begin{equation*} 
B(\theta_1,\eta \rho)\ , \ \ldots \ , \ B(\theta_G,\eta \rho) \ ,\ B(\alpha_1, r_1) \ , \ldots \ , \ B(\alpha_d, r_d)
\end{equation*}
are pairwise disjoint.  
Put $\hr = \min_{1 \le j \le d}(r_j)$ and take $\hB = G  \cdot (-\log_v(\eta \rho)) + (-\log_v(\hr))$. 

First suppose that $z \in B(\theta_h, \eta \rho)$ for some $h$, $1 \le h \le G$.  
Then $\theta_J = \theta_h$, and $\|z,\theta_j\|_v > \eta \rho$ for all $j = 1, \ldots, G$ with $j \ne h$.  
Furthermore, $\|z,\theta_j^*\|_v \le r < 1$ for $j = 1, \ldots, G$ and $\|z, \theta_K^*\|_v > \eta \rho$.  
Hence 
\begin{eqnarray*} 
D(z) & \le & (G-1) \cdot (-\log_v(\eta \rho)) + (-\log_v(\|z, \theta_J\|_v) + ( -\log_v(\eta \rho)) \\
    & \le & \hB + (-\log_v(\|z, \theta_J\|_v) .
\end{eqnarray*} 

Next suppose that $z \in B(\alpha_h, r_h)$ for some $h$, $1 \le h \le d$.  
In this case $\|z,\theta_j\|_v > \eta \rho$ for $j = 1, \ldots, G$.  
By our choice of $N_0$, at least four of the $\theta_k^*$ belong to $B(\alpha_h,r_h)$,
and at least one remains after $\theta_1^*, \ldots, \theta_G^*$ are moved.  This means that $\theta_J = \theta_\ell^*$
for some $\ell$ with $\theta_\ell^* \in B(\alpha_h, r_h)$, and that $\theta_K^* \in B(\alpha_h, r_h)$.  
By the definition of $K$, we have $\|z, \theta_J\|_v = \|z,\theta_\ell^*\|_v \ge \|z,\theta_K^*\|_v$.

If $\|z, \theta_J\|_v > \|z,\theta_K^*\|_v$, then $\theta_K^*$ was one of the roots moved out of $B(\alpha_h,r_h)$, 
and in particular $1 \le K \le G$.  Thus there is a term $\log_v(\|z,\theta_j^*\|_v)$ in second sum in (\ref{FDSum}) 
which cancels the term $-\log_v(\|z,\theta_K^*\|_v)$. (If $z = \theta_K^*$, 
then $-\log_v(\|z,\theta_K^*\|_v) = \infty$ but 
we have defined $D(z)$ by omitting these two terms from (\ref{FDSum}).)
This gives 
\begin{equation*}
D(z) \ \le \ G \cdot (-\log_v(\eta \rho)) \ \le \ \hB + (-\log_v(\|z,\theta_J\|)) \ .  
\end{equation*}  
On the other hand, if $\|z,\theta_J\|_v = \|z,\theta_K^*\|_v$, once more
\begin{equation*}
D(z) \ \le \ G \cdot (-\log_v(\eta \rho)) + (-\log_v(\|z,\theta_K^*\|) \ \le \ \hB + (-\log_v(\|z,\theta_J\|)) \ .
\end{equation*} 

Finally suppose $z \in B(a,r) \backslash \big( \bigcup_{j=1}^G B(\theta_j,\eta \rho) \cup \bigcup_{j=1}^d B(\alpha_j,r_j) \big)$.
Trivially $\|z,\theta_j\|_v > \eta \rho$ for $j = 1, \ldots, G$, and $\|z,\theta_K^*\| > \hr$. 
Since $\|z,\theta_J\|_v \le r < 1$, again we have 
\begin{equation*}
D(z) \ \le \ G \cdot (-\log_v(\eta \rho)) + (-\log_v(\hr)) \ \le \ \hB + (-\\log_v(\|z,\theta_J\|)) \ .
\end{equation*} 

This establishes (\ref{FDBound}), and completes the proof of condition (3) of Lemma \ref{ConstructionLemma}.  
Applying Lemma \ref{ConstructionLemma}, we obtain Proposition \ref{K_vReductionProp}A.  
\end{proof}

We have now completed the proof of Theorem \ref{CompactThm}.
\index{initial approximating functions $f_v(z)$!construction when $K_v$ nonarchimedean!the proof when $g(\cC_v) > 0$|)}

\section{ Corollaries to the Proof of Theorem $\ref{CompactThm}$}

In this section we note two consequences of the proof of Theorem \ref{CompactThm}  
which will be used in \S\ref{Chap11}.\ref{MovingLemmaSection}  
in the local patching construction for $K_v$-simple sets.
\index{$K_v$-simple!set} 
\index{patching argument!local}
\index{initial approximating functions $f_v(z)$!construction when $K_v$ nonarchimedean!consequences of the proof|(}   

\begin{definition} \label{MovePreparedDef}
Let $E_v \subset \cC_v(\CC_v)$ be a $K_v$-simple set\index{$K_v$-simple!set}  
with a $K_v$-simple decomposition\index{$K_v$-simple!decomposition}  
$E_v = \bigcup_{\ell=1}^D \big(B(a_\ell,r_\ell) \cap \cC_v(F_{w_\ell})\big)$. 
Let $H_v \subset E_v$ be a  $K_v$-simple set\index{$K_v$-simple!set!compatible with another set} 
compatible with $E_v$, with a $K_v$-simple decomposition\index{$K_v$-simple!decomposition} 
$H_v = \bigcup_{h=1}^N \big(B(\theta_h,\rho_h) \cap \cC_v(F_{u_h})\big)$ 
compatible with the 
decomposition\index{$K_v$-simple!decomposition!compatible with another decomposition}  
$E_v = \bigcup_{\ell=1}^D \big(B(a_\ell,r_\ell) \cap \cC_v(F_{w_\ell})\big)$.  
We will say the decomposition 
$H_v = \bigcup_{h=1}^N \big(B(\theta_h,\rho_h) \cap \cC_v(F_{u_h})\big)$
is {\em move-prepared}\index{move-prepared|ii} relative to $B(a_1,r_1), \ldots, B(a_D,r_D)$ if

 $(A)$ $g(\cC) = 0$, or

 $(B)$ $g = g(\cC) > 0$, and for each $\ell = 1, \ldots, D$ there are indices 
$h_{\ell 1}, \ldots, h_{\ell g}$ such that 

\qquad  $(1)$ 
  $B(\theta_{h_{\ell 1}},\rho_{h_{\ell 1}}), \ldots, B(\theta_{h_{\ell g}}, \rho_{h_{\ell g}})
           \subset B(a_{\ell},r_{\ell})$; 

\qquad $(2)$ there is a number $\rbar_{\ell}$ 
such that $\rho_{h_{\ell 1}}, \ldots, \rho_{h_{\ell g}} < \rbar_{\ell} < r_\ell$

\qquad \quad and $B(\theta_{h_{\ell 1}},\rbar_\ell), \ldots, B(\theta_{h_{\ell g}},\rbar_\ell)$
are pairwise disjoint and contained in $B(a_\ell,r_\ell)$;  

\qquad $(3)$ putting $\vtheta_{\ell} = (\theta_{h_{\ell 1}}, \ldots, \theta_{h_{\ell g}})$, 
the Abel map\index{Abel map} $J_{\vtheta_{\ell}} : \cC_v(\CC_v)^g \rightarrow \Jac(\cC_v)(\CC_v)$ 

\qquad \quad  is injective on $\prod_{j=1}^g B(\theta_{h_{\ell j}},\rbar_\ell)$, 
and $W_{\vtheta_{\ell}}(\rbar_{\ell}) 
:= \BFJ_{\vtheta_\ell} \big(\prod_{j=1}^g B(\theta_{h_{\ell j}},\rbar_\ell)\big)$ 

\qquad \quad  is an open subgroup of $\Jac(\cC_v)(\CC_v)$ with 
 the properties in Theorem $\ref{BKeyThm1A}$.\index{local action of the Jacobian} 

\qquad  We will call  
$B(\theta_{h_{\ell 1}},\rho_{h_{\ell 1}}), \ldots, B(\theta_{h_{\ell g}}, \rho_{h_{\ell g}})$
{\em distinguished balls} corresponding 
\index{distinguished balls}

\qquad  to $B(a_\ell,r_\ell)$ in the decomposition of $H_v$.  
\end{definition}

\begin{corollary} \label{CompactThmCor}
In Theorem $\ref{CompactThm}$, by choosing $N_v$ sufficiently large, 
we can arrange that the $K_v$-simple decomposition 
\index{$K_v$-simple!decomposition}   
$H_v = \bigcup_{h=1}^N B(\theta_h,\rho_h) \cap \cC_v(F_{u_h})$ of $H_v := f_v^{-1}(D(0,1)) \cap E_v$ 
is move-prepared relative to $B(a_1,r_1), \ldots, B(a_D,r_D)$.
\index{move-prepared}
\end{corollary} 

\begin{proof}
When $g(\cC_v) = 0$ there is nothing to show.
When $g = g(\cC_v) > 0$, the corollary follows by tracing through the proof of 
Theorem \ref{CompactThm}.  We note the key points in the argument, below.

The proof begins by using the decomposition 
$E_v = \bigcup_{\ell =1} \big(B(a_\ell,r_\ell) \cap \cC_v(F_{w_\ell})\big)$ to reduce 
to a single set of the form $E_{w,\ell} = B(a_\ell,r_\ell) \cap \cC_v(F_{w_\ell})$;  
the function $f_v$ is a scaled product of conjugates 
of functions $f_{w,\ell}$ for representatives of galois orbits of the balls $B(a_\ell,r_\ell)$ 
(see (\ref{NonarchFcnDef})). Since the Abel map\index{Abel map} is galois-equivariant, 
if the conditions in Definition
\ref{MovePreparedDef} are satisfied for some $B(a_\ell,r_\ell)$, they also hold for its conjugates.     
  
Proposition \ref{K_vReductionProp}A constructs $f_{w,\ell}$ for 
a single set $E_{w,\ell} = B(a_\ell,r_\ell) \cap \cC_v(F_{w_\ell})$. 
(For notational simplicity, the index $\ell$ is suppressed in its proof.)  
The first step in the proof (see (\ref{E_wContains}) and (\ref{tE_wSimple})) is to construct a subset 
\begin{equation*} 
\tE_{w,\ell} \ = \ \bigcup_{j=1}^{d_\ell} \big(\cC_v(F_{w_\ell}) \cap B(\alpha_{\ell j},r_{\ell j})\big) \ \subset \ E_{w,\ell}
\end{equation*} 
in which $r_{\ell 1} = \cdots = r_{\ell g} = \eta \rho < r_\ell$, 
such that $W_\valpha(\eta \rho) := \BFJ_{\valpha}\big(\prod_{j=1}^g B(\alpha_{\ell j},\eta \rho) \big)$ 
is an open subgroup of $\Jac(\CC_v)$ 
satisfying the conditions of Theorem \ref{BKeyThm1A}.\index{local action of the Jacobian}   
Put $\rbar_\ell = \eta \rho$.  In the construction of $f_{w,\ell}$, zeros are initially assigned to the 
cosets $\cC_v(F_{w_\ell}) \cap B(\alpha_{\ell j},r_{\ell j})$ 
in proportion to their weights under the $(\fX,\vs)$-equilibrium
distribution of $\tE_{w,\ell}$, giving a nonprincipal divisor $D_N^*$;  
then at most $3g$ zeros are moved to obtain a principal divisor $D_N$, 
which becomes the divisor of $f_{w,\ell}$. 
The corollary follows by noting that if $N_v$ (hence $N$)  
is sufficiently large, then $f_{w,\ell}$ 
has zeros in $\cC_v(F_{w,\ell}) \cap B(\alpha_{\ell j},\rbar_\ell)$
for each $j = 1, \ldots, g$,  
and the distance between these zeros is at most $\rbar_\ell$
so the corresponding balls $B(\theta_h,\rho_h)$ in the decomposition of $H_v$ have radii 
less than $\rbar_\ell$.  If we let $h_{\ell 1}, \ldots, h_{\ell g}$ be indices of zeros with
$\theta_{h_{\ell j}} \in B(\alpha_{\ell j}, r_{\ell})$, 
then by Theorem \ref{BKeyThm1A}(D)\index{local action of the Jacobian} 
\begin{equation} \label{WContainmentF}
W_{\vtheta_{\ell}}(\rbar_{\ell}) \ = \ 
\BFJ_{\vtheta_\ell} \big(\prod_{j=1}^g B(\theta_{h_{\ell j}},\rbar_\ell)\big) \ = \ 
\BFJ_{\valpha}\big(\prod_{j=1}^g B(\alpha_{\ell j},\eta \rho) \big) \ .
\end{equation}
When $\Char(K_v) = p > 0$,  the final step in the construction of $f_v$ 
replaces the function $f_v(z)$ described above with $\pi_v^{-N_v B/2} S_{p^B,\cO_v}(f(z))$
(see (\ref{FStirlingCompose}),
in order to assure that the leading coefficient $c_{v,i}$ belongs to $K_v(x_i)^{\sep}$ for each $i$.
\index{coefficients $A_{v,ij}$!leading} 
Since $x = 0$ is a root of the Stirling polynomial $S_{p^B,\cO_v}(x)$, the roots of the original $f_v(z)$
\index{Stirling polynomial!for $\cO_v$}
remain roots of the new one, and (\ref{WContainmentF}) still holds.  

Thus $H_v = f_v^{-1}(D(0,1)) \cap E_v$ is move-prepared relative to $B(a_1,r_1), \ldots, B(a_D,r_D)$.
\index{move-prepared}
\end{proof}

A second consequence of the proof of Theorem \ref{CompactThm} is

\begin{corollary} \label{EvHvExtraPtsCor} 
In Theorem $\ref{CompactThm}$,
by choosing $N_v$ appropriately large and divisible, 
we can arrange that 
$E_v = \bigcup_{\ell = 1}^D \big(B(a_\ell,r_\ell) \cap \cC_v(F_{w_\ell})\big)$
and $H_v = \bigcup_{h=1}^N \big(B(\theta_h,\rho_h) \cap \cC_v(F_{u_h})\big)$ 
have the property that there is a point 
$\wbar_\ell \in \big(B(a_\ell,r_\ell) \cap \cC_v(F_{w_\ell})\big) \backslash H_v$
for each $\ell = 1, \ldots, D$.
\end{corollary}  

\begin{proof}  This can be seen by tracing through the proof of Theorem \ref{CompactThm}.
However, a simple modification at the end of the proof gives the claim directly.  

Applying Theorem in its stated form,
given $\varepsilon > 0$ and $K_v$-simple decomposition 
\index{$K_v$-simple!decomposition}  
$E_v = \bigcup_{\ell = 1}^D \big(B(a_\ell,r_\ell) \cap \cC_v(F_{w_\ell})\big)$, 
there is a $K_v$-simple set $\tE_v \subset E_v$ such that for each $0 < \beta_v \in \QQ$
\index{$K_v$-simple!set}  
and each $K_v$-symmetric probability vector $\vs \in \cP^m(\QQ)$, 
\index{$K_v$-symmetric!probability vector}
assertions (A) and (B) in the Theorem hold.  

Fixing $0 < \beta_v < \QQ$ and $\vs$, write $\beta_v = \beta_v^{\prime} + \beta_v^{\prime \prime}$  
with $0 < \beta_v^{\prime}, \beta_v^{\prime \prime} \in \QQ$.  Applying the Theorem with 
$\beta_v$ replaced by $\beta_v^{\prime}$, there is an integer $N_v^{\prime} \ge 1$ such that 
for each positive integer $N$ divisible by $N_v^{\prime}$, 
there is an $(\fX,\vs)$-function\index{$(\fX,\vs)$-function} 
$f_v^{\prime} \in K_v(\cC_v)$ of degree $N$ such that assertion (B) holds for 
$\beta_v^{\prime}$: in particular, for each $x_i \in \fX$  
\begin{equation*}
\Lambda_{x_i}(f_v^{\prime},\vs) \ = \ \Lambda_{x_i}(\tE_v,\vs) + \beta_v^{\prime} \ ,
\end{equation*}
the zeros $\theta_1^{\prime}, \ldots, \theta_N^{\prime}$ of $f_v^{\prime}$
are distinct and belong to $E_v$, 
$(f_v^{\prime})^{-1}(D(0,1)) = \bigcup_{h = 1}^N B(\theta_h^{\prime},\rho_h^{\prime})$ where  
the balls $B(\theta_h^{\prime},\rho_h^{\prime})$ are pairwise disjoint and contained in 
$\bigcup_{\ell = 1}^D B(a_\ell,r_\ell)$,   
and the set $H_v^{\prime} := (f_v^{\prime})^{-1}(D(0,1)) \cap E_v$ has a $K_v$-simple decomposition 
\index{$K_v$-simple!decomposition}  
$H_v^{\prime} = \bigcup_{h=1}^N \big(B(\theta_h^{\prime},\rho_h^{\prime}) \cap \cC_v(F_{u_h^{\prime}}) \big)$ 
compatible with the decomposition 
$E_v = \bigcup_{\ell = 1}^D \big(B(a_\ell,r_\ell) \cap \cC_v(F_{w_\ell})\big)$.  
For each $h$, $\rho_h^{\prime}$ belongs to $|F_{u_h^{\prime}}^{\times}|_v$ and $f_v^{\prime}$
induces an $F_{u_h^{\prime}}$-rational scaled isometry from $B(\theta_h^{\prime},\rho_h^{\prime})$
\index{scaled isometry} 
onto $D(0,1)$.

Fix a positive integer $N_v^{\prime \prime}$ such that  
$N_v^{\prime \prime} \cdot \beta_v^{\prime \prime} \in \ZZ$, 
and put $N_v = N_v^{\prime}  N_v^{\prime \prime}$.  Given an integer $N$ divisible by $N_v$, 
let $f_v^{\prime}(z)$ be as above, and put 
\begin{equation*}
f_v(z) \ = \ \pi_v^{-N \beta_v^{\prime \prime}} \cdot f_v^{\prime}(z) 
\end{equation*} 
where $\pi_v$ is a uniformizer for the maximal ideal of $\cO_v$.  Clearly  
$\Lambda_{x_i}(f_v^{\prime},\vs) \ = \ \Lambda_{x_i}(\tE_v,\vs) + \beta_v$ for each $x_i \in \fX$.
The zeros $\theta_1, \ldots, \theta_N$ of $f_v(z)$ are the same as the zeros 
$\theta_1^{\prime}, \ldots, \theta_N^{\prime}$ of $f_v^{\prime}$, 
so they are distinct and belong to $E_v$.  

For each $h$, 
put $\rho_h = |\pi_v^{N \beta_v^{\prime \prime}}|_v \cdot \rho_h^{\prime} < \rho_h^{\prime}$ 
and put $F_{u_h} = F_{u_h^{\prime}}$.  
Since $\pi_v^{-N \beta_v^{\prime \prime}} \in K_v^{\times}$ 
and $\rho_h^{\prime} \in |F_{u_h^{\prime}}^{\times}|_v$,
we have $\rho_h \in |F_{u_h}^{\times}|_v$ 
and $B(\theta_h,\rho_h) \subsetneq B(\theta_h^{\prime},\rho_h^{\prime})$.  
Since $f_v^{\prime} : B(\theta_h^{\prime},\rho_h^{\prime}) \rightarrow D(0,1)$ is a scaled isometry,
\index{scaled isometry} 
it follows that 
\begin{equation*}
f_v^{-1}(D(0,1)) \ = \ \bigcup_{h=1}^N B(\theta_h,\rho_h) \ \subset \ 
\bigcup_{h=1}^N B(\theta_h^{\prime},\rho_h^{\prime}) \ = \ (f_v^{\prime})^{-1}(D(0,1)) \ ,
\end{equation*} 
and $H_v := f_v^{-1}(D(0,1)) \cap E_v$
has the $K_v$-simple decomposition 
\index{$K_v$-simple!decomposition}  
\begin{equation*}
H_v \ = \ \bigcup_{h=1}^N \big(B(\theta_h,\rho_h) \cap \cC_v(F_{u_h})\big)
\end{equation*}
compatible with the $K_v$-simple decomposition 
\index{$K_v$-simple!decomposition!compatible with another decomposition}   
$E_v = \bigcup_{\ell = 1}^D \big(B(a_\ell,r_\ell) \cap \cC_v(F_{w_\ell})\big)$. 

Now fix $\ell$, and take any $1 \le h \le N$ with
$B(\theta_h,\rho_h) \subset B(a_\ell,r_\ell)$.
Then $F_{u_h} = F_{u_h^{\prime}} = F_{w_\ell}$, 
$\theta_h = \theta_h^{\prime} \in B(a_\ell,r_\ell) \cap \cC_v(F_{w_\ell})$,
$B(\theta_h,\rho_h) \subsetneq B(\theta_h,\rho_h^{\prime}) \subset B(a_\ell,r_\ell)$,
and $\rho_h^{\prime} \in |F_{w_\ell}^{\times}|_v$.  
Since $f_v^{\prime} : B(\theta_h,\rho_h^{\prime}) \rightarrow D(0,1)$ 
is an $F_{w_\ell}$-rational scaled isometry, 
\index{scaled isometry}
there infinitely many points $\wbar  \in  \cC_v(F_{w_\ell})$ 
with $\|\theta_h,\wbar\|_v = \rho_h^{\prime}$.  Such points belong to 
$B(\theta_h,\rho_h^{\prime}) \cap \cC_v(F_{w_\ell})$ but not $B(\theta_h,\rho_h)$,
so they are not in $H_v$. 

Thus for each $\ell$ there is a point 
$\wbar_\ell \in \big(B(a_\ell,r_\ell) \cap \cC_v(F_{w_\ell})\big) \backslash H_v$.  
\end{proof} 
\index{initial approximating functions $f_v(z)$!construction when $K_v$ nonarchimedean!consequences of the proof|)}
\index{initial approximating functions $f_v(z)$!for nonarchimedean $K_v$-simple sets|)}

%% file: NewFSZChap7.tex
\chapter{ The Global Patching Construction} \label{Chap7}
 
In this section we will give the global patching construction for the proof of Theorem 4.2.  
\index{patching argument!global} 
This argument manages the patching process 
in such a way that the final patched function is $K$-rational.  
Part of the argument specifying the order in which the coefficients are patched, 
\index{coefficients $A_{v,ij}$}\index{initial patching functions|see{patching functions, initial $G_v^{(0)}(z)$}}
and the way the target coefficients are chosen.  
In expansions of functions, all coefficients are 
with respect to the $L$-rational basis, constructed in \S\ref{Chap3}.\ref{LRationalBasisSection}.
\index{basis!$L$-rational}    

The inputs to the global patching argument 
are the construction of the initial approximating functions, 
\index{patching argument!global}\index{initial approximating functions $f_v(z)$}
carried out in Chapters \ref{Chap5} and  \ref{Chap6} above, 
and the local patching constructions given in Chapters \ref{Chap8} -- \ref{Chap11} below.
\index{patching argument!local} 
The global and local patching constructions are largely independent;   
we have chosen to present the global argument first  
in order to provide the reader with an overview of the proof.  
Nonetheless, in order to understand some aspects of the global patching construction 
(in particular the need for patching the coefficients in bands 
\index{band!coefficients patched by bands}
and the reason for using different patching coefficient bounds for high, middle, 
and low-order coefficients), 
\index{coefficients $A_{v,ij}$!high-order}\index{coefficients $A_{v,ij}$!middle}\index{coefficients $A_{v,ij}$!low-order}
it is necessary to be acquainted with the local constructions. 
Therefore we encourage the reader to examine the local patching constructions
in parallel with the global one.  The patching argument for the 
non-archimedean $\RL$-domain case when $\Char(K)= 0$, given in Chapter \ref{Chap10}, is the easiest
\index{$\RL$-domain}  
and will shed light on the issues above.  

For the convenience of the reader, we restate Theorem \ref{aT1-B}.
The notion of a $K_v$-simple set is defined in Definition \ref{KvSimple}. 
\index{$K_v$-simple!set}  
\smallskip

\noindent{\bf Theorem \ref{aT1-B}} (FSZ with Local Rationality Conditions,
for $K_v$-simple sets). \label{Thm aT1-B-A}

{\em
Let $K$ be a global field, 
and let $\cC/K$ be a smooth, geometrically integral, projective curve.
Let $\fX = \{x_1, \ldots, x_m\} \subset \cC(\tK)$ 
be a finite set of points stable under $\Aut(\tK/K)$, and let
$\EE = \prod_v E_v \subset \prod_v \cC_v(\CC_v)$ be a 
an adelic set compatible with $\fX$,\index{compatible with $\fX$} such that each $E_v$ is stable under $\Aut_c(\CC_v/K_v)$.   
Let $S = S_K \subset \cM_K$ be a finite set of places $v \in \cM_K$ 
containing all archimedean $v$.      

Assume that $\gamma(\EE,\fX) > 1$, and that 

$(A)$ $E_v$ is $K_v$-simple for each $v \in S$,
\index{$K_v$-simple!set}  

$(B)$ $E_v$ is $\fX$-trivial for each $v \notin S$.
\index{$\fX$-trivial}  

\noindent{Then} 
there are infinitely many points $\alpha \in \cC(K^{\sep})$ 
such that for each $v \in \cM_K$, 
 the $\Aut(\tK/K)$-conjugates of $\alpha$ all belong to $E_v$. 
}

\smallskip
 
In Chapter \ref{Chap4} we have reduced Theorems 
\ref{aT1}, \ref{aT1-A1}, \ref{aT1-A}, \ref{FSZi}, \ref{FSZii}, \ref{aT1-B1}, and \ref{aT1-B2}
to Theorem \ref{aT1-B}.  In Chapters \ref{Chap5} and \ref{Chap6} we have constructed the 
initial approximating functions.\index{initial approximating functions $f_v(z)$} 

In this Chapter we prove Theorem \ref{aT1-B}, assuming 
the local patching constructions given in Chapters \ref{Chap8} -- \ref{Chap11}.  
\index{patching argument!local}   
In \S\ref{Chap7}.\ref{StrongApproxThmSection} -- \S\ref{Chap7}.\ref{SemiLocalSection} below 
we discuss some preliminaries:  
an adelic version of the Strong Approximation theorem,
\index{Strong Approximation theorem!adelic}
the existence of a dense set of subunits, and the semi-local theory. 
\index{subunit}
In \S\ref{Chap7}.\ref{Char0Section} we prove Theorem \ref{aT1-B} when $\Char(K) = 0$. 
First, we specify the patching parameters,
\index{patching parameters}    
then we construct the functions used to initiate the patching
process, and finally we carry out the global patching construction.  
\index{patching argument!global} 
In \S\ref{Chap7}.\ref{CharPSection} 
we prove Theorem \ref{aT1-B} when $\Char(K) = p > 0$.

We use the conventions concerning notation and absolute values from \S\ref{Chap3}.\ref{NotationSection}.
We assume familiarity with the theory of Green's functions  
from \S\ref{Chap3}.\ref{CompactGreenSection}, \S\ref{Chap3}.\ref{UpperGreenSection} 
\index{Green's function!of a compact set}
\index{Green's function!upper}
(or from \cite{RR1}, \S 3.2, \S 4.4), 
and with Green's matrices and the Cantor capacity
\index{capacity!Cantor capacity}  
from \S\ref{Chap3}.\ref{CantorCapacitySection} (or from \cite{RR1},\S 5.3).  

The Green's matrix and the Cantor capacity only depend on values of the Green's functions 
\index{Green's function}
\index{capacity!Cantor capacity}
outside $E_v$.  For a compact set $E_v$, if $x \notin E_v$, the upper Green's function $\Gbar(z,x;E_v)$ 
coincides with the Green's function $G(z,x;E_v)$ from \S\ref{Chap3}.\ref{CompactGreenSection} for all $z$,
and it coincides with the (lower) Green's function $G(z,x;E_v)$ studied in (\cite{RR1}) 
\index{Green's function!lower}
for all $z \notin E_v$.  Likewise, for an $\fX$-trivial set, or more generally for any algebraically
\index{$\fX$-trivial}
 capacitable set, the upper Green's function $\Gbar(z,x;E_v)$ coincides with the lower Green's function 
studied in (\cite{RR1}) for all $z \notin E_v$ (see \cite{RR1}, Theorem 4.4.4).
Hence for sets $\EE = \prod_v E_v$ meeting the conditions of Theorem \ref{aT1-B}, 
the upper Green's matrix $\Gammabar(\EE,\fX)$ and inner Cantor capacity $\gammabar(\EE,\fX)$ 
\index{capacity!inner Cantor capacity}
from \S\ref{Chap4}.\ref{CantorCapacitySection}   
coincide  with the Green's matrix $\Gamma(\EE,\fX)$ and Cantor capacity $\gamma(\EE,\fX)$
\index{Green's matrix!global} 
\index{capacity!Cantor capacity}
from (\cite{RR1},\ \S5.3).

For this reason, for the remainder of this chapter, we will drop the ``bar'' from 
$\Gbar(z,x;E_v)$, $\Gammabar(\EE,\fX)$, and $\gammabar(\EE,\fX)$ 
and simply write $G(z,x;E_v)$, $\Gamma(\EE,\fX)$ and $\gamma(\EE,\fX)$.
\index{Green's matrix!global}

\vskip .07 in
Let $S_K = S$ be the set of places of $K$ in Theorem \ref{aT1-B}, 
and let $\hS_K = \hS$ be the set of all places of $K$ 
where any of the following conditions holds:

\begin{enumerate} 
  
  \item  $v \in S_K$; in particular if 
  \begin{itemize} 
    \item  $v$ is archimedean;  or
    \item  $\cC$ has bad reduction at $v$ (with respect to the model $\fC$ 
                 determined by the given projective embedding of $\cC$); or
    \item  the points in $\fX$ do not specialize to distinct points $\pmod v$; or
    \item  $E_v$ is not $\fX$-trivial; 
\index{$\fX$-trivial}
  \end{itemize}
\vskip -.5in
\begin{equation} \label{hS_KList} \end{equation}
\vskip .05in  
  \item  $\cC$ has good reduction at $v$, and one or more of the uniformizing 
  \index{good reduction}
              parameters $g_{x_i}(z)$ fails to specialize to 
              a well-defined, nonconstant function 
              on the fibre $\fC \pmod v = \fC \times_{\cO_K} k_v$ 
              (in the classical terminology, ``has bad reduction at $v$'');                
  \item  $\cC$ has good reduction at $v$, and one or more of the basis functions 
  \index{good reduction}
              $\varphi_{ij}(z)$ and $\varphi_{\lambda}$ from 
              Section \ref{Chap3}.\ref{LRationalBasisSection} has bad reduction at $v$.            
\end{enumerate}
Note that although there are infinitely many basis functions, 
only finitely many places are affected by condition (3), 
because the basis functions belong to a multiplicatively finitely generated set.
\index{$L$-rational basis!multiplicatively finitely generated} 
Thus $\hS_K$ is finite.  

Put $L = K(\fX)$ and $\hS_L$ be the set of places of $L$ above $\hS_K$.  
Write $\EE_K = \EE$ and let $\EE_L = \prod_{w \in \cM_L} E_w$ 
be the set obtained from $\EE_K$ by base change, identifying $\CC_w$ with $\CC_v$ 
and putting $E_w = E_v$ if $w|v$.    
For each $w \notin \hS_L$, $E_w$ is $\fX$-trivial and the $g_{x_i}(z)$
\index{$\fX$-trivial}
have good reduction, so the local Green's matrix  $\Gamma(E_w,\fX)$ 
\index{Green's matrix!local}
\index{good reduction}
is the zero matrix.  Thus 
\begin{equation} \label{GreensMatrixSum}
\Gamma(\EE_K,\fX) \ = \ \frac{1}{[L:K]} \Gamma(\EE_L,\fX) \ = \ \frac{1}{[L:K]}
             \sum_{w \in \hS_L} \Gamma(E_w,\fX) \log(q_w) \ .
\end{equation}

In Theorem \ref{aT1-B}, our hypothesis that $\gamma(\EE,\fX) > 1$ 
is equivalent to $\Gamma(\EE,\fX)$ being negative definite.  
\index{Green's matrix!global}\index{Green's matrix!negative definite}
Let $\{\tE_v\}_{v \in \hS_K}$ be another collection of sets 
\index{Green's matrix!global}
for which $\widetilde{\EE} := \prod_{v \in \hS_K} \tE_v \times \prod_{v \notin \hS_K} E_v$
is $K$-rational and compatible with $\fX$.{compatible with $\fX$}  
By continuity, there are numbers $\varepsilon_v > 0$ for $v \in \hS_K$
such that $\Gamma(\widetilde{\EE},\fX)$ is also negative definite,\index{Green's matrix!negative definite} 
provided that for each $v \in \hS_K$ 
\begin{equation} \label{NegDefCloseness}
\left\{ \begin{array}{ll}
       |G(x_j,x_i;\tE_v)-G(x_j,x_i;E_v)| \ < \ \varepsilon_v &
                  \text{for all $i \ne j$ \ ,} \\
       |V_{x_i}(\tE_v)-V_{x_i}(E_v)| \ < \ \varepsilon_v &
                  \text{for all $i$ \ .}  
        \end{array}  \right. 
\end{equation}
\index{Green's function}
\index{Robin constant}

\section{ The Uniform Strong Approximation Theorem} \label{StrongApproxThmSection} 

Let $F$ be a global field.  Write $\cM_F$ for the set of all places of $F$,
and let $|x|_u$ be the absolute value associated to
$u \in \cM_F$, normalized as in \S{\ref{Chap3}}.\ref{NotationSection}.  
Let $\AA_F$ and $\JJ_F$ denote the adele ring and idele group of $F$,\index{adele ring}\index{idele group} 
respectively, and write $b = (b_u)_{u \in \cM_F}$
for an element of $\AA_F$ or $\JJ_F$.  Let
\begin{equation*}
\|b\|_F \ = \ \prod_{u \in \cM_F} |b_u|_u^{D_u}
\end{equation*}
denote the ``size'' of $b$.\index{size of an adele|ii}  Recall that $D_u=2$ 
if $u$ is archimedean and $K_u \cong \CC$, and $D_u=1$ otherwise.  

The following version of the Strong Approximation theorem is well known, 
but there seems to be no convenient reference for it.
\index{Strong Approximation theorem!adelic|ii}
   
\begin {lemma} \label{CStrong1}
Let $F$ be a global field.  There is a constant $B_F$ depending only on $F$
such that for each $b \in \JJ_F$  with  $\|b\|_F \ge B_F$,  
given any $c \in \AA_F$, there is an  $f \in F$ such that
\begin{equation*}
|f-c_u|_u \ \le |b_u|_u  \quad  \text{for all $u \in \cM_F$.}
\end{equation*}
\end{lemma}

\begin{proof}  By the Lemma on (\cite{Ca-F}, p.66), there is a constant $A(F) > 0$ 
such that for any $a \in \JJ_F$ with $\|a\|_F > A(F)$, 
there is a $\beta \in F^{\times}$ satisfying $|\beta|_u \le |a_u|_u$ for all $u \in \cM_F$.

Given any idele $y = (y_u)_{u \in \cM_F} \in \JJ_F$, put 
\begin{equation*} 
V(y) \ = \ \{x = (x_u)_{u \in \cM_F} \in \AA_F : |x|_u \le |y_u|_u \} \ .
\end{equation*}
By Corollary 1 of (\cite{Ca-F}, p.65), there is a $d \in \JJ_F$ such that $V(d)$
contains a fundamental domain for $\AA_F/F$:  that is, $\AA_F = V(d) + F$, 
where we view $F$ as embedded on the diagonal in $\AA_F$.  Let $D(F) = \|d\|_F$. 

Take $B_F = A(F) \cdot D(F)$.  Suppose $b \in \JJ_F$ is an idele with $\|b\|_F \ge B_F$.
Put $a = b \cdot d^{-1}$;  then $\|a\|_F \ge A(F)$.  Let $\beta \in F^{\times}$ be such that 
$|\beta|_u \le |a_u|_u$ for each $u \in \cM_F$;  then $|\beta d_u|_u \le |b_u|_v$ for each $u$.  
Since $V(d)$ contains a fundamental domain for $\AA_F$, so does $V(\beta d)$:
indeed,  
\begin{eqnarray*}
\AA_F & = & \beta \cdot \AA_F \ = \ \beta \cdot (V(d) + F) \\
   & = & \beta \cdot V(d) + \beta \cdot F \ = \ V(\beta d) + F \ .
\end{eqnarray*} 
Since $V(\beta d) \subseteq V(b)$, it follows that $\AA_F = V(b) + F$ as well.  

Now, take any adele $c = (c_u)_{u \in \cM_F} \in \AA_F$.  Let $x \in V(b)$ and $f \in F$ 
be such that $a = x + f$.  Then for each $u \in \cM_F$,
\begin{equation*}
|f - c_u|_u \ = \ |c_u - f|_u \ = \ |x_u|_u \ \le \ |b_u|_u 
\end{equation*}
as desired.
\end{proof}

\vskip .1 in 
Restricting to a finite set of places, we have 

\begin{corollary} \label{CStrong2}
Let $F$ be a global field, and let $\hS_F \subset \cM_F$ be a nonempty finite set of places.
Then there is a constant $B(\hS_F)$ with the following property.
For any set of numbers 
\begin{equation*}
\{ 0 < Q_u \in \RR : u \in \hS_F \}
\end{equation*}
such that $\prod_{u \in \hS_F} Q_u^{D_u} > B(\hS_F)$,
and any collection of elements $c_u \in F_u$ for $u \in \hS_F$,
there is an $f \in F$ satisfying 
\begin{equation*}
\left\{ \begin{array}{cl}
     |f-c_u|_u \le Q_u & \text{for all $u \in \hS_F$\ ,} \\
     |f|_u \le 1 &       \text{for all $u \notin \hS_F$\ .}
         \end{array} \right.
\end{equation*}
\end{corollary}

\begin{proof}  Let $B_F$ be the constant from Lemma \ref{CStrong1}, and put
\begin{equation*}
B(\hS_F) \ = \
   B_F \cdot (\prod_{\text{nonarchimedean $u \in \hS_F$}} q_u) \ ,
\end{equation*}
where $q_u$ is the order of the residue field at $u$.

Suppose $\prod_{u \in \hS_F} Q_u^{D_u} \ge B(\hS_F)$.
For each archimedean $u \in \hS_F$, 
there is a $b_u \in F_u$ with $|b_u|_u = Q_u$;
for each nonarchimedean $u \in \hS_F$ there is a $b_u \in F_u$ with
\begin{equation*}
q_u^{-1}Q_u \ < \ |b_u| \ \le \ Q_u \ .
\end{equation*}
For each $u \notin \hS_F$, put $b_u = 1$, and let 
$b = (b_u)_{u \in \cM_F} \in \JJ_F$.   Then $\|b\|_F \ge B_F$.
Let $c_u \in F_u$ for $u \in \hS_F$ be 
the given elements, and put $c_u = 0$ for $u \notin \hS_F$.   

By Lemma \ref{CStrong1} there is an $f \in F$ with $|f-c_u| \le |b_u|_u$
for all $u \in \cM_F$.
\end{proof}

\vskip .1 in
Now consider a global field $K$ and a nonempty finite set of places $\hS_K$ of $K$.
For any finite extension $F/K$, let $\hS_F$ be the set of places of $F$
above $\hS_K$.  The following extension of Corollary \ref{CStrong2}
will be used in the global patching process.
\index{patching argument!global} 

\begin{proposition} \label{CStrong3} 
{\bf $($Uniform Strong Approximation Theorem$)$} 
\index{Strong Approximation theorem!uniform|ii}
Let $K$ be a global field, and let $\hS_K$ be a nonempty finite
set of places of $K$.  Let $H/K$ be a finite normal extension. 
Then there is a constant $C_H(\hS_K) > 0$ with the following property.
Let $\{0 < Q_w \in \RR : w \in \hS_H\}$ be a $K$-symmetric set of numbers
\index{$K$-symmetric!set of numbers}
with $\prod_{w \in \hS_H} Q_w^{D_w} > C_H(\hS_K)$.  
Let $F$ be any field with $K \subseteq F \subseteq H$, and let 
$\{c_u \in F_u: u \in \hS_F\}$ be any set of elements. 
Then, regarding $c_u$ as an element of $H_w$ for each $w|u$,
there is an $f \in F$ satisfying
\begin{equation} \label{CFS3}
\left\{ \begin{array}{cl}
     |f-c_u|_w \le Q_w & \text{for each $w \in \hS_H$\ ,} \\
     |f|_w \le 1 &       \text{for each $w \notin \hS_H$\ .}
         \end{array} \right.
\end{equation}
\end{proposition}

\noindent{\bf Remark:}  If $H/K$ is a finite but not normal, 
the assumption that $\{Q_w\}_{w \in \hS_H}$ is $K$-symmetric 
\index{$K$-symmetric!set of numbers}
can be replaced by
the requirement that there is a set of numbers 
$\{ 0 < Q_v \in \RR : v \in \hS_K \}$ such that $Q_w^{D_w} = Q_v^{D_v [H_w:K_v]}$ 
whenever $w|v$.

\vskip .1 in
\begin{proof}  Define 
\begin{equation} \label{CFS1}
C_H(\hS_K) \ = \ \max_{K \subseteq F \subseteq H} (B(\hS_F)^{[H:F]}) \ .
\end{equation}
Let $\{Q_w\}_{w \in \hS_L}$, $F$, and $\{ c_u\}_{u \in \hS_F}$ 
be as in the Proposition.
For each $u \in \hS_F$ and each $w \in \hS_L$ with $w|u$, 
define $Q_u$ by $Q_u^{D_u[L_w:K_v]} = Q_w^{D_w}$.  
Note that $Q_u$ is well-defined since $\{Q_w\}_{w \in \hS_H}$ is $K$-symmetric. 
\index{$K$-symmetric!set of numbers}
By our normalization of absolute values,
if $x \in F_u$ is regarded as an element of $H_w$
then $|x|_w^{D_w} = (|x|_u^{D_u})^{[H_w:F_u]}$,
so $|x|_u \le Q_u$ iff $|x|_w \le Q_w$.

Since $\sum_{w|u} [H_w:F_u] = [H:F]$ we have
\begin{equation*}
(\prod_{u \in \hS_F} Q_u^{D_u})^{[H:F]}  
\ = \ \prod_{w \in \hS_L}  Q_w^{D_w} \ \ge \ C_L(\hS_K) \ \ge \ B(\hS_F)^{[H:F]} \ ,
\end{equation*}
so $\prod_{u \in \hS_F} Q_u^{D_u} \ge B(\hS_F)$.  
Let $f \in F$ be the element given by Lemma \ref{CStrong2}.
For each $u \in \hS_F$, we have $|f-c_u|_u \le Q_u$,  while for
each $u \notin \hS_F$, we have $|f|_u \le 1$.  Passing to the
extension $H/K$, for each $w \in \hS_H$
we have $|f - c_u|_w \le Q_w$ if $w|u$,
while for each $w \notin \hS_H$ we have $|f|_w \le 1$.
\end{proof}
  
 
\section{ $S$-units and $S$-subunits } \label{SUnitSection} 

Let $F$ be a global field, and let $\hS_F$ be a nonempty finite set of places
of $F$ containing all the archimedean places.      
Write $\hS_F = \hS_{F,\infty} \bigcup \hS_{F,0}$
where $\hS_{F,\infty}$ is the subset of archimedean places, 
and $\hS_{F,0}$ is the subset of nonarchimedean places;  
here, if $F$ is a function field, $\hS_{F,\infty}$ is empty.
The set of $\hS_F$-units is the group\index{$S$-unit|ii} 
\begin{equation*}
\cO_{F,\hS_F}^{\times}  \ = \ \{ f \in F^{\times} : 
            \ord_u(f) = 0 \ \text{for all $u \notin \hS_F$} \} \ .
\end{equation*}
By the $S$-unit theorem 
\begin{equation*}
\cO_{F,\hS_F}^{\times} \ \cong \ \ZZ/d \ZZ \times \ZZ^{\#(\hS_F)-1}
\end{equation*}
where $d = d_F$ is the number of roots of unity\index{roots of unity} in $F$.
Furthermore, the homomorphism $\log_{F,\hS_F} : \cO_{F,\hS_F}^{\times} \rightarrow \RR^{\#(S_F)}$
\begin{equation*}
\log_{F,\hS_F}(f) \ = \ (\log_u(|f|_u))_{u \in \hS_F}
\end{equation*}
maps $\cO_{F,\hS_F}^{\times}$ onto a $\ZZ$-lattice\index{lattice} which spans the hyperplane
\begin{equation*}
\cH_{\hS_F} \ = \ \{ \vt \in \RR^{\#(\hS_F)} : \sum_{u \in \hS_F} t_u \log(q_u) = 0 \}
    \ \subset \ \RR^{\#(\hS_F)} \ .
\end{equation*}
The kernel of $\log_{F,\hS}$ is the group of roots of unity $\mu_d$ in $F$.

Note that if $u \in \hS_{F,0}$ and $f \in F^{\times}$,
then $-\log_u(|f|_u) = \ord_u(f) \in \ZZ$.
The $S$-unit theorem\index{$S$-unit Theorem} therefore implies

\begin{proposition} \label{CPSUnit1} Let $F$ be a global field, 
and let $\hS_F$ be a finite set of places of $F$ containing $\hS_{F,\infty}$.  
Suppose $\vt \in \RR^{\#(\hS_F)}$ satisfies
\begin{equation*}
\sum_{u \in \hS_F} t_u \log(q_u) \ = \ 0 \ ,
\end{equation*}
with $t_u \in \QQ$ for each $u \in \hS_{F,0}$.
Then for each $\eta > 0$, there are an an integer $m_0 > 0$
and an $\hS_F$-unit $f_0 \in F^{\times}$ such that

$(1)$  $ \frac{1}{m_0} \log_u(|f_0|_u) = t_u$ for each $u \in \hS_{F,0}$;

$(2)$  $|\frac{1}{m_0} \log_u(|f_0|_u) - t_u| < \eta$ for each
$u \in \hS_{F,\infty}$.

\end{proposition}

When $F$ is a function field, this can be reformulated as follows:  

\begin{proposition} \label{CPSUnit2FF}
Let $F$ be a function field, 
and let $\hS_F$ be a finite set of places of $F$. 
Let $\{c_u \in F_u^{\times} : u \in \hS_F\}$ 
be a collection of elements such that 
\begin{equation*}
\sum_{u \in \hS_F} \log_u(|c_u|_u) \log(q_u) \ = \ 0 \ .
\end{equation*}
Then there are an $S_F$-unit $f \in \cO_{F,\hS_F}^{\times}$
and an integer $n_0 > 0$ such that $|c_u^{n_0}|_u = |f|_u$ for each $u \in \hS_{F,0}$.
\end{proposition} 

\begin{proof} 
Take $t_u = \log_u(|c_u|_u)$ for each $u \in \hS_F$, and let 
 $n_0$ and $f_0$ and be the integer $m_0$ and $\hS_F$ unit given by 
Proposition \ref{CPSUnit1}. 
\end{proof} 

When $F$ is a number field, there is a stronger version of Proposition \ref{CPSUnit2FF}   
using the concept of a subunit, introduced by Cantor\index{Cantor, David} (\cite{Can3}):
\index{subunit|ii}\index{$S$-subunit|ii} 

\begin{definition}  An $\hS_F$-subunit is a vector 
$\vec{\varepsilon} = (\varepsilon_u)_{u|\infty} \in \oplus_{u \in \hS_{F,\infty}} F_u^{\times}$  
for which there are an integer $n_0 > 0$ and an $\hS_F$-unit $f \in F^{\times}$ with
$\varepsilon_u^{n_0} = f$ for each $u \in \hS_{F,\infty}$.
\end{definition}

If $F$ is a number field, the group of roots of unity is dense in
the unit circle in $\CC^{\times}$, and the units $\{ \pm 1 \}$ represent both connected
components of $\RR^{\times}$.  Hence we have 

\begin{proposition} \label{CPSUnit2}
Let $F$ be a number field,
and let $\hS_F$ be a finite set of places of $F$ containing $\hS_{F,\infty}$. 
Let $\{c_u \in F_u^{\times} : u \in \hS_F\}$ 
be a collection of elements such that 
\begin{equation*}
\sum_{u \in \hS_F} \log_u(|c_u|_u) \log(q_u) \ = \ 0 \ .
\end{equation*}
Then for each $\delta > 0$, there are an $S_F$-unit $f \in \cO_{F,\hS_F}^{\times}$,
an integer $n_0 > 0$, and an $S_F$-subunit 
$\vec{\varepsilon} \in \oplus_{u \in \hS_{F,\infty}} F_u^{\times}$ such that

$(1)$ \  $|c_u^{n_0}|_u = |f|_u$ for each $u \in \hS_{F,0};$

$(2)$ \ $|\varepsilon_u - c_u|_u < \delta$ and $\varepsilon_u^{n_0} = f$ for each $u \in \hS_{F,\infty}$.
\end{proposition}

\begin{proof} Apply Proposition \ref{CPSUnit1}, taking $t_u = \log_u(|c_u|_u)$ for each $u \in \hS_F$.
Let $\eta > 0$ be small enough that for each $u \in \hS_{F,\infty}$,  
if $|x - t_u| < \eta$ then $|\exp(x) - |c_u|_u| < \delta/2$.
Let $m_0$ and $f_0$ and be the positive integer and $\hS_F$ unit given by 
Proposition \ref{CPSUnit1}.

For each $u \in \hS_{F,\infty}$ we have $||f_0|_u-|c_u|_u| < \delta/2$, 
and there is a root of unity $\omega_u \in F_u^{\times}$
such that $|\omega_u  f_0 - c_u| < \delta$;  put $\varepsilon_u = \omega_u f_0$, 
and let $\vec{\varepsilon} = (\varepsilon_u)_{u|\infty}$.  
Let $m_1$ be the least common multiple of the orders of the $\omega_u$;  
put $n_0 = m_0 m_1$ and take $f = f_0^{m_1}$. 
Clearly (1) and (2) hold for this $\vec{\varepsilon}$, $n_0$ and $f$.
\end{proof}  


\section{ The Semi-local Theory } \label{SemiLocalSection}
\index{semi-local theory|(} 

Let $K$ be a global field, and let $H/K$ be a finite extension.
For each place $v$ of $K$,  
there is a canonical isomorphism of topological algebras 
\begin{equation} \label{CFQQ1} 
H \otimes_K K_v \ \cong \ \oplus_{w|v} H_w \ .
\end{equation}
(This isomorphism holds even when $H/K$ is not separable; see (\cite{RR1}, p.321).)  
Under this isomorphism $K_v \cong K \otimes_K K_v$ is identified with the
set of diagonal elements $(\kappa_v,\ldots,\kappa_v)$, $\kappa_v \in K_v$.
More generally, for any field $F$ with $K \subseteq F \subseteq H$, the algebra    
$F \otimes_K K_v \cong \oplus_{u|v} F_u$ embeds in $H \otimes_K K_v$ in such a
way that $\oplus_{u|v} h_u \in \oplus_{u|v} F_u$ is sent to the
quasi-diagonal element\index{quasi-diagonal element} $\oplus_{w|v} f_w \in \oplus_{w|v} H_w$
where  $f_w = h_u$ for each $w|u$.  

When $H/K$ is normal, the group $\Aut(H/K)$ acts on $H \otimes_K K_v$
through its action on $H$: for each $\sigma \in \Aut(H/K)$,  $f \in H$, and $\kappa_v \in K_v$,
we have $\sigma(f \otimes_K \kappa_v) = \sigma(f) \otimes \kappa_v$.  
When this is interpreted on the right side of (\ref{CFQQ1}), it says that 
$\sigma$ induces a permutation $w \mapsto \sigma(w)$
of the places $w|v$, and a canonical isomorphism 
$\tau_{\sigma,w} : H_w \rightarrow H_{\sigma(w)}$ for each $w$.  
That is, the action of $\sigma$ on $\oplus_{w|v} H_w$ is gotten by 
applying $\tau_{\sigma,w}$ to the $w$-coordinate,  
while permuting the coordinates so the $w$-coordinate goes to the 
$\sigma(w)$-coordinate.  Furthermore, $\Aut(H/K)$ acts transitively on the places $w$ 
over $v$.  (When $H/K$ is galois this is well-known. When $H/K$ is merely normal, let $H^{\sep}$
be the separable closure of $K$ in $H$.  Then $\Aut(H/K) \cong \Gal(H^{\sep}/K)$ 
acts transitively on the places $w_0$ of $H^{\sep}$ lying over $v$, and the assertion 
follows because there is a unique place $w$ of $H$ over each $w_0$ of $H^{\sep}$.)

When $H/K$ is galois, $K_v$ is the sub-algebra fixed by $\Gal(H/K)$;
more generally, $F \otimes_K K_v$ is the sub-algebra
fixed by $\Gal(H/F)$, for each $F$ with $K \subseteq F \subseteq H$.

\smallskip

We will now apply these facts in the context of Theorem \ref{aT1-B}.  
Let $K$ be the global field in Theorem \ref{aT1-B}, and put $L = K(\fX)$. 
If $K$ is a number field, take $H = L$;  if $K$ is a function field, put $H = L^{\sep}$.
Then $H/K$ is galois.  Since $\cC/K$ is geometrically integral,
\begin{equation} \label{CFQQ1A}
H \otimes_K K_v(\cC) \ \cong \ \oplus_{w|v} H_w(\cC) \ .
\end{equation}
If $K \subseteq F \subseteq H$, the algebra 
$F \otimes_K K_v(\cC)$ embeds quasi-diagonally in $H \otimes_K K_v(\cC)$.

For each $\sigma \in \Gal(H/K)$ and each $w$, 
the isomorphism $\tau_{\sigma,w} : H_w \rightarrow H_{\sigma(w)}$ induces
an isomorphism $\widehat{\tau}_{\sigma,w} : H_w(\cC) \rightarrow H_{\sigma(w)}(\cC)$ fixing $K_v(\cC)$.   
As before, the action of $\Gal(H/K)$  on $H \otimes_K K_v(\cC)$ 
can be interpreted as applying $\widehat{\tau}_{\sigma,w}$ to the $w$-component of (\ref{CFQQ1A}), 
for each $w$, while permuting the coordinates so the $w$-component goes to the 
$\sigma(w)$-component.  $K_v(\cC)$ is the sub-algebra fixed by $\Gal(H/K)$;
more generally, $F \otimes_K K_v(\cC)$ is the sub-algebra
fixed by $\Gal(H/F)$, for each $F$ with $K \subseteq F \subseteq H$.

\smallskip

Let $J$ be the number from the construction of the $L$-rational 
and $L^{\sep}$ rational bases\index{basis!$L$-rational}\index{basis!$L^{\sep}$-rational} 
in \S\ref{Chap3}.\ref{LRationalBasisSection}, 
and let $\Lambda_0 = \dim_{\tK}(\tGamma(\sum_{i=1}^m J(x_i))$ be the number of low-order basis elements.
Given a probability vector $\vs \in \cP^m(\QQ)$ with positive coordinates, 
and an integer $N$ such that $N \vs \in \ZZ^m$, write $N_i = N s_i$ for $i = 1, \ldots, m$.
Assume that $N$ is large enough that $N_i \ge J$ for each $i$.
Suppose we are given an $(\fX,\vs)$-function\index{$(\fX,\vs)$-function} 
$\phi(z) \in K(\cC)$ of degree $N$.

\smallskip
If $K$ is a number field (that is, if $\Char(K) = 0$), 
we can expand $\phi(z)$ using the $L$-rational basis as
\index{basis!$L$-rational} 
\begin{equation} \label{CFDF1}
\phi(z) \ = \ \sum_{i=1}^m \sum_{j=0}^{N_i-J+1}
                             a_{ij} \varphi_{i,N_i-j}(z)
           + \sum_{\lambda=1}^{\Lambda_0} a_{\lambda} \varphi_{\lambda}(z) .
\end{equation}
with the $a_{ij}, a_{\lambda} \in L$.  
Since $\phi(z)$ is $K$-rational, for each $\sigma \in \Gal(F/K)$ we have $\sigma(\phi)(z) = \phi(z)$.
Applying $\sigma$ to (\ref{CFDF1}), and 
recalling that $\sigma(\varphi_{ij}) = \varphi_{\sigma(i),j}$ and
$\sigma(\varphi_{\lambda}) = \varphi_{\lambda}$, we find that
$\sigma(a_{ij}) = a_{\sigma(i),j}$ for all $i$, $j$,
and $\sigma(a_{\lambda}) = a_{\lambda}$ for all $\lambda$.
Thus, the $a_{ij}$ are $K$-symmetric relative to the
\index{$K$-symmetric!set of numbers}
action of $\Gal(L/K)$ on the $x_i$, and for each $\sigma \in \Gal(F/K(x_i))$
we have  $\sigma(a_{ij}) = a_{\sigma(i),j} = a_{ij}$, so $a_{ij}$
belongs to $K(x_i)$.  Likewise, since each $\varphi_{\lambda}$ is $K$-rational, 
each $a_{\lambda}$ belongs to $K$.  

Similarly, let $w$ be a place of $L$ with $w|v$.  
If we are given an $(\fX,\vs)$ function
$\phi_w(z) \in L_w(\cC)$ of degree $N$, we can write
\begin{equation} \label{CFDF2}
\phi_w(z) \ = \ \sum_{i=1}^m \sum_{j=0}^{N_i-J-1}
                             a_{w,ij} \varphi_{i,N_i-j}(z)
           + \sum_{\lambda=1}^{\Lambda} a_{w,\lambda} \varphi_{\lambda}(z)  
\end{equation}
where each $a_{w,ij} \in L_w$ and each $a_{w, \lambda} \in K_v$.  
In this context, we have:\index{semi-local theory!for number fields}  

\begin{proposition} \label{CPPr1}  Suppose $K$ is a number field,  
and let $v$ be a place of $K$. For each place $w$ of $L$ with $w|v$, 
let an $(\fX,\vs)$-function\index{$(\fX,\vs)$-function} $\phi_w(z) \in L_w(\cC)$ be given.  
Then the following are equivalent$:$

$(1)$  There is a $\phi_v(z) \in K_v(\cC)$
       such that $\phi_w(z) = \phi_v(z)$ for all $w|v$.

$(2)$  $\oplus_{w|v} \phi_w(z)$ is invariant under the 
action of $\Gal(L/K)$ on $L \otimes_K K_v(\cC)$.

$(3)$  If each $\phi_w(z)$ is expanded as in $(\ref{CFDF2})$, then

\quad $(a)$ for each $i$, $j$, each $w|v$, and each $\sigma \in \Gal(L/K)$,
\begin{equation*}
a_{\sigma(w),\sigma(i),j} \ = \ \tau_{\sigma,w}(a_{w,ij}) \ , \quad \text{and} 
\end{equation*}

\quad $(b)$ for each $\lambda$ there is an $a_{v,\lambda} \in K_v$ such that
      $a_{w,\lambda} = a_{v,\lambda}$ for all $w|v$.

\noindent{Under} these conditions, for each $i$ and $j$, if we write $F = K(x_i)$
then $\oplus_{w|v} a_{w,ij}$ belongs to $\oplus_{u|v} F_u$,
embedded semi-diagonally in $L \otimes_K K_v$.
\end{proposition}
 
\begin{proof}  The equivalence of (1) and (2)
follows from the description of the action of $\Gal(L/K)$ on $L \otimes_K K_v(\cC)$ given above.
For the equivalence of (2) and (3), note that if $\oplus_{w|v} \phi_w(z)$
is expanded as in (\ref{CFDF2}), with $a_{ij} = \oplus_{w|v} a_{w,ij}$ and
$a_{\lambda} = \oplus_{w|v} a_{w,\lambda}$, then
\begin{eqnarray*}
\lefteqn{\sigma\left(\sum_{i=1}^m \sum_{j=0}^{N_i-(2g+1)}
      a_{ij} \varphi_{i,N_i-j}(z)
  + \sum_{\lambda=1}^{\Lambda_0} a_{\lambda} \varphi_{\lambda}(z)\right)} \\
  & = & \sum_{i=1}^m \sum_{j=0}^{N_i-(2g+1)}
                \sigma(a_{ij}) \varphi_{\sigma(i),N_i-j}(z)
   + \sum_{\lambda=1}^{\Lambda_0} \sigma(a_{\lambda}) \varphi_{\lambda}(z) \ .
\end{eqnarray*}
Thus $\oplus_{w|v} \phi_w(z)$ is invariant under $\Gal(L/K)$ if and only if
$\sigma(a_{ij}) = a_{\sigma(i),j}$ and $\sigma(a_{\lambda}) = a_{\lambda}$
for all $\sigma$ and all $i$, $j$, and $\lambda$.  In view of the description
of the action of $\Gal(L/K)$ on $\oplus_{w|v} L_w$,
this holds if and only if
$\tau_{\sigma,w}(a_{w,ij}) = a_{\sigma(w),\sigma(i),j}$ for all
$\sigma$, $w$, $i$, $j$;
and $a_{w,\lambda} = a_{v,\lambda} \in K_v$ for all $w$, $\lambda$.
The assertion concerning the $F$-rationality of $\oplus_{w|v} a_{w,ij}$
follows from the discussion at the beginning of the section.
\end{proof}     

If $K$ is a function field (so $\Char(K) = p > 0$), 
we can expand $\phi(z)$ using the $L^{\sep}$-rational basis as
\index{basis!$L^{\sep}$-rational}  
\begin{equation} \label{CFDF1p}
\phi(z) \ = \ \sum_{i=1}^m \sum_{j=0}^{N_i-J+1}
                             \ta_{ij} \tphi_{i,N_i-j}(z)
           + \sum_{\lambda=1}^{\Lambda_0} \ta_{\lambda} \tphi_{\lambda}(z)  
\end{equation}
with the $\ta_{ij}, \ta_{\lambda} \in L^{\sep}$. Again, the $\ta_{ij}$ are $K$-symmetric relative 
\index{$K$-symmetric!set of numbers}
to the action of $\Gal(L^{\sep}/K)$, and each $\ta_{\lambda}$ belongs to $K$.  
If $w$ is a place of $L^{\sep}$ with $w|v$,   
and we are given an $(\fX,\vs)$ function
$\phi_w(z) \in L_w^{\sep}(\cC)$ of degree $N$, we can write
\begin{equation} \label{CFDF2p}
\phi_w(z) \ = \ \sum_{i=1}^m \sum_{j=0}^{N_i-J-1}
                             \ta_{w,ij} \tphi_{i,N_i-j}(z)
           + \sum_{\lambda=1}^{\Lambda} \ta_{w,\lambda} \tphi_{\lambda}(z)  
\end{equation}
with each $a_{w,ij} \in L_w^{\sep}$ and each $a_{w, \lambda} \in K_v$.  
In this case we have\index{semi-local theory!for function fields} 

\begin{proposition} \label{CPPr1p}  Suppose $K$ is a function field, 
and let $v$ be a place of $K$. For each place $w$ of $L^{\sep}$ with $w|v$, 
let an $(\fX,\vs)$-function\index{$(\fX,\vs)$-function} 
$\phi_w(z) \in L_w^{\sep}(\cC)$ be given.  
Then the following are equivalent$:$

$(1)$  There is a $\phi_v(z) \in K_v(\cC)$
       such that $\phi_w(z) = \phi_v(z)$ for all $w|v$.

$(2)$  $\oplus_{w|v} \phi_w(z)$ is invariant under the action of 
$\Gal(L^{\sep}/K)$ on $L^{\sep} \otimes_K K_v(\cC)$.

$(3)$  If each $\phi_w(z)$ is expanded as in $(\ref{CFDF2p})$, then

\quad $(a)$ for each $i$, $j$, each $w|v$, and each $\sigma \in \Gal(L^{\sep}/K)$,
\begin{equation*}
\ta_{\sigma(w),\sigma(i),j} \ = \ \tau_{\sigma,w}(\ta_{w,ij}) \ , \quad \text{and} 
\end{equation*}

\quad $(b)$ for each $\lambda$ there is an $\ta_{v,\lambda} \in K_v$ such that
      $\ta_{w,\lambda} = \ta_{v,\lambda}$ for all $w|v$.

\noindent{Under} these conditions, for each $i$ and $j$, if we write $F = K(x_i)^{\sep}$
then $\oplus_{w|v} \ta_{w,ij}$ belongs to $\oplus_{u|v} F_u$,
embedded semi-diagonally in $L^{\sep} \otimes_K K_v$.
\end{proposition}

\begin{proof}  The proof is similar to that of Proposition \ref{CPPr1}.
\end{proof}
\index{semi-local theory|)}


\section{ Proof of Theorem $\ref{aT1-B}$ when $\Char(K) = 0$ } \label{Char0Section} 
\index{patching argument!global!when $\Char(K) = 0$|(}

In this section we will prove Theorem \ref{aT1-B} when $\Char(K) = 0$. 
Let $K$, $L = K(\fX)$, $S_K$, and $\EE$ be as in Theorem \ref{aT1-B}, 
and let $\hS_K \supseteq S_K$ be the finite set of places of $K$ 
satisfying the conditions (\ref{hS_KList}) at the beginning this Chapter.
Let $\hS_L$ be the set of places of $L$ above $\hS_K$. 

The proof has three stages, and will occupy the rest of this section.
First, we choose the parameters governing the patching process.
\index{patching parameters}\index{initial approximating functions $f_v(z)$} 
Next, we construct a set of `initial approximating functions' $f_v(z)$ for $v \in \hS_K$,
\index{initial approximating functions $f_v(z)$}  
whose roots belong to $E_v$,   
and modify them to obtain `coherent approximating functions' $\phi_v(z)$ 
\index{coherent approximating functions $\phi_v(z)$} 
whose leading coefficients satisfy the conditions of Proposition \ref{CPSUnit2}.  
\index{coefficients $A_{v,ij}$!leading}   
By means of a degree-raising procedure,\index{degree-raising}
we use the coherent approximating functions to 
construct `initial patching functions' $G_v^{(0)}(z)$  
\index{patching functions, initial $G_v^{(0)}(z)$!construction of} 
\index{patching argument!global} 
whose whose roots also belong to $E_v$.  Finally, we patch the coefficients,  
\index{coefficients $A_{v,ij}$}
creating a sequence of $K_v$-symmetric functions 
\index{$K_v$-symmetric!set of functions}\index{patching functions, $G_v^{(k)}(z)$ for $1 \le k \le n$!constructed by patching}
$G_v^{(1)}(z), G_v^{(2)}(z), \cdots, G_v^{(n)}(z)$\index{patching functions, $G_v^{(k)}(z)$ for $1 \le k \le n$!constructed by patching} 
which have more and more coefficients in the global field $L$, 
but whose roots still belong to $E_v$.  
The final functions $G_v^{(n)}(z)$ have all their coefficients 
in $L$,\index{patching functions, $G_v^{(k)}(z)$ for $1 \le k \le n$!$G_v^{(n)}(z) = G^{(n)}(z)$ is independent of $v$} 
\index{coefficients $A_{v,ij}$}
and are $K_v$-rational but independent of $v$. 
Using Proposition \ref{CPPr1}, they can be seen to be $K$-rational. 
In this way we construct a function $G^{(n)}(z) \in K(\cC)$ whose roots belong to $E_v$ 
for each $v$.   

\begin{proof}[Proof of Theorem $\ref{aT1-B}$ when $\Char(K) =0$]
We begin by outlining the first two stages of the construction.
\index{global patching when $\Char(K) = 0$!outline of Stages 1 and 2}  

First we choose parameters $0 < h_v < r_v < R_v$ for $v \in \hS_K$, 
which control the amount of the freedom in the patching process.  
\index{patching parameters} 
Using the Green's matrix $\Gamma(\EE,\fX)$  
\index{Green's matrix!global}
we construct a $K$-symmetric probability vector $\vs$
\index{$K$-symmetric!probability vector}
with positive rational coefficients, and a positive integer $N$, 
which will be the common degree of the initial approximating functions.
\index{initial approximating functions $f_v(z)$} 
 
Using the approximation theorems from \S\ref{Chap5} and \S\ref{Chap6}, 
we construct the initial approximating functions $\{f_v(z)\}_{v \in \hS_K}$, 
\index{initial approximating functions $f_v(z)$}
which are $(\fX,\vs)$-functions\index{$(\fX,\vs)$-function} 
of common degree $N$ with roots belonging to $E_v$.
We then modify the $f_v(z)$ to obtain the coherent approximating functions\index{coherent approximating functions $\phi_v(z)$}
$\{\phi_v(z)\}_{v \in \hS_K}$.  The key properties of the $\phi_v(z)$ will be 
\begin{enumerate}

\item For each $v \in \hS_K$, $\phi_v(z) \in K_v(\cC)$ 
is an $(\fX,\vs)$ function of degree $N$, for which    
there are a constant $\kappa_v \in K_v^{\times}$ with $|\kappa_v|_v \ge 1$ 
and an initial approximating function $f_v(z) \in K_v(\cC)$ for $E_v$,  
\index{initial approximating functions $f_v(z)$}
such that 
\begin{equation*}
\phi_v(z) \ = \ \kappa_v f_v(z) \ . 
\end{equation*} 
Thus the $\phi_v(z)$ inherit the approximation properties of the $f_v(z)$, and their
roots are the same as those of the $f_v(z)$.  
  
\item  For each $w \in \hS_L$, put $\phi_w(z) = \phi_v(z)$ if $w|v$, and view $\phi_w(z)$
as an element of $L_w(\cC)$.  Although the $\phi_w(z)$ for $w|v$ are all the same, 
for distinct $w$ the points of $\fX$, which are their poles, 
are identified differently as points of $\cC_w(L_w)$.
Write $\tc_{w,i} = \lim_{z \rightarrow x_i} \phi_w(z) \cdot g_{x_i}(z)^{N s_i}$
for the leading coefficient of $\phi_w(z)$ at $x_i$. 
\index{coefficients $A_{v,ij}$!leading}

Then for each $i$, $\oplus_{w|\infty} \tc_{w,i}$ is an $\hS_L$-subunit and 
\index{subunit}
\begin{equation*}
\prod_{w \in \hS_L} |\tc_{w,i}|_w^{D_w} \ = \ 1 \ .
\end{equation*}
\end{enumerate} 

Our ability to achieve (1) 
uses that $\gamma(\EE,\fX) > 1$.  Our ability to achieve
(2) depends on the independent variability of the archimedean 
logarithmic leading coefficients, proved in Theorems \ref{CThm1} and \ref{RThm2}.
\index{logarithmic leading coefficients!independent variability of archimedean}
\index{independent variability!of logarithmic leading coefficients}

\medskip  
For the convenience of the reader, Theorem \ref{CCPX1} below 
summarizes Theorems \ref{CThm1}, \ref{RThm2}, \ref{RLThm}, \ref{CompactThm}  
and Corollaries \ref{CompactThmCor} and \ref{EvHvExtraPtsCor}
which construct the initial approximating functions.  In broad terms, those
\index{initial approximating functions $f_v(z)$}\label{`SymbolIndexInitialApproxf'}
theorems say that for each $K$-symmetric $\vs \in \cP^m(\QQ)$ and each $v \in \hS_K$, 
\index{$K$-symmetric!probability vector}
there are a set $\tE_v \subset E_v$ and an integer $N_v > 0$ 
such that for each $N > 0$ divisible by $N_v$, 
there is an $(\fX,\vs)$-function\index{$(\fX,\vs)$-function} $f_v(z) \in K_v(\cC)$ of degree $N$ 
such that $\frac{1}{N} \log_v(|f_v(z)|_v)$ approximates $\sum_{i=1}^m s_i G(z,x_i;\tE_v)$, 
and $f_v(z)$ has roots in $E_v$ and satisfies certain side conditions. 

Recall that if $f_v(z) \in K_v(\cC)$ is an $(\fX,\vs)$-function\index{$(\fX,\vs)$-function} 
of degree $N$,
then for each $x_i \in \fX$, the leading coefficient of $f_v$ at $x_i$ is 
\index{coefficients $A_{v,ij}$!leading}
\begin{equation*} 
c_{v,i} \ = \ (f_v \cdot g_{x_i}^{Ns_i})|_{x_i} \ = \ 
               \lim_{z \rightarrow x_i} f_v(z) \cdot g_{x_i}(z)^{N s_i} \ ,
\end{equation*} 
and that $\Lambda_{x_i}(f_v,\vs) = \frac{1}{N} \log_v(|c_{v,i}|_v)$. 
\index{global patching when $\Char(K) = 0$!Stage 1: Choices of sets and parameters!summary of the Initial Approximation theorems}  

\begin{theorem}  \label{CCPX1}  
{\bf $($Summary of the Initial Approximation Theorems$)$}
\index{initial approximation theorem!summary of initial approximation theorems!when $\Char(K_v) = 0$} 

Let $K$, $\EE$, and $\fX$ be as in Theorem $\ref{aT1-B}$. Then for each $v \in \hS_K$,  

\vskip .05 in
\noindent{$(A)$}  If $K_v \cong \CC$ $($so $E_v \subset \cC_v(\CC)$ is compact, 
$\CC$-simple,\index{$\CC$-simple set}
and disjoint from $\fX)$, let $U_v = E_v^0$ be the interior of $E_v$.   
Then for each $\varepsilon_v > 0$ there is a compact set $\tE_v$ contained in $E_v^0$  such that  

$(1)$ For each $x_i, x_j \in \fX$ with $x_i \ne x_j$, 
\begin{equation*} 
|V_{x_i}(\tE_v) - V_{x_i}(E_v)| \ < \ \varepsilon_v \ , \quad
|G(x_i,x_j;\tE_v) - G(x_i,x_j;E_v)| \ < \ \varepsilon_v \ ;
\end{equation*}
\index{Green's function!nonarchimedean}
\index{Robin constant!nonarchimedean}
   
$(2)$ There is a $\delta_v > 0$ with the property that for each $\vs \in \cP^m(\QQ)$, 
there is an integer $N_v > 0$ such that for 
each $\vbeta_v = {}^t(\beta_{v,1}, \ldots, \beta_{v,m}) \in [-\delta_v,\delta_v]^m$, 
and each positive integer  $N$ divisible by $N_v$, 
there is an $(\fX,\vs)$-function\index{$(\fX,\vs)$-function} $f_v(z) \in K_v(\cC)$ 
of degree $N$  satisfying
          
\quad $(a)$ $\{ z \in \cC_v(\CC) : |f_v(z)|_v \le 1 \} \subset E_v^0;$  

\quad $(b)$  For each $x_i \in \fX$, \ 
         $\frac{1}{N} \log_v(|c_{v,i}|_v) =  \Lambda_{x_i}(\tE_v,\vs) + \beta_{v,i}.$        

\vskip .05 in
\noindent{$(B)$}  If $K_v \cong \RR$ $($so $E_v \subset \cC_v(\CC)$ is compact, 
$\RR$-simple,\index{$\RR$-simple set} and disjoint from $\fX)$, 
let $E_v^0$ be the quasi-interior of $E_v$.  Then for each $\varepsilon_v > 0$, 
and each open set $U_v \subset \cC_v(\CC)$ which is stable under complex conjugation, 
bounded away from $\fX$, and satisfies $U_v \cap E_v = E_v^0$,  
there are a compact set $\tE_v$ contained in $E_v^0$ such that 

$(1)$ For each $x_i, x_j \in \fX$ with $x_i \ne x_j$, 
\begin{equation*} 
|V_{x_i}(\tE_v) - V_{x_i}(E_v)| \ < \ \varepsilon_v \ , \quad
|G(x_i,x_j;\tE_v) - G(x_i,x_j;E_v)| \ < \ \varepsilon_v \ ;
\end{equation*}
\index{Robin constant!nonarchimedean}
\index{Green's function!nonarchimedean}
   
$(2)$ For each $0 < \cR_v < 1$, there is a 
$\delta_v > 0$ with the property that for each $K_v$-symmetric $\vs \in \cP^m(\QQ)$, 
\index{$K_v$-symmetric!probability vector}
there is an integer $N_v > 0$  
such that for each $K_v$-symmetric 
\index{$K_v$-symmetric!vector}
$\vbeta_v = {}^t(\beta_{v,1}, \ldots, \beta_{v,m}) \in [-\delta_v,\delta_v]^m$ 
and each positive integer  $N$ divisible by $N_v$, 
there is an $(\fX,\vs)$-function\index{$(\fX,\vs)$-function} $f_v(z) \in K_v(\cC)$ 
of degree $N$  satisfying 

\quad $(a)$ $\{z \in \cC_v(\CC) : |f_v(z)| \le 1 \} \ \subset \ U_v$\ ,
all the zeros of $f_v(z)$ belong to $E_v^0$, and   
if $E_{v,i}$ is a component of $E_v$ contained in $\cC_v(\RR)$  
and $f_v(z)$ has $N_i$ zeros in $E_{v,i}$, 
then $f_v(z)$ oscillates $N_i$ times between $\pm \cR_v^N$ on $E_{v,i}$.  

\quad $(b)$  For each $x_i \in \fX$, \ 
      $\frac{1}{N} \log_v(|c_{v,i}|_v) = \Lambda_{x_i}(\tE_v,\vs) + \beta_{v,i}$.

\vskip .05 in
\noindent{$(C)$}  If $K_v$ is nonarchimedean and $v \in S_K$, $($so 
$E_v$ is compact, $K_v$-simple, and disjoint from $\fX)$, 
\index{$K_v$-simple!set}  
fix a $K_v$-simple decomposition 
\index{$K_v$-simple!decomposition}  
\begin{equation} \label{E_vSimple2}
E_v \ = \ \bigcup_{\ell=1}^{D_v} B(a_\ell,r_\ell) \cap \cC_v(F_{w_\ell}) \ .
\end{equation} 
and fix $\varepsilon_v > 0$.    
Then there is a $K_v$-simple set $\tE_v \subseteq E_v$  compatible with $E_v$
\index{$K_v$-simple!set!compatible with another set}
such that 

$(1)$ For each $x_i, x_j \in \fX$ with $x_i \ne x_j$, 
\begin{equation*} 
|V_{x_i}(\tE_v) - V_{x_i}(E_v)| \ < \ \varepsilon_v \ , \quad
|G(x_i,x_j;\tE_v) - G(x_i,x_j;E_v)| \ < \ \varepsilon_v \ ;
\end{equation*}
\index{Robin constant!nonarchimedean}
\index{Green's function!nonarchimedean}

$(2)$ For each $0 < \beta_v \in \QQ$  
and each $K_v$-symmetric $\vs \in \cP^m(\QQ)$, 
\index{$K_v$-symmetric!probability vector}
there is an integer $N_v \ge 1$  such that 
for each positive integer  $N$ divisible by $N_v$, 
there is an $(\fX,\vs)$-function\index{$(\fX,\vs)$-function} $f_v \in K_v(\cC_v)$ 
of degree $N$ satisfying
  
\quad $(a)$ The zeros $\theta_1, \ldots, \theta_N$ of $f_v$ are distinct and belong to $E_v$.

\quad $(b)$ $f_v^{-1}(D(0,1)) \subseteq \bigcup_{\ell = 1}^{D_v} B(a_{\ell},r_{\ell})$, 
and there is a decomposition $f_v^{-1}(D(0,1)) = \bigcup_{h=1}^N B(\theta_h,\rho_h)$, 
where the balls $B(\theta_h,\rho_h)$ are pairwise disjoint and isometrically parametrizable. 
\index{isometrically parametrizable ball} 
For each $h = 1, \ldots, N$, if $\ell = \ell(h)$ is such that 
$B(\theta_h,\rho_h) \subseteq B(a_{\ell},r_{\ell}),$ put $F_{u_h} = F_{w_\ell};$
then $\rho_h \in |F_{u_h}^{\times}|_v$
and $f_v$ induces an $F_{u_\ell}$-rational scaled isometry from $B(\theta_h,\rho_h)$ to $D(0,1)$,
\index{scaled isometry} 
with
\begin{equation*}
f_v\big( B(\theta_h,\rho_h) \cap \cC_v(F_{u_h})\big) \ = \ \cO_{F_{u_h}} \ , 
\end{equation*}      
such that  $|f_v(z_1)-f_v(z_2)|_v = (1/\rho_h) \|z_1,z_2\|_v$ for all $z_1, z_2 \in B(\theta_h, \rho_h)$. 

\quad $(c)$ The set $H_v := E_v \cap f_v^{-1}(D(0,1))$ is $K_v$-simple and compatible with $E_v$.
\index{$K_v$-simple!set!compatible with another set}
Indeed, 
\begin{equation} \label{E_vPBF2}
H_v \ = \ 
  \bigcup_{h=1}^N \big(B(\theta_h,\rho_h) \cap \cC_v(F_{u_h})\big) \ ,
\end{equation}
and $(\ref{E_vPBF2})$ is a 
$K_v$-simple decomposition of $H_v$ compatible with the 
\index{$K_v$-simple!decomposition!compatible with another decomposition}  
$K_v$-simple decomposition $(\ref{E_vSimple2})$ of $E_v$,
\index{$K_v$-simple!decomposition}  
which is move-prepared $($see Definition $\ref{MovePreparedDef})$ 
\index{move-prepared}
relative to the balls $B(a_1,r_1), \ldots, B(a_{D_v},r_{D_v})$.
For each $\ell$ there is a point 
$\wbar_\ell \in \big(B(a_\ell,r_\ell) \cap \cC_v(F_{w_\ell})\big) \backslash H_v$.  

\quad $(d)$  For each $x_i \in \fX$, \ 
 $\frac{1}{N} \log_v(|c_{v,i}|_v) \ = \ \Lambda_{x_i}(\tE_v,\vs) + \beta_v$.   

\vskip .05 in
\noindent{$(D)$}  If $K_v$ is nonarchimedean and $v \notin S_K$,
$($so $E_v$ is $\fX$-trivial and in particular is an $\RL$-domain disjoint from $\fX)$, 
\index{$\fX$-trivial}
\index{$\RL$-domain} 
put $\tE_v = E_v$.  

Then for each $K_v$-symmetric $\vs \in \cP^m(\QQ)$,
\index{$K_v$-symmetric!probability vector}  
there is an integer $N_v \ge 1$  such that for each positive integer  $N$ divisible by $N_v$, 
there is an $(\fX,\vs)$-function\index{$(\fX,\vs)$-function} $f_v \in K_v(\cC_v)$ 
of degree $N$ satisfying

$(a)$ $E_v = \tE_v  = \{z \in \cC_v(\CC_v) : |f_v(z)|_v \le 1 \};$

$(b)$ For each $x_i \in \fX$, \ $\frac{1}{N} \log_v(|c_{v,i}|_v) \ = \ \Lambda_{x_i}(\tE_v,\vs)$.
\end{theorem} 

\smallskip

Note that in (A) and (B) of Theorem \ref{CCPX1}, the number $\delta_v > 0$
depends on $\varepsilon_v$, $E_v$ and $U_v$, but the numbers $\beta_{v,i}$ for which 
\begin{equation*}
\frac{1}{N} \log_v(|c_{v,i}|_v) \ = \ \Lambda_{x_i}(\tE_v,\vs) + \beta_{v,i}
\end{equation*}
can be specified arbitrarily provided they are $K_v$-symmetric and satisfy
\index{$K_v$-symmetric!set of numbers}
$-\delta_v \le \beta_{v,i} \le \delta_v$ for each $i$.  
In (C) the number $0 < \beta_v \in \QQ$ is the same for all $i$.
In (D), the logarithmic leading coefficients match the $\Lambda_{x_i}(\tE_v,\vs)$ exactly.
\index{logarithmic leading coefficients}
For each $v$, the leading coefficients $c_{v,i}$ of $f_v$ are $K_v$-symmetric,
\index{coefficients $A_{v,ij}$!leading}
\index{$K_v$-symmetric!set of numbers} 
because $f_v(z)$ is $K_v$-rational and the $g_{x_i}(z)$ are $K_v$-symmetric.
\index{$K_v$-symmetric!set of functions}

We now turn to the details of the proof.

\medskip
{\bf Stage 1.  Choices of the sets and parameters.}
\index{patching parameters|ii}
\index{global patching when $\Char(K) = 0$!Stage 1: Choices of sets and parameters!the open sets $U_v$} 
We begin by making the choices that govern the patching process.  
Given a number field $F$ containing $K$, write $\hS_{F,\infty}$ for the set of archimedean places of $F$, 
and $\hS_{F,0}$ for the set of nonarchimedean places in $\hS_F$,
so $\hS_F = \hS_{F,\infty} \cup \hS_{F,0}$.  Similarly, write $S_F = S_{F,\infty} \cup S_{F,0}$.

\smallskip
{\bf The open sets $U_v$ for $v \in \hS_{K\infty}$.}  
For each $v \in \hS_{K,\infty}$ with $K_v \cong \CC$, 
$E_v$ is $\CC$-simple,\index{$\CC$-simple set} so it has finitely 
finitely many components, each of which is simply connected, has a piecewise smooth boundary\index{boundary!piecewise smooth} 
and is the closure of its interior.  Let  $U_v = E_v^0$ be its interior.
\index{$K_v$-simple!$\CC$-simple}\index{$\CC$-simple set}\index{simply connected}\index{closure of $\cC_v(\CC)$ interior}

For each $v \in \hS_{K,\infty}$ with $K_v \cong \RR$, 
$E_v$ is $\RR$-simple,\index{$\RR$-simple set} 
so it is stable under complex conjugation and has finitely many components, each of which is 
an interval of positive length contained in $\cC_v(\RR)$, or is disjoint from $\cC_v(\RR)$ 
and is simply connected, has a piecewise smooth boundary,\index{boundary!piecewise smooth}
\index{simply connected}\index{$K_v$-simple!$\RR$-simple}\index{$\RR$-simple set}
 and is the closure of its interior.\index{closure of $\cC_v(\CC)$ interior}  
Let $U_v \subset \cC_v(\CC)$ be an open set such that $U_v \cap E_v = E_v^0$, the quasi-interior of $E_v$.
We will choose $U_v$ so that it is stable under complex conjugation, bounded away from $\fX$, 
and has the same number of connected components as $E_v$. Thus, $U_v$ is the union of the interiors
of the components $E_{v,i}$ in $\cC_v(\CC) \backslash \cC_v(\RR)$, together with open sets $U_{v,i}$ 
such that $U_{v,i} \cap E_v$ is the real interior of $E_{v,i}$, 
for the $E_{v,i} \subset \cC_v(\RR)$. 

\smallskip
{\bf The $K_v$-simple decompositions of $E_v$ and the sets $U_v$, for $v \in S_{K,0}$.}
\index{$K_v$-simple!decomposition}
\index{global patching when $\Char(K) = 0$!Stage 1: Choices of sets and parameters!the $K_v$-simple decompositions}   
For each $v \in S_{K,0}$, the set $E_v$ is compact and $K_v$-simple 
\index{$K_v$-simple!set}  
(see Definition \ref{KvSimple}).  

Choose a $K_v$-simple decomposition    
\index{$K_v$-simple!decomposition}  
\begin{equation} \label{K_vSimpleDecompF}
E_v \ = \ \bigcup_{\ell=1}^{D_v} B(a_\ell,r_\ell) \cap \cC_v(F_{w_\ell}) \ .\
\end{equation}   
By refining this decomposition, if necessary, we can assume that 
$U_v := \bigcup_{\ell=1}^{D_v} B(a_\ell,r_\ell)$ is disjoint from $\fX$.   
Such a decomposition will be fixed for the rest of the construction.    

\smallskip
{\bf The sets $\tE_v$ for $v \in \hS_K$.} 
\index{global patching when $\Char(K) = 0$!Stage 1: Choices of sets and parameters!the sets $\tE_v$}  
By hypothesis, $\gamma(\EE,\fX) > 1$ in Theorem \ref{aT1-B}.  
This means that the Green's matrix  $\Gamma(\EE,\fX)$ is negative definite.
\index{Green's matrix!global}\index{Green's matrix!negative definite}  
Suppose $\tEE = \prod_{v \in \hS_K} \tE_v \times \prod_{v \notin \hS_K} E_v$ is another  
$K$-rational adelic set compatible with $\fX$.{compatible with $\fX$}  
By the discussion leading to (\ref{NegDefCloseness}), 
there are numbers $\varepsilon_v > 0$ for $v \in \hS_K$ such that 
such that $\Gamma(\widetilde{\EE},\fX)$ is also negative definite,\index{Green's matrix!negative definite} 
provided that for each $v \in \hS_K$ 
\begin{equation} \label{NegDefCloseness1}
\left\{ \begin{array}{ll}
       |G(x_j,x_i;\tE_v)-G(x_j,x_i;E_v)| \ < \ \varepsilon_v &
                  \text{for all $i \ne j$ \ ,} \\
       |V_{x_i}(\tE_v)-V_{x_i}(E_v)| \ < \ \varepsilon_v &
                  \text{for all $i$ \ .}  
        \end{array}  \right. 
\end{equation}
\index{Robin constant!nonarchimedean}
\index{Green's function!nonarchimedean}

For each $v \in \hS_K$, we will take $\tE_v \subseteq E_v$ to be the set given 
by Theorem \ref{CCPX1} for $E_v$, relative to the number $\varepsilon_v$ chosen above 
(and the set $U_v$, if $K_v \cong \RR$), satisfying (\ref{NegDefCloseness1}). 
Put $\tEE =  \prod_{v \in \hS_K} \tE_v \times \prod_{v \notin \hS_K} E_v$ 
with the $\tE_v$ chosen above, and let 
\begin{equation} \label{tV_KChoice} 
\tV_K \ := \ V(\tEE,\fX) \ = \ \val(\Gamma(\tEE,\fX)) 
\end{equation}
be the global Robin constant for $\tEE$ and $\fX$.  By construction, $\tV_K < 0$. 
\index{Robin constant!global} 

\smallskip
{\bf The local parameters $\eta_v$, $\cR_v$, $h_v$, $r_v$, and $R_v$.}
\index{global patching when $\Char(K) = 0$!Stage 1: Choices of sets and parameters!the local parameters $\eta_v$, $\cR_v$, $h_v$, $r_v$, and $R_v$}
\index{patching parameters!choice when $\Char(K) = 0$|ii}
Fix a collection of real numbers $\{\eta_v\}_{v \in \hS_K}$ 
with  $\eta_v > 0$ for each $v \in \hS_K$ and
$\eta_v \in \QQ$ for each $v \in \hS_{K,0}$, such that 
\begin{equation} \label{FPP1}
\sum_{v \in \hS_K} \eta_v \log(q_v) \ = \ |\tV_K| \ = \ - \tV_K \ .
\end{equation}  
\index{Robin constant!global}
The $\eta_v$ provide the freedom for adjustment needed  
in the construction of the initial approximating functions,
\index{initial approximating functions $f_v(z)$} 
and determine the scaling factors in passing from the initial approximating functions
to the coherent approximating functions. 
\index{initial approximating functions $f_v(z)$}\index{coherent approximating functions $\phi_v(z)$}

For each $v \in \hS_{K,\infty}$, fix a number $r_v$ with $1 < r_v < e^{\eta_v}$.
Then, choose a set of numbers $\{h_v\}_{v \in \hS_K}$  
with $\prod_{v \in \hS_K} h_v^{D_v} > 1$, such that 
\begin{equation} \label{FPP2}
\left\{ \begin{array}{ll} 1 < h_v < r_v & \text{if $v \in \hS_{K,\infty}$\ ,} \\
                          0 < h_v < 1 & \text{if $v \in \hS_{K,0}$\ .} \end{array} \right. 
\end{equation}  
Finally, for each $v \in S_{K,0}$, fix an $r_v$ with $h_v < r_v < 1$,
and for each $v \in \hS_{K,0} \backslash S_{K,0}$ put $r_v = 1$.
Note that $D_v = \log(q_v)$ for each archimedean $v$. 
Thus $0 < h_v < r_v$ for all $v$, and 
\begin{equation} \label{FPP3} 
1 \ < \ \prod_{v \in \hS_K} h_v^{D_v} \ < \ \prod_{v \in \hS_K} r_v^{D_v} \ .  
\end{equation} 
In the patching process, 
\index{patching parameters} 
the numbers $h_v$ will control how much the coefficients of the functions being patched 
\index{coefficients $A_{v,ij}$} 
can be changed, and the $r_v$ will be ``encroachment bounds'' which limit how close 
certain quantities can come to the $h_v$.
 
For each $v \in \cS_{K,\infty}$ with $K_v \cong \CC$, put $\hR_v = e^{\eta_v}$.
For each $v \in \cS_{K,\infty}$ with $K_v \cong \RR$, 
fix a number $0 < \cR_v < 1$ close enough to $1$ that 
\begin{equation} \label{FCCA1}
1 \ < \ r_v \ < \ \cR_v \cdot e^{\eta_v} \ < \ e^{\eta_v} \ . 
\end{equation} 
and put $\hR_v = \cR_v \cdot e^{\eta_v}$.  The number $\cR_v$ specifies the magnitude
of the oscillations of the initial approximating functions when $K_v \cong \RR$.
\index{initial approximating functions $f_v(z)$!when $K_v \cong \RR$}  

In either case, 
we can choose $R_v$ so that $r_v < R_v < \hR_v$;  thus 
\begin{equation} \label{FCCA2} 
1 \ < \ r_v  \ < \ R_v \ < \ \hR_v \ \le \ e^{\eta_v}  \ .
\end{equation} 
For each $v \in \hS_{K,0}$, put $R_v = q_v^{\eta_v}$.  
Then $0 < h_v < r_v < R_v$ for each $v \in \hS_K$, 
and $R_v \in |\CC_v^{\times}|_v$ for each $v \in \hS_{K,0}$.  

\smallskip
{\bf The numbers $\delta_v > 0$ for $v \in \hS_{K,\infty}$.}
\index{global patching when $\Char(K) = 0$!Stage 1: Choices of sets and parameters!the $\delta_v$ for $v \in \hS_{K,\infty}$}
If $K_{v} \cong \CC$, 
let $\delta_v$ be the number given by Theorem \ref{CCPX1}(A.2)  
for $E_v$ and $\tE_v$. If $K_v \cong \RR$, 
let $\delta_v$ be the number 
given by Theorem \ref{CCPX1}(B.2) for $E_v$ and $\tE_v$, 
relative to the number $\cR_v$ chosen in (\ref{FCCA1}).
For each $v \in \hS_{K,\infty}$, the number $\delta_v$ plays the role of a 
`radius of independent variability' for the logarithmic leading coefficients at $v$.
\index{logarithmic leading coefficients}
\index{logarithmic leading coefficients!independent variability of archimedean|ii}
\index{independent variability!of logarithmic leading coefficients|ii}   

\smallskip
{\bf The rational probability vector $\vs$.}
\index{global patching when $\Char(K) = 0$!Stage 1: Choices of sets and parameters!the probability vector $\vs$}
By construction, the Green's matrix $\Gamma(\tEE,\fX)$ is negative definite.\index{Green's matrix!negative definite} 
As above, put  $\tV_K = V(\tEE_K,\fX) < 0$.  
\index{Robin constant!global}
Let $\ts \in \cP^m(\RR)$ be the $K$-symmetric probability vector 
\index{$K$-symmetric!probability vector}
given by Proposition \ref{UniqueSVec} for which   
\begin{equation*}
\left( \begin{array}{c} \tV_K \\ \vdots \\ \tV_K \end{array} \right) 
       \ = \ \Gamma(\tEE_K,\fX) \ts \ . 
\end{equation*}
The entries of $\ts$ are positive, but they need not be rational.  
 
Fix an archimedean place $v_0$ of $K$, and let $\delta_{v_0}$ be the radius
of independent variability for the logarithmic leading coefficients at $v_0$, constructed above.
\index{logarithmic leading coefficients}
\index{logarithmic leading coefficients!independent variability of archimedean}
\index{independent variability!of logarithmic leading coefficients} 
By continuity, there is a $K$-symmetric probability vector
\index{$K$-symmetric!probability vector}
$\vs \in \cP^m(\QQ)$ close enough to $\ts$ that 
all its entries are positive, and such that for each $i = 1, \ldots, m$, 
the $i^{th}$ coordinate of $\Gamma(\tEE_K,\fX) \vs$ satisfies 
\begin{equation} \label{CGPF2} 
|\tV_K - (\Gamma(\tEE_K,\fX)\vs)_i| 
\ < \ \delta_{v_0} \log(q_{v_0}) \ .
\end{equation} 
\index{Robin constant!global}
This $\vs$ will be fixed for the rest of the construction.  

\smallskip
{\bf Stage 2. Construction of the Approximating Functions $f_v(z)$ and $\phi_v(z)$.} 

In this stage, we construct the initial approximating functions $f_v(z)$, 
\index{initial approximating functions $f_v(z)$}\index{coherent approximating functions $\phi_v(z)$|ii}
then modify them to obtain the coherent approximating functions $\phi_v(z)$. 
Our goal is to prove the following theorem:
\index{global patching when $\Char(K) = 0$!Stage 2: Constructing approximating functions!the Coherent Approximation theorem} 

\begin{theorem} \label{CTCX2}\label{`SymbolIndexCoherentApproxf'}
\index{coherent approximating functions $\phi_v(z)$!construction when $\Char(K) = 0$}  
Let $\cC$, $K$, $\EE$, $\fX$, and $S_K$ be as in Theorem $\ref{aT1-B}$, with $\Char(K) = 0$.
Let $\hS_K \supseteq S_K$ be the finite set of places satisfying conditions $(\ref{hS_KList})$.  
For each $v \in \hS_K$, let $\tE_v \subset E_v$ and $0 < h_v < r_v < R_v$ 
be the set and patching parameters constructed above. 
\index{patching parameters} 
For each $v \in \hS_{K,\infty}$, let $U_v \subset \cC_v(\CC)$ be the chosen set  
with $U_v \cap E_v = E_v^0$, and let $\delta_v > 0$ be the radius of independent variability
for the logarithmic leading coefficients.  For each $v \in S_{K,0}$ let
\index{logarithmic leading coefficients}
\index{logarithmic leading coefficients!independent variability of archimedean}
\index{independent variability!of logarithmic leading coefficients}  
$\bigcup_{\ell=1}^{D_v} B(a_\ell,r_\ell) \cap \cC_v(F_{w_\ell})$ 
be the chosen $K_v$-simple decomposition of $E_v$.  
\index{$K_v$-simple!decomposition}  
Let $\vs \in \cP_m(\QQ)$ be the chosen rational probability vector with positive coefficients, 
satisfying $(\ref{CGPF2})$.

Then there are a positive integer $N$ and 
$(\fX,\vs)$-functions\index{$(\fX,\vs)$-function} $\phi_v(z) \in K_v(\cC)$,
for $v \in \hS_K$, of common degree $N$, such that for each $v$ the zeros of $\phi_v(z)$ 
belong to $E_v$ and 

$(A)$ The $\phi_v(z)$ have the following mapping properties:  

\hskip .15in $(1)$  If  $K_v \cong \CC$, then 
  \begin{equation*}
         \{ z \in \cC_v(\CC_v) : |\phi_v(z)|_v \le R_v^N \} \ \subset \ U_v \ = \ E_v^0 \ .
  \end{equation*}

\hskip .15in $(2)$  If $K_v \cong \RR$, then 
  \begin{equation*}
      \{ z \in \cC_v(\CC_v) : |\phi_v(z)|_v \le 2R_v^N \} \ \subset \ U_v \ , 
  \end{equation*}
  
    \hskip .4in   and for each component $E_{v,j}$ of $E_v$ contained in $\cC_v(\RR)$, 
           if $\phi_v(z)$ has $\tau_j$ zeros 

    \hskip .4in  in $E_{v,j}$, then $\phi_v(z)$ oscillates $\tau_j$ times 
           between $\pm 2R_v^N$ on $E_{v,j}$.  
  
\hskip .15in $(3)$ If $K_v$ is nonarchimedean and $v \in S_K$, then 
\begin{equation} \label{CFBZ1}
r_v^N \ < \ q_v^{-1/(q_v-1)} \ < \ 1 \ , \quad \text{and}
\end{equation} 

    \hskip .3in $(a)$ the zeros $\theta_1, \ldots, \theta_N$ of $\phi_v(z)$ are distinct; 
%
          
    \hskip .3in $(b)$ $\phi_v^{-1}(D(0,1)) = \bigcup_{h=1}^N B(\theta_h,\rho_h)$, 
              where $B(\theta_1,\rho_1), \ldots, B(\theta_N,\rho_N)$ are pairwise 
             
    \hskip .5in  disjoint, isometrically parametrizable, 
\index{isometrically parametrizable ball}
             and contained in $\bigcup_{\ell = 1}^{D_v} B(a_{\ell},r_{\ell});$
           
    \hskip .3in $(c)$ $H_v := \phi_v^{-1}(D(0,1)) \cap E_v$ is $K_v$-simple, 
\index{$K_v$-simple!set}  
            with the $K_v$-simple decomposition 
\index{$K_v$-simple!decomposition}   
            \begin{equation*} 
            H_v \ = \ \bigcup_{k=1}^n \big( B(\theta_k,\rho_k) \cap \cC_v(F_{u_h}) \big)
            \end{equation*}

    \hskip .5in compatible with the $K_v$-simple decomposition 
    \index{$K_v$-simple!decomposition!compatible with another decomposition} 
           $\bigcup_{\ell = 1}^{D_v} \big( B(a_\ell,r_\ell) \cap \cC_v(F_{w_\ell}) \big)$ 
           of $E_v$, 
           
     \hskip .5in which is move-prepared relative to $B(a_1,r_1), \ldots, B(a_{D_v},r_{D_v})$. 
\index{move-prepared} 
           For each $\ell$,
     
     \hskip .5in there is a point 
       $\wbar_\ell \in \big(B(a_\ell,r_\ell) \cap \cC_v(F_{w_\ell})\big) \backslash H_v$.
           
  \hskip .3in $(d)$  For each $h = 1, \ldots, N$, $F_{u_h}/K_v$ is finite and separable. 
           If $\theta_h \in E_v \cap B(a_{\ell},r_{\ell})$, 
          
    \hskip .5in then $F_{u_h} = F_{w_\ell}$, $\rho_h \in |F_{w_\ell}^{\times}|_v$, 
           and $B(\theta_h,\rho_h) \subseteq B(a_\ell,r_\ell);$ and $\phi_v$ induces an
           
     \hskip .5in  
            $F_{u_h}$-rational scaled isometry from $B(\theta_h,\rho_h)$ onto $D(0,1)$
\index{scaled isometry} 
            with $\phi_v(\theta_h)= 0$,

    \hskip .5in  
           which takes $B(\theta_h,\rho_h) \cap \cC_v(F_{u_h})$ onto $\cO_{u_h}$. 
    
\hskip .15in $(4)$ If $K_v$ is nonarchimedean and $v \in \hS_K \backslash S_K$, then 
   \begin{equation*}
         E_v \ = \ \tE_v \ = \ \{ z \in \cC_v(\CC_v) : |\phi_v(z)|_v \le R_v^N \} \ .                  
   \end{equation*}
   
$(B)$  For each $w \in \hS_L$, put $\phi_w(z) = \phi_v(z)$ if $w|v$, 
and regard $\phi_w(z)$ as an element of $L_w(\cC)$.  
Each $x_i \in \fX$ is canonically embedded in \, $\cC_w(L_w);$  
let $\tc_{w,i} = \lim_{z \rightarrow x_i} \phi_w(z) \cdot g_{x_i}(z)^{N s_i}$
be the leading coefficient of $\phi_w(z)$ at $x_i$.  Then for each $i$, 
\index{coefficients $A_{v,ij}$!leading} 
\begin{equation*}
    \sum_{w \in \hS_L} \log_w(|\tc_{w,i}|_w) \log(q_w) \ = \ 0 \ . 
\end{equation*} 
Moreover, $\oplus_{w \in \hS_{L}} \tc_{w,i}$ is an $S_L$-subunit:  
\index{subunit}
there are an integer $n_0$ and a 
$K$-symmetric set of $\hS_L$-units $\mu_1, \ldots, \mu_m \in L$, 
\index{$K$-symmetric!set of numbers}
such that for each $i = 1, \ldots, m$,
\begin{equation} \label{Char0Coherence}
\left\{ \begin{array}{ll}
         \tc_{w,i}^{n_0} = \mu_i \ , & \text{if $w \in \hS_{L,\infty}$ \ ,} \\
         |\tc_{w,i}^{n_0}|_w = |\mu_i|_w \ , & \text{if $w \in \hS_{L,0}$ \ ,} 
        \end{array}
    \right.
\end{equation}
Necessarily $\mu_i \in K(x_i)$ for each $i$.
\end{theorem}

\begin{proof}
The proof has several steps, and consists of carefully choosing a 
compatible collection of initial approximating functions $f_v(z)$ 
\index{initial approximating functions $f_v(z)$}
of common degree $N$ in Theorem \ref{CCPX1}, scaling them, 
and then modifying their leading coefficients to satisfy (\ref{Char0Coherence}).   
\index{coefficients $A_{v,ij}$!leading}   

\vskip .1 in
{\bf The choice of $N$.}\index{patching parameters!choice when $\Char(K) = 0$|ii}
\index{global patching when $\Char(K) = 0$!Stage 2: Constructing approximating functions!the choice of $N$} 
For each archimedean $v$ with $K_v \cong \CC$, let $N_v > 0$ be the integer
given by Theorem \ref{CCPX1}(A.2) for $E_v$, $\varepsilon_v$, $\tE_v$, and $\vs$
as chosen above.  
For each archimedean $v$ with $K_v \cong \RR$, let $N_v > 0$ be the
integer given by Theorem \ref{CCPX1}(B.2) for 
$E_v$, $\varepsilon_v$, $\tE_v$, $U_v$, $\cR_v$, and $\vs$ as chosen above.  
For each nonarchimedean $v \in S_{K,0}$, put $\beta_v = \eta_v$ 
(where $0 < \eta_v \in \QQ$ is as in \ref{FPP1})
and let $N_v > 0$ be the integer given by Theorem \ref{CCPX1}(C.2) 
for $E_v$, $\tE_v$, $\vs$, $\beta_v$ and the $K_v$-simple decomposition 
\index{$K_v$-simple!decomposition}  
$E_v = \bigcup_{\ell = 1}^{D_v} \big(B(a_\ell,r_\ell) \cap \cC_v(F_{w_\ell})\big)$ chosen above.   
For each nonarchimedean $v \in \hS_{K,0} \backslash S_{K,0}$, 
let $N_v$ be as given by Theorem \ref{CCPX1}(D) for $E_v = \tE_v$ and $\vs$ as chosen above.

Let $N > 0$ be an integer which satisfies the following conditions: 
\begin{enumerate}
\item $N$ is divisible by $N_v$, for each $v \in \hS_K$;

\item $N$ is divisible by $J$, the number from the construction of the $L$-rational 
and $L^{\sep}$-rational bases in \S\ref{Chap3}.\ref{LRationalBasisSection};  

\item $N \cdot \eta_v \in \NN$ for each $v \in \hS_K$, where the $\eta_v$ are as in (\ref{FPP1});

\item $N$ is large enough that 
\begin{itemize}
\item $N s_i > J$ for each $i = 1, \ldots, m$; 

\item $1 < 2 R_v^N  < \hR_v^N$ for each archimedean $v$ with $K_v \cong \RR$; 

\item $r_v^N  < q_v^{-1/(q_v-1)} < 1$ for each  nonarchimedean $v \in S_{K,0}$. 
\end{itemize} 
\vskip -.5in
\begin{equation} \label{hS_KList1a} \end{equation} 
\vskip .35in
\end{enumerate} 
In particular (\ref{CFBZ1}) holds.

This $N$ will be fixed for the rest of the construction.
\index{patching parameters}

\medskip
{\bf The choice of the Initial Approximating Functions $f_v(z)$.} 
\index{initial approximating functions $f_v(z)$}
\index{coherent approximating functions $\phi_v(z)$!construction when $\Char(K)=0$!choice of the initial approximating functions}
\index{global patching when $\Char(K) = 0$!Stage 2: Constructing approximating functions!the Initial approximating functions} 
We will apply Theorem \ref{CCPX1} with the parameters chosen above.
Let $v_0$ be the archimedean place for which (\ref{CGPF2}) holds.
We first construct the $f_v(z)$ for $v \in \hS_K \backslash \{v_0\}$,
then choose $f_{v_0}(z)$ to `balance' their leading coefficients,
\index{coefficients $A_{v,ij}$!leading}   
so that (\ref{FJIG1}) below will hold.  
 
For each archimedean $v \ne v_0$, 
take $\vbeta_v = (\beta_{v,1}, \ldots, \vbeta_{v,m}) = \vORIG$ in 
Theorem \ref{CCPX1}, and let $f_v(z) \in K_v(\cC)$ 
be the $(\fX,\vs)$-function\index{$(\fX,\vs)$-function} of degree $N$ 
given by Theorem \ref{CCPX1}(A.2) if $K_v \cong \CC$,
or by Theorem \ref{CCPX1}(B.2) with $\cR_v$ as in (\ref{FCCA1}), 
if $K_v \cong \RR$.  
Thus the leading coefficients $c_{v,i}$ of $f_v(z)$ satisfy 
\index{coefficients $A_{v,ij}$!leading}   
$\frac{1}{N} \log_v(|c_{v,i}|_v) = \Lambda_{x_i}(\tE_v,\vs)$ 
for each $i$, for such $v$.       

For each $v \in S_{K,0}$, take $\beta_v = \eta_v$ as before 
(with $0 < \eta_v \in \QQ$ as in (\ref{FPP1}))
and let $f_v(z) \in K_v(\cC)$ be the 
$(\fX,\vs)$-function\index{$(\fX,\vs)$-function} of degree $N$ given by 
Theorem \ref{CCPX1}(C.2) with $\frac{1}{N} \log_v(|c_{v,i}|_v) = \Lambda_{x_i}(\tE_v,\vs) + \beta_v$ 
for each $i$.  For each $v \in \hS_{K,0} \backslash S_{K,0}$, let $f_v(z) \in K_v(\cC)$ 
be the $(\fX,\vs)$-function\index{$(\fX,\vs)$-function} of degree $N$ from Theorem \ref{CCPX1}(D),
with $\frac{1}{N} \log_v(|c_{v,i}|_v) = \Lambda_{x_i}(\tE_v,\vs)$ for each $i$. 

To construct $f_{v_0}(z)$, 
we must first specify $\vbeta_{v_0}$ in Theorem \ref{CCPX1}. 
For each $i$, define $\beta_{v_0,i}$ by  
\begin{equation} \label{FPIG}
\beta_{v_0,i} \log(q_{v_0}) \ = \ \tV_K - (\Gamma(\tEE_K,\fX) \vs)_i \ ;
\end{equation}\index{Robin constant!global}  
then $|\beta_{v_0,i}| < \delta_{v_0} $ by (\ref{CGPF2}). 
The vector $\vbeta_{v_0} := (\beta_{v_0,1}, \ldots, \beta_{v_0,m})$ is $K$-symmetric 
\index{$K$-symmetric!vector}
since $\Gamma(\tEE_K,\fX)$ and $\vs$ are;  this means that $\beta_{v_0,i} = \beta_{v_0,\sigma(i)}$
for each $\sigma \in \Gal(L/K)$.     
Let $f_{v_0}(z) \in K_{v_0}(\cC)$ 
be the $(\fX,\vs)$-function\index{$(\fX,\vs)$-function} of degree $N$
given by Theorem \ref{CCPX1}(A.2) for $\vbeta_{v_0}$ and $\tE_{v_0}$ 
if $K_{v_0} \cong \CC$, 
or by Theorem \ref{CCPX1}(B.2) with $\cR_{v_0}$ as in (\ref{FCCA1}), 
if $K_{v_0} \cong \RR$.  

Thus, for each $v \in \hS_K$ the leading coefficients $c_{v,i}$ of the $f_v(z)$ satisfy 
\index{coefficients $A_{v,ij}$!leading} 
\begin{equation} \label{Ff_vLeading}
\frac{1}{N} \log_v(|c_{v,i}|_v) \ = \ 
\left\{ \begin{array}{ll} 
       \Lambda_{x_i}(\tE_v,\vs)  
               &  \text{if $v \in \hS_K \backslash (S_{K_0} \cup \{v_0\})$,} \\
       \Lambda_{x_i}(\tE_{v_0},\vs) + \beta_{v_0,i}  & \text{if $v = v_0$,} \\
       \Lambda_{x_i}(\tE_v,\vs) + \eta_v & \text{if $v \in \hS_{K,0}$,}  
         \end{array} \right.                    
\end{equation} 
and each $f_v(z)$ has the mapping properties in  Theorem \ref{CCPX1}.

In particular, 
if $v \in S_{K,0}$, then  
$H_v := f_v^{-1}(D(0,1)) \cap E_v$ has a $K_v$-simple decomposition 
\index{$K_v$-simple!decomposition}  
$H_v = \bigcup_{h=1}^N \big(B(\theta_h,\rho_h) \cap \cC_v(F_{u_h})\big)$
compatible with the $K_v$-simple decomposition
\index{$K_v$-simple!decomposition!compatible with another decomposition} 
$E_v = \bigcup_{\ell = 1}^{D_v} \big(B(a_\ell,r_\ell) \cap \cC_v(F_{w_\ell})) \big)$,
which is move-prepared relative to $B(a_1,r_1), \ldots, B(a_{D_v},r_{D_v})$.
\index{move-prepared}  
Here $\theta_1,\ldots, \theta_N$ are the zeros of $f_v(z)$, 
$\rho_h \in |F_{u_h}^{\times}|_v$, and 
$f_v$ induces an $F_{u_h}$-rational scaled isometry from $B(\theta_h,\rho_h)$ to $D(0,1)$ 
\index{scaled isometry} 
which maps $B(\theta_h,\rho_h) \cap \cC_v(F_{u_h})$ onto $\cO_{u_h}$.  
For each $\ell = 1, \ldots, D_v$, there is a point 
$\wbar_\ell \in \big(B(a_\ell,r_\ell) \cap \cC_v(F_{w_\ell})) \big) \backslash H_v$.

\medskip
{\bf Preliminary choice of the Coherent Approximating Functions $\phi_v(z)$.} 
\index{coherent approximating functions $\phi_v(z)$!construction when $\Char(K)=0$!preliminary choice}
\index{global patching when $\Char(K) = 0$!Stage 2: Constructing approximating functions!preliminary choice of the Coherent approximating functions}  
First we will define functions $\phi_v(z) \in K_v(\cC)$ for $v \in \hS_K$, 
and then we will put $\phi_w(z) = \phi_v(z)$ for all $w|v$ 
and consider the leading coefficients of the collection 
\index{coefficients $A_{v,ij}$!leading} 
$\{\phi_w(z) \in L_w(\cC)\}_{w \in \hS_L}$.
For archimedean $v$, the $\phi_v(z)$ will be  modified later
to make the leading coefficients of the $\phi_w(z)$ subunits.
\index{subunit}
\index{coefficients $A_{v,ij}$!leading}   

If $v$ is archimedean, put $\kappa_v = e^{N \eta_v}$, where $\eta_v$ is as in \ref{FPP1}).
Recall from (\ref{FCCA2}) that $K_v \cong \CC$ we have $\hR_v = e^{\eta_v} > R_v$, 
while if $K_v \cong \RR$, we have $\hR_v = \cR_v \cdot e^{\eta_v} > R_v$.  
Thus $\kappa_v > R_v^N$ for each archimedean $v$, 
and $\kappa_v \cdot \cR_v^N = \hR_v^N > R_v^N$ if $K_v \cong \RR$.

If $v$ is nonarchimedean and $v \in S_{K,0}$, put $\kappa_v = 1$.  
If $v \in \hS_{K,0} \backslash S_{K,0}$, put $\kappa_v = \pi_v^{-N \eta_v}$, 
where again
$0 < \eta_v \in \QQ$  as in (\ref{FPP1}).  
Our choice of $N$ required that $N \eta_v \in \NN$,
so $\kappa_v \in K_v^{\times}$ and $|\kappa_v|_v = R_v^N > 1$.   

For each $v \in \hS_K$, define 
\begin{equation*} 
\phi_v(z) \ = \ \kappa_v f_v(z) \ \in \ K_v(\cC)  \ .
\end{equation*}

For each $v$ and each $i$, the leading coefficient $\tc_{v,i}$\label{`SymbolIndextcvi'}  
\index{coefficients $A_{v,ij}$!leading} 
of $\phi_v(z)$ at $x_i$ is given by $\tc_{v,i} = \kappa_v c_{v,i}$. 
By (\ref{Ff_vLeading}) and our choice of the  $\kappa_v$, it follows that 
\begin{equation} \label{FPhi_vLeading} 
\frac{1}{N} \log_v(|\tc_{v,i}|_v) \ = \ 
\left\{ \begin{array}{ll} 
                       \Lambda_{x_i}(\tE_v,\vs) + \eta_v & \text{if $v \ne v_0$} \\
             \Lambda_{x_i}(\tE_{v_0},\vs) + \eta_{v_0} + \beta_{v_0,i} & \text{if $v = v_0$}
                  \end{array} \right.  \ .        
\end{equation} 
Furthermore, the mapping properties of the $f_v(z)$ from Theorem \ref{CCPX1}, 
together with our choice of the $\kappa_v$, 
yield the following mapping properties for the $\phi_v(z)$.   
\begin{enumerate}

\item If $K_v \cong \CC$, 
then $\{z \in \cC_v(\CC) : |\phi_v(z)| \le \hR_v^N\} \subset E_v^0$;  

\item If $K_v \cong \RR$, 
then $\{z \in \cC_v(\CC) : |\phi_v(z)| \le \hR_v^N\} \subset U_v$,
and for each real component $E_{v,j}$ of $E_v$, if $\phi_v(z)$ has $\tau_j$
zeros in $E_{v,j}$ then $\phi_v(z)$ oscillates $\tau_j$ times between $\pm \hR_v^N$ on $E_{v,j}$.

\item If $v \in S_{K,0}$ then properties (a)-(d) in Theorem \ref{CTCX2}(A.3) hold for $\phi_v(z)$.
  Indeed, for $v \in S_{K,0}$ we have $\kappa_v = 1$, so $\phi_v(z) = f_v(z)$ and the mapping
  properties of $\phi_v(z)$ are inherited from those of $f_v(z)$.  
 
\item If $v \in \hS_{K,0} \backslash S_{K,0}$, then 
$E_v = \tE_v = \{z \in \cC_v(\CC_v) : |\phi_v(z)|_v \le R_v^N \}$.  
\end{enumerate} 

\smallskip
{\bf Coherence of the leading coefficients.} 
\index{coefficients $A_{v,ij}$!leading}
\index{coherent approximating functions $\phi_v(z)$!construction when $\Char(K)=0$!adjusting absolue values of leading coefficients}
\index{global patching when $\Char(K) = 0$!Stage 2: Constructing approximating functions!adjusting the leading coefficients}  
In order to view $\fX$ as a subset of $\cC_v(\CC_v)$, 
for each $v$ we have (non-canonically) fixed an embedding of $\tK$ in $\CC_v$ 
(see \S\ref{Chap3}.\ref{AssumptionsSection}),
leading to a distinguished choice of a place $w_v$\index{distinguished place $w_v$} of $L$ above $v$.  
Until now these choices have been a minor concern, 
since all constructions in the proof have been $K_v$-symmetric. 
\index{$K_v$-symmetric}
However, to properly understand the leading coefficients of the $\phi_v(z)$, 
\index{coefficients $A_{v,ij}$!leading} 
we must consider them over the fields $L_w$ for $w \in \hS_L$, 
since $\fX$ is canonically a subset of $\cC(L)$ 
and of $\cC_w(L_w)$ for each $w$, and the
uniformizer $g_{x_i}(z) \in L(\cC)$ is canonically an element of $L_w(\cC)$.  

For each $w \in \hS_L$, put $\phi_w(z) = \phi_v(z)$ if $w|v$, 
and view $\phi_w(z)$ as an element of $L_w(\cC)$.  Although the functions
$\phi_w(z)$ for $w|v$ are all the same, the points of $\fX$, which are their poles, 
are identified differently.  For each $i$ and $w$,   
let $\tc_{w,i} = \lim_{z \rightarrow x_i} \phi_w(z) \cdot g_{x_i}(z)^{Ns_i}$
be the leading coefficient of $\phi_w(z)$ at $x_i$.  
\index{coefficients $A_{v,ij}$!leading} 
Let $\sigma_w : L \hookrightarrow \CC_v$
be an embedding which induces the place $w$, and for each $i = 1, \ldots, m$ let 
$\sigma_w(i)$ be the index $j$ for which $\sigma_w(x_i) = x_j$ (where we identify $x_j$ with 
its image in $\cC_v(\CC_v)$ given by the fixed embedding of $\tK$ in $\CC_v$).  
Then $\tc_{w,i} = \tc_{v,\sigma_w(i)}$.    

The following proposition is the first step towards making the leading coefficients
\index{coefficients $A_{v,ij}$!leading} 
of the patching functions $\hS_L$-units: 
\index{patching functions, initial $G_v^{(0)}(z)$!leading coefficients of}  

\begin{proposition} \label{LeadingCoeffProp}
For each $i = 1, \ldots, m$
\begin{equation} \label{FJIG1} 
\sum_{w \in \hS_L} \log_w(|\tc_{w,i}|_w) \log(q_w) \ = \ 0 \ .
\end{equation} 
\end{proposition} 

To prove this we will need a lemma.  First note that by our normalizations of the 
absolute values on $L$ and $K$, if $w|v$ and $x \in \CC_v \cong \CC_w$, then 
\begin{equation} \label{FJIG2} 
\log_w(|x|_w) \log(q_w) \ = \ [L_w:K_v] \cdot \log_v(|x|_v) \log(q_v) \ .
\end{equation}
Recall that $\Gamma(\tEE,\fX) = \Gamma(\tEE_K,\fX) = \frac{1}{[L:K]} \Gamma(\tEE_L,\fX)$.   
Using (\ref{GreensMatrixSum}), this gives 
\begin{equation} \label{FJIG3}
[L:K] \cdot \Gamma(\tEE_K,\fX) \ = \ 
\sum_{w \in \hS_L} \Gamma(\tE_w,\fX) \log(q_w) \ ,
\end{equation}

\begin{lemma} \label{SumSigmaLem} For each $w \in \hS_L$ and each $i = 1, \ldots, m$,
the $i^{th}$ coordinate of $\Gamma(\tE_w,\fX) \vs$ satisfies  
\begin{equation} \label{FJIG4}
(\Gamma(\tE_w,\fX) \vs)_i \cdot \log(q_w) \ = \ 
[L_w:K_v] \cdot \Lambda_{\sigma_w(x_i)}(\tE_v,\vs) \log(q_v) \ .
\end{equation}  
\end{lemma}  

\begin{proof}  By definition, 
\begin{equation} \label{FJIG5}
(\Gamma(\tE_w,\fX) \vs)_i \cdot \log(q_w)  
\ = \ \sum_{j=1}^m \Gamma(\tE_w,\fX)_{ij} \log(q_w) s_j\ .
\end{equation}  
For each $i \ne j$, it follows from (\ref{FJIG2}) that  
\begin{eqnarray*} 
\Gamma(\tE_w,\fX)_{ij} \log(q_w) & = & G(x_i,x_j;\tE_w) \log(q_w) \\ 
& = & [L_w:K_v] \, G(\sigma_w(x_i),\sigma_w(x_j),\tE_v) \log(q_v) \ . 
\end{eqnarray*} 
\index{Green's function!nonarchimedean}
Similarly, since the uniformizers have been chosen in such a way that 
$\sigma_{w}(g_{x_i})(z) = g_{\sigma_w(x_i)}(z)$, for each $i$ 
\begin{equation*} 
\Gamma(\tE_w,\fX)_{ii} \log(q_w) \ = \ V_{x_i}(\tE_w) \log(q_w) 
\ = \ [L_w:K_v] V_{\sigma_w(x_i)}(\tE_v) \log(q_v) \ .
\end{equation*} 
\index{Robin constant!local}  
For compactness of notation, write 
\begin{equation*}
\tG(x_i,x_j;E_v) \ = \ \left\{ 
\begin{array}{ll} G(x_i,x_j;\tE_v) & \text{if\ $i \ne j$\ ,} \\
                  V_{x_i}(\tE_v)   & \text{if\ $i = j$\ .} 
\end{array} \right .
\end{equation*} 
\index{Green's function}
\index{Robin constant!local}
Since $\vs$ is $K$-symmetric, we have $s_j = s_{\sigma_w(j)}$ for each $j$.  Hence  
\index{$K$-symmetric!probability vector}    
\begin{eqnarray*}
(\Gamma(\tE_w,\fX) \vs)_i \cdot  \log(q_w) 
& = & \sum_{j=1}^m \Gamma(\tE_w,\fX)  s_j \log(q_w) \\ 
& = & [L_w:K_v] \cdot \sum_{j=1}^m  \tG(\sigma_w(x_i),\sigma_w(x_j);\tE_v) s_{\sigma_w(j)}  \log(q_v) \\ 
& = & [L_w:K_v] \cdot \Lambda_{\sigma_w(x_i)}(\tE_v,\vs) \log(q_v) 
\end{eqnarray*}
as desired.   
\end{proof} 

\vskip .1 in
\begin{proof}[Proof of Proposition \ref{LeadingCoeffProp}]
For each $v$, and each $w|v$, fix an embedding $\sigma_w : L \hookrightarrow \CC_v$
which induces the place $w$.  Since  $\tc_{w,i} = \tc_{v,\sigma_w(i)}$, 
it follows from (\ref{FPhi_vLeading}) that   
\begin{eqnarray}
\frac{1}{N} \sum_{w \in \hS_L} \log_w(|\tc_{w,i}|_w) \log(q_w) 
     & = & \sum_{v \in \hS_K} \sum_{w|v} [L_w:K_v] 
               \big(\frac{1}{N}\log_v(|\tc_{v,\sigma_w(x_i)}|_v)\big) \log(q_v) \notag \\
     & = & \sum_{v \in \hS_K} \sum_{w|v} [L_w:K_v] 
        \big(\Lambda_{\sigma_w(x_i)}(\tE_v,\vs) + \eta_v\big) \log(q_v) 
         \notag \\
     & & \qquad + \ \sum_{w|v_0} [L_w:K_v] \beta_{v_0,\sigma_w(i)} \log(q_{v_0}) \ .
         \label{BigSumFormula} 
\end{eqnarray} 
By Lemma \ref{SumSigmaLem} 
\begin{eqnarray}
\sum_{v \in \hS_K} \sum_{w|v} [L_w:K_v] \Lambda_{\sigma_w(x_i)}(\tE_v,\vs) \log(q_v) 
& = & \sum_{w \in \hS_L} (\Gamma(\tE_w,\fX) \vs)_i \log(q_w) \notag \\
& = & (\Gamma(\tEE_L,\fX) \vs)_i \ = \ [L:K]  \cdot (\Gamma(\tEE_K,\fX) \vs)_i \ .
         \label{FLambdaSum}
\end{eqnarray}  
By our choice of the $\eta_v$ in (\ref{FPP1}), 
\begin{equation}
\sum_{v \in \hS_K} \sum_{w|v} [L_w:K_v] \eta_v \log(q_v) 
\ = \ [L:K] \sum_{v \in \hS_K} \eta_v \log(q_v) \ = \ -[L:K] \cdot \tV_K \ .
\label{FEtaSum}
\index{Robin constant!global}
\end{equation}
Finally, since $\vbeta_{v_0}$ is $K$-symmetric, 
\index{$K$-symmetric!vector}
for each $w|v_0$ we have $\beta_{v_0,\sigma_w(i)} = \beta_{v_0,i}$, 
so by (\ref{FPIG})    
\begin{equation}  
\sum_{w|v_0} [L_w:K_v] \beta_{v_0,\sigma_w(i)} \log(q_{v_0})  
 \ = \ [L:K] \cdot \big(\tV_K - (\Gamma(\tEE_K,\fX) \vs)_i \big) \ . \label{FBetaSum}
\end{equation}
\index{Robin constant!global}
Combining (\ref{BigSumFormula}), (\ref{FLambdaSum}), (\ref{FEtaSum}) and (\ref{FBetaSum}) 
gives 
\begin{equation*}
\sum_{w \in \hS_L} \log_w(|\tc_{w,i}|_w) \log(q_w) \ = \ 0  
\end{equation*} 
as required.  
\end{proof} 

\vskip .1 in
{\bf Adjusting the leading coefficients to be $\hS_L$-subunits.}
\index{coherent approximating functions $\phi_v(z)$!construction when $\Char(K)=0$!making the leading coefficients $S$-subunits}
\index{subunit}\index{coefficients $A_{v,ij}$!leading}
\index{global patching when $\Char(K) = 0$!Stage 2: Constructing approximating functions!the Coherent approximating functions}
The final step in the proof of Theorem \ref{CTCX2} involves modifying the archimedean
$\phi_v(z)$ so that their leading coefficients become $\hS_L$-subunits.
\index{subunit}

By our choices of the $R_v$, $\hR_v$, and $N$ we have $R_v^N < \hR_v^N$ for each archimedean $v$, 
and $2R_v^N < \hR_v^N$ for each archimedean $v$ with $K_v \cong \RR$. 
For each $i = 1, \ldots, m$, put $N_i = Ns_i \in \NN$. 
Noting that our choice of $N$ has required that $N_i > J$,
let $\phi_{i,N_i}(z)$ be the corresponding function from the $L$-rational basis.  
\index{basis!$L$-rational} 
 
If $K_v \cong \CC$, the construction of $\phi_v(z)$ has arranged that  
\begin{equation*}
\{z \in \cC_v(\CC_v) : |\phi_v(z)| \le \hR_v^N  \} \ \subset  \ U_v \ ,
\end{equation*}  
Since $R_v^N < \hR_v^N$, by continuity there is a $\delta_v^{\prime} > 0$
such that if $\Delta_{v,1}, \ldots, \Delta_{v,m} \in \CC$ satisfy
$|\Delta_{v,i}| < \delta_v^{\prime}$ for each $i$, and if $\phi_v(z)$ is replaced by
$\widehat{\phi}_w(z) = \phi_w(z) + \sum_{i=1}^m \Delta_{w,i} \varphi_{i,N_i}(z)$, then
\begin{equation} \label{CFFZ1}
\{ z \in \cC_v(\CC_v) : |\widehat{\phi}_w(z)|_v \le R_v^N \} \ \subset  \ U_v \ .
\end{equation}

If $K_v \cong \RR$, the construction of $\phi_v(z)$ has arranged that
\begin{equation*}
\{z \in \cC_v(\CC_v) : |\phi_v(z)| \le 2 \hR_v^N  \} \ \subset  \ U_v \ ,
\end{equation*}  
and that for each component $E_{v,i}$ of $E_v$ contained in $\cC_v(\RR)$,
if $\phi_v(z)$ has $\tau_i$ zeros in $E_{v,i}$ then it oscillates $\tau_i$ times between
$\pm 2 \hR_v^N$ on $E_{v,i}$.   
Since $2 R_v^N < 2 \hR_v^N$, by continuity there is a $\delta_v^{\prime} > 0$ such that if  
$\Delta_{v,1}, \ldots, \Delta_{v,m} \in \CC$ are a $K_v$-symmetric set of numbers 
\index{$K_v$-symmetric!set of numbers}
with $|\Delta_{v,i}| < \delta_v^{\prime}$ for each $i$, 
and if $\phi_v(z)$ is replaced by
$\widehat{\phi}_v(z) =  \phi_v(z) + \sum_{i=1}^m \Delta_{v,i} \varphi_{i,N_i}(z)$, then
\begin{equation} \label{CFFZ2}
\{ z \in \cC_v(\CC_v) : |\widehat{\phi}_v(z)| \le 2R_v^N \} \ \subseteq \ U_v   
\end{equation}
and for each component $E_{v,i}$ contained in $\cC_v(\RR)$,
if $\phi_v(z)$ has $\tau_i$ zeros in $E_{v,i}$ then $\widehat{\phi}_v(z)$ 
oscillates $\tau_i$ times between $\pm 2 R_v^N$ on $E_{v,i}$.

Let $\delta^{\prime}$ be the minimum of the
$\delta_v^{\prime}$, for all $v \in \hS_{K,\infty}$.

\vskip .05 in
Fix $x_i \in \fX$, and put $F = K(x_i)$.  
Let $\hS_F$ be the set of places of $F$ above $\hS_K$.  
For each $v \in \hS_K$, since the $\phi_w(z) \in K_v(\cC)$
are the same for all $w|v$, Proposition \ref{CPPr1} tells us that
$\oplus_{w|v} \tc_{w,i} \in \oplus_{w|v} L_w$ actually belongs to 
$\oplus_{u|v} F_u$, embedded semi-diagonally in  $\oplus_{w|v} L_w$.
Write $\oplus_{u|v} \tc_{u,i}$ for the element of $\oplus_{u|v} F_u$
that induces it.  By (\ref{FJIG1})
\begin{equation*}
\sum_{u \in \hS_F} \log_u(|\tc_{u,i}|_u) \log(q_u)
   \ = \ \frac{1}{[L:F]} \sum_{w \in \hS_L} \log_w(|\tc_{w,i}|_w) \log(q_w)
   \ = \ 0 \ .
\end{equation*}
According to Proposition \ref{CPSUnit2FF} 
there are an $\hS_F$-unit $\mu_i \in F$, an integer $n_i$, and an $\hS_F$-subunit
\index{subunit}
$\oplus_{u \in \hS_{F,\infty}} \varepsilon_{u,i} \in \bigoplus_{u \in \hS_{F,\infty}} F_u^{\times}$
such that $\varepsilon_{u,i}^{n_i} = \mu_i$ for each $u \in \hS_{F,\infty}$, and 
\begin{equation*}  
\left\{ \begin{array}{ll} 
      |\tc_{u,i}-\varepsilon_{u,i}| < \delta^{\prime} \ , & \\
      |\tc_{u,i}^{n_i}|_u = |\mu_i|_u \ , & \text{for each  $u \in \hS_{F,0}$\ .}
        \end{array}  \right.
\end{equation*} 
Since $\mu_i$ is an $\hS_F$-unit (hence also an $\hS_L$-unit), 
$\sum_{w \in \hS_L} \log_u(|\mu_i|_w) \log(q_w) = 0$.  

For each archimedean $u$ and each $w|u$, 
put $\varepsilon_{w,i} = \varepsilon_{u,i}$.  
Then $\log_w(|\varepsilon_{w,i}|_w) = \frac{1}{n_i} \log_w(|\mu_i|_w)$
for each archimedean $w \in S_L$, and 
$\log_w(|c_{w,i}|_w) = \frac{1}{n_i} \log_w(|\mu_i|_w)$ for each nonarchimedean
$w \in S_L$.  It follows that
\begin{equation} \label{CFFF3}
\sum_{w \in S_{L,\infty}} \log_w(\varepsilon_{w,i}) \log(q_w)
+ \sum_{w \in S_{L,0}} \log_w(\tc_{w,i}) \log(q_w) \ = \ 0 \ .
\end{equation}

Now let $x_i$ vary. We next arrange for the $\mu_i$ and 
$\vec{\varepsilon}_{v,i} = \oplus_{w|v} \varepsilon_{w,i}$ to be $K$-symmetric.
\index{$K$-symmetric!vector}
By Proposition \ref{CPPr1}, the $\oplus_{w|v} \tc_{w,i}$ are $K$-symmetric. 
\index{$K$-symmetric!vector}
For each $\Gal(L/K)$-orbit $\fX_{\ell} \subset \fX$, 
fix an $x_i \in \fX_{\ell}$ and put $F = K(x_i)$ as before.  By construction,  
$\vec{\varepsilon}_{v,i} := \oplus_{w|v} \varepsilon_{w,i} \in L \otimes_K K_v$ 
belongs to $F \otimes_K K_v$ for each archimedean $v$.  For each $x_j \in \fX_{\ell}$,
choose $\sigma \in \Gal(L/K)$ with $\sigma(x_i) = x_j$, 
and replace $\mu_j$ with $\sigma(\mu_i)$,  
$\vec{\varepsilon}_{v,j}$ with $\sigma(\vec{\varepsilon}_{v,i})$.  Since
$\vec{\varepsilon}_{v,i} \in F \otimes_K K_v$, these objects are independent of
the choice of $\sigma$ with $\sigma(x_i) = x_j$, and are $K$-symmetric.
\index{$K$-symmetric!vector}  

After replacing the $\mu_i$ with powers of themselves, 
we can assume there is a number $n_0$ such that $n_i = n_0$, for all $i$. 
The numbers $\mu_1, \ldots, \mu_m$ form a $K$-symmetric system of $\hS_L$-units,
\index{$K$-symmetric!system of units}
and the $\varepsilon_{w,i}$ form a $K$-symmetric system of $\hS_L$-subunits.
\index{subunit}
\index{$K$-symmetric!system of subunits}
\index{subunit}
For each archimedean $w$ and each $i$, put
\begin{equation*}
\Delta_{w,i} = \varepsilon_{w,i} - c_{w,i}
\end{equation*}
and put 
$\widehat{\phi}_w(z) = \phi_w(z) + \sum_{i=1}^m \Delta_{w,i} \varphi_{i,N_i}(z)$.
Since $\oplus_{w|v} \Delta_{w,i} \varphi_{i,N_i} \in L \otimes_K K_v(\cC)$ 
and $\oplus_{w|v} \phi_w(z) \in L \otimes_K K_v(\cC)$ are $K$-symmetric,
\index{$K$-symmetric!vector} 
Proposition \ref{CPPr1} shows that $\widehat{\phi}_w(z)$ belongs to $K_v(\cC)$ for each $w|v$, 
and that $\widehat{\phi}_{w_1}(z) = \widehat{\phi}_{w_2}(z)$ 
for all $w_1, w_2|v$.   

Replace $\phi_w(z)$ with $\widehat{\phi}_w(z)$, for each archimedean $v$ and each $w|v$.
The leading coefficients of the new $\phi_w(z)$ 
\index{coefficients $A_{v,ij}$!leading} 
are the $\varepsilon_{w,i}$, so we can put $\phi_v(z) = \widehat{\phi}$, for any $w|v$.
By (\ref{CFFF3}), assertion (B) in the Theorem holds. 
Our construction has established assertions (A1) -- (A4), 
so the proof of Theorem \ref{CTCX2} is complete.
\end{proof}



{\bf Stage 3. The Patching Construction.}
\index{patching theorem!global patching construction when $\Char(K) = 0$|ii}
\index{global patching when $\Char(K) = 0$!Stage 3: The global patching construction!overview}

\smallskip
{\bf Overview.}  
\index{patching argument!global|ii} 
\index{patching argument!local|ii} 
The patching process has two parts, a global part and a local part. 
The global part concerns the way the patching coefficients 
\index{coefficients $A_{v,ij}$} 
are chosen, managing them so as to achieve global $K$-rationality for the final patched function.    
The local part is responsible for assuring $K_v$-rationality 
of the partially patched functions, and confining their roots to  $E_v$.  
\index{confinement argument} 

Although this description separates the roles of the global and local
parts of patching process, in fact the two interact,
and the coefficients are determined recursively, from highest to lowest order.  
Each local patching construction specifies certain parameters to the global patching process:  
the number of patching stages it considers high-order  
and bounds for the size of the patching coefficients it can handle.    
As patching is carried out, and high-order coefficients chosen by the global
\index{coefficients $A_{v,ij}$!high-order} 
process are achieved by the local process, lower-order coefficients 
are changed as a result.  The global process must take these changes 
into account in determining subsequent coefficients.  

\vskip .07 in
The patching process begins with the coherent approximating functions 
\index{patching argument!global}\index{coherent approximating functions $\phi_v(z)$}
$\{\phi_v(z)\}_{v \in \hS_K}$ given by Theorem \ref{CTCX2}.  
Its goal is to produce a function $G(z) \in K(\cC)$ independent of $v$, 
of much higher degree than the $\phi_v(z)$, whose zeros are points with the 
properties in Theorem \ref{aT1-B}.  

The first step is to compose each $\phi_v(z)$ with a 
``degree-raising polynomial'' $Q_{v,n}(x) \in K_v(x)$.\index{degree-raising|ii}  
The $Q_{v,n}(x)$ are monic, of common degree $n$.  
This allows the leading coefficients to be patched to become $\hS_L$-units,
\index{coefficients $A_{v,ij}$!leading}
and makes the degree large enough that certain analytic estimates are satisfied,
while keeping the roots in $E_v$.     

For each $v \in \hS_K$, the local patching process provides $Q_{v,n}(z)$.  
\index{patching argument!local} 
If $K_v \cong \CC$ or if $K_v$ is nonarchimedean and $E_v$ is an $\RL$-domain,
\index{$\RL$-domain} 
then $Q_{v,n}(x) = x^n$.  
If $K_v \cong \RR$, then $Q_{v,n}(x)$ is a composite of two Chebyshev polynomials.
\index{Chebyshev polynomial}   
If $K_v$ is nonarchimedean and $E_v \subset \cC(K_v)$, then $Q_{v,n}(x)$
is the Stirling polynomial of degree $n$ for the ring of integers $\cO_v$.  
  \index{Stirling polynomial!for $\cO_v$}
For appropriately large and divisible $n$, 
this yields the ``initial patching functions'' 
$G_v^{(0)}(z) = Q_{v,n}(\phi_v(z))$.\index{patching functions, initial $G_v^{(0)}(z)$!construction of} 
Although $N$ is likely quite large, 
$n$ should be thought of as astronomically larger than $N$.  
\index{$n$ astronomically larger than $N$}  

Each $G_v^{(0)}(z)$ can be expanded in terms of the 
$L$-rational basis functions\index{patching functions, $G_v^{(k)}(z)$ for $1 \le k \le n$!expansion of}
\index{basis!$L$-rational}  
$\varphi_{i,j}(z)$ and $\varphi_{\lambda}$, 
with $L_{w}$-rational coefficients for each $w|v$. 
\index{coefficients $A_{v,ij}$} 
For notational purposes, it will be useful to deem the 
basis functions $\varphi_{\lambda}(z)$ for $\lambda \le \Lambda_0$ 
and $\varphi_{i,j}(z)$ with $J < j \le N_i := Ns_i$, as being ``low-order'',
and list them as $\varphi_{\lambda}$, $\lambda = 1, \ldots, \Lambda$.\label{`SymbolIndexLambda'}   
Thus, for each $v$\label{`SymbolIndexAvij'}\index{patching functions, $G_v^{(k)}(z)$ for $1 \le k \le n$!expansion of} 
\begin{equation*}
G_v^{(0)}(z) \ = \
    \sum_{i=1}^m \sum_{j=0}^{(n-1)N_i-1} A_{v,ij} \varphi_{i,nN_i-j}(z)
       + \sum_{\lambda = 1}^{\Lambda} A_{v,\lambda} \varphi_{\lambda} \ .
\end{equation*}

The patching process initially adjusts the leading coefficients of 
\index{patching argument!global} 
\index{coefficients $A_{v,ij}$!leading}
the $G_v^{(0)}(z)$ to be global $S_L$-units, independent of $v$. 
Then, in stages,\index{patching functions, $G_v^{(k)}(z)$ for $1 \le k \le n$!constructed by patching} 
it inductively constructs functions $G_v^{(1)}(z), \ldots, G_v^{(n)}(z)$,\label{`SymbolIndexGvk'} 
where\index{compensating functions $\vartheta_{v,ij}^{(k)}(z)$}
\index{patching functions, $G_v^{(k)}(z)$ for $1 \le k \le n$!constructed by patching}
\begin{equation*}
G_v^{(k)}(z) \ =  \ G_v^{(k-1)}(z)
    + \sum_{i=1}^m \sum_{j=(k-1)N_i}^{kN_i-1} 
            \Delta_{v,ij}^{(k)} \vartheta_{v,ij}^{(k)}(z)
\end{equation*}
for $1 \le k \le n-1$ (see Proposition \ref{LocalPatch} for a more precise statement), 
and\index{patching functions, $G_v^{(k)}(z)$ for $1 \le k \le n$!constructed by patching} 
\begin{equation*}
G_v^{(n)}(z) \ = \ G_v^{(n-1)}(z)
     + \sum_{\lambda = 1}^{\Lambda}  \Delta_{v,\lambda}^{(n)} \varphi_{\lambda}(z) 
\end{equation*}
for $k = n$, in such a way that for each $v$ 
the coefficients $A_{v,ij}$, $A_{v,\lambda}$ 
\index{coefficients $A_{v,ij}$}
are changed into global $S_L$-integers $A_{ij}$, $A_{\lambda}$ independent of $v$.
This process is called ``patching'', because
\index{patching|ii} 
it pieces together a global function out of a collection of local ones.

The ``compensating functions'' $\vartheta_{v,ij}^{(k)}(z)$,\label{`SymbolIndexComp'} 
\index{compensating functions $\vartheta_{v,ij}^{(k)}(z)$|ii}
indexed by pairs $(i,j)$ with $(k-1)N_i \le j \le kN_i-1$ in ``bands'' for $k = 1, \ldots, n-1$, 
are determined by the local patching process and have poles are supported on $\fX$.
\index{patching argument!local} 
Each $\vartheta_{v,ij}^{(k)}(z)$
\index{compensating functions $\vartheta_{v,ij}^{(k)}(z)$!poles and leading coefficients of} 
has a pole of exact order $nN_i-j$ at $x_i$,  
and is chosen so that adding 
$\Delta_{v,ij}^{(k)} \vartheta_{v,ij}^{(k)}(z)$
\index{compensating functions $\vartheta_{v,ij}^{(k)}(z)$} 
to $G_v^{(k)}(z)$\index{patching functions, $G_v^{(k)}(z)$ for $1 \le k \le n$!constructed by patching} 
affects only poles of order $nN_i-j$ and below at $x_i$, 
and lower order poles outside the band, for $x_{i^{\prime}} \ne x_i$.  
It was Fekete and Szeg\"o's insight (\cite{F-SZ}) 
\index{Fekete, Michael}
\index{Szeg\"o, G\'abor}
that by using compensating functions more complicated 
\index{compensating functions $\vartheta_{v,ij}^{(k)}(z)$!more complicated than basis functions} 
than the basis functions $\varphi_{ij}(z)$,  
\index{basis!$L$-rational} 
one could control movement of the roots of the 
$G_v^{(k)}(z)$.\index{patching functions, $G_v^{(k)}(z)$ for $1 \le k \le n$!roots are confined to $E_v$}
For each $x_i$, the $\Delta_{v,ij}^{(k)}$ are chosen by ascending $j$ 
(decreasing order of the pole), so that coefficients $A_{ij}$
\index{coefficients $A_{v,ij}$} 
already patched are not changed in subsequent steps.

\vskip .1 in
The global patching process has two concerns.\index{patching argument!global} 

First, it must choose the $\Delta_{v,ij}^{(k)}$ in such a way that for each $(i,j)$, 
\begin{equation*}
A_{ij} \ = \ A_{v,ij} + \Delta_{v,ij}^{(k)}
\end{equation*}\label{`SymbolIndexDeltavij'}
is a global $\hS_L$ integer independent of $v$.  
This is accomplished by extending the base to $L$ 
and choosing the $A_{ij}$ via Proposition \ref{CStrong3},
\index{Strong Approximation theorem!uniform}
simultaneously patching the coefficients for
\index{coefficients $A_{v,ij}$} 
all $x_i$ belonging to a given galois orbit in $\fX$.

Second, it must impose conditions on the sizes of the
$\Delta_{v,ij}^{(k)}$ so that the local patching constructions can succeed.\index{patching argument!local}  
The choice of the $\Delta_{v,ij}^{(k)}$ involves tension between the global 
and local parts of the patching process.  The global part
\index{patching argument!tension between local and global}  
is charged with adjusting the $A_{v,ij}$ to make them algebraic numbers in $L$ independent of $v$.  
Doing so may require the $\Delta_{v,ij}^{(k)}$ to be fairly large.  
On the other hand, the local part is charged with assuring 
that the roots of $G_v^{(k)}(z)$\index{patching functions, $G_v^{(k)}(z)$ for $1 \le k \le n$!roots are confined to $E_v$} remain in $E_v$.  
For this, it is usually necessary that the $\Delta_{v,ij}^{(k)}$ be fairly small.  

If there is a bound $B_v > 0$ such that in the $k$-th stage of 
the local patching construction the $\Delta_{v,ij}^{(k)} \in L_v$ 
\index{patching argument!global} 
can be chosen arbitrarily, provided that $|\Delta_{v,ij}^{(k)}|_v \le B_v$ for all $i$, $j$, 
and the $\Delta_{v,ij}^{(k)}$ are $K_v$-symmetric, 
\index{$K_v$-symmetric!set of numbers}
we will say that the coefficients $A_{v,ij}$  for $(k-1)N_i \le j < kN_i$ can be 
\index{coefficients $A_{v,ij}$}
{\em sequentially patched with freedom $B_v$}.  Equivalently, we will say that 
{\em the $k$-th stage of the local patching process at $v$ can be carried out with freedom $B_v$}.
\index{patching argument!global}\index{freedom $B_v$ in patching|ii}   
 
As $k$ increases, there is greater and greater freedom in the patching process.  
However, for small $k$, balancing the demands of the global and local 
\index{patching argument!tension between local and global} 
patching constructions requires care.  The leading coefficients are the hardest to patch, 
\index{coefficients $A_{v,ij}$!leading} 
and they are controlled through the choice of $n$.  The high order coefficients are also
\index{coefficients $A_{v,ij}$!high-order} 
quite difficult to patch.  It turns out that there is a number $\kbar$,  
determined by the sets $E_v$ and initial approximating functions $f_v$ 
\index{initial approximating functions $f_v(z)$}
but fortunately independent of $n$, such that when $1 \le k \le \kbar$, 
the nonarchimedean $\Delta_{v,ij}^{(k)}$ must be very small.  
To compensate,  we must allow the archimedean $\Delta_{v,ij}^{(k)}$ 
to be quite large.  

The archimedean patching procedures accomplish this    
by exploiting a phenomenon of `magnification' introduced in (\cite{RR3}), 
by which small changes in $f_v(z)$\index{magnification argument|ii} 
create large changes in the leading coefficients of $G_v(z)$. 
\index{coefficients $A_{v,ij}$!leading}  
It is shown in Theorems \ref{Chap8}.\ref{DCPatch1} and \ref{Chap9}.\ref{DRPatch1}
that for any fixed $B_v > 0$, magnification enables us\index{magnification argument} 
to carry out the first $\kbar$ stages of the patching process with freedom $B_v$.  
\index{freedom $B_v$ in patching}  
 
In combination, the local and global patching processes determine 
\index{patching argument!global}\index{patching argument!local} 
$\kbar$, the number of patching stages deemed high order.
For appropriate numbers $B_v$, we will have  
\begin{equation*}
|\Delta_{v,ij}^{(k)}|_v \ \le \ 
   \left\{ \begin{array}{ll} B_v & \text{if $k \le \kbar$,} \\
                        h_v^{kN} & \text{if $k > \kbar$.} \end{array} \right. 
\end{equation*}
If these conditions are met, the local patching constructions will succeed.
\index{patching argument!local}   
On the other hand, for global target coefficients $A_{ij} \in L$ to exist, 
\index{coefficients $A_{v,ij}$!target} 
$(\prod_v B_v^{D_v})$ and $(\prod_v h_v^{D_v})^{\kbar}$
must be large enough that Proposition \ref{CStrong3} applies. 
\index{Strong Approximation theorem!uniform} 
Achieving this uses condition (\ref{FPP3}) that $\prod_v h_v^{D_v} > 1$,  
which ultimately depends on the fact that $\gamma(\EE,\fX) > 1$.     

The final patched functions $G_v^{(n)}(z)$\index{patching functions, $G_v^{(k)}(z)$ for $1 \le k \le n$!are $K_v$-rational} are $K_v$-rational 
but have all their coefficients in $L$.\index{coefficients $A_{v,ij}$} 
By Proposition \ref{CPPr1} there is a global function $G^{(n)}(z) \in K(\cC)$, independent of $v$,
such that $G^{(n)}(z) = G_v^{(n)}(z)$
\index{patching functions, $G_v^{(k)}(z)$ for $1 \le k \le n$!$G_v^{(n)}(z) = G^{(n)}(z)$ is independent of $v$} 
for each $v \in \hS_K$.  
The local patching constructions assure that its zeros belong to \index{patching argument!local} 
$E_v$ for all $v \in \hS_K$.  For each $v \notin \hS_K$, 
the coefficients of $G^{(n)}(z)$ are $\hS_L$-integers  
and its leading coefficients are $\hS_L$-units, so the fact that the
\index{coefficients $A_{v,ij}$!leading} 
basis functions $\varphi_{ij}(z)$ and $\varphi_{\lambda}$ have good reduction outside $\hS_K$,
\index{good reduction}
\index{basis!$L$-rational}  
combined with the fact that $E_v$ is $\fX$-trivial,  
\index{$\fX$-trivial}
show that $\{z \in \cC_v(\CC_v) : |G^{(n)}(z)|_v \le 1\} = E_v$. 
Thus the zeros of $G^{(n)}(z)$ belong to $E_v$ for all $v$.   

\medskip
{\bf Details.}
\index{global patching when $\Char(K) = 0$!Stage 3: The global patching construction!details}
We now give the details of the patching  construction.  
Let $K$, $S_K$, $\hS_K$, $\EE$, $\widetilde{\EE}$, 
\index{patching argument!global|ii} 
and the sets $U_v$ for $v \in \hS_{K,\infty}$ 
be as in Stage 1 of the proof.
Let $\vs \in \cP^m(\QQ)$ be the $K$-symmetric vector with positive rational
\index{$K$-symmetric!vector}
coefficients from (\ref{CGPF2}), and let $h_v$, $r_v$, and $R_v$  
be the local patching parameters from (\ref{FPP2}), with
\index{patching parameters}
$1 < h_v < r_v < R_v$ for archimedean $v$ 
and $0 < h_v < r_v \le 1 \le R_v$ for nonarchimedean $v$.  
Let the natural number $N$, 
the coherent $(\fX,\vs)$-functions\index{$(\fX,\vs)$-function}
 $\{\phi_v(z)\}_{v \in \hS_K}$\index{coherent approximating functions $\phi_v(z)$} 
of degree $N$, and the $\hS_K$-units $\mu_i$ from Theorem \ref{CTCX2}, 
be as Stage 2 of the proof.  

\smallskip
For each $v$, let $w_v$\index{distinguished place $w_v$} be the distinguished place of $L = K(\fX)$ over $v$, 
induced by the embedding  $\tK \hookrightarrow \CC_v$ chosen in \S\ref{Chap3}.\ref{AssumptionsSection}.   
This induces an embedding $L_{w_v} \hookrightarrow \CC_v$,
and allows us to identify $\fX$ with a subset of $\cC_v(\CC_v)$.  We will 
use these embeddings in comparing coefficients of functions over $K$ and over $L$.  
\index{coefficients $A_{v,ij}$}   

\smallskip
{\bf The order $\prec_N$.}  
\index{order!$\prec_N$|ii}\label{`SymbolIndexPrecN'}
\index{global patching when $\Char(K) = 0$!Stage 3: The global patching construction!the order $\prec_N$}
We will now define an ordering $\prec_N$ on the index set 
$\cI = \{(i,j) \in \ZZ^2 : 1 \le i \le m, 0 \le j < \infty\}$\label{`SymbolIndexI'}  
which specifies the sequence in which the coefficients are patched.  

Let $\Gal(L/K)$ act on $\cI$ in a $K$-symmetric way 
\index{$K$-symmetric!index set}
through its first coordinate, so $\sigma(i,j) = (\sigma(i),j)$ 
if $\sigma(x_i) = x_{\sigma(i)}$.      
With $\vs$ as in (\ref{CGPF2}) and $N$ as in Theorem \ref{CTCX2}, 
put $N_i = Ns_i$ for $i = 1, \ldots, m$. 
For each $(i,j) \in \cI$, we can uniquely write $j = (k-1)N_i + r$ 
with $k, r \in \ZZ$, $k \ge 1$ and $0 \le r < N_i$;  
put $k_N(i,j) = k$ and $r_N(i,j) = r$.   
Let the $\Gal(L/K)$-orbits in $\fX = \{x_1, \ldots, x_m\}$ 
be $\fX_1, \ldots, \fX_{m_1}$.  Without loss, 
we can assume the $x_i$ in a given $\fX_{\ell}$ have consecutive indices.  
If $x_i \in \fX_{\ell}$, put $\ell(i,j) = \ell$. 
 
Let $\prec_N$ be the total order on $\cI$ defined by 
\index{order!$\prec_N$|ii} 
$(i_1,j_1) \prec_N (i_2,j_2)$ iff
\begin{equation} \label{FPrec} 
\left\{ \begin{array}{l} 
    \text{$j_1 = j_2 = 0$ and $i_1 < i_2$, \quad or}  \\
    \text{$k_N(i_1,j_1) < k_N(i_2,j_2)$, \quad \ or} \\
    \text{$k_N(i_1,j_1) = k_N(i_2,j_2)$, $\max(j_1,j_2) \ge 1$
                              and $\ell(i_1,j_1) < \ell(i_2,j_2)$, 
                                         \quad or} \\                                  
    \text{$k_N(i_1,j_1) = k_N(i_2,j_2)$, $\ell(i_1,j_1) = \ell(i_2,j_2)$, 
                 and $j_1 <  j_2$, \quad or}  \\
     \text{$k_N(i_1,j_1) = k_N(i_2,j_2)$, $\ell(i_1,j_1) = \ell(i_2,j_2)$, 
                 $j_1 =  j_2 \ge 1$, and $i_1 < i_2.$}     
        \end{array} \right.        
\end{equation}
Write $(i_1,j_1) \preceq_N (i_2,j_2)$ 
iff $(i_1,j_1) \prec_N (i_2,j_2)$ or $(i_1,j_1) = (i_2,j_2)$.  
Define the ``bands'' of $\prec_N$, for $k =1, 2, \ldots $ by
\index{band!$\Band_N(k)$|ii}\label{`SymbolIndexBandN'}
\begin{equation} \label{FBand} 
\Band_N(k) \ = \ \{ (i,j) \in \cI : k_N(i,j) = k \} \ .
\end{equation} 
Note that the indices $(i,0)$ for the leading coefficients
\index{coefficients $A_{v,ij}$!leading} 
form the initial segment under $\prec_N$, and are contained in $\Band_N(1)$.
\index{order!$\prec_N$}
\index{band!$\Band_N(k)$} 
 
Let $\cong_N$ be the equivalence relation on $\cI$ 
defined by $(i_1,j_1) \cong_N (i_2,j_2)$ iff  
\begin{equation*} 
\left\{ \begin{array}{l} 
\text{$j_1 = j_2$, \ and \ } \\
\text{$x_{i_1}$, $x_{i_2}$ belong to the same 
                        galois orbit $\fX_{\ell}$. }
        \end{array} \right.
\end{equation*}   
Equivalently, $(i_1,j_1) \cong_N (i_2,j_2)$ iff $\sigma(\varphi_{i_1,j_1}) = \varphi_{i_2,j_2}$
for some $\sigma \in \Gal(L/K)$.  
Define the ``galois blocks'' of $\cI$ to be the equivalence classes for  $\cong_N$, 
and write\label{`SymbolIndexBlockN'}  
\begin{equation*} 
\Block(i,j) \ = \ \{ (i_1,j_1) \in \cI : (i_1,j_1) \cong_N (i,j) \} 
\ = \ \Gal(L/K)\big((i,j)\big) \ .
\end{equation*} 
\index{block, galois!$\Block(i,j)$|ii}

In patching, coefficients will be adjusted in $\prec_N$ order. 
\index{patching argument!global}\index{order!$\prec_N$}  
\index{coefficients $A_{v,ij}$} 
This means that the leading coefficients are modified first, 
\index{coefficients $A_{v,ij}$!leading} 
then the remaining coefficients are considered band by band. 
\index{band!coefficients patched by bands}
Within each band, they are considered block by block.
\index{block, galois!coefficients patched block by block}   
For each $i$, they are considered by increasing $j$.  The global patching process simultaneously
determines all the coefficients for a given block.
\index{patching argument!global}

\smallskip
{\bf Summary of the Local Patching Theorems.} 
\index{global patching when $\Char(K) = 0$!Stage 3: The global patching construction!summary of the Local patching theorems}
The global patching process interacts with the local patching processes
\index{patching argument!global}\index{patching argument!local|ii}  
to adjust the coefficients.  
The following Theorem summarizes the local patching constructions  
proved in Theorems \ref{DCPatch1}, \ref{DRPatch1}, \ref{DCPPatch1}, 
and \ref{DCPCPatch} below.
\index{confinement argument}   

\begin{theorem} \label{LocalPatch}\index{patching theorem!summary of local patching theorems!when $\Char(K_v) = 0$} 
Let $K$ be a number field.  Let $\cC/K$, $\EE$, $\fX$, and $S_K$ be as in Theorem $\ref{aT1-B}$.
Let $\hS_K \supseteq S_K$ be the finite set of places satisfying conditions $(\ref{hS_KList})$.  
For each $v \in \hS_K$, let $\tE_v \subset E_v$, and $0 < h_v < r_v < R_v$ 
be the set and patching parameters constructed in Stage $1$ of the proof. 
\index{patching parameters}
For each $v \in \hS_{K,\infty}$, let $U_v \subset \cC_v(\CC)$ be the chosen open set  
with $U_v \cap E_v = E_v^0$.  For each $v \in S_{K,0}$, let 
$\bigcup_{\ell=1}^{D_v} B(a_\ell,r_\ell) \cap \cC_v(F_{w_\ell})$ 
be the chosen $K_v$-simple decomposition of $E_v$. 
\index{$K_v$-simple!decomposition}  
Let the rational probability vector 
$\vs \in \cP^m(\QQ)$ be as in $(\ref{CGPF2})$, 
and let the natural number $N$ and the coherent approximating functions $\{\phi_v(z)\}_{v \in \hS_K}$ 
\index{coherent approximating functions $\phi_v(z)$} 
be those constructed in Theorem $\ref{CTCX2}$ in Stage $2$ of the proof.   

Then for each $v \in \hS_K$,
there is a constant $k_v > 0$ determined by the $E_v$, $U_v$, and $\phi_v(z)$,   
representing the minimal number of `high-order' stages  
in the local patching process for $K_v$.  \index{patching argument!local} 
Let $\kbar \ge k_v$ be a fixed integer. 
If $v$ is nonarchimedean, put $B_v = h_v^{\kbar N};$    
if $v$ is archimedean, let $B_v > 0$ be arbitrary.    
Then there is an integer $n_v > 0$, depending on $\kbar$ and $B_v$, such
that for each sufficiently large integer $n$ divisible by $n_v$, 
one can carry out the local patching process at $K_v$ as follows:

Put $G_v^{(0)}(z) = Q_{v,n}(\phi_v(z))$, where\index{patching functions, initial $G_v^{(0)}(z)$!construction of} 
\begin{equation*}
\left\{ \begin{array}{l}
   \text{If $K_v \cong \CC$, then $Q_{v,n}(x) =x^n;$} \\
   \text{If $K_v \cong \RR$, set $\hR_v = 2^{-1/n_vN} R_v$, 
         write $n = m_v n_v$, and let $T_{m,R}(x)$ be the} \\
   \text{\quad Chebyshev polynomial of degree $m$ for $[-2R,2R]$ 
\index{Chebyshev polynomial} 
         $($see $(\ref{DFSW2}))$.  Then} \\
   \text{ \qquad \qquad 
         $Q_{v,n}(x) \ = \ T_{m_v,\hR_v^{n_v N}}(T_{n_v,R_v^N}(x));$} \\ 
   \text{If $K_v$ is nonarchimedean and $v \in S_{K,0}$, 
      then $Q_{v,n}(x) = S_{n,v}(x)$} \\
   \text{\quad is the Stirling polynomial of degree $n$ for $\cO_v$
                   $($see $(\ref{FST1}));$} \\
   \text{If $K_v$ is nonarchimedean and $v \in \hS_{K,0} \backslash S_{K,0}$, 
      then $Q_{v,n}(x) = x^n$.}
   \end{array} \right.
\end{equation*}\index{Stirling polynomial!for $\cO_v$} 
For each $k$, $1 \le k \le n-1$, 
let $\{\Delta_{v,ij}^{(k)} \in \CC_v\}_{(i,j) \in \Band_N(k)}$ 
be a $K_v$-symmetric set of numbers,   
\index{$K_v$-symmetric!set of numbers}  
given recursively in $\prec_N$ order, subject to the conditions 
\index{order!$\prec_N$} 
that $\Delta_{v,i0}^{(1)} = 0$ for each archimedean $v$, and for all $(i,j)$  
\begin{equation} \label{FBound1}
|\Delta_{v,ij}^{(k)}|_v \ \le \ \left\{
      \begin{array}{ll} B_v & \text{if \ $k \le \kbar \ ,$} \\
                         h_v^{kN} & \text{if \ $k > \kbar \ .$}
      \end{array} \right.
\end{equation}
For $k = n$, let 
$\{\Delta_{v,\lambda}^{(n)} \in \CC_v\}_{1 \le \lambda \le \Lambda}$ 
be a $K_v$-symmetric set of numbers satisfying  
\index{$K_v$-symmetric!set of numbers}
\begin{equation} \label{FBound2} 
|\Delta_{v,\lambda}^{(n)}|_v \ \le \ h_v^{nN} \ .  
\end{equation} 

\vskip .05 in
\noindent{Then} for each $v \in \hS_K$, 
one can inductively construct $(\fX,\vs)$-functions\index{$(\fX,\vs)$-function}
$G_v^{(1)}(z), \ldots, G_v^{(n)}(z)$ in 
$K_v(\cC)$,\index{patching functions, $G_v^{(k)}(z)$ for $1 \le k \le n$!constructed by patching} 
of common degree $Nn$, such that 

\vskip .05 in
\noindent{$(A)$} For each $k$, $1 \le k \le n$, there are $K_v$-symmetric functions
\index{$K_v$-symmetric!set of functions} 
$\vartheta_{v,ij}^{(k)}(z) \in L_{w_v}(\cC)$,
\index{distinguished place $w_v$}\index{compensating functions $\vartheta_{v,ij}^{(k)}(z)$!are $K_v$-symmetric}
determined recursively in $\prec_N$ order, 
\index{order!$\prec_N$} 
and $(\fX,\vs)$-functions\index{$(\fX,\vs)$-function!$K_v$-rational} 
$\Theta_v^{(k)}(z) \in K_v(\cC)$ of degree at most 
$(n-k)N$, such that\index{compensating functions $\vartheta_{v,ij}^{(k)}(z)$}
\index{patching functions, $G_v^{(k)}(z)$ for $1 \le k \le n$!constructed by patching} 
\begin{eqnarray*} 
G_v^{(k)}(z) & = & G_v^{(k-1)}(z)
              +  \sum_{(i,j) \in \Band_N(k)}
                           \Delta_{v,ij}^{(k)} \vartheta_{v,ij}^{(k)}(z) 
                + \Theta_v^{(k)}(z)  
        \quad \text{for $k < n \ ,$} \\
G_v^{(n)}(z) & = & G_v^{(n-1)}(z) + \sum_{\lambda=1}^{\Lambda}
                   \Delta_{v,\lambda}^{(n)} \varphi_{\lambda}(z) \,                                                     
\end{eqnarray*}
and where for each $k < n$ and each $(i,j)$, 
if $\tc_{v,i}$ is the leading coefficient of $\phi_v(z)$ at $x_i$,   
\index{coefficients $A_{v,ij}$!leading} 

$(1)$ $\vartheta_{v,ij}^{(k)}(z)$
\index{compensating functions $\vartheta_{v,ij}^{(k)}(z)$!poles and leading coefficients of} 
has a pole of order $nN_i-j > (n-k-1)N_i$ at $x_i$ 
and leading coefficient $\tc_{v,i}^{n-k-1}$,  
\index{coefficients $A_{v,ij}$!leading} 
a pole of order at most $(n-k-1)N_{i^{\prime}}$ at each $x_{i^{\prime}} \ne x_i$, 
and no other poles;  

$(2)$ $\sum_{i=1}^m \sum_{j=(k-1)N_i}^{kN_i-1} 
  \Delta_{v,i^{\prime}j}^{(k)} \vartheta_{v,i^{\prime}j}^{(k)}(z)$
belongs to $K_v(\cC)\ ;$ 

$(3)$ $\Theta_v^{(k)}(z)$ is determined by the local patching process at $v$ 
after the coefficients in $\Band_N(k)$ have been modified by adding 
\index{coefficients $A_{v,ij}$}\index{band!$\Band_N(k)$} 
$\sum_{(i,j) \in \Band_N(k)} \Delta_{v,ij}^{(k)} 
\vartheta_{v,ij}^{(k)}(z)${compensating functions $\vartheta_{v,ij}^{(k)}(z)$} 
to $G_v^{(k)}(z);$\index{patching functions, $G_v^{(k)}(z)$} 
it has a pole of order at most $(n-k)N_i$ at each $x_i$ 
and no other poles, and may be the zero function.  

\vskip .05 in
\noindent{$(B)$} For each $k = 0, \ldots, n$,\index{patching functions, $G_v^{(k)}(z)$ for $1 \le k \le n$!roots are confined to $E_v$} 
\begin{equation*} 
\left\{ \begin{array}{l}
\text{If $K_v \cong \CC$, then  
  $\{ z \in \cC_v(\CC_v) : |G_v^{(k)}(z)|_v \le r_v^{nN} \} \subset U_v = E_v^0;$} \\
\text{If $K_v \cong \RR$, then} \\ 

\text{\quad $(1)$ 
    the zeros of $G_v^{(k)}(z)$\index{patching functions, $G_v^{(k)}(z)$ for $1 \le k \le n$!roots are confined to $E_v$} 
   all belong to $E_v^0$,\
    and for each component $E_{v,i}$ of $E_v$,} \\
\text{\qquad if $\phi_v(z)$ has $\tau_i$ zeros  in $E_{v,i}$, 
     then $G_v^{(k)}(z)$ has $T_i = n \tau_i$ zeros in $E_{v,i}$.} \\ 
\text{\quad $(2)$ $\{ z \in \cC_v(\CC_v) : 
         |G_v^{(k)}(z)|_v \le 2r_v^{nN}\} \ \subset \ U_v$,} \\ 
\text{\quad $(3)$ on each component $E_{v,i}$ contained in $\cC_v(\RR)$,} \\
\text{\qquad \qquad \qquad   
      $G_v^{(k)}(z)$ oscillates $T_i$ 
           times between $\pm 2 r_v^{nN}$ on  $E_{v,i}$.} \\
  
\text{If $K_v$ is nonarchimedean and $v \in S_{K,0}$, 
            then all the zeros 
    of $G_v^{(k)}(z)$\index{patching functions, $G_v^{(k)}(z)$ for $1 \le k \le n$!roots are confined to $E_v$} belong to $E_v$,} \\
 
\text{\quad and for $k = 0$ and $k = n$ they are distinct. When $k = n$,} \\

\text{\qquad \qquad \qquad 
$\{z \in \cC_v(\CC_v) : G_v^{(n)}(z) \in \cO_v \cap D(0,r_v^{nN}) \} \subset E_v.$} \\ 

\text{If $K_v$ is nonarchimedean and $v \in \hS_{K,0} \backslash S_{K,0}$, then} \\
\text{ \qquad \qquad \qquad 
$\{ z \in \cC_v(\CC_v) : |G_v^{(k)}(z)|_v \le R_v^{nN} \} = E_v.$ }
\end{array} \right. 
\end{equation*}
\end{theorem}

\vskip .1 in
\noindent{\bf Remark 1.} 
For almost all $v$ and $k$, we will have $\Theta_v^{(k)}(z) = 0$;
the only exception is for one value $k = k_1$ for each 
$v \in S_{K,0}$, where $\Theta_v^{(k_1)}(z)$ is chosen to `separate the roots'\index{separate the roots} 
of $G_v^{(k_1)}(z)$.\index{patching functions, $G_v^{(k)}(z)$ for $1 \le k \le n$!for nonarchimedean $K_v$-simple sets!roots are separated}  
See the discussion after Theorem 
\ref{DCPCPatch}, and Phase 3 in the proof of that theorem.  

\vskip .1 in
\noindent{\bf Remark 2.}  
Examining the proofs of the local patching theorems\index{patching argument!local}  
shows that from a local standpoint, 
the order in which the coefficients are received within a band 
\index{coefficients $A_{v,ij}$} 
\index{band!coefficients patched by bands}
is immaterial, provided that for each $x_i$, they are received by increasing $j$.  
For any such order, 
the same changes are produced in lower order coefficients within the band,
and the same functions $G_v^{(k)}(z)$ are obtained.\index{patching functions, $G_v^{(k)}(z)$}
For this reason, in the local process at each $v$, 
it is permissible to subdivide 
$\Gal(L/K)$-blocks into $\Gal^{c}(\CC_v/K_v)$-sub-blocks.   

In fact, for all bands for nonarchimedean $v$, and for the bands with $k > \kbar$
\index{band!coefficients patched by bands}
for archimedean $v$, the compensating functions $\vartheta_{v,ij}^{(k)}$ 
\index{compensating functions $\vartheta_{v,ij}^{(k)}(z)$!are $K_v$-symmetric}
are $K_v$-symmetric and are independent of the $\Delta_{v,ij}^{(k)}$.
\index{$K_v$-symmetric!set of functions}
For archimedean $v$ and bands with $k \le \kbar$, 
\index{band!coefficients patched by bands}
patching is carried out by a process called ``magnification''\index{magnification argument|ii} 
(see the proofs of Theorems \ref{DCPatch1} and \ref{DRPatch1}), 
and our description of the $\vartheta_{v,ij}^{(k)}$\index{compensating functions $\vartheta_{v,ij}^{(k)}(z)$}
in Theorem \ref{LocalPatch} is correct but artificial:    
rather, for each $\Gal^{c}(\CC_v/K_v)$-stable subset of the indices in a band,
the changes in those coefficients produce a canonical $K_v$-rational 
change in $G_v^{(k)}(z)$.\index{patching functions, $G_v^{(k)}(z)$ for $1 \le k \le n$!constructed by patching}  

\vskip .1 in
{\bf The choice of the parameters $\kbar$ and $B_v$.}
\index{global patching when $\Char(K) = 0$!Stage 3: The global patching construction!the choices of $\kbar$ and $B_v$}
\index{patching parameters} 
In Stage 1 of the construction we have chosen a collection of numbers $h_v$ for $v \in \hS_K$ 
such that $\prod_{v \in \hS_K} h_v^{D_v} > 1$.  
Likewise, for each $v \in \hS_K$, Theorem \ref{LocalPatch} provides a number $k_v$, 
the ``minimal number of stages considered high-order'' by the local patching process at $v$.
\index{patching argument!local} 

Let $\kbar$ be the smallest integer such that 
\begin{equation} \label{Fkbarchoice}
\left\{ \begin{array}{l} 
         \text{$\kbar \ \ge \ k_v$ \quad for each $v \in \hS_K$\ ,} \\
         \text{$(\prod_{v \in \hS_K} h_v^{D_v})^{\kbar N [L:K]} \ > \ C_L(\hS_K)$\ .} 
        \end{array} \right.
\end{equation}
where $C_L(\hS_K)$ is the constant from Proposition \ref{CStrong3}.
\index{Strong Approximation theorem!uniform}

\vskip .05 in 
For each nonarchimedean $v \in \hS_K$, the choice of $\kbar$
determines the constant $B_v = h_v^{\kbar N}$ 
in the local patching process (see Theorem \ref{LocalPatch}). 
For archimedean $v$, the constants $B_v$ in Theorem \ref{LocalPatch}
can be specified arbitrarily. Choose them large enough that
\begin{equation} \label{FBvchoice} 
(\prod_{v \in \hS_K} B_v^{D_v})^{[L:K]} \ > \ C_L(\hS_K) \ .
\end{equation}

Given $w \in \hS_L$, let $v$ be the place of $K$ under $w$.  
By our normalization of the absolute values 
in \S\ref{Chap3}.\ref{NotationSection}, $|x|_w^{D_w} = |x|_v^{D_v [L_w:K_v]}$ for each 
$x \in \CC_w \cong \CC_v$.  Define $h_w$ by $h_w^{D_w} = h_v^{D_v [L_w:K_v]}$,
and define $B_w$ by $B_w^{D_w} = B_v^{D_v [L_w:K_v]}$.  
Then  $|x|_w \le h_w^{kN}$ iff $|x|_v \le h_v^{kN}$, 
$|x|_w \le B_w$ iff $|x|_v \le B_v$, and 
\begin{equation*}
(\prod_{w \in \hS_L^{+}} h_w^{D_w})^{\kbar N} \ > \ C_L(\hS_K^{+}) \ , \quad  
\prod_{w \in \hS_L^{+}} B_w^{D_w} \ > \ C_L(\hS_K^{+}) \ .
\end{equation*}
  
\vskip .1 in
{\bf The choice of the initial patching functions.}\index{patching functions, initial $G_v^{(0)}(z)$|ii}
\index{global patching when $\Char(K) = 0$!Stage 3: The global patching construction!the choice of $n$} 
Theorem \ref{CTCX2} gives a degree $N$ and a collection of
coherent approximating functions $\phi_v(z) \in K_v(\cC)$ for $v \in \hS_K$, 
\index{initial approximating functions $f_v(z)$}
of common degree $N$.  For each $v \in \hS_K$, let $Q_{v,n}(x) \in K_v[x]$ be 
the monic degree-raising polynomial\index{degree-raising} of degree $n$ 
from Theorem \ref{LocalPatch}.  For suitable $n$, 
we will take the initial patching function at $v$ 
to be $G_v^{(0)}(z) = Q_{v,n}(\phi_v(z))$.\index{patching functions, initial $G_v^{(0)}(z)$!construction of} 
As explained above, our plan is to inductively
 construct functions $G_v^{(1)}(z), \ldots, G_v^{(n)}(z)$, 
making more and more coefficients global $\hS_L$-integers at each stage,
\index{coefficients $A_{v,ij}$}  
until finally the $G_v^{(n)}(z) = G^{(n)}(z)$ 
are $K$-rational and independent of $v$.\index{patching functions, $G_v^{(k)}(z)$ for $1 \le k \le n$!constructed by patching}       

At several places in the patching process,\index{patching argument!global}  
it is important to consider the $\phi_v(z)$ and $G_v^{(k)}(z)$
\index{patching functions, $G_v^{(k)}(z)$ for $1 \le k \le n$!viewed simultaneously over $K_v$ and $L_w$} 
over the fields $L_w$ with $w|v$,  
rather than over $K_v$.  This has already been seen in Theorem \ref{CTCX2}.   
However, our ultimate goal is to construct a $K$-rational function. 
Hence, the choices made in the local patching constructions must\index{patching argument!local}  
depend only on places $v$ of $K$, not on the places $w$ of $L$ with $w|v$. 

We resolve this by considering the $\phi_v(z)$ 
and $G_v^{(k)}(z)$\index{patching functions, $G_v^{(k)}(z)$ for $1 \le k \le n$!viewed simultaneously over $K_v$ and $L_w$} 
simultaneously over $K_v$ and the $L_w$.  
Viewing them over $L_w$ enables us examine their coefficients, 
\index{coefficients $A_{v,ij}$} 
which are canonically $L_w$-rational. Viewing them over $K_v$ assures that 
any choices in the local patching processes occur in the same way for all $w|v$. 

\vskip .1 in
{\bf The choice of $n$.}
The leading coefficients of the $G_v^{(0)}(z)$ 
are hardest to patch;\index{patching functions, initial $G_v^{(0)}(z)$!leading coefficients of} 
\index{coefficients $A_{v,ij}$!leading} 
we must make them $\hS_L$-units.  The key to this is our choice of $n$. 
\index{patching functions, initial $G_v^{(0)}(z)$!making the leading coefficients $S$-units}  

Given $v \in \hS_K$, put $\phi_w(z) = \phi_v(z)$  
for each $w \in \hS_L$ with $w|v$, viewing the $\phi_w(z)$ 
as functions in $L_w(\cC)$.   
By Theorem \ref{CTCX2} the leading coefficients $\tc_{w,i}$ of the $\phi_w(z)$ 
\index{coefficients $A_{v,ij}$!leading} 
have the property that there are an integer $n_0$, 
and a $K$-symmetric system of $\hS_L$-units $\mu_i$, such that for each $i$  
\index{$K$-symmetric!system of units}
\begin{equation} \label{FLead}
 \left\{  \begin{array}{cl}
\tc_{w,i}^{n_0} \ = \ \mu_i & \text{for each archimedean $w \in \hS_L$,} \\
|\tc_{w,i}^{n_0}|_w \ = \ |\mu_i|_w & 
             \text{for each nonarchimedean $w \in \hS_L$.} 
   \end{array} \right. 
\end{equation}

For each nonarchimedean $v \in \hS_K$, 
all the fields $L_w$ for $w|v$ are isomorphic, 
and by the structure of the group of units $O_w^{\times}$ there is an integer 
$n_v^{\prime} > 0$ such that for each  $x \in \cO_w^{\times}$, and
each integer $n^{\prime}$ divisible by $n_v^{\prime}$, 
\begin{equation} \label{FCloseUp} 
|x^{n^{\prime}} - 1|_v \cdot \max_{1 \le i \le m}(|\tc_{v,i}|_v^2) \ \le \ B_v \ .
\end{equation}
Let $n_1$ be the least common multiple of the  $n_v^{\prime}$. 

For each $v \in \hS_K$, Theorem \ref{LocalPatch} provides a number $n_v$
\index{patching parameters} 
such that the local patching process\index{patching argument!local}  at $v$ 
will preserve the properties of the roots of $G_v^{(0)}(z)$,
\index{patching functions, $G_v^{(k)}(z)$ for $1 \le k \le n$!roots are confined to $E_v$} 
provided $n_v|n$ and $n$ is sufficiently large, 
and the $\Delta_{v,ij}^{(k)}$, $\Delta_{v,\lambda}^{(n)}$ are $K_v$-symmetric and 
\index{$K_v$-symmetric!set of functions}
satisfy the size constraints (\ref{FBound1}), (\ref{FBound2}) 
relative to $h_v$, $\kbar$ and the $B_v$ chosen above. 
Let $n_2$ be the least common multiple of the $n_v$ for $v \in \hS_K$.  

Finally, let $n$ be a positive integer such that  
\begin{equation} \label{FnChoice}
 n_0 n_1 n_2 | n \ .
\end{equation}
By Theorem \ref{LocalPatch} there is an $n_3$ such that if $n \ge n_3$,
then for each $v \in \hS_K$ the local patching process\index{patching argument!local}  can 
be successfully completed.

Until last step in the proof, $n \ge n_3$ will be a fixed integer 
satisfying (\ref{FnChoice}).

\medskip
{\bf Patching the Leading Coefficients.}\index{patching!leading coefficients|ii}
\index{global patching when $\Char(K) = 0$!Stage 3: The global patching construction!patching the leading coefficients}
Given such an $n$, for each $v \in \hS_K$ 
put $G_w^{(0)}(z) = G_v^{(0)}(z)$\index{patching functions, $G_v^{(k)}(z)$ for $1 \le k \le n$!viewed simultaneously over $K_v$ and $L_w$}
for all $w|v$, and expand 
\begin{equation*}
G_w^{(0)}(z) \ = \
    \sum_{i=1}^m \sum_{j=0}^{(n-1)N_i-1} A_{w,ij} \varphi_{i,nN_i-j}(z)
       + \sum_{\lambda = 1}^{\Lambda} A_{w,\lambda} \varphi_{\lambda} \ , 
\end{equation*}
with the $A_{w,ij}, A_{w,\lambda} \in L_w$.
Since each $Q_{v,n}(z)$ is monic, 
for each $w \in \hS_L$ the leading coefficient of $G_w^{(0)}(z)$ at $x_i$ is
\index{coefficients $A_{v,ij}$!leading}  
\begin{equation*} 
A_{w,i0} \ = \ \tc_{w,i}^n \ .
\end{equation*}

For each archimedean $v$, and all $w|v$, by (\ref{FLead})  
$\tc_{w,i}^n =  \mu_i^{n/n_0}$ is a global $\hS_L$-unit independent of $w$
and $v$.  By construction the $\mu_i$ are $K$-symmetric.  
\index{$K$-symmetric!system of units}
In the local patching process at $v$,\index{patching argument!local} 
take $\Delta_{v,i0}^{(1)} = 0$ for $i = 1, \ldots, m$.  
Trivially the $\Delta_{v,i0}^{(1)}$ are $K_v$-symmetric, with 
\index{$K_v$-symmetric!set of numbers}
$|\Delta_{v,i0}^{(1)}|_v \le B_v$ for each $i$.  
Put $\hG_v(z) = G_v^{(0)}(z)$.\index{patching functions, initial $G_v^{(0)}(z)$!leading coefficients of}
\index{patching functions, initial $G_v^{(0)}(z)$!making the leading coefficients $S$-units}

For nonarchimedean $v  \in \hS_K$, we claim that we can adjust the leading coefficients  
\index{coefficients $A_{v,ij}$!leading} 
to be $\mu_i^{n/n_0}$, as well.  This depends on the fact that $n_0 n_2|n$.  

Let $v \in \hS_{K,0}$. For each $w|v$, since $|\tc_{w,i}^{n_0}|_w = |\mu_i|_w$, 
it follows that $\mu_i/\tc_{w,i}^{n_0} \in \cO_w^{\times}$.   If we put 
$\Delta_{w,i0}^{(1)} \ = \ \tc_{w,i}^2 \big(\frac{\mu_i^{n/n_0}}{\tc_{w,i}^n} - 1\big)$, 
then by (\ref{FCloseUp}),
\begin{equation}
|\Delta_{w,i0}^{(1)}|_v \ = \ 
  |\tc_{w,i}^2|_v \cdot \Big|\left(\frac{\mu_i}{\tc_{w,i}^{n_0}}\right)^{n/n_0} - 1 \Big|_v \ \le \  B_v \ .
\end{equation}
Moreover
\begin{equation} \label{FChange} 
\mu_i^{n/n_0} \ = \ A_{w,i0} + \Delta_{w,i0}^{(1)} \tc_{w,i}^{n-2} \ .
\end{equation}
As will be seen below, this is what is needed for   
the local patching constructions in Theorem \ref{LocalPatch} 
to change the leading coefficients to $\mu_i^{n/n_0}$.
\index{coefficients $A_{v,ij}$!leading}   

However, for the local patching process\index{patching argument!local} at $v$
we need changes $\Delta_{v,i0}^{(1)}$ independent of $w$, 
not the $\Delta_{w,i0}^{(1)}$ which \'a priori could depend on $w$.  
We will now present a ``see-saw'' argument\index{see-saw argument|ii} using Proposition \ref{CPPr1}  
which shows that for all $w|v$, the functions 
$\sum_{i=1}^m \Delta_{w,i0}^{(1)} \varphi_{i,nN_i}(z) \in L_w(\cC)$
belong to $K_v(\cC)$, and are independent of $w$, so we can take
$\Delta_{v,i0}^{(1)} = \Delta_{w_v,i0}^{(1)}$
for the distinguished place $w_v$\index{distinguished place $w_v$} induced by the embedding 
$\tK \hookrightarrow \CC_v$ used to identify $\fX$ with a subset of $\cC_v(\CC_v)$.   
A similar argument applies at later steps of the patching process, 
and in the future we will omit some details.  

The $G_w^{(0)}(z)$ with $w|v$ are all the same and belong to $K_v(\cC)$, so 
\begin{equation*} 
\oplus_{w|v} G_w^{(0)}(z) \ \in \ \oplus_{w|v} L_w(\cC) 
          \ \cong \ L \otimes_K K_v(\cC)
\end{equation*}           
is $\Gal(L/K)$-invariant in the sense of \S\ref{Chap7}.\ref{SemiLocalSection}.  
By Proposition \ref{CPPr1}, if we put $F = K(x_i)$, 
then $\oplus_{w|v} A_{w,i0}$ belongs to $\oplus_{u|v} F_u$ (embedded
semi-diagonally in $\oplus_{w|v} L_w \cong L \otimes_K K_v$), and for each
$\sigma \in \Gal(L/K)$, 
\begin{equation*} 
\sigma(\oplus_{w|v} A_{w,i0}) = \oplus_{w|v} A_{w,\sigma(i) 0} \ .   
\end{equation*}
By a similar argument, 
$\sigma(\oplus_{w|v} \tc_{w,i}) = \oplus_{w|v} \tc_{w,\sigma(i)}$.

On the other hand, by Theorem \ref{CTCX2}, $\mu_i \in K(x_i)$ and 
$\sigma(\mu_i) = \mu_{\sigma(i)}$.  Hence, viewing $\mu_i$ 
as embedded semi-diagonally in $\oplus_{w|v} L_w$, we see that 
\begin{equation*}
\oplus_{w|v} \Delta_{w,i0}^{(1)} 
\ = \ \oplus_{w|v} (\mu_i^{n/n_0} - A_{w,i0})/\tc_{w,i}^{n-2}
\end{equation*}
also satisfies 
$\sigma(\oplus_{w|v} \Delta_{w,i0}^{(1)}) 
= \oplus_{w|v} \Delta_{w,\sigma(i) 0}^{(1)}$
for each $\sigma \in \Gal(L/K)$. 

For the basis functions we have $\sigma(\varphi_{i,j}) = \varphi_{\sigma(i),j}$
\index{basis!$L$-rational} 
by construction.  Thus for each galois orbit $\fX_{\ell}$
\begin{equation*}
\oplus_{w|v} 
  \big( \sum_{x_i \in \fX_{\ell}} \Delta_{w,i0}^{(1)} \varphi_{i,nN_i}(z) \big) 
         \ \in \ L \otimes K_v(\cC) 
\end{equation*}
is $\Gal(L/K)$-invariant.  Applying Proposition \ref{CPPr1} in reverse, 
there is a function
$H_{v,\ell}(z) = \sum_{x_i \in \fX_{\ell}} \Delta_{v,i0}^{(1)} \varphi_{i,nN_i} \in K_v(\cC)$ 
such that  
\begin{equation*}
 \sum_{x_i \in \fX_{\ell}} \Delta_{w,i0}^{(1)} \varphi_{i,nN_i}(z) 
        \ = \ H_{v,\ell}(z)
\end{equation*}
for each $w|v$.  Thus the $\Delta_{v,i0}^{(1)} := \Delta_{w_v,i0}^{(1)}$\index{distinguished place $w_v$} 
are well-defined and $K_v$-symmetric.   
\index{$K_v$-symmetric!set of numbers} 
 
Patch $G_v^{(0)}(z)$ by setting\index{patching functions, $G_v^{(k)}(z)$ for $1 \le k \le n$!constructed by patching}
\begin{equation}  \label{CFPWQ}
\hG_v(z) \ = \ G_v^{(0)}(z) + \sum_{i=1}^m \Delta_{v,i0}^{(1)} \vartheta_{v,i0}^{(1)}(z)
\end{equation}
where the $\vartheta_{v,i0}^{(1)}(z)$ are the compensating functions
\index{compensating functions $\vartheta_{v,ij}^{(k)}(z)$!poles and leading coefficients of}
from Theorem \ref{LocalPatch}.  The leading coefficient 
\index{coefficients $A_{v,ij}$!leading} 
of $\vartheta_{v,i0}^{(1)}(z)$ at $x_i$ is $\tc_{v,i}^{n-2} = \tc_{w_v,i}^{n-2}$,\index{distinguished place $w_v$} 
so by (\ref{FChange}) this changes the leading coefficient of $G_w^{(0)}(z)$ 
at $x_i$ to $\mu_i^{n/n_0}$, for each $w$ and $i$.  
(The lower-order coefficients are changed as well, 
but they will be dealt with in subsequent patching steps.)\index{patching argument!global} 
By Theorem \ref{LocalPatch}.A.2, $\hG_v(z)$ is $K_v$-rational.        

\vskip .1 in
{\bf Patching the High Order Coefficients.}\index{patching!high-order coefficients}
\index{global patching when $\Char(K) = 0$!Stage 3: The global patching construction!patching the high-order coefficients}
Next we patch the remaining coefficients for the stage $k = 1$
\index{coefficients $A_{v,ij}$!high-order} 
and inductively carry out the patching process for stages  $k = 2, \ldots, \kbar$.  

Suppose that for some $k$, 
we have constructed functions $G_v^{(k-1)}(z) \in K_v(\cC)$, 
$v \in \hS_K$.\index{patching functions, $G_v^{(k)}(z)$ for $1 \le k \le n$!constructed by patching} 
In the $k^{th}$ stage we patch the coefficients with indices in $\Band_N(k)$
\index{band!$\Band_N(k)$} 
by increasing $\prec_N$ order, patching all the coefficients\index{patching argument!global}
\index{order!$\prec_N$} 
in a given block at once.
\index{block, galois!coefficients patched block by block} 
  
Suppose that after patching a certain number of blocks, 
\index{block, galois!coefficients patched block by block}
we have obtained functions $\hG_v(z) \in K_v(\cC)$ for $v \in \hS_K$.  
(When $k =1$, we view patching the high order coefficients 
\index{coefficients $A_{v,ij}$!high-order} 
as taking the initial step in passing 
from $G_v^{(0)}(z)$ to $G_v^{(1)}(z)$.)\index{patching functions, $G_v^{(k)}(z)$ for $1 \le k \le n$!constructed by patching} 
To lighten notation, we update the coefficients after each step:
for each $v$, write 
\begin{equation*}
\hG_v(z) \ = \
\sum_{i=1}^m \sum_{j=0}^{(n-1)N_i} \, A_{v,ij} \varphi_{i,nN_i-j}(z)
   + \sum_{\lambda = 1}^{\Lambda} A_{v,\lambda} \, \varphi_{\lambda}   
\end{equation*}
with the $A_{v,ij}, A_{v,\lambda} \in L_{w_v}$.\index{distinguished place $w_v$}  
For each $w|v$, put $\hG_w(z) = \hG_v(z)$ 
and regard $\hG_w(z)$ as belonging to $L_w(\cC)$.  Expand
\begin{equation*}
\hG_w(z) \ = \
\sum_{i=1}^m \sum_{j=0}^{(n-1)N_i} A_{w,ij} \, \varphi_{i,nN_i-j}(z)
   + \sum_{\lambda = 1}^{\Lambda} \, A_{w,\lambda} \varphi_{\lambda}  
\end{equation*}
where the $A_{w,ij}$ and $A_{w,\lambda}$ belong to $L_w$.  

Let $(i_0,j_0) \in \Band_N(k)$ be the least index for which the coefficients
have not been patched.   Thus, for each $(i,j) \prec_N (i_0,j_0)$ there 
is an $A_{ij} \in L$ such that $A_{w,ij} = A_{ij}$ for all $w$.  
To patch the coefficients for the indices $(i,j) \in \Block(i_0,j_0)$, 
we first determine a target value $A_{i_0 j_0} \in K(x_{i_0})$ 
for the $A_{w,i_0 j _0}$, $w \in \hS_L$, and then, to preserve galois equivariance, 
we define the target values for the other $(i,j)$ in $\Block(i_0,j_0)$ 
\index{block, galois!$\Block(i,j)$}
by requiring that if $\sigma \in \Gal(L/K)$ 
is such that $\sigma(i_0) = i$, then $A_{i,j} = \sigma(A_{i_0 j_0})$.  
This is well-defined, since if $\sigma_1, \sigma_2 \in \Gal(L/K)$ are such that 
$\sigma_1(i_0) = \sigma_2(i_0)$, then $\sigma_2^{-1}\sigma_1$ fixes $x_{i_0}$,
and so since $A_{i_0 j_0} \in K(x_{i_0})$ 
we have $\sigma_1(A_{i_0 j_0}) = \sigma_2(A_{i_0 j_0})$.   
     
Consider the vector
\begin{equation*}
\vec{A}_{L,i_0 j_0} \ := \ \oplus_{w \in \hS_L} A_{w,i_0 j_0} 
                  \ \in \ \oplus_{w \in \hS_L} L_w \ .
\end{equation*}
Put $F = K(x_{i_0})$.  
For each $v \in \hS_K$, the $\hG_w(z) \in K_v(\cC)$ are the same for all $w|v$, 
so Proposition \ref{CPPr1} tells us that $\vec{A}_{L,i_0,j_0}$ 
belongs to $\oplus_{u \in \hS_F} F_u$, embedded semi-diagonally in 
$\oplus_{w \in \hS_L} L_w$.

For each $w \in \hS_L$, put
\begin{equation}  
Q_w \ = \ B_w \cdot |\tc_{w,i_0}^{n-k-1}|_w \ ,  \label{CFGH1}
\end{equation}
where $\tc_{w,i_0}$ is the leading coefficient of $\phi_w(z)$ at $x_{i_0}$.
\index{coefficients $A_{v,ij}$!leading} 
Theorem \ref{CTCX2} has arranged that 
$\prod_{w \in \hS_L} |\tc_{w,i_0}|_w^{D_w} = 1$, so   
\begin{equation*}
\prod_{w \in \hS_L} Q_w^{D_w} 
\ = \ \prod_{w \in \hS_L} B_w^{D_w} \ > \ C_L(\hS_K) \ .
\end{equation*}
Note that $B_w$ depends only on the place $v$ of $K$ below $w$,
while $|\tc_{w,i_0}|_w$ depends only on the place $u$ of $F$ below $w$,
since the $\phi_w(z) \in K_v(\cC)$ are the same for all $w|v$.
Hence $Q_w$ depends only on the place $u$ below $w$.  
Similarly the coefficients $A_{w,i_0 j_0}$ with $w|u$ belong to $F_u$ 
and depend only on $u$.  

Thus we can apply Proposition \ref{CStrong3} to the elements
\index{Strong Approximation theorem!uniform}
$c_u = A_{w,i_0 j_0} \in F_u$, and to the $Q_w$.
By Proposition \ref{CStrong3}, there is an $A_{i_0 j_0} \in K(x_{i_0})$ such that
\index{Strong Approximation theorem!uniform} 
\begin{equation*}
\left\{ \begin{array}{cl}
     |A_{i_0 j_0} - A_{w,i_0 j_0}|_w \le Q_w & \text{for each $w \in \hS_L$\ ,} \\
     |A_{i_0 j_0}|_w \le 1         & \text{for each $w \notin \hS_L$\ .}
        \end{array} \right.
\end{equation*}
This $A_{i_0 j_0}$ will be the target in patching the $A_{w,i_0 j_0}$.
For each $(i,j_0) \in \Block(i_0,j_0)$, choose a $\sigma \in \Gal(L/K)$
with $\sigma(x_{i_0}) = x_i$, and put $A_{i j_0} = \sigma(A_{i_0 j_0})$.
Since $A_{i_0,j_0} \in K(x_{i_0})$, the $A_{i j_0}$ are well-defined    
and satisfy $\sigma(A_{i j_0}) = A_{\sigma(i) j_0}$ for all $\sigma \in \Gal(L/K)$.

Put $\Delta_{w,ij}^{(k)} = (A_{ij} - A_{w,ij})/\tc_{w,i}^{n-k-1}$
for each $(i,j) \in \Block(i_0,j_0)$ and each $w \in \hS_L$.  Thus 
\begin{equation}  \label{FGas1}
|\Delta_{w,ij}^{(k)}|_w \ \le \ Q_w/|\tc_{w,i}^{n-k-1}|_w \ = \ B_w  
\end{equation}
and 
\begin{equation} \label{FGas2}  
     A_{ij} \ = \ A_{w,ij} + \Delta_{w,ij}^{(k)} \tc_{w,i}^{n-k-1} \ .
\end{equation}
By construction the $\oplus_{w|v} \Delta_{w,ij}^{(k)}$ are equivariant under  
$\Gal(L/K)$.

Let $\fX_{\ell}$ be the galois orbit of $x_{i_0}$. 
By a see-saw argument like the one used in patching the leading coefficients,
\index{see-saw argument} 
\index{coefficients $A_{v,ij}$} 
for each $v \in \hS_K$ the functions 
\begin{equation*}
 \sum_{x_i \in \fX_{\ell}} \Delta_{w,i j_0}^{(k)} 
              \varphi_{i,j_0}(z) 
\ \in \ L_w(\cC) 
\end{equation*}
are independent of $w|v$ and belong to $K_v(\cC)$. 

Define the patching coefficients 
\index{coefficients $A_{v,ij}$!patching}\index{patching coefficients} 
by $\Delta_{v,ij_0}^{(k)} = \Delta_{w_v,ij_0}^{(k)}$,\index{distinguished place $w_v$}  
for each $(i,j_0) \in \Block(i_0,j_0)$ and each $v \in \hS_K$. 

\vskip .1 in
For each $v \in \hS_K$, the local patching construction\index{patching argument!local} for $v$
produces functions $\vartheta_{v,ij_0}(z) \in L_{w_v}(\cC)$\index{distinguished place $w_v$}
such that $\vartheta_{v,ij}(z)$
\index{compensating functions $\vartheta_{v,ij}^{(k)}(z)$!poles and leading coefficients of}
has a pole of order $nN_i-j_0$ at $x_i$,  
with leading coefficient $\tc_{v,i}^{n-k-1} = \tc_{w_v,i}^{n-k-1}$ at $x_i$,\index{distinguished place $w_v$} 
\index{coefficients $A_{v,ij}$!leading} 
and poles of order $\le (n-k-1)N_{i^{\prime}}$ for all $i^{\prime} \ne i$.  
Theorem \ref{LocalPatch}.A2 shows that
$\sum_{x_i \in \fX_{\ell}} \Delta_{v,i j_0}^{(k)} \vartheta_{v,i j_0}(z) 
\in K_v(\cC)$.  Replace  $\hG_v(z)$ by 
\begin{equation*} 
\ckG_v(z) \ = \ 
\hG_v(z) + \ \sum_{x_i \in \fX_{\ell}} \Delta_{v,i j_0}^{(k)} \vartheta_{v,i,j_0}(z) \ .
\end{equation*} 
  
By (\ref{FGas2}), for each $w|v$ and each $(i,j_0) \in \Block(i_0,j_0)$,  
this changes the coefficient $A_{w,ij_0}$ of $\hG_w(z)$ to $A_{ij_0}$,   
and leaves the coefficients preceding $(i_0,j_0)$ unchanged. 
Hence, the induction can continue. 
 
When all the coefficients in $\Band_N(k)$ have been patched, 
\index{band!coefficients patched by bands}\index{band!$\Band_N(k)$}   
the local patching process at $v$ determines 
an $(\fX,\vs)$-function\index{$(\fX,\vs)$-function!$K_v$-rational} 
$\Theta_v^{(k)}(z) \in K_v(\cC)$ with a pole of order at most $(n-k)N_i$ at each $x_i$.
Using the current $\hG_v(z)$, 
we set $G_v^{(k)}(z) = \hG_v(z) + \Theta_v^{(k)}(z)$\index{patching functions, $G_v^{(k)}(z)$ for $1 \le k \le n$!constructed by patching} 
and replace $k$ by $k+1$.

\vskip .1 in
{\bf Patching the Middle Coefficients.}
\index{coefficients $A_{v,ij}$!middle}\index{patching!middle coefficients}
\index{global patching when $\Char(K) = 0$!Stage 3: The global patching construction!patching the middle coefficients}  
In this stage we carry out the patching process for $k = \kbar+1, \ldots, n-1$.
The construction is the same as for the high 
order coefficients, except that in place of (\ref{CFGH1}) we take
\begin{equation} \label{CFGH3}  
Q_w\ = \ h_w^{kN} \cdot |\tc_{w,i}^{n-k-1}|_w \ .
\end{equation}
The construction succeeds because for each $k \ge \kbar$, 
\begin{equation} \label{CFGH5}
\big(\prod_{w \in \hS_L} h_w^{D_w}\big)^{k N} \ > \ C_L(\hS_K) .
\end{equation}

\medskip
{\bf Patching the Low Order Coefficients.}
\index{coefficients $A_{v,ij}$!low-order}\index{patching!low-order coefficients}
\index{global patching when $\Char(K) = 0$!Stage 3: The global patching construction!patching the low-order coefficients} 
The final stage of the global patching process deals with the coefficients
$A_{v,\lambda}$ in the functions\index{patching functions, $G_v^{(k)}(z)$ for $1 \le k \le n$!constructed by patching}  
\begin{equation*}
G_v^{(n-1)}(z) \ = \
\sum_{i=1}^m \sum_{j=0}^{(n-1)N_i} A_{v,ij} \varphi_{i,nN_i-j}(z)
   + \sum_{\lambda = 1}^{\Lambda} A_{v,\lambda} \varphi_{\lambda} \ .
\end{equation*}
All the $A_{v,\lambda}$ will be patched simultaneously.  

For each $w|v$, put $G_w^{(n-1)}(z) = G_v^{(n-1)}(z)$ 
and expand\index{patching functions, $G_v^{(k)}(z)$ for $1 \le k \le n$!viewed simultaneously over $K_v$ and $L_w$} 
\begin{equation*}
G_w^{(n-1)}(z) \ = \
\sum_{i=1}^m \sum_{j=0}^{(n-1)N_i} A_{w,ij} \varphi_{i,nN_i-j}(z)
   + \sum_{\lambda = 1}^{\Lambda} A_{w,\lambda} \varphi_{\lambda} \ .
\end{equation*}

By construction, for each $v \in \hS_K$, 
the $\oplus_{w|v} G_w^{(n-1)}(z) \in L \otimes_K K_v(\cC)$
are $\Gal(L/K)$ invariant.  By Proposition \ref{CPPr1} this means that 
for each $\lambda$, the coefficient vector $\oplus_{w|v} A_{w,\lambda}$ 
has the same galois-equivariance properties as $\varphi_{\lambda}(z)$.  
In particular, if $K \subset F_{\lambda} \subset L$ is the smallest 
field of rationality for $\varphi_{\lambda}(z)$, 
then $\oplus_{w|v} A_{w,\lambda} \in F_{\lambda} \otimes_K K_v$.  

Since 
\begin{equation*}
(\prod_{w \in \hS_L} h_w^{D_w})^{nN} \ > \ C_L(\hS_K) \ ,
\end{equation*}
taking $Q_w = h_w^{nN}$ in Proposition \ref{CStrong3} we can find
\index{Strong Approximation theorem!uniform} 
an $A_{\lambda} \in F_{\lambda}$ such that 
\begin{equation*}
\left\{ \begin{array}{ll}
         |A_{\lambda} - A_{w,\lambda}|_w \le h_w^{nN} &
                                \text{for all $w \in \hS_L$\ ,} \\
         |A_{\lambda}|_w \le 1 & \text{for all $w \notin \hS_L$ \ .}
        \end{array} \right.
\end{equation*}
By working with representatives of galois orbits as before,
we can arrange that for each $\sigma \in \Gal(L/K)$ we have
$\sigma(A_{\lambda}) = A_{\lambda^{\prime}}$
if $\sigma(\varphi_{\lambda}) = \varphi_{\lambda^{\prime}}$.  

Put 
\begin{equation*}
\Delta_{w,\lambda} \ = \ A_{\lambda} - A_{w,\lambda} 
\end{equation*}                                                 
for each $w$ and $\lambda$, and put 
\begin{equation*}
H_w(z) \ = \ \sum_{\lambda} \Delta_{w,\lambda} \varphi_{\lambda}(z) \ .
\end{equation*}
Then $\oplus_{w|v} H_w(z) \in L \otimes_K K_v(\cC)$ 
is stable under $\Gal(L/K)$, for each $v \in \hS_K$.  It follows that the
$H_w(z)$ belong to $K_v(\cC)$ and are the same for all $w|v$.  
Let $H_v(z) = H_{w_v}(z)$,\index{distinguished place $w_v$} and expand 
$H_v(z) = \sum_{\lambda=1}^{\Lambda} \Delta_{v,\lambda} \varphi_{\lambda}(z)$.
Then \begin{equation*}
|\Delta_{v,\lambda}|_v \ \le \ h_v^{nN} 
\end{equation*}
for each $v$ and $\lambda$.  

Patch $G_v^{(n-1)}(z)$ by setting\index{patching functions, $G_v^{(k)}(z)$ for $1 \le k \le n$!constructed by patching}
\begin{equation*}
G_v^{(n)}(z) \ = \ G_v^{(n-1)}(z) + H_v(z) \
\end{equation*}
This replaces the low-order coefficients of the $G_v^{(n)}(z)$ 
with\index{patching functions, $G_v^{(k)}(z)$ for $1 \le k \le n$!constructed by patching} 
\index{coefficients $A_{v,ij}$!low-order} 
the $A_{\lambda}$.  

\vskip .1 in
{\bf Conclusion of the Patching Argument.}\index{patching argument!conclusion of global}
\index{global patching when $\Char(K) = 0$!conclusion of the patching argument}
The patching process has now arranged that the $G_v^{(n)}(z) \in K_v(\cC)$
\index{patching functions, $G_v^{(k)}(z)$ for $1 \le k \le n$!$G_v^{(n)}(z) = G^{(n)}(z)$ is independent of $v$}
for $v \in \hS_K$ all coincide with a single function  $G^{(n)}(z)$, 
whose coefficients belong to $L$.  
Fix any $v$, and put $G_w(z) = G(z)$ for all $w|v$;  
then $\oplus_{w|v} G_w(z) \in \oplus_{w|v} L(\cC) \cong L \otimes_K K(\cC)$ 
is invariant under $\Gal(L/K)$, so by Proposition \ref{CPPr1} it belongs to $K(\cC)$.   

For each $v \in \hS_K$, our restrictions on the magnitudes of the 
$\Delta_{v,ij}^{(k)}$ and the $\Delta_{v,\lambda}$ assure that the conclusions
of Theorem \ref{LocalPatch} apply.  Thus 
\begin{equation*} 
\left\{ \begin{array}{l}
\text{If $K_v \cong \CC$, then  
  $\{ z \in \cC_v(\CC_v) : |G^{(n)}(z)|_v \le r_v^{Nn} \} \subset U_v = E_v^0;$} \\
\text{If $K_v \cong \RR$, then} \\ 

\text{\quad $(1)$ the zeros of $G^{(n)}(z)$ all belong to $E_v^0$,\
    and for each component $E_{v,i}$ of $E_v$,} \\
\text{\qquad if $\phi_v(z)$ has $\tau_i$ zeros  in $E_{v,i}$, 
     then $G^{(n)}(z)$ has $T_i = n \tau_i$ zeros in $E_{v,i}$.} \\ 
\text{\quad $(2)$ $\{ z \in \cC_v(\CC_v) : 
         |G^{(n)}(z)|_v \le 2r_v^{nN}\} \ \subset \ U_v$,} \\ 
\text{\quad $(3)$ for each component $E_{v,i}$ contained in $\cC_v(\RR)$, then} \\
\text{\qquad \qquad  
      $G^{(n)}(z)$ oscillates $T_i$ 
           times between $\pm 2 r_v^{nN}$ on  $E_{v,i}$.} \\
  
\text{If $K_v$ is nonarchimedean and $v \in S_{K,0}$, then the zeros of $G^{(n)}(z)$ 
are distinct} \\

\text{\qquad \qquad and belong to $E_v$, and\index{patching functions, $G_v^{(k)}(z)$ for $1 \le k \le n$!roots are confined to $E_v$}
      $\{z \in \cC_v(\CC_v) : G_v^{(n)}(z) \in \cO_v \cap D(0,r_v^{nN}) \} \subset E_v$. } \\

\text{If $K_v$ is nonarchimedean and $v \in \hS_{K,0} \backslash S_{K,0}$, then} \\
\text{ \qquad \qquad
$\{ z \in \cC_v(\CC_v) : |G^{(n)}(z)|_v \le R_v^{Nn} \} = E_v.$ }
\end{array} \right. 
\end{equation*}
On the other hand, for each $v \notin \hS_K$,  
our construction has arranged that in the expansion 
\begin{equation*}
G^{(n)}(z) \ = \ 
\sum_{i=1}^m \sum_{j=0}^{(n-1)N_i} A_{ij} \varphi_{i,nN_i-j}(z)
   + \sum_{\lambda = 1}^{\Lambda} A_{\lambda} \varphi_{\lambda} \ ,
\end{equation*}
all the coefficients belong to $\hcO_v$ and the leading coefficients
belong to $\hcO_v^{\times}$.  Our choice of $\hS_K$
assures that $\cC_v$ and the functions $\varphi_{ij}(z)$
and $\varphi_{\lambda}(z)$ all have good reduction at $v$, 
\index{good reduction}
and the $x_i$ specialize to distinct points $\pmod{v}$.  
Hence $G^{(n)}(z) \pmod{v}$ is a nonconstant function 
with a pole of order $nN_i > 0$ at each $x_i$.
  
Thus for each $v \notin \hS_K$,
\begin{equation*}
\{ z \in \cC_v(\CC_v) : |G^{(n)}(z)|_v \le 1 \} 
  \ = \ \cC_v(\CC_v) \backslash (\bigcup_{i=1}^m B(x_i,1)^{-}) \ = \ E_v \ .
\end{equation*}
  
{\bf Construction of the points in Theorem \ref{aT1-B}.}
\index{global patching when $\Char(K) = 0$!constructing the points in Theorem \ref{aT1-B}}
The patching argument\index{patching argument!conclusion of global} holds for each integer $n > n_3$ 
divisible by $n_0 n_1 n_2$.  For any such $n$, 
the zeros of $G^{(n)}(z)$ satisfy the conditions of the Theorem.
If there are any archimedean $v \in \hS_K$ with $K_v \cong \RR$ 
such that some component $E_{v,i}$ of $E_v$ is contained in $\cC_v(\RR)$, 
or if there are any nonarchimedean $v \in \hS_K$, the construction shows that 
the zeros of $G^{(n)}(z)$ are distinct.  Letting $n \rightarrow \infty$, 
we obtain infinitely many points satisfying the conditions of the Theorem.  

However, if there are no such $v$, then since $\prod_{v \in \hS_K} r_v^{Nn}$
grows arbitrarily large as $n \rightarrow \infty$, 
the number of $\hS_K$-integers $\kappa \in K$ satisfying
$|\kappa|_v \le r_v^{Nn}$ for all $v \in \hS_K$ also becomes arbitrarily large.
For any such $\kappa$, the roots of $G^{(n)}(z) = \kappa$ are points
satisfying the conditions of the Theorem. 
Hence there are infinitely many such points.

This completes the proof of Theorem \ref{aT1-B} when $\Char(K) = 0$. 
\end{proof} 
\index{patching argument!global!when $\Char(K) = 0$|)}

\section{ Proof of Theorem $\ref{aT1-B}$ when $\Char(K) = p > 0$}  
\label{CharPSection}
\index{patching argument!global!when $\Char(K) = p > 0$|(}

When $\Char(K) = p > 0$, the proof of $\ref{aT1-B}$ is similar to that 
when $\Char(K) = 0$, but because there are no archimedean 
places and all the residue fields are lie over the same prime field $\FF_p$, 
many of the details are simpler.    
On the other hand, there are some complications which arise from the fact that $L/K$ may 
be inseparable.  For this reason we carry out the patching process using
the $L^{\sep}$-rational basis rather than the $L$-rational basis.
\index{basis!$L^{\sep}$-rational} 
\index{basis!$L$-rational}   

\medskip 
\begin{proof}[Proof of Theorem $\ref{aT1-B}$ when $\Char(K) = p > 0$] { \ }

Let $K$ be a function field, 
and let $\cC/K$, $\fX$, $\EE_K = \EE = \prod_v E_v$, and $S_K$ be as in Theorem \ref{aT1-B}.  
The overall structure of the proof is similar to that when $\Char(K) = 0$. 

\smallskip
{\bf Stage 1.  Choices of the sets and parameters.}\index{patching parameters!choice when $\Char(K) = p > 0$|ii}
We begin by making the choices governing the patching process: 
\index{patching parameters}

\smallskip
{\bf The place $v_0$.}
\index{global patching when $\Char(K) = p > 0$!Stage 1: Choices of sets and parameters!the place $v_0$}    
Let $\hS_K$ be the finite set of places of $K$ containing $S_K$ and 
satisfying the conditions in (\ref{hS_KList}).  Fix a place $v_0 \in \cM_K \backslash \hS_K$, 
which will play the role of a place ``at $\infty$''. Put  
\begin{equation*}
\hS_K^{+} \ = \ \hS_K \cup \{v_0\} \ .
\end{equation*}
By our choice of $\hS_K$,
the curve $\cC$, the uniformizing parameters\index{uniformizing parameter!normalizes $L$-rational basis}
$g_{x_i}(z)$, and the basis functions $\varphi_{ij}(z), \tphi_{ij}(z)$, 
\index{basis!$L$-rational} 
$\varphi_{\lambda}$, and $\tphi_{\lambda}$ all have good reduction at $v_0$,
\index{good reduction}
and the set $E_{v_0}$ is $\fX$-trivial.  
\index{$\fX$-trivial}

{\bf Summary of the Initial Approximation Theorems.}
\index{global patching when $\Char(K) = p > 0$!Stage 1: Choices of sets and parameters!summary of the initial approximation theorems}   
We will only need the initial approximation theorems concerning $K_v$-simple sets 
\index{$K_v$-simple!set}  
and $\RL$-domains.  Theorem \ref{CCPX1p} below 
\index{$\RL$-domain} 
summarizes Theorems \ref{RLThm}, \ref{CompactThm} and Corollaries \ref{CompactThmCor}
and \ref{EvHvExtraPtsCor} in the context of function fields.  
The main difference from the corresponding results 
when $\Char(K) = 0$ is that we can require that the 
leading coefficients belong to $K_v(x_i)^{\sep}$, not just $K_v(x_i)$.  
\index{coefficients $A_{v,ij}$!leading} 

\begin{theorem}  \label{CCPX1p}\index{initial approximation theorem!summary of initial approximation theorems!when $\Char(K_v) = p > 0$}  
Let $K$ be a function field with $\Char(K) = p > 0$, 
and let $\EE$, and $\fX$ be as in Theorem $\ref{aT1-B}$. 
Then for each place $v$ of $K$,  

\smallskip
\noindent{\ \ $(A)$}  If $v \in S_K$ 
$($so $E_v$ is compact,  $K_v$-simple, and disjoint from $\fX)$, fix a $K_v$-simple decomposition 
\index{$K_v$-simple!set}  
\begin{equation} \label{E_vSimple2p}
E_v \ = \ \bigcup_{\ell=1}^{D_v} B(a_\ell,r_\ell) \cap \cC_v(F_{w_\ell}) \ .
\end{equation} 
and fix $\varepsilon_v > 0$.    
Then there is a compact, $K_v$-simple set $\tE_v \subseteq E_v$  compatible with $E_v$
\index{$K_v$-simple!set!compatible with another set}  
such that 

$(1)$ For each $x_i, x_j \in \fX$ with $x_i \ne x_j$, 
\index{Robin constant!local}
\index{Green's function}
\begin{equation*} 
|V_{x_i}(\tE_v) - V_{x_i}(E_v)| \ < \ \varepsilon_v \ , \quad
|G(x_i,x_j;\tE_v) - G(x_i,x_j;E_v)| \ < \ \varepsilon_v \ ;
\end{equation*}

$(2)$ For each $0 < \beta_v \in \QQ$  
and each $K_v$-symmetric $\vs \in \cP^m(\QQ)$, 
\index{$K_v$-symmetric!probability vector}
there is an integer $N_v \ge 1$  such that 
for each positive integer  $N$ divisible by $N_v$, 
there is an $(\fX,\vs)$-function\index{$(\fX,\vs)$-function!$K_v$-rational} 
$f_v \in K_v(\cC_v)$ of degree $N$ satisfying
  
\quad $(a)$ The zeros $\theta_1, \ldots, \theta_N$ of $f_v$ are distinct and belong to $E_v$.

\quad $(b)$ $f_v^{-1}(D(0,1)) \subseteq \bigcup_{\ell = 1}^{D_v} B(a_{\ell},r_{\ell})$, 
and there is a decomposition $f_v^{-1}(D(0,1)) = \bigcup_{h=1}^N B(\theta_h,\rho_h)$, 
where the balls $B(\theta_h,\rho_h)$ are pairwise disjoint and isometrically parametrizable. 
\index{isometrically parametrizable ball} 
For each $h = 1, \ldots, N$, if $\ell = \ell(h)$ is such that 
$B(\theta_h,\rho_h) \subseteq B(a_{\ell},r_{\ell}),$ put $F_{u_h} = F_{w_\ell}$.
Then $\rho_h \in |F_{u_h}^{\times}|_v$
and $f_v$ induces an $F_{u_h}$-rational scaled isometry from $B(\theta_h,\rho_h)$ to $D(0,1)$,
\index{scaled isometry}  
with
\begin{equation*}
f_v\big( B(\theta_h,\rho_h) \cap \cC_v(F_{u_h})\big) \ = \ \cO_{F_{u_h}} \ , 
\end{equation*}      
such that  $|f_v(z_1)-f_v(z_2)|_v = (1/\rho_h) \|z_1,z_2\|_v$ for all $z_1, z_2 \in B(\theta_h, \rho_h)$. 

\quad $(c)$ The set $H_v := E_v \cap f_v^{-1}(D(0,1))$ is $K_v$-simple and compatible with $E_v$.
\index{$K_v$-simple!set!compatible with another set}
Indeed, 
\begin{equation} \label{E_vPBF2p}
H_v \ = \ 
  \bigcup_{h=1}^N \big(B(\theta_h,\rho_h) \cap \cC_v(F_{u_h})\big) 
\end{equation}
is a $K_v$-simple decomposition of $H_v$ compatible with the 
\index{$K_v$-simple!decomposition!compatible with another decomposition} 
$K_v$-simple decomposition $(\ref{E_vSimple2p})$ of $E_v$, 
which is move-prepared $($see Definition $\ref{MovePreparedDef})$ relative to 
\index{move-prepared}
$B(a_1,r_1), \ldots, B(a_{D_v},r_{D_v})$.  Moreover, 
for each $\ell =1, \ldots, D_v$, 
there is a point $\wbar_\ell \in \big(B(a_\ell,r_\ell) \cap \cC_v(F_{w_\ell})) \big) \backslash H_v$.

\quad $(d)$ For each $x_i \in \fX$, \ the leading coefficient 
\index{coefficients $A_{v,ij}$!leading}\label{`SymbolIndexcvi'}
$c_{v,i} = \lim_{z \rightarrow x_i} f_v(z) \cdot g_{x_i}(z)^{Ns_i}$ 
belongs to $K_v(x_i)^{\sep}$, and   
$\frac{1}{N} \log_v(|c_{v,i}|_v) \ = \ \Lambda_{x_i}(\tE_v,\vs) + \beta_v$.   
 
\vskip .05 in
\noindent{\ \ $(B)$}  If $v \notin S_K$,
$($so $E_v$ is $\fX$-trivial and in particular is an $\RL$-domain disjoint from $\fX)$,
\index{$\fX$-trivial} 
\index{$\RL$-domain} 
put $\tE_v = E_v$.  Then for each $K_v$-symmetric $\vs \in \cP^m(\QQ)$, 
\index{$K_v$-symmetric!probability vector} 
there is an integer $N_v \ge 1$  such that for each positive integer  $N$ divisible by $N_v$, 
there is an $(\fX,\vs)$-function\index{$(\fX,\vs)$-function!$K_v$-rational}
 $f_v \in K_v(\cC_v)$ of degree $N$ such that 

\quad $(a)$ $E_v = \tE_v  = \{z \in \cC_v(\CC_v) : |f_v(z)|_v \le 1 \};$

\quad $(b)$ For each $x_i \in \fX$, \ the leading coefficient 
\index{coefficients $A_{v,ij}$!leading} 
$c_{v,i} = \lim_{z \rightarrow x_i} f_v(z) \cdot g_{x_i}(z)^{Ns_i}$ 
belongs to $K_v(x_i)^{\sep}$, and   
$\frac{1}{N} \log_v(|c_{v,i}|_v) \ = \ \Lambda_{x_i}(\tE_v,\vs) + \beta_v$.   
\end{theorem} 
 
\smallskip
{\bf The $K_v$-simple decompositions of $E_v$  and sets $U_v$, for $v \in S_K$.}
\index{$K_v$-simple!decomposition}
\index{global patching when $\Char(K) = p > 0$!Stage 1: Choices of sets and parameters!the $K_v$-simple decompositions}    
For each $v \in S_K$, the set $E_v$ is compact and $K_v$-simple 
\index{$K_v$-simple!set}  
(see Definition \ref{KvSimple}).
Choose a $K_v$-simple decomposition    
\begin{equation} \label{K_vSimpleDecompFp}
E_v \ = \ \bigcup_{\ell=1}^{D_v} B(a_\ell,r_\ell) \cap \cC_v(F_{w_\ell}) \ .\
\end{equation}   
By refining this decomposition, if necessary, we can assume that 
$U_v := \bigcup_{\ell=1}^{D_v} B(a_\ell,r_\ell)$ is disjoint from $\fX$.   
This decomposition will be fixed for the rest of the construction.    

\smallskip
{\bf The sets $\tE_v$ for $v \in \hS_K^{+}$.}
\index{global patching when $\Char(K) = p > 0$!Stage 1: Choices of sets and parameters!the sets $\tE_v$}    
By hypothesis, $\gamma(\EE,\fX) > 1$ in Theorem \ref{aT1-B}.  
This means that the Green's matrix  $\Gamma(\EE,\fX)$ is negative definite.\index{Green's matrix!negative definite}  
\index{Green's matrix!global}
Suppose $\tEE = \prod_{v \in \hS_K^{+}} \tE_v \times \prod_{v \notin \hS_K^{+}} E_v$ 
is another $K$-rational adelic set compatible with $\fX$.{compatible with $\fX$}  
By the discussion leading to (\ref{NegDefCloseness}), 
there are numbers $\varepsilon_v > 0$ for $v \in \hS_K^{+}$ such that 
such that $\Gamma(\widetilde{\EE},\fX)$ is also negative definite,\index{Green's matrix!negative definite} 
provided that for each $v \in \hS_K^{+}$ 
\begin{equation} \label{NegDefCloseness1p}
\left\{ \begin{array}{ll}
       |G(x_j,x_i;\tE_v)-G(x_j,x_i;E_v)| \ < \ \varepsilon_v &
                  \text{for all $i \ne j$ \ ,} \\
       |V_{x_i}(\tE_v)-V_{x_i}(E_v)| \ < \ \varepsilon_v &
                  \text{for all $i$ \ .}  
        \end{array}  \right. 
\end{equation}
\index{Green's function} 
\index{Robin constant}

For each $v \in S_K$, we will take $\tE_v \subseteq E_v$ to be the set given 
by Theorem \ref{CCPX1p} for $E_v$, relative to the number $\varepsilon_v$ chosen above 
satisfying (\ref{NegDefCloseness1p}).  For each $v \in \hS_K^{+} \backslash S_K$ 
the set $E_v$ is an $\RL$-domain, and we will take $\tE_v = E_v$ 
\index{$\RL$-domain} 
as in Theorem \ref{CCPX1p}.
Put $\tEE_K =  \prod_{v \in \hS_K^{+}} \tE_v \times \prod_{v \notin \hS_K^{+}} E_v$ 
with the sets $\tE_v$ just chosen, and let 
\begin{equation} \label{tV_KChoicep} 
\tV_K \ := \ V(\tEE_K,\fX) \ = \ \val(\Gamma(\tEE_K,\fX)) 
\end{equation}
be the global Robin constant for $\tEE_K$ and $\fX$.  By construction, $\tV_K < 0$.  
\index{Robin constant!global}

\smallskip
{\bf The rational probability vector $\vs$.}\index{patching parameters!choice when $\Char(K) = p > 0$|ii}
\index{global patching when $\Char(K) = p > 0$!Stage 1: Choices of sets and parameters!the probability vector $\vs$}  
By construction, the Green's matrix $\Gamma(\tEE_K,\fX)$ is $K$-symmetric 
and negative definite.\index{Green's matrix!negative definite}
\index{$K$-symmetric!matrix} 
However, in a major simplification from the case when $\Char(K) = 0$, 
there is a matrix $\Gamma_0 \in M_m(\QQ)$ such that 
$\Gamma(\tEE_K,\fX) = \Gamma_0 \cdot \log(p)$. This means that the unique probability vector $\vs$
for which the components of $\Gamma(\tEE_K,\fX) \vs$ are equal has rational coordinates, 
and that $\tV_K \in \QQ \cdot \log(p)$. 
\index{Robin constant!global}   

To see this, note that by (\ref{GreensMatrixSum}) we have 
\begin{equation} 
\Gamma(\tEE_K,\fX) \ = \ \frac{1}{[L:K]} \Gamma(\tEE_L,\fX) \ = \ \frac{1}{[L:K]}
             \sum_{w \in \hS_L^{+}} \Gamma(\tE_w,\fX) \log(q_w) \ .
\end{equation}
\index{Green's matrix!global}
Since each $\tE_w$ is either $L_w$-simple or is an $\RL$-domain, 
\index{$\RL$-domain} 
Proposition \ref{GreenRationalProp} shows that the entries $G(x_i,x_j;\tE_w)$ and $V_{x_i}(\tE_w)$
\index{Green's function}  
\index{Robin constant!local}
in the local Green's matrices belong to $\QQ$. 
For each $w$ there is a natural number $f_w$ such that $q_w = p^{f_w}$, 
so $\Gamma(\tEE_K,\fX) = \Gamma_0 \cdot \log(p)$ with $\Gamma_0 \in M_m(\QQ)$ as claimed.
Clearly $\Gamma_0$ is $K$-symmetric and negative definite.\index{Green's matrix!negative definite}
\index{$K$-symmetric!matrix}

By the same argument as in the proof of Proposition \ref{UniqueSVec}, 
the vector $\vs^{\prime}$ defined by 
\begin{equation*} 
\vs^{\prime} \ = 
\ -\Gamma_0^{-1} \left( \begin{array}{c} 1 \\ \vdots \\ 1 \end{array} \right) 
\end{equation*}  
is $K$-symmetric with positive entries, and since $\Gamma_0 \in M_m(\QQ)$ its entries are rational. 
\index{$K$-symmetric!matrix} 
By suitably scaling $\vs^{\prime}$ we arrive at a $K$-symmetric probability vector $\vs \in \cP^m(\QQ)$
\index{$K$-symmetric!probability vector}
for which the entries of $\Gamma_0 \vs$ are equal.  
Evidently we have 
\begin{equation} \label{GreasyVEq}
\Gamma(\tEE_K,\fX) \, \vs \ = \ \Gamma_0 \, \vs \cdot \log(p) 
\ = \ \left( \begin{array}{c} \tV \\ \vdots \\ \hV \end{array} \right) 
\end{equation}  
for some $\tV \in \QQ \cdot \log(p)$.  By the minimax property defining $\tV_K = \val(\Gamma(\tEE_K,\fX))$
\index{minimax property} 
\index{Robin constant!global}
it must be that $\tV = \tV_K$, and so $\tV_K \in \QQ \cdot \log(p)$.   

This $\vs$ will be fixed for the rest of the construction.

\smallskip
{\bf The local parameters $\eta_v$, $h_v$, $r_v$, and $R_v$.}
\index{patching parameters!choice when $\Char(K) = p > 0$|ii}
\index{global patching when $\Char(K) = p > 0$!Stage 1: Choices of sets and parameters!the parameters $\eta_v$, $h_v$, $r_v$, and $R_v$}   
Since $\tV_K \in \QQ \cdot \log(p)$ and $\tV_K < 0$, 
\index{Robin constant!global}
and since $\log(q_v) = f_v \log(p)$ for each $v$, 
we can choose a collection of numbers $\{\eta_v\}_{v \in \hS_K^{+}}$ 
with $0 < \eta_v \in \QQ$ for each $v$, such that  
\begin{equation} \label{FPP1p}
\sum_{v \in \hS_K^{+}} \eta_v \log(q_v) \ = \ |\tV_K| \ = \ - \tV_K \ .
\end{equation} 
The $\eta_v$ provide the freedom for adjustment needed 
in the construction of the initial approximating functions, 
\index{initial approximating functions $f_v(z)$}
and determine the scaling factors in passing from the initial approximating functions
to the coherent approximating functions. 
\index{coherent approximating functions $\phi_v(z)$}
 
For the place $v_0$, fix a number $r_{v_0}$ such that $1 < r_{v_0} < e^{\eta_{v_0}}$.
Then, choose a set of numbers $\{h_v\}_{v \in \hS_K^{+}}$  
with $\prod_{v \in \hS_K^{+}} h_v > 1$, such that 
\begin{equation} \label{FPP2p}
\left\{ \begin{array}{ll} 1 < h_{v_0} < r_{v_0} & \text{if $v = v_0$\ ,} \\
                          0 < h_v < 1 & \text{if $v \in \hS_K$\ .} \end{array} \right. 
\end{equation}  
Finally, for each $v \in S_K$, fix an $r_v$ with $h_v < r_v < 1$,
and for each $v \in \hS_K \backslash S_K$ put $r_v = 1$.
For each $v \in \hS_K^{+}$, put $R_v = q_v^{\eta_v}$.  
Then $0 < h_v < r_v < R_v$ for each $v$, and  
\begin{equation} \label{FPP3p} 
1 \ < \ \prod_{v \in \hS_K^{+}} h_v \ < \ \prod_{v \in \hS_K^{+}} r_v \ .  
\end{equation} 
Furthermore, $R_v \in |\CC_v^{\times}|_v$ for each $v \in \hS_K^{+}$.  

As in the proof when $\Char(K) = 0$, 
the numbers $h_v$ control how much the Laurent coefficients of the patching functions 
\index{patching functions, $G_v^{(k)}(z)$}
can be changed, and the $r_v$ are ``encroachment bounds'' which limit how close 
certain quantities can come to the $h_v$.

\medskip
{\bf Stage 2. Construction of the Coherent Approximating Functions $\phi_v(z)$.}
\index{coherent approximating functions $\phi_v(z)$|ii}
\index{global patching when $\Char(K) = p > 0$!Stage 2: Constructing approximating functions!the Coherent Approximation theorem}    
We will now construct the coherent approximating functions $\phi_v(z)$, 
modifying the initial approximating functions $f_v(z)$ given by Theorem \ref{CCPX1p}
\index{initial approximating functions $f_v(z)$}
for the sets and parameters chosen above.  Let $J$ be the number 
from the construction of the $L$-rational and $L^{\sep}$-rational bases in 
\S\ref{Chap3}.\ref{LRationalBasisSection}.   

\begin{theorem} \label{CTCX2p}\index{coherent approximating functions $\phi_v(z)$!construction when $\Char(K) = p > 0$}  
Let $\cC$, $K$, $\EE$, $\fX$, and $S_K$ be as in Theorem $\ref{aT1-B}$, 
where $\Char(K) = p > 0$.
Let $\hS_K^{+} \supseteq S_K$ be the finite set of places constructed above.
For each $v \in \hS_K^{+}$, let $\tE_v \subset E_v$ and $0 < h_v < r_v < R_v$ 
be the set and patching parameters constructed above.
\index{patching parameters} 
For each $v \in S_K$, let 
$\bigcup_{\ell=1}^{D_v} B(a_\ell,r_\ell) \cap \cC_v(F_{w_\ell})$ 
be the $K_v$-simple decomposition of $E_v$ chosen above.  
\index{$K_v$-simple!decomposition}    
Let $\vs \in \cP_m(\QQ)$ be the rational probability vector with positive coefficients  
such that $\tV_K = \Gamma(\tEE_K,\fX) \vs$.
\index{Robin constant!global}

Then there are a positive integer $N$ 
and $(\fX,\vs)$-functions\index{$(\fX,\vs)$-function!$K_v$-rational} $\phi_v(z) \in K_v(\cC)$
for $v \in \hS_K^{+}$, of common degree $N$, such that $N_i := Ns_i$ belongs to $\NN$ 
and is divisible by $J$, for each $i =1, \ldots, m$, and 

$(A)$ The $\phi_v(z)$ have the following properties:  

\hskip .15in $(1)$ If $v \in S_K$, then 
\begin{equation} \label{CFBZ1p}
r_v^N \ < \ q_v^{-1/(q_v-1)} \ < \ 1 \ , \quad \text{and}
\end{equation} 

    \hskip .3in $(a)$ the zeros $\theta_1, \ldots, \theta_N$ of $\phi_v(z)$ are distinct
             and belong to $E_v$;

          
    \hskip .3in $(b)$ $\phi_v^{-1}(D(0,1)) = \bigcup_{h=1}^N B(\theta_h,\rho_h)$, 
             where the balls $B(\theta_1,\rho_1), \ldots, B(\theta_N,\rho_N)$ are 
             
    \hskip .5in pairwise disjoint, isometrically parametrizable, 
\index{isometrically parametrizable ball}    
             and contained in $\bigcup_{\ell = 1}^{D_v} B(a_{\ell},r_{\ell});$
           
    \hskip .3in $(c)$ $H_v := \phi_v^{-1}(D(0,1)) \cap E_v$ is $K_v$-simple, 
    \index{$K_v$-simple!decomposition}  
    \index{$K_v$-simple!set}  
            with the $K_v$-simple decomposition 
            \begin{equation*} 
            H_v \ = \ \bigcup_{h=1}^N \big( B(\theta_h,\rho_h) \cap \cC_v(F_{u_h}) \big)
            \end{equation*}

    \hskip .5in compatible with the $K_v$-simple decomposition
\index{$K_v$-simple!decomposition!compatible with another decomposition}   
           $\bigcup_{\ell = 1}^{D_v} \big( B(a_\ell,r_\ell) \cap \cC_v(F_{w_\ell}) \big)$ 
           of $E_v$,
           
    \hskip .5in which is move-prepared relative to $B(a_1,r_1), \ldots, B(a_{D_v},r_{D_v})$.
    \index{move-prepared}
            For each $\ell $,
    
    \hskip .5in  there is a point 
           $\wbar_\ell \in \big(B(a_\ell,r_\ell) \cap \cC_v(F_{w_\ell})) \big) \backslash H_v$.
           
   \hskip .3in $(d)$  For each $h = 1, \ldots, N$, $F_{u_h}/K_v$ is finite and separable. 
           If $\theta_h \in E_v \cap B(a_{\ell},r_{\ell})$, 
          
    \hskip .5in then $F_{u_h} = F_{w_\ell}$, $\rho_h \in |F_{w_\ell}^{\times}|_v$, 
           and $B(\theta_h,\rho_h) \subseteq B(a_\ell,r_\ell);$ and $\phi_v$ induces an
           
     \hskip .5in  
            $F_{u_h}$-rational scaled isometry from $B(\theta_h,\rho_h)$ onto $D(0,1)$
 \index{scaled isometry} 
            with $\phi_v(\theta_h)= 0$,

    \hskip .5in  
           which takes $B(\theta_h,\rho_h) \cap \cC_v(F_{u_h})$ onto $\cO_{u_h}$. 
    
\hskip .15in $(2)$ If $v \in \hS_K^{+} \backslash S_K$, then 
   \begin{equation*}
     E_v \ = \ \tE_v \ = \ \{ z \in \cC_v(\CC_v) : |\phi_v(z)|_v \le R_v^N \} \ . 
   \end{equation*}
   
$(B)$  For each $w \in \hS_L$, put $\phi_w(z) = \phi_v(z)$ if $w|v$, 
and regard $\phi_w(z)$ as an element of $L_w(\cC)$.  
Each $x_i \in \fX$ is canonically embedded in \, $\cC_w(L_w);$  
let $\tc_{w,i} = \lim_{z \rightarrow x_i} \phi_w(z) \cdot g_{x_i}(z)^{N s_i}$
be the leading coefficient of $\phi_w(z)$ at $x_i$.  Then for each $i$ 
\index{coefficients $A_{v,ij}$!leading} 
\begin{equation} \label{LeadingCoherenceCond}
    \sum_{w \in \hS_L^{+}} \log_w(|\tc_{w,i}|_w) \log(q_w) \ = \ 0 \ . 
\end{equation} 
Moreover, there are a positive integer $n_0$ and a 
$K$-symmetric set of $\hS_L^{+}$-units $\mu_1, \ldots, \mu_m \in L$,
\index{$K$-symmetric!system of units}
with $\mu_i \in K(x_i)^{\sep}$ for each $i$, 
such that $|\tc_{w,i}^{n_0}|_w = |\mu_i|_w$ for each $w \in \hS_L^{+}$ 
and each $i = 1, \ldots, m$.
\end{theorem}

\begin{proof}
Because there are no archimedean places, the proof is simpler than when $\Char(K) = 0$.  
It consists of choosing a collection of initial approximating functions 
\index{initial approximating functions $f_v(z)$}
$f_v(z)$ of common degree $N$ using Theorem \ref{CCPX1p}, then scaling them 
so their leading coefficients satisfy (\ref{LeadingCoherenceCond}). 
\index{coefficients $A_{v,ij}$!leading} 

\smallskip
{\bf The choice of $N$.}\index{patching parameters!choice when $\Char(K) = p > 0$|ii}
\index{global patching when $\Char(K) = p > 0$!Stage 2: Constructing approximating functions!the choice of $N$}  
For each $v \in S_K$, put $\beta_v = \eta_v$ 
(where $0 < \eta_v \in \QQ$ is the number from (\ref{FPP1p}))
and let $N_v > 0$ be the integer given by Theorem \ref{CCPX1p}(A.2) 
for $E_v$, $\tE_v$, $\vs$, $\beta_v$ and the $K_v$-simple decomposition 
\index{$K_v$-simple!decomposition}  
$E_v = \bigcup_{\ell = 1}^{D_v} \big(B(a_\ell,r_\ell) \cap \cC_v(F_{w_\ell})\big)$ chosen above.   
For each $v \in \hS_K^{+} \backslash S_K$, 
let $N_v$ be as given by Theorem \ref{CCPX1p} for $E_v = \tE_v$ and $\vs$ as chosen above.

Fix an integer $N > 0$ be an integer which satisfies the following conditions: 
\begin{enumerate}
\item $N$ is divisible by $N_v$, for each $v \in \hS_K^{+}$;

\item  $N_i := N s_i$ belongs to $\NN$  
and is divisible by $J$ for each $i = 1, \ldots, m$;     

\item $N \cdot \eta_v \in \NN$ for each $v \in \hS_K^{+}$, 
where the $\eta_v$ are as in (\ref{FPP1p});

\item $N$ is large enough that 
\begin{itemize}
 \item $N s_i > J$ for each $i = 1, \ldots, m$; 
 \item $r_v^N  < q_v^{-1/(q_v-1)} < 1$, for each $v \in S_K$.
 \end{itemize} 
 \vskip -.35in
\begin{equation} \label{hS_KList1} \end{equation}  
 \vskip .2in
\end{enumerate} 
In particular (\ref{CFBZ1p}) holds.

\smallskip
{\bf The choice of the Initial Approximating Functions $f_v(z)$.}
\index{global patching when $\Char(K) = p > 0$!Stage 2: Constructing approximating functions!the Initial approximating functions}   
\index{coherent approximating functions $\phi_v(z)$!construction when $\Char(K) = p > 0$!choice of the initial approximating functions}
We will apply Theorem \ref{CCPX1p} with the parameters chosen above.
For each $v \in S_K$, take $\beta_v = \eta_v$ (with $0 < \eta_v \in \QQ$ as in (\ref{FPP1p}))
and let $f_v(z) \in K_v(\cC)$ 
be the $(\fX,\vs)$-function\index{$(\fX,\vs)$-function!$K_v$-rational} 
of degree $N$ given by Theorem
\ref{CCPX1p}(A.2) with $\frac{1}{N} \log_v(|c_{v,i}|_v) = \Lambda_{x_i}(\tE_v,\vs) + \beta_v$
and $c_{v,i} \in K_v(x_i)^{\sep}$  for each $i$.  
For each $v \in \hS_K^{+} \backslash S_K$, let $f_v(z) \in K_v(\cC)$ 
be the $(\fX,\vs)$-function\index{$(\fX,\vs)$-function!$K_v$-rational} 
of degree $N$ from Theorem \ref{CCPX1p}(B),
with $\frac{1}{N} \log_v(|c_{v,i}|_v) = \Lambda_{x_i}(\tE_v,\vs)$ for each $i$. 

Each $f_v(z)$ has the mapping properties from Theorem \ref{CCPX1p}.  In particular, 
if $v \in S_K$, then  
$H_v := f_v^{-1}(D(0,1)) \cap E_v$ has a $K_v$-simple decomposition 
\index{$K_v$-simple!decomposition}  
$H_v = \bigcup_{h=1}^N \big(B(\theta_h,\rho_h) \cap \cC_v(F_{u_h})\big)$
compatible with the $K_v$-simple decomposition 
\index{$K_v$-simple!decomposition!compatible with another decomposition} 
$E_v = \bigcup_{\ell = 1}^{D_v} \big(B(a_\ell,r_\ell) \cap \cC_v(F_{w_\ell})) \big)$,
which is move-prepared relative to $B(a_1,r_1), \ldots, B(a_{D_v},r_{D_v})$. 
\index{move-prepared} 
Here $\theta_1,\ldots, \theta_N$ are the zeros of $f_v(z)$, 
$\rho_h \in |F_{u_h}^{\times}|_v$, and 
$f_v$ induces an $F_{u_h}$-rational scaled isometry from $B(\theta_h,\rho_h)$ to $D(0,1)$ 
\index{scaled isometry} 
which maps $B(\theta_h,\rho_h) \cap \cC_v(F_{u_h})$ onto $\cO_{u_h}$.  
For each $\ell = 1, \ldots, D_v$, there is a point 
$\wbar_\ell \in \big(B(a_\ell,r_\ell) \cap \cC_v(F_{w_\ell})) \big) \backslash H_v$.

\smallskip
{\bf The choice of the Coherent Approximating Functions $\phi_v(z)$.}  
\index{coherent approximating functions $\phi_v(z)$!construction when $\Char(K) = p > 0$!choice of the coherent approximating functions} 
\index{coherent approximating functions $\phi_v(z)$|ii}
\index{global patching when $\Char(K) = p > 0$!Stage 2: Constructing approximating functions!the Coherent approximating functions}
If $v \in S_K$, put $\kappa_v = 1$.  
If $v \in \hS_K^{+} \backslash S_K$, put $\kappa_v = \pi_v^{-N \eta_v}$, 
where $0 < \eta_v \in \QQ$ is as in (\ref{FPP1p}).  
Our choice of $N$ required that $N \eta_v \in \NN$,
so $\kappa_v \in K_v^{\times}$ and $|\kappa_v|_v = R_v^N > 1$.   
For each $v \in \hS_K^{+}$, put 
\begin{equation*}
\phi_v(z) \ = \ \kappa_v f_v(z) \ \in \ K_v(\cC) \ .
\end{equation*}
For each $v$ and each $i$, the leading coefficient $\tc_{v,i}$
of $\phi_v(z)$ at $x_i$ is given by $\tc_{v,i} = \kappa_v c_{v,i}$, 
so $\tc_{v,i} \in K_v(x_i)^{\sep}$.   
By our choices of the $\beta_v$ and $\kappa_v$, for each $v \in \hS_K^{+}$ we have 
\begin{equation} \label{FPhi_vLeadingp} 
\frac{1}{N} \log_v(|\tc_{v,i}|_v) \ = \ 
                       \Lambda_{x_i}(\tE_v,\vs) + \eta_v  \ .                         
\end{equation} 
Furthermore, the mapping properties of the $f_v(z)$ from Theorem \ref{CCPX1p}, 
together with our choice of the $\kappa_v$, 
yield the following mapping properties for the $\phi_v(z)$.   
\begin{enumerate}
\item If $v \in S_K$ then properties (a)-(d) in Theorem \ref{CTCX2p}(A.1) hold for $\phi_v(z)$
  and $H_v$.  Indeed, since $\kappa_v = 1$ for $v \in S_K$, then $\phi_v(z) = f_v(z)$ 
  so the mapping properties of $\phi_v(z)$ are inherited from those of $f_v(z)$.
  

\item If $v \in \hS_K^{+} \backslash S_K$, then 
$E_v = \tE_v = \{z \in \cC_v(\CC_v) : |\phi_v(z)|_v \le R_v^N \}$.  
\end{enumerate} 

\smallskip
{\bf Coherence of the leading coefficients.} 
\index{coefficients $A_{v,ij}$!leading}
\index{global patching when $\Char(K) = p > 0$!Stage 2: Constructing approximating functions!adjusting the leading coefficients}
\index{coherent approximating functions $\phi_v(z)$!construction when $\Char(K) = p > 0$!making the leading coefficients $S$-subunits}
To understand the leading coefficients of the $\phi_v(z)$, 
we must consider them over the fields $L_w$ for $w \in \hS_L^{+}$, 
since $\fX$ is canonically a subset of $\cC(L)$ 
and of $\cC_w(L_w)$ for each $w$, and the
uniformizer $g_{x_i}(z) \in L(\cC)$ is canonically an element of $L_w(\cC)$.  

For each $w \in \hS_L$, put $\phi_w(z) = \phi_v(z)$ if $w|v$, 
and view $\phi_w(z)$ as an element of $L_w(\cC)$.  Although the functions
$\phi_w(z)$ for $w|v$ are all the same, the points of $\fX$, which are their poles, 
are identified differently.  For each $i$ and $w$,   
let $\tc_{w,i} = \lim_{z \rightarrow x_i} \phi_w(z) \cdot g_{x_i}(z)^{Ns_i}$
be the leading coefficient of $\phi_w(z)$ at $x_i$.  
\index{coefficients $A_{v,ij}$!leading} 
Let $\sigma_w : L \hookrightarrow \CC_v$
be an embedding which induces the place $w$, and for each $i = 1, \ldots, m$ let 
$\sigma_w(i)$ be the index $j$ for which $\sigma_w(x_i) = x_j$ (where we identify $x_j$ with 
its image in $\cC_v(\CC_v)$ given by the fixed embedding of $\tK$ in $\CC_v$).  
Then $\tc_{w,i} = \tc_{v,\sigma_w(i)}$.    

Recall that $\Gamma(\tEE_K,\fX) = \frac{1}{[L:K]} \Gamma(\tEE_L,\fX)$.   
It follows from (\ref{GreensMatrixSum}) that
\begin{equation} \label{FJIG3p}
[L:K] \cdot \Gamma(\tEE_K,\fX) \ = \ 
\sum_{w \in \hS_L^{+}} \Gamma(\tE_w,\fX) \log(q_w) \ .
\end{equation} 
Just as when $\Char(K) = 0$, Lemma \ref{SumSigmaLem} shows that   
for each $w \in \hS_L^{+}$, and each $i = 1, \ldots, m$,
the $i^{th}$ coordinate of $\Gamma(\tE_w,\fX) \vs$ satisfies  
\begin{equation} \label{FJIG4p}
(\Gamma(\tE_w,\fX) \vs)_i \cdot \log(q_w) \ = \ 
[L_w:K_v] \cdot \Lambda_{\sigma_w(x_i)}(\tE_v,\vs) \log(q_v) \ .
\end{equation}  

We can now prove (\ref{LeadingCoherenceCond}). Since  $\tc_{w,i} = \tc_{v,\sigma_w(i)}$, 
it follows from (\ref{FPhi_vLeading}) that   
\begin{eqnarray}
\frac{1}{N} \sum_{w \in \hS_L^{+}} \log_w(|\tc_{w,i}|_w) \log(q_w) 
     & = & \sum_{v \in \hS_K^{+}} \sum_{w|v} [L_w:K_v] 
               \big(\frac{1}{N}\log_v(|\tc_{v,\sigma_w(x_i)}|_v)\big) \log(q_v) \notag \\
     & = & \sum_{v \in \hS_K^{+}} \sum_{w|v} [L_w:K_v] 
        \big(\Lambda_{\sigma_w(x_i)}(\tE_v,\vs) + \eta_v\big) \log(q_v) \ .
         \label{BigSumFormulap} 
\end{eqnarray} 
By Lemma \ref{SumSigmaLem} and our choice of $\vs$ in (\ref{GreasyVEq}), 
\begin{eqnarray}
\sum_{v \in \hS_K^{+}} \sum_{w|v} [L_w:K_v] \Lambda_{\sigma_w(x_i)}(\tE_v,\vs) \log(q_v) 
& = & \sum_{w \in \hS_L^{+}} (\Gamma(\tE_w,\fX) \vs)_i \log(q_w) \notag \\
& = & (\Gamma(\tEE_L,\fX) \vs)_i \ = \ [L:K]  \cdot \tV_K \ .
         \label{FLambdaSump}
\end{eqnarray}  
\index{Robin constant!global}
By our choice of the $\eta_v$ in (\ref{FPP1p}), 
\begin{equation}
\sum_{v \in \hS_K^{+}} \sum_{w|v} [L_w:K_v] \eta_v \log(q_v) 
\ = \ [L:K] \sum_{v \in \hS_K^{+}} \eta_v \log(q_v) \ = \ -[L:K] \cdot \tV_K \ .
\label{FEtaSump}
\end{equation}
\index{Robin constant!global}
Combining (\ref{BigSumFormulap}), (\ref{FLambdaSump}), and (\ref{FEtaSump}) gives 
\begin{equation*}
\sum_{w \in \hS_L^{+}} \log_w(|\tc_{w,i}|_w) \log(q_w) \ = \ 0 \ ,
\end{equation*} 
which is (\ref{LeadingCoherenceCond}). 

\smallskip
The final assertion in Theorem \ref{CTCX2p} concerns the existence of a $K$-symmetric system of 
$S_L^{+}$-units $\mu_1, \ldots, \mu_m$ with $\mu_i \in K(x_i)^{\sep}$ for each $i$,
\index{$K$-symmetric!system of units}
and a positive integer $n_0$ such that  
$|\tc_{w,i}^{n_0}|_w = |\mu_i|_w$ for each $w \in \hS_L^{+}$ and each $i = 1, \ldots, m$.

To show this, note that by our choice of $N$ (see (\ref{hS_KList1})) 
$N_i := N s_i$ is an integer divisible by $J$ for each $i =1, \ldots, m$.  By the
construction of the $L$-rational and $L^{\sep}$-rational bases in 
\S\ref{Chap3}.\ref{LRationalBasisSection}, this means that the basis functions 
\index{basis!$L$-rational} 
$\varphi_{i,N_i} = \tphi_{i,N_i}$ for each $i$.  For each $v \in \hS_K^{+}$, 
the leading coefficients $\tc_{v,i}$ of $\phi_v(z)$ are $K_v$-symmetric 
\index{coefficients $A_{v,ij}$!leading} 
\index{$K_v$-symmetric!set of numbers}
since $\phi_v$ is $K_v$-rational, with each $\tc_{v,i} \in K_v(x_i)^{\sep}$
by construction.  The $L^{\sep}$-rational basis is $K_v$-symmetric
\index{basis!$L^{\sep}$-rational} 
\index{$K_v$-symmetric!set of functions}
by construction, so the function
\begin{equation*}
\phi_v^0(z) \ := \ \sum_{i=1}^m \tc_{v,i} \varphi_{i,N_i}(z) 
       \ = \ \sum_{i=1}^m \tc_{v,i} \tphi_{i,N_i}(z) 
\end{equation*}
consisting of the leading terms of $\phi_v(z)$, is $K_v$-rational.  

Consider the field $H = L^{\sep}$.  Since $L/L^{\sep}$ is purely inseparable, 
for each place of $w_0$ of $H$ there is a unique place $w$ of $L$ with $w|w_0$,
and $w$ is totally ramified over $w_0$ with ramification index $[L:K]^{\insep}$.    
Since $H/K$ is galois, the group $\Aut(L/K) \cong \Gal(H/K)$ 
acts transitively on the places $w|v$ of $L$, and the places $w_0|v$ of $H$.  
For each $w_0|v$, put $\phi_{w_0}^0(z) = \phi_v^0(z)$, 
regarding $\phi_{w_0}^0(z) \in K_v(\cC)$ as an element of $H_{w_0}(\cC)$.   
Write $\tc_{w_0,i}$, $i=1,\ldots, m$ for its leading coefficients; 
\index{coefficients $A_{v,ij}$!leading} 
thus if $w$ is the place of $L$ over $w_0$, 
then $\tc_{w_0,i} = \tc_{w,i} = \tc_{v,\sigma_w(i)}$.  

First fix $x_i \in \fX$, and put $F = K(x_i)^{\sep}$.  
Let $\hS_F^{+}$ be the set of places of $F$ above $\hS_K^{+}$.  
For each $v \in \hS_K^{+}$, since the $\phi_{w_0}^0(z) \in K_v(\cC)$
are the same for all $w_0|v$, applying Proposition \ref{CPPr1p} 
to $\oplus_{w_0|v} \phi_{w_0}^0(z)$ tells us that
$\oplus_{w|v} \tc_{w,i} = \oplus_{w_0|v} \tc_{w_0,i} \in \oplus_{w_0|v} H_{w_0}$ 
actually belongs to $\oplus_{u|v} F_u$, embedded semi-diagonally in  $\oplus_{w_0|v} H_{w_0}$.
Write $\oplus_{u|v} \tc_{u,i}$ for the element of $\oplus_{u|v} F_u$
that induces it.  Then by (\ref{LeadingCoherenceCond}), 
(and the fact that $[L: K] = \sum_{w|v} [L_w:K_v]$,
even when $\Char(K) = p > 0$; see (\cite{RR1}, p.321))
\begin{eqnarray*}
\sum_{u \in \hS_F^{+}} \log_u(|\tc_{u,i}|_u) \log(q_u)
 & = & \frac{1}{[H:F]} \sum_{w_0 \in \hS_H^{+}} \log_{w_0}(|\tc_{w_0,i}|_{w_0}) \log(q_{w_0}) \\ 
   & = & \frac{1}{[L:F]} \sum_{w \in \hS_L^{+}} \log_w(|\tc_{w,i}|_w) \log(q_w)
   \ = \ 0 \ .
\end{eqnarray*}
By Proposition \ref{CPSUnit2FF} 
there are an $\hS_F^{+}$-unit $\mu_i \in F$, and an integer $n_i$ such that 
\begin{equation*}     
      |\tc_{u,i}^{n_i}|_u \ = \ |\mu_i|_u 
\end{equation*} 
for each  $u \in \hS_F^{+}$.

Now let $x_i$ vary;  we will arrange for the $\mu_i$  to be $K$-symmetric.
\index{$K$-symmetric!system of units}
By Proposition \ref{CPPr1p}, the $\oplus_{w_0|v} \tc_{w_0,i}$ are $K$-symmetric. 
\index{$K$-symmetric!vector}
For each $\Aut(L/K)$-orbit $\fX_{\ell} \subset \fX$, 
fix an $x_{i_\ell} \in \fX_{\ell}$. 
For each $x_j \in \fX_{\ell}$, choose $\sigma \in \Aut(L/K)$ with $\sigma(x_{i_\ell}) = x_j$, 
and replace $\mu_j$ with $\sigma(\mu_{i_\ell})$.  These $\mu_j$ are independent of
the choice of $\sigma$ with $\sigma(x_{i_\ell}) = x_j$, and are $K$-symmetric. 
\index{$K$-symmetric!system of units} 
After further replacing $\mu_1, \ldots, \mu_m$ with appropriate powers of themselves, 
we can assume there is a number $n_0$ such that $n_i = n_0$, for all $i$. 
Thus $\mu_1, \ldots, \mu_m$ form a $K$-symmetric system of $\hS_L^{+}$-units,
\index{$K$-symmetric!system of units}
with $\mu_i \in K(x_i)^{\sep}$ for each $i$, and $|\tc_{w,i}^{n_0}|_w = |\mu_i|_w$
for each $w \in \hS_L^{+}$.   

This completes the proof of Theorem \ref{CTCX2p}.
\end{proof}

{\bf Stage 3. The Patching Construction.}
\index{patching theorem!global patching construction when $\Char(K) = p > 0$|ii}
\index{global patching when $\Char(K) = p > 0$!Stage 3: The global patching construction!overview}    
In the function field case,\index{patching argument!global|ii}
there are several differences in the patching argument from the number field case. 

One complication arises from the fact that $L/K$ may be inseparable.
In order to preserve the $K_v$-rationality of the patching 
functions $G_v^{(k)}(z)$,\index{patching functions, $G_v^{(k)}(z)$ for $1 \le k \le n$!are $K_v$-rational}
it is helpful to expand them\index{patching functions, $G_v^{(k)}(z)$ for $1 \le k \le n$!expansion of}   
in terms of the $L^{\sep}$-rational basis $\{\tphi_{ij}, \tphi_{\lambda}\}$ 
\index{basis!$L^{\sep}$-rational} 
rather than the $L$-rational basis.  
\index{basis!$L$-rational} 
However, the several the basis functions $\tphi_{ij}$ can contribute to 
poles of the $G_v^{(k)}(z)$ of the same order.\index{patching functions, $G_v^{(k)}(z)$ for $1 \le k \le n$!expansion of} 
To deal with this, instead of patching the coefficients of the $\tphi_{ij}$ in $\prec_N$ order, 
\index{order!$\prec_N$}
we patch all the coefficients in a band simultaneously.
\index{band!coefficients patched by bands!simultaneously when $\Char(K) = p > 0$} 
 
In addition, in the global patching construction\index{patching argument!global} 
we cannot use the same method for patching the high order coefficients 
\index{patching!high-order coefficients}
\index{coefficients $A_{v,ij}$!high-order} 
as when $\Char(K) = 0$, because there are no archimedean places where a `magnification argument'
\index{magnification argument} 
can apply.  Instead, by the choice of $n$, in the local patching construction
\index{patching argument!local}  
we arrange that all the high order coefficients (apart from the leading coefficient) are $0$, 
so they do not need to be patched. This is possible because $\Char(K) = p$.

The following theorem summarizes the local patching constructions
\index{patching argument!local|ii} 
proved in Theorems \ref{DCPPatch1p} and \ref{DCPCPatch1p} below.  
After stating the theorem, we compare the patching constructions 
when $\Char(K) = 0$ and $\Char(K) = p > 0$, and 
establish some estimates for the coefficients  
needed to carry out the patching process in bands.  
We then choose the parameters $\kbar$ and $n$, 
\index{band!coefficients patched by bands}\index{confinement argument} 
and give the details of the patching process. 
\index{global patching when $\Char(K) = p > 0$!Stage 3: The global patching construction!summary of the Local patching theorems}

\begin{theorem} \label{LocalPatchp}\index{patching theorem!summary of local patching theorems!when $\Char(K_v) = p > 0$} 
Let $K$ be a function field.  Let $\cC/K$, $\EE$, $\fX$, and $S_K$ be as in Theorem $\ref{aT1-B}$.
Let $\hS_K^{+} \supseteq S_K$ be the finite set of places satisfying conditions $(\ref{hS_KList})$,
together with $v_0$.  
For each $v \in \hS_K^{+}$, let $\tE_v \subset E_v$, and $0 < h_v < r_v < R_v$ 
be the set and patching parameters constructed in Stage $1$ above.
\index{patching parameters} 
For each $v \in S_K$, let 
$\bigcup_{\ell=1}^{D_v} B(a_\ell,r_\ell) \cap \cC_v(F_{w_\ell})$ 
be the chosen $K_v$-simple decomposition of $E_v$; 
\index{$K_v$-simple!decomposition}  
by construction, $U_v = \bigcup_{\ell=1}^{D_v} B(a_\ell,r_\ell)$ is disjoint from $\fX$.
Let the rational probability vector 
$\vs \in \cP^m(\QQ)$ be as in $(\ref{GreasyVEq})$, 
and let the natural number $N$ and the coherent approximating functions $\{\phi_v(z)\}_{v \in \hS_K^{+}}$ 
\index{coherent approximating functions $\phi_v(z)$} 
be the ones constructed in Theorem $\ref{CTCX2p}$.  
Put $N_i = Ns_i$ for $i=1, \ldots, m$.
By construction $N_i \in \NN$ and $J|N_i$ for each $i$, 
and for each $v \in \hS_K^{+}$ and each $i$, the leading coefficient 
\index{coefficients $A_{v,ij}$!leading} 
$\tc_{v,i} = \lim_{z \rightarrow x_i} \phi_v(z) \cdot g_{x_i}(z)^{N_i}$ 
belongs to $K_v(x_i)^{\sep}$.     

For each $v \in \hS_K^{+}$,
Theorem $\ref{DCPPatch1p}$ or $\ref{DCPCPatch1p}$ 
provides a number $k_v > 0$ determined by $E_v$ and $\phi_v(z)$,   
representing the minimal number of `high-order' stages  
in the local patching process at $v$.
\index{patching argument!local}  
Let $\kbar \ge k_v$ be a fixed integer.
 
Then for each $v \in \hS_K^{+}$, there are an integer $n_v > 0$ and a number $0 < B_v < 1$, 
depending on $\kbar$, $E_v$, and $\phi_v(z)$, such
that for each sufficiently large integer $n$ divisible by $n_v$, 
the local patching process at $v$ can be carried out as follows:

Put $G_v^{(0)}(z) = Q_{v,n}(\phi_v(z))$,\index{patching functions, initial $G_v^{(0)}(z)$!construction of} where 
\begin{equation*}
\left\{ \begin{array}{l}
   \text{If $v \in S_K$, then $Q_{v,n}(x) = S_{n,v}(x)$} \\
   \text{\quad is the Stirling polynomial of degree $n$ for $\cO_v$
                   $($see $(\ref{FST1}));$} \\
   \text{If $v \in \hS_K^{+} \backslash S_K$, then $Q_{v,n}(x) = x^n$.}
   \end{array} \right.
\end{equation*}\index{Stirling polynomial!for $\cO_v$} 
Then all the zeros of $G_v^{(0)}(z)$ belong to $E_v$, 
and if $v \in S_K$ they are distinct.\index{patching functions, initial $G_v^{(0)}(z)$!roots are confined to $E_v$} 
For each $x_i$, the leading coefficient of $G_v^{(0)}(z)$ 
at $x_i$ is $\tc_{v,i}^n$\index{patching functions, initial $G_v^{(0)}(z)$!leading coefficients of} 
\index{coefficients $A_{v,ij}$!leading} 
and when $G_v^{(0)}(z)$ is expanded in terms of the 
$L$-rational basis as\index{patching functions, initial $G_v^{(0)}(z)$!expansion of} 
\index{basis!$L$-rational} 
\begin{equation*} 
 G_v^{(0)}(z) \ = \ \sum_{i=1}^m \sum_{j=0}^{(n-1)N_i-1} A_{v,ij} \varphi_{i,nN_i-j}(z) 
   + \sum_{\lambda=1}^{\Lambda} A_{v,\lambda} \varphi_{\lambda} \ ,
\end{equation*} 
then $A_{v,ij} = 0$ for all $(i,j)$ with $1 \le j < \kbar N_i$.  
\index{patching functions, initial $G_v^{(0)}(z)$!when $\Char(K) = p > 0$!high-order coefficients are $0$}

For each $k$, $1 \le k \le n-1$, 
let $\{\Delta_{v,ij}^{(k)} \in L_{w_v}\}_{(i,j) \in \Band_N(k)}$\index{distinguished place $w_v$} 
be a $K_v$-symmetric set of numbers satisfying
\index{$K_v$-symmetric!set of numbers} 
\begin{equation} \label{FBound1p}
\left\{ \begin{array}{ll}
\text{$|\Delta_{v,i0}^{(1)}|_v \le B_v $ 
            and $\Delta_{v,ij}^{(1)} = 0$ for $j = 1, \ldots, N_i-1$,}  
               & \text{if \ $k = 1$\ ,} \\
\text{$\Delta_{v,ij}^{(k)} = 0$ for $j = (k-1)N_i, \ldots, kN_i - 1$,} 
               & \text{if \ $k = 2, \ldots, \kbar$\ , } \\
|\Delta_{v,ij}^{(k)}|_v  \le  h_v^{kN},  
               & \text{if \ $k = \kbar+1, \ldots, n-1$\ ,}
      \end{array} \right.
\end{equation}
such that $\Delta_{v,i0}^{(1)} \in K_v(x_i)^{\sep}$ for each $i$ 
and such that for each $k = \kbar + 1, \ldots, n-1$ 
\begin{equation} \label{DFCCV1Bp}
\Delta_{v,k}(z) \ := \ 
\sum_{i=1}^m \sum_{j=(k-1)N_i}^{kN_i-1} \Delta_{v,ij}^{(k)} \cdot   \varphi_{i,(k+1)N_i - j} 
 \ \in \ K_v(\cC) \ .
\end{equation} 
For $k = n$, let 
$\{\Delta_{v,\lambda}^{(n)} \in L_{w_v}\}_{1 \le \lambda \le \Lambda}$\index{distinguished place $w_v$} 
be a $K_v$-symmetric set of numbers such that 
\index{$K_v$-symmetric!set of numbers}
\begin{equation} \label{FBound2p} 
|\Delta_{v,\lambda}^{(n)}|_v \ \le \ h_v^{nN}   
\end{equation} 
for each $\lambda$, and 
\begin{equation} \label{DFCCF2Bp} 
\Delta_{v,n}(z) \ := \ 
\sum_{\lambda = 1}^{\Lambda} \Delta_{v,\lambda}^{(n)} \cdot \varphi_{\lambda} \ \in K_v(\cC) \ . 
\end{equation} 

Then one can inductively construct 
 $(\fX,\vs)$-functions\index{$(\fX,\vs)$-function!$K_v$-rational}
\index{patching functions, $G_v^{(k)}(z)$ for $1 \le k \le n$!constructed by patching} 
$G_v^{(1)}(z), \ldots, G_v^{(n)}(z)$ in $K_v(\cC)$, 
of common degree $Nn$, such that$:$

\smallskip
\noindent{$(A)$} For each $k = 1, \ldots, n$, $G_v^{(k)}(z)$ 
is obtained from $G_v^{(k-1)}(z)$\index{patching functions, $G_v^{(k)}(z)$ for $1 \le k \le n$!constructed by patching}
as follows$:$  

$(1)$  When $k = 1$, the local patching process at $v$ provides a
\index{patching argument!local}
$K_v$-symmetric set of functions 
\index{$K_v$-symmetric!set of functions}
$\ttheta_{v,10}^{(1)}(z), \ldots, \ttheta_{v,m0}^{(1)}(z) \in L_{w_v}^{\sep}(\cC)$\index{distinguished place $w_v$}
such that 
\begin{equation*}
G_v^{(1)}(z) \ = \ 
    G_v^{(0)}(z) + \sum_{i=1}^m \Delta_{v,i0}^{(1)} \cdot \ttheta_{v,i0}^{(1)}(z) \ , 
\end{equation*}
where for each $i = 1, \ldots, m$,
$\ttheta_{v,i0}^{(1)}(z) \in K_v(x_i)^{\sep}(\cC)$ has the form 
\begin{equation*}
\ttheta_{v,i0}^{(1)}(z) \ = \ \tc_{v,i}^{n} \varphi_{i,nN_i}(z) + \tTheta_{v,i0}^{(1)}(z) 
\end{equation*}   
for an $(\fX,\vs)$-function\index{$(\fX,\vs)$-function!$K_v$-rational} $\tTheta_{v,i0}^{(1)}(z)$
with a pole of order at most $(n-\kbar)N_{i^\prime}$ at each $x_{i^\prime}$.  
Thus, in passing from $G_v^{(0)}(z)$ to $G_v^{(1)}(z)$, 
each of the leading coefficients\index{patching functions, $G_v^{(k)}(z)$ for $1 \le k \le n$!leading coefficients of}
$A_{v,i0} = \tc_{v,i}^n$ is replaced with $\tc_{v,i}^n + \Delta_{v,i0}^{(1)} \cdot \tc_{v,i}^n$,
and the coefficients $A_{v,ij}$ for $1 \le j < \kbar N_i$ remain $0$.   

$(2)$ For $k = 2, \ldots, \kbar$, 
we have $G_v^{(k)}(z) = G_v^{(k-1)}(z)$.\index{patching functions, $G_v^{(k)}(z)$} 

$(3)$ For $k = \kbar + 1, \ldots, n-1 $, 
we have\index{patching functions, $G_v^{(k)}(z)$ for $1 \le k \le n$!constructed by patching} 
\begin{equation} \label{GvkOmegaModification} 
G_v^{(k)}(z) \ = \ G_v^{(k-1)}(z)
     + \Delta_{v,k}(z) \cdot F_{v,k}(z) + \Theta_v^{(k)}(z) \ ,    
\end{equation}
where 

\quad $(a)$ $\Delta_{v,k}(z) = \sum_{(i,j) \in \Band_N(k)} \Delta_{v,ij}^{(k)} \varphi_{i,(k+1)N_i-j}(z)$
 belongs to $K_v(\cC)$ by $(\ref{DFCCV1Bp});$

\quad $(b)$ $F_{v,k}(z) \in K_v(\cC)$ 
is an $(\fX,\vs)$-function\index{$(\fX,\vs)$-function!$K_v$-rational} 
determined by the local patching process\index{patching argument!local}

\qquad using $G_v^{(k-1)}(z)$,\index{patching functions, $G_v^{(k)}(z)$ for $1 \le k \le n$!roots are confined to $E_v$}
whose roots belong to $E_v$. For each $x_i$, it has a pole of order 

\qquad $(n-k-1)N_i$ at $x_i$, and its leading coefficient 
\index{coefficients $A_{v,ij}$!leading} 
$d_{v,i} = \lim_{z \rightarrow x_i} F_{v,k}(z) \cdot g_{x_i}(z)^{(n-k-1)N_i}$

\qquad has absolute value $|d_{v,i}|_v = |\tc_{v,i}|_v^{n-k-1}$. 

\quad $(c)$ $\Theta_v^{(k)}(z) \in K_v(\cC)$ 
is an $(\fX,\vs)$-function\index{$(\fX,\vs)$-function!$K_v$-rational} 
determined by the local patching process 

\qquad after the coefficients in $\Band_N(k)$
\index{band!coefficients patched by bands!simultaneously when $\Char(K) = p > 0$} 
\index{band!$\Band_N(k)$}  have been modified;  it has  a pole of order 

\qquad at most $(n-k)N_i$ at each $x_i$ and no other poles, and may be the zero function.  

$(4)$ For $k = n$\index{patching functions, $G_v^{(k)}(z)$ for $1 \le k \le n$!constructed by patching}                              
\begin{equation*}
G_v^{(n)}(z) \ = \ G_v^{(n-1)}(z) + \sum_{\lambda=1}^{\Lambda}
                   \Delta_{v,\lambda}^{(n)} \cdot \varphi_{\lambda}(z) \ . 
\end{equation*}

\noindent{$(B)$ For} each $v \in \hS_K^{+}$ 
and each $k = 1, \ldots, n$,\index{patching functions, $G_v^{(k)}(z)$ for $1 \le k \le n$!constructed by patching} 
\begin{equation*} 
\left\{ \begin{array}{l}
\text{If $v \in S_K$, then all the zeros of $G_v^{(k)}(z)$ belong to $E_v$,} \\
 
\text{\quad and for $k = 0$ and $k = n$ they are distinct. When $k = n$,} \\

\text{\qquad \qquad \qquad 
$\{z \in \cC_v(\CC_v) : G_v^{(n)}(z) \in \cO_v \cap D(0,r_v^{Nn}) \} \subset E_v.$} \\ 

\text{If $v \in \hS_K^{+} \backslash S_K$, 
then all the zeros of $G_v^{(k)}(z)$ belong to $E_v$, and } \\
\text{ \qquad \qquad \qquad $\{ z \in \cC_v(\CC_v) : |G_v^{(k)}(z)|_v \le R_v^{Nn} \} = E_v.$ }
\end{array} \right. 
\end{equation*}
\end{theorem}

\noindent{\bf Remark.} 
As when $\Char(K) = 0$, we will have $\Theta_v^{(k)}(z) = 0$ except for 
one value $k = k_1$ for each $v \in S_K$, 
where $\Theta_v^{(k_1)}(z)$ is chosen to `separate the roots'\index{separate the roots} of $G_v^{(k_1)}(z)$.
\index{patching functions, $G_v^{(k)}(z)$ for $1 \le k \le n$!for nonarchimedean $K_v$-simple sets!roots are distinct}  
See the discussion after Theorem \ref{DCPCPatch1p}, 
and Phase 3 in the proof of that theorem. 

\medskip
\index{global patching when $\Char(K) = p > 0$!Stage 3: The global patching construction!comparison with the case $\Char(K) = 0$}
The underlying patching constructions when $\Char(K) = 0$ and $\Char(K) = p > 0$
\index{patching argument!comparing Characteristic $0$ and $p > 0$}
are the same:  we expand\index{patching functions, $G_v^{(k)}(z)$ for $1 \le k \le n$!expansion of} 
\begin{equation*}
G_v^{(k-1)}(z) \ = \ \sum_{i=1}^m \sum_{j=0}^{(n-1)N_i} A_{v,ij} \varphi_{i,nN_i-j}(z) 
            + \sum_{\lambda = 1}^{\Lambda} A_{v,\lambda} \varphi_{\lambda}
\end{equation*}
and we modify the coefficients of $G_v^{(k-1)}(z)$\index{patching functions, $G_v^{(k)}(z)$ for $1 \le k \le n$!constructed by patching} 
\index{coefficients $A_{v,ij}$}\index{compensating functions $\vartheta_{v,ij}^{(k)}(z)$} 
\index{band!coefficients patched by bands!simultaneously when $\Char(K) = p > 0$}  
in $\Band_N(k)$\index{band!$\Band_N(k)$}  by setting
\begin{equation} \label{LModification}
G_v^{(k)}(z)\ = \ G_v^{(k-1)}(z) + \sum_{i=1}^m \sum_{j=(k-1)N_i}^{kN_i-1} \Delta_{v,ij}^{(k)} \vartheta_{v,ij}^{(k)}(z) \ .
\end{equation}
Here $\vartheta_{v,ij}^{(k)}(z)$
\index{compensating functions $\vartheta_{v,ij}^{(k)}(z)$!poles and leading coefficients of} 
has a pole of order $nN_i-j$ at $x_i$ and a pole of order at most $(n-k)N_{i^{\prime}}$ at $x_{i^{\prime}}$ 
for $i^{\prime} \ne i$.  Examining the local patching constructions 
\index{patching argument!local}
for nonarchimedean $v$ shows that 
\begin{equation} \label{varThetaDef}
\vartheta_{v,ij}^{(k)}(z) \ = \ \varphi_{i,r}(z) \cdot F_{v,k}(z) 
\end{equation}\index{compensating functions $\vartheta_{v,ij}^{(k)}(z)$}
where $F_{v,k}(z)$ is $K_v$-rational, with a pole of order $(n-k-1)N_i$ at each $x_i$,
and $nN_i - j = (n-k-1)N_i + r$, so $N_i + 1 \le r = (k+1)N_i - j \le 2N_i$.   
(When  $E_v$ is an $\RL$-domain, $F_{v,k}(z) = \phi_v(z)^{n-k-1}$;  
\index{$\RL$-domain} 
it is more complicated when $E_v$ is $K_v$-simple.)  Since the $\varphi_{ir}$ are $K_v$-symmetric,
\index{$K_v$-simple!set}  
\index{$K_v$-symmetric!set of functions}
the $\vartheta_{v,ij}^{(k)}(z)$\index{compensating functions $\vartheta_{v,ij}^{(k)}(z)$!are $K_v$-symmetric} 
are $K_v$-symmetric.
\index{$K_v$-symmetric!set of functions}
Rewriting (\ref{LModification}) using  (\ref{varThetaDef}) 
gives\index{patching functions, $G_v^{(k)}(z)$ for $1 \le k \le n$!constructed by patching} 
\begin{eqnarray} 
G_v^{(k)}(z) & = & G_v^{(k-1)}(z) + \Big(\sum_{i=1}^m \sum_{j=(k-1)N_i}^{kN_i-1} 
          \Delta_{v,ij}^{(k)} \varphi_{i,(k+1)N_i-j}(z)\Big) \cdot F_{v,k}(z) \notag \\
             & = & G_v^{(k-1)}(z) + \Delta_{v,k}(z) \cdot F_{v,k}(z) \ , \label{PrelimGvkModification} 
\end{eqnarray} 
which is (\ref{GvkOmegaModification}) before the addition of $\Theta_v^{(k)}(z)$.  

However, the patching constructions when $\Char(K) = 0$ and $\Char(K) = p > 0$ have different aims. 
\index{patching argument!comparing Characteristic $0$ and $p > 0$} 
When $\Char(K) = 0$, the patching construction modifies the coefficients $A_{v,ij}$ one by one 
in $\prec_N$ order, making them global numbers $A_{ij} \in L$ \index{order!$\prec_N$} 
which depend on the numbers chosen in earlier patching steps.  
Since each $\varphi_{ij}$ in the $L$-rational basis has a pole of different order,
\index{basis!$L$-rational} 
patching the coefficients of lower degree basis functions does not change coefficients patched earlier,
and since $L_{w_v}/K_v$\index{distinguished place $w_v$} is separable, the modification term 
$\sum_{i=1}^m \sum_{j=(k-1)N_i}^{kN_i-1} \Delta_{v,ij}^{(k)} \vartheta_{v,ij}^{(k)}(z)$  
is $K_v$-rational if the $\Delta_{v,ij}^{(k)}$ are $K_v$-symmetric. 
\index{$K_v$-symmetric!set of numbers}\index{compensating functions $\vartheta_{v,ij}^{(k)}(z)$}
 
When $\Char(K) = p > 0$, since $L_{w_v}/K_v$ may be inseparable,\index{distinguished place $w_v$} 
we cannot simply use galois equivariance of the $\Delta_{v,ij}^{(k)}$ 
to deduce the $K_v$-rationality of the 
modification term.\index{patching functions, $G_v^{(k)}(z)$ for $1 \le k \le n$!are $K_v$-rational}  
Instead, we expand $G_v^{(k-1)}(z)$ and $\Delta_{v,k}(z) F_{v,k}(z)$ 
using the $L^{\sep}$-rational basis, writing 
\index{basis!$L^{\sep}$-rational}\label{`SymbolIndextAvij'} 
\begin{eqnarray}
G_v^{(k-1)}(z) & = & \sum_{i=1}^m \sum_{j=0}^{(n-1)N_i} \tA_{v,ij} \tphi_{i,nN_i-j}(z) 
            + \sum_{\lambda = 1}^{\Lambda} \tA_{v,\lambda} \tphi_{\lambda} \ ,\label{FGv1X} \\
\Delta_{v,k}(z) F_{v,k}(z) & = & 
\sum_{i=1}^m \sum_{j=(k-1)N_i}^{kN_i-1} \tdelta_{v,ij} \tphi_{i,nN_i-j}(z)
\ + \ \text{lower degree terms} \ ,    \label{LSepModification}        
\end{eqnarray} 
and  patch the coefficients of $G_v^{(k-1)}(z)$\index{patching functions, $G_v^{(k)}(z)$ for $1 \le k \le n$!constructed by patching} 
in $\Band_N(k)$\index{band!coefficients patched by bands!simultaneously when $\Char(K) = p > 0$} 
\index{band!$\Band_N(k)$}  
choosing the $\Delta_{v,ij}^{(k)}$ so as to modify the coefficients
\index{coefficients $A_{v,ij}$} 
$\tA_{v,ij}$ relative to the $L^{\sep}$-rational basis.
\index{basis!$L^{\sep}$-rational}   
By Proposition \ref{FPushBound} below, 
$L_{w_v}^{\sep}$-rationality\index{distinguished place $w_v$} and galois equivariance for the $\tdelta_{v,ij}$
implies the $K_v$-rationality of $\Delta_{v,k}(z) F_{v,k}(z)$.

Since several terms of the $L^{\sep}$-rational basis can have poles of the same order,
\index{basis!$L^{\sep}$-rational} 
this requires us to choose all the patching coefficients in a band simultaneously.
\index{band!coefficients patched by bands}\index{patching coefficients}
\index{band!coefficients patched by bands!simultaneously when $\Char(K) = p > 0$}      
By the construction of the $L$-rational and $L^{\sep}$-rational bases in 
\S\ref{Chap3}.\ref{LRationalBasisSection}, 
the transition matrix from the $L$-rational basis to the $L^{\sep}$-rational basis
\index{basis!$L$-rational} 
\index{basis!$L^{\sep}$-rational} 
is block-diagonal, with blocks of size $J$.  
Since we have required that $J|N_i$ in Theorem \ref{LocalPatchp},
we can modify the coefficients $\tA_{v,ij}$ from $\Band_N(k)$\index{band!$\Band_N(k)$}
\index{band!coefficients patched by bands!simultaneously when $\Char(K) = p > 0$}    
by patching all the coefficients $A_{v,ij}$ from $\Band_N(k)$\index{band!$\Band_N(k)$}  at once,
and the patching modifications from later bands do not affect the 
\index{band!coefficients patched by bands}
modifications made in earlier ones. 
     
The patching process must also keep the roots of 
the $G_v^{(k)}(z)$ in $E_v$.\index{patching functions, $G_v^{(k)}(z)$ for $1 \le k \le n$!roots are confined to $E_v$}
Doing so requires an analysis of the relationship between patching modifications of the 
form (\ref{LModification}) and those of the form (\ref{LSepModification}). 
This is carried out in Proposition \ref{FPushBound}. 
We study the growth rates of the coefficients, 
\index{coefficients $A_{v,ij}$!growth rate} 
and show that by using modifications of the form (\ref{PrelimGvkModification}) 
we can independently vary the $L^{\sep}$-rational coefficients of the $G_v^{(k)}(z)$
\index{coefficients $A_{v,ij}$!$L^{\sep}$-rational}\index{patching functions, $G_v^{(k)}(z)$}  
in a given band, in a uniform way independent of $k$. 
\index{band!coefficients patched by bands!simultaneously when $\Char(K) = p > 0$}  
The uniformity ultimately depends on the fact that each $\phi_v(z)$ has its zeros in $E_v$ 
and its poles in $\fX$, 
and $E_v$ is bounded away from $\fX$. 

\smallskip

To motivate the formulation of Proposition \ref{FPushBound}, note that 
since $(k-1)N_i \le j \le kN_i-1$ in (\ref{PrelimGvkModification}), 
if we replace $j$ by $s = j-(k-1)N_i$ and write $\Delta_{v,is} = \Delta_{v,ij}^{(k)}$, 
then (\ref{PrelimGvkModification}) becomes\index{patching functions, $G_v^{(k)}(z)$ for $1 \le k \le n$!constructed by patching} 
\begin{equation} \label{LSepModification2}
G_v^{(k)}(z) \ = \ G_v^{(k-1)}(z) + 
\Big(\sum_{i=1}^m \sum_{s=0}^{N_i-1} \Delta_{v,is} \cdot \varphi_{i,2N_i-s}(z)\Big) F_{v,k}(z) \ .
\end{equation}
We will use Proposition \ref{FPushBound} again in the proof of Theorem \ref{DCPCPatch1p},
so we state it in more generality than is needed for Theorem \ref{LocalPatchp},  
letting $\Char(K)$ be arbitrary and letting $\ell \ge 1$ bands be patched at once. 
\index{band!coefficients patched by bands!simultaneously when $\Char(K) = p > 0$} 
(In the proof of Theorem \ref{LocalPatchp} we will take $\ell = 1$.)
\index{global patching when $\Char(K) = p > 0$!Stage 3: The global patching construction!the patching by Blocks theorem}


\begin{proposition}  \label{FPushBound} 
Let $\Char(K)$ be arbitrary, and let $v$ be a nonarchimedean place of $K$.  
Then there are numbers $\Lambda_v, \tUpsilon_v > 0$, 
depending only on $E_v$, $\fX$, 
the choice of the $L$-rational and $L^{\sep}$-rational bases,
and the projective embedding of $\cC_v$, 
with the following property.  

Let $r > 0$ be small enough that 

\quad $(1)$  $r < \min_{i \ne j} (\|x_i,x_j\|_v);$ 

\quad $(2)$ each of the balls $B(x_i,r)$ is isometrically parametrizable and disjoint from $E_v;$ 
\index{isometrically parametrizable ball}

\quad $(3)$ for each $i$, none of the $\varphi_{ij}(z)$ has a zero in $B(x_i,r).$

\noindent{Put} $\varpi_v = \min(1,\Lambda_v \cdot r)$, 
and let $\ell$, $k$ be integers with $\ell \ge 1$ and $1 \le k \le n-1$.  

Let $\vs \in \cP^m(\QQ)$
be a positive rational probability vector;  
let $N$ be a positive integer such that $N_i = Ns_i \in \NN$ and $J|N_i$, 
for each $i = 1, \ldots, m$.  
Let $F_v(z) \in \CC_v(\cC)$ be an $(\fX,\vs)$-function\index{$(\fX,\vs)$-function} 
which has a pole of order $(n-k-1)N_i$ at $x_i$, for each $i$, 
and whose zeros belong to $E_v$. 
Let $0 \ne d_{v,i} = \lim_{z \rightarrow x_i} F_v(z) \cdot g_{x_i}(z)^{(n-k-1)N_i}$ 
be its leading coefficient at $x_i$.
\index{coefficients $A_{v,ij}$!leading} 

Given $\vDelta = (\Delta_{v,is})_{1 \le i \le m, 0 \le s < \ell N_i} \in \CC_v^{\ell N}$, 
let $\Delta_v(z) \in \CC_v(\cC)$ be the $(\fX,\vs)$-function\index{$(\fX,\vs)$-function}
\begin{equation} \label{Fh_vF} 
\Delta_v(z) \ = \ \sum_{i=1}^m \sum_{s=0}^{\ell N_i-1} \Delta_{v,is} \cdot \varphi_{i,(\ell+1)N_i-s}(z) \ .
\end{equation}  
Expand $\Delta_v(z)F_v(z)$ in terms of the $L^{\sep}$-rational basis as 
\index{basis!$L^{\sep}$-rational} 
\begin{equation*}
\Delta_v(z) F_v(z)
\ = \ \sum_{i=1}^m \sum_{s=0}^{\ell N_i-1}  \tdelta_{v,is} \cdot \tphi_{i,(k+\ell)N_i - s}(z) 
                  \ + \ \text{lower order terms \ .} \notag                   
\end{equation*}
and write $\tdelta = (\tdelta_{v,is})_{1 \le i \le m, 0 \le s < \ell N_i} \in \CC_v^{\ell N}$.
Let $\Phi^{\sep}_{F_v} :\CC_v^{\ell N} \rightarrow \CC_v^{\ell N}$
be the linear map defined by $\Phi^{\sep}_{F_v}(\vDelta) = \tdelta$.  

Then $\Phi^{\sep}_{F_v}$
is an isomorphism and for each $\rho > 0$  
\begin{equation}
\Phi_{F_v}^{\sep}\big(\bigoplus_{i=1}^m D(0,\rho)^{\ell N_i}\big) \ \supseteq \
\bigoplus_{i=1}^m D(0,\tUpsilon_v \varpi_v^{\ell N} |d_{v,i}|_v \rho)^{\ell N_i}  \label{FBallC2Q}
\end{equation} 
\noindent{and} 
\begin{equation} 
\Phi^{\sep}_{F_v}\big(\bigoplus_{i=1}^m \bigoplus_{s=0}^{\ell N_i-1} D(0,\varpi_v^{-s} \rho) \big) \ \supseteq \ 
\bigoplus_{i=1}^m \bigoplus_{s=0}^{\ell N_i-1} D(0, \tUpsilon_v \varpi_v^{-s} |d_{v,i}|_v \cdot \rho ) \ .
\label{FBallC2} 
\end{equation}

Moreover, if $F_v(z)$ is $K_v$-rational,  
then for each $K_v$-symmetric $\tdelta \in \big(L_{w_v}^{\sep}\big)^{\ell N}$,\index{distinguished place $w_v$}
\index{$K_v$-symmetric!vector} 
the unique solution to $\Phi^{\sep}_{F_v}(\vDelta) = \tdelta$ belongs to $L_{w_v}^{\ell N}$
and is $K_v$-symmetric, and the corresponding function 
\index{$K_v$-symmetric!vector}
$\Delta_v(z) = \sum_{i=1}^m \sum_{s=0}^{\ell N_i-1} \Delta_{v,is} \cdot \varphi_{i,(\ell+1)N_i-s}(z)$  
is $K_v$-rational. 
\end{proposition}

The proof of Proposition \ref{FPushBound} 
will be given in \S\ref{Chap7}.\ref{FPushBoundProofSec} below.
We now choose the patching parameters $\kbar$ and $n$ in Theorem \ref{LocalPatchp}, 
\index{patching parameters}
and give the details of the global patching process.\index{patching argument!global|ii} 

\medskip
{\bf The choice of $\kbar$.}
\index{global patching when $\Char(K) = p > 0$!Stage 3: The global patching construction!the choice of $\kbar$}
In Stage 1 we have chosen a collection of numbers $h_v$ 
for $v \in \hS_K^{+}$ such that $\prod_{v \in \hS_K^{+}} h_v > 1$.  
Likewise, for each $v \in \hS_K^{+}$, Proposition \ref{LocalPatch} provides a number $k_v$  
(the minimal number of stages considered high-order by the patching process at $v$).
Finally, for each $v \in \hS_K^{+}$, fix a number $r = r_v > 0$ satisfying the conditions of 
Proposition \ref{FPushBound} and small enough that each ball $B(x_i,r_v)$ is disjoint from $U_v$;  
then Proposition \ref{FPushBound} provides numbers $\tUpsilon_v > 0$ and $0 < \varpi_v \le 1$  
(the comparison constants for the transition between the $L$-rational and $L^{\sep}$-rational bases). 
Put $\tUpsilon = \prod_{v \in \hS_K^{+}} (\tUpsilon_v \varpi_v^N) $ 
and $h = \prod_{v \in \hS_K^{+}} h_v > 1$. 
Finally, put $H = L^{\sep}$ and let $C_H(\hS_K^{+})$ be the constant from Proposition \ref{CStrong3}.
\index{Strong Approximation theorem!uniform} 
 
Let $\kbar$ be the smallest integer such that 
\begin{equation} \label{Fkbarchoicep}
\left\{ \begin{array}{l} 
         \text{$\kbar \ \ge \ k_v$ \quad for each $v \in \hS_K^{+}$\ ,} \\
         \text{$ (\tUpsilon \cdot h^{\kbar N})^{[H:K]} \ > \ C_H(\hS_K^{+})$\ .} 
        \end{array} \right.
\end{equation}

\medskip
{\bf The choice of $n$.}
\index{global patching when $\Char(K) = p > 0$!Stage 3: The global patching construction!the choice of $n$}
\index{patching parameters|ii}
As in the patching construction when $\Char(K) = 0$, 
for suitable $n$ we will take the initial patching functions to be
\index{patching functions, initial $G_v^{(0)}(z)$|ii}
$G_v^{(0)}(z) = Q_{v,n}(\phi_v(z))$, where $Q_{v,n}(x) \in \cO_v(z)$ is the monic polynomial
of degree $n$ given by Theorem \ref{LocalPatchp}.  

One consideration in choosing $n$ is to facilitate 
patching the leading coefficients of the $G_v^{(0)}(z)$ 
to be $\hS_L^{+}$-units.\index{patching functions, initial $G_v^{(0)}(z)$!leading coefficients of}
\index{patching!leading coefficients}
Given $v \in \hS_K$, put $\phi_w(z) = \phi_v(z)$  
for each $w \in \hS_L$ with $w|v$, viewing the $\phi_w(z)$ 
as functions in $L_w(\cC)$.   
By Theorem \ref{CTCX2p} the leading coefficients $\tc_{w,i}$ of the $\phi_w(z)$ 
\index{coefficients $A_{v,ij}$!leading} 
have the property that there are an integer $n_0$, 
and a $K$-symmetric system of $\hS_L$-units $\mu_i$, such that for each $i$ 
\index{$K$-symmetric!system of units}
and each $w \in \hS_L^{+}$  
\begin{equation} \label{FLeadp}
|\tc_{w,i}^{n_0}|_w \ = \ |\mu_i|_w \ . 
\end{equation}
For each $v \in \hS_K^{+}$, let $0 < B_v < 1$ be the number 
from Theorem \ref{LocalPatchp} controlling the freedom in patching the leading coefficients.
\index{coefficients $A_{v,ij}$!leading} 
All the fields $L_w$ for $w|v$ are isomorphic, 
and by the structure of the group of units $O_w^{\times}$ there is an integer 
$n_v^{\prime} > 0$ such that for each  $x \in \cO_w^{\times}$, and
each integer $n^{\prime}$ divisible by $n_v^{\prime}$, 
\begin{equation} \label{FCloseUpp} 
|x^{n^{\prime}} - 1|_v  \ \le \ B_v \ .
\end{equation}
Let $n_1$ be the least common multiple of the  $n_v^{\prime}$.  

Another consideration in choosing $n$ is to assure that various analytic estimates
are satisfied.  For each $v \in \hS_K^{+}$, Theorem \ref{LocalPatchp} provides a number $n_v$
such that the local patching process at $v$ \index{patching argument!local}
will keep the roots of $G_v^{(0)}(z)$ 
in $E_v$,\index{patching functions, initial $G_v^{(0)}(z)$!roots are confined to $E_v$} 
provided $n_v|n$ and $n$ is sufficiently large.  
Let $n_2$ be the least common multiple of the $n_v$ for $v \in \hS_K^{+}$. 

Finally, let $n$ be a positive integer such that  
\begin{equation} \label{FnChoicep}
 n_0 n_1 n_2 | n 
\end{equation}
and is large enough that 
\begin{equation} \label{Fnbigenoughp}
(\prod_{v \in \hS_H^{+}} h_v)^{nN [H:K]} \ > \ C_H(\hS_K^{+}) \ .
\end{equation} 
By Theorem \ref{LocalPatchp} there is an $n_3$ such that if $n \ge n_3$,
then (\ref{Fnbigenoughp}) holds and for each $v \in \hS_K^{+}$ the local patching process can 
be successfully completed.\index{patching argument!local}

Until last step in the proof, $n \ge n_3$ will be a fixed integer 
satisfying (\ref{FnChoicep}).  

\medskip
{\bf The order $\prec_N$.}
\index{global patching when $\Char(K) = p > 0$!Stage 3: The global patching construction!the order $\prec_N$}  
The index set $\cI = \{(i,j) \in \ZZ^2 : 1 \le i \le m, 0 \le j < \infty\}$,
the order $\prec_N$, the bands \index{order!$\prec_N$|ii} 
\index{band!coefficients patched by bands}\index{band!$\Band_N(k)$|ii} 
$\Band_N(k) = \{(i,j) \in \cI : 1 \le i \le m, (k-1)N_i \le j \le kN_i-1\}$, 
and the galois blocks\index{block, galois!$\Block(i,j)$|ii} 
$\Block((i_0,j_0)) = \{ (i,j_0) : \text{ $\sigma(x_{i_0}) = x_i$ for some $\sigma \in \Aut(L/K)$} \}$
will be the same as those when $\Char(K) = 0$ (see (\ref{FPrec})).  
%
%

\medskip
{\bf Patching the Leading Coefficients.}
\index{patching!leading coefficients}
\index{global patching when $\Char(K) = p > 0$!Stage 3: The global patching construction!patching the leading coefficients}
For each $v \in \hS_K^{+}$, the initial patching function 
\index{patching functions, initial $G_v^{(0)}(z)$!construction of}
$G_v^{(0)}(z) = Q_{v,n}(\phi_v(z)) \in K_v(\cC)$ 
is an $(\fX,\vs)$-function\index{$(\fX,\vs)$-function} of degree $nN$,
with a pole of order $nN_i$ and leading coefficient $\tc_{v,i}^n$ at each $x_i$.
\index{coefficients $A_{v,ij}$!leading}   
Let $\mu_1, \ldots, \mu_m \in L^{\sep}$ and $n_0 \ge 1$ 
be the $K$-symmetric set of $\hS_L^{+}$-units
\index{$K$-symmetric!system of units}
and integer constructed in Stage 2 above.
\index{patching functions, initial $G_v^{(0)}(z)$!making the leading coefficients $S$-units}   

Fix $v \in \hS_K^{+}$, and view $\mu_1, \ldots, \mu_m$ as embedded in $L_{w_v}^{\sep}$;\index{distinguished place $w_v$}   
thus $|\tc_{v,i}^{n_0}|_v = |\mu_i|_v$.  For each $i = 1, \ldots, m$ put 
\begin{equation*} 
\Delta_{v,i0}^{(1)} \ = \ (\mu_i^{n/n_0}/\tc_{v,i}^n) - 1 
              \ = \ (\mu_i/\tc_{v,i}^{n_0})^{n/n_0} - 1 \ .
\end{equation*} 
Since $n_0 n_1|n$, (\ref{FCloseUpp}) shows that 
\begin{equation} \label{FHighBound}
|\Delta_{v,i0}^{(1)}|_v \ \le \ B_v \ .
\end{equation}  
Since $\tc_{v,i} \in L_{w_v}^{\sep}$,\index{distinguished place $w_v$} we have $\Delta_{v,i0}^{(1)} \in L_{w_v}^{\sep}$.
Since the $\mu_i$ and $\tc_{v,i}$ are $K_v$-symmetric, the $\Delta_{v,i0}^{(1)}$
\index{$K_v$-symmetric!set of numbers}
are $K_v$-symmetric as well.  We will take the $\Delta_{v,ij}^{(1)}$  
\index{$K_v$-symmetric!set of numbers}
for $j = 1, \ldots, N_i-1$, to be $0$.  

By Theorem \ref{LocalPatchp}, 
when $G_v^{(0)}(z)$ is expanded using the $L^{\sep}$-rational basis 
as\index{patching functions, initial $G_v^{(0)}(z)$!expansion of}
\index{basis!$L^{\sep}$-rational} 
\begin{equation*}
G_v^{(0)}(z) \ = \ \sum_{i=1}^m \sum_{j=1}^{(n-1)N_i} \tA_{v,ij} \tphi_{ij}(z) 
         + \sum_{\lambda=1}^{\Lambda} \tA_{v,\lambda} \tphi_{\lambda} \ ,
\end{equation*}
then for each $i$ we have $\tA_{v,i0} = \tc_{v,i}^n$ and $A_{v,ij} = 0$
for $j = 1, \ldots, \kbar N_i -1$.  
The local patching construction at $v$ provides a $K_v$-symmetric set of functions 
\index{$K_v$-symmetric!set of functions}\index{patching argument!local}
$\ttheta_{v,i0}^{(1)}(z) \in L_{w_v}^{\sep}(\cC)$\index{distinguished place $w_v$} for  $i =1, \ldots, m$, such that for each $i$
there is an $(\fX,\vs)$-function\index{$(\fX,\vs)$-function} 
$\tTheta_{v,i}(z) \in L_{w_v}^{\sep}(\cC)$\index{distinguished place $w_v$} 
with a pole of order at most $(n-\kbar)N_{i^{\prime}}$
for each $i^{\prime}$, for which  
\begin{equation*} 
\ttheta_{v,i0}^{(1)}(z) \ = \ \tc_{v,i}^n \tphi_{i,nN_i}(z) + \tTheta_{v,i}(z) \ .
\end{equation*}  
(See (\ref{FphiQm1p1}) and (\ref{FphiQm1p1Q}); 
note that $\tphi_{i,nN_i} = \varphi_{i,nN_i}$ since $J|N_i$.) 
Thus if we put\index{patching functions, $G_v^{(k)}(z)$ for $1 \le k \le n$!constructed by patching} 
\begin{equation*}
G_v^{(1)}(z) \ = \ G_v^{(0)}(z) + \sum_{i=1}^m \Delta_{v,i0}^{(1)} \cdot \ttheta_{v,i0}^{(1)}(z) \ ,
\end{equation*} 
then the leading coefficient of $G_v^{(1)}(z)$ at $x_i$ 
becomes\index{patching functions, $G_v^{(k)}(z)$} 
\index{coefficients $A_{v,ij}$!leading} 
\begin{equation*}
\tA_{v,i0} + \Delta_{v,i0}^{(1)} \cdot \tc_{v,i}^n \ = \ 
\tc_{v,i}^n + \big((\mu_i^{n/n_0}/\tc_{v,i}^n)-1) \cdot \tc_{v,i}^n \ = \ \mu_i^{n/n_0} 
\end{equation*} 
while the coefficients $\tA_{v,ij}$ for $j = 1, \ldots, \kbar N_i -1$ remain $0$. 
Since the $\Delta_{v,i0}^{(1)}$ and $\vartheta_{v,i0}^{(1)}(z)$ are are $K_v$-symmetric
\index{$K_v$-symmetric!set of functions}
and defined over $L_{w_v}^{\sep}$,\index{distinguished place $w_v$} $G_v^{(1)}(z)$ belongs to $K_v(\cC)$.\index{patching functions, $G_v^{(k)}(z)$}   

By construction all the zeros of $G_v^{(0)}(z)$\index{patching functions, $G_v^{(k)}(z)$ for $1 \le k \le n$!roots are confined to $E_v$} 
belong to $E_v$, and since (\ref{FHighBound}) 
holds, the local patching process assures that the 
zeros of $G_v^{(1)}(z)$\index{patching functions, $G_v^{(k)}(z)$ for $1 \le k \le n$!roots are confined to $E_v$} belong to $E_v$ as well.
\index{patching argument!local} 

\smallskip
{\bf Patching the High Order Coefficients.}  
\index{coefficients $A_{v,ij}$!high-order}\index{patching!high-order coefficients}
\index{global patching when $\Char(K) = p > 0$!Stage 3: The global patching construction!patching the high-order coefficients} 
No patching is needed for $k = 2, \ldots, \kbar$.
Since the coefficients $\tA_{v,ij} = 0$ 
for $j=1, \ldots, \kbar N_i - 1$, they are already independent of $v$ and belong to $L^{\sep}$.
Hence for $k = 2, \ldots, \kbar$ we can take $\Delta_{v,ij}^{(k)} = 0$ for all $i,j$
and set $G_v^{(k)}(z) = G_v^{(k-1)}(z) = G_v^{(1)}(z)$.\index{patching functions, $G_v^{(k)}(z)$ for $1 \le k \le n$!constructed by patching} 

\smallskip
{\bf Patching the Middle Coefficients.}  For each $k = \kbar+1, \ldots, n-1$ we first choose the target 
\index{coefficients $A_{v,ij}$!middle}\index{patching!middle coefficients}
\index{global patching when $\Char(K) = p > 0$!Stage 3: The global patching construction!patching the middle coefficients}
coefficients using  Proposition \ref{CPPr1p} and a see-saw argument.\index{see-saw argument}
We then choose the patching coefficients $\Delta_{v,ij}^{(k)}$ using Proposition \ref{FPushBound},
and patch using Theorem \ref{LocalPatchp}.

As before, let $H = L^{\sep}$.  Then $L/H$ is purely inseparable and $H/K$ is galois,
with $\Aut(L/K) \cong \Gal(H/K)$.  
For each place $w_0$ of $H$ there is a unique place $w$ of $L$ lying over $w_0$, 
and $w/w_0$ is totally ramified.  In the discussion below, we will work primarily with $H$, 
and to simplify notation we will write $w$ both for places of $L$ and $H$.  
Let $\hS_H^{+}$ be the set of places $w$ of $H$ over places $v \in \hS_K^{+}$.

Suppose that for some $k > \kbar$, we have completed the patching process through stage $k-1$,
and have constructed functions 
$G_v^{(k-1)}(z) \in K_v(\cC)$\index{patching functions, $G_v^{(k)}(z)$ for $1 \le k \le n$!constructed by patching} 
with the properties in Theorem \ref{LocalPatchp}.
For each $v \in \hS_K^{+}$, and each $w \in \hS_H^{+}$ with $w|v$, put $G_w^{(k-1)}(z) = G_v^{(k-1)}(z)$
\index{patching functions, $G_v^{(k)}(z)$ for $1 \le k \le n$!viewed simultaneously over $K_v$ and $L_w$} 
and expand $G_w^{(k-1)}(z)$ using the $L^{\sep}$-rational basis as 
\index{basis!$L^{\sep}$-rational} 
\begin{equation*}
G_w^{(k-1)}(z) \ = \
    \sum_{i=1}^m \sum_{j=0}^{(n-1)N_i-1} \tA_{w,ij} \tphi_{i,nN_i-j}(z)
       + \sum_{\lambda = 1}^{\Lambda} \tA_{w,\lambda} \tphi_{\lambda} \ . 
\end{equation*}
Since $G_w^{(k-1)}(z)$ is rational over $K_v$ and the $\tphi_{i,nN_i-j}$ and $\tphi_{\lambda}$ 
are rational over $H = L^{\sep}$, the $\tA_{w,ij}$ and $\tA_{w,\lambda}$ belong to $H_w$ 
and are $K_v$-symmetric.
\index{$K_v$-symmetric!set of functions} 
Since $\tphi_{i,nN_i-j}$ is rational over $K(x_i)^{\sep}$,  
by galois equivariance $\tA_{w,ij}$ in fact belongs to $K_v(x_i)^{\sep} \subset H_w$.  
 
\smallskip
We will now choose the target coefficients $A_{ij} \in K(x_i)^{\sep}$. 
\index{coefficients $A_{v,ij}$!target}   
Let $(i_0,j_0) \in \Band_N(k)$ is the least index under $\prec_N$ for which the target coefficient
\index{order!$\prec_N$}\index{band!$\Band_N(k)$} 
has not been chosen.  Because of the way the $L^{\sep}$-rational basis
\index{basis!$L^{\sep}$-rational} 
was constructed,  $\Block_N((i_0,j_0)) \subset \Band_N(k)$
\index{block, galois!$\Block(i,j)$}\index{band!$\Band_N(k)$} 
is the set of indices $(i,j_0) \in \cI$ for which there is 
a $\sigma \in \Gal(L^{\sep}/K)$ such that $\sigma(\tphi_{i_0,nN_{i_0}-j_0}) = \tphi_{i,nN_i-j_0}$.
We first determine the target coefficient $A_{i_0,j_0} \in K(x_{i_0})^{\sep}$ and then
define the target coefficients for the other $(i,j)$ in $\Block((i_0,j_0))$\index{block, galois!$\Block(i,j)$} 
so as to preserve galois equivariance.

Consider the vector
\begin{equation*}
\vec{A}_{H,i_0 j_0} \ := \ \oplus_{w \in \hS_H^{+}} \tA_{w,i_0 j_0} 
                  \ \in \ \bigoplus_{w \in \hS_H^{+}} H_w  
\end{equation*}
and let $F = K(x_{i_0})^{\sep}$.  
For each $v \in \hS_K$, the functions $G_w^{(k-1)}(z) \in K_v(\cC)$ 
are the same for all places $w$ of $H$ with $w|v$, 
and the $\tA_{w,i_0,j_0}$ belong to $H_w^{\sep}$, 
so Proposition \ref{CPPr1p} tells us that $\vec{A}_{H,i_0,j_0}$ 
belongs to $\oplus_{u \in \hS_F} F_u$, embedded semi-diagonally in 
$\oplus_{w \in \hS_H} H_w$.

\smallskip
Since $k > \kbar$, our choice of $\kbar$ (see \ref{Fkbarchoicep}) assures that  
\begin{equation*} 
\prod_{v \in \hS_K^{+}} (\tUpsilon_v \varpi_v^N \cdot h_v^{k N})^{[H:K]} \ > \ C_H(\hS_K^{+}) \ .
\end{equation*} 
Recalling that for each $w \in \hS_H^{+}$ we have $|x|_w = |x|_v^{[H_w:K_v]}$,  
put $\tUpsilon_w = \tUpsilon_v^{[H_w:K_v]}$, $b_w = \varpi_v^{[H_w:K_v]}$, and  $h_w = h_v^{[H_w:K_v]}$. 
Since $\sum_{w|v} [H_w:K_v] = [H:K]$ for each $v$, it follows that 
\begin{equation} \label{HLevelProd}
\prod_{w \in \hS_H^{+}} (\tUpsilon_w b_w^N \cdot h_w^{k N}) \ > \ C_H(\hS_K^{+})  \ .
\end{equation}  

For each each $w \in \hS_H^{+}$, put
\begin{equation}  
Q_w \ = \ \tUpsilon_w b_w^N \cdot h_w^{k N} \cdot |\tc_{w,i_0}^{n-k-1}|_w \ ,  \label{CFGH1p}
\end{equation}
where $\tc_{w,i_0}$ is the leading coefficient of $\phi_w(z)$ at $x_{i_0}$.
By (\ref{LeadingCoherenceCond}) we have 
$\prod_{w \in \hS_H^{+}} |\tc_{w,i_0}|_w = 1$, so  (\ref{HLevelProd}) gives    
\begin{equation} \label{QBigEnuf}
\prod_{w \in \hS_H^{+}} Q_w \ > \ C_H(S_H^{+}) \ . 
\end{equation}
Note that $\tUpsilon_w$, $b_w$ and $h_w$ depend only on the place $v$ of $K$ below $w$,
while $|\tc_{w,i_0}|_w$ depends only on the place $u$ of $F$ below $w$,
since the $\phi_w(z) \in K_v(\cC)$ are the same for all $w|v$.
Hence $Q_w$ depends only on the place $u$ below $w$.  
Similarly the coefficients $\tA_{w,i_0 j_0}$ with $w|u$ belong to $F_u$ 
and depend only on $u$.  

By (\ref{QBigEnuf}) we can apply Proposition \ref{CStrong3} to the elements
\index{Strong Approximation theorem!uniform} 
$c_u = \tA_{w,i_0 j_0} \in F_u$, and to the $Q_w$.
(Note that in the function field case, there are no archimedean places, 
so the exponents $D_w$ in Proposition \ref{CStrong3} are all $1$.)
\index{Strong Approximation theorem!uniform} 
By Proposition \ref{CStrong3}, there is an $\tA_{i_0 j_0} \in K(x_{i_0})^{\sep}$ such that
\begin{equation*}
\left\{ \begin{array}{cl}
     |\tA_{i_0 j_0} - \tA_{w,i_0 j_0}|_w \le Q_w & \text{for each $w \in \hS_H^{+}$\ ,} \\
     |\tA_{i_0 j_0}|_w \le 1         & \text{for each $w \notin \hS_H^{+}$\ .}
        \end{array} \right.
\end{equation*}
This $\tA_{i_0 j_0}$ will be the target in patching the $\tA_{w,i_0 j_0}$.
For each $(i,j_0) \in \Block((i_0,j_0))$, take $\sigma \in \Gal(H/K)$
with $\sigma(\tphi_{i_0,nN_{i_0}-j_0}) = \tphi_{i,nN_i-j_0}$, and put $\tA_{i j} = \sigma(\tA_{i_0 j_0})$.
Since $\tA_{i_0,j_0} \in K(x_{i_0})^{\sep}$, the $\tA_{i j}$ are well-defined    
and galois equivariant.

Repeat this process until target coefficients $\tA_{ij}$ have been chosen for all $(i,j) \in \Block_N(k)$.
Since the $\tA_{ij}$ belong to $L^{\sep}$ and are $K$-symmetric, the function
\index{$K$-symmetric!set of numbers}
\begin{equation*}
H^{(k)}(z) \ := \ \sum_{(i,j) \in \Band_N(k)} \tA_{ij} \tphi_{ij}(z)
\end{equation*}
is $K$-rational.   

\medskip
We next choose the patching coefficients $\Delta_{v,ij}^{(k)}$ 
\index{patching coefficients}
so as to replace the part of the $L^{\sep}$-rational 
expansions of the $G_v^{(k-1)}(z)$\index{patching functions, $G_v^{(k)}(z)$ for $1 \le k \le n$!expansion of} 
coming from $\Band_N(k)$\index{band!$\Band_N(k)$}   with $H^{(k)}(z)$. 

Fix $v \in \hS_K^{+}$, and let $F_{v,k}(z) \in K_v(\cC)$ 
be the $(\fX,\vs)$-function\index{$(\fX,\vs)$-function} of degree $(n-k-1)N$
provided by Theorem \ref{LocalPatchp}.  View $\fX$ and the $\tA_{ij}$ as embedded in $L_{w_v}$,\index{distinguished place $w_v$}  
and let $d_{v,i} \in L_{w_v}$\index{distinguished place $w_v$} 
be the leading coefficient of $F_{v,k}(z)$ at $x_i$.  By hypothesis 
$|d_{v,i}|_v = |\tc_{v,i}^{n-k-1}|_v$.   

For each $(i,j) \in \Band_N(k)$, put 
$\tdelta_{v,ij}^{(k)} = \tA_{ij} - \tA_{w_v,ij} \in K_v(x_i)^{\sep} \subseteq H_{w_v}$,\index{distinguished place $w_v$} 
and let $Q_v = Q_{w_v}^{1/[H_{w_v}:K_v]} = \tUpsilon_v \varpi_v^N  h_v^{kN} \cdot |\tc_{v,i}^{n-k-1}|_v$. 
Since $|\tdelta_{v,ij}^{(k)} |_{w_v} = |\tA_{ij} - \tA_{w_v,ij}|_{w_v} \le Q_{w_v}$, we have 
\begin{equation}  \label{FGas1p}
|\tdelta_{v,ij}^{(k)}|_v \ \le \ Q_v \ = \ \tUpsilon_v \varpi_v^N  |d_{v,i}|_v \cdot h_v^{kN} \ ,  
\end{equation}
and 
\begin{equation} \label{FGas2p}  
    \tA_{ij} \ = \ \tA_{v,ij} + \tdelta_{v,ij}^{(k)} \ .
\end{equation}

The $\tA_{v,ij}$  are $K_v$-symmetric since 
\index{$K_v$-symmetric!set of numbers}\index{patching functions, $G_v^{(k)}(z)$ for $1 \le k \le n$!are $K_v$-rational}
$G_v^{(k)}(z)$ and $F_{v,k}(z)$ are $K_v$-rational, and the $\tA_{ij}$
are $K_v$-symmetric since they are $K$-symmetric.  Hence the $\tdelta_{v,ij}^{(k)}$ are $K_v$-symmetric.
\index{$K$-symmetric!set of numbers}
\index{$K_v$-symmetric!set of numbers}
Put\index{distinguished place $w_v$} 
\begin{equation*}
\tdelta_v^{(k)} \ = \ (\tdelta_{v,ij}^{(k)})_{(i,j) \in \Band_N(k)} \ \in \ (L_{w_v}^{\sep})^N  
\end{equation*}
and let $\Phi_{F_{v,k}}^{\sep} : \CC_v^N \rightarrow \CC_v^N$ 
be the map from Proposition \ref{FPushBound}.  
By Proposition \ref{FPushBound} there a unique 
$\vDelta_v^{(k)} = (\Delta_{v,ij}^{(k)})_{(i,j) \in \Band_N(k)} \in \CC_v^N$ such that  
$\Phi_{F_{v,k}}^{\sep}(\vDelta_v^{(k)}) = \tdelta_v^{(k)}$.  Using (\ref{FGas1p}), 
and applying (\ref{FBallC2Q}) of Proposition \ref{FPushBound} with $\rho = h_v^{kN}$ and $\ell = 1$, 
we see that  
\begin{equation} \label{FGas3p} 
|\Delta_{v,ij}^{(k)}|_v \ \le \ h_v^{kN} 
\end{equation} 
for all $(i,j)$.  Furthermore, Proposition \ref{FPushBound} tells us that  
\begin{equation*}
\Delta_{v,k}(z) \ := \ \sum_{(i,j) \in \Band_N(k)} \Delta_{v,ij}^{(k)} \varphi_{ij}(z)
\end{equation*} 
is $K_v$-rational.  By the definition of the map $\Phi_{F_{v,k}}^{\sep}$, 
\begin{equation} \label{FGas4p}
\Delta_{v,k}(z) \cdot F_{v,k}(z) \ = \ 
\sum_{(i,j) \in \Band_N(k)} \tdelta_{v,ij}^{(k)} \cdot \tphi_{ij}(z) 
\ + \ \text{terms of lower order} \ .  
\end{equation} 

\smallskip
If we patch $G_v^{(k-1)}(z)$  
using the $\Delta_{v,ij}^{(k)}$,\index{patching functions, $G_v^{(k)}(z)$ for $1 \le k \le n$!constructed by patching} 
then Theorem \ref{LocalPatchp} provides a $K_v$-rational 
$(\fX,\vs)$-function\index{$(\fX,\vs)$-function} $\tTheta_v^{(k)}(z)$  with a pole of 
order at most $(n-k)N_i$ at each $x_i$, such that\index{patching functions, $G_v^{(k)}(z)$ for $1 \le k \le n$!constructed by patching} 
\begin{equation*}
G_v^{(k)}(z)  
\ = \ G_v^{(k-1)} + \Delta_{v,k}(z) \cdot F_{v,k}(z) + \tTheta_v^{(k)}(z) \ .
\end{equation*} 
By (\ref{FGas2p}) and (\ref{FGas4p}), for each $(i,j) \in \Band_N(k)$ 
the coefficient of $\tphi_{ij}$ in the $L^{\sep}$-rational 
expansion of $G_v^{(k)}(z)$\index{patching functions, $G_v^{(k)}(z)$ for $1 \le k \le n$!constructed by patching}
becomes $\tA_{ij}$.  Since (\ref{FGas3p}) holds, the roots of $G_v^{(k)}(z)$ belong to $E_v$.  
Finally, since $G_v^{(k)}(z)$, $\Delta_{v,k}(z) F_{v,k}(z)$,
\index{patching functions, $G_v^{(k)}(z)$ for $1 \le k \le n$!roots are confined to $E_v$}
and $\tTheta_v^{(k)}(z)$ are $K_v$-rational, so is $G_v^{(k)}(z)$.
\index{patching functions, $G_v^{(k)}(z)$ for $1 \le k \le n$!are $K_v$-rational}   

\medskip
{\bf Patching the Low Order Coefficients.}
\index{coefficients $A_{v,ij}$!low-order}\index{patching!low-order coefficients}
\index{global patching when $\Char(K) = p > 0$!Stage 3: The global patching construction!patching the low-order coefficients} 
The final stage of the global patching process deals with the coefficients
$\tA_{v,\lambda}$ in the expansions\index{patching functions, $G_v^{(k)}(z)$ for $1 \le k \le n$!expansion of}
\begin{equation*}
G_v^{(n-1)}(z) \ = \
\sum_{i=1}^m \sum_{j=0}^{(n-1)N_i} \tA_{v,ij} \tphi_{i,nN_i-j}(z)
   + \sum_{\lambda = 1}^{\Lambda} \tA_{v,\lambda} \tphi_{\lambda} \ .
\end{equation*}

\smallskip
For each $v \in \hS_K^{+}$, all the $\tA_{v,\lambda}$ will be patched simultaneously.  
As before, we use a see-saw argument.\index{see-saw argument}   Let $H = L^{\sep}$.  For each $w$ of $H$ with $w|v$, 
put $G_w^{(n-1)}(z) = G_v^{(n-1)}(z)$ 
and expand\index{patching functions, $G_v^{(k)}(z)$ for $1 \le k \le n$!viewed simultaneously over $K_v$ and $L_w$} 
\begin{equation*}
G_w^{(n-1)}(z) \ = \
\sum_{i=1}^m \sum_{j=0}^{(n-1)N_i} A_{w,ij} \tphi_{i,nN_i-j}(z)
   + \sum_{\lambda = 1}^{\Lambda} A_{w,\lambda} \tphi_{\lambda} \ .
\end{equation*}

By construction, for each $v \in \hS_K^{+}$, 
the vector $\oplus_{w|v} G_w^{(n-1)}(z) \in H \otimes_K K_v(\cC)$
is $\Gal(H/K)$ invariant.  By Proposition \ref{CPPr1} this means that 
for each $\lambda$, the coefficient vector $\oplus_{w|v} \tA_{w,\lambda}$ 
has the same galois-equivariance properties as $\varphi_{\lambda}(z)$.  
In particular, if $K \subset F_{\lambda} \subset L$ is the smallest 
field of rationality for $\varphi_{\lambda}(z)$, 
then $\oplus_{w|v} A_{w,\lambda} \in F_{\lambda} \otimes_K K_v$.  

By (\ref{Fnbigenoughp}) in our choice of $n$, we have  
\begin{equation*}
(\prod_{w \in \hS_L} h_w)^{nN} \ > \ C_L(\hS_K^{+}) \ ,
\end{equation*}
so taking $Q_w = h_w^{nN}$ in Proposition \ref{CStrong3} we can find
\index{Strong Approximation theorem!uniform} 
an $\tA_{\lambda} \in F_{\lambda}$ such that 
\begin{equation*}
\left\{ \begin{array}{ll}
         |\tA_{\lambda} - \tA_{w,\lambda}|_w \le h_w^{nN} &
                                \text{for all $w \in \hS_L$\ ,} \\
         |\tA_{\lambda}|_w \le 1 & \text{for all $w \notin \hS_L$ \ .}
        \end{array} \right.
\end{equation*}
By working with representatives of galois orbits as before,
we can arrange that for each $\sigma \in \Gal(H/K)$ we have
$\sigma(\tA_{\lambda}) = \tA_{\lambda^{\prime}}$
if $\sigma(\varphi_{\lambda}) = \varphi_{\lambda^{\prime}}$.  

Put 
\begin{equation*}
\tDelta_{w,\lambda} \ = \ \tA_{\lambda} - \tA_{w,\lambda} 
\end{equation*}                                                 
for each $w$ and $\lambda$, and put 
\begin{equation*}
\Delta_w^{(n)}(z) \ = \ \sum_{\lambda} \tDelta_{w,\lambda} \tphi_{\lambda}(z) \ .
\end{equation*}
Then $\oplus_{w|v} \Delta_w(z) \in L \otimes_K K_v(\cC)$ 
is stable under $\Gal(L/K)$, for each $v \in \hS_K$.  It follows that the
$\Delta_w(z)$ belong to $K_v(\cC)$ and are the same for all $w|v$.  
Put $\Delta_{v,n}(z) = H_{w_v}^{(n)}(z)$,\index{distinguished place $w_v$} and expand 
$\Delta_{v,n}(z) = \sum_{\lambda=1}^{\Lambda} \Delta_{v,\lambda} \tphi_{\lambda}(z)$.
Then \begin{equation*}
|\tDelta_{v,\lambda}|_v \ \le \ h_v^{nN} 
\end{equation*}
for each $v$ and $\lambda$.  

Patch $G_v^{(n-1)}(z)$ by setting\index{patching functions, $G_v^{(k)}(z)$ for $1 \le k \le n$!constructed by patching}
\begin{equation*}
G_v^{(n)}(z) \ = \ G_v^{(n-1)}(z) + \Delta_{v,n}(z) \
\end{equation*}
This replaces the low-order coefficients of the $G_v^{(n)}(z)$ with\index{patching functions, $G_v^{(k)}(z)$} 
the $\tA_{\lambda}$.  

\vskip .1 in
{\bf Conclusion of the Patching Argument.}\index{patching argument!conclusion of global}
\index{global patching when $\Char(K) = p > 0$!conclusion of the patching argument}
The patching process has now arranged that the $G_v^{(n)}(z) \in K_v(\cC)$
\index{patching functions, $G_v^{(k)}(z)$ for $1 \le k \le n$!$G_v^{(n)}(z) = G^{(n)}(z)$ is independent of $v$} 
for $v \in \hS_K$ all coincide with a single function  $G^{(n)}(z)$, 
whose coefficients relative to the $L^{\sep}$-rational basis belong to $L^{\sep}$. 
\index{basis!$L^{\sep}$-rational}  
Fix any $v$, and put $G_w(z) = G(z)$ for all places $w$ of $L^{\sep}$ with $w|v$;  
then $\oplus_{w|v} G_w(z) \in \oplus_{w|v} L^{\sep}(\cC) \cong L^{\sep} \otimes_K K(\cC)$ 
is invariant under $\Gal(L^{\sep}/K)$, so by Proposition \ref{CPPr1p} it belongs to $K(\cC)$.   

For each $v \in \hS_K$, our restrictions on the magnitudes of the 
$\Delta_{v,ij}^{(k)}$ and the $\Delta_{v,\lambda}$ assure that the conclusions
of Proposition \ref{LocalPatch} apply.  Thus 
\begin{equation*} 
\left\{ \begin{array}{l}
\text{If $v \in S_K$, so $E_v$ is $K_v$-simple, 
\index{$K_v$-simple!set}  
        then the zeros of $G^{(n)}(z)$ are distinct and belong to $E_v$\ .} \\
\text{If $v \in \hS_K^{+} \backslash S_K$,
      then $\{ z \in \cC_v(\CC_v) : |G^{(n)}(z)|_v \le R_v^{Nn} \} = E_v\ .$ }
\end{array} \right. 
\end{equation*}
On the other hand, for each $v \notin \hS_K^{+}$,  
our construction has arranged that in the expansion 
\begin{equation*}
G^{(n)}(z) \ = \ 
\sum_{i=1}^m \sum_{j=0}^{(n-1)N_i} A_{ij} \varphi_{i,nN_i-j}(z)
   + \sum_{\lambda = 1}^{\Lambda} A_{\lambda} \varphi_{\lambda} \ ,
\end{equation*}
all the coefficients belong to $\hcO_v$ and the leading coefficients
belong to $\hcO_v^{\times}$.  Our choice of $\hS_K^{+}$
assures that $\cC_v$ and the functions $\varphi_{ij}(z)$
and $\varphi_{\lambda}(z)$ all have good reduction at $v$, 
\index{good reduction}
and the $x_i$ specialize to distinct points $\pmod{v}$.  
Hence $G^{(n)}(z) \pmod{v}$ is a nonconstant function 
with a pole of order $nN_i > 0$ at each $x_i$.
It follows that for each $v \notin \hS_K^{+}$,
\begin{equation*}
\{ z \in \cC_v(\CC_v) : |G^{(n)}(z)|_v \le 1 \} 
  \ = \ \cC_v(\CC_v) \backslash (\bigcup_{i=1}^m B(x_i,1)^{-}) \ = \ E_v \ .
\end{equation*}
  
{\bf Construction of the points in Theorem \ref{aT1-B}.}
\index{patching argument!conclusion of global}
\index{global patching when $\Char(K) = p > 0$!constructing the points in Theorem \ref{aT1-B}}  
The patching argument holds for each integer $n > n_3$ 
divisible by $n_0 n_1 n_2$.  For any such $n$, 
the zeros of $G^{(n)}(z)$ satisfy the conditions of the Theorem.
If there are any $v$ for which $E_v$ is $K_v$-simple, the construction shows that 
\index{$K_v$-simple!set}  
the zeros of $G^{(n)}(z)$ are distinct, and letting $n \rightarrow \infty$ 
we obtain the points in the Theorem.  

However, if there are no such $v$, then since $\prod_{v \in \hS_K^{+}} r_v^{Nn}$
grows arbitrarily large as $n \rightarrow \infty$, 
the number of $\hS_K^{+}$-integers $\kappa \in K$ satisfying
$|\kappa|_v \le r_v^{Nn}$ for all $v \in \hS_K^{+}$ also becomes arbitrarily large.
For any such $\kappa$, the roots of $G^{(n)}(z) = \kappa$ are points
satisfying the conditions of the Theorem.\index{patching functions, $G_v^{(k)}(z)$ for $1 \le k \le n$!roots are confined to $E_v$}  
Hence there are infinitely many such points.

This completes the proof of Theorem \ref{aT1-B} when $\Char(K) = p > 0$.
\end{proof}
\index{patching argument!global!when $\Char(K) = p > 0$|)}

\section{ Proof of Proposition $\ref{FPushBound}$} \label{FPushBoundProofSec}
 
Fix integers $\ell \ge 1$ and $1 \le k \le n-1$.  
Let $F_v(z) \in \CC_v(\cC)$ be an $(\fX,\vs)$-function\index{$(\fX,\vs)$-function} 
with a pole of order $(n-k-1)N_i$ and leading coefficient $d_{v,i} \ne 0$ at each $x_i$, 
whose zeros all belong to $E_v$.  In proving Proposition \ref{FPushBound}, it will be useful
to introduce a scaled version of the $L$-rational basis, consisting of the basis functions 
\index{basis!scaled $L$-rational|ii} 
$\{d_{v,i} \varphi_{ij}, \varphi_\lambda\}$.  
Write  
\begin{eqnarray} 
\Delta_v  & = & \sum_{i=1}^m \sum_{s=0}^{\ell N_i-1} \Delta_{v,is} \cdot \varphi_{i,(\ell+1) N_i-s} \ , 
\label{omegavdef} \\
\Delta_v F_v 
& = & \sum_{i=1}^m \sum_{s=0}^{\ell N_i-1}  
          \delta_{v,is} \cdot d_{v,i} \varphi_{i,(n-k+\ell)N_i - s} 
                 \ + \ \text{lower order terms}  \ , \label{omegavFvexp}
\end{eqnarray} 
put 
$\vDelta = (\Delta_{v,is})_{1 \le i \le m , 0 \le s  < N_i}$ 
$\vdelta  =  (\delta_{v,is})_{1 \le i \le m , 0 \le s < \ell N_i}$, 
and let $\Phi_{F_v} : \CC_v^{\ell N} \rightarrow \CC_v^{\ell N}$ be the linear map defined by 
\begin{equation} \label{Fdelta_vs1B}
\Phi_{F_v}(\vDelta) \ = \ \vdelta \ ,
\end{equation} 
which takes the coefficients of $\Delta_v$ 
to the high-order coefficients of $\Delta_v F_v$.
Note that $\Phi_{F_v}$ decomposes as direct sum of maps 
$\Phi_{F_v,i} : \CC_v^{\ell N_i} \rightarrow \CC_v^{\ell N_i}$, 
since only the terms in $\Delta_v$ involving $\varphi_{i,(\ell+1)N_i-s}$ for $s = 0, \ldots, \ell N_i-1$
can contribute to poles of $\Delta_v F_v$ with order greater than $(n-k)N_i$ at $x_i$.  

\smallskip
The maps $\Phi_{F_v}$, $\Phi_{F_v,i}$ have an intrinsic interpretation as follows.  
Put 
\begin{equation*}
V \ = \ \bigoplus_{i=1}^m \bigoplus_{s=0}^{\ell N_i-1} \CC_v \varphi_{i,(\ell+1) N_i-s} \ , \qquad 
W \ = \ \bigoplus_{i=1}^m \bigoplus_{s=0}^{\ell N_i - 1} \CC_v d_{v,i} \varphi_{i,(n-k+\ell )N_i-s} \ .
\end{equation*}
Then $\Phi_{F_v}$ is the map on coordinates associated to a 
linear transformation $\Phi_{F_v}^0 : V \rightarrow W$, defined as follows. 
For any divisor $D$ on $\cC_v(\CC_v)$, 
let $\Gamma_{\CC_v}(D) = \{f \in \CC_v(\cC) : \div(f) + D \ge 0\}$. 
Then for $D = \sum_{i=1}^m N_i(x_i)$, 
the inclusion of $W$ into $\Gamma_{\CC_v}((n-k+\ell)D )$ induces an isomorphism
\begin{equation*} 
\iota : W \ \cong \ \Gamma_{\CC_v}\big((n-k+\ell)D \big) / \Gamma_{\CC_v}\big((n-k)D\big) \ , 
\end{equation*} 
and for each function $\Delta_v(z) \in V$ 
\begin{equation*}
\Phi_{F_v}^0(\Delta_v) \ = \ \iota^{-1}\big(\, \Delta_v  F_v \!\!\!\! \pmod{\Gamma_{\CC_v}((n-k)D)} \, \big) \ .
\end{equation*} 
Similarly, for each $i$, put 
\begin{equation*}
V_i \ = \ \bigoplus_{s=0}^{\ell N_i-1} \CC_v \varphi_{i,(\ell+1)N_i-s} \ , \qquad 
W_i \ = \ \bigoplus_{s=0}^{\ell N_i - 1} \CC_v d_{v,i} \varphi_{i,(n-k+\ell)N_i-s} \ ,
\end{equation*} 
Then $\Phi_{F_v,i}$ is the map on coordinates associated to a map 
$\Phi_{F_v,i}^0 : V_i \rightarrow W_i$, defined as follows.   
The inclusion of $W_i$ into $\Gamma_{\CC_v}((n-k)D + \ell N_i (x_i))$ induces an isomorphism
\begin{equation*} 
\iota_i : W_i \ \cong \ \Gamma_{\CC_v}\big((n-k)D + \ell N_i(x_i)\big) / \Gamma_{\CC_v}\big((n-k)D\big) \ , 
\end{equation*} 
and for each function $\Delta_v(z) \in V_i$ 
\begin{equation*}
\Phi_{F_v,i}^0(\Delta_v) 
\ = \ \iota_i^{-1}\big(\, \Delta_v  F_v \!\!\!\! \pmod{\Gamma_{\CC_v}((n-k)D)} \, \big) \ .
\end{equation*}  

For each $i$, since $J|N_i$ the functions in $\{\varphi_{i,(\ell+1)N_i-s}\}_{0 \le s < \ell N_i}$
and $\{\tphi_{i,(\ell+1)N_i-s}\}_{0 \le s < \ell N_i}$ are $K(x_i)$-linear combinations of each other, 
and each set forms a basis for $V_i$.  Similarly, the functions in 
$\{d_{v,i} \varphi_{i,(n-k+\ell)N_i-s}\}_{0 \le s < \ell N_i}$
and $\{\tphi_{i,(n-k+\ell)N_i-s}\}_{0 \le s < \ell N_i}$ 
are $K(x_i)$-linear combinations of each other, 
and each set forms a basis for $W_i$.  

\smallskip
The map $\Phi_{F_v}^{\sep}$ in Proposition \ref{FPushBound} is the coordinate map associated to 
$\Phi_{F_v}^0$ using the $L$-rational basis on the source and the $L^{\sep}$-rational basis on the target:
\index{basis!$L$-rational} 
\index{basis!$L^{\sep}$-rational} 
if we write 
\begin{equation*}
\Delta_v F_v 
\ = \ \sum_{i=1}^m \sum_{s=0}^{\ell N_i-1}  
          \tdelta_{v,is} \cdot \tphi_{i,(n-k+\ell)N_i - s} 
                 \ + \ \text{lower order terms}  \ , 
\end{equation*}
and put $\tdelta = (\tdelta_{v,is})_{1 \le i \le m, 0 \le s < \ell N_i}$
then $\Phi_{F_v}^{\sep}(\vDelta) = \tdelta$.  Clearly $\Phi_{F_v}^{\sep}$ decomposes as a 
direct sum of maps $\Phi_{F_v,i}^{\sep} : \CC_v^{\ell N_i} \rightarrow \CC_v^{\ell N_i}$ 
associated to the $\Phi_{F_v,i}^0 : V_i \rightarrow W_i$.

\smallskip
Before proving Proposition \ref{FPushBound}, we will need two lemmas.  
Recall that the $L$-rational basis is multiplicatively generated by finitely many functions. 
\index{basis!$L$-rational}\index{$L$-rational basis!multiplicatively finitely generated}
This means that collectively, the basis functions $\varphi_{ij}(z)$ have only finitely many distinct zeros. 

\begin{lemma} \label{ExpansionCoeffBound} 
Let $K_v$ be nonarchimedean.  
Then there is a constant $\Lambda_v > 0$, depending only on $\fX$, 
the choice of the uniformizing parameters\index{uniformizing parameter!normalizes $L$-rational basis}
 $g_{x_i}(z)$, and the projective embedding of $\cC_v$, 
with the following property$:$  

Let $r > 0$ be small enough that 

\quad $(1)$  $r < \min_{i \ne j} (\|x_i,x_j\|_v);$ 

\quad $(2)$ each of the balls $B(x_i,r)$ is isometrically parametrizable$;$
\index{isometrically parametrizable ball}

\quad $(3)$ for each $i$, none of the $\varphi_{ij}(z)$ has a zero in $B(x_i,r).$

\noindent{Put} $\varpi_v = \min(1,\Lambda_v \cdot r)$, 
and let $\ell$, $k$ be integers with $\ell \ge 1$ and $1 \le k \le n-1$.  

Suppose $F_v(z) \in \CC_v(\cC)$ is an $(\fX,\vs)$-function\index{$(\fX,\vs)$-function} 
which has a pole of order $(n-k-1)N_i$ and leading coefficient $d_{v,i} \ne 0$ at $x_i$, for each $i$. 
Assume $F_v(z)$ has no zeros in $\bigcup_{i=1}^m B(x_i,r)$.  
Then for each $i$, and each integer $0 \le s < \ell N_i$,   
when we expand $\varphi_{i,(\ell+1)N_i - s}(z) F_v(z)$ using the scaled $L$-rational basis as 
\index{basis!scaled $L$-rational}
\begin{eqnarray} 
\varphi_{i,(\ell+1)N_i-s} \cdot F_v
& = & \sum_{t=0}^{\ell N_i - s-1} C_i(s,t) \cdot d_{v,i} \varphi_{i,(n-k+\ell)N_i - s - t} \notag \\
& & \qquad  \ + \ \text{terms with poles of order $\le (n-k)N_{i^\prime}$ 
                 at each $x_{i^\prime}$ \ ,} \label{FCCF0}
\end{eqnarray}
we have $|C_i(s,t)|_v \le 1 /\varpi_v^t$ for each $t$.  
\end{lemma}


\begin{proof}
For each $i$, let $g_{x_i}(z)$  be the uniformizing parameter\index{uniformizing parameter!normalizes $L$-rational basis} 
used to normalize the basis functions $\varphi_{ij}(z)$;  
\index{basis!$L$-rational} 
thus $\lim_{z \rightarrow x_i} F_v(z) \cdot g_{x_i}(z)^{(n-k-1)N_i} = d_{v,i}$ 
and $\lim_{z \rightarrow x_i} \varphi_{ij}(z) \cdot g_{x_i}(z)^j = 1$ for each $j > N_i$.  
Let $\varrho_i : D(0,r) \rightarrow B(x_i,r)$ be an isometric parametrization
\index{isometric parametrization} 
with $\varrho_i(0) = x_i$, and put $b_{v,i} = \lim_{Z \rightarrow 0} g_{x_i}(\varrho_i(Z))/Z$. 

Take $\Lambda_v = \min_{1 \le i \le m} (|b_{v,i}|_v)$, and put 
\begin{equation} \label{varpivDef}
\varpi_v \ = \ \min(1,\Lambda_v r) \ = \ \min(1,r |b_{v,1}|_v, \ldots, r |b_{v,m}|_v) \ .  
\end{equation}  
To establish the bounds in the Lemma, first fix $i$.  
By abuse of notation, write $F_v(Z)$ for $F_v(\varrho_i(Z))$  
and $\varphi_{ij}(Z)$ for $\varphi_{ij}(\varrho_i(Z))$.  
For compactness of notation, temporarily write $h = n-k-1$. Then 
\begin{equation*}
\lim_{Z \rightarrow 0} F_v(Z) \cdot Z^{h N_i} 
\ = \ \lim_{Z \rightarrow 0} \Big(F_v(\varrho_i(Z))\cdot g_{x_i}(\varrho_i(Z))^{hN_i}\Big) 
               \cdot \Big(\frac{Z}{g_{x_i}(\varrho_i(Z))}\Big)^{hN_i}  
\ = \ d_{v,i} b_{v,i}^{-hN_i} \ , 
\end{equation*}
so $F_v(Z)$ has a Laurent expansion of the form 
\begin{equation} \label{FCCF1} 
F_v(Z) \ = \ d_{v,i} b_{v,i}^{-(n-k-1)N_i} \cdot Z^{-(n-k-1)N_i} \cdot (1 + \sum_{j=1}^\infty f_j Z^j ) \ .
\end{equation}
Since $F_v(Z)$ has no zeros in $D(0,r)$, the theory of Newton Polygons\index{Newton Polygon} shows that $|f_j|_v < 1/r^j$
for each $j \ge 1$ (see Lemma \ref{BLem1} and the discussion before it).  
Similarly, for each $0 \le s < \ell N_i$
\begin{equation} \label{FCCF3}
\varphi_{i,(\ell + 1)N_i - s}(Z) \cdot  F_v(Z) \ = \ 
d_{v,i} b_{v,i}^{-(n-k+\ell)N_i + s} \cdot Z^{-(n-k+\ell ) N_i + s} 
\cdot (1 + \sum_{j=1}^{\infty} f_{is,j} Z^j  ) 
\end{equation}
with $|f_{is,j}|_v < 1/r^j$ for each $j \ge 1$, and for each $0 \le t \le \ell N_i - s - 1$ 
\begin{equation} \label{FCCF2}
\varphi_{i, (n-k+\ell)N_i - s - t}(Z) 
\ = \ b_{v,i}^{-(n-k+\ell)N_i + s + t} Z^{-(n-k+\ell )N_i + s + t} \cdot (1 + \sum_{j=1}^{\infty} c_{i,s+t,j} Z^j )
\end{equation}
with $|c_{i,s+t,j}|_v < 1/r^j$ for each $j \ge 1$. 

To prove the lemma, fix $0 \le s < \ell N_i$, insert the expansions (\ref{FCCF2}) into (\ref{FCCF0}) 
and compare the coefficients of the resulting series with those in (\ref{FCCF3}). 
Comparing the coefficients of $Z^{-(n-k+\ell)N_i + s}$ we see that when $t = 0$
\begin{equation*} 
C_i(s,0) \cdot b_{v,i}^{-(n-k+\ell) N_i + s} \ = \ b_{v,i}^{-(n-k+\ell)N_i + s} \ ,
\end{equation*} 
so $C_i(s,0) = 1$; trivially, $|C_i(s,0)|_v \le 1/(r|b_{v,i}|_v)^0$.  
Inductively, take $1 \le t \le \ell N_i - s - 1$ 
and assume that $|C_i(s,j)|_v \le 1/(r|b_{v,i}|_v)^j$ for $0 \le j \le t-1$.  
Comparing the coefficients of $Z^{-(n-k+\ell)N_i + s + t}$ we find that 
\begin{equation*}  
 C_i(s,t) \cdot b_{v,i}^{-(n-k+\ell) N_i + s + t} 
+ \sum_{j = 0}^{t-1} C_i(s,j) \cdot b_{v,i}^{-(n-k+\ell)N_i + s + j} c_{i,s + j, t-j} 
\ = \ b_{v,i}^{-(n-k+\ell)N_i + s} \cdot f_{is,t} \ ,
\end{equation*} 
or equivalently   
\begin{equation*}  
C_i(s,t) \ = \  b_{v,i}^{-t} \cdot f_{is,t} - 
\sum_{j = 0}^{t-1} C_i(s,j) \cdot b_{v,i}^{j-t } c_{i, s + j, t-j} \ .
\end{equation*}
Since $|f_{i s, t}|_v \le 1/r^t$, 
it follows that $|b_{v,i}^{-t} f_{i s,t}|_v \le 1/(r|b_{v,i}|_v))^t$. 
Similarly, for $0 \le j \le t-1$ we have $|c_{i,s + j, t-j}|_v \le 1/r^{t-j}$, 
so by induction, for each such $j$ 
\begin{equation*} 
|C_i(s,j) \cdot b_{v,i}^{j-t} c_{i,s + t, t-j} |_v \ \le \ 
      1/(r |b_{v,i}|_v)^j \cdot |b_{v,i}|_v^{j-t} \cdot 1/r^{t-j} \ = \ 1/(r|b_{v,i}|_v)^t \ .
\end{equation*}   
By the ultrametric inequality $|C_i(s,t)|_v \le 1/(r|b_{v,i}|_v)^t$,
and the induction can continue.

\smallskip  
Now let $i$ vary.  For each $i$ we have $r |b_{v,i}|_v \ge \Lambda_v r \ge \varpi_v$, 
and the Lemma follows.
\end{proof} 

The following lemma gives bounds for the entries of a the inverse of a 
unipotent lower triangular matrix,\index{lower triangular matrix|ii} given a suitable bound for the entries in each subdiagonal. 
The indices $i, j, k, \ell$ in the lemma are unrelated to $i, j, k$ and $\ell$ as used elsewhere.   

\begin{lemma} \label{InverseCoeffBound}
Let $C \in M_k(\CC_v)$ be a lower triangular matrix whose diagonal elements are $1$,
and write it as 
\begin{equation} \label{FMatrixF}
C \ = \ \left( \begin{array}{ccccc} 
            1     &    0      &   0     & \cdots &    0   \\
         C_{2,1}  &    1      &   0     & \cdots &    0   \\
         C_{3,1}  &  C_{3,2}  &   1     & \cdots &    0   \\
           \vdots & \vdots    &  \vdots & \ddots & \vdots \\
         C_{k,1}  &  C_{k,2}  & C_{k,3} & \ldots &    1       
               \end{array} \right) \ .
\end{equation} 
Let $\varpi_v > 0$ be such that for each $h = 1, \ldots, k$,
the elements $C_{i,j}$ belonging to the $\ell^{th}$ subdiagonal $($i.e. those with $i-j = \ell)$, 
satisfy $|C_{i,j}|_v \le 1/\varpi_v^\ell$. Then 
\begin{equation} \label{FInverseF}
C^{-1} \ = \ \left( \begin{array}{ccccc} 
            1     &    0      &   0     & \cdots &    0   \\
         c_{2,1}  &    1      &   0     & \cdots &    0   \\
         c_{3,1}  &  c_{3,2}  &   1     & \cdots &    0   \\
           \vdots & \vdots    &  \vdots & \ddots & \vdots \\
         c_{k,1}  &  c_{k,2}  & c_{k,3} & \ldots &    1       
               \end{array} \right) \ .
\end{equation}
and for all $c_{i,j}$ in the $\ell^{th}$ subdiagonal, we have $|c_{i,j}|_v \le 1/\varpi_v^\ell$.
\end{lemma}  

\begin{proof} Clearly $C^{-1}$ exists and has the form (\ref{FInverseF}).  
To show that $|c_{i,j}|_v \le 1/\varpi_v^{i-j}$, we use induction on $\ell = i-j$.
When $i-j = 1$, the $(i,j)$ term in $C \cdot C^{-1} = I$ is 
\begin{equation*} 
c_{i,j}  +  C_{i,j} \ = \ 0 \ ,
\end{equation*}
so $|c_{i,j}|_v = |C_{i,j}|_v \le 1/\varpi_v$.  Now suppose $i - j = \ell > 1$.  
The $(i,j)$ term in  $C \cdot C^{-1} = I$ is 
\begin{equation*} 
C_{i,j} + C_{i,j+1} c_{j+1,j} + \cdots + C_{i,i-1}c_{i-1,j} + c_{i,j} \ = \ 0 \ .
\end{equation*}    
Assuming that $|c_{i^{\prime},j^{\prime}}|_v \le 1/\varpi_v^{1^{\prime}-j^{\prime}}$ 
for all $(i^{\prime},j^{\prime})$ with $i^{\prime}-j^{\prime} < \ell$, 
and using our hypothesis on $C$, the ultrametric inequality gives  
\begin{equation*} 
|c_{i,j}|_v \ \le \ 
\max\big( |C_{i,j+1}|_v |c_{j+1,j}|_v, \cdots,  |C_{i,i-1}|_v |c_{i-1,j}|_v, |C_{i,j}|_v  \big)
\ \le \ 1/\varpi_v^{i-j} 
\end{equation*}
as desired.
\end{proof}

We can now prove Proposition \ref{FPushBound}.

\begin{proof}[Proof of Proposition \ref{FPushBound}] 
Let $F_v(z) \in \CC_v(\cC)$ be an $(\fX,\vs)$-function\index{$(\fX,\vs)$-function} 
with a pole of order 
$(n-k-1)N_i$ and leading coefficient $d_{v,i} \ne 0$ at each $x_i$,  
whose zeros all belong to $E_v$.  Writing  
\begin{eqnarray} 
\qquad \Delta_v(z)  & = & 
 \sum_{i=1}^m \sum_{s=0}^{\ell N_i-1} \Delta_{v,is} \cdot \varphi_{i,(\ell+1)N_i-s}(z) 
\ , 
\label{omegavF2}\\
\qquad \Delta_v(z) F_v(z) 
& = & \sum_{i=1}^m \sum_{s=0}^{\ell N_i-1}  
          \delta_{v,is} \cdot d_{v,i} \varphi_{i,(n-k+\ell)N_i - s} (z)
                 \ + \ \text{lower order terms}  \ , \label{omegavFvLscaled}
\end{eqnarray} 
put 
\begin{equation*}
\vDelta \ = \ (\Delta_{v,is})_{\substack{1 \le i \le m \\ 0 \le s < \ell N_i}} \ , \qquad
\vdelta \ = \ (\delta_{v,is})_{\substack{1 \le i \le m \\ 0 \le s < \ell N_i}} \ , 
\end{equation*} 
and define $\Phi_{F_v} : \CC_v^{\ell N} \rightarrow \CC_v^{\ell N}$ by 
$\Phi_{F_v}(\vDelta) = \vdelta$ as in (\ref{Fdelta_vs1B}). 
If we write $\vDelta_i = (\Delta_{v,is})_{0 \le s < \ell N_i}$
and $\vdelta_i = (\delta_{v,it})_{0 \le t < \ell N_i}$, then 
$\Phi_{F_v}$ is a direct sum of the maps $\Phi_{F_v,i} : \CC_v^{\ell N_i} \rightarrow \CC_v^{\ell N_i}$
with $\Phi_{F_v,i}(\vDelta_i) = \vdelta_i$ for each $i$.  

Let $ r > 0$ be small enough that 

\smallskip

\quad $(1)$ $r < \min_{i \ne j}( \|x_i,x_j\|_v)$;

\quad $(2)$ each of the balls $B(x_i,r)$ is isometrically parametrizable and disjoint from $E_v$;
\index{isometrically parametrizable ball}

\quad $(3)$ for each $i$, none of the $\varphi_{ij}$ has a zero in $B(x_i,r)$. 

\smallskip

\noindent{Let} $\Lambda_v > 0$ and  $\varpi_v = \min(1,\Lambda_v \cdot r) > 0$ 
be as in Lemma \ref{ExpansionCoeffBound}.

We begin by showing that $\Phi_{F_v}$ 
is an isomorphism and that for each $\rho > 0$ 
\begin{equation} \label{FBallC1}
\Phi_{F_v}\big(\bigoplus_{i=1}^m \bigoplus_{s=0}^{\ell N_i-1} D(0,\varpi_v^{-s} \rho) \big) \ \supseteq \ 
\bigoplus_{i=1}^m \bigoplus_{t=0}^{\ell N_i-1} D(0,\varpi_v^{-t} \cdot \rho ) \ .
\end{equation} 
For this, it is enough to show that each $\Phi_{F_v,i}$ is an isomorphism, and that  
\begin{equation} \label{FSupset} 
\Phi_{F_v,i} \big( \bigoplus_{s=0}^{\ell N_i-1} D(0,\varpi_v^{-s}\rho) \big) \ \supseteq \ 
\bigoplus_{s=0}^{\ell N_i-1} D(0,\varpi_v^{-s} \cdot \rho ) \ .
\end{equation}


\smallskip
Fix $i$.  As in (\ref{FCCF0}), for each $0 \le t < \ell N_i$ the product  
$\varphi_{i,(\ell+1)N_i -t} \cdot F_v$ can be expanded using the scaled $L$-rational basis as 
\index{basis!scaled $L$-rational}   
\begin{equation*}
d_{v,i} \varphi_{i, (n-k+\ell)N_i - s}
      + \sum_{t=1}^{\ell N_i - s -1} C_i(s,t) \cdot  d_{v,i} \varphi_{i,(n-k+\ell)N_i - s - t} 
                 \ + \ \text{terms not contributing to $\Phi_{F_v,i}$ \ .} 
\end{equation*} 
This means that the matrix of $\Phi_{F_v,i}$ is 
\begin{equation} \label{FCiMatrix}
C_i \ = \  \left(  \begin{array}{ccccc} 
                1    &   0        &   0       & \cdots &    0    \\
            C_i(0,1) &   1        &   0       & \cdots &    0   \\
            C_i(0,2) &  C_i(1,1)  &   1       & \cdots &    0  \\
             \vdots  &  \vdots    &   \vdots        & \ddots &  \vdots  \\
     C_i(0,\ell N_i - 1)  &  C_i(1,\ell N_i - 2) & \ldots  & C_i(\ell N_i - 2,1)  &  1       
               \end{array} \right) \ ;
\end{equation}  
in particular, $\Phi_{F_v,i}$ is nonsingular. 
By Lemma \ref{ExpansionCoeffBound} we have $C_i(s,t) \varpi_v^{-t}$ for all $s, t$.  
Hence by Lemma \ref{InverseCoeffBound}, 
\begin{equation*} 
C_i^{-1} \ = \  \left(  \begin{array}{ccccc} 
                1    &   0        &   0       & \cdots &    0    \\
            c_i(0,1) &   1        &   0       & \cdots &    0   \\
            c_i(0,2) &  c_i(1,1)  &   1       & \cdots &    0  \\
             \vdots  &  \vdots    &   \vdots        & \ddots &  \vdots  \\
     c_i(0,\ell N_i - 1)  &  c_i(1,\ell N_i - 2) & \ldots  & c_i(\ell N_i - 2,1)  &  1       
               \end{array} \right) \ , 
\end{equation*} 
with $c_i(s,t) \le \varpi_v^{-t}$ for all $s, t$.  

Fix $\rho > 0$, and assume that $|\delta_{v,it}|_v \le \varpi_v^{-s} \rho$ for $0 \le t < \ell N_i$.
Since $\vDelta_i = \Phi_{F_v}^{-1}(\vdelta_i) = C_i^{-1} \vdelta_i$, this means that for each $s$ 
\begin{equation} \label{FPhiInvF}
\Delta_{v,is} \ = \ \big(\sum_{t=0}^{s-1} c_i(t,s-t) \delta_{i,t} \big) +  \delta_{i,s} \ , 
\end{equation} 
and the ultrametric inequality shows that $|\Delta_{v,is}|_v \le \varpi_v^{-s} \rho$. 
 This proves (\ref{FSupset}) and (\ref{FBallC1}).

\smallskip
We can now prove (\ref{FBallC2Q}) and (\ref{FBallC2}) in Proposition \ref{FPushBound}.  
Write $\Delta_v(z) F_v(z)$ using the $L$-rational basis as 
\index{basis!$L$-rational} 
\begin{equation} \label{omegavFvL} 
\Delta_v(z) F_v(z) \ = \ 
\sum_{i=1}^m \sum_{s=0}^{\ell N_i - 1} \hdelta_{v,is} \varphi_{i,(n-k+\ell)N_i - s}(z)  
\ + \ \text{lower order terms} \ .
\end{equation} 
Comparing (\ref{omegavFvL}) and (\ref{omegavFvLscaled}) 
shows that $|\hdelta_{v,is}|_v = |d_{v,is}|_v |\delta_{v,is}|_v$
for all $i, s$.  Next, expand $\Delta_v(z) F_v(z)$ in terms of the $L^{\sep}$-rational basis as 
\index{basis!$L^{\sep}$-rational}  
\begin{equation} \label{omegavFvLswp} 
\Delta_v(z) F_v(z) \ = \ 
\sum_{i=1}^m \sum_{s=0}^{\ell N_i - 1} \tdelta_{v,is} \tphi_{i,(n-k+\ell)N_i - s}(z)  
\ + \ \text{lower order terms} \ .
\end{equation}

By Proposition \ref{TransitionProp}(C), for each $i = 1, \ldots, m$  
there is an invertible $J \times J$ matrix $\tcB_i$
which expresses each set of $J$ consecutive basis 
elements $\{\varphi_{ij}\}_{h J + 1 \le j \le (h+1)J}$ of the $L$-rational basis 
\index{basis!$L$-rational} 
in terms of the corresponding set  $\{\tphi_{ij}\}_{h J + 1 \le j \le (h+1)J}$ from the 
$L^{\sep}$-rational basis.  
\index{basis!$L^{\sep}$-rational} 
Since $J|N_i$ for each $i$, it follows that there is a constant $\tUpsilon_v > 0$
such that if $|\tdelta_{v,is}|_v \le \tUpsilon_v \varpi_v^{-s} |d_{v,i}|_v \cdot \rho$ for all $i, s$, 
then $|\hdelta_{v,is}|_v \le \varpi_v^{-s} |d_{v,i}|_v \cdot \rho$ for all $i, s$.  This in turn means 
$|\delta_{v,is}|_v \le \varpi_v^{-s}  \rho$ for all $i, s$.   
By the discussion above, $|\Delta_{v,is}|_v \le \varpi_v^{-s} \rho$ for all $i, s$.   

Since $\Phi_{F_v}^{\sep}$ is the coordinate map for $\Phi_{F_v}^0$ using the $L$-rational basis 
\index{basis!$L$-rational} 
on the source and the $L^{\sep}$-rational basis on the target, we see that
\index{basis!$L^{\sep}$-rational} 
\begin{equation} \label{FBallC2A}
\Phi^{\sep}_{F_v}\big(\bigoplus_{i=1}^m \bigoplus_{s=0}^{\ell N_i-1} D(0,\varpi_v^{-s} \rho) \big) \ \supseteq \ 
\bigoplus_{i=1}^m \bigoplus_{s=0}^{\ell N_i-1} D(0, \tUpsilon_v \varpi_v^{-s} |d_{v,i}|_v \cdot \rho ) \ ,
\end{equation} 
which is (\ref{FBallC2}). To show (\ref{FBallC2Q}), 
replace $\rho$ with $\varpi_v^{\ell N} \rho$ in (\ref{FBallC2A}).
Since $0 < \varpi_v \le 1$, 
\begin{eqnarray*}
\Phi^{\sep}_{F_v}\big(\bigoplus_{i=1}^m \bigoplus_{s=0}^{\ell N_i-1} D(0,\varpi_v^{\ell N-s} \rho) \big)  
& \supseteq &
\bigoplus_{i=1}^m \bigoplus_{s=0}^{\ell N_i-1} D(0, \tUpsilon_v \varpi_v^{\ell N-s} |d_{v,i}|_v \cdot \rho ) \\
& \supseteq & 
\bigoplus_{i=1}^m D(0, \tUpsilon_v \varpi_v^{\ell N} |d_{v,i}|_v \cdot \rho )^{\ell N_i} \ . 
\end{eqnarray*}
This yields (\ref{FBallC2Q}) since $D(0,\rho)^{\ell N} \supseteq 
\bigoplus_{i=1}^m \bigoplus_{s=0}^{\ell N_i-1} D(0,\varpi_v^{\ell N-s} \rho)$.

\smallskip
Finally, we prove the rationality assertions in Proposition \ref{FPushBound}. 

Assume that $F_v(z)$ is $K_v$-rational.  Using the $L^{\sep}$-rational basis, we can write 
\index{basis!$L^{\sep}$-rational} 
\begin{eqnarray} 
\Delta_v(z)  & = & \sum_{i=1}^m \sum_{s=0}^{\ell N_i-1} \tDelta_{v,is} \cdot \tphi_{i,(\ell+1)N_i-s}(z) \ , 
\label{omegavF2t} \\
\Delta_v(z) F_v(z) 
& = & \sum_{i=1}^m \sum_{s=0}^{\ell N_i-1}  
          \tdelta_{v,is} \cdot \tphi_{i,(n-k+\ell)N_i - s}(z) 
                 \ + \ \text{lower order terms}  \ , \label{omegavFvLscaledt}
\end{eqnarray} 
Put 
\begin{equation*}
\tDelta \ = \  (\tDelta_{v,is})_{\substack{1 \le i \le m \\ 0 \le s < \ell N_i}} \ , \qquad 
\tdelta \ = \ (\tdelta_{v,is})_{\substack{1 \le i \le m\\ 0 \le s < \ell N_i}} \ , 
\end{equation*} 
and for each $i = 1, \ldots, m$ put
\begin{equation*}
\tDelta_i \ = \ (\tDelta_{v,is})_{ 0 \le s < \ell N_i} \ , \qquad 
\tdelta_i \ = \ (\tdelta_{v,is})_{ 0 \le s < \ell N_i} \ . 
\end{equation*} 

Fix $i$.  Then for each $0 \le t < \ell N_i$ the product  
$\tphi_{i,(\ell+1)N_i -t} \cdot F_v$ is rational over $K_v(x_i)^{\sep}$ and 
can be expanded using the $L^{\sep}$-rational basis as  
\index{basis!$L^{\sep}$-rational}  
\begin{equation*}
   \tphi_{i,(\ell+1)N_i -t} \cdot F_v  \ = \ \sum_{s=0}^{\ell N_i -1} \tC_i(s,t) \cdot  \varphi_{i,(n-k+\ell)N_i - s} 
                 \ + \ \text{terms not contributing to $\Phi_{F_v,i}^0$ \ .} 
\end{equation*} 
Since each $\tphi_{i,(\ell+1)N_i -t}$ is rational over $K_v(x_i)^{\sep}$, by galois equivariance
each $\tC_i(s,t)$ belongs to $K_v(x_i)^{\sep}$.  Thus the matrix for $\Phi_{F_v,i}^0$ 
using the $L^{\sep}$-rational basis on the source and the target is 
\index{basis!$L^{\sep}$-rational} 
\begin{equation*}
\tC_i \ = \ \left( \begin{array}{ccc}
                             \tC_i(0,0)       &  \cdots    &    \tC_i(0,\ell N_i - 1) \\
                             \vdots           &  \ddots    &          \vdots          \\
                         \tC_i(\ell N_i-1,0)  &   \cdots   &    \tC_i(\ell N_i - 1,\ell N_i - 1) 
                   \end{array} \right)  \ \in \ GL_{\ell N_i}(K_v(x_i)^{\sep}) \ .  
\end{equation*} 
Because the $L^{\sep}$-rational basis is $K_v$-symmetric, 
\index{basis!$L^{\sep}$-rational} 
\index{$K_v$-symmetric!set of functions}
the collection of matrices $\{\tC_1, \ldots, \tC_m\}$
is $K_v$-symmetric.  

\smallskip
Suppose  $\tdelta$ belongs to $(L_{w_v}^{\sep})^{\ell N}$\index{distinguished place $w_v$} and is $K_v$-symmetric. 
\index{$K_v$-symmetric!vector} 
Then $\tdelta_i$ belongs to $(K_v(x_i)^{\sep})^{\ell N_i}$ for each $i$, 
and the set of vectors $\{\tdelta_1, \ldots, \tdelta_m\}$ is $K_v$-symmetric. 
\index{$K_v$-symmetric!set of vectors} 
It follows that $\tDelta_i = \tC_i^{-1} \tdelta_i$ belongs to $(K_v(x_i)^{\sep})^{\ell N_i}$ for each $i$,
and the set of vectors  $\{\tDelta_1, \ldots, \tDelta_m\}$ is $K_v$-symmetric. 
\index{$K_v$-symmetric!set of vectors} 
Since 
\begin{equation*}
\Delta_v(z) \ = \ \sum_{i=1}^m \sum_{s=0}^{\ell N_i-1} \tDelta_{v,is} \cdot \tphi_{i,(\ell+1)N_i-s}(z)\ , 
\end{equation*}
where the $\tDelta_{v,is}$ and $\tphi_{i,(\ell+1)N_i-s}$ 
are $K_v$-symmetric and rational over $L_{w_v}^{\sep}$,\index{distinguished place $w_v$}
\index{$K_v$-symmetric!set of functions}
it follows that $\Delta_v$ is $K_v$-rational.

If we re-express $\Delta_v$ in terms of the $L$-rational basis as 
\index{basis!$L$-rational} 
\begin{equation*} 
\Delta_v(z) \ = \ \sum_{i=1}^m \sum_{s=1}^{\ell N_i} \Delta_{v,is} \cdot \varphi_{i,N_i+s}(z) \ ,
\end{equation*} 
then the associated vector $\vDelta = (\Delta_{v,is})_{1 \le i \le m, 1 \le s \le \ell N_i}$
is the unique solution to $\Phi^{\sep}_{F_v}(\vDelta) = \tdelta$ in $\CC_v^{\ell N}$.   
Since $\Delta_v$ is $K_v$-rational, necessarily $\vDelta$ belongs to $L_{w_v}^{\ell N}$\index{distinguished place $w_v$} 
and is $K_v$-symmetric. 
\index{$K_v$-symmetric!vector}  

This completes the proof of Proposition \ref{FPushBound}.  
\end{proof}
 

%% file: NewFSZChap8.tex
\chapter{ The Local Patching Construction when $K_v \cong \CC$ } 
\label{Chap8}

In this section we give the confinement argument for Theorem \ref{aT1-B} when
\index{confinement argument}\index{patching construction!for the case when $K_v \cong \CC$|(}
\index{patching argument!local!for the case when $K_v \cong \CC$|(}    

$K_v \cong \CC$.  Write $\CC$ for $\CC_v$
and $|x|$ for $|x|_v$.  Let $w_v$ be the distinguished place\index{distinguished place $w_v$} 
of $L = K(\fX)$ determined by the embedding $\tK \hookrightarrow \CC$
used to identify $\fX$ with a subset of $\cC_v(\CC)$, 
and identify $L_{w_v} = K_v = \CC$.  

      At several places in this section, we assert that certain objects are $K_v$-symmetric.
\index{$K_v$-symmetric}   
Since $\Aut_c(\CC_v/K_v) = \Aut(\CC/\CC)$ is trivial, this is a vacuous condition. 
However, we include it for compatibility with the results stated 
in Chapters \ref{Chap9} -- \ref{Chap11}.    

\smallskip
Following the construction of the coherent approximating functions in Theorem \ref{CTCX2},
\index{coherent approximating functions $\phi_v(z)$} 
we begin with the following data:   
\begin{enumerate}
  \item A $K_v$-symmetric probability vector $\vs \in \cP^{m}(\QQ)$ 
\index{$K_v$-symmetric!probability vector}
           with positive rational coefficients.

  \item A $\CC$-simple set $E_v \subset \cC_v(\CC) \backslash \fX$:  
           thus, $E_v$ is compact and nonempty  
           with finitely many connected components, each of which is simply connected\index{simply connected},  
           has a piecewise smooth boundary,\index{boundary!piecewise smooth} 
           and is the closure of its interior $E_v^0$.\index{closure of $\cC_v(\CC)$ interior} 
           
  \item Parameters $h_v$, $r_v$, $R_v$ with $1 < h_v < r_v < R_v$,
          which govern the freedom in the patching process.
          \index{patching parameters}\index{patching argument!local}
          
  \item A number $N$  and an $(\fX,\vs)$-function $\phi_v(z) \in K_v(\cC)$ 
           of degree $N$ such that 
            \begin{equation*}
              \{ z \in \cC_v(\CC) : |\phi_v(z)| \le R_v^N \} \ \subset \ E_v^0 \ .      
            \end{equation*} 
            
  \item An order $\prec_N$ on the index set\index{order!$\prec_N$}
              $\cI = \{(i,j) \in \ZZ^2 : 1 \le i \le m, 0 \le j\}$ 
           determined by $N$ and $\vs$ as in (\ref{FPrec}), 
           which gives the sequence in which coefficients are patched. 
 \index{coefficients $A_{v,ij}$} 
\end{enumerate} 

We will use the $L$-rational basis  $\{\varphi_{ij}, \varphi_{\lambda}\}$ 
\index{basis!$L$-rational}
from \S\ref{Chap3}.\ref{LRationalBasisSection} to expand all functions, 
and $\Lambda = \dim_K(\Gamma(\sum_{i=1}^m N_i(x_i)))$ 
will be the number of low-order basis elements, 
as in the global patching process.\index{patching argument!global}
The order $\prec_N$\index{order!$\prec_N$} respects the $N$-bands (\ref{FBand}),\index{band!$\Band_N(k)$} 
and for each $x_i \in \fX$,  
specifies the terms to be patched in decreasing pole order.

\begin{theorem}  \label{DCPatch1}\index{patching theorem!for the case when $K_v \cong \CC$} 
Suppose $K_v \cong \CC$. 
Let $E_v \subset \cC_v(\CC) \backslash \fX$ be $\CC$-simple, 
with interior $E_0$. 
Let $\vs \in \cP^m(\QQ)$ be a $K_v$-symmetric probability vector with positive
\index{$K_v$-symmetric!probability vector} 
rational coefficients, let $1 < h_v < r_v < R_v$ be numbers, 

Let $\phi_v(z) \in K_v(\cC)$ be an $(\fX,\vs)$-function of degree $N$ satisfying 
\begin{equation*}
\{ z \in \cC_v(\CC) : |\phi_v(z)| \le R_v^{N} \} \ \subset \ E_v^0 \ .
\end{equation*}  
Let $N_i = Ns_i$ for each $x_i$, 
and let $\tc_{v,i} = \lim_{z \rightarrow x_i} \phi_v(z) \cdot g_{x_i}(z)^{N_i}$
be the leading coefficient of $\phi_v(z)$ at $x_i$. Put
\index{coefficients $A_{v,ij}$!leading}  
\begin{equation*}
M_v \ = \ \max( \max_{\substack{ 1 \le i \le m \\ N_i < j \le 2N_i}} \|\varphi_{ij}\|_{E_v },
    \max_{1 \le \lambda \le \Lambda} \|\varphi_{\lambda}\|_{E_v } ) \ .
\end{equation*} 
Let  $k_v > 0$  be the least integer such that
\begin{equation} \label{DCPFF0}
\frac{2 N M_v}{1 - (h_v/r_v)^N}
                   \cdot \Big(\frac{h_v}{r_v}\Big)^{k_v N}  \ < \ \frac{1}{4} \ , 
\end{equation}
and let $\kbar \ge k_v$ be a fixed integer.  
Let $B_v > 0$  be an arbitrary constant.  Then there is an integer $n_v$ 
depending on $\phi_v(z)$, $\kbar$, $B_v$, $r_v$, and $R_v$, such
that for each sufficiently large integer $n$ divisible by $n_v$,  
one can carry out the local patching process at $K_v$ as follows:
\index{patching argument!local}  

Put $G_v^{(0)}(z) = \phi_v(z)^n$.\index{patching functions, initial $G_v^{(0)}(z)$!construction of}  
For each $k = 1, \ldots, n-1$, 
let $\{\Delta_{v,ij}^{(k)} \in L_{w_v}\}_{(i,j) \in \Band_N(k)}$\index{distinguished place $w_v$} 
be a $K_v$-symmetric set of numbers given in $\prec_N$ order,
\index{order!$\prec_N$}\index{band!$\Band_N(k)$}
\index{$K_v$-symmetric!set of numbers} 
subject to the conditions that for each $i$, we have $\Delta_{v,i0}^{(1)} = 0$ and for each $j > 0$   
\begin{equation} \label{FBound1A}
|\Delta_{v,ij}^{(k)}| \ \le \ \left\{
      \begin{array}{ll} B_v & \text{if \ $k \le \kbar \ ,$} \\
                         h_v^{kN} & \text{if \ $k > \kbar \ .$}
      \end{array} \right.
\end{equation}
For $k = n$, let 
$\{\Delta_{v,\lambda}^{(n)} \in L_{w_v}\}_{1 \le \lambda \le \Lambda}$\index{distinguished place $w_v$} 
be an arbitrary $K_v$-symmetric set of numbers satisfying 
\index{$K_v$-symmetric!set of numbers} 
\begin{equation} \label{FBound2A} 
|\Delta_{v,\lambda}^{(n)}| \ \le \ h_v^{nN} \ .  
\end{equation} 

\vskip .05 in
\noindent{Then} one can inductively construct $(\fX,\vs)$-functions 
$G_v^{(1)}(z), \ldots, G_v^{(n)}(z) \in K_v(\cC)$,\index{patching functions, $G_v^{(k)}(z)$ for $1 \le k \le n$!constructed by patching} 
of common degree $nN$, having leading coefficient $\tc_{v,i}^n$ at each $x_i$,
\index{coefficients $A_{v,ij}$!leading}
and satisfying  

\vskip .05 in
\noindent{$(A)$} For each $k = 1, \ldots, n$, 
there are functions\index{compensating functions $\vartheta_{v,ij}^{(k)}(z)$} 
$\vartheta_{v,ij}^{(k)}(z) \in L_{w_v}(\cC)$,\index{distinguished place $w_v$} 
determined recursively in $\prec_N$ order\index{order!$\prec_N$}, such that 
\begin{eqnarray*} 
G_v^{(k)}(z) & = & G_v^{(k-1)}(z)
              +  \sum_{(i,j) \in \Band_N(k)}
                           \Delta_{v,ij}^{(k)} \vartheta_{v,ij}^{(k)}(z)  
        \quad \text{for $k < n \ ,$} \\
G_v^{(n)}(z) & = & G_v^{(n-1)}(z) + \sum_{\lambda=1}^{\Lambda}
                       \Delta_{v,\lambda}^{(n)} \varphi_{\lambda}(z) \ ,                      
\end{eqnarray*}
and where for each $(i,j)$,\index{band!$\Band_N(k)$}  

$(1)$ $\vartheta_{v,ij}^{(k)}(z)$
\index{compensating functions $\vartheta_{v,ij}^{(k)}(z)$!poles and leading coefficients of} 
has a pole of order $nN_i-j > (n-k-1)N_i$ at $x_i$ 
and leading coefficient $\tc_{v,i}^{n-k-1}$,  
\index{coefficients $A_{v,ij}$!leading}
a pole of order at most $(n-k-1)N_{i^{\prime}}$ at each $x_{i^{\prime}} \ne x_i$, 
and no other poles $;$  

$(2)$ $\sum_{(i^{\prime},j) \in \Aut_c(\CC_v/K_v)(i,j)} 
  \Delta_{v,i^{\prime}j}^{(k)} \vartheta_{v,i^{\prime}j}^{(k)}(z)$
belongs to $K_v(\cC)$ $;$ 

\vskip .05 in
\noindent{$(B)$} For each $k = 1, \ldots, n$, \quad  
$\{ z \in \cC_v(\CC) : |G_v^{(k)}(z)| \le r_v^{Nn} \} \subset E_v^0$.
\end{theorem}


\noindent{\bf Remark.} 
A key aspect of Theorem \ref{DCPatch1} is that by choosing $n$ appropriately, 
the freedom $B_v$ in patching the coefficients for $k \le \kbar$ can be made arbitrarily large.
\index{patching argument!freedom $B_v$ in patching}\index{freedom $B_v$ in patching}
\index{coefficients $A_{v,ij}$} 
The patching procedure accomplishes this by 
exploiting a phenomenon of `magnification' introduced in (\cite{RR3}). 
\index{magnification argument|ii}  
It first raises $\phi_v(z)$ to a power $n_v > \kbar$, 
so that $F_v(z) = \phi_v(z)^{n_v}$ has enough coefficients to adjust independently.
\index{coefficients $A_{v,ij}$}   
It then varies those coefficients `infinitesimally', preserving the analytic properties of $F_v(z)$.
Finally it raises the modified $F_v(z)$ to a further power $m_v$ with $n = m_v n_v$, 
creating large changes in the 
coefficients of $G_v^{(0)}(z)$.\index{patching functions, initial $G_v^{(0)}(z)$}

\medskip
\begin{proof}[Proof of Theorem \ref{DCPatch1}] 
Let $\kbar \ge k_v$ and $B_v > 0$ be as in the Theorem.  
Let  $n_v > 0$ be an integer large enough that
\begin{eqnarray}
&& n_v \ > \ \kbar \quad \text{and}  \label{DFX1}  \\
&& \hR_v \ := \ 2^{-1/(Nn_v)} R_v \ > \ r_v \ , \label{DFX2}
\end{eqnarray}
and suppose $n$ is a multiple of $n_v$, say $n = m_v n_v$ for an appropriate 
integer $m_v$.  

The construction will be carried out in three phases.  

\medskip
\noindent{\bf Phase 1.  Patching the high-order coefficients.}
\index{coefficients $A_{v,ij}$!high-order}
\index{patching!high-order coefficients}
\index{local patching for $\CC$-simple sets!Phase 1: high-order coefficients|(}

In this phase we carry out the patching for stages $k = 1, \ldots, \kbar$.

Using the basis functions 
\index{basis!$L$-rational}
$\varphi_{i,k}$ and $\varphi_{\lambda}$ we can write 
\begin{equation} \label{DFE1}
G_v^{(0)}(z) \ = \ \phi_v(z)^n \ = \ \sum_{i=1}^{m} \sum_{j=0}^{(n-1)N_i-1}
                             A_{v,ij} \varphi_{i,nN_i-j}(z)
           + \sum_{\lambda=1}^{\Lambda} A_{v,\lambda} \varphi_{\lambda}(z) 
                \ . 
\end{equation}
Here $A_{v,i0} = \tc_{v,i}^n$, for each $i$.  
The coefficients $A_{v,ij}$ with $j = 1, \ldots, \kbar N_i - 1$ 
\index{coefficients $A_{v,ij}$!high-order}
will be deemed ``high order''.

Put $F_v(z) = \phi_v(z)^{n_v}$.   
Then $G_v^{(0)}(z) = \phi_v(z)^n = (F_v(z))^{m_v}$.
\index{patching functions, initial $G_v^{(0)}(z)$!construction of} 
We will patch the high-order coefficients of $G_v^{(0)}(z)$
\index{patching functions, initial $G_v^{(0)}(z)$!for archimedean sets $E_v$!patched by magnification}
\index{coefficients $A_{v,ij}$!high-order}
by sequentially adjusting corresponding coefficients of $F_v(z)$.  
As will be seen, a small change in the latter produces a large change 
in the former.
\index{magnification argument|ii}  
Write 
\begin{equation} \label{DFE2} 
F_v(z)  
   \ = \ \sum_{i=1}^{m}
        \sum_{j=0}^{(n_v-1)N_i-1} a_{v,ij} \varphi_{i,n_v N_i-j}(z)
    + \sum_{\lambda=1}^\Lambda a_{v,\lambda} \varphi_{\lambda}(z)
                  \ . \notag
\end{equation}  
To adjust $G_v^{(0)}(z)$,
\index{patching functions, initial $G_v^{(0)}(z)$!for archimedean sets $E_v$!patched by magnification}  
we will replace $F_v(z)$ with
\begin{equation} \label{DFE3}
\hF_v(z) \ = \ \phi_v(z)^{n_v} 
  + \sum_{i=1}^{m} \sum_{j=0}^{\kbar N_i-1} \eta_{v,ij} \varphi_{i,n_v N_i-j}(z) 
\end{equation}
for appropriately chosen $\eta_{v,ij}$.  

We will take $\eta_{v,i0} = 0$ for each $i$, since $\Delta_{v,i0}^{(1)} = 0$.    
The remaining $\eta_{v,ij}$  will be determined recursively, 
in terms of the $\Delta_{v,ij}^{(k)}$, in $\prec_N$ order\index{order!$\prec_N$};  
in particular, for each $x_i$ the $\eta_{v,ij}$ will be determined
in order of increasing $j$.  As $F_v(z)$ is changed stepwise to $\hF_v(z)$, 
then $G_v^{(0)}(z) = (F_v(z))^{m_v}$
\index{patching functions, initial $G_v^{(0)}(z)$!for archimedean sets $E_v$!patched by magnification}  
is changed stepwise to
$G_v^{(\kbar)}(z) = \hF_v(z)^{m_v}$,\index{patching functions, $G_v^{(k)}(z)$ for $1 \le k \le n$!constructed by patching} passing through     
$G_v^{(1)}(z)$, $G_v^{(2)}(z)$, $\ldots$, $G_v^{(\kbar-1)}(z)$\index{patching functions, $G_v^{(k)}(z)$} at 
intermediate steps.

It will be useful to consider what happens as each
$\eta_{v,ij}$ is varied in turn.  
Suppose $\ckF_v(z)$ is a function obtained  
at one of the intermediate steps, and at the next step $\ckF_v(z)$ is replaced 
by $\ckF_v^{\prime}(z) = \ckF_v(z) + \eta_{v,ij} \varphi_{i,n_v N_i-j}(z)$.
Let $k$ be such that $(k-1)N_i \le j < kN_i$.  
When $\ckF_v^{\prime}(z)^{m_v}$ is expanded using the binomial theorem,\index{binomial theorem} 
the result is\index{compensating functions $\vartheta_{v,ij}^{(k)}(z)$} 
\begin{equation} \label{DXYZ2}
\ckF_v^{\prime}(z)^{m_v} 
   \ = \ \ckF_v(z)^{m_v} 
         + (m_v \eta_{v,ij}/\tc_{v,i}^{n_v-k-1}) \cdot \vartheta_{v,ij}^{(k)}(z)
\end{equation}
where\index{compensating functions $\vartheta_{v,ij}^{(k)}(z)$!construction of} 
\begin{eqnarray} \label{DFRH1}
\vartheta_{v,ij}^{(k)}(z) & = & 
        \tc_{v,i}^{n_v-k-1} \Big( \varphi_{i,n_v N_i - j}(z) \ckF_v(z)^{m_v-1}  \\
   & & \quad + \sum_{t = 2}^{m_v} \frac{1}{m_v} \binom{m_v}{t} 
        \eta_{v,ij}^{t-1} 
             \varphi_{i,n_vN_i - j}(z)^t \ckF_v(z)^{m_v-t} \Big) \ . \notag
\end{eqnarray}
The first term on the right has a pole of order $nN_i - j$ at $x_i$, 
while all the other terms have poles of lower order at $x_i$.  
Likewise, for each $x_{i^{\prime}} \ne x_i$, the first term 
has a pole of order $(n-n_v)N _{i^{\prime}}$ at $x_{i^{\prime}}$ and all the other
terms have poles of lower order.

Since we have required that the $\Delta_{v,i0}^{(k)} = 0$, 
the leading coefficient of each $\ckF_v(z)$ at $x_i$ is $\tc_{v,i}^{n_v}$, 
\index{coefficients $A_{v,ij}$!leading}
the same as that of $\phi_v(z)^{n_v}$. 
It follows that each\index{compensating functions $\vartheta_{v,ij}^{(k)}(z)$!poles and leading coefficients of} 
$\vartheta_{v,ij}^{(k)}(z)$ 
has leading coefficient $\tc_{v,i}^{n-k-1}$ 
and meets the conditions of the Theorem. 

\vskip .1 in
It is essentially trivial by continuity  
(but will be rigorously proved in Lemma \ref{DLemCP1} below),
that there is an $\epsilon_v > 0$ such that if $|\eta_{v,ij}| < \epsilon_v$ 
for each $i, j$ then 
\begin{eqnarray} 
\{ z : |\hF_v(z)| \le \hR_v^{n_v N} \}
 & \subset & \{ z : |\phi_v(z)^{n_v}| \le R_v^{n_v N} \} \notag \\
 & = & \{ z  : |\phi_v(z)| \le R_v^N \}  \ \subset \ E_v^0 \ . \label{DFGA1}
\end{eqnarray} 

\vskip .1 in
The numbers $n_v$ and  $\tc_{v,i} \ne 0$ are fixed.  
Assuming the existence of such an $\epsilon_v$, let $B_v$ be the number\index{freedom $B_v$ in patching|ii} 
in the statement of Theorem \ref{DCPatch1}.   
Suppose $n$ (and hence $m_v = n/n_v$) is large enough that
\begin{equation} \label{DFE5}
\frac{B_v}{m_v} \cdot  \max(1,\max_{1 \le i \le m}(|\tc_{v,i}|))^{n_v} 
\ < \ \epsilon_v \ .
\end{equation}
For all $(i, j)$ with $1 \le j < \kbar N_i$,  
the $\Delta_{v,ij}^{(k)} \in \CC$ satisfy 
$|\Delta_{v,ij}^{(k)}| \le B_v$.  Hence, taking 
\begin{equation*}
\eta_{v,ij} \ = \ \frac{1}{m_v} \Delta_{v,ij}^{(k)} \tc_{v,i}^{n_v-k-1}
\end{equation*}
we have $|\eta_{v,ij}| < \epsilon_v$ and (\ref{DFGA1}) holds.
On the other hand, (\ref{DXYZ2}) becomes\index{compensating functions $\vartheta_{v,ij}^{(k)}(z)$}
\begin{equation*}
\ckF^{\prime}_v(z)^{m_v} \ = \ \ckF_v(z)^{m_v} + \Delta_{v,ij}^{(k)} \cdot \vartheta_{v,ij}^{(k)}(z)
\end{equation*}
for the chosen $\Delta_{v,ij}^{(k)}$.
 
In summary, small changes  
$\eta_{v,ij}$ in the coefficients of $\phi_v(z)^{m_v}$ are ``magnified'' to  
\index{coefficients $A_{v,ij}$}
\index{magnification argument|ii}
large changes $\Delta_{v,ij}^{(k)}$ 
in the coefficients of $G_v^{(0)}(z)$.
\index{patching functions, initial $G_v^{(0)}(z)$!for archimedean sets $E_v$!patched by magnification}
If these changes are carried out in $\prec_N$ order\index{order!$\prec_N$}, 
then at appropriate steps in the construction 
we obtain functions\index{compensating functions $\vartheta_{v,ij}^{(k)}(z)$!construction of}
\begin{equation*} 
G_v^{(k)}(z) \ = \ G_v^{(k-1)}(z) + \sum_{i=1}^{m} \sum_{j=(k-1)N_i}^{kN_i-1}
                           \Delta_{v,ij}^{(k)} \vartheta_{v,ij}^{(k)}(z) 
\end{equation*}                           
for $k = 1, \ldots, \kbar$.  
Note that the leading coefficients 
of the $G_v^{(k)}(z)$\index{patching functions, $G_v^{(k)}(z)$ for $1 \le k \le n$!leading coefficients of} are never changed. 
\index{coefficients $A_{v,ij}$!leading} 
Since the leading coefficient of 
$G_v^{(0)}(z)$\index{patching functions, initial $G_v^{(0)}(z)$!leading coefficients of} at $x_i$ is $\tc_{v,i}^{n}$, 
the same is true for each $G_v^{(k)}(z)$.\index{patching functions, $G_v^{(k)}(z)$ for $1 \le k \le n$!leading coefficients of}  
Similarly, the leading coefficient of $\hF_v(z)$ 
at $x_i$ is $\tc_{v,i}^{n_v}$.  
 
\vskip .1 in
To rigorously prove the existence of an $\epsilon_v$ for which (\ref{DFGA1}) holds, 
we must use some information about the $\varphi_{ij}(z)$.

Given a function $F(z) \in \CC(\cC)$ and a number $R > 0$, write
\begin{eqnarray*}
W_R  & = & \{ z \in \cC_v(\CC) : |F(z)| \le R \} \ , \\
V_R  & = & \{ z \in \cC_v(\CC) : |F(z)| \ge R \}  
\end{eqnarray*} 
(regarding the poles of $F(z)$ as belonging to $V_R$), and put
\begin{equation*}
\Gamma_R  \ = \ \{ z \in \cC_v(\CC) : |F(z)| = R \} \ .
\end{equation*}
Then $\Gamma_R$ is the common boundary of $W_R$ and $V_R$.

\begin{lemma} \label{DLemCP1}
Let $F(z) \in \CC(\cC_v)$ have polar divisor $\div(F)_{\infty}$.
Suppose that $H(z) \in \CC(\cC_v)$
has polar divisor $\div(H)_{\infty} \le \div(F)_{\infty}$,  
and for some $\delta < 1$ we have $|H(z)| < \delta \cdot R$ on $\Gamma_R$.
Then
\begin{equation*}
\{ z \in \cC_v(\CC) : |F(z) + H(z)| \le (1-\delta) \cdot R \}
                     \ \subset \ W_R \ .
\end{equation*}
\end{lemma}

\begin{proof}
Consider $G(z) = H(z)/F(z)$.  On $\Gamma_R$ we have $|G(z)| < \delta$.
The hypothesis on the poles implies that $G(z)$ extends
to a function holomorphic in $V_R$.  By the Maximum Modulus Principle 
\index{Maximum principle!for holomorphic functions}
$|G(z)| < \delta$ on $V_R$, so $|H(z)| < \delta \cdot |F(z)|$ on $V_R$.
It follows that $|F(z) + H(z)| > (1-\delta) |F(z)| \ge (1-\delta)R$
on $V_R$, so
$\{ z \in \cC_v(\CC) : |F(z) + H(z)| \le (1-\delta) \cdot R \} \subset  W_R$.
\end{proof}

\vskip .1 in
To obtain (\ref{DFGA1}), put 
$\hM_v = \max_{1 \le i \le m} \Big( \max{1 \le j \le \kbar N_i} \|\varphi_{i,n_v N_i-j}\|_{E_v} \Big)$  
and let $\epsilon_v >0$ be small enough that
\begin{equation*}
\epsilon_v \cdot \kbar N \hM_v \ < \ R_v^{n_v N} - \hR_v^{n_v N} \ .
\end{equation*}     
Apply Lemma \ref{DLemCP1} with $F(z) = \phi_v(z)^{n_v}$ and $R = R_v^{n_v N}$, 
taking $\delta = 1 - (\hR_v/R_v)^{n_v N}$.  
By hypothesis, we have $W_R \subset E_v^0$. 
Take 
\begin{equation*}
H(z) \ = \ \sum_{i=1}^{m} \sum_{j=1}^{\kbar N_i-1} 
                        \eta_{v,ij} \varphi_{i,n_vN_i-j}(z)
\end{equation*}
where $|\eta_{v,ij}| \le \epsilon_v$ for each $(i,j)$.  Then
$|H(z)| \le \epsilon_v \cdot \kbar N \hM_v < \delta R$
on $\Gamma_R = \{ z \in \cC_v(z) : |\phi_v(z)^{n_v}| = R_v^{n_vN} \}$,
and $\hF_v(z) = F(z) + H(z)$, while $\hR_v^{n_v N} = (1-\delta) R$, 
so (\ref{DFGA1}) follows from the Lemma.  
                       
\vskip .1 in
Let $\hGamma_v$ denote the level curve
$\{ z : |\hF_v(z)| = \hR_v^{n_vN} \}$.
By (\ref{DFGA1}), 
\begin{equation} \label{DFGA2}
\hGamma_v \ \subset \ \{ z \in \cC_v(\CC) : |\phi_v(z)| \le R_v^N \} \ \subset \ E_v^0 \ .
\end{equation}  
The function $\hF_v(z)$ and the curve $\hGamma_v$ will play a
key role in the rest of the construction.\index{local patching for $\CC$-simple sets!Phase 1: high-order coefficients|)}
 
\vskip .1 in
\noindent{\bf Phase 2.  Patching the middle coefficients.}
\index{patching!middle coefficients}
\index{coefficients $A_{v,ij}$!middle}
\index{local patching for $\CC$-simple sets!Phase 2: middle coefficients|(}

In this phase we carry out the patching process for $k = \kbar+1, \ldots, n-1$.
For each $k$ we begin with a function $G_v^{(k-1)}(z)$,\index{patching functions, $G_v^{(k)}(z)$ for $1 \le k \le n$!constructed by patching} 
and we modify the
coefficients with $(k-1)N_i \le j < kN_i$, for each $i$. 
For each such  $j$ we can uniquely write
\begin{equation*}
nN_i -j \ = \ r_{ij} + (n-k-1)N_i, \quad 
\text{with $N_i < r_{ij} \le 2N_i$.}
\end{equation*}  
We can then write
\begin{equation*}
n-k-1 \ = \ \ell_1 + \ell_2 n_v, \quad 
            \text{with \ $0 \le \ell_1 < n_v$, \  $0 \le \ell_2 < m_v$,}
\end{equation*}
so $nN_i-j = r_{ij} + \ell_1 N_i + \ell_2 n_v N_i$.  
ut \index{compensating functions $\vartheta_{v,ij}^{(k)}(z)$!construction of}
\begin{equation*}
\vartheta_{v,ij}^{(k)}(z) \ = \ 
 \varphi_{i,r_{ij}}(z) \phi_v(z)^{\ell_1} \hF_v(z)^{\ell_2} \ .
\end{equation*}
Then $\vartheta_{v,ij}^{(k)}(z)$
\index{compensating functions $\vartheta_{v,ij}^{(k)}(z)$!poles and leading coefficients of} 
has a pole of exact order $nN_i - j$ at $x_i$, 
with leading coefficient $\tc_{v,i}^{n-k-1}$. 
\index{coefficients $A_{v,ij}$!leading} 
Its poles at the $x_{i^{\prime}} \ne x_i$
are of order at most $(n-k-1)N_{i^{\prime}}$, so it meets the conditions
of the theorem.

Modifying the coefficients stepwise in $\prec_N$ order\index{order!$\prec_N$}, we put 
\index{coefficients $A_{v,ij}$!middle}\index{compensating functions $\vartheta_{v,ij}^{(k)}(z)$}
\begin{equation} \label{DFGG1} 
G_v^{(k)}(z) \ = \ G_v^{(k-1)}(z) + \sum_{i=1}^{m} \sum_{j=(k-1)N_i + 1}^{kN_i}
                           \Delta_{v,ij}^{(k)} \vartheta_{v,ij}^{(k)}(z) 
\end{equation}    
where $|\Delta_{v,ij}^{(k)}| \le h_v^{kN}$.

We now seek a bound for $|\vartheta_{v,ij}^{(k)}(z)|$
\index{compensating functions $\vartheta_{v,ij}^{(k)}(z)$!bounds for} 
on the level curve $\hGamma_v$.
By definition, $|\hF_v(z)^{\ell_2}| = \hR_v^{Nn_v\ell_2}$ on $\hGamma_v$.
Since $\hR_v = 2^{-1/(Nn_v)}R_v$ it follows from (\ref{DFGA2}) that  on $\hGamma_v$, 
for $0 \le \ell_1 < n_v$,  
\begin{equation} \label{DFGA3}
|\phi_v(z)^{\ell_1}| \ \le \ 2 \hR_v^{N\ell_1}
\end{equation} 
Finally, 
$| \varphi_{i,r_{ij}}(z)| \le M_v$ on $\hGamma_v$ for all
$N_i + 1 \le r_{ij} \le 2N_i$ (in fact this holds for all $z \in E_v$).  
Hence\index{compensating functions $\vartheta_{v,ij}^{(k)}(z)$!bounds for}
\begin{equation*}
|\vartheta_{v,ij}^{(k)}(z)| \ \le \ 
| \varphi_{i,r_{ij}}| \cdot |\phi_v(z)^{\ell_1}| \cdot |\hF_v(z)^{\ell_2}| 
\ \le \ 2 M_v \hR_v^{N(n-k-1)} \ . 
\end{equation*}
Since there are $N$ terms in the sum (\ref{DFGG1}), on $\hGamma_v$
\begin{equation} \label{DFXY1}
|G_v^{(k)}(z)| \ \le \ |G_v^{(k-1)}(z)| +
                              N h_v^{kN} \cdot 2 M_v \hR_v^{N(n-k-1)} \ .
\end{equation}
\index{local patching for $\CC$-simple sets!Phase 2: middle coefficients|)}

\vskip .1 in                               
\noindent{\bf Phase 3.  Patching the low-order coefficients.}
\index{patching!low-order coefficients}
\index{coefficients $A_{v,ij}$!low-order}
\index{local patching for $\CC$-simple sets!Phase 3: low-order coefficients|(}

In the final step we take 
\begin{equation*}
G_v^{(n)}(z) \ = \ G_v^{(n-1)}(z)
          + \sum_{\lambda = 1}^{\Lambda} \Delta_{v,\lambda}^{(n)} \varphi_{\lambda}
\end{equation*}
where $|\Delta_{v,\lambda}^{(n)}| \le h_v^{nN}$ for each $\lambda$.  

Since $\Lambda \le N$, and each $|\varphi_{\lambda}(z)| \le M_v$ on $\hGamma_v$
(indeed on all of $E_v$), on $\hGamma_v$
\begin{equation} \label{DFXY2}
|G_v^{(n)}(z)| \ \le \ |G_v^{(n-1)}(z)| + N M_v h_v^{nN} \ .
\end{equation}
\index{local patching for $\CC$-simple sets!Phase 3: low-order coefficients|)}

\medskip
To complete the proof, we must show that if $n$ is sufficiently large then 
part (B) of Theorem \ref{DCPatch1} holds. 
Assume that $n$ is large enough that
\begin{equation} \label{DFNN1}
N M_v \left(\frac{h_v}{r_v}\right)^{nN} \ < \ \frac{1}{4} \ ,
\end{equation}
and recall that $\kbar$ satisfies 
\begin{equation} \label{DFX3} 
\frac{2 N M_v}{1 - (h_v/r_v)^N} 
                \cdot \Big(\frac{h_v}{r_v}\Big)^{\kbar N}  \ < \ \frac{1}{4} \ .
\end{equation}

Consider the total change on $\hGamma_v$ in passing from
$G_v^{(\kbar)}(z) = \hF_v(z)^{m_v}$ to $G_v^{(n)}(z)$.\index{patching functions, $G_v^{(k)}(z)$ for $1 \le k \le n$!constructed by patching} 
By (\ref{DFXY1}) and (\ref{DFXY2}), for each $z \in \hGamma_v$,
\begin{equation} \label{DFZQ1}
|G_v^{(k)}(z) - \hF_v(z)^{m_v}| \ \le \ N M_v h_v^{nN} +
            \frac{2 N M_v}{\hR_v^N} \cdot \hR_v^{nN} 
              \cdot \sum_{k= \kbar}^{n-1} \frac{h_v^{kN}}{\hR_v^{kN}} 
\end{equation}
Since $\hR_v > r_v > 1$, 
by inserting (\ref{DFNN1}), and (\ref{DFX3}) in (\ref{DFZQ1}), 
we find that on $\hGamma_v$  
\begin{equation} \label{DFX4} 
|G_v^{(n)}(z) - \hF_v(z)^{m_v}| \ < \ \frac{1}{2} \hR_v^{nN} \ .
\end{equation}
As $|\hF_v(z)^{m_v}| = \hR_v^{nN}$ on $\hGamma_v$, by applying
Lemma \ref{DLemCP1} with $F(z) = \hF_v(z)^{m_v}$ and
$H(z) = G_v^{(n)}(z) - \hF_v(z)^{m_v}$, taking $\delta = \frac{1}{2}$,
we see that 
\begin{equation*} 
\{ z \in \cC_v(\CC) : |G_v^{(n)}(z)| \le \frac{1}{2} \hR_v^{nN} \}
 \ \subset \ \{ z \in \cC_v(\CC) : |\hF_v(z)^{m_v}| \le  \hR_v^{nN} \}
\end{equation*}
which is contained in $E_v^0$.  

Finally, if $n$ is also large enough that 
\begin{equation} \label{DFNN2}
\frac{1}{2} \hR_v^{nN} \ > \ r_v^{nN} \ ,
\end{equation} 
then 
\begin{equation}
\{ z \in \cC_v(\CC) : |G_v^{(n)}(z)| \le r_v^{nN} \} \ \subset \ E_v^0 \ .
\end{equation}
A similar argument shows that\index{patching functions, $G_v^{(k)}(z)$ for $1 \le k \le n$!roots are confined to $E_v$}  
$\{ z \in \cC_v(\CC) : |G_v^{(k)}(z)| \le r_v^{nN} \} \subset E_v^0$ for each $k = 1, \ldots, n$.

In summary if $n > \kbar$ is divisible by $n_v$ and large enough that conditions
(\ref{DFE5}), (\ref{DFNN1}) and (\ref{DFNN2}) hold, the construction succeeds.
\end{proof}
\index{patching construction!for the case when $K_v \cong \CC$|)}
\index{patching argument!local!for the case when $K_v \cong \CC$|)}

%% file: NewFSZChap9.tex
\chapter{ The Local Patching Construction when $K_v \cong \RR$ } 
\label{Chap9}

\vskip .1 in
In this section we give the confinement argument\index{confinement argument} for Theorem \ref{aT1-B} when  
$K_v \cong \RR$. \index{patching construction!for the case when $K_v \cong \RR$|(}
\index{patching argument!local!for the case when $K_v \cong \RR$|(} 
Write $\CC_v$ for $\CC$ and $| \ |_v$ for $| \ |$.  
Let $w_v$ be the distinguished place\index{distinguished place $w_v$} 
of $L = K(\fX)$ determined by the embedding $\tK \hookrightarrow \CC_v$
used to identify $\fX$ with a subset of $\cC_v(\CC_v)$. 
Identify $K_v$ with $\RR$, and $L_{w_v}$ with $\RR$ or $\CC$ as appropriate. 

Following the construction of the coherent approximating functions in Theorem \ref{CTCX2},
\index{coherent approximating functions $\phi_v(z)$} 
we begin with the following data:   
\begin{enumerate}
  \item A $K_v$-symmetric probability vector $\vs \in \cP^{m}(\QQ)$ 
\index{$K_v$-symmetric!probability vector}
           with positive rational coefficients.
  \item An $\RR$-simple set $E_v$:  in particular $E_v$ is nonempty and compact,
            stable under complex conjugation, and is a union of finitely 
            many pairwise disjoint, nonempty compact sets $E_{v,1}, \ldots, E_{v,\ell}$
            such that each $E_{v,i}$ is either 

            (a) a closed interval of positive length 
            contained in $\cC_v(\RR)$, or 
            
            (b) is disjoint from $\cC_v(\RR)$, simply connected,\index{simply connected} 
            has a piecewise smooth boundary,\index{boundary!piecewise smooth} 
             and is the closure of its $\cC_v(\CC)$-interior.\index{closure of $\cC_v(\CC)$ interior}
           
   \item Let $E_v^0$ be the quasi-interior of $E_v$,
\index{quasi-interior} 
           the union of the real interiors of the components $E_{v,i} \subset \cC_v(\RR)$  
           and the complex interiors of the components $E_{v,i} \subset \cC_v(\CC) \backslash \cC_v(\RR)$;  
           then we are given a $\cC_v(\CC)$-open set $U_v$ such that 

             (a) $U_v \cap E_v = E_v^0$, 
             
             (b) the components of $U_v$ are simply connected, and 

             (c) the closure $\Ubar_v$ is disjoint from $\fX$.  
           
  \item Parameters $h_v$, $r_v$, $R_v$, with $1 < h_v < r_v < R_v$,
          which govern the freedom in the patching process.
          \index{patching parameters}
          
  \item A number $N$  and an $(\fX,\vs)$-function $\phi_v(z) \in K_v(\cC)$ 
           of degree $N$ whose zeros all belong to $E_v^0$, and which has the following properties:
           
           (a) $\phi_v^{-1}(D(0,2R_v^{N})) \subset U_v$.   
           
           (b) For each component $E_{v,i} \subset \cC_v(\RR)$, if $\phi_v(z)$ has $\tau_j$ 
               zeros in $E_{v,j}$, then $\phi_v(z)$ oscillates $\tau_j$ times between $\pm 2 R_v^N$ 
               on $U_v \cap E_{v,j}$.  
        
  \item Put $N_i = Ns_i$ for each $i$, and write
           $\tc_{v,i} = \lim_{z \rightarrow x_i} \phi_v(z) \cdot g_{x_i}(z)^{N_i}$
           for the leading coefficient of $\phi_v(z)$ at $x_i$;   
   \index{coefficients $A_{v,ij}$!leading}
           then we are given an order $\prec_N$\index{order!$\prec_N$} on the index set 
              $\cI = \{(i,j) \in \ZZ^2 : 1 \le i \le m, 0 \le j\}$ 
           determined by $N$ and $\vs$ as in (\ref{FPrec}), 
           which gives the sequence in which coefficients are patched.
\index{coefficients $A_{v,ij}$}   
           We will use the $L$-rational basis  
\index{basis!$L$-rational}
           $\{\varphi_{ij}, \varphi_{\lambda}\}$ 
           from \S\ref{Chap3}.\ref{LRationalBasisSection} to expand all functions, and 
           $\Lambda = \dim_K(\Gamma(\sum_{i=1}^m N_i(x_i)))$ 
           will be the number of low-order basis elements, 
           as in \S\ref{Chap7}.\ref{Char0Section}.  
           The order $\prec_N$\index{order!$\prec_N$} respects the $N$-bands (\ref{FBand}),\index{band!$\Band_N(k)$} 
           and for each $x_i \in \fX$,  
           specifies the terms to be patched in decreasing pole order. 
\end{enumerate} 

\smallskip
Let $D(0,R) = \{z \in \CC : |z| \le R\}$ be the filled disc, 
and let $E(a,b) = \{x+iy \in \CC : x^2/a^2 + y^2/b^2 \le 1\}$ 
be the filled ellipse.\index{filled ellipse|ii}\label{`SymbolIndexFilledEllipse'}
Write $C(0,R) = \partial D(0,R)$ and $\partial E(a,b)$ for their boundaries, 
and note that if $a > b$ then $D(0,b) \subset E(a,b) \subset D(0,a)$.   

\smallskip
Let $T_n(z)$ be the Chebyshev polynomial of degree $n$ for the interval
\index{Chebyshev polynomial|ii} 
\label{`SymbolIndexChebPoly'} 
$[-2,2]$, defined by $T_n(2 \cos(\theta)) \ = \ 2 \cos(n \theta)$.
Equivalently, $T_n(z)$ is the unique polynomial of degree $n$ for which $T_n(z+1/z) = z^n + 1/z^n$.  
For each $R > 0$, let $T_{n,R}(z) = R^n T_n(z/R)$ 
\label{`SymbolIndexChebPolyR'}
be the Chebyshev polynomial for the interval $[-2R,2R]$.
Then $T_{n,R}(z)$ is monic of degree $n$ with coefficients in $\RR$, 
and as noted in (\cite{Rob1}), 
\begin{equation} \label{DFSW2}
T_{n,R}(z) 
           \ = \  \ z^n + \sum_{k=1}^{\lfloor n/2 \rfloor} 
                           \frac{n}{k} \binom{n-k-1}{k-1} R^{2k} z^{n-2k} \ .
\end{equation}
Furthermore, $T_{n,R}(z)$ has the following mapping properties:\index{Chebyshev polynomial!mapping properties|ii}
  
First, $T_{n,R}([-2R,2R]) = [-2R^n,2R^n]$ and $T_{n,R}^{-1}([-2R^n,2R^n]) = [-2R,2R]$.  
These facts follow from the identity $T_n(2 \cos(\theta)) = 2 \cos(n \theta)$, 
which means that $T_{n,R}$ oscillates $n$ times between $\pm 2R^n$ on $[-2R,2R]$.  
Second, for each $t > R$, 
\begin{eqnarray*} 
T_{n,R}\big(E(t + \frac{R^2}{t}, t - \frac{R^2}{t})\big) 
\ = \ E(t^n + \frac{R^{2n}}{t^n}, t^n - \frac{R^{2n}}{t^n})  \qquad \\
\text{and} \quad \big(T_{n,R}\big)^{-1}\big(E(t^n + \frac{R^{2n}}{t^n}, t^n - \frac{R^{2n}}{t^n})\big) 
\ = \ E(t + \frac{R^2}{t}, t - \frac{R^2}{t}) \ .
\end{eqnarray*}
Indeed, $T_{n,R}$ gives an $n$-to-$1$ map from $E(t + R^2/t,t-R^2/t)$
onto $E(t^n + R^{2n}/t^n,t^n-R^{2n}/t^n)$ (counting multiplicities).   
This follows from the definition of $T_{n,R}$ and the commutativity of the diagram  
\begin{equation*}
\begin{array}{ccc}
 C(0,t/R) & \stackrel{z^n}{\longrightarrow} & C(0,t^n/R^n) \\
 \ \ \downarrow \ z + \frac{1}{z} & & \ \ \downarrow \ z + \frac{1}{z}  \\
      \partial E(t/R+R/t, t/R-R/t) & \stackrel{T_n(z)}{\longrightarrow} &
      \partial E(t^n/R^n+R^n/t^n, t^n/R^n-R^n/t^n) \ .
 \end{array}
\end{equation*}

\vskip .1 in
\begin{theorem}  \label{DRPatch1}\index{patching theorem!for the case when $K_v \cong \RR$}
Suppose $K_v \cong \RR$. 
Let $E_v \subset \cC_v(\RR) \backslash \fX$ be a $K_v$-simple set. 
\index{$K_v$-simple!set}
Let $U_v$ be an open set in $\cC_v(\CC_v)$ with $U_v \cap E_v = E_v^0$,  
whose components are simply connected, and whose closure $\Ubar_v$ is disjoint from $\fX$.  
Let $\vs \in \cP^m(\QQ)$ be a $K_v$-symmetric probability
\index{$K_v$-symmetric!probability vector}
vector with positive rational coefficients, and let $1 < h_v < r_v < R_v$ be numbers. 

Let $\phi_v(z) \in K_v(\cC)$ be an $(\fX,\vs)$-function of degree $N$ 
whose zeros belong to $E_v^0$, satisfying  

         \quad $(1)$ $\varphi^{-1}(D(0,2R_v^{N})) \subset U_v$,   
           
         \quad $(2)$ For each component $E_{v,i}$ contained in $\cC_v(\RR)$, if $\phi_v(z)$ has $\tau_i$ 
               zeros in $E_{v,i}$, then $\phi_v(z)$ oscillates $\tau_i$ times between $\pm 2 R_v^N$ 
               on $E_{v,i}$. 
               
Let $\tc_{v,i}$ be the leading coefficient of $\phi_v(z)$ at $x_i$, and put  
\index{coefficients $A_{v,ij}$!leading}
\begin{equation*}
M_v \ = \ \max( 
     \max_{\substack{ 1 \le i \le m \\ N_i < j \le 2N_i}} \|\varphi_{ij}\|_{\Ubar_v },
    \max_{1 \le \lambda \le \Lambda} \|\varphi_{\lambda}\|_{\Ubar_v } ) \ .
\end{equation*} 
Let  $k_v > 0$  be the least integer such that
\begin{equation} \label{DRPFF0}
\frac{16 N M_v}{1 - (h_v/r_v)^N}
                   \cdot \big(\frac{h_v}{r_v}\big)^{k_v N}  \ < \ \frac{1}{4} \ , 
\end{equation}
and let $\kbar \ge k_v$ be a fixed integer.  Let $B_v > 0$  be an arbitrary constant.  
Then there is an integer $n_v$, depending on $\phi_v(z)$, $\kbar$, $B_v$, $r_v$, and $R_v$, 
such that for each sufficiently large integer $n$ divisible by $n_v$, 
one can carry out the local patching process at $K_v$ as follows: 

Write $n = m_v n_v$. For suitable $\hR_1, \hR_2$ with $r_v < \hR_2 < \hR_1 < R_v$, 
put\index{patching functions, initial $G_v^{(0)}(z)$!construction of}
\begin{equation*}
G_v^{(0)}(z) \ = \ T_{m_v,\hR_2^{n_v N}}(T_{n_v,\hR_1^N}(\phi_v(z))) \ .
\end{equation*}    
For each $k$, $1 \le k < n$, 
let $\{\Delta_{v,ij}^{(k)} \in \CC_v\}_{(i,j) \in \Band_N(k)}$ 
be an arbitrary $K_v$-symmetric set of numbers given recursively in $\prec_N$ order\index{order!$\prec_N$},
\index{$K_v$-symmetric!set of numbers}\index{band!$\Band_N(k)$} 
subject to the conditions that for each $i$, we have $\Delta_{v,i0}^{(1)} = 0$ and for each $j > 0$  
\begin{equation} \label{DFRCV1}
|\Delta_{v,ij}^{(k)}|_v \ \le \ \left\{
      \begin{array}{ll} B_v & \text{if \ $k \le \kbar \ ,$} \\
                         h_v^{kN} & \text{if \ $k > \kbar \ .$}
      \end{array} \right.
\end{equation}
For $k = n$, let 
$\{\Delta_{v,\lambda}^{(n)} \in \CC_v\}_{1 \le \lambda \le \Lambda}$ 
be an arbitrary $K_v$-symmetric set of numbers satisfying  
\index{$K_v$-symmetric!set of numbers}
\begin{equation} \label{DFRF2}
|\Delta_{v,\lambda}^{(n)}|_v \ \le \ h_v^{nN} \ .  
\end{equation} 

\vskip .05 in
\noindent{Then} one can inductively construct $(\fX,\vs)$-functions 
$G_v^{(1)}(z), \ldots, G_v^{(n)}(z)$ in $K_v(\cC)$,\index{patching functions, $G_v^{(k)}(z)$ for $1 \le k \le n$!constructed by patching} 
of common degree $nN$, having the following properties:  

\vskip .05 in
\noindent{$(A)$} For each $k = 1, \ldots, n$, there are functions 
$\vartheta_{v,ij}^{(k)}(z) \in L_{w_v}(\cC)$,\index{distinguished place $w_v$}
\index{compensating functions $\vartheta_{v,ij}^{(k)}(z)$} 
determined recursively in $\prec_N$ order\index{order!$\prec_N$}, 
such that\index{patching functions, $G_v^{(k)}(z)$ for $1 \le k \le n$!constructed by patching} 
\begin{eqnarray*} 
G_v^{(k)}(z) & = & G_v^{(k-1)}(z)
              +  \sum_{(i,j) \in \Band_N(k)}
                           \Delta_{v,ij}^{(k)} \vartheta_{v,ij}^{(k)}(z)  
        \quad \text{for $k < n \ ,$} \\
G_v^{(n)}(z) & = & G_v^{(n-1)}(z) + \sum_{\lambda=1}^{\Lambda}
                           \Delta_{v,\lambda}^{(n)} \varphi_{\lambda}(z) \ ,                                                     
\end{eqnarray*}
and where for each $(i,j)$,\index{band!$\Band_N(k)$} 

\quad $(1)$ $\vartheta_{v,ij}^{(k)}(z)$
\index{compensating functions $\vartheta_{v,ij}^{(k)}(z)$!poles and leading coefficients of} 
has a pole of order $nN_i-j > (n-k-1)N_i$ at $x_i$ 
and leading coefficient $\tc_{v,i}^{n-k-1}$,  
\index{coefficients $A_{v,ij}$!leading}
a pole of order at most $(n-k-1)N_{i^{\prime}}$ at each $x_{i^{\prime}} \ne x_i$, 
and no other poles;  

\quad $(2)$ $\sum_{(i^{\prime},j) \in \Aut_c(\CC_v/K_v)(i,j)} 
  \Delta_{v,i^{\prime}j}^{(k)} \vartheta_{v,i^{\prime}j}^{(k)}(z)$
belongs to $K_v(\cC);$\index{compensating functions $\vartheta_{v,ij}^{(k)}(z)$!are $K_v$-symmetric}  

\vskip .05 in
\noindent{$(B)$} For each $k = 1, \ldots, n$, 

\quad $(1)$ the zeros of $G_v^{(k)}(z)$\index{patching functions, $G_v^{(k)}(z)$ for $1 \le k \le n$!roots are confined to $E_v$}  
all belong to $E_v^0$, and for each component $E_{v,i}$
of $E_v$, if $\phi_v(z)$ has $\tau_i$ zeros in $E_{v,i}$, then $G_v^{(k)}(z)$ has $T_i = n \tau_i$ zeros
in $E_{v,i}$.  

\quad $(2)$ $\{ z \in \cC_v(\CC_v) : |G_v^{(k)}(z)|_v \le 2r_v^{nN}\} \ \subset \ U_v$, and 

\quad $(3)$ for each component $E_{v,i}$ contained in $\cC_v(\RR)$, $G_v^{(k)}(z)$
\index{patching functions, $G_v^{(k)}(z)$ for $1 \le k \le n$!for archimedean sets $E_v$!oscillate on real components of $E_v$}  
oscillates $T_i$ times between $\pm 2 r_v^{nN}$ on  $E_{v,i}$.  
\end{theorem}

\vskip .1 in
\noindent{\bf Remark.}  As in the patching construction when $K_v \cong \CC$,
a key feature of Theorem \ref{DRPatch1} is that by choosing $n$ appropriately, 
the freedom $B_v$ in patching the coefficients for $k \le \kbar$ can be made arbitrarily large.
\index{patching argument!freedom $B_v$ in patching}   
\index{coefficients $A_{v,ij}$}
Again this is accomplished by using `magnification'.\index{magnification argument} The degree of $\phi_v(z)$
is raised by a two-stage composition with Chebyshev polynomials.
\index{Chebyshev polynomial} 

The argument confining the roots 
of the $G_v^{(k)}(z)$\index{patching functions, $G_v^{(k)}(z)$ for $1 \le k \le n$!roots are confined to $E_v$}  
to $E_v$ has two parts. 
One part, which goes back to Fekete and Szeg\"o (\cite{F-SZ}) 
\index{Fekete, Michael}
\index{Szeg\"o, G\'abor} 
and uses the Maximum Modulus principle, confines the roots to $U_v$
\index{Maximum principle!for holomorphic functions}
and shows that the number of roots in each component of $U_v$ is preserved.
Since $U_v \cap E_v = E_v^0$, if $E_{v,i}$ is a component of $E_v$ which is disjoint from $\cC_v(\RR)$, 
this means that roots in $E_{v,i}^0$ must remain there. 
The other part, which goes back to Robinson (\cite{Rob1}),  
\index{Robinson, Raphael} 
is based on oscillation properties of Chebyshev polynomials and the intermediate value theorem. 
\index{Chebyshev polynomial} 
It shows that the number of roots in each component $E_{v,i}$ contained in $\cC_v(\RR)$ is preserved.  
The mapping properties of Chebyshev polynomials discussed before the statement of the Theorem  
\index{Chebyshev polynomial} 
enable to us carry out both confinement arguments simultaneously.

\vskip .1 in

\begin{proof}[Proof of Theorem \ref{DRPatch1}]  
Let $\kbar \ge k_v$ and $B_v > 0$ be as in the Theorem.  
Choose $n_v \in \NN$ large enough that $n_v > \kbar$ and $8^{-1/(n_v N)} R_v  >  r_v$.  
Put $\hR_1 = 2^{-1/(n_v N)} R_v$ and $\hR_2 = 8^{-1/(n_v N)} R_v$, 
so that $r_v < \hR_2 < \hR_1 < R_v$ 
and $2\hR_1^{n_v N} = R_v^{n_vN}$, $4\hR_2^{n_v N} = \hR_1^{n_vN}$.  
Set 
\begin{equation*}
F_v(z) \ = \ T_{n_v,\hR_1^N}(\phi_v(z)) \ .  
\end{equation*}
Let $n$ be a multiple of $n_v$, write $n = m_v n_v$, 
and put\index{patching functions, initial $G_v^{(0)}(z)$!construction of}  
\begin{equation*}
G_v^{(0)}(z) \ = \ T_{m_v,\hR_2^{n_v N}}\big(T_{n_v,\hR_1^N}(\phi_v(z))\big) 
\ = \ T_{m_v,\hR_2^{n_v N}}\big(F_v(z)\big) \ .
\end{equation*} 

We will begin by investigating the 
mapping properties of  $F_v(z)$ and 
$G_v^{(0)}(z)$.\index{patching functions, initial $G_v^{(0)}(z)$!mapping properties of}  
We first show that all the zeros of $F_v(z)$ belong to $E_v^0$, 
and that for each component $E_{v,i}$ of $E_v$, if $\phi_v$ has $\tau_i$ zeros in $E_{v,i}$,
then $F_v(z)$ has $n_v \tau_i$ zeros in $E_{v,i}$ (counted with multiplicities).  

Let $t_1$ be the largest real root of 
\begin{equation*}
t_1^{n_v N} + \frac{\hR_1^{2 n_v N}}{t_1^{n_v N}} \ = \ 4\hR_1^{n_v N} \ , 
\end{equation*} 
so that  $t_1^{n_v N} = (2 + \sqrt{3}) \hR_1^{n_v N} = ((2+\sqrt{3})/2) R_v^{n_v N}$.  Put 
\begin{eqnarray*}
a_1  & = & t_1^{N} + \frac{\hR_1^{2N}}{t_1^{N}} \ , 
\qquad  A_1 \ = \ t_1^{n_v N} + \frac{\hR_1^{2n_v N}}{t_1^{n_v N}}  
                 \ = \ 4 \hR_1^{n_v N} \ = \ 2 R_v^{n_v N} ,  \\
b_1 & = & t_1^{N} - \frac{\hR_1^{2N}}{t_1^{N}}  \ , 
\qquad  B_1 \ = \ t_1^{n_v N} - \frac{\hR_1^{2n_v N}}{t_1^{n_v N}} \ = \ 2 \sqrt{3} \hR_1^{n_v N}   
                 \ = \ \sqrt{3} R_v^{n_v N} \ ;
\end{eqnarray*}
then 
\begin{equation*}
T_{n_v,\hR_1^N}(E(a_1,b_1)) \ = \ E(A_1,B_1) \ .
\end{equation*} 
Here $a_1 = g(1/n_v) \cdot R_v^N$ where $g(x) = ((2 + \sqrt{3})/2)^x + (2(2+\sqrt{3}))^{-x}$. 
Using Calculus, one sees that $g(x) < 2$ for $0 < x < 1$, so $a_1 < 2 R_v^N$.  It follows that  
\begin{eqnarray} 
& & \qquad \qquad E(a_1,b_1) \ \subset \ D(0,2R_v^N) \ ,  \label{FContainment1} \\
& & D(0,\sqrt{3} R_v^{n_v N}) \ \subset \ E(A_1,B_1) \ \subset \ D(0,2 R_v^{n_v N}) \ . 
    \label{FContainment2}
\end{eqnarray} 
As $2\hR_1^{n_v N} = R_v^{n_v N}$,  (\ref{FContainment2}) shows that
\begin{equation} \label{FContainment4} 
D(0,2\hR_1^{n_v N}) \ \subset \ E(A_1,B_1) \ .
\end{equation} 

Since $\phi_v^{-1}(D(0,2R_v^N)) \subset U_v$
and $T_{n_v,\hR_1^N}^{-1}(E(A_1,B_1)) = E(a_1,b_1)$, (\ref{FContainment1}) gives 
\begin{equation*}
F_v^{-1}(E(A_1,B_1)) \ \subset \ U_v \ .
\end{equation*} 
For each component $E_{v,i}$ of $E_v$ contained in $\cC_v(\RR)$, 
the function $\phi_v$ is real-valued and oscillates $\tau_i$ times between $\pm 2 R_v^N$ on $E_{v,i}$. 
Since $[-2 \hR_1^N, 2 \hR_1^N] \subset [-2 R_v^N, 2 R_v^N]$ and  $T_{n_v,\hR_1^N}$ 
oscillates $n_v$ times between $\pm 2 \hR_1^{n_v N}$ on $[-2\hR_1^N,2\hR_1^N]$,  
it follows that $F_v(z)$ oscillates $n_v \tau_i$ times between $\pm 2 \hR_1^{n_v N}$ 
on $E_{v,i}$.  

If $E_{v,i}$ is a component of $E_v$ disjoint from $\cC_v(\RR)$, 
then $U_v \cap E_{v,i} = E_{v,i}^0$.  
Since $T_{n_v,\hR_1^N}$ has $n_v$ zeros in $E(a_1,b_1)$,
and $E_{v,i}$ is simply connected with a piecewise smooth boundary,\index{boundary!piecewise smooth}
the Argument Principle\index{Argument Principle} shows that $F_v(z)$ has $n_v \tau_i$ zeros in $E_{v,i}^0$.  
On the other hand, if $E_{v,i}$ is a component contained in $\cC_v(\RR)$, 
then by the discussion above $F_v(z)$ has at least $n_v \tau_i$ zeros in $E_{v,i}^0$. 
Since $\sum_i n_v \tau_i= n_v N$ and $F_v(z)$ has degree $n_v N$, these zeros   
account for all the zeros of $F_v(z)$.  
Thus all the zeros of $F_v(z)$ belong to $E_v^0$, 
and $F_v(z)$ has exactly $n_v \tau_i$ zeros in each $E_{v,i}$.

\smallskip
Next, let $t_2$ be the largest real root of 
\begin{equation*}
t_2^{n_v N} + \frac{\hR_2^{2 n_v N}}{t_2^{n_v N}} \ = 4 \hR_2^{n_v N} \ , 
\end{equation*} 
so that $t_2^{n_v N} = (2 + \sqrt{3}) \hR_2^{n_v N}$.  If we put  
\begin{eqnarray*}
a_2 & = & t_2^{n_v N} + \frac{\hR_2^{2 n_v N}}{t_2^{n_v N}} 
\ = \ 4 \hR_2^{n_v N} \ = \ \hR_1^{n_v N} \ , \\
b_2 & = & t_2^{n_v N} - \frac{\hR_2^{2 n_v N}}{t_2^{n_v N}} \ = \ 2\sqrt{3} \hR_2^{n_v N} 
\ = \ \frac{\sqrt{3}}{2} \hR_1^{n_v N} \ , 
\end{eqnarray*}
and
\begin{eqnarray*}
A_2 & = & t_2^{m_v n_v N} + \frac{\hR_2^{2m_v n_v N}}{t_2^{m_v n_v N}} 
\ = \ \big(1 + (2-\sqrt{3})^{2 m_v} \big) t_2^{nN} \ , \\
B_2 & = & t_2^{m_v n_v N} - \frac{\hR_2^{2m_v n_v N}}{t_2^{m_v n_v N}} 
\ = \ \big(1 - (2-\sqrt{3})^{2 m_v} \big) t_2^{nN} \ , \\
\end{eqnarray*}
then $T_{m_v,\hR_2^{n_v N}}(E(a_2,b_2)) = E(A_2,B_2)$.  


Since $T_{m_v,\hR_2^{n_v N}}^{-1}(E(A_2,B_2)) \ = \ E(a_2,b_2)$ 
and $E(a_2,b_2) \subset D(0,\hR_1^{n_v N}) \subset E(A_1,B_1)$, 
it follows that\index{patching functions, initial $G_v^{(0)}(z)$!mapping properties of} 
\begin{equation} \label{FContainment5}
\big(G_v^{(0)}\big)^{-1}(E(A_2,B_2)) \ \subset \ F_v^{-1}(E(A_1,B_1)) \ \subset \ U_v \ .
\end{equation} 
Since $[-2\hR_2^{n_vN},2\hR_2^{n_vN}] \subset E(a_2,b_2)$ 
and $T_{m_v,\hR_2^{n_vN}}$ oscillates $m_v$ times between $\pm 2\hR_2^{nN}$ 
on $[-2\hR_2^{n_vN},2\hR_2^{n_vN}]$, an argument similar to the one for $F_v(z)$ 
shows that $G_v^{(0)}(z)$\index{patching functions, initial $G_v^{(0)}(z)$!mapping properties of} 
oscillates $n \tau_i$ times between $\pm 2 \hR_2^{nN}$ on each 
$E_{v,i}$ contained in $\cC_v(\RR)$, 
that all the zeros of $G_v^{(0)}(z)$\index{patching functions, initial $G_v^{(0)}(z)$!roots are confined to $E_v$} 
belong to $E_v^0$, and that $G_v^{(0)}(z)$ has $T_i = n \tau_i$ zeros in each $E_{v,i}$ (counting multiplicities).

\medskip
We now turn to the patching construction, which will be carried out in three phases. 

\medskip 
\noindent{\bf Phase 1. Patching the high-order coefficients.}
\index{coefficients $A_{v,ij}$!high-order}\index{patching!high-order coefficients}
\index{local patching for $\RR$-simple sets!Phase 1: high-order coefficients|(}  

In this phase we carry out the patching for stages $k = 1, \ldots, \kbar$.  
The fact that $E(a_2,b_2) \subset D(0,\hR_1^{n_v N})$, while $D(0,2\hR_1^{n_v N}) \subset E(A_1,B_1)$, 
gives us freedom to adjust $F_v(z)$ while maintaining (\ref{FContainment5}), and this is 
the basis for the magnification argument.  Write $\hT_{\ell}(z)$ for $T_{\ell,\hR_2^{n_v N}}(z)$. 
\index{magnification argument|ii}


%
%
Using the basis functions $\varphi_{ij}(z)$, $\varphi_{\lambda}(z)$
\index{basis!$L$-rational}
we can write\index{patching functions, initial $G_v^{(0)}(z)$!expansion of} 
\begin{equation} \label{DRFE1}
G_v^{(0)}(z)  = \ \sum_{i=1}^{m} \sum_{j=0}^{(n-1)N_i-1}
                             A_{v,ij} \varphi_{i,nN_i-j}(z)
           + \sum_{\lambda=1}^{\Lambda} A_{v,\lambda} \varphi_{\lambda}(z) 
                \ . 
\end{equation}
Here $A_{v,i0} = \tc_{v,i}^n$, for each $i$.  
The coefficients $A_{v,ij}$ with $j = 0, \ldots, \kbar N_i - 1$ 
\index{coefficients $A_{v,ij}$!high-order}
will be deemed ``high order''.  They will be patched a magnification argument
\index{magnification argument}
similar to the one when $K_v \cong \CC$:  we will sequentially 
modify the coefficients of $F_v(z)$, changing it from
\index{coefficients $A_{v,ij}$}
\begin{equation*}
F_v(z) \ = \ \sum_{i=1}^{m} \sum_{j=0}^{(n_v-1)N_i-1}
                         a_{v,ij} \varphi_{i,n_vN_i-j}(z)
          + \sum_{\lambda=1}^{\Lambda} a_{v,\lambda} \varphi_{\lambda}(z) 
\end{equation*}
to
\begin{equation*}
\hF_v(z) \ = \ F_v(z) 
 + \sum_{i=1}^{m} \sum_{j=1}^{\kbar N_i-1} \eta_{v,ij} \varphi_{i,n_vN_i-j}(z) \ ,  
\end{equation*}
thereby stepwise changing $G_v^{(0)}(z) = \hT_{m_v}(F_v(z))$ 
to $\hG_v(z) = \hT_{m_v}(\hF_v(z))$.
\index{patching functions, initial $G_v^{(0)}(z)$!for archimedean sets $E_v$!patched by magnification}  
We will require the $\eta_{v,ij}$ to be $K_v$-symmetric, so $\hF_v(z)$ is $K_v$-rational.
\index{$K_v$-symmetric!set of numbers}  

We claim there there is an $\epsilon_v > 0$  such that if the $\eta_{v,ij}$ 
are $K_v$-symmetric 
\index{$K_v$-symmetric!set of numbers}
and each $|\eta_{v,ij}|_v < \epsilon_v$, then $\hF_v$ oscillates $n_v \tau_i$ times
between $\pm \hR_1^{n_v N}$ on each $E_{v,i}$ contained in $\cC_v(\RR)$, and  
\begin{equation} \label{FContainment6} 
\hF_v^{-1}\big(D(0,\hR_1^{n_v N})\big) \ \subset U_v \ .
\end{equation}

To see this, put 
$\hM_v = \max_{1 \le i \le m} 
         \Big( \max_{1 \le j \le \kbar N_i} \|\varphi_{i,n_v N_i-j}\|_{\Ubar_v} \Big)$ 
and take $\epsilon_v$ small enough that $\epsilon_v \cdot \kbar N \hM_v < \hR_1^{n_v N}$. Write 
\begin{equation*}
H_v(z) \ = \ \sum_{i=1}^{m} \sum_{j=1}^{\kbar N_i-1} 
                        \eta_{v,ij} \varphi_{i,n_vN_i-j}(z)\ , 
\end{equation*}
so that $\hF_v(z) = F_v(z) + H_v(z)$.  
Since the $\eta_{v,ij}$ and $\varphi_{i,n_v N_i - j}(z)$ are $K_v$-symmetric, 
\index{$K_v$-symmetric!set of functions}
$H_v(z)$ is $K_v$-symmetric and in particular is real-valued on $\cC_v(\RR)$.
\index{$K_v$-symmetric}  
At each $z$ where $|F_v(z)|_v = 2 \hR_1^{n_v N}$ 
we have $|H_v(z)|_v < \varepsilon \cdot \kbar N \hM_v < \hR_1^{n_v N}$, 
and so $|\hF_v(z)|_v > \hR_1^{n_v N}$.
If $E_{v,i}$ is a component of $E_v$ contained in $\cC_v(\RR)$,    
then since $F_v(z)$ oscillates $n_v \tau_i$ times between $\pm 2 \hR_1^{n_v N}$ on $E_{v,i}$, 
it follows that $\hF_v(z)$ oscillates $n_v \tau_i$ times between $\pm\hR_1^{n_v N}$ on $E_{v,i}$.
It remains to show (\ref{FContainment6}).  For this, 
put $\Gamma = \{z \in \cC_v(\CC) : |F_v(z)|_v = 2 \hR_1^{n_v N}\}$, 
and apply Lemma \ref{DLemCP1} to $F_v(z)$ and $H_v(z)$, taking $R = 2 \hR_1^{n_v N}$ and $\delta = 1/2$.
Since $F_v^{-1}\big(D(0,2 \hR_1^{n_v N})\big) \subset U_v$, 
we conclude that $\hF_v^{-1}\big(D(0,\hR_1^{n_v N})\big) \subset U_v$.

Since $[-2 \hR_2^{n_v N}, 2 \hR_2^{n_v N}] \subset E(a_2,b_2) \subset D(0,\hR_1^{n_v N})$, 
the same argument as for $G_v^{(0)}(z)$\index{patching functions, initial $G_v^{(0)}(z)$!mapping properties of} 
shows that $\hG_v(z) = T_{m_v,\hR_2^{n_v N}}(\hF_v(z))$ 
oscillates $n \tau_i$ times between $\pm 2 \hR_2^{nN}$ on each 
$E_{v,i}$ contained in $\cC_v(\RR)$, that all the zeros of $\hG_v(z)$ belong to $E_v^0$, 
that $\hG_v(z)$ has $T_i = n \tau_i$ zeros in each $E_{v,i}$ (counting multiplicities), 
and that 
\begin{equation} \label{FContainment7}
\hG_v^{-1}\big(E(A_2,B_2)\big) \ \subset U_v \ .
\end{equation}    

\smallskip
Let $B_v > 0$ be the number in the statement of Theorem \ref{DRPatch1}.  
We now show that by choosing the $n$ and the $\eta_{v,ij}$ appropriately, 
we can achieve freedom $B_v$ in  patching the high order coefficients.
\index{patching argument!freedom $B_v$ in patching}  
\index{freedom $B_v$ in patching|ii}
\index{coefficients $A_{v,ij}$!high-order}  
That is, the $\Delta_{v,ij}^{(k)}$  with  $1 \le j < \kbar N_i$ can be can be specified arbitrarily, 
provided that they are $K_v$-symmetric and satisfy $|\Delta_{v,ij}^{(k)}|_v \le B_v$.  
\index{$K_v$-symmetric!set of numbers}

 By (\ref{DFSW2}),\index{patching functions, initial $G_v^{(0)}(z)$!construction of}
\begin{eqnarray*}
G_v^{(0)}(z) & = & \hT_{m_v}(F_v(z)) \\
& = & F_v(z)^{m_v} + \sum_{k=1}^{\lfloor m_v/2 \rfloor}
 (-1)^{k} \frac{m_v}{k} \binom{m_v-k-1}{k-1} \hR_2^{2k n_v N} F_v(z)^{m_v - 2k} \ .
\end{eqnarray*}
When the right side is expanded in terms of the basis functions, only the
\index{basis!$L$-rational}
pure power $F_v^{m_v}(z)$ can contribute to the coefficients of
\index{coefficients $A_{v,ij}$!high-order}
$\varphi_{i,nN_i-j}(z)$ with $j < 2 n_v N_i$; 
in particular this holds for $j < \kbar N_i$.  

Since only $F_v(z)^{m_v}$ contributes to the high order
coefficients, essentially the same argument applies here as in the
\index{coefficients $A_{v,ij}$!high-order}
patching construction over $\CC$.  By sequentially
adjusting the numbers $\eta_{v,ij}$ we will modify the
corresponding coefficients in the expansion of $\hT_{m_v}(F_v(z))$.
By (\ref{DXYZ2}), the change in  $A_{v,ij}$ induced by replacing $a_{v,ij}$ with
$a_{v,ij} + \eta_{v,ij}$ is
\begin{equation*}
\Delta_{v,ij}^{(k)} \ = \ m_v \eta_{v,ij}/ \tc_{v,ij}^{n_v - k + 1} \ .
\end{equation*}
Conversely, if a desired change $\Delta_{v,ij}^{(k)}$ is given, then taking
\begin{equation} \label{FEtaChoice}
\eta_{v,ij} = \frac{ \tc_{v,ij}^{n_v-k-1} \Delta_{v,ij}^{(k)}} {m_v}
\end{equation}
will produce that change.
\index{magnification argument|ii}

Henceforth we will assume $n$  is large enough that (with $m_v = n/n_v$) 
\begin{equation} \label{DRPF1}
\frac{B_v}{m_v} \cdot  \max(1,\max_{1 \le i \le m}(|\tc_{v,i}|_v))^{n_v} 
\ < \ \epsilon_v \ .
\end{equation}
If $|\Delta_{v,ij}^{(k)}|_v < B_v$, and $\eta_{v,ij}$ is defined by (\ref{FEtaChoice}), 
then $|\eta_{v,ij}|_v < \epsilon_v$  
so our discussion about the mapping properties of $\hF_v(z)$ applies.  

\vskip .1 in
For $i = 1,\ldots, m$ and $j = 0$, since $\Delta_{v,i0}^{(1)} = 0$ we have $\eta_{v,i0} = 0$.    
The remaining $\eta_{v,ij}$  will be determined recursively, 
in terms of the $\Delta_{v,ij}^{(k)}$, in $\prec_N$ order\index{order!$\prec_N$};  
in particular, for each $x_i$ the $\eta_{v,ij}$ will be determined
in order of increasing $j$.  As $F_v(z)$ is changed stepwise to $\hF_v(z)$, 
then $G_v^{(0)}(z) = \hT_{m_v}(F_v(z))$ 
is changed stepwise to
\index{patching functions, initial $G_v^{(0)}(z)$!for archimedean sets $E_v$!patched by magnification}
$G_v^{(\kbar)}(z) = \hT_{m_v}(\hF_v(z))$,
\index{patching functions, $G_v^{(k)}(z)$ for $1 \le k \le n$!for archimedean sets $E_v$!patched by magnification}  
passing through     
$G_v^{(1)}(z)$, $G_v^{(2)}(z)$, $\ldots$, $G_v^{(\kbar-1)}(z)$ at 
intermediate steps.

Consider what happens as each $\eta_{v,ij}$ is varied in turn.  
Suppose $\ckF_v(z)$ is a function obtained  
at one of the intermediate steps, and at the next step $\ckF_v(z)$ is replaced 
by $\ckF_v^{\prime}(z) = \ckF_v(z) + \eta_{v,ij} \varphi_{i,n_v N_i-j}(z)$.

If $x_i \in \cC_v(\RR)$, 
then $\varphi_{i,n_vN_i-j} \in K_v(\cC)$ and $\tc_{v,i} \in \RR$.   
By hypothesis, $\Delta_{v,ij}^{(k)} \in \RR$, so $\eta_{v,ij} \in \RR$ and
\begin{equation*}
\ckF_v^{\prime}(z) \ := \ \ckF_v(z) + \eta_{v,ij} \varphi_{i,n_vN_i-j}(z)
\end{equation*}
is $K_v$-symmetric. 
\index{$K_v$-symmetric}
Considering an expansion similar to (\ref{DFRH1}) 
one sees that\index{compensating functions $\vartheta_{v,ij}^{(k)}(z)$}
\begin{equation*}
\hT_{m_v}(\ckF_v^{\prime}(z)) 
\ = \ \hT_{m_v}(\ckF_v(z)) + \Delta_{v,ij}^{(k)} \vartheta_{v,ij}(z) 
\end{equation*}
for a $K_v$-rational (hence $K_v$-symmetric) function
\index{compensating functions $\vartheta_{v,ij}^{(k)}(z)$!are $K_v$-symmetric}
$\vartheta_{v,ij}(z)$ meeting\index{$K_v$-symmetric!function}
the conditions of the theorem.  

If $x_i \in \cC_v(\CC) \backslash \cC_v(\RR)$, 
let $\overline{z}$ denote the complex conjugate of $z$, 
and let $\ibar$ be the index such that $x_{\ibar} = \overline{x_i}$.  
As observed in the remark after Theorem \ref{CTCX2}, we can assume without loss
that $(i,j)$ and $(\ibar,j)$ are consecutive indices under $\prec_N$\index{order!$\prec_N$}.  
By our hypothesis of $K_v$-symmetry, 
$\Delta_{v,\ibar j}^{(k)} = \overline{\Delta_{v,i j}^{(k)}}$, 
and $\tc_{v,\ibar} = \overline{\tc_{v,i}}$, so 
$\eta_{v,\ibar j} = \overline{\eta_{v,ij}}$ and $N_{\ibar} = N_i$.
Hence, after two steps 
\begin{equation*}
\ckF_v^{\prime \prime}(z) \ = \ \ckF_v(z) + \eta_{v,ij} \varphi_{i,n_v N_i-j}(z) 
         + \eta_{v,\ibar j} \varphi_{\ibar,n_v N_i-j}(z)
\end{equation*}
is $K_v$-symmetric.  
\index{$K_v$-symmetric}
After grouping the terms in the multinomial
expansions\index{multinomial theorem} appropriately, 
one sees that\index{compensating functions $\vartheta_{v,ij}^{(k)}(z)$!construction of} 
\begin{equation*}
\hT_{m_v}(\ckF_v^{\prime \prime}(z)) \ = \ 
       \hT_{m_v}(\ckF_v(z)) + \Delta_{v,ij}^{(k)} \vartheta_{v,ij}^{(k)}(z)
                       + \Delta_{v,\ibar j}^{(k)} \vartheta_{v,\ibar j}^{(k)}(z)
\end{equation*}
where $\vartheta_{v,ij}^{(k)}(z), \vartheta_{v,\ibar j}^{(k)}(z) \in L_{w_v}(\cC)$\index{distinguished place $w_v$} 
meet the conditions of the theorem.\index{compensating functions $\vartheta_{v,ij}^{(k)}(z)$!are $K_v$-symmetric}

\medskip 

The order $\prec_N$\index{order!$\prec_N$} specifies the $\eta_{v,ij}$ 
in ``$N$-bands'' with $(k-1)N_i \le j < kN_i$, for $k = 1, 2, \ldots, \kbar$.\index{band!$\Band_N(k)$}
When we have completed patching the $k^{th}$ band, 
we obtain a function
\index{patching functions, $G_v^{(k)}(z)$ for $1 \le k \le n$!for archimedean sets $E_v$!oscillate on real components of $E_v$}  
$G_v^{(k)}(z)$ which oscillates $n\tau_i$ times between $\pm 2 \hR_2^{nN}$ on each real component $E_{v,i}^0$,
and satisfies $\{z : G_v^{(k)}(z) \in E(A_2,B_2) \} \subset U_v$.  
Since the $\Delta_{v,ij}^{(k)}$ are $K_v$-symmetric, 
\index{$K_v$-symmetric!function}
each $G_v^{(k)}(z)$ will be $K_v$-rational.\index{patching functions, $G_v^{(k)}(z)$ for $1 \le k \le n$!are $K_v$-rational}   

\medskip
The final function $\hF_v(z)$ thus obtained, for which
\index{patching functions, $G_v^{(k)}(z)$ for $1 \le k \le n$!for archimedean sets $E_v$!patched by magnification}    
\begin{equation*}
G_v^{(\kbar)}(z) \ = \ \hT_{m_v}(\hF_v(z)) \ ,
\end{equation*}
will play an important role in the rest of the argument.  
Note that the leading coefficient of $\hF_v(z)$ at
\index{coefficients $A_{v,ij}$!leading}
$x_i$ is $\tc_{v,i}^{n_v}$.\index{local patching for $\RR$-simple sets!Phase 1: high-order coefficients|)}

\medskip
\noindent{\bf Phase 2.  Patching the middle coefficients.}
\index{coefficients $A_{v,ij}$!middle}\index{patching!middle coefficients}
\index{local patching for $\RR$-simple sets!Phase 2: middle coefficients|(}  

In this phase we will construct
functions $G_v^{(k)}(z)$ for $k = \kbar+1, \ldots, n-1$, 
setting\index{patching functions, $G_v^{(k)}(z)$ for $1 \le k \le n$!constructed by patching} 
\begin{equation} \label{FXCM1}
G_v^{(k)}(z) \ = \  G_v^{(k-1)}(z) + \sum_{i=1}^{m} \sum_{j=(k-1)N_i}^{kN_i-1}
                        \Delta_{v,ij}^{(k)} \vartheta_{v,ij}(z)
\end{equation}
for the given $\Delta_{v,ij}^{(k)}$ and appropriate functions 
$\vartheta_{v,ij}(z) \in L_{w_v}(\cC)$,\index{distinguished place $w_v$}
\index{compensating functions $\vartheta_{v,ij}^{(k)}(z)$} 
adjoining the terms in $\prec_N$\index{order!$\prec_N$} order.  
The conditions of theorem require that 
the $\Delta_{v,ij}^{(k)} \in \CC$ be $K_v$-symmetric; 
\index{$K_v$-symmetric!set of numbers}
they also satisfy 
\begin{equation*}
|\Delta_{v,ij}^{(k)}|_v \ \le \ h_v^{kN} \ .
\end{equation*}
The $\vartheta_{v,ij}(z)$ will be $\Aut_c(\CC_v/K_v)$-equivariant,
\index{compensating functions $\vartheta_{v,ij}^{(k)}(z)$!are $K_v$-symmetric} 
so for each $(i,j)$ 
\begin{equation*}
\sum_{(i^{\prime},j) \in \Aut_c(\CC_v/K_v)(i,j)} 
       \Delta_{v,ij}^{(k)} \vartheta_{v,ij}(z) \ \in \ K_v(\cC) \ ;
\end{equation*}        
consequently $G_v^{(k)}(z) \in K_v(\cC)$ 
as well.\index{patching functions, $G_v^{(k)}(z)$ for $1 \le k \le n$!are $K_v$-rational} 

Fix $k$, and write
\begin{equation*}
n-k-1 \ = \ \ell_1 + n_v \ell_2, 
    \quad \text{with $0 \le \ell_1 < n_v$, $0 \le \ell_2 < m_v$.}
\end{equation*}
For each $i$, and each $j$ with $(k-1)N_i \le j < k N_i$,
we can uniquely write
\begin{equation*}
nN_i-j = r_{ij} + (n-k-1) N_i \quad \text{where $N_i < r_{ij} \le 2N_i$}
\end{equation*}
so $nN_i-j = r_{ij} + \ell_1 N_i + \ell_2 n_v N_i$.  
Put\index{compensating functions $\vartheta_{v,ij}^{(k)}(z)$!construction of} 
\begin{equation*}
\vartheta_{v,ij}(z) \ = \  \varphi_{i,r_{ij}}(z) 
       \cdot T_{\ell_1,R_1^N}(\phi_v(z)) \cdot \hT_{\ell_2}(\hF_v(z)) \ .
\end{equation*}      
Then $\vartheta_{v,ij}(z)$
\index{compensating functions $\vartheta_{v,ij}^{(k)}(z)$!poles and leading coefficients of} 
has a pole of order $nN_i-j$ at $x_i$, and leading
coefficient $\tc_{v,i}^{n-k-1}$.  Its poles at the $x_{i^{\prime}} \ne x_i$ are 
\index{coefficients $A_{v,ij}$!leading}
of order at most $(n-k-1)N_{i^{\prime}}$, 
so it meets the requirements of the global patching process. 
\index{patching argument!global}   
Since the $\varphi_{ij}(z)$ are $\Aut_c(\CC_v/K_v)$-equivariant, 
and $\phi_v(z)$ and $\hF_v(z)$ are $K_v$-rational, 
the $\vartheta_{v,ij}(z)$\index{compensating functions $\vartheta_{v,ij}^{(k)}(z)$!are $K_v$-symmetric} 
are $\Aut_c(\CC_v/K_v)$-equivariant. 

\smallskip
Define\index{patching functions, $G_v^{(k)}(z)$ for $1 \le k \le n$!mapping properties of} 
\begin{eqnarray*}
& & \qquad \hE_1 \ = \ \{ z \in E_v \cap \cC_v(\RR) : \phi_v(z) \in [-2\hR_1^{n_v N},2\hR_1^{n_v N}] \} \ , \\
& & \hE_2 \ = \ \{ z \in E_v \cap \cC_v(\RR) : \hF_v(z) \in [-2\hR_2^{n_v N},2\hR_2^{n_v N}] \} \\    
& & \qquad \qquad \qquad \qquad 
       \ = \ \{ z \in \cC_v(\RR) : G_v^{(\kbar)}(z) \in [-2\hR_2^{n N},2\hR_2^{n N}] \} \ . 
\end{eqnarray*}
By construction, $\hE_2 \subset \hE_1 \subset E_v \cap \cC_v(\RR)$.  
Likewise, put\index{patching functions, $G_v^{(k)}(z)$ for $1 \le k \le n$!mapping properties of}  
\begin{eqnarray*}
& & \quad \hW_1 \ = \ \phi_v^{-1}(E(a_1,b_1)) \ = \ F_v^{-1}(E(A_1,B_1)) \ , \\ 
& & \hW_2 \ = \ \hF_v^{-1}(E(a_2,b_2)) \ = \ \big(G_v^{(\kbar)}\big)^{{-1}}(E(A_2,B_2)) \ .    
\end{eqnarray*}
Then $\hW_2 \subset \hW_1 \subset U_v$.   

We will now bound the change\index{patching functions, $G_v^{(k)}(z)$ for $1 \le k \le n$!constructed by patching}  
in $|G_v^{(k)}(z)-G_v^{(k-1)}(z)|_v$ on $\hE_2$ and $\hW_2$.  
We first bound\index{patching functions, $G_v^{(k)}(z)$}  $\|G_v^{(k)}-G_v^{(k-1)}\|_{\hE_2}$,  
and we begin by showing that for each $(i,j)$ 
occurring in (\ref{FXCM1}),\index{compensating functions $\vartheta_{v,ij}^{(k)}(z)$!bounds for} 
\begin{equation} \label{DRPF3A}
\|\vartheta_{v,ij}\|_{\hE_2} \le 16 M_v \hR_2^{(n-k-1)N} \ .
\end{equation}
To see this, note that by the definition of $M_v$,
we have $| \varphi_{i,r_{ij}}(z)|_v \ \le \ M_v$
for all $z \in E_v \cap \cC_v(\RR)$.  Also, since $\ell_1 < n_v$, 
for each $z \in \hE_1$, 
\begin{equation*}
|T_{\ell_1,R_1^N}(\phi_v(z))|_v \ \le \ 2 \hR_1^{\ell_1 N} 
                  \ \le \ 8 \hR_2^{\ell_1 N} \ .
\end{equation*}
by the definitions of $\hR_1$ and $\hR_2$.  
Finally, by the properties of Chebyshev polynomials, for each $z \in \hE_2$.
\index{Chebyshev polynomial} 
\begin{equation*}
|\hT_{\ell_2}(\hF_v(z))|_v \ \le \ 2 \hR_2^{\ell_2 n_v N} 
\end{equation*}
Combining these, and using that $\hE_2 \subset \hE_1 \subset E_v \cap \cC_v(\RR)$, 
gives (\ref{DRPF3A}).

Since each $|\Delta_{v,ij}^{(k)}|_v \le h_v^{kN}$ and $1 < r_v < \hR_2$, 
it follows that\index{compensating functions $\vartheta_{v,ij}^{(k)}(z)$}
\index{patching functions, $G_v^{(k)}(z)$ for $1 \le k \le n$!roots are confined to $E_v$}  
\begin{eqnarray}
 \|G_v^{(k)} - G_v^{(k-1)}\|_{\hE_2}
  & \le & \sum_{i=1}^{m} \sum_{j=(k-1)N_i}^{kN_i-1} 
                    |\Delta_{v,ij}^{(k)}|_v  \|\vartheta_{v,ij}(z)\|_{\hE_2} \notag \\
  & \le & N \cdot h_v^{kN} \cdot 16 M_v \hR_2^{(n-k-1)N}  
  \ = \ \frac{16 N M_v}{\hR_2^N} \cdot \hR_2^{nN} \cdot 
                \Big(\frac{h_v}{\hR_2}\Big)^{kN} \label{DRPF4} \\
  & < & 16 N M_v \cdot \hR_2^{nN} \cdot \left(\frac{h_v}{r_v}\right)^{kN} \ .  \notag
\end{eqnarray}

We next bound\index{patching functions, $G_v^{(k)}(z)$}   $\|G_v^{(k)} - G_v^{(k-1)}\|_{\hW_2}$.  
We claim that for each $(i,j)$ in (\ref{FXCM1}), 
\begin{equation} \label{DRPF3AB}
\|\vartheta_{v,ij}\|_{\hW_2} \le 5 M_v t_2^{(n-k-1)N} \ .
\end{equation}
To see this,\index{compensating functions $\vartheta_{v,ij}^{(k)}(z)$!bounds for} 
recall that $| \varphi_{i,r_{ij}}(z)|_v \ \le \ M_v$
for all $z \in \Ubar_v$.  Also, note that $T_{\ell_1,\hR_1^N}$ maps $E(a_1,b_1)$ to $E(A^{(\ell_1)},B^{(\ell_1)})$ where $B^{(\ell_1)} < A^{(\ell_1)}$ and 
\begin{equation*}  
A^{(\ell_1)} \ = \ t_1^{\ell_1 N} +  \frac{\hR_1^{2 \ell_1 N}}{t_1^{\ell_1 N}}  
  \ = \ \big(1 + (2-\sqrt{3})^{2 \ell_1}\big) \cdot t_1^{\ell_1 N} 
   \ < \ 1.1 \cdot t_1^{\ell_1 N} \ .
\end{equation*} 
Since  $t_1^{n_v N} = (2 + \sqrt{3}) \hR_1^{n_v N}$,
$\hR_1^{n_v N} = 4 \hR_2^{n_v N}$ and $t_2^{n_v N} = (2 + \sqrt{3}) \hR_2^{n_v N}$, 
we have $t_1^{n_v N} = 4 t_2^{n_v N}$ and 
for $0 \le \ell_1 < n_v$ it follows that $t_1^{\ell_1 N} < 4 t_2^{\ell_1 N}$.   
This means that $E(A^{(\ell_1)},B^{(\ell_1)}) \subset D(0,1.1 \cdot 4t_2^{\ell_1 N})$, 
so for each $z \in \hW_1 = \phi_v^{-1}(E(a_1,b_1))$, 
\begin{equation} \label{FZZBB1}
|T_{\ell_1,\hR_1^N}(\phi_v(z))|_v \ \le \ 1.1 \cdot 4 t_2^{\ell_1 N} \ .
\end{equation}
Similarly, since $\hT_{\ell_2} = T_{\ell_2, \hR_2^{n_v N}}$ 
maps $E(a_2,b_2)$ to $E(\tA^{(\ell_2)},\tB^{(\ell_2)})$
where $\tB^{(\ell_2)} < \tA^{(\ell_2)}$ and 
\begin{equation*}  
\tA^{(\ell_1)} \ = \ t_2^{\ell_2 n_v N} + \frac{\hR_2^{2 \ell_2 n_v N}}{t_2^{\ell_2 n_v N}} 
\ = \ \big(1 + (2-\sqrt{3})^{2 \ell_2}\big) \cdot t_2^{\ell_2 n_v N} 
\ < \ 1.1 \cdot t_2^{\ell_1 n_v N} \ ,
\end{equation*}
for each $z \in \hW_2 = \hF_v^{-1}(E(a_2,b_2))$ we have  
\begin{equation} \label{FZZBB2}
|\hT_{\ell_2}(\hF_v(z))|_v \ \le \ 1.1 \cdot t_2^{\ell_2 n_v N} \ .
\end{equation}
Since $\hW_2 \subset \hW_1 \subset \Ubar_v$ and $4\cdot(1.1)^2 < 5$, 
combining (\ref{FZZBB1}) and (\ref{FZZBB2}) gives (\ref{DRPF3AB}).

Since each $|\Delta_{v,ij}^{(k)}|_v \le h_v^{kN}$ and $1 < r_v < \hR_2 < t_2$, 
it follows that\index{compensating functions $\vartheta_{v,ij}^{(k)}(z)$}
\index{patching functions, $G_v^{(k)}(z)$ for $1 \le k \le n$!roots are confined to $E_v$}  
\begin{eqnarray}
 \|G_v^{(k)} - G_v^{(k-1)}\|_{\hW_2}
  & \le & \sum_{i=1}^{m} \sum_{j=(k-1)N_i}^{kN_i-1} |\Delta_{v,ij}^{(k)}|_v
                  \|\vartheta_{v,ij}(z)\|_{\hW_2} \notag \\
  & \le & N \cdot h_v^{kN} \cdot 5 M_v t_2^{(n-k-1)N}  
  \ = \ \frac{5 N M_v}{t_2^{N}} \cdot t_2^{nN} 
             \cdot \Big(\frac{h_v}{t_2}\Big)^{kN} \label{DRPF4B} \\
  & < & 5 N M_v \cdot t_2^{nN} \cdot \left(\frac{h_v}{r_v}\right)^{kN} \ .  \notag
\end{eqnarray}\index{local patching for $\RR$-simple sets!Phase 2: middle coefficients|)}

\medskip
\noindent{\bf Phase 3.  Patching the low-order coefficients.}
\index{coefficients $A_{v,ij}$!low-order}\index{patching!low-order coefficients}
\index{local patching for $\RR$-simple sets!Phase 3: low-order coefficients|(}  

In the final step we take\index{patching functions, $G_v^{(k)}(z)$ for $1 \le k \le n$!constructed by patching} 
\begin{equation*}
G_v^{(n)}(z) \ = \ G_v^{(n-1)}(z)
          + \sum_{\lambda = 1}^{\Lambda} \Delta_{v,\lambda}^{(n)} \varphi_{\lambda}
\end{equation*}
with $K_v$-symmetric $\Delta_{v,\lambda}^{(n)}$ 
\index{$K_v$-symmetric!set of numbers}
satisfying $|\Delta_{v,\lambda}^{(n)}|_v \le h_v^{nN}$ for each $\lambda$.  
Since $\Lambda \le N$, and each $|\varphi_{\lambda}|_v \le M_v$ on $\Ubar_v$, 
it follows that on $\hE_2$\index{patching functions, $G_v^{(k)}(z)$ for $1 \le k \le n$!roots are confined to $E_v$} 
\begin{equation} \label{DRFP6}
\|G_v^{(n)}-G_v^{(n-1)}\|_{\hE_2}  
\ \le \ N M_v h_v^{nN} \ < \ 16 N M_v \cdot \hR_2^{nN} \cdot \Big(\frac{h_v}{r_v}\Big)^{nN}\ , 
\end{equation}
while on $\hW_2$\index{patching functions, $G_v^{(k)}(z)$ for $1 \le k \le n$!roots are confined to $E_v$}  
\begin{equation} \label{DRFP6B}
\|G_v^{(n)}-G_v^{(n-1)}\|_{\hW_2}  
\ \le \ N M_v h_v^{nN} \ < \ 5 N M_v \cdot t_2^{nN} \cdot \Big(\frac{h_v}{r_v}\Big)^{nN}\ .
\end{equation}\index{local patching for $\RR$-simple sets!Phase 3: low-order coefficients|)}  

\vskip .1 in
To conclude the proof, we show that if $n$ is sufficiently large, 
then $G_v^{(n)}(z)$\index{patching functions, $G_v^{(k)}(z)$ for $1 \le k \le n$!mapping properties of}  
has the mapping properties in part (B) of the Theorem.  
A similar argument applies\index{patching functions, $G_v^{(k)}(z)$}  to $G_v^{(k)}(z)$ for each $k = 1, \ldots, n$.    

Consider the total change in passing from\index{patching functions, $G_v^{(k)}(z)$ for $1 \le k \le n$!mapping properties of} 
$G_v^{(\kbar)}(z) = \hT_{m_v}(\hF_v(z))$ to $G_v^{(n)}(z)$.  By (\ref{DRPF4})
and (\ref{DRFP6}), for each $z \in \hE_2$,
\begin{equation*}
|G_v^{(n)}(z) - G_v^{(\kbar)}(z)|_v \ \le \
            16 N M_v \cdot \hR_2^{nN} 
              \cdot \sum_{k=\kbar}^{n} \left(\frac{h_v}{r_v}\right)^{kN} \ < \ 
   \frac{16 N M_v}{1 - (h_v/r_v)^{N}} 
                \cdot \Big(\frac{h_v}{r_v}\Big)^{\kbar N} \cdot \hR_2^{nN} \ .
\end{equation*}
Since $\kbar \ge k_v$, assumption (\ref{DRPFF0}) in Theorem \ref{DRPatch1} shows that 
\index{patching functions, $G_v^{(k)}(z)$ for $1 \le k \le n$!constructed by patching} 
\begin{equation} \label{DRPF9}
\|G_v^{(n)} - G_v^{(\kbar)}\|_{\hE_2} \ < \ \frac{1}{4} \hR_2^{nN} \ .
\end{equation}
Similarly, on $\hW_2$\index{patching functions, $G_v^{(k)}(z)$}  
\begin{equation} \label{DRPF9B}
\|G_v^{(n)} - G_v^{(\kbar)}\|_{\hW_2} \ < \ \frac{1}{4} t_2^{nN} \ .
\end{equation}

We first show that\index{patching functions, $G_v^{(k)}(z)$ for $1 \le k \le n$!mapping properties of}  
\begin{equation} \label{FWEnuf}
\{z \in \cC_v(\CC) : |G_v^{(n)}(z)|_v \le \frac{1}{2}t_2^{nN} \}  \ \subset \ U_v \ .
\end{equation} 
Let $\hGamma_2 = \partial \hW_2 = \{z \in \cC_v(\CC) : G_v^{(\kbar)}(z) \in \partial E(A_2,B_2)\}$.  
Since $A_2 > B_2$ and 
\begin{equation*} 
B_2 \ = \ (1-(2-\sqrt{3})^{2 m_v}) t_2^{nN} \ > \ 0.9 t_2^{nN} \ ,
\end{equation*}
for each $z \in \hGamma_2$ we have $|G_v^{(\kbar)}(z)|_v > 0.9 t_2^{nN}$.
\index{patching functions, $G_v^{(k)}(z)$ for $1 \le k \le n$!roots are confined to $E_v$}   
By (\ref{DRPF9B}), 
$|G_v^{(n)}(z)-G_v^{(\kbar)}(z)|_v < 0.25 t_2^{nN}$.
Applying Lemma \ref{DLemCP1} with $F(z) = G_v^{(\kbar)}(z)$  
and $H(z) = G_v^{(n)}(z)-G_v^{(\kbar)}(z)$, we see that 
\begin{equation*}
\{z \in \cC_v(\CC) : |G_v^{(n)}(z)|_v \le 0.65 t_2^{nN} \} \ \subset \ \hW_2 \ ,
\end{equation*} 
which yields (\ref{FWEnuf}).  

If $E_{v,i}$ is a component of $E_v$ which is disjoint from $\cC_v(\RR)$, 
then $G_v^{(\kbar)}(z)$\index{patching functions, $G_v^{(k)}(z)$ for $1 \le k \le n$!roots are confined to $E_v$} 
has $n \tau_i$ zeros in $E_{v,i}^0$.  Put $\hW_{2,i} = \hW_2 \cap E_{v,i}$.  
Since $\hW_2 \subset U_v$ and $U_v \cap E_{v,i} = E_{v,i}^0$, 
it follows that $\hW_{2,i} \subset E_{v,i}^0$.  
Applying Rouch\'e's theorem to $G_v^{(\kbar)}(z)$ and $G_v^{(n)}(z)$ on $\partial \hW_{2,i}$,
\index{Rouch\'e's theorem}\index{patching functions, $G_v^{(k)}(z)$} 
we conclude that $G_v^{(n)}(z)$ has $n \tau_i$ zeros in $E_{v,i}^0$ as well. 

We next show that if $E_{v,i}$ is a component of $E_v$ contained in $\cC_v(\RR)$, 
then $G_v^{(n)}(z)$
\index{patching functions, $G_v^{(k)}(z)$ for $1 \le k \le n$!for archimedean sets $E_v$!oscillate on real components of $E_v$}  
oscillates $n \tau_i$ times between $\pm (7/4) \hR_2^{nN}$ on $E_{v,i}$. 
Recall that $G_v^{(\kbar)}(z) = \hT_{m_v}(\hF_v(z))$\index{patching functions, $G_v^{(k)}(z)$}  
oscillates $n \tau_i$ times  
between $\pm 2 \hR_2^{nN}$ on $E_{v,i}$. Equation (\ref{DRPF9}) shows that
at each  $z_i \in E_{v,i}$ where $G_v^{(\kbar)}(z_i) = \pm 2 \hR_2^{nN}$,
then $G_v^{(n)}(z_i)$ has the same sign as $G_v^{(\kbar)}(z_i)$ and
$|G_v^{(n)}(z_i)|_v \ge (7/4) \hR_2^{nN}$.
Hence $G_v^{(n)}(z)$ oscillates $n \tau_i$ times\index{patching functions, $G_v^{(k)}(z)$ for $1 \le k \le n$!roots are confined to $E_v$} 
between $\pm (7/4) \hR_2^{nN}$ on $E_{v,i}$, and in particular it has $n \tau_i$ zeros 
in $E_{v,i}^0$.  

Since $G_v^{(n)}(z)$ has degree $nN$ and $\sum n \tau_i = nN$, all the zeros of $G_v^{(n)}(z)$
lie in $E_v^0$.\index{patching functions, $G_v^{(k)}(z)$ for $1 \le k \le n$!roots are confined to $E_v$}   
Since $t_2 > \hR_2 > r_v$, for all sufficiently large $n$ we have 
$(1/2) t_2^{nN} > 2 r_v^{nN}$ and $(7/4) \hR_2^{nN} > 2 r_v^{nN}$.  For such $n$, 
$G_v^{(n)}(z)$ oscillates $n \tau_i$ times between $\pm 2 r_v^{nN}$ on each $E_{v,i}$
contained in $\cC_v(\RR)$, and 
\begin{equation*}  
\{z \in \cC_v(\CC) : |G_v^{(n)}(z)|_v \le 2 r_v^{nN} \} \ \subset \ U_v \ .
\end{equation*}
Thus the construction succeeds for any integer $n$ divisible by $n_v$ which is  
large enough that condition (\ref{DRPF1}) holds and $n > \kbar$, $(1/2) t_2^{nN} > 2 r_v^{nN}$, 
and $(7/4) \hR_2^{nN} > 2 r_v^{nN}$.  
\end{proof}
\index{patching construction!for the case when $K_v \cong \RR$|)}
\index{patching argument!local!for the case when $K_v \cong \RR$|}   
 

%% file: NewFSZChap10.tex
\chapter{The Local Patching Construction for Nonarchimedean $\RL$-domains} 
\index{$\RL$-domain}
\label{Chap10}

In this section we give the confinement argument for Theorem \ref{aT1-B} when $K_v$ is nonarchimedean, 
and $E_v$ is a $K_v$-rational $RL$-domain.
\index{confinement argument}  \index{$\RL$-domain}
\index{patching argument!local!for nonarchimedean $\RL$-domains|(}  

Let $q_v$ be the order of the residue field of $K_v$. 
Let $w_v$ be the distinguished place\index{distinguished place $w_v$} of $L = K(\fX)$ 
determined by the embedding $\tK \hookrightarrow \CC_v$
used to identify $\fX$ with a subset of $\cC_v(\CC_v)$, 
and view $L_{w_v}$ as a subset of $\CC_v$.  Following the construction
of the coherent approximating functions in Theorems \ref{CTCX2} and \ref{CTCX2p}, 
\index{coherent approximating functions $\phi_v(z)$}
we begin with the following data: 
    
\begin{enumerate}
  \item A $K_v$-symmetric probability vector $\vs \in \cP^{m}(\QQ)$ 
\index{$K_v$-symmetric!probability vector}
           with positive rational coefficients;  
      
  \item A number $N$, a number $R_v \in |\CC_v^{\times}|_v$,
           and an $(\fX,\vs)$-function $\phi_v(z) \in K_v(\cC)$ of degree $N$ such that 
            \begin{equation*}
              \{ z \in \cC_v(\CC_v) : |\phi_v(z)|_v \le R_v^N \} \ = \ E_v \ .      
            \end{equation*}  
        If $\Char(K_v) = p > 0$, we will assume that the number $J$ from the construction 
        of the $L$-rational and $L^{\sep}$-rational bases 
        in \S\ref{Chap3}.\ref{LRationalBasisSection}\index{basis!$L$-rational}\index{basis!$L^{\sep}$-rational}
        divides $N_i := Ns_i$, 
        and that the leading coefficient\index{coefficients $A_{v,ij}$!leading}
        $\tc_{v,i}$ of $\phi_v(z)$ at $x_i$ belongs to $K_v(x_i)^{\sep}$, for each $i = 1, \ldots, m$.

   \item Parameters $h_v$ and $r_v$ such that $0 < h_v <r_v \le R_v$, 
              which govern the freedom in the patching process;
            \index{patching argument!freedom in}\index{patching parameters}              
           
  \item  An order $\prec_N$\index{order!$\prec_N$} on the index set 
              $\cI = \{(i,j) \in \ZZ^2 : 1 \le i \le m, 0 \le j\}$ 
           determined by $N$ and $\vs$ as in (\ref{FPrec}), 
           which gives the sequence in which coefficients are patched.
\index{coefficients $A_{v,ij}$}
\end{enumerate} 

\smallskip
We will use the $L$-rational basis  
\index{basis!$L$-rational}
$\{\varphi_{ij}, \varphi_{\lambda}\}$ 
from \S\ref{Chap3}.\ref{LRationalBasisSection} to expand functions, 
and $\Lambda = \dim_K(\Gamma(\sum_{i=1}^m N_i(x_i)))$ 
will be the number of low-order basis elements, 
as in \S\ref{Chap7}.\ref{Char0Section}.  
The order $\prec_N$\index{order!$\prec_N$} respects the $N$-bands (\ref{FBand}),\index{band!$\Band_N(k)$} 
and for each $x_i \in \fX$, specifies the terms to be patched in decreasing pole order. 

\medskip

When $\Char(K_v) = 0$, we will need the following the following patching theorem.
\index{patching construction!for nonarchimedean $\RL$-domains!when $\Char(K_v) = 0$}    

 

\begin{theorem}  \label{DCPPatch1}\index{patching theorem!for nonarchimedean $\RL$-domains!when $\Char(K_v) = 0$} 
Suppose $K_v$ is nonarchimedean, and that $\Char(K_v) = 0$.
Let $\vs \in \cP^m(\QQ)$ be a $K_v$-symmetric probability vector with positive
\index{$K_v$-symmetric!probability vector} 
rational coefficients. Suppose $R_v \in |\CC_v^{\times}|_v$, let $0 < h_v < r_v \le R_v$,  
and let $\phi_v(z) \in K_v(\cC)$ be an $(\fX,\vs)$-function of degree $N$ satisfying 
\begin{equation*}
\{ z \in \cC_v(\CC_v) : |\phi_v(z)|_v \le R_v^N \} \ = \ E_v  \ . 
\end{equation*} 

For each $i = 1, \ldots, m$ let $N_i = Ns_i$  
and let $\tc_{v,i} = \lim_{z \rightarrow x_i} \phi_v(z) \cdot g_{x_i}(z)^{N_i}$
be the leading coefficient of $\phi_v(z)$ at $x_i$.  Put 
\index{coefficients $A_{v,ij}$!leading}
\begin{equation*}
M_v \ = \ \max( \max_{\substack{ 1 \le i \le m \\ N_i < j \le 2N_i}} \|\varphi_{ij}\|_{E_v },
    \max_{1 \le \lambda \le \Lambda} \|\varphi_{\lambda}\|_{E_v } ) \ .
\end{equation*} 
Let  $k_v \ge 1$  be the least integer such that
\begin{equation} \label{DCPF2C}
\Big(\frac{h_v}{r_v}\Big)^{k_v N} \cdot \frac{M_v}{\min(1,R_v^N)} \ < \ 1 \ ,
\end{equation}
and let $\kbar \ge k_v$ be a fixed integer.  Put $n_v = 1$. 
Then there is a number $0 < B_v < 1$ depending on $\kbar$, $E_v$, and $\phi_v(z)$ 
such that for each sufficiently large integer $n$,  
one can carry out the local patching process at $K_v$ as follows $:$  

Put $G_v^{(0)}(z) = \phi_v(z)^n$.\index{patching functions, initial $G_v^{(0)}(z)$!construction of}
For each $k$, $1 \le k < n$, 
let $\{\Delta_{v,ij}^{(k)} \in L_{w_v}\}_{(i,j) \in \Band_N(k)}$\index{distinguished place $w_v$} 
be a $K_v$-symmetric set of numbers given recursively in $\prec_N$ order,\index{order!$\prec_N$}\index{band!$\Band_N(k)$}
\index{$K_v$-symmetric!set of numbers} 
such that for each $(i,j)$   
\begin{equation} \label{DFCCV1}
|\Delta_{v,ij}^{(k)}|_v \ \le \ \left\{
      \begin{array}{ll} B_v & \text{if \ $k \le \kbar \ ,$} \\
                         h_v^{kN} & \text{if \ $k > \kbar \ .$}
      \end{array} \right.
\end{equation}
For $k = n$, let 
$\{\Delta_{v,\lambda}^{(n)} \in L_{w_v}\}_{1 \le \lambda \le \Lambda}$\index{distinguished place $w_v$} 
be a $K_v$-symmetric set of numbers such that for each $\lambda$ 
\index{$K_v$-symmetric!set of numbers}
\begin{equation} \label{DFCCF2} 
|\Delta_{v,\lambda}^{(n)}|_v \ \le \ h_v^{nN} \ .  
\end{equation}

\vskip .05 in
Then one can inductively construct 
$(\fX,\vs)$-functions $G_v^{(1)}(z), \ldots, G_v^{(n)}(z) \in K_v(\cC)$,
\index{patching functions, $G_v^{(k)}(z)$ for $1 \le k \le n$!constructed by patching}
of common degree $nN$, having the following properties $:$  
 
\noindent{$(A)$} For each $k = 1, \ldots, n-1$,\index{patching functions, $G_v^{(k)}(z)$ for $1 \le k \le n$!constructed by patching} 
\begin{equation*} 
G_v^{(k)}(z) \ = \ G_v^{(k-1)}(z) +  \sum_{(i,j) \in \Band_N(k)} 
                           \Delta_{v,ij}^{(k)} \cdot \vartheta_{v,ij}^{(k)}(z)
\end{equation*}                           
where $\vartheta_{v,ij}^{(k)}(z) = \varphi_{i,(k+1)N_i - j} \cdot \phi_v(z)^{n-k-1}$,
\index{band!$\Band_N(k)$}\index{compensating functions $\vartheta_{v,ij}^{(k)}(z)$!construction of}  
and for $k = n$\index{patching functions, $G_v^{(k)}(z)$ for $1 \le k \le n$!constructed by patching}                            
\begin{equation*}                           
G_v^{(n)}(z) \ = \ G_v^{(n-1)}(z) + \sum_{\lambda=1}^{\Lambda}
                   \Delta_{v,\lambda}^{(n)} \cdot \varphi_{\lambda}(z) \ .  
\end{equation*}
In particular

$(1)$ Each $\vartheta_{v,ij}^{(k)}(z)$ belongs to $K_v(x_i)(\cC)$,
\index{compensating functions $\vartheta_{v,ij}^{(k)}(z)$!poles and leading coefficients of}  
has a pole of order $nN_i-j > (n-k-1)N_i$ at $x_i$ with leading coefficient $\tc_{v,i}^{n-k-1}$, 
\index{coefficients $A_{v,ij}$!leading} 
has poles of order at most $(n-k-1)N_{i^{\prime}}$ at each $x_{i^{\prime}} \ne x_i$, 
and has no other poles $;$  

$(2)$ The $\vartheta_{v,ij}^{(k)}(z)$ are $K_v$-symmetric $;$
\index{compensating functions $\vartheta_{v,ij}^{(k)}(z)$!are $K_v$-symmetric} 
\index{$K_v$-symmetric!set of functions}




$(3)$  In passing from $G_v^{(0)}(z)$\index{patching functions, initial $G_v^{(0)}(z)$!leading coefficients of} 
to $G_v^{(1)}(z)$,\index{patching functions, $G_v^{(k)}(z)$ for $1 \le k \le n$!leading coefficients of} each of the leading coefficients
\index{coefficients $A_{v,ij}$!leading}
$A_{v,i0} = \tc_{v,i}^n$ is replaced with $\tc_{v,i}^n + \Delta_{v,i0}^{(1)} \cdot \tc_{v,i}^{n-2}$. 

\vskip .05 in
\noindent{$(B)$} For each $k = 0, \ldots, n$,\index{patching functions, $G_v^{(k)}(z)$ for $1 \le k \le n$!roots are confined to $E_v$} 
\quad  
$\{ z \in \cC_v(\CC_v) : |G_v^{(n)}(z)|_v \le R_v^{Nn} \} \ = \ E_v \ .$ 
\end{theorem}

When $\Char(K_v) = p > 0$, we will use the following patching theorem instead.
The $K_v$-rationality assumptions in the theorem are addressed by the global patching process.
\index{patching construction!for nonarchimedean $\RL$-domains!when $\Char(K_v) = p > 0$}    

\begin{theorem}  \label{DCPPatch1p}\index{patching theorem!for nonarchimedean $\RL$-domains!when $\Char(K_v) = p > 0$}  
Suppose $K_v$ is nonarchimedean and $\Char(K_v) = p > 0$. 
Let $\vs \in \cP^m(\QQ)$ be a $K_v$-symmetric probability vector with positive 
\index{$K_v$-symmetric!probability vector}
rational coefficients. Suppose $R_v \in |\CC_v^{\times}|_v$, let $0 < h_v < r_v \le R_v$,  
and let $\phi_v(z) \in K_v(\cC)$ be an $(\fX,\vs)$-function of degree $N$ satisfying 
\begin{equation*}
\{ z \in \cC_v(\CC_v) : |\phi_v(z)|_v \le R_v^N \} \ = \ E_v \ . 
\end{equation*} 

Let $N_i = Ns_i$ for each $x_i$, 
and let $\tc_{v,i} = \lim_{z \rightarrow x_i} \phi_v(z) \cdot g_{x_i}(z)^{N_i}$
be the leading coefficient of $\phi_v(z)$ at $x_i$.  
\index{coefficients $A_{v,ij}$!leading} 
Assume that $J|N_i$ and $\tc_{v,i} \in K_v(x_i)^{\sep}$, for each $i = 1, \ldots, m$.    
Put 
\begin{equation*}
M_v \ = \  \max\big( 
    \max_{1 \le i \le m} \max_{N_i < j \le 2N_i} \| \varphi_{ij}\|_{E_v },
    \max_{1 \le \lambda \le \Lambda} \|\varphi_{\lambda}\|_{E_v } \big) \ . 
\end{equation*} 

Let  $k_v > 0$  be the least integer such that
\begin{equation} \label{DCPF2Cp}
\Big(\frac{h_v}{r_v}\Big)^{k_v N} \cdot \frac{M_v}{\min(1,R_v^N)} \ < \ 1 \ , 
\end{equation}
and let $\kbar \ge k_v$ be a fixed integer.   
Then there are an integer $n_v \ge 1$ and a number $0 < B_v < 1$ 
depending on $\kbar$, $E_v$, and $\phi_v(z)$,
such that for each sufficiently large integer $n$ divisible by $n_v$,  
one can carry out the local patching process at $K_v$ as follows:  

Put $G_v^{(0)}(z) = \phi_v(z)^n$.\index{patching functions, initial $G_v^{(0)}(z)$!construction of}  
Then the leading coefficient 
of $G_v^{(0)}(z)$\index{patching functions, initial $G_v^{(0)}(z)$!leading coefficients of} 
at each $x_i$ is $\tc_{v,i}^n$,
\index{coefficients $A_{v,ij}$!leading} 
$\{ z \in \cC_v(\CC_v) : |G_v^{(0)}(z)|_v \le R_v^{nN}\} = E_v$,  
and when $G_v^{(0)}(z)$\index{patching functions, initial $G_v^{(0)}(z)$!expansion of} 
is expanded in terms of the $L$-rational basis as 
\index{basis!$L$-rational}
\begin{equation*} 
 G_v^{(0)}(z) \ = \ \sum_{i=1}^m \sum_{j=0}^{(n-1)N_i-1} A_{v,ij} \varphi_{i,nN_i-j}(z) 
   + \sum_{\lambda=1}^{\Lambda} A_{v,\lambda} \varphi_{\lambda} \ ,
\end{equation*} 
we have $A_{v,ij} = 0$ for all $(i,j)$ with $1 \le j < \kbar N_i$. 
\index{patching functions, initial $G_v^{(0)}(z)$!when $\Char(K) = p > 0$!high-order coefficients are $0$}  

For each $k$, $1 \le k \le n-1$, 
let $\{\Delta_{v,ij}^{(k)} \in L_{w_v}\}_{(i,j) \in \Band_N(k)}$\index{distinguished place $w_v$} 
be a $K_v$-symmetric set of numbers satisfying 
\index{$K_v$-symmetric!set of numbers}\index{band!$\Band_N(k)$}
\begin{equation} \label{DFCCV1p}
\left\{ \begin{array}{ll}
\text{$|\Delta_{v,i0}^{(1)}|_v \le B_v $ 
            and $\Delta_{v,ij}^{(1)} = 0$ for $j = 1, \ldots, N_i-1$,}  
               & \text{if \ $k = 1$\ ,} \\
\text{$\Delta_{v,ij}^{(k)} = 0$ for $j = (k-1)N_i, \ldots, kN_i - 1$,} 
               & \text{if \ $k = 2, \ldots, \kbar$\ , } \\
|\Delta_{v,ij}^{(k)}|_v  \le  h_v^{kN},  
               & \text{if \ $k = \kbar+1, \ldots, n-1$\ ,}
      \end{array} \right.
\end{equation}
such that $\Delta_{v,i0}^{(1)} \in K_v(x_i)^{\sep}$ for each $i$ 
and such that for each $k = \kbar + 1, \ldots, n-1$ 
\begin{equation} \label{DFCCV1B}
\sum_{i=1}^m \sum_{j=(k-1)N_i}^{kN_i-1} \Delta_{v,ij}^{(k)} \cdot   \varphi_{i,(k+1)N_i - j} 
 \ \in \ K_v(\cC) \ .
\end{equation} 
For $k = n$, let 
$\{\Delta_{v,\lambda}^{(n)} \in L_{w_v}\}_{1 \le \lambda \le \Lambda}$\index{distinguished place $w_v$} 
be a $K_v$-symmetric set of numbers such that
\index{$K_v$-symmetric!set of numbers} 
\begin{equation} \label{DFCCF2np} 
|\Delta_{v,\lambda}^{(n)}|_v \ \le \ h_v^{nN}   
\end{equation} 
for each $\lambda$, and 
\begin{equation} \label{DFCCF2B} 
\sum_{\lambda = 1}^{\Lambda} \Delta_{v,\lambda}^{(n)} \cdot \varphi_{\lambda} \ \in K_v(\cC) \ . 
\end{equation} 

Then one can inductively construct\index{patching functions, $G_v^{(k)}(z)$ for $1 \le k \le n$!constructed by patching} 
 $(\fX,\vs)$-functions $G_v^{(1)}(z), \ldots, G_v^{(n)}(z)$ in $K_v(\cC)$, 
of common degree $Nn$, such that$:$  
 
\smallskip
\noindent{$(A)$} For each $k = 1, \ldots, n$, $G_v^{(k)}(z)$
\index{patching functions, $G_v^{(k)}(z)$ for $1 \le k \le n$!constructed by patching} 
is obtained from $G_v^{(k-1)}(z)$ as follows$:$  

$(1)$  When $k = 1$, there is a $K_v$-symmetric set of functions
\index{$K_v$-symmetric!set of functions} 
$\ttheta_{v,10}^{(1)}(z), \ldots, \ttheta_{v,m0}^{(1)}(z) \in L_{w_v}^{\sep}(\cC)$\index{distinguished place $w_v$}
such that\index{patching functions, initial $G_v^{(0)}(z)$}
\index{patching functions, $G_v^{(k)}(z)$ for $1 \le k \le n$!constructed by patching} 
\begin{equation*}
G_v^{(1)}(z) \ = \ 
    G_v^{(0)}(z) + \sum_{i=1}^m \Delta_{v,i0}^{(1)} \cdot \ttheta_{v,i0}^{(1)}(z) \ , 
\end{equation*}
where for each $i = 1, \ldots, m$,
$\ttheta_{v,i0}^{(1)}(z) \in K_v(x_i)^{\sep}(\cC)$ has the form 
\begin{equation*}
\ttheta_{v,i0}^{(1)}(z) \ = \ \tc_{v,i}^{n} \varphi_{i,nN_i}(z) + \tTheta_{v,i0}^{(1)}(z) 
\end{equation*}   
for an $(\fX,\vs)$-function $\tTheta_{v,i0}^{(1)}(z)$
with a pole of order at most $(n-\kbar)N_i$ at each $x_i$.  
Thus, in passing from $G_v^{(0)}(z)$\index{patching functions, initial $G_v^{(0)}(z)$} 
to $G_v^{(1)}(z)$,\index{patching functions, $G_v^{(k)}(z)$ for $1 \le k \le n$!leading coefficients of} 
each of the leading coefficients
\index{coefficients $A_{v,ij}$!leading} 
$A_{v,i0} = \tc_{v,i}^n$ is replaced with $\tc_{v,i}^n + \Delta_{v,i0}^{(1)} \cdot \tc_{v,i}^n$,
and the coefficients $A_{v,ij}$ for $1 \le j < \kbar N_i$ remain $0$.
\index{coefficients $A_{v,ij}$}    

$(2)$ For $k = 2, \ldots, \kbar$, we have $G_v^{(k)}(z) = G_v^{(k-1)}(z)$.
\index{patching functions, $G_v^{(k)}(z)$ for $1 \le k \le n$!constructed by patching} 

$(3)$ For $k = \kbar + 1, \ldots, n-1 $, we have\index{patching functions, $G_v^{(k)}(z)$ for $1 \le k \le n$!constructed by patching} 
\begin{equation} \label{GvkOmegaModificationRL} 
G_v^{(k)}(z) \ = \ G_v^{(k-1)}(z)
     + \omega_v^{(k)}(z) \cdot F_{v,k}(z) + \Theta_v^{(k)}(z) \ ,    
\end{equation}
where 

\quad $(a)$ $\omega_v^{(k)}(z) = \sum_{(i,j) \in \Band_N(k)} \Delta_{v,ij}^{(k)} \varphi_{i,(k+1)N_i-j}(z)$, 
which belongs to $K_v(\cC)$ by $(\ref{DFCCV1B});$\index{band!$\Band_N(k)$}

\quad $(b)$ $F_{v,k}(z) = \phi_v(z)^{n-k-1}$ is a $K_v$-rational $(\fX,\vs)$-function whose roots belong to $E_v$.

\qquad  For each $x_i$, it has a pole of order $(n-k-1)N_i$ at $x_i$, and its leading coefficient 
\index{coefficients $A_{v,ij}$!leading} 

\qquad $d_{v,i} = \lim_{z \rightarrow x_i} F_{v,k}(z) \cdot g_{x_i}(z)^{(n-k-1)N_i}$
has absolute value $|d_{v,i}|_v = |\tc_{v,i}|_v^{n-k-1}$. 

\quad $(c)$ $\Theta_v^{(k)}(z) \in K_v(\cC)$ is an $(\fX,\vs)$-function 
determined by the local patching process 

\qquad at $v$  after the coefficients in $\Band_N(k)$ have been modified; it has a pole of order
\index{coefficients $A_{v,ij}$}\index{band!$\Band_N(k)$}

\qquad at most $(n-k)N_i$ at each $x_i$ and no other poles, and may be the zero function.  

$(4)$ For $k = n$\index{patching functions, $G_v^{(k)}(z)$ for $1 \le k \le n$!constructed by patching}                              
\begin{equation*}
G_v^{(n)}(z) \ = \ G_v^{(n-1)}(z) + \sum_{\lambda=1}^{\Lambda}
                   \Delta_{v,\lambda}^{(n)} \cdot \varphi_{\lambda}(z) \ . 
\end{equation*}

\vskip .05 in
\noindent{$(B)$} For each $k = 1, \ldots, n$,\index{patching functions, $G_v^{(k)}(z)$ for $1 \le k \le n$!roots are confined to $E_v$} 
\begin{equation*}   
\{ z \in \cC_v(\CC_v) : |G_v^{(n)}(z)|_v \le R_v^{nN} \} \ = \ E_v \ .
\end{equation*}  
\end{theorem}

\smallskip 
Theorems \ref{DCPPatch1} and \ref{DCPPatch1p} will be proved together.  There are some 
differences in the way the leading and high-order coefficients are treated, 
but the underlying\index{coefficients $A_{v,ij}$!leading}\index{coefficients $A_{v,ij}$!high-order}  
patching constructions for the middle and low-order coefficients are the same. 
\index{patching construction!for nonarchimedean $\RL$-domains!differences when $\Char(K) = 0$ and $\Char(K) = p$}  
\index{patching argument!for nonarchimedean $\RL$-domains|(}   
\index{coefficients $A_{v,ij}$!middle}\index{coefficients $A_{v,ij}$!low-order}
\index{compensating functions $\vartheta_{v,ij}^{(k)}(z)$!construction of}   
In Theorem \ref{DCPPatch1} the compensating functions are
$\vartheta_{v,ij}^{(k)}(z) = \varphi_{i,(k+1)N_i - j}(z) \cdot \phi_v(z)^{n-k-1}$, 
while in Theorem \ref{DCPPatch1p} we have 
$\omega_v^{(k)}(z) \ = \ \sum_{(i,j) \in \Band_N(k)} \Delta_{v,ij}^{(k)}  \varphi_{i,(k+1)N_i - j}(z)$ and  
$F_{v,k}(z) = \phi_v(z)^{n-k-1}$,\index{band!$\Band_N(k)$} 
so as noted after Theorem \ref{LocalPatchp},\index{compensating functions $\vartheta_{v,ij}^{(k)}(z)$}   
\begin{equation*}
\sum_{(i,j) \in \Band_N(k)} \Delta_{v,ij}^{(k)} \vartheta_{v,ij}^{(k)}(z) 
\ = \ \omega_v^{(k)}(z) \cdot F_{v,k}(z) \ .
\end{equation*}\index{band!$\Band_N(k)$}  
 
To prove Theorems \ref{DCPPatch1} and \ref{DCPPatch1p}, 
we will need the following nonarchimedean analogue of Lemma \ref{DLemCP1},
which is valid both when $\Char(K_v) = 0$ and when $\Char(K_v) = p > 0$:   

\begin{lemma} \label{DCPL1} 
Let $F(z) \in \CC_v(\cC_v)$ be a nonconstant rational function, and let $R > 0$
be an element of the value group of $\CC_v^{\times}$.  Put
\begin{equation*}
U \ = \ \{ z \in \cC_v(\CC_v) : |F(z)|_v \le R \} \ .
\end{equation*}
Suppose $H(z) \in \CC_v(\cC_v)$ is a function
such that $|H(z)|_v < R$ for all $z \in U$,  and whose polar divisor
satisfies $\div_{\infty}(H) \le \div_{\infty}(F)$.  Then 
\begin{equation*} 
\{ z \in \cC_v(\CC_v) : |F(z) + H(z)|_v \le R \} \ = \ U \ .
\end{equation*}
\end{lemma} 

Lemma \ref{DCPL1} depends on the following nonarchimedean Maximum Modulus Principle:
\index{Maximum principle!nonarchimedean!for $\RL$-domains}\index{distinguished boundary|ii}

\begin{proposition}[Maximum Principle with Distinguished Boundary] \label{DCPP1}
Let $f(z) \in \CC_v(\cC_v)$ be a nonconstant rational function,
and let $R > 0$ belong to the value group of $\CC_v^{\times}$. Put 
\begin{eqnarray*}
U & = & \{ z \in \cC_v(\CC_v) : |f(z)|_v \le R \} \ ,\\
\partial U(f) & = & \{ z \in \cC_v(\CC_v) : |f(z)|_v = R \} \ .
\end{eqnarray*}
Let $g(z) \in \CC_v(\cC_v)$ be a function with no poles in $U$.  
Then $|g(z)|_v$ achieves its maximum for $z \in U$ at a point 
$z_0 \in \partial U(f)$.  
\end{proposition}

\begin{proof} See (\cite{RR1}, Theorem 1.4.2, p.51).
\end{proof} 

\begin{proof}[Proof of Lemma \ref{DCPL1}]  By the ultrametric inequality, if $z \in U$ then 
$|F(z)+H(z)|_v \le R$, so
\begin{equation*}  
U \ \subseteq \ \{ z \in \cC_v(\CC_v) : |F(z) + H(z)|_v \le R \} \ .
\end{equation*}
To establish the reverse containment, put
\begin{equation*}
V \ = \ \{ z \in \cC_v(\CC_v) : |(1/F)(z)|_v \le 1/R \}
\end{equation*}
regarding $1/F$ as a rational function whose value is $0$ on the poles of $F$.  
The distinguished boundaries of $U$ and $V$ satisfy 
\begin{eqnarray*}
\partial U(F) & = & \{ z \in \cC_v(\CC_v) : |F(z)|_v = R \} \\
       & = & \{ z \in \cC_v(\CC_v) : |1/F(z)|_v = 1/R \} \ = \partial V(1/F) \ .
\end{eqnarray*} 
Put $G(z) = H(z)/F(z) \in \CC_v(\cC)$.  By hypothesis, 
$G(z)$ has no poles in $V$, and $|G(z)|_v < 1$ for all
$z \in \partial V(1/F) = \partial U(F) \subset U$.  By Proposition \ref{DCPP1}
$|G(z)|_v < 1$ for all $z \in V$.  Equivalently $|H(z)|_v < |F(z)|_v$ for
all $z \in V$ which are not poles of $F$, so by the ultrametric inequality 
$|F(z)+H(z)|_v = |F(z)|_v$ for such $z$.  In particular, for $z \notin U$
we have $|F(z)+H(z)|_v = |F(z)|_v > R$ so
\begin{equation*}  
\{ z \in \cC_v(\CC_v) : |F(z) + H(z)|_v \le R \} \ \subseteq \ U \ .
\end{equation*}
\end{proof}
 
\begin{proof}[Proof of Theorems \ref{DCPPatch1} and \ref{DCPPatch1p}]

The patching construction will be carried out in three phases.
The proofs differ only in their treatment of the high-order coefficients.  
\index{coefficients $A_{v,ij}$!high-order} 

\medskip
\noindent{\bf Phase 1. Patching the high-order coefficients.}
\index{coefficients $A_{v,ij}$!high-order}\index{patching!high-order coefficients}
    
In this phase we carry out the patching process for stages $k = 1, \ldots, \kbar$.

\smallskip
\index{patching!high-order coefficients!for $\RL$-domains when $\Char(K_v) = 0$}
\index{local patching for nonarch $\RL$-domains!Phase 1: high-order coefficients!when $\Char(K_v) = 0$|(} 
First assume $\Char(K_v) = 0$.  
Let the $(\fX,\vs)$-function $\phi_v(z)$ of degree $N$, and the numbers 
$k_v$, $M_v$, $0 < h_v < r_v \le R_v$, and $\kbar \ge k_v$ be as in Theorem \ref{DCPPatch1}.  
Take $n_v = 1$, and put 
\begin{equation*} 
B_v \ = \ 
\frac{  \min(1,R_v)^{(\kbar+1) N} } {2 M_v } \ .
\end{equation*}   
Assume $n > \kbar$, 
and let $G_v^{(0)}(z) = \phi_v(z)^n$.\index{patching functions, initial $G_v^{(0)}(z)$!construction of} 

For each $k = 1, \ldots, \kbar$, 
we begin with a $K_v$-rational $(\fX,\vs)$-function $G_v^{(k-1)}(z)$
\index{patching functions, $G_v^{(k)}(z)$ for $1 \le k \le n$!roots are confined to $E_v$}
\index{patching functions, $G_v^{(k)}(z)$ for $1 \le k \le n$!are $K_v$-rational} 
satisfying
\begin{equation*}
 \{ z \in \cC_v(\CC_v) : |G_v^{(k-1)}(z)|_v \le 1 \} \ = \ E_v \ .
\end{equation*}
Expand\index{patching functions, $G_v^{(k)}(z)$ for $1 \le k \le n$!expansion of} 
\begin{equation*}
G_v^{(k-1)}(z) \ = \
     \sum_{i=1}^m \sum_{j=0}^{(n-1)N_i-1} A_{v,ij} \varphi_{i,nN_i-j}(z) 
   + \sum_{\lambda=1}^{\Lambda} A_{v,\lambda} \varphi_{\lambda} \ .
\end{equation*}   
We will patch the coefficients
\index{coefficients $A_{v,ij}$!high-order} 
$A_{v,ij}$ with $(k-1)N_i \le j < kN_i$ in $\prec_N$\index{order!$\prec_N$} order.
Given $(i,j)$ we can uniquely write 
\begin{equation*}
nN_i-j \ = \ r_{ij} + (n-k-1)N_i, \quad \text{with $N_i + 1 \le r_{ij} \le 2N_i$\ ;}
\end{equation*}
thus $r_{ij} = (k+1)N_i - j$.  Put\index{compensating functions $\vartheta_{v,ij}^{(k)}(z)$!construction of}  
\begin{equation} \label{FvTChoice}
\vartheta_{v,ij}^{(k)}(z) \ = \   \varphi_{i,r_{ij}} \cdot \phi_v(z)^{n-k-1} 
                    \ = \   \varphi_{i,(k+1)N_i - j} \cdot \phi_v(z)^{n-k-1}\ . 
\end{equation}
Then $\vartheta_{v,ij}^{(k)}(z)$
\index{compensating functions $\vartheta_{v,ij}^{(k)}(z)$!poles and leading coefficients of}  
has pole of order $nN_i-j$ at $x_i$, with leading coefficient
\index{coefficients $A_{v,ij}$!leading} 
$\tc_{v,i}^{n-k-1}$, and a pole of order at most $(n-k-1)N_{i^{\prime}}$ at each
$x_{i^{\prime}} \ne x_i$.  
As a collection, the $\vartheta_{v,ij}^{(k)}(z)$ 
are $K_v$-symmetric.\index{compensating functions $\vartheta_{v,ij}^{(k)}(z)$!are $K_v$-symmetric} 
\index{$K_v$-symmetric!set of functions}
By the definition of $M_v$, for each $(i,j)$ 
\begin{equation} \label{DCPF7A}
\|\vartheta_{v,ij}^{(k)}(z)\|_{E_v} \ \le \ M_v R_v^{(n-k-1)N} \ .
\end{equation}

Since $1 \le k \le \kbar$, the definition of $B_v$ shows that $B_v M_v < R_v^{(k+1)N}$.
By hypothesis, $\{\Delta_{v,ij}^{(k)} \in L_{w_v}\}_{(i,j) \in \Band_N(k)}$
\index{band!$\Band_N(k)$}\index{distinguished place $w_v$} 
is a $K_v$-symmetric collection of numbers such that $|\Delta_{v,ij}^{(k)}|_v \le B_v$
\index{$K_v$-symmetric!set of functions} 
for each $(i,j)$.  Thus\index{compensating functions $\vartheta_{v,ij}^{(k)}(z)$!bounds for}  
\begin{equation} \label{FM23} 
\|\Delta_{v,ij}^{(k)} \vartheta_{v,ij}^{(k)}(z) \|_{E_v} \ \le \ B_v M_v R_v^{(n-k-1)N} \ < \ R_v^{nN} \ .
\end{equation} 

Put\index{patching functions, $G_v^{(k)}(z)$ for $1 \le k \le n$!constructed by patching} 
\begin{equation} \label{DCPF8}
G_v^{(k)}(z) \ = \ G_v^{(k-1)}(z)
     + \sum_{i=1}^m \sum_{j=(k-1)N_i}^{kN_i-1} \Delta_{v,ij}^{(k)} \cdot \vartheta_{v,ij}(z) \ .
\end{equation}

Let $H(z)$ denote the sum on the right in (\ref{DCPF8}). 
$H(z)$ is $K_v$-rational,
since $\Char(K_v) = 0$ and the $\tc_{v,i}$, $\Delta_{v,ij}^{(k)}$ 
and $\vartheta_{v,ij}(z)$\index{compensating functions $\vartheta_{v,ij}^{(k)}(z)$!are $K_v$-symmetric}  
are defined over $L_{w_v}$\index{distinguished place $w_v$} and are $K_v$-symmetric.
\index{$K_v$-symmetric!set of functions}  
It follows that $G_v^{(k)}(z)$\index{patching functions, $G_v^{(k)}(z)$ for $1 \le k \le n$!are $K_v$-rational}  
is a $K_v$-rational $(\fX,\vs)$-function of degree $nN$.  
Clearly $\div_{\infty}(H) \le \div_{\infty}(G_v^{(k)})$.\index{patching functions, $G_v^{(k)}(z)$}
By (\ref{FM23}) and the ultrametric inequality, we have $|H(z)|_v < 1$ for each $z \in E_v$.  
Hence by Lemma \ref{DCPL1},\index{patching functions, $G_v^{(k)}(z)$ for $1 \le k \le n$!roots are confined to $E_v$}  
\begin{equation*}
\{ z \in \cC_v(\CC_v) : |G_v^{(k)}(z)|_w \le R_v^{nN} \} \ = \ E_v  \ .
\end{equation*}
This completes the patching process for the high-order coefficients when $\Char(K_v) = 0$.
\index{coefficients $A_{v,ij}$!high-order}
\index{local patching for nonarch $\RL$-domains!Phase 1: high-order coefficients!when $\Char(K_v) = 0$|)}  

\medskip 
Next assume $\Char(K_v) = p > 0$.
\index{patching!high-order coefficients!for $\RL$-domains when $\Char(K_v) = p > 0$}
\index{local patching for nonarch $\RL$-domains!Phase 1: high-order coefficients!when $\Char(K_v) = p > 0$|(}
\index{patching functions, initial $G_v^{(0)}(z)$!when $\Char(K) = p > 0$!high-order coefficients are $0$|(} 

Let the $(\fX,\vs)$-function $\phi_v(z)$ of degree $N$, and the numbers 
$k_v$, $M_v$, $0 < h_v < r_v \le R_v$, and $\kbar \ge k_v$ be as in Theorem \ref{DCPPatch1p}.   
Let $J = p^A \ge \max(2g+1, \max_i([K_v(x_i):K_v]^{\insep}))$
be the number from \S\ref{Chap3}.\ref{LRationalBasisSection}  
in the construction of the $L$-rational and $L^{\sep}$-rational bases.  
By assumption, $J|N_i$ for each $i$.
The leading coefficient of $\phi_v(z)$ at $x_i$ is
\index{coefficients $A_{v,ij}$!leading}  
$\tc_{v,i} = \lim_{z \rightarrow x_i} \phi_v(z) \cdot g_{x_i}(z)^{N_i}$;  
by hypothesis, $\tc_{v,i} \in K_v(x_i)^{\sep}$ for each $i$.  

\smallskip
We will choose $n_v$ and $B_v$ differently from the way they were chosen when $\Char(K_v) = 0$. 
Recall that $k_v$ is the least integer for which 
\begin{equation*}
\big(\frac{h_v}{R_v}\big)^{k_v} \cdot M_v \ < \ 1 \ ,
\end{equation*}
and that $\kbar \ge k_v$ is a a fixed integer (specified by the global patching process).
\index{patching argument!global}  
Let $n_v = p^r$ be the least power of $p$ for which 
\begin{equation} \label{F1pNvDef}
p^r \ \ge \ \max(\kbar N_1, \ldots, \kbar N_m)\ ,
\end{equation}
and let   
\begin{equation} \label{F1pBvDef1}
B_v \ = \ 
\frac{R_v^{n_v N}}{2 \max_i \big(|\tc_{v,i}^{n_v}|_v \|\varphi_{i,n_v N_i}\|_{E_v} \big)} \ . 
\end{equation}

\smallskip 
We now show that for all sufficiently large $n$ divisible by $n_v$, 
we can carry out the patching process imposing the conditions in the Theorem.
In fact, we will see that we can take $n = n_v \cdot Q = p^r \cdot Q$ 
for any integer $Q > \max(3,\kbar)$.

Given such an $n$, 
put $G_v^{(0)}(z) = \phi_v(z)^n$.\index{patching functions, initial $G_v^{(0)}(z)$!construction of}  
We first show that $G_v^{(0)}(z)$\index{patching functions, initial $G_v^{(0)}(z)$} 
has the properties in the Theorem.
It is clear that 
$E_v = \{z \in \cC_v(\CC_v) : |G_v^{(0)}(z)|_v \le R_v^{nN}\}$.
\index{patching functions, initial $G_v^{(0)}(z)$!mapping properties of}
For each $i$, the leading coefficient of 
$G_v^{(0)}(z)$ at $x_i$ is $\tc_{v,i}^n$,\index{patching functions, initial $G_v^{(0)}(z)$!leading coefficients of} 
\index{coefficients $A_{v,ij}$!leading} 
which belongs to $K_v(x_i)^{\sep}$ since $\tc_{v,i} \in K_v(x_i)^{\sep}$.     
If we  expand $\phi_v(z)^Q$ using the $L$-rational basis as 
\index{basis!$L$-rational}
\begin{equation*}
\phi_v(z)^Q  \ = \ \sum_{i=1}^m \sum_{j=0}^{(Q-1)N_i - 1} B_{v,ij} \varphi_{i,QN_i -j}(z) 
                \ + \ \sum_{\lambda = 1}^{\Lambda} B_{\lambda} \varphi_{\lambda}(z) \ ,
\end{equation*} 
then since $\Char(K_v) = p > 0$ and $n_v = p^r$ 
it follows that\index{patching functions, initial $G_v^{(0)}(z)$!expansion of} 
\begin{equation} \label{FGv0exp1}
G_v^{(0)}(z) \ = \ (\phi_v(z)^Q)^{n_v} 
  \ = \ \sum_{i=1}^m \sum_{j=0}^{(Q-1)N_i - 1} B_{v,ij}^{n_v} \varphi_{i,QN_i -j}(z)^{n_v} 
                \ + \ \sum_{\lambda = 1}^{\Lambda} B_{\lambda}^{n_v} \varphi_{\lambda}^{n_v} \ .
\end{equation}
Since $J|N_i$ for each $i$, Proposition \ref{TransitionProp}(B) shows that 
\begin{equation*} 
\varphi_{i,QN_i}(z)^{n_v} \ = \ \varphi_{i,nN_i}(z) 
\end{equation*}
belongs to the $L$-rational basis.
\index{basis!$L$-rational}  
On the other hand if $D$ is the divisor $\sum_{i=1}^m N_i(x_i)$, then since $n_v \ge \kbar N_i$
for each $i$,  and  since $Q \ge 2$, 
all other terms in the expansion (\ref{FGv0exp1}) belong to $\Gamma_{\CC_v}((n-\kbar)D)$.  
Since $J|N_i$ for each $i$, this means that 
when we expand $G_v^{(0)}(z)$\index{patching functions, initial $G_v^{(0)}(z)$!expansion of} 
in terms of the $L$-rational basis as\index{basis!$L$-rational}
\begin{equation} \label{FtGExp1} 
 G_v^{(0)}(z) \ = \ \sum_{i=1}^m \sum_{j=0}^{(n-1)N_i-1} A_{v,ij} \varphi_{i,nN_i-j}(z) 
   + \sum_{\lambda=1}^{\Lambda} A_{v,\lambda} \varphi_{\lambda} \ ,
\end{equation} 
then $A_{v,ij} = 0$ for all $(i,j)$ with $1 \le j < \kbar N_i$.  
\index{patching functions, initial $G_v^{(0)}(z)$!when $\Char(K) = p > 0$!high-order coefficients are $0$|)}

\smallskip
We next carry out the patching process for the stage $k = 1$.   
We want to modify the leading coefficients $A_{v,i0}$  
\index{coefficients $A_{v,ij}$!leading} 
and leave the remaining-order high coefficients $A_{v,ij}$ for $1 \le j < \kbar N_i$
\index{coefficients $A_{v,ij}$!high-order} 
(which are $0$) unchanged.   
By assumption, we are given a $K_v$-symmetric set of numbers 
\index{$K_v$-symmetric!set of numbers}
$\{\Delta_{v,ij}^{(1)}  \in L_{w_v}\}_{1 \le i \le m, 0 \le j < N_i}$,\index{distinguished place $w_v$} 
with $|\Delta_{v,i0}^{(1)}|_v \le B_v$ for each $i$,  
and $\Delta_{v,ij}^{(1)} = 0$ for all $j \ge 1$ and all $i$, 
such that $\Delta_{v,i0}^{(1)}$ belongs to $K_v(x_i)^{\sep}$ for each $i$, 
and we wish to replace $A_{v,i0} = \tc_{v,i}^n$ with $\tc_{v,i}^n + \Delta_{v,i0} \tc_{v,i}^n$ 
in (\ref{FtGExp1}). 

Recall that $n = n_v Q$.  We claim that setting
\begin{equation} \label{FBunttheta}
\ttheta_{v,i0}^{(1)}(z) \ = \ \tc_{v,i}^{n_v} \varphi_{i,n_v N_i} \cdot \phi_v(z)^{n_v(Q-1)} 
\end{equation} 
in Theorem \ref{DCPPatch1p} for each $i = 1, \ldots, m$,
and then putting\index{patching functions, $G_v^{(k)}(z)$ for $1 \le k \le n$!constructed by patching} 
\begin{equation} \label{FLCSum1}
G_v^{(1)}(z) \ = \ G_v^{(0)}(z) \ + \ 
 \sum_{i=1}^m \Delta_{v,i0}^{(1)} \, \ttheta_{v,i0}^{(1)}(z) \ ,
\end{equation} 
accomplishes what we need.  
Let $H(z)$ denote the sum on the right side of (\ref{FLCSum1}). 

First, adding $H(z)$ to $G_v^{(0)}(z)$\index{patching functions, $G_v^{(k)}(z)$ for $1 \le k \le n$!constructed by patching} 
adds $\Delta_{v,i0}^{(1)} \tc_{v,i}^n$ to $A_{v,i0}$, for each $i$. 
This follows from the fact that 
$\tc_{v,i}^{n_v} \varphi_{i,n_v N_i} \cdot \phi_v(z)^{n_v(Q-1)}$ has a pole of order $nN_i$ at $x_i$ 
with leading coefficient $\tc_{v,i}^n$, and at each $x_{i^{\prime}} \ne x_i$
\index{coefficients $A_{v,ij}$!leading}  
its pole has order less than $(n-\kbar)N_{i^{\prime}}$. 
 
Second, adding $H(z)$ to $G_v^{(0)}(z)$\index{patching functions, initial $G_v^{(0)}(z)$} 
leaves $A_{v,ij} = 0$ for $1 \le j < \kbar N_i$.  
This follows by considering an expansion of  $\varphi_{i,n_v N_i}  \cdot \big(\phi_v(z)^{(Q-1)}\big)^{n_v}$ 
like the one in (\ref{FGv0exp1}):  if we write $\phi_v(z)^{Q-1}$  as  
\begin{equation*}
\phi_v(z)^{Q-1}  \ = \ \sum_{\ell=1}^m \sum_{j=0}^{(Q-2)N_\ell - 1} 
       C_{v,\ell j} \varphi_{\ell,QN_\ell -j}(z) 
                \ + \ \sum_{\lambda = 1}^{\Lambda} C_{\lambda} \varphi_{\lambda}(z) \ ,
\end{equation*} 
then 
\begin{eqnarray}
& & \ttheta_{v,i0}^{(1)}(z) \ = \ \varphi_{i,n_v N_i}  \cdot \big( \phi_v(z)^{Q-1} \big)^{n_v} \notag \\
& & \qquad  \ = \ \sum_{\ell=1}^m \sum_{j=0}^{(Q-2)N_\ell - 1} 
        C_{v,\ell j}^{n_v} \cdot \varphi_{i,n_v N_i} \cdot \varphi_{\ell,(Q-1)N_\ell -j}(z)^{n_v} 
                \ + \ \sum_{\lambda = 1}^{\Lambda} C_{\lambda}^{n_v} \cdot 
                   \varphi_{i,n_v N_i} \cdot \varphi_{\lambda}^{n_v} \ , \label{FphiQm1p1}
\end{eqnarray} 
and since $n_v \ge \kbar N_i$ and $Q \ge 3$, all the terms in (\ref{FphiQm1p1}) besides 
the one with $(\ell,j) = (i,0)$ belong to $\Gamma_{\CC_v}((n-\kbar)D)$.  
On the other hand, since $J|N_i$, Proposition \ref{TransitionProp}(B) shows that 
that term coincides with $\tc_{v,i}^{n_v(Q-1)} \varphi_{i,nN_i}(z)$.  
Note that $\varphi_{i,n_v N_i} = \tphi_{i,n_v N_i}$ is rational over $K_v(x_i)^{\sep}$ 
by Proposition \ref{TransitionProp}(B), since $J|N_i$.  Since the 
$\tc_{v,i}$ and $\varphi_{i,n_v N_i}$ are $K_v$-symmetric, and $\phi_v(z)^{n_v(Q-1)}$
\index{$K_v$-symmetric!set of functions}
is $K_v$-rational, the $\ttheta_{v,i0}^{(1)}(z)$ are $K_v$-symmetric.  
\index{$K_v$-symmetric!set of functions}
This discussion also shows that 
each $\ttheta_{v,i0}^{(1)}(z)$ belongs to $K_v(x_i)^{\sep}(\cC)$ and has the form 
$\ttheta_{v,i0}^{(1)}(z) = \tc_{v,i}^n \varphi_{i,nN_i}(z) + \tTheta_{v,i0}(1)(z)$
for an $(\fX,\vs)$-function $\tTheta_{v,i0}^{(1)}(z)$ with poles of order at most $(n-\kbar)N_{i^\prime}$ 
at each $x_{i^\prime}$, as asserted in Theorem \ref{DCPPatch1p}.

Third, $G_v^{(1)}(z)$ is $K_v$-rational.\index{patching functions, $G_v^{(k)}(z)$ for $1 \le k \le n$!are $K_v$-rational} 
Indeed, the  $\tc_{v,i}^{n_v}$ and $\ttheta_{v,i0}^{(1)}(z)$ 
are $L_{w_v}^{\sep}$-rational\index{distinguished place $w_v$} and $K_v$-symmetric, so $H(z)$ is $K_v$-rational.
\index{$K_v$-symmetric!set of functions}  
Since $G_v^{(0)}(z)$ is $K_v$-rational,\index{patching functions, initial $G_v^{(0)}(z)$!are $K_v$-rational} 
so is $G_v^{(1)}(z)$. 

Finally, $E_v = \{z \in \cC_v(\CC_v) : |G_v^{(1)}(z)|_v \le R_v^{nN}\}$.
\index{patching functions, $G_v^{(k)}(z)$ for $1 \le k \le n$!roots are confined to $E_v$}  
To see this, note that our choice  of $B_v$ in (\ref{F1pBvDef1}), and the fact that 
$|\Delta_{v,i0}^{(1)}|_v \le B_v$ for each $i$, means that 
$\|H(z)\|_{E_v} \le \frac{1}{2} R_v^{nN}$.  
Hence the claim follows by applying Lemma \ref{DCPL1}
 to $F(z) = G_v^{(0)}(z)$ and $H(z)$,\index{patching functions, initial $G_v^{(0)}(z)$!mapping properties of} 
taking $R = R_v^{nN}$.

\smallskip
For $k = 2, \ldots, \kbar$, we have $\Delta_{v,ij}^{(k)} = 0$ for all $(i,j)$,  
and we take $G_v^{(k)}(z) = G_v^{(k-1)}(z)$.\index{patching functions, $G_v^{(k)}(z)$ for $1 \le k \le n$!constructed by patching}
\index{local patching for nonarch $\RL$-domains!Phase 1: high-order coefficients!when $\Char(K_v) = p > 0$|)}    

\vskip .1 in
\noindent{\bf Phase 2. Patching the middle coefficients.}
\index{coefficients $A_{v,ij}$!middle}\index{patching!middle coefficients}
\index{local patching for nonarch $\RL$-domains!Phase 2: middle coefficients|(}   

In this phase we carry out the patching process for $k = \kbar+1, \ldots, n-1$.  
The construction is the same regardless of $\Char(K_v)$, 
and coincides with the one in Phase 1 when $\Char(K_v) = 0$, 
except that for each $k$, 
instead of  $|\Delta_{v,ij}^{(k)}|_v \le B_v$ we have 
\begin{equation*}
|\Delta_{v,ij}^{(k)}|_v \le h_v^{kN} \ .
\end{equation*}
Since $k > k_v$, if we take 
$\vartheta_{v,ij}^{(k)}(z) \ = \   \varphi_{i,(k+1)N_i - j} \cdot \phi_v(z)^{n-k-1}$,  
then by condition (\ref{DCPF2C}) (resp. condition (\ref{DCPF2Cp}))
\begin{equation*}
\|\Delta_{v,ij}^{(k)} \vartheta_{v,ij}(z) \|_{U_v} \ \le \ h_v^{kN} M_v R_v^{(n-k-1)N} 
\ \le \ \Big(\frac{h_v}{r_v}\Big)^{kN} \frac{M_v}{R_v^N} \cdot R_v^{nN} \ < \ R_v^{nN} \ .
\end{equation*} 
Hence if\index{compensating functions $\vartheta_{v,ij}^{(k)}(z)$}
\index{patching functions, $G_v^{(k)}(z)$ for $1 \le k \le n$!constructed by patching}  
\begin{equation} \label{FMidSum}
G_v^{(k)}(z) \ = \ G_v^{(k-1)}(z)
     + \sum_{i=1}^m \sum_{j=(k-1)N_i}^{kN_i-1} \Delta_{v,ij}^{(k)} \cdot \vartheta_{v,ij}(z) \ ,
\end{equation}
then as before, by Lemma \ref{DCPL1}\index{patching functions, $G_v^{(k)}(z)$ for $1 \le k \le n$!roots are confined to $E_v$}  
\begin{equation*}
E_v \ = \ \{ z \in \cC_v(\CC_v) : |G_v^{(k)}(z)|_v \le R_v^{nN} \}  \ .
\end{equation*}

When $\Char(K_v) = 0$, the sum on the right in (\ref{FMidSum}) is $K_v$-rational 
for the same reasons of $K_v$-symmetry as in Phase 1.  
When $\Char(K_v) = p > 0$, it can be written as 
\begin{equation*}
\Big( \sum_{i=1}^m \sum_{j=(k-1)N_i}^{kN_i-1} \Delta_{v,ij}^{(k)} \cdot \varphi_{v,(k+1)N_i - j}(z) \Big) \cdot \phi_v(z)^{n-k-1} 
\end{equation*} 
which is $K_v$-rational by assumption (\ref{DFCCV1B}).  
Thus $G_v^{(k)}(z)$ is $K_v$-rational.\index{patching functions, $G_v^{(k)}(z)$ for $1 \le k \le n$!are $K_v$-rational}
\index{local patching for nonarch $\RL$-domains!Phase 2: middle coefficients|)}    

\medskip
\noindent{\bf Phase 3. Patching the low-order coefficients.}
\index{coefficients $A_{v,ij}$!low-order}\index{patching!low-order coefficients}
\index{local patching for nonarch $\RL$-domains!Phase 3: low-order coefficients|(}   

In the final step we take\index{patching functions, $G_v^{(k)}(z)$ for $1 \le k \le n$!constructed by patching}
\begin{equation*}
G_v^{(n)}(z) \ = \ G_v^{(n-1)}(z)
          + \sum_{\lambda = 1}^{\Lambda} \Delta_{v,\lambda}^{(n)} \varphi_{\lambda}
\end{equation*}
with $K_v$-symmetric $\Delta_{v,\lambda} \in L_{w_v}$\index{distinguished place $w_v$} 
\index{$K_v$-symmetric!set of numbers}
satisfying $|\Delta_{v,\lambda}|_v \le h_v^{nN}$ for each $\lambda$.  
When $\Char(K_v) = 0$, the sum $\sum_{\lambda = 1}^{\Lambda} \Delta_{v,\lambda}^{(n)} \varphi_{\lambda}$ 
is $K_v$-rational by the $K_v$-symmetry of the $\Delta_{v,\lambda}^{(n)}$ 
and the $\varphi_{\lambda}$. 
When $\Char(K_v) = p > 0$, it is $K_v$-rational by assumption (\ref{DFCCF2B}).
Thus $G_v^{(n)}(z)$ is $K_v$-rational.\index{patching functions, $G_v^{(k)}(z)$ for $1 \le k \le n$!are $K_v$-rational}

Since $n > k_v$, condition (\ref{DCPF2C}) (resp. condition (\ref{DCPF2Cp}))
shows that for each $\lambda$
\begin{equation*}
\|\Delta_{v,ij}^{(n)} \varphi_{\lambda}(z) \|_{E_v} 
   \ \le \ h_v^{nN} \cdot M_v \ \le \ (h_v/r_v)^{nN} M_v \cdot R_v^{nN} \ < \ R_v^{nN} \ . 
\end{equation*}
Hence by Lemma \ref{DCPL1}, as before,\index{patching functions, $G_v^{(k)}(z)$ for $1 \le k \le n$!roots are confined to $E_v$}
\begin{equation*}
E_v \ = \ \{ z \in \cC(\CC_v) : |G_v^{(n)}(z)|_v \le R_v^{nN} \} \ .
\end{equation*}
This completes the proof.  
\end{proof}\index{local patching for nonarch $\RL$-domains!Phase 3: low-order coefficients|)}  

\index{patching argument!local!for nonarchimedean $\RL$-domains|)}

%% file: NewFSZChap11.tex
\chapter{The Local Patching Construction for Nonarchimedean $K_v$-simple Sets }
\label{Chap11}

In this section we give the confinement argument\index{confinement argument} for Theorem \ref{aT1-B} when  
$K_v$ is nonarchimedean, and $E_v \subset \cC_v(\CC_v)$ is $K_v$-simple, hence compact. 
\index{$K_v$-simple!set}\index{patching argument!local!for nonarchimedean $K_v$-simple sets|(}
This construction is the most intricate of the four confinement arguments, 
and uses results from Appendices \ref{AppC} and \ref{AppD}.  It breaks into two cases,
when $\Char(K_v) = 0$ and when $\Char(K_v) = p > 0$. 

Let $q_v$ be the order of the residue field of $K_v$. 
Let $w_v$ be the distinguished place\index{distinguished place $w_v$} of $L = K(\fX)$ 
determined by the embedding $\tK \hookrightarrow \CC_v$
used to identify $\fX$ with a subset of $\cC_v(\CC_v)$, 
and view $L_{w_v}$ as a subset of $\CC_v$.  Following the construction
of the coherent approximating functions in Theorems \ref{CTCX2} and \ref{CTCX2p},
\index{coherent approximating functions $\phi_v(z)$} 
we begin with the following data: 
    
\begin{enumerate}
  \item A $K_v$-symmetric probability vector $\vs \in \cP^{m}(\QQ)$
\index{$K_v$-symmetric!probability vector} 
           with positive rational coefficients.   

  \item A compact, $K_v$-simple set $E_v \subset \cC_v(\CC_v)$ equipped with a $K_v$-simple decomposition 
  \index{$K_v$-simple!set}
  \index{$K_v$-simple!decomposition}
           $E_v = \bigcup_{\ell = 1}^{D_v} \big( B(a_\ell,r_\ell) \cap \cC_v(F_{w_\ell}) \big)$
           such that $U_v := \bigcup_{\ell = 1}^{D_v}  B(a_\ell,r_\ell)$ is disjoint from $\fX$.
           (To ease notation, we henceforth write $D$ for $D_v$.) In particular  

            $(a)$ the balls $B(a_1,r_1), \ldots, B(a_D,r_D)$ 
           are pairwise disjoint and isometrically parametrizable; 
           each $B(a_\ell,r_\ell)$ is disjoint from $\fX$.
           
            $(b)$ the collection of balls $\{B(a_\ell,r_{\ell})\}_{1 \le \ell \le D}$
           is stable under $\Aut_c(\CC_v/K_v)$.   For each $\sigma \in \Aut_c(\CC_v/K_v)$,   
           if $\sigma(B(a_j,r_j)) = B(a_k,r_k)$, then $\sigma(F_{w_j}) = F_{w_k}$.
           For each $\ell$, $F_{w_\ell}$ is a finite separable extension of $K_v$, 
           $a_{\ell} \in \cC_v(F_{w_\ell})$, $r_{\ell} \in |F_{w_\ell}^{\times}|_v$, and 
           $B(a_\ell,r_\ell)$ has exactly $[F_{w_\ell}:K_v]$ conjugates under $\Aut_c(\CC_v/K_v)$.

  \item A number $N$  and an $(\fX,\vs)$-function $\phi_v(z) \in K_v(\cC)$ 
           of degree $N$ such that 

          $(a)$ the zeros $\theta_1, \ldots, \theta_N$ of $\phi_v(z)$ are distinct and belong to $E_v$;  
          
          $(b)$ $\phi_v^{-1}(D(0,1)) = \bigcup_{h=1}^N B(\theta_h,\rho_h)$, 
             where $B(\theta_1,\rho_1), \ldots, B(\theta_N,\rho_N)$ 
             are pairwise disjoint, isometrically parametrizable, 
             and contained in $\bigcup_{\ell = 1}^D B(a_{\ell},r_{\ell})$. 
           
          $(c)$ $H_v := \phi_v^{-1}(D(0,1)) \cap E_v$ is $K_v$-simple, 
          \index{$K_v$-simple!set}
          \index{$K_v$-simple!decomposition}
            with a $K_v$-simple decomposition 
            \begin{equation*}             
            H_v \ = \ \bigcup_{h=1}^N \big( B(\theta_h,\rho_h) \cap \cC_v(F_{u_h}) \big)
            \end{equation*}         
           which is compatible with the $K_v$-simple decomposition
 \index{$K_v$-simple!decomposition!compatible with another decomposition}  
           $\bigcup_{\ell = 1}^D \big( B(a_\ell,r_\ell) \cap \cC_v(F_{w_\ell}) \big)$ of $E_v$, 
           and move-prepared (see Definition \ref{MovePreparedDef}) relative to 
  \index{move-prepared} 
           $B(a_1,r_1), \ldots, B(a_D,r_D)$. For each $\ell = 1, \ldots, D$, there is a point 
            $\wbar_\ell \in \big(B(a_\ell,r_\ell) \cap \cC_v(F_{w_\ell})) \big) \backslash H_v$.
            
           
         $(d)$ For each $h = 1, \ldots, N$, $\phi_v$ induces 
           an $F_{u_h}$-rational scaled isometry from $B(\theta_h,\rho_h)$ onto $D(0,1)$ 
           with $\phi_v(\theta_h)= 0$,
           which takes $B(\theta_h,\rho_h) \cap \cC_v(F_{u_h})$ onto $\cO_{u_h}$.

   \item  Put $N_i = Ns_i \in \NN$ for each $i = 1, \ldots, m$.
       
        If $\Char(K_v) = p > 0$, the  number $J$ from the construction of the $L$-rational 
        and $L^{\sep}$-rational bases\index{basis!$L$-rational}\index{basis!$L^{\sep}$-rational}
        in \S\ref{Chap3}.\ref{LRationalBasisSection} divides $N_i$, 
        and the leading coefficient 
   \index{coefficients $A_{v,ij}$!leading} 
        $\tc_{v,i} = \lim_{z \rightarrow x_i} \phi_v(z) \cdot g_{x_i}(z)^{N_i}$ of $\phi_v(z)$ at $x_i$ 
        belongs to $K_v(x_i)^{\sep}$, for each $i = 1, \ldots, m$.  
            
   \item Parameters $0 < h_v < r_v < 1$ satisfying $h_v^N < r_v^N  < q_v^{-1/(q_v-1)} < 1$.\index{patching parameters} 
           
    \item  An order $\prec_N$\index{order!$\prec_N$} on the index set 
              $\cI = \{(i,j) \in \ZZ^2 : 1 \le i \le m, 0 \le j\}$ 
           determined by $N$ and $\vs$ as in (\ref{FPrec}), 
           which gives the sequence in which coefficients are patched.
      \index{coefficients $A_{v,ij}$} 
\end{enumerate} 

\smallskip
We will use the $L$-rational basis  
\index{basis!$L$-rational}
$\{\varphi_{ij}, \varphi_{\lambda}\}$ 
from \S\ref{Chap3}.\ref{LRationalBasisSection} to expand functions, 
and $\Lambda = \dim_K(\Gamma(\sum_{i=1}^m N_i(x_i)))$ 
will be the number of low-order basis elements, 
as in \S\ref{Chap7}.\ref{Char0Section}.  
The order $\prec_N$\index{order!$\prec_N$} respects the $N$-bands\index{band!$\Band_N(k)$} (\ref{FBand}), 
and for each $x_i \in \fX$, specifies the terms to be patched in decreasing pole order. 

As in Definition \ref{BGD1}, let the Stirling polynomial of degree $n$ for $\cO_v$ be 
\index{Stirling polynomial!for $\cO_v$}
\begin{equation} \label{FStirling} 
S_{n,v}(x)  \ = \ \prod_{j=0}^{n-1} (x-\psi_v(j)) \ ,
\end{equation}
where $\{\psi_v(j)\}_{0 \le j < \infty}$ is the basic well-distributed
sequence in $\cO_v$.  

\smallskip
When $\Char(K_v) = 0$, we will prove the following patching theorem.

\begin{theorem}  \label{DCPCPatch}\index{patching theorem!for nonarchimedean $K_v$-simple sets!when $\Char(K_v) = 0$} 
Suppose $K_v$ is nonarchimedean, and $\Char(K_v) = 0$.   
Let $E_v$ be $K_v$-simple, with a $K_v$-simple decomposition 
\index{$K_v$-simple!set}
\index{$K_v$-simple!decomposition}
$E_v = \bigcup_{\ell = 1}^D \big( B(a_\ell,r_\ell) \cap \cC_v(F_{w_\ell}) \big)$  
such that $U_v := \bigcup_{\ell=1}^D  B(a_{\ell},r_{\ell})$ is disjoint from $\fX$. 
Let $\vs \in \cP^m(\QQ)$ be a $K_v$-symmetric probability vector with positive coefficients,
\index{$K_v$-symmetric!probability vector} 
and let $\phi_v(z) \in K_v(\cC)$ be an $(\fX,\vs)$-function of degree $N$ 
with distinct zeros $\theta_1, \ldots, \theta_N \in E_v$, such that

\vskip .05 in
  $(1)$ $\phi_v^{-1}(D(0,1)) = \bigcup_{k=1}^N B(\theta_k,\rho_k)$, 
             where the balls $B(\theta_1,\rho_1), \ldots, B(\theta_N,\rho_N)$ 
             are pairwise disjoint, isometrically parametrizable, 
             and contained in $\bigcup_{\ell = 1}^D B(a_{\ell},r_{\ell})$ $;$ 
         
  $(2)$ The set $H_v := \phi_v^{-1}(D(0,1)) \cap E_v$ is $K_v$-simple, 
  \index{$K_v$-simple!decomposition}
  \index{$K_v$-simple!set}
            with a $K_v$-simple decomposition 
            \begin{equation*} 
            H_v \ = \ \bigcup_{h=1}^N \big( B(\theta_h,\rho_h) \cap \cC_v(F_{u_h}) \big)           
            \end{equation*}\index{$K_v$-simple!decomposition!compatible with another decomposition}  
           compatible with the $K_v$-simple decomposition 
           $\bigcup_{\ell = 1}^D \big( B(a_\ell,r_\ell) \cap \cC_v(F_{w_\ell}) \big)$ of $E_v$,
           and is move-prepared relative to $B(a_1,r_1), \ldots, B(a_D,r_D)$.
     \index{move-prepared} 
           In particular, if $\theta_h \in B(a_{\ell},r_{\ell}) \cap \cC_v(F_{w_\ell})$, 
           then $F_{u_h} = F_{w_\ell}$, $\rho_h \in |F_{w_\ell}^{\times}|_v$, 
           and $B(\theta_h,\rho_h) \subseteq B(a_{\ell},r_{\ell})$.
           For each $\ell = 1, \ldots, D$, there is a point 
            $\wbar_\ell \in \big(B(a_\ell,r_\ell) \cap \cC_v(F_{w_\ell})) \big) \backslash H_v$.
           
  $(3)$ For each $h = 1, \ldots, N$, $\phi_v$ induces 
           an $F_{u_h}$-rational scaled isometry from $B(\theta_h,\rho_h)$ onto $D(0,1)$,
           which takes $B(\theta_h,\rho_h) \cap \cC_v(F_{u_h})$ onto $\cO_{u_h}$. 

\vskip .05 in
Let $0 < h_v < r_v < 1$ be numbers satisfying 
\begin{equation} \label{hvrvIneq1}
h_v^N \ < \ r_v^N \ < \ q_v^{-1/(q_v-1)} \ < \ 1 \ .
\end{equation} 
Put 
\begin{equation} \label{CompactMvDef}
M_v \ = \ \max( \max_{\substack{ 1 \le i \le m \\ N_i < j \le 2N_i}} \|\varphi_{ij}\|_{U_v },
    \max_{1 \le \lambda \le \Lambda} \|\varphi_{\lambda}\|_{U_v } ) \ ,
\end{equation} 
and let $k_v > 0$ be the least integer such that for all $k \ge k_v$,
\begin{equation} \label{DFCCZ1}
h_v^{kN} \cdot M_v \ < \ q_v^{-\frac{k+1}{q_v - 1} - \log_v(k+1)} \ . 
\end{equation}
Let $\kbar \ge k_v$ be a fixed integer, and put $B_v = h_v^{\kbar N}$.
Then there is an integer $n_v \ge 1$ such that for each sufficiently large integer $n$
divisible by $n_v$, the local patching process at $K_v$ can be carried out as follows: 

\smallskip
Put  $G_v^{(0)}(z) = S_{n,v}(\phi_v(z))$.\index{patching functions, initial $G_v^{(0)}(z)$}   

For each $k$, $1 \le k < n$, 
let $\{\Delta_{v,ij}^{(k)} \in \CC_v\}_{(i,j) \in \Band_N(k)}$\index{band!$\Band_N(k)$} 
be an arbitrary $K_v$-symmetric set of numbers given recursively in $\prec_N$ order,\index{order!$\prec_N$} 
\index{$K_v$-symmetric!set of numbers}
subject to the condition that for each $(i,j)$   
\begin{equation} \label{DFCCV1A}
|\Delta_{v,ij}^{(k)}|_v \ \le \ \left\{
      \begin{array}{ll} B_v & \text{if \ $k \le \kbar \ ,$} \\
                         h_v^{kN} & \text{if \ $k > \kbar \ .$}
      \end{array} \right.
\end{equation}
For $k = n$, let 
$\{\Delta_{v,\lambda}^{(n)} \in \CC_v\}_{1 \le \lambda \le \Lambda}$ 
be an arbitrary $K_v$-symmetric set of numbers satisfying
\index{$K_v$-symmetric!set of numbers}  
\begin{equation} \label{DFCCF2A} 
|\Delta_{v,\lambda}^{(n)}|_v \ \le \ h_v^{nN} \ .  
\end{equation} 

\vskip .05 in
\noindent{Then} one can inductively construct $(\fX,\vs)$-functions 
$G_v^{(1)}(z), \ldots, G_v^{(n)}(z)$ in $K_v(\cC)$,\index{patching functions, $G_v^{(k)}(z)$ for $1 \le k \le n$!constructed by patching} 
of common degree $Nn$, such that 

\smallskip
\noindent{$(A)$} For each $k = 1, \ldots, n-1$,\index{compensating functions 
$\vartheta_{v,ij}^{(k)}(z)$}\index{patching functions, $G_v^{(k)}(z)$ for $1 \le k \le n$!constructed by patching}  
\begin{equation*} 
G_v^{(k)}(z) \ = \ G_v^{(k-1)}(z) +  \sum_{(i,j) \in \Band_N(k)} 
                           \Delta_{v,ij}^{(k)} \cdot \vartheta_{v,ij}^{(k)}(z) \ + \ \Theta_v^{(k)}(z) \ ,
\end{equation*}\index{band!$\Band_N(k)$}                           
where $\vartheta_{v,ij}^{(k)}(z) = \varphi_{i,(k+1)N_i - j} \cdot F_{v,k}(z)$ 
with\index{compensating functions $\vartheta_{v,ij}^{(k)}(z)$!construction of} 
an $(\fX,\vs)$ function $F_{v,k}(z) \in K_v(\cC)$ 
of degree $(n-k-1)$ independent of $(i,j)$ whose roots belong to $E_v$, 
and $\Theta_v^{(k)}(z) \in K_v(\cC)$ has  a pole of order at most $(n-k)N_i$ at each $x_i$ 
and no other poles.  For each $i$, the leading coefficient $d_{v,i}$ of $F_k(z)$ at $x_i$
\index{coefficients $A_{v,ij}$!leading} 
belongs to $K_v(x_i)$ and has absolute value $|d_{v,i}|_v = |\tc_{v,i}|_v^{n-k-1}$.  For $k = n$ 
\begin{equation*}                           
G_v^{(n)}(z) \ = \ G_v^{(n-1)}(z) + \sum_{\lambda=1}^{\Lambda}
                   \Delta_{v,\lambda}^{(n)} \cdot \varphi_{\lambda}(z) \ .  
\end{equation*}
In particular

$(1)$ Each $\vartheta_{v,ij}^{(k)}(z)$ belongs to $K_v(x_i)(\cC)$,
\index{compensating functions $\vartheta_{v,ij}^{(k)}(z)$!poles and leading coefficients of}  
has a pole of order $nN_i-j > (n-k-1)N_i$ at $x_i$ with leading coefficient $\tc_{v,i}^{n-k-1}$,  
has poles of order at most $(n-k-1)N_{i^{\prime}}$ at each $x_{i^{\prime}} \ne x_i$, 
and has no other poles $;$  

$(2)$ The $\vartheta_{v,ij}^{(k)}(z)$ are $K_v$-symmetric $;
$\index{compensating functions $\vartheta_{v,ij}^{(k)}(z)$!are $K_v$-symmetric} 
\index{$K_v$-symmetric!set of functions}

$(3)$  In passing from $G_v^{(0)}(z)$ to 
$G_v^{(1)}(z)$,\index{patching functions, $G_v^{(k)}(z)$ for $1 \le k \le n$!leading coefficients of}  
each of the leading coefficients\index{patching functions, initial $G_v^{(0)}(z)$!leading coefficients of} 
\index{coefficients $A_{v,ij}$!leading} 
$\tc_{v,i}^n$ of $G_v^{(0)}(z)$ is replaced with 
$\tc_{v,i}^n + \Delta_{v,i0}^{(1)} \cdot \tc_{v,i}^{n-2}$.

\vskip .05 in
\noindent{$(B)$} For each $k = 0, 1, \ldots, n$, 
the zeros of $G_v^{(k)}(z)$ belong to $E_v$,
\index{patching functions, $G_v^{(k)}(z)$ for $1 \le k \le n$!for nonarchimedean $K_v$-simple sets!roots are distinct}  
and for $k = 0$ and $k = n$ they are distinct.  Furthermore 
\begin{equation*}
\{ z \in \cC_v(\CC_v) : G_v^{(n)}(z) \in \cO_v 
    \text{ and } |G_v^{(n)}(z)|_v \le r_v^{Nn} \} \ \subseteq \ E_v \ .
\end{equation*} 
\end{theorem}

When $\Char(K_v) = p > 0$, we have the following patching theorem.
The $K_v$-rationality assumptions (\ref{DFCCV1BB}), (\ref{DFCCF2BB})
in the theorem are addressed by the global patching process.  

\begin{theorem}  \label{DCPCPatch1p}\index{patching theorem!for nonarchimedean $K_v$-simple sets!when $\Char(K_v) = p > 0$} 
Suppose $K_v$ is nonarchimedean, and $\Char(K_v) = p > 0$.   
Let $E_v$ be $K_v$-simple, with a $K_v$-simple decomposition 
\index{$K_v$-simple!set}
\index{$K_v$-simple!decomposition}
$E_v = \bigcup_{\ell = 1}^D \big( B(a_\ell,r_\ell) \cap \cC_v(F_{w_\ell}) \big)$  
such that $U_v := \bigcup_{\ell=1}^D  B(a_{\ell},r_{\ell})$ is disjoint from $\fX$. 
Let $\vs \in \cP^m(\QQ)$ be a $K_v$-symmetric probability vector with positive coefficients,
\index{$K_v$-symmetric!probability vector} 
and let $\phi_v(z) \in K_v(\cC)$ be an $(\fX,\vs)$-function of degree $N$ 
with distinct zeros $\theta_1, \ldots, \theta_N \in E_v$, such that

\vskip .05 in
  $(1)$ $\phi_v^{-1}(D(0,1)) = \bigcup_{h=1}^N B(\theta_h,\rho_h)$, 
             where the balls $B(\theta_1,\rho_1), \ldots, B(\theta_N,\rho_N)$ 
             are pairwise disjoint, isometrically parametrizable, 
             and contained in $\bigcup_{\ell = 1}^D B(a_{\ell},r_{\ell})$ $;$ 
         
  $(2)$ The set $H_v := \phi_v^{-1}(D(0,1)) \cap E_v$ is $K_v$-simple, 
  \index{$K_v$-simple!set}
  \index{$K_v$-simple!decomposition}
            with a $K_v$-simple decomposition 
            \begin{equation*} 
            H_v \ = \ \bigcup_{h=1}^N \big( B(\theta_h,\rho_h) \cap \cC_v(F_{u_h}) \big)           
            \end{equation*}
           compatible with the $K_v$-simple decomposition 
    \index{$K_v$-simple!decomposition!compatible with another decomposition}  
           $\bigcup_{\ell = 1}^D \big( B(a_\ell,r_\ell) \cap \cC_v(F_{w_\ell}) \big)$ of $E_v$,
           and is move-prepared relative to $B(a_1,r_1), \ldots, B(a_D,r_D)$.
    \index{move-prepared} 
           In particular, if $\theta_h \in B(a_{\ell},r_{\ell}) \cap \cC_v(F_{w_\ell})$, 
           then $F_{u_h} = F_{w_\ell}$, $\rho_h \in |F_{w_\ell}^{\times}|_v$, 
           and $B(\theta_h,\rho_h) \subseteq B(a_{\ell},r_{\ell})$.
           For each $\ell = 1, \ldots, D$, there is a point 
            $\wbar_\ell \in \big(B(a_\ell,r_\ell) \cap \cC_v(F_{w_\ell})) \big) \backslash H_v$.
           
  $(3)$ For each $h = 1, \ldots, N$, $\phi_v$ induces 
           an $F_{u_h}$-rational scaled isometry from $B(\theta_h,\rho_h)$ onto $D(0,1)$,
           which takes $B(\theta_h,\rho_h) \cap \cC_v(F_{u_h})$ onto $\cO_{u_h}$. 
 
\vskip .05 in
Let $0 < h_v < r_v < 1$ be numbers satisfying 
\begin{equation} \label{hvrvIneq2}
h_v^N \ < \ r_v^N \ < \ q_v^{-1/(q_v-1)} \ < \ 1 \ ,
\end{equation}
and let 
\begin{equation} \label{CompactMvDefp}
M_v \ = \ \max( \max_{\substack{ 1 \le i \le m \\ N_i < j \le 2N_i}} \|\varphi_{ij}\|_{U_v },
    \max_{1 \le \lambda \le \Lambda} \|\varphi_{\lambda}\|_{U_v } ) \ .
\end{equation} 
Let $k_v > 0$ be the least integer such that for all $k \ge k_v$,
\begin{equation} \label{DFCCZ1p}
h_v^{kN} \cdot M_v \ < \ q_v^{-\frac{k+1}{q_v - 1} - \log_v(k+1)} \ . 
\end{equation}
Let $\kbar \ge k_v$ be a fixed integer. 
Then there are an integer $n_v \ge 1$ and a number $0 < B_v < 1$ 
such that for each sufficiently large integer $n$
divisible by $n_v$, the local patching process at $K_v$ can be carried out as follows: 
     
\smallskip
Put  $G_v^{(0)}(z) = S_{n,v}(\phi_v(z))$.\index{patching functions, initial $G_v^{(0)}(z)$!construction of}   

Then the leading coefficient\index{patching functions, initial $G_v^{(0)}(z)$!leading coefficients of} 
of $G_v^{(0)}(z)$ at $x_i$ is $\tc_{v,i}^n$,\index{coefficients $A_{v,ij}$!leading} 
the zeros of $G_v^{(0)}(z)$
\index{patching functions, initial $G_v^{(0)}(z)$!for nonarchimedean $K_v$-simple sets!roots are distinct}  
are distinct and belong to $E_v$, 
and when $G_v^{(0)}(z)$\index{patching functions, initial $G_v^{(0)}(z)$!roots are confined to $E_v$}  
is expanded in terms of the $L$-rational basis as 
\index{basis!$L$-rational}\index{patching functions, initial $G_v^{(0)}(z)$!expansion of}
\begin{equation*} 
 G_v^{(0)}(z) \ = \ \sum_{i=1}^m \sum_{j=0}^{(n-1)N_i-1} A_{v,ij} \varphi_{i,nN_i-j}(z) 
   + \sum_{\lambda=1}^{\Lambda} A_{v,\lambda} \varphi_{\lambda} \ ,
\end{equation*} 
we have $A_{v,ij} = 0$ for all $(i,j)$ with $1 \le j < \kbar N_i$.
\index{patching functions, initial $G_v^{(0)}(z)$!when $\Char(K) = p > 0$!high-order coefficients are $0$}  

For each $k$, $1 \le k \le n-1$, 
let $\{\Delta_{v,ij}^{(k)} \in L_{w_v}\}_{(i,j) \in \Band_N(k)}$\index{band!$\Band_N(k)$}\index{distinguished place $w_v$} 
be a $K_v$-symmetric set of numbers satisfying 
\index{$K_v$-symmetric!set of numbers}
\begin{equation} \label{DFCCV1Ap}
\left\{ \begin{array}{ll}
\text{$|\Delta_{v,i0}^{(1)}|_v \le B_v $ 
            and $\Delta_{v,ij}^{(1)} = 0$ for $j = 1, \ldots, N_i-1$,}  
               & \text{if \ $k = 1$\ ,} \\
\text{$\Delta_{v,ij}^{(k)} = 0$ for $j = (k-1)N_i, \ldots, kN_i - 1$,} 
               & \text{if \ $k = 2, \ldots, \kbar$\ , } \\
|\Delta_{v,ij}^{(k)}|_v  \le  h_v^{kN},  
               & \text{if \ $k = \kbar+1, \ldots, n-1$\ ,}
      \end{array} \right.
\end{equation}
such that $\Delta_{v,i0}^{(1)} \in K_v(x_i)^{\sep}$ for each $i$ 
and such that for each $k = \kbar + 1, \ldots, n-1$ 
\begin{equation} \label{DFCCV1BB}
\Delta_{v,k}(z) \ := \ 
\sum_{i=1}^m \sum_{j=(k-1)N_i}^{kN_i-1} \Delta_{v,ij}^{(k)} \cdot   \varphi_{i,(k+1)N_i - j} 
 \ \in \ K_v(\cC) \ .
\end{equation} 
For $k = n$, let 
$\{\Delta_{v,\lambda}^{(n)} \in L_{w_v}\}_{1 \le \lambda \le \Lambda}$\index{distinguished place $w_v$} 
be a $K_v$-symmetric set of numbers such that 
\index{$K_v$-symmetric!set of numbers}
\begin{equation} \label{DFCCF2AA}  
|\Delta_{v,\lambda}^{(n)}|_v \ \le \ h_v^{nN}   
\end{equation} 
for each $\lambda$, and 
\begin{equation} \label{DFCCF2BB}  
\Delta_{v,n}(z) \ := \ 
\sum_{\lambda = 1}^{\Lambda} \Delta_{v,\lambda}^{(n)} \cdot \varphi_{\lambda} \ \in K_v(\cC) \ . 
\end{equation} 

Then one can inductively construct\index{patching functions, $G_v^{(k)}(z)$ for $1 \le k \le n$!constructed by patching}  
 $(\fX,\vs)$-functions $G_v^{(1)}(z), \ldots, G_v^{(n)}(z)$ in $K_v(\cC)$, 
of common degree $Nn$, such that$:$  
 
\smallskip
\noindent{$(A)$} For each $k = 1, \ldots, n$, $G_v^{(k)}(z)$ is obtained from $G_v^{(k-1)}(z)$
as follows$:$\index{patching functions, $G_v^{(k)}(z)$}   

$(1)$  When $k = 1$, there is a $K_v$-symmetric set of functions
\index{$K_v$-symmetric!set of functions} 
$\ttheta_{v,10}^{(1)}(z), \ldots, \ttheta_{v,m0}^{(1)}(z) \in L_{w_v}^{\sep}(\cC)$\index{distinguished place $w_v$}
such that\index{patching functions, $G_v^{(k)}(z)$ for $1 \le k \le n$!constructed by patching}  
\begin{equation*}
G_v^{(1)}(z) \ = \ 
    G_v^{(0)}(z) + \sum_{i=1}^m \Delta_{v,i0}^{(1)} \cdot \ttheta_{v,i0}^{(1)}(z) \ , 
\end{equation*}
where for each $i = 1, \ldots, m$,
$\ttheta_{v,i0}^{(1)}(z) \in K_v(x_i)^{\sep}(\cC)$ has the form 
\begin{equation*}
\ttheta_{v,i0}^{(1)}(z) \ = \ \tc_{v,i}^{n} \varphi_{i,nN_i}(z) + \tTheta_{v,i0}^{(1)}(z) 
\end{equation*}   
for an $(\fX,\vs)$-function $\tTheta_{v,i0}^{(1)}(z)$
with a pole of order at most $(n-\kbar)N_{i^{\prime}}$ at each $x_{i^\prime}$.  
Thus, in passing from $G_v^{(0)}(z)$ 
to $G_v^{(1)}(z)$,\index{patching functions, $G_v^{(k)}(z)$ for $1 \le k \le n$!leading coefficients of}  
each of the leading coefficients\index{patching functions, initial $G_v^{(0)}(z)$!leading coefficients of}  
$A_{v,i0} = \tc_{v,i}^n$\index{coefficients $A_{v,ij}$!leading} is replaced with $\tc_{v,i}^n + \Delta_{v,i0}^{(1)} \cdot \tc_{v,i}^n$,
and the coefficients $A_{v,ij}$ for $1 \le j < \kbar N_i$ remain $0$.   
\index{coefficients $A_{v,ij}$}  

$(2)$ For $k = 2, \ldots, \kbar$, we have 
$G_v^{(k)}(z) = G_v^{(k-1)}(z)$.\index{patching functions, $G_v^{(k)}(z)$ for $1 \le k \le n$!constructed by patching} 

$(3)$ For $k = \kbar + 1, \ldots, n-1 $, we have\index{patching functions, $G_v^{(k)}(z)$ for $1 \le k \le n$!constructed by patching} 
\begin{equation} \label{GvkOmegaModificationp} 
G_v^{(k)}(z) \ = \ G_v^{(k-1)}(z)
     + \Delta_{v,k}(z) \cdot F_{v,k}(z) + \Theta_v^{(k)}(z) \ ,    
\end{equation}
where 

\quad $(a)$ $\Delta_{v,k}(z) = \sum_{(i,j) \in \Band_N(k)} \Delta_{v,ij}^{(k)} \varphi_{i,(k+1)N_i-j}(z)$ belongs to $K_v(\cC)$ by $(\ref{DFCCV1BB});$\index{band!$\Band_N(k)$}

\quad $(b)$ $F_{v,k}(z) \in K_v(\cC)$ is an $(\fX,\vs)$-function 
determined by the local patching process 

\qquad using $G_v^{(k-1)}(z)$,\index{patching functions, $G_v^{(k)}(z)$} 
whose roots belong to $E_v$. For each $x_i$, it has a pole of order 

\qquad $(n-k-1)N_i$ at $x_i$, and the leading coefficient 
\index{coefficients $A_{v,ij}$!leading} 
$d_{v,i} = \lim_{z \rightarrow x_i} F_{v,k}(z) \cdot g_{x_i}(z)^{(n-k-1)N_i}$

\qquad has absolute value $|d_{v,i}|_v = |\tc_{v,i}|_v^{n-k-1}$. 

\quad $(c)$ $\Theta_v^{(k)}(z) \in K_v(\cC)$ is an $(\fX,\vs)$-function 
determined by the local patching process 

\qquad after the coefficients in $\Band_N(k)$ have been modified;  it has  a pole of order 
\index{coefficients $A_{v,ij}$}\index{band!$\Band_N(k)$}

\qquad at most $(n-k)N_i$ at each $x_i$ and no other poles, and may be the zero function.  

$(4)$ For $k = n$\index{patching functions, $G_v^{(k)}(z)$ for $1 \le k \le n$!constructed by patching}                              
\begin{equation*}
G_v^{(n)}(z) \ = \ G_v^{(n-1)}(z) + \Delta_{v,n}(z) \ . 
\end{equation*}
\smallskip
\noindent{$(B)$} For each $k = 1, \ldots, n$, 
the zeros of $G_v^{(k)}(z)$ belong to $E_v$,
\index{patching functions, $G_v^{(k)}(z)$ for $1 \le k \le n$!roots are confined to $E_v$} 
\index{patching functions, $G_v^{(k)}(z)$ for $1 \le k \le n$!for nonarchimedean $K_v$-simple sets!roots are distinct} 
and for $k = n$ they are distinct.  Furthermore   
\begin{equation*}
\{ z \in \cC_v(\CC_v) : G_v^{(n)}(z) \in \cO_v 
 \text{ and } |G_v^{(n)}(z)|_v \le r_v^{Nn} \} \ \subseteq \ E_v \ . 
\end{equation*} 
\end{theorem}

\smallskip
\noindent{\bf Remark.} In both theorems, 
we will have $\Theta_v^{(k)}(z) = 0$ for all but one value $k = k_1$. 
For $k = k_1$, 
after computing\index{compensating functions $\vartheta_{v,ij}^{(k)}(z)$}
\index{patching functions, $G_v^{(k)}(z)$ for $1 \le k \le n$!constructed by patching} 
\begin{equation*}
G_v^{(k_1)}(z) \ = \ G_v^{(k_1-1)}(z)
              +  \sum_{(i,j) \in \Band_N(k_1)}
                           \Delta_{v,ij}^{(k_1)} \vartheta_{v,ij}^{(k_1)}(z) 
\end{equation*}
\index{band!$\Band_N(k)$}
we `move the roots of $G_v^{(k_1)}(z)$ apart',
\index{patching functions, $G_v^{(k)}(z)$ for $1 \le k \le n$!for nonarchimedean $K_v$-simple sets!roots are separated} 
constructing an $(\fX,\vs)$-function 
\index{move roots apart}
$\hG_v^{(k_1)}(z)$ whose patched coefficients 
\index{move roots apart} 
are the same as those of $G_v^{(k_1)}(z)$,
\index{patching functions, $G_v^{(k)}(z)$ for $1 \le k \le n$!for nonarchimedean $K_v$-simple sets!roots are separated} 
but whose roots are well-separated, 
and replace $G_v^{(k_1)}(z)$ by $\hG_v^{(k_1)}(z)$.\index{patching functions, $G_v^{(k)}(z)$}  
The function
$\Theta_v^{(k_1)}(z)= \hG_v^{(k_1)}(z)-G_v^{(k_1)}(z)$ 
is chosen to accomplish this change: see Phase 3 in the proof,\index{separate the roots} 
 given in \S\ref{Chap11}.\ref{NonArchPatchingProof} below. 

\smallskip 
Theorems \ref{DCPCPatch} and \ref{DCPCPatch1p} will be proved together.  There are some 
differences in the way the leading and high-order coefficients are treated, but the underlying
\index{patching construction!for nonarchimedean $K_v$-simple sets!differences when $\Char(K) = 0$ and $\Char(K) = p$}   
\index{coefficients $A_{v,ij}$!leading}\index{coefficients $A_{v,ij}$!high-order} 
patching constructions for the middle and low-order coefficients are the same. 
\index{coefficients $A_{v,ij}$!middle}
\index{coefficients $A_{v,ij}$!low-order}

\smallskip
Given a positive integer $n$, 
we begin the patching process by composing $\phi_v(z)$ with the
Stirling polynomial $S_{n,v}(x)$.  This yields the $K_v$-rational function
\index{Stirling polynomial!for $\cO_v$}\index{patching functions, initial $G_v^{(0)}(z)$!construction of} 
\begin{equation*} 
G_v^{(0)}(z) \ := \ S_{n,v}(\phi_v(z)) \ = \ \prod_{j=0}^{n-1} (\phi_v(z)-\psi_v(j)) \ .
\end{equation*}
An important observation is that 
$S_{n,v}(\phi_v(z))$ is highly factorized, and by taking products 
$\prod_{j \in \cS} (\phi_v(z)-\psi_v(j))$
corresponding to subsets $\cS \subset \{0, \ldots, n-1\}$ 
we can easily obtain $K_v$-rational functions 
dividing $G_v^{(0)}(z)$.
\index{patching functions, initial $G_v^{(0)}(z)$!for nonarchimedean $K_v$-simple sets!are highly factorized} 
On the other hand, by restricting $S_{n,v}(\phi_v(z))$ to one of  
the isometrically parametrizable balls $B(\theta_h,\rho_h)$ and composing it
with a parametrization, we obtain a power series which behaves
much like $S_{n,v}(x)$.  The key to the construction is the interaction 
between these global and local ways of 
viewing $G_v^{(0)}(z)$.\index{patching functions, initial $G_v^{(0)}(z)$}   

Before we can give the proof, we must develop some machinery.


\section{ The Patching Lemmas } \label{PatchingLemmasSection}

In this section, we consider aspects of the construction involving power series. 
\index{patching construction!for nonarchimedean $K_v$-simple sets!basic patching lemmas}   

For each zero $\theta_h$ of $\phi_v(z)$, let 
$\sigma_h : D(0,\rho_h) \rightarrow B(\theta_h,\rho_h)$
be an $F_{u_h}$-rational isometric parametrization with $\sigma_h(0) = \theta_h$.
Let $d_h \in F_{u_h}^{\times}$ be such that $|d_h|_v = \rho_h$, and define
$\hsigma_h : D(0,1) \rightarrow B(\theta_h,\rho_h)$ by
\begin{equation*}
\hsigma_h(Z) \ = \ \sigma_h(d_hZ) \ .
\end{equation*}
Put $\hPhi_h(Z) = \phi_v(\hsigma_h(Z))$.
Thus $\hPhi_h(Z)$ is a power series converging on $D(0,1)$  
which induces an $F_{u_h}$-rational distance-preserving isomorphism
from $D(0,1)$ to itself and takes $\cO_{u_h}$ to $\cO_{u_h}$.  
It satisfies $\hPhi_h(0) = 0$, and $|\hPhi_h^{\prime}(0)|_v = 1$.  
After replacing $d_h$ by $\mu_h d_h$ 
for an appropriate $\mu_h \in \cO_{u_h}^{\times}$, if necessary, 
we can assume that $\hPhi_h^{\prime}(0) = 1$.  Hence we can expand 
\begin{equation*}
\hPhi_h(Z) \ = \ Z + \sum_{i=2}^{\infty} C_{hi} Z^i \in \cO_{u_h}[[Z]] \ ,
\end{equation*}
where $\ord_v(C_{hi}) > 0$ for each $i \ge 2$,  
and $\ord_v(C_{hi}) \rightarrow \infty$ as $i \rightarrow \infty$.  
Let $\tPhi_h(Z)$ be the inverse power series to $\hPhi_h(Z)$, 
so $\tPhi_h(\hPhi_h(Z)) = \hPhi_h(\tPhi_h(Z)) = Z$.  

The restriction of $S_{n,v}(\phi_v(z))$ to $B(\theta_h,\rho_h)$ 
corresponds to the function $S_{n,v}(\hPhi_h(Z))$ on $D(0,1)$.  
There is a $1-1$ correspondence between the zeros 
$\theta_{hj}$ of $S_{n,v}(\phi_v(z))$ in $B(\theta_h,\rho_h)$ 
and the zeros $\alpha_{hj}$ of $S_{n,v}(\hPhi_h(Z))$ in $D(0,1)$,
given by $\theta_{hj} = \hsigma_h(\alpha_{hj})$;  the $\theta_{hj}$
belong to $\cC_v(F_{u_h})$, and the $\alpha_{hj}$ belong to $\cO_{u_h}$. 
Since the zeros of $S_{n,v}(z)$ are the $\psi_v(j)$ for $j = 0, \ldots, n-1$ we have
$\hPhi_h(\alpha_{hj}) = \psi_v(j)$ (or equivalently, $\tPhi_h(\psi_v(j)) = \alpha_{hj}$)
for each $j$, $0 \le j \le n-1$.
Thus the power series
\begin{equation*}
G_h^{(0)}(Z) \ := \ S_{n,v}(\hPhi_h(Z)) \ = \ \prod_{j=0}^{n-1} (\hPhi_h(Z) - \psi_v(j))
\end{equation*}
can be factored as 
\begin{equation*}
G_h^{(0)}(Z) \ = \ \prod_{j=0}^{n-1} (Z - \alpha_{hj}) \cdot \cH_h(Z)
\end{equation*}
where $\cH_h(Z)$ is an invertible power series in $\cO_{u_h}[Z]$.
Write $q = q_v$.  By (\ref{BGF1}), 
since $\hPhi_h(Z)$ is distance-preserving, 
for all $0 \le j, k < n$ with $j \ne k$ 
\begin{equation} 
\ord_v(\alpha_{h j}-\alpha_{h k}) \ = \ \ord_v(\psi_v(j) - \psi_v(k))
      \ = \ \val_q(|j - k|) \ . \label{FReg1} 
\end{equation}

We will now consider how to modify the functions 
$G_h^{(0)}(Z)$, while keeping their zeros in $\cO_{u_h}$.
To do this, we establish a series of lemmas concerning power series,  
analogous to those proved for polynomials in (\cite{RR3}).
For the rest of this section, $F_u/K_v$ will denote a finite extension in $\CC_v$,
with ring of integers $\cO_u$.

\begin{definition} \label{DDefCPC1}  Let $F_u/K_v$ be a finite extension.
Suppose $\cS = \{j_1, j_1+1, \ldots, j_1+\ell-1\}$ 
is a sequence of $\ell$ consecutive non-negative integers, and that  
$\hPhi : D(0,1) \rightarrow D(0,1)$ is a distance-preserving automorphism 
defined by an $F_u$-rational power series, so $\hPhi(\cO_u) = \cO_u$.
\index{regular sequence!$\psi_v$-regular sequence|ii}
A {\em $\psi_v$-regular sequence
of length $\ell$ in $\cO_u$ attached to $\cS$ and $\hPhi(z)$} is a sequence
$\{\alpha_j\}_{j \in \cS} \subset \cO_u$ such that
\begin{equation*}
\ord_v(\hPhi(\alpha_{j}) - \psi_v(j)) \ \ge \ \log_v(\ell)
\end{equation*}
(or equivalently, $\ord_v\big(\alpha_j-\tPhi(\psi_v(j))\big) \ge \log_v(\ell)$), for each $j \in \cS$.
\end{definition}

In particular, for each $h$, $\{\alpha_{hj}\}_{0 \le j < n}$ is a $\psi_v$-regular
sequence\index{regular sequence!$\psi_v$-regular sequence} of length $n$ in $\cO_{u_h}$ relative to $\hPhi_h(Z)$. 
In the applications, we will have $\hPhi(Z) = \hPhi_h(Z)$,
and $\cS$ will be a subsequence of $\{0,1, \ldots, n-1\}$.
However, often the precise power series defining the $\psi_v$-regular sequence in $\cO_u$ 
\index{regular sequence!$\psi_v$-regular sequence}
is not important (generally all that is used is the fact the power series is an isometry), 
and we will frequently speak of a $\psi_v$-regular sequence of length $\ell$,
\index{regular sequence!$\psi_v$-regular sequence} 
or just a $\psi_v$-regular sequence, if $\cS$, $\hPhi(Z)$, and $\cO_u$  
are understood from context.  Note that each subsequence of a $\psi_v$-regular sequence
\index{regular sequence!$\psi_v$-regular sequence}
consisting of consecutive elements, is itself a $\psi_v$-regular sequence.
\index{regular sequence!$\psi_v$-regular sequence}

For the rest of this section, 
we will assume that $\hPhi(Z) \in \cO_u[[Z]]$ as in Definition \ref{DDefCPC1} 
has been fixed; let $\tPhi(Z)$ denote its inverse. 
If $\{\alpha_j\}_{j \in \cS}$ is a
$\psi_v$-regular sequence of length $\ell$ in $\cO_u$ for $\hPhi(Z)$, then for each $k \ne j \in \cS$,
\index{regular sequence!$\psi_v$-regular sequence}
\begin{equation} \label{DFXXY1}
\ord_v(\alpha_{k}-\alpha_j) \ = \ \val_q(|k-j|) 
\end{equation}
because
\begin{equation*}
\ord_v(\alpha_{k}-\tPhi(\psi_v(k)))  \ \ge \ \log_v(\ell) \ , \quad
\ord_v(\alpha_{j}-\tPhi(\psi_v(j))) \ \ge \ \log_v(\ell) \ ,
\end{equation*}
while by formula (\ref{BGF1}), since $\tPhi(Z)$ is distance-preserving, 
\begin{eqnarray*}
\ord_v(\tPhi(\psi_v(k))-\tPhi(\psi_v(j))) 
& = & \ord_v(\psi_v(k)-\psi_v(j)) \\
& = & \val_q(|k-j|) \ < \ \log_v(\ell) \ .
\end{eqnarray*}

\begin{lemma} \label{DLemCPC2}
Let $\{ \alpha_j \}_{j_1 \le j < j_1+\ell}$ be a $\psi_v$-regular sequence of
\index{regular sequence!$\psi_v$-regular sequence}
length $\ell$ in $\cO_u$.  Given $z \in D(0,1)$, let $J$ be an index for which
$\ord_v(z-\alpha_J)$ is maximal.  
Then for each $j \ne J$ with $j_1 \le j < j_1+\ell$, we have
\begin{equation*}
\ord_v(z-\alpha_j) \ \le \ \ord_v(\alpha_J - \alpha_j) \ = \ \val_q(|J-j|) \ .
\end{equation*}
\end{lemma}

\begin{proof}  Fix $j \ne J$.  
By hypothesis $\ord_v(z-\alpha_j) \le \ord_v(z-\alpha_J)$.  By the ultrametric inequality, 
\begin{equation*}
\ord_v(\alpha_J-\alpha_j)
    \ \ge \ \min(\ord_v(z-\alpha_J),\ord_v(z-\alpha_j))  
    \  =  \ \ord_v(z-\alpha_j) \ . 
\end{equation*}
By (\ref{DFXXY1}), $\ord_v(\alpha_J-\alpha_j) = \val_q(|J-j|)$, 
so we obtain the result.
\end{proof}

\vskip .1 in
Recall from (\ref{FactorialFormula}) that 
\begin{equation} \label{DFWSD1}
\sum_{k=1}^\ell \val_q(k) \ = \
     \frac{\ell}{q-1} - \frac{1}{q-1} \sum_{j \ge 0} d_j(\ell) \ .
\end{equation}
The following lemma generalizes the parts of Proposition \ref{BGProp1} we will need:  

\begin{lemma} \label{DLemCPC3}
Let $\{ \alpha_j \}_{j \in \cS}$ be a $\psi_v$-regular sequence of
\index{regular sequence!$\psi_v$-regular sequence}
length $\ell$ in $\cO_u$.  Put $P_{\cS}(Z) = \prod_{j \in \cS} (Z-\alpha_j)$.  Then

$(A)$  For each $J \in \cS$
\begin{equation*}
        \ord_v\big(\prod_{\substack{ j \in \cS \\ j \ne J }}
                  (\alpha_J-\alpha_j)\big) \ < \ \frac{\ell}{q-1} \ .
\end{equation*}

$(B)$  For each $z \in D(0,1)$, if $J \in \cS$ is such that
$\ord_v(z-\alpha_J)$ is maximal, then
\begin{equation*}
\ord_v(P_{\cS}(z)) \ < \ \frac{\ell}{q-1} + \ord_v(z-\alpha_J) \ .
\end{equation*}
\end{lemma}

\begin{proof}  Suppose $\cS = \{j_1, \ldots, j_1+\ell-1\}$.  
To prove $(A)$,  fix $J \in \cS$ and recall that if $j \in \cS$ and $j \ne J$ then
$\ord_v(\alpha_J-\alpha_j) = \val_q(|J-j|)$.  Hence
\begin{eqnarray} 
\lefteqn{\ord_v(\prod_{\substack{ j \in \cS \\ j \ne J }} (\alpha_J-\alpha_j)) 
 \ = \ \sum_{j=j_1}^{J-1} \val_q(|J-j|)
               + \sum_{j=J+1}^{j_1+\ell-1} \val_q(|J-j|)} \qquad \qquad & & \notag \\
 & = & \sum_{k=1}^{J-j_1} \val_q(k) + \sum_{k=1}^{j_1+\ell-J-1} \val_q(k) 
                    \label{DFCXV1} \\
 & = & \frac{(J-j_1) - \sum_{j \ge 0} d_j(J-j_1)}{q-1}  \notag \\
 &   & \qquad + \ \frac{(j_1+\ell-J-1) - \sum_{j \ge 0} d_j(j_1+\ell-J-1)}{q-1} \notag \\
 & = & \frac{\ell}{q-1}
      - \frac{ 1 + \sum_{j \ge 0} d_j(J-j_1) + \sum_{j \ge 0} d_j(j_1+\ell-J-1)}{q-1}
      \ < \ \frac{\ell}{q-1} \ . \notag
\end{eqnarray}

For $(B)$, let $J \in \cS$ be such that $\ord_v(z-\alpha_J)$ is maximal,
or equivalently $|z-\alpha_J|_v$ is minimal.
By Lemma \ref{DLemCPC2}
\begin{equation*}
\ord_v(P_{\cS}(z)) \ \le \ \ord_v(z-\alpha_J)
                      + \ord_v(\prod_{j \ne J} (\alpha_J - \alpha_j)) \ .
\end{equation*}
Applying part $A)$, we get
\begin{equation*} 
\ord_v(P_{\cS}(z)) \ < \ \frac{\ell}{q-1} + \ord_v(z-\alpha_J) \ ,
\end{equation*}
as required.
\end{proof}

\vskip .1 in
We now come to the basic lemma governing the patching process.  
If $\cQ(Z) \in \CC_v[[Z]]$ converges in $D(0,1)$, write
\begin{equation*}
\|\cQ\|_{D(0,1)} \ = \ \sup_{Z \in D(0,1)} |\cQ(Z)|_v \ .
\end{equation*}
for its $\sup$ norm relative to the absolute value $|x|_v$.  The lemma will 
be applied to functions of the form 
$\cQ(Z) = \cQ_{\cS,h}(Z) = \prod_{j \in \cS} (\hPhi_h(Z)-\psi_v(j))$, 
for appropriate (usually short) subsequences $\cS \subset \{0,1, \ldots, n-1\}$.  

\begin{lemma}[\bf Basic Patching Lemma] \label{DLemCPC4}
\index{patching construction!for nonarchimedean $K_v$-simple sets!basic patching lemmas}
\index{Basic Patching Lemma|ii} 
Let $F_u/K_v$ be a finite extension in $\CC_v$, and 
let  $\cQ(Z) \in F_u[[Z]]$ be a power series converging on $D(0,1)$,
with $\|\cQ\|_{D(0,1)} = 1$.  Suppose that $\cQ(Z)$ has exactly $\ell$ roots
in $D(0,1)$, and that these roots form a $\psi_v$-regular sequence 
\index{regular sequence!$\psi_v$-regular sequence}
$\{\alpha_j\}_{j \in \cS}$ of length $\ell$ in $\cO_u$ with respect to $\hPhi(Z)$.  
Let $M \ge \log_v(\ell)$ be given.  Then for any power series 
$\Delta(Z) \in F_u[[Z]]$ converging on $D(0,1)$, with $\sup$ norm
\begin{equation*}
\|\Delta\|_{D(0,1)} \ \le \ q_v^{-\frac{\ell}{q-1} -M} \ ,   
\end{equation*}
the roots $\alpha_j^*$ of \, $\cQ^*(Z) = \cQ(Z) + \Delta(Z)$ 
again form a $\psi_v$-regular sequence of length $\ell$ in $\cO_u$
\index{regular sequence!$\psi_v$-regular sequence}
with respect to $\hPhi(Z)$. 
They can be uniquely labeled in such a way that
\begin{equation*}
\ord_v(\alpha_j^* - \alpha_j) \ > \ M
\end{equation*}
for each $j \in \cS$.
\end{lemma}

\begin{proof}
If $\cQ(Z) = \sum_{i=0}^{\infty} B_i Z^i$, then since $\|\cQ\|_{D(0,1)} = 1$
it follows that $\ord_v(B_i) \ge 0$ for all $i$ and $\ord_v(B_i) = 0$
for some $i$.  Since $\cQ(Z)$ converges in the closed disc $D(0,1)$,
\begin{equation*}
\lim_{i \rightarrow \infty} \ord_v(B_i) \ = \ \infty \ .
\end{equation*}
Hence there is a largest index $K$ for which $\ord_v(B_K)= 0$, 
and the theory of Newton Polygons\index{Newton Polygon} shows that $K = \ell$ is the number of roots
of $\cQ(Z)$ in $D(0,1)$ (see Lemma \ref{BLem1}).  

Similarly, if $\Delta(Z) = \sum_{i=0}^{\infty} \Delta_i Z^i$,
the fact that $\Delta(Z)$ converges in $D(0,1)$, with
$\|\Delta\|_{D(0,1)} \le q_v^{-\frac{\ell}{q-1}-M}$, tells us that 
$\ord_v(\Delta_i) \ge \frac{\ell}{q-1} + M$ for all $i$.

Now consider $\cQ^*(Z) = \cQ(Z) + \Delta(Z)$.
Writing $\cQ^*(Z) = \sum_{i=0}^{\infty} C_i Z^i$, we have
$C_i = B_i + \Delta_i$, so $\ord_v(C_i) \ge 0$ for all $i$,
$\ord_v(C_\ell) = 0$, and $\ord_v(C_i) > 0$ for all $i > \ell$.  
By the theory of Newton Polygons\index{Newton Polygon}, $\cQ^{*}(Z)$ also has exactly $\ell$
roots in $D(0,1)$.

\vskip .1 in
By Lemma \ref{BLem1} 
\begin{equation*}
\cQ(Z) \ = \  B \cdot \cP(Z) \cdot \cH(Z)
\end{equation*}
where $B \in F_u$ is a constant, $\cP(Z) \in F_u(Z)$ is the polynomial
\begin{equation*}
\cP(Z) \ = \ \cP_{\cS}(Z) \ = \ \prod_{j \in \cS } (Z - \alpha_j),
\end{equation*}
and $\cH(Z) \in \cO_u[[Z]]$ is an invertible power series with constant term $1$.

Since $\cQ(Z)$ has only finitely many roots $\alpha_j$,
there is a point $z_0 \in D(0,1)$ with $|z_0 - \alpha_j|_v = 1$ for all $j$.
At each such point $|\cQ(z_0)|_v =  \|\cQ\|_{D(0,1)} = 1$.
Using $|\cP(z_0)|_v = |\cH(z_0)|_v = 1$, we find that $|B|_v = 1$.

\vskip .1 in
Now fix a root $\alpha_J$, and expand $\cQ(Z)$ and $\Delta(Z)$
as power series in $Z - \alpha_J$:
\begin{eqnarray*}
\cQ(Z)      & = & B_0^{(J)} + B_1^{(J)} (Z-\alpha_J)
                       + B_2^{(J)} (Z-\alpha_J)^2 + \cdots \ , \\
\Delta(Z) & = & \Delta_0^{(J)} + \Delta_{v,1}^{(J)} (Z-\alpha_J) +
                         \Delta_2^{(J)} (Z-\alpha_J)^2 + \cdots \ .
\end{eqnarray*}
We will use the theory of Newton Polygons\index{Newton Polygon} to show that $\cQ^*(Z)$ has a
unique root $\alpha_J^* \in \cO_u$ satisfying $\ord_v(\alpha_J^*-\alpha_J) > M$.

By Lemma \ref{BLem1}, the initial part of the Newton Polygon of $\cQ(Z)$
\index{Newton Polygon}
(expanded about $\alpha_J$) coincides with that of the polynomial $\cP(Z)$,
while its remaining sides have slope $\ge 0$.  Here
\begin{equation*}
\cP(Z) \ = \ \prod_{j \in \cS} ((Z-\alpha_J) - (\alpha_j-\alpha_J)) 
     \ = \ \sum_{i=1}^\ell A_i^{(J)} (Z-\alpha_J)^i \ .
\end{equation*}
Up to sign, the coefficients $A_i^{(J)}$ are elementary symmetric polynomials 
in the $\alpha_j-\alpha_J$.  In particular $A_0^{(J)} = 0$ and
\begin{equation*}
A_1^{(J)} \ = \ \pm \prod_{j \ne J} (\alpha_j - \alpha_J) \ .
\end{equation*}
Since $\{\alpha_j\}_{j \in \cS}$ 
is a $\psi_v$-regular sequence of length $\ell$,
\index{regular sequence!$\psi_v$-regular sequence}
Lemma \ref{DLemCPC3} tells us that
\begin{equation*}
\ord_v(A_1^{(J)}) \ < \ \frac{\ell}{q-1} \ .
\end{equation*}
For each $i \ge 2$, considering the expansions of $A_1^{(J)}$ and $A_i^{(J)}$ 
as elementary symmetric functions in the $\alpha_j-\alpha_J$ gives   
\begin{eqnarray*}
A_i^{(J)} \ = \ \pm A_1^{(J)} \cdot
   \left( \sum_{\substack{ 0 \le j_1 < \cdots < j_{i-1} < \ell \\
                         \text{each $j_{\ell} \ne J$} }}
 \frac{1}{(\alpha_{j_1}-\alpha_J) \cdots (\alpha_{j_{i-1}}-\alpha_J)}
              \right) \ .
\end{eqnarray*}
For each $j \ne J$, we have
$\ord_v(\alpha_j-\alpha_J) = \val_q(|j-J|) < \log_v(\ell)$.  Hence
\begin{equation*}
\ord_v(A_i^{(J)}) \ > \ \ord_v(A_1^{(J)}) - (i-1) \log_v(\ell) \ .
\end{equation*}
Returning to the Newton polygon of $\cQ(Z)$, we see that $B_0^{(J)} = 0$,\index{Newton Polygon} 
\begin{equation*}
\ord_v(B_1^{(J)}) \ = \ \ord_v(A_1^{(J)}) \ < \ \frac{\ell}{q-1} \ ,
\end{equation*}
and for each $i \ge 2$
\begin{equation*}
\ord_v(B_i^{(J)}) \ > \ \ord_v(B_1) - (i-1) \log_v(\ell) \ .
\end{equation*}
(This holds trivially for $i \ge \ell$.)

For the power series $\Delta(Z)$, elementary estimates show that for each $i$ 
\begin{equation*}
\ord_v(\Delta_i^{(J)}) \ \ge \ \frac{\ell}{q-1} + M \ .
\end{equation*}
In particular $\ord_v(\Delta_{v,1}^{(J)}) \ge \frac{\ell}{q-1} > \ord_v(B_1)$,
and $\ord_v(\Delta_i^{(J)}) > \ord_v(B_1) - (i-1) \log_v(\ell)$
for each $i \ge 2$.

\vskip .1 in
Now consider the Newton polygon of $\cQ^*(Z) = \cQ(Z) + \Delta(Z)$,\index{Newton Polygon}
expanded about $\alpha_J$:
\begin{equation*}
\cQ^*(Z) \ = \ C_0^{(J)} + C_1^{(J)} (Z-\alpha_J)
                     + C_2^{(J)} (Z-\alpha_J)^2 + \cdots \ .
\end{equation*}
By the discussion above,
\begin{eqnarray*}
\ord_v(C_0^{(J)}) & = & \ord_v(\Delta_0^{(J)}) \ \ge \ \frac{\ell}{q-1} + M \ , \\
\ord_v(C_1^{(J)}) & = & \ord_v(B_1^{(J)}) \ < \ \frac{\ell}{q-1} \ ,
\end{eqnarray*}
and for each $i \ge 2$,
\begin{equation*}
\ord_v(C_i^{(J)}) \ > \ \ord_v(B_1^{(J)}) - (i-1) \log_v(\ell) \ .
\end{equation*}
Since $M \ge \log_v(\ell)$, the Newton polygon\index{Newton Polygon} of $\cQ^*(Z)$
has a break at $i = 1$, and its initial segment has slope $< -M$.
Hence $\cQ^*(Z)$ has a unique root $\alpha_J^*$ satisfying
\begin{equation*}
\ord_v(\alpha_J^* - \alpha_J) \ > \ M \ .
\end{equation*}
Since $\cQ^*(Z) \in F_u[[Z]]$, the theory of Newton Polygons\index{Newton Polygon} shows that $Z-\alpha_J^*$ 
is a linear factor in the Weierstrass factorization of $\cQ^*(Z)$ over $F_u$
\index{Weierstrass Preparation Theorem}  
(see Proposition \ref{BProp2}(B)).  
Thus $\alpha_J^* \in \cO_u$.

This applies for each $J$.  Since $M \ge \log_v(\ell)$, the $\ell$ roots
$\{\alpha_j^*\}_{j \in \cS}$ 
are distinct and form a $\psi_v$-regular sequence of length $\ell$ in $\cO_u$
\index{regular sequence!$\psi_v$-regular sequence}
attached to $\cS$.
\end{proof}

\vskip .1 in
Lemma \ref{DLemCPC4} will be applied roughly as follows:  

We begin with the function
\index{patching functions, initial $G_v^{(0)}(z)$!for nonarchimedean $K_v$-simple sets!are highly factorized}   
$G_v^{(0)}(z) = S_{n,v}(\phi_v(z)) = \prod_{j=0}^{n-1} (\phi_v(z) - \psi_v(j))$.
The early stages of the patching process seek to preserve this factorization 
as much as possible.  Suppose that at the beginning of the $k^{th}$ stage, 
$k \ge 1$, there is a sequence of $k+1$ consecutive integers 
\begin{equation*}
\cS_k \ = \ \{j_1, j_1+1, \ldots, j_1+k\} \ \subset \ \{0,1, \ldots, n-1 \} \ ,
\end{equation*}
such that $G_v^{(k-1)}(z) = \prod_{j \in \cS_k} (\phi_v(z)-\psi_v(j)) \cdot F_{v,k}(z)$
\index{patching functions, $G_v^{(k)}(z)$ for $1 \le k \le n$!factorization of}
for some $(\fX,\vs)$-function $F_{v,k}(z) \in K_v(\cC)$.  
Expand\index{patching functions, $G_v^{(k)}(z)$ for $1 \le k \le n$!expansion of}
\begin{equation*}
G_v^{(k-1)}(z)
   \ = \ \sum_{i=1}^m \sum_{j=0}^{(n-1)N_i} A_{v,ij} \varphi_{i,nN_i-j}(z)
        + \sum_{\lambda=1}^{\Lambda} A_{v,\lambda} \varphi_{\lambda}(z) \ ,
\end{equation*}
where the $A_{v,ij}, A_{v,\lambda} \in L_{w_v}$.\index{distinguished place $w_v$}
We must patch the coefficients $A_{v,ij}$ in the range $(k-1)N_i \le j < k N_i$.
\index{coefficients $A_{v,ij}$}

Given $(i,j)$ with $(k-1)N_i \le j < kN_i$, 
write $nN_i - j = (n-k-1)N_i + r_{ij}$ where $N_i < r_{ij} \le 2N_i$  
(so $r_{ij} = (k+1)N_i - j$), and put
\begin{equation*}
\vartheta_{v,ij}^{(k)}(z) \ = \ \varphi_{i,r_{ij}}(z) \cdot F_{v,k}(z) \ .
\end{equation*}
Then $\vartheta_{v,ij}^{(k)}(z)$
\index{compensating functions $\vartheta_{v,ij}^{(k)}(z)$!poles and leading coefficients of}  
has a pole of order $nN_i-j$ at $x_i$ and a pole
of order at most $(n-k-1) N_{i^{\prime}}$ at each $x_{i^{\prime}} \ne x_i$.  
By construction, the $\vartheta_{v,ij}^{(k)}(z)$
\index{compensating functions $\vartheta_{v,ij}^{(k)}(z)$!are $K_v$-symmetric}  are $K_v$-symmetric.  
\index{$K_v$-symmetric!set of functions}

Let $\{\Delta_{v,ij}^{(k)}\}$ for $1 \le i \le m$, $(k-1)N_i+1 \le j < kN_i$
be a $K_v$ symmetric set of numbers belonging to $L_{w_v}$,\index{distinguished place $w_v$} 
as given by the global patching process.
\index{patching argument!global} 
Patch $G_v^{(k-1)}(z)$ by setting\index{patching functions, $G_v^{(k)}(z)$ for $1 \le k \le n$!constructed by patching} 
\begin{equation*}
G_v^{(k)}(z) \ = \ 
G_v^{(k-1)}(z) + \sum_{i=1}^m \sum_{j=(k-1)N_i}^{kN_i-1} 
                \Delta_{v,ij}^{(k)} \vartheta_{v,ij}^{(k)}(z) \ .
\end{equation*}
Then the coefficients $A_{v,ij}$ with $j < (k-1)N_i$ are unchanged, 
the\index{coefficients $A_{v,ij}$} 
coefficients $A_{v,ij}$ with $(k-1)N_i < j \le kN_i$ are adjusted
by the $\Delta_{v,ij}$, and the coefficients $A_{v,ij}$ with $j \ge kN_i$
are changed in complicated ways that are unimportant to us.  

If $\Char(K_v) = 0$, then since the $\Delta_{v,ij}^{(k)}$ are $K_v$-symmetric, the sum 
\index{$K_v$-symmetric!set of numbers}\index{compensating functions $\vartheta_{v,ij}^{(k)}(z)$}
\begin{equation*}
\sum_{i=1}^m \sum_{j=(k-1)N_i}^{kN_i-1} 
                \Delta_{v,ij}^{(k)} \vartheta_{v,ij}^{(k)}(z) 
\end{equation*} 
is $K_v$-rational.  If $\Char(K_v) = p > 0$, then\index{compensating functions $\vartheta_{v,ij}^{(k)}(z)$}  
\begin{equation*}
\sum_{i=1}^m \sum_{j=(k-1)N_i}^{kN_i-1} 
                \Delta_{v,ij}^{(k)} \vartheta_{v,ij}^{(k)}(z) 
           \ = \ \Big( \sum_{i=1}^m \sum_{j=(k-1)N_i}^{kN_i-1} 
                \Delta_{v,ij}^{(k)} \varphi_{i,(k+1)N_i-j}(z) \Big) \cdot F_{v,k}(z)
\end{equation*} 
is $K_v$-rational by assumption (\ref{DFCCV1BB}) and the fact that $F_{v,k}(z)$ is $K_v$-rational.  
It follows that $G_v^{(k)}(z)$ is $K_v$-rational.\index{patching functions, $G_v^{(k)}(z)$ for $1 \le k \le n$!are $K_v$-rational}  

Put $G_k(z) = \prod_{j \in \cS_k} (\phi_v(z)-\psi_v(j))$, 
so $G_v^{(k-1)}(z) = G_k(z) \cdot F_{v,k}(z)$.  
By our choice of the $\vartheta_{v,ij}^{(k)}(z)$,
\index{compensating functions $\vartheta_{v,ij}^{(k)}(z)$}
\index{patching functions, $G_v^{(k)}(z)$ for $1 \le k \le n$!constructed by patching} 
\begin{equation*}
G_v^{(k)}(z) \ = \ \left( G_k(z) +
   \sum_{i=1}^m \sum_{j=(k-1)N_i + 1}^{kN_i} \Delta_{v,ij}
                      \cdot \varphi_{i,r_{ij}}(z)
     \right) \cdot F_{v,k}(z) \ .
\end{equation*}
Thus the changes in the roots of $G_v^{(k)}(z)$\index{patching functions, $G_v^{(k)}(z)$ for $1 \le k \le n$!roots are confined to $E_v$} 
have been localized to $G_k(z)$.  
Put $\Delta_{v,k}(z) = \sum_{i=1}^m \sum_{j=(k-1)N_i}^{kN_i-1} \Delta_{v,ij}
                      \cdot \varphi_{i,r_{ij}}(z)$, 
and let $G_k^*(z) = G_k(z) + \Delta_{v,k}(z)$.    

Observe that $G_k(z) = \prod_{j \in \cS_k} (\phi_v(z) - \psi_v(j))$ 
has $\sup$ norm $1$ on each ball $B(\theta_h,\rho_h)$. 
Fixing $h$, when we compose $G_k(z)$ with the parametrization
$\hsigma_h : D(0,1) \rightarrow B(\theta_h,\rho_h)$,
we obtain a function $\cQ_{k,h}(Z)$ whose roots in
$D(0,1)$ form a $\psi_v$-regular sequence of length $k+1$ in $\cO_{u_h}$
\index{regular sequence!$\psi_v$-regular sequence}
attached to $\cS_k$.
Taking $\Delta_{k,h}(Z) = \Delta_{v,k}(\hsigma_h(Z))$, 
we can apply Lemma \ref{DLemCPC4} to 
$G_{k,h}^*(Z) = \cQ_{k,h}(Z) + \Delta_{k,h}(Z)$.
If the $|\Delta_{v,ij}^{(k)}|_v$ are small enough, 
the roots of $\cQ_{k,h}^*(Z)$ will form a $\psi_v$-regular sequence of length
\index{regular sequence!$\psi_v$-regular sequence}
$k+1$ in $\cO_{u_h}$.  Since $\hsigma_h(Z)$ is $F_{u_h}$-rational,
the roots of $G_v^{(k)}(z)$ in $B(a_h,\rho_h)$ belong to $\cC(F_{u_h})$.

\vskip .1 in
The actual patching argument will be more complicated, 
because we eventually we will run out of ``new'' subsequences 
$\cS_k$ to use in patching, and we must deal with roots which 
have been previously patched.  

In the latter steps of the construction, we will use the following lemma:    

\begin{lemma} \label{DLemC6} {\bf (Refined Patching Lemma)}
\index{patching construction!for nonarchimedean $K_v$-simple sets!refined patching lemma}
\index{Refined Patching Lemma|ii} 
Let  $\cQ(Z) \in F_u[[Z]]$ be a power series converging on $D(0,1)$,
with $\sup$ norm $\|\cQ\|_{D(0,1)} = 1$.  Suppose the roots $\{\alpha_j\}$
of $\cQ(Z)$ in $D(0,1)$ belong to $\cO_u$ and can be partitioned into $r$ disjoint
$\psi_v$-regular sequences in $\cO_u$ attached to index sets 
\index{regular sequence!$\psi_v$-regular sequence}
$\cS_1, \ldots, \cS_r$ of respective lengths $\ell_1, \ldots, \ell_r$.
Put $\ell = \sum_{k=1}^r \ell_k$.  Suppose further that there is a bound
$T \ge \max_i(\log_v(\ell_i))$ such that
\begin{equation*}
\ord_v(\alpha_j-\alpha_{k}) \ \le \ T
\end{equation*}
for all $j \ne k$.

Then for any $M \ge T$, and any power series $\Delta(Z) \in F_u[[Z]]$ 
converging on $D(0,1)$ which satisfies 
\begin{equation*}
\|\Delta\|_{D(0,1)} \ \le \ q_v^{-\frac{\ell}{q-1} - (r-1)T - M} \ ,
\end{equation*}
the roots $\{\alpha_j^*\}$ of $\cQ^*(Z) = \cQ(Z) + \Delta(Z)$ in $D(0,1)$ 
again form a union of $\psi_v$-regular
sequences in $\cO_u$ attached to $\cS_1, \ldots, \cS_r$.  They can 
uniquely be labeled in such a way that
\begin{equation*}
\ord_v(\alpha_j^* - \alpha_j) \ > \ M
\end{equation*}
for each $j \in \bigcup_{i=1}^r S_i$.
\end{lemma}

\begin{proof}
By the Weierstrass Preparation Theorem (Lemma \ref{BLem1}), we can write
\index{Weierstrass Preparation Theorem}  
\begin{equation*}
\cQ(Z) \ = \  B \cdot \cP(Z) \cdot \cH(Z)
\end{equation*}
where $B \in F_u$ is a constant, $\cP(Z) \in \cO_u[Z]$ is the polynomial
\begin{equation*}
\cP(Z) \ = \ \prod_{k=1}^r \prod_{j \in \cS_k} (Z - \alpha_j) \ ,
\end{equation*}
and $\cH(Z) \in \cO_u[[Z]]$ is an invertible power series with constant term $1$.
As before, $|B|_v = 1$. We are concerned with the roots of 
\begin{equation*}
\cQ^*(Z) \ = \ \cQ(Z) + \Delta(Z) \ .
\end{equation*} 
As in the proof of the Basic Patching Lemma, the theory of Newton Polygons
\index{Newton Polygon}\index{Basic Patching Lemma}
shows that $\cQ(Z)$ and $\cQ^*(Z)$ both have exactly $\ell$ roots in $D(0,1)$.

\vskip .1 in
Fix a root $\alpha_J$, and expand $\cQ(Z)$ and $\Delta(Z)$
as power series in $Z - \alpha_J$:
\begin{eqnarray*}
\cQ(Z)      & = & B_0^{(J)} + B_1^{(J)} (Z-\alpha_J)
                       + B_2^{(J)} (Z-\alpha_J)^2 + \cdots \ , \\
\Delta(Z) & = & \Delta_0^{(J)} + \Delta_{v,1}^{(J)} (Z-\alpha_J) +
                         \Delta_2^{(J)} (Z-\alpha_J)^2 + \cdots \ .
\end{eqnarray*}

The initial part of the Newton Polygon\index{Newton Polygon} of $\cQ(Z)$
(expanded about $\alpha_J$) coincides with that of $\cP(Z)$,
while its remaining sides have slope $> 0$.  Here
\begin{eqnarray*}
\cP(Z) \ = \ \prod_j ((Z-\alpha_J) - (\alpha_j-\alpha_J)) 
     \ = \ \sum_{i=0}^\ell A_i^{(J)} (Z-\alpha_J)^i \ .
\end{eqnarray*}
The coefficients $A_i^{(J)}$ are  symmetric polynomials in 
\index{coefficients $A_{v,ij}$} 
the $\alpha_j-\alpha_J$.  In particular $A_0^{(J)} = 0$ and
\begin{equation*}
A_1^{(J)} \ = \ \pm \prod_{j \ne J} (\alpha_j - \alpha_J) \ .
\end{equation*}
Suppose $J \in \cS_i$.  By part (A) of Lemma \ref{DLemCPC3},
\begin{equation*}
\ord_v(\prod_{\substack{ j \in \cS_i \\ j \ne J }} (\alpha_J-\alpha_j))
    \ < \ \frac{\ell_i}{q-1} \ .
\end{equation*}
For each $k \ne i$, by part (B) of Lemma \ref{DLemCPC3}, 
\begin{equation*}
\ord_v(\prod_{ j \in \cS_k } (\alpha_J-\alpha_j))
    \ < \ \frac{\ell_k}{q-1} + T
\end{equation*}
since $\max_{j \in \cS_k}(\ord_v(\alpha_J - \alpha_j)) \le T$.
Summing these, we see that
\begin{equation*}
\ord_v(A_1^{(J)}) \ < \ \frac{\ell}{q-1} + (r-1)T \ .
\end{equation*}

For each $i \ge 2$
\begin{eqnarray*}
A_i^{(J)} & = & \pm A_1^{(J)} \cdot
   \left( \sum_{\substack{ j_1 < \cdots < j_{i-1}  \\ \text{each $j_k \ne J$} }}
 \frac{1}{(\alpha_{j_1}-\alpha_J) \cdots (\alpha_{j_{i-1}}-\alpha_J)}
              \right) 
\end{eqnarray*}
just as in the Basic Patching Lemma.\index{Basic Patching Lemma}  We have
$\ord_v(\alpha_j-\alpha_J) \le T$ for each $j \ne J$.  Hence
\begin{equation*}
\ord_v(A_i^{(J)}) \ \ge \ \ord_v(A_1^{(J)}) - (i-1) T \ .
\end{equation*}
Returning to the Newton polygon\index{Newton Polygon} of $\cQ(Z)$, we see that  $B_0^{(J)} = 0$, 
\begin{equation*}
\ord_v(B_1^{(J)}) \ = \ \ord_v(A_1^{(J)}) \ < \ \frac{\ell}{q-1} + (r-1)T \ ,
\end{equation*}
and for each $i \ge 2$
\begin{equation*}
\ord_v(B_i^{(J)}) \ \ge \ \ord_v(B_1^{(J)}) - (i-1)T \ .
\end{equation*}
(Note that $T > 0$ unless $\ell = 0$, in which case there is nothing to prove.)

For the power series $\Delta(Z)$, elementary estimates give
\begin{equation*}
\ord_v(\Delta_i^{(J)}) \ \ge \ \frac{\ell}{q-1} + (r-1)T + M
\end{equation*}
for each $i$.
In particular $\ord_v(\Delta_{v,1}^{(J)}) > \frac{\ell}{q-1} + (r-1)T$,
and $\ord_v(\Delta_i^{(J)}) > \ord_v(B_1) - (i-1)T$
for each $i \ge 2$.

\vskip .1 in
Now consider the Newton polygon\index{Newton Polygon} of $\cQ^*(Z) = \cQ(Z) + \Delta(Z)$,
expanded about $\alpha_J$:
\begin{equation*}
\cQ^*(Z) \ = \ C_0^{(J)} + C_1^{(J)} (Z-\alpha_J)
                     + C_2^{(J)} (Z-\alpha_J)^2 + \cdots \ .
\end{equation*}
By the discussion above,
\begin{eqnarray*}
\ord_v(C_0^{(J)}) & = & \ord_v(\Delta_0^{(J)})
               \ \ge \ \frac{k}{q-1} + (r-1)T + M \ , \\
\ord_v(C_1^{(J)}) & = & \ord_v(B_1^{(J)}) \ < \ \frac{k}{q-1} + (r-1)T \ ,
\end{eqnarray*}
and for each $i \ge 2$,
\begin{equation*}
\ord_v(C_i^{(J)}) \ \ge \ \ord_v(B_1^{(J)}) - (i-1) T \ .
\end{equation*}
Since $M \ge T$, the Newton polygon\index{Newton Polygon} of $\cQ^*(Z)$
has a break at $i = 1$, and its initial segment has slope $< -M$.
Hence $\cQ^*(Z)$ has a unique root $\alpha_J^*$ satisfying
\begin{equation*}
\ord_v(\alpha_J^* - \alpha_J) \ > \ M \ .
\end{equation*}
Since $\cQ^*(Z) \in F_u[[Z]]$, the theory of Newton Polygons shows that $Z-\alpha_J^*$ 
\index{Newton Polygon}
is a linear factor in the Weierstrass factorization of $\cQ^*(Z)$ over $F_u$ 
\index{Weierstrass Preparation Theorem}  
(see Proposition \ref{BProp2}(B)). 
Thus $\alpha_J^* \in \cO_u$.

This applies for each $J$.  Since $M \ge T$, the roots $\alpha_j^*$
are distinct and form a union of $\psi_v$-regular sequences in $\cO_u$
\index{regular sequence!$\psi_v$-regular sequence}
attached to $\cS_1, \ldots, \cS_r$.  
\end{proof}

\section{ Stirling Polynomials when $\Char(K_v) = p > 0$.} \label{StirlingPolypSection} 
\index{Stirling polynomial!for $\cO_v$!when $\Char(K_v) = p > 0$}

In this section, assume that $\Char(K_v) = p > 0$. 
 
In Proposition \ref{SnCharpForm} below, 
we will show that by requiring that $n$ be divisible by a sufficiently high power of $p$,
we can make arbitrarily many high order coefficients of the Stirling polynomial 
\index{Stirling polynomial}
$S_{n,v}(z) = \prod_{j=0}^{n-1} (z-\psi_v(j))$ be $0$.
\index{Stirling polynomial!for $\cO_v$!high-order coefficients vanish}   
This fact plays a key role in the degree-raising argument in
the proof of Theorem \ref{DCPCPatch1p} in \S\ref{Chap11}.\ref{NonArchPatchingProof}.

We begin with a lemma concerning homogeneous products of linear forms over a finite field.  

\begin{lemma} \label{QTLem}  Let $\FF_q$ be the finite field with $q$ elements, 
let $r \ge 1$ be an integer, and put 
\begin{equation} \label{QTDef} 
Q_r(z;T_1, \ldots, T_r) \ = \ 
\prod_{a_1, \ldots, a_r \in \FF_q} 
                      \big( z - \sum_{i=1}^r a_i T_i \big) \ .
\end{equation} 
Then $Q_r$ has the form 
\begin{equation} \label{QTForm} 
Q_r(z;T_1, \ldots, T_r) 
\ = \ z^{q^r} + \sum_{\ell=1}^r P_{r,\ell}(T_1, \ldots, T_r) \cdot z^{q^{r-\ell}} 
\end{equation} 
where $P_{r,\ell}(T_1, \ldots,T_r)$ is a homogeneous 
polynomial of degree $q^r - q^{r-\ell}$ in $\FF_q[T_1, \ldots, T_r]$ 
which is symmetric in $T_1, \ldots, T_r$, for each $\ell=1, \ldots, r$. 
\end{lemma}  

\begin{proof} Since $Q_r$ is symmetric in $T_1, \ldots, T_r$ and is homogeneous of degree $q^r$, 
if $Q_r$ has the form (\ref{QTForm}), necessarily the $P_{r,\ell}(T_1, \ldots, T_r)$ 
are symmetric in $T_1, \ldots, T_r$ and homogeneous of degree $q^r-q^{r-\ell}$.  

We now prove (\ref{QTForm}) by induction on $r$.  When $r = 1$, 
\begin{eqnarray} 
Q_1(z;T) & = & \prod_{a\in \FF_q} (z - a T) \ = \ T^q \cdot \prod_{a \in \FF_q} (\frac{z}{T} - a) 
                      \notag  \\
      & = & T^q \cdot \big( \big(\frac{z}{T}\big)^q - \frac{z}{T} \big) \ = \ z^q - T^{q-1} \cdot z \ .
                      \label{FQT1Form}
\end{eqnarray} 
Now suppose that (\ref{QTForm}) holds for some $r$.  Then for $r + 1$,  
\begin{eqnarray*}
Q_{r+1}(z;T_1, \ldots, T_{r+1}) & = & \prod_{a_1 \in \FF_q} 
        \Big( \prod_{a_2 \in \FF_q} \cdots \prod_{a_{r+1} \in \FF_q} 
                                 \Big((z - a_1 T_1) - \sum_{i=2}^{r+1} a_i T_i \Big) \Big) \\
       & = & \prod_{a_1 \in \FF_q} Q_r(z-a_1T_1; T_2, \ldots, T_{r+1}) \ .
\end{eqnarray*}  
Using (\ref{QTForm}) for  $Q_r(z-a_1T_1; T_2, \ldots, T_{r+1})$, 
noting that $(z-a_1 T_1)^q = z^q- a_1 T_1^q$, and applying (\ref{FQT1Form}), we see that
$Q_{r+1}(z;T_1, \ldots, T_{r+1})$ has the form 
\begin{eqnarray*}
Q_{r+1}(z;T_1, \ldots, T_{r+1})                                 
    & = & \prod_{a_1 \in \FF_q} \Big( \big(z^{q^r} - a_1 T_1^{q^r}\big) + 
        \sum_{\ell=1}^r \big(z^{q^{r-\ell}}-a_1 T_1^{q^{r-\ell}}\big) 
                                   P_{r,\ell}(T_2,\ldots,T_{r+1}) \Big) \\
    & = & \prod_{a_1 \in \FF_q} \Big( Q_r(z;T_2, \ldots, T_{r+1}) 
                                 - a_1 Q_r(T_1;T_2, \ldots, T_r) \Big) \\
    & = & Q_r(z;T_2, \ldots, T_{r+1})^q 
               - Q_r(z;T_2, \ldots, T_{r+1})\cdot Q_r(T_1;T_2, \ldots, T_r)^{q-1} \ .
\end{eqnarray*} 
Using (\ref{QTForm}) for $Q_r(z;T_2, \ldots, T_{r+1})$, 
then expanding the $q$-th power and collecting terms, we obtain
\begin{eqnarray*}
Q_{r+1}(z;T_1, \ldots, T_{r+1})  
   \ = \ z^{q^{r+1}} + \sum_{\ell=1}^{r+1}  
                        P_{r+1,\ell}(T_1, \ldots, T_{r+1}) \cdot z^{q^{r+1-\ell}} \ .            
\end{eqnarray*}
\vskip -.35 in
\end{proof} 

\vskip .25 in
Now let $\pi_v$ be a uniformizer for the 
maximal ideal of $\cO_v$.
Let $q = q_v = p^{f_v}$ be the order of the residue field $\cO_v/\pi_v \cO_v$.
By the structure theory of local fields in positive characteristic,   
$\cO_v \cong \FF_q[[\pi_v]]$ and $K_v \cong \FF_q((\pi_v))$.  The following proposition uses the
fact that in the basic well-distributed sequence $\{\psi_v(k)\}_{0 \le k < \infty}$ for $\cO_v$,
the representatives $\psi_v(0), \ldots, \psi_v(q-1)$ for $\cO_v/\pi_v \cO_v$ are the 
Teichm\"uller representatives, the elements of $\FF_q$.\index{Teichm\"uller representatives}    

\begin{proposition} \label{SnCharpForm} 
Let $K_v \cong \FF_q((\pi_v))$ be a local field of characteristic $p > 0$.  If $q^r | n $
for some $r > 0$, then $S_{n,v}(z)$ can be expanded as  
\begin{equation} \label{SnTerms} 
S_{n,v}(z) \ = \ z^n + \sum_{j=q^r-q^{r-1}}^n C_j \, z^{n-j} 
\end{equation} 
with each $C_j \in \cO_v$.
\end{proposition} 

\begin{proof}
First suppose $n = q^r$ for some $r > 0$.  
The numbers $\psi_v(k)$, for $k = 0, \ldots, q^r-1$,
run over all possible sums $a_0 + a_1 \pi_v + \cdots + a_{r-1} \pi_v^{r-1}$   
with $a_0, \ldots, a_{r-1} \in \FF_q$.  By Lemma \ref{QTLem}, 
\begin{equation} 
S_{q^r,v}(z) \ = \ Q_r(z;1,\pi_v, \ldots, \pi_v^{r-1}) 
\ = \ z^{q^r} + A_r z^{q^{r-1}} + \text{\rm terms of lower degree} \ , 
\end{equation}
where $A_r = P_{r,1}(1,\pi_v, \ldots, \pi_v^{r-1})$. 
Now suppose $n = q^r \ell$.  If $0 \le k < \ell$ and we write $k = j q^r + s$ 
with $0 \le j < \ell$, $0 \le s < q^r$,
then $\psi_v(k) = \psi_v(j) \pi_v^r + \psi_v(s)$.  It follows that 
\begin{eqnarray}
S_{n,v}(z) 
& = & \prod_{j=0}^{\ell-1} \prod_{s = 0}^{q^r-1} (z - \psi_v(j) \pi_v^r - \psi_v(s))
\ = \  \prod_{j=0}^{\ell-1} S_{q^r,v}(z - \psi_v(j) \pi_v^r) \notag \\
& = &  \prod_{j=0}^{\ell-1} \Big((z - \psi_v(j) \pi_v^r)^{q^r} 
          + A_r \cdot (z - \psi_v(j) \pi_v^r)^{q^{r-1}} 
              + \text{\rm terms of lower degree} \Big) \ . \label{SqrF2} 
\end{eqnarray}  
Since $(z-\psi_v(j) \pi_v^r)^{q} = z^q - \psi_v(j)^q \pi_v^{qr}$, upon multiplying out (\ref{SqrF2}) 
we see that 
\begin{equation*} 
S_{n,v}(z) 
\ = \ (z^{q^r})^\ell + \ell \cdot (z^{q^r})^{\ell-1} (A_r z^{q^{r-1}}) 
+ \text{\rm terms of lower degree} \ ,
\end{equation*} 
which yields (\ref{SnTerms}).
\end{proof}


\section{ Proof of Theorems $\ref{DCPCPatch}$ and $\ref{DCPCPatch1p}$ } \label{NonArchPatchingProof} 

In this section we prove Theorems \ref{DCPCPatch} and \ref{DCPCPatch1p}.
The construction is a generalization of those in (\cite{RR2}) and (\cite{RR3}).  

\vskip .05 in
We begin the construction with\index{patching functions, initial $G_v^{(0)}(z)$!construction of}  
\begin{eqnarray*}
G_v^{(0)}(z) & = & 
S_{n,v}(\phi_v(z)) \ = \ \prod_{j=0}^{n-1} (\phi_v(z)-\psi_v(j)) \\
& = & \sum_{i=1}^m \sum_{j=0}^{(n-1)N_i-1} A_{v,ij} \varphi_{i,nN_i-j}
      + \sum_{\lambda=1}^{\Lambda} A_{v,\lambda} \varphi_{\lambda} \ , 
\end{eqnarray*}
whose leading coefficient at $x_i$ is $A_{v,i0} = \tc_{v,i}^n$.  
\index{coefficients $A_{v,ij}$!leading}

The roots $\{\theta_{hj}\}_{1 \le h \le N, 0 \le j < n}$ 
of $G_v^{(0)}(z)$
\index{patching functions, initial $G_v^{(0)}(z)$!for nonarchimedean $K_v$-simple sets!roots are distinct}  
belong to $H_v := E_v \cap \big(\bigcup_{h=1}^N B(\theta_h,\rho_h)\big)$,
and are distinct.  
Indeed, as noted at the beginning of \S\ref{Chap11}.\ref{PatchingLemmasSection}, 
for each $h$ there is a $1-1$ correspondence between the roots 
$\theta_{hj}$ of $S_{n,v}(\phi_v(z))$ in $B(\theta_h,\rho_h)$ 
and the zeros $\alpha_{hj}$ of $S_{n,v}(\hPhi_h(Z))$ in $D(0,1)$,
given by $\theta_{hj} = \hsigma_h(\alpha_{hj})$ for $j = 0, \ldots, n-1$.  
Moreover, the $\alpha_{hj}$ belong to $\cO_{u_h}$,
and form a $\psi_v$-regular sequence of length $n$ in $\cO_{u_h}$. 
\index{regular sequence!$\psi_v$-regular sequence}
Put
\begin{equation} \label{FtUvDef} 
U_v^0 \ := \ \bigcup_{h=1}^N B(\theta_h,\rho_h) \ \subset \ U_v  
\ = \ \bigcup_{\ell = 1}^D B(a_\ell,r_\ell) \ .
\end{equation} 

For suitable $n$, the patching process will inductively construct $K_v$-rational 
$(\fX,\vs)$-functions $G_v^{(k)}(z)$, for $k = 1, \ldots, n$,
\index{patching functions, $G_v^{(k)}(z)$ for $1 \le k \le n$!constructed by patching} 
whose roots belong to $H_v$ for all $k$.  
For each $k$, there will be a natural $1-1$ correspondence between\index{roots!natural $1-1$ correspondence|ii}
the roots of $G_v^{(k-1)}(z)$ and $G_v^{(k)}(z)$, 
and the roots of $G_v^{(n)}(z)$ will be distinct.
\index{patching functions, $G_v^{(k)}(z)$ for $1 \le k \le n$!for nonarchimedean $K_v$-simple sets!roots are distinct}  
The conditions on $n$ needed for the construction to succeed will be noted as they arise;  
they will all require that $n$ be sufficiently large, 
or that it be divisible by a certain integer.  Only a finite number of conditions will be imposed,
so there is an $n_v \ge 1$ such that the construction will succeed if $n$ 
is sufficiently large and divisible by $n_v$.  

\smallskip
Recall that  
\begin{equation} \label{FBD} 
M_v \ = \
   \max \big(\max_i \max_{N_i < j \le 2N_i} (\|\varphi_{i,j}\|_{U_v}),
       \max_{\lambda} (\|\varphi_{\lambda}\|_{U_v})\big) \ .
\end{equation} 
We have defined $k_v$ to be the smallest integer such that for all $k \ge k_v$,
\begin{equation} \label{FkvBound}
h_v^{Nk} \cdot M_v \ < \ q_v^{-(\frac{k+1}{q_v - 1} + \log_v(k+1))} \ .
\end{equation}
The global patching process specifies a number $\kbar \ge k_v$, 
the number of bands of coefficients considered ``high-order''.
\index{band!high-order}\index{coefficients $A_{v,ij}$!high-order}  

\vskip .1 in
\noindent{\bf Phase 1.  Patching the leading and high-order coefficients, for $k = 1, \ldots, \kbar$.}
\index{coefficients $A_{v,ij}$!leading}
\index{coefficients $A_{v,ij}$!high-order}

The patching constructions for the leading and high-order coefficients 
are different when $\Char(K_v) = 0$ and when $\Char(K_v) = p > 0$.  

\smallskip
{\bf Case A.}  Suppose $\Char(K_v) = 0$.
\index{local patching for nonarch $K_v$-simple sets!Phase 1: leading and high-order coefficients!when $\Char(K_v) = 0$|(}  

In this case, the leading coefficient is patched along with 
\index{coefficients $A_{v,ij}$!leading} 
the other high order coefficients.  The bound for the high order patching 
\index{coefficients $A_{v,ij}$!high-order}
coefficients in Theorem \ref{DCPCPatch} is $B_v = h_v^{\kbar N}$.  
Since $0 < h_v < 1$ and $\kbar \ge k_v$, it follows from (\ref{FkvBound}) that 
\begin{equation} \label{FBM}
B_v \cdot M_v \ = \ h_v^{\kbar N} \cdot M_v \ \le \ q_v^{-\frac{\kbar+1}{q_v - 1} - \log_v(\kbar+1)} 
\ \le \ q_v^{-\frac{k+1}{q_v - 1} - \log_v(k+1)} 
\end{equation}
for each $k = 1, \ldots, \kbar$. 

\smallskip
When $k = 1$, we are given a $K_v$-symmetric set of numbers 
\index{$K_v$-symmetric!set of numbers}
$\{\Delta_{v,ij}^{(1)} \in L_{w_v}\}_{(i,j) \in \Band_N(1)}$,\index{distinguished place $w_v$} 
determined recursively in $\prec_N$ order,\index{order!$\prec_N$}\index{band!$\Band_N(k)$}
such that for each $(i,j)$,  
\begin{equation}
|\Delta_{v,ij}^{(1)}|_v \ \le \ B_v \ .   \label{FBUG1}
\end{equation}
Put $\cS_1 = \{0, 1 \}$, and let
\begin{equation*}
 P_1(x) \ = \ \prod_{j \in \cS_1} (x-\psi_v(j)), \qquad
\hP_1(x) \ = \ \prod_{j=2}^{n-1} (x-\psi_v(j))
\end{equation*}
so $S_{n,v}(x) = P_1(x) \cdot \hP_1(x)$.  
Set $Q_1(z) = P_1(\phi_v(z))$; then
\index{patching functions, initial $G_v^{(0)}(z)$!for nonarchimedean $K_v$-simple sets!are highly factorized}  
\begin{equation*}
G_v^{(0)}(z) \ = \ Q_1(z) \cdot \hP_1(\phi(z)) \ .
\end{equation*} 

For each $(i,j) \in \Band_N(1)$, 
put\index{band!$\Band_N(k)$}\index{compensating functions $\vartheta_{v,ij}^{(k)}(z)$!construction of}
\begin{equation*}
\vartheta_{v,ij}^{(1)}(z) \ = \ \varphi_{i,2N_i-j}(z) \cdot \hP_1(\phi_v(z)) \ .
\end{equation*}
Thus $\vartheta_{v,ij}^{(1)}(z)$ has a pole of order $nN_i-j$ at $x_i$ with
leading coefficient $\tc_{v,i}^{n-2}$, and a pole of order $(n-2)N_{i^{\prime}}$
\index{coefficients $A_{v,ij}$!leading}
at each $x_{i^{\prime}} \ne x_i$.  
By construction the $\vartheta_{v,ij}^{(1)}(z)$ are $K_v$-symmetric.
\index{compensating functions $\vartheta_{v,ij}^{(k)}(z)$!poles and leading coefficients of} 
\index{$K_v$-symmetric!set of functions} 
Since the $\Delta_{v,ij}^{(1)}$ are $K_v$-symmetric as well, for each $(i,j)$
\index{$K_v$-symmetric!set of numbers}
\begin{equation*}
\sum_{x_{i^{\prime}} \in \Aut_c(\CC_v/K_v)(x_i)} 
    \Delta_{v,i^{\prime}j}^{(1)} \vartheta_{v,i^{\prime}j}^{(1)} \ \in \ K_v(\cC) \ .
\end{equation*}
This assures that when the global patching process determines the patching coefficients 
\index{coefficients $A_{v,ij}$!patching}
$\Delta_{v,ij}^{(k)}(z)$ recursively in $\prec_N$ order,\index{order!$\prec_N$} 
the partially patched function $G_v^{(0)}(z)$ is 
$K_v$-rational\index{patching functions, initial $G_v^{(0)}(z)$!are $K_v$-rational}  
after each $\Aut(\tK/K)$-orbit of coefficients is chosen. 
\index{coefficients $A_{v,ij}$}  

Patch $G_v^{(0)}(z)$\index{patching functions, initial $G_v^{(0)}(z)$} 
 by setting\index{compensating functions $\vartheta_{v,ij}^{(k)}(z)$}
\index{patching functions, $G_v^{(k)}(z)$ for $1 \le k \le n$!constructed by patching}
\begin{equation*}
G_v^{(1)}(z) \ = \ G_v^{(0)}(z) 
+ \sum_{i=1}^m \sum_{j=0}^{N_i-1} \Delta_{v,ij}^{(1)} \cdot \vartheta_{v,ij}^{(1)}(z) \ .
\end{equation*}
Then $G_v^{(1)}(z)$ is a $K_v$-rational $(\fX,\vs)$-function of degree $Nn$.
\index{patching functions, $G_v^{(k)}(z)$ for $1 \le k \le n$!are $K_v$-rational}  

Now consider how the roots $\{\theta_{hj}\}$ change 
in passing from
\index{patching functions, $G_v^{(k)}(z)$ for $1 \le k \le n$!for nonarchimedean $K_v$-simple sets!movement of roots}  
to $G_v^{(0)}(z)$ to $G_v^{(1)}(z)$.  If we write
$\Delta_{v,1}(z) = \sum_{i=1}^m \sum_{j=0}^{N_i-1}
                 \Delta_{v,ij}^{(1)} \cdot \varphi_{i,2N_i-j}(z)$,
and put 
\begin{equation*}
G_1^*(z) \ = \ Q_1(z) + \Delta_{v,1}(z) \ ,
\end{equation*}
then by our choice of the $\vartheta_{v,ij}^{(1)}(z)$\index{compensating functions $\vartheta_{v,ij}^{(k)}(z)$}  
\begin{equation*}
G_v^{(1)}(z) \ = \ G_1^*(z) \cdot \hP_1(\phi_v(z)) \ .
\end{equation*}  
Hence the roots $\theta_{hj}$ with $j \ge 2$ are all preserved.
 
For each ball $B(\theta_h,\rho_h)$,
let $\hsigma_h : D(0,1) \rightarrow B(\theta_h,\rho_h)$
be the $F_{u_h}$-rational parametrization 
from \S\ref{Chap11}.\ref{PatchingLemmasSection};   
recall that $\theta_{hj} = \hsigma_h(\alpha_{hj})$.
By abuse of language, we will refer to both the $\theta_{hj}$ 
and the $\alpha_{hj}$ as ``roots''.\index{roots} 
By (\ref{FBD}), (\ref{FBM}) and (\ref{FBUG1}),  
\begin{equation*}
\|\Delta_{v,1}(z)\|_{U_v^0} \ \le \ B_v \cdot M_v
                    \ \le \ q^{-(2/(q-1) + \log_v(2))} \ .
\end{equation*}
Put $\cQ_{1,h}(Z) = Q_1(\hsigma_h(Z)) = P_1(\hPhi_h(Z)))$ and  
$\Delta_{1,h}(Z) = \Delta_{v,1}(\hsigma_h(Z))$.
The roots $\{\alpha_{h0}, \alpha_{h1}\}$ of $\cQ_{1,h}(Z)$ 
form a $\psi_v$-regular sequence of length $2$ in $\cO_{u_h}$, and
\index{regular sequence!$\psi_v$-regular sequence}
$\|\Delta_{1,h}\|_{D(0,1)} \le q^{-(2/(q-1) + \log_v(2))}$.  
By Lemma \ref{DLemCPC4} the roots of 
$\cQ_{1,h}^*(Z) = \cQ_{1,h}(Z) + \Delta_{1,h}(Z)$
form a $\psi_v$-regular sequence $\{ \alpha_{h0}^*, \alpha_{h1}^*\}$
\index{regular sequence!$\psi_v$-regular sequence}
of length $2$ in $\cO_w$, with
\begin{equation*}
\ord_v(\alpha_{h0}^*-\alpha_{h0}) > \log_v(2), \quad
\ord_v(\alpha_{h1}^*-\alpha_{h1}) > \log_v(2) \ .
\end{equation*} 
Thus the roots of
\index{patching functions, $G_v^{(k)}(z)$ for $1 \le k \le n$!for nonarchimedean $K_v$-simple sets!movement of roots}
$G_v^{(1)}(\hsigma_h(Z))$ in $D(0,1)$ are
\begin{equation*}
\{ \alpha_{h0}^*, \alpha_{h1}^*,
           \alpha_{h2}, \ldots, \alpha_{h,n-1}\} \ ,
\end{equation*} 
a union of a $\psi_v$-regular sequence of length $2$ 
\index{regular sequence!$\psi_v$-regular sequence}
and the remaining $n-2$ elements of the original $\psi_v$-regular sequence of length $n$.
\index{regular sequence!$\psi_v$-regular sequence}
Transferring this back to $G_v^{(1)}(z)$,
\index{patching functions, $G_v^{(k)}(z)$ for $1 \le k \le n$!for nonarchimedean $K_v$-simple sets!movement of roots}
put $\theta_{h0}^* = \hsigma_h(\alpha_{h0}^*)$, 
$\theta_{h1}^* = \hsigma_h(\alpha_{h1}^*)$.  
Then $\theta_{h0}^*, \theta_{h1}^* \in 
B(\theta_h,\rho_h) \cap \cC(F_{u_h}) \subset \tE_v$.
The roots of $G_v^{(1)}(z)$ in $B(\theta_h,\rho_h)$ are
\index{patching functions, $G_v^{(k)}(z)$ for $1 \le k \le n$!for nonarchimedean $K_v$-simple sets!movement of roots} 
\begin{equation*}
\{ \theta_{h0}^*, \theta_{h1}^*,
           \theta_{h2}, \ldots, \theta_{h,n-1}\} \ .
\end{equation*} 

Let $\varepsilon_{v,i} = 1 + \Delta_{v,i0}^{(1)}$.  Since the $\Delta_{v,i0}^{(1)}$
are $K_v$-symmetric, with $|\Delta_{v,i0}^{(1)}|_v < 1$ for each $i$, 
\index{$K_v$-symmetric!set of numbers}
the $\varepsilon_{v,i}$ form a $K_v$-symmetric system of units in 
\index{$K_v$-symmetric!system of units}
$\cO_{w_v}^{\times}$.\index{distinguished place $w_v$}  The leading coefficient $A_{v,i0} = \tc_{v,i}^n$
\index{coefficients $A_{v,ij}$!leading}\index{patching functions, $G_v^{(k)}(z)$ for $1 \le k \le n$!leading coefficients of} 
of $G_v^{(0)}(z)$ at $x_i$ is changed to $\varepsilon_{v,i} \tc_{v,i}^n$ 
in $G_v^{(1)}(z)$.  (From the global patching process, 
we know that $\varepsilon_{v,i} \tc_{v,i}^N = \mu_i^{n/n_0}$ where the $\mu_i$
are the $\hS$-units from Theorem \ref{CTCX2};  
however, from a local standpoint, this is irrelevant.)  The leading coefficient
\index{coefficients $A_{v,ij}$!leading}
of $Q_1^*(z)$ at $x_i$ is $\varepsilon_{v,i} $.  

\smallskip
Next we construct the $G_v^{(k)}(z)$ for $k = 2, \ldots, \kbar$.
\index{patching functions, $G_v^{(k)}(z)$ for $1 \le k \le n$!constructed by patching}
Each $G_v^{(k)}(z)$ will be a $K_v$-rational $(\fX,\vs)$-function
\index{patching functions, $G_v^{(k)}(z)$ for $1 \le k \le n$!are $K_v$-rational}
of degree $nN$, having a pole of order $nN_i$ and
leading coefficient $\varepsilon_{v,i} \tc_{v,i}^n$ 
\index{coefficients $A_{v,ij}$!leading}
at $x_i$, and with a factorization into $K_v$-rational $(\fX,\vs)$-functions 
of the form\index{patching functions, $G_v^{(k)}(z)$ for $1 \le k \le n$!factorization of} 
\begin{equation*}
G_v^{(k)}(z) \ = \ Q_1^*(z) Q_2^*(z) \cdots Q_k^*(z) \cdot \hP_k(\phi_v(z)) \ .
\end{equation*}
The properties of the $Q_k^*(z)$ will be discussed below.

\vskip .1 in
Recall that $\cS_1 = \{0,1\}$.  
Put $\cS_2 = \{2,3,4\}$, $\cS_3 = \{5,6,7,8\}$, and so on, through $\cS_{\kbar}$, 
where $\cS_k$ consists of the next $k+1$ integers in $\{0,1, \ldots, n-1\}$ 
after $\cS_{k-1}$.  Explicitly, 
\begin{equation} \label{DFFT1}
\cS_k \ = \ \{ j_k, \ldots, j_k+k \}
\end{equation}
where $j_k = (k^2 + k -2)/2$.  For this to be possible we need
\begin{equation} \label{DFXYZ1}
n \ \ge \ (\kbar^2+3\kbar-2)/2
\end{equation}
which we henceforth assume.  

In the $k^{th}$ step, the roots corresponding to $\cS_k$ will be moved.\index{roots!move roots}  Put
\begin{equation} \label{DFFQ1}
P_k(x) = \prod_{j \in \cS_k} (x-\psi_v(j)), \quad
\hP_k(x) = \prod_{0 \le j \le n-1, \ j \notin \cS_1 \cup \ldots \cS_k} (x-\psi_v(j)) \ ,
\end{equation}
so that $S_{n,v}(x) = P_1(x) \cdots P_k(x) \cdot \hP_k(x)$ and
$\hP_{k-1}(x) = P_k(x) \cdot \hP_k(x)$.  Set 
\begin{equation}
Q_k(z) \ = \ P_k(\phi_v(z)) \ . \label{DFFR2}
\end{equation} 
For each $k$, $2 \le k \le \kbar$, the function $Q_k^*(z)$ 
will have the following properties:
\begin{enumerate}
\item $Q_k^*(z)$ is obtained by perturbing $Q_k(z)$; 
\item $Q_k^*(z)$ is a $K_v$-rational $(\fX,\vs)$-function
of degree $(k+1)N$, with a pole of order $(k+1)N_i$ and leading
coefficient $\tc_{v,i}^{k+1}$ at $x_i$;
\index{coefficients $A_{v,ij}$!leading}
\item For each ball $B(\theta_h,\rho_h)$, 
the roots $\{ \alpha_{h,j_k }^*, \ldots, \alpha_{h,j_k + k}^* \}$ 
of $Q_k^*(\hsigma_h(Z))$ in $D(0,1)$ form a $\psi_v$-regular sequence
\index{regular sequence!$\psi_v$-regular sequence}
$\{\alpha_{hj}\}_{j \in \cS_k}$ of length $k+1$ in $\cO_{u_h}$ 
relative to $\hPhi_h(Z)$, and for each $j \in \cS_k$
\begin{equation*}
\ord_v(\alpha_{hj}^*-\alpha_{hj})
        \ > \ \frac{k+1}{q-1} + \log_v(k+1);
\end{equation*} 
\item For each ball $B(\theta_h,\rho_h)$, 
$\|Q_k^*\|_{B(\theta_h,\rho_h)} = 1$. 
\end{enumerate}

\vskip .1 in
Inductively suppose $G_v^{(k-1)}(z)$ has been constructed, 
with\index{patching functions, $G_v^{(k)}(z)$ for $1 \le k \le n$!factorization of}   
\begin{eqnarray}
G_v^{(k-1)}(z) & = & Q_1^*(z) \cdots Q_{k-1}^*(z) \cdot \hP_{k-1}(\phi_v(z)) 
                                   \notag \\
      & = & Q_1^*(z) \cdots Q_{k-1}^*(z) \cdot Q_k(z) \cdot \hP_k(\phi_v(z)) \ .
                                   \label{FGFact} 
\end{eqnarray}
Expand\index{patching functions, $G_v^{(k)}(z)$ for $1 \le k \le n$!expansion of}
\begin{equation*}
G_v^{(k-1)}(z) \ = \
      \sum_{i=1}^m \sum_{j=0}^{(n-1)N_i-1} A_{v,ij} \varphi_{i,nN_i-j}(z)
            + \sum_{\lambda = 1}^{\Lambda} A_{\lambda} \varphi_{\lambda}(z) \ .
\end{equation*}
The $A_{v,ij}$ and $A_{v,\lambda}$ belong to $L_{w_v}$\index{distinguished place $w_v$} and are $K_v$-symmetric.
\index{$K_v$-symmetric!set of functions}

We are given a $K_v$-symmetric set 
\index{$K_v$-symmetric!set of numbers}
of numbers $\{\Delta_{v,ij}^{(k)} \in L_{w_v}\}_{(i,j) \in \Band_N(k)}$,\index{band!$\Band_N(k)$}\index{distinguished place $w_v$} 
determined recursively in $\prec_N$ order,\index{order!$\prec_N$} such that for each $(i,j)$ 
\begin{equation} \label{DFCY2}
|\Delta_{v,ij}^{(k)}|_v \ \le \ B_v \ .
\end{equation}
 
For each $(i,j) \in \Band_N(k)$,\index{band!$\Band_N(k)$} we have $(k-1)N_i \le j < kN_i$. 
Let $r_{ij} = (k+1)N_i -j$, so $N_i < r_{ij} \le 2N_i$ 
and $nN_i - j = (n-k-1)N_i + r_{ij}$, 
and put\index{compensating functions $\vartheta_{v,ij}^{(k)}(z)$!construction of} 
\begin{equation} \label{FTheta}
\vartheta_{v,ij}^{(k)}(z) \ = \  \varepsilon_{v,i}^{-1}  
        \cdot \varphi_{i,r_{ij}}(z) \cdot  \prod_{\ell=1}^{k-1} Q_{\ell}^*(z) 
            \cdot \hP_k(\phi_v(z)) \ .
\end{equation}
Then $\vartheta_{v,ij}^{(k)}(z)$ has a 
pole of order $nN_i-j$ at $x_i$.
\index{compensating functions $\vartheta_{v,ij}^{(k)}(z)$!poles and leading coefficients of}  
Its leading coefficient at $x_i$ is $\tc_{v,i}^{n-k-1}$, 
\index{coefficients $A_{v,ij}$!leading}
because the leading coefficient of $Q_1^*(z)$ is $\varepsilon_{v,i}$, 
\index{coefficients $A_{v,ij}$!leading}
while for $2 \le \ell \le k-1$ the leading
coefficient of $Q_{\ell}^*(z)$ is $\tc_{v,i}^{\ell+1}$.  
\index{coefficients $A_{v,ij}$!leading}
It has a pole of order at most $(n-k-1)N_{i^{\prime}}$ 
at each $x_{i^{\prime}} \ne x_i$.

Since the $\tc_{v,i}$, $\varepsilon_{v,i}$, and $\varphi_{ij}(z)$
are $K_v$-symmetric, and $\hP_k(\phi_v(z))$ is $K_v$-rational, 
\index{$K_v$-symmetric!set of functions}
the $\vartheta_{v,ij}^{(k)}(z)$ are $K_v$-symmetric.
\index{compensating functions $\vartheta_{v,ij}^{(k)}(z)$!are $K_v$-symmetric}  
Thus for each $(i,j)$,
\begin{equation*}
\sum_{x_{i^{\prime}} \in \Aut_c(\CC_v/K_v)(x_i)} 
   \Delta_{v,i^{\prime}j}^{(k)} \vartheta_{v,i^{\prime}j}^{(k)} \ \in \ K_v(\cC) \ .
\end{equation*}
We patch $G_v^{(k-1)}(z)$ by setting\index{compensating functions $\vartheta_{v,ij}^{(k)}(z)$}
\index{patching functions, $G_v^{(k)}(z)$ for $1 \le k \le n$!constructed by patching}
\begin{equation} \label{DFFZ4A}
G_v^{(k)}(z) \ = \ G_v^{(k-1)}(z)
   + \sum_{i=1}^m \sum_{j=(k-1)N_i}^{kN_i-1} \Delta_{v,ij}^{(k)} 
                 \vartheta_{v,ij}^{(k)}(z) \ .
\end{equation}
The coefficients $A_{v,ij}$ 
\index{coefficients $A_{v,ij}$}
with $j < (k-1)N_i$ remain unchanged.
In particular, $G_v^{(k)}(z)$ 
has the same leading coefficients\index{patching functions, $G_v^{(k)}(z)$ for $1 \le k \le n$!leading coefficients of} 
$\varepsilon_{v,i} \tc_{v,i}^n$ as $G_v^{(1)}(z)$.  
\index{coefficients $A_{v,ij}$!leading}
Each $G_v^{(k)}(z)$\index{patching functions, $G_v^{(k)}(z)$ for $1 \le k \le n$!are $K_v$-rational} 
is a $K_v$-rational $(\fX,\vs)$-function 
of degree $nN$, so its roots are $K_v$-symmetric.  
\index{$K_v$-symmetric!set of roots}  

\vskip .1 in
Put
\begin{eqnarray}
\Delta_{v,k}(z) & = &  \sum_{i=1}^m \sum_{j=(k-1)N_i + 1}^{k N_i}
   \Delta_{v,ij}^{(k)} \cdot 
         \varepsilon_{v,i}^{-1}  \varphi_{i,r_{ij}}(z) \ , 
                 \label{DFFZ3}\\
Q_k^*(z)    & = & Q_k(z) + \Delta_{v,k}(z)  \label{DFFZ4} \ .
\end{eqnarray}
Then $\Delta_{v,k}(z)$ and $Q_k^*(z)$ are $K_v$ rational.  Since $Q_k(z)$ has a pole 
of order $(k+1)N_i$ at each $x_i$, with leading coefficient $\tc_{v,i}^{k+1}$,
\index{coefficients $A_{v,ij}$!leading}
while $\Delta_{v,k}(z)$ has a pole of order at most $2N_i$ at $x_i$, 
$Q_k^*(z)$ has the same leading coefficients as $Q_k(z)$. 
\index{coefficients $A_{v,ij}$!leading}

By (\ref{FGFact}) and the definition of 
the $\vartheta_{v,ij}^{(k)}(z)$,\index{compensating functions $\vartheta_{v,ij}^{(k)}(z)$} 
changes in the roots of $G_v^{(k-1)}(z)$
\index{patching functions, $G_v^{(k)}(z)$ for $1 \le k \le n$!for nonarchimedean $K_v$-simple sets!movement of roots} 
are localized to the factor $Q_k(z) = P_k(\phi_v(z))$:\index{patching functions, $G_v^{(k)}(z)$ for $1 \le k \le n$!factorization of}      
\begin{equation}
G_v^{(k)}(z) \ = \ Q_1^*(z) Q_2^*(z) \cdots Q_k^*(z)
                            \cdot \hP_k(\phi_v(z)) \ . \label{DFFZ5}
\end{equation}
Furthermore, by (\ref{FBD}), (\ref{FBM}), and (\ref{DFCY2}) 
\begin{equation*}
\|\Delta_{v,k}(z)\|_{U_v^0} \ \le \ B_v \cdot M_v
                    \ \le \ q^{-\frac{k+1}{q-1} - \log_v(k+1)} \ .
\end{equation*}
In particular, for each ball $B(\theta_h,\rho_h)$,
we have $\|\Delta_{v,k}(z)\|_{B(\theta_h,\rho_h)} < 1$ 
while $\|Q_k(z)\|_{B(\theta_h,\rho_h)} = 1$, 
so $\|Q_k^*(z)\|_{B(\theta_h,\rho_h)} = 1$.

Let $\hsigma_h : D(0,1) \rightarrow B(\theta_h,\rho_h)$
be the chosen $F_{u_h}$-rational parametrization;   
recall that $\theta_{hj} = \hsigma_h(\alpha_{hj})$.
Apply Lemma \ref{DLemCPC4} to $\cQ_{k,h}(Z) = Q_k(\hsigma_h(Z)) = P_k(\hPhi_h(Z)))$ 
and $\Delta_{k,h}(Z) = \Delta_{v,k}(\hsigma_h(Z))$.
The roots of $\cQ_{k,h}(Z)$ in $D(0,1)$ are $\{\alpha_{hj}\}_{j \in \cS_k}$, 
which is a $\psi_v$-regular sequence of length $k+1$ in
\index{regular sequence!$\psi_v$-regular sequence}
$\cO_{u_h}$ in attached to $\cS_k$, and 
$\|\Delta_{k,h}\|_{D(0,1)} \le q^{-(k+1)/(q-1) - \log_v(k+1)}$.
By Lemma \ref{DLemCPC4} the roots
$\{\alpha_{hj}^*\}_{j \in \cS_k}$  of 
$\cQ_{k,h}^*(Z) = \cQ_{k,h}(Z) + \Delta_{k,h}(Z)$
form a $\psi_v$-regular sequence of length $k+1$ in $\cO_{u_h}$, with
\index{regular sequence!$\psi_v$-regular sequence}
\begin{equation} \label{DFCCD1}
\ord_v(\alpha_{hj}^*-\alpha_{hj}) \ > \ \log_v(k+1)
\end{equation}
for each $j \in \cS_k$.

Thus the roots of
$G_v^{(k)}(\hsigma_h(Z))$ in $D(0,1)$ are
\index{patching functions, $G_v^{(k)}(z)$ for $1 \le k \le n$!for nonarchimedean $K_v$-simple sets!movement of roots}
\begin{eqnarray*}
& & \{ \alpha_{h0}^*, \alpha_{h1}^*, \
      \alpha_{h2}^*, \alpha_{h3}^*, \alpha_{h4}^*, \ 
      \ldots, \alpha_{h,j_k+k}^*, \ \ 
        \alpha_{h,j_{k+1}}, \ldots, \alpha_{h,n-1} \} \ ,
\end{eqnarray*}
which is a union of $k$ $\psi_v$-regular sequences of lengths $2, 3, 4, \ldots, k+1$,
\index{regular sequence!$\psi_v$-regular sequence}
together with the remainder of the original $\psi_v$-regular sequence of length $n$.
\index{regular sequence!$\psi_v$-regular sequence}
Transferring this back to $G_v^{(k)}(z)$,
\index{patching functions, $G_v^{(k)}(z)$ for $1 \le k \le n$!for nonarchimedean $K_v$-simple sets!movement of roots}
put $\theta_{hj}^* = \hsigma_h(\alpha_{hj}^*)$ for each $j \in \cS_k$.  
The roots of $G_v^{(k)}(z)$ in $B(\theta_h,\rho_h)$ are
\index{patching functions, $G_v^{(k)}(z)$ for $1 \le k \le n$!for nonarchimedean $K_v$-simple sets!movement of roots} 
\begin{eqnarray*}
& & \{ \theta_{h0}^*, \theta_{h1}^*, \
      \theta_{h2}^*, \theta_{h3}^*, \theta_{h4}^*, \ 
    \ldots, \theta_{h,j_k+k}^*, \ \ 
        \theta_{h,j_{k+1}}, \ldots, \theta_{h,n-1} \} \ , 
\end{eqnarray*}
and they all belong to $\cC_v(F_{u_h}) \cap B(\theta_h,\rho_h)$.
\index{local patching for nonarch $K_v$-simple sets!Phase 1: leading and high-order coefficients!when $\Char(K_v) = 0$|)}   

\smallskip
{\bf Case B.}  Suppose $\Char(K_v) = p > 0$. \
\index{local patching for nonarch $K_v$-simple sets!Phase 1: leading and high-order coefficients!when $\Char(K_v) = p > 0$|(}
\index{patching functions, initial $G_v^{(0)}(z)$!when $\Char(K) = p > 0$!high-order coefficients are $0$|(}

In this case the goals of the high-order patching process are different.
First, we need to choose $n_v$ so that if $n_v|n$, then 
when $G_v^{(0)}(z) = S_{n,v}(\phi_v(z))$ is expanded\index{patching functions, initial $G_v^{(0)}(z)$!expansion of}  
using the $L$-rational basis as
\index{basis!$L$-rational} 
\begin{equation} \label{Gv0ExpA}
G_v^{(0)}(z) \ = \ \sum_{i=1}^m \sum_{j=0}^{(n-1)N_i-1} A_{v,ij} \varphi_{i,nN_i-j}(z) 
           + \sum_{\lambda = 1}^{\Lambda} A_{v,\lambda} \varphi_{\lambda}(z) \ ,
\end{equation} 
we have $A_{v,ij} = 0$ for all $(i,j)$ with $1 \le i \le m$, $1 \le j < \kbar N_i$.
Second, when we patch the leading coefficients $A_{v,i0} = \tc_{v,i}^n$, 
\index{coefficients $A_{v,ij}$!leading}
we must carry out the patching process in such a way that 
the $A_{v,ij}$ with $1 \le i \le m$, $1 \le j < \kbar N_i$ remain $0$.   

To accomplish these goals, we will choose $n_v$ and $B_v$ differently from how they were chosen 
when $\Char(K_v) = 0$.  Recall that $q = q_v = p^{f_v}$ is the order of the residue field of $\cO_v$, 
and that $k_v$ is the least integer for such that for all $k \ge k_v$, 
\begin{equation} \label{FkvBound1}
h_v^{Nk} \cdot M_v \ < \ q^{-(\frac{k+1}{q - 1} + \log_v(k+1))} \ .
\end{equation}
Recall also that $J = p^A$ is the least power of $p$ such that 
\begin{equation*} 
p^A \ \ge \ \max\big(2g+1,\max_{1 \le i \le m}([K(x_i):K]^{\insep})\big) \ ,
\end{equation*} 
(see \S\ref{Chap3}.\ref{LRationalBasisSection})
and that we have chosen $N$ in such a way that $J|N_i$ for each $i$.  
By the construction of the $L$-rational and $L^{\sep}$-rational bases, 
this means that $\varphi_{i,nN_i}(z) = \tphi_{i,nN_i}(z)$
is $L^{\sep}$-rational for each $i =1, \ldots, m$.  

Let $\kbar \ge k_v$ is a fixed integer specified by the global patching process.  
We will take $n_v = q^r$ to be the least power of $q$ such that
\begin{equation} \label{qwCBig}
\left\{ \begin{array}{l}
         q^r \ \ge \ \max(\kbar N_1, \ldots, \kbar N_m) \ , \\
         q^r - q^{r-1} \ \ge \ \kbar \ , 
        \end{array} \right.
\end{equation}
and take   
\begin{equation} \label{F1pBvDef1Q}
B_v \ = \ 
\min\left(\frac{1}{2}\, ,\frac{q^{-\big(\frac{n_v}{q-1} + \log_v(n_v)\big)}}
{ \max_i \big(|\tc_{v,i}^{n_v}|_v \|\varphi_{i, n_v N_i}\|_{U_v} \big)}\right) \ . 
\end{equation}
Note that $J|n_v$, since $n_v \ge \kbar N_1 \ge J$ 
and both $n_v$ and $J$ are powers of $p$.    

We will that for see all sufficiently large $n$ divisible by $n_v$, 
we can carry out the patching process described in Theorem \ref{DCPCPatch1p}.

\smallskip
Given $n$ divisible by $n_v$, 
write $n = n_v Q = q^r \cdot Q$,\index{patching functions, initial $G_v^{(0)}(z)$!construction of}  
and put $G_v^{(0)}(z) = S_{n,v}(\phi_v(z))$.  
We first show when $G_v^{(0)}(z)$ 
is expanded as in (\ref{Gv0ExpA}),\index{patching functions, initial $G_v^{(0)}(z)$!expansion of} 
then the high-order coefficients $A_{v,ij}$ are $0$  for $1 \le i \le m$, 
\index{coefficients $A_{v,ij}$!high-order}
$1 \le j < \kbar N_i$, as required.  For this, we apply Proposition \ref{SnCharpForm},
which says that $S_{n,v}(x)$ has the form  
\begin{equation} \label{FSnTwist}
S_{n,v}(x) \ = \ x^n + \sum_{j = q^r - q^{r-1}}^n C_j x^{n-j} 
\end{equation} 
with each $C_j \in \cO_v$. 

For each $i$, the leading coefficient of $G_v^{(0)}(z)$ at $x_i$ is $\tc_{v,i}^n$, 
\index{coefficients $A_{v,ij}$!leading}\index{patching functions, initial $G_v^{(0)}(z)$!leading coefficients of} 
which belongs to $K_v(x_i)^{\sep}$ since $\tc_{v,i} \in K_v(x_i)^{\sep}$ 
by the construction of $\phi_v(z)$.     
If we expand $\phi_v(z)^Q$ using the $L$-rational basis as 
\index{basis!$L$-rational}
\begin{equation*}
\phi_v(z)^Q  \ = \ \sum_{i=1}^m \sum_{j=0}^{(Q-1)N_i - 1} B_{v,ij} \varphi_{i,QN_i -j}(z) 
                \ + \ \sum_{\lambda = 1}^{\Lambda} B_{\lambda} \varphi_{\lambda}(z) \ ,
\end{equation*} 
then since $\Char(K_v) = p > 0$ and $n_v$ is a power of $p$, it follows that
\begin{equation*}
\phi_v(z)^n  \ = \ (\phi_v(z)^Q)^{n_v} \ = \ 
 \sum_{i=1}^m \sum_{j=0}^{(Q-1)N_i - 1} B_{v,ij}^{n_v} \varphi_{i,QN_i -j}(z)^{n_v}  
                \ + \ \sum_{\lambda = 1}^{\Lambda} B_{\lambda}^{n_v} \varphi_{\lambda}(z)^{n_v} \ .
\end{equation*} 
Here for each $i$, since $J|N_i$, Proposition \ref{TransitionProp}(B) shows that 
$\varphi_{i,Q N_i}^{n_v} = \varphi_{i,nN_i}$, and since the leading coefficient of $\phi_v(z)$ 
\index{coefficients $A_{v,ij}$!leading}
is $\tc_{v,i}$, we have $B_{v,ij}^{n_v} = \tc_{v,i}^n$.  
Each term $B_{v,ij}^{n_v} \varphi_{i,QN_i -j}(z)^{n_v}$ for $j \ge 1$ has a pole of order
at most $nN_i - n_v$ at $x_i$ and no other poles, 
and the $B_{\lambda}^{n_v} \varphi_{\lambda}(z)^{n_v}$ have poles
of order at most $n_v N_{i^\prime}$ at each $x_{i^\prime}$.  
Thus if $n \ge n_v + 1$, which we henceforth assume, then 
\begin{equation} \label{FRump5}
\phi_v(z)^n \ = \ \sum_{i=1}^m \tc_{v,i}^n \varphi_{i,nN_i}(z) 
\ + \ \text{terms with poles of order $\le (n-n_v)N_i$ at each $x_i$\ .} 
\end{equation}
Now consider\index{patching functions, initial $G_v^{(0)}(z)$!expansion of}  
\begin{equation} \label{Gv0Form}
G_v^{(0)}(z) \ = \ S_{n,v}(\phi_v(z)) 
\ = \ \phi_v(z)^n + \sum_{j = q^r - q^{r-1}}^n C_j \phi_v(z)^{n-j} \ . 
\end{equation}
By (\ref{qwCBig}), (\ref{FSnTwist}) and (\ref{FRump5}), 
when $G_v^{(0)}(z)$ is expanded as in (\ref{Gv0ExpA}) 
the coefficients $A_{v,ij}$ for $1 \le i \le m$, $1 \le j < \kbar N_i$  are all $0$.  
\index{band!$\Band_N(k)$}\index{coefficients $A_{v,ij}$!high-order}
\index{patching functions, initial $G_v^{(0)}(z)$!when $\Char(K) = p > 0$!high-order coefficients are $0$|)}

\smallskip
We next carry out the patching process for the stage $k = 1$.   
We want to modify the leading coefficients $A_{v,i0}$ 
\index{coefficients $A_{v,ij}$!leading} 
and leave the remaining-order high coefficients $A_{v,ij}$ for $1 \le j < \kbar N_i$
\index{coefficients $A_{v,ij}$!high-order}
(which are $0$) unchanged.   
By assumption, we are given a $K_v$-symmetric set of numbers 
\index{$K_v$-symmetric!set of numbers}
$\{\Delta_{v,i0}^{(1)}  \in L_{w_v}\}_{1 \le i \le m, 0 \le j < N_i}$,\index{distinguished place $w_v$} 
with $|\Delta_{v,i0}^{(1)}|_v \le B_v$ for each $i$,  
and $\Delta_{v,ij}^{(1)} = 0$ for all $j \ge 1$ and all $i$, 
such that $\Delta_{v,i0}^{(1)}$ belongs to $K_v(x_i)^{\sep}$ for each $i$, 
and we wish to replace $A_{v,i0} = \tc_{v,i}^n$ with $\tc_{v,i}^n + \Delta_{v,i0} \tc_{v,i}^n$ 
in (\ref{Gv0ExpA}). 

We claim that taking
\begin{equation} \label{FzphiQm1p1Q}
\ttheta_{v,i0}^{(1)}(z) \ = \ \tc_{v,i}^{n_v} \varphi_{i,n_vN_i}(z) \cdot S_{n-n_v,v}(\phi_v(z))
\end{equation}
in Theorem \ref{DCPCPatch1p} for each $i=1, \ldots, m$, 
and then putting\index{patching functions, $G_v^{(k)}(z)$ for $1 \le k \le n$!constructed by patching}   
\begin{equation} \label{FLCSum1QQ}
G_v^{(1)}(z) \ = \ G_v^{(0)}(z) \ + \ 
 \sum_{i=1}^m \Delta_{v,i0}^{(1)} \, \ttheta_{v,i0}^{(1)}(z) \ ,
\end{equation} 
accomplishes what we need.  
Let $H(z)$ denote the sum on the right side of (\ref{FLCSum1QQ}). 

First, adding $H(z)$ to $G_v^{(0)}(z)$\index{patching functions, initial $G_v^{(0)}(z)$}  
adds $\Delta_{v,i0}^{(1)} \tc_{v,i}^n$ to $A_{v,i0} = \tc_{v,i}^n$, for each $i$. 
This follows from the fact that 
$\tc_{v,i}^{n_v} \varphi_{i,n_vN_i}(z) \cdot S_{n-n_v,v}(\phi_v(z))$ has a pole of order $nN_i$ at $x_i$ 
with leading coefficient $\tc_{v,i}^n$, and at each $x_{i^{\prime}} \ne x_i$ 
\index{coefficients $A_{v,ij}$!leading}
its pole has order less than $(n-\kbar)N_{i^{\prime}}$.  

Put $\varepsilon_{v,i} = 1 + \Delta_{v,i0}^{(1)}$; 
then the leading coefficient of $G_v^{(1)}(z)$\index{patching functions, $G_v^{(k)}(z)$ for $1 \le k \le n$!leading coefficients of} 
at $x_i$ is $\varepsilon_{v,i} \tc_{v,i}^n$.  
\index{coefficients $A_{v,ij}$!leading} 
Since the $\Delta_{v,i0}^{(1)}$ are $K_v$-symmetric, 
\index{$K_v$-symmetric!set of numbers}
with $\Delta_{v,i0}^{(1)} \in K_v(x_i)^{\sep}$ and $|\Delta_{v,i0}^{(1)}|_v \le B_v < 1$ 
for each $i$, the $\varepsilon_{v,i}$ form a $K_v$-symmetric system of $\cO_{w_v}^{\sep}$-units\index{distinguished place $w_v$} 
\index{$K_v$-symmetric!system of units}
with $\varepsilon_{v,i} \in K_v(x_i)^{\sep}$ for each $i$.   
Since $\tc_{v,i}$ belongs to $K_v(x_i)^{\sep}$, so does $\varepsilon_{v,i} \tc_{v,i}^n$.   

Second, adding $H(z)$ to $G_v^{(0)}(z)$\index{patching functions, initial $G_v^{(0)}(z)$}  
leaves $A_{v,ij} = 0$ for $1 \le i \le m$, $1 \le j < \kbar N_i$.  
To see this, note that if $n = n_vQ$ then $n - n_v = n_v(Q-1)$. 
Since $n_v|(n - n_v)$, Proposition \ref{SnCharpForm} and an argument 
like the one which gave (\ref{Gv0Form}) show that 
\begin{equation} \label{FRump6}
S_{n-n_v,v}(\phi_v(z)) 
\ = \ \phi_v(z)^{n-n_v} + \sum_{j = q^r - q^{r-1}}^{n-n_v} C_j^{\prime} \phi_v(z)^{n-n_v-j} 
\end{equation}
for certain $C_j^{\prime} \in \cO_v$.
Likewise, by an argument similar to the
one which gave (\ref{FRump5}), if $n \ge n_v + 2$ (which we henceforth assume), then  
\begin{eqnarray} 
\phi_v(z)^{n-n_v} & = & \sum_{i=1}^m \tc_{v,i}^{n-n_v} \varphi_{i,(n-n_v)N_i}(z) \label{FRump7} \\
& & \ + \ \text{terms with poles of order $\le (n-2n_v)N_i$ at each $x_i$\ .} \notag
\end{eqnarray}
Finally, since $J|N_i$ for each $i$, Proposition \ref{TransitionProp}(B) shows that for
each basis function $\varphi_{ij}(z)$ we have 
\index{basis!$L$-rational}
\begin{equation}  \label{FRump8}
\varphi_{i,n_v N_i}(z) \cdot \varphi_{ij}(z) \ = \ \varphi_{i,n_vN_i + j}(z) \ .
\end{equation} 
Combining (\ref{FRump6}), (\ref{FRump7}) and (\ref{FRump8}), 
and using that $q^r-q^{r-1} \ge \kbar$, we see that 
\begin{eqnarray}
\ttheta_{v,i0}^{(1)}(z) & = & \tc_{v,i}^{n_v} \varphi_{i,n_v N_i}(z) \cdot S_{n-n_v,v}(\phi_v(z)) \notag \\
                        & = & \tc_{v,i}^n \varphi_{i,nN_i}(z) 
           \ + \ \text{terms with poles of order $\le (n-\kbar)N_{i^{\prime}}$ at each $x_i^{\prime}$\ .} 
           \label{FphiQm1p1Q}
\end{eqnarray} 
Thus $\ttheta_{v,i0}^{(1)}(z) = \tc_{v,i}^n \varphi_{i,nN_i}(z) + \tTheta_{v,i0}^{(1)}(z)$
for an $(\fX,\vs)$-function $\tTheta_{v,i0}^{(1)}(z)$ with poles of order at most $(n-\kbar)N_{i^\prime}$ 
at each $x_{i^\prime}$, as asserted in Theorem \ref{DCPCPatch1p}.  

Third, $\ttheta_{v,i0}^{(1)}(z)$ is rational over $K_v(x_i)^{\sep}$ for each $i$,
and the  $\ttheta_{v,i0}^{(1)}(z)$ are $K_v$-symmetric. 
\index{$K_v$-symmetric!set of functions}
To see this, note that $\varphi_{i,n_v N_i}(z) = \tphi_{i,n_vN_i}(z)$ is rational over $K_v(x_i)^{\sep}$ 
\index{$K_v$-symmetric!set of functions}
by Proposition \ref{TransitionProp}(B), since $J|N_i$.  
Furthermore, $\tc_{v,i} \in K_v(x_i)^{\sep}$ by hypothesis, 
and $S_{n-n_v,v}(\phi_v(z))$ is $K_v$-rational, 
so $\ttheta_{v,i0}^{(1)}(z)$ is $K_v(x_i)^{\sep}$-rational.   
The $\ttheta_{v,i0}^{(1)}(z)$ 
are $K_v$-symmetric since the $\tc_{v,i}$ and $\varphi_{i,n_v N_i}(z)$ are $K_v$-symmetric
\index{$K_v$-symmetric!set of functions}
and $S_{n-n_v,v}(\phi_v(z))$ is $K_v$-rational.  

Fourth, $G_v^{(1)}(z)$ is $K_v$-rational.\index{patching functions, $G_v^{(k)}(z)$ for $1 \le k \le n$!are $K_v$-rational} 
To see this, note that $\tc_{v,i}^{n_v}$ belongs to $K_v(x_i)^{\sep}$,  
and the $\ttheta_{v,i0}^{(1)}(z)$ are rational over $L_{w_v}$-rational\index{distinguished place $w_v$} and $K_v$-symmetric
\index{$K_v$-symmetric!set of functions}
as remarked above.  It follows that $H(z)$ is $K_v$-rational.  
Since $G_v^{(0)}(z)$ is $K_v$-rational,\index{patching functions, initial $G_v^{(0)}(z)$}  
$G_v^{(1)}(z)$ is $K_v$-rational as well.\index{patching functions, $G_v^{(k)}(z)$ for $1 \le k \le n$!are $K_v$-rational} 

\smallskip
We will now show that the roots of $G_v^{(1)}(z)$ belong to $E_v$.
\index{patching functions, $G_v^{(k)}(z)$ for $1 \le k \le n$!roots are confined to $E_v$}  
Let 
\begin{equation*} 
\cS_1 \ = \ \{n-n_v, n-n_v + 1, \ldots, n-1\} \ , 
\end{equation*} 
and put   
$P_1(x) = \prod_{j \in \cS_1}(x - \psi_v(j))$, 
$\hP_1(x) = S_{n-n_v,v}(x) = \prod_{j=0}^{n-n_v-1} (x - \psi_v(j))$,   
so  
\begin{equation*}
S_{n,v}(x) \ = \ P_1(x) \cdot S_{n-n_v,v}(x) \ = \ P_1(x) \cdot \hP_1(x) \ .
\end{equation*} 
Since $G_v^{(0)}(z) = S_{n,v}(\phi_v(z))$,\index{patching functions, $G_v^{(k)}(z)$ for $1 \le k \le n$!expansion of}  
using (\ref{FzphiQm1p1Q}) we can write (\ref{FLCSum1QQ}) as 
\begin{eqnarray*} 
G_v^{(1)}(z) & = & G_v^{(0)}(z) \ + \ 
\Big(\sum_{i=1}^m \Delta_{v,i0}^{(1)} \cdot 
       \tc_{v,i}^{n_v} \varphi_{i,n_v N_i}(z) \Big) \cdot S_{n-n_v,v}(\phi_v(z)) \\ 
& = & \Big(P_1(\phi_v(z)) + \Delta_{v,1}(z)\Big) \cdot \hP_1(\phi_v(z)) \\
& = & Q_1^*(z) \cdot  \hP_1(\phi_v(z)) \ ,
\end{eqnarray*} 
where $\Delta_{v,1}(z) = \sum_{i=1}^m \Delta_{v,i0}^{(1)} \cdot \tc_{v,i}^{n_v} \varphi_{i,n_v N_i}(z)$
and $Q_1^*(z) = P_1(\phi_v(z)) + \Delta_{v,1}(z)$.  
Since $|\Delta_{v,i0}^{(1)}|_v \le B_v$ for each $i$, 
by (\ref{F1pBvDef1Q}) and the ultrametric inequality we have 
\begin{equation*}
\|\Delta_{v,1}(z)\|_{U_v} \ \le \ q^{-\frac{n_v}{q-1} - \log_v(n_v)} \ .
\end{equation*}
In particular, for each ball $B(\theta_h,\rho_h)$,
we have $\|\Delta_{v,1}(z)\|_{B(\theta_h,\rho_h)} \le q^{-n_v/(q-1) - \log_v(n_v)} < 1$. 
Since $\|P_1(\phi_v(z))\|_{B(\theta_h,\rho_h)} = 1$, it follows that 
$\|Q_1^*(z)\|_{B(\theta_h,\rho_h)} = 1$.  

For each $h = 1, \ldots, N$, let $\hsigma_h : D(0,1) \rightarrow B(\theta_h,\rho_h)$
be the $F_{u_h}$-rational parametrization chosen at the beginning of 
\S\ref{Chap11}.\ref{PatchingLemmasSection}; recall that
the roots of $S_{n,v}(\phi_v(z))$ in $B(\theta_h,\rho_h)$ are $\theta_{hj}$ for $j = 0, \ldots, n-1$
and that $\theta_{hj} = \hsigma_h(\alpha_{hj})$ where $\alpha_{hj} = \tPhi_h(\psi_v(j)) \in \cO_{u_h}$. 
 
Just as when $\Char(K_v) =0$, we can apply Lemma \ref{DLemCPC4} to $\cQ_{1,h}(Z) = Q_1(\hsigma_h(Z))$ 
and $\Delta_{1,h}(Z) = \Delta_{v,1}(\hsigma_h(Z))$.
The roots of $\cQ_{1,h}(Z)$ in $D(0,1)$ are $\{\alpha_{hj}\}_{j \in \cS_1}$, 
which is a $\psi_v$-regular sequence of length $n_v$ in
\index{regular sequence!$\psi_v$-regular sequence}
$\cO_{u_h}$ in attached to $\cS_1$.
By Lemma \ref{DLemCPC4} the roots
$\{\alpha_{hj}^*\}_{j \in \cS_1}$  of 
$\cQ_{1,h}^*(Z) = \cQ_{1,h}(Z) + \Delta_{1,h}(Z)$
form a $\psi_v$-regular sequence of length $n_v$ in $\cO_{u_h}$, with
\index{regular sequence!$\psi_v$-regular sequence}
\begin{equation} \label{DFCCD1b}
\ord_v(\alpha_{hj}^*-\alpha_{hj}) \ > \ \log_v(n_v)
\end{equation}
for each $j \in \cS_1$.

Thus the roots of
$G_v^{(1)}(\hsigma_h(Z))$ in $D(0,1)$ are
\index{patching functions, $G_v^{(k)}(z)$ for $1 \le k \le n$!for nonarchimedean $K_v$-simple sets!movement of roots}
\begin{equation*}
 \{ \alpha_{h0}, \alpha_{h1}, \ldots, \alpha_{h,n-n_v-1}, \ 
           \alpha_{h,n-n_v}^{*}, \ldots, \alpha_{h,n-1}^* \} \ ,
\end{equation*}
which is the union of 
the initial part of the original $\psi_v$-regular sequence of length $n$
\index{regular sequence!$\psi_v$-regular sequence}
and a $\psi_v$-regular sequence of length $n_v$.  
Transferring this back to $G_v^{(k)}(z)$,
\index{patching functions, $G_v^{(k)}(z)$ for $1 \le k \le n$!for nonarchimedean $K_v$-simple sets!movement of roots}
put $\theta_{hj}^* = \hsigma_h(\alpha_{hj}^*)$ for each $j \in \cS_1$.  
Then the roots of $G_v^{(1)}(z)$ in $B(\theta_h,\rho_h)$ are
\index{patching functions, $G_v^{(k)}(z)$ for $1 \le k \le n$!for nonarchimedean $K_v$-simple sets!movement of roots} 
\begin{equation*}
 \{ \theta_{h0}, \theta_{h1}, \ldots, \theta_{h,n-n_v-1}, \ 
           \theta_{h,n-n_v}^{*}, \ldots, \theta_{h,n-1}^* \} \ ,
\end{equation*}
and they all belong to $\cC_v(F_{u_h}) \cap B(\theta_h,\rho_h) \subseteq E_v$.  
 
\smallskip
For $k = 2, \ldots, \kbar$, we have $\Delta_{v,ij}^{(k)} = 0$ for all $(i,j)$,  
and we take $G_v^{(k)}(z) = G_v^{(k-1)}(z)$.\index{patching functions, $G_v^{(k)}(z)$ for $1 \le k \le n$!constructed by patching}  
For notational compatibility with Case A,
for each $k =2, \ldots, \kbar$ let $\cS_k$ be the empty set and put $Q_k^*(z) = 1$, 
$\hP_k(z) = \hP_1(z)$. Thus for each $k$\index{patching functions, $G_v^{(k)}(z)$ for $1 \le k \le n$!factorization of}  
\begin{equation*} 
G_v^{(k)}(z) \ = \ Q_1^*(z)  Q_2^*(z) \cdots Q_{k}^*(z) \cdot \hP_k(\phi_v(z)) \ .
\end{equation*} 
Note that the leading coefficient of $Q_1^*(z)$ at $x_i$ is $\varepsilon_{v,i} \tc_{v,i}^{n_v}$.  
\index{coefficients $A_{v,ij}$!leading}
\index{local patching for nonarch $K_v$-simple sets!Phase 1: leading and high-order coefficients!when $\Char(K_v) = p > 0$|)}

\vskip .1 in
\noindent{\bf Phase 2.  Patching for  $k=\kbar+1, \ldots, k_1$.}
\index{local patching for nonarch $K_v$-simple sets!Phase 2: carry on|(}  

Since $\kbar \ge k_v$, for each $k \ge \kbar+1$ we have 
\begin{equation*} 
h_v^{kN}\cdot M_v \ < \ q^{-\frac{k+1}{q-1} - \log_v(k+1)} \ .
\end{equation*} 
Let $k_1$ be the least integer such that for all $k \ge k_1$,
\begin{equation} \label{FPatchBoundE}
h_v^{kN}\cdot M_v \ < \ q^{-\frac{k+1}{q-1} - \log_v(n)} \ .
\end{equation}
If $n$ is large enough that 
\begin{equation} \label{FMAN} 
q^{-\frac{\kbar+1}{q-1} - \log_v(n)} \ < \ h_v^{\kbar N} \cdot M_v \ ,
\end{equation}
which we henceforth assume, then $k_1 > \kbar$.  
Since $h_v^N < q^{-1/(q-1)}$, 
there is a constant $A_1 > 0$ such that 
\begin{equation*}
k_1 \ \le \ A_1 \log_v(n) \ ,
\end{equation*}
so if $n$ is sufficiently large, then $n > k_1$ which we also assume.  

The purpose of Phase 2 is simply to ``carry on'' until $k$ is large enough that
(\ref{FPatchBoundE}) holds, at which point Lemma \ref{DLemCPC4} 
preserves the position of the roots $\alpha_{hj}$ 
within balls of size $q_v^{-\lceil \log_v(n) \rceil}$.  
Equivalently, in patching steps for $k > k_1$, 
each pair of patched and unpatched roots will satisfy  
\begin{equation*}
\ord_v(\alpha_{hj}^*-\alpha_{hj}) \ \ge \ \log_v(n) \ .
\end{equation*}
The patching process in Phase 2 is the same
as the one for steps  $k \ge 2$ in Case A of Phase 1, 
except that in place of (\ref{DFCY2}) we require that
each $\Delta_{v,ij}^{(k)}$ satisfy
\begin{equation} \label{DFCY3}
|\Delta_{v,ij}^{(k)}|_v \ \le \ h_v^{kN} .
\end{equation}
When $\Char(K_v) = 0$, for each $k = \kbar+1, \ldots, k_1$ we put $\cS_k = \{j_k, \ldots, j_k+k\}$
as in (\ref{DFFT1}).
When $\Char(K_v) = p > 0$, we put $\cS_{\kbar + 1} = \{0,1, \ldots, \kbar\}$   
and for $k = \kbar + 2, \ldots, k_1$ we let $\cS_k$ 
be the set consisting of the next $k+1$ integers after $\cS_{k-1}$.
These $\cS_k$ will be subsequences of $\{0,1, \ldots, n-1\}$ if $n$ is large enough that
\begin{equation} \label{DFXC2}
n \ > \ n_v + ( A_1 \log_v(n) + 1)(\ A_1 \log_v(n) + 2)/2 \ , 
\end{equation}
which we henceforth assume.

Let the $P_k(z)$, $\hP_k(z)$, $Q_k(z)$, and $Q_k^*(z)$ 
be as in (\ref{DFFQ1}), (\ref{DFFR2}), and (\ref{DFFZ4}).
Inductively suppose $G_v^{(k-1)}(z)$\index{patching functions, $G_v^{(k)}(z)$ for $1 \le k \le n$!factorization of} 
has been constructed with   
\begin{eqnarray}
G_v^{(k-1)}(z) & = & Q_1^*(z) \cdots Q_{k-1}^*(z) \cdot \hP_{k-1}(\phi_v(z)) 
                                   \notag \\
      & = & Q_1^*(z) \cdots Q_{k-1}^*(z) \cdot  P_k(\phi_v(z)) \cdot \hP_k(\phi_v(z))  
                                   \label{FGFactB} 
\end{eqnarray}  
When $\Char(K_v) = 0$  
we patch $G_v^{(k-1)}(z)$\index{patching functions, $G_v^{(k)}(z)$ for $1 \le k \le n$!constructed by patching} 
by setting\index{compensating functions $\vartheta_{v,ij}^{(k)}(z)$} 
\begin{eqnarray*} 
G_v^{(k)}(z) & = & G_v^{(k-1)}(z) + \sum_{i=1}^m \sum_{j=(k-1)N_i}^{kN_i-1} 
\Delta_{v,ij}^{(k)} \vartheta_{v,ij}^{(k)}(z) \\
& = & Q_1^*(z) \cdots Q_{k-1}^*(z) \cdot \big(P_k(\phi_v(z)) + \Delta_{v,k}(z) \big) \cdot \hP_k(\phi_v(z)) \ ,
\end{eqnarray*}
where the compensating functions are\index{compensating functions $\vartheta_{v,ij}^{(k)}(z)$!construction of}
\begin{equation*}
\vartheta_{v,ij}^{(k)}(z) \ = \ \varepsilon_{v,i}^{-1} \varphi_{i,(k+1)N_i - j}(z) \cdot 
Q_1^*(z) \cdots Q_{k-1}^*(z) \cdot \hP_k(\phi_v(z)) \ , 
\end{equation*}  
and where as in (\ref{DFFZ3}),  
\begin{equation*}
\Delta_{v,k}(z) \ = \ \sum_{i=1}^m \sum_{j=(k-1)N_i}^{kN_i-1} 
   \Delta_{v,ij}^{(k)} \cdot \varepsilon_{v,i}^{-1} \varphi_{i,(k+1)N_i-j}(z) \ .
\end{equation*}
It is easy to see that the $\vartheta_{v,ij}^{(k)}(z)$
\index{compensating functions $\vartheta_{v,ij}^{(k)}(z)$!poles and leading coefficients of}  
have the properties required by
Theorem \ref{DCPCPatch}. Including the factor of $\varepsilon_{v,i}^{-1}$ 
makes the leading coefficient of $\vartheta_{v,ij}^{(k)}(z)$
\index{coefficients $A_{v,ij}$!leading}
at $x_i$ be $\tc_{v,i}^{n-k-1}$. 

When $\Char(K_v) = p > 0$, we put 
\begin{equation*}
\Delta_{v,k}(z) \ = \
 \sum_{i=1}^m \sum_{j=(k-1)N_i}^{kN_i-1} \Delta_{v,ij}^{(k)} \cdot \varphi_{i,(k+1)N_i-j}(z) \ ,
\end{equation*} 
which is $K_v$-rational by the hypotheses of Theorem \ref{DCPCPatch1p},
and let 
\begin{equation*} 
F_{v,k}(z) \ = \ Q_1^*(z) \cdots Q_{k-1}^*(z) \cdot \hP_k(\phi_v(z)) \ .
\end{equation*}
 We patch $G_v^{(k-1)}(z)$\index{patching functions, $G_v^{(k)}(z)$ for $1 \le k \le n$!constructed by patching} by setting 
\begin{eqnarray*} 
G_v^{(k)}(z) & = & G_v^{(k-1)}(z) + \Delta_{v,k}(z) \cdot F_{v,k}(z) \\
& = &  Q_1^*(z) \cdots Q_{k-1}^*(z) \cdot \big(P_k(\phi_v(z)) + 
\Delta_{v,k}(z) \big) \cdot  \hP_k(\phi_v(z)) 
\end{eqnarray*}
By induction $F_{v,k}(z)$ is $K_v$-rational, its roots belong in $E_v$, 
and it has a pole of order $(n-k-1)N_i$
with leading coefficient $d_{v,i} = \varepsilon_{v,i} \tc_{v,i}^{n-k-1}$ at each $x_i$.
\index{coefficients $A_{v,ij}$!leading}  
In particular $|d_{v,i}|_v = |\tc_{v,i}^{n-k-1}|_v$.
Thus $F_{v,k}(z)$ satisfies the conditions of Theorem \ref{DCPCPatch1p}.

The two patching constructions are the same except for a 
minor difference in the choice of $\Delta_{v,k}(z)$, 
and in both cases 
\begin{equation*} 
\|\Delta_{v,k}\|_{U_v^0} \ \le \ 
\max_{i,j}(\Delta_{|v,ij}|_v) \cdot \max_{i,j} (\|\varphi_{(k+1)N_i - j}\|_{U_v^0}) 
\ \le \ h_v^{kN} \cdot M_v \ < \ q^{-\frac{k}{q-1} - \log_v(k+1)} \ .
\end{equation*} 
The same argument as in Phase 1, using Lemma \ref{DLemCPC4},
shows that $Q_k^*(z) = P_k(\phi_v(z)) + \Delta_{v,k}(z)$ is $K_v$-rational and 
has its roots in $E_v$.  For each $h = 1, \ldots, N$, 
the roots of $Q_k^*(z)$ in $B(\theta_h,\rho_h)$  
are $\theta_{hj}^* = \hsigma_h(\alpha_{hj}^*)$, for $j \in \cS_k$.
Since for each $h$ 
\begin{equation*}
\|\Delta_{v,k}\|_{B(\theta,\rho_h)} \ < \ q^{-\frac{k+1}{q-1} - \log_v(k+1)} \ ,
\end{equation*}
the $\alpha_{hj}^*$ form a $\psi_v$-regular sequence of length $k+1$ in $\cO_{u_h}$,
\index{regular sequence!$\psi_v$-regular sequence}  
and for each $j \in \cS_k$ 
\begin{equation} \label{DFCCD2}
\ord_v(\alpha_{hj}^*-\alpha_{hj}) \ > \ \log_v(k+1) \ .
\end{equation} 
After patching, we have a factorization\index{patching functions, $G_v^{(k)}(z)$ for $1 \le k \le n$!factorization of}
\begin{equation*}
G_v^{(k)}(z) \ = \ Q_1^*(z) \cdots Q_k^*(z) \cdot \hP_k(\phi_v(z)) \ ,
\end{equation*}
so the induction in (\ref{FGFactB}) can continue.
\index{local patching for nonarch $K_v$-simple sets!Phase 2: carry on|)}   

\vskip .1 in
\noindent{\bf Phase 3.  Moving the roots apart.}
\index{move roots apart}\index{separate the roots} 
\index{local patching for nonarch $K_v$-simple sets!Phase 3: move roots apart|(} 

At this point we have obtained a $K_v$-rational 
$(\fX,\vs)$-function\index{patching functions, $G_v^{(k)}(z)$ for $1 \le k \le n$!are $K_v$-rational}
$G_v^{(k_1)}(z)$ of degree $nN$ whose coefficients $A_{v,ij}$ with 
$1 \le i \le m$, $0 \le j < k_1 N_i$ have been patched, 
\index{coefficients $A_{v,ij}$!high-order}
and which has the factorization\index{patching functions, $G_v^{(k)}(z)$ for $1 \le k \le n$!factorization of}
\begin{equation*}
G_v^{(k_1)}(z) \ = \
   Q_1^*(z) \cdots  Q_{k_1}^*(z) \cdot \hP_{k_1}(\phi_v(z)) \ .
\end{equation*}
Its zeros belong to $E_v$, 
and it has $n$ roots in each ball $B(\theta_h,\rho_h)$:  when $\Char(K_v) = 0$, these are
\begin{equation*}
 \{ \theta_{h0}^*, \theta_{h1}^*, \
      \theta_{h2}^*, \theta_{h3}^*, \theta_{h4}^*, \ \ldots, \
        \theta_{h,j_{k_1}+k_1}^*, \
          \theta_{h,j_{k_1}+k_1+1}, \ldots, \theta_{h,n-1}\} \ .
\end{equation*}
When the zeros in $B(\theta_h,\rho_h)$ are pulled back to $D(0,1)$ using $\hsigma_h$, they become 
\begin{equation*}
 \{ \alpha_{h0}^*, \alpha_{h1}^*, \
      \alpha_{h2}^*, \alpha_{h3}^*, \alpha_{h4}^*, \ \ldots, \
        \alpha_{h,j_{k_1}+k_1}^*, \
          \alpha_{h,j_{k_1}+k_1+1}, \ldots, \alpha_{h,n-1}\} \ ,
\end{equation*}
a union of the $\psi_v$-regular sequences in $\cO_{u_h}$ corresponding to
\index{regular sequence!$\psi_v$-regular sequence} 
$\cS_1, \ldots, \cS_{k_1}$, together with the unpatched roots.   
For each $k = 1, \ldots, k_1$
the ``patched'' roots $\{\alpha_{hj}^*\}_{j \in \cS_k}$ satisfy 
\begin{equation*}
\ord_v(\alpha_{hj}^* - \alpha_{hj}) \ > \ \log_v(\#(\cS_k))
\end{equation*}
relative to the original ``unpatched'' roots $\alpha_{hj}$ for $j \in \cS_k$.
\index{regular sequence!$\psi_v$-regular sequence}

\vskip .1 in
Although the patched roots corresponding to $\cS_{k}$ 
\index{patched roots} 
are well-separated from each other, 
they may have come close to (or even coincide with) roots corresponding to
some other $\cS_{\ell}$ or to unpatched roots.
\index{unpatched roots} 
The purpose of Phase 3 is to move the patched roots 
$\alpha_{hj}^*$ to new points $\alpha_{hj}^{**}$ which are well-separated from each other and
\index{well-separated}\index{separate the roots}
the unpatched roots, while preserving the coefficients of 
\index{coefficients $A_{v,ij}$!high-order}
$G_v^{(k_1)}(z)$ that have already been patched.\index{patching functions, $G_v^{(k)}(z)$}

We accomplish this in two steps.
First, we move the patched roots away from any roots they have come  
too near to.  This changes the high order coefficients.  
\index{coefficients $A_{v,ij}$!high-order} 
\index{coefficients $A_{v,ij}$!restoring} 
Second, we restore the high order coefficients to their
original values, by re-patching\index{re-patch} using a new sequence of roots. 
We allow these new roots only to move in such a way that they stay 
well-separated from the other roots, 
which in turn limits how far the 
original $\alpha_{hj}^*$ can be moved.  
A computation shows that for a suitable constant $A_4$, 
the roots $\alpha_{hj}^{**}$ of the re-patched\index{re-patch} function $\hG_v^{(k_1)}(z)$
can be required to satisfy $|\alpha_{hj}^{**}-\alpha_{h \ell}^{**}|_v \ge n^{-A_4}$  
(or equivalently $\ord_v(\alpha_{hj}^{**}-\alpha_{h \ell}^{**}) \le A_4 \log_v(n)$)
for all $h$ and all $j \ne \ell$.  
We then replace $G_v^{(k_1)}(z)$ 
 with $\hG_v^{(k_1)}(z)$.\index{patching functions, $G_v^{(k)}(z)$ for $1 \le k \le n$!constructed by patching} 

\vskip .1 in
There are several obstacles to carrying this out.  
Moving the roots $\alpha_{hj}^*$ generally results in
a non-principal divisor.  To compensate, we choose a collection of roots
in ``good position'' which are well-separated from the
patched roots, and move them to regain a principal divisor.  (It is here
that we use the assumption that $H_v$ is `move-prepared'.)
\index{move-prepared} 
In Lemma \ref{MovingLemma1} below, 
we construct a function $Y(z) \in K_v(\cC)$ whose
zeros are the moved roots, whose poles are the original roots, 
and which is very close to $1$ outside $U_v$. 
Multiplying $G_v^{(k_1)}(z)$ by $Y(z)$\index{patching functions, $G_v^{(k)}(z)$}
yields a function $\Gbar_v^{(k_1)}(z)$ 
with the desired new roots.
By standard estimates for Laurent coefficients,\index{Laurent expansion} 
the amount the high-order coefficients are changed in passing from 
\index{coefficients $A_{v,ij}$!high-order} 
$G_v^{(k_1)}(z)$ to $\Gbar_v^{(k_1)}(z)$
depends on how close $Y(z)$ is to $1$ at the points $x_i$, 
which in turn depends on how far the roots were moved.

To restore the high-order coefficients to their previous values, 
\index{coefficients $A_{v,ij}$!high-order} 
\index{coefficients $A_{v,ij}$!restoring}
we apply the basic patching lemma but use a carefully chosen, 
previously unpatched sequence of roots to absorb the resulting movement 
in the roots.  We prove an estimate giving the ``cost'' 
of independently adjusting each coefficient.  This 
tells us how far the $\alpha_{hj}^*$ can be moved. 

\vskip .1 in
We now give the details of the construction, postponing the proofs of 
the three Moving Lemmas \ref{MovingLemma1}, \ref{MovingLemma2}, 
and \ref{MovingLemma3} to \S\ref{Chap11}.\ref{MovingLemmaSection}.  Put
\begin{equation*}
\delta_n = q_v^{-\lceil \log_v(n) \rceil} \ .
\end{equation*}
For each $\alpha \in D(0,1)$, we will call the disc $D(\alpha,\delta_n)$ 
the ``$\delta_n$-coset'' of $\alpha$,\index{$\delta_n$-coset|ii}
and  we will refer to elements of $\{0,1,\ldots, n-1\}$ as ``indices''.
\index{indices|ii}
 
Let  $\cS^{\triangle} = \cS_1 \cup \cS_2 \cup \cdots \cup \cS_{k_1}$ 
be the collection of indices which have already been used in the patching process; 
we will call it the set of ``patched'' indices.  We will call its complement 
$\cS^{\Diamond} = \{0,1, \ldots, n-1\} \backslash \cS^{\triangle}$ the set of 
``unpatched'' indices.\index{indices!unpatched|ii} 
Let $\cS^{\dagger} \subset \cS^{\Diamond}$ be the collection of unpatched indices
which Steps 1 and 2 of the patching process have ``endangered'':
the set of indices $\ell \in \cS^{\Diamond}$ 
for which some patched root\index{roots!patched|ii} $\alpha_{hj}^*$ has moved too close to an   
unpatched root $\alpha_{h \ell}$:\index{roots!unpatched|ii} 
\begin{eqnarray*}
\cS^{\dagger} & = &
  \{ \ell \in \cS^{\Diamond}
       : \ord_v(\alpha_{hj}^*-\alpha_{h \ell}) \ge \log_v(n) 
                      \text{\ for some $j \in \cS_1 \cup \cS_2 \cup \cdots \cup \cS_{k_1}$} \\
 & & \qquad \qquad \qquad  \text{and some $h$, $1 \le h \le N$} \} \ .
\end{eqnarray*}
This set of indices must be ``protected'' until later in the patching process;
we will call the corresponding roots ``endangered''.\index{roots!endangered|ii} 
Finally, put
\begin{equation*}
\cS^{\heartsuit} \ = \ \cS^{\Diamond} \backslash \cS^{\dagger} \ .
\end{equation*}
This is the set of indices which are ``safe'' to use in re-patching:\index{re-patch}\index{indices!safe|ii}  
for each $\ell \in \cS^{\heartsuit}$, there is no $h$ for which
any $\alpha_{hj}^*$ belongs to the $\delta_n$-coset\index{$\delta_n$-coset} of
$\alpha_{h \ell}$.\index{roots!safe|ii} 

We have 
$\#(\cS^{\triangle}) \le A_2 \cdot (\log_v(n))^2$  
for an appropriate constant $A_2$.  
Since each patched root\index{roots!patched} $\alpha_{hj}^*$ can belong to the
$\delta_n$-coset of at\index{$\delta_n$-coset} most one unpatched root\index{roots!unpatched} $\alpha_{h \ell}$,
it follows that 
\begin{equation} \label{DFLIP1}
\#(\cS^{\triangle} \cup \cS^{\dagger}) 
        \ \le \ (N+1) A_2 (\log_v(n))^2 \ .
\end{equation}  
We view $\cS^{\triangle} \cup \cS^{\dagger}$ 
as a collection of marked indices which partitions its complement in 
$\{0,1, \ldots, n-1 \}$ (the set of ``safe'' indices $\cS^{\heartsuit}$)\index{indices!safe} 
into a collection of sequences of consecutive integers.  
Let $\cS^0$ be the longest such sequence 
(to be specific, the first one, if there are two or more of the same length). 
The Pigeon-hole Principle\index{Pigeon-hole Principle} shows that by taking $n$ sufficiently large,
we can make $S^0$ arbitrarily long.  If 
\begin{equation} \label{DFHWC1}
n \ \ge \ ((N+1) A_2 (\log_v(n))^2 + 1) \cdot (A_1 \log_v(n) + 2) \ , 
\end{equation}
which we will henceforth assume, then    
\begin{equation} \label{FK10} 
\#(\cS^0) \ \ge \ A_1 \log_v(n) + 2 \ \ge \ k_1 + 2 \ .
\end{equation}
 
We will call $S^0$ the ``long safe sequence''.\index{long safe sequence|ii} Write $\cS^0 = \{j_0, j_1, \ldots, j_L\}$.  
The index $j_0$ will be used to provide the roots\index{roots!in good position} in ``good position''
needed to recover principal divisors when the $\alpha_{hj}^*$ are moved in
the first step of the process.  Let 
\begin{equation*}
\cS^0[k_1] \ = \ \{j_1, \ldots, j_{k_1+1}\} 
\end{equation*}
be the subsequence of $\cS^0$ consisting of the next $k_1 + 1$ 
integers.  This is the sequence of indices that will be used for ``re-patching''\index{re-patch}  
in the second step.  

Our first lemma describes the properties of the divisor $\cD$ and the function $Y(z)$. 
It relates $\varepsilon$, the distance we can move the roots,\index{roots!move roots}   
to  $C_2 \varepsilon$, a bound for the size of $|Y(z)-1|_v$.  
 

\begin{lemma} \label{MovingLemma1} {\bf (First Moving Lemma)}
\index{First Moving Lemma|ii}

Let $\delta_n = q_v^{-\lceil \log_v(n) \rceil}$,  $\cS^{\triangle}$, $j_0$, and 
$\{\theta_{hj}^* \in \cC_v(F_{u_h})\}_{1 \le h \le N, j \in S^{\triangle}}$ be as above.  
Then there are an $\varepsilon_1 > 0$ and constants $C_1, C_2 \ge 1$  
$($depending on $\phi_v(z)$, $E_v$, $H_v$, and their $K_v$-simple decompositions, 
\index{$K_v$-simple!set}
\index{$K_v$-simple!decomposition}
but not on $n)$, with the following property$:$ 

For each $0 < \varepsilon < \varepsilon_1$ such that 
\begin{equation} \label{FEPSCOND1}
C_1  \varepsilon \ \le \ \delta_n \cdot \min_{1 \le h \le N}(\rho_h) \ ,
\end{equation}  
given any $K_v$-symmetric set 
\index{$K_v$-symmetric!set of numbers}
$\{\theta_{hj}^{**} \in \cC_v(F_{u_h}) \cap B(\theta_h,\rho_h) \}_{1 \le h \le N, j \in S^{\triangle}}$
with $\|\theta_{hj}^{**},\theta_{hj}^*\|_v < \varepsilon$ for all $(h,j)$,  
there is a $K_v$-symmetric collection of points 
\index{$K_v$-symmetric!set of points}
$\{\theta_{h,j_0}^{**} \in \cC_v(F_{u_h}) \cap B(\theta_h,\rho_h) \}_{1 \le h \le N}$ 
satisfying  
\begin{equation*}  
\|\theta_{h,j_0}^{**}, \theta_{h,j_0}\|_v \ \le \ C_1 \varepsilon \ \le \ \delta_n \rho_h 
\end{equation*}
for each $h$, such that 

$(A)$ The divisor 
\begin{equation*} 
\cD \ = \ \sum_{j \in S^{\triangle}} 
                \sum_{h=1}^N ((\theta_{hj}^{**}) - (\theta_{hj}^*))
        + \sum_{h=1}^N ((\theta_{h,j_0}^{**})-(\theta_{h,j_0}))
\end{equation*}
is $K_v$-rational and principal;

$(B)$ Writing $U_v^0 = \bigcup_{h=1}^N B(\theta_h,\rho_h) \subset U_v$ as before, we have  
\begin{enumerate}
\item $|Y(z)|_v = 1$ for all $z \in \cC_v(\CC_v) \backslash U_v^0$;
\item $|Y(z)-1|_v \le C_2 \varepsilon$ for all $z \in \cC_v(\CC_v) \backslash U_v^0$. 
\end{enumerate} 
\end{lemma} 

For now, let $\varepsilon > 0$ and 
$\{\theta_{hj}^{**}\}_{1 \le h \le N, j \in \cS^{\triangle} \bigcup \{j_0\}}$
be any number and collection of points satisfying the conditions of Lemma \ref{MovingLemma1}, 
and let $\cD$ and $Y(z)$ be the corresponding divisor and function.
We will explain the rest of the construction,
then make the final choice of $\varepsilon$ and the $\theta_{hj}^{**}$ at the end. 

\smallskip 
Put $\Gbar_v^{(k_1)}(z) = Y(z) G_v^{(k_1)}(z)$.
We first consider how the coefficients of $G_v^{(k_1)}(z)$ 
change in\index{patching functions, $G_v^{(k)}(z)$ for $1 \le k \le n$!constructed by patching}
\index{coefficients $A_{v,ij}$!high-order} 
passing from $G_v^{(k_1)}(z)$ to $\Gbar_v^{(k_1)}(z)$.  
Write\index{patching functions, $G_v^{(k)}(z)$ for $1 \le k \le n$!constructed by patching}
\begin{eqnarray}
G_v^{(k_1)}(z) \ = \
      \sum_{i=1}^m \sum_{j=0}^{(n-1)N_i-1} A_{v,ij} \varphi_{i,nN_i-j}(z)
        + \sum_{\lambda = 1}^{\Lambda} A_{\lambda} \varphi_{\lambda}(z) \ ,
             \label{DGMN1} \\
\Gbar_v^{(k_1)}(z) \ = \
      \sum_{i=1}^m \sum_{j=0}^{(n-1)N_i-1} \Abar_{v,ij} \varphi_{i,nN_i-j}(z)
        + \sum_{\lambda = 1}^{\Lambda} \Abar_{\lambda} \varphi_{\lambda}(z) \ .
             \label{DGMN2}
\end{eqnarray} 
Because $G_v^{(k_1)}(z)$ and $\Gbar_v^{(k_1)}(z)$ 
are $K_v$-rational,\index{patching functions, $G_v^{(k)}(z)$ for $1 \le k \le n$!are $K_v$-rational}
the $A_{v,ij}$ and $\Abar_{v,ij}$ belong to $L_{w_v}$\index{distinguished place $w_v$} and are $K_v$-symmetric.
\index{$K_v$-symmetric!set of functions}

\begin{lemma} \label{MovingLemma2} {\bf (Second Moving Lemma)}
\index{Second Moving Lemma|ii}
There are constants $\varepsilon_2 > 0$ and $C_3, C_4 \ge 1$ 
$($depending on $E_v$, $\fX$, 
the choices of the $L$-rational and $L^{\sep}$-rational bases, the uniformizers $g_{x_i}(z)$, 
and the projective embedding of $\cC_v),$
such that if $\varepsilon$ and $Y(z)$ are as in Lemma $\ref{MovingLemma1}$,
and in addition $\varepsilon < \varepsilon_2$ and $n$
is sufficiently large, then for all $i = 1, \ldots, m$
and all $0 \le j < k_1 N_i$, we have  
\begin{equation*} 
|\Abar_{v,ij}-A_{v,ij}|_v  \ \le \ C_3 C_4^j (|\tc_{v,i}|_v)^n \cdot \varepsilon \ . 
\end{equation*} 
\end{lemma} 

\smallskip
We next ask about the largest change in the coefficients
\index{coefficients $A_{v,ij}$!high-order} 
$A_{v,ij}$  we can correct for, by re-patching.\index{re-patch}
Write
\begin{equation*}
\Gbar_v^{(k_1)}(z) \ = \ Y(z) \cdot G_v^{(k_1)}(z) \ = \
    Y(z) Q_1^*(z) \cdots Q_{k_1}^*(z) \cdot \hP_{k_1}(\phi_v(z)) \ .
\end{equation*}
Multiplying $G_v^{(k_1)}(z)$ by $Y(z)$ 
moves the $\theta_{hj}^*$\index{patching functions, $G_v^{(k)}(z)$ for $1 \le k \le n$!constructed by patching}
for $j \in \cS^{\triangle}$ to the $\theta_{hj}^{**}$, 
moves the $\theta_{h j_0}$ to the $\theta_{h j_0}^{**}$, 
and leaves all other zeros unchanged.
Recalling that $x-\psi_v(j_0)$ is a factor of $\hP_{k_1}(x)$, write 
\begin{equation*}
\hP_{k_1}(x)  \ = \ (x-\psi_v(j_0)) \cdot \hP_{k_1,j_0}(x)  
\end{equation*}
and put
\begin{equation*}
\Qbar_{k_1}(z) \ = \ Y(z) Q_1^*(z) \cdots Q_{k_1}^*(z) (\phi_v(z)-\psi_v(j_0))
\end{equation*}
so that
\begin{equation} \label{DMUG1}
\Gbar_v^{(k_1)}(z) = \Qbar_{k_1}(z) \cdot \hP_{k_1,j_0,}(\phi_v(z)) \ .
\end{equation}

The function $\Qbar_{k_1}(z)$ accounts 
for the change in passing from\index{patching functions, initial $G_v^{(0)}(z)$} 
$G_v^{(0)}(z)$ to $\Gbar_v^{(k_1)}(z)$.\index{patching functions, $G_v^{(k)}(z)$ for $1 \le k \le n$!constructed by patching}   
It has the following properties.

First, $\Qbar_{k_1}(z)$ is a $K_v$-rational $(\fX,\vs)$-function,   
and it extends to a function defined and finite on $B(\theta_h,\rho_h)$, 
for each $1 \le h \le N$.  Indeed, the zeros of $\Qbar_{k_1}(z)$
are those of $Y(z)$, and the poles of $\Qbar_{k_1}(z)$ are those of 
$Q_1^*(z) \cdots Q_{k_1}^*(z) (\phi_v(z)-\psi_v(j_0))$: 
the poles of $Y(z)$ cancel with the zeros of $Q_1^*(z) \cdots Q_{k_1}^*(z) (\phi_v(z)-\psi_v(j_0))$.     

Second, $\|\Qbar_{k_1}\|_{B(\theta_h,\rho_h)} = 1$ for each $1 \le h \le N$. 
To see this, fix $h$ and restrict $\Qbar_{k_1}(z)$ to $B(\theta_h,\rho_h)$.
Let $\{\eta_{h \ell}\}_{1 \le \ell \le T_h}$ be a list of 
the zeros and poles of $Y(z)$ and the $Q_k^*(z)$ in $B(\theta_h,\rho_h)$, 
and write $\eta_{h \ell} = \hsigma_h(\tau_{h \ell})$ for each $h, \ell$.    
For all $z \in B(\theta_{h},\rho_h) \backslash \big( \bigcup_{\ell = 1}^{T_h} 
B(\eta_{h \ell},\rho_h)^-\big)$ we have $|\Qbar_{k_1}(z)|_v = 1$,
since this is true for each of the factors in its definition.
Pulling this back to  $D(0,1)$, we see that $\Qbar_{k_1}(\hsigma_h(Z))$
is a power series converging in $D(0,1)$, with absolute value $1$ except on
the finitely many subdiscs $D(\tau_{h \ell},1)^-$.
By the Maximum Modulus principle for power series,
$\|\Qbar_{k_1}\|_{B(\theta_h,\rho_h)} = 1$.
\index{Maximum principle!nonarchimedean!for power series}

Third, when the zeros
of $\Qbar_{k_1}(z)$ in $B(\theta_h,\rho_h)$ are pulled back to $D(0,1)$ using $\hsigma_h(Z)$,
they form a union of $\psi_v$-regular sequences of lengths
\index{regular sequence!$\psi_v$-regular sequence}
$1, \#(\cS_1), \dots, \#(\cS_{k_1})$ attached to the sets $ \{j_0\}, \cS_1, \ldots, \cS_{k_1}$.
This holds by construction, since the zeros of
$\phi_v(z)-\psi_v(j_0)$ and $Q_1^*(z), \ldots, Q_{k_1}^*(z)$ have this property, 
and multiplying by $Y(z)$ moves each root\index{roots!move roots} only by an amount
which preserves its position in its $\psi_v$-regular sequence.
\index{regular sequence!$\psi_v$-regular sequence}

Fourth, $\Qbar_{k_1}(z)$ has degree $Nt$, where $t = \#(\cS_1 \cup \cdots \cup \cS_{k+1})$.  
For each $x_i \in \fX$, 
the leading coefficient of $\Qbar_{k_1}(z)$ at $x_i$ has the
\index{coefficients $A_{v,ij}$!leading} 
form $\mu_{v,i} \cdot \tc_{v,i}^t$ where $\mu_{v,i} \in \cO_{w_v}^{\times}$.\index{distinguished place $w_v$}  
Indeed, the leading coefficient of 
\index{coefficients $A_{v,ij}$!leading} 
$(\phi_v(z)-\psi_v(j_0)) Q_1^*(z) \dots Q_{k_1}^*(z)$ 
at $x_i$ is $\varepsilon_{v,i} \cdot \tc_{v,i}^t$,   
while $Y(x_i) \in \cO_{w_v}^{\times}$\index{distinguished place $w_v$} by Lemma \ref{MovingLemma1}(B2).

\smallskip
Our plan is to replace $\Gbar_v^{(k_1)}(z)$ with a new function $\hG_v^{(k_1)}(z)$
with the same high-order coefficients 
as the original $G_v^{(k_1)}(z)$,
\index{patching functions, $G_v^{(k)}(z)$ for $1 \le k \le n$!for nonarchimedean $K_v$-simple sets!roots are separated}
\index{coefficients $A_{v,ij}$!high-order} 
and whose zeros are well-separated from each other.
To do this, we will use the basic patching lemma
via the sequence of  `safe' indices\index{indices!safe} $\cS^0[k_1]$ of length $k_1 + 1$.

Put
\begin{equation*}
P_{k_1}^{0}(x) = \prod_{j \in \cS^0[k_1]} (x-\psi_v(j))\ , \quad
\hP_{k_1}^{0}(x) =
      \prod_{j \in \cS^{\Diamond} \backslash (\{j_0\} \cup \cS^0[k_1])} 
                           (x-\psi_v(j)) \ ,
\end{equation*}
noting that
$\hP_{k_1,j_0}(x) = P_{k_1}^{0}(x) \cdot \hP_{k_1}^{0}(x)$.
Then write
\begin{equation} \label{DMUG3}
\Fbar_{v,k_1}(z) \ = \ \Qbar_{k_1}(z) \cdot \hP_{k_1}^{0}(\phi_v(z)) \ ,
\end{equation}
so that 
\begin{equation} \label{DMUG2} 
\Gbar_v^{(k_1)}(z) \ = \ P_{k_1}^{0}(\phi_v(z)) \cdot \Fbar_{v,k_1}(z) \ .
\end{equation} 
Thus $\Fbar_{v,k_1}(z)$ is a $K_v$-rational $(\fX,\vs)$-function of degree
$N \cdot (n-k_1-1)$, whose roots form a union of $\psi_v$-regular sequences accounting 
\index{regular sequence!$\psi_v$-regular sequence}
for the indices\index{indices} in $\{0, \ldots, n-1\} \backslash \cS^0[k_1]$.  
It has a pole of order $(n-k_1-1)N_i$ at each $x_i$,
with leading coefficient $\mu_{v,i} \tc_{v,i}^{n-k_1-1}$. 
\index{coefficients $A_{v,ij}$!leading} 
Furthermore $\|\Fbar_{v,k_1}\|_{U_v^0} \le 1$, since
$\|\Qbar_{k_1}\|_{U_v^0} = 1$ and $\|\hP_{k_1}^{0}\|_{U_v^0} = 1$.

\begin{lemma} \label{MovingLemma3} {\bf (Third Moving Lemma)}
\index{Third Moving Lemma|ii}  
There are constants $\varepsilon_3 > 0$ and $C_6, C_7 \ge 1$ 
$($depending only on $\phi_v(z)$, $E_v$, $H_v$, their $K_v$-simple decompositions 
\index{$K_v$-simple!set}
\index{$K_v$-simple!decomposition}
$\bigcup_{\ell = 1}^D \big(B(a_\ell,r_\ell) \cap \cC_v(F_{w_\ell})\big)$ 
and $\bigcup_{h=1}^N \big(B(\theta_h,\rho_h) \cap \cC_v(F_{u_h})\big)$, 
the choices of the $L$-rational and $L^{\sep}$-rational bases, the uniformizers $g_{x_i}(z)$, 
and the projective embedding of $\cC_v)$, 
such that if $0 < \varepsilon < \varepsilon_3$,
then there is a $K_v$-rational $(\fX,\vs)$-function $\Deltabar_{v,k_1}(z)$ of the form
\begin{equation*}
\Deltabar_{v,k_1}(z) \ = \ 
\sum_{i=1}^m \sum_{j=0}^{k_1 N_i-1} \Deltabar_{v,ij} \varphi_{i,(k_1+1)N_i-j}(z) \ ,
\end{equation*}
and for which 
\begin{equation*}
\|\Deltabar_{v,ij}(z)\|_{U_v^0} \ \le \  C_6 C_7^{k_1} \varepsilon  \ ,
\end{equation*}
such that 
when $\Gbar_v^{(k_1)}(z)$ from Lemma $\ref{MovingLemma2}$ i
s replaced with\index{patching functions, $G_v^{(k)}(z)$ for $1 \le k \le n$!constructed by patching} 
\begin{equation} \label{FGRepatch}
\hG_v^{(k_1)}(z) \ = \ \Gbar_v^{(k_1)}(z) + \Deltabar_{v,k_1}(z)\Fbar_{v,k_1}(z)\ ,
\end{equation}
then for each $(i,j)$ with $1 \le i \le m$, $0 \le j < k_1 N_i$,
the coefficient $\Abar_{v,ij}$ of $\Gbar_v^{(k_1)}(z)$ is restored to the coefficient 
\index{coefficients $A_{v,ij}$!high-order} 
\index{coefficients $A_{v,ij}$!restoring}\index{patching functions, $G_v^{(k)}(z)$ for $1 \le k \le n$!constructed by patching}
$A_{v,ij}$ of $G_v^{(k_1)}(z)$ in $\hG_v^{(k_1)}(z)$.
\end{lemma} 


Finally, we consider the amount of movement in the roots\index{roots!move roots}
 caused by (\ref{FGRepatch}).\index{separate the roots}  
Since 
\begin{equation*}
\hG_{k_1}(z) \ = \ \big(P_{k_1}^{0}(\phi_v(z)) + \Deltabar_{v,k_1}(z)\big)\cdot \Fbar_{v,k_1}(z) \ ,
\end{equation*} 
the movement is isolated to the roots of $P_{k_1}^{0}(\phi_v(z))$. 

Write  $Q_{k_1}^{0}(z) = P_{k_1}^{0}(\phi_v(z))$, put 
\begin{equation} \label{Q0Stark1Def} 
Q_{k_1}^{0*}(z) \ := \ Q_{k_1}^{0}(z) + \Deltabar_{v,k_1}(z) \ ,   
\end{equation}  
and consider $Q_{k_1}^{0}(z)$ and $Q_{k_1}^{0*}(z)$  on each ball $B(\theta_h,\rho_h)$. 
Pull them back to $D(0,1)$ using $\hsigma_h(Z)$,
and apply the Basic Patching Lemma (Lemma \ref{DLemCPC4}).\index{Basic Patching Lemma}  The roots of
$Q_{k_1}^{0}(\hsigma_h(Z)))$ in $D(0,1)$ form a $\psi_v$-regular sequence
\index{regular sequence!$\psi_v$-regular sequence}
of length $k_1+1$ in $\cO_{u_h}$ attached to $\cS^0[k_1]$,
namely $\{\alpha_{hj}\}_{j \in \cS^0[k_1]}$.
We have chosen $\cS^0[k_1]$ to be a `safe' sequence, which means
these $\alpha_{hj}$ are the only roots within their $\delta_n$-cosets.\index{$\delta_n$-coset}  

If we can arrange that the roots\index{roots!move roots} only move within their $\delta_n$-cosets,\index{$\delta_n$-coset}
they will remain separated from all the other roots.\index{separate the roots}   
To assure this it is enough to require
\begin{equation} \label{DGEST5}
C_6 C_7^{k_1} \varepsilon \ \le \ q^{-\frac{k_1+1}{q-1} - \log_v(n)}  \ .
\end{equation}
Under this condition, Lemma \ref{DLemCPC4} shows that for each $h$ the roots
$\{\alpha_{hj}^*\}_{j \in \cS^0[k_1]}$ of 
$\Qbar_{k_1}^{0*}(\hsigma_h(Z)))$
form a $\psi_v$-regular sequence of length $k_1+1$ in $\cO_{u_h}$, satisfying
\index{regular sequence!$\psi_v$-regular sequence}
\begin{equation*}
\ord_v(\alpha_{hj}^*-\alpha_{hj}) \ \ge \ \log_v(n)
\end{equation*}
for each $j \in \cS^0[k_1]$.

\vskip .1 in
We will now specify  $\varepsilon$ and the $\theta_{hj}^{**}$.  

We want $\varepsilon$ to be as large as possible.  
For the construction to succeed, we must have  
$\varepsilon < \min(\varepsilon_1,\varepsilon_2,\varepsilon_3)$, 
and (\ref{FEPSCOND1}) and (\ref{DGEST5}) must hold:  
\begin{equation*}
C_1 \varepsilon \le q^{-\lceil \log_v(n) \rceil} \cdot \min_{1 \le h \le N}(\rho_h), 
\qquad 
C_6  C_7^{k_1} \varepsilon \ \le \ q^{-\frac{k_1+1}{q-1} - \log_v(n)} \ .
\end{equation*} 
Since $k_1 \le A_1 \log_v(n)$, 
then for an appropriate constant $A_3$ we can choose $\varepsilon$ so that 
\begin{equation} \label{DGEST6}
-\log_v(\varepsilon) \ = \ A_3 \log_v(n) \ ,
\end{equation}
provided $n$ is sufficiently large.

We next choose the  $\theta_{hj}^{**}$ for $1 \le h \le N$, $j \in \cS^{\triangle}$.   

Fixing $h$, let $\hsigma_h : D(0,1) \rightarrow B(\theta_h,\rho_h)$ be the parametrization 
used before;  specifying the $\theta_{hj}^{**}$ is equivalent to specifying numbers 
$\alpha_{hj}^{**} \in \cO_{u_h}$ 
such that $\hsigma_h(\alpha_{hj}^{**}) = \theta_{hj}^{**}$.
Recall that $\|\hsigma_h(x),\hsigma_h(y)\|_v = \rho_h |z-x|_v$ for all 
$x, y \in D(0,1)$.  Put $\rhobar = \max_{1 \le h \le N}(\rho_h)$ 
and let $\varepsilon_0 = \varepsilon/\rhobar$, 
noting that $\varepsilon_0 \le \delta_n$ by (\ref{FEPSCOND1}).  
In Lemma \ref{MovingLemma1} we can move the $\alpha_{hj}^*$ 
for $j \in \cS^{\triangle}$
to arbitrary points $\alpha_{hj}^{**} \in \cO_{u_h}$
such that $|\alpha_{hj}^{**}-\alpha_{hj}^*|_v  \le  \varepsilon_0$
(provided the collection $\{\alpha_{hj}^{**}\}$ is $K_v$-symmetric),
\index{$K_v$-symmetric!set of points}
while only moving the $\alpha_{h j_0}$ within their $\delta_n$-cosets.\index{$\delta_n$-coset}  
However, we need to choose the $\alpha_{hj}^{**}$ in such a way that they 
become well-separated from each other and from the unpatched roots.\index{roots!unpatched}  

For a given $\alpha \in \cO_{u_h}$, consider the $\varepsilon_0$-coset of $\alpha$ in $\cO_{u_h}$,  
\begin{equation*}
D_{u_h}(\alpha,\varepsilon_0) \ := \ \cO_{u_h} \cap D(\alpha,\varepsilon_0)
      \ = \ \{ z \in \cO_{u_h} : |z-\alpha|_v \le \varepsilon_0 \} \ .
\end{equation*}
The roots $\alpha_{hj}^*$ (which correspond to the indices\index{indices} 
$j \in \cS^{\triangle} = \cS_1 \cup \cdots \cup \cS_{k_1}$),
form a union of $\psi_v$-regular subsequences of respective lengths
$\#(\cS_1),  \ldots, \#(\cS_{k_1})$.  At most one $\alpha_{hj}^*$ from each subsequence
can belong to $D_{u_h}(\alpha,\varepsilon_0)$.  Since the original sequence
$\{\alpha_{hj}\}_{0 \le j \le n-1}$ was a $\psi_v$-regular sequence of length $n$ in $\cO_{u_h}$ 
\index{regular sequence!$\psi_v$-regular sequence}
there is at most one unpatched\index{roots!unpatched} root $\alpha_{h j^{\prime}}$ 
in $D_{u_h}(\alpha,\varepsilon_0)$.  Thus $D_{u_h}(\alpha,\varepsilon_0)$ contains
at most $k_1+1$ roots.  Furthermore, 
if $\alpha_{h,j_0}^{**} \in D_{u_h}(\alpha,\varepsilon_0)$, 
then since $\varepsilon_0 \le \delta_n$,
our choice of $j_0$ means that $D_{u_h}(\alpha,\varepsilon_0)$ does not contain
any of the $\alpha_{hj}^*$.

Put $\delta_0 = q^{-\lceil \log_v(k_1 + 1) \rceil }$. 
Since $\log_v(k_1) \le \log_v(A_1 \log_v(n))$, 
it follows from (\ref{DGEST6}) that
\begin{equation} \label{DGEST8}
-\log_v(\delta_0 \varepsilon_0) \  =  \ -\log_v(\delta_0 \varepsilon/\rhobar) 
                              \ \le \ A_4 \log_v(n) 
\end{equation}
for an appropriate constant $A_4$ .
There are at least  $q^{\lceil \log_v(k_1 + 1) \rceil } \ge k_1 + 1$
distinct $\delta_0 \varepsilon_0$-cosets\index{$\delta_n$-coset} 
$D_{u_h}(\beta,\delta_0 \varepsilon_0) \subset D_{u_h}(\alpha,\delta_0)$
with $\beta \in \cO_{u_h}$.  
By simultaneously adjusting all the $\alpha_{hj}^*$
belonging to $D_{u_h}(\alpha,\varepsilon_0)$ 
we can choose new roots\index{roots!move roots} 
$\alpha_{hj}^{**} \in D_{u_h}(\alpha,\varepsilon_0) = \cO_{u_h} \cap D(\alpha,\varepsilon_0)$ 
which are separated from each other and from the unpatched root $\alpha_{h j^{\prime}}$\index{roots!unpatched}
(if it exists), by a distance at least $\delta_0 \varepsilon_0$.  
Do this for each $h$ and each $D_{u_h}(\alpha,\varepsilon_0)$, making the choices  
for different $h$ in a galois-equivariant way. 
 
It follows that for each $h$, and each
$j \in \cS_1 \cup \cdots \cup \cS_{k_1}$, we can choose the $\alpha_{hj}^{**} \in \cO_{u_h}$
so that $|\alpha_{hj}^{**}-\alpha_{hj}^*|_v \le \varepsilon_0$ and 
\begin{enumerate}
\item for each $\ell \in \cS_1 \cup \cdots \cup \cS_{k_1}$ with $\ell \ne j$,
\begin{equation*}
     \ord_v(\alpha_{hj}^{**}-\alpha_{h \ell}^{**})
            \ \le \ A_4 \log_v(n) \ ,
\end{equation*}
\item for each $\ell \in \{j_0\} \cup \cS^0[k_1]$ 
\begin{equation*} 
      \ord_v(\alpha_{hj}^{**}-\alpha_{h \ell}^*) \ \le \ A_4 \log_v(n) \ , 
\end{equation*}   
\item for each $\ell \in \{0,1,\ldots,n-1\} \backslash 
           (\{j_0\} \cup \cS_1 \cup \cdots \cup \cS_{k_1} \cup \cS^0[k_1])$
\begin{equation*}
     \ord_v(\alpha_{hj}^{**}-\alpha_{h \ell}) \ \le \ A_4 \log_v(n) \ .
\end{equation*} 
\end{enumerate}
Property (2) holds because  
$\{j_0\} \bigcup \cS_{k_1}^0 \subset \cS^{\heartsuit}$.
Fix such a choice of the $\alpha_{hj}^{**}$, 
and define $\hG_v^{(k_1)}(z)$\index{patching functions, $G_v^{(k)}(z)$ for $1 \le k \le n$!constructed by patching}
by means of 
Lemma \ref{MovingLemma3}.

\vskip .1 in
Given  $\alpha, \alpha^{\prime} \in \CC_v$, we will say that 
$\alpha$ and $\alpha^{\prime}$ are 
``logarithmically separated by at least $T$''  if\index{logarithmically separated|ii}
\begin{equation*}
\ord_v(\alpha^{\prime}-\alpha) \ \le \ T  
\end{equation*}
or equivalently, if $|\alpha^{\prime}-\alpha|_v \ge q_v^{-T}$. 

To ease the notation in subsequent steps,
relabel the roots $\alpha_{hj}^{**}$ as $\alpha_{hj}^*$, 
and replace $G_v^{(k_1)}(z)$ with $\hG_v^{(k_1)}(z)$.
\index{patching functions, $G_v^{(k)}(z)$ for $1 \le k \le n$!for nonarchimedean $K_v$-simple sets!roots are separated}
In the statements of Theorems \ref{DCPCPatch} and \ref{DCPCPatch1p}, 
this is accounted for by adding\index{patching functions, $G_v^{(k)}(z)$ for $1 \le k \le n$!constructed by patching} 
$\Theta_v^{(k_1)}(z) := \hG_v^{(k_1)}(z) - G_v^{(k_1)}(z)$ to $G_v^{(k_1)}(z)$.  
Thus, the new function $G_v^{(k_1)}(z)$ has the same high-order coefficients\index{patching functions, $G_v^{(k)}(z)$} 
\index{coefficients $A_{v,ij}$!high-order} 
as the old one, and its roots\index{separate the roots}  
are logarithmically separated\index{logarithmically separated} from each other by at least $A_4 \log_v(n)$. 
\index{local patching for nonarch $K_v$-simple sets!Phase 3: move roots apart|)}

\vskip .1 in
\noindent{\bf Phase 4.  Patching with the long safe sequence for $k = k_1 + 1, \ldots, k_2$.}
\index{local patching for nonarch $K_v$-simple sets!Phase 4: using the long safe sequence|(}  

We have now arrived at a function $G_v^{(k_1)}(z)$
\index{patching functions, $G_v^{(k)}(z)$ for $1 \le k \le n$!for nonarchimedean $K_v$-simple sets!movement of roots} 
with ``patched'' roots $\alpha_{hj}^*$ for\index{roots!patched} 
$j \in \cS_1 \cup \cdots \cup \cS_{k_1} \cup \{j_0\} \cup \cS^0[k_1]$,
and `unpatched' roots $\alpha_{hj}$ for all other $j$.  By the construction\index{roots!unpatched}
in Phase 3, for all $h$, and all $j \ne \ell$, the roots
satisfy
\begin{equation*}
\ord_v(\alpha_{hj}^* - \alpha_{h \ell}^*), \
\ord_v(\alpha_{hj}^* - \alpha_{h \ell}),  \
\ord_v(\alpha_{hj} - \alpha_{h \ell}) \ \le A_4 \log_v(n) \ .
\end{equation*}
as appropriate for each $j, \ell$. 

The number $k_1$ has the property that for all $k \ge k_1$
\begin{equation} \label{FK1}
h_v^{kN} \ \le \ q^{-\frac{k}{q-1}-\log_v(n)} \ .
\end{equation}
The purpose of Phase 4 is to carry on the patching process
until $h_v^{kN}$ is so much smaller than  $q^{-k/(q-1)}$ 
that the ``endangered'' roots $\alpha_{hj}$ for\index{roots!endangered}
$j \in \cS^{\dagger}$ can be included in the patching process: 
this will allow us to apply the Refined Patching Lemma\index{Refined Patching Lemma} 
(Lemma \ref{DLemC6}) in Phase 5.

Let $k_2$ be the least integer
such that for all $k \ge k_2$,
\begin{equation} \label{FK2} 
h_v^{kN} \cdot M_v 
     \ \le \ q^{-\frac{k+1}{q-1} - 3 A_4 \log_v(n)} \ .
\end{equation}
Thus, for an appropriate constant $A_5$, 
\begin{equation} \label{DFROG1}
k_2 \ \le \ A_5 \log_v(n) \ .
\end{equation} 
 
As in Phase 3, we will use the ``long safe sequence''\index{long safe sequence} of roots\index{roots!safe} $\cS^0$ 
in patching.\index{roots!long safe sequence}  We will now impose a condition on $n$ which means 
that $\cS^0$ is actually much longer than was required by (\ref{FK10}).
By (\ref{FK1}), patched roots only move within their\index{roots!patched} 
$\delta_n$-cosets,\index{$\delta_n$-coset} so they maintain their position within a $\psi_v$-regular 
sequence of length $n$.  This means that instead of choosing a new 
$\psi_v$-regular subsequence of length $k$ to use in patching at each step, 
we can simply extend the previous one.  

\vskip .1 in 
Recall (\ref{DFHWC1}).  If $n$ is large enough that 
\begin{equation} \label{DBBB1} 
n \ \ge \ ((N+1) A_2 (\log_v(n))^2 + 1) \cdot (A_5 \log_v(n) + 2) \ .
\end{equation}
which we will henceforth assume, then by the Pigeon-hole Principle,\index{Pigeon-hole Principle} 
the long safe sequence satisfies\index{roots!long safe sequence|ii}\index{long safe sequence|ii}  
\begin{equation*}
\#(\cS^0) \ \ge \ A_5 \log_v(n) + 2 \ \ge \ k_2 + 2 \ .
\end{equation*}
Recall that we write $\cS^0 = \{j_0, j_1, \ldots, j_L \}$.
For each $ k = k_1+1, \ldots, k_2$, put 
\begin{equation*}
\cS^0[k] \ = \ \{j_1, j_2, \ldots, j_{k+1}\} \ .
\end{equation*} 
Also recall that  $\cS^{\Diamond} = \{0,1,\ldots,n-1\} \backslash
    (\cS_1 \cup \cdots \cup \cS_{k_1})$, and that\index{patching functions, $G_v^{(k)}(z)$ for $1 \le k \le n$!factorization of}  
\begin{equation*}
G_v^{(k_1)}(z) 
     \ = \ \Qbar_{k_1}(z) \cdot Q_{k_1}^{0*}(z) \cdot \hP_{k_1}^{0}(\phi_v(z)) \ .
\end{equation*}

For $k = k_1 + 1, \ldots, k_2$ we will patch\index{patching functions, $G_v^{(k)}(z)$ for $1 \le k \le n$!constructed by patching} 
$G_v^{(k-1)}(z)$ to $G_v^{(k)}(z)$ as follows.  
Noting that $\cS^0[k] = \cS^0[k-1] \cup \{j_{k+1}\}$,  
define 
\begin{equation*}
\hP_k^0(x) \ = \ 
        \prod_{j \in \cS^{\Diamond} \backslash (\{j_0\} \cup \cS^0[k])} 
              (x - \psi_v(j)) \ = \ \hP_{k-1}^0(x)/(x-\psi_v(j_{k+1})) \ .
\end{equation*}
Then $\hP_{k-1}^0(x) = (x-\psi_v(j_{k+1})) \cdot \hP_k^0(x)$, and if we set 
\begin{equation*}
Q_k^0(z) \ = \ Q_{k-1}^{0*}(z) \cdot (\phi_v(z) - \psi_v(j_{k+1})) \ , 
\end{equation*}
then\index{patching functions, $G_v^{(k)}(z)$ for $1 \le k \le n$!factorization of}
\begin{equation} \label{DBBC2}
G_v^{(k-1)}(z) \ = \ \Qbar_{k_1}(z) \cdot Q_k^0(z) \cdot \hP_k^0(\phi_v(z)) \ . 
\end{equation}

By construction, when the roots of $Q_k^0(z)$ in  
$B(\theta_h,\rho_h)$ are pulled back to $D(0,1)$ using $\hsigma_h(Z)$,
they form a $\psi_v$-regular sequence of length $k+1$ in $\cO_{u_h}$.  
\index{regular sequence!$\psi_v$-regular sequence}
For notational simplicity, we will relabel these roots
(the $\alpha_{hj}^*$ for $j \in \cS^0[k-1]$, 
together with $\alpha_{h,j_{k+1}}$) as 
$\{\dot{\alpha}_{hj}\}_{j \in \cS^0[k]}$.  

\smallskip
When $\Char(K_v) = 0$, we are given a $K_v$-symmetric collection of numbers 
\index{$K_v$-symmetric!set of numbers}\index{band!$\Band_N(k)$}
$\{\Delta_{v,ij} \in L_{w_v}\}_{(i,j) \in \Band_N(k)}$\index{distinguished place $w_v$} 
determined recursively in $\prec_N$ order,\index{order!$\prec_N$} and   
we patch $G_v^{(k-1)}(z)$ 
by setting\index{compensating functions $\vartheta_{v,ij}^{(k)}(z)$}
\index{patching functions, $G_v^{(k)}(z)$ for $1 \le k \le n$!constructed by patching}
\begin{eqnarray*} 
G_v^{(k)}(z) & = & G_v^{(k-1)}(z) + \sum_{i=1}^m \sum_{j=(k-1)N_i}^{kN_i-1} 
\Delta_{v,ij}^{(k)} \vartheta_{v,ij}^{(k)}(z) \\
& = & Q_1^*(z) \cdots Q_{k-1}^*(z) \cdot \big(P_k(\phi_v(z)) + \Delta_{v,k}(z) \big) \cdot \hP_k(\phi_v(z)) \ ,
\end{eqnarray*}
where the compensating functions are\index{compensating functions $\vartheta_{v,ij}^{(k)}(z)$!construction of}  
\begin{equation*}
\vartheta_{v,ij}^{(k)}(z) \ = \ \varepsilon_{v,i}^{-1} \varphi_{i,(k+1)N_i - j}(z) \cdot 
Q_1^*(z) \cdots Q_{k-1}^*(z) \cdot \hP_k(\phi_v(z)) \ , 
\end{equation*}  
and where as in (\ref{DFFZ3}),  
\begin{equation*}
\Delta_{v,k}(z) \ = \ \sum_{i=1}^m \sum_{j=(k-1)N_i}^{kN_i-1} 
   \Delta_{v,ij}^{(k)} \cdot \varepsilon_{v,i}^{-1} \varphi_{i,(k+1)N_i-j}(z) \ .
\end{equation*}

We claim that the leading coefficient of $\vartheta_{v,ij}^{(k)}(z)$
\index{compensating functions $\vartheta_{v,ij}^{(k)}(z)$!poles and leading coefficients of} 
\index{coefficients $A_{v,ij}$!leading} 
at $x_i$ is $\tc_{v,i}^{n-k-1}$. 
To see this, note that the leading coefficient of  
$G_v^{(k-1)}(z)$\index{patching functions, $G_v^{(k)}(z)$ for $1 \le k \le n$!leading coefficients of} at $x_i$ 
is $A_{v,i0} = \varepsilon_{v,i} \tc_{v,i}^n$,   
while the leading coefficient of $\hP_k^0(\phi_v(z))$ at $x_i$ 
is $\tc_{v,i}^{k+1}$.  By (\ref{DBBC2}), the leading 
coefficient of $ \Qbar_{k_1}(z) \cdot \hP_k(\phi_v(z))$ 
\index{coefficients $A_{v,ij}$!leading} 
is $\varepsilon_{v,i} \tc_{v,i}^{n-k-1}$, which gives what we want. 
Note also that\index{compensating functions $\vartheta_{v,ij}^{(k)}(z)$!bounds for} 
\begin{equation*}  
\|\vartheta_{v,ij}^{(k)}\|_{U_v^0} \ \le \ M_v \ .
\end{equation*} 
Indeed,  
$\|\Qbar_{k_1}(z)\|_{U_v^0} \le 1$ and  $\|\hP_k^0(\phi_v(z))\|_{U_v^0} \le 1$, 
while $\| \varphi_{i,r_{ij}}\|_{U_v^0} \le M_v$. 

Clearly the $\vartheta_{v,ij}^{(k)}(z)$ are $K_v$-symmetric.
\index{compensating functions $\vartheta_{v,ij}^{(k)}(z)$!are $K_v$-symmetric}  
\index{$K_v$-symmetric!set of functions} 
It follows that $\Delta_{v,k}(z)$ and $G_v^{(k)}(z)$ 
are $K_v$-rational,\index{patching functions, $G_v^{(k)}(z)$ for $1 \le k \le n$!are $K_v$-rational} 
and for each $(i,j)$ 
\begin{equation*}
\sum_{x_{i^{\prime}} \in \Aut_c(\CC_v/K_v)(x_i)} 
   \Delta_{v,i^{\prime}j}^{(k)} \vartheta_{v,i^{\prime}j}^{(k)}    
\end{equation*}
is $K_v$-rational.

\smallskip
When $\Char(K_v) = p > 0$, let  
\begin{equation*} 
F_{v,k}(z) \ = \ Q_1^*(z) \cdots Q_{k-1}^*(z) \cdot \hP_k(\phi_v(z)) \ .
\end{equation*}
By arguments similar to those above, $F_{v,k}(z)$ is $K_v$-rational, its roots belong in $E_v$, 
and it has a pole of order $(n-k-1)N_i$
with leading coefficient $d_{v,i} = \varepsilon_{v,i} \tc_{v,i}^{n-k-1}$ at each $x_i$.
\index{coefficients $A_{v,ij}$!leading}   
In particular $|d_{v,i}|_v = |\tc_{v,i}^{n-k-1}|_v$.
Thus $F_{v,k}(z)$ satisfies the conditions of Theorem \ref{DCPCPatch1p}.

By the hypotheses of Theorem \ref{DCPCPatch1p} we are given a $K_v$-symmetric collection of numbers 
\index{$K_v$-symmetric!set of functions}\index{band!$\Band_N(k)$}
$\{\Delta_{v,ij} \in L_{w_v}\}_{(i,j) \in \Band_N(k)}$\index{distinguished place $w_v$} such that the function 
\begin{equation*}
\Delta_{v,k}(z) \ = \
 \sum_{i=1}^m \sum_{j=(k-1)N_i}^{kN_i-1} \Delta_{v,ij}^{(k)} \cdot \varphi_{i,(k+1)N_i-j}(z) \ ,
\end{equation*} 
is $K_v$-rational.  
We patch $G_v^{(k-1)}(z)$ by setting\index{patching functions, $G_v^{(k)}(z)$ for $1 \le k \le n$!constructed by patching} 
\begin{eqnarray*} 
G_v^{(k)}(z) & = & G_v^{(k-1)}(z) + \Delta_{v,k}(z) \cdot F_{v,k}(z) \\
& = &  Q_1^*(z) \cdots Q_{k-1}^*(z) \cdot \big(P_k(\phi_v(z)) + \Delta_{v,k}(z) \big) \cdot  \hP_k(\phi_v(z)) 
\end{eqnarray*}

\smallskip
The two patching constructions differ only in the choice of $\Delta_{v,k}(z)$,
and in both cases  
$\|\Delta_{v,k}\|_{U_v^0} \le h_v^{kN} \cdot M_v < q^{-(k/(q+1) - \log_v(n)}$.   
Write 
\begin{equation*}
Q_k^{0*}(z) \ = \ Q_k^0(z) + \Delta_{v,k}(z)   \ .
\end{equation*}
The change in passing from $G_v^{(k-1)}(z)$ to $G_v^{(k)}(z)$
\index{patching functions, $G_v^{(k)}(z)$ for $1 \le k \le n$!for nonarchimedean $K_v$-simple sets!movement of roots} 
is localized to $G_k^{0*}(z)$, and  
\begin{equation*}
G_v^{(k)}(z) \ = \Qbar_{k_1}(z) \cdot Q_k^{0*}(z)
                            \cdot \hP_k(\phi_v(z)) \ . 
\end{equation*}

By Lemma \ref{DLemCPC4}, when the roots of $Q_k^{0*}(z)$ 
in  $B(\theta_h,\rho_h)$ are pulled back to $D(0,1)$ 
using $\hsigma_h(Z)$, they form a $\psi_v$-regular sequence $\{\alpha_{hj}^*\}$
\index{regular sequence!$\psi_v$-regular sequence}
of length $k+1$ in $\cO_{u_h}$ attached to $\cS^0[k]$.  Since $k > k_1$ 
we have $h_v^{kN} \cdot M_v \le q^{-\frac{k+1}{q-1} -  \log_v(n)}$, which means that 
\begin{equation*}
\ord_v(\alpha_{hj}^*-\dot{\alpha}_{hj}) \ \ge \ \log_v(n) 
\end{equation*}
for each $j$. The fact that $\cS^0$ consists only of ``safe'' indices\index{indices!safe} means
that the $\alpha_{hj}^*$ for $j \in \cS^0[k]$ remain the only roots within 
their $\delta_n$-cosets,\index{$\delta_n$-coset} and have not moved nearer to any of the other roots.\index{roots} 
Thus
\begin{equation*}
\{\alpha_{hj}^*\}_{j \in \cS^0[k]} 
     \cup \{\alpha_{hj}\}_{\cS^0 \backslash \cS^0[k]}
\end{equation*}
is again a $\psi_v$-regular sequence attached to $\cS^0$,  
\index{regular sequence!$\psi_v$-regular sequence}
and the induction can continue.
\index{local patching for nonarch $K_v$-simple sets!Phase 4: using the long safe sequence|)}    

\vskip .1 in 
\noindent{\bf Phase 5. Patching using the remaining unpatched indices.}\index{indices!unpatched} 
\index{local patching for nonarch $K_v$-simple sets!Phase 5: patch unpatched indices|(}  

\vskip .1 in
At this point  
we have constructed a $K_v$-rational $(\fX,\vs)$-function\index{patching functions, $G_v^{(k)}(z)$ for $1 \le k \le n$!factorization of}
\begin{equation*}
G_v^{(k_2)}(z) \ = \ \Qbar_{k_1}(z) \cdot Q_{k_2}^{0*}(z)
                            \cdot \hP_{k_2}^0(\phi_v(z)) 
\end{equation*} 
whose coefficients $A_{v,ij}$ have been patched for all $(i,j)$ 
\index{coefficients $A_{v,ij}$!middle} 
with  $1 \le i \le m$, $0 \le j < k_2 N_i$.  
When the roots of $G_v^{(k_2)}(z)$
\index{patching functions, $G_v^{(k)}(z)$ for $1 \le k \le n$!for nonarchimedean $K_v$-simple sets!movement of roots} 
in $B(\theta_h,\rho_h)$ are pulled back to $D(0,1)$
using $\hsigma_h(Z)$, they form a union of $\psi_v$-regular sequences
\index{regular sequence!$\psi_v$-regular sequence}
consisting of patched roots $\alpha_{hj}^*$ corresponding to the sets\index{roots!patched} 
$\{j_0\}$, $\cS_1, \ldots, \cS_{k_1}$,  and $\cS^0[k_2]$, 
together with unpatched roots $\alpha_{hj}$ for $j$ in the set\index{roots!unpatched}
\begin{equation*}
\cS^{\Diamond} \backslash (\{j_0\} \cup \cS^0[k]) 
       \ = \ \{0,1,\ldots,n-1\} \backslash
         ( \{j_0\} \cup \cS_1 \cup \cdots \cup \cS_{k_1} \cup \cS^0[k_2] ) \ .
\end{equation*}
As the roots of $G_v^{(k_1)}(z)$ were 
logarithmically separated by\index{roots!logarithmically separated}\index{logarithmically separated}
\index{patching functions, $G_v^{(k)}(z)$ for $1 \le k \le n$!for nonarchimedean $K_v$-simple sets!roots are separated}
at least $A_4 \log_v(n)$ and since Phase 4 only patched using ``safe'' roots,\index{roots!safe}
the roots of $G_v^{(k_2)}(z)$ remain logarithmically separated by\index{logarithmically separated}
at least $A_4 \log_v(n)$.
\index{patching functions, $G_v^{(k)}(z)$ for $1 \le k \le n$!for nonarchimedean $K_v$-simple sets!roots are separated}   

\smallskip

Note that $\cS^{\Diamond} = \{0,1, \ldots, n-1\} \backslash (\cS_1 \cup \cdots \cup \cS_{k_1})$
is a sequence of consecutive indices,\index{indices!consecutive} both when $\Char(K_v) = 0$ and when $\Char(K_v) = p > 0$. 
Since $\{j_0\} \cup\cS^0[k_2] \subset \cS^0 \subset \cS^{\Diamond}$ 
is also a sequence of consecutive indices,\index{indices!consecutive} 
its complement in $\cS^{\Diamond}$ consists of at most two sequences of consecutive indices.  
Recalling that $\{j_0\} \cup \cS^0[k_2] = \{j_0, j_1, \ldots, j_{k_2+1}\}$,
put $k_3 = \#(\cS^{\Diamond})-3$ and list the elements of 
$\cS^{\Diamond} \backslash (\{j_0\} \cup \cS^0[k_2])$ in increasing order 
as $\{j_{k_2+2}, \ldots, j_{k_3+1} \}$.  For each $k$ with 
$k_2 < k \le k_3$, put 
\begin{equation*}
\cS^{\Diamond}[k] \ = \ \cS^0[k_2] \cup \{j_{k_2+2}, \ldots, j_{k+1} \} \ .
\end{equation*} 
By the discussion above, $\cS^{\Diamond}[k]$ 
is a union of most $3$ subsequences of consecutive indices.\index{indices!consecutive}  

Recall from \S\ref{Chap11}.\ref{PatchingLemmasSection} the 
\vskip .1 in 

\noindent{\bf Lemma $\text{\bf \ref{DLemC6}}$. (Refined Patching Lemma)}\index{Refined Patching Lemma} 
{\it Let  $\cQ(Z) \in K_v[[Z]]$ be a power series converging on $D(0,1)$,
with $\sup$ norm $\|\cQ\|_{D(0,1)} = 1$.  Suppose the roots $\{\alpha_j\}$
of $\cQ(Z)$ in $D(0,1)$ can be partitioned into $r$ disjoint
$\psi_v$-regular sequences in $\cO_w$ attached to index sets 
\index{regular sequence!$\psi_v$-regular sequence}
$\cS_1, \ldots, \cS_r$ of respective lengths $\ell_1, \ldots, \ell_r$.
Put $\ell = \sum_{k=1}^r \ell_k$.  Suppose further that there is a bound
$T \ge \max_i(\log_v(\ell_i))$ such that
\begin{equation*}
\ord_v(\alpha_j-\alpha_k) \ \le \ T
\end{equation*}
for all $j \ne k$.

Then for any $M \ge T$, and any power series $\Delta(Z) \in K_v[[z]]$ 
converging on $D(0,1)$ with  
\begin{equation*}
\|\Delta\|_{D(0,1)} \ \le \ q^{-\frac{\ell}{q-1} - (r-1)T - M} \ .
\end{equation*}
the roots $\{\alpha_j^*\}$ of $\cQ^*(Z) = \cQ(Z) + \Delta(Z)$ in $D(0,1)$ 
again form a union of $\psi_v$-regular
sequences in $\cO_w$ attached to $\cS_1, \ldots, \cS_r$.  They can 
uniquely be labelled in such a way that
\begin{equation*}
\ord_v(\alpha_j^* - \alpha_j) \ > \ M
\end{equation*}
for each $j$.
}
\vskip .1 in

The number $k_2$ was chosen so that  
if $r = 3$ and $T = A_4 \log_v(n)$, then for all $k \ge k_2$,
\begin{equation} \label{DMONK}
h_v^{kN} \cdot M_v 
     \ \le \ q^{-\frac{k+1}{q-1} - rT} \ .
\end{equation}
This means that when we 
apply the Refined Patching Lemma\index{Refined Patching Lemma} using at most $3$ sequences of roots, 
all roots\index{roots!logarithmically separated} 
will remain logarithmically separated\index{logarithmically separated}
 by at least $A_4 \log_v(n)$.   
   
\vskip .1 in 
When $k = k_2$ write
\begin{equation*}
Q_{k_2}^{\Diamond *}(z) \ = \ Q_{k_2}^{0*}(z) \ , 
     \qquad  \hP_{k_2}^{\Diamond}(z) \ = \ \hP_{k_2}^{0}(z) \ ,
\end{equation*}
so that with this notation\index{patching functions, $G_v^{(k)}(z)$ for $1 \le k \le n$!factorization of}
\begin{equation*}
G_v^{(k_2)}(z) \ = \ \Qbar_{k_1}(z) \cdot Q_{k_2}^{\Diamond *}(z) 
                          \cdot \hP_{k_2}^{\Diamond}(\phi_v(z)) \ .
\end{equation*}    
For $k = k_2+1, \ldots, k_3$,  inductively suppose that\index{patching functions, $G_v^{(k)}(z)$ for $1 \le k \le n$!factorization of} 
\begin{equation*}
G_v^{(k-1)}(z) \ = \ \Qbar_{k_1}(z) \cdot Q_{k-1}^{\Diamond *}(z) 
                       \cdot \hP_{k-1}^{\Diamond}(\phi_v(z)) 
\end{equation*} 
where the roots of $Q_{k-1}^{\Diamond *}(z)$ 
correspond to $\cS^{\Diamond}[k-1]$.  Put 
\begin{equation*}    
\hP_k^{\Diamond}(z) 
 \ = \ \prod_{j \in \cS^{\Diamond} \backslash 
                  (\{j_0\} \cup \cS^{\Diamond}[k])} (z-\psi_v(j)) 
 \ = \ \hP_{k-1}^{\Diamond}(z)/(z-\psi_v(j_{k+1}))                               
\end{equation*} 
and put 
\begin{equation*}
Q_k^{\Diamond}(z)  \ = \ Q_{k-1}^{\Diamond *}(z) \cdot (\phi_v(z)-\psi_v(j_{k+1}))
\end{equation*}
so that\index{patching functions, $G_v^{(k)}(z)$ for $1 \le k \le n$!factorization of} 
\begin{equation*}
G_v^{(k-1)}(z) \ = \  \Qbar_{k_1}(z) \cdot Q_k^{\Diamond}(z)
                             \cdot \hP_k^{\Diamond}(\phi_v(z)) \ .
\end{equation*}

\smallskip
The patching argument in Phase 5 is very similar to that in Phase 4.  

When $\Char(K_v) = 0$, we are given a $K_v$-symmetric collection of numbers 
\index{$K_v$-symmetric!set of numbers}\index{band!$\Band_N(k)$}
$\{\Delta_{v,ij} \in L_{w_v}\}_{(i,j) \in \Band_N(k)}$\index{distinguished place $w_v$} 
determined recursively in $\prec_N$ order,\index{order!$\prec_N$} and   
we patch $G_v^{(k-1)}(z)$\index{patching functions, $G_v^{(k)}(z)$ for $1 \le k \le n$!constructed by patching} 
by setting\index{compensating functions $\vartheta_{v,ij}^{(k)}(z)$}  
\begin{eqnarray} 
G_v^{(k)}(z) & = & G_v^{(k-1)}(z) + \sum_{i=1}^m \sum_{j=(k-1)N_i}^{kN_i-1} 
\Delta_{v,ij}^{(k)} \vartheta_{v,ij}^{(k)}(z) \label{DBLB2} \\
& = & \Qbar_{k_1}(z) \cdot \big(P_k(\phi_v(z)) + \Delta_{v,k}(z) \big) \cdot \hP_k^{\Diamond}(\phi_v(z)) \ ,
\notag 
\end{eqnarray}
where the compensating functions are\index{compensating functions $\vartheta_{v,ij}^{(k)}(z)$!construction of}  
\begin{equation} \label{DBLB1}
\vartheta_{v,ij}^{(k)}(z) \ = \ 
    \varepsilon_{v,i}^{-1}  \varphi_{i,(k+1)N_i-j}(z) \cdot 
            \Qbar_{k_1}(z) \cdot \hP_k^{\Diamond}(\phi_v(z))
\end{equation} 
and where as in (\ref{DFFZ3}),  
\begin{equation*}
\Delta_{v,k}(z) \ = \ \sum_{i=1}^m \sum_{j=(k-1)N_i}^{kN_i-1} 
   \Delta_{v,ij}^{(k)} \cdot \varepsilon_{v,i}^{-1} \varphi_{i,(k+1)N_i-j}(z) \ .
\end{equation*}

As in Phase 4, the $\vartheta_{v,ij}^{(k)}(z)$ 
are $K_v$-symmetric.\index{compensating functions $\vartheta_{v,ij}^{(k)}(z)$!are $K_v$-symmetric}
\index{$K_v$-symmetric!set of functions} 
Each $\vartheta_{v,ij}^{(k)}(z)$ 
\index{compensating functions $\vartheta_{v,ij}^{(k)}(z)$!poles and leading coefficients of}
has a pole of order $nN_i-j$ at $x_i$ with leading
coefficient $\tc_{v,i}^{n-k-1}$, and poles of order at most $(n-k-1)N_{i^{\prime}}$ 
\index{coefficients $A_{v,ij}$!middle} 
at each $x_{i^{\prime}} \ne x_i$;   
furthermore  $\|\vartheta_{v,ij}^{(k)}\|_{U_v^0} \ \le \ M_v$.  
$G_v^{(k)}(z)$\index{patching functions, $G_v^{(k)}(z)$ for $1 \le k \le n$!are $K_v$-rational} is $K_v$-rational by the $K_v$-symmetry 
of the $\Delta_{v,ij}^{(k)}$ and $\vartheta_{v,ij}^{(k)}(z)$,
\index{compensating functions $\vartheta_{v,ij}^{(k)}(z)$!are $K_v$-symmetric} 
and for each $(i,j)$ 
\begin{equation*}
\sum_{x_{i^{\prime}} \in \Aut_c(\CC_v/K_v)(x_i)} 
   \Delta_{v,i^{\prime}j}^{(k)} \vartheta_{v,i^{\prime}j}^{(k)} 
\end{equation*}
is $K_v$-rational.  

\smallskip
When $\Char(K_v) = p > 0$, let  
\begin{equation*} 
F_{v,k}(z) \ = \ \Qbar_{k_1}(z) \cdot \hP_k^{\Diamond}(\phi_v(z)) \ .
\end{equation*}
By arguments similar to those before, $F_{v,k}(z)$ is $K_v$-rational, its roots belong in $E_v$, 
and it has a pole of order $(n-k-1)N_i$
with leading coefficient $d_{v,i} = \varepsilon_{v,i} \tc_{v,i}^{n-k-1}$ at each $x_i$.
\index{coefficients $A_{v,ij}$!leading}   
In particular $|d_{v,i}|_v = |\tc_{v,i}^{n-k-1}|_v$.
Thus $F_{v,k}(z)$ satisfies the conditions of Theorem \ref{DCPCPatch1p}.

By the hypotheses of Theorem \ref{DCPCPatch1p} we are given a $K_v$-symmetric collection of numbers
\index{$K_v$-symmetric!set of numbers}\index{band!$\Band_N(k)$} 
$\{\Delta_{v,ij} \in L_{w_v}\}_{(i,j) \in \Band_N(k)}$\index{distinguished place $w_v$} such that the function 
\begin{equation*}
\Delta_{v,k}(z) \ = \
 \sum_{i=1}^m \sum_{j=(k-1)N_i}^{kN_i-1} \Delta_{v,ij}^{(k)} \cdot \varphi_{i,(k+1)N_i-j}(z) \ ,
\end{equation*} 
is $K_v$-rational.  
We patch $G_v^{(k-1)}(z)$ by setting\index{patching functions, $G_v^{(k)}(z)$ for $1 \le k \le n$!constructed by patching} 
\begin{eqnarray*} 
G_v^{(k)}(z) & = & G_v^{(k-1)}(z) + \Delta_{v,k}(z) \cdot F_{v,k}(z) \\
& = &  Q_1^*(z) \cdots Q_{k-1}^*(z) \cdot \big(P_k(\phi_v(z)) + \Delta_{v,k}(z) \big) \cdot  \hP_k(\phi_v(z)) 
\end{eqnarray*}

\smallskip
The two patching constructions differ only in the choice of $\Delta_{v,k}(z)$, and in both cases  
$\|\Delta_{v,k}\|_{U_v^0} \le h_v^{kN} \cdot M_v \le q^{-\frac{k+1}{q-1} - 3 A_4 \log_v(n)}$.  Let 
\begin{equation*}
Q_k^{\Diamond *}(z) \ = \ Q_k^{\Diamond}(z) + \Delta_{v,k}(z)   \ .
\end{equation*}
In passing from $G_v^{(k-1)}(z)$ to $G_v^{(k)}(z)$,
\index{patching functions, $G_v^{(k)}(z)$ for $1 \le k \le n$!for nonarchimedean $K_v$-simple sets!movement of roots} 
the change is isolated 
in the factor $Q_k^{\Diamond}(z)$, and we have\index{patching functions, $G_v^{(k)}(z)$ for $1 \le k \le n$!factorization of}  
\begin{equation*}
G_v^{(k)}(z) \ = \Qbar_{k_1}(z) \cdot Q_k^{\Diamond *}(z)
                            \cdot \hP_k^{\Diamond}(\phi_v(z)) \ . 
\end{equation*}

When the roots of $Q_k^{\Diamond *}(z)$ 
in $B(\theta_h,\rho_h)$ are pulled back to $D(0,1)$ 
using  $\hsigma_h(Z)$, then Lemma \ref{DLemC6},\index{Refined Patching Lemma} applied with $r \le 3$ and 
\begin{equation*}
M \ = \ -\log(h_v^{kN} \cdot M_v) -\frac{k+1}{q-1} - (r-1)) \ \ge \ T,
\end{equation*}
shows they form a union of at most three $\psi_v$-regular sequences in $\cO_{u_h}$ 
\index{regular sequence!$\psi_v$-regular sequence}
attached to $\cS^{\Diamond}[k]$ and  
\begin{equation*}
\ord_v(\alpha_{hj}^*-\alpha_{hj}) \ \ge \ T \ = \ A_4 \log_v(n) 
\end{equation*}
for each $j$.  Hence the roots of $Q_k^{\Diamond *}(z)$
have not moved closer\index{roots!move roots} 
to any of the other roots of $G_v^{(k-1)}(z)$,
\index{patching functions, $G_v^{(k)}(z)$ for $1 \le k \le n$!for nonarchimedean $K_v$-simple sets!roots are separated} 
and the induction can continue.
\index{local patching for nonarch $K_v$-simple sets!Phase 5: patch unpatched indices|)}      

\vskip .1 in
\noindent{\bf Phase 6.  Completing the patching process}
\index{local patching for nonarch $K_v$-simple sets!Phase 6: complete the patching|(}

\vskip .05 in
We have now obtained a function $G_v^{(k_3)}(z)$
\index{patching functions, $G_v^{(k)}(z)$ for $1 \le k \le n$!for nonarchimedean $K_v$-simple sets!movement of roots} 
whose roots have all been patched.\index{roots!patched}  When the roots in each ball 
$B(\theta_h,\rho_h)$ are pulled back to $D(0,1)$ using $\hsigma_h(Z)$,
they form a union of at most $r := k_1 + 4$  $\psi_v$-regular sequences
\index{regular sequence!$\psi_v$-regular sequence}
in $\cO_{u_h}$, with total length $n$.  
These roots are logarithmically separated\index{logarithmically separated} 
from each other by\index{roots!logarithmically separated}
at least $T = A_4 \log_v(n)$.

We must now include all the roots in the patching process.  
To be able to apply Lemma \ref{DLemC6},\index{Refined Patching Lemma} we need that for all $k > k_3$
\begin{equation} \label{DFRUG}
h_v^{kN} \cdot M_v \ \le \ q_v^{-\frac{n}{q-1} - (k_1+4) T} \ .
\end{equation}
However $k_3$ is quite large:  
$k_3  = \#(\cS^{\Diamond})-3 \ \ge \ n - A_2 (\log_v(n))^2 - 3$. 
Thus (\ref{DFRUG}) will hold if 
\begin{equation} \label{FnCond} 
n \ \ge \ A_6 \cdot (\log_v(n))^2
\end{equation}
for a suitable constant $A_6$, which we henceforth assume.  
 
\smallskip 
When $\Char(K_v) = 0$, for $k = k_3 + 1, \ldots, n-1$ we patch as follows.
For each $(i,j)$ with $1 \le i \le m$, $(k-1)N_i \le j < kN_i$, 
put\index{compensating functions $\vartheta_{v,ij}^{(k)}(z)$!construction of}
\begin{equation} \label{DBL3}
\vartheta_{v,ij}^{(k)}(z) \ = \  \varphi_{i,(k+1)N_i-j}(z) 
           \cdot \prod_{j=1}^{n-k-1}(\phi_v(z)-\psi_v(j)) \ .             
\end{equation}  
The $\vartheta_{v,ij}^{(k)}$ are $K_v$-symmetric.
\index{compensating functions $\vartheta_{v,ij}^{(k)}(z)$!are $K_v$-symmetric}
\index{$K_v$-symmetric!set of functions}   
It is easy to see that each $\vartheta_{v,ij}^{(k)}$ 
has a pole of order $nN_i - j$ at $x_i$
with leading coefficient $\tc_{v,i}^{n-k-1}$, and has a pole of order 
\index{coefficients $A_{v,ij}$!leading} 
at most $(n-k-1)N_{i^{\prime}}$ at each $x_{i^{\prime}} \ne x_i$, 
with  $\|\vartheta_{v,ij}^{(k)}(z)\|_{U_v^0} \le M_v$.

By Theorem \ref{DCPCPatch} we are given a $K_v$-symmetric 
\index{$K_v$-symmetric!set of numbers}\index{band!$\Band_N(k)$}
set of numbers $\{\Delta_{v,ij}^{(k)} \in L_{w_v}\}_{(i,j) \in \Band_N(k)}$,\index{distinguished place $w_v$}
determined recursively in $\prec_N$ order,\index{order!$\prec_N$} 
such that $|\Delta_{v,ij}^{(k)}|_v \le h_v^{kN}$ for each $i, j$.  
Put\index{compensating functions $\vartheta_{v,ij}^{(k)}(z)$}
\index{patching functions, $G_v^{(k)}(z)$ for $1 \le k \le n$!constructed by patching}  
\begin{equation*} 
G_v^{(k)}(z) \ = \ G_v^{(k-1)}(z)
          + \sum_{i=1}^m \sum_{j=(k-1)N_i+1}^{kN_i} \Delta_{v,ij}^{(k)} \vartheta_{v,ij}^{(k)}(z) \ .
\end{equation*}
Here $G_v^{(k)}(z)$ is $K_v$-rational\index{patching functions, $G_v^{(k)}(z)$ for $1 \le k \le n$!are $K_v$-rational} 
by the $K_v$-symmetry\index{compensating functions $\vartheta_{v,ij}^{(k)}(z)$!are $K_v$-symmetric}
of the $\Delta_{v,ij}^{(k)}$ and $\vartheta_{v,ij}^{(k)}(z)$, 
and for each $(i,j)$ 
\begin{equation*}
\sum_{x_{i^{\prime}} \in \Aut_c(\CC_v/K_v)(x_i)} 
   \Delta_{v,i^{\prime}j}^{(k)} \vartheta_{v,i^{\prime}j}^{(k)} \ \in \ K_v(\cC) \ .
\end{equation*}
Furthermore $h_v^{kN} M_v \le q^{-\frac{n}{q-1} - rT}$, 
so Lemma \ref{DLemC6}\index{Refined Patching Lemma} shows that 
the roots of $G_v^{(k)}(z)$\index{patching functions, $G_v^{(k)}(z)$ for $1 \le k \le n$!roots are confined to $E_v$}
 belong to $E_v$ and have the 
same separation properties as those of $G_v^{(k-1)}(z)$.
\index{patching functions, $G_v^{(k)}(z)$ for $1 \le k \le n$!for nonarchimedean $K_v$-simple sets!roots are separated}   
 
When $k = n$, we are given a $K_v$-symmetric set of numbers 
\index{$K_v$-symmetric!set of numbers}
$\{\Delta_{v,\lambda}\}_{1 \le \lambda \le \Lambda}$ 
with $|\Delta_{v,\lambda}|_v \le h_v^{nN}$ 
for each $\lambda$.  Put\index{patching functions, $G_v^{(k)}(z)$ for $1 \le k \le n$!constructed by patching}
\begin{equation*}
G_v^{(n)}(z) \ = \ G_v^{(n-1)}(z)
       + \sum_{\lambda = 1}^{\Lambda} \Delta_{v,\lambda} \varphi_{\lambda} \ .
\end{equation*}
Clearly $G_v^{(n)}(z)$ is $K_v$-rational.\index{patching functions, $G_v^{(k)}(z)$ for $1 \le k \le n$!are $K_v$-rational} 
Since $\|\varphi_{\lambda}\|_{U_v^0} \le M_v$ for each $\lambda$, 
Lemma \ref{DLemC6}\index{Refined Patching Lemma} shows that the roots of $G_v^{(n)}(z)$ belong to $E_v$ and have the 
same separation properties as those of $G_v^{(n-1)}(z)$.
\index{patching functions, $G_v^{(k)}(z)$ for $1 \le k \le n$!for nonarchimedean $K_v$-simple sets!roots are separated} 
In particular, they are distinct.

\smallskip
When $\Char(K_v) = p > 0$, for $k = k_3 + 1, \ldots, n-1$ we patch as follows.
For each $k$, put 
\begin{equation} \label{DBL3P}
F_{v,k}(z) \ = \ \prod_{j=1}^{n-k-1}(\phi_v(z)-\psi_v(j)) \ .             
\end{equation} 
Then  $F_{v,k}(z)$ is $K_v$-rational, its roots belong in $E_v$, 
and it has a pole of order $(n-k-1)N_i$
with leading coefficient $d_{v,i} = \tc_{v,i}^{n-k-1}$ at each $x_i$.  
\index{coefficients $A_{v,ij}$!leading} 
In particular $|d_{v,i}|_v = |\tc_{v,i}^{n-k-1}|_v$.
Thus $F_{v,k}(z)$ satisfies the conditions of Theorem \ref{DCPCPatch1p}.

By the hypotheses of Theorem \ref{DCPCPatch1p} we are given a $K_v$-symmetric collection of numbers 
\index{$K_v$-symmetric!set of numbers}\index{band!$\Band_N(k)$}
$\{\Delta_{v,ij} \in L_{w_v}\}_{(i,j) \in \Band_N(k)}$\index{distinguished place $w_v$} such that the function 
\begin{equation*}
\Delta_{v,k}(z) \ = \
 \sum_{i=1}^m \sum_{j=(k-1)N_i}^{kN_i-1} \Delta_{v,ij}^{(k)} \cdot \varphi_{i,(k+1)N_i-j}(z) \ ,
\end{equation*} 
is $K_v$-rational.  
We patch $G_v^{(k-1)}(z)$ by setting\index{patching functions, $G_v^{(k)}(z)$ for $1 \le k \le n$!constructed by patching} 
\begin{equation*} 
G_v^{(k)}(z) \ = \ G_v^{(k-1)}(z) + \Delta_{v,k}(z) \cdot F_{v,k}(z) 
\end{equation*}
Here $G_v^{(k)}(z)$\index{patching functions, $G_v^{(k)}(z)$ for $1 \le k \le n$!are $K_v$-rational} is $K_v$-rational 
since $\Delta_{v,k}(z)$ and $F_{v,k}(z)$ are $K_v$-rational.  
Since $\|\Delta_{v,k}\|_{U_v^0} \le h_v^{kN} M_v \le q^{-\frac{n}{q-1} - rT}$ 
and $\|F_{v,k}\|_{U_v^0} \le 1$,  
by Lemma \ref{DLemC6}\index{Refined Patching Lemma} the roots of $G_v^{(k)}(z)$ belong to $E_v$ and have the 
same separation properties as those of $G_v^{(k-1)}(z)$.
\index{patching functions, $G_v^{(k)}(z)$ for $1 \le k \le n$!for nonarchimedean $K_v$-simple sets!roots are separated} 

When $k = n$, by the hypotheses of Theorem \ref{DCPCPatch1p} 
we are given a $K_v$-symmetric set of numbers 
\index{$K_v$-symmetric!set of numbers}
$\{\Delta_{v,\lambda}\}_{1 \le \lambda \le \Lambda}$ 
with $|\Delta_{v,\lambda}|_v \le h_v^{nN}$ 
for each $\lambda$, such that 
\begin{equation*} 
\Delta_{v,n}(z) \ = \ \sum_{\lambda = 1}^{\Lambda} \Delta_{v,\lambda} \varphi_{\lambda}(z) 
\end{equation*}
is $K_v$-rational.  Put\index{patching functions, $G_v^{(k)}(z)$ for $1 \le k \le n$!constructed by patching} 
\begin{equation*}
G_v^{(n)}(z) \ = \ G_v^{(n-1)}(z) + \Delta_{v,n}(z) \ .
\end{equation*}
Clearly $G_v^{(n)}(z)$ is $K_v$-rational.\index{patching functions, $G_v^{(k)}(z)$ for $1 \le k \le n$!are $K_v$-rational} 
Since $\|\varphi_{\lambda}\|_{U_v^0} \le M_v$ for each $\lambda$, 
Lemma \ref{DLemC6}\index{Refined Patching Lemma} shows that the roots of $G_v^{(n)}(z)$ belong to $E_v$ and have the 
same separation properties as those of $G_v^{(n-1)}(z)$.
\index{patching functions, $G_v^{(k)}(z)$ for $1 \le k \le n$!for nonarchimedean $K_v$-simple sets!roots are separated}  
In particular, they are distinct.

\smallskip 
The final assertion in Theorems \ref{DCPCPatch} and \ref{DCPCPatch1p} is that 
\begin{equation*}
\{z \in \cC_v(\CC_v) : G_v^{(n)}(z) \in \cO_v \cap D(0,r_v^{nN} \} \ \subseteq \ E_v \ .
\end{equation*} 
To assure this, we must assume that\index{patching functions, $G_v^{(k)}(z)$ for $1 \le k \le n$!mapping properties of} 
\begin{equation} \label{rvnNIneq}
r_v^{nN} \ < \ q_v^{-\frac{n}{q-1} - (k_1 + 4)T} \ , 
\end{equation} 
which holds for all sufficiently large $n$ since $r_v^N < q_v^{-1/(q-1)}$ 
by (\ref{hvrvIneq1}) and (\ref{hvrvIneq2}).

Given (\ref{rvnNIneq}), for each $h = 1, \ldots, N$ 
when we restrict $G_v^{(n)}(z)$ to $B(\theta_h,\rho_h)$\index{patching functions, $G_v^{(k)}(z)$ for $1 \le k \le n$!mapping properties of}
and pull it back to $D(0,1)$ using $\hsigma_h(Z)$, 
the Refined Patching Lemma (Lemma \ref{DLemC6})\index{Refined Patching Lemma}
shows that for each $\kappa_v \in \cO_v$ with $|\kappa_v|_v \le r_v^{nN}$, the function
$G_v^{(n)}(\hsigma_h(Z)) - \kappa_v$ has $n$ distinct roots in $\cO_{u_h}$.
Correspondingly $G_v^{(n)}(z) = \kappa_v$ has $n$ distinct solutions in\index{patching functions, $G_v^{(k)}(z)$} 
$\cC_v(F_{u_h}) \cap B(\theta_h,\rho_h)$.  Since there are $N$ balls $B(\theta_h,\rho_h)$ 
and $G_v^{(n)}(z)$ has degree $nN$, this accounts for all the solutions to  $G_v^{(n)}(z) = \kappa_v$
in $\cC_v(\CC_v)$.  It follows that\index{patching functions, $G_v^{(k)}(z)$ for $1 \le k \le n$!mapping properties of} 
\begin{equation*}
\{z \in \cC_v(\CC_v) : G_v^{(n)}(z) \in \cO_v \cap D(0,r_v^{nN}) \} 
\ \subset \ \bigcup_{h=1}^N \cC_v(F_{u_h}) \cap B(\theta_h,\rho_h) \ \subseteq \ E_v \ .
\end{equation*}  
\index{local patching for nonarch $K_v$-simple sets!Phase 6: complete the patching|)}

\smallskip

This completes the proof of Theorems \ref{DCPCPatch} and \ref{DCPCPatch1p}, 
subject to the proofs of the three Moving Lemmas below.   


\section{ Proofs of the Moving Lemmas } \label{MovingLemmaSection} 

In this section we give the proofs of the three Moving Lemmas.  
\index{patching construction!for nonarchimedean $K_v$-simple sets!proofs of the moving lemmas|(}
Our notation 
and assumptions are the same as in \S\ref{Chap11}.\ref{NonArchPatchingProof}.  For the convenience
of the reader, before giving the proofs we restate the lemmas making the  
hypotheses more explicit.  
\vskip .1 in

\noindent{\bf Lemma $\text{\bf \ref{MovingLemma1}}$. (First Moving Lemma)}\index{First Moving Lemma|ii} 
{\it Let $E_v$ and $\phi_v(z)$ be as in Theorems $\ref{DCPCPatch}$ and $\ref{DCPCPatch1p}:$ 
$E_v$ has the $K_v$-simple decomposition 
\index{$K_v$-simple!set}
\index{$K_v$-simple!decomposition}
$E_v = \bigcup_{\ell = 1}^D \big( B(a_\ell,r_\ell) \cap \cC_v(F_{w_\ell}) \big)$ 
such that $U_v := \bigcup_{\ell=1}^D  B(a_{\ell},r_{\ell})$ is disjoint from $\fX$, 
and $H_v := \phi_v^{-1}(D(0,1))$ has the 
$K_v$-simple decomposition $H_v = \bigcup_{h=1}^N \big(B(\theta_h,\rho_h) \cap \cC_v(F_{u_h})\big)$
\index{$K_v$-simple!set}\index{$K_v$-simple!decomposition}
\index{$K_v$-simple!decomposition!compatible with another decomposition}  
which is compatible with $\bigcup_{\ell = 1}^D \big( B(a_\ell,r_\ell) \cap \cC_v(F_{w_\ell}) \big)$
and move-prepared with respect to $B(a_1,r_1), \ldots, B(a_D,r_D)$. 
\index{move-prepared} 
For each $\ell$, there is a point 
$\wbar_\ell \in \big(\cC_v(F_{w_\ell}) \cap B(a_\ell,r_\ell)\big) \backslash H_v$.
   
Let $S^{\triangle} = \cS_1 \cup \ldots \cup \cS_{k_1}$, 
so the set of  patched roots\index{roots!patched} from Phases $1$ and $2$ is 
$\{\theta_{hj}^*\}_{1 \le h \le N, j \in S^{\triangle}}$.  
Let $j_0$ be as in Phase $3$, so 
$\{\theta_{h j_0} \in \cC_v(F_{u_h}) \cap B(\theta_h, \rho_h)\}_{1 \le h \le N}$ 
is a set of ``safe'' roots.\index{roots!safe}  
 
Then there are constants $\varepsilon_1 > 0$ and $C_1, C_2 \ge 1$ 
$($depending on $\phi_v(z)$, $E_v$, $H_v$, and their $K_v$-simple decompositions$)$, 
\index{$K_v$-simple!set}
\index{$K_v$-simple!decomposition}
with the following property$:$ 

Put $\delta_n = q_v^{-\lceil \log_v(n) \rceil}$. 
Then for each $0 < \varepsilon < \varepsilon_1$ small enough that 
\begin{equation*} 
C_1  \varepsilon \ \le \ \delta_n \cdot \min_{1 \le h \le N}(\rho_h) \ ,
\end{equation*}  
given any $K_v$-symmetric set 
\index{$K_v$-symmetric!set of numbers}
$\{\theta_{hj}^{**} \in \cC_v(F_{u_h}) \cap B(\theta_h,\rho_h) \}_{1 \le h \le N, j \in S^{\triangle}}$
with $\|\theta_{hj}^{**},\theta_{hj}^*\|_v < \varepsilon$ for all $(h,j)$,  
there is a $K_v$-symmetric collection of points 
\index{$K_v$-symmetric!set of points}
$\{\theta_{h,j_0}^{**} \in \cC_v(F_{u_h}) \cap B(\theta_h,\rho_h) \}_{1 \le h \le N}$ 
satisfying  
\begin{equation*}  
\|\theta_{h,j_0}^{**}, \theta_{h,j_0}\|_v \ \le \ C_1 \varepsilon \ \le \ \delta_n \rho_h 
\end{equation*}
for each $h$, such that 

$(A)$ The divisor 
\begin{equation} \label{cDForm}
\cD \ = \ \sum_{j \in S^{\triangle}} 
                \sum_{h=1}^N ((\theta_{hj}^{**}) - (\theta_{hj}^*))
        + \sum_{h=1}^N ((\theta_{h,j_0}^{**})-(\theta_{h,j_0}))
\end{equation} 
is $K_v$-rational and principal;

$(B)$ Writing $U_v^0 = \bigcup_{h=1}^N B(\theta_h,\rho_h) \subset U_v$ as before, we have 
\begin{enumerate}
\item $|Y(z)|_v = 1$ for all $z \in \cC_v(\CC_v) \backslash U_v^0$;
\item $|Y(z)-1|_v \le C_2 \varepsilon$ for all $z \in \cC_v(\CC_v) \backslash U_v^0$. 
\end{enumerate} 
}  
\medskip
\noindent{\bf Remark.}  Although the statement of this lemma is rather technical, 
it is a deep result which depends on the theory of the Universal Function\index{universal function}
developed in Appendix \ref{AppC} and the local action of the Jacobian\index{local action of the Jacobian} 
studied in Appendix \ref{AppD}.
It is the key to the local patching construction for nonarchimedean $K_v$-simple sets. 
\index{$K_v$-simple!set} 

\begin{proof}  Consider the $K_v$-simple decompositions 
\index{$K_v$-simple!decomposition}
$E_v = \bigcup_{\ell = 1}^{D} \big(B(a_\ell,r_\ell) \cap \cC_v(F_{w_\ell})\big)$
and $H_v = \bigcup_{h=1}^N \big(B(\theta_h,\rho_h) \cap \cC_v(F_{u_h})\big)$. 
Since these decompositions 
are compatible,\index{$K_v$-simple!decomposition!compatible with another decomposition}   
we have $F_{u_h} = F_{w_\ell}$ for each $h$ and $\ell$
such that $B(\theta_h,\rho_h) \subset B(a_\ell,r_\ell)$.

\smallskip
We first reduce the Lemma to a similar assertion for a single ball $B(a_\ell,r_\ell)$.  
For each $\ell = 1, \ldots, D$, 
let $\cI_{\ell}$ be the set of indices\index{indices} $1 \le h \le N$ 
such that $B(\theta_h,\rho_h) \subseteq B(a_\ell,r_\ell)$.
Suppose that for each $\ell$ there are constants 
$\varepsilon_1^{(\ell)} > 0$ and $C_1^{(\ell)}, C_2^{(\ell)} \ge 1$
such that if $0 < \varepsilon < \varepsilon_1^{(\ell)}$ and 
\begin{equation*} 
C_1^{(\ell)} \cdot \varepsilon \ \le \ \delta_n \cdot \min_{h \in \cI_\ell}(\rho_h) \ ,
\end{equation*} 
then given any set of points 
$\{\theta_{hj}^{**} \in \cC_v(F_{w_\ell}) \cap B(\theta_h,\rho_h)\}_{h \in \cI_\ell, j \in \cS^{\triangle}}$
satisfying $\|\theta_{hj}^{**},\theta_{hj}\|_v \le \varepsilon$ for all $h, j$, 
there are points
$\{\theta_{h, j_0}^{**} \in \cC_v(F_{w_\ell}) \cap B(\theta_h,\rho_h)\}_{h \in \cI_{\ell}}$
with $\|\theta_{h,j_0}^{**},\theta_{h,j_0}\|_v \le C_1^{(\ell)} \varepsilon \le \delta_n \rho_h$
such that the $F_{w_\ell}$-rational divisor
\begin{equation} \label{YellProps}
\cD_\ell \ = \ \sum_{h \in \cI_\ell, j \in S^{\triangle}} \big((\theta_{hj}^{**})-(\theta_{hj}^{*})\big) 
                      \ + \ \sum_{h \in \cI_\ell} \big((\theta_{h, j_0}^{**})-(\theta_{h, j_0})\big)  
\end{equation}
is principal, and if we put $U_{v,\ell}^0 = \bigcup_{h \in \cI_\ell} B(\theta_h,\rho_h)$, 
then there is an $F_{w_\ell}$-rational function $Y_\ell(z)$ with divisor $\cD_\ell$
such that 

\vskip .03in

\quad $(1_\ell)$ \ $|Y_\ell(z)|_v = 1$ for all $z \in \cC_v(\CC_v) \backslash U_{v,\ell}^0$\ ;

\quad $(2_\ell)$ \ $|Y_\ell(z)-1|_v \le C_2^{(\ell)} \varepsilon$ 
for all $z \in \cC_v(\CC_v) \backslash U_{v,\ell}^0$ \ . 

\vskip .03in

Put $\varepsilon_1 = \min_\ell(\varepsilon_1^{(\ell)})$, $C_1 = \max_\ell(C_1^{(\ell)})$
and $C_2 = \max_\ell(C_2^{(\ell)})$.  Let $\rhobar = \min_{1 \le h \le N}(\rho_h)$. 
Take $0 < \varepsilon \le \varepsilon_1$ small enough that 
$C_1 \varepsilon \le \delta_n \rhobar$, and let
$\{\theta_{hj}^{**} \in \cC_v(F_{u_h}) \cap B(\theta_h,\rho_h) \}_{1 \le h \le N, j \in S^{\triangle}}$
be a $K_v$-symmetric set of points satisfying $\|\theta_{hj}^{**},\theta_{hj}^*\|_v < \varepsilon$
\index{$K_v$-symmetric!set of points}
for all $h, j$.   

We now construct the $K_v$-symmetric set of points 
\index{$K_v$-symmetric!set of points}
$\{\theta_{h,j_0}^{**} \in \cC_v(F_{u_h}) \cap B(\theta_h,\rho_h) \}_{1 \le h \le N}$, 
the divisor $\cD$, and the function $Y(z)$ in the Lemma, by using galois equivariance:
we keep the divisors $\cD_\ell$ and functions $Y_\ell(z)$ for a set of representatives 
of the galois orbits for the balls $B(a_\ell,r_\ell)$, then throw away the others 
and replace them with the galois conjugates for the representatives.  
 
To be precise, write $\tK_v^{\sep}$ for the maximal separable extension of $K_v$.  
Since $E_v = \bigcup_{\ell=1}^{D} \big(\cC_v(F_{w_\ell}) \cap B(a_\ell,r_\ell) \big)$ 
is a $K_v$-simple decomposition, $F_{w_\ell}/K_v$ is separable for each $\ell$, and 
\index{$K_v$-simple!decomposition}
the orbit of $B(a_\ell,r_\ell)$ under $\Gal(\tK_v^{\sep}/K_v)$ 
has exactly $d_\ell := [F_{w_\ell}:K_v]$ elements; this means there is an action 
of $\Gal(\tK_v^{\sep}/K_v)$ on the index set $\{\ell \in \NN: 1 \le \ell \le D\}$ such that 
$B(a_{\sigma(\ell)},r_{\sigma(\ell)}) = \sigma(B(a_\ell,r_\ell))$ 
and $F_{w_{\sigma(\ell)}} = \sigma(F_{w_\ell})$ for each $\ell$ and each $\sigma \in \Gal(\tK_v^{\sep}/K_v)$.
Similarly, since $H_v = \bigcup_{h=1}^N \big(\cC_v(F_{u_h}) \cap B(\theta_h,\rho_h)\big)$ 
is a $K_v$-simple decomposition, there is an action of 
\index{$K_v$-simple!decomposition}
of $\Gal(\tK_v^{\sep}/K_v)$ on $\{h \in \NN: 1 \le h \le N\}$ such that 
$B(\theta_{\sigma(h)},\rho_{\sigma(h)}) = \sigma(B(\theta_h,\rho_h))$ 
and $F_{u_{\sigma(h)}} = \sigma(F_{u_h})$ for each $h$ and $\sigma$.  
The fact that the $\theta_{hj}^{**}$ are $K_v$-symmetric implies
\index{$K_v$-symmetric!set of functions} 
that $\theta_{\sigma(h),j}^{**} = \sigma(\theta_{hj}^{**})$ for all $h$, $\sigma$.
The compatibility of the decompositions of $H_v$ and $E_v$ means that 
$\cI_{\sigma(\ell)} = \{\sigma(h) : h \in \cI_\ell\}$ for each $\sigma$, and that 
for each $h$ such that $B(\theta_h,\rho_h) \subseteq B(a_\ell,r_\ell)$, 
we have $\sigma(h) = h $ if and only if $\sigma(\ell) = \ell$.


Let $\cL = \{\ell_1, \ldots, \ell_r\}$ be a set of representatives for the distinct galois orbits
of the balls $B(a_1,r_1), \ldots, B(a_D,r_D)$.     
For each $\ell_k \in \cL$, we have $F_{u_h} = F_{w_{\ell_k}}$ for all $h \in \cI_{\ell_k}$. Let 
\begin{equation*} 
\{\theta_{hj}^{**} \in \cC_v(F_{w_{\ell_k}}) \cap B(\theta_h,\rho_h)\}_{h \in \cI_{\ell_k}, j \in \cS^{\triangle}} 
\end{equation*} 
be the corresponding subset of  
$\{\theta_{hj}^{**} \in \cC_v(F_{u_h}) \cap B(\theta_h,\rho_h) \}_{1 \le h \le N, j \in S^{\triangle}}$.  
By hypothesis, there is a collection of points 
$\{\theta_{h, j_0}^{**} \in \cC_v(F_{w_{\ell_k}}) \cap B(\theta_h,\rho_h)\}_{h \in \cI_{\ell_k}}$
with $\|\theta_{h,j_0}^{**},\theta_{h,j_0}\|_v \le C_1^{(\ell_k)} \varepsilon \le \delta_n \rho_h$ 
such that the $F_{w_{\ell_k}}$-rational divisor 
\begin{equation*}
\cD_{\ell_k} \ = \ \sum_{h \in \cI_{\ell_k}, j \in S^{\triangle}} \big((\theta_{hj}^{**})-(\theta_{hj}^{*})\big) 
                      \ + \ \sum_{h \in \cI_{\ell_k}} \big((\theta_{h, j_0}^{**})-(\theta_{h, j_0})\big)  
\end{equation*}
is principal;  let $Y_{\ell_k}(z) \in F_{w_{\ell_k}}(\cC)$ be the corresponding function.
For an arbitrary $1 \le \ell \le D$ there are an $\ell_k \in \cL$ and a $\sigma \in \Gal(\tK_v^{\sep}/K_v)$
such that $\ell = \sigma(\ell_k)$;  redefine 
\begin{equation*}
D_\ell = \sigma(D_{\ell_k}) \ , \qquad  Y_{\ell}(z) = \sigma(Y_{\ell_k})(z) 
\end{equation*}
and redefine the $\theta_{h,j_0}^{**}$ for $h \in \cI_\ell$ by putting 
$\theta_{\sigma(h),j_0}^{**} = \sigma(\theta_{h,j_0}^{**}) \in \cC_v(F_{w_\ell}) \cap B(a_\ell,r_\ell)$ 
for each $h \in \cI_{\ell_k}$.  
By the discussion above, all of these are well-defined, the set $\{\theta_{h,j_0}^{**}\}_{1 \le h \le N}$
is $K_v$-symmetric, and for all $\ell$ and $\sigma$
\index{$K_v$-symmetric!set of functions}
we have $D_{\sigma(\ell)} = \sigma(D_\ell)$, $Y_{\sigma(\ell)}(z) = \sigma(Y_\ell)(z)$.
 
Finally, define 
\begin{equation*} 
\cD = \sum_{\ell =1}^D \cD_\ell \ , \qquad 
Y(z) = \prod_{\ell =1}^D Y_\ell(z)  \ .
\end{equation*} 
Since $\cD$ is $\tK_v^{\sep}$-rational and is fixed by $\Gal(\tK_v^{\sep}/K_v)$, it is $K_v$-rational;
similarly $Y(z)$ is $K_v$-rational.  Clearly $\div(Y(z)) = \cD$.  
By construction, $\{\theta_{h,j_0}^{**}\}_{1 \le h \le N}$
is $K_v$-symmetric, and by galois equivariance,
\index{$K_v$-symmetric!set of functions}  
$\theta_{h,j_0}^{**} \in \cC_v(F_{u_h}) \cap B(\theta_h,\rho_h)$ and $\|\theta_{h,j_0}^{**},\theta_{h,j_0}\|_v \le C_1^{(\ell)} \varepsilon \le \delta_n \rho_h$ 
for each $h$. 

Since the original sets 
$\{\theta_{hj}^{**} \}_{1 \le h \le N, j \in S^{\triangle}}$,
$\{\theta_{hj} \}_{1 \le h \le N, j \in S^{\triangle}}$,
and $\{\theta_{h,j_0}^{**}\}_{1 \le h \le N}$
were $K_v$-symmetric, $\cD$ has the form (\ref{cDForm}).
\index{$K_v$-symmetric!set of functions}    
Clearly $U_{v,\ell}^0 \subset B(a_\ell,r_\ell) \subseteq U_v$ for each $\ell$.  
For each $z \notin U_{v,\ell}^0$, we have $|Y_\ell(z)|_v = 1$, 
so for each $z \notin U_v^0 = \bigcup_{\ell=1}^D U_{v,\ell}^0$ 
we have $|Y(z)|_v = 1$, and (B1) in the Lemma holds;  
similarly, for each $z \notin U_{v,\ell}^0$, 
we have $|Y_\ell(z)-1|_v \le C_2 \varepsilon$, so since 
\begin{equation*}
Y(z) - 1 \ = \ \sum_{\ell=1}^D (Y_\ell(z)-1) \cdot \left(\prod_{k = \ell+1}^D Y_k(z) \right) \ ,
\end{equation*}
the ultrametric inequality shows that for each $z \notin U_v^0$ we have $|Y(z)-1|_v \le C_2 \varepsilon$, 
and (B2) holds.

\medskip
Now fix $\ell$;  we will construct $\cD_\ell$ and $Y_\ell(z)$ for $B(a_\ell,r_\ell)$, 
and show they satisfy properties $(1_\ell)$ and $(2_\ell)$.  
The proof has two steps:  
first we use the local action of the Jacobian, from Appendix D, 
to construct the principal divisor $\cD_\ell$; 
then we use the theory of the Universal Function,\index{universal function} 
from Appendix C, to construct $Y_\ell(z)$. 
Put $E_{v,\ell} = E_v \cap B(a_\ell,r_\ell) = \cC_v(F_{w_\ell}) \cap B(a_\ell,r_\ell)$ 
and put $H_{v,\ell} = H_v \cap B(a_\ell,r_\ell) \subset E_{v,\ell}$.  As noted above, 
we have $F_{u_h} = F_{w_\ell}$ for each $h \in \cI_\ell$,   
so $H_{v,\ell} = \bigcup_{h \in \cI_\ell} \big(\cC_v(F_{w_\ell}) \cap B(\theta_h,\rho_h)) 
= \cC_v(F_{w,\ell}) \cap U_{v,\ell}^0$.
By hypothesis, there is a point 
$\wbar_\ell \in \big(\cC_v(F_{w_\ell}) \cap B(a_\ell,r_\ell)\big) \backslash H_{v,\ell}$;  
clearly $\wbar_\ell \notin U_{v,\ell}^0$.  

\smallskip
We begin by constructing $\cD_\ell$.    

First assume that $g = g(C_v) > 0$. By hypothesis, the $K_v$-simple 
\index{$K_v$-simple!decomposition}
decomposition $H_v = \bigcup_{h=1}^N \big(\cC_v(F_{u_h}) \cap B(\theta_h,\rho_h)\big)$ is 
move-prepared relative to $B(a_1,r_1), \ldots, B(a_D,r_D)$. For simplicity,
\index{move-prepared}   
relabel the roots $\theta_1, \ldots, \theta_N$ of $\phi_v(z)$ so that
$B(\theta_1,\rho_1), \ldots, B(\theta_M,\rho_M)$ are contained in $B(a_\ell,r_\ell)$; 
thus $\cI_\ell = \{1,\ldots,M\}$.  Suppose also that $B(\theta_1,\rho_1), 
\ldots, B(\theta_g,\rho_g)$ 
are the distinguished balls corresponding to $B(a_\ell,r_\ell)$ in the 
\index{distinguished balls} 
definition of move-preparedness (Definition \ref{MovePreparedDef}).  This means there is
\index{move-prepared} 
a number $\rbar_\ell$ with $\rho_1, \ldots, \rho_g < \rbar_\ell < r_\ell$ such that    
$B(\theta_1,\rbar_\ell), \ldots, B(\theta_g,\rbar_\ell)$ are pairwise disjoint
and contained in $B(a_\ell,r_\ell)$, and if we put $\vtheta_\ell = (\theta_1, \ldots, \theta_g)$ 
then 
\begin{equation*}
W_{\vtheta_\ell}(\rbar_\ell) \ := \ \BFJ_{\vtheta_\ell} \Big(\prod_{h=1}^g B(\theta_h,\rbar_\ell) \Big) 
\end{equation*} 
is an open subgroup of $\Jac(\cC_v)(\CC_v)$ 
satisfying the properties in Theorem \ref{BKeyThm1A}.\index{local action of the Jacobian}

Since $E_{v,\ell}$ is compact, by Proposition \ref{BPropD2New}\index{Lipschitz continuity!of the Abel map}
there are constants $\varepsilon_0^{(\ell)} > 0$, $C_0^{(\ell)} > 0$ such that if $0 < \varepsilon \le \varepsilon_0^{(\ell)}$, 
then for all $x, y \in E_{v,\ell}$ with $\|x,y\|_v \le \varepsilon$, 
the divisor class $\BFj_{x}(y) = [(y)-(x)]$ 
belongs to $W_{\vtheta_\ell}( C_0^{(\ell)} \varepsilon \cdot \rbar_\ell)$.  

In the Lemma, we will take 
\begin{equation*}
\varepsilon_1^{(\ell)} = \varepsilon_0^{(\ell)} \ , \qquad C_1^{(\ell)} = \max(1,C_0^{(\ell)}\rbar_\ell) \ .
\end{equation*}
Put $\rhobar_\ell = \min_{h \in \cI_\ell} (\rho_h)$, and let 
$0 < \varepsilon \le \varepsilon_1^{(\ell)}$ be small enough that 
$C_1^{(\ell)} \varepsilon \le \delta_n \rhobar_\ell$.
Thus $B(\theta_{hj},C_1^{(\ell)} \varepsilon) \subseteq B(\theta_{hj},\delta_n \rho_h)$ 
for all $1 \le h \le M$, $1 \le j \le n$. 
The balls $B(\theta_{hj},\delta_n \rho_h)$ are pairwise disjoint
and isometrically parametrizable, hence the same is true for the balls 
$B(\theta_{hj},C_1^{(\ell)} \varepsilon)$.  Using the ``safe'' index $j_0$, 
put $\vtheta_{\ell,j_0} = (\theta_{1,j_0}, \ldots, \theta_{g,j_0})$.  
By Theorem \ref{BKeyThm1A}(D),\index{local action of the Jacobian}
\begin{equation*}
W_\ell(C_1^{(\ell)} \varepsilon) \ := \ W_{\vtheta_{\ell,j_0}}(C_1^{(\ell)} \varepsilon) 
\ = \ \BFJ_{\vtheta_{\ell,j_0}} \Big(\prod_{h=1}^g B(\theta_{h,j_0},C_1^{(\ell)} \varepsilon) \Big) 
\end{equation*} 
is an open subgroup of $W_{\vtheta_\ell}(\rbar_\ell)$.
By our choice of $C_1^{(\ell)}$ it contains  $W_{\vtheta_\ell}( C_0^{(\ell)} \varepsilon \cdot \rbar_\ell)$.

Suppose we are given a set of points  
$\{\theta_{hj}^{**} \in \cC_v(F_{w_\ell}) \cap B(\theta_h,\rho_h)\}_{1 \le h \le M, j \in \cS^{\triangle}}$ 
with $\|\theta_{hj}^{**},\theta_{hj}^*\|_v \le \varepsilon$ for each $(h,j)$.
Since $\theta_{hj}^*$ belongs to $B(\theta_h,\rho_h)$ and $\varepsilon \le \rho_h$,
we have $\theta_{hj}^{**} \in B(\theta_h,\rho_h)$ as well.  Thus, 
$\theta_{hj}^{**} \in \cC_v(F_{w_\ell}) \cap B(\theta_h,\rho_h) \subset H_{v,\ell}$ 
for each $h, j$. 
 
Using the action $\ap$ of the group $W_{\vtheta_{\ell,j_0}}(C_1^{(\ell)} \varepsilon)$ on 
$\prod_{h=1}^g B(\theta_{h,j_0},C_1^{(\ell)} \varepsilon)$  
from Theorem \ref{BKeyThm1A},\index{local action of the Jacobian} we will construct points 
$\theta_{h,j_0}^{**} \in \cC_v(F_{w_\ell}) \cap B(\theta_{h,j_0},C_1^{(\ell)} \varepsilon)$, 
for $h = 1, \ldots, g$, such that
\begin{equation*} 
\cD_\ell \ := \ \sum_{j \in S^{\triangle}} 
                \sum_{h \in \cI_\ell} ((\theta_{hj}^{**}) - (\theta_{hj}^*))
        + \sum_{h=1}^g ((\theta_{h,j_0}^{**})-(\theta_{h,j_0}))
\end{equation*}
is $F_{w_\ell}$-rational and principal.  Consider the divisor class  
\begin{equation*}
x \ = \ \sum_{j \in S^{\triangle}} 
         \sum_{h \in \cI_\ell} [(\theta_{hj}^{**}) - (\theta_{hj}^*))] \ \in \ \Jac(\cC_v)(F_{w_\ell}) \ .
\end{equation*} 
As noted above, the $\theta_{hj}^{**}$ and $\theta_{hj}^*$ belong to $H_{v,\ell} \subset U_{v,\ell}^0$.  
By our choice of $\varepsilon$, we have 
$[(\theta_{hj}^{**}) - (\theta_{hj}^*))] \in W_{\vtheta_\ell}(C_0^{\ell} \varepsilon \cdot \rbar_\ell) 
 \subseteq W_{\vtheta_{\ell,j_0}}(C_1^{(\ell)} \varepsilon)$ for all $h, j$.  
Since $W_{\vtheta_{\ell,j_0}}(C_1^{(\ell)} \varepsilon)$ is a group, it follows that 
$x \in W_{\vtheta_{\ell,j_0}}(C_1^{(\ell)} \varepsilon) \cap \Jac(\cC_v)(F_{w_\ell})$.    
Define 
\begin{equation*}   
(\theta_{1,j_0}^{**}, \ldots, \theta_{g,j_0}^{**}) 
\ = \ (-x) \ \ap \ (\theta_{1,j_0}, \ldots, \theta_{g,j_0})
\ \in \ \prod_{h=1}^g B(\theta_{h,j_0},C_1^{(\ell)} \varepsilon)  \ .
\end{equation*}  
By Theorem \ref{BKeyThm1A}(C),\index{local action of the Jacobian} 
\begin{equation*}
 \sum_{h=1}^g [(\theta_{h,j_0}^{**})-(\theta_{h,j_0})] \ = \ 
\crJ_{\vtheta_{\ell,j_0}}((\theta_{1,j_0}^{**}, \ldots, \theta_{g,j_0}^{**})) 
\ = \ -x \ ,
\end{equation*}
so $\cD_\ell$ is principal.  By Theorem \ref{BKeyThm1A}(E),\index{local action of the Jacobian} 
the action $\ap$ preserves $F_{w_\ell}$
rationality, so each $\theta_{h,j_0}^{**} \in \cC_v(F_{w_\ell})$, 
and $\cD_\ell$ is $F_{w_\ell}$-rational.  Finally, our choice of $\varepsilon$
required that $C_1^{(\ell)} \varepsilon \le \delta_n \rho_h$, 
so $\theta_{h,j_0}^{**} \in \cC_v(F_{w_\ell}) \cap B(\theta_{h,j_0},\delta_n \rho_h)$ 
for each $h = 1, \ldots, g$.

For $h = g+1, \ldots, M$, put $\theta_{h,j_0}^{**} = \theta_{h,j_0}$.    

\smallskip
We next construct $Y_\ell(z)$.  
For this, it will be useful to relabel the $\theta_{hj}^*$ and $\theta_{hj}^{**}$ 
by gathering them in groups of size $g$.   For simplicity, first assume that $g$ 
divides $M \cdot  \#(\cS^{\triangle})$. Put $T = M \cdot \#(\cS^{\triangle})/g$ and write
\begin{equation*}
\{\theta_{hj}^*\}_{1 \le h \le M, j \in \cS^{\triangle}} 
\ = \ \{c_k^{(t)}\}_{1 \le k \le g, 1 \le t \le T} \ .
\end{equation*}
Using the same correspondence between indices,\index{indices} write 
\begin{equation*}
\{\theta_{hj}^{**}\}_{1 \le h \le M, j \in \cS^{\triangle}} 
\ = \ \{\hc_k^{(t)}\}_{1 \le k \le g, 1 \le t \le T} \ .  
\end{equation*} 
Put $\vc^{(t)} = (c_1^{(t)}, \ldots, c_g^{(t)})$,   
$\hc^{(t)} = (\hc_1^{(t)}, \ldots, \hc_g^{(t)})$.  Clearly
\begin{equation} \label{BreakDownF}
 \sum_{j \in S^{\triangle}} \sum_{h=1}^M  [(\theta_{hj}^{**}) - (\theta_{hj}^*))]
\ = \ \sum_{t = 1}^T \sum_{k=1}^g [(\hc_k^{(t)})-(c_k^{(t)})] \ . 
\end{equation} 

If $g$ does not divide $M \cdot \#(\cS^{\triangle})$,
put $T = \lceil M \cdot \#(\cS^{\triangle})/g\rceil$  
and set $r = T \cdot g - M \cdot \#(\cS^{\triangle})$.
Fix an element $j_1 \in \cS^{\triangle}$ 
and augment the lists $\{\theta_{hj}^{*}\}$ and $\{\theta_{hj}^{**}\}$ by adjoining $r$ copies of
$\theta_{1,j_1}^*$ at the end of each, then break down the lists into groups of size $g$ as before.  
In this way the final vectors $\vc^{(T)}$, $ \hc^{(T)}$ have their last $r$ components equal to 
$\theta_{1,j_1}^*$, and (\ref{BreakDownF}) still holds.  

Put $\vd^{(0)} 
= \vtheta_{\ell,j_0} = (\theta_{1,j_0}, \ldots, \theta_{g,j_0})$ and write 
$\vd^{(0)} = (d_1^{(0)}, \ldots, d_g^{(0)})$.  We will inductively construct vectors 
$\vd^{(t)} = (d_1^{(t)}, \ldots, d_g^{(t)}) 
 \in \prod_{h=1}^g \big( \cC_v(F_{w_\ell}) \cap B(\theta_{h,j_0},\cC_1^{(\ell)} \varepsilon) \big)$ 
such that for each $t = 1, \ldots, T$, the divisor 
\begin{equation*} 
\cD_\ell^{(t)} \ = \ \sum_{h=1}^g ((\hc_h^{(t)}) - (c_h^{(t)}))
                 + \sum_{h=1}^g ((d_h^{(t)}) - (d_h^{(t-1)}))
\end{equation*}
is $F_{w_\ell}$-rational and principal, and such that
\begin{equation*} 
\cD_\ell \ = \ \sum_{t = 1}^T \cD_\ell^{(t)} \ .
\end{equation*}  
Suppose $\vd^{(t-1)}$ has been constructed.  Put 
\begin{equation*}
x^{(t)} \ = \ \sum_{j=1}^g [(\hc_j^{(t)}) - (c_j^{(t)})]
\ \in \ \Jac(\cC_v)(F_{w_\ell}) \cap W_{\vtheta_{\ell,j_0}}(C_1^{(\ell)}
 \varepsilon) \ . 
\end{equation*}
Using the operation $\ap$ of the group $W_{\vtheta_{\ell,j_0}}(C_1^{(\ell)} \varepsilon)$ on 
$\prod_{h=1}^g B(\theta_{h,j_0},C_1^{(\ell)} \varepsilon)$  
in Theorem \ref{BKeyThm1A},\index{local action of the Jacobian} define 
\begin{equation*}   
\vd^{(t)} \ = \ (-x^{(t)}) \ \ap \ \vd^{(t-1)} 
\ \in \ \prod_{h=1}^g B(\theta_{h,j_0},C_1^{(\ell)} \varepsilon) \ . 
\end{equation*}  
By Theorem \ref{BKeyThm1A}(C),\index{local action of the Jacobian} 
\begin{equation*} 
\sum_{j=1}^g [(d_j^{(t)}) - (d_j^{(t-1)})] \ = \ -x^{(t)} \ = \ -\sum_{j=1}^g [(\hc_j^{(t)}) - (c_j^{(t)})]
\end{equation*} 
so $\cD^{(t)}$ is principal. 
Since $C_1^{(\ell)} \varepsilon \le \delta_n \rhobar_\ell$, 
Theorem \ref{BKeyThm1A}(E)\index{local action of the Jacobian} 
shows that $\vd^{(t)}$ belongs to 
$\prod_{h=1}^g \big( \cC_v(F_{w_\ell}) \cap B(\theta_{h,j_0},\delta_n \rho_h) \big)$, 
and $\cD^{(t)}$ is $F_{w_\ell}$-rational. 

Since $x^{(1)} + \ldots + x^{(T)} = x$, the fact that $\ap$ is an action 
 assures that when $t = T$, we have
$\vd^{(T)} = (\theta_{1,j_0}^{**}, \ldots, \theta_{g,j_0}^{**})$ with the points $\theta_{h,j_0}^{**}$ 
constructed earlier.  Thus the divisor class 
\begin{equation*} 
 \sum_{t=1}^T  \Big(\sum_{j=1}^g [(d_j^{(t)}) - (d_j^{(t-1)})]\Big)   
\end{equation*} 
telescopes to $\sum_{h=1}^g [(\theta_{h,j_0}^{**})-(\theta_{h,j_0})]$, 
and  $\cD_\ell = \sum_{t=1}^T \cD_\ell^{(t)}$ as claimed.  

\vskip .05 in
If $g(\cC) = 0$, we can again assume the roots of $\phi_v(z)$ 
are labelled so that $\cI_\ell = \{1, \ldots, M\}$.  The divisor 
$\cD_\ell := \sum_{j \in \cS^{\triangle}} \sum_{h=1}^M (\theta_{hj}^{**}) - (\theta_{hj}^*)$
is already principal, so we can take $\theta_{h,j_0}^{**} = \theta_{h,j_0}$ for 
each $h = 1, \ldots, M$.  For compatibility with the notation above,
put $T = M \cdot \#(\cS^{\triangle})$, 
and relabel the sets $\{\theta_{hj}\}_{1 \le h \le M, j \in \cS^{\triangle}}$, 
$\{\theta_{hj}^{**}\}_{1 \le h \le M, j \in \cS^{\triangle}}$, 
as $\{c^{(t)}\}_{1 \le t \le T}$, $\{\hc^{(t)}\}_{1 \le t \le T}$, 
respectively.  For each $t = 1, \ldots, T$, 
put $\cD_\ell^{(t)} = (\hc^{(t)})-(c^{(t)})$.  Then each $\cD_\ell^{(t)}$
is $F_{w_\ell}$-rational, and $\cD_\ell = \sum_{t=1}^T \cD^{(t)}$. 

\smallskip
We can now construct $Y_\ell(z)$.  Recall that 
\begin{equation*}
C_1^{(\ell)} \varepsilon \ \le \ \delta_n \cdot \rhobar_\ell 
\ = \ \delta_n \cdot \min_{h \in \cI_\ell} \rho_h\ .
\end{equation*} 
For each $(h,j)$ with   $1 \le h \le M$ and $j \in \cS^{\triangle}$ 
we have $\|\theta_{hj}^{**},\theta_{hj}^{*}\|_v \le \varepsilon$
and in particular 
$\theta_{hj}^{**}  \in  B(\theta_{hj},\delta_n \rho_h)$ since $C_1^{(\ell)} \ge 1$;  
while for $j = j_0$ and $h = 1, \ldots, M$, we have 
$\|\theta_{h,j_0}^{**},\theta_{h,j_0}\|_v \le C_1^{(\ell)}\varepsilon$, hence 
$\theta_{h,j_0}^{**} \in  B(\theta_{h,j_0},\delta_n \rho_h)$. 

We now apply Theorem \ref{BMT2} of Appendix C to the set $H_{v,\ell}$, 
taking $d = \max(1,2g)$, $r = \rhobar_\ell$, 
and replacing $\varepsilon$ in the Theorem with $C_1^{(\ell)} \varepsilon \le \delta_n \rhobar_\ell$.
Let $D(H_{v,\ell},d)$ be the constant from the Theorem,
and take
\begin{equation*}
C_2^{(\ell)} \ = \ C_1^{(\ell)} \cdot \frac{D(H_{v,\ell},d)}{(\rhobar_\ell)^d} 
\end{equation*} 
in the Lemma.  Let $\wbar_\ell \in \big(\cC_v(F_{w_\ell}) \cap B(a_\ell,r_\ell)\big) \backslash H_v
= \big(\cC_v(F_{w_\ell}) \cap B(a_\ell,r_\ell)\big) \backslash U_{v,\ell}^0$ 
be the point from the statement of the Lemma.  
For each $p \in U_{v,\ell}^0$ we have $\|p,\wbar_\ell\|_v > r = \rhobar_\ell $. 

For each $t = 1, \ldots, T$, by specializing the Universal Function\index{universal function}
$f(z,w;\vp,\vq)$ of degree $d$ in Theorem \ref{UniversalFcn} 
of Appendix C, taking $w = \wbar_\ell$, and letting $\vp$ (resp. $\vq$) be 
vectors consisting of the zeros (resp. poles) of $\cD_\ell^{(t)}$,
we obtain a function $Y_\ell^{(t)}(z)$
for which $\div(Y_\ell^{(t)}) = \cD_\ell^{(t)}$ and $Y_\ell^{(t)}(\wbar_\ell) = 1$.   
Each $Y_\ell^{(t)}(z)$ is $F_{w_\ell}$-rational,
since $\cD_\ell^{(t)}$ is $F_{w_\ell}$-rational and $\wbar_\ell \in \cC_v(F_{w_\ell})$.  

The sets $(\bigcup_{j=1}^d B(p_j,r_j)^-) \cup (\bigcup_{j=1}^d B(q_j,r_j)^-))$
and $(\bigcup_{j=1}^d B(p_j,r)^-) \cup (\bigcup_{j=1}^d B(q_j,r)^-)$ from Theorem \ref{BMT2}(A,B)
are both contained in $U_{v,\ell}^0 = \bigcup_{h \in \cI_\ell} B(\theta_h,\rho_h)$.  
Hence for each $t = 1, \ldots, T,$

\smallskip
$(1_{\ell,t})$ \ \ $|Y_\ell^{(t)}(z)|_v = 1$ for all $z \in \cC_v(\CC_v) \backslash U_{v,\ell}^0$\ ;

$(2_{\ell,t})$ \ \ $|Y_\ell^{(t)}(z)-1|_v \le C_2^{(\ell)} \varepsilon$ 
for all $z \in \cC_v(\CC_v) \backslash U_{v,\ell}^0$ \ . 
\smallskip

\noindent{Put} $Y_\ell(z) = \prod_{t=1}^T Y^{(t)}(z) \in F_{w_\ell}(\cC_v)$.  
Then $\div(Y_\ell) = \cD_\ell$, 
and $|Y_\ell(z)|_v = 1$ for all $z \notin U_{v,\ell}^0$. Thus assertion $(1_\ell)$ holds.
Since   
\begin{equation*}
Y_\ell(z) - 1 \ = \ \sum_{t=1}^T \big(Y_\ell^{(t)}(z)-1\big) \cdot \prod_{s=t+1}^T Y_\ell^{(s)}(z) 
\end{equation*}
and $(1_{\ell,t})$ and $(2_{\ell,t})$ above hold for all $t$, 
the ultrametric inequality shows that $|Y_\ell(z)-1|_v \le C_2^{(\ell)} \varepsilon$ 
for all $z \notin U_{v,\ell}^0$. Thus assertion $(2_\ell)$  holds. 
\end{proof} 

\vskip .1 in

\noindent{\bf Lemma $\text{\bf \ref{MovingLemma2}}$. (Second Moving Lemma)}
\index{Second Moving Lemma|ii} 
{\it There are constants $\varepsilon_2 > 0$ and $C_3, C_4 \ge 1$ 
$($depending on $E_v$, $\fX$, 
the choices of the $L$-rational and $L^{\sep}$-rational bases, the uniformizers $g_{x_i}(z)$, 
and the projective embedding of $\cC_v),$
such that if $0 < \varepsilon < \varepsilon_1$ and $Y(z)$ 
are as in Lemma $\ref{MovingLemma1}$,
and in addition $\varepsilon < \varepsilon_2$ and $n$
is sufficiently large, then when we expand\index{patching functions, $G_v^{(k)}(z)$ for $1 \le k \le n$!constructed by patching} 
\begin{eqnarray*}
G_v^{(k_1)}(z) \ = \
      \sum_{i=1}^m \sum_{j=0}^{(n-1)N_i-1} A_{v,ij} \varphi_{i,nN_i-j}(z)
        + \sum_{\lambda = 1}^{\Lambda} A_{\lambda} \varphi_{\lambda}(z) \ , \qquad \qquad \quad \\
\Gbar_v^{(k_1)}(z) \ := \ Y(z) \cdot G_v^{(k_1)}(z) \ = \ 
      \sum_{i=1}^m \sum_{j=0}^{(n-1)N_i-1} \Abar_{v,ij} \varphi_{i,nN_i-j}(z)
        + \sum_{\lambda = 1}^{\Lambda} \Abar_{\lambda} \varphi_{\lambda}(z) \ ,
\end{eqnarray*}
for all $i = 1, \ldots, m$ and all $0 \le j < k_1 N_i$, we have  
\begin{equation*} 
 |\Abar_{v,ij}-A_{v,ij}|_v 
   \ \le \ C_3 C_4^j (|\tc_{v,i}|_v)^n \varepsilon \ .
\end{equation*} 
} 

\vskip .03 in 
For the proof, we will need the following lemma concerning power series,
which is closely related to Lemma \ref{ExpansionCoeffBound}:  

\begin{lemma} \label{BLemCPC5}
Let $r$ belong to the value group of $\CC_v$, 
and let $b \in \CC_v^{\times}$.  Suppose that for each $j \ge J_0$,
$\Phi_j(Z) = b^{-j} Z^{-j}\big(1+\sum_{\ell=1}^{\infty} 
C_{\ell}^{(j)} Z^{\ell}\big) \in \CC_v((Z))$ is a Laurent series\index{Laurent expansion} 
with leading coefficient $b^{-j}$, which converges in $D(0,r)\backslash\{0\}$ and has no zeros there.
\index{coefficients $A_{v,ij}$!leading}  
Let 
\begin{equation*}
G(Z) \ = \ Z^{-M}\big(g_0 + \sum_{j=1}^{\infty} g_j Z^{j}\big) \in \CC_v((Z))
\end{equation*} 
be another Laurent series\index{Laurent expansion} converging in $D(0,r) \backslash \{ 0 \}$
and having no zeros there, with $g_0 \ne 0$.  

We can uniquely expand $G(Z)$ as a linear combination of the $\Phi_j(Z)$ 
and a residual series in $Z$, writing
\begin{equation} \label{BFGL1}
G(Z) \ = \ \sum_{j=0}^{M-J_0} B_j \Phi_{M-j}(Z)
              + Z^{-J_0+1} \big(\sum_{j = 0}^{\infty} B_j^{\prime} Z^j \big) \ .
\end{equation}
Suppose $Y(Z) = \sum_{\ell=0}^{\infty} h_{\ell} Z^{\ell}$ is a power series
converging in $D(0,r)$, and there is an $\varepsilon$ with $0 < \varepsilon < 1$
such that $|Y(Z)-1|_v \le \varepsilon$ for all $Z \in D(0,r)$.

Consider the product $Y(Z)G(Z) = \sum_{j=0}^{\infty} \gbar_j Z^{-M+j}$.
If we expand
\begin{equation} \label{BFGL2}
Y(Z)G(Z) \ = \  \sum_{j=0}^{M-J_0} \Bbar_j \Phi_{M-j}(Z)
              + Z^{-J_0+1}\big(\sum_{j = 0}^{\infty} \Bbar_j^{\prime} Z^j \big)
\end{equation}
then $|\Bbar_j-B_j|_v \le \varepsilon \cdot |g_0|_v \cdot (r |b|_v)^{-j}$
        for each $j$ in the range $0 \le j \le M-J_0$.
\end{lemma}

\begin{proof}
After replacing $Z$ by $bZ$, and $D(0,r)$ by $D(0,r|b|_v)$, we can
assume without loss that $b=1$.  In particular, we can assume that
each $\Phi_j(Z)$ has leading coefficient $1$.
\index{coefficients $A_{v,ij}$!leading} 

Under this hypothesis, we will first show that for each $j \ge 0$,
\begin{equation*}
|\gbar_j - g_j|_v \ \le \ \varepsilon \cdot \frac{|g_0|_v}{r^j} \ .
\end{equation*}
Multiplying $G(Z)$ by $Z^M$, we obtain a power series
converging in $D(0,r)$, having no zeros in $D(0,r)$, 
whose Taylor coefficients are the $g_j$.  The theory of Newton Polygons
\index{Newton Polygon}
shows that $|g_j|_v \le |g_0|/r^j$ for all $j$ 
(see Lemma \ref{BLem1} and the discussion preceding it; in fact, strict inequality holds when $j \ge 1$).
On the other hand, by the Maximum Modulus Principle for power series,
\index{Maximum principle!nonarchimedean!for power series}
since $|Y(Z)-1|_v \le \varepsilon$ for all $Z \in D(0,r)$,
we have $|h_0-1|_v \le \varepsilon$
and $|h_{\ell}|_v \le \varepsilon/r^{\ell}$ for all $\ell \ge 1$.

In the product $Y(Z)G(Z)$ we have $\gbar_j = \sum_{k=0}^j h_k g_{j-k}$
for each $j$.  Hence
\begin{eqnarray*}
|\gbar_j - g_j|_v & = & |(h_0-1) g_j  + h_1 g_{j-1} + \ldots + h_j g_0|_v \\
   & \le & \max(|h_0-1|_v |g_j|_v, |h_1|_v |g_{j-1}|_v , \ldots,
                            |h_j|_v |g_0|_v ) \\
   & \le & \max(\varepsilon \cdot |g_0|_v/r^j,
                \varepsilon/r \cdot |g_0|_v/r^{j-1} ,  \ldots ,
                \varepsilon/r^j \cdot |g_0|_v) \\
   & = & \varepsilon \cdot |g_0|_v/r^j \ .
\end{eqnarray*}

Now consider the expansions (\ref{BFGL1}) and (\ref{BFGL2}).  Clearly
$B_0 = g_0$ and $\Bbar_0 = \gbar_0$, so $|\Bbar_0-B_0|_v \le \varepsilon |B_0|_v$.
Since the $\Phi_j(Z)$ have no zeros in $D(0,r)$ and have leading
coefficient $1$, the theory of Newton Polygons shows that
\index{coefficients $A_{v,ij}$!leading}\index{Newton Polygon}
$|C_{\ell}^{(j)}|_v \le 1/r^{\ell}$ for each $j$ and $\ell$, as before.

Suppose inductively that for some $J \le M-J_0$, we have shown that
$|\Bbar_j-B_j|_v \le \varepsilon |B_0|_v/r^j$ for all $0 \le j < J$.
Using (\ref{BFGL1}) and (\ref{BFGL2}) we have
\begin{eqnarray*}
G(Z) - \sum_{j=0}^{J-1} B_j \Phi_{M-j}(Z) 
          & = & \sum_{k=0}^{\infty} \delta_k Z^{-M+J+k} \ , \\
Y(Z)G(Z) - \sum_{j=0}^{J-1} \Bbar_j \Phi_{M-j}(Z)
          & = & \sum_{k=0}^{\infty} \deltabar_k Z^{-M+J+k} \ ,
\end{eqnarray*}
for certain numbers $\delta_k, \deltabar_k \in \CC_v$.
Inserting the Laurent expansions\index{Laurent expansion} for $G(Z)$,
$Y(Z) G(Z)$ and the $\Phi_{M-j}(Z)$, we see that
for each $k$
\begin{eqnarray*}
\delta_k & = & g_{J+k} - B_0 C_{k+J}^{(M)} - B_1 C_{k+J-1}^{(M-1)} - \cdots
                         - B_{J-1} C_{k+1}^{(M-J+1)} \ , \\
\deltabar_k & = & \gbar_{J+k} - \Bbar_0 d_{J+k}^{(M)} - \cdots
                         - \Bbar_{J-1} C_{k+1}^{(M-J+1)} \ .
\end{eqnarray*}
By the ultrametric inequality and the estimates above, 
\begin{eqnarray*}
|\deltabar_k-\delta_k|_v
        & \le & \max(\varepsilon |B_0|_v/r^{J+k},
                           \varepsilon |B_0|_v \cdot 1/r^{J+k}, \ldots,
                  \varepsilon |B_0|_v/r^{J-1} \cdot 1/r^{k+1} ) \\
        & = & \varepsilon |B_0|_v/r^{J+k} \ .
\end{eqnarray*}
When $k = 0$, the fact that $\Phi_{M-J}(Z)$ has leading term $Z^{-M+J}$
shows that $B_J = \delta_0$ and $\Bbar_J = \deltabar_0$.  Hence
$|\Bbar_j-B_j|_v \le \varepsilon |B_0|/r^J$ and the induction can continue.

When $J = M-J_0$, the induction stops because there is no 
function $\Phi_{J_0-1}(Z)$.
\end{proof}

\begin{proof}[Proof of Lemma \ref{MovingLemma2}:]  
Since $\fX$ is disjoint from $E_v$, 
and since the basis functions $\varphi_{ij}(z)$ and $\varphi_{\lambda}$ 
\index{basis!$L$-rational}
belong to a multiplicatively finitely generated set,
there is a radius $r > 0$ in the value group of $\CC_v^{\times}$
such that 

\quad $(1)$  $r < \min_{i \ne j} (\|x_i,x_j\|_v);$ 

\quad $(2)$ each of the balls $B(x_i,r)$ is isometrically parametrizable and disjoint from $E_v$; 

\quad $(3)$ for each $i$, none of the $\varphi_{ij}(z)$ has a zero in $B(x_i,r)$;.

Fixing $x_i \in \fX$, let $\varrho_i : D(0,r) \rightarrow B(x_i,r)$
be an $L_{w_0}$-rational isometric parametrization of $B(x_i,r)$
with $\varrho_i(0) = x_i$.  If $Z$ is the coordinate on $D(0,r)$, 
then we can expand $G_v^{(k_1)}(z)$, $Y(z)$,\index{patching functions, $G_v^{(k)}(z)$ for $1 \le k \le n$!expansion of} 
and the $\varphi_{ij}(z)$
as Laurent series\index{Laurent expansion} in $Z$ converging in $D(0,r) \backslash \{0\}$, 
putting $G(Z) = G_v^{(k_1)}(\varrho_i(Z))$,\index{patching functions, $G_v^{(k)}(z)$ for $1 \le k \le n$!expansion of} 
$H(Z) = Y(\varrho_i(Z))$, and $\Phi_j(z) = \varphi_{ij}(\varrho_i(Z))$.

With respect to the uniformizer $g_{x_i}(z)$ the $\varphi_{ij}(z)$ are monic.
That is,
\begin{equation*}
\lim_{z \rightarrow x_i} \varphi_{ij}(z) \cdot g_{x_i}(z)^j \ = \ 1 \ .
\end{equation*}
When $g_{x_i}(z)$ is expanded in terms of $Z$, its leading coefficient
\index{coefficients $A_{v,ij}$!leading} 
will be some $b_{v,i} \in K_v(x_i)^{\times}$:
\begin{equation*}
b_{v,i} = \lim_{Z \rightarrow 0} \frac{g_{x_i}(\varrho_i(Z))}{Z} \ .
\end{equation*}
It follows that $\varphi_{ij}(\varrho_i(Z))$ has the leading term 
$(b_{v,i})^{-j} Z^{-j}$.

We now apply Lemma \ref{BLemCPC5} with $\Phi_j(Z)$, $G(Z)$ and $H(Z)$ 
as above, taking $b = b_{v,i}$ and $J_0 = k_1 N_i$. 
Because only the $\varphi_{ij}(z)$ with $j \ge (n-k_1) N_i$
have poles of order $(n-k_1) N_i$ or more at $x_i$, we can write
\begin{eqnarray*}
G(Z) & = & \sum_{j=0}^{k_1 N_i - 1 } A_{v,ij} \varphi_{i,nN_i-j}(\varrho_i(Z))
         + \sum_{j=(n-k_1)N_i}^{\infty} A_{ij}^{\prime} Z^{-nN_i + j} \\
Y(Z)G(Z) & = & \sum_{j=0}^{k_1 N_i - 1} \Abar_{v,ij} \varphi_{i,nN_i-j}(\varrho_i(Z))
         + \sum_{j=(n-k_1)N_i}^{\infty} \Abar_{ij}^{\prime} Z^{-nN_i + j}
\end{eqnarray*}
where the $A_{v,ij}$ and $\Abar_{v,ij}$ are the
same as in (\ref{DGMN1}), (\ref{DGMN2}).  
Thus $g_0 = A_{v,i0}$ in Lemma \ref{BLemCPC5}.
Recall that $|A_{v,i0}|_v = |\tc_{v,i}|_v^n$ and that
$|Y(Z)-1|_v = |Y(z)-1|_v \le C_2 \varepsilon$ for all $Z \in D(0,r)$, 
where $C_2$ is the constant from Lemma \ref{MovingLemma1}.  
Put $C_3 = C_2$ in Lemma \ref{MovingLemma2}.  By Lemma \ref{BLemCPC5}, 
\begin{equation} \label{DGEST1}
|A_{v,ij}-\Abar_{v,ij}|_v \ \le \
   C_3 \varepsilon \cdot \frac{(|\tc_{v,i}|_v)^n}{(r|b_{v,i}|_v)^j}
\end{equation}
for each $j = 0, \ldots, k_1 N_i$.  

Letting $x_i$ vary, Lemma \ref{MovingLemma2} holds with
\begin{equation} \label{C4Def} 
C_4 \ = \ \frac{1}{\min(1,r|b_{v,1}|_v,\ldots,r|b_{v,m}|_v)} \ .
\end{equation}
\vskip -.35in
\end{proof}

\vskip .25 in
\noindent{\bf Lemma $\text{\bf \ref{MovingLemma3}}$. (Third Moving Lemma)}
\index{Third Moving Lemma|ii} 
{\it  There are constants $\varepsilon_3 > 0$ and $C_6, C_7 \ge 1$ 
$($depending only on $\phi_v(z)$, $E_v$, $H_v$, their $K_v$-simple decompositions 
\index{$K_v$-simple!set}
\index{$K_v$-simple!decomposition}
$\bigcup_{\ell = 1}^D \big(B(a_\ell,r_\ell) \cap \cC_v(F_{w_\ell})\big)$ 
and $\bigcup_{h=1}^N \big(B(\theta_h,\rho_h) \cap \cC_v(F_{u_h})\big)$, 
the choices of the $L$-rational and $L^{\sep}$-rational bases, the uniformizers $g_{x_i}(z)$, 
and the projective embedding of $\cC_v)$, 
such that if $0 < \varepsilon < \varepsilon_3$,
and if $\Fbar_{v,k_1}(z)$ is as in $(\ref{DMUG3})$ and $U_v^0 = \bigcup_{h=1}^N B(\theta_h,\rho_h)$, 
then there is a $K_v$-rational $(\fX,\vs)$-function $\Deltabar_{v,k_1}(z)$ of the form
\begin{equation*}
\Deltabar_{v,k_1}(z) \ = \ 
\sum_{i=1}^m \sum_{j=0}^{k_1 N_i-1} \Deltabar_{v,ij} \varphi_{i,(k_1+1)N_i-j}(z) \ ,
\end{equation*}
satisfying 
\begin{equation*}
\|\Deltabar_{v,ij}\|_{U_v^0} \ \le \  C_6 C_7^{k_1} \varepsilon  \ ,
\end{equation*}
such that 
when $\Gbar_v^{(k_1)}(z)$ from Lemma $\ref{MovingLemma2}$ is replaced with 
\begin{equation*} 
\hG_v^{(k_1)}(z) \ = \ \Gbar_v^{(k_1)}(z) + \Deltabar_{v,k_1}(z)\Fbar_{v,k_1}(z)\ ,
\end{equation*}
then for each $(i,j)$ with $1 \le i \le m$, $0 \le j < k_1 N_i$,
the coefficient $\Abar_{v,ij}$ of $\Gbar_v^{(k_1)}(z)$ is restored to the coefficient 
\index{coefficients $A_{v,ij}$!restoring}\index{patching functions, $G_v^{(k)}(z)$ for $1 \le k \le n$!expansion of}
$A_{v,ij}$ of $G_v^{(k_1)}(z)$ in $\hG_v^{(k_1)}(z)$.
}
\vskip .1 in

\begin{proof}  This is a consequence of Proposition \ref{FPushBound}.  
Using the $L$-rational basis, expand  
\index{basis!$L$-rational}\index{patching functions, $G_v^{(k)}(z)$ for $1 \le k \le n$!expansion of}  
\begin{eqnarray*} 
G_v^{(k_1)}(z) & = & \sum_{i=1}^m \sum_{j=0}^{(n-1)N_i-1} A_{v,ij} \varphi_{i,nN_i-j}(z) 
+ \sum_{\lambda=1}^\Lambda A_\lambda \varphi_\lambda(z) \ , \\
\Gbar_v^{(k_1)}(z) & = & \sum_{i=1}^m \sum_{j=0}^{(n-1)N_i-1} \Abar_{v,ij} \varphi_{i,nN_i-j}(z) 
+ \sum_{\lambda=1}^\Lambda \Abar_\lambda \varphi_\lambda(z) \ ,
\end{eqnarray*}
and using the $L^{\sep}$-rational basis, write
\index{basis!$L^{\sep}$-rational}\index{patching functions, $G_v^{(k)}(z)$ for $1 \le k \le n$!expansion of}
\begin{eqnarray*} 
G_v^{(k_1)}(z) & = & \sum_{i=1}^m \sum_{j=0}^{(n-1)N_i-1} \ta_{v,ij} \tphi_{i,nN_i-j}(z) 
+ \sum_{\lambda=1}^\Lambda \ta_\lambda \tphi_\lambda(z) \ , \\
\Gbar_v^{(k_1)}(z) & = & \sum_{i=1}^m \sum_{j=0}^{(n-1)N_i-1} \abar_{v,ij} \tphi_{i,nN_i-j}(z) 
+ \sum_{\lambda=1}^\Lambda \abar_\lambda \varphi_\lambda(z) \ .
\end{eqnarray*}

By Proposition \ref{TransitionProp}(C), the transition matrix from the $L$-rational basis
\index{basis!$L$-rational}
\index{$L$-rational basis!transition matrix block diagonal}
to the $L^{\sep}$-rational basis is block diagonal with blocks of size $J$, and for a given $i$
\index{basis!$L^{\sep}$-rational}
the same $J \times J$ matrix $\cB_{i,jk}$ occurs for each block.  
By Lemma \ref{MovingLemma2} 
we have $|A_{v,ij} - \Abar_{v,ij}|_v \le C_3 C_4^j(|\tc_{v,i}|_v)^n \varepsilon$ 
for all $1 \le i \le m$, $0 \le j < k_1 N_i$. 
Since $J|N_i$ for each $i$, there is a constant $C_5$ such that 
\begin{equation} \label{aabarbound}
|\ta_{v,ij} - \abar_{v,ij}|_v \ \le \ C_5 C_4^j (|\tc_{v,i}|_v)^n \varepsilon 
\end{equation} 
for all $1 \le i \le m$, $0 \le j < k_1 N_i$.
Indeed, putting $B_v = \max_{1 \le i \le m, 1 \le j,k \le J} |\cB_{i,jk}|_v$,  
we can take $C_5 = B_v C_3 C_4^{J-1}$.  

For each $i,j$ put $\tdelta_{v,ij} = \ta_{v,ij} - \abar_{v,ij}$.  
Since $G_v^{(k_1)}(z)$\index{patching functions, $G_v^{(k)}(z)$ for $1 \le k \le n$!are $K_v$-rational}
and $\Gbar_v^{(k_1)}(z)$ are $K_v$-rational, the $\tdelta_{v,ij}$ 
belong to $L_{w_v}^{\sep}$\index{distinguished place $w_v$} and are $K_v$-symmetric
 (Proposition \ref{CoeffRationalityCor}). 
\index{$K_v$-symmetric!set of numbers} 

We now apply Proposition \ref{FPushBound} taking $\ell = k_1$, $F_v(z) = \Fbar_{v,k_1}(z)$, and 
\begin{equation*}
\tdelta \ = \ (\tdelta_{v,ij})_{1 \le i \le m, 0 \le j < k_1 N_i} \ .
\end{equation*}   
We will take $r$ in Proposition \ref{FPushBound} to be the same number 
as in the proof of Lemma \ref{MovingLemma2}.  Comparing (\ref{varpivDef}) and (\ref{C4Def})
shows that if $\varpi_v$ is the constant from Proposition \ref{FPushBound} 
and $C_4$ is the constant from Lemma \ref{MovingLemma2}, then $C_4 = \varpi_v^{-1}$.
Letting $\tUpsilon_v$ be the constant from Proposition \ref{FPushBound},
and recalling from the discussion after (\ref{DMUG2}) that 
the leading coefficient $d_{v,i}$ of $\Fbar_{v,k_1}(z)$ at $x_i$ has absolute value 
\index{coefficients $A_{v,ij}$!leading} 
$|d_{v,i}|_v = |\tc_{v,i}|_v^{n-k_1-1}$, 
we will take $\rho$ in Proposition \ref{FPushBound} to be
\begin{equation} \label{rhoChoice}
\rho \ = \ \frac{C_5 |\tc_{v,i}|_v^{k_1+1}}{\tUpsilon_v} \cdot \varepsilon \ .
\end{equation}
Then for all $(i,j)$ with $1 \le i \le m, 0 \le j < k_1 N_i$ we have 
\begin{equation} \label{tdeltaBound}
|\tdelta_{v,ij}|_v \ = \ 
|\ta_{v,ij} - \abar_{v,ij}|_v \ \le \ C_5 C_4^j (|\tc_{v,i}|_v)^n \varepsilon \ = \ 
\tUpsilon_v \varpi_v^{-j} |d_{v,i}|_v \rho \ .  
\end{equation} 

By Proposition \ref{FPushBound} there is a unique 
$\vDelta = (\Deltabar_{v,is})_{1 \le i \le m, 0 \le s < k_1 N_i} \in (L_{w_v}^{\sep})^{k_1 N}$\index{distinguished place $w_v$}  
for which
\begin{equation*}
\Phi_{\Fbar_{v,k_1}}^{\sep}(\vDelta) \ = \ \tdelta \ ; 
\end{equation*}
moreover the $\Deltabar_{v,ij}$ are $K_v$-symmetric and 
\index{$K_v$-symmetric!set of numbers}
\begin{equation*} 
\Deltabar_{v,k_1}(z) \ := \ \sum_{i=1}^m \sum_{j=0}^{k_1 N_i-1} \Deltabar_{v,ij} \varphi_{(k_1+1)N_i -j}(z) 
\end{equation*}
is $K_v$-rational.  The fact that $\Phi_{\Fbar_{v,k_1}}^{\sep}(\vDelta) = \tdelta$
and each $\tdelta_{v,ij} = \ta_{v,ij}-\abar_{v,ij}$ means that 
\begin{equation*}
\Deltabar_{v,k_1}(z) \Fbar_{v,k_1}(z)
\ = \ \sum_{i=1}^m \sum_{j=0}^{\ell N_i-1}  (\ta_{v,ij}-\abar_{v,ij}) \cdot \tphi_{i,(k+\ell)N_i - j}(z) 
                  \ + \ \text{lower order terms \ .}                   
\end{equation*}
Consequently, when we expand $\Deltabar_{v,k_1}(z) \Fbar_{v,k_1}(z)$ using the $L$-rational basis,
\index{basis!$L$-rational}
we have  
\begin{equation*}
\Deltabar_{v,k_1}(z) \Fbar_{v,k_1}(z)
\ = \ \sum_{i=1}^m \sum_{j=0}^{\ell N_i-1}  (A_{v,ij}-\Abar_{v,ij}) \cdot \varphi_{i,(k+\ell)N_i - j}(z) 
                  \ + \ \text{lower order terms \ .}                    
\end{equation*}
This means that when $\Gbar_v^{(k_1)}(z)$\index{patching functions, $G_v^{(k)}(z)$ for $1 \le k \le n$!constructed by patching} 
is replaced with 
\begin{equation*} 
\hG_v^{(k_1)}(z) \ = \ \Gbar_v^{(k_1)}(z) + \Deltabar_{v,k_1}(z)\Fbar_{v,k_1}(z)\ ,
\end{equation*}
then for each $(i,j)$ with $1 \le i \le m$, $0 \le j < k_1 N_i$,
the coefficient $\Abar_{v,ij}$ of $\Gbar_v^{(k_1)}(z)$ is changed 
\index{coefficients $A_{v,ij}$} 
to $A_{v,ij}$ in $\hG_v^{(k_1)}(z)$.\index{patching functions, $G_v^{(k)}(z)$ for $1 \le k \le n$!constructed by patching}

Finally, by (\ref{rhoChoice}) and (\ref{tdeltaBound}), 
and by (\ref{FBallC2}) of Proposition \ref{FPushBound}, for each  
$(i,j)$ with $1 \le i \le m, 0 \le j < k_1 N_i$ we have 
\begin{equation*}
|\Deltabar_{v,ij}|_v \ \le \ \varpi_v^{-j} \rho 
\ = \ \frac{C_5}{\tUpsilon_v} C_4^j \cdot (|\tc_{v,i}|_v)^{k_1+1} \cdot \varepsilon \ .
\end{equation*}
On the other hand, by Proposition \ref{GBL} there is a constant $C_v^0$ 
such that $\|\varphi_{ij}(z)\|_{U_v^0} \le (C_v^0)^j$ for all $i$ and $j$.  
Without loss, we can assume that $C_v^0 \ge 1$.  It follows that 
\begin{eqnarray*}
\|\Deltabar_{v,k_1}(z)\|_{U_v^0} & \le & \max_{1 \le i \le m, 0 \le j < k_1 N_i}
\big(|\Deltabar_{v,ij}|_v \cdot \|\varphi_{i,(k_1+1)N_i -j} \|_{U_v^0} \big) \\
& \le & \max_{1 \le i \le m, 0 \le j < k_1 N_i}
\Big( \frac{C_5}{\tUpsilon_v} C_4^j \cdot (|\tc_{v,i}|_v)^{k_1+1} 
     \cdot (C_v^0)^{(k_1+1)N_i-j} \cdot \varepsilon \Big) 
\ \le \ C_6 C_7^{k_1} \cdot \varepsilon \ ,
\end{eqnarray*}
where 
\begin{equation*}
C_6 \ = \ \max \Big( 1, \  
   \frac{C_5}{\tUpsilon_v} \cdot \max_i(|\tc_{v,i}|_v) \cdot (C_v^0)^{\max_i(N_i)} \Big) 
\end{equation*}
and 
\begin{equation*}
C_7 \ = \ \max \Big( 1, \ 
  \max_i(|\tc_{v,i}|_v) \cdot \max \big(C_v^0,C_4\big)^{\max_i(N_i)} \Big) \ .
\end{equation*} 
This completes the proof. 
\end{proof}
\index{patching construction!for nonarchimedean $K_v$-simple sets!proofs of the moving lemmas|)}
\index{patching argument!local!for nonarchimedean $K_v$-simple sets|)}

%% file: NewFSZAppA.tex
\chapter{$(\fX,\vs)$-Potential theory} \label{AppA}  

In this appendix we study potential theory for the 
\index{potential theory!$(\fX,\vs)$|ii}\index{canonical distance!$[z,w]_{\fX,\vs}$!potential theory for|ii}
\index{$(\fX,\vs)$-potential theory|ii} 
$(\fX,\vs)$-canonical distance.  
\index{canonical distance!$[z,w]_{\fX,\vs}$} 
In section \ref{AppA}.\ref{XSPotTheorySection}
we discuss the basic facts of $(\fX,\vs)$-potential theory for compact sets
\index{potential theory!$(\fX,\vs)$}   
concerning potential functions, 
equilibrium distributions, the transfinite diameter, and the Chebyshev constant. 
\index{potential function!$(\fX,\vs)$}
\index{equilibrium potential!$(\fX,\vs)$}
\index{equilibrium distribution!$(\fX,\vs)$}
\index{transfinite diameter!$(\fX,\vs)$} 
\index{Chebyshev constant!$(\fX,\vs)$} 
In section \ref{AppA}.\ref{ArchMassBoundsSection}, 
which concerns the archimedean case, 
we derive bounds for the mass the $(\fX,\vs)$-equilibrium distribution 
of $H$ can give to ``small'' subsets of $H$.  
In section \ref{AppA}.\ref{NonArchMuXSSection}, which concerns the nonarchimedean case,
we determine the $(\fX,\vs)$-equilibrium distributions for 
a class of well-behaved sets.  

\vskip .1 in
Fix a projective embedding of $\cC/K$. Given a place $v$ of $K$, 
let $\|z,w\|_v$ be the corresponding spherical metric on $\cC_v(\CC_v)$.
\index{spherical metric}  
As in \S\ref{Chap3}.\ref{AssumptionsSection}, if $v$ is nonarchimedean let 
$q_v$ be the order of the residue field of $K_v$, and let $\log_v(x)$ be
the logarithm to the base $q_v$.  If $v$ is archimedean, 
put $\log_v(x) = \ln(x)$.  

Let $\fX = \{x_1, \ldots, x_m\} \subset \cC(\tilde{K})$ 
be the finite, galois-stable set of points from \S\ref{Chap3}.\ref{AssumptionsSection}.  
For each $x_i \in \fX$, let $g_{x_i}(z) \in K(\cC)$ be the uniformizer at 
$x_i$ chosen in \S\ref{Chap3}.\ref{AssumptionsSection}, 
and let the canonical distance $[z,w]_{x_i}$ be normalized so
\index{canonical distance!$[z,w]_{\zeta}$!normalization of}
that for each $w \ne x_i$, 
\begin{equation*}  
\lim_{z \rightarrow x_i} 
      \left([z,w]_{x_i} \cdot |g_{x_i}(z)|_v \right) \ = \ 1 \ .
\end{equation*}
Let $\vs = (s_1, \ldots, s_m) \in \cP^m$ be a probability vector.  
As in \S\ref{Chap3}.\ref{CanonicalDistanceSection}, we define the $(\fX,\vs)$-canonical distance by 
\index{canonical distance!$[z,w]_{\fX,\vs}$}
\begin{equation*}
[z,w]_{\fX,\vs} \ = \ \prod_{i=1}^n ([z,w]_{x_i})^{s_i} \ .
\end{equation*}


\section{ $(\fX,\vs)$-Potential Theory for Compact Sets } \label{XSPotTheorySection} 
  
Let $H \subset \cC_v(\CC_v) \backslash \fX$ be a compact set. 
In this section we will define analogues 
of the classical logarithmic capacity,\index{capacity!logarithmic}\index{capacity!$(\fX,\vs)$}  
transfinite diameter,\index{transfinite diameter!$(\fX,\vs)$} 
Chebyshev constant,\index{Chebyshev constant!$(\fX,\vs)$}\index{Chebyshev constant}  
potential functions,\index{potential function!$(\fX,\vs)$}\index{equilibrium potential!$(\fX,\vs)$}
and Green's functions,\index{Green's function!$(\fX,\vs)$}   
relative to the kernel $[z,w]_{\fX,\vs}$.\index{potential theory!$(\fX,\vs)$}
   
We will study these objects and 
their relation with the corresponding objects when 
$\fX$ consists of a single point. 
The proofs of all the results below are classical and 
(with minor modifications) are the same as 
those in (\cite{RR1}, \S3 and \S4), so for the most part we only sketch them.  

\vskip .1 in
We first define the $(\fX,\vs)$-capacity. \index{capacity!$(\fX,\vs)$|ii} 
For any probability measure $\nu$ supported on $H$, 
the $(\fX,\vs)$-energy  is
\begin{equation} \label{FVF0} 
I_{\fX,\vs}(\nu) 
\ = \ \iint_{H \times H} -\log_v([z,w]_{\fX,\vs}) \, d\nu(z) d\nu(w) \ .
\end{equation} 
and the $(\fX,\vs)$-potential function is\index{potential function!$(\fX,\vs)$}
\index{equilibrium potential!$(\fX,\vs)$|ii}
\begin{equation} \label{FVC0} 
u_{\fX,\vs}(z,\nu) 
\ = \ \int_H -\log_v([z,w]_{\fX,\vs}) \, d\nu(w) \ .
\end{equation} 
The $(\fX,\vs)$-Robin constant is defined by
\index{Robin constant!$(\fX,\vs)$|ii}
\begin{equation} \label{FVF1}
V_{\fX,\vs}(H) \ = \ 
\inf_{\substack{ \text{probability measures} \\ \text{$\nu$ supported on $H$} }} I_{\fX,\vs}(\nu) \ , 
\end{equation}
and the $(\fX,\vs)$-capacity is given by
\index{capacity!$(\fX,\vs)$|ii}
\begin{equation} \label{FCF1} 
\gamma_{\fX,\vs}(H) \ = \ 
\left\{ \begin{array}{ll} e^{-V_{\fX,\vs}(H)} & \text{if $v$ is archimedean,} \\
                         q_v^{-V_{\fX,\vs}(H)}& \text{if $v$ is nonarchimedean.} 
         \end{array} \right. 
\end{equation}

We next define the $(\fX,\vs)$-transfinite diameter.\index{transfinite diameter!$(\fX,\vs)$|ii}  
For $N = 2, 3, ...$ let 
\begin{equation} \label{FTD1}
d_N(H) \ = \ \sup_{q_1, \ldots, q_N \in H} 
         (\prod^N_{\substack{ i, j = 1 \\ i \ne j }} 
                      [q_i,q_j]_{\fX,\vs})^{1/N^2}  \ ;
\end{equation}
then the $(\fX,\vs)$-transfinite diameter is                    
\begin{equation} \label{FTD2}
d_{\fX,\vs}(H) \ = \ \lim_{N \rightarrow \infty} d_N(H) \ .
\end{equation}
The existence of the limit 
follows by a classical argument, given in (\cite{RR1}, p.150 and pp.203-204)
for the kernel $[z,w]_{\zeta}$.   
There the exponent $1/N^2$ in (\ref{FTD1}) is replaced by $1/N(N-1)$, 
and the $d_N(H)$ are shown to be monotonically decreasing.  
Our modification to the exponent does not affect the convergence in (\ref{FTD2}),  
or the value of the limit.  

Finally, we define the restricted $(\fX,\vs)$-Chebyshev constant $\CH^*_{\fX,\vs}(H)$.
\index{Chebyshev constant!$(\fX,\vs)$} 
Given points $a_1, \ldots, a_N \in \cC_v(\CC) \backslash \fX$, 
consider the $(\fX,\vs)$-pseudopolynomial (see \S\ref{Chap3}.\ref{XSfunctionSection})
\index{pseudopolynomial!$(\fX,\vs)$}
\begin{equation*} 
P(z;a_1, \ldots, a_N) 
\ = \ \prod_{i=1}^N [z,a_i]_{\fX,\vs} \ .
\end{equation*} 
Writing $\|P\|_H  = \sup_{z \in H} P(z)$, put  
\begin{equation} \label{FCH0}
\CH^*_N(H) \ = \ \inf_{a_1, \ldots, a_N \in H} (\|P(z;a_1, \ldots, a_N)\|_H)^{1/N} \ ,
\end{equation}  
and then define
\begin{equation} \label{FCH1} 
\CH^*_{\fX,\vs}(H) \ = \ \lim_{N \rightarrow \infty} \CH^*_N(H) \ .
\end{equation}
The existence of the limit in (\ref{FCH1}) follows from arguments similar to those
in (\cite{RR1}, p.151 and pp.203-304).  
We call $\CH^*_{\fX,\vs}(H)$ the restricted Chebyshev constant  
\index{Chebyshev constant!restricted} 
since the points $a_1, \ldots, a_N$ are required to be in $H$; lifting that restriction, 
it is also possible to define an unrestricted Chebyshev constant $\CH_{\fX,\vs}(H)$, whose value turns out 
to be the same as the restricted one. 

The following theorems summarize the main facts concerning these objects:  

\begin{theorem} \label{ATE10A}
Let  $H \subset \cC_v(\CC_v) \backslash \fX$ be compact.  
Then for each probability vector $\vs \in \cP^m$, 
\begin{equation} \label{FRev1} 
 \gamma_{\fX,\vs}(H) \ = \ d_{\fX,\vs}(H) \ = \ \CH^*_{\fX,\vs}(H) \ ,
\end{equation}
and these quantities are $0$ if and only if $H$ has capacity $0$
\index{capacity $= 0$}
in the sense of Definition $\ref{FCapDef}$.
\end{theorem}
\index{Chebyshev constant!$(\fX,\vs)$} 
\index{transfinite diameter!$(\fX,\vs)$} 
\index{capacity!$(\fX,\vs)$} 

\begin{proof} The proofs are analogous to those of \cite{RR1}, Theorems 3.1.18 and 4.1.19.
\end{proof}  

\begin{theorem} \label{ATE10B}  Let  $H \subset \cC_v(\CC_v) \backslash \fX$ be compact,
with positive capacity.  Then for each probability vector $\vs \in \cP^m$,
\index{capacity $> 0$} 

$(A)$  If $H$ has positive capacity, then there
\index{capacity $> 0$}
is a unique probability measure $\mu = \mu_{\fX,\vs}$ on $H$,
called the $(\fX,\vs)$-equilibrium distribution of $H$, for which
\index{equilibrium distribution!$(\fX,\vs)$|ii} 
\index{Robin constant!$(\fX,\vs)$}
\begin{equation*}
V_{\fX,\vs}(H) \ = \ I_{\fX,\vs}(\mu) \ ;
\end{equation*}

$(B)$  For this measure $\mu_{\fX,\vs}$, the potential function
\index{equilibrium distribution!$(\fX,\vs)$} 
\index{equilibrium potential!$(\fX,\vs)$|ii}
\begin{equation*}
u_{\fX,\vs}(z,H) \ := \ u_{\fX,\vs}(z,\mu_{\fX,\vs})  
\ = \ \int_{H} -\log_v([z,w]_{\fX,\vs}) \, d\mu_{\fX,\vs}(w)
\end{equation*}
satisfies $u_{\fX,\vs}(z,H) \le V_{\fX,\vs}(H)$ for all $z \in \cC_v(\CC_v)$,
\index{Robin constant!$(\fX,\vs)$}
with  $u_{\fX,\vs}(z,H) = V_{\fX,\vs}(H)$
\index{equilibrium potential!takes constant value a.e. on $E_v$} 
for all $z \in H$ except possibly an $F_\sigma$-set $e_{\fX,\vs} \subset H$ of inner capacity $0$.
\index{capacity $= 0$} 
Moreover, $u_{\fX,\vs}(z,H)$ is continuous on $\cC_v(\CC_v) \backslash e_{\fX,\vs}$.  
 
In the archimedean case, $u_{\fX,\vs}(z,H) < V_{\fX,\vs}(H)$ on each component of 
\index{equilibrium potential!$(\fX,\vs)$}
$\cC_v(\CC) \backslash H$ which contains a point $x_i \in \fX$ with $s_i > 0$, and 
$u_{\fX,\vs}(z,H) = V_{\fX,\vs}(H)$ on all other components of $\cC_v(\CC) \backslash H$.
The exceptional set $e_{\fX,\vs}$ is contained in $\partial H_{\fX,\vs}$,
\index{Robin constant!$(\fX,\vs)$}\index{exceptional set} 
the part of the boundary of $H$ shared by the components of $\cC_v(\CC) \backslash H$ 
on which $u_{\fX,\vs}(z,H) < 0$, and $H$ and $\partial H_{\fX,\vs}$ have the same capacity,\index{capacity}
potential function, and equilibrium distribution with respect to $[z,w]_{\fX,\vs}$.   
\index{equilibrium distribution!$(\fX,\vs)$}\index{equilibrium potential!$(\fX,\vs)$}
\index{equilibrium potential!$(\fX,\vs)$}
Furthermore, $u_{\fX,\vs}(z,H)$ 
is superharmonic\index{superharmonic} on $\cC_v(\CC) \backslash \fX$, 
subharmonic\index{subharmonic} on $\cC_v(\CC) \backslash H$,   
and harmonic on  $\cC_v(\CC) \backslash (H \bigcup \fX)$.  At each $x_i \in \fX$, 
$u_{\fX,\vs}(z,H) + s_i \log(|z-x_i|)$ extends to a function harmonic
in a neighborhood of $x_i$. 

In the nonarchimedean case, $u_{\fX,\vs}(z,H) < V_{\fX,\vs}(H)$ for all $z \in \cC_v(\CC_v) \backslash H$,
\index{Robin constant!$(\fX,\vs)$}
and $u_{\fX,\vs}(z,H) + s_i \log_v(|z-x_i|_v)$ has a finite limit at each $x_i \in \fX$ .  
\index{equilibrium potential!$(\fX,\vs)$} 
\end{theorem}

\noindent{\bf Remark.} 
When $H$ is clear from the context, we will often write $u_{\fX,\vs}(z)$
for $u_{\fX,\vs}(z,H)$.  
\index{equilibrium potential!$(\fX,\vs)$}

\vskip .1 in
\begin{proof}  In the classical theory, the assertions in Theorems \ref{ATE10A} and \ref{ATE10B}
are the main consequences of Maria's Theorem and Frostman's theorem. 
\index{Maria's Theorem}\index{Frostman's Theorem}
They are established for the kernel $[z,w]_{\zeta}$ 
in Theorems 3.1.6, 3.1.7, 3.1.12, 
and 3.1.18 of (\cite{RR1}, \S 3.1) in the archimedean case, 
and in Theorems 4.1.11, 4.1.19, and 4.1.22 of (\cite{RR1}, \S 4.1)
in the nonarchimedean case.  

In the archimedean case, the continuity/harmonicity properties of 
$\log([z,w]_{\fX,\vs})$ shown in Proposition \ref{APropA2}, 
together with the Maximum principle for harmonic functions,
\index{Maximum principle!for harmonic functions} 
allow the proofs in (\cite{RR1}, \S3.1) 
to be carried over for $[z,w]_{\fX,\vs}$.
The property of the canonical distance needed for those proofs 
\index{canonical distance!archimedean!comparable with absolute value}
is that for each $q \in H$, and each disc $D(q,r)$ with $\zeta \notin \overline{D(q,r)}$, 
if we fix a coordinate chart on $D(q,r)$, then 
there is a constant $C$ (depending on $\zeta$ and the choice of coordinates) 
such that for all $z \ne w \in D(q,r)$ 
\begin{equation*}
-\log([z,w]_{\zeta}) - C \ \le \ -\log(|z-w|) \ \le \ -\log([z,w]_{\zeta}) - C 
\end{equation*} 
(see \cite{RR1}, p.139).  By Proposition \ref{APropA2} this holds for 
$[z,w]_{\fX,\vs}$.    

In the nonarchimedean case, the proofs given in 
(\cite{RR1}, \S4.1) 
use two properties of the canonical distance.  
First, for each $q \in H$, each $\zeta \in \fX$, 
and each isometrically parametrizable ball $B(q,r)$ disjoint from $\fX$
(see Definition \ref{IsoParamDef}), there
is a constant $C = C_{q,\zeta}$ such that $[z,w]_{\zeta} = C \|z,w\|_v$
\index{canonical distance!nonarchimedean!comparable with chordal distance}
for all $z, w \in B(q,r)$.
Since $[z,w]_{\fX,\vs}$ is a weighted product of the $[z,w]_{x_i}$, 
with the weights summing to $1$, 
Proposition \ref{APropA2} shows that this property holds for $[z,w]_{\fX,\vs}$ as well. 
Second, for each pair of 
points $w,\zeta \in \cC_v(\CC)$, 
and each isometrically parametrizable ball $B(w,r)$ not containing $\zeta$,
there is a function $f(z) \in \CC_v(\cC)$ of degree $N$ say, 
having all its zeros in $B(w,r)$ and having poles only at $\zeta$, such that
\begin{equation*}
           [z,w]_{\zeta} \ = \ (|f(z)|_v)^{1/N} 
\end{equation*} 
for all $z \notin B(w,r)$ (see \cite{RR1}, Proposition 2.1.6). 
For $[z,w]_{\fX,\vs}$ the analogue
of this is that if $B(w,r)$ is an isometrically parametrizable ball disjoint
from $\fX$, then there are functions $f_i(z) \in \CC_v(\cC)$ of degree $N_i$
say, having all their zeros in $B(w,r)$ and such that $f_i(z)$ has poles
only at $x_i$, for which 
\begin{equation*}
 -\log_v([z,w]_{\fX,\vs}) \ = \ 
              -\sum_{i=1}^m s_i \cdot \frac{1}{N_i} \log_v(|f_i(z)|_v) 
\end{equation*}
for all $z \notin B(w,r)$. 
\end{proof} 

\vskip .1 in
The following proposition can often be used to show that the exceptional set $e_{\fX,\vs}$ in
\index{exceptional set}
Theorem \ref{ATE10B} is empty.  By an arc, we mean a homeomorphic image of the segment $[0,1]$.
\index{exceptional set}\index{arc|ii}   

\begin{proposition} \label{ACE11}  Let all assumptions be as in Theorem $\ref{ATE10B}$.

$(A)$ If $K_v$ is archimedean, then $u_{\fX,\vs}(z,H) = V_{\fX,\vs}(H)$ at each point $z_0 \in H$ for which
\index{Robin constant!$(\fX,\vs)$}
\index{equilibrium potential!$(\fX,\vs)$}
there is an arc $A \subset H$ with $z_0 \in A$.
\index{equilibrium potential!takes constant value a.e. on $E_v$}
\index{arc}

$(B)$ If $K_v$ is nonarchimedean with $\Char(K_v) = 0$, and if $p$ is the rational prime lying under $v$, 
then $u_{\fX,\vs}(z,H) = V_{\fX,\vs}(H)$ at each point $z_0 \in H$
\index{Robin constant!$(\fX,\vs)$} 
for which, for some $r > 0$, there is an isometric parametrization 
$f_{z_0}: D(0,r) \rightarrow B(z_0,r) \subset \cC_v(\CC_v)$ with $f_{z_0}(0) = z_0$, 
such that $f(\ZZ_p \cap D(0,r)) \subset H$.  If $K_v$ is nonarchimedean with $\Char(K_v) = p > 0$, 
the analogous assertion holds with $\ZZ_p$ replaced by $\FF_p[[T]]$. 
\end{proposition}

\begin{proof} In the archimedean case, 
this is a classical consequence of the existence of a ``barrier''.\index{barrier}  
The proof is given (\cite{RR1}, Theorem 3.1.9) when $\fX$ is a single point, 
and the argument, which is purely local, carries over unchanged in the general case.  

\smallskip 
In the nonarchimedean case, the proof uses the monotonicity of upper Green's functions of compact sets
\index{Green's function!upper $(\fX,\vs)$} 
(\cite{RR1}, Proposition 4.4.1(A)).  We give the argument when $\Char(K_v) = 0$;  
the proof when $\Char(K_v) > 0$ is similar.  
After shrinking $r$ if necessary, we can assume that $B(z_0,r) \cap \fX = \phi$.  
By Proposition \ref{APropA2}.(B1), there is a constant $C > 0$ such that 
$[z,w]_{\fX,\vs} = C \|z,w\|_v$ for all $z, w \in B(z_0,r)$.  
Since $f_{z_0}$ is an isometric parametrization,
we have $[f_{z_0}(x),f_{z_0}(y)]_{\fX,\vs} = C|x-y|_v$ for all $x, y \in D(0,r)$.  

Since a set $e \subset D(0,r)$ has positive inner capacity if and only if it supports a probability measure
\index{capacity $> 0$} 
$\nu$ for which $\iint -\log_v(|x-y|_v)\, d\nu(x) d\nu(y) < \infty$, 
by pushing forward or pulling back appropriate measures one sees that 
$f_{z_0}$ takes sets of positive inner capacity to sets of positive inner capacity, 
\index{capacity $> 0$}
and sets of inner capacity $0$ to sets of inner capacity $0$. 
\index{capacity $= 0$} 

Put $H_0 = H \cap B(z_0,r)$, and let $\nu_0 = \mu_{\fX,\vs}|_{H_0}$. 
\index{equilibrium distribution!$(\fX,\vs)$|ii}  
Since $H_0$ contains $f_{z_0}(\ZZ_p \cap D(0,r))$,
it has positive capacity.  This means that $\nu_0(H_0) > 0$.  
\index{capacity $> 0$}
Let $u(z,\nu_0) = \int -\log_v(\|z,w\|_v) \, d\nu_0(w)$.  
\index{potential function}
Since $[z,w]_{\fX,\vx}$ is constant on 
pairwise disjoint isometrically parametrizable balls in $\cC_v(\CC_v) \backslash \fX$, there 
is a constant $D$ such that for all $z \in B(z_0,r)$,  
\begin{equation} \label{FUDFormula} 
u_{\fX,\vs}(z,H) \ = \ D + u(z,\nu_0) \ .
\end{equation} 
\index{equilibrium potential!$(\fX,\vs)$}

Pull back $[z,w]_{\fX,\vs}$, $u_{\fX,\vx}(z,H)$, $H_0$ and $\nu_0$ to $D(0,r)$ using $f_{z_0}$.
Let $E = f_{z_0}^{-1}(H_0)$ and put $\nu = (1/\nu_0(H_0)) f_{z_0}^*(\nu_0)$; then $\nu$ is a probability
measure supported on $E$.  Consider the potential function 
\begin{equation*}
u_{\infty}(x,\nu) \ := \ \int -\log_v(|x-y|_v)\, d\nu(y) 
\end{equation*}
\index{potential function}
\index{equilibrium potential}
on $\PP^1(\CC_v)$.  Since $\mu_{\fX,\vs}$ is the equilibrium measure of $H$, 
\index{equilibrium distribution!$(\fX,\vs)$}
we have $u_{\fX,\vs}(z,H) \le  V_{\fX,\vs}(H)$ for all $z$,
\index{equilibrium potential!$(\fX,\vs)$}
\index{Robin constant!$(\fX,\vs)$}
with equality on $H$ except on a set of inner capacity $0$.
\index{equilibrium potential!takes constant value a.e. on $E_v$}
\index{capacity $= 0$}  
By (\ref{FUDFormula}) there is a constant $V$ 
such that  $u_{\infty}(z,\nu_0) \le V$ on $E$, 
with equality except on a set of inner capacity $0$. 
\index{potential function} 
\index{capacity $= 0$}
It follows from (\cite{RR1}, Proposition 4.1.23)   
that $\nu$ is the equilibrium measure of $E$.  This means that $G(x,\infty;H)  := V - u_{\infty}(x,\nu)$
\index{Green's function}
\index{equilibrium distribution}
is the upper Green's function of $H$ with respect to $\infty$.  Since $\ZZ_p \cap D(0,r) \subseteq H$, 
\index{Green's function!upper $(\fX,\vs)$} 
by (\cite{RR1}, Proposition 4.4.1(A)) we have  
\begin{equation*} 
G(x,\infty;\ZZ_p \cap D(0,r)) \ \ge \ G(x,\infty;H) 
\end{equation*} 
for all $x \in \CC_v$.  The explicit computation in  
Proposition \ref{OwProp} gives $G(0,\infty;\ZZ_p \cap D(0,r)) = 0$,   
so $G(0,\infty;H) = 0$ as well.  
This means that $u_{\infty}(0,\nu) = V$, and hence by (\ref{FUDFormula}) that 
\index{equilibrium potential!$(\fX,\vs)$}
$u_{\fX,\vs}(z_0,H) = V_{\fX,\vs}(H)$.
\index{Robin constant!$(\fX,\vs)$}
\end{proof}

\vskip .1 in 
The following proposition sometimes lets us determine the equilibrium distribution. 
\index{equilibrium distribution!$(\fX,\vs)$} 
\index{equilibrium distribution!determining|ii}

\begin{proposition} \label{IDProp}
Let $\mu_0$ be a probability measure on $H$ for which 
there is a constant $V < \infty$ such that 
the potential function $u_{\fX,\vs}(z,\mu_0)$ equals $V$ on $H$,
\index{equilibrium potential!$(\fX,\vs)$}
except possibly on a set of inner capacity $0$. 
\index{equilibrium potential!takes constant value a.e. on $E_v$}
\index{capacity $= 0$} 
Then $V_{\fX,\vs}(H) = V$, and $\mu_{\fX,\vs} = \mu_0$.   
\index{Robin constant!$(\fX,\vs)$}
\index{equilibrium distribution!$(\fX,\vs)$|ii} 
\end{proposition} 

\begin{proof} By the same argument as in (\cite{RR1}, Lemmas 3.1.4 and 4.17),
a positive measure $\nu$ for which $I_{\fX,\vs}(\nu) < \infty$
cannot charge sets of inner capacity $0$.
\index{capacity $= 0$}  
In particular, this applies
to $\mu_0$ and $\mu_{\fX,\vs}$.  Hence by Theorem \ref{ATE10B}(B) and Fubini-Tonelli, 
\index{equilibrium distribution!$(\fX,\vs)$} 
\index{equilibrium potential!$(\fX,\vs)$}
\index{Fubini-Tonelli theorem} 
\begin{eqnarray*}
V  & = & \int_{H} u_{\fX,\vs}(z,\mu_0) \, d\mu_{\fX,\vs}(z) 
   \ = \ \iint_{H \times H} -\log_v([z,w]_{\fX,\vs}) 
               \, d\mu_0(w) d\mu_{\fX,\vs}(z) \\
   & = & \int_{H} u_{\fX,\vs}(w,\mu_{\fX,\vs}) \, d\mu_0(w) \ = \ V_{\fX,\vs}(H) \ .
\end{eqnarray*}
\index{Robin constant!$(\fX,\vs)$}
Consequently 
\begin{eqnarray*}
I_{\fX,\vs}(\mu_0) & = & \iint_{H \times H} -\log_v([z,w]_{\fX,\vs}) 
               \, d\mu_0(w) d\mu_0(z)  \\
& = & \int_H u_{\fX,\vs}(z,\mu_0) \, d\mu_0(z) \ = \ V_{\fX,\vs}(H) \ .
\end{eqnarray*}
\index{equilibrium potential!$(\fX,\vs)$}
\index{Robin constant!$(\fX,\vs)$}
Since $\mu_{\fX,\vs}$ is the unique probability measure minimizing the $(\fX,\vs)$-energy integral
\index{energy integral!$(\fX,\vs)$} 
(Theorem \ref{ATE10B}(A)), it follows that $\mu_{\fX,\vs} = \mu_0$.
\index{equilibrium distribution!$(\fX,\vs)$}
\index{equilibrium distribution!determining}
\end{proof}  

\vskip .1 in
We define the Green's function $G_{\fX,\vs}(z;H)$ to be 
\begin{equation*}
G_{\fX,\vs}(z;H) \ = \ V_{\fX,\vs}(H) - u_{\fX,\vs}(z) \ .
\end{equation*}
\index{equilibrium potential!$(\fX,\vs)$}
\index{Robin constant!$(\fX,\vs)$}
\index{Green's function!$(\fX,\vs)$|ii} 
For each $x_i \in \fX$, and each positive measure $\nu$ 
supported on $\cC_v(\CC_v) \backslash \{x_i\}$, put 
\begin{equation*}
u_{x_i}(z,\nu) 
      \ = \ \int_{\cC_v(\CC_v)} -\log_v([z,w]_{x_i}) \, d\nu(w) \ .
\end{equation*}
\index{potential function}
Let $\mu_i$ be the equilibrium distribution of $H$ with 
\index{equilibrium distribution}
respect to $[z,w]_{x_i}$, and write  
\begin{equation*}
u_{x_i}(z) \ = \ u_{x_i}(z,\mu_i) 
               \ = \ \int_{H} -\log_v([z,w]_{x_i}) \, d\mu_i(w) \ .
\end{equation*}
\index{equilibrium potential}
\index{potential function}
As in \S\ref{Chap3}.\ref{UpperGreenSection}, 
for each $x_i \in \fX$ the Green's function $G(z,x_i;H)$ is defined by  
\begin{equation*}
G(z,x_i;H) \ = \ V_{x_i}(H) - u_{x_i}(z) \ . 
\end{equation*}
\index{Robin constant}
\index{Green's function}
\index{equilibrium potential} 

We can express   $\mu_{\fX,\vs}$ and $G_{\fX,\vs}(z;H)$
\index{equilibrium distribution!$(\fX,\vs)$}
\index{Green's function!$(\fX,\vs)$}
in terms of the $\mu_i$ and $G(z,x_i;H)$: 
\index{Green's function} 

\begin{proposition}  \label{BPropF1}
Suppose $H \subset \cC_v(\CC_v) \backslash \fX$ is compact and has positive capacity.  
\index{capacity $> 0$}
Then \begin{eqnarray}
\mu_{\fX,\vs} & = & \sum_{i=1}^m s_i \mu_i \ , \label{AFE82} \\
G_{\fX,\vs}(z;H) & = & \sum_{i=1}^m s_i G(z,x_i;H)  \ , 
                                 \label{BFFF7B} 
\end{eqnarray}
\index{Green's function}
\index{Green's function!$(\fX,\vs)$|ii}
\index{equilibrium distribution!$(\fX,\vs)$|ii}
\index{equilibrium distribution}
and
\begin{equation}  \label{BFFF8}
V_{\fX,\vs}(H) \ = \ 
  \sum_{i,j,k=1}^m s_i s_j s_k \iint_{H \times H} 
             -\log_v([z,w]_{x_i}) \, d\mu_j(z) d\mu_k(w) \ .
\end{equation}  
\index{Robin constant!$(\fX,\vs)$|ii}
\index{equilibrium distribution}
In particular, $V_{\fX,\vs}(H)$ is a continuous function of $\vs$.     
\end{proposition}

\begin{proof} {(Archimedean Case.)} 
Formula (\ref{BFFF7B}) follows from the strong form of the Maximum 
principle for harmonic functions, 
\index{Maximum principle!for harmonic functions!strong form} 
applied to $G_{\fX,\vs}(z,H) - \sum s_i G(z,x_i;H)$ 
\index{Green's function}
\index{Green's function!$(\fX,\vs)$}
on $\cC_v(\cC) \backslash (H \bigcup \fX)$
(see \cite{RR1}, Proposition 3.1.1).  
Formula (\ref{AFE82}) is a consequence of (\ref{BFFF7B}), 
since by the Riesz decomposition theorem (see \cite{RR1}, Theorem 3.1.11)
\index{Riesz Decomposition theorem}
a measure can be recovered from its potential function.  In modern terminology, 
\index{equilibrium potential}
applying the $dd^c$ operator on $\cC_v(\CC) \backslash \fX$, one has  
$\mu_{\fX,\vs} = dd^c(-u_{\fX,\vs}(z))$ and
\index{equilibrium distribution!$(\fX,\vs)$}
\index{equilibrium potential!$(\fX,\vs)$}
$\mu_i = dd^c(-u_{x_i}(z)) = dd^c(G(z,x_i;H))$ for each $i$. 
\index{Green's function}

For the assertion about $V_{\fX,\vs}(H)$, note that 
\index{Robin constant!$(\fX,\vs)$}
\begin{eqnarray*}
V_{\fX,\vs}(H) 
   &=& \iint_{H \times H} -\log([z,w]_{\fX,\vs}) \, 
                           d\mu_{\fX,\vs}(z) d\mu_{\fX,\vs}(w) \\
   &=& \sum_{i,j,k = 1}^m \, s_i s_j s_k \left(
    \iint_{H \times H} -\log([z,w]_{x_i} \, d\mu_j(z) d\mu_k(w) 
    \right)  \ .
\end{eqnarray*}
\index{equilibrium distribution!$(\fX,\vs)$}
\index{equilibrium distribution}
        
{(Nonarchimedean Case.)} 
The proof is more complicated in this case, because of the absence of a
Laplacian operator on $\cC_v(\CC_v)$.   
(Actually, a suitable Laplacian has been defined in the  
context of Berkovich Spaces,\index{Berkovich!analytic space} 
for $\PP^1$ by Baker and Rumely (\cite{B-R}),\index{Baker, Matthew}\index{Rumely, Robert} 
and for curves of arbitrary genus by Thuillier (\cite{Th}).  
\index{Thuillier, Amaury} 
However, introducing that theory would take us too far afield.) 
Instead, we use the approximability of Green's functions by algebraic functions,
\index{Green's function}
the $[z,w]_{\fX,\vs}$-factorization of pseudopolynomials 
\index{pseudopolynomial!$(\fX,\vs)$}
(\ref{FGHB}), and the nonarchimedean Maximum modulus principle.   
 
First suppose $\vs \in \cP^m \bigcap \QQ^m$.  
Choose decreasing sequences of numbers $r_n > 0$ and $\varepsilon_n > 0$ with 
$\lim_{n \rightarrow \infty} r_n 
      = \lim_{n \rightarrow \infty} \varepsilon_n = 0$. 
For each $n$, put $W_n = \{ z \in \cC_v(\CC_v) : 
\|z,w\|_v \le r_n \ \text{for some $w \in H$} \}$.  Thus, the $W_n$ form a
decreasing sequence of neighborhoods of $H$ with  
$\bigcap_{n=1}^{\infty} W_n= H$. 

For each $i$ and $n$, Proposition 4.1.5 of (\cite{RR1}) provides a function 
$f^{(n)}_i(z) \in \CC_v(\cC)$, with poles only at $x_i$ and normalized so
that $|f^{(n)}_i(z)|_v \cdot |g_{x_i}(z)^N|_v = 1$ if $N = \deg(f_i^{(n)})$,  
such that for all $z \in \cC_v(\CC_v) \backslash (W_n \cup \{x_i\})$, 
\begin{equation} \label{BFF1A}
\left| u_{x_i}(z) - \Big( \frac{-1}{N} \log_v|f^{(n)}_i(z)|_v \Big) \right|_v 
         \ < \ \varepsilon_n \ .
\end{equation}
\index{potential function}
By raising the $f^{(n)}_i(z)$ to appropriate powers, we can assume without loss
that for a given $n$, they have common degree $N_n$.  Let the zeros of 
$f^{(n)}_i(z)$ (with multiplicity) be  $a_{ij}^{(n)}$, 
$j = 1, \ldots, N_n$.   Let $\nu^{(n)}_i$ be the probability measure 
which gives weight $1/N_n$ to each $a_{ij}^{(n)}$.  
By the construction in (\cite{RR1}, Proposition 4.1.5), for each $i$ the 
$\nu^{(n)}_i$ converge weakly to $\mu_i$.  
\index{equilibrium distribution}

The normalization of the $f^{(n)}_i(z)$ means that for each $n$ and $i$
\begin{equation*}
|f^{(n)}_i(z)|_v \ = \ \prod_{j=1}^{N_n} [z,a_{ij}^{(n)}]_{x_i} \ ,
\end{equation*}
and so 
\begin{equation} \label{BFFF1A}
u_{x_i}(z,\nu^{(n)}_i) 
    \ = \ \frac{1}{N_n} \cdot \sum_{j=1}^{N_n} -\log_v([z,a_{ij}^{(n)}]_{x_i})  
    \ = \ - \frac{1}{N_n} \log_v(|f^{(n)}_i(z)|_v) \ .
\end{equation}
\index{potential function}
In particular, for each $z \notin H$, by (\ref{BFF1A}), 
\begin{equation} \label{BFFF1B}
\lim_{n \rightarrow \infty} u_{x_i}(z,\nu^{(n)}_i) 
        \ = \ u_{x_i}(z,\mu_i) \ = \ u_{x_i}(z) \ .
\end{equation}        

Let $M$ be a common denominator for the $s_i$ and for each $n$ put 
\begin{equation} \label{BFFF2}
F^{(n)}(z) \ = \ \prod_{i=1}^m f^{(n)}_i(z)^{Ms_i} \ .
\end{equation}
Then $F^{(n)}(z)$ is an $(\fX,\vs)$-function in the sense of 
Definition \ref{FXvsFunctionDef}.  
Let the zeros  $F^{(n)}(z)$, listed with multiplicities, be 
$a^{(n)}_1, \ldots, a^{(n)}_{MN_n}$;
these are of course just the $a^{(n)}_{ij}$, repeated certain 
numbers of times.  By formula (\ref{FGHB})
there is a constant $C_n$ such that for all $z \in \cC_v(\CC_v)$ 
\begin{equation} \label{BFFF3}
|F^{(n)}(z)|_v 
    \ = \ C_n \cdot \prod_{j=1}^{MN_n} [z,a^{(n)}_j]_{\fX,\vs} \ .
\end{equation}
Let $\omega^{(n)}$ be the probability measure which gives mass $1/(MN_n)$ to 
each of the points $a^{(n)}_j$;  clearly 
\begin{equation*}
\omega^{(n)} \ = \ \sum_{i=1}^m s_i \nu^{(n)}_i \ .
\end{equation*} 
Combining (\ref{BFFF1A}), (\ref{BFFF2}) and (\ref{BFFF3}), we see that 
\begin{eqnarray}
u_{\fX,\vs}(z,\omega^{(n)}) 
& = & \int -\log_v([z,w]_{\fX,\vs}) \, d\omega^{(n)}(w) \notag \\
& = & \frac{1}{MN_n} \log_v(C_n) - \frac{1}{MN_n} \log_v(|F^{(n)}(z)|_v 
                           \notag \\
& = & \frac{1}{MN_n} \log_v(C_n) + 
    \sum_{i=1}^m s_i \left(-\frac{1}{N_n} \log_v(|f^{(n)}_i(z)|_v \right) 
                           \notag \\
& = & \frac{1}{MN_n} \log_v(C_n) + \sum_{i=1} s_i u_{x_i}(z,\nu^{(n)}_i) 
                           \ . \label{BFFF4} 
\end{eqnarray}
\index{equilibrium potential!$(\fX,\vs)$}

Since the $\nu^{(n)}_i$ converge weakly to the $\mu_i$, the $\omega^{(n)}$ 
converge weakly to $\omega = \sum s_i \mu_i$.  Hence, for all $z \notin H$, 
\index{equilibrium distribution}
\begin{equation} \label{BFFF5}
\lim_{n \rightarrow \infty} u_{\fX,\vs}(z,\omega^{(n)}) 
                    \ = \ u_{\fX,\vs}(z,\omega) \ .
\end{equation}
\index{equilibrium potential!$(\fX,\vs)$}
Comparing (\ref{BFFF1B}), (\ref{BFFF5}) and (\ref{BFFF4}) it follows that 
$C := \lim_{n \rightarrow \infty} \frac{1}{NM_n} \log_v(C_n)$ exists, and that 
for all $z \notin H$
\begin{equation} \label{FEQQ}
u_{\fX,\vs}(z,\omega) \ = \ C + \sum_{i=1} s_i u_{x_i}(z) \ .
\end{equation}
Using (\cite{RR1}, Lemma 4.1.3) and its $(\fX,\vs)$-analogue, 
we conclude that $\omega = \sum s_i \mu_i$.  Hence (\ref{FEQQ}) holds for
$z \in H$ as well.  
\index{equilibrium distribution}

By (\cite{RR1}, Theorem 4.1.11), for each $i$ there is a set $e_i \subset H$ 
of capacity $0$ such that $u_{x_i}(z) = V_{x_i}(H)$ for all
\index{equilibrium potential!takes constant value a.e. on $E_v$}
\index{Robin constant}
\index{capacity $= 0$} 
$z \in H \backslash e_i$.  Moreover, since 
$u_{x_i}(z) \le V_{x_i}(H)$ for all $z$, (\cite{RR1}, Lemma 4.1.9)
\index{Robin constant}
shows that $u_{x_i}(z)$ is continuous at each $z \in H \backslash e_i$.
\index{equilibrium potential}
Put $e = \bigcup_{i=1}^m e_i$.  
By (\cite{RR1}, Lemma 4.1.9),  for all $z \in H \backslash e$, 
\begin{eqnarray*}
u_{\fX,\vs}(z,\omega) 
    & = & \lim_{\substack{ w \rightarrow z \\ w \notin H }} 
                   ( \ C + \sum_{i=1} s_i u_{x_i}(w) ) \\  
    & = &  \ C + \sum_{i=1} s_i V_{x_i}(H) \ .
\end{eqnarray*}   
\index{equilibrium potential!$(\fX,\vs)$}  
\index{Robin constant}                        
Moreover, by (\cite{RR1}, Corollary 4.1.15) $e$ has inner capacity $0$. 
\index{capacity $= 0$} 
Thus $u_{\fX,\vs}(z,\omega)$ is constant on $H$, except on a set of
inner capacity $0$.   
\index{capacity $= 0$} 
\index{equilibrium potential!takes constant value a.e. on $E_v$}

By Proposition \ref{IDProp}, $\omega$ is the equilibrium 
distribution $\mu_{\fX,\vs}$ of $H$.  That is, 
\index{equilibrium distribution!$(\fX,\vs)$}  
\index{equilibrium distribution!determining}  
\begin{equation} \label{BFFF7A}
\mu_{\fX,\vs} \ = \ \sum_{i=1}^m s_i \mu_i \ .
\end{equation}
From this (\ref{BFFF7B}) and (\ref{BFFF8}) follow at once.   
         
Now consider the general case where possibly $\vs \notin \QQ^m$.  
Let $f(\vs)$ be the function given by the right side of (\ref{BFFF8}). 
Fix $\vs$, and choose a sequence of probability vectors $\vs^{(n)} \in \QQ^m$ 
approaching $\vs$.  Let $\mu_0 := \mu_{\fX,\vs}$ be the $(\fX,\vs)$-equilibrium 
\index{equilibrium distribution!$(\fX,\vs)$}
distribution of $H$.  Then
\begin{eqnarray}
V_{\fX,\vs}(H) & = & \iint_{H \times H} -\log_v([z,w]_{\fX,\vs}) \, 
                            d\mu_0(z) d\mu_0(w) \notag \\ 
                 & = & \sum_{i=1}^m s_i 
                         \iint_{H \times H} -\log_v([z,w]_{x_i}) \, 
                            d\mu_0(z) d\mu_0(w) \ . \label{BFFF9}
\end{eqnarray}
\index{Robin constant!$(\fX,\vs)$}
Suppose $V_{\fX,\vs}(H) < f(\vs)$.  By the continuity $f$ and 
of the right side of (\ref{BFFF9}) in $\vs$, for sufficiently large $n$
\begin{equation*}
\iint_{H \times H} -\log_v([z,w]_{\fX,\vs^{(n)}}) \, 
                            d\mu_0(z) d\mu_0(w)  \ < \ f(\vs^{(n)}) \ .
\end{equation*}
However, this contradicts that $f(\vs^{(n)}) = V_{\fX,\vs^{(n)}}(H)$ is the
\index{Robin constant!$(\fX,\vs)$}
minimal value of the energy integral for  $[z,w]_{\fX,\vs^{(n)}}$.\index{energy integral!$(\fX,\vs)$}  

Consequently $V_{\fX,\vs}(H) \ge f(\vs)$.  But $\sum s_i \mu_i$ is 
a probability measure for which the $(\fX,\vs)$-energy integral equals $f(\vs)$.\index{energy integral!$(\fX,\vs)$} 
Hence this is the minimal value of the energy integral, and by the 
uniqueness of the $(\fX,\vs)$-equilibrium distribution
\index{equilibrium distribution!$(\fX,\vs)$}
\index{equilibrium distribution}
$\mu_{\fX,\vs} \ = \ \sum_{i=1}^m s_i \mu_i$.  
\end{proof} 


\section{ Mass Bounds in the Archimedean case }
\label{ArchMassBoundsSection} 

Throughout this section we assume $K_v$ is archimedean, and we identify
$\CC_v$ with $\CC$.  The results proved here will be used in Theorem \ref{RThm2},
the construction of the initial local approximating functions when $K_v \cong \RR$.
\index{initial approximating functions $f_v(z)$!archimedean}

Suppose $H \subset \cC_v(\CC) \backslash \fX$ can be decomposed 
as $H = H_1 \cup e$, where $H_1$ and $e$ are closed and disjoint;  
we think of $e$ as being ``small''.  The case of interest is when 
$K_v \cong \RR$, $H \subset \cC_v(\RR)$, and $e$ is a short interval.  

If  $\mu_{\fX,\vs}$ is the $(\fX,\vs)$-equilibrium distribution of $H$, 
\index{equilibrium distribution!$(\fX,\vs)$}
we seek upper and lower bounds for $\mu_{\fX,\vs}(e)$ in terms of 
the Robin constants and Green's functions of $H_1$ and $e$.
\index{Robin constant}  
\index{Green's function}
\index{mass bounds|ii}
    
\begin{lemma} \label{bL3} 
Let $H \subset \cC_v(\CC) \backslash \fX$ be compact with positive capacity.
\index{capacity $> 0$} 
Let $M$ be a constant such that for each $x_i \in \fX$, 
\begin{equation*}
\max_{z,w \in H} [z,w]_{x_i} \ < \ M \ ;
\end{equation*} 
put $C = \log(M)$.

Suppose $H = H_1 \bigcup e$ where $H_1$ and $e$ are closed and disjoint.
Given $\vs \in \cP^m$, let $\mu_{\fX,\vs}$ be the 
$(\fX,\vs)$-equilibrium distribution of $H$.  Then 
\index{equilibrium distribution!$(\fX,\vs)$}
\begin{equation}   \label{AFJ122} 
\mu_{\fX,\vs}(e)  \ \le \ 
       \frac{V_{\fX,\vs}(H) + C}{V_{\fX,\vs}(e) + C} 
\end{equation}
\index{Robin constant!$(\fX,\vs)$}
\end{lemma}             

\begin{proof}  
By our hypothesis, $[z,w]_{\fX,\vs} < M$ for all $z, w \in H$.  

If $M = 1$, by the same 
argument as in the proof of (\cite{RR1}, formula (15), p.148) one obtains 
\begin{equation*}
\mu_{\fX,\vs}(e) \ \le \ \frac{V_{\fX,\vs}(H)}{V_{\fX,\vs}(e)} \ .
\end{equation*}
\index{Robin constant!$(\fX,\vs)$}
In the general case, if we 
renormalize $[z,w]_{\fX,\vs}$ by replacing it with
$\frac{1}{M} [z,w]_{\fX,\vs}$,   
then for each compact set $X \subset H$, 
$V_{\fX,\vs}(X)$ is replaced by $V_{\fX,\vs}(X) + \log(M)$. 
\index{Robin constant!$(\fX,\vs)$}
The result follows.  
\end{proof} 

\medskip
To obtain a lower bound, we need information about the potential-theoretic
\index{potential theoretic separation}
separation between $H_1$ and $e$.  

\begin{lemma} \label{ALJ29}  
Let $H = H_1 \bigcup e \subset \cC_v(\CC) \backslash \fX$, and the constants 
$M$ and $C$, be as in Lemma $\ref{bL3}$.  Let  $m > 0$ be such that for each 
$x_i \in \fX$, and all $z \in e$, 
\begin{equation} \label{AFJ123}
G(z,x_i;H_1) \ \ge \ m \ .
\end{equation}
\index{Green's function}

Fixing $\vs$, let $\mu_{\fX,\vs}$ be the $(\fX,\vs)$-equilibrium distribution of $H$.
\index{equilibrium distribution!$(\fX,\vs)$}   
\index{Robin constant!$(\fX,\vs)$}
\index{Robin constant}
Suppose  $V_{\fX,\vs}(e) \ge V_{\fX,\vs}(H_1)$.  Then  
\begin{equation}                \label{AFJ124}
\mu_{\fX,\vs}(e) \ \ge \ 
     \frac{m^2}{2 (V_{\fX,\vs}(H_1) + C)(V_{\fX,\vs}(e) + C + 2m)} \ .
\end{equation}
\end{lemma} 

\begin{proof}
If $e$ has capacity $0$, the result is trivial since 
\index{capacity $= 0$} 
$\mu_{\fX,\vs}(e) = 0$ and $V_{\fX,\vs}(e) = \infty$.  Hence
\index{Robin constant!$(\fX,\vs)$} 
without loss we can assume that $e$ has positive capacity.  
\index{capacity $> 0$} 
First suppose  $M = 1$.  

Write $\mu$ (resp. $\mu_1$, resp. $\mu_2$) 
for the $(\fX,\vs)$-equilibrium distribution of $H$ (resp. $H_1$, resp. $e$). 
\index{equilibrium distribution!$(\fX,\vs)$}
Put $V_1 = V_{\fX,\vs}(H_1) = I_{\fX,\vs}(\mu_1)$,
\index{Robin constant!$(\fX,\vs)$}
$V_2 = V_{\fX,\vs}(e) = I_{\fX,\vs}(\mu_2)$, and let    
\begin{equation*}
I_{\fX,\vs}(\mu_1,\mu_2) 
\ = \ \iint - \log([z,w]_{\fX,\vs}) \, d\mu_1(z) d\mu_2(w) \ . 
\end{equation*}
Then for each $0 \le t \le 1$,  
\begin{eqnarray}
V_{\fX,\vs}(H) & = & I_{\fX,\vs}(\mu) \ \le \ I_{\fX,\vs}((1-t)\mu_1+t\mu_2)  \notag   \\
   & = &  (1-t)^2 I_{\fX,\vs}(\mu_1) + 2t(1-t) I_{\fX,\vs}(\mu_1,\mu_2) 
                + t^2 I_{\fX,\vs}(\mu_2) \ .  \label{AFJ125}
\end{eqnarray} 
\index{Robin constant!$(\fX,\vs)$}
By our hypothesis, 
$G_{\fX,\vs}(z;H_1) = \sum_{i=1}^m s_i G(z,x_i;H_1) \ge m$ on $e$; 
\index{Green's function}
\index{Green's function!$(\fX,\vs)$}
hence the potential function 
\index{potential function!$(\fX,\vs)$}
\begin{equation*}
u_{\fX,\vs}(z;H_1) 
     \ = \  V_1 - G_{\fX,\vs}(z;H_1) 
\end{equation*}
satisfies $u_{\fX,\vs}(z;H_1) \ \le \  V_1 - m$ on $e$. 
Thus
\begin{equation*}
I_{\fX,\vs}(\mu_1,\mu_2) \ = \ \int_e u_{\fX,\vs}(z;H_1) d\mu_2(z) 
               \ \le \ V_1 - m \ . 
\end{equation*}               
Since $I_{\fX,\vs}(\mu_1) = V_1$ and $I_{\fX,\vs}(\mu_2) = V_2$, (\ref{AFJ125}) gives
\begin{equation}   \label{AFJ126} 
V_{\fX,\vs}(H) \ \le \ (V_2 + 2m -V_1)t^2 -2mt + V_1 \ .
\end{equation}     
\index{Robin constant!$(\fX,\vs)$}         
The minimum of the right side occurs at $t = m/(V_2 + 2m - V_1)$, 
which lies in the interval $[0,1]$ because of our assumption that $V_2 \ge V_1$.
Inserting this in (\ref{AFJ126}), we get
\begin{equation} \label{AFJ127}                   
V_{\fX,\vs}(H) \ \le \ V_1 - \frac{m^2}{V_2 + 2m - V_1} \ .
\end{equation}
\index{Robin constant!$(\fX,\vs)$}

Put $\beta = \mu_{\fX,\vs}(e)$.  Because $-\log([z,w]_{\fX,\vs}) \ge 0$ on $H$, we have  
\begin{eqnarray*}
V_{\fX,\vs}(H) 
    &=& \iint_{H \times H} -\log([z,w]_{\fX,\vs}) \,d\mu(z) d\mu(w) \\
    &\ge& \iint_{H_1 \times H_1} -\log([z,w]_{\fX,\vs}) \,d\mu(z) d\mu(w)
    \ = \ I_{\fX,\vs}(\mu|_{H_1}) \ .
\end{eqnarray*}
\index{Robin constant!$(\fX,\vs)$}
Since  $\frac{1}{1-\beta} \mu|_{H_1}$ is a probability measure on $H_1$, 
upon dividing by $(1-\beta)^2$, we get  
\begin{equation}
\frac{V_{\fX,\vs}(H)}{(1-\beta)^2} 
\ \ge \ I_{\fX,\vs}(\frac{1}{1-\beta} \mu|_{H_1}) \ \ge \ I_{\fX,\vs}(\mu_1) 
\ = \ V_1 \ .  \label{AFJ128}
\end{equation}
\index{Robin constant!$(\fX,\vs)$}
Combining (\ref{AFJ127}) and (\ref{AFJ128}) gives 
\begin{equation*}
(1-\beta)^2 \ \le \ 1 - \frac{m^2}{V_1(V_2 + 2m - V_1)} 
\ \le \ 1 - \frac{m^2}{V_1(V_2 + 2m)} \ .
\end{equation*}
Taking square roots, we see that 
\begin{equation*}
1 - \beta \ \le \ 1 -  \frac{m^2}{2V_1(V_2 + 2m)} \ ,
\end{equation*} 
which is equivalent to (\ref{AFJ124}).  

The general case follows upon scaling $[z,w]_{\fX,\vs}$ by $1/M$.  
\end{proof}

\vskip .1 in
Fix a local coordinate patch $U \subset \cC_v(\CC) \backslash \fX$,  
with coordinate function $z$ say.  We can describe subsets of $U$, 
such as intervals or discs, in terms of the coordinate function $z$. 
Our last result concerns the behavior of the $(\fX,\vs)$-Robin constant
\index{Robin constant!$(\fX,\vs)$}
of an interval, as its length goes to $0$.    

\begin{lemma} \label{bL2}  Let $U \subset \cC_v(\CC) \backslash \fX$
be a local coordinate patch, and let $I$ be a compact subset of $U$.  
Then there is a constant $A$ depending only $I$ 
$($and the choice of the local coordinate function $z)$
 such that for any interval $e_a(h) = [a-h,a+h] \subset I$, 
and any probability vector $\vs \in \cP^m$, we have 
\begin{equation} \label{bF21}
-\log(h) - A \ \le \ V_{\fX,\vs}(e_a(h)) \ \le \ -\log(h) + A \ .
\end{equation}
\index{Robin constant!$(\fX,\vs)$}
\end{lemma}

\begin{proof}  There is a constant $A_0$ such that for each 
$x_i \in \fX$, and all $z, w \in I$ with $z \ne w$,  
$-\log(|z-w|) - A_0  \le -\log( [z,w]_{x_i})  \le  -\log(|z-w|) + A_0$.
 Hence for all $z, w \in I$ with $z \ne w$, and all $\vs \in \cP^m$.  
\begin{equation*}
-\log(|z-w|) - A_0 \ \le \ -\log( [z,w]_{\fX,\vs}) 
                   \ \le \ -\log(|z-w|) + A_0 \ .
\end{equation*} 

Fix $\vs$, and let $\mu = \mu_{\fX,\vs}$ be the 
$(\fX,\vs)$-equilibrium distribution of $e_a(h)$.  
\index{equilibrium distribution!$(\fX,\vs)$}
Also, let $\mu_0$ be the equilibrium
distribution of $e_a(h)$ considered as a subset of $\CC$, via the local coordinate
function $z$.   Then by the energy minimizing property of $\mu_0$ and the
fact that the classical Robin constant of a segment of length $L$ is $-\log(L/4)$,  
\index{Robin constant} 
\begin{eqnarray*} 
V_{\fX,\vs}(e_a(h)) 
  &=& \iint_{e_a(h) \times e_a(h)} -\log( [z,w]_{\fX,\vs}) \, d\mu(z) d\mu(w) \\
  &\ge& \iint_{e_a(h) \times e_a(h)} (-\log(|z-w|) - A_0)  \, d\mu(z) d\mu(w) \\
  &\ge& \iint_{e_a(h) \times e_a(h)} -\log(|z-w|) \, d\mu_0(z) d\mu_0(w)  - A_0 \\
  &=& -\log(h/2)  - A_0  \ .
\end{eqnarray*}
\index{Robin constant!$(\fX,\vs)$}
Similarly, using the energy minimizing property of $\mu_{\fX,\vs}$, 
\index{equilibrium distribution!$(\fX,\vs)$}
\begin{eqnarray*}
-\log(h/2)  + A_0  
   &=& \iint_{e_a(h) \times e_a(h)} (-\log(|z-w|) + A_0) \, d\mu_0(z) d\mu_0(w) \\ 
   &\ge& \iint_{e_a(h) \times e_a(h)} -\log( [z,w]_{\fX,\vs}) \, d\mu(z) d\mu(w) \\
   &=& V_{\fX,\vs}(e_a(h)) \ .
\end{eqnarray*}
\index{Robin constant!$(\fX,\vs)$}
Putting $A = A_0 + \log(2)$, we obtain the result.  
\end{proof} 


\section{ Description of $\mu_{\fX,\vs}$ in the Nonarchimedean Case }
\label{NonArchMuXSSection} 

Throughout this section we assume that $K_v$ is nonarchimedean.  
Our goal is to determine $\mu_{\fX,\vs}$ for a class of well-behaved compact sets.
\index{equilibrium distribution!$(\fX,\vs)$}
\index{equilibrium distribution!determining}
The results proved here will be used in Theorem \ref{DCPCPatch},  
the construction of the initial local approximating functions 
\index{initial approximating functions $f_v(z)$!nonarchimedean}
in the nonarchimedean compact case.

For the remainder of this section, $F_w \subset \CC_v$ will be a fixed  
finite extension of $K_v$, with ramification index $e = e_{w/v}$ 
and residue degree $f = f_{w/v}$.  

\vskip .1 in
We begin by considering the special case when $\cC_v = \PP^1/K_v$
and $\zeta = \infty$.  Identify $\PP_v^1(\CC_v) \backslash \{\infty\}$ 
with $\CC_v$ and 
normalize the canonical distance so that $[z,x]_{\infty} = |z-x|_v$.
\index{canonical distance!$[z,w]_{\zeta}$!normalization of}
The equilibrium distribution and potential function can be determined 
\index{equilibrium distribution!$(\fX,\vs)$}
\index{potential function!$(\fX,\vs)$}
\index{equilibrium potential!$(\fX,\vs)$}
explicitly when $H$ is a coset of $\cO_w$.  The following result is 
a mild generalization of (\cite{RR1}, Example 4.1.24):   

\begin{lemma} \label{BFLem2}  
Let $F_w/K_v$ be a finite extension with ramification index $e$
and residue degree $f$.  Suppose $H = a + b \cO_w$, where $a \in \CC_v$ 
and $b \in \CC_v^{\times}$.  Then the equilibrium distribution $\mu$ 
\index{equilibrium distribution}
of $H$ relative to $[x,y]_{\infty} = |x-y|_v$ is the pushforward of   
additive Haar measure\index{Haar measure} on $\cO_w$ by the affine map $x = a+bz$,   
and if $|b|_v = r$, the potential function 
$u_{\infty}(x) = u_{\infty}(x,H)$ is given by 
\index{potential function}\index{equilibrium potential}
\begin{equation} \label{BFFM4A}
u_{\infty}(x) \ = \ \left\{
   \begin{array}{ll} -\log_v(r) + \frac{1}{e(q_v^f-1)} 
                    & \text{for $x \in H$,} \\
                     -\log_v(|x-a|_v) & \text{for $x \notin D(a,r)$\ .}
    \end{array} \right. 
\end{equation}
\end{lemma}

\begin{proof}
After a change of coordinates, we can assume without loss that $H = \cO_w$. 
Since $H$ and $|x-y|_v$ are invariant under translation by elements 
of $\cO_w$, the uniqueness of the equilibrium distribution 
\index{equilibrium distribution}
shows that it must be translation-invariant as well.  
It follows that $\mu$ is the additive Haar measure $\mu_w$ on $\cO_w$.  
\index{Haar measure}

Let $\pi_w$ be a generator for the maximal ideal $\fm_w$ of $\cO_w$.  
We can compute $u_{\infty}(0)$ directly:  
\index{potential function}
\index{equilibrium potential}
\index{equilibrium potential!determining}
\begin{eqnarray}
u_{\infty}(0) & = & \sum_{\ell = 0}^{\infty}
      \int_{\fm_w^{\ell} \backslash \fm_w^{\ell+1}} -\log(|0-y|_v) \, d\mu(y) \\
      & = & \sum_{\ell = 0}^{\infty} \ell \cdot \frac{1}{e} 
                \left(\frac{1}{q_v^{f \ell}} - \frac{1}{q_v^{f(\ell+1)}} \right)
      \ = \ \frac{1}{e(q_v^f - 1)} \ .
\end{eqnarray}
By translation invariance, $u_{\infty}(x) = u_{\infty}(0)$ for all $x \in \cO_v$. 
  
For $x \notin D(0,1)$, the ultrametric inequality gives 
\begin{eqnarray*} 
u_{\infty}(x) & = & \int_{\cO_v} -\log(|x-y|_v) \, d\mu(y)  \\
        & = & \int_{\cO_v} -\log(|x|_v) \, d\mu(y) \ = \ -\log_v(|x|_v) 
\end{eqnarray*}
It is not hard to give a formula for $u_{\infty}(x)$ 
when $x \in D(0,1) \backslash \cO_w$ (see \cite{RR1}, Example 4.1.26), 
but we will not need this.   
\end{proof}

\vskip .1 in
Now let $\cC_v/K_v$ be arbitrary.
Suppose $a \in \cC_v(F_w)$, and let $B(a,r)$ be an isometrically
parametrizable ball disjoint from $\fX$, whose radius $r$ belongs
to $|F_w^{\times}|_v$.   Take $H = \cC_v(F_w) \bigcap B(a,r)$.  
By Theorem \ref{IsoParamThm}, 
there is an $F_w$-rational isometric parametrization
$\Lambda : D(0,r) \rightarrow B(a,r)$ with $\Lambda(0) = a$, 
and if $b \in F_w^{\times}$ is such that $|b|_v = r$, 
then  $\Lambda(b \cO_w) = H$.   
  
Fix $\vs \in \cP^m$.  By Proposition \ref{APropA2}
there is a constant $C_a(\vs)$, which belongs to the value group of
$\CC_v^{\times}$ if $\vs \in \cP^m \cap \QQ^m$, 
such that $[z,w]_{\fX,\vs} = C_a(\vs) \|z,w\|_v$
for all $z,w \in B(a,r)$.   
Using Lemma \ref{BFLem2}, we obtain:  

\begin{corollary} \label{BFCor1}
Suppose  $a \in \cC_v(F_w)$, $r \in |F_w^{\times}|_v$, 
and that $B(a,r) \subset \cC_v(\CC_v)$
is an isometrically parametrizable ball disjoint from $\fX$.  
Let $\vs \in \cP^m$ be arbitrary, and take  $H = \cC_v(F_w) \bigcap B(a,r)$.

Then the $(\fX,\vs)$-equilibrium distribution of $H$ 
\index{equilibrium distribution!$(\fX,\vs)$}
is the pushforward $\Lambda_*(\mu_w)$
of additive Haar measure on $b \cO_w$, normalized to have total mass $1$, 
\index{Haar measure} 
and the $(\fX,\vs)$-potential function of $H$ satisfies
\index{equilibrium potential!$(\fX,\vs)$}
\begin{equation} \label{BFFM4B}
u_{\fX,\vs}(z,H) \ = \ \left\{
   \begin{array}{ll}  - \log_v(C_a(\vs) \cdot r) + \frac{1}{e(q_v^f - 1)}
                          & \text{for all $z \in H$,} \\
                      -\log_v([z,a]_{\fX,\vs}) & \text{for all $z \notin B(a,r)$.}
   \end{array} \right. 
\end{equation}
\end{corollary}

\begin{proof} Write $z = \Lambda(x)$, $w = \Lambda(y)$
for $z, w \in B(a,r)$ and $x, y \in D(0,r)$.  Then $\|z,w\|_v = |x-y|_v$, and 
$[z,w]_{\fX,\vs} = C_a(\vs)|x-y|_v$.

Write $\mu_w$ for the additive Haar measure on $b \cO_w = F_w \cap D(0,r)$, 
\index{Haar measure} 
\index{equilibrium distribution}
normalized to have total mass $1$, and put $\mu_0 = \Lambda_*(\mu_w)$.  
By Lemma \ref{BFLem2}, for each $z \in H$
\begin{eqnarray*}
u_{\fX,\vs}(z,\mu_0) & = & \int_H -\log_v([z,w]_{\fX,\vs}) \, d\mu_0(w) \\ 
   & = & \int_{b \cO_w} -\log_v(C_a(\vs) |x-y|_v) \, d\mu_w(y) \\
     & = & -\log_v(C_a(\vs)) - \log_v(r) + \frac{1}{e(q_v^f-1)} \ .
\end{eqnarray*}  
\index{equilibrium potential!$(\fX,\vs)$}
In particular, $u_{\fX,\vs}(z,\mu_0)$ is constant on $H$. 
\index{equilibrium potential!takes constant value a.e. on $E_v$}
By Proposition \ref{IDProp}, it follows that $\mu_{\fX,\vs} = \mu_0$.  

If $z \notin B(a,r)$ then Proposition \ref{APropA2}
gives $[z,w]_{\fX,\vs} = [z,a]_{\fX,\vs}$ 
for all $w \in B(a,r)$, so
\begin{equation*}
u_{\fX,\vs}(z,\mu_0) \ = \ \int_H -\log_v([z,w]_{\fX,\vs}) \, d\mu_0(w) 
     \ = \  -\log_v([z,a]_{\fX,\vs}) \ .
\end{equation*}
\index{equilibrium potential!$(\fX,\vs)$}
\end{proof} 

\vskip .1 in
Note that in Corollary \ref{BFCor1}, the equilibrium measure $\mu_{\fX,\vs}$
\index{equilibrium distribution!$(\fX,\vs)$}
of $H$ is independent of $\vs$.  This is a general phenomenon for compact subsets
of isometrically parametrizable balls:  

\begin{lemma} \label{BFLem1}  
Let $B(a,r) \subset \cC_v(\CC_v)$ be an isometrically parametrizable ball
disjoint from $\fX$, and let $H \subset B(a,r)$ be compact with positive capacity. 
\index{capacity $> 0$} 
Then the equilibrium distribution $\mu_{\fX,\vs}$ of $H$ is 
\index{equilibrium distribution!$(\fX,\vs)$}
a probability measure $\mu^*$ independent of $\vs$.
\end{lemma}

\begin{proof} For a given $\vs$, the equilibrium distribution $\mu_{\fX,\vs}$
\index{equilibrium distribution!$(\fX,\vs)$}
is the unique probability measure $\mu$ supported on $H$ which minimizes the
energy integral\index{energy integral!$(\fX,\vs)$} 
\begin{equation*}
I_{\fX,\vs}(\mu) 
\ = \ \iint_{H \times H} -\log_v([z,w]_{\fX,\vs}) \, d\mu(z) d\mu(w) \ .
\end{equation*}
Fix an isometric parametrization of $B(a,r)$.  
By Proposition \ref{APropA2}, for each $\vs$ there is 
a constant $C_a(\vs)$ such that for all $z,w \in B(a,r)$, 
\begin{equation*}
-\log_v([z,w]_{\fX,\vs}) \ = \ -\log_v(\|z,w\|_v) - \log_v(C_a(\vs)) \ .
\end{equation*}
Hence the same measure $\mu^*$ minimizes the energy integral,\index{energy integral!$(\fX,\vs)$}  for all $\vs$.
\end{proof}  

\vskip .1 in
Now let $B(a_{\ell},r_{\ell})$ for $\ell = 1, \ldots, D$ be isometrically
parametrizable balls in $\cC_v(\CC_v)$, disjoint from each other and from $\fX$.
Suppose $H = \bigcup_{\ell=1}^D H_{\ell}$, 
where $H_{\ell} \subset B(a_{\ell},r_{\ell})$ 
is compact and has positive capacity for each $\ell$.
\index{capacity $> 0$}   
Let $\mu_{\ell}^*$ be the $(\fX,\vs)$-equilibrium
\index{equilibrium distribution!$(\fX,\vs)$}
distribution of $H_{\ell}$, which is independent of $\vs$ by Lemma \ref{BFLem1}.  
Let 
\begin{equation*}
u_{\fX,\vs}(z,H_{\ell}) 
\ = \ \int_{H_{\ell}} -\log([z,w]_{\fX,\vs}) \, d\mu_{\ell}^*(w)
\end{equation*} 
be the $(\fX,\vs)$-potential function of $H_{\ell}$.
\index{equilibrium potential!$(\fX,\vs)$}

\begin{proposition} \label{BFProp2}
Let $H = \bigcup_{\ell=1}^D H_{\ell}$ be as above.  For each $\vs$,
there are weights $w_{\ell}(\vs) > 0$ 
with $\sum_{\ell=1}^D w_{\ell}(\vs) = 1$ such that the equilibrium distribution
$\mu_{\fX,\vs}$ of $H$ satisfies 
\index{equilibrium distribution!$(\fX,\vs)$}
\begin{equation} \label{BFFM1}
\mu_{\fX,\vs} \ = \ \sum_{\ell=1}^D w_{\ell}(\vs) \, \mu_{\ell}^* \ ,
\end{equation}
and the potential function is given by 
\index{equilibrium potential!$(\fX,\vs)$}
\begin{equation} \label{BFFM2}
u_{\fX,\vs}(z,H) \ = \ \sum_{\ell=1}^D w_{\ell}(\vs) \, u_{\fX,\vs}(z,H_{\ell}) \ .
\end{equation}
\end{proposition}

\begin{proof}
Let $w_{\ell}(\vs) = \mu_{\fX,\vs}(H_{\ell})$.  
Since $H_{\ell}$ has positive capacity,
\index{capacity $> 0$} 
necessarily $w_{\ell}(\vs) > 0$ (see \cite{RR1}, Lemma 4.1.7).  
Since the balls $B(a_{\ell},r_{\ell})$ are pairwise disjoint and do not meet $\fX$,
Proposition \ref{APropA2} shows that $[z,w]_{\fX,\vs}$ is constant 
for $z \in B(a_{\ell},r_{\ell})$, $w \in B(a_k,r_k)$, if $\ell \ne k$.  
Hence the same arguments as in (\cite{RR1}, Proposition 4.1.27) 
yield (\ref{BFFM1}) and (\ref{BFFM2}).  
\end{proof}

\vskip .1 in
\noindent{\bf Remark.}
Using Proposition \ref{BPropF1}, one sees that there are constants $W_{i \ell} > 0$
such that $w_{\ell}(\vs) = \sum_{i=1}^m s_i W_{i \ell}$ 
for all $\vs$ and all $\ell$.

\begin{theorem} \label{BFThm2A} 
Suppose $B(a_1,r_1), \ldots, B(a_D,r_D)$  
are pairwise disjoint isometrically parametrizable balls in $\cC_v(\CC_v)$,
whose union is disjoint from $\fX$.  For each $\ell$, let $H_{\ell} \subset B(a_{\ell},r_{\ell})$
be a compact set of positive capacity,  
\index{capacity $> 0$} 
and let $H = \bigcup_{\ell=1}^D H_{\ell}$.

Given $\vs \in \cP^m$, 
let $V_{\fX,\vs}(H_{\ell})$ be the $(\fX,\vs)$-Robin constant of $H_{\ell}$.
\index{Robin constant!$(\fX,\vs)$}
Let $\mu_{\ell}^*$ be the $(\fX,\vs)$-equilibrium distribution of $H_{\ell}$ 
\index{equilibrium distribution!$(\fX,\vs)$}
$($which is independent of $\vs$, by Lemma $\ref{BFLem1})$.   
 
Then the Robin constant $V = V_{\fX,\vs}(H)$ 
\index{Robin constant!$(\fX,\vs)$}
and the weights $w_{\ell} = w_{\ell}(\vs) = \mu_{\fX,\vs}(H_{\ell}) > 0$
such that $\mu_{\fX,\vs} = \sum_{\ell=1}^D w_{\ell}(\vs) \, \mu_{\ell}^*$ 
$($given by Proposition $\ref{BFProp2})$   
are uniquely determined by the $D+1$ linear equations
\begin{equation} \label{BFFM6A}
\left\{ \begin{array}{ccl}
        1  & = & 0 \cdot V + \sum_{\ell=1}^D w_{\ell} \ ,     \\   
0 & = & V + w_j \cdot \big(-V_{\fX,\vs}(H_j)\big) 
           + \displaystyle{\sum^D_{\substack{ \ell=1 \\ \ell \ne j }}} 
            w_{\ell}  \cdot \log_v \big([a_{\ell},a_j]_{\fX,\vs}\big)  \\
 &  & \qquad \qquad \text{for $j = 1, \ldots, D$.}  
          \end{array} \right.
\end{equation}
\index{Robin constant!$(\fX,\vs)$}
If $\vs \in \cP^m \cap \QQ^m$ and $V_{\fX,\vs}(H_{\ell}) \in \QQ$ for each $\ell$,
then $V_{\fX,\vs}(H)$ and the $w_{\ell}(\vs)$ belong to $\QQ$.
\end{theorem}
 
\begin{proof}
By Theorem \ref{ATE10B} and Proposition \ref{IDProp}, $u_{\fX,\vs}(z,H)$
\index{equilibrium potential!$(\fX,\vs)$}
\index{equilibrium potential!takes constant value a.e. on $E_v$}
takes the constant value $V$ on $H$,  
except possibly on an exceptional set $e \subset H$ of inner capacity $0$.
\index{exceptional set}
\index{capacity $= 0$} 
Similarly, each $u_{\fX,\vs}(z,H_j)$ 
takes the constant value $V_{\fX,\vs}(H_j)$ on $H_j$, 
\index{Robin constant!$(\fX,\vs)$}
except possibly on an exceptional set $e_j \subset H_j$ of inner capacity $0$.
\index{exceptional set}
\index{capacity $= 0$}   
Since $H_j$ has positive capacity, $H_j \backslash (e \cup e_j)$ is nonempty.
\index{capacity $> 0$} 
For each $j$, let $a_j^* \in H_j$ be a point where $u_{\fX,\vs}(a_j^*,H) = V$
and $u_{\fX,\vs}(a_j^*,H_j) = V_{\fX,\vs}(H_j)$.
\index{Robin constant!$(\fX,\vs)$}
The first equation in (\ref{BFFM6A}) follows from Proposition \ref{BFProp2}.
Using Corollary \ref{BFCor1} and evaluating $u_{\fX,\vs}(H)$ at each $a_j^*$, 
we obtain the last $D$ equations in (\ref{BFFM6A}) with the $a_{\ell}, a_j$ 
\index{equilibrium potential!determining}
replaced by the $a_{\ell}^*, a_j^*$.
However, if $\ell \ne j$, then $\log_v([x,y]_{\fX,\vs})$ 
is constant for $(x,y) \in B(a_{\ell},r_{\ell}) \times B(a_j,r_j)$ by Proposition \ref{APropA2}. 
Hence $\log_v([a_{\ell}^*,a_j^*]_{\fX,\vs}) = \log_v([a_{\ell},a_j]_{\fX,\vs})$. 

Conversely, we claim that the system of linear equations (\ref{BFFM6A}) 
in the variables $V$ and $w_{\ell}$ is nonsingular.  
To see this, first note that the values 
$V = V_{\fX,\vs}(H)$ and $w_{\ell} = \mu_{\fX,\vs}(H_{\ell})$ 
\index{Robin constant!$(\fX,\vs)$}
provide one solution to this system, with positive $w_{\ell}$. 
On the other hand, any solution to the system, with positive $w_{\ell}$, 
determines a probability measure on $H$ having the properties of the equilibrium distribution. 
\index{equilibrium distribution!$(\fX,\vs)$}
If the system were singular, there would be other solutions arbitrarily close to the one given 
above, contradicting the uniqueness of the equilibrium distribution.  Thus the equations
(\ref{BFFM6A}) uniquely determine $V$ and the $w_{\ell}$.

If $\vs \in \cP^m \cap \QQ^m$ and the $V_{\fX,\vs}(H_j) \in \QQ$, 
\index{Robin constant!$(\fX,\vs)$}
then the coefficients of the linear equations are rational, since for all $\ell \ne j$ we have
$\log_v([a_{\ell},a_j]_{\fX,\vs}) = \sum_{i=1}^m s_i \log_v([a_{\ell},a_j]_{x_i}) \in \QQ$.    
Hence $V$ and the $w_{\ell}$ must be rational as well.
\end{proof} 
 
\vskip .1 in
We now apply the preceding results to sets $H = \bigcup_{\ell=1}^D H_{\ell}$ of a special form: 

\begin{corollary} \label{BFCor2}  
Let $B(a_1,r_1), \ldots, B(a_D,r_D) \subset \cC_v(\CC_v)$  
be pairwise disjoint isometrically parametrizable balls 
whose union is disjoint from $\fX$. For each $\ell$, let $F_{w_\ell}/K_v$ be a finite 
extension in $\CC_v$, with residue degree $f_{\ell} = f_{w_{\ell}/v}$ and ramification index
$e_{\ell} = e_{w_{\ell}/v}$.  Assume that each $a_{\ell} \in \cC_v(F_{w_{\ell}})$ 
and each $r_{\ell} \in |F_{w_{\ell}}^{\times}|_v$, 
and put $H_{\ell} = \cC_v(F_{w_\ell}) \bigcap B(a_{\ell},r_{\ell})$. 
Let $H = \bigcup_{\ell=1}^D H_{\ell}$.

Given $\vs \in \cP^m$, 
let $C_{a_{\ell}}(\vs)$
be the constant such that $[z,w]_{\fX,\vs} = C_{a_{\ell}}(\vs) \cdot \|z,w\|_v$
for $z,w \in B(a_{\ell},r_{\ell})$.  
Let $\mu_{\ell}^*$ be the $(\fX,\vs)$-equilibrium distribution of $H_{\ell}$
\index{equilibrium distribution!$(\fX,\vs)$} 
$($which is independent of $\vs$, by Lemma $\ref{BFLem1}$, and is given by
a pushforward of additive Haar measure\index{Haar measure} on $F_w \cap D(0,r_{\ell})$, normalized to have
mass $1$, by Corollary $\ref{BFCor1})$.  
 
Then the Robin constant $V = V_{\fX,\vs}(H)$ 
\index{Robin constant!$(\fX,\vs)$}
and the weights $w_{\ell} = w_{\ell}(\vs) = \mu_{\fX,\vs}(H_{\ell}) > 0$
such that $\mu_{\fX,\vs} = \sum_{\ell=1}^D w_{\ell} \, \mu_{\ell}^*$ 
$($given by Proposition $\ref{BFProp2})$   
are uniquely determined by the $D+1$ linear equations
\begin{equation} \label{BFFM6}
\left\{ \begin{array}{ccl}
        1  & = & 0 \cdot V + \sum_{\ell=1}^D w_{\ell} \ ,     \\   
0 & = & V + w_j \cdot \big(\log_v(C_{a_j}(\vs) \cdot r_j) - \frac{1}{e_j(q_v^{f_j}-1)}\big)
           + \displaystyle{\sum^D_{\substack{ \ell=1 \\ \ell \ne j }}} 
            w_{\ell}  \cdot \log_v \big([a_{\ell},a_j]_{\fX,\vs}\big)  \\
 &  & \qquad \qquad \text{for $j = 1, \ldots, D$.}  
          \end{array} \right.
\end{equation}
If $\vs \in \cP^m \cap \QQ^m$, then $V_{\fX,\vs}(H)$ and the $w_{\ell}(\vs)$ belong to $\QQ$.
\index{Robin constant!$(\fX,\vs)$}
\index{Robin constant!nonarchimedean!takes on rational values|ii}
\end{corollary}

\begin{proof}
This follows by combining Theorem \ref{BFThm2A}, Proposition \ref{BFProp2}, and Corollary \ref{BFCor1}.
Note that if $\vs \in \cP^m \cap \QQ^m$, then for each $j$ we have 
$\log_v(C_{a_j}(\vs)) \in \QQ$ by Proposition \ref{APropA2}(B2), 
so $V_{\fX,\vs}(H_{j}) = \frac{1}{e_j(q_v^{f_j}-1)} - \log_v(C_{a_j}(\vs) \cdot r_j) \in \QQ$.  
\index{Robin constant!$(\fX,\vs)$}
\index{Robin constant!nonarchimedean!takes on rational values}
\end{proof}

As a consequence, we show that for sets $H$ of form in Corollary \ref{BFCor2}, 
then for each $\zeta \notin H$ the Green's function $G(x,\zeta;H)$ 
\index{Green's function!nonarchimedean!takes on rational values}
and the Robin constant $V_{\zeta}(H)$ take on rational values.
\index{Robin constant!nonarchimedean!takes on rational values}

\begin{corollary} \label{BFCor3}  
Let $B(a_1,r_1), \ldots, B(a_D,r_D) \subset \cC_v(\CC_v)$  
be pairwise disjoint isometrically parametrizable balls 
whose union is disjoint from $\fX$. For each $\ell$, let $F_{w_\ell}/K_v$ be a finite 
extension in $\CC_v$.  Assume that each $a_{\ell} \in \cC_v(F_{w_{\ell}})$ 
and each $r_{\ell} \in |F_{w_{\ell}}^{\times}|_v$, 
and put $H_{\ell} = \cC_v(F_{w_\ell}) \bigcap B(a_{\ell},r_{\ell})$. 

Let $H = \bigcup_{\ell=1}^D H_{\ell}$.
Then for each $\zeta \in \cC_v(\CC_v) \backslash H$, we have $V_\zeta(H) \in \QQ$, 
and for each $x \in \cC_v(\CC_v) \backslash \{ \zeta \}$ we have $G(x,\zeta;H) \in \QQ$.
\index{Robin constant!nonarchimedean!takes on rational values}
\index{Green's function!nonarchimedean!takes on rational values}
\end{corollary}

\begin{proof}  We apply the preceding results, taking $\fX = \{\zeta\}$ and $\vs = (1)$. 
Fix a uniformizing parameter\index{uniformizing parameter!normalizes canonical distance} 
$g_{\zeta}(z) \in \CC_v(\cC)$ and normalize the canonical distance
\index{canonical distance!$[z,w]_{\zeta}$!normalization of}
by $\lim_{x \rightarrow \zeta} [x,y]_{\zeta} \cdot |g_{\zeta}(z)|_v = 1$ 
as usual.  By Proposition \ref{BFProp2} for each $z$ we have 
\begin{equation} \label{PotFFormula} 
u_{\zeta}(z,H) \ = \ \sum_{\ell = 1}^D w_{\ell} \cdot u_{\zeta}(z,H_{\ell}) \ .
\end{equation} 
\index{equilibrium potential}
By Corollary \ref{BFCor2}, under our hypotheses on $H$,
the Robin constant $V_{\zeta}(H) = V_{\fX,\vs}(H)$ 
\index{Robin constant!$(\fX,\vs)$}
and the weights $w_1, \ldots, w_D$  belong to $\QQ$. 
By the definition of the Green's function, for each $z$ 
\index{Green's function} 
\begin{equation*}
G(z,\zeta;H) \ = \ V_{\zeta}(H) - u_{\zeta}(z,H) \ .
\end{equation*}
\index{Robin constant}
\index{equilibrium potential}

First, suppose $x \in H$.  
The set $H$ satisfies the hypotheses of Proposition \ref{ACE11}, 
so we have $u_{\zeta}(x,H) = V_{\zeta}(H)$ and $G(x,\zeta;H) = 0 \in \QQ$.  
\index{Robin constant}
\index{Green's function!nonarchimedean!takes on rational values}
Next, take $x \in \cC_v(\CC_v) \backslash \Big(\big( \bigcup_{\ell=1}^D
             B(a_{\ell},r_{\ell}) \big) \cup \{\zeta\} \Big)$.   
By Corollary \ref{BFCor1}, for each $z \notin B(a_{\ell},r_{\ell})$ we have 
$u_{\zeta}(z,H_{\ell}) = -\log_v([z,a_\ell]_{\zeta})$. Inserting this in (\ref{PotFFormula}) gives
\begin{equation*} 
u_{\zeta}(x,H) \ = \  - \sum_{\ell = 1}^D w_{\ell} \cdot \log_v([x,a_{\ell}]_{\zeta}) \ .
\end{equation*} 
By Proposition \ref{APropA2}(B1) we have $\log_v([x,a_{\ell}]_{\zeta}) \in \QQ$\
for each $\ell$, so $u_{\zeta}(x,H) \in \QQ$.   
Hence $G(x,\zeta;H) = V_{\zeta}(H) -u_{\zeta}(x,H) \in \QQ$.
\index{Robin constant!nonarchimedean!takes on rational values|ii}
\index{Green's function!nonarchimedean!takes on rational values|ii}
Finally, let $x \in \cC_v(\CC_v) \backslash  \Big( H \cup \{\zeta\} \Big)$ be arbitrary.  
By replacing the cover $B(a_1,r_1), \ldots, B(a_D,r_D)$ of $H$ with a finer cover, 
we can arrange that $x \notin \bigcup_{\ell=1}^D B(a_{\ell},r_{\ell})$.  Applying 
the argument above to this new cover, we see that $G(x,\zeta;H) \in \QQ$.
\index{Green's function!nonarchimedean!takes on rational values}
\end{proof}

\vskip .1 in
We close with a proposition which says that deleting small balls 
from a set $H$ of positive capacity does not significantly change 
\index{capacity $> 0$} 
its $(\fX,\vs)$-Robin constant or potential function.
\index{Robin constant!$(\fX,\vs)$}
\index{potential function!$(\fX,\vs)$}
\index{equilibrium potential!$(\fX,\vs)$}

\begin{proposition} \label{BPropF3}
Let $H \subset \cC_v(\CC_v) \backslash \fX$ be a compact set of positive capacity,
\index{capacity $> 0$}  
and fix $q_1, \ldots, q_d \in H$.   For each $r > 0$, write 
\begin{equation*}
\ckH(r) \ = \ H \backslash \left( \bigcup_{\ell=1}^d B(q_{\ell},r) \right) \ .
\end{equation*}
Let $W$ be a neighborhood of $H$.  Then for each $\varepsilon > 0$,
there is an $R > 0$ such that, uniformly for compact sets $H^{\prime}$ 
with $\ckH(R) \subseteq H^{\prime} \subseteq H$ and for $\vs \in \cP^m$, 
\begin{equation*}
V_{\fX,\vs}(H) \ \le \ V_{\fX,\vs}(H^{\prime}) 
\ < \ V_{\fX,\vs}(H) + \varepsilon \ ,
\end{equation*}
\index{Robin constant!$(\fX,\vs)$}
and for all $z \in \cC_v(\CC_v) \backslash (\fX \cup W)$,
\begin{equation*}
|u_{\fX,\vs}(z,H) - u_{\fX,\vs}(z,H^{\prime})| \ < \ \varepsilon \ .  
\end{equation*}
\index{equilibrium potential!$(\fX,\vs)$}
\end{proposition}

\begin{proof}  After shrinking $W$ if necessary, we can assume that $W$ has the form 
$W = \bigcup_{\ell=1}^D B(a_{\ell},r_{\ell})$, where the balls $B(a_{\ell},r_{\ell})$
are isometrically parametrizable, pairwise disjoint, and do not meet $\fX$.

We will first prove analogous assertions for the $V_{x_i}(H^{\prime})$
\index{Robin constant}
and $u_{x_i}(z,H^{\prime})$, and then consider the situation for general $\vs$.
\index{equilibrium potential}
We begin by showing that for each $x_i \in \fX$, 
\begin{equation} \label{FLimF} 
\lim_{r \rightarrow 0} V_{x_i}(\ckH(r))  \ = \ V_{x_i}(H) \ .
\end{equation}
A finite set has capacity $0$, so for each $x_i$, by 
\index{capacity $= 0$} 
(\cite{RR1}, Corollary 4.1.15),  
\begin{equation*}
\gamma_{x_i}(H) 
  \ = \ \gammabar_{x_j}(H \backslash \{q_1, \ldots, q_d\})  
  \ := \ \sup_{\substack{ \text{compact} 
                  \\ A \subset H \backslash \{q_1, \ldots, q_d\} }} 
                 \gamma_{x_i}(A) \ .
\end{equation*}
As $r \rightarrow 0$, the $\ckH(r)$ form an increasing sequence of compact 
sets whose union is $H \backslash \{q_1, \ldots, q_d\}$.  Hence 
\begin{equation*}
\lim_{r \rightarrow 0} \gamma_{x_i}(\ckH(r)) 
     \ = \ \gammabar_{x_i}(H \backslash \{q_1, \ldots, q_d\})
     \ = \ \gamma_{x_i}(H)
\end{equation*}
which is equivalent to (\ref{FLimF}).  As the $\ckH(r)$ increase, 
the $V_{x_i}(\ckH(r))$ decrease.
\index{Robin constant!monotonicity of}  

We next show weak convergence of the equilibrium distributions.  Again fix $x_i$,
\index{equilibrium distribution!$(\fX,\vs)$}
and let $\mu_i$ be the equilibrium distribution of $H$ with respect to $x_i$.
Take a sequence 
$r_1 > r_2 > \cdots > 0$ with $r_j \rightarrow 0$, and put $H_j = \ckH(r_j)$. 
We can assume  $r_1$ is small enough that each $H_j$ has positive capacity.
\index{capacity $> 0$} 
Let $\mu_i^{(j)}$ be the equilibrium distribution of $H_j$ with respect to $x_i$.
As shown above, the Robin constants $V_{x_i}(H_j)$ decrease monotonically to $V_{x_i}(H)$.
\index{Robin constant!monotonicity of}
\index{Robin constant!monotonicity of}
Since the equilibrium measure of $H$ with respect to $x_i$ is unique (\cite{RR1}, Theorem 4.1.22), 
\index{equilibrium distribution}
the argument on (\cite{RR1}, p.190) shows 
that the $\mu_i^{(j)}$ converge weakly to $\mu_i$.  

Since $[z,w]_{x_i}$ is constant for $z, w$ belonging to disjoint 
isometrically parametrizable balls in $\cC_v(\CC_v) \backslash \{x_i\}$, for each 
$z \in \cC_v(\CC_v) \backslash (\{x_i\} \cup W)$ we have 
\begin{eqnarray*}
u_{x_i}(z,H) & = & 
\sum_{\ell=1}^D -\log_v([z,a_{\ell}]_{x_i}) \cdot \mu_i(B(a_{\ell},\rho_{\ell})) \ , \\
u_{x_i}(z,\ckH(r_j)) & = & 
\sum_{\ell=1}^D -\log_v([z,a_{\ell}]_{x_i}) \cdot \mu_i^{(j)}(B(a_{\ell},\rho_{\ell})) \ . \\
\end{eqnarray*}
\index{equilibrium potential}
\index{equilibrium distribution}
Let $g_{x_i}(z)$ be the uniformizing parameter\index{uniformizing parameter!normalizes canonical distance} 
which determines the normalization of $[z,w]_{x_i}$.  
By the same argument as in the proof of (\cite{RR1}, Proposition 4.1.5), 
there is an isometrically parametrizable ball $B(x_i,\delta)$ such that 
$u_{x_i}(z,H) = u_{x_i}(z,\ckH(r_j)) = \log_v(|g_{x_i}(z)|_v)$ for all 
\index{isometrically parametrizable ball}
\index{equilibrium potential}
$z \in B(x_i,\delta) \backslash \{x_i\}$  and all $j$.  Furthermore, 
by (\cite{RR1}, Proposition 4.1.1) there is a finite bound $B$ 
such that $|\log_v([z,w]_{x_i})| \le B$ for each $w \in W$ and each 
$z \in \cC_v(\CC_v) \backslash (W \cup B(x_i,\delta))$.  
Since the $\mu_i^{(j)}$ converge weakly to $\mu_i$ as $j \rightarrow \infty$, 
it follows that the $u_{x_i}(z,\ckH(r_j))$ 
converge uniformly to $u_{x_i}(z,H)$ on $\cC_v(\CC_v) \backslash (\{x_i\} \cup W)$. 
\index{equilibrium distribution}
\index{equilibrium potential}

The argument above applies for each $x_i$.  
Given $\varepsilon > 0$, take $j$ large enough that for all $i = 1, \ldots, m$,   we have
$|V_{x_i}(H)-V_{x_i}(\ckH(r_j))| < \varepsilon/3$ 
\index{Robin constant}
\index{equilibrium potential}
and $|u_{x_i}(z,H) - u_{x_i}(\ckH(r_j))| < \varepsilon/3$
for all $z \in \cC_v(\CC_v) \backslash (\{x_i\} \cup W)$. Put $R_0 = r_j$.  

Let $H^{\prime}$ be a compact set with $\ckH(R_0) \subseteq H^{\prime} \subseteq H$.
By the monotonicity of the Robin constant, for each $x_i$ 
\index{Robin constant!monotonicity of}
\begin{equation} \label{FV00} 
V_{x_i}(H) \ \le \ V_{x_i}(H^{\prime}) \ \le \ V_{x_i}(\ckH(R_0)) \ < \ V_{x_i}(H) + \varepsilon/3 \ .
\end{equation}
We claim as well that  $|u_{x_i}(z,H) - u_{x_i}(z,H^{\prime})| < \varepsilon$ outside $W$.
\index{equilibrium potential}
This follows from a 3-$\varepsilon$'s argument, 
using monotonicity of the upper Green's functions.  
\index{Green's function!upper} 
\index{Green's function!monotonic}
Recall that $G(z,x_i,H) = V_{x_i}(H)-u_{x_i}(z,H)$,
\index{Robin constant} 
\index{Green's function}
with similar equalities for $G(z,x_i,H^{\prime})$ and $G(z,x_i\ckH(R_0))$.
By (\cite{RR1}, Proposition 4.4.1), for all $z \in \cC_v(\CC_v) \backslash \{x_i\}$ 
\begin{equation*}
G(z,x_i;H) \ \le \ G(z,x_i;H^{\prime}) \ \le \ G(z,x_i;\ckH(R_0)) \ .
\end{equation*}
It follows that for all $z \in \cC_v(\CC_v) \backslash (\{x_i\} \cup W)$, 
\begin{eqnarray}
\lefteqn{|u_{x_i}(z,H)- u_{x_i}(z,H^{\prime})|} \quad &  & \notag \\
    & = & |(V_{x_i}(H) - G(z,x_i;H)) -
             (V_{x_i}(H^{\prime})-G(z,x_i;H^{\prime})|  \label{FV01} \\
    & \le & |V_{x_i}(H)-V_{x_i}(H^{\prime})| + |G(z,x_i;H)-G(z,x_i;H^{\prime})|  \notag \\ 
    & \le & |V_{x_i}(H)-V_{x_i}(\ckH(R_0))| + |G(z,x_i;H)-G(z,x_i;\ckH(R_0))| \notag \\  
    & \le & 2|V_{x_i}(H)-V_{x_i}(\ckH(R_0))| + |u_{x_i}(z,H)-u_{x_i}(z,\ckH(R_0))| 
    \ < \ \varepsilon \ . \notag 
\end{eqnarray}
\index{Green's function}
\index{Robin constant}
\index{equilibrium potential}
      
\vskip .05 in
Now fix $\vs \in \cP^m$.  Recall that $W = \bigcup_{\ell=1}^D B(a_\ell,r_{\ell})$. 
For  each $\ell = 1, \ldots, D$,
put $H_{\ell} = H \cap B(a_{\ell},r_{\ell})$.  We can assume without loss that
each $H_{\ell}$ has positive capacity, since removing a set of capacity $0$
\index{capacity $> 0$} 
\index{capacity $= 0$} 
from $H$ does not change its Robin constant or potential function 
\index{Robin constant} 
\index{equilibrium potential}
with respect to any $x_i$ (\cite{RR1}, Corollary 4.1.15).  Applying the first 
part of the proof to each $H_{\ell}$ with respect to its cover 
$W_{\ell} := B(a_{\ell},r_{\ell})$,
there is an $R > 0$ such that for all $\ell$, 
if $H^{\prime}_{\ell}$ is a compact set with 
$\ckH_{\ell}(R) \subset H^{\prime}_{\ell} \subset H_{\ell}$, then for each $x_i$ 
and for all $z \in \cC_v(\CC_v) \backslash (\{x_i\} \cup B(a_{\ell},r_{\ell}))$ 
\begin{eqnarray}
V_{x_i}(H_{\ell}) \ \le \ V_{x_i}(H^{\prime}_{\ell}) 
\ < \ V_{x_i}(H_{\ell}) + \varepsilon \ , \quad \text{and}  \label{FCA1} \\
|u_{x_i}(z,H_{\ell}) - u_{x_i}(z,H^{\prime}_{\ell})| \ < \ \varepsilon \ . 
\qquad \label{FCA2}
\end{eqnarray}
\index{Robin constant}
\index{equilibrium potential}
We can assume $R \le R_0$, so (\ref{FV00}) and (\ref{FV01}) continue to hold.

Fix a compact set $H^{\prime}$ with $\ckH(R) \subseteq H^{\prime} \subseteq H$,
and for each $\ell$ put $H^{\prime}_{\ell} = H^{\prime} \cap B(a_{\ell},r_{\ell})$.
As above, we can assume that each $H^{\prime}_{\ell}$ has positive capacity. 
\index{capacity $> 0$} 
Fix $\vs \in \cP^m$.  We will now show that 
\begin{equation} \label{FVN1} 
V_{\fX,\vs}(H) \ \le \ V_{\fX,\vs}(H^{\prime}) 
               \ < \ V_{\fX,\vs}(H) + \varepsilon \ .
\end{equation} 
\index{Robin constant!$(\fX,\vs)$}

By Lemma \ref{BFLem1}, for each $\ell$ 
the equilibrium distribution $\mu_{\ell}^*$ of
\index{equilibrium distribution!$(\fX,\vs)$}
$H_{\ell}$ is independent of $x_i$, so the expression (\ref{BFFF8})
for $V_{\fX,\vs}(H_{\ell})$ in Proposition \ref{BPropF1} simplifies to
\index{Robin constant!$(\fX,\vs)$} 
\index{Robin constant}
\begin{equation} \label{FVC1}
V_{\fX,\vs}(H_{\ell}) \ = \ \sum_{i=1}^m s_i V_{x_i}(H_{\ell}) \ .
\end{equation}
Similarly the equilibrium distribution $\mu_{\ell}^{\prime *}$ of 
$H^{\prime}_{\ell}$ is independent of $x_i$, so   
\begin{equation} \label{FVC2}
V_{\fX,\vs}(H^{\prime}_{\ell}) \ = \ \sum_{i=1}^m s_i V_{x_i}(H^{\prime}_{\ell}) \ .
\end{equation}
\index{Robin constant!$(\fX,\vs)$}
\index{Robin constant}
It follows from (\ref{FCA1}), (\ref{FVC1}), and (\ref{FVC2}) that 
\begin{equation} \label{FBBX1} 
V_{\fX,\vs}(H_{\ell}) \ \le \ V_{\fX,\vs}(H^{\prime}_{\ell}) \ < \ 
V_{\fX,\vs}(H_{\ell}) + \varepsilon \ .
\end{equation} 

By Proposition \ref{BFProp2}, there are weights $w_{k}(\vs) > 0$ 
with $\sum_{k=1}^D w_{k}(\vs) = 1$ 
such that the $(\fX,\vs)$-equilibrium distribution of $H$
\index{equilibrium distribution!$(\fX,\vs)$}
is given by $\mu_{\fX,\vs} = \sum_k w_k(\vs) \mu_k^*$.  
\index{equilibrium distribution}

Consider the probability measure 
$\nu := \sum_k w_{\ell}(\vs) \mu_k^{\prime *}$
on $H^{\prime}$.  Although $\nu$ may not be the equilibrium measure 
\index{equilibrium distribution!$(\fX,\vs)$}
$\mu_{\fX,\vs}^{\prime}$ of $H^{\prime}$, we will show that 
\begin{equation} \label{FNX1} 
u_{\fX,\vs}(z,\nu) \ < \ V_{\fX,\vs}(H) + \varepsilon \ 
\end{equation} 
\index{equilibrium potential!$(\fX,\vs)$}
\index{Robin constant!$(\fX,\vs)$}
for all $z \in \bigcup_{\ell=1}^D B(a_{\ell},r_{\ell})$.  
By the definition of $V_{\fX,\vs}(H^{\prime})$, it follows that 
\index{Robin constant!$(\fX,\vs)$}
\begin{equation*} 
V_{\fX,\vs}(H^{\prime}) \ \le \ I_{\fX,\vs}(\nu) \ = \ 
\int_{H^{\prime}} u_{\fX,\vs}(z,\nu) \, d\nu(z) 
\ < \ V_{\fX,\vs}(H) + \varepsilon \ .
\end{equation*}
\index{equilibrium potential!$(\fX,\vs)$}
This will yield (\ref{FVN1}) since the lower bound there is trivial. 

By Lemma \ref{BFLem1}, for each $\ell = 1, \ldots, D$ 
\begin{equation*}
u_{\fX,\vs}(z,H_{\ell}) \ = \ \sum_{i=1}^m s_i u_{x_i}(z,H_{\ell})\ , \ \quad  
u_{\fX,\vs}(z,H^{\prime}_{\ell}) 
\ = \ \sum_{i=1}^m s_i u_{x_i}(z,H^{\prime}_{\ell}) \ .
\end{equation*}
\index{equilibrium potential!$(\fX,\vs)$} 
Hence by (\ref{FCA2}), 
for each $z \in \cC_v(\CC_v) \backslash (\fX \cup B(a_{\ell},r_{\ell}))$
\begin{equation} \label{FBBX2}
|u_{\fX,\vs}(z,H_{\ell})-u_{\fX,\vs}(z,H^{\prime}_{\ell})| \ < \ \varepsilon \ .
\end{equation} 
Fix $\ell$ and note that 
$u_{\fX,\vs}(z,H^{\prime}_{\ell}) \le V_{\fX,\vs}(H^{\prime}_{\ell})$
\index{Robin constant!$(\fX,\vs)$}
\index{equilibrium potential!$(\fX,\vs)$}
for all $z$, hence in particular for $z \in B(a_{\ell},r_{\ell})$.  
Thus for each $k \ne \ell$, 
$u_{\fX,\vs}(z,H_k)$ and $u_{\fX,\vs}(z,H^{\prime}_k)$ are constant
on $B(a_{\ell},r_{\ell})$, since $[z,w]_{\fX,\vs}$ is constant on 
pairwise disjoint isometrically parametrizable balls not meeting $\fX$.  
This means that for any $z \in B(a_{\ell},r_{\ell})$
\begin{equation} \label{FBBX3}
V_{\fX,\vs}(z,H) \ = \ w_{\ell}(\vs) V_{\fX,\vs}(H_{\ell}) + 
         \sum_{k \ne \ell} w_k(\vs) u_{\fX,\vs}(z,H_k) \ .
\end{equation} 
\index{Robin constant!$(\fX,\vs)$}
\index{equilibrium potential!$(\fX,\vs)$}
On the other hand,  by the definition of $\nu$ and 
(\ref{FBBX1}), (\ref{FBBX2}), and (\ref{FBBX3}), on $B(a_{\ell},r_{\ell})$, 
\begin{eqnarray*}
u_{\fX,\vs}(z,H^{\prime}) & = & \sum_{k=1}^D w_k(\vs) u_{\fX,\vs}(z,H^{\prime}_k) \\
         & \le & w_\ell(\vs) (V_{\fX,\vs}(H) + \varepsilon) 
                 + \sum_{k \ne \ell} w_k(\vs) (u_{\fX,\vs}(z,H_k) + \varepsilon) \\
         & = & V_{\fX,\vs}(H) + \varepsilon \ ,
\end{eqnarray*}
\index{Robin constant!$(\fX,\vs)$}
\index{equilibrium potential!$(\fX,\vs)$}
which gives (\ref{FNX1}) as $\ell$ varies.  

Finally, we show that $|u_{\fX,\vs}(z,H) - u_{\fX,\vs}(z,H^{\prime})| < 3\varepsilon$
for all $z \notin \fX \cup W$.  Indeed for such $z$, by (\ref{FV00}),
(\ref{FV01}), (\ref{FVN1}), and Proposition \ref{BPropF1},  
\begin{eqnarray*}
\lefteqn{|u_{\fX,\vs}(z,H) - u_{\fX,\vs}(z,H^{\prime})|} \qquad &  & \\
  & = & |(V_{\fX,\vs}(H)- G_{\fX,\vs}(z,H)) - 
                 (V_{\fX,\vs}(H^{\prime})- G_{\fX,\vs}(z,H^{\prime}))| \\
  & \le & |V_{\fX,\vs}(H)-V_{\fX,\vs}(H^{\prime})| 
             + \sum_{i=1}^m s_i |G(z,x_i;H)-G(z,x_i;H^{\prime})|   \\
  & \le & |V_{\fX,\vs}(H)-V_{\fX,\vs}(H^{\prime})| 
                  + \sum_{i=1}^m s_i |V_{x_i}(H)-V_{x_i}(H^{\prime})| \\
  &  & \qquad \qquad + \sum_{i=1}^m s_i |u_{x_i}(z,H)-u_{x_i}(z,H^{\prime})| 
  \ < \ 3 \varepsilon  \ .               
\end{eqnarray*} 
\index{equilibrium potential!$(\fX,\vs)$}
\index{Robin constant!$(\fX,\vs)$}
\index{Green's function}
Since $\varepsilon > 0$ is arbitrary, this proves the Proposition.  
\end{proof} 
 

%% file: NewFSZAppB.tex
\chapter{The Construction of Oscillating Pseudopolynomials} 
\label{AppB}

The purpose of this appendix is to construct $(\fX,\vs)$-pseudopolynomials
\index{pseudopolynomial!$(\fX,\vs)$} 
with all their roots in $H$, and having large oscillations on $H \cap \cC_v(\RR)$, 
in the case when $K_v \cong \RR$.
This is accomplished in Theorem \ref{bT3}, 
providing the potential-theoretic input to Theorem \ref{RThm2}.   
\index{potential theory}

\vskip .1 in
Throughout this appendix we assume that $K_v$ is archimedean.   
We specialize to the case $K_v \cong \RR$ only in section 
\ref{AppB}.\ref{OscillatingPseudoPolySection}.  
Let  $\cC_v/K_v$ be a smooth, connected projective curve, 
so $\cC_v(\CC)$ is a Riemann surface. 
\index{Riemann surface}     

We keep the notation and assumptions of \S\ref{Chap3}.\ref{AssumptionsSection}.  
Let $\fX = \{x_1, \ldots, x_m\} \subset  \cC_v(\tK) \subset \cC_v(\CC)$ 
be a finite, $\Aut(\tK/K)$-stable set of points.  
As usual, $H \subset \cC_v(\CC) \backslash \fX$ will be compact.   
Let the canonical distances $[z,w]_{x_i}$ be normalized so that
\index{canonical distance!$[z,w]_{\zeta}$}
\index{canonical distance!$[z,w]_{\zeta}$!normalization of}
$\lim_{z \rightarrow x_i} [z,w]_{x_i} \cdot g_{x_i}(z) = 1$, 
with the $g_{x_i}(z)$ as in \S\ref{Chap3}.\ref{AssumptionsSection}.  
As in Chapter \ref{Chap5}, if $\vs = (s_1, \ldots, s_m) \in \cP^m$ 
is a $K_v$-symmetric probability vector, 
\index{$K_v$-symmetric!probability vector}
the $(\fX,\vs)$-canonical distance is defined by
\index{canonical distance!$[z,w]_{\fX,\vs}$}
\begin{equation*}
[z,w]_{\fX,\vs} \ = \ \prod_{i=1}^{m} ([z,w]_{x_i})^{s_i} \ .
\end{equation*}  
Recall that an $(\fX,\vs)$-pseudopolynomial is a function of the form
\index{pseudopolynomial!$(\fX,\vs)$|ii}
\begin{equation*}
P(z) \ = \ \prod_{i=1}^N [z,\alpha_i]_{\fX,\vs}
\end{equation*}
where each $\alpha_i \in \cC_v(\CC) \backslash \fX$.  For simplicity, we will
often just speak of pseudopolynomials, rather than $(\fX,\vs)$-pseudopolynomials. 
\index{pseudopolynomial!$(\fX,\vs)$}

\vskip .1 in
The motivation for the construction is the classical fact that the 
Chebyshev polynomial $\tP_N(z)$ of degree $N$ for an interval $[a,b]$ oscillates 
\index{Chebyshev polynomial} 
$N$ times between $\pm M^N$,  where $M^N$ is the $\sup$ norm  $\|\tP_N\|_{[a,b]}$.  
Thus $|\tP_N(z)|$ varies $N$ times from $M^N$ to $0$ to $M^N$ on $[a,b]$.  
One could hope to prove the same thing for an $(\fX,\vs)$-Chebyshev 
\index{pseudopolynomial!$(\fX,\vs)$}
pseudopolynomial (that is, a pseudopolynomial having minimal $\sup$ norm on $H$),
but this seems difficult because while $[z,w]_{\fX,\vs}$ is locally 
well-understood, it is globally quite mysterious.  Furthermore, for disconnected
sets $H$, Chebyshev pseudopolynomials need not have all their roots in $H$.  
\index{pseudopolynomial!Chebyshev}

    Instead, we consider {\em restricted} Chebyshev pseudopolynomials for $H$, 
\index{Chebyshev pseudopolynomial!restricted} 
which by definition have all their roots in $H$.  
The result which makes everything work is Proposition \ref{bP1}, 
which asserts that restricted Chebyshev pseudopolynomials for sufficiently short intervals 
\index{pseudopolynomial!restricted}
(``short'' is made precise in Definition \ref{ShortnessDef}) have an 
\index{short@`short' interval}  
oscillation property like that of classical Chebyshev polynomials.  

    This suggests that when $H = \bigcup_{\ell = 1}^D H_{\ell}$ is a disjoint 
union of short intervals, we could obtain the function we want by taking 
\index{short@`short' interval}  
the product of the Chebyshev pseudopolynomials for the intervals $H_{\ell}$. 
\index{pseudopolynomial!Chebyshev}
\index{Chebyshev pseudopolynomial} 
However, this would only be partially right, because the 
pseudopolynomials from intervals $H_j$ with $j \ne \ell$ would affect the magnitude of the
\index{pseudopolynomial}
oscillations on $H_{\ell}$.  Instead, we need to take the product of appropriate
``weighted'' Cheybshev pseudopolynomials for the sets $H_{\ell}$, 
\index{pseudopolynomial!weighted Chebyshev}
\index{Chebyshev pseudopolynomial!weighted} 
incorporating the background function coming from the other $H_j$ from the very start. 
This leads to the topic of ``weighted potential theory'', 
\index{potential theory!weighted} 
or ``potential theory in the presence of an external field'', 
a subject which goes back at least as far as Gauss. 
\index{potential theory}

We refer the reader to the book by Saff and Totik (\cite{ST}) 
\index{Saff, Ed} 
\index{Totik, Vilmos} 
for an exposition of weighted potential theory for subsets of $\CC$. 
\index{potential theory!weighted} 
In the classical unweighted case, it is known that under appropriate hypotheses on $H$,
the discrete probability measures associated to Cheybshev polynomials converge weakly
\index{Chebyshev polynomial}
to the equilibrium distribution of $H$.  A variant of this for weighted Chebyshev polynomials 
\index{Chebyshev polynomial!weighted}
is given in (\cite{ST}, Theorem III.4.2).  Thus, one expects that the correct weight
function to use in defining the weighted $(\fX,\vs)$-Chebyshev pseudopolynomials 
for the\index{pseudopolynomial!weighted Chebyshev}\index{Chebyshev pseudopolynomial!weighted} 
sets $H_{\ell}$, should be the exponential of the part of the potential function 
\index{equilibrium potential!$(\fX,\vs)$}
for $H$ coming from $H \backslash H_{\ell}$.  

This appendix works out that idea.  To reach the goal, 
it is necessary to develop a considerable amount of machinery.  Fortunately, most of 
the arguments are direct adaptations of classical proofs to 
the $(\fX,\vs)$-context on curves.  
Sections \ref{AppB}.\ref{WeightedXSCapacitySection}--\ref{AppB}.\ref{Weight0CaseSection} below 
``turn the crank'' of a standard potential-theoretic machine,
\index{potential theory} 
developing the background needed to prove three key facts: 
Theorems \ref{APF15}, \ref{ATF18} and \ref{ATF19}.
Theorem \ref{APF15} says that in the case of interest to us, 
the weighted equilibrium distribution for each $H_{\ell}$ 
coincides with the restriction to $H_{\ell}$ 
of the unweighted $(\fX,\vs)$-equilibrium distribution for $H$.
Theorem \ref{ATF18} shows that the weighted Chebyshev constant for $H_{\ell}$ coincides with the
\index{Chebyshev constant!weighted} 
unweighted $(\fX,\vs)$-capacity of $H$.  And, Theorem \ref{ATF19} proves the convergence of the
\index{capacity!$(\fX,\vs)$}
weighted discrete Chebyshev measures to the equilibrium distribution.  
\index{Chebyshev measure} 
If one accepts these results, 
the construction of oscillating pseudopolynomials can be carried out quickly.  
Proposition \ref{bP1} in section \ref{AppB}.\ref{ChebShortIntervalSection} establishes the
oscillation property of weighted Chebyshev pseudopolynomials for short intervals, 
\index{pseudopolynomial!weighted Chebyshev}\index{short@`short' interval}\index{Chebyshev pseudopolynomial!weighted}  
and the main Theorem \ref{bT3} in section \ref{AppB}.\ref{OscillatingPseudoPolySection} assembles the pieces. 
\index{capacity!weighted $(\fX,\vs)$}
We now briefly outline the contents of sections 
\ref{AppB}.\ref{WeightedXSCapacitySection}--\ref{AppB}.\ref{Weight0CaseSection}.  Sections 
\ref{AppB}.\ref{WeightedXSCapacitySection}, \ref{AppB}.\ref{WeightedChebyshevSection} 
and \ref{AppB}.\ref{WeightedTDSection} introduce the weighted $(\fX,\vs)$-capacity, 
\index{capacity!weighted $(\fX,\vs)$}
weighted $(\fX,\vs)$-Chebyshev constant, and weighted $(\fX,\vs)$-transfinite diameter and prove
\index{transfinite diameter!weighted $(\fX,\vs)$}
\index{Chebyshev constant!weighted $(\fX,\vs)$} 
their existence.  Our definition of these quantities is somewhat different from that
in \cite{ST}.  The detailed convergence proof for the $(\fX,\vs)$-Chebyshev constant is needed
\index{Chebyshev constant!$(\fX,\vs)$}
because a strong notion of convergence for the finite-level Chebyshev constants to the
asymptotic one is required in Theorem \ref{bT3}.  
A peculiar aspect of classical potential theory is the ``rock-paper-scissors'' nature of the
\index{rock-paper-scissors argument}
\index{potential theory!classical}
capacity, Chebyshev constant, and transfinite diameter:  although all three are equal, 
\index{Chebyshev constant} 
\index{capacity}
\index{transfinite diameter} 
this is seen only after one proves inequalities between them in a cyclic manner.  
Section \ref{AppB}.\ref{Weight0CaseSection}
notes that when the weight function is trivial, the weighted objects constructed
in sections \ref{AppB}.\ref{WeightedXSCapacitySection}--\ref{AppB}.\ref{WeightedTDSection} 
coincide with the unweighted objects studied 
in Appendix \ref{AppA}.  
This observation has the consequence that in the case of interest, 
the weighted equilibrium distribution is unique (Theorem \ref{APF15}), 
a fact we are unable to establish in general. 

In our path through this material,   
we prove only what is needed to prove the main Theorem \ref{bT3}.  
Many standard facts,    
which could be established if we were developing the theory for its own sake, 
are not touched upon.  For a more complete treatment of the theory in 
classical case, see \cite{ST}.   


\section{ Weighted $(\fX,\vs)$-Capacity Theory }
       \label{WeightedXSCapacitySection}
\index{capacity!weighted $(\fX,\vs)$|ii}


In this section we introduce the weighted $(\fX,\vs)$-capacity.
\index{capacity!weighted $(\fX,\vs)$|ii}   

\medskip
{\bf Motivation.}   
Consider the energy minimization problem in the definition 
of the logarithmic capacity for compact sets $H \subset \CC$: 
\index{capacity!logarithmic} 
for each probability measure $\nu$ supported on $H$, 
put $I(\nu) = \iint -\log|z-w| \, d\mu(z) d\mu(w)$;  then the 
Robin constant is defined by 
\index{Robin constant}
\begin{equation} \label{BFA1}
V_{\infty}(H) \ = \ \inf_{\nu} I(\nu) 
\end{equation}
and the capacity is $\gamma(H) = e^{-V_{\infty}(H)}$. 
\index{capacity}  

If $\gamma(H) > 0$, there is a unique probability measure $\mu$
on $H$, called the equilibrium distribution, for which 
\begin{equation} \label{AFA2} 
\iint_{H \times H} -\log(|z-w|) \, d\mu(z) d\mu(w) \ = \ V_{\infty}(H) \ .
\end{equation} 
\index{Robin constant}
The equilibrium potential of $H$ is then defined by 
\index{equilibrium potential}
\begin{equation*}
u_{\infty}(z,H) \ = \ \int_{H} -\log(|z-w|) \, d\mu(w) \ .
\end{equation*} 

Now suppose $H = H_1 \bigcup H_2$, where
the $H_i$ are closed and disjoint, and put $\mu_1 = \mu \vert_{H_1}$,
$\mu_2 = \mu \vert_{H_2}$.  Let 
\begin{equation*}
\widehat{u}(z) \ = \ \int_{H_2} -\log(|z-w|) \, d\mu_2(w) 
\end{equation*}
be the part of the equilibrium potential coming from $H_2$. Then 
\index{equilibrium potential}
\begin{eqnarray}
V_{\infty}(H) & = & \iint_{H \times H} -\log(|z-w|) \, d\mu(z) d\mu(w)    \notag \\ 
    &=& \iint_{H_1 \times H_1} -\log(|z-w|) \, d\mu_1(z) d\mu_1(w) 
           + 2 \int_{H_1} \widehat{u}(z) d\mu_1(z)          \notag         \\
    & & \qquad \qquad + 
       \iint_{H_2 \times H_2} -\log(|z-w|) \, d\mu_2(z) d\mu_2(w) \ . 
\label{AFA3}
\end{eqnarray}
\index{Robin constant}
If one thinks of $\mu_2$ on $H_2$ as given, then the minimization 
problem (\ref{BFA1}) can be viewed as asking for a positive measure  
on $H_1$ of mass $\sigma = \mu(H_1)$, which minimizes the sum of the first two terms in 
(\ref{AFA3}).  We now generalize this situation, replacing $H_1$ by the
full set $H$ and allowing an arbitrary weight function $W(z)$.    

\vskip .1 in
To simplify notation, throughout the discussion below 
we fix $\fX$ and $\vs$ and suppress the $(\fX,\vs)$-dependence 
in all quantities, though that dependence is present.

\medskip
{\bf Definition of the Weighted $(\fX,\vs)$-Capacity.}
\index{capacity!weighted $(\fX,\vs)$|ii} 
Let $H$ be a compact subset of $\cC_v(\CC) \backslash \fX$, 
and let $W(z) : \cC_v(\CC) \rightarrow [0,\infty]$ be a function which 
is positive, bounded and continuous on a neighborhood of $H$.  

Fix a number $\sigma > 0$, and put 
\begin{equation*}
\widehat{u}(z) \ = \ -\log(W(z)) \ .
\end{equation*}
Given a positive Borel measure $\nu$ on $H$ with total mass $\sigma$, 
define the energy 
\begin{equation*}
I_{\sigma}(\nu,W) \ = \ \iint_{H \times H} -\log([z,w]_{\fX,\vs}) \, d\nu(z) d\nu(w)  
        + 2 \int_H \widehat{u}(z) \, d\nu(z) \ .
\end{equation*}
Let the weighted Robin constant be
\index{Robin constant!weighted $(\fX,\vs)$|ii}
\begin{equation} \label{FMu1} 
V_{\sigma}(H,W) \ = \ \inf_{\substack{ \nu \ge 0 \\ \nu(H) = \sigma }}
                I_{\sigma}(\nu,W) \ ,
\end{equation}
where the $\inf$ is taken over positive Borel measures on $H$ with 
total mass $\sigma$.  Then, define
the weighted capacity
\begin{equation}
\gamma_{\sigma}(H,W) \ = \ e^{-V_{\sigma}(H,W)} \ . \label{FNQ3A} 
\end{equation}

Because $W(z)$ is finite and bounded away from $0$,  
and $[z,w]_{\fX,\vs}/\|z,w\|_v$ is bounded for $z \ne w \in H$,  
it is easy to see that $\gamma_{\sigma}(H,W) > 0$
if and only if $H$ has positive capacity in the sense of Definition
\index{capacity $> 0$}
\ref{FCapDef}.

\medskip
If $\nu$ is a positive Borel measure on $H$, define the potential function 
\begin{eqnarray*}
u_{\nu}(z) & = & \int_H -\log([z,w]_{\fX,\vs}) \, d\nu(w)  \\
           & = & \lim_{t \rightarrow \infty}
          \int_H -\log^{(t)}([z,w]_{\fX,\vs}) \, d\nu(w)
\end{eqnarray*}     
\index{potential function}
where $-\log^{(t)}([z,w]_{\fX,\vs}) := \max(t,-\log([z,w]_{\fX,\vs})$. 
Then $u_{\nu}(z)$ is harmonic in $\cC_v(\CC) \backslash (H \cup \fX)$ and  
superharmonic\index{superharmonic} in $\cC_v(\CC) \backslash \fX$.  For any 
$w \notin \fX$, $u_{\nu}(z) + r \log([z,w]_{\fX,\vs})$ extends to 
a function harmonic in $\cC_v(\CC) \backslash (H \cup \{w\})$.
\index{potential function!properties of}

\medskip
Suppose $H$ has positive capacity.  
\index{capacity $> 0$} 
In this section we show the existence of a measure $\mu = \mu_{\sigma,H,W}$ 
achieving $V_{\sigma}(H,W)$, which we call an equilibrium distribution. 
\index{Robin constant!weighted $(\fX,\vs)$} 
We will prove uniqueness later, and only in the cases of interest to us 
(Theorems \ref{ATE8} and \ref{APF15}).   
However, it seems likely that the equilibrium distribution is always unique 
(this holds in the classical case:  see \cite{ST}, p.27).  

The proof of the existence of $\mu$ is standard.  
Take a sequence of measures $\nu_k \ge 0$ on $H$, with $\nu_k(H) = \sigma$ for all
$k$, such that 
\begin{equation*}
\lim_{k \rightarrow \infty} I_{\sigma}(\nu_k,W) \ = \ V_{\sigma}(H,W) \ .
\end{equation*} 
\index{Robin constant!weighted $(\fX,\vs)$}
After passing to a subsequence, we can assume the $\nu_k$ 
converge weakly to a measure $\mu$ on $H$, which is necessarily positive 
and satisfies $\mu(H) = \sigma$.  We claim that $I_{\sigma}(\mu,W) = V_{\sigma}(H,W)$.  

Tautologically $I_{\sigma}(\mu,W) \ge V_{\sigma}(H,W)$.  For the reverse inequality, 
\index{Robin constant!weighted $(\fX,\vs)$}
note that 
\begin{eqnarray}
I_{\sigma}(\mu,W) &=& \lim_{t \rightarrow \infty}
         \left( \iint_{H \times H} -\log^{(t)}([z,w]_{\fX,\vs}) \, d\mu(z) d\mu(w) 
                  + 2 \int_H \widehat{u}(z) d\mu(z) \right)  \notag \\
   &=&  \lim_{t \rightarrow \infty} \lim_{k \rightarrow \infty} 
     \left( \iint_{H \times H} -\log^{(t)}([z,w]_{\fX,\vs}) \, d\nu_k(z) d\nu_k(w) 
          + 2 \int_H \widehat{u}(z) d\nu_k(z) \right)  \label{AFA5} \\ 
   &\le&  \liminf_{k \rightarrow \infty} \lim_{t \rightarrow \infty} 
      \left( \iint_{H \times H} -\log^{(t)}([z,w]_{\fX,\vs}) \, d\nu_k(z) d\nu_k(w) 
          + 2 \int_H \widehat{u}(z) d\nu_k(z) \right)  \label{AFA6} \\
   &=& \liminf_{k \rightarrow \infty} I_{\sigma}(\nu_k,W) \ = \ V_{\sigma}(H,W)  \ . \notag
\end{eqnarray}
\index{Robin constant!weighted $(\fX,\vs)$}
The second equality follows from the weak convergence of the $\nu_k$ and
the continuity of $-\log^{(t)}([z,w]_{\fX,\vs})$ on $H$.  
The interchange of limits between (\ref{AFA5}) and (\ref{AFA6}) is valid 
because the kernels $-\log^{(t)}([z,w]_{\fX,\vs})$ are increasing with $t$.  

Let $H_{\mu}^* \subset H$ be the carrier\index{carrier|ii} of $\mu$:
\begin{equation*} 
H_{\mu}^* \ = \ \Big\{z \in H : \mu(B(z,\delta) \bigcap H) > 0 
\text{\ for each $\delta > 0$} \Big\} \ .
\end{equation*}  
where the balls $B(a,\delta) = \{z \in \cC_v(\CC) : \|z,a\|_v \le \delta\}$ 
are computed relative to a fixed spherical metric on $\cC_v(\CC)$.\index{spherical metric}  
Note that $H_{\mu}^*$ is closed.     

We now show that $u_{\mu}(z)$ satisfies an analogue of Frostman's theorem
(compare \cite{ST}, Theorem 1.3, p.27):  
\index{Frostman's Theorem} 
\index{equilibrium potential!takes constant value a.e. on $E_v$|ii}

\begin{theorem} \label{ATA1}  Suppose $H$ has positive capacity, and 
\index{capacity $> 0$}
let $\mu$ be any equilibrium distribution for $H$ relative to $W(z)$.  
Then there exist a constant $\cV_{\mu}$ and an 
$F_{\sigma}$ set $e_{\mu}$ of inner capacity $0$ contained in $H$  
\index{capacity $= 0$}
such that for all $z \in H_{\mu}^* \backslash e_{\mu}$  
\begin{equation*}
              u_{\mu}(z) + \widehat{u}(z) \ = \ \cV_{\mu} \ ,
\end{equation*}
and for all $z \in H \backslash e_{\mu}$
\begin{equation*}
              u_{\mu}(z) + \widehat{u}(z) \ \ge \ \cV_{\mu} \ .
\end{equation*}
Moreover, 
\begin{equation} \label{AFA8}
\sigma \cdot \cV_{\mu} + \int_H \widehat{u}(z) \, d\mu(z) \ = \ V_{\sigma}(H,W) \ . 
\end{equation}  
\index{potential function}
\end{theorem}

\begin{proof} 
Write $f(z) = u_{\mu}(z) + \widehat{u}(z)$ and set 
\begin{equation} \label{AFA9} 
\cV_{\mu} \ := \ \sup_{z \in H_{\mu}^*} f(z) \ .
\end{equation} 
We will now show that $f(z) \ge \cV_{\mu}$ for `almost all' $z \in H$.  Put 
\begin{equation*}
e_n \ = \ \{z \in H : f(z) \le \cV_{\mu} - \frac{1}{n} \}  
              \quad \text{for $n = 1, 2, 3 , \ldots$} \ ; 
\end{equation*}
and 
\begin{equation*}              
e_{\mu} \ = \ \{z \in H : f(z) < \cV_{\mu} \} \ ,             
\end{equation*}  
so $e_{\mu} = \bigcup_{n=1}^{\infty} e_n$.  Here, each $e_n$ is closed, 
since $u_{\mu}(z)$ is superharmonic\index{superharmonic} and hence lower semi-continuous.
\index{potential function}\index{potential function!is superharmonic}
\index{potential function!is lower semi-continuous}
\index{semi-continuous!potential function is lower semi-continuous}\index{superharmonic}
We will show that each $e_n$ has capacity $0$, which implies that $e_{\mu}$ is
\index{capacity $= 0$}
an $F_{\sigma}$ set of inner capacity $0$.  
\index{capacity $= 0$}

Suppose to the contrary that some $e_n$ had positive capacity.  Then
\index{capacity $> 0$}
there would be a probability measure $\eta$ supported on $e_n$ such that 
\begin{equation*}
I(\eta) \ = \ \int_{H \times H} -\log([z,w]_{\fX,\vs}) \, d\eta(z) d\eta(w) \ < \ \infty .
\end{equation*}
By the definition of $\cV_{\mu}$, there is a  
$q \in H^*$ where $f(q) > \cV_{\mu} - \frac{1}{2n}$.  By the lower 
semi-continuity of $u_{\mu}(z)$, there is also a $\delta > 0$ such that
\index{potential function!is lower semi-continuous} 
\index{semi-continuous!potential function is lower semi-continuous}
for all $z \in B(q,2\delta)$  
\begin{equation*}
f(z) \ > \ \cV_{\mu} - \frac{1}{2n} \ .
\end{equation*}
Put $A = H \bigcap B(q,\delta)$, noting that $A$ and $e_n$ 
are closed and disjoint, hence bounded 
away from each other.  Since $q \in H^*$, 
\begin{equation*}
\mu(A) \ > \ 0 \ .
\end{equation*}
Write $a = \mu(A)$, and put 
\begin{equation*}
\Delta \ = \ \left\{ 
      \begin{array}{ll} a \cdot \eta & \text{on $e_n$ \ ,} \\
                       -\mu & \text{on $A$ \ .} \end{array} \right.
\end{equation*}                       
Then $\Delta(H) = 0$, and it is easy to see that 
$I(\Delta) = \iint -\log([z,w]_{\fX,\vs}) \, d\Delta(z) d\Delta(w)$ is finite:
when it is expanded as a sum of three integrals, the diagonal terms are 
finite by hypothesis and the cross term is finite because $A$ and $e_n$
are bounded apart.    

For $0 \le t \le 1$ put 
\begin{equation*}
\mu_t \ = \ \mu + t \Delta \ . 
\end{equation*}
Then $\mu_t$ is a positive measure of total mass $\sigma$ supported on $H$.  
Comparing the integrals $I_{\sigma}(\mu_t,W)$ and $V_{\sigma}(\mu,W)$, one finds 
\index{Robin constant!weighted $(\fX,\vs)$}
\begin{equation*}
I_{\sigma}(\mu_t,W) - I_{\sigma}(\mu,W) \ = \  
  2 t \int_{H} u_{\mu}(z) + \widehat{u}(z) \, d\Delta(z) \ + \ t^2 I(\Delta) \ .
\end{equation*}
By construction 
\index{potential function}
\begin{eqnarray*}
\int_{H} u_{\mu}(z) + \widehat{u}(z) \, d\Delta(z) 
     &=& a \cdot \int_{e_n} f(z) \, d\eta(z) 
              - \int_{A} f(z) \, d\mu(z) \\
     &\le& a \cdot (\cV_{\mu} - \frac{1}{n}) 
                   - a \cdot (\cV_{\mu} - \frac{1}{2n}) 
     \ = \ - \frac{a}{2n} \ .
\end{eqnarray*}
Thus $I_{\sigma}(\mu_t,W) - I_{\sigma}(\mu,W) < 0$ for sufficiently small $t > 0$, 
contradicting the minimality of $I_{\sigma}(\mu,W)$.  It follows that 
$f(z) = u_{\mu}(z) + \widehat{u}(z) \ge \cV_{\mu}$ on $H \backslash e_{\mu}$.  

By (\ref{AFA9}), $f(z) \le \cV_{\mu}$ on $H_{\mu}^*$.  Thus $f(z) = \cV_{\mu}$ 
on $H_{\mu}^* \backslash e_{\mu}$.  Necessarily $\mu(e_{\mu}) = 0$, since 
$I(\mu) < \infty$ while $e_{\mu}$ has inner capacity $0$.  Hence 
\index{capacity $= 0$}
\begin{eqnarray*}
V_{\sigma}(H,W) &=& \iint_{H \times H} -\log([z,w]_{\fX,\vs}) \, d\mu(w) d\mu(z) 
                  + 2 \int_H \widehat{u}(z) \, d\mu(z) \\
         &=& \int_H u_{\mu}(z) + \widehat{u}(z) \, d\mu(z) + \int_H \widehat{u}(z) \, d\mu(z) \\
         &=& \sigma \cdot \cV_{\mu} + \int_H \widehat{u}(z) \, d\mu(z) \ .
\end{eqnarray*} 
\index{Robin constant!weighted $(\fX,\vs)$}
\index{potential function} 
\end{proof}


\section{ The Weighted Cheybshev Constant } \label{WeightedChebyshevSection} 
\index{Chebyshev constant!weighted|ii}

{\bf Motivation.} Consider the classical restricted Chebyshev constant, 
\index{Chebyshev constant!restricted} 
which is defined for compact sets $H \subset \CC$ by 
first putting
\begin{equation} \label{BFCH1}  
\CH^*_N(H) \ = \ \inf_{\substack{ \text{ monic $P(z) \in \CC[z]$,} \\
                               \text{ degree $N$,} \\
                               \text{ roots in $H$} }} (\|P(z)\|_E)^{1/N}
\end{equation} 
for $N = 1, 2, 3, \ldots$, 
and then setting $\CH(H) = \lim_{N \rightarrow \infty} \CH^*_N(H)$. 
 
By the compactness of $H$, for each $N$ there is a  
polynomial $P_N(z) = \prod_{i=1}^N (z-\alpha_i)$ 
which achieves the $\inf$  in (\ref{BFCH1});  
it is called a (restricted) Chebyshev polynomial.
\index{Chebyshev polynomial!restricted}   
It is known that as $N \rightarrow \infty$, the discrete measures 
$\omega_N = \frac{1}{N} \sum \delta_{\alpha_i}(z)$ associated to the $P_N(z)$ 
converge weakly to the equilibrium measure $\mu$ of $H$.
\index{equilibrium distribution}  
  
Now suppose $H = H_1 \bigcup H_2$, where the $H_i$ are closed, 
nonempty, and disjoint.  Put 
\begin{equation*}
\widehat{u}(z) \ = \ \int_{H_2} -\log(|z-w|) \, d\mu(z)   
\end{equation*}   
and set $W(z) = \exp(-\widehat{u}(z))$.   
Label the roots of $P_N(z)$ so  
$\alpha_1, \ldots, \alpha_{n} \in H_1$ 
and $\alpha_{n+1}, \ldots, \alpha_N \in H_2$.  
For large $N$, one has  
$\prod_{i = n + 1}^N |z-a_i| \cong \exp(-N\widehat{u}(z))$ 
outside $H_2$.  Thus, on $H_1$  
\begin{equation}
|P_N(z)| \ \cong \ \prod_{i=1}^{n} |z-\alpha_i| \cdot W(z)^N \ . 
                                    \label{AFB11}
\end{equation}           
If one thinks of $W(z)$ as given, then the minimization problem (\ref{BFCH1})   
can be thought of as varying $\alpha_1, \ldots, \alpha_n$ over $H_1$ 
to achieve the minimum in (\ref{BFCH1}).   

\medskip
{\bf Definition of the weighted $(\fX,\vs)$-Chebyshev Constant.}
Let $H$ be a compact subset of $\cC_v(\CC) \backslash \fX$, and let
$W(z) : \cC_v(\CC) \rightarrow [0,\infty]$ be a function which is positive, 
bounded and continuous on a neighborhood of $H$.  
Let $[z,w]_{\fX,\vs}$ be as before.
By a weighted $(\fX,\vs)$-pseudopolynomial of 
\index{pseudopolynomial!weighted $(\fX,\vs)$}
bidegree $(n,N)$ (relative to $W(z)$) we mean a function of the form
\begin{equation*}
P(z) \ = \ P_{(n,N)}(z,W) \ = \ 
\prod_{i = 1}^{n} [z,\alpha_i]_{\fX,\vs} \cdot W(z)^N  \ .
\end{equation*}
The $\alpha_i$ will be called the roots of $P(z)$.   

To simplify notation,
for the rest of this section we will generally omit explicit mention of the 
$(\fX,\vs)$-dependence in the quantities discussed.    

Since $[z,w]_{\fX,\vs}$ and $W(z)$ are continuous on $H$, and $H$ is compact, 
among the weighted pseudopolynomials $P_{(n,N)}(z,W)$ of bidegree $(n,N)$
\index{pseudopolynomial!weighted $(\fX,\vs)$}
whose roots belong to $H$, there is at least one with minimal minimal sup norm 
on $H$ (it need not be unique).  Fix one, and write $\tP_{(n,N)}(z,W)$ for it.    
We will call it a weighted, restricted Chebyshev pseudopolynomial.  
\index{pseudopolynomial!weighted Chebyshev}
\index{Chebyshev pseudopolynomial!weighted}

Put
\begin{eqnarray}
\CH^*_{(n,N)} & = & \CH^*_{(n,N)}(H,W) 
     \ := \ \min_{\substack{ P_{(n,N)} \\ \text{roots in $H$} }} 
                       (\|P_{(n,N)}(z,W)\|_H)^{1/N} \notag \\
           & = & \left(\|\tP_{(n,N)}\|_H\right)^{1/N} \ .  \label{FMe1} 
\end{eqnarray}
For each $\sigma > 0$, 
define the weighted Chebyshev constant $\CH^*_{\sigma}(H,W)$ by
\index{Chebyshev constant!weighted} 
\begin{equation}
\CH^*_{\sigma}(H,W) \ = \ 
   \lim_{\substack{ n/N \rightarrow \sigma \\ N \rightarrow \infty }}
               \CH^*_{(n,N)} \label{AFB12}
\end{equation}
provided the limit exists, that is, if for each $\varepsilon > 0$, 
there are a $\delta > 0$ and an $N_0$ such that if
$|\frac{n}{N} - \sigma| < \delta$ and $N \ge N_0$,
then $|\CH^*_{(n,N)} - \CH^*_{\sigma}(H,W)| < \varepsilon$.  

\begin{theorem} \label{ATB1} For each $\sigma > 0$, $\CH^*_{\sigma}(H,W)$ exists.
The function $g(r) = \CH^*_r(H,W)$ is continuous for $r > 0$.  
If $g(r_0) = 0$ for one $r_0$, then $g(r) \equiv 0$.  
Moreover, $g(r) \equiv 0$ if and only if $H$ has capacity $0$. 
\index{capacity $= 0$} 
\end{theorem}

Recall that $H$ has capacity $0$ iff for some
\index{capacity $= 0$} 
$\zeta \in \cC_v(\CC) \backslash H$, we have $\gamma_{\zeta}(H) = 0$ 
in the sense of (\cite{RR1}, \S3.1);   
this holds for one $\zeta \notin H$ iff it holds for all $\zeta \notin H$.  

\vskip .1 in
\begin{proof}
First fix $0 < r \in \QQ$, and take $(n,N)$ with $n/N = r$.  Given a
rational number $\kappa > 0$ such that $\kappa n, \kappa N \in \ZZ$, we can write
\begin{eqnarray*}
\kappa n &=& k n + a \qquad \text{where $k = \lfloor \kappa  \rfloor \in \ZZ$
                                 and $0 \le a < n$ \ ,}  \\
\kappa N &=& kN + A \qquad \text{where $0 \le A = a/r < N$ \ .}
\end{eqnarray*}
By the compactness of $H$ and our assumptions on $W(z)$, 
there is a constant $C > 1$ such that $1/C \le W(z) \le C$ for all $z \in H$.
Put $D = \max(1,\max_{z,w \in H} [z,w]_{\fX,\vs})$. 
If $\tP_{(n,N)}(z,W) = \prod_{i = 1}^{n} [z,\alpha_i]_{\fX,\vs} \cdot W(z)^N$, then
\begin{eqnarray*}
(\CH^*_{(\kappa n,\kappa N)})^{\kappa N}
  &\le& \max_{z \in H} \left( \prod_{i = 1}^{n} [z,\alpha_i]_{\fX,\vs}^k \cdot
       \prod_{i=1}^{a} [z,\alpha_i]_{\fX,\vs} \cdot W(z)^{\kappa N} \right) \\
  &=& \max_{z \in H}
      \left( \prod_{i = 1}^{n} [z,\alpha_i]_{\fX,\vs} \cdot W(z)^N \right)^k
        \cdot \prod_{i=1}^{a} [z,\alpha_i]_{\fX,\vs} \cdot W(z)^{A} \\
  &\le& (\CH^*_{(n,N)})^{kN} \cdot D^{a} \cdot C^{A} \ .
\end{eqnarray*}
Thus
\begin{equation*}
\CH^*_{(\kappa n,\kappa N)} \ \le \ (\CH^*_{(n,N)})^{\frac{k}{\kappa }} \cdot D^{\frac{a}{\kappa N}}
         \cdot C^{\frac{A}{\kappa N}} \ .
\end{equation*}
Letting $\kappa  \rightarrow \infty$ gives
\begin{equation*}
\CH^*_{(n,N)} \ \ge \ \limsup_{\kappa  \rightarrow \infty} \CH^*_{(\kappa n,\kappa N)} \ , 
\end{equation*}
then replacing $(n,N)$ by $(\kappa n,\kappa N)$ and again letting
$\kappa  \rightarrow \infty$ gives
\begin{equation*}
\liminf_{\kappa  \rightarrow \infty} \CH^*_{(\kappa n,\kappa N)} \ \ge \
       \limsup_{\kappa  \rightarrow \infty} \CH^*_{(\kappa n,\kappa N)} \ .
\end{equation*}
Consequently, as $\kappa $ passes through all rational numbers with  
$\kappa n,\kappa N \in \ZZ$, the limit 
\begin{equation} \label{Fg0}
g(r) \ := \ \lim_{\kappa  \rightarrow \infty} \CH^*_{(\kappa n,\kappa N)}
\end{equation}
exists, and the same limit is obtained when $\kappa $ passes through {\emph{any}}
sequence of rationals with $\kappa  \rightarrow \infty$ and $\kappa n, \kappa N \in \ZZ$.

We now compare $g(r)$ and $g(s)$, when $0 < r, s$ are rationals
with $r \ne s$.  Fix positive integers $n, N$ with $r = n/N$
and $m, M$ with $s = m/M$.  Let $\lambda $ be a positive integer such
that $\lambda n \ge m$.  Write
\begin{equation*}
\lambda n = \ell m + b \qquad 
      \text{where $\ell = \left\lfloor \frac{\lambda n}{m} \right\rfloor$
                                 and $0 \le b < m$ \ .}  
\end{equation*}
If $\tP_{(m,M)}(z,W) = \prod_{i=1}^{m} [z,\alpha_i] \cdot W(z)^M$, we have
\begin{eqnarray*}
(\CH^*_{(\lambda n,\lambda N)})^{\lambda N}
  &\le& \max_{z \in H} \left( \prod_{i = 1}^{m} [z,\alpha_i]_{\fX,\vs}^{\ell} \cdot
       \prod_{i=1}^{b} [z,\alpha_i]_{\fX,\vs} \cdot W(z)^{\lambda N} \right) \\
  &=& \max_{z \in H}
      \left( \prod_{i = 1}^{m} [z,\alpha_i]_{\fX,\vs} \cdot W(z)^M \right)^{\ell}
        \cdot \prod_{i=1}^{b} [z,\alpha_i]_{\fX,\vs} \cdot W(z)^{\lambda N-{\ell}M} \\
  &\le& (\CH^*_{(m,M)})^{{\ell}M} \cdot D^{b} \cdot C^{|\lambda N-{\ell}M|} \ .
\end{eqnarray*}
The absolute value appears in  the last term because we use 
$W(z) \le C$ to obtain the final inequality  if $\lambda N-{\ell}M > 0$, 
and $W(z) > 1/C$ if $\lambda N-{\ell}M < 0$. Consequently
\begin{equation} \label{AFB13}
(\CH^*_{(\lambda n,\lambda N)})^{\frac{\lambda N}{\lambda n}}
\ \le \ CH_{(m,M)}^{\frac{{\ell}M}{{\ell}m + b}} \cdot D^{\frac{b}{{\ell}m + b}}
                 \cdot C^{| \frac{\lambda N}{\lambda n} - \frac{{\ell}M}{{\ell}m + b}|} \ .
\end{equation}
Letting $\lambda  \rightarrow \infty$ (so ${\ell} \rightarrow \infty$ as well), we get
\begin{equation} \label{AFB14}
g(r)^{1/r} \ \le \ (\CH^*_{(m,M)})^{\frac{M}{m}}
               \cdot C^{|\frac{1}{r} - \frac{M}{m}|} \ ;
\end{equation}
then replacing $(m,M)$ by $(\kappa m,\kappa M)$ 
and letting $\kappa \rightarrow \infty$,
\begin{equation} \label{AFB15}
g(r)^{1/r} \ \le \ g(s)^{1/s} \cdot C^{|\frac{1}{r} - \frac{1}{s}|} \ .
\end{equation}
Reversing the roles of $r$ and $s$, we similarly obtain
\begin{equation} \label{AFB16}
g(s)^{1/s} \ \le \ g(r)^{1/r} \cdot C^{|\frac{1}{s} - \frac{1}{r}|} \ .
\end{equation}
From (\ref{AFB15}) and (\ref{AFB16}) it follows that if $g(r_0) = 0$ for
one $r_0$, then $g(r) \equiv 0$.

Suppose $g(r) \ne 0$.  Taking logarithms in (\ref{AFB15}) and (\ref{AFB16})
we see that
\begin{equation} \label{ABug1}
\left| \frac{1}{r} \log(g(r)) - \frac{1}{s} \log(g(s)) \right|
 \ \le \ \left| \frac{1}{r} - \frac{1}{s} \right| \log(C) \ .
\end{equation}
Thus if $\{r_i\}_{1 \le i < \infty}$ is a Cauchy
sequence of rationals converging to some real $r > 0$, then
$\{g(r_i)\}$ is also a Cauchy sequence;  given two Cauchy sequences
converging to the same $r$, the sequences of values $g(r_i)$ converge to
the same value.   Thus, we can extend $g(r)$ to a function
defined on all positive reals.  It is easy to see that
the extended function satisfies (\ref{ABug1})
for all real $r, s > 0$ and hence is continuous.  

Trivially $g(r)$ can be extended continuously to all $r > 0$
if $g(r) \equiv 0$ on $\QQ$.

We will now show that the limit (\ref{AFB12}) exists, and that 
\begin{equation}
\lim_{\substack{ n/N \rightarrow r \\ N \rightarrow \infty }}
               \CH^*_{(n,N)} \ = \ g(r) \ . \label{AFB17}
\end{equation}
First suppose $g(r) \ne 0$.  Then (\ref{AFB17}) is equivalent to
\begin{equation}
\lim_{\substack{ n/N \rightarrow r \\ N \rightarrow \infty }}
          \frac{N}{n} \log(\CH^*_{(n,N)}) \ = \ \frac{1}{r} \log(g(r)) \ . 
                    \label{AFB18}
\end{equation}
To prove (\ref{AFB18}), fix $\varepsilon > 0$ and take $r_1, r_2 \in \QQ$
with $0 < r_1 < r < r_2$ such that
\begin{eqnarray}
\left| \frac{1}{r_1} \log(g(r_1)) - \frac{1}{r} \log(g(r)) \right|
  &<& \frac{\varepsilon}{5} \ , \label{AFB19} \\
\text{and} \qquad \left| \frac{1}{r_1} - \frac{1}{r_2} \right| \log(C) 
                    \ < \ \frac{\varepsilon}{5} \ .
                                \label{AFB20}
\end{eqnarray}
By (\ref{Fg0}), there is a pair $(m_0,M_0)$ with $m_0/M_0 = r_1$, such that
\begin{equation} \label{Fg1}
|\frac{M_0}{m_0} \log(\CH^*_{(m_0,M_0)}) - \frac{1}{r_1}\log(g(r_1))| \ < \ 
                   \frac{\varepsilon}{5} \ .
\end{equation}
Fix an integer $N_0 \ge M_0$ large enough that for each integer 
$\ell \ge \lfloor N_0/M_0 \rfloor$,  
\begin{equation} 
 \frac{1}{\ell} \log(D) \ < \ \frac{\varepsilon}{5} \ ,  \label{Fg2} 
\end{equation}
\begin{equation}
|\frac{\ell M_0}{\ell m_0 + m_0} - \frac{M_0}{m_0}| 
     \cdot |\log(\CH^*_{(m_0,M_0)})| \ < \ \frac{\varepsilon}{5} \ ,  \label{Ft3}
\end{equation}
\begin{equation} 
\frac{\ell M_0} {\ell m_0 + m_0} \ >  \ \frac{1}{r_2}  \ , \label{Ft4}
\end{equation}
Consider any pair $(n,N)$ with $N \ge N_0$ and $r_1 < \frac{n}{N} < r_2$.  

Applying (\ref{AFB14}), with $r$ replaced by $r_1$  
and $(m,M)$ replaced by $(n,N)$, then taking logarithms, we get 
\begin{equation} \label{AFB21}
\frac{N}{n} \log(\CH^*_{(n,N)}) \ \ge \ \frac{1}{r_1} \log(g(r_1))
        - \left| \frac{1}{r_1} - \frac{N}{n} \right| \log(C) \ .
\end{equation}
Combining (\ref{AFB21}), (\ref{AFB19}) and (\ref{AFB20}) gives
\begin{equation} \label{AFB22}
\frac{N}{n} \log(\CH^*_{(n,N)})
     \ > \ \frac{1}{r} \log(g(r)) - \frac{2 \varepsilon}{5} \ .
\end{equation}

For the opposite equality, note that since $m_0/M_0 = r_1 < n/N < r_2$ and
$N \ge N_0> M_0$, we have $n \ge m_0$.  Apply (\ref{AFB13}),   
taking  $(m,M) = (m_0,M_0)$ and $\lambda  = 1$, so that   
$n = \ell m_0 + b$, 
with $\ell = \lfloor n/m_0 \rfloor \ge \lfloor N_0/M_0 \rfloor$ and $0 \le b < m_0$;
this yields    
\begin{eqnarray} \label{AFB23}
\lefteqn{ \frac{N}{n} \log(\CH^*_{(n,N)}) } \\ 
&\le& \frac{{\ell} M_0}{{\ell} m_0 + b} \log(\CH^*_{(m_0,M_0)})
        +  \frac{b}{{\ell} m_0 + b} \log(D) 
      + \left| \frac{N}{n} - \frac{{\ell}M_0}{{\ell}m_0 + b} \right| \log(C) \ .
                      \notag
\end{eqnarray} 
Using (\ref{Fg2}) we see that 
\begin{equation} \label{Fg3}
\frac{b}{\ell m_0 + b} \log(D) \ \le \ \frac{1}{\ell} \log(D) 
       \ < \ \frac{\varepsilon}{5} \ .  
\end{equation} 
Likewise $1/r_1 > N/n > 1/r_2$, 
and $1/r_1 > \ell M_0/(\ell m_0 + b) > 1/r_2$ by (\ref{Ft4}) ,
so by (\ref{AFB20})
\begin{equation} \label{Ft5} 
\left| \frac{N}{n} - \frac{\ell M_0}{ \ell m_0 + b} \right| \log(C) 
\ < \ \frac{\varepsilon}{5}
\end{equation}  
Combining (\ref{AFB23}), (\ref{AFB19}),  (\ref{AFB20}),  (\ref{Ft3}),
 (\ref{Fg3}), and (\ref{Ft5}) gives 
\begin{equation*}
\frac{N}{n} \log(\CH^*_{(n,N)})
    \  < \ \frac{1}{r} \log(g(r)) +  \varepsilon \ .
\end{equation*}  

In the case where $g(r) \equiv 0$, first note that if $H$ is finite,
then $\CH^*_{(n,N)} = 0$ whenever $n \ge \#(H)$, so (\ref{AFB17}) is trivial.
If $H$ is infinite, then  $\CH^*_{(n,N)} > 0$ for all $(n,N)$.  Take $0 < \varepsilon < 1$. 
Fix $r_1, r_2 \in \QQ$ with $0 < r_1 < r < r_2$ such that
\begin{equation} \label{AFB28}
 C^{| \frac{1}{r_1} - \frac{1}{r_2} |} \ < \ 2
\end{equation}
and fix $(m_0,M_0)$ such that $m_0/M_0 = r_1$ 
and $(\CH^*_{(m_0,M_0)})^{M_0/m_0} < \varepsilon/5$.  Let $N_0$ be  
large enough that for each integer $\ell \ge \lfloor N_0/M_0\rfloor$ we have
\begin{equation} \label{fD1} 
(\CH^*_{(m_0,M_0)})^{\frac{\ell M_0}{\ell m_0+ m_0}} \ < \ \frac{\varepsilon}{4} \ ,
\end{equation} 
\begin{equation} \label{fD2}
D^{\frac{1}{\ell}} \ < \ 2 \ , \qquad
\frac{\ell M_0}{\ell m_0+ m_0} \ > \ \frac{1}{r_2} \ .
\end{equation} 
Consider any pair $(n,N)$ with $r_1 < n/N < r_2$ and $N \ge N_0$.  
Again apply (\ref{AFB13}), 
taking  $(m,M) = (m_0,M_0)$ and $\lambda  = 1$, 
so that  $n = \ell m_0 + b$, 
with $\ell = \lfloor n/m_0 \rfloor \ge \lfloor N_0/M_0 \rfloor$ 
and $0 \le b < m_0$; this gives 
\begin{equation} \label{AFB29}
(\CH^*_{(n,N)})^{\frac{N}{n}}
\ \le \ (\CH^*_{(m_0,M_0)})^{\frac{\ell M_0}{\ell  m_0 + b}} 
             \cdot D^{\frac{b}{\ell m_0 + b}}
                 \cdot C^{| \frac{N}{n} - \frac{\ell M_0}{\ell m_0 + b}|} \ .
\end{equation}
By (\ref{AFB28}) and  (\ref{fD2}), 
\begin{equation} \label{AFB31}
D^{\frac{b}{\ell m_0 + b}}
          \cdot C^{| \frac{N}{n} - \frac{\ell M_0}{\ell m_0 + b}|} \ < \ 4 \ .
\end{equation}
Combining (\ref{AFB29}), (\ref{fD1}) and (\ref{AFB31}) gives 
\begin{equation*}
(\CH^*_{(n,N)})^{\frac{N}{n}} \ < \ \varepsilon \ .
\end{equation*}

Finally, let us show that $g(r) \equiv 0$ if and only if $H$ has capacity $0$.
\index{capacity $= 0$}
If instead of $[z,w]_{\fX,\vs}$ and $W(z)$ we had used another distance function
\begin{equation*}
[z,w]_{\fX,\vs}^* \ = \ \prod_{i = 1}^{m^*} ([z,w]_{x_i^*})^{s_i^*}
\end{equation*}
and another continuous, positive background function $W^*(z)$, then
there would be a constant $A > 0 $ such that
\begin{eqnarray*}
& & 1/A \ < \ [z,w]_{\fX,\vs}^*/[z,w]_{\fX,\vs} \ < \ A \ , \\
& & \ \ \ 1/A \ < \ W^*(z)/W(z)  \ < \ A \ ,
\end{eqnarray*}
for all $z \ne w \in H$.  Let $\CH^*_{(n,N)}(N,W^*)$ and $\CH^*_{r}(H,W^*)$
denote the
weighted Chebyshev constants computed relative to $[z,w]_{\fX,\vs}^*$ and $W^*(z)$.
\index{Chebyshev constant!weighted} 
Taking $\alpha_1, \ldots, \alpha_{n} \in H$
such that $\tP_{(n,N)}(z,W) = \prod_{i = 1}^{n} [z,\alpha_i]_{\fX,\vs} \cdot W(z)^N$,
we get
\begin{equation*}
\CH^*_{(n,N)}(H,W^*) \ \ge \ 
      \CH^*_{(n,N)}\cdot A^{- \frac{n}{N} - 1 } .
\end{equation*}
A similar inequality holds with $\CH^*_{(n,N)}(H,W^*)$ and $\CH^*_{(n,N)}$ reversed, and 
it follows that $\CH^*_{r}(H,W^*) \ne 0$ iff $\CH^*_{r}(H,W) \ne 0$.

  Taking $r = 1$, $W^*(z) = 1$, and
$[z,w]_{\fX,\vs}^* = [z,w]_{\zeta}$ for a fixed $\zeta \in \cC_v(\CC) \backslash H$, 
we have $\CH^*_1(H,W^*) = \gamma_{\zeta}(H)$ by (\cite{RR1}, Theorem 3.1.18).  
Thus $g(r) \not \equiv 0$ if and only if for some (hence any) $\zeta \notin H$, 
we have $\gamma_{\zeta}(H) > 0$.  
\end{proof}


\section{ The Weighted Transfinite Diameter } \label{WeightedTDSection}  

{\bf Motivation.}
Consider the classical transfinite diameter for compact sets $H \subset \CC$, 
\index{transfinite diameter}
defined by first putting 
\begin{equation} \label{AFC32}
d_N(H) \ = \ \left( \sup_{z_1, \ldots,z_N \in H}
     \prod^N_{\substack{ i, j = 1 \\ i \ne j }} |z_i-z_j|  \right)^{1/N^2} \ ,
\end{equation}
and then setting $d_{\infty}(H) \ = \ \lim_{N \rightarrow \infty} d_N(H)$.
(Usually the exponent $1/N^2$ in (\ref{AFC32}) is replaced by $1/N(N-1)$, 
but this does not affect the value of the limit.)

Since $H$ is compact, for each $N$ 
there exist points $\alpha_1, \ldots, \alpha_N \in H$ 
realizing the supremum in (\ref{AFC32}): 
$d_N(H) = (\prod^N_{i, j = 1, i \ne j} |\alpha_i - \alpha_j|)^{1/N^2}$.
The collection $\{\alpha_1, \ldots, \alpha_n\}$ is called a set of Fekete points. 
\index{Fekete!points|ii}  
 
Now suppose $H = H_1 \bigcup H_2$, where the $H_i$ are closed, nonempty,
and disjoint.  We can write
\begin{equation} \label{AFC33}
d_N(H) 
\ = \ \prod_{\substack{ \alpha_i, \alpha_j \in H_1 \\ i \ne j }} |\alpha_i-\alpha_j| 
\cdot \prod_{\alpha_i \in H_1} \prod_{\alpha_j \in H_2} |\alpha_i-\alpha_j|^2
\cdot \prod_{\substack{ \alpha_i, \alpha_j \in H_2 \\ i \ne j }} |\alpha_i-\alpha_j| \ .
\end{equation}
Let $\nu_N = \sum_{i=1}^N \frac{1}{N} \delta_{\alpha_i}(z)$ be 
the probability measure equally supported on the $\alpha_i$.  It is known that
the $\nu_N$ converge weakly to the equilibrium
distribution $\mu$ of $H$.  As before, put 
\index{equilibrium distribution}
\begin{equation*}
\widehat{u}(z) \ = \ \int_{H_2} -\log([z,w]_{\fX,\vs}) \, d\mu(w)
\end{equation*} 
and set $W(z) = \exp(-\widehat{u}(z))$. Then 
\begin{equation*}
\prod_{\alpha_i \in H_1} \prod_{\alpha_j \in H_2} (|\alpha_i-\alpha_j|)^2
\ \cong \ \prod_{\alpha_i \in H_1} W(\alpha_i)^{2N} \ .
\end{equation*}

Label the Fekete points so that $\alpha_1, \ldots, \alpha_n \in H_1$ 
\index{Fekete!points} 
and $\alpha_{n+1}, \ldots, \alpha_N \in H_2$.   
If one is interested in $H_1$ and thinks of $\alpha_{n+1}, \ldots, \alpha_N$ as given, 
then the maximization problem (\ref{AFC32}) 
can be thought of as varying $z_1, \ldots, z_n$ over $H_1$ 
so as to maximize
\begin{equation*}
\prod^{n}_{\substack{ i, j = 1 \\ i \ne j }} |z_i-z_j| 
\cdot \prod_{i=1}^{n_1} W(z_i)^{2N} \ . 
\end{equation*}
Note that as $N \rightarrow \infty$, then $n/N \rightarrow \sigma := \mu(H_1)$.

\medskip
{\bf Definition of the Weighted Transfinite Diameter $d_{\sigma}(H,W)$.}
Let $H \subset \cC_v(\CC) \backslash \fX$ be compact, 
and let $W(z) : \cC_v(\CC) \rightarrow [0,\infty]$ be a function which 
is positive, bounded and  continuous on a neighborhood of $H$.  

For any positive integers $n$ and $N$, define 
\begin{equation} \label{AFC34}
d_{(n,N)} \ = \ d_{(n,N)}(H,W) \ =  \
         \max_{\substack{ z_1, \ldots, z_{n} \in H }} 
    \left( \prod^{n}_{\substack{ i, j = 1 \\ i \ne j }} [z_i,z_j]_{\fX,\vs} 
             \cdot \prod_{i=1}^{n} W(z_i)^{2N} \right)^{1/N^2} \ .
\end{equation}
Then, for each $0 < \sigma \in \RR$ let the weighted 
$(\fX,\vs)$-transfinite diameter be 
\index{transfinite diameter!weighted $(\fX,\vs)$}
\begin{equation} \label{AFC35}
d_{\sigma}(H,W) \ = \ \lim_{\substack{ n/N \rightarrow \sigma \\
                           N  \rightarrow \infty }}
                                 d_{(n,N)} 
\end{equation}
provided the limit exists.  This is understood to mean that for each
$\varepsilon > 0$, there exist a $\delta > 0$ and an $N_0$ such that if
$|\frac{n}{N} - r| < \delta$ and $N \ge N_0$,
then $|d_{(n,N)}  - d_{\sigma}(H,W)| < \varepsilon$.  
For notational simplicity we suppress the $(\fX,\vs)$-dependence 
in $d_{(n,N)} $ and $d_{\sigma}(H,W)$.  

A set of points $\alpha_1, \ldots, \alpha_n \in H$ which achieve the maximum value in (\ref{AFC34})
will be called a set of $(n,N)$-Fekete points for $H$ relative to the weight $W(z)$.
\index{Fekete!points!$(n,N)$}  
 
\begin{theorem} \label{ATC3} For each $\sigma > 0$, $d_{\sigma}(H,W)$ exists.
The function $f(r) = d_r(H,W)$ is continuous for $r > 0$;  and if
$f(r_0) = 0$ for one $r_0$, then $f(r) \equiv 0$.  Moreover, $f(r) \equiv 0$ 
if and only if $H$ has capacity $0$.  
\index{capacity $= 0$}
\end{theorem}

\begin{proof}
First consider what happens to $d_{(n,N)} $ when
$(n,N)$ is replaced by $(\lambda n,\lambda N)$ for some rational $\lambda  \ge 1$  for
which $\lambda n$ and $\lambda N$ are integers.
If $\alpha_1, \ldots, \alpha_{\lambda n} \in H$ are chosen to maximize 
(\ref{AFC34}) for $d_{(\lambda n,\lambda N)} $, then
\begin{equation}
d_{(\lambda n,\lambda N)}^{(\lambda N)^2} \ = \ 
  \prod^{\lambda n}_{\substack{ i, j = 1 \\ i \ne j }} [\alpha_i,\alpha_j]_{\fX,\vs} 
             \cdot \prod_{i=1}^{\lambda n} W(\alpha_i)^{2\lambda N} \ . \label{AFC36}
\end{equation}
For any $n$-element subset $\alpha_{k_1}, \ldots, \alpha_{k_{n}}$ we have 
\begin{equation*}
d_{(n,N)}^{N^2} \ \ge \ 
  \prod^{n}_{\substack{ i, j = 1 \\ i \ne j }} [\alpha_{k_i},\alpha_{k_j}]_{\fX,\vs} 
             \cdot \prod_{i=1}^{n} W(\alpha_{k_i})^{2N} \ .
\end{equation*}
There are $\binom{\lambda n}{n}$ such subsets; each pair $\{\alpha_i,\alpha_j\}$ belongs to 
$\binom{\lambda n - 2}{n - 2}$ subsets, and each $\alpha_i$ belongs to 
$\binom{\lambda n - 1}{n - 1}$ subsets.  If $C > 1$ is such that 
$1/C \le W(z) \le C$ on $H$, then 
\begin{eqnarray}
(d_{(n,N)}^{N^2})^{\binom{\lambda n}{n}}  
 &\ge& \left( \prod^{\lambda n}_{\substack{ i, j = 1 \\ i \ne j }} 
        [\alpha_i,\alpha_j]_{\fX,\vs} \right)^{ {\binom{\lambda n - 2}{n - 2}} \cdot 
  \left( \prod_{i=1}^{\lambda n} W(\alpha_i)^{2N} \right)^{\binom{\lambda n - 1}{n - 1} } }
       \notag \\
 &=& \left( d_{(\lambda n,\lambda N)}  \right)^{(\lambda N)^2 {\binom{\lambda n - 2}{n - 2} } } 
    \cdot \left( \prod_{i=1}^{\lambda n} W(\alpha_i))\right)^{ 
   2N { \binom{\lambda n - 1}{n - 1} } - 2\lambda N {\binom{\lambda n - 2}{n - 2} } } 
          \notag \\ 
 &\ge& \left( d_{(\lambda n,\lambda N)} \right)^{(\lambda N)^2 {\binom{\lambda n - 2}{n - 2} } }  
    \cdot C^{- \lambda n \left(
    2N{ \binom{\lambda n - 1}{n - 1} - 2\lambda N { \binom{\lambda n - 2}{n - 2 }}  } \right) }
       \ . \label{AFC37} 
\end{eqnarray}
Simplifying exponents, we find  
\begin{equation} \label{AFC38}
d_{(n,N)}^{\frac{\lambda n - 1}{\lambda (n - 1)}} 
       \cdot C^{\frac{2n(\lambda -1)}{\lambda N(n-1)}} \ \ge \ d_{(\lambda n,\lambda N)}  \ .
\end{equation}

This holds for any $(n,N)$ and $\lambda $ satisfying the conditions above.
Given $0 < \sigma \in \QQ$, take $n$, $N$ with 
$\sigma = n/N$.  Letting $\lambda  \rightarrow \infty$ in (\ref{AFC38}) (where $\lambda $
runs over all rationals such that $\lambda n, \lambda N \in \ZZ$), we find
\begin{equation} \label{AFC39}
d_{(n,N)}^{\frac{n}{n-1}} \cdot C^{\frac{2n}{N(n - 1)}} 
\ \ge \ \limsup_{\lambda  \rightarrow \infty}\, d_{(\lambda n,\lambda N)}  \ .
\end{equation}
Replacing $(n,N)$ by $(\kappa n,\kappa N)$ on the left side of (\ref{AFC39}) 
and again letting $\kappa  \rightarrow \infty$ gives
\begin{equation*}
\liminf_{\kappa  \rightarrow \infty} \, d_{(\kappa n,\kappa N)}  \ \ge \
       \limsup_{\lambda  \rightarrow \infty} \, d_{(\lambda n,\lambda N)}  \ .
\end{equation*}
Hence 
\begin{equation} \label{FGo1} 
f(r) \ := \ \lim_{\lambda  \rightarrow \infty} d_{(\lambda n,\lambda N)} 
\end{equation}
is well-defined;  moreover, the same limit is obtained when $\lambda $ passes
through any sequence of values such that $\lambda  \rightarrow \infty$ and
$\lambda n, \lambda N$ are integers.

We now seek to compare $f(r)$ and $f(s)$, when $0 < r, s$ are rationals
with $r \ne s$.  Fix positive integers $n, N$ with $r = n/N$
and $m, M$ with $s = m/M$.  Let $\lambda $ be a positive integer such
that $\lambda n \ge m$.  As before, let $\alpha_1, \ldots, \alpha_{\lambda n} \in H$
realize the maximum
in  (\ref{AFC34}) for $d_{(\lambda n,\lambda N)} $.  Then as in (\ref{AFC37}), 
\begin{eqnarray*}
 d_{(m,M)}^{ M^2 { \binom{\lambda n}{m} } }
   & \ge & \left( d_{(\lambda n,\lambda N)} \right)^{(\lambda N)^2 { \binom{\lambda n - 2}{m - 2} } }
    \cdot \left( \prod_{i = 1}^{\lambda n} W( \alpha_i ) \right)^{ 
              \left( 2M { \binom{\lambda n - 1}{m - 1} } 
                     - 2\lambda N { \binom{\lambda n - 2}{m - 2} } \right) } \\
  & \ge & \left( d_{(\lambda n,\lambda N)} \right) ^{ (\lambda N)^2 {\binom{\lambda n - 2}{m - 2}} }
   \cdot C^{ - | \lambda n \left( 2M  \binom{\lambda n - 1}{m - 1} 
                   - 2 \lambda N \binom{\lambda n - 2}{m - 2} \right) |  } \ .
\end{eqnarray*}
The absolute value in the last term occurs since we use 
$W(z) > 1/C$ to obtain
the final inequality if the exponent in the previous line is positive, and 
$W(z) < C$ if it is negative.  Simplifying exponents gives
\begin{equation} \label{AFC40}
d_{(m,M)}^{\frac{M}{m} \cdot \frac{M}{m - 1} }
 \ \ge \ d_{(\lambda n,\lambda N)}^{\frac{\lambda N}{\lambda n} \cdot \frac{\lambda N}{\lambda n - 1} }
    \cdot C^{-2 | \frac{M}{m - 1}
                             - \frac{\lambda N}{\lambda n - 1} | } \ .
\end{equation}
Letting $\lambda  \rightarrow \infty$ in (\ref{AFC40}) gives 
\begin{equation} \label{AFC41}
d_{(m,M)}^{\frac{M}{m} \cdot \frac{M}{m - 1} }
 \ \ge \ f(r)^{1/r^2} \cdot C^{-2 | \frac{M}{m - 1} - \frac{1}{r} | } \ ;
\end{equation}
then, replacing $(m,M)$ by $(\kappa m,\kappa M)$ and letting $\kappa  \rightarrow \infty$
in (\ref{AFC41}), we find 
\begin{equation} \label{AFC42}
f(s)^{1/s^2}
 \ \ge \ f(r)^{1/r^2} \cdot C^{-2 | \frac{1}{s} - \frac{1}{r} | } \ .
\end{equation}
Interchanging the role of $r$ and $s$ in (\ref{AFC42}), also
\begin{equation} \label{AFC43}
f(r)^{1/r^2}
 \ \ge \ f(s)^{1/s^2} \cdot C^{-2 | \frac{1}{r} - \frac{1}{s} | } \ .
\end{equation}

It follows from (\ref{AFC42}) and (\ref{AFC43}) that $f(r) \equiv 0$ if and
only if $f(r_0) = 0$ for one $r_0$.  Suppose $f(r) \not \equiv 0$.  
Taking logarithms in (\ref{AFC42}) and (\ref{AFC43}), we find 
\begin{equation} \label{AFC44}
\left| \frac{1}{s^2} \log(f(s)) - \frac{1}{r^2} \log(f(r)) \right|
 \ \le \ 2 \left| \frac{1}{s} - \frac{1}{r} \right| \log(C) \ .
\end{equation}

From (\ref{AFC44}) we see that if $\{r_i\}_{1 \le i < \infty}$ is a Cauchy
sequence of rationals converging to some $r > 0$, then
$\{f(r_i)\}$ is also a Cauchy sequence;  and given two Cauchy sequences
converging to the same $r$, the sequences of values $f(r_i)$ converge to
the same value.   Thus, we can extend $f(r)$ to a function 
defined on all positive reals, which is easily seen to satisfy (\ref{AFC44}) 
for all real $r, s > 0$ and hence is continuous.  
Note that trivially $f(r)$ can be extended by continuity if $f(r) \equiv 0$ on $\QQ$.

We will now show that the limit (\ref{AFC35}) exists.   
Fix $0 < r \in \RR$.  We claim that 
\begin{equation} \label{AFC45}
\lim_{\substack{ n/N \rightarrow r \\
            N  \rightarrow \infty }} d_{(n,N)}  \ =  \ f(r) \ .
\end{equation}
To show this, first suppose $f(r) \ne 0$.
Then (\ref{AFC45}) is equivalent to
\begin{equation} \label{AFC46}
\lim_{\substack{ n/N \rightarrow r \\
            N  \rightarrow \infty }}
        \frac{N}{n}\frac{N}{n - 1} \log\left(d_{(n,N)}  \right)
             \ =  \ \frac{1}{r^2} \log(f(r)) \ .
\end{equation}
Given $\varepsilon > 0$, fix $r_1, r_2 \in \QQ$ with
$0 < r_1 < r < r_2$ such that
\begin{equation} \label{AFC47}
\left| \frac{1}{r_1^2} \log(f(r_1)) - \frac{1}{r^2} \log(f(r)) \right|
       \ < \ \frac{\varepsilon}{6}  \ , 
\end{equation}
\begin{equation} \label{AFC48} 
 2 \left| \frac{1}{r_1} - \frac{1}{r_2} \right| \log(C)
       \ < \ \frac{\varepsilon}{6} \ .  
\end{equation}
Then, using (\ref{FGo1}) and (\ref{AFC48}), 
take $(n_0,N_0)$ with $n_0/N_0 = r_1$ and $N_0$ large enough that  
\begin{equation} \label{FFor1} 
\left|\frac{N_0}{n_0} \frac{N_0}{n_0-1} \log(d_{(n_0,N_0)})
                    - \frac{1}{r_1^2} \log(f(r_1))\right| 
              \ < \ \frac{\varepsilon}{6} \ , 
\end{equation} 
\begin{equation} \label{FFor2} 
2\left|\frac{N_0}{n_0-1} - \frac{1}{r_1}\right| \log(C) 
          \ < \ \frac{\varepsilon}{6} \ .
\end{equation}

Consider a pair $(n,N)$ with $r_1 < n/N < r_2$ and $N \ge N_0$.
First, replacing $(m,M)$ by $(n,N)$ and $r$ by $r_1$ in (\ref{AFC41}), we have 
\begin{equation} \label{FFor3} 
\frac{N}{n}\frac{N}{n-1} \log(d_{(n,N)}) 
 \ \ge \ \frac{1}{r_1^2} \log(f(r_1))  
  - 2 \left|\frac{N}{n-1}-\frac{1}{r_1}\right| \log(C) \ .
\end{equation} 
Note that since $r_1 < n/N < r_2$ and $N \ge N_0$, also 
\begin{equation*}
\frac{N_0}{n_0-1} \ \ge \ \frac{N}{n-1}  \ > \ \frac{1}{r_2} \ .
\end{equation*}
Hence (\ref{AFC48}) and (\ref{FFor2}) give
\begin{eqnarray}  
2 \left| \frac{N}{n-1} - \frac{1}{r_1} \right|\log(C) 
 & < & \left(2 \left| \frac{N}{n-1} - \frac{1}{r_2} \right| 
            + 2 \left| \frac{1}{r_1} - \frac{1}{r_2}\right| \right)\log(C) \notag \\
 & < & \left( 2 \left| \frac{N_0}{n_0-1} - \frac{1}{r_2} \right|   
            + 2 \left| \frac{1}{r_1} - \frac{1}{r_2}\right| \right)\log(C) \notag \\
 & < & \left( 2 \left| \frac{N_0}{n_0-1} - \frac{1}{r_1} \right| 
            + 2 \left| \frac{1}{r_1} - \frac{1}{r_2}\right|    
            + 2 \left| \frac{1}{r_1} - \frac{1}{r_2}\right| \right)\log(C)            
 \ < \ \frac{3\varepsilon}{6} \ . \label{FFor4} 
\end{eqnarray} 
Combining (\ref{FFor3}), (\ref{AFC47}) and  (\ref{FFor4}) gives 
\begin{equation} \label{FOut1}
\frac{N}{n}\frac{N}{n-1} \log(d_{(n,N)}) 
 \ \ge \ \frac{1}{r^2} \log(f(r)) - \frac{4 \varepsilon}{6} \ . 
\end{equation} 

For the opposite inequality, replace $(m,M)$ by $(n_0,N_0)$ in (\ref{AFC40}), 
and take $\lambda = 1$ (which is permissible since $n \ge n_0$ 
under our hypotheses).  This yields 
\begin{equation} 
\frac{N}{n} \frac{N}{n-1} \log(d_{(n,N)}) \ \le \ 
     \frac{N_0}{n_0} \frac{N_0}{n_0-1} \log(d_{(n_0,N_0)})
     + 2 \left| \frac{N}{n-1} - \frac{N_0}{n_0-1} \right| \log(C) \ . \label{FFor5}
\end{equation} 
Combining (\ref{FFor5}), (\ref{AFC47}), (\ref{FFor1}), 
(\ref{FFor2}), and (\ref{FFor4}) gives 
\begin{equation} \label{FFor6}
\frac{N}{n} \frac{N}{n-1} \log(d_{(n,N)})
 \ < \ \frac{1}{r^2} \log(f(r)) + \varepsilon \ .
\end{equation}

In the case where $f(r) \equiv 0$, we need
only to show that 
\begin{equation} \label{AFC55}
\lim_{\substack{ n/N \rightarrow r \\ N \rightarrow \infty }}
      d_{(n,N)}^{ \frac{N}{n} \cdot \frac{N}{n - 1} } \ = \ 0 \ .  
\end{equation} 
If $H$ is finite, then 
$d_{(n,N)}  = 0$ whenever $n > \#(H)$, and hence (\ref{AFC55}) holds
trivially.  If $H$ is infinite, then each $d_{(n,N)}  > 0$.   
Fix $\varepsilon > 0$ and take $r_1,r_2 \in \QQ$ with $0 < r_1 < r < r_2$, 
such that 
\begin{equation} \label{AFC56}
C^{2 | \frac{1}{r_1} - \frac{1}{r_2} | } \ \le \ 2 \ . 
\end{equation} 
Take $(n_0,N_0)$ with $n_0/N_0 = r_1$.  Since
$\lim_{\lambda  \rightarrow \infty} d_{(\lambda n_0,\lambda N_0)}  = 0$ 
we can assume $N_0$ is large enough that 
\begin{equation*}
d_{(n_0,N_0)}^{\frac{N_0}{n_0} \cdot \frac{N_0}{n_0 - 1}} 
        \ < \ \varepsilon/3 \ .
\end{equation*}
In view of (\ref{AFC56}) we can also assume that $N_0$ is large enough 
that for all $(n,N)$ with  $r_1 < \frac{n}{N} < r_2$ and $N \ge N_0$, then
\begin{equation*}
C^{2 \left( \frac{N_0}{n_0 - 1} - \frac{N}{n - 1} \right)} \ \le \ 3 \ . 
\end{equation*} 
Taking $\lambda  = 1$ and $(m,M) = (n_0,N_0)$ in (\ref{AFC40}), for all such $(n,N)$ we have 
\begin{equation*}   
d_{(n,N)}^{\frac{N}{n} \cdot \frac{N}{n - 1}} \ < \ \varepsilon \ .
\end{equation*}

Finally, let us show that $f(r) \equiv 0$ if and only if $H$ has capacity $0$.
\index{capacity $= 0$}
If instead of $[z,w]_{\fX,\vs}$ and $W(z)$ we had used another distance function 
\begin{equation*}
[z,w]_{\fX,\vs}^* \ = \ \prod_{i = 1}^{m^*} ([z,w]_{x_i^*})^{s_i^*}
\end{equation*}
and another continuous, positive background function $W^*(z)$, then 
there is a constant $A > 0 $ such that 
\begin{eqnarray*}
& & 1/A \ < \ [z,w]_{\fX,\vs}^*/[z,w]_{\fX,\vs} \ < \ A \ , \\
& & \ \ \ 1/A \ < \ W^*(z)/W(z) \ < \ A \ ,
\end{eqnarray*}
for all $z \ne w \in H$.  Let $d^*_r(H,W^*)$ denote the 
weighted transfinite diameter computed relative to $[z,w]_{\fX,\vs}^*$ 
\index{transfinite diameter!weighted $(\fX,\vs)$}
and $W^*(z)$, and write $d^*_{(n,N)}(H,W^*)$ for the finite terms in its definition.   
Taking points $\alpha_1, \ldots, \alpha_{n} \in H$ 
which realize the maximum in (\ref{AFC34}), we see that 
\begin{equation*}
d^*_{(n,N)}(H,W^*) \ \ge \ d_{(n,N)}(H,W)  
 \cdot A^{- \left( \frac{n}{N} \frac{n - 1}{N} + \frac{2n}{N} \right) } . 
\end{equation*}
Passing to a limit as $N \rightarrow \infty$ and $n/N \rightarrow r$ gives 
\begin{equation*}
d^*_r(H,W^*) \ \ge \ d_{r}(H,W)  \cdot A^{-(r^2 + 2r)} 
\end{equation*}
A similar inequality holds with  
$d^*_r(H,W^*)$ and $d_{r}(H,W)$ reversed.  
Thus  $d_{r}(H,W) \ne 0$ if and only if $d^*_r(H,W^*) \ne 0$.  

Taking $r = 1$, $W^*(z) = 1$, and 
$[z,w]_{\fX,\vs}^* = [z,w]_{\zeta}$ 
for a fixed $\zeta \in \cC_v(\CC) \backslash H$, 
by (\cite{RR1}, Theorem 3.1.18) we have 
$d^*_1(H,W^*) = \gamma_{\zeta}(H)$.  Thus $f(r) \not \equiv 0$ if and only
if for some (hence any) $\zeta \notin H$, $\gamma_{\zeta}(H) > 0$.  
\end{proof}


\section{ Comparisons }   \label{CapacityComparisonSection}

In this section we will compare the weighted Chebyshev constant, 
\index{Chebyshev constant!weighted $(\fX,\vs)$}
the weighted transfinite diameter, and the weighted capacity 
\index{transfinite diameter!weighted $(\fX,\vs)$}
\index{capacity!weighted $(\fX,\vs)$}
(for fixed $(\fX,\vs)$). 

Fix a compact set $H \subset \cC_v(\CC) \backslash \fX$, 
and let $W(z) : \cC_v(\CC) \rightarrow [0,\infty]$ be a function which is 
positive and continuous on a neigbhorhood of $H$.
Put $\widehat{u}(z) = -\log(W(z))$.  

\medskip
{\bf The Weighted Capacity and the Weighted Transfinite Diameter.}
\index{capacity!weighted $(\fX,\vs)$}
We first prove an inequality between the weighted Robin constant
\index{Robin constant!weighted $(\fX,\vs)$}
and the weighted transfinite diameter. 
\index{transfinite diameter!weighted $(\fX,\vs)$} 

\begin{proposition} \label{APD4}
For each $\sigma > 0$,
\begin{equation}
V_{\sigma}(H,W) \ \le \ -\log(d_{\sigma}(H,W))  \ .  \label{AFD57}
\end{equation}
\index{Robin constant!weighted $(\fX,\vs)$}
\end{proposition}

Before giving the proof, we will need a lemma.  
Let $U$ be a neighborhood of $H$, bounded away from $\fX$, 
on which $\widehat{u}(z)$ is continuous and bounded.  
  
   For any set $F$ with $F \subset U$, define the inner weighted 
Robin constant and inner weighted capacity  by 
\index{Robin constant!weighted $(\fX,\vs)$}
\index{capacity!weighted $(\fX,\vs)$}
\begin{eqnarray*}
\Vbar_{\sigma}(F,W)
    &=& \inf_{\substack{ K \subset F \\ \text{$K$ compact} }} 
          V_{\sigma}(K,W) \ , \\
\gammabar_{\sigma}(F,W) &=& 
      \sup_{\substack{ K \subset F \\ \text{$K$ compact} }} 
                          \gamma_{\sigma}(K,W) \ \, = \ \,  \exp(-\Vbar_{\sigma}(F,W)) \ .
\end{eqnarray*}
\index{Robin constant!weighted $(\fX,\vs)$}
Note that for any compact $K \subset U$, we have $\Vbar_{\sigma}(K,W) = V_{\sigma}(K,W)$, and that   
if $F_1 \subset F_2 \subset U$, then trivially 
$\Vbar_{\sigma}(F_1,W) \ge \Vbar_{\sigma}(F_2,W)$.  

\begin{lemma} \label{ALD5}
Let $H \subset U$ and $W(z)$ be as above.  For any $\sigma > 0$, 
and any $\varepsilon > 0$,  there is a neighborhood $\widetilde{U}$ of $H$ contained
in $U$ such that 
\begin{equation*}
V_{\sigma}(H,W) \ \ge \ \Vbar_{\sigma}(\widetilde{U},W) \ \ge \ V_{\sigma}(H,W) - \varepsilon \ .
\end{equation*}
\index{Robin constant!weighted $(\fX,\vs)$}
\end{lemma}

\begin{proof}
Take a sequence of open sets 
$ U_1 \supset U_2 \supset \ldots \supset H$ whose
closures are compact and satisfy $\overline{U}_i \subset U$, with  
$\bigcap_{k = 1}^{\infty} \overline{U}_k = H$.   
We claim that 
\begin{equation} \label{TCXM1}
\lim_{k \rightarrow \infty} V_{\sigma}(\overline{U}_k,W)\ = \ V_{\sigma}(H,W) \ .
\index{Robin constant!weighted $(\fX,\vs)$}
\end{equation}
To see this, for each $k$ 
let $\mu_k$ be an equilibrium distribution for 
\index{equilibrium distribution!$(\fX,\vs)$}
$\overline{U}_k$ and $W$ with mass $\sigma$, so $I_{\sigma}(\mu_k,W)) = V_{\sigma}(\overline{U}_k,W)$.
After passing to a subsequence if necessary, we can assume that the $\mu_k$ 
converge weakly to a measure $\mu_{\infty}$. 
Clearly $\mu_{\infty}$ is positive, supported on $H$, and has total mass $\sigma$.  
Then 
\begin{eqnarray}
 V_{\sigma}(H,W) & \le & I_{\sigma}(\mu_{\infty},W)  \notag \\ 
   & = & \ \lim_{t \rightarrow \infty}
         \left( \iint_{H \times H} 
              -\log^{(t)}([z,w]_{\fX,\vs}) \, d\mu_{\infty}(z) d\mu_{\infty}(w) 
                       + 2 \int_H \widehat{u}(z) d\mu_{\infty}(z) \right)  \notag \\
   & = & \lim_{t \rightarrow \infty} \lim_{k \rightarrow \infty} 
     \left( \iint_{\overline{U}_k \times \overline{U}_k} 
         -\log^{(t)}([z,w]_{\fX,\vs}) \, d\mu_k(z) d\mu_k(w) 
          + 2 \ \int_{\overline{U}_k} \widehat{u}(z) d\mu_k(z) \right)  \label{AFD58} \\ 
   & \le & \ \liminf_{k \rightarrow \infty} \lim_{t \rightarrow \infty} 
      \left( \iint_{\overline{U}_k \times \overline{U}_k} 
      -\log^{(t)}([z,w]_{\fX,\vs}) \, d\mu_k(z) d\mu_k(w) 
          + 2 \int_{\overline{U}_k} \widehat{u}(z) d\mu_k(z) \right)  \label{AFD59} \\
   & = &   \liminf_{k \rightarrow \infty} \ I_{\sigma}(\mu_k,W) 
         \ = \ \liminf_{k \rightarrow \infty} \ V_{\sigma,H}(\overline{U}_k,W) \ .
                  \label{AFD60}
\end{eqnarray}      
\index{Robin constant!weighted $(\fX,\vs)$}
\index{equilibrium distribution} 
The interchange of limits in (\ref{AFD58}), (\ref{AFD59}) 
is valid because the kernels $-\log^{(t)}([z,w]_{\fX,\vs})$ are increasing with $t$.

     If $H$ has capacity $0$, then $V_{\sigma}(H,W) = \infty$, so 
\index{Robin constant!weighted $(\fX,\vs)$}
\index{capacity $= 0$}
$\lim_{k \rightarrow \infty} V_{\sigma}(\overline{U}_k,W) = \infty$.  If 
$H$ has positive capacity, then $V_{\sigma}(H,W)$ is finite.
\index{capacity $> 0$}  
Since $V_{\sigma}(H,W) \ge V_{\sigma}(\overline{U}_k,W)$ for all $k$, (\ref{AFD60}) gives
\index{Robin constant!weighted $(\fX,\vs)$}
$\lim_{k \rightarrow \infty} V_{\sigma}(\overline{U}_k,W) = V_{\sigma}(H,W)$.  
In either case (\ref{TCXM1}) holds, and we obtain the Lemma 
by taking $\widetilde{U} = U_k$ for sufficiently large $k$.  
\end{proof}  

\medskip
{\bf Proof of Proposition \ref{APD4}:  }  
By Theorem \ref{ATC3} and the remarks after (\ref{FNQ3A}),  
\begin{equation*}
V_{\sigma}(H,W)  =  \infty 
\quad \text{iff} \quad -\log(d_{\sigma}(H,W)) =  \infty 
\quad \text{iff} \quad \text{$H$ has capacity $0$,} 
\end{equation*} 
\index{capacity $= 0$}
\index{Robin constant!weighted $(\fX,\vs)$}
so we can assume that $V_{\sigma}(H,W)$ is finite, $d_{\sigma}(H,W) > 0$, 
and $H$ has positive capacity.
\index{capacity $> 0$}

Fix $\sigma$, and fix $\varepsilon > 0$.  Let $U$ be a
neighborhood of $H$ on which $\widehat{u}(z)$ is continuous and bounded.  
By Lemma \ref{ALD5},
there is a neighborhood $\widetilde{U}$ of $H$, whose closure is contained in $U$,
such that $\Vbar_{\sigma}(\widetilde{U},W) > V_{\sigma}(H,W) - \varepsilon$.  
\index{Robin constant!weighted $(\fX,\vs)$}
We can assume that $H$ is
covered by a finite number of local coordinate patches,  
each of which contained in $\widetilde{U}$.
Then there is an $R_0 > 0$ such that for each $x \in H$, 
there is a coordinate patch which contains the closed disc
$D(x,R_0) = \{ z : |z-w| \le R_0 \}$ relative to that coordinate.  
For each $x \in H$ fix such a coordinate patch, and given $0 < R \le R_0$,  
write $D(x,R)$ for the corresponding closed disc of radius $R$.  

To prove (\ref{AFD57}) we use a familiar construction involving
``smearing out'' Fekete measures.  
\index{Fekete!measure}
Take a sequence of pairs $(n_k,N_k)$ with
$n_k/N_k \rightarrow \sigma$  and  $N_k \rightarrow \infty$.
Given $k$, let $\alpha_1, \ldots, \alpha_{n_k} \in H$ be points where the
maximum in the definition of $d_{(n_k,N_k)}(H,W)$ is achieved; let 
\begin{equation} \label{AFD67}
\nu_k \ =  \ \sum_{i=1}^{n_k} \frac{1}{N_k} \delta_{\alpha_i}(z)
\end{equation}
be the associated measure of mass $n_k/N_k$ on $H$; 
it will be called a Fekete measure.  After passing to a
\index{Fekete!measure|ii}
subsequence, if necessary, we can assume that the $\nu_k$ converge
weakly to a measure $\nu$ on $H$.

Fix $k$, and write $(n,N) = (n_k,N_k)$,
$\nu_k = \sum_{i=1}^{n} \frac{1}{N} \delta_{\alpha_i}(z)$.
Without loss, we can assume that $k$ (hence $N$) is large enough that
$\frac{1}{\sqrt{\pi N}} \le R_0$.  Let $dm_i$ be the measure which coincides
with Lebesgue measure on the disc $D_i = D(\alpha_i,\frac{1}{\sqrt{\pi N}})$,
and is $0$ outside that disc.  Thus, $dm_i$ has total mass $\frac{1}{N}$.
Put $F_k = \bigcup_{i=1}^{n} D_i$, and let
\begin{equation*}
\widetilde{\nu}_k \ = \ \frac{r}{n/N} \sum_{i=1}^{n} dm_i  \ .
\end{equation*}
Then $F_k \subset \widetilde{U}$, and $\widetilde{\nu}_k$
is a positive measure of mass $\sigma$ on $F_k$. It follows that
\begin{eqnarray} \label{AFD61}
V_{\sigma}(H,W) - \varepsilon &\le& \Vbar_{\sigma}(\widetilde{U},W) 
            \ \le \ V_{\sigma}(F_k,W) \ \le \ I_{\sigma}(\widetilde{\nu}_k,W)  \\
&=& \frac{\sigma}{n/N} \sum_{i, j = 1}^{n}
    \iint_{D_i \times D_j} -\log([z,w]_{\fX,\vs}) \, dm_i(z) dm_j(w) \notag \\
& & \qquad \qquad
       + \frac{\sigma}{n/N} \cdot 2 \int_{F_k} \widehat{u}(z) \, d\nu_k(z) \ . \notag
\end{eqnarray}
\index{Robin constant!weighted $(\fX,\vs)$}
For each fixed $z$, the function $-\log([z,w]_{\fX,\vs})$ 
is superharmonic\index{superharmonic} in $w$.
Hence, using polar coordinates in $D_j$, we have
\begin{eqnarray*}
\int_{D_j} -\log([z,w]_{\fX,\vs}) \, dm_j(w)
     &=& \int_0^{\frac{1}{\sqrt{\pi N}}} \int_0^{2 \pi}
             -\log([z,\alpha_j + te^{i\theta} ]_{\fX,\vs}) \, d\theta \, t \, dt \\
     &\le& -\frac{1}{N} \log([z,\alpha_j]_{\fX,\vs}) \ .
\end{eqnarray*}
If $i \ne j$, then since $-\frac{1}{N} \log([z,\alpha_j]_{\fX,\vs})$
is superharmonic\index{superharmonic} in $w$, this gives
\begin{eqnarray}
\iint_{D_i \times D_j} -\log([z,w]_{\fX,\vs}) \, dm_i(z) dm_j(w)
     &\le& \int_{D_i} -\frac{1}{N} \log([z,\alpha_j]_{\fX,\vs}) \, dm_i(z) \notag \\
     &\le& -\frac{1}{N^2} \log([\alpha_i,\alpha_j]_{\fX,\vs}) \ . \label{AFD62}
\end{eqnarray}
If $i = j$, then since $-\log([z,w]_{\fX,\vs}) = -\log(|z-w|) + \eta(z,w)$
for a $\cC^{\infty}$ function $\eta(z,w)$, we see that
\begin{eqnarray}
& &\iint_{D_i \times D_i}  -\log([z,w]_{\fX,\vs}) dm_i(z) dm_i(w)                                
     \ \le\ \int_{D_i} -\frac{1}{N} \log([z,\alpha_i]_{\fX,\vs}) \, dm_i(z)  \notag \\
& & \quad = \ \frac{1}{N} \left( \int_0^{\frac{1}{\sqrt{\pi N}}} \int_0^{2 \pi}
    -\log(t) \, d\theta \, t\,  dt + \int_{D_i} \eta(z,\alpha_i) \, dm_i(z) \right) \notag \\
& & \quad = \ O(\frac{1}{N^2} \log(N)) \ . \label{AFD63}
\end{eqnarray}
Thus
\begin{eqnarray}
 I_{\fX,\vs}(\widetilde{\nu}_k)&=& \sum_{i, j = 1}^{n}
    \iint_{D_i \times D_j} -\log([z,w]_{\fX,\vs}) \, dm_i(z) dm_j(w)   \notag \\
 &\le& \frac{1}{N^2} \sum^{n}_{\substack{ i,j = 1 \\ i \ne j }}
    -\log([\alpha_i,\alpha_j]_{\fX,\vs}) + O(\frac{\log(N)}{N})  \ . \label{AFD64}
\end{eqnarray}

On the other hand, $W(z)$, and hence $\widehat{u}(z)$,
is continuous on the closure of $\widetilde{U}$.  
Hence as $k \rightarrow \infty$
\begin{equation} \label{AFD65}
\left| \int \widehat{u}(z) d\widetilde{\nu}_k(z) - \int \widehat{u}(z) d\nu_k \right|
  \ = \ o(1) \ .
\end{equation}
Inserting (\ref{AFD64}) and (\ref{AFD65}) in (\ref{AFD61}), we obtain
\begin{equation} \label{AFD66}
V_{\sigma}(H,W) - \varepsilon 
     \ \le \ -\log\left(d_{(n_k,N_k)}(H,W)\right) + o(1) \ .
\end{equation}
\index{Robin constant!weighted $(\fX,\vs)$}
Passing to the limit as $k \rightarrow \infty$, and
using that $\varepsilon > 0$ is arbitrary gives 
\begin{equation*}
V_{\sigma}(H,W) \ \le \ -\log(d_{\sigma}(H,W)) \ . 
\end{equation*}
\index{Robin constant!weighted $(\fX,\vs)$}
\hfill $\blacksquare$

\medskip
\noindent{\bf The Weighted Transfinite Diameter and Chebyshev Constant.}
\index{Chebyshev constant!weighted}

Fix $\sigma > 0$, and consider a sequence of pairs of positive integers 
$(n_k,N_k)$ with $N_k \rightarrow \infty$ and $n_k /N_k \rightarrow \sigma$.  
For each  $k$, let $\nu_k$ be the associated Fekete measure, 
\index{Fekete!measure}
as in (\ref{AFD67}).

\begin{proposition}  \label{APD6}
With the notation above, if $\nu$ is any weak limit 
of the Fekete measures $\nu_k$, then
\index{Fekete!measure}
\begin{equation}
-\log(d_{\sigma}(H,W)) 
 \ \le \ \sigma \cdot (-\log(\CH^*_{\sigma}(H,W))) + \int_{H} \widehat{u}(z) d\nu(z) \ .  \label{AFD68}
\end{equation}
\end{proposition}

\begin{proof}
If $H$ has capacity $0$, 
\index{capacity $= 0$}
both $d_{\sigma}(H,W)$ and $\CH^*_{\sigma}(H,W)$ are $0$, and the 
inequality is trivial.  Hence we can assume $H$ has positive capacity,
\index{capacity $> 0$}
so $d_{\sigma}(H,W)$ and $\CH^*_{\sigma}(H,W)$ are positive. 

Fix $k$ and write $(n, N) = (n_k,N_k)$.  As before, let
$\alpha_1, \ldots, \alpha_{n} \in H$ be points where $d_{(n,N)}(H,W)$ is achieved.  
Fixing $\alpha_2, \ldots \alpha_{n}$, consider
\begin{equation} \label{AFD69}
F(z) \ := \ \left( \prod_{i = 2}^{n} [z,\alpha_i]_{\fX,\vs} \cdot W(z)^N \right)^2
   \cdot \prod_{\substack{ 2 \le i, j \le n \\ i \ne j }} [\alpha_i,\alpha_j]_{\fX,\vs} 
             \cdot \prod_{i=2}^{n} W(\alpha_i)^{2N} \ . 
\end{equation}
By hypothesis $F(z)$ takes its maximum value on $H$ when $z = \alpha_1$.  
Thus the pseudopolynomial 
\begin{equation*}
P(z)\ = \ \prod_{i = 2}^{n} [z,\alpha_i]_{\fX,\vs} \cdot W(z)^N
\end{equation*}
\index{pseudopolynomial!weighted $(\fX,\vs)$}
achieves its maximum on $H$ at $z = \alpha_1$.  By definition, if
\begin{equation*}
\tP_{(n-1, N)}(z) \ = \ \prod_{i = 2}^{n} [z,\beta_i]_{\fX,\vs} W(z)^N
\end{equation*}
is the pseudopolynomial with minimal sup norm on $H$, then 
\index{pseudopolynomial!weighted $(\fX,\vs)$}
\begin{equation*}
(\CH^*_{(n-1,N)})^N \ := \ \|\tP_{(n -1, N)}\|_{H} \ \le \ P(\alpha_1) \ .   
\end{equation*}
Thus $(\CH^*_{n -1,N})^N \le  \prod_{i = 2}^{n} [\alpha_1,\alpha_i]_{\fX,\vs} \cdot W(\alpha_1)^N.$ 
The same argument applies for each $\alpha_k$, so 
\begin{eqnarray*}
d_{(n,N)}(H,W)^{N^2} 
     &=& \prod_{i = 1}^{n} 
          \left( \prod^{n}_{\substack{ j = 1 \\ j \ne i }} [\alpha_i,\alpha_j]_{\fX,\vs} 
                    \cdot W(\alpha_i)^N \right) \cdot 
          \left( \prod_{i = 1}^{n} W(\alpha_i)^N \right) \\
     &\ge& (\CH^*_{(n - 1,N)})^{Nn} \cdot 
                    \prod_{i = 1}^{n} W(\alpha_i)^N \ .
\end{eqnarray*}
Taking logarithms, dividing by $N^2$,
and recalling that all quantities involved depend on $k$, we get 
\begin{eqnarray} 
- \log(d_{(n_k,N_k)}(H,W))
    \ \le \ \frac{n_k}{N_k}  
           \left( - \log(\CH^*_{(n_k -1, N_k)}) \right) 
   +  \sum_{i = 1}^{n_k} \widehat{u}(\alpha_i) \, \frac{1}{N_k} \ . \notag \\
                                                                 \label{AFD70} 
\end{eqnarray}

Now let $ k\rightarrow \infty$.  
The fractions $(n_k -1)/N_k$ converge to $\sigma$, 
so by Theorem \ref{ATB1} the weighted Chebyshev constant satisfies
\index{Chebyshev constant!weighted} 
\begin{equation*}
\CH^*_{\sigma}(H,W) \ = \ \lim_{k \rightarrow \infty} \CH^*_{(n_k -1,N_k)}
\end{equation*}
After passing to a subsequence, if necessary, we can assume that the measures 
$\nu_k = \sum_{i = 1}^{n_k} \frac{1}{N_k} \delta_{\alpha_i}(z)$
converge weakly to $\nu$.  Hence by (\ref{AFD70}) and Theorem \ref{ATC3}, 
\begin{eqnarray}\label{AFD71}
-\log(d_{\sigma}(H,W)) &=& \lim_{k \rightarrow \infty}
            -\log(d_{(n_k,N_k)}(H,W))  \\
  & \le & \sigma \cdot (-\log(\CH^*_{\sigma}(H,W)) + \int_{H} \widehat{u}(z) d\nu(z) \ . \notag
\end{eqnarray}
\end{proof}
 

\section{ Particular cases of interest}
   \label{Weight0CaseSection}  

In this section we consider the `classical' case where the
weight function $W(z)$ is trivial, and specialize to the situation of interest for our applications.  
As in the previous sections, we assume $(\fX,\vs)$ has been fixed.
 
\begin{theorem} \label{ATE8}  
Let $H \subset \cC_v(\CC) \backslash \fX$ be compact. 
Assume that $\sigma = 1$ and $W(z) \equiv 1$.
Then $\gamma_1(H,W) = \gamma_{\fX,\vs}(H)$, 
$\CH^*_1(H,W) = \CH^*_{\fX,\vs}(H)$, and $d_1(H,W) = d_{\fX,\vs}(H)$, so  
\begin{equation*}
\gamma_1(H,W)  \ = \ d_1(H,W) \ = \ \CH^*_1(H,W) \ = \ \gamma_{\fX,\vs}(H) \ .
\end{equation*}
Furthermore, if $H$ has positive capacity, then the equilibrium distribution $\mu$ for $H$
\index{capacity $> 0$}
relative to $W(z)$ is unique and coincides with $\mu_{\fX,\vs}$, and the constant
$\cV_{\mu}$ in Theorem $\ref{ATA1}$ coincides with $V_{\fX,\vs}(H)$.  
\index{Robin constant!$(\fX,\vs)$}
\end{theorem}

\begin{proof}  For the first assertion, it suffices to note that the 
optimization problems (\ref{FMu1}), (\ref{FVF1}) defining 
$V_1(H,W)$ and $V_{\fX,\vs}(H)$ are the same; 
\index{Robin constant!weighted $(\fX,\vs)$}
\index{Robin constant!$(\fX,\vs)$} 
for each $N$, the optimization problems 
(\ref{FMe1}), (\ref{FCH0}) defining $\CH^*_{(N,N)}(H,W)$ and $\CH^*_N(H)$ are the same;
and for each $N$, 
the optimization problems (\ref{AFC34}), (\ref{FTD1}) defining $d_{(N,N)}(H,W)$ 
and $d_N(H)$ are the same.

The second assertion follows from Theorem \ref{ATE10A}. 
The final assertions follow from the fact that $V_1(H,W) = V_{\fX,\vs}(H)$,
\index{Robin constant!weighted $(\fX,\vs)$}
\index{Robin constant!$(\fX,\vs)$}
and the uniqueness of $\mu_{\fX,\vs}$ in Theorem \ref{ATE10B}.  
\end{proof} 

\vskip .1 in 

Consider a sequence of pairs $(n_k,N_k)$ with 
$n_k/N_k \rightarrow 1$ and $N_k \rightarrow \infty$.  
For a given $k$, 
let $\alpha_1, \ldots, \alpha_{n_k}$ achieve the maximum in the definition
of $d_{(n_k,N_k)}(H,W)$ for $W(z) \equiv 1$.  Let   
\begin{equation*}
\nu_k \ = \ \sum_{i=1}^{n_k} \frac{1}{N_k} \delta_{\alpha_i}(z)
\end{equation*}
be the associated Fekete measure of mass $n_k/N_k$ on $H$.
\index{Fekete!measure} 

We will now see that, just as in the classical case, 
the Fekete measures converge weakly to the equilibrium 
\index{Fekete!measure} 
distribution $\mu_{\fX,\vs}$ of $H$. 
\index{equilibrium distribution!$(\fX,\vs)$} 

\begin{corollary} \label{ACE14} 
Let $H \subset \cC_v(\CC) \backslash \fX$ be compact with positive capacity,
\index{capacity $> 0$}
and suppose $W(z) \equiv 1$.
Then for any sequence of pairs $(n_k,N_k)$ 
with $n_k/N_k \rightarrow 1$ and $N_k \rightarrow \infty$,
the associated sequence of Fekete measures 
\index{Fekete!measure}  
$\{\nu_k\}$ converges weakly to $\mu_{\fX,\vs}$.  
\end{corollary}

\begin{proof}  Let $\nu_{\infty}$ be a weak limit of a subsequence of the
$\nu_k$.  Upon passing to that subsequence, we see that 
\begin{equation*}
V_{\fX,\vs}(H) \ = \ \lim_{k \rightarrow \infty} 
                               - \log(d_{(n_k,N_k)}(H,W)) \ .
\end{equation*} 
\index{Robin constant!$(\fX,\vs)$}
Let $\widetilde{\nu}_k$ be a sequence of measures `smearing out' the $\nu_k$ 
as in the proof of Proposition \ref{APD4};  
then the $\widetilde{\nu}_k$ also converge weakly to $\nu_{\infty}$.  
By the argument leading to (\ref{AFD66}), as $k \rightarrow \infty$,   
\begin{equation*}
I_{\fX,\vs}(\widetilde{\nu}_k) \ = \ I_{n_k/N_k}(\widetilde{\nu}_k,W) 
\ \le \ - \log(d_{(n_k,N_k)}(H,W)) + o(1) \ . 
\end{equation*}                            
On the other hand, by (\ref{AFD61}), for any $\varepsilon > 0$ and 
sufficently large $k$  
\begin{equation*}
V_{\fX,\vs}(H) - \varepsilon \ \le \ I_{\fX,\vs}(\widetilde{\nu_k}) \ .
\end{equation*}
\index{Robin constant!$(\fX,\vs)$}
Thus, by Theorem \ref{ATE8}
\begin{eqnarray*}
V_{\fX,\vs}(H) 
   &\le& I_{\fX,\vs}(\nu_{\infty}) 
         \ \le \ \liminf_{k \rightarrow \infty} I_{n_k/N_k}(\widetilde{\nu}_k,H)  \\ 
   &\le& \liminf_{k \rightarrow \infty} -\log(d_{(n_k,N_k)}(H,W)) 
          \ = \ V_{\fX,\vs}(H) \ .
\end{eqnarray*}
\index{Robin constant!$(\fX,\vs)$}
Since $\mu_{\fX,\vs}$ is the unique probability measure on $H$ 
which minimizes the energy integral (Theorem \ref{ATE10B}),\index{energy integral!$(\fX,\vs)$}  
it follows that $\nu_{\infty} = \mu_{\fX,\vs}$.
\end{proof}


\medskip

We now come to the case of interest for our application.

Suppose $H \subset \cC_v(\CC) \backslash \fX$ is compact, 
and $H_\ell$ is a connected component of $H$.
We will assume that $H_\ell$ has positive capacity.  
\index{capacity $> 0$}
Let $\mu = \mu_{\fX,\vs}$ be the equilibrium distribution of $H$.  Put
\index{equilibrium distribution!$(\fX,\vs)$}
\begin{eqnarray}
\sigma_{\ell} & = & \mu_{\fX,\vs}(H_\ell) \ , \notag \\
\widehat{u}_\ell(z) &=& \int_{H \backslash H_\ell} -\log([z,w]_{\fX,\vs}) \, d\mu_{\fX,\vs}(w) 
                        \ , \label{Fu1} \\
W_\ell(z) &=& \exp(-\widehat{u}_\ell(z)) \ . \notag
\end{eqnarray}
(We take $\widehat{u}_\ell(z) \equiv 0$ and $W_{\ell}(z) \equiv 1$ if $H \backslash H_{\ell}$ is empty.)
We will use the results of previous sections to study weighted potential 
\index{potential theory!weighted}
theory for $H_\ell$, $\sigma_{\ell}$, and $W_\ell(z)$.  The problem is a somewhat subtle 
one:  to show that in the weighted case for $H_\ell$, with the weight $W_\ell(z)$ coming
from the unweighted case for $H \backslash H_\ell$, the extremal objects for $H_\ell$ are the 
restrictions of the global unweighted objects for $H$.  

The key results are Theorems \ref{APF15}, \ref{ATF18} and \ref{ATF19}.

\begin{theorem}  \label{APF15} 
Let $H$, $H_\ell$, $\sigma_\ell = \mu_{\fX,\vs}(H_\ell)$
and $W_\ell(z)$ be as in $(\ref{Fu1})$.  
Then the equilibrium distribution $\mu_{\sigma_{\ell},H_\ell,W_\ell}$ 
\index{equilibrium distribution!weighted $(\fX,\vs)$}
of $H_\ell$ relative to $W_\ell(z)$ 
with mass $\sigma_\ell$ is unique, and is given by 
\begin{equation}
\mu_{\sigma_\ell,H_\ell,W_\ell} \ = \ \mu_{\fX,\vs} \vert_{H_\ell} \ .\label{AFF83}
\end{equation}
\end{theorem}
   
\begin{proof} Put $\mu_\ell = \mu_{\fX,\vs} \vert_{H_\ell}$, 
$\dot{\mu}_\ell = \mu_{\fX,\vs} \vert_{H \backslash H_\ell}$, and let $\widetilde{\mu}_\ell$ be any  
equilibrium distribution for $H_\ell$ relative to $W_\ell(z)$ with mass $\sigma_\ell$. 
\index{equilibrium distribution!weighted $(\fX,\vs)$} 
By definition, $\widetilde{\mu}_\ell$ minimizes the weighted energy 
\begin{equation*}
I_{\sigma_\ell}(\nu,W_\ell) \ = \ I_{\fX,\vs}(\nu) + 2 \int_H \widehat{u}_{\ell}(z) \, d\nu(z)  
\end{equation*}
among all positive measures of mass $\sigma_\ell$ on $H_\ell$,
so $I_{\sigma_\ell}(\widetilde{\mu}_\ell,W_\ell) = V_{\sigma_\ell}(H_\ell,W_\ell)$.  Thus, 
\index{Robin constant!weighted $(\fX,\vs)$}
\begin{equation}
I_{\sigma_\ell}(\mu_\ell,W_\ell) \ \ge \ I_{\sigma_\ell}(\widetilde{\mu}_\ell,W_\ell) \ .  \label{AFF84}
\end{equation}  

     On the other hand, by Theorem \ref{ATE10B}, 
$\mu := \mu_{\fX,\vs}$ is the unique probability measure minimizing 
\begin{equation*}
I_{\fX,\vs}(\mu) \ = \ \iint_{H \times H} -\log([z,w]_{\fX,\vs}) \, d\mu(z) d\mu(w) \ .
\end{equation*}
Thus, if we put $\widetilde{\mu} = \widetilde{\mu}_\ell + \dot{\mu}_\ell$, 
then $I_{\fX,\vs}(\widetilde{\mu}) \ge I_{\fX,\vs}(\mu)$.  Expanding this, we see that 
\begin{eqnarray*} 
I_{\fX,\vs}(\widetilde{\mu}_\ell) 
         + 2 \int_{H_\ell} \widehat{u}_\ell(z) d\widetilde{\mu}_\ell(z) + I_{\fX,\vs}(\dot{\mu}_\ell) 
             \qquad \qquad \qquad \qquad \qquad \qquad  \\
     \qquad \qquad \qquad \qquad \qquad
  \ge \ I_{\fX,\vs}(\mu_\ell) +  2 \int_{H_\ell} \widehat{u}_\ell(z) d\mu_\ell(z) + I_{\fX,\vs}(\dot{\mu}_\ell) \ ,
\end{eqnarray*}
which implies that $I_{\sigma_\ell}(\widetilde{\mu}_\ell,W_\ell) \ge I_{\sigma_\ell}(\mu_\ell,W_\ell)$.
Combining this with (\ref{AFF84}) gives $I_{\sigma_\ell}(\widetilde{\mu}_\ell,W_\ell) = I_{\sigma_\ell}(\mu_\ell,W_\ell)$.  
Consequently
\begin{equation*}
I_{\fX,\vs}(\widetilde{\mu}) \ = \ I_{\sigma_\ell}(\widetilde{\mu}_\ell,W_\ell) + I_{\fX,\vs}(\dot{\mu}_\ell)  
            \ = \ I_{\sigma_\ell}(\mu_\ell,W_\ell) + I_{\fX,\vs}(\dot{\mu}_\ell) \ = \ I_{\fX,\vs}(\mu) \ .
\end{equation*}
From the uniqueness of $\mu$, we conclude that 
$\widetilde{\mu} = \mu_{\fX,\vs}$, hence $\widetilde{\mu}_\ell = \mu_\ell$. 
\end{proof}             

\vskip .1 in
Given measures $\nu_1$, $\nu_2$, write
\begin{equation*}
I_{\fX,\vs}(\nu_1,\nu_2) \ = \ \iint -\log([z,w]_{\fX,\vs}) \, d\nu_1(z) d\nu_2(z) \ , 
\end{equation*}
provided the integral is defined.  Also, given a measure $\nu$, put 
\begin{equation*}
\widehat{I}_{\fX,\vs}(\nu) = \iint_{ H \times H \backslash \{\text{diagonal}\} } 
                          -\log([z,w]_{\fX,\vs}) \, d\nu(z) d\nu(w) \ ,
\end{equation*} 
when the integral is defined.  
Then $I_{\fX,\vs}(\nu_1,\nu_2)$ is symmetric, 
and is bilinear when all relevant terms are defined and finite.
If $\nu$ is a positive measure, then $I_{\fX,\vs}(\nu) = I_{\fX,\vs}(\nu,\nu)$,
and if in addition $\nu$ does not charge the diagonal, then $I_{\fX,\vs}(\nu) = \widehat{I}_{\fX,\vs}(\nu)$.

\medskip
We will now show that in the situation of (\ref{Fu1}),
the weighted transfinite diameter and the weighted capacity coincide.
\index{transfinite diameter!weighted $(\fX,\vs)$}
\index{capacity!weighted $(\fX,\vs)$} 
As in the proof of Proposition \ref{APF15}, put
\begin{equation}  \label{AFF85}
   \mu_\ell = \mu_{\fX,\vs} \vert_{H_\ell} \ , \qquad
   \dot{\mu}_\ell = \mu_{\fX,\vs} \vert_{H \backslash H_\ell} \ .
\end{equation} 
Note that $I_{\fX,\vs}(\mu_\ell,\dot{\mu}_\ell) = \int_{H_\ell} \widehat{u}_\ell(z) \, d\mu_\ell(z)$. 

\begin{theorem} \label{ATF16}
Let $H$, $H_\ell$, $\sigma_\ell = \mu_{\fX,\vs}(H_\ell)$,
and $W_\ell(z)$ be as in Theorem $\ref{APF15}$.  Then
\begin{equation} \label{AFF86}
-\log(d_{\sigma_\ell}(H_\ell,W_\ell)) \ = \ V_{\sigma_\ell}(H_\ell,W_\ell) 
      \ = \ \sigma_\ell \cdot V_{\fX,\vs}(H) + \int_{H_\ell} \widehat{u}_\ell(z) \, d\mu_\ell(z) \ .
\end{equation}
\index{Robin constant!weighted $(\fX,\vs)$}
\index{Robin constant!$(\fX,\vs)$}
\end{theorem}

\begin{proof}
The second equality in (\ref{AFF86}) is easy:  by Proposition \ref{APF15} 
\begin{eqnarray*}
V_{\sigma_\ell}(H_\ell,W_\ell) &=& I_{\sigma_\ell}(\mu_\ell,W_\ell) \\ 
             &=& \int_{H_\ell \times H_\ell} -\log([z,w]_{\fX,\vs}) \, d\mu_\ell(z) d\mu_\ell(w) 
                       + 2 \int_{H_\ell} \widehat{u}_\ell(z) \, d\mu_\ell(z) \\                  
             &=& I_{\fX,\vs}(\mu_\ell,\mu_\ell) + 2 I_{\fX,\vs}(\mu_\ell,\dot{\mu}_\ell)
                     \ =  \ I_{\fX,\vs}(\mu_{\fX,\vs},\mu_\ell) + I_{\fX,\vs}(\mu_\ell,\dot{\mu}_\ell) \\ 
             &=& \int_{H_\ell} u_{\fX,\vs}(z) \, d\mu_\ell(z) 
                     + I_{\fX,\vs}(\mu_\ell,\dot{\mu}_\ell) \\
             &=& \sigma_\ell \cdot V_{\fX,\vs}(H) 
                          + \int_{H_\ell} \widehat{u}_\ell(z) \, d\mu_\ell(z) \ .
\end{eqnarray*} 
\index{Robin constant!weighted $(\fX,\vs)$}
where the last inequality holds since $u_{\fX,\vs}(z)$ takes the constant value $V_{\fX,\vs}(H)$
\index{equilibrium potential!$(\fX,\vs)$}
\index{equilibrium potential!takes constant value a.e. on $E_v$}
\index{Robin constant!$(\fX,\vs)$}
on $H$, except possibly on a set of inner capacity $0$;  
\index{capacity $= 0$}
and the exceptional set necessarily has $\mu_\ell$-measure $0$.

We now turn to the first equality.  By Proposition \ref{APD4}, 
$V_{\sigma_\ell}(H_\ell,W_\ell) \le  -\log(d_{\sigma_\ell}(H_\ell,W_\ell))$ so we need only show the reverse
inequality.   
\index{Robin constant!weighted $(\fX,\vs)$}
For this, we use the fact that 
the unweighted Fekete measures for $H$ converge weakly to the unweighted equilibrium distribution of $H$,
\index{Fekete!measure} 
\index{equilibrium distribution!$(\fX,\vs)$}
as shown in Corollary \ref{ACE14}.  
Given an integer $N > 0$, let $\alpha_\ell^*, \ldots, \alpha_N^*$ be points 
maximizing the transfinite diameter $d_{(N,N)}(H,W)$ for 
$W(z) \equiv 1$.  Label them so that $\alpha_\ell^*, \ldots, \alpha_{n}^* \in H_\ell$  
and $\alpha_{n + 1}^*, \ldots, \alpha_N^* \in H \backslash H_\ell$, and put 
\begin{equation*}
D_{n,N}^{(1)} \ = \ \left( \prod^{n}_{\substack{ i, j = 1 \\ i \ne j }} 
             [\alpha_i^*,\alpha_j^*]_{\fX,\vs} \cdot 
           \prod_{i=1}^n W_\ell(\alpha_i^*)^{2N} \right)^{1/{N^2}} \ .
\end{equation*}                   
By the definition of the weighted transfinite diameter,
\index{transfinite diameter!weighted $(\fX,\vs)$}
$D_{n,N}^{(1)} \le  d_{(n,N)}(H_\ell,W_\ell)$.  Also put
\begin{equation*}
\nu^{(N)} \ = \ \sum_{i=1}^N \frac{1}{N} \delta_{\alpha_i^*}(z)  
\end{equation*}
and put 
$\nu_\ell^{(N)} = \nu^{(N)} \vert_{H_\ell}$, $\dot{\nu}_\ell^{(N)} = \nu^{(N)} \vert_{H \backslash H_\ell}$;
then 
\begin{equation*}
-\log(D_{n,N}^{(1)}) \ = \ \widehat{I}_{\fX,\vs}(\nu_\ell^{(N)}) + 2 I_{\fX,\vs}(\nu_\ell^{(N)},\dot{\nu}_\ell^{(N)}) \ .
\end{equation*} 

It follows easily by weak convergence that 
\begin{equation} \label{AFF88}
\lim_{N \rightarrow \infty} I_{\fX,\vs}(\nu_\ell^{(N)},\dot{\nu}_\ell^{(N)}) \ = \ I_{\fX,\vs}(\mu_\ell,\dot{\mu}_\ell) \ .
\end{equation}
We will show below that 
\begin{equation} \label{AFF87}
\lim_{N \rightarrow \infty} \widehat{I}_{\fX,\vs}(\nu_\ell^{(N)}) \ = \ I_{\fX,\vs}(\mu_\ell) \ .
\end{equation}
Granting (\ref{AFF87}), by Proposition \ref{APD4} we then have 
\begin{eqnarray*}
I_{\fX,\vs}(\mu_\ell) &+& 2 I_{\fX,\vs}(\mu_\ell,\dot{\mu}_\ell) 
                \ = \  V_{\sigma_\ell}(H_\ell,W_\ell) \ \le \ -\log(d_{\sigma_\ell}(H_\ell,W_\ell)) \\
    &=& \lim_{N \rightarrow \infty} -\log(d_{(n,N)}(H_\ell,W_\ell)) 
       \ \le \ \limsup_{N \rightarrow \infty} -\log(D_N^{(1)}) \\
    &=& \limsup_{N \rightarrow \infty} \widehat{I}_{\fX,\vs}(\nu_\ell^{(N)})
                  + 2 I_{\fX,\vs}(\nu_\ell^{(n)},\dot{\nu}_\ell^{(N)})
        \ = \ I_{\fX,\vs}(\mu_\ell) + 2 I_{\fX,\vs}(\mu_\ell,\dot{\mu}_\ell) \ ,
\end{eqnarray*}
\index{Robin constant!weighted $(\fX,\vs)$}
so equalities hold throughout, yielding the theorem.  

To prove (\ref{AFF87}), note that when $W(z) \equiv 1$, by the definition
of $d_{(N,N)}(H,W)$ we have that for each $N$,  
\begin{equation*}
-\log(d_{(N,N)}(H,W)) \ = \ 
 \widehat{I}_{\fX,\vs}(\nu_\ell^{(N)}) + 2 I_{\fX,\vs}(\nu_\ell^{(N)},\dot{\nu}_\ell^{(N)}) + \widehat{I}_{\fX,\vs}(\dot{\nu}_\ell^{(N)}) \ .
\end{equation*}  
Passing to the limit as $N \rightarrow \infty$ and using Theorem \ref{ATE8} and 
Corollary \ref{ACE14}, we have 
\begin{eqnarray*}
 I_{\fX,\vs}(\mu_\ell) &+& 2 I_{\fX,\vs}(\mu_\ell,\dot{\mu}_\ell) \ + \ I_{\fX,\vs}(\dot{\mu}_\ell) 
                     \ = \ I_{\fX,\vs}(\mu) \ = \ V_{\fX,\vs}(H)  \\
  &=& \lim_{N \rightarrow \infty} -\log(d_{(N,N)}(H,W)) \\
  &=& \lim_{N \rightarrow \infty}  
     \left( \widehat{I}_{\fX,\vs}(\nu_\ell^{(N)}) + 2 I_{\fX,\vs}(\nu_\ell^{(N)},\dot{\nu}_\ell^{(N)}) 
                  + \widehat{I}_{\fX,\vs}(\nu_\ell^{(N)}) \right) \\
  &=& 2 I_{\fX,\vs}(\mu_\ell,\dot{\mu}_\ell) + \lim_{N \rightarrow \infty} 
            \left(\widehat{I}_{\fX,\vs}(\nu_\ell^{(N)}) + \widehat{I}_{\fX,\vs}(\nu_\ell^{(N)}) \right)
\end{eqnarray*}
\index{Robin constant!$(\fX,\vs)$}
and hence
\begin{equation} \label{AFF89}
 I_{\fX,\vs}(\mu_\ell) +  I_{\fX,\vs}(\dot{\mu}_\ell) \ = \ \lim_{N \rightarrow \infty} 
            \left( \widehat{I}_{\fX,\vs}(\nu_\ell^{(N)}) + \widehat{I}_{\fX,\vs}(\dot{\nu}_\ell^{(N)}) \right) \ .
\end{equation}
On the other hand, consider the `smearing out' $\widetilde{\nu}_\ell^{(N)}$ of 
$\nu_\ell^{(N)}$ as in Proposition \ref{APD4}.  
The same argument which gave (\ref{AFD64}) gives
\begin{equation*} 
I_{\fX,\vs}(\widetilde{\nu}_\ell^{(N)}) \ \le \ \widehat{I}_{\fX,\vs}(\nu_\ell^{(N)}) 
             + O \left( \frac{\log(N)}{N} \right) \ .
\end{equation*}
Since the $\nu_\ell^{(N)}$ convege weakly to $\mu_\ell$,
so do the $\widetilde{\nu}_\ell^{(N)}$, and therefore
\begin{equation} \label{AFF90}
I_{\fX,\vs}(\mu_\ell) \ = \ \lim_{N \rightarrow \infty} I_{\fX,\vs}(\widetilde{\nu}_\ell^{(N)})  
     \ \le \ \liminf_{N \rightarrow \infty} \widehat{I}_{\fX,\vs}(\nu_\ell^{(N)}) \ .
\end{equation}      
Symmetrically,                       
\begin{eqnarray} \label{AFF91}
I_{\fX,\vs}(\dot{\mu}_\ell) &\le& \liminf_{N \rightarrow \infty} \widehat{I}_{\fX,\vs}(\dot{\nu}_\ell^{(N)}) \ .
\end{eqnarray} 
Combining (\ref{AFF89}), (\ref{AFF90}) and (\ref{AFF91}) gives (\ref{AFF87}).
\end{proof}                                 
                  
\medskip
Just as in the unweighted case, it follows that 
the Fekete measures for the weighted transfinite diameter associated to 
\index{transfinite diameter!weighted $(\fX,\vs)$}
\index{Fekete!measure} 
$H_\ell$ and $W_\ell(z)$ converge weakly to the equilibrium distribution. 
\index{equilibrium distribution!weighted $(\fX,\vs)$}
More precisely, consider a sequence of pairs $(n_k,N_k)$ with 
$n_k/N_k \rightarrow \sigma_\ell$ and $N_k \rightarrow \infty$.  For a given $k$, 
let $\alpha_1, \ldots, \alpha_{n_k}$ achieve the maximum in the definition
of $d_{(n_k,N_k)}(H_\ell,W_\ell)$, and let    
\begin{equation*}
\nu_k \ = \ \sum_{i=1}^{n_k} \frac{1}{N_k} \delta_{\alpha_i}(z)
\end{equation*}
be the associated Fekete measure of mass $n_k/N_k$ on $H_\ell$.
\index{Fekete!measure} 

\begin{corollary} \label{ACF17} 
Let $H$, $H_\ell$, $\sigma_\ell = \mu_{\fX,\vs}(H_\ell)$,
and $W_\ell(z)$ be as in Theorem $\ref{APF15}$.  Then for any sequence of pairs $(n_k,N_k)$ with 
$n_k/N_k \rightarrow \sigma_\ell$ and $N_k \rightarrow \infty$,   
the corresponding sequence of Fekete measures 
\index{Fekete!measure} 
$\{\nu_k\}$ for $H_\ell$ relative to $W_\ell(z)$
converges weakly to $\mu_{\sigma_\ell,H_\ell,W_\ell} = \mu_{\fX,\vs} \vert_{H_\ell}$.
\end{corollary}

\begin{proof} The proof is the same as that of Corollary \ref{ACE14}, using
Theorems \ref{APF15} and \ref{ATF16}.
\end{proof} 

\medskip
We now determine the weighted Chebyshev constant, in the situation of Theorem \ref{APF15}.
\index{Chebyshev constant!weighted} 

\begin{theorem} \label{ATF18}
Let $H$, $H_\ell$, $\sigma_\ell = \mu_{\fX,\vs}(H_\ell)$,
and $W_\ell(z)$ be as in Theorem $\ref{APF15}$.  Then
\begin{equation} \label{AFF92}
  \CH^*_{\sigma_\ell}(H_\ell,W_\ell) \ = \ \gamma_{\fX,\vs}(H) \ .
\end{equation}
\end{theorem}

\begin{proof}
Since the Fekete measures converge weakly to the weighted equilibrium 
\index{Fekete!measure} 
distribution $\mu_\ell$, Proposition \ref{APD6} gives
\index{equilibrium distribution!weighted $(\fX,\vs)$}
\begin{equation*} 
 -\log(d_{\sigma_\ell}(H_\ell,W_\ell))    
  \ \le \ \sigma_\ell \cdot (-\log(\CH^*_{\sigma_\ell}(H_\ell,W_\ell)))
           + \int_{H_\ell} \widehat{u}_\ell(z) \, d\mu_\ell(z) \ .
\end{equation*} 
Comparing this with (\ref{AFF86}) gives 
\begin{equation}
V_{\fX,\vs}(H) \ \le \ -\log(\CH^*_{\sigma_\ell}(H_\ell,W_\ell)) \ . \label{AFF93}
\end{equation}
\index{Robin constant!$(\fX,\vs)$}

We will now show that equality holds in (\ref{AFF93}).  If not,   
then there would be an $\varepsilon > 0$ such that 
\begin{equation}
V_{\fX,\vs}(H) + \varepsilon \ < \ -\log(\CH^*_{\sigma_\ell}(H_\ell,W_\ell)) \ . \label{AFF94}
\end{equation} 
\index{Robin constant!$(\fX,\vs)$}
Take a sequence of pairs $(n_k,N_k)$ with 
$n_k/N_k \rightarrow \sigma_\ell$ and $N_k \rightarrow \infty$.  
Given $k$, let 
$\tP_k(z) := \tP_{(n_k,N_k)}(z,W_\ell)$ be the corresponding Chebyshev pseudopolynomial
\index{pseudopolynomial!weighted Chebyshev}
\index{Chebyshev pseudopolynomial!weighted} 
for $H_\ell$ relative to $W_\ell(z)$;  
put $\CH^*_k = \CH^*_{(n_k,N_k)} = (\|\tP_k\|_{H_\ell})^{1/N_k}$.  Also, let $\omega_k$ 
be the usual discrete measure of mass $n_k/N_k$ supported on the roots of 
$\tP_k(z)$.  

Our assumption (\ref{AFF94}) implies that on $H_\ell$, for each sufficiently large $k$,
\index{potential function} 
\begin{eqnarray*}
V_{\fX,\vs}(H) + 3\varepsilon/4  &<&  -\log(\CH^*_k) \\
      &\le& -\frac{1}{N_k} \log(\tP_k(z)) \ = \ u_{\omega_k}(z) + \widehat{u}_\ell(z) \ .
\end{eqnarray*} 
\index{Robin constant!$(\fX,\vs)$}
Let
\begin{equation*}
\widetilde{\omega}_k \ = \ \frac{\sigma_\ell}{n_k/N_k} \omega_k \ ,
\end{equation*}
renormalizing $\omega_k$ to have mass $\sigma_\ell$.  Then for all sufficiently large
$k$, we have 
\begin{equation}  \label{AFF95}        
V_{\fX,\vs}(H) + \varepsilon/2 \ \le \ u_{\widetilde{\omega}_k}(z) + \widehat{u}_\ell(z)
\end{equation}
\index{Robin constant!$(\fX,\vs)$}
on $H_\ell$.   But $\widehat{u}_\ell(z) = u_{\dot{\mu}_\ell}(z)$ and by Theorem \ref{ATE10B},
\index{potential function} 
\begin{equation} \label{AFF96}
u_{\mu_\ell}(z) + u_{\dot{\mu}_\ell}(z) \ = \ u_{\fX,\vs}(z) \ \le \ V_{\fX,\vs}(H) 
\end{equation}
\index{Robin constant!$(\fX,\vs)$}
on $H$.  By (\ref{AFF95}) and (\ref{AFF96}), 
\begin{equation}
u_{\mu_\ell}(z) + \varepsilon/2 \ \le \ u_{\widetilde{\omega}_k}(z)  \label{AFF97}
\end{equation}
on $H_\ell$.  However, $u_{\widetilde{\omega}_k}(z) - u_{\mu_\ell}(z)$ extends
to a function harmonic on $\cC_v(\CC) \backslash H_\ell$. By the Maximum principle,
\index{Maximum principle!for harmonic functions}
it follows that (\ref{AFF97}) holds throughout $\cC_v(\CC)$.  

We can obtain a contradiction from this as follows.  
Put $\nu_k = \widetilde{\omega}_k + \dot{\mu}_\ell$;  then $\nu_k$ and $\mu = \mu_\ell + \dot{\mu}_\ell$ 
\index{potential function}
are probability measures supported on $H$.  By (\ref{AFF96}), 
$u_{\nu_k}(z) = u_{\widetilde{\omega}_k}(z) + u_{\dot{\mu}_\ell}(z) 
\ge u_{\mu}(z) + \varepsilon/2$ on $H$.   
Hence by 
Since $\int_H u_{\mu}(z) \, d\mu(z) = I_{\fX,\vs}(\mu) = V_{\fX,\vs}(H)$, 
\index{Robin constant!$(\fX,\vs)$}
it follows from the Fubini-Tonelli theorem and the fact that 
\index{Fubini-Tonelli theorem} 
$u_{\mu}(z) \le V_{\fX,\vs})$ for all $z \in H$ (Theorem \ref{ATE10B}), that 
\index{equilibrium potential!$(\fX,\vs)$}
\begin{eqnarray*}
V_{\fX,\vs}(H) + \varepsilon/2 &\le& \int_H u_{\nu_k} \, d\mu(z) \\ 
       &=& \int_H u_{\mu}(z) \, d\nu_k(z) \ \le \ V_{\fX,\vs}(H) \ .
\end{eqnarray*}  
\index{Robin constant!$(\fX,\vs)$}
This contradiction shows that $V_{\fX,\vs}(H) = -\log(\CH^*_{\sigma_\ell}(H_\ell,W_\ell))$,
which is equivalent to the assertion in the theorem.  
\end{proof}       

\medskip  
Finally, we show that under appropriate hypotheses,
the discrete measures attached to Chebyshev pseudopolynomials 
\index{Chebyshev pseudopolynomial!weighted} 
\index{pseudopolynomial!weighted Chebyshev}
\index{Chebyshev measure!converge weakly to $\mu_{\fX,\vs}$|ii} 
converge weakly to the equilibrium distribution.
\index{equilibrium distribution}      
As before, take a sequence of pairs $(n_k,N_k)$ with 
$n_k/N_k \rightarrow \sigma_\ell$ and $N_k \rightarrow \infty$. Let the Chebyshev measure
\index{Chebyshev measure}  
$\omega_k$ be the discrete measure of mass $n_k/N_k$ supported equally on the roots of 
$\tP_{(n_k,N_k)}(z)$ for $H_\ell$ relative to $W_\ell(z)$.    

\begin{theorem} \label{ATF19} 
Let $H$, $H_\ell$, $\sigma_\ell = \mu_{\fX,\vs}(H_\ell)$,
and $W_\ell(z)$ be as in Theorem $\ref{APF15}$.  Assume also that $\cC_v(\CC) \backslash H_\ell$
is connected, and that $H_\ell$ has empty interior.
Then for any sequence of pairs $(n_k,N_k)$ with 
$n_k/N_k \rightarrow \sigma_\ell$ and $N_k \rightarrow \infty$,
the corresponding sequence of Chebyshev measures
\index{Chebyshev measure} 
$\{\omega_k\}$ for $H_\ell$ relative to $W_\ell(z)$ 
converges weakly to the equilibrium distribution
$\mu_{\sigma_\ell,H_\ell,W_\ell} = \mu_{\fX,\vs} \vert_{H_\ell}$.
\end{theorem}
\index{equilibrium distribution!$(\fX,\vs)$}

\begin{proof}  
Recall that $\mu_\ell = \mu_{\fX,\vs}|_{H_\ell}$, $\dot{\mu}_\ell = \mu_{\fX,\vs}|_{H_2}$.
Let the numbers $\CH^*_k$, the Chebyshev pseudopolynomials $\tP_k(z)$, and
\index{pseudopolynomial!weighted Chebyshev}
\index{Chebyshev pseudopolynomial!weighted} 
the measures $\widetilde{\omega}_k$, 
$\nu_k = \widetilde{\omega}_k + \dot{\mu}_\ell$ be as in the proof of 
Theorem \ref{ATF18}.  Let $\omega$ be a weak limit of the $\omega_k$;  
after passing to a subsequence, if necessary, we can assume that the 
full sequence converges weakly to it.  Clearly $\omega$ is also the weak limit
of the $\widetilde{\omega}_k$.  

Fix $\varepsilon > 0$.  
Since $\CH^*_{\sigma_\ell}(H_\ell,W_\ell) = \CH^*_{\fX,\vs}(H)$ by Theorem \ref{ATF18}, 
an argument similar to the one in the proof of Theorem \ref{ATF18}
shows that for sufficiently large $k$ 
\index{potential function}
\begin{equation*} 
V_{\fX,\vs}(H) - \varepsilon \ < \ u_{\nu_k}(z) 
         \ = \ u_{\widetilde{\omega}_k}(z) + u_{\dot{\mu}_\ell}(z) \ .
\end{equation*}    
\index{Robin constant!$(\fX,\vs)$}      
on $H_\ell$.  On the other hand,
\begin{equation*}
u_{\mu_\ell}(z) + u_{\dot{\mu}_\ell}(z) \ = \ u_{\fX,\vs}(z) \ \le \ V_{\fX,\vs}(H) \ .
\end{equation*}  
Subtracting, we see that 
\index{potential function}
\begin{equation} \label{AFF98}
u_{\widetilde{\omega}_k}(z) \ \ge \ u_{\mu_\ell}(z) - \varepsilon 
\end{equation}
on $H_\ell$.  Since $u_{\widetilde{\omega}_k}(z) - u_{\mu_\ell}(z)$ extends to a function
harmonic in all of $\cC_v(\CC) \backslash H_\ell$, the Maximum principle for
harmonic functions shows (\ref{AFF98}) holds in all of 
\index{Maximum principle!for harmonic functions}
$\cC_v(\CC)$.  

As the $u_{\widetilde{\omega}_k}(z)$ converge uniformly to $u_{\omega}(z)$ 
\index{potential function}
on compact subsets of $\cC_v(\CC) \backslash H_\ell$ and $\varepsilon > 0$ is
arbitrary, it follows from (\ref{AFF98}) that 
\begin{equation} \label{AFF99} 
          u_{\omega}(z) \ \ge \ u_{\mu_\ell}(z) 
\end{equation}          
for all $z \in \cC_v(\CC) \backslash H_\ell$.  

Suppose equality in (\ref{AFF99}) failed to hold for some 
$z_0 \in \cC_v(\CC) \backslash H_\ell$.  Since $\cC_v(\CC) \backslash H_\ell$
is connected, the Maximum principle implies that
\index{Maximum principle!for harmonic functions} 
$u_{\omega}(z) > u_{\mu_\ell}(z)$ for all  $z \in \cC_v(\CC) \backslash H_\ell$.  
In particular, there would be a $\delta > 0$ such that 
\begin{equation*}
          u_{\omega}(z) \ \ge \ u_{\mu_\ell}(z) + \delta
\end{equation*}  
on $H_2$.  Now put $\nu = \omega + \dot{\mu}_\ell$;  then
\index{potential function}
\begin{eqnarray*} 
u_{\nu}(z) &=& u_{\omega}(z) + u_{\dot{\mu}_\ell}(z) \\
        &\ge&  u_{\mu_\ell}(z) + u_{\dot{\mu}_\ell}(z) \ = \ u_{\fX,\vs}(z) 
\end{eqnarray*}
on $\cC_v(\CC) \backslash H_\ell$, with $u_{\nu} \ge u_{\fX,\vs}(z) + \delta$
on $H_2$.  

Let $e = e_{\fX,\vs}$ be the exceptional subset of $H$ inner capacity $0$ 
\index{capacity $= 0$}
where $u_{\fX,\vs}(z) < V_{\fX,\vs}(H)$, given by Theorem \ref{ATE10B}. 
\index{equilibrium potential!$(\fX,\vs)$}
\index{equilibrium potential!takes constant value a.e. on $E_v$}
\index{Robin constant!$(\fX,\vs)$} 
We claim that $u_{\nu}(z) \ge V_{\fX,\vs}(H)$ on $H \backslash e$.  
To see this, for each $\eta > 0$ let
\begin{equation*}
U_{\eta} \ = \ \{ z \in \cC_v(\CC) : u_{\fX,\vs}(z) > V_{\fX,\vs}(H) - \eta \} \ .
\end{equation*}
\index{Robin constant!$(\fX,\vs)$}
This is an open set, since $u_{\fX,\vs}(z)$ is lower semi-continuous.
\index{equilibrium potential!$(\fX,\vs)$}
\index{equilibrium potential!is lower semi-continuous}
\index{semi-continuous!potential function is lower semi-continuous}

Furthermore, at each point of $H \backslash e$,
$u_{\fX,\vs}(z)$ is continuous and equal to $V_{\fX,\vs}(H)$ 
\index{Robin constant!$(\fX,\vs)$}
(Theorem \ref{ATE10B})
so apart from a subset of inner capacity $0$,
\index{capacity $= 0$}
$\partial U_{\eta}$ is contained in $\cC_v(\CC) \backslash H$.  In particular
on $\partial U_{\eta} \backslash e$ we have
\begin{equation*}
u_{\nu}(z) \ \ge \ u_{\fX,\vs}(z) \ = \ V_{\fX,\vs}(H) - \eta \ .
\end{equation*}   
\index{equilibrium potential!$(\fX,\vs)$}
\index{Robin constant!$(\fX,\vs)$}
\index{potential function!is superharmonic}
Since $u_{\nu}(z)$ is superharmonic\index{superharmonic} and bounded from below on $U_{\eta}$, 
the strong form of the Maximum principle (see \cite{RR1}, Proposition 3.1.1) shows that 
\index{Maximum principle!for harmonic functions!strong form}
$u_{\nu}(z) \ge V_{\fX,\vs}(H) - \eta$ on $U_{\eta}$.  Since
\index{Robin constant!$(\fX,\vs)$}
$\eta > 0$ is arbitrary, and $H \backslash e \subset U_{\eta}$,
it follows that $u_{\nu}(z) \ge V_{\fX,\vs}(H)$  on $H \backslash e$.
\index{potential function}

Because $u_{\fX,\vs}(z) = V_{\fX,\vs}(H)$ on $(H \backslash H_\ell) \backslash e$, 
\index{Robin constant!$(\fX,\vs)$}
\index{equilibrium potential!$(\fX,\vs)$}
we have $u_{\nu}(z) \ge V_{\fX,\vs}(H) + \delta$  
on $(H \backslash H_\ell) \backslash e$.  However, a set of inner capacity $0$
\index{capacity $= 0$}
necessarily has $\mu$-measure $0$ (\cite{RR1}, Lemma 3.1.4).  By Fubini-Tonelli, 
\index{Fubini-Tonelli theorem}
\begin{eqnarray} 
V_{\fX,\vs}(H) + \mu_{\fX,\vs}(H \backslash H_\ell) \cdot \delta  &\le& \int_H u_{\nu}(z) d\mu_{\fX,\vs}(z) \notag \\
  &=&  \int_{H} u_{\fX,\vs}(z) d\nu(z) \ \le \ V_{\fX,\vs}(H) \ . \label{AFF100}
\end{eqnarray} 
\index{equilibrium potential!$(\fX,\vs)$}
\index{Robin constant!$(\fX,\vs)$}
We claim that $\mu_{\fX,\vs}(H \backslash H_\ell) > 0$.  
Otherwise $u_{\fX,\vs}$ would be supported on $H_\ell$, 
and the fact that $\cC_v(\CC) \backslash H_\ell$ is connected and contains $\fX$ would mean 
that $u_{\fX,\vs}(z) < V_{\fX,\vs}(H)$ on $\cC_v(\CC) \backslash H_\ell$. 
\index{equilibrium potential!$(\fX,\vs)$} 
\index{Robin constant!$(\fX,\vs)$}
However, $u_{\fX,\vs}(z) = V_{\fX,\vs}(H)$ for all $z \in H$ except possibly a set of inner
capacity $0$, which contradicts that 
\index{equilibrium potential!takes constant value a.e. on $E_v$}
\index{capacity $= 0$} 
$H \backslash H_\ell \subset \cC_v(\CC) \backslash H_\ell$ has positive capacity.
\index{capacity $> 0$} 
Thus (\ref{AFF100}) is impossible.  

We conclude that in (\ref{AFF99}), we have $u_{\nu}(z) = u_{\fX,\vs}(z)$ for all $z \notin H_\ell$.
\index{equilibrium potential!$(\fX,\vs)$}
Moreover, we have shown that
$u_{\nu}(z) \ge V_{\fX,\vs}(H)$ for all $z \in H \backslash e$.
\index{Robin constant!$(\fX,\vs)$}
We now claim that $u_{\nu}(z) \le V_{\fX,\vs}(H)$ for all $z$.  Suppose
to the contrary that $u_{\nu}(z_0) > V_{\fX,\vs}(H)$ for some $z_0$.  Then
since $u_{\nu}(z)$ is lower semi-continuous,
\index{potential function!is lower semi-continuous}
\index{semi-continuous!potential function is lower semi-continuous}
\begin{equation*}
U \ := \ \{ z \in \cC_v(\CC) : u_{\nu}(z) > V_{\fX,\vs}(H) \}
\end{equation*}\index{Robin constant!$(\fX,\vs)$}
would be a nonempty open set.  Since $H_\ell$ has no interior, $U$ contains
points of $\cC_v(\CC) \backslash H_\ell$.  However, at these points
$u_{\nu}(z) = u_{\fX,\vs}(z) \le V_{\fX,\vs}(H)$, 
contradicting the definition of $U$.
\index{Robin constant!$(\fX,\vs)$}

It follows that $u_{\nu}(z)$ coincides with $u_{\fX,\vs}(z)$ 
\index{equilibrium potential!$(\fX,\vs)$}
\index{equilibrium potential!takes constant value a.e. on $E_v$}
except possibly on the exceptional set $e$ of inner capacity $0$.\index{capacity $= 0$}  
However, a superharmonic\index{superharmonic}\index{potential function!is superharmonic}
function which is bounded below is determined by its values
on the complement any set of inner capacity $0$
\index{capacity $= 0$}
(\cite{Ts}, Theorem III.28, p.78).  
Thus $u_{\nu}(z) = u_{\fX,\vs}(z)$ for all $z$.
\index{equilibrium potential!$(\fX,\vs)$}
By the Riesz Decomposition theorem (\cite{RR1}, Theorem 3.1.11)
\index{Riesz Decomposition theorem}
we can recover the measure from the potential function, 
\index{equilibrium potential!$(\fX,\vs)$}
so $\nu = \mu_{\fX,\vs}$, 
and hence $\omega = \mu_\ell$.
\end{proof}

\medskip
\noindent{\bf Remark.}  Some hypotheses on $H_\ell$ are necessary 
in Theorem \ref{ATF19};  it is not always true that the Chebyshev
\index{Chebyshev measure} 
measures converge to the equilibrium distribution.
\index{equilibrium distribution!$(\fX,\vs)$}  
For example, in the classical case consider the
unit disc $H_\ell = D(0,1) \subset \CC$, with weight $W(z) \equiv 1$,  
relative to $\fX = \{\infty\}$. 
The Chebyshev polynomial of degree $n$ for $H_\ell$ is $z^n$, 
\index{Chebyshev polynomial} 
with roots only at the origin.  However, the equilibrium 
distribution of $H_\ell$ is the uniform measure supported on $\partial H_\ell = C(0,1)$. 
\index{equilibrium distribution!$(\fX,\vs)$}


\section{ Chebyshev Pseudopolynomials for short intervals }
\label{ChebShortIntervalSection}
 
In this section we will show that when $H$ is a ``sufficiently short'' interval, 
\index{short@`short' interval}  
the weighted $(\fX,\vs)$-Chebyshev pseudopolynomials for $H$ have oscillation
\index{Chebyshev pseudopolynomial!weighted $(\fX,\vs)$} 
\index{pseudopolynomial!weighted $(\fX,\vs)$}
properties like those of classical Chebyshev polynomials. 
\index{Chebyshev polynomial}
The notion of ``shortness'' depends on the location of $H$ relative to $\fX$
and the existence of a suitable system of coordinates, 
but is independent of the choice of the weight function.
\index{short@`short' interval}   

The motivating case is case when
$K_v \cong \RR$ and $H \subset \cC_v(\RR) \backslash \fX$ is a closed interval. 
However, since any analytic arc becomes an interval in suitable coordinates,
\index{arc!analytic}
the results apply more generally.

\vskip .1 in
Fix a local coordinate patch $U \subset \cC_v(\CC)$,
with coordinate function $z$ say.  Thus, $z$ gives a holomorphic 
isomorphism between $U$ and a simply connected\index{simply connected} open set $z(U) \subset \CC$.  
We can decompose the coordinate function into its real and imaginary parts, 
$z = u + i v$, and speak of the real and imaginary coordinates of points in $U$.  
For us, the case of interest is when $z(U) \cap \RR$ is nonempty;  
to simplify notation, we will assume that is the situation, 
and that $U$ and $z$ have been chosen so that $v = 0$ on $z^{-1}(z(U) \cap \RR)$.  

By a real interval $H \subset U$, we mean
a set of the form $H = z^{-1}([a,b])$ where $[a,b] \subset z(U) \cap \RR$.
By abuse of notation, we will simply write $H = [a,b]$.  
Similarly, we can speak of a disc $D(t,r) \subset U$.
Using the coordinate function to identify $U$ with $z(U) \subset \CC$,  
we can speak of translating a point $p \in U$ by a number $c \in \CC$:  
$p \mapsto p+c$, provided that both points involved belong to $U$.  
(Formally, if $p \in U$, then $p+c$ means the point $z^{-1}(z(p) + c)$, 
if $z(p) + c \in z(U)$.)    

Recall that $\fX = \{x_1, \ldots, x_m\}$.  
By Proposition \ref{APropA2}, relative to the given local coordinate,  
for each $x_j \in \fX$ there is a $\cC^{\infty}$ function $\eta_j(z,w)$ 
on $(U \backslash \ \{x_j\}) \times (U \backslash \{x_j\})$ 
(which is harmonic in each variable separately) 
such that for all $z, w \in U \backslash \{x_j\}$ 
\begin{equation*}
-\log([z,w]_{x_j}) \ = \ -\log(|z-w|) + \eta_j(z,w) \ .
\end{equation*}
Writing $z = u_1 + i v_1$, $w = u_2 + i v_2$, we can speak of the partial 
derivatives of $\eta_j(z,w)$ relative to $u_1$, $v_1$, $u_2$ and $v_2$.  

\begin{lemma} \label{bL1}  Let $U \subset \cC_v(\CC)$ be a local
coordinate patch, and let $H = [a,b] \subset U$  be a real interval 
disjoint from $\fX$.

Let $C$ be a bound such that, uniformly for all $x_j \in \fX$ and  
all $z, w \in H$,  
\begin{equation} \label{FAX1}
\left|\frac{\partial \eta_j}{\partial u_2}(z,w)\right| \ \le \ C \ , \qquad 
\left|\frac{\partial^2 \eta_j}{\partial u_2^2}(z,w)\right| \le \ C \ .  
\end{equation}
Then for each $x_j \in \fX$, 
all $z, c, d \in H$, and each $0 < \varepsilon \in \RR$ 
such that $c \pm \varepsilon, d \pm \varepsilon$ belong to $H$, 
  
$(1)$  $|\eta_j(z,c-\varepsilon) - \eta_j(z,c) + \eta_j(z,d+\varepsilon)- \eta_j(z,d)| 
          \ < \ C |d-c| \varepsilon + 2 C \varepsilon^2$ \ ;
           
$(2)$  $| \eta_j(z,c-\varepsilon) + \eta_j(z,c+\varepsilon) - 2\eta(z,c)| 
          \ < \ 2C \varepsilon^2$ \ ;

$(3)$  $|\eta_j(z,c+\varepsilon) - \eta_j(z,c)| \ < \ C  \varepsilon$ \ ;

$(4)$  $|\eta_j(z,d-\varepsilon) - \eta_j(z,d)| \ < \ C \varepsilon$ \ .
\end{lemma}

\begin{proof}
All of these follow from the Mean Value Theorem.  
\index{Mean Value theorem}
For example, we prove $(1)$:
for appropriate $c^* \in (c-\varepsilon,c)$, $d^* \in (d,d+\varepsilon)$ and 
$e^* \in (c^*,d^*)$
\begin{eqnarray*}
\lefteqn{ \left| \eta(z,c-\varepsilon) - \eta(z,c) 
     + \eta(z,d+\varepsilon)- \eta(z,d) \right| } \\
   &=& \left| -\varepsilon(\frac{\partial \eta}{\partial u_2}(z,c^*))
          + \varepsilon(\frac{\partial \eta}{\partial u_2}(z,d^*)) \right|
   \ = \ \left| \varepsilon \cdot (d^*-c^*) 
              (\frac{\partial^2 \eta}{\partial u_2^2}(z,e^*)) \right| \\
   &<& \varepsilon (|d-c| + 2\varepsilon) \cdot C \ .
\end{eqnarray*} 
Clearly $(2)$ is a special case of $(1)$, and $(3)$ and $(4)$ are easy.    
\end{proof}

\vskip .1 in
Given a probability vector $\vs = (s_1, \ldots, s_m)$,  
on $U \backslash \fX$ we have 
\begin{equation*}
-\log([z,w]_{\fX,\vs}) \ = \ -\log(|z-w|) + \eta(z,w) \ .
\end{equation*}
where $\eta(z,w) = \eta_{\fX,\vs}(z,w) = \sum_{j=1}^m s_j \eta_j(z,w)$.  
By the triangle inequality and the fact that $\sum s_j = 1$, 
the bounds in Lemma \ref{bL1} hold with $\eta_j(z)$ replaced by $\eta(z)$. 

\vskip .1 in
 
\begin{definition} \label{ShortnessDef} Let $U \subset \cC_v(\CC)$ be
a coordinate patch, with coordinate function $z$.  A real interval
$H = [a,b] \subset U$ is {\em short} $($relative to $\fX$ and the coordinate function $z$ on $U)$ 
if it is disjoint from $\fX$ and $|b-a| < \min(1/C(H,\fX), 1/\sqrt{2C(H,\fX)})$, where 
\begin{equation} \label{FShort} 
C(H,\fX) \ = \ \max_{x_j \in \fX} \max_{z, w \in H} 
      \left(\left|\frac{\partial \eta_j}{\partial u_2}(z,w)\right|, 
           \left|\frac{\partial^2 \eta_j}{\partial u_2^2}(z,w)\right|\right) \ . 
\end{equation}     
\end{definition} 
\index{short@`short' interval|ii}  

While this notion of ``shortness'' is ugly, it is easy to apply.    
In practice, we will be given a coordinate patch $U$ whose closure
$\overline{U}$ is disjoint from $\fX$.  Defining $C(U)$ as in 
(\ref{FShort}) with $H$ replaced by $\overline{U}$, 
one sees that any real interval $[a,b] \subset U$ with  
$|b-a| < \max(1/C(U), 1/\sqrt{2C(U)})$  is ``short''. 
\index{short@`short' interval}   

We will now show that restricted, 
weighted $(\fX,\vs)$-Chebyshev pseudopolynomials for 
\index{Chebyshev pseudopolynomial!weighted $(\fX,\vs)$}
\index{pseudopolynomial!weighted Chebyshev}
short real intervals behave like classical Chebyshev polynomials.
\index{short@`short' interval}  
\index{Chebyshev polynomial}

\begin{proposition} \label{bP1} 
Let $U \subset \cC_v(\CC)$ be a coordinate patch with coordinate function $z$.
Suppose $H = [a,b] \subset U \backslash \fX$ is a short real interval. 
\index{short@`short' interval}  
 
Fix $\vs \in \cP^m$, and let $W(z)$ be weight function which is 
continuous, positive and bounded on a neighborhood of $H$.  Then for  
any pair $(n,N)$ with $n \ge 1$, each Chebyshev pseudopolynomial
\index{pseudopolynomial!weighted Chebyshev}
\index{Chebyshev pseudopolynomial!weighted $(\fX,\vs)$} 
$\tP_{(n,N)}(z,W) = \prod_{i=1}^n [z,\alpha_i]_{\fX,\vs} \cdot W(z)^N$ 
for $H$ relative to $W(z)$ has the following properties.  Assume
the roots $\alpha_i$ are labeled in increasing order.  

$(1)$  $\tP_{(n,N)}(z,W)$ has distinct roots 
           which lie in the interior of $[a,b]$.   

$(2)$  If $\alpha_1$ is the leftmost root of $\tP_{(n,N)}(z,W)$, there is 
    a point $\alpha_0 \in [a,\alpha_1]$ where  
    \begin{equation*}
       \tP_{(n,N)}(\alpha_0,W) \ = \ \|\tP_{(n,N)}\|_H \ . 
    \end{equation*}

$(3)$   For each pair of consecutive roots $\alpha_i, \alpha_{i+1}$ 
    there is a point $\alpha_i \in (\alpha_i,\alpha_{i+1})$ where 
    \begin{equation*}      
       \tP_{(n,N)}(\alpha_i,W) \ = \ \|\tP_{(n,N)}\|_H \ .
    \end{equation*}
    
$(4)$ If $\alpha_{n}$ is the rightmost root, 
     there is a point $\alpha_{n} \in (\alpha_{n},b]$ where 
    \begin{equation*}
       \tP_{(n,N)}(\alpha_{n},W) \ = \ \|\tP_{(n,N)}\|_H \ .
    \end{equation*} 
\end{proposition}

\medskip
\begin{proof} 
Write $C = C(H,\fX)$, where $C(H,\fX)$ is as in Definition \ref{ShortnessDef}.  
The hypothesis that $H$ is ``short'' implies that 
\index{short@`short' interval}  
\begin{equation}  \label{FShApp} 
\frac{1}{b-a} >  C \ , \qquad \frac{1}{(b-a)^2} > 2 C \ .
\end{equation}

We will now show that in the presence of the bounds (\ref{FShApp}),
the classical arguments concerning oscillation properties 
of Chebyshev polynomials carry over. 
\index{Chebyshev polynomial}
The weight function plays a negligible role; 
it cancels out in all the ratios below.
   
Fix $(n,N)$.  We already know that $\tP_{(n,N)}(z,W)$ exists;  
the problem is to show it has the properties above.  
Write $\tP(z) = \tP_{(n,N)}(z,W)$, $M = \|\tP(z)\|_H$.

First, suppose $\tP(z)$ had a root at $z = a$.  Then  
$\tP(a) = 0$ and so there would be an interval $[a,a+\delta]$ 
on which $\tP(z) < M$.  We claim that by replacing the factor $[z,a]_{\fX,\vs}$ in $\tP(z)$
by $[z,a+\varepsilon]_{\fX,\vs}$ for an appropriately small $\varepsilon$, we could
reduce $\|\tP\|_H$.  Let $P_{\varepsilon}(z)$ be the pseudopolynomial thus obtained.
\index{pseudopolynomial!weighted Chebyshev}
By the continuity of $[z,w]_{\fX,\vs}$, for sufficiently small $\varepsilon >0$ 
we would still have $P_{\varepsilon}(z) < M$ for $z \in [a,a+\delta]$.  
For $z \in (a+\delta,b]$, 
\begin{equation*}
\frac{P_{\varepsilon}(z)}{\tP(z)} \ = \ \frac{[z,a+\varepsilon]_{\fX,\vs}}{[z,a]_{\fX,\vs}} 
\end{equation*}
and 
\begin{equation*}
-\log\left(\frac{[z,a+\varepsilon]_{\fX,\vs}}{[z,a]_{\fX,\vs}}\right) \ = \ 
   -\log\left(\frac{|z-a-\varepsilon|}{|z-a|}\right) 
          + \eta(z,a+\varepsilon) - \eta(z,a) \ .
\end{equation*}
By Lemma \ref{bL1}, for all $z \in H$,
$|\eta(z,a+\varepsilon) - \eta(z,a)| < C \varepsilon$.  On the other hand,
for $z \in (a+\delta,b]$ and $0 < \varepsilon < \delta$,
\begin{equation*}
 -\log\left(\frac{|z-a-\varepsilon|}{|z-a|}\right) 
    \ = \ -\log\left(|1 - \frac{\varepsilon}{z-a} |\right) 
    \ > \ \frac{\varepsilon}{|z-a|} \ .
\end{equation*}
Since $[a,b]$ is short enough that $1/(b-a) > C$, also $1/(z-a) > C$ and so  
\index{short@`short' interval}  
\begin{equation*}
-\log\left(\frac{[z,a+\varepsilon]_{\fX,\vs}}{[z,a]_{\fX,\vs}}\right) \ > \ 0 \ ,
\end{equation*}
Thus $[z,a+\varepsilon]_{\fX,\vs} < [z,a]_{\fX,\vs}$ for all $z \in [a+\delta,b]$.  
It follows that $\|P_{\varepsilon}(z)\|_H < \|\tP\|_H$,
contradicting the minimality of $\|\tP\|_H$.  A similar argument shows
$\tP(z)$ cannot have a root at $z = b$.

Let $\alpha_1$ be the leftmost root of $\tP(z)$.
If $\tP(z)$ did not achieve its maximum in $[a,\alpha_1]$, 
an argument like the one above shows 
we could reduce $\|\tP\|_H$ by moving $\alpha_1$ slightly to the right.
For similar reasons, if $\alpha_{n}$ is the rightmost root and 
$\tP(z)$ did not achieve its maximum
in $[\alpha_{n},b]$, we could reduce $\|\tP\|_H$ by moving $\alpha_{n}$ to the left.

Next, suppose $\tP(z)$ had a double root at $\alpha_i$, say.  As shown above, 
$\alpha_i \in (a,b)$.  Since $\tP(\alpha_i) = 0$, there is a $\delta > 0$ such that 
$\tP(z) < M$ for $z \in [\alpha_i-\delta,\alpha_i+\delta]$.
If we define $P_{\varepsilon}(z)$ by 
replacing the factor $[z,\alpha_i]_{\fX,\vs}^2$ in $\tP(z)$ by
$[z,\alpha_i-\varepsilon]_{\fX,\vs}[z,\alpha_i+\varepsilon]_{\fX,\vs}$, 
then by the continuity of $[z,w]_{\fX,\vs}$, for sufficiently small $\varepsilon$
we will have $P_{\varepsilon}(z) < M$ for $z \in [\alpha_i-\delta,\alpha_i+\delta]$.

For $z \in H \backslash [\alpha_i-\delta,\alpha_i+\delta]$,  
\begin{equation*}
\frac{P_{\varepsilon}(z)}{\tP(z)} \ = \ 
           \frac{[z,\alpha_i-\varepsilon]_{\fX,\vs}[z,\alpha_i+\varepsilon]_{\fX,\vs}}{[z,\alpha_i]_{\fX,\vs}^2}
\end{equation*}
and  
\begin{eqnarray*}
\lefteqn{-\log\left(\frac{[z,\alpha_i-\varepsilon]_{\fX,\vs}[z,\alpha_i+\varepsilon]_{\fX,\vs}}
                          {[z,\alpha_i]_{\fX,\vs}^2}\right)} & & \\
& = &  -\log\left(\frac{(z-\alpha_i+\varepsilon)(z-\alpha_i-\varepsilon)}
             {(z-\alpha_i)^2}\right) 
 + \eta(z,\alpha_i-\varepsilon)  
                 + \eta(z,\alpha_i+\varepsilon) - 2\eta(z,\alpha_i) \ .
\end{eqnarray*}
We can assume $0 < \varepsilon < \delta$, so since $|z-\alpha_i| > \delta$,
\begin{equation*}
 -\log\left(\frac{(z-\alpha_i+\varepsilon)(z-\alpha_i-\varepsilon)}{(z-\alpha_i)^2}\right)
    \ = \ -\log\left(1 - \frac{\varepsilon^2}{(z-\alpha_i)^2} \right) 
    \ > \ \frac{\varepsilon^2}{(z-\alpha_i)^2} \ .
\end{equation*}
On the other hand, by Lemma \ref{bL1},    
$|\eta(z,\alpha_i-\varepsilon) +\eta(z,\alpha_i+\varepsilon) - 2\eta(z,\alpha_i)|
< 2C \varepsilon^2$.  As $[a,b]$ is short enough that $1/(b-a)^2 > 2C$,
\index{short@`short' interval}  
then $\varepsilon^2/(z-\alpha_i)^2 > 2C\varepsilon^2$ for
$z \in H \backslash [\alpha_i-\varepsilon,\alpha_i+\varepsilon]$, so  
\begin{equation*}
-\log\left(\frac{[z,\alpha_i-\varepsilon]_{\fX,\vs}[z,\alpha_i+\varepsilon]_{\fX,\vs}} 
                  {[z,\alpha_i]_{\fX,\vs}^2}\right) \ > \  0 \ \ .
\end{equation*}
Hence $[z,\alpha_i-\varepsilon]_{\fX,\vs}[z,\alpha_i+\varepsilon]_{\fX,\vs} < [z,\alpha_i]_{\fX,\vs}^2$ for all 
$z \notin [\alpha_i-\delta,\alpha_i+\delta]$, which implies that   
$\|P_{\varepsilon}(z)\|_H < \|\tP\|_H$.  
This contradicts the minimality of $\|\tP\|_H$.  

Finally, suppose $\tP(z)$ did not take on the value $M$ between two
consecutive roots $\alpha_i, \alpha_{i+1}$.  There is a $\delta > 0$
such that $\tP(z) < M$ for all $z \in [\alpha_i-\delta,\alpha_{i+1}+\delta]$.  
Define $P_{\varepsilon}(z)$ by replacing the product 
$[z,\alpha_i]_{\fX,\vs}[z,\alpha_{i+1}]_{\fX,\vs}$ in $\tP(z)$
with $[z,\alpha_i-\varepsilon]_{\fX,\vs}[z,\alpha_{i+1}+\varepsilon]_{\fX,\vs}$.
By the continuity of $[z,w]_{\fX,\vs}$, for sufficiently small $\varepsilon > 0$ we
have $P_{\varepsilon}(z) < M$ for $z \in [\alpha_i-\delta,\alpha_{i+1}+\delta]$.       
For $z \ne \alpha_i, \alpha_{i+1}$ we have 
\begin{equation*}
\frac{P_{\varepsilon}(z)}{\tP(z)}\ = \ 
           \frac{[z,\alpha_i-\varepsilon]_{\fX,\vs}[z,\alpha_{i+1}+\varepsilon]_{\fX,\vs}}
                            {[z,\alpha_i]_{\fX,\vs}[z,\alpha_{i+1}]_{\fX,\vs}} \ .
\end{equation*}
Furthermore, for $z \in [a,b] $ but $z \notin [\alpha_i-\delta,\alpha_{i+1}+\delta]$,  
\begin{eqnarray*}
 \lefteqn{ -\log\left( \frac{[z,\alpha_i-\varepsilon]_{\fX,\vs}[z,\alpha_{i+1}+\varepsilon]_{\fX,\vs}}
         {[z,\alpha_i]_{\fX,\vs}[z,\alpha_{i+1}]_{\fX,\vs}} \right) 
 \ = \ -\log\left(\frac{(z-\alpha_i+\varepsilon)(z-\alpha_{i+1}-\varepsilon}
                         {(z-\alpha_i)(z-\alpha_{i+1})}\right) }\\ 
  & & \qquad \qquad \qquad \qquad + \ \eta(z,\alpha_i-\varepsilon) 
      + \eta(z,\alpha_{i+1}+\varepsilon) - \eta(z,\alpha_i) - \eta(z,\alpha_{i+1})\ .
\end{eqnarray*}
Here
$|\eta(z,\alpha_i-\varepsilon) + \eta(z,\alpha_{i+1} 
          + \varepsilon) - \eta(z,\alpha_i) - \eta(z,\alpha_{i+1})|
< C (\alpha_{i+1}-\alpha_i) \varepsilon + 2 C \varepsilon^2$
by Lemma \ref{bL1}.  On the other hand, for small enough 
$\varepsilon$ and for $z \notin [\alpha_i-\delta,\alpha_{i+1}+\delta]$,
\begin{eqnarray*}
  -\log\left(\frac{(z-\alpha_i+\varepsilon)(z-\alpha_{i+1}-\varepsilon)}
                      {(z-\alpha_i)(z-\alpha_{i+1})}\right)
    &=& -\log\left(1 - \frac{\varepsilon(\alpha_{i+1}-\alpha_i) +
                   \varepsilon^2}{(z-\alpha_i)(z-\alpha_{i+1})} \right) \\
    &>& \frac{\varepsilon(\alpha_{i+1}-\alpha_i)}{(z-\alpha_i)(z-\alpha_{i+1})} \ .
\end{eqnarray*}
As $1/(b-a)^2 > 2C$, then $\frac{\varepsilon (\alpha_{i+1}-\alpha_i)}{(z-\alpha_i)(z-\alpha_{i+1})} 
                 \ > \ 2C(\alpha_{i+1}-\alpha_i)\varepsilon$.
Thus for sufficiently small $\varepsilon > 0$
\begin{equation*}
-\log\left( \frac{[z,\alpha_i-\varepsilon]_{\fX,\vs}[z,\alpha_{i+1}+\varepsilon]_{\fX,\vs}}
                   {[z,\alpha_i]_{\fX,\vs}[z,\alpha_{i+1}]_{\fX,\vs}} \right) \ > \ 0 \ ,
\end{equation*}
whence
$[z,\alpha_i-\varepsilon]_{\fX,\vs}[z,\alpha_{i+1}+\varepsilon]_{\fX,\vs} < [z,\alpha_i]_{\fX,\vs}[z,\alpha_{i+1}]_{\fX,\vs}$ 
for all $z$ as above. This shows that $\|P_{\varepsilon}(z)\|_H < \|\tP\|_H$ 
for small enough $\varepsilon$, once more contradicting the 
minimality of $\|\tP\|_H$.
\end{proof}

\vskip .1 in
\noindent{\bf Remark.}  Because of the presence of the weight function,
$\tP_{(n,N)}(z,W)$  need not take on its maximum value at $a$ and $b$
(as holds classically).  However, $\tP_{(n,N)}(z,W)$ does 
``vary $n$ times from $M$ to $0$ to $M$'' on $H$,
which is enough for the application in the next section.    


\section{ Oscillating Pseudopolynomials.}
   \label{OscillatingPseudoPolySection} 

In this section we specialize to the case $K_v \cong \RR$. 
Thus, $\cC_v(\RR)$ has meaning, 
and there is an action of complex conjugation $z \mapsto \zbar$ on $\cC_v(\CC)$.  

We will assume that $\fX$ and $H \subset \cC_v(\CC) \backslash \fX$ are stable under complex conjugation, 
that $H$ is compact, and that $H$ has finitely many connected components $H_1, \ldots, H_D$,   
where no $H_{\ell}$ is reduced to a point, and each $H_{\ell}$ is simply connected.
Under these hypotheses $\cC_v(\CC) \backslash H_{\ell}$ is connected for each $\ell$, 
and $\cC_v(\CC) \backslash H$ is connected.  Since $H$ is stable 
under complex conjugation, for each $\ell$ there is an index $\ellbar$ 
for which $\overline{H_{\ell}} = H_{\ellbar}$\, ;  possibly $H_{\ell} = H_{\ellbar}$. 
We will say that a component $H_{\ell}$ is a ``short interval'' if it satisfies the following condition: 
\index{short@`short' interval|ii}       
\begin{eqnarray} \label{ShortIntDef} 
\centerline{\text{$H_{\ell} = [a_{\ell},b_{\ell}] \subset \cC_v(\RR) \backslash \fX$
is short relative to $\fX$ and a suitable  }} \\
\centerline{\text{local coordinate function $z_{\ell}$, 
in the sense of Definition \ref{ShortnessDef}.}} \notag 
\end{eqnarray}
  
Recall that a probability vector $\vs \in \cP^m$ is $K_v$-symmetric if $s_j = s_k$ 
\index{$K_v$-symmetric!probability vector}
whenever $\overline{x_{j}} = x_k$.  We will say that a vector 
$\vec{n} = (n, \ldots, n_D) \in \NN^D$ is $K_v$-symmetric 
\index{$K_v$-symmetric!vector}
if $n_{j} = n_{k}$ whenever $\overline{H_j} = H_{k}$ (i.e., when $\jbar = k$).   

\smallskip
Let $\vs \in \cP^m$ be a $K_v$-symmetric probability vector.
\index{$K_v$-symmetric!probability vector}
Let  $\vec{n} \in \NN^D$ be $K_v$-symmetric, and put $N = \sum_{\ell = 1}^D n_{\ell}$.
\index{$K_v$-symmetric!vector}
For each such $\vn$ we will construct an $(\fX,\vs)$-pseudopolynomial 
\index{pseudopolynomial!$(\fX,\vs)$}
$P_{\vn}(z) = P_{(\fX,\vs),\vn}(z)$ whose roots belong to $H$, 
which satisfies $P_{\vn}(z) = P_{\vn}(\zbar)$ for all $z$, 
which has large oscillations on the sets $H_{\ell}$ which are short intervals, 
\index{short@`short' interval}  
and whose normalized logarithm $(-1/N) \log(P_{\vn}(z))$ approximates $u_{\fX,\vs}(z,H)$ outside 
\index{equilibrium potential!$(\fX,\vs)$}
a neighborhood of $H$.  

Most of the roots of $P_{\vn}(z)$ will be roots of the weighted Chebyshev polynomials, 
\index{Chebyshev polynomial!weighted $(\fX,\vs)$}
or weighted Fekete points, for the sets $H_{\ell}$. 
\index{Fekete!points}   
Some care is needed to assure that $P_{\vn}(z) = P_{\vn}(\zbar)$.  

Let $\mu_{\fX,\vs}$ be the equilibrium 
distribution of $H$ relative to $\fX$ and $\vs$.  For each $\ell$ put  
\index{equilibrium distribution!$(\fX,\vs)$}
\begin{equation*}
\widehat{u}_{\ell}(z) \ = \ \int_{H \backslash H_{\ell}} -\log([z,w]_{\fX,\vs}) 
                      \,  d\mu_{\fX,\vs}(w)  
\end{equation*}
and let $W_{\ell}(z) = \exp(-\widehat{u}_{\ell}(z))$.  
Since $\fX$ and $H$ are stable under complex conjugation,
and $\vs$ is $K_v$-symmetric, we have  $\widehat{u}_{\ellbar}(z) = \widehat{u}_{\ell}(\zbar)$ and  
\index{$K_v$-symmetric!probability vector}
$W_{\ellbar}(z) = W_{\ell}(\zbar)$.

If $H_{\ell}$ is a short interval, let 
\index{short@`short' interval}  
\begin{equation} \label{fSnug0} 
\tP_{\ell,(n_{\ell},N)}(z) \ = \ 
\prod_{i=1}^{n_{\ell}} [z,\alpha_{{\ell},i}]_{\fX,\vs} \cdot W_{\ell}(z)^N
\end{equation}
be the weighted Chebyshev pseudopolynomial for $H_{\ell}$ 
\index{pseudopolynomial!weighted Chebyshev}
\index{Chebyshev pseudopolynomial!weighted $(\fX,\vs)$}
with weight $W_{\ell}(z)$.  We are interested in its roots 
$\alpha_{\ell,1}, \ldots, \alpha_{\ell,n_{\ell}}$, which belong to $H_{\ell}$.

If $H_{\ell}$ is not a short interval and $H_\ell \ne H_\ellbar$, 
\index{short@`short' interval}  
let $\alpha_{\ell,1}, \ldots, \alpha_{\ell,n_{\ell}} \in H_{\ell}$ 
be a set of $(n_{\ell},N)$-Fekete points for $H_{\ell}$ relative to the weight $W_{\ell}(z)$,\index{Fekete!points} 
that is, a set of points achieving the maximum value 
\begin{equation*} 
 d_{(n_{\ell},N)}(H_{\ell},W_{\ell}) \ =  \
    \left( \prod^{n_{\ell}}_{\substack{ i, j = 1 \\ i \ne j }} [\alpha_{\ell,i},\alpha_{\ell,j}]_{\fX,\vs} 
             \cdot \prod_{i=1}^{n_{\ell}} W(z_i)^{2N} \right)^{1/N^2} 
\end{equation*} 
in (\ref{AFC34}).
We take the Fekete points $\alpha_{\ellbar,i}$\index{Fekete!points}
for $H_{\ellbar}$ to be the conjugates of the $\alpha_{\ell,i}$ for $H_{\ell}$.

Finally, suppose $H_{\ell}$ is not a short interval, but $H_{\ell} = H_{\ellbar}$.  
\index{short@`short' interval}  
We first show\footnote{The author thanks Will Kazez for this argument.}
\index{Kazez, William} 
that $H_{\ell}$ contains a point $\beta_{\ell}$ fixed by complex conjugation, that is, 
a point in $\cC_v(\RR)$.  Let $S$ be $\cC_v(\CC)$, viewed as a topological surface,
and let $S_0$ be the quotient of $S$ under the action of complex conjugation.  Write $\pi : S \rightarrow S_0$
for the quotient map.   Let $\dot{S} = \cC_v(\CC) \backslash (\cC_v(\RR) \cup \fX)$ 
and put $\dot{S}_0 = \pi(\dot{S}) \subset S_0$.  
Then  $S_0$ is a compact, connected (possibly non-orientable) surface with boundary, 
$\dot{S}_0$ consists of the interior of $S_0$ with a finite number of points removed,
and $\dot{S}$ is a $2$-to-$1$ unramified cover of $\dot{S}_0$.   
Suppose $H_{\ell} \cap \cC_v(\RR) = \phi$. Then $H_{\ell} \subset \dot{S}$.
Choose a point $P \in H_{\ell}$, and let $\Pbar \in H_{\ell}$ be its image under complex conjugation. 
By hypothesis $P \ne \Pbar$.  
Since $H_{\ell}$ is simply connected, and in particular path connected,
\index{simply connected} 
there is a path $\alpha$ from $P$ to $\Pbar$ in $H_{\ell}$.  
Let $\overline{\alpha}$ be the conjugate path; then the concatenation 
$\overline{\alpha} * \alpha$ is a loop in $H_{\ell}$ 
(by a loop, we mean a continuous image of the unit circle).  Since $H_{\ell}$ is simply connected,
$\overline{\alpha} * \alpha$ is homotopic in $H_{\ell}$ to a point.     
Put $\alpha_0 = \pi(\alpha) = \pi(\overline{\alpha}) \subset \dot{S}_0$.  Then 
$\alpha_0*\alpha_0$ is homotopic in $\dot{S}_0$ to a point.  
However, $\pi(P) = \pi(\Pbar)$, so $\alpha_0$ itself is a loop,   
and the fundamental group of a surface with at least one puncture  
is a free group and in particular is torsion-free.
This means $\alpha_0$ is homotopic in $\dot{S}_0$ to a point;  let $\sigma$ be such a homotopy.   
Since $\dot{S}$ is an unramified cover of $\dot{S}_0$, 
we can lift $\sigma$ to a homotopy of loops in $\dot{S}$  
whose initial element is the pre-image $\alpha$ of $\alpha_0$.    
Thus $\alpha$ is a loop, and so $P = \Pbar$, a contradiction.   

Still assuming $H_{\ell}$ not a short interval but $H_{\ell} = H_{\ellbar}$, 
\index{short@`short' interval}   
put $m_{\ell} = \lfloor n_{\ell}/2 \rfloor$, $M = \lfloor N/2 \rfloor$, 
and let $\alpha_{\ell,1}, \ldots, \alpha_{\ell,m_{\ell}} \in H_{\ell}$ 
be a set of $(m_{\ell},M)$-Fekete points for $H_{\ell}$ relative to $W_{\ell}(z)$. 
\index{Fekete!points}
If $n = 2m_{\ell}$ is even, put $\alpha_{\ell,i} = \overline{\alpha_{\ell,i-m_{\ell}}}$ 
for $i = m_\ell+1, \ldots, 2m_{\ell}$.  If $n_{\ell} = 2m_{\ell}+1$ is odd, 
define $\alpha_{\ell,m_{\ell}+1}, \ldots, \alpha_{\ell,2m_{\ell}}$ as above and put 
$\alpha_{\ell,2m_{\ell}+1} = \beta_{\ell}$.
Now define 
\begin{equation} \label{bF3} 
P_{\vec{n}}(z) \ := \ \prod_{{\ell}=1}^D \, \prod_{i=1}^{n_{\ell}} \, 
                    [z,\alpha_{{\ell},i}]_{\fX,\vs} \ .
\end{equation}   
Then $P_{\vec{n}}(z)$ is an $(\fX,\vs)$-pseudopolynomial of degree $N$ for $H$,
\index{pseudopolynomial!$(\fX,\vs)$}
satisfying $P_{\vn}(z) = P_{\vn}(\zbar)$ for all $z$.   
Let 
\begin{equation*}
\cM_{\vn} \ = \ \Big(\|P_{\vn}(z)\|_H\Big)^{1/N} \ . 
\end{equation*}
If there are components $H_{\ell}$ which are short intervals, 
\index{short@`short' interval}  
let $\rho_{\vn} > 0$ be the largest number such that $P_{\vn}(z)$ varies
$n_{\ell}$ times from $(\rho_{\vn})^N$ to $0$ to $(\rho_{\vn})^N$ 
each of those components;  otherwise put $\rho_{\vn} = \cM_{\vn}$. 

Let 
\begin{equation*}
\omega_{\vn} \ = \ \frac{1}{N} \sum_{\ell=1}^D 
             \sum_{i=1}^{n_{\ell}} \delta_{\alpha_{\ell,i}}(z)
\end{equation*} 
be the discrete measure of mass $1$ associated to $P_{\vn}(z)$.  
By construction, it is stable under complex conjugation.  

For each $\ell = 1, \ldots, D$, put $\sigma_{\ell} = \mu_{\fX,\vs}(H_{\ell}) > 0$, 
and let $\vec{\sigma} = (\sigma_1, \ldots, \sigma_D)$.  
Then $\vsigma$ is $K_v$-symmetric, i.e. $\sigma_{\ellbar} = \sigma_{\ell}$ for each $\ell$. 
\index{$K_v$-symmetric!vector}
Given a sequence
of $K_v$-symmetric vectors $\vn_k = (n_{k,1}, \ldots, n_{k,D})\in \NN^D$ for $k = 1, 2, 3, \ldots$,
\index{$K_v$-symmetric!vector}
put $N_k = \sum_{\ell = 1}^D n_{k,\ell}$ for each $k$.  

\begin{proposition} \label{bP2}  Suppose $K_v \cong \RR$.  
Assume that $\fX$ and $H \subset \cC_v(\CC) \backslash \fX$ 
are stable under complex conjugation, and that $H$ is compact 
and has finitely many connected components $H_1, \ldots, H_D$,   
where no $H_{\ell}$ is reduced to a point, and each $H_{\ell}$ is simply connected.
\index{closure of $\cC_v(\CC)$ interior}\index{simply connected}
Let $\vs \in \cP^m$ be a $K_v$-symmetric probability vector. 
\index{$K_v$-symmetric!probability vector}  
Then for each $K_v$-symmetric vector $\vn \in \NN^D$, 
\index{$K_v$-symmetric!vector}
the $(\fX,\vs)$-pseudopolynomial $P_{\vn}(z)$
\index{pseudopolynomial!$(\fX,\vs)$}
has all its roots belong to $H$, with has exactly $n_{\ell}$ 
roots in each $H_{\ell}$. 
It satisfies $P_{\vn}(z) = P_{\vn}(\zbar)$ for all $z \in \cC_v(\CC)$.  
In addition, for each $H_{\ell} = [a_\ell,b_\ell]$ which is a short interval, 
\index{short@`short' interval}  
the roots of $P_{\vn}(z)$ in $H_{\ell}$ are distinct,  
and $P_{\vn}(z)$ varies\, 
$n_\ell$ times from $\rho_{\vn}^N$ to $0$ to $\rho_{\vn}^N$ 
on $H_{\ell}$, where  $N = \deg(P_{\vn_k}) = \sum_{\ell=1}^d n_{\ell}$.   

For any sequence of $K_v$-symmetric vectors\, $\vn_k \in \NN^D$ for which $N_k \rightarrow \infty$
\index{$K_v$-symmetric!vector}
and $\vn_k/N_k \rightarrow \vec{\sigma}$,   

$(1)$ The discrete measures $\omega_{\vn_k}$ associated to the $P_{\vn_k}$
converge weakly to the equilibrium distribution $\mu_{\fX,\vs}$ of $H\, ;$  
\index{equilibrium distribution!$(\fX,\vs)$}

$(2)$ If $U$ is any open neighborhood of $H$, then as $k \rightarrow \infty$
the functions $- \frac{1}{N_k} \log(|P_{\vn_k}(z)|)$ converge uniformly to the equilibrium potential 
\index{equilibrium potential!$(\fX,\vs)$}
$u_{\fX,\vs}(z,H)$ on $\cC_v(\CC) \backslash (U \cup \fX)$.

$(3)$ $\lim_{k \rightarrow \infty} \cM_{\vn_k} 
          = \lim_{k \rightarrow \infty} \rho_{\vn_k} = \gamma_{\fX,\vs}(H) \, .$

\end{proposition}

\begin{proof} By construction, $P_{\vn}(z) = P_{\vn}(\zbar)$ for all $z \in \cC_v(\CC)$,
the roots of $P_{\vn}(z)$ belong to $H$, and $P_{\vn}$ has exactly $n_{\ell}$ roots in $H_{\ell}$.  
For each $H_{\ell}$ which is a short interval,  
\index{short@`short' interval}  
the roots of $P_{\vn}(z)$ in $H_{\ell}$ are distinct by Proposition \ref{bP1}, 
and $P_{\vn}(z)$ varies $n_{\ell}$ times from $\rho_{\vn}$ to $0$ to $\rho_{\vn}$ on $H_{\ell}$
by the definition of $\rho_{\vn}$.    
 
Now consider a sequence $\{\vn_k\}_{k \in \NN}$ with $N_k \rightarrow \infty$ 
and $n_{k,\ell}/N_k \rightarrow \sigma_{\ell}$ for each $\ell$. Put   
\begin{equation*} 
\omega_{\vec{n_k},{\ell}} 
\ = \ \sum_{i=1}^{n_{\ell}} \frac{1}{N_k} \delta_{\alpha_{{\ell},i}}(z) \ . 
\end{equation*}
By  Corollary \ref{ACF17} and Theorem \ref{ATF19},  for each $\ell$, as $k \rightarrow \infty$,
the measures $\omega_{\vn_k,{\ell}}$ converge weakly to $\mu_{\fX,\vs}|_{H_{\ell}}$.
(The presence of the points $\beta_{\ell}$ do not affect this.)  
Hence the $\omega_{\vn_k} = \sum_{\ell=1}^d \omega_{\vn_k,\ell}$ 
converge weakly to $\mu_{\fX,\vs}$.  

This implies that outside any neighborhood $U$ of $H$, 
the potential functions $u_{\fX,\vs}(z,\omega_{\vn_k})$ converge uniformly to $u_{\fX,\vs}(z,H)$.
\index{equilibrium potential!$(\fX,\vs)$}
Indeed, since $\partial U$ and $H$ are compact and disjoint,  
$-\log([z,w]_{\fX,\vs})$ is uniformly continuous on $\partial U \times H$.   
Since the $\omega_{\vn_k}$ converge weakly to 
$\mu_{\fX,\vs}$, as $k \rightarrow \infty$ 
\begin{equation*}
u_{\fX,\vs}(z,\omega_{\vn_k}) \ = \ \int_{H} -\log([z,w]_{\fX,\vs}) \, d\omega_{\vn_k}(w)
\ = \ -\frac{1}{N_k} \log(P_{\vn_k}(z)) 
\end{equation*}
converges uniformly to\index{equilibrium potential!$(\fX,\vs)$} 
$u_{\fX,\vs}(z,H) = \int_H -\log([z,w]_{\fX,\vs}) \, d\mu_{\fX,\vs}(w)$
on $\partial U$.  However, for each $k$, by Theorem \ref{ATE10B} the function 
$u_{\fX,\vs}(z,\omega_{\vn_k}) - u_{\fX,\vs}(z,H)$ on $\cC_v(\CC) \backslash (U \cup \fX)$
extends to a function harmonic on a neighborhood of each $x_i \in \fX$.  
By the Maximum principle, $u_{\fX,\vs}(z,\omega_{\vn_k})= -\frac{1}{N_k} \log(P_{\vn_k}(z))$ 
converges uniformly to $u_{\fX,\vs}(z,H)$
on $\cC_v(\CC) \backslash (U \cup \fX)$. 
\index{Maximum principle!for harmonic functions}

\smallskip 
We will next show that $\lim_{k \rightarrow \infty} \cM_{n_k} = \gamma_{\fX,\vs}(H)$,
or equivalently, that 
\begin{equation} \label{FSnug}  
\lim_{k \rightarrow \infty} -\log(\cM_{\vn_k}) \ = \ V_{\fX,\vs}(H) \ .
\end{equation}
\index{Robin constant!$(\fX,\vs)$}
Recall that a Hausdorff space $X$ 
is simply connected if and only if any two points are joined by an arc\index{simply connected} 
\index{arc} 
and every loop in $X$ is homotopic in $X$ to a point.  Since each component of $H$ is simply connected and no
\index{simply connected}
component is a point, it follows from Theorem \ref{ATE10B} and Proposition \ref{ACE11} that 
the potential function $u_{\fX,\vs}(z,H) $ is continuous,
\index{equilibrium potential!$(\fX,\vs)$}  
with $u_{\fX,\vs}(z,H) = V_{\fX,\vs}(H)$ for all $z \in H$.
\index{Robin constant!$(\fX,\vs)$}
Since $\cC_v(\CC) \backslash H$ is connected, in addition $u_{\fX,\vs}(z,H) < V_{\fX,\vs}(H)$ 
for all $z \notin H$.

We first claim that $-\log(\cM_{n_k}) \le V_{\fX,\vs}(H)$ for each $k$.   
\index{Robin constant!$(\fX,\vs)$}
For each $z \in H$, $u_{\fX,\vs}(z,\omega_{\vn_k})$ satisfies 
\begin{equation*} 
u_{\fX,\vs}(z,\omega_{\vn_k}) \ = \ -\frac{1}{N_k} \log(P_{\vn_k}(z)) \ \ge \ -\log(\cM_{\vn_k}) \ .
\end{equation*} 
\index{equilibrium potential!$(\fX,\vs)$}
By Fubini-Tonelli, it follows that 
\index{Fubini-Tonelli theorem}
\begin{equation*} 
V_{\fX,\vs}(H) \ = \ \int_H u_{\fX,\vs}(z,H) \, d\omega_{\vn_k}(z) 
               \ = \ \int_H u_{\fX,\vs}(z,\omega_{\vn_k}) \, d\mu_{\fX,\vs}(z) \ \ge \ -\log(\cM_{\vn_k})
\end{equation*} 
\index{Robin constant!$(\fX,\vs)$}
as asserted.

Now fix $\varepsilon > 0$, and put 
\begin{equation*}
U_\varepsilon \ = \ \{z \in \cC_v(\CC) : u_{\fX,\vs}(z,H) > V_{\fX,\vs}(H) - \varepsilon \} \ . 
\end{equation*}
\index{equilibrium potential!$(\fX,\vs)$}
\index{Robin constant!$(\fX,\vs)$}
By the properties of $u_{\fX,\vs}(z,H)$ noted above, $U_\varepsilon$ is open and 
$H \subset U_\varepsilon $.  Since the $\omega_{\vn_k}$ 
converge weakly to $\mu_{\fX,\vs}$, the functions $u_{\fX,\vs}(z,\omega_{\vn_k})$ converge 
uniformly to $u_{\fX,\vs}(z,H)$ on $\cC_v(\CC) \backslash (U_\varepsilon \cup \fX)$. 
\index{equilibrium potential!$(\fX,\vs)$}
Thus, there is a $k_0$ such that for all $k \ge k_0$ and 
all $z \in \cC_v(\CC) \backslash (U_\varepsilon \cup \fX)$ we have 
\begin{equation*} 
|u_{\fX,\vs}(z,H) - u_{\fX,\vs}(z,\omega_{\vn_k})| \ < \ \varepsilon  
\end{equation*} 
\index{equilibrium potential!$(\fX,\vs)$}
This means that for $k \ge k_0$ and $z \in \partial U_\varepsilon$,
\begin{equation} \label{f7OmegaIneq}
u_{\fX,\vs}(z,\omega_{\vn_k}) \ > \ V_{\fX,\vs}(z) - 2\varepsilon \ .
\end{equation}\index{Robin constant!$(\fX,\vs)$}
However, $u_{\fX,\vs}(z,\omega_{\vn_k})$ is superharmonic\index{superharmonic} on $U_\varepsilon$
\index{equilibrium potential!is superharmonic}
\index{Maximum principle!for superharmonic functions}
so by the Maximum principle for superharmonic functions,\index{superharmonic}
(\ref{f7OmegaIneq}) holds throughout $U_\varepsilon$.
In particular, $-\log(\cM_{\vn_k}) > V_{\fX,\vs}(z) - 2\varepsilon$.  
Since $\varepsilon > 0$ is arbitrary, (\ref{FSnug}) holds.\index{Robin constant!$(\fX,\vs)$}



Finally, we show that $\lim_{k \rightarrow \infty} \rho_{\vn_k} = \gamma_{\fX,\vs}(H)$. 
Recall that $W_{\ell}(z) = \exp(-\widehat{u}_{\ell}(z))$ where 
\begin{equation*}
\widehat{u}_{\ell}(z) \ = \ \int_{H \backslash H_{\ell}} -\log([z,w]_{\fX,\vs}) \,  d\mu_{\fX,\vs}(w) \ .
\end{equation*}  
Write $\tP_{k,\ell}(z)$ for the weighted 
Chebyshev pseudopolynomial $\tP_{\ell,(n_{k,\ell},N_k)}(z,W_{\ell})$ 
as in (\ref{fSnug0}),\index{pseudopolynomial!weighted Chebyshev}\index{Chebyshev pseudopolynomial!weighted $(\fX,\vs)$}
and write $\CH^*_{k,\ell}$ for the weighted Chebyshev constant\index{Chebyshev constant!weighted}  
$\CH^*_{(n_{k,\ell},N_k)}(H_{\ell},W_{\ell}) = \|\tP_{k,\ell}\|_{H_{\ell}}^{1/N_k}$.
Also, for each $k$ and $\ell$, put 
\begin{equation*} 
\widehat{u}_{k,\ell}(z) 
      \ := \ \int_{H \backslash H_{\ell}} -\log([z,w]_{\fX,\vs}) \,  d\omega_{\vec{n_k}}(w)  \ .
\end{equation*}
By Theorem \ref{ATB1} and Theorem \ref{ATF18}, as $k \rightarrow \infty$  
and $\vn_k/N_k \rightarrow \vec{\sigma}$, the $\CH^*_{k,\ell}$ 
converge to $\CH^*_{\sigma_{\ell}}(H_{\ell},W_{\ell}) = \gamma_{\fX,\vs}(H)$.  

Fix $\varepsilon > 0$.  Since the measures $\omega_{\vn_k}$ converge weakly to $\mu_{\fX,\vs}$, if $k$ is 
sufficiently large, then by Theorem \ref{ATF18}, for each $\ell$     
\begin{equation}
|\widehat{u}_{k,\ell}(z) - \widehat{u}_{\ell}(z)| 
       <  \varepsilon \quad \text{on \ $H_{\ell}$}  
                                        \label{bF4}
\end{equation}
and 
\begin{equation}
e^{-\varepsilon} \cdot \gamma_{\fX,\vs}(H) \ < \ \CH^*_{k,\ell}  \ .  \label{bF5}
\end{equation}
Put $\rho = e^{-2\varepsilon} \gamma_{\fX,\vs}(H)$. 
  
Suppose $H_{\ell}$ is a short interval.   
\index{short@`short' interval}  
Since $\|\tP_{k,\ell}\|_{H_{\ell}} = (\CH^*_{k,\ell})^{N_k}$,  
Proposition \ref{bP1} shows that $\tP_{k,\ell}(z)$ oscillates $n_{k,\ell}$ times from 
$(\CH^*_{k,\ell})^{N_k}$ to $0$ to $(\CH^*_{k,\ell})^{N_k}$ on $H_{\ell}$.  
At each point in $H_{\ell}$ where 
$\tP_{k,\ell}(z) = (\CH^*_{\ell,\vn})^{N_k}$, we have  
\begin{equation*}
P_{\vn_k}(z) \ > \ (e^{-\varepsilon}  \CH^*_{k,\ell})^{N_k} \ > \ \rho^{N_k}  
\end{equation*} 
and at each point where $\tP_{k,\ell}(z) = 0$, also $P_{\vn_k}(z) = 0$. 
 
Thus $P_{\vn_k}(z)$ oscillates $n_{k,\ell}$ times from 
$\rho^{N_k}$ to $0$ to $\rho^{N_k}$ on each $H_{\ell}$ which is a short interval,
\index{short@`short' interval}   
so $\rho_{\vn_k} \ge \rho = e^{-2\varepsilon} \gamma_{\fX,\vs}(H)$.  
Trivially $\rho_{\vn_k} \le \cM_{\vn_k}$, 
so since $\lim_{k \rightarrow \infty} \cM_{\vn_k} = \gamma_{\fX,\vs}(H)$, and since $\varepsilon > 0$
is arbitrary, assertion (3) in the Proposition follows.
\end{proof} 
 
\noindent{\bf Remark.}  In constructing the functions $P_{\vn}(z)$, 
the use of the Chebyshev points for the components which are short intervals is essential.
\index{short@`short' interval}  
\index{Chebyshev points}   
However, the use of the Fekete points for the other components is
\index{Fekete!points}  
a matter of convenience;  all that is needed is that the associated discrete measures be stable under 
complex conjugation, and converge weakly to $\mu_{\fX,\vs}$.  

 On some curves $\cC_v/\RR$ there exist connected sets 
$H \subset \cC_v(\CC) \backslash \cC_v(\RR)$ which are stable under complex conjugation, 
but are not simply connected.\index{simply connected}  
For example, consider an elliptic curve $\cE/\RR$ 
such as the one defined by $y^2 = x^3-1$, 
for which the real locus $\cE(\RR)$ is 
homeomorphic to a loop.  Such a curve which has two real $2$-torsion points $P_0 = O$ and  
$P_1$, and two conjugate non-real $2$-torsion points $P_2$, $P_3$.  
The set $H = P_2 + \cE(\RR)$  
(where addition is under the group law on $\cE$) is stable under complex conjugation, 
and complex conjugation acts on $H$ via translation by $P_1 = P_3-P_2$, so it has no fixed points in $H$.

      In Proposition \ref{bP2}, one can replace the hypothesis that the components $H_{\ell}$ 
are simply connected\index{simply connected} with two additional hypotheses:  first, that
each connected component of $\cC_v(\CC) \backslash H$ must contain at least one $x_i$ 
for which $s_i > 0$; and second, that for each component $H_{\ell}$ which is not a short interval but 
\index{short@`short' interval}  
satisfies $H_{\ell} = H_{\ellbar}$, the number $n_{\ell}$ is even.

\smallskip

\vskip .1 in
We can reformulate Proposition \ref{bP2} in a useful way, as follows.  
For each $\vn$, put
\begin{equation} \label{FQ1}
Q_{\vn}(z) \ = \ \frac{1}{\cM_{\vn}^N} P_{\vn}(z) \ , 
\end{equation}
and put $R_{\vn} = \rho_{\vn}/\cM_{\vn} \le 1$.  Then $\|Q_{\vn}\|_H = 1$, 
$Q_{\vn}(z) = Q_{\vn}(\zbar)$ for all $z$, 
and $Q_{\vn}(z)$ varies $N$ times from $R_{\vn}^N$ to $0$ to $R_{\vn}^N$ 
on each $H_{\ell}$ which is a short interval.    
\index{short@`short' interval}  

Recall that the Green's function of $H$ is 
$G_{\fX,\vs}(z,H) = V_{\fX,\vs}(H) - u_{\fX,\vs}(z,H)$. 
\index{Robin constant!$(\fX,\vs)$} 
\index{Green's function!$(\fX,\vs)$} 
\index{equilibrium potential!$(\fX,\vs)$}
As in \S\ref{Chap5}.\ref{ArchApproxThmSection}, 
define the logarithmic leading coefficient of $H$ at $x_i$ by\index{logarithmic leading coefficients}
\begin{eqnarray} 
\Lambda_{x_i}(H,\vs) & = & 
\lim_{z \rightarrow x_i} G_{\fX,\vs}(z;H) + s_i \log(|g_{x_i}(z)|) 
\label{FQ3} \\
& = & V_{x_i}(H) + \sum_{j \ne i} G(x_i,x_j;H) \ . \notag 
\end{eqnarray}  
\index{Robin constant}
\index{Green's function}
Likewise, as in \S\ref{Chap5}.\ref{OutlineSection}, 
for each $x_i \in \fX$,  
define the logarithmic leading coefficient of $Q_{\vn}(z)$ at $x_i$ by 
\index{logarithmic leading coefficients!of $Q_{\vn}(z)$}
\begin{equation} \label{FQ2}
\Lambda_{x_i}(Q_{\vn},\vs) \ = \ 
\lim_{z \rightarrow x_i} \frac{1}{N} \log(Q_{\vn}(z)) + s_i \log(|g_{x_i}(z)|) \ .
\end{equation}  
Recalling that for any sequence $\{\vn_k\}_{k \in \NN}$ with $N_k \rightarrow \infty$
and $\vn_k/N_k \rightarrow \vec{\sigma}$ we have 
\begin{equation*}
\lim_{k \rightarrow \infty} \cM_{\vn_k} \ = \ \gamma_{\fX,\vs}(H) 
\ = \ e^{-V_{\fX,\vs}(H)} \ ,  
\end{equation*} 
\index{Robin constant!$(\fX,\vs)$}
by Proposition \ref{bP2} we have:    

\begin{theorem} \label{bT3} Suppose $K_v \cong \RR$.  
Assume that $\fX$ and $H \subset \cC_v(\CC) \backslash \fX$ are stable under complex conjugation,  
that $H$ is compact, 
and that $H$ has finitely many connected components $H_1, \ldots, H_D$,   
where no $H_{\ell}$ is reduced to a point and each $H_{\ell}$ is simply connected.\index{simply connected}  
  
Fix a $K_v$-symmetric probability vector $\vs \in \cP^m$.
\index{$K_v$-symmetric!probability vector}  
For each $\ell = 1, \ldots, D$ put $\sigma_{\ell} = \mu_{\fX,\vs}(H_{\ell})$ 
and let $\vsigma = (\sigma_1, \ldots, \sigma_D)$. 
For each $K_v$-symmetric vector $\vn \in \NN^D$ 
\index{$K_v$-symmetric!vector}
write $N = N_{\vn_k} = \sum_{\ell} n_{\ell}$.
Then the collection of $(\fX,\vs)$-pseudopolynomials 
\index{pseudopolynomial!$(\fX,\vs)$}
\begin{equation*}
Q_{\vn}(z) \ = \ \frac{1}{\cM_{\vn}^N} \prod_{{\ell}=1}^D \, \prod_{i=1}^{n_{\ell}} \, 
                    [z,\alpha_{{\ell},i}]_{\fX,\vs}  \ ,
\end{equation*}
constructed above for $\vn \in \NN^D$, has the following properties: 

\vskip .03 in
$(1)$ For each $\vn$, $Q_{\vn}$ satisfies $\|Q_{\vn}\|_H = 1$,  
with $Q_{\vn}(z) = Q_{\vn}(\zbar)$ for all $z \in \cC_v(\CC)$.
The roots of $Q_{\vn}$ all belong to $H$, with $n_{\ell}$ roots in each $H_{\ell}$.
Put $R_{\vn} = \rho_{\vn}/\cM_{\vn}$, so $0 < R_{\vn} \le 1$.  
For each $H_{\ell}$ which is a short interval 
\index{short@`short' interval}  
in the sense of $(\ref{ShortIntDef})$, 
the roots of $Q_{\vn}$ in $H_{\ell}$ are distinct, 
and $Q_{\vn}$ varies $n_{\ell}$ times from $R_{\vn}^N$ to $0$ to $R_{\vn}^N$ 
on $H_{\ell}$. 

\vskip .03 in
$(2)$ Let  $\{\vn_k\}_{k \in \NN}$ be any sequence  
with $N_{\vn_k} \rightarrow \infty$ and $\vn_k/N_{\vn_k} \rightarrow \vec{\sigma}$.  
Then $\lim_{k \rightarrow \infty} R_{\vn_k}  = 1$, and 
the discrete measures $\omega_{\vn_k}$ associated to the $Q_{\vn_k}$
converge weakly to the equilibrium distribution $\mu_{\fX,\vs}$ of $H$.  
\index{equilibrium distribution!$(\fX,\vs)$} 
For each neighborhood $U$ of $H$,     
the functions $\frac{1}{N_{\vn_k}}\log(Q_{\vn_k}(z))$ 
converge uniformly to $G_{\fX,\vs}(z,H)$ 
\index{Green's function!$(\fX,\vs)$}
on $\cC_v(\CC_v) \backslash (U \cup \fX)$,
and for each $x_i \in \fX$,   
\begin{equation*} 
\lim_{k \rightarrow \infty} \Lambda_{x_i}(Q_{\vn_k},\vs) 
\ = \ \Lambda_{x_i}(H,\vs) \ .    
\end{equation*}  
\end{theorem}  
 

%% file: NewFSZAppC.tex
\chapter{The Universal Function} \label{AppC}  

\index{universal function|(} 
In this Appendix, we construct a parametrization of rational functions
of degree $d$ on a curve by their zeros and poles.  
\index{universal function|ii}
We then establish a $v$-adic bound for how much 
a rational function changes when its zeros and poles are moved slightly.  
This is used in  \S\ref{Chap11}.\ref{NonArchPatchingProof}, in Step 4 of 
the patching process in the nonarchimedean compact case, 
\index{patching argument!for nonarchimedean $K_v$-simple sets}
with $d = \max(1,2g)$, in moving the roots of the partially patched 
function away from each other.  


\vskip .1 in
Let $F$ be a field, and let $\Fbar$ be a fixed algebraic closure of $F$.   
When $\cC = \PP^1/F$, let $z, w, p, q$ be  independent variables 
and consider the crossratio\index{crossratio}
\begin{equation*}
\chi(z,w;p,q) \ = \ \frac{(z-p)}{(z-q)} \frac{(w-q)}{(w-p)} \ ,
\end{equation*}    
which extends to a rational function on  $(\PP^1)^4$.  
Now specialize $p, q$ to $\PP^1(\Fbar)$, and take 
$w \in \PP^1(\Fbar)$ distinct from $p, q$.  Then $f_{w,p,q}(z) = \chi(z,w;p,q)$ 
is a rational function on $\PP^1$ with divisor $(p)-(q)$,
normalized by the condition that $f_{w,p,q}(w) = 1$.  

More generally, for arbitrary  
$p_1, \ldots, p_d, q_1, \ldots, q_d \in \PP^1(\Fbar)$, 
if $w \in \PP^1(\Fbar)$ is distinct from the $p_i$ and $q_i$, then 
\begin{equation} \label{FP1Case} 
f_{w,\vp,\vq}(z) \ := \  f(z,w;\vp,\vq) \ := \ \prod_{i=1}^d \chi(z,w;p_i,q_i)
\end{equation}
is the unique rational function on $\PP^1$ with divisor 
$\div_{\cC}(f_{w,\vp,\vq}) = \sum (p_i) - \sum (q_i)$,
for which $f_{w,\vp,\vq}(w) = 1$.  Conversely, 
for any nonconstant rational function $h(z) \in \Fbar(z)$, 
there is a point $w \in \PP^1(\Fbar)$ where $h(w) = 1$.  
In this way, we obtain a parametrization of 
all rational functions of degree $d$ on $\PP^1$  
by means of their zeros and poles and a normalizing point.    

\vskip .1 in
We will now show the existence of similar parametrizations for arbitrary curves $\cC/F$.
 
\vskip .1 in
Let $\cC/F$  be a smooth, projective, geometrically integral curve of genus $g  > 0$. 
Given $D \in \Div^d_{\cC/F}(\Fbar)$, let $[D] \in \Pic^d_{\cC/F}(\Fbar)$
be the linear equivalence class of $D$.
Likewise, given $\vp = (p_1, \ldots, p_d) \in \cC(\Fbar)^d$,
write $[\vp] \in \Pic^d_{\cC/F}(\Fbar)$ for the linear equivalence class of
$\sum (p_i)$.  Let $\Jac(\cC)/F$ be the Jacobian of $\cC$, 
\index{Jacobian variety}  
\label{`SymbolIndexJacobian'}
and let $\Phi : \cC^d \times \cC^d \rightarrow \Jac(\cC)$ 
be the $F$-rational map defined by $\Phi(\vp,\vq) = [\vp]-[\vq]$ 
for $\vp, \vq \in \cC^d(\Fbar)$.  Put $Y = \Phi^{-1}(0)$.
Then $Y$ is the $F$-rational subvariety of $\cC^d \times \cC^d$ for which 
\begin{equation*}
Y(\Fbar) \ = \ \{ (\vp,\vq) \in \cC(\Fbar)^d \times \cC(\Fbar)^d :   
                     \text{$\sum (p_i) - \sum (q_i)$ is principal} \} \ .
\end{equation*} 
We will construct a ``universal rational function'' 
\index{universal function}
$f(z,w;\vp,\vq)$ on $\cC \times \cC \times Y$ 
which parametrizes normalized rational functions of degree $d$ 
on $\cC$ in the sense above.  
For this it will be necessary to assume that $d \ge 2g-1$.  

We will then specialize to the case where $F$ is a nonarchimedean local field $K_v$,
and study the continuity properties of $f(z,w;\vp,\vq)$.    

\begin{theorem}  \label{UniversalFcn} 

Let $F$ be a field, and let $\cC/F$ be a smooth, projective,
geometrically integral curve of genus $g \ge 0$.
Fix $d \ge \max(1, 2g-1)$, and let 
$Y$ be the $F$-rational subvariety of $\cC^d \times \cC^d$ for which 
\begin{equation*} 
Y(\Fbar) \ = \ \{ (\vp,\vq) \in \cC^d(\Fbar) \times \cC^d(\Fbar) : 
 \sum (p_i) \sim \sum (q_i) \} \ .
\end{equation*}  
Then there is an $F$-rational function\footnote{The author thanks Robert Varley for pointing out 
\index{Varley, Robert} 
that $f(z,w;\vp,\vq)$ can be defined even when $\cC(F) = \phi$.} 
$f(z,w;\vp,\vq)$ on $\cC \times \cC \times Y$
uniquely defined by the property that for each $(\vp,\vq) \in Y(\Fbar)$ 
and each $w \in \cC(\Fbar)$ distinct from the $p_i$, $q_i$, 
the function $f_{w,\vp,\vq}(z) := f(z,w;\vp,\vq) \in \Fbar(\cC)$
satisfies $\div_{\cC}(f_{w,\vp,\vq}) = \sum (p_i) - \sum (q_i)$ 
and $f_{w,\vp,\vq}(w) = 1$.  

Moreover, if $F \subseteq M \subseteq \Fbar$ is an extension such that 
$w \in \cC(M)$ and $\sum (p_i) - \sum (q_i)$ is rational over $M$,
then $f_{w,\vp,\vq}(z) \in M(\cC)$.  
 \end{theorem} 

\vskip .1 in 
The proof uses the theory of the Picard scheme, due to Grothendieck and Mumford. 
\index{Picard scheme} 
\index{Grothendieck, Alexander}
\index{Mumford, David} 
The part of the theory we need goes back to Weil and Matsusaka. 
\index{Weil, Andr\'e}
\index{Matsusaka, Teruhisa}  

\vskip .1 in
A modern reference for this is Kleiman (\cite{Klei});\index{Kleiman, Steven}  
see also Milne (\cite{Milne}).\index{Milne, James} 
We follow Kleiman's notation.\index{Kleiman, Steven} 
Given a separated map of locally Noetherian schemes $F : X \rightarrow S$, 
let $\Pic_{X/S}$ be the relative Picard functor, defined by
\begin{equation*}
\Pic_{X/S}(T) \ = \ \Pic(X_T)/\Pic(T)
\end{equation*}
for any locally Noetherian scheme $T/S$ (\cite{Klei}, Definition 2.2);  
here $X_T = X \times_S T$. 
 
If $F : X \rightarrow S$ is projective and flat, with reduced, connected geometric
fibres, and if $X(S) \ne \phi$, then $\Pic_{X/S}$
is represented by a commutative group scheme $\BPic_{X/S}$ 
which is separated and locally of finite
type over $S$ (see \cite{Klei}, Theorems 2.5 and 4.8, and Exercise 3.11).  

The scheme $\BPic_{X/S}$ commutes with base change:  
for any locally Noetherian scheme $S^{\prime}/S$,  
$\BPic_{X_{S^{\prime}}/S^{\prime}}$ exists and 
equals $\BPic_{X/S} \times_S S^{\prime}$ (\cite{Klei}, Exercise 4.4).  
Points of $\BPic_{X/S}$ correspond
in a natural bijective way to classes of invertible sheaves on the 
fibres of $X/S$ (\cite{Klei}, Exercise 4.5). 
There is an invertible sheaf $\cP$ on $X \times \BPic_{X/S}$,  
called a Poincar\'e sheaf,\index{Poincar\'e sheaf} such that for any locally Noetherian scheme $T/S$,
and any invertible sheaf $\cL$ on $X_T$, there exists a unique $S$-map 
$h : T \rightarrow \BPic_{X/S}$ such that for some invertible sheaf $\cN$ on $T$,
\begin{equation*}
\cL \ \cong \ (1 \times h)^*\cP \otimes f_T^*(\cN)
\end{equation*} 
(see \cite{Klei}, Exercise 4.3).  In general, the Poincar\'e sheaf\index{Poincar\'e sheaf} is not unique. 

Let $\Div_{X/S}$ be the functor defined by 
\begin{equation*} 
\Div_{X/S}(T)  \ = \ \{ \text{relative effective divisors on $X_T/T$} \} \ ;
\end{equation*} 
see (\cite{Klei}, \S3) for details.  
It is represented by an open subscheme $\BDiv_{X/S}$ of the Hilbert scheme
\index{Hilbert scheme} 
$\BHilb_{X/S}$ (\cite{Klei}, Theorem 3.7). 
By (\cite{Klei}, \S3.10 and Exercise 4.7), 
there is a coherent $\cO_{\BPic_{X/S}}$-module $\cQ$ for which  
$\BDiv_{X/S} \cong \BP(\cQ)$.
There is a natural map of functors 
$A_{X/S}(T) : \Div_{X/S}(T) \rightarrow \Pic_{X/S}(T)$,
called the ``Abel map'',\index{Abel map|ii} which sends a relative 
effective divisor $D$ on $X_T/T$ to the sheaf $\cO_{X_T}(D)$,   
and there is a corresponding Abel map of $S$-schemes 
$\BA_{X/S} : \BDiv_{X/S} \rightarrow \BPic_{X/S}$.  

\medskip
If $S = \Spec(F)$ where $F$ is a field, and $X = \cC$ is a smooth, projective, 
geometrically integral curve of genus $g > 0$ with $\cC(F) \ne \phi$, 
the connected component of the identity, $\BPic_{X/S}^0$,  
is an abelian variety of dimension $g$ (\cite{Klei}, Exercise 5.23)  
which is $F$-isomorphic to the Jacobian $\Jac(\cC)/F$. 
\index{Jacobian variety}   
%
In this setting, 
we will write $\BPic_{\cC/F}$ for $\BPic_{X/S}$
and $\BDiv_{\cC/F}$ for $\BDiv_{X/S}$.  
Here $\BPic_{\cC/F}$ is a disjoint union of open subschemes $\BPic^d_{\cC/F}$ 
representing invertible sheaves $\cL$ of degree $d$,
and $\BPic^d_{\cC/F}$ is a $\BPic^0_{\cC/F}$-torsor for each $d$ 
(\cite{Klei}, Exercise 6.21). Similarly $\BDiv_{\cC/F}$ 
is a disjoint union of open subschemes $\BDiv^d_{\cC/F}$ for $d \ge 1$, 
and for each $d$ the Abel map\index{Abel map} takes $\BDiv^d_{\cC/F}$ to $\BPic^d_{\cC/F}$.       
Let $\cC^d$ be the $d$-fold product of $\cC$ with itself,
and let $\Sym^{(d)}(\cC)$ be the $d$-fold symmetric product;  it is smooth, 
since $\cC$ is a curve (\cite{Milne}, Proposition 3.2, p.94). 
\index{Milne, James} 
For each $d \ge 1$ there is a canonical surjective morphism 
$\alpha_d : \cC^d \rightarrow \BDiv^d_{\cC/F}$ given on $\Fbar$-points
by $\alpha_d(p_1, \ldots, p_d) = \sum (p_i)$.  
It induces an isomorphism $\Sym^d(\cC) \cong \BDiv^d_{\cC/F}$ 
(see \cite{Klei}, Remark 3.9).

\vskip .1 in
\begin{proof}[Proof of Theorem \ref{UniversalFcn}]  
We first carry out the construction over $\Fbar$.  
Make a base change to $\Fbar$, and until further notice, 
replace $\cC$ by $\cCbar = \cC_{\Fbar}$ and $Y$ by $\Ybar = Y_{\Fbar}$.   
We will denote the function in the theorem constructed over $\Fbar$ by $\fbar(z,w;\vp,\vq)$.  

Clearly $\cCbar(\Fbar) \ne 0$.  
If $\cC$ has genus $0$, then $\cCbar \cong \PP^1/\Fbar$, 
so we can construct $\fbar(z,w;\vp,\vq)$ by using
the cross-ratio as in (\ref{FP1Case}).  
In the argument below we will assume that $g > 0$. 

Note that for each $(\vp,\vq) \in \Ybar(\Fbar)$, 
and each $w \in \cCbar(\Fbar)$ distinct from the $p_i, q_i$, 
there is a unique function $\fbar_{w,\vp,\vq}(z) \in \Fbar(\cCbar)$ 
for which $\div_{\cCbar}(\fbar_{w,\vp,\vq}) = \sum (p_i) - \sum (q_i)$ 
and $\fbar_{w,\vp,\vq}(w) = 1$.  We must show that these functions 
glue to give a globally defined rational function 
on $\cCbar \times \cCbar \times \Ybar$.  

Put $\cZbar = \cCbar \times \BPic_{\cCbar/\Fbar}^d$.
If $p_2 : \cZbar \rightarrow \BPic_{\cCbar/\Fbar}^d$
is the projection on the second factor, 
and $y = [D] \in \BPic_{\cCbar/\Fbar}^d(\Fbar)$,
write $\cZbar_y$ for the fibre $p_2^*(y) \cong \cCbar$.
Let $\cPbar$ be a Poincar\'e\index{Poincar\'e sheaf} sheaf on $\cZbar$.  Then if
$i_y : \cCbar \hookrightarrow \cZbar$ is the inclusion
$i_y(P) = (P,y) \in \cZbar$, we have 
\begin{equation*}
\cPbar_y \ := \ i_y^*(\cPbar) \ \cong \ \cPbar|_{\cZbar_{y}} 
\ \cong \ \cO_{\cCbar}(D) \ .
\end{equation*}

Since $d \ge 2g-1$, the Riemann-Roch theorem shows that 
\index{Riemann-Roch theorem}
$\dim(H^0(\cCbar,\cO_{\cCbar}(D))) = d-g+1$ for all 
$D \in \Div_{\cCbar/\Fbar}^d(\Fbar)$.  
The projection $p_2 : \cZbar \rightarrow \BPic_{\cCbar/\Fbar}^d$ 
is a flat, projective morphism of
Noetherian schemes, and $\cPbar$ is flat over $\cO_{\cZbar}$,
hence also flat over $\cO_{\BPic_{\cCbar/\Fbar}^d}$.
By Grauert's Theorem\index{Grauert's theorem} (see [H], p.288), $(p_2)_*(\cPbar)$
is locally free of rank $d-g+1$ over $\cO_{\BPic_{\cCbar/\Fbar}^d}$,
and for each $y = [D] \in \BPic_{\cCbar/\Fbar}^d(\Fbar)$, the natural map
\begin{equation} \label{BFP1}
(p_2)_*(\cPbar) \otimes k(y) 
     \ \rightarrow \ H^0(\cZbar_y,\cPbar_y)
     \ \cong \ H^{0}(\cCbar,\cO_{\cCbar}(D))
\end{equation}
is an isomorphism.

Via the isomorphism $\Sym^{(d)}(\cCbar) \cong \BDiv_{\cCbar/\Fbar}^d$, 
we can identify $\Sym^{(d)}(\cCbar)(\Fbar)$ with  
the set of effective divisors
$D = \sum_{i=1}^d (p_i)$ of degree $d$ supported on $\cCbar(\Fbar)$.
Let $Q : \cCbar^d \rightarrow \Sym^{(d)}(\cCbar)$
be the quotient map, and let
$P = \BA_{\cCbar/\Fbar} : \Sym^{(d)}(\cCbar) \rightarrow \BPic_{\cCbar/\Fbar}^d$ 
be the Abel map,\index{Abel map} so that 
$P(D) = [D]$.  Then $Q$ is finite of degree $d!$, and the fibres of $P$ are
isomorphic to $\PP^{d-g}$.  Indeed, if $y = [D] \in \BPic_{\cCbar/\Fbar}^d(\Fbar)$,
\begin{eqnarray*}
P^{-1}(y) & = & \{ D^{\prime} \in \Sym^{(d)}(\cCbar)(\Fbar) : D^{\prime} \sim D \} \\
            & = & \{ D^{\prime} \in \Div_{\cCbar/\Fbar}^d(\Fbar) :
                          D^{\prime} \ge 0, D^{\prime} \sim D \} \\
            & \cong & \Proj(H^0(\cCbar,\cO_\cCbar(D))) \ .
\end{eqnarray*}

Let $U \subset \BPic_{\cCbar/\Fbar}^d$ be an affine subset small enough that
$(p_2)_*(\cPbar)$ is free over $\cO_{\BPic_{\cCbar/\Fbar}^d}(U)$.  
For each $y_0 \in \BPic_{\cCbar/\Fbar}^d(\Fbar)$, there is such an affine containing $y_0$.

Put $\cZbar|_U = p_2^{-1}(U) \cong \cCbar \times U \subset \cZbar$, and let 
\begin{eqnarray*}
\Sym^{(d)}(\cCbar)|_U \ = \  P^{-1}(U) \ \subset \ \Sym^{(d)}(\cCbar) \ , 
                 \qquad \\
\cCbar^d|_U \ = \ (P \circ Q)^{-1}(U) 
       \ = \ Q^{-1}(\Sym^{(d)}(\cCbar)|_U) \ \subset \ \cCbar^d \ . 
\end{eqnarray*} 
Then $\Sym^{(d)}(\cCbar)|_U  \cong  \PP^{d-g} \times U$.  
Let $\cF_0, \ldots, \cF_{d-g}$ be a
set of $\cO_{\BPic_{\cCbar/\Fbar}^d}(U)$-generators for 
$H^0(\cZbar|_U,\cPbar) \cong H^0(U,(p_2)_*(\cPbar))$.  
Then for each $y = [D] \in U(\Fbar)$,
the sections $i_{y}^*(\cF_j) = \cF_j|_{\cZbar_{y}}$ form a basis for
$H^0(\cCbar,\cO_{\cCbar}(D))$.  
Recall that $\Sym^{(d)}(\cCbar) \cong \BDiv^d_{\cCbar/\Fbar} \subset \BHilb_{\cCbar/\Fbar}$, 
and let $((a_0:\cdots:a_{d-g}),y)$ vary over $\PP^{d-g} \times U$.  
By the universal property of the Hilbert scheme, the flat family of divisors 
\index{Hilbert scheme}
\begin{equation} \label{Forward} 
\widetilde{\mu}((a_0:\cdots:a_{d-g}),y)
     \ = \ \div_{\cCbar}(i_{y}^*(\sum a_j \cF_j)) \ =: \ \sum_{i=1}^d (p_i) 
\end{equation}
corresponds to a morphism $\mu : \PP^{d-g} \times U \rightarrow \Sym^{(d)}(\cCbar)|_U$.  
However, this morphism is simply a realization of the fibration of $\Sym^{(d)}(\cCbar)|_U$
as a trivial $\PP^{d-g}$-bundle over $U$, and hence is an isomorphism. 
%
%
  
Let $\lambda : \Sym^{(d)}(\cCbar)|_U \rightarrow \PP^{d-g} \times U$ 
be the isomorphism inverse to $\mu$.  The surjective morphism
\begin{equation*}
\varphi_U = \lambda \circ Q: \cCbar^d|_U
                      \rightarrow \PP^{d-g} \times U 
\end{equation*}
provides a means of parametrizing sections of $\cPbar$ by their zeros: 
 using the homogeneous coordinates $a_0, \ldots, a_{d-g}$ on $\PP^{d-g}$, 
given $\vp \in \cCbar|_U^d(\Fbar)$, write
\begin{equation*}
\varphi_U(\vp) \ = \ \big((a_0(\vp): \cdots :a_{d-g}(\vp)),[\vp]\big) 
               \ \in \ (\PP^{d-g} \times U)(\Fbar) \ .
\end{equation*}
Then by (\ref{Forward}), tautologically 
\begin{equation} \label{FDF} 
\div_{\cCbar}(i_{[\vp]}^*(\sum a_j(\vp) \cF_j)) \ = \ \mu(\varphi_{U}(\vp)) 
 \ = \ \sum (p_i) \ .
\end{equation} 

For each pair $(j, k)$ with  $(j,k) \in \{0, \ldots, d-g\}$, consider the rational
function on $\cCbar^d$ defined by  $h_{j,k}(\vp) = a_j(\vp)/a_k(\vp)$ 
on $\cCbar^d|_U(\Fbar)$.
Its domain includes all points in $U$ where $a_k(\vp) \ne 0$.  Put   
\begin{equation} \label{BFLM1A}
\cG_{k,U}(z,\vp) \ = \ \sum_{j=0}^{d-g} h_{j,k}(\vp) \cF_j(z,[\vp])   
          \ \in \ \Fbar(\cCbar^d) \otimes_{\Fbar} H^0(\cZbar|_U,\cPbar) \ . 
\end{equation}
Formula (\ref{FDF}) shows that
for each $\vp \in \cCbar^d|_U(\Fbar)$ with $a_k(\vp) \ne 0$, the pullback
\begin{equation*}
G_{k,\vp}(z) \ = \ i_{[\vp]}^*(\cG_{k,U}(z,\vp))
\end{equation*}  
is a section of $H^0(\cCbar,\cO_{\cCbar}([\vp]))$ satisfying 
\begin{equation} \label{BFLM1B} 
\div_{\cCbar}(G_{k,\vp}(z)) \ = \ \sum_{i=1}^d (p_i)  \ .
\end{equation}

Put $\Ybar|_U = \Ybar \cap (\cCbar^d|_U \times \cCbar^d|_U)$.  
For each $(k,\ell)$ with $0 \le k, \ell \le d-g$, let 
\begin{equation} \label{FfDef1}
\fbar_{k,\ell,U}(z,w;\vp,\vq) 
      \ = \   \frac{\cG_{k,U}(z,\vp)}{\cG_{\ell,U}(z,\vq)} 
                    \frac{\cG_{\ell,U}(w,\vq)}{\cG_{k,U}(w,\vp)} \ .      
\end{equation}
This a rational function on $\cCbar \times \cCbar \times \Ybar$, 
defined at least for $\vp, \vq \in \cCbar^d|_U(\Fbar)$ 
where $a_k(\vp), a_{\ell}(\vq) \ne 0$,
and for $z, w \in \cCbar(\Fbar)$ distinct from the $p_i$, $q_i$.     

Suppose $\vp \in \cCbar^d|_U(\Fbar)$.  
If $\vq \in \cCbar^d(\Fbar)$ satisfies $\sum (q_i) \sim \sum (p_i)$, 
that is, if $[\vq] = [\vp] \in U$ (so in particular $\vq \in \cCbar^d|_U(\Fbar)$), 
there is an index $\ell$ such that $a_{\ell}(\vq) \ne 0$. 
For each $w \in \cCbar(\Fbar)$ distinct from the $p_i, q_i$,  
\begin{equation*}
\fbar_{w,\vp,\vq}(z)  \ = \ i_{w,\vp,\vq}^*(\fbar_{k,\ell,U}) 
  \ = \   \frac{G_{k,\vp}(z)}{G_{\ell,\vq}(z)} 
                   \frac{G_{\ell,\vq}(w)}{G_{k,\vp}(w)} \ \in \ \Fbar(\cCbar) 
\end{equation*}
is the unique rational function on $\cCbar$ for which  
$\div_{\cCbar}(\fbar_{w,\vp,\vq}) = \sum (p_i) - \sum (q_i)$
and $\fbar_{w,\vp,\vq}(w)$ = 1.  Hence as  $U$ and $k, \ell$  vary, 
the $\fbar_{k,\ell,U}$ glue to give the 
desired function $\fbar(z,w;\vp,\vq)$. 

\vskip .1 in

Now consider the field of definition of $\fbar(z,w;\vp,\vq)$.  
Recall that a field $M$ is pseudo-algebraically closed\index{pseudo-algebraically closed field}
 (PAC) if every absolutely irreducible 
variety $V/M$ has an $M$-rational point.  It is well-known that for any field $F$, the separable closure 
$F^{\sep}$ is PAC (see \cite{FJar}, p.130, and \cite{L}, p.76).  
Since $\cC(F^{\sep})$ is nonempty, 
there is a finite separable extension $\hF/F$ for which $\cC(\hF)$ is nonempty.  
Thus the theory of the Picard scheme applies over $\hF$. 
\index{Picard scheme} 
Let $\hf(z,w;\vp,\vq)$ be the function $\fbar(z,w;\vp,\vq)$ constructed above, 
regarding $\hF$ as the ground field.  

We will first show that $\hf(z,w;\vp,\vq)$ is $\hF$-rational.  
Put $\hcC = \cC_{\hF}$, $\hY = Y_{\hF}$, $\hcZ = \cZ_{\hF}$.  
If $g = 0$, then $\hcC \cong \PP^1_{\hF}$, 
$\hY \cong (\PP^1_{\hF})^d \times (\PP^1_{\hF})^d$,
and $\hcZ \cong \PP^1_{\hF} \times \PP^1_{\hF} \times (\PP^1_{\hF})^d \times (\PP^1_{\hF})^d$.  
In this case $\hf(z,w;\vp,\vq)$ is defined using the crossratio\index{crossratio}
 and is $\hF$-rational by construction.
If $g > 0$, then $\hcZ = \hcC \times \BPic_{\hcC/\hF}^d$; 
let $\hcP$ be a Poincar\'e sheaf\index{Poincar\'e sheaf} on $\hcZ$.  Then $\cZbar = (\hcZ)_{\Fbar}$ 
and we can take $\cPbar = (\hcP)_{\Fbar}$.
The varieties $\hcC$ and $\hY$ are $\hF$-rational, 
so $\hcC \times \hcC \times \hY$ is $\hF$-rational.  
For each $y_0 \in \BPic_{\hcC/\hF}^d(\Fbar)$ there is an $\hF$-rational 
affine neighborhood $U$ of $y_0$ such that
$(p_2)_*(\hcP)$ is free over $\cO_{\BPic_{\hcC/\hF}^d}(U)$, 
so we can assume the affines $U$ in the construction above are $\hF$-rational.  
The invertible sheaf $\hcP$ is $\hF$-rational, so the sections 
$\hcF_0, \ldots, \hcF_d$ can be chosen to be $\hF$-rational, 
and then the maps $\varphi_U$ and the functions $h_{j,k}(\vp)$ 
will be $\hF$-rational.  Hence $\hf(z,w;\vp,\vq)$ is $\hF$-rational.

Recall that $\Aut(\Fbar/F) \cong \Gal(F^{\sep}/F)$.  Given $\sigma \in \Aut(\Fbar/F)$, 
let $\hcC^{\sigma}$ (resp. $\hY^{\sigma}$) be the conjugate variety to $\hcC$ (resp. $\hY$).  
Similarly, put $\hf^{\sigma} = \sigma \circ f \circ \sigma^{-1}$;
it is a function on $\hcC^{\sigma} \times \hcC^{\sigma} \times \hY^{\sigma}$.  
It has properties analogous to those of $\hf$:  
for each $(w, \vp, \vq) \in \hcC^{\sigma} \times \hY^{\sigma}$ with $w \ne p_i,q_i$ for all $i$, 
if we put $\hf^{\sigma}_{w,\vp,\vq}(z) = \hf^{\sigma}(z,w;\vp,\vq)$ then 
 $\div_{\hcC^{\sigma}}(\hf^{\sigma}_{w,\vp,\vq}) = \sum (p_i)-\sum (q_i)$ 
and $\hf^{\sigma}_{w,\vp,\vq}(w) = 1$.  
Regarding $\cC$ and $\hC$ as projective varieties cut out by $\hF$-rational equations, 
fix an $\hF$-rational isomorphism $\gamma : \cC \rightarrow \hC$.  By abuse of notation, 
we will denote the induced isomorphisms $\cC^{d} \rightarrow \hcC^{d}$ 
and  $\cC \times \cC \times Y \rightarrow \hcC \times \hcC \times \hY$  
by $\gamma$ as well. 

\smallskip
Define $f(z,w;\vp,\vq)$ on $\cC \times \cC \times Y$ by $f = \hf \circ \gamma$. 
We claim that $f$ is $F$-rational.  
It suffices to show that $f^{\sigma} = f$ for all $\sigma \in \Aut(\Fbar/F)$.
For this, note that by the defining properties of $\hf$, for each $(\vp,\vq) \in Y(\Fbar)$ 
and each $w \in \cC(\Fbar)$ distinct from the $p_i, q_i$, the function $f_{w,\vp,\vq} \in \Fbar(\cC)$
given by $f_{w,\vp,\vq}(z) = f(z,w;\vp,\vq)$ satisfies $f_{w,\vp,\vq}(w) = 1$ and 
$\div_{\cC}(f_{w,\vp,\vq}) = \sum (p_i) - \sum (q_i)$.  Indeed, 
$f_{w,\vp,\vq}(w) = \hf(\gamma(w),\gamma(w);\gamma(\vp),\gamma(\vq)) = 1$ and 
\begin{eqnarray*} 
\div_{\cC}(f_{w,\vp,\vq}) 
& = & \gamma^{-1}\big(\div_{\hC}(\hf_{\gamma(w),\gamma(\vp),\gamma(\vq)})\big) \\
& = & \gamma^{-1}\big(\sum (\gamma(p_i)) - \sum (\gamma(q_i))\big)
                          \ = \ \sum (p_i) - \sum (q_i) \ .
\end{eqnarray*} 
These properties uniquely determine $f_{w,\vp,\vq}$.  However, for each $\sigma$, 
an analous computation using the representation $f^{\sigma} = \hf^{\sigma} \circ \gamma^{\sigma}$ 
shows that $(f^{\sigma})_{w,\vp,\vq}$ has the same properties.  
Hence $f_{w,\vp,\vq} = (f^{\sigma})_{w,\vp,\vq}$.  Since $f(z,w;\vp,\vq) = f^{\sigma}(z,w;\vp,\vq)$ 
for a Zariski-dense set of points $(z,w,\vp,\vq)$, it follows that $f = f^{\sigma}$.                         

\smallskip 
The final assertion in the theorem is that for any subfield $F \subset M \subset \Fbar$,
if $\sum (p_i) - \sum (q_i)$ is $M$-rational and $w \in \cC(M)$, then 
$f_{w,\vp,\vq}$ is $M$-rational.  This is clear, since there exists a  
$g \in M(\cC)$ for which $\div_{\cC}(g) = \sum (p_i) - \sum (q_i)$, and then the function 
$g(w)^{-1} \cdot g$ is $M$-rational and has the 
properties that characterize $f_{w,\vp,\vq}$.  
\end{proof}    
 
\noindent{\bf Remark.}  Robert Varley has given an alternate construction of $f(z,w;\vp,\vq)$
\index{Varley, Robert} 
which does not use the theory of the Picard scheme, but only requires Grauert's theorem.
\index{Picard scheme}\index{Grauert's theorem}   
His construction applies even when the degree $d$ is not in the 
stable range $d \ge \max(1,2g-1)$.  

We now sketch this construction\footnote{The author thanks Varley\index{Varley, Robert} 
for permission to include his construction here.}  
As before, let $F$ be a field, and let $\cC/F$ be a smooth, projective, geometrically integral 
curve of genus $g \ge 0$.  Fix $d \ge 1$ and let $\tY \subset Y$ 
be a reduced, irreducible, locally closed  $F$-rational subvariety. 
(Note that if $d$ is in the stable range, then $Y$ is irreducible, 
but in general it may have more than one component.)   
We claim there is a function $\widetilde{f}(z,w;\vp,\vq) \in F(\cC \times \cC \times \tY)$ 
with the properties in Theorem \ref{UniversalFcn}.  

Let $U \subset \cC \times \tY$ be 
the open $F$-rational subvariety for which  
\begin{equation*} 
U(\Fbar) \ = \ \{ (w,\vp,\vq) \in \cC(\Fbar) \times \tY(\Fbar) : \text{$w$ is distinct 
from $p_1, \ldots, p_d, q_1, \ldots, q_d$} \} \ .
\end{equation*}  
Define sections 
$\varphi, \sigma_i, \tau_i : U \rightarrow \cC \times U \hookrightarrow \cC \times \cC \times \tY$ by 
\begin{equation*} 
\varphi(w,\vp,\vq) = (w,w,\vp,\vq), \quad \sigma_i(w,\vp,\vq) = (p_i,w,\vp,\vq), \quad 
\tau_i(w,\vp,\vq) = (q_i,w,\vp,\vq) 
\end{equation*}
for $i = 1, \ldots, d$.  Let $W$ be the subvariety $\varphi(U) \subset \cC \times U$, 
and let $\cD$ and $\cE$ be the Cartier divisors on $\cC \times U$ 
corresponding to the $F$-rational Weil divisors $\sum_{i=1}^d \sigma_i(U)$, $\sum_{i=1}^d \tau_i(U)$
\index{Weil!divisor} 
respectively.  
Consider the line bundles $\cO_{\cC \times U}(\cE - \cD)$ and $\cO_W(\cE - \cD)$, 
the restriction map $r : \cO_{\cC \times U}(\cE - \cD) \rightarrow \cO_W(\cE - \cD)$, 
the projection $p_2 : \cC \times U \rightarrow U$, and the direct images 
$(p_2)_*(\cO_{\cC \times U}(\cE - \cD))$ and $(p_2)_*(\cO_W(\cE - \cD))$.
Then one can show that
\begin{enumerate}

\item $H^0(U,(p_2)_*(\cO_{\cC \times U}(\cE - \cD)))$ embeds naturally 
in the function field $F(\cC \times \cC \times \tY)$;

\item $(p_2)_*(\cO_{\cC \times U}(\cE - \cD))) \cong \cO_U$;  

\item the function $\widetilde{f}(z,w;\vp,\vq) \in F(\cC \times \cC \times \tY)$ 
corresponding to the canonical section $1 \in H^0(U,\cO_U)$ has the desired properties.
\end{enumerate} 

Indeed, (1) follows from a standard interpretation of sections of 
$\cO_{\cC \times U}(\cE - \cD)$ as elements of the function field 
$F(\cC \times U) = F(\cC \times \cC \times \tY)$. 
Assertion (2) follows from two subclaims:  first, 
$(p_2)_*(\cO_{\cC \times U}(\cE - \cD)) \cong (p_2)_*(\cO_W(\cE - \cD))$, 
and second, $(p_2)_*(\cO_W(\cE - \cD)) \cong \cO_U$.  
For the first, note that the fibres of $(p_2)_*(\cO_{\cC \times U}(\cE - \cD))$ are one-dimensional, 
so by Grauert's theorem\index{Grauert's theorem} $(p_2)_*(\cO_{\cC \times U}(\cE - \cD))$
is a line bundle on $U$. 
It maps pointwise nontrivially into $(p_2)_*(\cO_W(\cE - \cD))$,
which is also a line bundle on $U$.
For the second, note that $\cO_W(\cE) \cong \cO_W$ and $\cO_W(\cD) \cong \cO_W$
since $\cD$ and $\cE$ are disjoint from $W$. 
However $(p_2)_*\cO_W \cong \cO_U$ since $\varphi : U \rightarrow W \subset \cC \times U$ is a section,
and so $(p_2)_*(\cO_W(\cE - \cD)) \cong \cO_U$.

\vskip .1 in 
We now specialize to the case $F = K_v$, where $K_v$ is a nonarchimedean local field.

Our main result is the following.  Let $\cC_v/K_v$ be a smooth, projective,
geometrically integral curve of genus $g \ge 0$. 
Fix a spherical metric $\|z,w\|_v$ on $\cC_v(\CC_v)$,\index{spherical metric}   
and recall that for each $p \in \cC_v(\CC_v)$ and each $r > 0$, 
\begin{eqnarray*}
B(p,r)^- & = & \{z \in \cC_v(\CC_v) : \|z,p\|_v < r \} \ , \\
B(p,r) \  & = & \{ z \in \cC_v(\CC_v) : \|z,p\|_v \le r \} \ .
\end{eqnarray*}
We will show that if $\vp, \vq \in \cC_v(\CC_v)^d$ 
are sufficiently near each other, and $[\vp] = [\vq]$,  
then $f(z,w;\vp,\vq)$ is close to $1$ outside 
fixed balls containing the $p_j$ and $q_j$.    

\begin{theorem} \label{BMT2}  Let $K_v$ be a nonarchimedean local field.
Suppose $\cC_v/K_v$ is a smooth, projective, 
geometrically integral curve of genus $g \ge 0$, and let $E_v \subset \cC_v(\CC_v)$ be compact.
Then for each $d \ge \max(1,2g)$, 
there are a radius $r_0 = r_0(E_v,d) > 0$ and a constant $D = D(E_v,d) > 0$ 
with the following property:    

Given $0 < \varepsilon < r \le r_0$, 
suppose  $\vp, \vq \in E_v^d$ are such that $\sum (p_j) \sim \sum (q_j)$ 
and $\|p_j,q_j\|_v \le \varepsilon$ for each $j = 1, \ldots, d$.  
Put $r_j = \|p_j,q_j\|_v$.  Then 

$(A)$ for all $z, w$ in
 $\cC_v(\CC_v) \backslash 
       \big((\bigcup_{j=1}^d B(p_j,r_j)^-) 
            \cup (\bigcup_{j=1}^d B(q_j,r_j)^-)\big)$ 
we have  
\begin{equation*} 
|f(z,w;\vp,\vq)|_v \ = \ 1 \ ; 
\end{equation*} 

$(B)$ for all $z, w$ in
 $\cC_v(\CC_v) \backslash 
       \big((\bigcup_{j=1}^d B(p_j,r)^-) 
            \cup (\bigcup_{j=1}^d B(q_j,r)^-)\big)$ 
we have
\begin{equation*} 
|f(z,w;\vp,\vq) - 1|_v  \ \le \ \frac{D}{r^d} \varepsilon \ .
\end{equation*} 
\end{theorem}

We will use Theorem \ref{BMT2} in the proof of Lemma \ref{MovingLemma1}, the ``First Moving Lemma'' 
in the patching process in the nonarchimedean compact case.  In our application 
\index{patching argument!for nonarchimedean $K_v$-simple sets}
$r$ will be fixed, and the important factor governing $|f(z,w;\vp,\vq) - 1|_v$ 
will be $\max_j(\|p_j,q_j\|_v)$. 

\begin{proof}[Proof of Theorem \ref{BMT2} when $g = 0$.]  
Fixing an isomorphism of $\cCbar_v/\CC_v$ with $\PP^1/\CC_v$,  
we can assume that $f(z,w;\vp,\vq) = \prod_{j=1}^d \chi(z,w;p_j,q_j)$.
Since any two spherical metrics on $\cC_v(\CC_v)$ 
are comparable\index{spherical metric!from different embeddings comparable}  
 (\cite{RR1}, Theorem 1.1.1)
we can assume that $\|x,y\|_v$ is the standard metric on $\PP^1(\CC_v)$ given for $x, y \in \CC_v$ by 
\begin{equation} \label{FWAX1}
\|x,y\|_v \ = \ \frac{|x-y|_v}{\max(1,|x|_v) \max(1,|y|_v)} \ .
\end{equation}
By simple algebraic manipulations, one sees that that for $z,w, p, q \in \CC_v$ 
\begin{equation} \label{FWAX2} 
|\chi(z,w;p,q)|_v \ = \ \left| \frac{(z-p)(w-q)}{(z-q)(w-p)} \right|_v 
                \ = \ \frac{\|z,p\|_v \|w,q\|_v}{\|z,q\|_v \|w,p\|_v} \ ,
\end{equation}
and that 
\begin{equation} \label{FWAX2A} 
\chi(z,w;p,q) - 1 \ = \ \frac{(z-w)(p-q)}{(z-q)(w-p)}\ , \quad 
|\chi(z,w;p,q) - 1|_v \ = \ \frac{\|z,w\|_v \|p,q\|_v}{\|z,q\|_v \|w,p\|_v} \ .
\end{equation} 
By continuity, these formulas extend to $z, w, p, q \in \PP^1(\CC_v)$. 
Furthermore, by a telescoping argument 
\begin{equation} \label{FWAX3} 
f(z,w;\vp,\vq) - 1 \ = \ \sum_{j=1}^d \Big(\chi(z,w;p_j,q_j) - 1\Big) 
\cdot \left( \prod_{k=j+1}^d \chi(z,w;p_k,q_k) \right) \ . 
\end{equation} 
 
Take $r_0 = 1$ and $D = 1$.  
 
If $z, w \in \PP^1(\CC_v)  \backslash \big((\bigcup_{j=1}^d B(p_j,r_j)^-) 
            \cup (\bigcup_{j=1}^d B(q_j,r_j)^-)\big)$ 
then since $\|x,y\|_v$ satisfies the ultrametric inequality and $\|p_j,q_j\| = r_j$, 
we have $\|z,p_j\|_v = \|z,q_j\|_v$ and $\|w,p_j\|_v = \|w,q_j\|_v$ for each $j$.  
From (\ref{FWAX2}) it follows that 
\begin{equation*} 
|f(z,w,\vp,\vq)|_v \ = \ \prod_{j=1}^d |\chi(z,w;p_j,q_j)|_v \ = \ 1 \ ,
\end{equation*}
which is assertion (A).  

If $z, w \in \PP^1(\CC_v)  \backslash \big((\bigcup_{j=1}^d B(p_j,r)^-) 
            \cup (\bigcup_{j=1}^d B(q_j,r)^-)\big)$, then 
since $\|x,y\|_v$ satisfies the ultrametric inequality, and $\|z,q_j\|_v, \|w,p_j\|_v \ge r$ 
while $\|p_j,q_j\|_v = r_j \le r$, it follows that 
\begin{equation*}
\|z,w\|_v \ \le \ \max(\|z,q_j\|, \|q_j,p_j\|_v, \|w,p_j\|_v) \ = \ 
\max(\|z,q_j\|, \|w,p_j\|_v) \ .
\end{equation*} 
It then follows from (\ref{FWAX2A}) that 
\begin{equation} \label{FWAX4}
 |\chi(z,w;p_j,q_j) - 1|_v \ = \  \frac{\|z,w\|_v \|p_j,q_j\|_v}{\|z,q_j\|_v \|w,p_j\|_v} \
\le \ \frac{\|p_j,q_j\|_v}{r}  \ .
\end{equation}
Hence by (\ref{FWAX2}), (\ref{FWAX3}), (\ref{FWAX4}) and the ultrametric inequality, 
\begin{equation*} 
|f(z,w;\vp,\vq)-1|_v \ \le \ \frac{1}{r} \cdot \max_j(\|p_j,q_j\|_v) 
\end{equation*} 
which is stronger than the inequality claimed in (B).
\end{proof} 

Note that in the proof when $g = 0$, we did not use anything about the compact set $E_v$, 
and the bound in Theorem \ref{BMT2}(B) holds for all $\vp,\vq \in \cC_v(\CC_v)^d$
with each $\|p_j,q_j\|_v \le r$.  It seems likely that this remains true when $g > 0$ as well. 
For our application we only need the bound when $\vp,\vq \in E_v^d$, so we have not pursued it.      

\vskip .1 in
For the remainder of this Appendix, we will assume that $g > 0$.
In this case the proof of Theorem \ref{BMT2} requires much more machinery.  
The idea is to first locally control the functions $G_{k,\cF,U}(z,\vp)$
in the factorization (\ref{FNEW1}) below, using power series expansions, 
and then extend that control to all of $\cC_v(\CC_v)$ using various forms of 
the Maximum modulus principle and the theory of the canonical distance.   
\index{Maximum principle!nonarchimedean}\index{canonical distance!factorization property}  

Before giving the proof, we will need several technical lemmas.  We wish to apply the theory of 
rigid analysis, so we work with $\cCbar_v$ rather than $\cC_v$.  
The first three lemmas prepare the way to use the Maximum modulus principle of rigid analysis 
on $\cCbar_v \times \cCbar_v^d$. 
\index{Maximum principle!nonarchimedean!from Rigid analysis} 


\begin{lemma} \label{BLemA1}
Let $p_1, \ldots, p_d \in \cCbar_v(\CC_v)$.  Suppose $0 < r < 1$ belongs 
to value group of $\CC_v^{\times}$ and is small enough that each ball $B(p_j,r)$ 
is isometrically parametrizable.  Then for each
$\zeta \in \ \cCbar_v(\CC_v) \backslash (\bigcup_{j=1}^d B(p_j,r))$,
there is a function $g(z) \in \CC_v(\cCbar_v)$ with poles only at $\zeta$,
such that
\begin{eqnarray}
\bigcup_{j=1}^d B(p_j,r)\ & = & \{ z \in \cCbar_v(\CC_v) : |g(z)|_v \le 1 \} \ ,
        \label{BFGG1} \\
\bigcup_{j=1}^d B(p_j,r)^- & = & \{ z \in \cCbar_v(\CC_v) : |g(z)|_v < 1 \} \ .
        \label{BFGG2}
\end{eqnarray}
\end{lemma}

\begin{proof}
If the balls $B(p_j,r)$ are pairwise disjoint, the result follows from
(\cite{RR1}, Theorem 4.2.16) and its proof.  In the general case, note that
since $\|z,w\|_v$ satisfies the ultrametric inequality, then any two balls
$B(p_i,r)$ and $B(p_j,r)$ either coincide, or are disjoint.
For each $\ell = 1, \ldots, d$ we can represent $\bigcup_j B(p_j,r)$ as a 
disjoint union of a subset of the balls $B(p_j,r)$, in such a way that 
$B(p_{\ell},r)$ occurs in the representation.  Let $g_{\ell}(z)$ 
be the function obtained for this representation, and put 
$g(z) = \prod_{\ell=1}^d g_{\ell}(z)$.  Then (\ref{BFGG1}) and 
(\ref{BFGG2}) hold for this $g(z)$.
\end{proof}

\vskip .1 in
Recall that a subset 
$W \subset \cCbar_v(\CC_v)$ is called an $\RL$-domain  
\index{$\RL$-domain|ii}
(``Rational Lemniscate domain''; see \cite{RR1}, p.220) 
if there is nonconstant function $h \in \CC_v(\cCbar_v)$ for which
\begin{equation*}
W  \ = \ \{ z \in \cCbar_v(\CC_v) : |h(z)|_v \le 1 \} \ .
\end{equation*}
In that case  
\begin{equation*}
\partial W \ = \ 
\partial W(h) \ = \ \{ z \in \cCbar_v(\CC_v) : |h(z)|_v = 1 \} 
\end{equation*}
is called its boundary (with respect to $h$).  By (\cite{RR1}, Corollary 4.2.14),  
a finite intersection (or union) of $RL$-domains is again an $RL$-domain.  However,
\index{$\RL$-domain} 
that Corollary does not explicitly give a boundary.     

\begin{lemma} \label{BLemA2}
Under the hypotheses of Lemma $\ref{BLemA1}$,

$(A)$  $\cCbar_v(\CC_v) \backslash (\bigcup_{j=1}^d B(p_j,r)^-)$ is an 
$RL$-domain with boundary
\index{$\RL$-domain}
\begin{equation*}
(\bigcup_{j=1}^d B(p_j,r)) 
        \backslash (\bigcup_{i=1}^d B(p_j,r)^-) \ ; 
\end{equation*}         

$(B)$  $(\bigcup_{j=1}^d B(p_j,r)) \backslash 
(\bigcup_{j=1}^d B(p_j,r)^-)$ is an $RL$-domain.
\index{$\RL$-domain}

$(C)$  Each isometrically parametrizable ball $B(p_j,r_j)$ 
with $r_j$ in the value group of $\CC_v^{\times}$, is an $RL$-domain 
\index{$\RL$-domain}
with boundary $B(p_j,r_j) \backslash B(p_j,r_j)^-$.  
\end{lemma}

\begin{proof}
Let $g(z) \in \CC_v(\cCbar_v)$ be the function from Lemma \ref{BLemA1}.
Then the $RL$-domain 
\index{$\RL$-domain}
\begin{equation*}
\{ z \in \cCbar_v(\CC_v) : |1/g(z)|_v \le 1 \} 
\ = \ \cCbar_v(\CC_v) \backslash (\bigcup_{j=1}^d B(p_j,r)^-)  
\end{equation*}
has boundary 
\begin{equation*}
\{ z \in \cCbar_v(\CC_v) : |1/g(z)|_v = 1 \} \ = \
               (\bigcup_{j=1}^d B(p_i,r)) \backslash (\bigcup_{j=1}^d B(p_j,r)^-)
\end{equation*}
which proves (A).  For part (B), note that 
\begin{equation*}
(\bigcup_{j=1}^d B(p_j,r)) 
     \backslash (\bigcup_{j=1}^d B(p_j,r)^-) \ = \
          \{ z \in \cCbar_v(\CC_v) : |1/g(z)|_v \le 1, |g(z)|_v \le 1 \}
\end{equation*}
and apply (\cite{RR1}, Corollary 4.2.14).  
Part (C) follows by applying Lemma \ref{BLemA1}
to each $B(p_j,r_j)$ by itself.
\end{proof}

\vskip .1 in
Recall that there is a faithful functor from the category of varieties
over $\CC_v$ to the category of rigid analytic spaces\index{rigid analytic space} over $\CC_v$ 
(see \cite{BGR}, p.363).  If $\Xbar_v/\CC_v$ is a variety, 
we will say that a subset of $\Xbar_v(\CC_v)$ is a 
{\em affinoid domain} if its image 
\index{affinoid!domain|ii}
in the rigid analytic space $\Xbar_v^{an}$ associated to $\Xbar_v$ is an 
admissible affinoid in the sense of rigid analysis (see \cite{BGR}, p.277, p.357):
\index{affinoid!admissible} 
essentially, if its image under the functor above is the underlying point set of 
$\Sp(T)$ for some Tate algebra $T = \CC_{v}\<<z_1, \ldots, z_k\>> /\cI$ associated to 
an affine subset of $\Xbar_v$. 
Here $\CC_{v}\<<z_1, \ldots, z_k\>>$ is the ring of power series converging on the unit polydisc,
and $\cI$ is a finitely generated ideal. 
   
An $\RL$-domain on a curve is an affinoid domain in the sense above: 
\index{$\RL$-domain}
\index{affinoid!domain|ii}

\begin{lemma} \label{RLeqAffinoid}  
If $\cCbar_v/\CC_v$ is a curve, and $W \subset \cCbar_v(\CC_v)$ is an $RL$-domain, 
\index{$\RL$-domain}
then $W$ is an affinoid domain.  
\index{affinoid!domain|ii} 
\end{lemma} 

\begin{proof}  Let $h(z) \in \CC_v(\cCbar_v)$ be a nonconstant function for which 
$W = \{z \in \cCbar_v(\CC_v) : |h(z)|_v \le 1 \}$.  
The function field $\CC_v(\cCbar_v)$ is finite over $\CC_v(h)$.  
If $\Char(\CC_v) = 0$, then $\CC_v(\cCbar_v)/\CC_v(h)$ is separably algebraic.  
If $\Char(\CC_v) = p > 0$, then since $\CC_v$ is algebraically closed, 
$\CC_v(\cCbar_v)$ is separably generated. 
After replacing $h$ by $h^{1/p^m}$ for some $m$ we can assume that 
$\CC_v(\cCbar_v)/\CC_v(h)$ is separably algebraic.  
By the primitive element theorem, 
there is a function $G(z) \in \CC_v(\cCbar_v)$ for which $\CC_v(\cCbar_v) = \CC_v(h,G)$.  
Let $f(x,h) = x^d + a_1(h) x^{d-1} + \ldots + a_d(h)$ 
be the minimal polynomial of $G$ over $\CC_v(h)$, 
and let $\pi_v \in \CC_v^{\times}$ satisfy $\ord_v(\pi_v) > 0$.  

After multiplying $G$ by a power of the product of denominators of the rational functions $a_i(h)$ and an 
appropriate power of $\pi_v$, we can assume that the $a_i(h)$ are polynomials in 
$h$ with coefficients in $\hcO_v$.  By the ultrametric inequality, 
for each $z \in W$ we must  have $|G(z)|_v \le 1$.   
Hence the map which sends $z$ to $(x,y) = (G(z),h(z))$  
induces an isomorphism of $W$ with the underlying point set of $\Sp(\CC_v\<<x,y\>>/(f(x,y))$. 
Examining the construction in (\cite{BGR}, Example 2, p.363), one sees that 
this map is the one realizing $W$ as an affinoid domain under the functor above.  
\index{affinoid!domain|ii}
\end{proof}  

\begin{corollary} \label{AffinoidCor} 
Let $p_1, \ldots, p_d \in \cCbar_v(\CC_v)$; suppose that $r, r_1, \ldots, r_d$ belong to the 
value group of $\CC_v^{\times}$ and are small enough that each ball $B(p_j,r)$ and $B(p_j,r_j)$ 
is isometrically parametrizable.  Then  
\begin{eqnarray*}
W & := & ((\bigcup_{j=1}^d B(p_j,r))
               \backslash (\bigcup_{j=1}^d B(p_j,r)^-))
                          \times \prod_{j=1}^d B(p_j,r_j) \\
  &   & \qquad \qquad \subset \ \cCbar_v(\CC_v) \times \cCbar_v^d(\CC_v)
\end{eqnarray*}
is an affinoid domain in $\cCbar_v \times \cCbar_v^d$.  
\index{affinoid!domain}
\end{corollary}

\begin{proof}  A product of admissible affinoids is an admissible affinoid
\index{affinoid!admissible} 
(\cite{BGR}, \S9.3.5). 
\end{proof}   

\vskip .1 in
In order to study $f(z,w;\vp,\vq)$, it will be useful to reformulate 
(\ref{FfDef1}) using functions rather than sections of a line bundle. 
We begin with the following lemma, keeping the notations 
in the proof of Theorem \ref{UniversalFcn}.  

\begin{lemma} \label{NVSec} Let $\cC/F$ be a smooth, projective, geometrically 
integral curve of genus $g > 0$.  Suppose $d \ge \max(1,2g)$.  
Let $U \subset \BPic_{\cCbar/\Fbar}^d$ 
be an affine subset over which $(p_2)_*(\cPbar)$ 
is free.  Then for any given $\vp, \vq \in \cCbar^d|_U(\Fbar)$ 
with $[\vp] = [\vq]$, and any finite set of points 
$z_1, \ldots, z_k \in \cCbar(\Fbar)$, 
there is a section $\cF \in H^0(\cZbar|_U,\cPbar)$ for which the support of  
\begin{equation*}
\div_{\cC}(i_{[\vp]}^*(\cF)) \ = \ \sum_{i=1}^d (p_i^{\prime})
\end{equation*}
is disjoint from the $p_i$, $q_i$, and $z_i$.  
\end{lemma}

\begin{proof}  Put $D = \sum (p_i)$.  
Our assumption on $d$ assures there is movement
in the linear system on $\cCbar$ associated to $D$.  
Consider the $\Fbar$-vector space 
$\Gamma(D) = \{h(z) \in \Fbar(\cC) : \div_{\cC}(h) \ge D \}$.   
Since $d \ge \max(1,2g)$, the Riemann-Roch theorem shows that 
\index{Riemann-Roch theorem} 
$\dim_{\Fbar}\big(\Gamma(D)\big) = d-g+1$, 
while for each $p \in \cCbar(\Fbar)$  
\begin{equation*}  
\dim_{\Fbar}\big(\Gamma(D-(p))\big) \ = \ d-g \ .
\end{equation*} 
The set 
\begin{equation*} 
\Gamma^{\prime}(D) \ := \ \Gamma(D) \backslash 
      \Big( \bigcup_{i=1}^d \Gamma(D-(p_i)) \cup \bigcup_{i=1}^d \Gamma(D-(q_i))
                     \cup \bigcup_{i=1}^k \Gamma(D-(z_i)) \Big) 
\end{equation*}
is nonempty because $\Fbar$ is infinite, 
and each of the finitely many subspaces removed
is a proper $\Fbar$-subspace.  For any function $h \in \Gamma^{\prime}(D)$, 
the polar divisor of $h$ is precisely $\sum (p_i)$, and the zeros of $h$
are distinct from the $p_i$, $q_i$ and $z_i$. Fix such an $h$ and write 
\begin{equation*}
\div_{\cC}(h) \ = \ \sum_{i=1}^d (p_i^{\prime}) - \sum_{i=1}^d (p_i) \ .
\end{equation*}
Then $D^{\prime} = \sum (p_i^{\prime})$ is linearly equivalent to 
$\sum (p_i)$ and $\sum (q_i)$, and the $p_i^{\prime}$ are distinct from 
the $p_i$, $q_i$, and $z_i$.      

Now let $\cF_0, \ldots, \cF_{d-g}$ be generators for 
$H^0(\cZbar|_U,\cPbar)$ over $\cO_{\BPic_{\cCbar/\Fbar}^d(U)}$.  
Let $\cH$ be the $\Fbar$-vector space generated by the $\cF_i$.  
By our assumptions, the map 
\begin{equation*} 
i_{[\vp]}^* : \cH \rightarrow H^0(\cCbar,\cO_{\cCbar}(\sum (p_i)))
\end{equation*}
is an isomorphism. Thus there is an $\cF \in \cH$ with  
$\div_{\cC}(i_{[\vp]}^*(\cF)) = \sum (p_i^{\prime})$.  
\end{proof} 

Henceforth, assume $d \ge \max(1,2g)$.  
Given an affine $U \subset \cZbar$ and a basis of sections 
$\cF_0, \ldots, \cF_{d-g} \in H^0(\cZbar_U,\cPbar)$ as above, 
consider the sections $\cG_{k,U}(z,\vp)$ defined in (\ref{BFLM1A}).   
Let $0 \ne \cF \in H^0(\cZbar_U,\cPbar)$ be an arbitrary section and put
\begin{equation} \label{BFLM1}
G_{k,\cF,U}(z,\vp) \ = \ \frac{\cG_{k,U}(z,\vp)}{\cF(z,[\vp])} 
\ = \ \sum_{j=0}^{d-g} \frac{a_j(\vp)}{a_k(\vp)} \cdot
                         \frac{\cF_j(z,[\vp])}{\cF(z,[\vp])} \ . 
\end{equation}
Then $G_{k,\cF,U}$ is an $\Fbar$-rational function on $\cC \times \cC^d$, 
defined at least for $(z,\vp) \in \cC(\Fbar) \times \cC^d|_U(\Fbar)$ where
$a_k(\vp) \ne 0$ and $\cF(z,[\vp]) \ne 0$.  The important point is that $\cF$ 
depends on $\vp$ only through $[\vp]$, so for each $\vp$ the polar divisor 
of $G_{k,\cF,U}(\vp)$  depends only on $[\vp]$.

Fix $\vp, \vq \in \cC^d|_U(\Fbar)$ 
with $[\vp] = [\vq]$, fix $z_0 \in \cC(\Fbar)$, and fix $w \in \cC(\Fbar)$  
distinct from the $p_i$ and $q_i$.  
Choose $k$ with $a_k(\vp) \ne 0$ and $\ell$ with $a_{\ell}(\vq) \ne 0$.  
By Lemma \ref{NVSec}, there is an $\cF$
for which the support of $\div_{\cC}(i_{[\vp]}^*(\cF))$ is disjoint from 
$z_0$, $w$, and the $p_i$ and $q_i$.  
Since $\cF(z,[\vp])$ depends on $\vp$\, only through $[\vp]$, 
we have  $\cF(z,[\vp]) = \cF(z,[\vq])$ for all $z$.  
It follows from (\ref{BFLM1B}) that 
\begin{equation*} 
\div_{\cC}\Big(\frac{G_{k,\cF,U}(z,\vp)}{G_{\ell,\cF,U}(z,\vq)}\Big)  
\ = \ \sum (p_i) - \sum (q_i) \ .
\end{equation*}  
Thus, for an appropriate choice of $U$, $k$, $\ell$ and $\cF$,
we can represent $f_{w,\vp,\vq}(z) = f(z,w;\vp,\vq)$ 
in a neighborhood of $z_0$ by   
\begin{equation} \label{FNEW1}
f(z,w;\vp,\vq) \ = \ 
\frac{ G_{k,\cF,U}(z,\vp)}{ G_{\ell,\cF,U}(z,\vq)} \cdot 
            \frac{G_{\ell,\cF,U}(w,\vq)}{G_{k,\cF,U}(w,\vp)} \ .  
\end{equation} 
This means that to understand $f(z,w;\vp,\vq)$, it suffices to understand the 
$G_{k,\cF,U}(z,\vp)$, which is simpler to do because only the zeros of 
$G_{k,\cF,U}(z,\vp)$ are controlled.   

\begin{lemma} \label{BLemA3}
Let $E_v \subset \cCbar_v(\CC_v)$ be compact, and let $d \ge \max(1,2g)$. 
Then there are a radius $R = R(E_v,d) > 0$ 
in the value group of $\CC_v^{\times}$,
and a number $B = B(E_v,d) > 0$, with the following properties$:$

There are finitely many affine subsets $U_i \subset \BPic_{\cCbar_v/\CC_v}^d$
such that $(p_2)_*(\cPbar)$ is free over $U_i$, 
with functions $G_i(z,\vp) = G_{k_i,\cF_i,U_i}(z,\vp)$ as in $(\ref{BFLM1})$, 
such that for each $\vp = (p_1, \ldots, p_d) \in E_v^d$, the balls $B(p_j,R)$, 
for $j = 1, \ldots, d$, are isometrically parametrizable, 
and there is some $i$ for which 
$(\bigcup_{j=1}^d B(p_j,R)) \times \prod_{j=1}^d B(p_j,R) 
\subset \cCbar_v(\CC_v) \times \cCbar_v^d |_{U_i}(\CC_v)$ and 

$(A)$ $|G_i(z,\vq)|_v \le 1$ \ for all 
             $(z,\vq)$ in $(\bigcup_{j=1}^d B(p_j,R)) 
                            \times \prod_{j=1}^d B(p_j,R)$;

$(B)$ $|G_i(z,\vq)|_v \ge B$ for all $(z,\vq)$  in
\begin{equation*}
 (\bigcup_{j=1}^d B(p_j,R)) \backslash
     (\bigcup_{j=1}^d B(p_j,R)^-) \times \prod_{j=1}^d B(p_j,\frac{1}{2}R) \ .
\end{equation*}     
\end{lemma}

\begin{proof}
By Theorem \ref{IsoParamThm} there is a number $0 < R_0 \le 1$
(depending on the spherical metric $\|z,w\|_v$)\index{spherical metric}  such that
each ball $B(a,r)$ with $a \in \cCbar_v(\CC_v)$ and $0 < r \le R_0$ is isometrically
parametrizable.

Fix $\vp \in E_v^d$.  Choose an affine set 
$U \subset \BPic_{\cCbar_v/\CC_v}^d$ for which $\vp \in \cCbar_v^d|_U(\CC_v)$
and which is small enough that $(p_2)_*(\cPbar)$
is free over $U$.  Then for each sufficiently small $r > 0$ we will have 
$\prod_{j=1}^d B(p_j,r) \subset \cCbar_v^d|_U(\CC_v)$
and all of the balls $B(p_j,r)$ will be isometrically parametrizable.

By Lemma \ref{NVSec} there is a section $\cF = \sum c_j \cF_j$ of
$(p_2)_*(\cPbar)(U)$ with coefficients $c_j \in \CC_v$
for which $\div_{\cCbar_v}(i_{[\vp]}^*(\cF)) = \sum (p_i^{\prime})$ 
is coprime to $\sum (p_i)$.
Thus, we can find $k$ and $\cF$ so that 
if $G_{k,\cF,U}$ is as in (\ref{BFLM1})), then 
$i_{[\vp]}^*(G_{k,\cF,U})$ has polar divisor $\sum (p_i^{\prime})$ with 
support disjoint from $\{p_1, \ldots, p_d\}$.
Since $G_{k,\cF,U}(p_j,\vp) = 0$ for each $j = 1, \ldots, d$, 
by continuity there is an $r > 0$ such that
$(\bigcup_{j=1}^d B(p_j,r)) \times \prod_{j=1}^d B(p_j,r)$
is contained in 
\begin{equation*}
 \{ (z,\vq) \in \cCbar_v(\CC_v) \times \cCbar_v^d|_{U}(\CC_v)
                              : |G_{k,\cF,U}(z,\vq)|_v \le 1 \}
\end{equation*}
Without loss we can assume $r \le R_0$, 
so the balls $B(p_j,r)$ are isometrically parametrizable.

By compactness, there are a finite number of points $\vp^{(i)}$ and radii
$r^{(i)}$  such that the sets  
$\prod_{j=1}^d B(p_j^{(i)},r^{(i)})$
cover $E_v^d$.  Let the $U_i$ and $G_i=G_{k_i,\cF_i,U_i}(z,\vp)$ be the
corresponding affine sets and functions, and let $R = R(E_v,d)$  be
the minimum of the $r^{(i)}$.  After shrinking $R$ if necessary, we can
assume $R$ belongs to the value group of $\CC_v^\times$.  
Then for any $\vp \in E_v^d$, there is some $\vp^{(i)}$ for which
\begin{equation*}
\prod_{j=1}^d B(p_j,R) \ \subset \
            \prod_{j=1}^d B(p_j^{(i)},r^{(i)})
\end{equation*}
and so (A) holds for this $R$.

Again fix $\vp \in E_v^d$, and choose $U = U_i$ and $G(z,\vq) = G_i(z,\vq)$
so that (A) holds. Fix $1/2 < C < 1$ in the value group of $\CC_v^{\times}$.
By Lemma \ref{BLemA2}, the set
\begin{eqnarray*}
W & := & ((\bigcup_{j=1}^d B(p_j,R))
               \backslash (\bigcup_{j=1}^d B(p_j,R)^-))
                          \times \prod_{j=1}^d B(p_j,C R) \\
  & \subset & \cCbar_v(\CC_v) \times \cCbar_v^d|_{U_i}(\CC_v)
\end{eqnarray*}
is an affinoid domain.  Moreover,
\index{affinoid!domain}
for each $\vq = (q_1, \ldots, q_d) \in \prod_{j=1}^d B(p_j,cR)$,
the function $G_{\vq}(z) = G(z,\vq)$ has zeros only at $q_1, \ldots, q_d$,
and so in particular it does not vanish on  
$((\bigcup_{j=1}^d B(p_j,R))\backslash (\bigcup_{j=1}^d B(p_j,R)^-))$.  
It follows that $1/G(z,\vq)$ is a rigid analytic function\index{rigid analytic function} on $W$. By the Maximum
Modulus principle of rigid analysis (see \cite{BGR}, p.237),
\index{Maximum principle!nonarchimedean!from Rigid analysis}
there is a number $B(\vp)$ with $0 < B(\vp) < 1$ such that
for all $(z,\vq) \in W$
\begin{equation*}
\left| \frac{1}{G(z,\vq)} \right|_v \ \le \ \frac{1}{B(\vp)} \ ;
\end{equation*}
equivalently, $|G(z,\vq)|_v \ge B(\vp)$.

Again by compactness, there are finitely many points $\vp^{(\ell)} \in E_v^d$
such that the sets $\prod_{j=1}^d B(p_j^{(\ell)},C R)$ cover $E_v^d$.  Moreover,
for any two points $\vp, \vq \in E_v^d$, the sets
$\prod_{j=1}^d B(p_j,C R)$ and $\prod_{j=1}^d B(q_j,C R)$
either coincide or are disjoint.  Hence, if  $B$ is the minimum of
the corresponding numbers $B(\vp^{(\ell)})$, then part (B) holds for all $\vp \in E_v^d$,
with this $B$.
\end{proof}

\vskip .1 in
The following lemma uses power series to obtain uniform local control of $|G(z,\vp)|_v$.

\begin{lemma} \label{BLemA4}
With the notation and hypotheses of Lemma $\ref{BLemA3}$, 
write $B = B(E_v,d)$ and $R = R(E_v,d)$.
Fix $\vp \in E_v^d$, and take $U = U_i$ and $G = G_i$
so that the assertions of Lemma $\ref{BLemA3}$ hold for $\vp$ with
respect to $U$ and $G$.  Put 
\begin{equation*}
r_0 \ = \ r_0(E_v,d) \ := \ \min(1/2, B(E_v,d)) \cdot R(E_v,d) \ < \ R \ = \ R(E_v,d) \ . 
\end{equation*} 
Then 

$(A)$  For each $\ell = 1, \ldots, d$, 
put $\cM_{\ell} = \cM_{\ell}(\vp) = \max_{z \in B(p_{\ell},R)} |G(z,\vp)|_v$.  
Then for each $\vq \in \prod_{j=1}^d B(p_j,r_0)$, and each $\ell$, 
we have 
\begin{equation*}
\max_{z \in B(q_{\ell},R)} |G(z,\vq)|_v \ = \ \cM_{\ell} \ .
\end{equation*}

$(B)$ For each $\ell = 1, \ldots, d$, there is a constant $C_{\ell} = C_{\ell}(\vp)$ 
with the following property:  For each $\vq \in \prod_{j=1}^d B(p_j,r_0)$, 
each $\ell$, and each $z \in B(q_{\ell},R)$,  
\begin{equation*}
|G(z,\vq)|_v \ = \ C_{\ell} \cdot \prod_{q_j \in B(q_{\ell},R)} \|z,q_j\|_v \ .
\end{equation*}
\end{lemma}

\begin{proof}
Note that if  $\vq \in \prod_{j=1}^d B(p_j,r_0)$, 
then $B(q_{\ell},R) = B(p_{\ell},R)$ for each $\ell = 1, \ldots, d$. 

For part (A), fix $\ell$, and note that 
$\cM_{\ell} = \max_{z \in B(p_{\ell},R)} |G(z,\vp)|_v \ge B(E_v,d)$.
Let $z_{\ell} \in B(p_{\ell},R)$ be a point where
$|G(z_{\ell},\vp)|_v = \cM_{\ell}$.
Choose isometric parametrizations of the balls
$B(p_{\ell},R)$ and $B(p_1,R), \ldots, B(p_d,R)$
in terms of local coordinate functions $Z, P_1, \ldots, P_d$ on
$D(0,R) = \{z \in \CC_v : |z|_v \le 1\}$, in such a way that
$z_{\ell} = Z(0)$, and $p_j = P_j(0)$.
Let $Q_1, \ldots, Q_d \in D(0,R)$ be such that $q_j = P_j(Q_j)$.
For each $j = 1, \ldots, d$, since $q_j \in B(p_j,r_0)$, 
the definition of isometric parametrizability shows that $|Q_j|_v = \|p_j,q_j\|_v \le r_0$.

Using these parametrizations, on $B(p_{\ell},R) \times \prod_{j=1}^d B(p_j,R)$ 
we can expand $G$ as a power series
\begin{equation*}
\cG(Z,\vP) \ = \ \sum_{i,k} a_{i,k} Z^i \vP^k \ .
\end{equation*}
Here $a_{0,0} = G(z_{\ell},\vp)$, so $|a_{0,0}| = \cM_{\ell}$.  Moreover,
since $\cG(Z,\vP)$ converges on $D(0,R) \times D(0,R)^d$ and
$|\cG(Z,\vP)|_v \le 1$ for all $(Z,\vP) \in D(0,R) \times D(0,R)^d$, we have 
\begin{equation*}
|a_{i,k}|_v \ \le \ \frac{1}{R^{i+|k|}}
\end{equation*}
for all $i$, $k$.  Consequently, for each $(i,k)$ with $|k| > 0$,
\begin{eqnarray*}
|a_{i,k} Z^i \vQ^k|_v
   & \le & \frac{1}{R^{i+|k|}} \cdot R^i \cdot r_0^{|k|}
                    \ = \ \Big(\frac{r_0}{R}\Big)^{|k|} \\
   & \le & \frac{r_0}{R} \ \le \ B \ < \  \cM_{\ell} \ = \ |a_{0,0}|_v \ .
\end{eqnarray*}

For each $z \in B(q_{\ell},R)$ we can write
\begin{equation*}
G(z,\vq) \ = \ \sum_{i=0}^{\infty} a_{i,0} Z^i +
            \sum_{i=0}^{\infty} \sum_{|k| > 0} a_{i,k} Z^i \vQ^k 
         \ = \ G(z,\vp) +
            \sum_{i=0}^{\infty} \sum_{|k| > 0} a_{i,k} Z^i \vQ^k \ ,
\end{equation*}
so 
\begin{equation*}
|G(z,\vq)|_v \ \le \ \max(|G(z,\vp)|_v,
         \max_i (\max_{|k|> 0} (|a_{i,k} Z^i \vQ^k|_v))) \ \le \ \cM_{\ell} \ .
\end{equation*}
On the other hand, when $z = z_{\ell}$ 
\begin{eqnarray*}
|G(z_{\ell},\vq)|_v 
    & = & |\cG(0,\vQ)|_v
      \ = \ |a_{0,0} + \sum_{i=1}^{\infty} \sum_{|k|>0} a_{0,k} \vQ^k|_v \\
    & = & |a_{0,0}|_v \ = \ \cM_{\ell} \ .
\end{eqnarray*}
This proves (A).  

For part (B), write 
\begin{equation*} 
\div_{\cCbar_v}(G(z,\vp))  =  \sum_{j=1}^d (p_j) - \sum_{j=1}^d (\delta_j) , \quad
\div_{\cCbar_v}(G(z,\vq))  =  \sum_{j=1}^d (q_j) - \sum_{j=1}^d (\Delta_j) \ .
\end{equation*} 
By the definition of $R$ in Lemma \ref{BLemA3}, we have
$|G(z,\vp)|_v \le 1$ and $|G(z,\vq)|_v \le 1$ on
$\cD := \bigcup_{j=1}^d B(p_j,R)$, so the $\delta_j$ and $\Delta_j$
lie outside $\cD$.  Fix a point $\zeta \in \cCbar_v(\CC_v) \backslash \cD$, distinct
from the $\delta_j$ and $\Delta_j$, and consider the canonical distance
$[z,w]_{\zeta}$.  By the factorization property of the canonical distance 
\index{canonical distance!factorization property}\index{canonical distance!$[z,w]_{\zeta}$}
(see \S\ref{Chap3}.\ref{CanonicalDistanceSection}), 
there are constants $C$ and $D$ such that for all $z \in \cCbar_v(\CC_v)$,
\begin{equation*}
|G(z,\vp)|_v =  C \cdot \frac{ \prod_{j=1}^d [z,p_j]_{\zeta} }
                            { \prod_{j=1}^d [z,\delta_j]_{\zeta} }, \qquad
|G(z,\vq)|_v =  D \cdot \frac{ \prod_{j=1}^d [z,q_j]_{\zeta} }
                            { \prod_{j=1}^d [z,\Delta_j]_{\zeta} } \ .
\end{equation*}
By Proposition \ref{APropA2}(B.2) (applied with $\fX = \{\zeta\}$), for each 
isometrically parametrizable ball $B(a,r_a)$ not containing $\zeta$,
there is a constant $c_a$ such that $[z,w]_{\zeta} = c_a \|z,w\|_v$
for all $z,w \in D(a,r_a)$.  By Proposition \ref{APropA2}(B.1),
if $B(a,r_a)$ and $B(b,r_b)$ are disjoint isometrically
parametrizable balls not containing $\zeta$, then $[z,w]_{\zeta}$ is constant
for $z \in B(a,r_a)$ and $w \in B(b,r_b)$.  It follows that there
are constants $C_{\ell}$ and $D_{\ell}$ such that for all
$z \in B(p_{\ell},R) = B(q_{\ell},R)$,
\begin{equation*}
|G(z,\vp)|_v
   =  C_{\ell} \cdot \prod_{p_j \in B(p_{\ell},R)} \|z,p_j\|_v \ , \quad
|G(z,\vq)|_v
   =  D_{\ell} \cdot \prod_{q_j \in B(q_{\ell},R)} \|z,p_j\|_v \ .
\end{equation*}
Clearly $|G(z,\vp)|_v$ achieves its maximum $\cM_{\ell}$ at a point 
$z_{\ell} \in B(p_{\ell},R)$ if and only if $\|z_{\ell},p_j\|_v = R$ 
for all $p_j \in B(p_{\ell},R)$, and then $\cM_{\ell}(\vp) = C_{\ell} R^{m_{\ell}}$ where 
$m_{\ell}$ is the number of points $p_j$ (or $q_j$) in $B(p_{\ell},R)$. 
Similarly $|G(z_{\ell},\vq)|_v$ 
achieves its maximum $\cM_{\ell} = D_{\ell} R^{m_{\ell}}$ 
on $B(q_{\ell},R) = B(p_{\ell},R)$ if and only
$\|z_{\ell},q_j\|_v = R$ for all $q_j \in B(q_{\ell},R)$.  
Since there are infinitely many  $z \in B(p_{\ell},R)$ 
satisfying both conditions simultaneously, we must have $C_{\ell} = D_{\ell}$.  
\end{proof}

\vskip .1 in 
Observe that in the notation of Lemma \ref{BLemA4}, 
for each $\vq \in \prod_{j=1}^d B(p_j,r_0)$ and each $\ell = 1, \ldots, d$, 
we have $\cM_{\ell}(\vq) = \cM_{\ell}(\vp) < 1$ and $C_{\ell}(\vq) = C_{\ell}(\vp)$.  

\begin{lemma} \label{BLemA5}
With the notation and hypotheses of Lemmas $\ref{BLemA3}$ and $\ref{BLemA4}$,
there is a constant $C = C(E_v,d)$ such that for each $\vp \in E_v^d$,
if $0 < r \le r_0$ belongs to the value group of $\CC_v$, and if $U = U_i$ and
$G = G_i$ are chosen for $\vp$ as in Lemma $\ref{BLemA4}$, then
\begin{equation*}
|G(z,\vp)|_v \ \ge \ C \cdot r^d
\end{equation*}
for all $z \in (\bigcup_{j=1}^d B(p_j,r)) \backslash
                    (\bigcup_{j=1}^d B(p_j,r)^-)$.
\end{lemma}

\begin{proof}
We can cover $E_v^d$ with finitely many sets of the form
$\prod_{j=1}^d B(p_j,r_0)$.  For each of these sets, Lemma \ref{BLemA4}
gives constants $C_{\ell} = C_{\ell}(\vp)$ such that for all 
$\vq \in \prod_{j=1}^d B(p_j,r_0)$, and all $z \in B(q_{\ell},r)$
\begin{equation*}
|G(z,\vq)|_v = C_{\ell} \cdot \prod_{q_j \in B(q_{\ell},r)} \|z,q_j\|_v \ .
\end{equation*}
Let $C$ be the minimum of these constants, for all the representative sets  
and all $\ell$. 

If $z \in B(p_{\ell},r) \backslash (\bigcup_{j=1}^d B(p_j,r)^-)$,
then $\|z,p_j\|_v = r$ for all $j$.  There are at most $d$ points $p_j$
in $B(p_{\ell},r)$, so $|G(z,\vp)|_v \ge C r^d$.
\end{proof}

\vskip .1 in
The lemma below uses the Maximum Modulus principle for $RL$-domains 
\index{Maximum principle!nonarchimedean!for $\RL$-domains}
\index{$\RL$-domain}
to control $|\frac{G(z,\vp)}{G(z,\vq)}-1|_v$ outside
$\bigcup_{j=1}^d B(q_j,r_j)^-$.

\begin{lemma} \label{BLemA6}
With the notation and hypotheses of Lemmas $\ref{BLemA3}$ and $\ref{BLemA4}$, 
there is a constant $D = D(E_v,d)$ with the following property. 
Let $0 < r \le r_0$ belong to the value group of $\CC_v^{\times}$.  
Suppose $\vp, \vq \in E_v^d$ are such that $\max_j (\|p_j,q_j\|_v) \le r$, 
and $\sum (p_j) \sim \sum (q_j)$.  Put $r_j = \|p_j,q_j\|_v$. 
Take $U = U_i$ and $G = G_i$ as in Lemma \ref{BLemA3}.  Then  for each
$z \in \cCbar_v(\CC_v) \backslash (\bigcup_{j=1}^d B(p_j,r)^-)$,
\begin{equation*}
\left| \frac{G(z,\vq)}{G(z,\vp)} - 1 \right|_v
       \ \le \ \frac{D}{r^d} \cdot \max_j(r_j) \ .
\end{equation*}
\end{lemma}

\begin{proof}
Suppose $G(z,\vp) = G_{k,\cF,U}(z,\vp)$ corresponds to a section $\cF$ and an affine set $U$ 
as in (\ref{BFLM1}). Fix $\vp$, $\vq$ as in the Lemma.
Noting that the polar divisors $s_{\cF}([\vp])$ and $s_{\cF}([\vq])$
of $G(z,\vp)$ and $G(z,\vq)$ depend only on the class
$[\vp] = [\vq] \in \Pic_{\cCbar_v/\CC_v}^d$, we have
\begin{equation*}
s_{\cF}([\vp]) \ = \ s_{\cF}([\vq]) \ = \ \sum_{j=1}^d (\Delta_j) \ .
\end{equation*}
Hence $G(z,\vq)/G(z,\vp)$ has poles only at the points $p_j \in \supp(\vp)$.

By Lemma \ref{BLemA1} and the Maximum Modulus Principle for $\RL$-domains 
\index{Maximum principle!nonarchimedean!for $\RL$-domains}
with boundary (see \cite{RR1}, Theorem 1.4.2, p. 51), 
it suffices to establish the bound in the Lemma for each
$z_0 \in (\bigcup_{j=1}^d B(p_j,r))
            \backslash (\bigcup_{j=1}^d B(p_j,r)^-)$.

Fix $\ell$, fix
$z_0 \in B(p_{\ell},r) \backslash (\bigcup_{j=1}^d B(p_j,r)^-)$,
and introduce local coordinate functions $Z,\vP$ on 
$B(p_{\ell},R)$ and the $B(p_j,R)$ as in Lemma \ref{BLemA4} 
so that $z_0 = Z(0)$, $p_j = P_j(0)$, 
and $q_j = P_j(Q_j)$.  On $B(p_{\ell},R) \times \prod_j B(p_j,R)$, 
expand $G_{k,\cF,U}$ as a power series
\begin{equation*}
\cG(Z,\vP) \ = \ \sum_{i,k} a_{i,k} Z^i \vP^k
\end{equation*}
where $|a_{i,k}|_v \le 1/R^{i+|k|}$ for all $i, k$.
Then $G(z_0,\vp) = \cG(0,\vec{0}) = a_{0,0}$, and
$G(z_0,\vq) = \cG(0,\vQ) = a_{0,0} + \sum_{|k| > 0} a_{0,k} \vQ^k$.  Hence
\begin{equation*}
  \frac{G(z_0,\vq)}{G(z_0,\vp)} - 1
      \ = \ \sum_{|k| > 0} \frac{a_{0,k}}{a_{0,0}} \vQ^k \ .
\end{equation*}
By Lemma \ref{BLemA5}, $|a_{0,0}|_v \ge C r^d$.  Hence, term by term,
\begin{eqnarray*}
\left| \frac{a_{0,k}}{a_{0,0}} \vQ^k \right|_v
   & \le & \frac{1}{Cr^d} \cdot \frac{1}{R^{|k|}} \cdot (\max_j(r_j))^{|k|} \\
   & \le & \frac{1}{Cr^d} \max(\frac{r_j}{R})
          \ = \ \frac{1}{CRr^d} \max(r_j) \ .
\end{eqnarray*}
This gives the result, with $D = 1/(C \cdot R)$.
\end{proof}

\vskip .1 in 

\begin{lemma} \label{BLemA7}
With the notation and hypotheses of Lemmas $\ref{BLemA3}$ and $\ref{BLemA4}$,
let $\vp, \vq \in E_v^d$ be such that $\max_j \|p_j,q_j\|_v \le r_0$,
and assume $\sum (p_j) \sim \sum (q_j)$.
Put $\|p_j,q_j\|_v = r_j$ for $j = 1, \ldots, d$. Then 
\begin{equation*}
\left| \frac{G(z,\vq)}{G(z,\vp)} \right|_v \ = \ 1 
\end{equation*}
for all $z \in \cCbar_v(\CC_v) \backslash
((\bigcup_{j=1}^d B(p_j,r_j)^-) \cup (\bigcup_{j=1}^d B(q_j,r_j)^-))$. 
\end{lemma}

\begin{proof}
As in Lemma \ref{BLemA6},  $G(z,\vp)$ and $G(z,\vq)$
have common polar divisor $\sum (\Delta_j)$.
Fix $\zeta \in \cCbar_v(\CC_v)\backslash ((\bigcup_{j=1}^d B(p_j,r_j)^-) 
\cup (\bigcup_{j=1}^d B(q_j,r_j)^-))$, distinct from the $\Delta_j$.
By the theory of the canonical distance there are constants
\index{canonical distance!factorization property}
$C_{\vp}$ and $C_{\vq}$ such that for all $z \in \cCbar_v(\CC_v)$,
\begin{equation*}
G(z,\vp) \ = \ C_{\vp} \cdot \frac{ \prod_{j=1}^d [z,p_j]_{\zeta} }
                            { \prod_{j=1}^d [z,\Delta_j]_{\zeta} }, \qquad 
G(z,\vq) \ = \ C_{\vq} \cdot \frac{ \prod_{j=1}^d [z,q_j]_{\zeta} }
                            { \prod_{j=1}^d [z,\Delta_j]_{\zeta} } \ .
\end{equation*}
Hence
\begin{equation*}
\left| \frac{G(z,\vq)}{G(z,\vp)} \right|_v \ = \
\frac{C_{\vp}}{C_{\vq}} \cdot \frac{ \prod_{j=1}^d [z,q_j]_{\zeta} }
                       { \prod_{j=1}^d [z,p_j]_{\zeta} } \ .
\end{equation*}
As noted in the proof of Lemma \ref{BLemA4},
$[z,w]_{\zeta}$ is constant for $z$ and $w$
belonging to isometrically parametrizable balls disjoint from each other
and from $\zeta$, while on a given isometrically parametrizable ball disjoint
from $\zeta$ it is a constant multiple of $\|z,w\|_v$.
Hence for $z \notin (\bigcup_{j=1}^d B(p_j,r_j)^-)
\cup (\bigcup_{j=1}^d B(q_j,r_j)^-) \cup \{\zeta\}$, we have
$[z,p_j]_{\zeta} = [z,q_j]_{\zeta}$ for each $j$,
and it follows that
\begin{equation*}
\left| \frac{G(z,\vq)}{G(z,\vp)} \right|_v \ = \ \frac{C_{\vp}}{C_{\vq}} \ .
\end{equation*}
This holds for $z = \zeta$ as well, by continuity, if we view
$G(z,\vq)/G(z,\vp)$ as a rational function with divisor 
$\sum (q_j) - \sum (p_j)$.
By the proof of Lemma \ref{BLemA4}(B) there are points
$z_{\ell} \notin \bigcup_{j=1}^d B(p_j,r_j)^- \cup
\bigcup_{j=1}^d B(q_j,r_j)^-$
where $|G(z_{\ell},\vp)|_v = |G(z_{\ell},\vq)|_v$, so $C_{\vp}/C_{\vq} = 1$,
and the result follows.
\end{proof}

\vskip .1 in
We can now prove the main theorem.

\vskip .1 in
\begin{proof}[Proof of Theorem \ref{BMT2} when $g > 0$.]  Let $r_0$ be as in Lemma \ref{BLemA4}.   
Given $\vp, \vq \in E_v^d$ with $\sum (p_m) \sim \sum (q_j)$ and 
$\max_j(\|p_j,q_j\|_v) \le r \le r_0$, choose $U = U_i$ and $G = G_i$ as in 
Lemma \ref{BLemA3}.  Then
\begin{equation*}
f(z,w;\vp,\vq) \ = \ \frac{G(z,\vp)}{G(z,\vq)} \cdot \frac{G(w,\vq)}{G(w,\vp)}
\end{equation*}
so (A) follows from Lemma \ref{BLemA7}.  Furthermore  
\begin{equation*}
f(z,w;\vp,\vq) - 1 \ = \
  \frac{G(z,\vp)}{G(z,\vq)} \cdot \left( \frac{G(w,\vq)}{G(w,\vp)} - 1 \right)
      + \left( \frac{G(z,\vp)}{G(z,\vq)} - 1 \right) \ .
\end{equation*}
Let $D$ be the constant from Lemma \ref{BLemA6}.  Then by Lemma \ref{BLemA6}, 
for all $z, w \notin 
 ((\bigcup_{j=1}^d B(p_j,r)^-) \cup (\bigcup_{j=1}^d B(q_j,r)^-))$
\begin{equation*}
\left| \frac{G(z,\vp)}{G(z,\vq)} - 1 \right|_v
                         \ \le \ \frac{D}{r^d} \max_j(r_j) \ , \qquad 
\left| \frac{G(w,\vq)}{G(w,\vp)} - 1 \right|_v
                         \ \le \ \frac{D}{r^d} \max_j(r_j) \ ,
\end{equation*}
while by Lemma \ref{BLemA7}
\begin{equation*}
\left|  \frac{G(z,\vp)}{G(z,\vq)} \right|_v \ = \ 1 \ .
\end{equation*}
Combining these gives (B).
\end{proof}
\index{universal function|)}

%% file: NewFSZAppD.tex
\chapter{ The Local Action of the Jacobian}
\label{AppD}

\index{Jacobian variety}\index{local action of the Jacobian|(}
\vskip .1 in
Let $K_v$ be a nonarchimedean local field, 
and suppose $\cC_v/K_v$ is a smooth, connected, projective curve of genus $g > 0$. 
Write $\Jac(\cC_v)$ for the Jacobian of $\cC_v$ over $K_v$. 
\index{Jacobian variety}

In this Appendix we will show that 
for a dense set of points $\va \in \cC_v(\CC_v)^g$, there is an action 
of a neighborhood of the origin in $\Jac(\cC_v)(\CC_v)$ 
\index{Jacobian variety}
on a sufficiently small neighborhood of $\va$ in $\cC_v(\CC_v)^g$, 
which makes that neighborhood into a principal homogeneous space.\index{principal homogeneous space}
This action is used in \S\ref{Chap6}.\ref{NonArchProofSection} in    
the construction of the initial local approximating 
\index{initial approximating functions $f_v(z)$!nonarchimedean}
functions in the nonarchimedean compact case, 
and in \S\ref{Chap11}.\ref{NonArchPatchingProof},
\index{patching argument!for nonarchimedean $K_v$-simple sets} 
the patching process in the nonarchimedean compact case, 
in moving the roots of the partially patched function away from each other.

\vskip .1 in
Fix a spherical metric $\|x,y\|_v$ on $\cC_v(\CC_v)$.\index{spherical metric!on curve}    
We will be working simultaneously with balls in $\cC_v(\CC_v)$,  
$\cC_v(\CC_v)^g$, and $\Jac(\cC_v)(\CC_v)$, so we will write $B_{\cC_v}(a,r)$ for the ball
\index{Jacobian variety}
$\{ z \in \cC_v(\CC_v) : \|z,a\|_{\cC_v,v} \le r \}$  
simply denoted $B(a,r)$ elsewhere.  As usual, we put $D(0,r) = \{z \in \CC_v : |z|_v \le r \}$
and $D(0,r)^- = \{z \in \CC_v : |z|_v < r \}$.
Given a point $\va = (a_1, \ldots, a_g) \in \cC_v(\CC_v)^g$, and radii  $r_1, \ldots, r_g > 0$, 
we write $\BCg(\va,\vr) := \prod_{i=1}^g B_{\cC_v}(a_i,r_i)$. 
We also put $D(\vORIG,\vr) = \prod_{i=1}^g D(0,r_i)$, and if $r_1 = \cdots = r_g = R$ we write 
$D(\vORIG,R) = \prod_{i=1}^g D(0,R)$, $D(\vORIG,R)^- = \prod_{i=1}^g D(0,R)^-$. 

If $r_1\ldots, r_g > 0$ are small enough, then by Theorem \ref{IsoParamThm} 
(proved in \cite{RR1}, Theorem 1.2.3)
each ball $B_{\cC_v}(a_i,r_i)$ can be isometrically parametrized by power series, that is, 
there is an analytic isomorphism $\varphi_i :  D(0,r_i) \rightarrow B_{\cC_v}(a_i,r_i)$
defined by convergent power series 
(which are $F_u$-rational provided $a_1, \ldots, a_g \in \cC_v(F_u)$, 
where $K_v \subseteq F_u \subseteq \CC_v$ and $F_u$ is complete), such that
$\|\varphi_i(x),\varphi_i(y)\|_{\cC_v,v} = |x-y|_v$ for all $x,y \in D(0,r_i)$.  It follows that
\begin{equation} \label{FFdef} 
\Phi_{\va} \ := \
   (\varphi_1,\ldots,\varphi_g) :  D(\vORIG,\vr) \rightarrow \BCg(\va,\vr) 
\end{equation}
is an analytic isomorphism.  The construction in (\cite{RR1}, Theorem 1.2.3) shows that the maps 
$\varphi_i$ can be chosen in such a way that $\varphi_i^{-1} : B(a_i,r_i) \rightarrow D(0,r_i)$ is projection 
on one of the coordinates, followed  by a translation, and we will assume that that is the case.  

\smallskip
The Jacobian $\Jac(\cC_v)/K_v$ is an abelian variety characterized by the property that  
\index{Jacobian variety}
$\Jac(\cC_v) \times_{K_v} \Spec(F_u)$ becomes isomorphic to $\BPic^0_{\cC_u/F_u}$ 
over any extension $F_u/K_v$ such that $\cC_v(F_u) \ne \phi$. 
For each $\va = (a_1, \ldots, a_g) \in \cC_v(\CC_v)^g$, there is a morphism 
$\crJ_{\va} : \cC_v^g \rightarrow \Jac(\cC_v)$, defined over $K_v(a_1, \ldots, a_g)$,
\index{Jacobian variety}
which takes $\vx = (x_1, \ldots, x_g)$ to the linear equivalence class of the divisor 
$\sum (x_i) - \sum (a_i)$.  
It induces a birational morphism from $\Sym^{(g)}(\cC_v)$ onto $\Jac(\cC_v)$. 
\index{Jacobian variety}
This was the idea behind Weil's algebraic construction 
\index{Weil, Andr\'e}
of $\Jac(\cC_v)$:  using the Riemann-Roch theorem, he showed that there was a birational,
\index{Jacobian variety}
\index{Riemann-Roch theorem}   
commutative law of composition defined on an open subset of $\Sym^{(g)}(\cC_v)$ (a `group chunk'),
\index{group chunk}  
which could be extended to an addition law on an abelian variety.  Later, he showed that 
every abelian variety is projective, and Matsusaka showed that $\Jac(\cC_v)$ and its group law
\index{Jacobian variety}
\index{Matsusaka, Teruhisa}
were defined over $K_v$.  A modern account of this theory can be found in (\cite{Milne}).
\index{Milne, James} 

Since $\Jac(\cC_v)$ is smooth and projective,  each point of $\Jac(\cC_v)(\CC_v)$ 
\index{Jacobian variety}
has a neighborhood in the $v$-topology which is isometrically parametrizable by power series 
(Theorem \ref{IsoParamThm}).  If $x \in \Jac(\cC_v)(F_u)$, 
where $K_v \subseteq F_u \subseteq \CC_v$ and $F_u$ is complete, 
then those power series are defined over $F_u$.  It follows that each point of $\Jac(\cC_v)(F_u)$ 
has a neighborhood in $\Jac(\cC_v)(F_u)$ analytically isomorphic to $\cO_u^g$.

\smallskip
Write $\JacNer(\cC_v)$ for the N\'eron model of $\Jac(\cC_v)$.  
\index{Jacobian variety}
\index{N\'eron model|ii} 
\label{`SymbolIndexNeron'}
The N\'eron model (see \cite{Artin}, \cite{BLR})
is a smooth, separated group scheme of finite type over $\Spec(\cO_v)$ 
whose generic fibre is isomorphic to $\Jac(\cC_v)$,  
characterized by the property that each point of $\Jac(\cC_v)(K_v)$ 
extends to a section of $\JacNer(\cC_v)/\Spec(\cO_v)$.  
By (\cite{BLR}, Theorem 1, p.153), $\JacNer(\cC_v)$ is quasi-projective.   
Let it be embedded in $\PP^N/\cO_v$, for an appropriate $N$.  
Fix a corresponding system of homogeneous coordinates on $\PP^N/K_v$. 
We will identify $\Jac(\cC_v)$ with generic fibre of $\JacNer(\cC_v)$, 
\index{Jacobian variety}
viewing it as locally cut out of $\PP^N$ by the equations defining $\JacNer(\cC_v)$.   

Let $\|x,y\|_{J,v}$ be the induced spherical metric on $\Jac(\cC_v)(\CC_v)$,
\index{Jacobian variety}\index{spherical metric!on Jacobian}  
and let $O$ be the origin of $\Jac(\cC_v)$.  
Then the ball $B_J(O,1)^- := \{z \in \Jac(\cC_v)(\CC_v) : \|z,O\|_{J,v} < 1\}$ is a subgroup.
Since $O$ is nonsingular on the special fibre of the N\'eron model, 
$B_J(O,1)^-$ can be isometrically parametrized by parametrized by power series  
converging on $D(\vORIG,1)^-$, taking $\vORIG$ to $O$ 
(Theorem \ref{IsoParamThm}). Let 
\begin{equation} \label{JPsiDef}
\Psi : D(\vORIG,1)^- \ \rightarrow \ B_J(O,1)^- \ .
\end{equation}
be such an isometric parametrization.  By the construction in (\cite{RR1}, Theorem 1.2.3) 
we can assume $\Psi$ has been chosen in such a way that $\Psi^{-1} : B_J(O,1)^- \rightarrow D(\vORIG,1)^-$ 
is projection on some of the coordinates.
Pulling the group action back to $D(\vORIG,1)^-$ using $\Psi$ 
yields the formal group of $\Jac(\cC_v)$ over $\cO_v$. 
\index{Jacobian variety}
\index{formal group}
 
Writing $\vX = (X_1, \ldots, X_g)$ and $\vY = (Y_1, \ldots, Y_g)$, let $S(\vX,\vY) \in \cO_v[[\vX,\vY]]^g$ 
and $M(\vX) \in \cO_v[[\vX]]^g$ be the vectors of power series defining addition and negation 
in the formal group.
\index{formal group}  
Since  $S(\vX,\vORIG) = \vX$, $S(\vX,\vY) = S(\vY,\vX)$, and $S(\vX,M(\vX)) = \vORIG$, 
modulo terms of degree $\ge 2$ we have 
\begin{equation} \label{SumMinusSeries} 
S(\vX,\vY)  \equiv  \vX + \vY , \qquad  M(\vX)  \equiv  -\vX \ .  
\end{equation}  
The following facts are well known:

\begin{proposition} \label{FormalGroupProperties}  
Let $p$ be the residue characteristic of $K_v$.  Then 

$(A)$ For each $0< r < 1$ the ball 
\begin{equation*}
B_J(O,r) \ := \ \{ z \in \Jac(\cC_v)(\CC_v): \|z,O\|_{J,v} \le r \}
\end{equation*}
is an open subgroup of $B_J(O,1)^-$.  

$(B)$ $B_J(O,1)^-$ is a topological pro-$p$-group.\index{pro-$p$-group}

$(C)$ There is an $R > 0$ such that $B_J(O,R)$ is torsion-free.  

\end{proposition}  

\begin{proof}  For the convenience of the reader we recall the proofs.
It suffices to prove the assertions for $D(\vORIG,1)^-$ with the group law defined by $S(\vX,\vY)$.  
Let $[2](\vX) = S(\vX,\vX)$ and inductively put $[n](\vX) = S(\vX,[n-1](\vX))$ for $n = 3, 4, \ldots$.    

For (A), fix $0 < r < 1$ and suppose $\vx,\vy \in D(\vORIG,r)$.  
It follows easily from (\ref{SumMinusSeries}) that $S(\vx,\vy)$ and $M(\vx,\vy)$ belong to $D(\vORIG,r)$.  

For (B), note that by (\ref{SumMinusSeries}) we have $[p](\vX) \equiv p\vX$
modulo terms of degree $\ge 2$.  If $\vx \in D(\vORIG,1)^-$, then $\vx \in D(\vORIG,r)$ for some $r < 1$, 
and each term in the series defining $[p](\vx)$ has absolute value at most $R = \max(|p|_v r, r^2)$, 
so $[p](\vx) \in D(\vORIG,R)$.  Iterating this we see that 
\begin{equation*} 
\lim_{k \rightarrow \infty} [p^k](\vx) \ = \ \vORIG \ .
\end{equation*} 

For (C), note that by (A) and (B) above, $B_J(O,1)^-$ can only have $p$-power torsion.  
By the general theory of abelian varieties, $\Jac(\cC_v)(\CC_v)$ has at most $p^{2g}-1$ elements of order $p$,
\index{Jacobian variety}
so the same is true for $B_J(O,1)^-$.   If $R > 0$ is small enough, then $B_J(O,R)$ 
has no elements of order $p$, and hence is torsion free.  
\end{proof} 

\noindent{\bf Remark.}  
When $\Char(K_v) = 0$, it follows from the existence of the $v$-adic `logarithm map' 
(see \cite{LALG}, Corollary 4, p.LG5.36)
that there is a subgroup $B_J(0,R)$ analytically isomorphic to the additive group $\hcO_v^g$.  
However, when $\Char(K_v) = p > 0$, no logarithm map exists, 
and no subgroup $B_J(0,R)$ can be isomorphic to $\hcO_v^g$ with the additive structure induced from $\CC_v$, 
since $B_J(0,r)$ is torsion-free for all small $r$, while $\hcO_v^g$ is purely $p$-torsion.  
Nonetheless, by considering the form of the power series $S(\vX,\vY)$ and $M(\vX)$, one sees easily 
that for any $0 \ne \pi \in \hcO_v$, there is an $R_0 > 0$ such that if $0 < R \le R_0$
and $R \in |\CC_v^{\times}|_v$, then $B_J(O,R)/B_J(O,|\pi|_v R)$ is isomorphic to $\hcO_v^g/\pi \hcO_v^g$.    

\section{  The Local Action of the Jacobian on $\cC_v^g$ } \label{LocalActionThmSection} 

Write $\cCbar_v = \cC_v \times_{K_v} \Spec(\CC_v)$.  
Then $\Jac(\cCbar_v) \cong \Jac(\cC_v) \times_{K_v} \Spec(\CC_v)$.  We will identify
\index{Jacobian variety}
$\cC_v(\CC_v)$ with $\cCbar_v(\CC_v)$, and $\Jac(\cC_v)(\CC_v)$ with $\Jac(\cCbar_v)(\CC_v)$.   

Let $\Div_{\cCbar_v/\CC_v}(\CC_v)$ be the divisor group of $\cCbar_v$,
and let $\sim$ denote the relation of linear equivalence for divisors.
 Let $\Pic_{\cCbar_v/\CC_v}(\CC_v) = \Div_{\cCbar_v/\CC_v}(\CC_v)/\sim$ 
be the relative Picard group (see the discussion after Theorem \ref{UniversalFcn} 
\index{Picard group!relative}  
in Appendix C),  and let $\BPic_{\cCbar_v/\CC_v}$ be the associated Picard scheme. 
\index{Picard scheme}  
Let $\BPic_{\cCbar_v/\CC_v}^\nu$ be its degree $\nu$ component, regarded as an algebraic variety.  
Then $\BPic_{\cCbar_v/\CC_v}^{0} \cong \Jac(\cCbar_v)$ as a  group scheme, 
\index{Jacobian variety}
and $\BPic_{\cCbar_v/\CC_v}^\nu(\CC_v)$ is a principal homogeneous space\index{principal homogeneous space} 
for $\Jac(\cCbar_v)(\CC_v) = \Jac(\cC_v)(\CC_v)$. 

On the product $\cCbar_v^g$
we have the cycle class map \ $[\ ] : \cCbar_v^g \rightarrow \BPic_{\cCbar_v/\CC_v}^{g}$,
defined on $\cC_v(\CC_v)^g$ by 
\begin{equation*}
[\vx] \ = \ [(x_1, \ldots, x_g)] \ = \ ((x_1) + \ldots + (x_g))/\sim \ .
\end{equation*}
In the notation of Appendix C, $[\vx] = P \circ Q(\vx)$ where
$Q : \cCbar_v^g \rightarrow \Sym^{(g)}(\cCbar_v) \cong \BDiv^g_{\cCbar_v/\CC_v}$ 
is the quotient by the symmetric group $S_g$, and
$P : \Sym^{(g)}(\cCbar_v) \rightarrow \BPic_{\cCbar_v/\CC_v}^{g}$ 
is the Abel map\index{Abel map} $P(\sum_{i=1}^g (x_i))) = [\sum_{i=1}^g (x_i)]$.  
The morphism $Q$ is flat and finite of degree $g!$,
and the morphism $P$ is a birational isomorphism.

Let $\Jp$ and $\Jm$ denote addition and subtraction
under the group law on $\BPic_{\cCbar_v/\CC_v}(\CC_v)$, and by restriction, on $\Jac(\cC_v)(\CC_v)$.    
For each $\va \in \cC_v(\CC_v)^g$,  
\index{Jacobian variety}
let $\crJ_{\va} : \cC_v(\CC_v)^g \rightarrow \Jac(\cC_v)(\CC_v)$ be the map 
\begin{equation} \label{LambdaDef}
\crJ_{\va}(\vx) \ = \ [\vx] \Jm [\va] \ .
\end{equation} 

Now suppose $F_u/K_v$ is a separable finite extension, 
and let $H_w$ be the galois closure of $F_u$ over $K_v$.
Put $d = [F_u:K_v]$, and let $\sigma_1, \ldots, \sigma_d$ be the distinct embeddings of $F_u$ into $H_w$.
Extend each $\sigma_i$ to an automorphism of $H_w$. 
Since the addition law $\Jp$ in $\Jac(\cC_v)$ is defined over $K_v$, for each $x \in \Jac(\cC_v)(F_u)$
the trace  $\Tr_{F_u/K_v} : \Jac(\cC_v)(F_u) \rightarrow \Jac(\cC_v)(K_v)$ is given by  
\index{Jacobian variety}
\begin{equation*} 
Tr_{F_u/K_v}(x) \ = \ \sigma_1(x) \Jp \sigma_2(x) \Jp \cdots \Jp \sigma_d(x) \ .
\end{equation*}  

\smallskip
Our main result is as follows:\index{local action of the Jacobian|ii}   

\begin{theorem}  \label{BKeyThm1}
Let $K_v$ be a nonarchimedean local field, 
and let $\cC_v/K_v$ be a smooth, projective, 
geometrically integral curve of genus $g > 0$.  
Then the points $\va = (a_1, \ldots, a_g) \in \cC_v(\CC_v)^g$ such that 
$\crJ_{\va} : \cC_v(\CC_v)^g \rightarrow \Jac(\cC_v)(\CC_v)$
is nonsingular at $\va$ are dense in $\cC_v(\CC_v)^g$ for the $v$-topology.
If $F_u/K_v$ is a finite extension and $\cC_v(F_u)$ is nonempty, 
they are dense in $\cC_v(F_u)^g$.    
  
Fix such an $\va;$  then $a_1, \ldots, a_g$ are distinct, 
and for each $0 < \eta < 1$, there is a number $0 < R < 1$ 
$($depending on $\va$ and $\eta)$ such that the balls $B_{\cC_v}(a_1,R), \ldots, B_{\cC_v}(a_g,R)$ 
are pairwise disjoint and isometrically parametrizable,
and for each $\vr = (r_1, \ldots, r_g)$ satisfying 
\begin{equation} \label{FBCond} 
0 < r_1, \ldots, r_g \le R \quad \text{and} \quad \eta \cdot \max(r_i) \le \min(r_i) \ ,
\end{equation}

$(A)$ {\rm $($Subgroup$)$} The map $\crJ_{\va} : \cC_v(\CC_v)^g \rightarrow \Jac(\cC_v)(\CC_v)$ 
is injective on $\prod_{i=1}^g B_{\cC_v}(a_i,r_i)$, 
and the image $W_{\va}(\vr) := \crJ_{\va}\big(\prod_{i=1}^g B_{\cC_v}(a_i,r_i)\big)$ 
is an open subgroup of $\Jac(\cC_v)(\CC_v)$. 
\index{Jacobian variety} 

$(B)$ {\rm $($Limited Distortion$)$} 
For $i = 1, \ldots g$, let $\varphi_i : D(0,R) \rightarrow B(a_i,R)$
be isometric parametrizations with $\varphi_i(0) = a_i$.
Given $0 < r_1, \ldots, r_g \le R$, 
let $\Phi_{\va} = (\varphi_1,\ldots,\varphi_g) : \prod_{i=1}^g D(0,r_i) \rightarrow \prod_{i=1}^g B(a_i,r_i)$ 
be the associated map. Let $\Psi : D(\vORIG,1)^- \ \rightarrow \ B_J(O,1)^-$ 
be the isometric parametrization inducing the formal group,  
\index{formal group}
and let $L_{\va} : \CC_v^g \rightarrow \CC_v^g$ be the linear map 
$(\Psi^{-1} \circ \crJ_{\va} \circ \Phi_{\va})^{\prime}(\vORIG)$. 

Then $W_{\va}(\vr) = \Psi(L_{\va}(D(\vORIG,\vr)))$.  Giving $D(\vORIG,\vr)$ its structure as an additive
subgroup of $\CC_v^g$, the map $\Psi \circ L_{\va}$ induces an isomorphism of groups 
\begin{equation*} 
D(\vORIG,\vr)/D(\vORIG,\eta \vr) \ \cong \ W_{\va}(\vr)/W_{\va}(\eta \vr) 
\end{equation*} 
with the property that for each $\vx \in D(\vORIG,\vr)$, 
\begin{equation*} 
\crJ_{\va}(\Phi_{\va}(\vx)) \ \equiv \ \Psi(L_{\va}(\vx))  \pmod{W_{\va}(\eta \vr)} \ . 
\end{equation*} 

$(C)$ {\rm $($Action$)$} There is an action $(w,\vx) \mapsto w \ap \vx$
of $W_{\va}(\vr)$ on $\prod_{i=1}^g B_{\cC_v}(a_i,r_i)$ 
which makes $\prod_{i=1}^g B_{\cC_v}(a_i,r_i)$ into a principal homogeneous space\index{principal homogeneous space} 
for $W_{\va}(\vr)$. 
It is defined by $w \ap \vx = \crJ_{\va}^{-1}(w \Jp \crJ_{\va}(\vx))$ 
if we restrict the domain of $\crJ_{\va}$ to $\prod_{i=1}^g B_{\cC_v}(a_i,r_i)$, 
and has the property that for each $w \in W_{\va}(\vr)$   
and each $\vx \in \prod_{i=1}^g B_{\cC_v}(a_i,r_i)$,
\begin{equation} \label{FNat1}
[w \ap \vx] \ = \ w \Jp [\vx] \ .
\end{equation}

$(D)$ {\rm $($Uniformity$)$} For each $\vb \in \prod_{i=1}^g B_{\cC_v}(a_i,r_i)$, 
\begin{equation} \label{FCenter1}  
W_{\va}(\eta \vr) \, \ap \, \vb  \ = \ \prod_{i=1}^g B_{\cC_v}(b_i,\eta r_i) 
\quad \text{and} \quad 
\crJ_{\vb}\Big(\prod_{i=1}^g B_{\cC_v}(b_i,\eta r_i)\Big) \ = \ W_{\va}(\eta \vr)\ . 
\end{equation} 

$(E)$ {\rm $($Rationality$)$} If $F_u/K_v$ is a finite extension, and $\va \in \cC_v(F_u)^g$,  then   
\begin{eqnarray}
\crJ_{\va}(\prod_{i=1}^g (B_{\cC_v}(a_i,r_i) \cap \cC_v(F_u))) \ = \ 
W_{\va}(\vr) \cap \Jac(\cC_v)(F_u) \ , \label{FRat1A} \\
(W_{\va}(\vr) \cap \Jac(\cC_v)(F_u)) \, \ap \, \va 
\ = \ \prod_{i=1}^g \big(B_{\cC_v}(a_i,r_i) \cap \cC_v(F_u)\big) \ . \label{FRat2A}
\end{eqnarray}

$(F)$ {\rm $($Trace$)$} If $F_u/K_v$ is finite and separable, 
there is a constant $C = C(F_u,\va) > 0$, depending on $F_u$ and $\va$ but not on $\vr$, such that 
if\, $r = \min_i(r_i)$ then 
\begin{equation} \label{FTraceContains} 
B_J(O,C r) \cap \Jac(\cC_v)(K_v) \ \subseteq \ \Tr_{F_u/K_v}\big(W_{\va}(\vr) \cap \Jac(\cC_v)(F_u)\big)  \ .
\end{equation} 
\end{theorem}

For the proof of the Fekete-Szeg\"o theorem we will need one more property of the subgroups 
\index{Fekete, Michael}\index{Szeg\"o, G\'abor} 
$W_{\va}(\vr)$, asserting Lipschitz continuity\index{Lipschitz continuity!of the Abel map|ii}
 of the Abel maps\index{Abel map} $\BFj_x(z) = [(z)-(x)]$ 
on compact sets $E_v \subset \cC_v(\CC_v)$. 

\begin{proposition} \label{BPropD2New} Let $K_v$ be a nonarchimedean local field,
and let $\cC_v/K_v$ be a smooth, projective, geometrically integral curve of genus 
$g = g(\cC_v) > 0$.  Let $E_v \subset \cC_v(\CC_v)$ be compact. 
Suppose that $\va \in E_v^g$ is a point such that 
$\BFJ_{\va} : \cC_v(\CC_v)^g \rightarrow \Jac(\cC_v)(\CC_v)$ is nonsingular at $\va$, 
and let $0 < \eta < 1$, $0 < R < 1$ and $0 < r_1, \ldots, r_g \le R$ be numbers satisfying  
the conditions $(\ref{FBCond})$,  so $W_{\va}(\vr) = \BFJ_{\va}\big(\prod_{i=1}^g B_{\cC_v}(a_i,r_i)\big)$ 
is an open subgroup of $\Jac(\cC_v)(\CC_v)$ with the properties in Theorem $\ref{BKeyThm1}$.
\index{Jacobian variety}

Then there are constants $\varepsilon_0, C_0 > 0$ $($depending on $\va$ and $E_v)$
such that if $0 < \varepsilon \le \varepsilon_0$, 
then for all $x, z \in E_v$ with $\|x,z\|_v \le \varepsilon$, 
the divisor class $\BFj_{x}(z) = [(z)-(x)]$ belongs to $W_{\va}( C_0 \varepsilon \cdot \vr)$.  
\end{proposition}

The proofs, given in \S\ref{AppD}.\ref{LocalActionSection}, 
involve expanding the maps in question in terms of power series, 
and applying properties of power series proved below.  


\section{  Lemmas on Power Series in Several Variables }
   \label{PowerSeriesLemmasSection} 
 
In this section we recall some facts about power series in several variables.
All the results are standard.  
Fix $0 < d \in \NN$;  in the application we will take $d = g$.

Given variables $X_1, \ldots, X_d$  
and natural numbers $k_1, \ldots, k_d$,
we write $\vX = (X_1, \ldots, X_d)$, $k = (k_1,\ldots,k_d)$, 
and $\vX^{k} = X_1^{k_1} \cdots X_d^{k_d}$.  
Put $\vORIG = (0,\ldots, 0) \in \CC_v^d$. 
If  $\va = (a_1, \ldots, a_d) \in \CC_v^d$ and $\vr = (r_1, \ldots, r_d) \in \RR^d$, 
with $r_1, \ldots, r_d > 0$, we write 
\begin{equation*}
D(\va,\vr) \ = \ \prod_{i=1}^d D(a_i,r_i) \ ,
\qquad D(\va,\vr)^- = \prod_i D(a_i,r_i)^- \ . 
\end{equation*} 
If $r > 0$ we put $D(\vORIG,r) = D(0,r)^d$ and $D(\vORIG,r)^- = (D(0,1)^-)^d$.
The norm $|\vx|_v = \max_i(|x_i|_v)$ induces a metric $|\vx-\vy|_v$ on $\CC_v^d$ 
and on each polydisc $D(\va,\vr)$, $D(\va,\vr)^-$.       

\vskip .1 in 
First, recall that a power series $g(\vX) \in \CC_v[[\vX]]$  
which converges on a polydisc $D(\vORIG,\vr)$   
is determined by its values on $D(\vORIG,\vr) \cap K_v^d$:  

\begin{lemma} \label{BVLem} 
Suppose $g(\vX) \in \CC_v[[\vX]]$ converges on $D(\vORIG,\vr)$, with 
$g(\va) = 0$ for each $\va \in D(\vORIG,\vr) \cap K_v^d$.  
Then $g(\vX)$ is the zero power series.
\end{lemma}

\begin{proof}
If $g(\vX) \not \equiv 0$,  
after making a $K_v$-rational change of variables and 
shrinking the $r_i$ if necessary, we can apply the Weierstrass Preparation Theorem
\index{Weierstrass Preparation Theorem}  
(see \cite{BGR}, Theorem 1, p.201, and Proposition 2, p.205).
This means we can factor $g(\vX)$ as $G(\vX) h(\vX)$ where $h(\vX)$ 
is an invertible power series converging in $D(\vORIG,\vr)$ and 
\begin{equation*}
G(\vX) = X_d^M + \sum_{i=1}^M c_i(X_1, \ldots, X_{d-1}) X_d^{M-i}
\end{equation*}
is a monic polynomial in $X_d$ with coefficients  
$c_i(X_1, \ldots, X_{d-1}) \in \CC_v[[X_1, \ldots, X_{d-1}]]$ which converge in 
$D(\vORIG,(r_1, \ldots, r_{d-1}))$.  Since $g(\vX)$ vanishes on 
$D(\vORIG,\vr) \cap K_v^d$, so does $G(\vX)$.  Fixing 
$a_1, \ldots, a_{d-1} \in K_v$ with $|a_i|_v \le r_i$ for each $i$, 
we see that $G(a_1, \ldots, a_{d-1}, X_d)$ is a monic polynomial in $X_d$ 
with infinitely many roots.  

This is impossible, so $g(\vX) \equiv 0$.
\end{proof}

\vskip .1 in
Suppose $h(\vX) = \sum a_{k} \vX^{k} \in \CC_v[[\vX]]$  
converges on the polydisc $D(\vORIG,1)$.
It follows from the Maximum Modulus Principle (see \cite{BGR}, p.201) that
\index{Maximum principle!nonarchimedean!for power series}
\begin{equation}  \label{BCAF1}
\|h\|_{D(\vORIG,1)} \ = \ \max_{k} (|a_{k}|_v) \ .
\end{equation} 
Here $\|h\|_{D(\vORIG,1)} = \sup_{\vz \in D(\vORIG,1)}(|h(\vz)|_v)$ is the
$\sup$ norm of $h$, and $\max_{k} (|a_{k}|_v)$ is the so-called `Gauss norm'. 
\index{Gauss norm}
If $h(\vX)$ converges on the polydisc $D(\vORIG,\vr)$, an analogous result holds:
\begin{equation} \label{BCAF2}
\|h\|_{D(\vORIG,\vr)} \ = \ \max_{k} (|a_{k}|_v \vr^{k}) \ .
\end{equation}   
 
Uniform convergence of values implies convergence
in the Gauss norm:  
\index{Gauss norm}

\begin{lemma}  \label{BCL1}  Fix a polydisc $D(\vORIG,\vr)$.  
Let $g^{(\ell)}(\vX) \in \CC_v[[\vX]]$, for $\ell = 1, 2, \ldots$,
be power series for which there is a function $G(\vX)$ on $D(\vORIG,\vr)$ 
such that uniformly for $\va \in D(\vORIG,\vr)$,
the values $g^{(\ell)}(\va)$ converge to $G(\va)$.  
Then there is a power series $g(\vX) \in \CC_v[[\vX]]$ converging on $D(\vORIG,\vr)$
such that the coefficients of the $g^{(\ell)}(\vX)$  converge to the coefficients
of $g(\vX)$, and  $g(\va) = \lim_{\ell \rightarrow \infty} g^{(\ell)}(\va) = G(\va)$
for each $\va \in D(\vORIG,\vr)$.
\end{lemma}

\begin{proof} 
Write $g^{(\ell)}(\vX) = \sum_k b_{k}^{(\ell)} \vX^{k}$.
By the Maximum Modulus principle, applied to $D(\vORIG,\vr)$,
\index{Maximum principle!nonarchimedean!for power series}
for each $k \in \NN^d$ the coefficients $b_{k}^{(\ell)}$ converge
to a number $b_{k}$, and if $g(\vX) = \sum_{k} b_{k} \vX^{k}$
then in the Gauss norm for $D(\vORIG,\vr)$ the $g^{(\ell)}(\vX)$ converge to 
\index{Gauss norm}
$g(\vX)$. But convergence in the Gauss norm implies convergence 
of values on $D(\vORIG,\vr)$.  
\end{proof}

\vskip .1 in 
Recall that $\hcO_v$ denotes the ring of integers of $\CC_v$.  

\begin{proposition} \label{BCProp1}  Let $H(\vX) = (h_1(\vX), \ldots, h_d(\vX)) \in \hcO_v[[\vX]]^d$ 
be such that $H(\vX) \equiv \vX$ modulo terms of degree $\ge 2$.  Then 

$(A)$ $H$ induces an isometry from $D(\vORIG,1)^-$ onto $D(\vORIG,1)^-$. 

$(B)$  For each polydisc $D(\vORIG,\vr) \subset D(\vORIG,1)^-$ satisfying $\max_i(r_i)^2 \le \min_i(r_i)$, 
$H$ induces an isometry from $D(\vORIG,\vr)$ onto $D(\vORIG,\vr)$, and 
if $F_u \subseteq \CC_v$ is a complete field such that $H(\vX)$ is rational over $F_u$,
then $H(D(\vORIG,\vr) \cap F_u^d) = D(\vORIG,\vr) \cap F_u^d$.
%
%
\end{proposition}  

\begin{proof}
First, note that the form of $H(\vX)$ shows that $H$ converges on $D(\vORIG,1)^-$ and that 
$H(D(\vORIG,1)^-) \subset D(\vORIG,1)^-$.

Next we claim that $H$ preserves distances on $D(\vORIG,1)^-$.
To see this, fix $\vp, \vq \in D(\vORIG,1)^-$ with $\vp \ne \vq$.
Write $h_i(\vX) = \sum_{k} a_{i,k} \vX^{k}$,
and put  $|k| = k_1 + \ldots + k_d$.
For each $i$ we have
\begin{equation*}
h_i(\vp) - h_i(\vq)
\ = \ (p_i - q_i) + \sum_{|k| \ge 2} a_{i,k} (\vp^{k} - \vq^{k}) \ .
\end{equation*}  
Using the ultrametric property of $\CC_v$, 
it is easy to see that for each $k$
\begin{equation} \label{BCF1}
|a_{i,k} (\vp^{k} - \vq^{k})|_v \ \le \
     (\max_j |p_j -q_j|_v) \cdot \big(\max(|\vp|_v,|\vq|_v)\big)^{|k|-1} \ .
\end{equation}
Consequently $|H(\vp)-H(\vq)|_v \le \max_i |p_i-q_i|_v = |\vp-\vq|_v$.  

Now take $i$ so that $|p_i - q_i|_v$ is maximal.
Noting that $\max(|\vp|_v,|\vq|_v) < 1$, it follows from (\ref{BCF1}) that
$|p_i - q_i|_v > |a_{i,k} (\vp^{k} - \vq^{k})|_v$
for all $k$ with $|k| \ge 2$.
Consequently $|h_i(\vp) - h_i(\vq)|_v = |p_i - q_i|_v \ne 0$,
so  $|H(\vp) - H(\vq)|_v = |\vp-\vq|_v$.  In particular, $H$ is $1 - 1$
on $D(\vORIG,1)^-$.  

Next we will show that $H(\vX)$ has a right inverse
$G(\vX) = (g_1(\vX), \ldots, g_d(\vX))$ belonging to $\hcO_v[[\vX\]]$, 
such that $G(\vX) \equiv \vX$ modulo terms of degree $\ge 2$.
To do this, we apply Newton's method to power series,
and use convergence of values on $D(\vORIG,1)^-$ to deduce convergence of coefficients of
the power series.   

Let $J_H(\vX) = (\frac{\partial h_i}{\partial z_j}(\vX))$
be the Jacobian matrix of $H(\vX)$, computed using formal partial
derivatives of series.  By hypothesis, $J_H(\vX) \equiv I$
modulo terms of degree $\ge 1$, so \ $\det(J_H(\vX)) \in \hcO_v[[\vX]]$
is a power series with constant term $1$.
Thus, the formal geometric series for $\det(J_H(\vX))^{-1}$ belongs to
$\hcO_v[[\vX]]$ and converges for all $\vz \in D(\vORIG,1)^-$, 
and $J_H(\vX)^{-1} =  \det(J_H(\vX))^{-1} \cdot \Adj(J_H(\vX))$
has components given by power series in $\hcO_v[[\vX]]$
which also converge for all $\vz \in D(\vORIG,1)^-$.

Define a sequence of maps
$G^{(\ell)}(\vX) = (g_1^{(\ell)}(\vX), \ldots, g_d^{(\ell)}(\vX))$ 
with coordinate functions belonging to $\hcO_v[[\vX]]$, 
by setting $G^{(0)}(\vX) = \vX$ and putting
\begin{equation}  \label{BCF3}
G^{(\ell+1)}(\vX) \ = \
    G^{(\ell)}(\vX)
       - J_H(G^{(\ell)}(\vX))^{-1} \cdot (H(G^{(\ell)}(\vX)) - \vX)
\end{equation}
for each $\ell$.  Inductively one sees that $G^{(\ell)}(\vX) \equiv \vX$
modulo terms of degree $\ge 2$, so the substitutions make sense formally and
all the component functions $g_i^{(\ell)}(\vX)$ belong to $\hcO_v[[\vX]]$.
In particular they converge on $D(\vORIG,1)^-$.   

Next fix $\vq \in D(\vORIG,1)^-$ and put $r = |\vq|_v < 1$.
The usual sequence of Newton iterates $\vp_0, \vp_1, \ldots \in D(\vORIG,r)$
converging to a solution of $H(\vp) = \vq$ is defined by setting $\vp_0 = \vq$
and putting
\begin{equation} \label{BCF4}
\vp_{\ell+1} \ = \
      \vp_{\ell} - J_H(\vp_{\ell})^{-1} \cdot (H(\vp_{\ell}) - \vq) \ .
\end{equation}
By the form of $H(\vX)$, clearly $|H(\vp_0) - \vq|_v \le r^2$.  Assume
inductively that $\vp_{\ell} \in D(\vORIG,r)$ and
$|H(\vp_{\ell}) - \vq |_v \le r^{\ell+2}$.  By (\ref{BCF4}) we have
$\vp_{\ell+1} \in D(\vORIG,r)$.  Expanding $H(\vp_{\ell+1})$ 
and using (\ref{BCF4}) we find that
\begin{equation*}
|H(\vp_{\ell+1}) - \vq|_v \ \le \ r^{2(\ell+2)} \ \le \ r^{(\ell+1) + 2} \ .
\end{equation*}
From (\ref{BCF4}) we see that the $\vp_{\ell}$ converge to a vector 
$\vp \in D(\vORIG,r)$ such that $H(\vp) = \vq$.

Comparing (\ref{BCF3}) and (\ref{BCF4}) shows that
that $\vp_{\ell} = G^{(\ell)}(\vq)$  for each $\ell$,
that is, the values of the $G^{(\ell)}(\vq)$ converge 
for each $\vq \in D(\vORIG,1)^-$.  Moreover by the ultrametric inequality,
for each $r < 1$ the convergence is uniform on $D(\vORIG,r)$, with
$|G^{(\ell+1)}(\vq) - G^{(\ell)}(\vq)|_v \le r^{\ell+2}$ for all $\ell$
and all $\vq \in D(\vORIG,r)$.  Hence,
Lemma \ref{BCL1} produces a function $G(\vX) = (g_1(\vX), \ldots, g_d(\vX))$
with components $g_i(\vX) \in \hcO_v[[\vX]]$, such that $H(G(\vq)) = \vq$
for all $\vq \in D(\vORIG,1)^-$.  The fact that each $G^{(\ell)}(\vX) \equiv \vX$
modulo terms of degree $\ge 2$ means that $G(\vX)$ has this property as well.
From this, it follows that $G(D(\vORIG,1)^-) \subset D(\vORIG,1)^-$.
Hence, $H$ gives a surjection from $D(\vORIG,1)^-$ onto $D(\vORIG,1)^-$.
To see that $G(\vX)$ is a left inverse to $H(\vX)$ as well as a right
inverse, note that if $\vp \in D(\vORIG,1)^-$ and $\vq = H(\vp)$
then
\begin{equation*}
H(\vp) \ = \ \vq \ = H(G(\vq)) \ ;
\end{equation*}
however, since $H(\vz)$ is $1-1$ on $D(\vORIG,1)^-$, necessarily $\vp = G(\vq)$.
That is, $\vp = G(H(\vp))$  for all $\vp \in D(\vORIG,1)^-$.  

The fact that $H$ and $G$ define inverse functions on $D(\vORIG,1)^-$ means that 
$G(H(\vX))-\vX$ and $H(G(\vX))-\vX$
are identically equal to $\vORIG$ on $D(\vORIG,1)^-$, 
and then Lemma \ref{BVLem} 
shows their component power series are identically $0$.  
Hence, $H$ and $G$ are formal inverses.  

If $D(\vORIG,\vr) \subset D(\vORIG,1)^-$ is a polydisc satisfying $\max_i(r_i)^2 \le \min_i(r_i)$, 
then the fact that $H(\vX) \in \hcO_v[[\vX]]^d$ and $H(\vX) \equiv \vX$ modulo terms of degree $\ge 2$ 
shows that $H$ maps $D(\vORIG,\vr)$ into itself. However, $G(\vX)$ has the same form,
so $G$ maps $D(\vORIG,\vr)$ into itself as well.  Since $H$ and $G$ are inverses, 
and $H$ preserves distances, $H$ induces an isometry from $D(\vORIG,\vr)$ onto itself.  

Finally, if $H(\vX)$ is rational over $F_u$ for some subfield
$F_u \subset \CC_v$, the formulas (\ref{BCF3}) show that the $G^{(\ell)}(\vX)$
are rational over $F_u$, and if $F_u$ is complete then $G(\vX)$ is also rational over $F_u$.
It follows that $H(D(\vORIG,\vr) \cap F_u^d) \subseteq D(\vORIG,\vr) \cap F_u^d$
and $G(D(\vORIG,\vr) \cap F_u^d) \subseteq D(\vORIG,\vr) \cap F_u^d$,
so $H(D(\vORIG,\vr) \cap F_u^d) = D(\vORIG,\vr) \cap F_u^d$
\end{proof} 

\noindent{\bf Remark.}  Under the hypotheses of Proposition \ref{BCProp1}, 
in general $H(\vX)$ will not converge on the closed unit polydisc $D(\vORIG,1)$, 
and even if it does, it will not in general be $1-1$ and 
$G(\vX)$ will not converge on $D(\vORIG,1)$.  
For example, in one variable, if $H(X) = X - X^2$ then $H(0) = H(1) = 0$,
while the power series expansion for its inverse 
$G(X) = (1-\sqrt{1+4X}))/2$ only converges on $D(0,1)^-$.

\medskip

For the proof of Proposition \ref{BPropD2New}, we will need the following lemma.

\begin{lemma} \label{DiagonalZeroLemma} Let $R \in |\CC_v^{\times}|_v$, and suppose 
$h : D(0,R) \times D(0,R) \rightarrow D(0,r)$
is a map defined by a power series $h(X,Y) = \sum_{j,k=0}^\infty c_{jk}X^j Y^k \in \CC_v[[X,Y]]$ 
which converges on $D(0,R) \times D(0,R)$, and satisfies $h(x,x) = 0$ for all $x \in D(0,R)$.  
Then for all $x, y \in D(0,R)$, 
\begin{equation*} 
|h(x,y)|_v \ \le \ \frac{r}{R} \cdot |x-y|_v \ .
\end{equation*} 
\end{lemma}  

\begin{proof}  By the Maximum Modulus Principle for power series,
\index{Maximum principle!nonarchimedean!for power series} 
we have $|c_{jk}|_v \le r/R^{j+k}$ for all $j, k$, and 
\begin{equation} \label{FcBound}
\lim_{j, k \rightarrow \infty} |c_{jk}|_v \cdot R^{j+k} \ = \ 0 
\end{equation}
since $h(X,Y)$ converges on $D(0,R) \times D(0,R)$.  Fix $x, y \in D(0,R)$.  Then 
\begin{eqnarray}
h(x,y) & = & h(x,y) - h(x,x) \ = \ \sum_{j,k=0}^{\infty} c_{jk} x^j (y^k-x^k) \notag \\ 
& = & (y-x) \cdot \sum_{j,k=0}^{\infty} c_{jk} x^j \big(\sum_{\ell = 0}^{k-1} x^{\ell} y^{k-1-\ell}\big) \ .
\label{hxySum}
\end{eqnarray}
By our estimate for $|c_{ij}|_v$ and the ultrametric inequality, 
each term in the sum on the right in (\ref{hxySum}) 
has absolute value at most $r/R$, and by (\ref{FcBound}) the sum converges and has absolute value 
at most $r/R$.  Thus $|h(x,y)|_v \le |x-y|_v \cdot r/R$.  
\end{proof} 


\section{  Proof of the Local Action Theorem } \label{LocalActionSection}
\index{Jacobian variety} 

In this section we prove Theorem \ref{BKeyThm1} and Proposition \ref{BPropD2New}.  
We use the notation established prior to the statement of the Theorem, 
and begin with four lemmas.  

\smallskip
In our first lemma, we show that certain subsets of the formal group of $\Jac(\cC_v)$ are subgroups.
\index{Jacobian variety}
\index{formal group}
We denote a group structure by a triple consisting of the underlying set 
and two vectors of power series, representing addition and negation, which converge on the set.  

\begin{lemma} \label{SMGroupLemma}
Let $\cD(\vORIG,1)^- := (D(\vORIG,1)^-,S(\vX,\vY),M(\vX))$ be the formal group of $\Jac(\cC_v)$.
\index{formal group} 
Let $L : \CC_v^g \rightarrow \CC_v^g$ be a nonsingular linear map,    
and fix $0 < \eta < 1$.  Then there is an $R_1 > 0$ such that for each $\vr = (r_1, \ldots, r_g)$ 
satisfying $0 < r_1, \ldots, r_g \le R_1$ and $\eta \cdot \max_i(r_i) \le \min_i(r_i)$, 
and each $0 < \lambda \le 1$,

$(A)$ $\cD_L(\vORIG,\lambda\vr) := (L(D(\vORIG,\lambda\vr)),S(\vX,\vY),M(\vX))$ 
is a subgroup of $\cD(\vORIG,1)^-;$  

$(B)$ If $D(\vORIG,\vr)$ is given its structure as an additive subgroup of $\CC_v^g$, 
the map of sets $L : D(\vORIG,\vr) \rightarrow L(D(\vORIG,\vr))$ induces an isomorphism of groups   
\begin{equation*} 
D(\vORIG,\vr)/D(\vORIG,\eta \vr) \ \cong \ \cD_L(\vORIG,\vr)/\cD_L(\vORIG,\eta \vr) \ .
\end{equation*} 
\end{lemma}  

\begin{proof}  Let $(\alpha_{ij}), (\beta_{ij}) \in M_g(\CC_v)$ be the 
matrices associated to $L^{-1}$ and $L$, respectively, 
and put $A = \max_{i,j} |\alpha_{ij}|_v$, $B = \max_{i,j} |\beta_{ij}|_v$.  Let    
\begin{equation*} 
\tS(\vX,\vY) \ = \ L^{-1}(S(L(\vX),L(\vY))) \ , \qquad  \tM(\vX) = L^{-1}(M(L(\vX))) \ .
\end{equation*} 
Then $\cD_L(\vORIG,\vr)$ is a group if and only if 
$\widetilde{\cD}(\vORIG,\vr) := (D(\vORIG,\vr),\tS(\vX,\vY),\tM(\vX))$ is a group,
and if they are groups, the map $L : D(\vORIG,\vr) \rightarrow L(D(\vORIG,\vr))$ 
induces isomorphism between them.  
Note that if $\vr$ satisfies the conditions in the Lemma, then so does $\lambda \vr$
for each $0 < \lambda \le 1$.   Hence it suffices to consider the $\widetilde{\cD}(\vORIG,\vr)$.

Write $\tS(\vX,\vY) = (\tS_1(\vX,\vY), \ldots, \tS_g(\vX,\vY))$,  
$\tM(\vX) = (\tM_1(\vX), \ldots, \tM_g(\vX))$, and for each $i$, expand 
\begin{equation*} 
\tS_i(\vX,\vY) = X_i + Y_i + \sum_{|k|+|\ell| \ge 2} c_{i,k,\ell} \vX^k \vY^{\ell} \ , \qquad 
\tM_i(\vX) = -X_i + \sum_{|k| \ge 2} d_{i,k} \vX^k \ .
\end{equation*} 
Since the power series defining $S(\vX,\vY)$ and $M(\vX)$ have coefficients in $\hcO_v$, 
it is easy to see that for each $i = 1, \ldots, g$ and all $k, \ell \in \NN^g$ with $|k|+|\ell| \ge 2$
(resp. all $k \in \NN^g$ with $|k| \ge 2$), we have 
$|c_{i,k,\ell}|_v \ \le \ AB^{|k|+|\ell|}$ and $|d_{i,k}|_v \le AB^{|k|}$.     
Hence for all $\vx, \vy \in D(\vORIG,\vr)$ 
\begin{equation*}
|c_{i,k,\ell} \vx^k \vy^k|_v \ \le \ \ A \big(B \max_i(r_i)\big)^{|k|+|\ell|} \ , \qquad  
|d_{i,k} \vx^k|_v \ \le \ \ A \big(B \max_i(r_i)\big)^{|k|} \ .
\end{equation*}  

If $\max_i(r_i) < 1/B$ then $\tS(\vX,\vY)$ and 
$\tM(\vX)$ converge for $\vx, \vy \in D(\vORIG,\vr)$.   
If also $\max_i(r_i) \le \eta/(AB^2)$ and $\eta \max_i(r_i) \le \min_i(r_i)$, 
then for all $|k|+|\ell| \ge 2$, 
\begin{equation*}  
|c_{i,k,\ell} \vx^k \vy^\ell|_v \, \le \,  A \big(B \max_i(r_i)\big)^{|k|+|\ell|} \, \le \, A B^2 \max_i(r_i)^2 
\, \le \, \eta \max_i(r_i) \, \le \, \min_i(r_i) \ .
\end{equation*} 
Similarly for all $|k| \ge 2$, we have $|d_{i,k} \vx^k|_v \le \min_i(r_i)$, 
and so $\tS(\vX,\vY)$ and $\tM(\vX)$ map $D(\vORIG,\vr)$ into itself. 
Since $\tS(\tS(\vx,\tM(\vy)),\vy)=\vx$ and $\tM(\tM(\vx)) = \vx$ for all $\vx, \vy \in D(\vORIG,\vr)$, 
they are surjective.      
Thus $\widetilde{\cD}(\vORIG,\vr)$ is a group.  

Finally, if $\max(r_i) \le \eta^2 / A B^2$ then by an argument similar to the one above, 
for all $\vx,\vy \in D(\vORIG,\vr)$ and all $k, \ell$ with $|k|+|\ell| \ge 2$, 
one has $|c_{i,k,\ell} \vx^k \vy^{\ell}|_v \le \eta \min_i(r_i)$. 
Likewise, for all $\vx \in D(\vORIG,\vr)$ and all $k$ with $|k| \ge 2$, one has 
$|d_{i,k} \vx^k| \le \eta \min_i(r_i)$.  Thus  
\begin{equation*} 
\tS(\vx,\vy) = \vx + \vy + \delta_1(\vx,\vy) \ , \qquad 
\tM(\vx) = -\vx + \delta_2(\vx) \ ,
\end{equation*} 
where $\delta_1(\vx,\vy)$ and $\delta_2(\vx)$ belong to $D(\vORIG,\eta \vr)$.  This means that 
if $D(\vORIG,\vr)$ is viewed as an additive subgroup of $\CC_v^g$, then  
$\widetilde{\cD}(\vORIG,\vr)/\widetilde{\cD}(\vORIG,\eta \vr) \ \cong \ 
D(\vORIG,\vr)/D(\vORIG,\eta \vr)$. 
 
Put $R_1 = \frac{1}{2} \min(1/B,\eta^2/(AB^2))$.
Then if $0 < r_1, \ldots, r_g \le R_1$ and $\eta \max(r_i) \le \min(r_i)$, 
both $\widetilde{\cD}(\vORIG,\vr)$ and  $\widetilde{\cD}(\vORIG,\eta \vr)$ are groups,
and 
\begin{equation*}
D(\vORIG,\vr)/D(\vORIG,\eta \vr) \ \cong \ 
\widetilde{\cD}(\vORIG,\vr)/\widetilde{\cD}(\vORIG,\eta \vr) \ \cong \ 
\cD_L(\vORIG,\vr)/\cD_L(\vORIG,\eta \vr) \ .
\end{equation*} 
This yields the result.  
\end{proof} 

Our next lemma shows that if $H(\vX) : \CC_v^g \rightarrow \CC_v^g$ 
is a map defined by convergent power series, 
whose derivative $L = H^{\prime}(\vORIG)$
is nonsingular, then for suitable polydiscs $D(\vORIG,\vr)$ the image $H(D(\vORIG,\vr))$ 
coincides with $L(D(\vORIG,\vr))$.  

\begin{lemma} \label{ConvergenceLemma}  
Let $H(\vX) \in \CC_v[[\vX]]^g$ 
be a vector of power series which converges on $D(\vORIG,r)$  
for some $r > 0$, and maps $D(\vORIG,r)$ into $D(0,1)^-$.  
Assume the derivative $L = H^{\prime}(\vORIG) : \CC_v^g \rightarrow \CC_v^g$ 
is nonsingular.  

Then for each $0 < \eta < 1$, there is an $R_2$ with $0 < R_2 \le r$ 
such that for each $D(\vORIG,\vr)$ 
with $0 < r_1, \ldots, r_g \le R_2$ and $\eta \cdot \max_i(r_i) \le \min_i(r_i)$, 

$(A)$ $H(\vX)$ gives an analytic isomorphism from $D(\vORIG,\vr)$ onto $L(D(\vORIG,\vr))$, 
with the property that for each $\vz \in D(\vORIG,\vr)$,  
\begin{equation*}
H(\vz) \ \equiv \ L(\vz) \pmod{L(D(\vORIG,\eta \vr))} \ . 
\end{equation*}    

$(B)$  For each complete field $F_u$ such that $K_v \subseteq F_u \subseteq \CC_v$ 
and $H(\vX)$ is rational over $F_u$,  
\begin{equation*}
H(D(\vORIG,\vr) \cap F_u^g) \ = \ L(D(\vORIG,\vr)) \cap F_u^g \ .
\end{equation*}  
\end{lemma}

\begin{proof}  
Let $(\alpha_{ij}) \in M_g(\CC_v)$ be the matrix associated to $L^{-1}$,
and put $A = \max_{i,j} (|\alpha_{ij}|_v)$.   
Choose $\pi \in K_v^{\times}$ so that $|\pi|_v < \min(r,r^2/A)$, 
and put $B = |\pi|_v$.  

Our plan is to apply Proposition \ref{BCProp1} to 
\begin{equation} \label{HLFormula} 
\widetilde{H}(\vX) \ := \ L^{-1}(\frac{1}{\pi} H(\pi \vX)) \ .
\end{equation}
By construction, $\widetilde{H}^{\prime}(\vORIG) = \id$.  Write 
\begin{equation*} 
H(\vX) \ = \ L(\vX) + \sum_{|k| \ge 2} \vc_k \vX^k \ , \qquad \   
\widetilde{H}(\vX) \ = \ \vX + \sum_{|k| \ge 2} \tc_k \vX^k \ ,
\end{equation*} 
where $\vc_k = (c_{1,k}, \ldots, c_{g,k})$ and $\tc_k = (\tc_{1,k}, \ldots, \tc_{g,k})$.  
Since $H(\vX)$ converges on $D(\vORIG,r)$ and maps $D(\vORIG,r)$ into $D(\vORIG,1)^-$,
we have $|c_{i,k}|_v \le 1/r^{|k|}$ for all $i, k$.  
Since $B/r \le 1$ and $B \le r^2/A$,  for each $i$ and each $k$ with $|k| \ge 2$ we have  
\begin{equation*} 
|\tc_{i,k}|_v \ \le \ \frac{A}{B} \big(\frac{B}{r}\big)^{|k|} 
\ < \ \frac{A}{B} \big(\frac{B}{r}\big)^2 \ \le \ 1 \ .
\end{equation*} 

Thus $\widetilde{H}(\vX) \in \hcO_v[[\vX]]^g$ and $\widetilde{H}(\vX) \equiv \vX$ 
modulo terms of degree $\ge 2$.
By Proposition \ref{BCProp1}, if $0 < s_1, \ldots, s_d < 1$ and $\max(s_i)^2 \le \min(s_i)$, 
then $\widetilde{H}(D(\vORIG,\vs)) = D(\vORIG,\vs)$.  
Furthermore, if $F_u$ is a complete field with $K_v \subset F_u \subset \CC_v$,
and if $H(\vX)$ is rational over $F_u$, then $\widetilde{H}(\vX)$ is rational over $F_u$,    
so $\widetilde{H}(D(\vORIG,\vs) \cap F_u^g) = D(\vORIG,\vs) \cap F_u^g$. 

Given $0 < \eta < 1$, 
put $R_2 = \eta ^2 B$.  Then $R_2 < r$.  

Suppose $0 < r_1, \ldots, r_g \le R_2$, with $\eta \cdot \max_i(r_i) \le \min_i(r_i)$.  
We first show that $H(D(\vORIG,\vr)) = L(D(\vORIG,\vr))$. 
Put $\vs = (1/B) \vr$.  Then $0 < s_i < 1$ for each $i$, 
since $\max_i(r_i) \le R_2$ and $R_2/B = \eta^2 < 1$, and  
\begin{equation*} 
\max_i(s_i)^2 \ \le \ \frac{R_2}{B^2} \max_i(r_i) 
  \ \le \ \frac{\eta}{B} \big( \eta \max_i(r_i) \big) \ < \ \frac{1}{B} \min_i(r_i) 
              \ = \ \min_i(s_i) \ .
\end{equation*} 
Since $H(\pi \vX) = \pi L(\widetilde{H}(\vX))$ 
and $D(\vORIG,\vr) = D(\vORIG,B\vs) = \pi D(\vORIG,\vs)$, it follows that 
$H(D(\vORIG,\vr)) = L(D(\vORIG,\vr))$ and $H(D(\vORIG,\vr) \cap F_u^g) = L(D(\vORIG,\vr)) \cap F_u^g$. 

Third, we show that  $H(\vx) \equiv L(\vx) \pmod{L(D(\vORIG,\eta \vr))}$ 
for each $\vx \in D(\vORIG,\vr)$.  Note that $H(\vX) = L( \pi \widetilde{H}(\frac{1}{\pi} \vX))$,
and that if $\vx \in D(\vORIG,\vr)$ then $\frac{\vx}{\pi} \in D(\vORIG,1)^-$.
By the inequalities $\max_i(r_i) \le R_2$ and $ R_2 \le \eta^2 B < B$, 
and our estimate  $|\tc_{i,k}|_v \le 1$ above,
for each $i$ and each $k$ with $|k| \ge 2$ we have  
\begin{eqnarray*}
| \pi \tc_{i,k} \big(\frac{\vx}{\pi}\big)^k|_v & \le & B \big(\frac{\max_i(r_i)}{B} \big)^{|k|}
\ < \ B \big(\frac{\max_i(r_i)}{B} \big)^2 \ \le \ \frac{R_2}{B} \max_i(r_i) \\
& \le & \eta^2 \max_i(r_i) \ \le \ \eta \min_i(r_i) \ .
\end{eqnarray*} 
Furthermore, $\widetilde{H}(\vX) \equiv \vX$ modulo terms of degree $\ge 2$, so  
$\pi \widetilde{H}(\frac{\vx}{\pi}) \equiv \vx \pmod{D(\vORIG,\eta \vr)}$.  
Since $H(\vx) = L(\pi \widetilde{H}(\frac{\vx}{\pi}))$,  
it follows that $H(\vx) \equiv L(\vx) \pmod{L(D(\vORIG,\eta \vr))}$. 
\end{proof} 

Our third lemma shows that for a dense set of points $\va \in \cC_v(\CC_v)^g$, 
the map $[\ ] : \cCbar_v^g \rightarrow \BPic_{\cCbar_v/\CC_v}^{g}$ is nonsingular at $\va$.

\begin{lemma} \label{DensityLemma} 
There is a non-empty Zariski-open subset $U \subset \cCbar_v^g$ such that for 
each $\va \in U(\CC_v)$, the map $[\ ] : \cCbar_v^g \rightarrow \BPic_{\cCbar_v/\CC_v}^{g}$ 
is nonsingular at $\va$.  Moreover, for each complete field $F_u$ with $K_v \subseteq F_u \subseteq \CC_v$
such that $\cC_v(F_u)$ is nonempty, $U(F_u)$ is dense in $\CC_v(F_u)^g$ for the $v$-topology. 
\end{lemma} 

\begin{proof}
Fix $\va \in \cC_v^g(\CC_v)$.  
The map $[\ ] : \cCbar_v^g \rightarrow \BPic_{\cCbar_v/\CC_v}^{g}$, 
where $[\vx]$ is the linear equivalence class of the divisor $(x_1) + \cdots + (x_g)$, factors as 
\begin{equation*}
[\ ] : \cCbar_v^g \ \stackrel{Q}{\longrightarrow} \ \Sym^g(\cCbar_v) \ 
\stackrel{P}{\longrightarrow} \ \BPic^g_{\cCbar_v/\CC_v} 
\end{equation*} 
where $Q : \cCbar_v^g \rightarrow \Sym^g(\cCbar_v)$ is the quotient, and  
$P : \Sym^g(\cCbar_v) \rightarrow \ \BPic^g_{\cCbar_v/\CC_v}$ is the Abel map.\index{Abel map}  

First consider $Q : \cCbar_v^g \rightarrow \Sym^{(g)}(\cCbar_v)$.
Let $S_g$ be the symmetric group on the letters $\{1, \ldots, g\}$.  For each $\pi \in S_g$, 
write $\pi(\vx) = (x_{\pi(1)}, \ldots, x_{\pi(g)})$. 
If $\vx \in \cC_v^g(\CC_v)$ has distinct coordinates, then $Q(\vx)$ has $g!$ distinct preimages, 
namely the points $\pi(\vx)$ for $\pi \in S_g$.  Since $Q$ is a finite, flat morphism of degree $g!$, 
it is \'etale at $\vx$.  In particular, it is nonsingular at $\vx$.
Thus, $Q$ is nonsingular in the complement of the generalized diagonal 
$X = \{\vx = (x_1, \ldots, x_g) : x_i = x_j \ \text{for some $i \ne j$}\}$.  

The morphism $P : \Sym^{(g)}(\cCbar_v) \rightarrow \ \BPic^g_{\cCbar_v/\CC_v}$ is
a birational isomorphism, and so it is nonsingular outside a proper Zariski-closed 
set $Y \subset \Sym^{(g)}(\cCbar_v)$.  Since $Q$ is dominant, $Q^{-1}(Y) \subset \cCbar_v^g$  
is a proper, Zariski-closed subset of $\cCbar_v^g$. 
Thus $P \circ Q$ is nonsingular on the nonempty, Zariski-dense set 
$U = \cCbar_v^g \backslash (X \cup Q^{-1}(Y))$.  

We now show that if $F_u$ is a complete field with $K_v \subseteq F_u \subseteq \CC_v$
for which $\cC_v(F_u)$ is nonempty, then $U(F_u)$ is dense in $\cC_v^g(F_u)$ for the $v$-topology.
Suppose to the contrary that there was a point $\va \in \cC_v^g(F_u)$ such that 
 $\BCg(\va,\vr) \cap U(F_u)$ is empty, for some $r_1, \ldots, r_g > 0$. Then $\va \in X \cup Q^{-1}(Y)$.
Let $h$ be one of the equations cutting out $X \cup Q^{-1}(Y)$ in a neighborhood of $\va$.   
After shrinking the $r_i$,
if necessary, we can assume that each ball $B(a_i,r_i)$ is isometrically parametrizable.  
Let $\Phi_{\va} : D(\vORIG,\vr) \rightarrow \BCg(\va,\vr)$ 
be a product of $F_u$-rational isometric parametrization maps as in (\ref{FFdef}). 
By composing $h$ with $ \Phi_{\va}$, 
we would obtain a power series $h \circ \Phi_{\va}(\vz) \in \CC_v[[\vX]]$ 
converging on $D(\vORIG,\vr)$ and vanishing identically on 
$D(\vORIG,\vr) \cap F_u^g$.  
By Lemma \ref{BVLem}, $h \circ \Phi_{\va} \equiv 0$, so 
$h$ would vanish on all of $\BCg(\va,\vr) = \Phi_{\va}(D(\vORIG,\vr))$.  
This is impossible by what has been shown above.  
\end{proof} 

\smallskip
Our final lemma concerns the trace $\Tr_{F_u/K_v} : \Jac(\cC_v)(F_u) \rightarrow \Jac(\cC_v)(K_v)$ 
\index{Jacobian variety}
when $F_u/K_v$ is finite and separable. 
Let $H_w$ be the galois closure of $K_v$ in $\CC_v$.
The subgroup $B_{J}(O,1)^- \cap \Jac(\cC_v)(H_w)$ is stable under $\Gal(H_w/K_v)$, so we can pull back 
the trace $\Tr_{F_u/K_v}$ to the formal group $\cD(\vORIG,1)^-$.  
\index{formal group}
Let $\vX_1 = (X_{11}, \ldots, X_{1g}), \, \cdots, $
$\vX_d = (X_{d1}, \ldots,X_{dg})$ be $d$ vectors of independent variables, and put  
\begin{equation*} 
S^{(d)}(\vX_1, \cdots, \vX_d) \ = \ S(\vX_1, S(\vX_2, \cdots, S(\vX_{d-1},\vX_d ))) 
\ \in \ O_v[[\vX_1, \ldots, \vX_d]] \ .
\end{equation*} 
Since $S(\vX,\vY)$ is commutative and associative, 
for each permutation $\pi$  on the letters $\{1, \ldots, d\}$  
we have $S^{(d)}(\vX_{\pi(1)}, \cdots, \vX_{\pi(d)}) = S^{(d)}(\vX_1, \cdots, \vX_d)$. 
Let $\sigma_1, \ldots, \sigma_d$ be the distinct embeddings of $F_u$ into $H_w$, and 
extend each $\sigma_i$ to an automorphism of $H_w$.   
Define $\Tr_{F_u/K_v,S} : D(\vORIG,1)^- \cap F_u^g \rightarrow D(\vORIG,1)^- \cap K_v^g$ by 
\begin{equation*}
\Tr_{F_u/K_v,S}(\vx) \ = \ S^{(d)}(\sigma_1(\vx), \cdots, \sigma_d(\vx)) \ .
\end{equation*} 
Then $\Tr_{F_u/K_v}(\Psi(\vx)) = \Psi(\Tr_{F_u/K_v,S}(\vx))$ for each $\vx \in D(\vORIG,1)^- \cap F_u^g$.

\begin{lemma} \label{TraceBall}
Let $F_u/K_v$ be a separable finite extension.  Then there are constants $R_u$ and $C_u$,
depending on $F_u$, with $0 < R_u < 1$ and $0 < C_u \le 1$, such that for each $r$ with $0 < r \le R_u$
we have 
\begin{equation} \label{FTraceF} 
\Tr_{F_u/K_v}(B_J(O,r) \cap \Jac(\cC_v)(F_u)) \ \supseteq \ B_J(O,C_u r) \cap \Jac(\cC_v)(K_v) \ .  
\end{equation}
\index{Jacobian variety}
   
\end{lemma} 

\begin{proof}  We keep the notation above.  
Since $\Psi : D(\vORIG,1)^- \rightarrow B_J(O,1)^-$ 
is a $K_v$-rational isometric parametrization, for each $0 < r < 1$ 
it takes $D(\vORIG,r) \cap F_u^g$ to $B_J(0,r) \cap \Jac(\cC_v)(F_u)$ 
and $D(\vORIG,r) \cap K_v^g$ to $B_J(0,r) \cap \Jac(\cC_v)(K_v)$.
\index{Jacobian variety}
Thus it suffices to prove the assertion corresponding to (\ref{FTraceF}) 
for $\Tr_{F_u/K_v,S}$.   

The idea is that one can describe the image of $\Tr_{F_u/K_v,S}$ 
by restricting $\Tr_{F_u/K_v,S}$ to a carefully chosen $K_v$-rational subspace.  
Since $F_u/K_v$ is separable, there is a $0 \ne \beta \in \cO_v$ such that $\Tr_{F_u/K_v}(\cO_u) = \beta \cO_v$.
Let $\alpha \in \cO_u$ be such that $\Tr_{F_u/K_v}(\alpha) = \beta$;  
then $|\alpha|_v = 1$, $|\beta|_v \le 1$.  
Define
\begin{equation*}   
T(\vX) \ = \ S^{(d)}(\sigma_1(\alpha) \vX, \cdots, \sigma_d(\alpha) \vX) \ .
\end{equation*} 
Then for each $\vx \in D(\vORIG,1)^- \cap K_v^g$, we have $\alpha \vx \in D(\vORIG,1)^- \cap F_u^g$, and 
\begin{equation} \label{FTrMotive}
\Tr_{F_u/K_v,S}(\alpha \vx) \ = \ T(\vx) \ .
\end{equation}   

\`A priori $T(\vX) \in \cO_w[[\vX]]$.  However, for each $\sigma \in \Gal(H_w/K_v)$ there is a 
permutation $\pi = \pi_{\sigma}$ of $\{1, \ldots, d\}$ such that $\sigma  \sigma_i = \sigma_{\pi(i)}$, 
so (writing $\sigma T$ for the power series obtained by letting $\sigma$ act on the coefficients of $T$)
\begin{eqnarray*} 
(\sigma T)(\vX) & = & S^{(d)}(\sigma (\sigma_1(\alpha)) \vX, \cdots, \sigma (\sigma_d(\alpha)) \vX) \\
& = & S^{(d)}(\sigma_{\pi(1)}(\alpha) \vX, \cdots, \sigma_{\pi(d)}(\alpha) \vX)
\ = \ T(\vX) 
\end{eqnarray*} 
Thus $T(\vX) \in \cO_v[[\vX]]$.  Since $S(\vX,\vY) \equiv \vX + \vY$ modulo terms of degree $\ge 2$,
it follows that $T(\vX) \equiv \Tr_{F_u/K_v}(\alpha) \vX = \beta \vX$ modulo terms of degree $\ge 2$.  

Put $\tT(\vX) = \frac{1}{\beta^2} T(\beta \vX)$.  
Then $\tT(\vX) \in \cO_v[[\vX]]$ and $\tT(\vX) \equiv \vX$ modulo terms of degree $\ge 2$.  
By Proposition \ref{BCProp1}, for each $0 < s < 1$ we have 
$\tT(D(\vORIG,s) \cap K_v^g) = D(\vORIG,s) \cap K_v^g$.
Suppose $0 < r < |\beta|_v$, and put $s = r/|\beta|_v$.  Then $0 < s < 1$. 
Since $T(\beta \vX) = \beta^2 \tT(\vX)$ and $\beta \in K_v^{\times}$, it follows that 
\begin{eqnarray} 
T(D(\vORIG,r) \cap K_v^g)  
& = & T(\beta (D(\vORIG,s) \cap K_v^g)) \ = \ \beta^2 \tT(D(\vORIG,s) \cap K_v^g) \notag \\
& = &  \beta^2(D(\vORIG,s) \cap K_v^g) \ = \ D(\vORIG,|\beta|_v r) \cap K_v^g \ . \label{FTrC} 
\end{eqnarray}  
Set $R_u = \frac{1}{2}|\beta|_v$ and $C_u = |\beta|_v$.  
By (\ref{FTrMotive}) and (\ref{FTrC}), if $0 < r \le R_u$, then  
\begin{equation*} 
\Tr_{F_u/K_v,S}(D(\vORIG,r) \cap F_u^g) \ \supseteq \ 
T(D(\vORIG,r) \cap K_v^g) \ = \ D(\vORIG,C_u r) \cap K_v^g \ . 
\end{equation*} 
This yields (\ref{FTraceF}).
\end{proof} 

\begin{proof}[Proof of Theorem \ref{BKeyThm1}.]  
Fix a point $\va = (a_1, \ldots, a_g) \in \cC_v^g(\CC_v)$ 
for which $[\ ] : \cCbar_v^g \rightarrow \BPic_{\cCbar_v/\CC_v}^{g}$ is nonsingular at $\va$.  
By Lemma \ref{DensityLemma}, such $\va$ are dense in $\cC_v^g(\CC_v)$ for the $v$-topology, 
and if $F_u$ is a complete field with $K_v \subseteq F_u \subseteq \CC_v$ 
such that $\cC_v(F_u)$ is nonempty, they are dense in $\cC_v^g(F_u)$. 

The associated map $\crJ_{\va} : \cC_v(\CC_v)^g \rightarrow \Jac(\cC_v)$ given by 
$\crJ_{\va}(\vx) = [\vx] \Jm [\va]$ is nonsingular at $\va$, since the subtraction morphism  
$\Jm [\va]$ in $\BPic(\cCbar_v/\CC_v)$ is an isomorphism.  
For a suitably small $r > 0$, let $\Phi_{\va} : D(\vORIG,r) \rightarrow \BCg(\va,r)$ 
be the isometric parametrization from (\ref{FFdef}).  Without loss, we can assume 
$r$ is small enough that $\crJ_{\va}(\Phi_{\va}(D(\vORIG,r))) \subset B_J(O,1)^-$.  
Let $\Psi : D(\vORIG,1)^- \rightarrow B_J(O,1)^-$
be the isometric parametrization from (\ref{JPsiDef}), and put 
\begin{equation*}    
H_{\va}(\vX) \ = \ \Psi^{-1} \circ \crJ_{\va} \circ \Phi_{\va}(\vX) \ .
\end{equation*} 
Then $H_{\va}(\vX) \in \CC_v[[\vX]]^g$, and  $H_{\va}(\vX)$ converges on $D(\vORIG,r)$.  
If $\va \in \cC_v(F_u)^g$, 
with $F_u$ as above,  then $H_{\va}(\vX) \in F_u[[\vX]]^g$. 
Put $L_{\va} = H_{\va}^{\prime}(\vORIG)$, viewing it as a linear map 
$L_{\va} : \CC_v^g \rightarrow \CC_v^g$.  By our discussion, $L_{\va}$ is invertible, 
and if $H_{\va}(\vX) \in F_u[[\vX]]$ with $F_u$ as above, then $L_{\va} \in \GL_g(F_u)$.  

Fix $0 < \eta < 1$, and let $R_1, R_2 > 0$ be the numbers given by Lemmas \ref{SMGroupLemma} 
and \ref{ConvergenceLemma} applied to $H_{\va}(\vX)$, $L_{\va}$ and $\eta$.  
Put $R = \min(1/2,R_1,R_2)$. 
Suppose $0 < r_1, \ldots, r_g \le R$ and $\eta \max_i(r_i) \le \min(r_i)$.  
By Lemma \ref{ConvergenceLemma}, $H_{\va}$ induces an analytic isomorphism from $D(\vORIG,\vr))$ 
onto $L_{\va}(D(\vORIG,\vr))$.  In particular, $H_{\va}$ is injective on $D(\vORIG,\vr)$, 
which means (since $\Phi_{\va}$ and $\Psi$ are isomorphisms and $\crJ_{\va}(\vx) = [\vx] \Jm [\va]$) 
that the balls $B_{\cC_v}(a_1,r_1), \dots, B_{\cC_v}(a_g,r_g)$ must be pairwise disjoint, 
and especially that the points $a_1, \ldots, a_g$ must be distinct.  

By Lemma \ref{SMGroupLemma}, 
$L_{\va}(D(\vORIG,\vr))$ (equipped with the addition law $S(\vX,\vY)$ and negation $M(\vX)$)
is a subgroup of the formal group of $\Jac(\cC_v)$,
\index{Jacobian variety}\index{formal group}
denoted by $\cD_{L_{\va}}(\vORIG,\vr)$.    
By Lemma \ref{ConvergenceLemma}, for each $\vz \in D(\vORIG,\vr)$ 
\begin{equation} \label{FHCong1} 
H_{\va}(\vz) \ \equiv \ L_{\va}(\vz) \pmod{L_{\va}(D(\vORIG,\eta \vr))} \ ,  
\end{equation} 
and the map $L_{\va} : D(\vORIG,\vr) \rightarrow L_{\va}(D(\vORIG,\vr))$ induces an isomorphism 
from the additive group $D(\vORIG,\vr)/D(\vORIG,\eta \vr)$ 
onto $\cD_{L_\va}(\vORIG,\vr)/\cD_{L_\va}(\vORIG,\eta \vr)$.
  
\smallskip
We now define an action $\ap_0$ of $\cD_{L_\va}(\vORIG,\vr)$ 
on the polydisc $D(\vORIG,\vr)$ 
by setting, for $\vw \in L_{\va}(D(\vORIG,\vr))$ and $\vz \in D(\vORIG,\vr)$,  
\begin{equation*} 
\vw \ap_0 \vz \ = \ H_{\va}^{-1}(S(\vw,H_{\va}(\vz))) \ .
\end{equation*} 
It is easy to check that for all $\vw_1, \vw_2 \in L_{\va}(D(\vORIG,\vr))$ 
we have $\vw_1 \ap_0 (\vw_2 \ap_0 \vz) = S(\vw_1,\vw_2) \ap_0 \vz$.  
Thus, $\cD_{L_\va}(\vORIG,\vr)$ acts on $D(\vORIG,\vr)$.  
The action is transitive, since if $\vz_1, \vz_2 \in D(\vORIG,\vr)$ 
and $w = S(H_{\va}(\vz_2),M(H_{\va}(\vz_1)))$ then $w \ap_0 \vz_1 = \vz_2$.  
It is simple, since if $w_1 \ap \vz = w_2 \ap \vz$ then $S(w_1,H_{\va}(\vz)) = S(w_2,H_{\va}(\vz))$,
and hence $w_1 = w_2$.  Thus, $D(\vORIG,\vr)$ is a principal homogeneous space\index{principal homogeneous space} for 
$\cD_{L_\va}(\vORIG,\vr)$ under $\ap_0$.   

If in addition $\va \in \cC_v(F_u)^g$ where $F_u$ is as above, then by Lemma \ref{ConvergenceLemma}, 
$H_{\va}$ maps $D(\vORIG,\vr) \cap F_u^g$ isomorphically onto $L_{\va}(D(\vORIG,\vr)) \cap F_u^g$.  
Since $S(\vX,\vY)$ and $M(\vX)$ are rational over $K_v$, 
$\cD_{L_\va}(\vORIG,\vr) \cap F_u^g$ acts 
simply and transitively on $D(\vORIG,\vr) \cap F_u^g$.  


\smallskip
We can now prove parts (A) and (B) of Theorem \ref{BKeyThm1}.  Given $\va$ and $\vr$ as above, 
define $W_{\va}(\vr) \subset \Jac(\cC_v)(\CC_v)$ by $W_{\va}(\vr) = \crJ_{\va}(\BCg(\va,\vr))$. 
\index{Jacobian variety}  
Since $\Psi \circ H_{\va} = \crJ_{\va} \circ \Phi_{\va}$, we also have  
\begin{equation*} 
W_{\va}(\vr)  \ = \ \crJ_{\va}(\Phi_{\va}(D(\vORIG,\vr)))  \\
             \ = \ \Psi(H_{\va}(D(\vORIG,\vr))) \ = \ \Psi(\cD_{L_{\va}}(\vORIG,\vr))  \ .
\end{equation*} 
Since $\cD_{L_{\va}}(\vORIG,\vr))$ is a subgroup of the formal group $(D(\vORIG,1)^-,S(\vX,\vY),M(\vX))$,
\index{formal group}
it follows that $W_{\va}(\vr)$ (equipped with the group operations $\Jp$, $\Jm$) 
is a subgroup of $\Jac(\cC_v)(\CC_v)$.  Furthermore, by Lemma \ref{SMGroupLemma} 
the map $\Psi \circ L_{\va}$ induces an isomorphism
\begin{equation} \label{FDEquiv}
D(\vORIG,\vr)/D(\vORIG,\eta \vr) \ \cong \ W_{\va}(\vr)/W_{\va}(\eta \vr) \ .
\end{equation}  
For each $\vz \in D(\vORIG,\vr)$, we have 
$L_{\va}(\vz) \equiv H_{\va}(\vz)  \pmod{L_{\va}(D(\vORIG,\eta \vr))}$, 
and $H_{\va} = \Psi^{-1} \circ \crJ_{\va} \circ \Phi_{\va}$. 
Hence (\ref{FHCong1})  shows that 
\begin{equation} \label{FCong3}
\crJ_\va(F_\va(\vz)) \ \equiv \ \Psi(L_\va(\vz))  \pmod{W_{\va}(\eta \vr)} \ .
\end{equation}

To prove (C), we define the action $\ap$ of $W_{\va}(\vr)$ on 
$\BCg(\va,\vr) = \prod_{i=1}^g B_{\cC_v}(a_i,r_i)$ by 
restricting the domain of $\crJ_{\va}$ to $\BCg(\va,\vr)$ and setting, 
for $w \in W_{\va}(\vr)$ and $\vx \in \BCg(\va,\vr)$, 
\begin{equation} \label{FNuz} 
w \ap \vx \ = \ \crJ_{\va}^{-1}(w \Jp \crJ_{\va}(\vx))) \ .
\end{equation}
Applying $\crJ_{\va}$ to (\ref{FNuz}) shows that $[w \ap \vx] \Jm [\va] = w \Jp ([\vx] \Jm [\va])$, 
or equivalently that 
\begin{equation} \label{FAct}
[w \ap \vx] = w \Jp [\vx] \ .
\end{equation} 

Tracing through the definitions, one sees that if $w = \Psi(\vw)$ and $\vx = \Phi_{\va}(\vz)$
with $\vw \in L_{\va}(D(\vORIG,\vr))$ and $\vz \in D(\vORIG,\vr)$, 
then $w \ap \vx = \Phi_{\va}(\vw \ap_0 \vz)$. 
Thus, $\ap$ is the pushforward of the action $\ap_0$ to $\cC_v^g$ and $\Jac(\cC_v)$, 
\index{Jacobian variety}  
using the maps $\Phi_{\va}$ and $\Psi$. 
Since $D(\vORIG,\vr)$ is a principal homogeneous space\index{principal homogeneous space} for $\cD_{L_\va}(\vORIG,\vr))$ under $\ap_0$,
it follows that $\prod_{i=1}^g B_{\cC_v}(a_i,r_i)$ is a 
principal homogeneous space\index{principal homogeneous space} for $W_{\va}(\vr)$ 
under $\ap$. 

To prove (D), given $\vb \in \prod_{i=1}^g B_{\cC_v}(a_i,r_i)$, write $\vb = \Phi_{\va}(\vbeta)$
with $\vbeta \in D(\vORIG,\vr)$.  Since the component functions of $\Phi_{\va}$ are isometric parametrization
of the balls $B_{\cC_v}(a_i,r_i)$, for each $\vz \in D(\vORIG,\vr)$ 
we will have $\Phi_{\va}(\vz) \in \prod_{i=1}^g B_{\cC_v}(b_i,\eta r_i)$ if and only if 
$\vz = \vbeta + \vdelta$ with $\vdelta \in D(\vORIG,\eta \vr)$. 
Suppose $\Phi_{\va}(\vz) \in \prod_{i=1}^g B_{\cC_v}(b_i,\eta r_i)$.   
Since $W_{\va}(\eta \vr) = \Psi(L_{\va}(D(\vORIG,\eta \vr)))$, 
by (\ref{FDEquiv}) and (\ref{FCong3})
\begin{equation*} 
\crJ_{\va}(\Phi_{\va}(\vz)) 
\ \equiv \ \Psi(L_{\va}(\vb) + L_{\va}(\vdelta)) 
\ \equiv \ \Psi(L_{\va}(\vb)) 
\ \equiv \ \crJ_{\va}(\vb) \pmod{W_{\va}(\eta \vr)} \ .
\end{equation*} 
Hence $\Phi_{\va}(\vz) \in W_{\va}(\eta \vr) \, \ap \, \vb$. 
Conversely, if $\Phi_{\va}(\vz) \in W_{\va}(\eta \vr) \, \ap \, \vb$, then  
\begin{equation*} 
\Psi(L_{\va}(\vz)) \equiv \crJ_{\va}(\Phi_{\va}(\vz)) 
\ \equiv \ \crJ_{\va}(\vb) \ \ \equiv \ \Psi(L_{\va}(\vbeta)) \pmod{W_{\va}(\eta \vr)} \ ,
\end{equation*} 
which means that $\vz - \vbeta \in D(\vORIG,\eta \vr)$, and hence that 
$\Phi_{\va}(\vz) \in \prod_{i=1}^g B_{\cC_v}(b_i,\eta r_i)$.
Thus $W_{\va}(\eta \vr) \, \ap \, \vb = \prod_{i=1}^g B_{\cC_v}(b_i,\eta r_i)$.
It follows immediately that $J_{\vb}(\prod_{i=1}^g B_{\cC_v}(b_i,\eta r_i)) = W_{\va}(\eta \vr)$.

\smallskip
Next, we prove (E).  Let $F_u/K_v$ be a finite extension in $\CC_v$.
If $\va \in \cC_v(F_u)^g$, 
then $\Phi_{\va}$, $L_{\va}$, $\crJ_{\va}$ and $\Psi$ are $F_u$-rational.  
Since $\Phi_{\va}$ and $\Psi$ are isometric parametrizations, this means that 
$\Phi_{\va}(D(\vORIG,\vr) \cap F_u^g) = \prod_{i=1}^g \big(B_{\cC_v}(a_i,r_i) \cap \cC_v(F_u)\big)$ 
and $\Psi\big(L_{\va}(D(\vORIG,\vr) \cap F_u^g)\big) = W_{\va}(\vr) \cap \Jac(\cC_v)(F_u)$.  
Since $\crJ_{\va}$ is $F_u$-rational, 
$\prod_{i=1}^g \big(B_{\cC_v}(a_i,r_i) \cap \cC_v(F_u)\big)$ is a 
principal homogeneous space\index{principal homogeneous space} for 
$W_{\va}(\vr) \cap \Jac(\cC_v)(F_u)$ under $\ap$.   This implies (\ref{FRat1A}) and (\ref{FRat2A}).
\index{Jacobian variety}   

To prove (F), assume in addition that $F_u/K_v$ is separable.  
Let $(\alpha_{ij}) \in M_g(F_u)$ be the matrix of the linear map $L_{\va}^{-1}$, 
and let $A = \max_{ij}(|\alpha_{ij}|_v)$ be its operator norm.  
Write $r = \min_i(r_i)$, and put $s = (1/A)r$. Then $L_{\va}^{-1}(D(\vORIG,s)) \subseteq D(\vORIG,\vr)$, so  
$D(\vORIG,s) \subseteq L_{\va}(D(\vORIG,\vr))$, hence $B_J(O,s)$ is contained in $W_{\va}(\vr)$.
Let $0 < R_u, C_u \le 1$ be the constants from Lemma \ref{TraceBall}.   
If $s \le R_u$, then by Lemma \ref{TraceBall}, $\Tr_{F_u/K_v}\big(W_{\va}(\vr) \cap \Jac(\cC_v)(F_u)\big)$ 
contains $B_J(O,C_u s) \cap \Jac(\cC_v)(K_v)$.  
If $s > R_u$, it contains $B_J(0, C_u R_u) \cap \Jac(\cC_v)(K_v)$. 
\index{Jacobian variety}  
Put $C = C_u \min(1/A, R_u)$. Since $ r \le 1$, in either case 
\begin{equation*} 
\Tr_{F_u/K_v}\big(W_{\va}(\vr) \cap \Jac(\cC_v)(F_u)\big) \ \supseteq \ B_J(O,C r) \cap \Jac(\cC_v)(K_v)\ .
\end{equation*} 

This completes the proof. 
\end{proof}

Finally, we give the proof of Proposition \ref{BPropD2New}.\index{Lipschitz continuity!of the Abel map|ii}

\begin{proof}[Proof of Proposition \ref{BPropD2New}.]
Let $W_{\va}(\vr) = \BFJ_{\va}\big(\prod_{i=1}^g B_{\cC_v}(a_i,r_i)\big) = L_{\va}(D(\vORIG,\vr))$ 
be an open subgroup of $\Jac(\cC_v)(\CC_v)$ with the properties in Theorem \ref{BKeyThm1},
\index{Jacobian variety}   
and let $E_v \subset \cC_v(\CC_v)$ be a nonempty compact subset.  

Next consider the map $\BFj_x(y) = [(y)-(x)]$, 
which takes $\cC_v(\CC_v) \times \CC_v(\CC_v)$ to $\Jac(\cC_v)(\CC_v)$.  
We claim that for each  $\tau \in E_v$, there are an isometrically parametrizable ball 
$B_{\cC_v}(\tau,R_\tau)$ with $\BFj_\tau(B_{\cC_v}(\tau,R_\tau)) \subseteq W_{\va}(\vr)$, 
such that if $0 < \varepsilon \le R_\tau$ 
then for all $x, y \in B_{\cC_v}(\tau,R_\tau)$ with $\|x,y\|_v \le \varepsilon$, 
we have $\BFj_x(y) \in W_{\va}((\varepsilon/R_\tau) \cdot \vr)$.  

To see this,
let $\varphi_\tau : D(0,R_\tau) \rightarrow B_{\cC_v}(\tau,R_\tau)$ be an isometric parametrization.
Without loss we can assume that $R_\tau \in |\CC_v^{\times}|_v$.  Let $g_\tau(X)$ be the composite
of the sequence of maps
\begin{equation*}
D(0,R_\tau) \stackrel{\varphi_\tau}{\longrightarrow} B_{\cC_v}(\tau,R_\tau) \stackrel{\BFj_\tau}{\longrightarrow} 
W_{\va}(\vr) \stackrel{\Psi^{-1}}{\longrightarrow} D(\vORIG,1)^- \ \subset \ \CC_v^g \ .
\end{equation*} 
Here $\varphi_\tau$ is an analytic map defined by convergent power series in $\CC_v[[X]]$, 
$\BFj_\tau$ is an algebraic morphism, and $\Psi^{-1}$ is the inverse of the isometric parametrization 
$\Psi : D(\vORIG,1)^- \rightarrow  B_J(O,1)^-$; it is given by projection on $g$ of the coordinates 
(see \ref{JPsiDef}).  Thus $g_\tau(X) = (g_{\tau,1}(X), \ldots, g_{\tau,g}(X))$ is a map whose coordinate
functions are defined by convergent power series in $\CC_v[[X]]$. 

Since $\BFj_x(y) = \BFj_\tau(y) \Jm \BFj_\tau(x)$ on $B(\tau,R_\tau)$,
and since the group operations $\Jp$, $\Jm$ on $B_J(O,1)^-$ 
correspond to $S(\vX,\vY)$, $M(\vX)$ in the formal group $D(\vORIG,1)^-$,
\index{formal group} 
when $\BFj_x(y)$ is pulled back to $D(0,R_\tau) \times D(0,R_\tau)$ using $\varphi_\tau$, 
it is represented by the power series map $S(g_\tau(X),M(g_\tau(Y))$.  By hypothesis,
$\BFj_\tau(B(\tau,R_\tau)$ is contained in $W_{\va}(\vr) = L_{\va}(D(0,\vr))$, 
where $L_{\va} = (\Psi^{-1} \circ \BFJ_{\va} \circ \Phi_{\va})^{\prime}(0) : \CC_v^g \rightarrow \CC_v^g$
is a nonsingular linear map.  Since $W_{\va}(\vr)$ is a group
the image of $S(g_\tau(X),M(g_\tau(Y))$ is contained in $L_{\va}(D(0,\vr))$.

Now define $h_\tau(X,Y) : D(0,R_{\tau}) \times D(0,R_{\tau}) \rightarrow \CC_v^g$ by 
\begin{equation*}
h_{\tau}(X,Y) \ = \ L_{\va}^{-1}(S(g_\tau(X),M(g_\tau(Y))) 
\ = \ (h_{\tau,1}(X,Y), \ldots, h_{\tau,g}(X,Y))\ .
\end{equation*} 
The image of $h_{\tau}(X,Y)$ is contained in $D(0,\vr) = \prod_{j=1}^g D(0,r_j)$, 
so each $h_{\tau,j}(X,Y) \in \CC_[[X,Y]]$ 
is a map defined by power series converging on $D(0,R_\tau) \times D(0,R_\tau)$, 
with image contained in $D(0,r_j)$.  Since $\BFj_x(x) = 0$ for each $x$, it follows that
$h_{\tau,j}(x,x) = 0$ for all $x \in D(0,R_\tau)$.  It follows from Lemma \ref{BPropD2New} that
$|h_{\tau,j}(x,y)|_v \le |x-y|_C \cdot (r_j/R_\tau)$ for all $x, y \in D(0,R_\tau)$.  
In particular, if $0 < \varepsilon \le R_\tau$ and $x, y \in D(0,R_\tau)$ satisfy 
$|x-y|_v \le \varepsilon$, then $h_\tau(x,y) \in D(0,(\varepsilon/R_\tau) \vr)$.
Since $\varphi_\tau$ is an isometric parametrization, 
and since $L_{\va}(D(0,(\varepsilon/R_\tau) \vr) = W_{\va}((\varepsilon/R_\tau) \vr)$, 
for all $x, y \in B(\tau,R_{\tau})$ 
with $\|x,y\|_v \le \varepsilon$ we have $\BFj_x(y) \in W_{\va}((\varepsilon/R_\tau) \cdot \vr)$.

To conclude the proof, consider the nonempty compact set $E_v \subset \cC_v(\CC_v)$. 
Cover $E_v$ with finitely many balls $B(\tau_1,R_{\tau_1}), \ldots, B(\tau_n,R_{\tau_n})$
satisfying the conditions above.  Put $\varepsilon_0 = \min_{1 \le j \le n}(R_{\tau_j})$, 
and put $C_0 = 1/\varepsilon_0$.  Then for each $0 < \varepsilon \le \varepsilon_0$,
and all $x, y \in E_v$ with $\|x,y\|_v \le \varepsilon$, 
we have $\BFj_x(y) \in W_{\va}(C_0 \varepsilon \cdot \vr)$.
\end{proof}
\index{local action of the Jacobian|)}

%% file: NewFSZ.bbl
\newpage
\begin{thebibliography}{999}

\bibitem{Ahl} L. Ahlfors, Complex Analysis, McGraw-Hill, 1979.

\bibitem{Akh} N.I. Akhiezer, Elements of the Theory of Elliptic Functions, 
AMS Translations of Mathematical Monographs 79, AMS, Providence, 1990.   

\bibitem{Am} Y. Amice, Les nombres p-adiques, Presses Universitaires de France,
   Paris, 1975.

\bibitem{Ar} E. Artin, Algebraic Numbers and Function Fields,
   Gordon and Breach, New York, 1967
   
\bibitem{Artin} M.~Artin, N\'eron Models, chapter VIII in:  G. Cornell, J. Silverman, eds., 
Arithmetic Geometry, Springer-Verlag, New York, 1986. 
   
\bibitem{Aut} P. Autissier, Pointes entiers et th\'eor\`emes de Bertini arithm\'etiques, 
Th\`ese de l'Universit\'e Paris XI (D\'ecembre 2001). 

\bibitem{B-R} M. Baker and R. Rumely, Potential Theory and Dynamics on the 
Berkovich Projective Line, Mathematical Surveys and Monographs 159, AMS, Providence, 2010. 

\bibitem{BIR} M. Baker, S. Ih, and R. Rumely, A finiteness property of torsion points,
Algebra and Number Theoyr 2 (2008), 217-248. 

\bibitem{BY} M. L. Balinski and H. P. Young, The quota method of apportionment,
American Mathematical Monthly 82 (1975), 701-730.

\bibitem{Berk} V. Berkovich, Spectral theory and analytic geometry over non-Archimedean fields, 
Mathematical Surveys and Monographs 33, AMS, Providence, 1990.

\bibitem{BGR} S. Bosch, U. G\"untzer, R. Remmert, Non-Archimedean Analysis,
   Grundlehren der mathematischen Wissenschaften 261, Springer-Verlag,
   New York, 1984

\bibitem{BLR} S. Bosch, W. L\"utkebohmert, M. Raynaud, N\'eron Models, 
Ergebnisse der Mathematik Band 21, Springer-Verlag, New York, 1990.

\bibitem{C-H} R. Courant and D. Hilbert, Methods of Mathematical Physics (Volume 1),
Interscience, New York, 1953. 

\bibitem{Can1} D. Cantor, On approximation by polynomials with algebraic
integer coefficients, in: Proc. Symposia in Pure Math 12  
(W. Leveque and E. Strauss, eds.), AMS Providence, R.I., 1969.

\bibitem{Can2} D. Cantor, On certain algebraic integers and approximation by
rational functions with integral coefficients, Pacific J. Math 67 (1976), 
323-338.

\bibitem{Can3} D. Cantor, On an extension of the definition of transfinite
diameter and some applications, J. Reine Angew. Math 316 (1980), 160-207.

\bibitem{C-R} D. Cantor and P. Roquette, On diophantine equations over
the ring of all algebraic integers, J. Number Theory 18 (1984), 1-26.

\bibitem{Ca-F} J. W. Cassels and A. Fr\"olich, Algebraic Number Theory,
Academic Press, London, 1967. 

\bibitem{CR} T. Chinburg and R. Rumely, The capacity pairing,
J. Reine Angew. Math 434 (1993), 1-44.

\bibitem{PClark} P. Clark, 
Rational points on Atkin-Lehner Quotients of Shimura Curves,  
Doctoral Thesis, Harvard University, 2003. 

\bibitem{Cremona} J. Cremona, 
Algorithms for Modular Elliptic Curves, 
Cambridge University Press, Cambridge, 1992. 

\bibitem{Fall} Th. Falliero, S\'eries d'Eisenstein Hyperboliques et Capacit\'e 
d'Unions d'Intervalles, Ph.D. Thesis, University of Bordeaux I, 1998.  

\bibitem{FS} Th. Falliero and A. Sebbar, Capacit\'e de trois intervalles 
et Fonctions th\^{e}ta de genre 2, J. Math. Pures. Appl 80 (2001), 409--443. 

\bibitem{FS2} Th. Falliero and A. Sebbar, Capacities and Jacobi Matrices, 
Proceedings of the Edinburgh Mathematical Society 46 (2003), 719-745.

\bibitem{F-SZ} M. Fekete and G. Szeg\"o, On algebraic equations with 
integral coefficients whose roots belong to a given point set,
Math. Zeit. 63 (1955), 158-172.

\bibitem{Fies} K.-H. Fieseler, Zariski's Main Theorem f\"ur affinoide Kurven, 
          Math. Ann. 251 (1980), 97-110. 

\bibitem{FJar} M. D. Fried and M. Jarden, Field Arithmetic, Ergebnisse der Mathematik, 3 Folge, Band 11,
Springer-Verlag, Berlin-Heidelberg-New York, 1986.

\bibitem{G} F. R. Gantmacher, The Theory of Matrices, vol. II, Chelsea, 1974.

\bibitem{Gro} A. Grothendieck and J. Dieudonn\'e, El\'ements de Geometrie Alg\'ebrique, 
IHES Puplcations Math\'ematiques 17, 29, 24, 28, 32, Presses Universitaire de France, Paris, 1961-1967.

\bibitem{G-M-P} B. Green, M. Matignon and F. Pop, On the local Skolem property,
J. Reine Angew. Math. 458 (1995), 183-199.

\bibitem{Klei} S. Kleiman, The Picard scheme,
pp. 235--321 in: Fundamental Algebraic Geometry, Mathematical Surveys and Monographs 123,
AMS, Providence, 2005. 

\bibitem{Kl}  M. Klimek, Pluripotential theory, London Mathematical Society
Monographs (New Series) 6, Oxford Science Publications, 1991.

\bibitem{Kob}  H. Kober, Dictionary of Conformal Representations, 
Dover, New York, 1952.

\bibitem{LANT}  S. Lang, Algebraic Number Theory, Addison-Wesley, Reading, 1968.

\bibitem{L}  S. Lang, Introduction to algebraic geometry, Addison-Wesley, 
Reading;  $4^{th}$ printing, 1973.

\bibitem{LE} S. Lang, Elliptic Functions, Addison-Wesley, Reading, 1973.  

\bibitem{Mat} H. Matsumura, Commutative Algebra (2nd edition), 
Lecture notes in Mathematics Series, Benjamin/Cummings Publishing Company, Reading, 1980.

\bibitem{MB1} L. Moret-Bailly, Groupes de Picard et probl\`emes de  Skolem I, 
Ann. Scient. \'Ec. Norm. Sup (4 ser.) 22 (1989), 161-179.

\bibitem{MB2} L. Moret-Bailly, Groupes de Picard et probl\`emes de  Skolem II, 
Ann. Scient. \'Ec. Norm. Sup (4 ser.) 22 (1989), 181-194.

\bibitem{MB3} L. Moret-Bailly, Probl\'emes de Skolem sur les champs algebraiques, 
Compositio Math. 125 (2001), 1-30.  

\bibitem{MR} B. Mazur, Modular Curves and the Eisenstein Ideal, 
IHES Publ. Math 47 (1977), 33-186:  Appendix by B. Mazur and M. Rapoport, 
``Behavior of the N\'eron model of the Jacobian of $X_0(N)$ at bad primes''.   

\bibitem{McC}  W. McCallum, The degenerate fibre of the Fermat Curve, 
in:  N. Koblitz, ed., Number theory related to Fermat's Last Theorem, 
Progress in Math. 26, Birkha\"user, Boston (1982), 57--70.

\bibitem{Milne} J.S.~Milne, Jacobian Varieties, chapter VII in:  G. Cornell, J. Silverman, eds., 
Arithmetic Geometry, Springer-Verlag, New York, 1986. 


\bibitem{Neh} Z. Nehari, Conformal Mapping, Dover, New York, 1952.

\bibitem{DPark} D. Park, The Fekete-Szeg\"o theorem with splitting conditions
on the projective line of positive characteristic $p$, Doctoral thesis,
University of Georgia, 2006.

\bibitem{P} G. P\'olya, \"Uber ganzwertige Polynome in algebraischen
Zahlk\"orpern, J. Reine Angew. Math. 149 (1919), 97-116.

\bibitem{P-S} A. Prestel and J. Schmid, Decidability of the rings of real algebraic
and $p$-adic algebraic integers, J. Reine Angew. Math. 414 (1991), 141-148.

\bibitem{Rob1} R. Robinson, Conjugate algebraic numbers in real point sets, 
Math Zeit. 84 (1964), 415-427.

\bibitem{Rob2} R. Robinson, Intervals containing infinitely many 
conjugate sets of algebraic units, Ann Math (2) 80 (1968), 532-543.  

\bibitem{R-G-P} P. Roquette, B. Green and F. Pop, On Rumely's local-global principle,
Jahresber. Deutch. Math.-Verein 97 (1995), 43-74.

\bibitem{RR1} R. Rumely, Capacity theory on algebraic curves, Lecture Notes
in Mathematics 1378, Springer-Verlag, Berlin-Heidelberg-New York, 1989.

\bibitem{RR2} R. Rumely, The Fekete-Szeg\"o theorem with splitting conditions,
part I, Acta Arithmetica 93 (2000), 99-116. 

\bibitem{RR3} R. Rumely, The Fekete-Szeg\"o theorem with splitting conditions, 
part II, Acta Arithmetica 103 (2002), 347--410.

\bibitem{RR4} R. Rumely, On the relation between Cantor's capacity and the 
sectional capacity, Duke Mathematical Journal 70 (1993), 517-574.

\bibitem{RR5} R. Rumely, Arithmetic over the ring of all algebraic integers,
J. Reine Ange. Math. 368 (1986), 127-133.  

\bibitem{RR6} R. Rumely, An intersection theory for curves with analytic contributions 
from nonarchimedean places,  in:  Canadian Mathematical Society Conference Proceedings 15
(1995), AMS, Providence, 325-357.

\bibitem{ST} E. Saff and V. Totik, Logarithmic potentials with external fields,
Grundlehren  316, Springer-Verlag, Berlin-Heidelberg-New York, 1991.  

\bibitem{SS} M. Schiffer and D. Spencer, Functionals of Finite Riemann Surfaces,
Pinceton University Press, Princeton, 1954.  

\bibitem{LALG} J.-P. Serre, Lie Algebras and Lie Groups, 1964 Lectures given at 
Harvard University, W. A. Benjamin, New York, 1965.  

\bibitem{ShimuraTheta} G. Shimura, Theta functions with complex multiplication,
Duke Mathematical Journal 43 (1976), 673--696.

\bibitem{Sil2} J. Silverman, Advanced Topics in the Arithmetic of Elliptic Curves,
Graduate Texts in Mathematics 151, Springer-Verlag, New York, 1994.  

\bibitem{SM} N. Steinmetz, Rational Iteration, Walter de Gruyter, New York, 1993.

\bibitem{Tam} A. Tamagawa, Unramified Skolem Problems and Unramified Bertini Theorems in Positive Characteristic, 
Documenta Mathematica, Extra Volume Kato (2003), 789--831.  

\bibitem{Th} A. Thuillier,   
Th\'eorie du potentiel sur les courbes en g\'eometrie analytique
non archimje\'edienne.  Applications \`a la th\'eorie d'Arakelov.
Ph.D. Thesis, University of Rennes, 2005.  Preprint available at
{http://tel.ccsd.cnrs.fr/documents/archives0/00/01/09/90/index.html}.  

\bibitem{Ts}  M. Tsuji, Potential Theory in Modern Function Theory, Maruzen,
Tokyo, 1959 (reprinted: Chelsea, New York, 1975).

\bibitem{vdD} L. van den Dries, 
Elimination theory for the ring of all algebraic integers,
J. Reine Angew. Math. 388 (1988), 189-205. 

\bibitem{WW} E.T. Whittaker and G.N. Watson, 
A Course of Modern Analysis ($4^{th}$) edition, 
Cambridge University Press, 1927;  reprinted 1973.

\bibitem{Wid} H. Widom,  Extremal polynomials associated with a system 
of curves in the complex plane, Advances in Math. 3 (1969), 127--232.

%

\end{thebibliography}
